\documentclass[11pt]{book}             
\usepackage{amssymb}
\usepackage{amsfonts}
\usepackage{amsmath}
\usepackage{geometry}
\usepackage{bm}
\numberwithin{equation}{section}
\usepackage[toc,page]{appendix}
\usepackage{rotating}

\newtheorem{theorem}{Theorem}[section]
\newtheorem{definition}[theorem]{Definition}
\newtheorem{lemma}[theorem]{Lemma}

\newtheorem{remark}[theorem]{Remark}

\newenvironment{proof}[1][Proof]{\begin{trivlist}
\item[\hskip \labelsep {\bfseries #1}]}{\end{trivlist}}
\newcommand{\qed}{\nobreak \ifvmode \relax \else
      \ifdim\lastskip<1.5em \hskip-\lastskip
      \hskip1.5em plus0em minus0.5em \fi \nobreak
      \vrule height0.75em width0.5em depth0.25em\fi}

\raggedright                            

\title{\bf Special functions and reversible three-term recurrence formula (R3TRF)}    
\author{\textsc{Yoon Seok Choun} \footnote{The first series was `Special functions and three-term recurrence formula (3TRF).' It is available as arXiv. This book is a second series.} \footnote{Yoon.Choun@baruch.cuny.edu; ychoun@gmail.com}\\
  Baruch College, The City University of New York,\\
  Natural Science Department, A506,\\
  17 Lexington Avenue,\\
  New York, NY 10010
} 
\date{\today}                           

\begin{document}
\frontmatter                            
\maketitle 
\tableofcontents                        
\mainmatter 
\chapter*{Preface}

For the past 350 years, we only have constructed the power series expansion in closed forms using the two term recurrence relation in linear differential equation. However, currently the analytic solution and its integral forms for more than three term cases are unknown.

For the three term case, we assume that its solution is
\begin{equation}
y(x)= \sum_{n=0}^{\infty } c_n x^{n+\lambda }\nonumber
\end{equation}
On the above, $\lambda $ is an indicial root, $y(x)$ is its analytic solution for the linear differential equation and $x$ is an independent variable for its solution. When we substitute the above equation into any linear ordinary differential equation having three different coefficients in their power series expansion, the recurrence relation for frobenius solution is
\begin{equation}
c_{n+1}=A_n \;c_n +B_n \;c_{n-1} \hspace{1cm};n\geq 1
\nonumber
\end{equation}
where
\begin{equation}
c_1= A_0 \;c_0
\nonumber
\end{equation}
where $c_1= A_0 \;c_0$ and $\lambda $ is an indicial root. On the above, $A_n$ and $B_n$ are themselves polynomials of degree $m$: for the second-order ODEs, a numerator and a denominator of $A_n$ are usually equal or less than polynomials of degrees 2. Likewise, a numerator and a denominator of $B_n$ are also equal or less than polynomials of degrees 2
   
For the case of three term recurrence in its power series expansion, all other differential equations having no analytic solution in closed forms can be described as in the above recurrence relation equation.\footnote{For more details, look at arXiv:1303.0806 in ``Special functions and three term recurrence formula (3TRF).''}  

In the first series ``Special functions and three term recurrence formula (3TRF)'', I showed how to construct the power series expansions in closed forms, its integral representations and its generating functions for linear ordinary differential equations (infinite series and polynomial which makes $B_n$ term terminated) such as Heun, Mathieu, Lame, Grand Confluent Hypergeometric (GCH) equations.\footnote{These 4 linear differential equations have three different coefficients in their power series expansion. Biconfluent Heun equation is the special case of GCH equation. Look at arXiv:1303.0813 and arXiv:1303.0819 in ``Special functions and three term recurrence formula (3TRF).''} There are three types polynomial in their three term recurrence relation: (1) polynomial which makes $B_n$ term terminated, (2) polynomial which makes $A_n$ term terminated, (3) polynomial which makes $A_n$ and $B_n$ terms terminated at same time. However, by using 3TRF, there are no way to construct the power series in closed forms indeed its integral representation for the type 2 and 3 polynomials. 

In this series natural numbers $\mathbb{N}_{0}$ means $\{ 0,1,2,3,\cdots \}$. Pochhammer symbol $(x)_n$ is used to represent the rising factorial: $(x)_n = \frac{\Gamma (x+n)}{\Gamma (x)}$. And I will show how to obtain the formula for the type 2 polynomial and infinite series of its power series expansion. In the future series I will show you for the type 3 case; the power series, integral formalism and generating function such as Heun, Confluent Heun, Lame and GCH equations will be constructed analytically.    
\section*{Structure of book}

In chapter 1, I will generalize the three term recurrence relation in linear differential equation in a backward for the infinite series and polynomial which makes $A_n$ term terminated including all higher terms of $B_n$'s.

In chapters 2--9, I will apply reversible three term recurrence formula to (1) the power series expansion in closed forms, (2) its integral representation and (3) generating functions for Heun, Confluent Heun, GCH, Lame and Mathieu equations that consist of three term recursion relation for the infinite series and polynomial which makes $A_n$ term terminated including all higher terms of $B_n$'s. The short descriptions of all 9 chapters are as follows.

\vspace{1mm}

\textbf{Chapter 1.} ``Generalization of the reversible three-term recurrence formula and its applications''-- generalize the three term recurrence relation in linear differential equation in a backward. Obtain the exact solution of the three term recurrence relation for infinite series and polynomial which makes $A_n$ term terminated
\vspace{1mm}

\textbf{Chapter 2.} ``Heun function using reversible three-term recurrence formula''-- apply the reversible three term recurrence formula to the power series expansion in closed forms and its integral forms of Heun functions (infinite series and polynomial which makes $A_n$ term terminated) including all higher terms of $B_n$'s. 
\vspace{1mm}

\textbf{Chapter 3.} ``The generating function for Heun polynomial using reversible three-term recurrence formula''-- apply the  reversible three term recurrence formula and derive generating function for Heun polynomial which makes $A_n$ term terminated including all higher terms of $B_n$'s. 
\vspace{1mm}

\textbf{Chapter 4.} ``Confluent Heun function using three-term recurrence formula''-- (1) apply the three term recurrence formula to the power series expansion in closed forms and its integral forms of the Confluent Heun functions (infinite series and polynomial which makes $B_n$ term terminated) including all higher terms of $A_n$'s. (2) apply the three term recurrence formula and derive the generating function for the Confluent Heun polynomial which makes $B_n$ term terminated including all higher terms of $A_n$'s.
\vspace{1mm}

\textbf{Chapter 5.} ``Confluent Heun function using reversible three-term recurrence formula''-- (1) apply the reversible three term recurrence formula to the power series expansion in closed forms and its integral forms of the Confluent Heun functions (infinite series and polynomial which makes $A_n$ term terminated) including all higher terms of $B_n$'s. (2) apply three term recurrence formula and derive the generating function for the Confluent Heun polynomial which makes $A_n$ term terminated including all higher terms of $B_n$'s.
\vspace{1mm}

\textbf{Chapter 6.} ``Grand Confluent Hypergeometric function using reversible three-term recurrence formula''-- (1) apply the reversible three term recurrence formula and formulate the analytic solution for the power series expansion and its integral forms of grand confluent hypergeometric function (infinite series and polynomial which makes $A_n$ term terminated) including all higher terms of $B_n$'s. Replacing $\mu $ and $\varepsilon \omega $ by 1 and $-q$ transforms the grand confluent hypergeometric function into the  Biconfluent Heun function. 
(2) apply the reversible three term recurrence formula and construct the generating function for grand confluent hypergeometric  polynomial which makes $A_n$ term terminated including all higher terms of $B_n$'s. 
\vspace{1mm}

\textbf{Chapter 7.} ``Mathieu function using reversible three-term recurrence formula''-- apply the reversible three term recurrence formula, and analyze the power series expansion in closed forms of Mathieu function and its integral forms (infinite series and polynomial which makes $A_n$ term terminated) including all higher terms of $B_n$'s. Construct the generating function for Mathieu polynomial which makes $A_n$ term terminated. 
\vspace{1mm}

\textbf{Chapter 8.} ``Lame function in the algebraic form using reversible three-term recurrence formula''-- (1) apply the reversible three term recurrence formula, and analyze the power series expansion in closed forms and integral forms of Lame function in the algebraic form (infinite series and polynomial which makes $A_n$ term terminated) including all higher terms of $B_n$'s. (2) apply the reversible three term recurrence formula and construct the generating function for Lame polynomial which makes $A_n$ term terminated in the algebraic form including all higher terms of $B_n$'s. 
\vspace{1mm}

\textbf{Chapter 9.} ``Lame function in Weierstrass's form using reversible three-term recurrence formula''-- (1) apply the  reversible three term recurrence formula, and analyze the power series expansion in closed forms of Lame function in Weierstrass's form and its integral forms (infinite series and polynomial which makes $A_n$ term terminated) including all higher terms of $B_n$'s. (2) apply the reversible three term recurrence formula and construct the generating function for Lame polynomial which makes $A_n$ term terminated in Weierstrass's form including all higher terms of $B_n$'s. 

\section*{Acknowledgements}
I thank Bogdan Nicolescu.  The endless discussions I had with him were of great joy.  
\chapter{Generalization of the reversible three-term recurrence formula (R3TRF) and its applications}
\chaptermark{R3TRF}
In ``Generalization of the three-term recurrence formula and its applications''\cite{Chou2012}, I generalize the three term recurrence relation in a linear differential equation for infinite series and polynomial which makes $B_n$ term terminated including all higher terms of $A_n$'s.\footnote{`` higher terms of $A_n$'s'' means at least two terms of $A_n$'s.}  

In this chapter, I generalize the three term recurrence relation in a linear ordinary differential equation in a backward for infinite series and polynomial which makes $A_n$ term terminated including all higher terms of $B_n$'s.\footnote{`` higher terms of $B_n$'s'' means at least two terms of $B_n$'s.}

The rest of chapters in this series show how to apply reversible three-term recurrence formula (R3TRF) into linear differential equations having three term recursion relations and go on producing analytic solutions of some of the well known special functions including: Heun, Confluent Heun, GCH (Biconfluent Heun), Lame and Mathieu equations. 
\section{\label{sec:level1}Introduction}

The history of linear differential equations is over 350 years. By using Frobenius method and putting the power series expansion into linear differential equations, the recursive relation of coefficients starts to appear. There can be between two and infinity number of coefficients in the recurrence relation in the power series expansion. During this period mathematicians developed analytic solutions of only two term recursion relation in closed forms. Currently the analytic solution of three term recurrence relation is unknown.

In ``Generalization of the three-term recurrence formula and its applications''\cite{Chou2012}, I show how to obtain the Frobenius solutions for infinite series and polynomial which makes $B_n$ term terminated including all higher terms of $A_n$'s using three term recurrence formula (3TRF) in a linear ordinary differential equation having a recursive relation between a 3-term. The sequence $c_n$ consists of combinations $A_n$ and $B_n$. I observed the term inside parenthesis of sequence $c_n$ which includes zero term of $A_n$'s, one term of $A_n$'s, two terms of $A_n$'s, three terms of $A_n$'s, etc in order to construct general expressions of the power series in closed forms for infinite series and polynomial which makes $B_n$term terminated. 

Any problems having three term recurrence relation in their linear differential equations in modern physics (E \& M, Newtonian mechanics, quantum mechanic, QCD, supersymmetric field theories, string theories, general relativity, etc)  mostly involve the quantum mechanical energy in $A_n$ or $B_n$ term. In general, the energy is quantized in quantum mechanical points of a view. The power series expansion of a eigen-function is the polynomial. In Ref.\cite{Chou2012} I construct the Frobenius solution in closed forms of three term recurrence relation in a linear differential equation for polynomial which makes $B_n$ term terminated. Also we need to derive the analytic solution of three term recurrence relation in a linear differential equation for polynomial which makes $A_n$ term terminated. Because in some modern physics problems the eigenvalue (quantum mechanical energy) is involved in $A_n$ term: Mathieu equation, Heun equation and its confluent forms (black hole problems \cite{Suzu1998,Suzu1999}) , Lame equation in the Weierstrass's form \cite{Lame1837,Erde1955,Hobs1931,Whit1952}, etc. 
 
In this chapter I construct three term recurrence relation in a linear ordinary differential equation in a backward to obtain the analytic solution of the three term recurrence relation for infinite series and polynomial which makes $A_n$ term terminated including all higher terms of $B_n$'s.
I designate the mathematical formulas of infinite series and polynomial which makes $A_n$ terminated of a linear differential equation having a 3-term recurrence relation between successive coefficients in its power series expansion as ``reversible three-term recurrence formula (R3TRF)'' that will be expressed below. 

In chapters 2--9, using the R3TRF, I show how to obtain analytic solutions of infinite series and polynomial which makes $A_n$ term terminated (1) in the power series expansion, (2) in the integral formalism and (3) in the generating function for any linear ordinary differential equations having three term recurrence relation. 

\section[Mathieu equation, Frobenius method and reversible three term recurrence relation]{Mathieu equation, Frobenius method and reversible three term recurrence relation
  \sectionmark{Mathieu equation, Frobenius method and R3TRF}}
\sectionmark{Mathieu equation, Frobenius method and R3TRF}
In some rare cases, there are no such solutions of the power series expansion for polynomial which makes $B_n$ term terminated in a linear ordinary differential equation having a 3-term recurrence relation. One of examples is Mathieu differential equation in the algebraic form.
Mathieu functions\cite{Guti2003}, is an example of three term recurrence relation appears in physical problems involving elliptical shapes\cite{Troe1973} or periodic potentials, were introduced by Mathieu (1868)\cite{Math1868} when he investigated the vibrating elliptical drumhead. 

Mathieu ordinary differential equation is of Fuchsian types with the two regular and one irregular singularities. In contrast, Heun equation of Fuchsian types has the three regular and one irregular singularities \cite{Heun1889}. Heun equation has the four kind of confluent forms: (1) Confluent Heun (two regular and one irregular singularities) \cite{Ronv1995,Deca1978}, (2) Doubly confluent Heun (two irregular singularities), (3) Biconfluent Heun (one regular and one irregular singularities) \cite{Maro1979}, (4) Triconfluent Heun equations (one irregular singularity).  Mathieu equation in algebraic forms is also derived from the Confluent Heun equation by changing all coefficients $\delta =\gamma =\frac{1}{2}$, $\beta =0$, $\alpha =\frac{q}{0}$ and $q=\frac{\lambda +2q}{4}$: see chapters 4 and 5.

Mathieu equation is 
\begin{equation}
 \frac{d^2{y}}{d{z}^2} + \left( \lambda - 2q\;\mathrm{\cos}2z \right) y = 0 \label{eq:1}
\end{equation}
where $\lambda $ and $q$ are parameters. This is an equation with periodic-function coefficient. Mathieu equation also can be described in algebraic forms putting $x=\mathrm{\cos}^2z$:
\begin{equation}
 4x (1-x ) \frac{d^2{y}}{d{x}^2} + 2( 1-2x) \frac{d{y}}{d{x}} + ( \lambda + 2 q - 4 q x ) y = 0\label{eq:2}
\end{equation}
This equation has two regular singularities: $x=0$ and $x=1$; the other singularity $x=\infty $ is irregular. Assume that its solution is
\begin{equation}
y(x)= \sum_{n=0}^{\infty } c_n x^{n+\nu }\hspace{1cm}\mbox{where}\;\nu =\mbox{indicial}\;\mbox{root}\label{eq:3}
\end{equation}
Plug (\ref{eq:3}) into (\ref{eq:2}):
\begin{equation}
c_{n+1}=A_n \;c_n +B_n \;c_{n-1} \hspace{1cm};n\geq 1\label{eq:4}
\end{equation}
where
\begin{subequations}
\begin{equation}
A_n = \frac{4(n+\nu )^2-(\lambda +2q)}{2(n+1+\nu )(2(n+\nu )+1)}\label{eq:5a}
\end{equation}
\begin{equation}
B_n = \frac{4q}{2(n+1+\nu )(2(n+\nu )+1)}\label{eq:5b}
\end{equation}
\begin{equation}
c_1= A_0 \;c_0 \label{eq:5c}
\end{equation}
\end{subequations}
We have two indicial roots which are $\nu = 0$ and $\frac{1}{2} $. As we see (\ref{eq:5b}), there is no way to make $B_n$ term terminate at certain index value of $n$, because the numerator of (\ref{eq:5b}) is just consist of constant $q$ parameter.\footnote{Whenever index $n$ increase in (\ref{eq:5b}), the $B_n$ term never terminated with a fixed constant parameter $q$.} So there are only two kind of power series expansions: infinite series and the polynomial which makes $A_n$ term terminated. 

In ``The power series expansion of Mathieu function and its integral formalsim''\cite{Chou2013a}, I construct the analytic solution of Mathieu equation for infinite series by using 3TRF \cite{Chou2012}.
In this chapter I generalize the three term recurrence relation in a linear differential equation in a backward in order to construct the analytic solution for polynomial which makes $A_n$ term terminated.

\section{\label{sec:level3}Infinite series}
Assume that
\begin{equation}
c_1= A_0 \;c_0
\label{eq:5}
\end{equation}
(\ref{eq:5}) is a necessary boundary condition. The three term recurrence relation in all differential equations having no analytic solution in closed forms follow (\ref{eq:5}).
\begin{equation}
\prod _{n=a_i}^{a_i-1} B_n =1 \hspace{1cm} \mathrm{where}\;\mathrm{a_i}\;\mathrm{is}\;\mathrm{positive}\;\mathrm{integer}\;\mathrm{including}\;0
\label{eq:6}
\end{equation}
(\ref{eq:6}) is also a necessary condition. Every differential equations having no analytic solution in closed forms also take satisfied with (\ref{eq:6}). 

My definition of $B_{i,j,k,l}$ refer to $B_{i}B_{j}B_{k}B_{l}$. Also, $A_{i,j,k,l}$ refer to $A_{i}A_{j}A_{k}A_{l}$. For $n=0,1,2,3,\cdots $, (\ref{eq:4}) gives
\begin{equation}
\begin{tabular}{  l  }
  \vspace{2 mm}
  $c_0$ \\
  \vspace{2 mm}
  $c_1 = A_0 c_0 $ \\
  \vspace{2 mm}
  $c_2 = (A_{0,1}+B_1) c_0 $ \\
  \vspace{2 mm}
  $c_3 = (A_{0,1,2}+ B_1 A_2+ B_2 A_0) c_0 $\\
  \vspace{2 mm}
  $c_4 = (A_{0,1,2,3}+B_1 A_{2,3} + B_2 A_{0,3} + B_3 A_{0,1} + B_{1,3}) c_0 $ \\
  \vspace{2 mm}
  $c_5 = (A_{0,1,2,3,4}+ B_1 A_{2,3,4}+ B_2 A_{0,3,4} + B_3 A_{0,1,4} + B_4 A_{0,1,2} $  \\
  \vspace{2 mm}
  \hspace{0.8 cm} $ + B_{1,3} A_4 + B_{1,4} A_2 + B_{2,4} A_0 ) c_0 $ \\     
  \vspace{2 mm}
  $c_6 = (A_{0,1,2,3,4,5}+ B_1 A_{2,3,4,5}+ B_2 A_{0,3,4,5}+ B_3 A_{0,1,4,5} + B_4 A_{0,1,2,5}+ B_5 A_{0,1,2,3} $ \\
  \vspace{2 mm}
  \hspace{0.8 cm} $+ B_{1,3} A_{4,5} + B_{1,4} A_{2,5} + B_{2,4} A_{0,5} + B_{1,5} A_{2,3} + B_{2,5} A_{0,3} + B_{3,5} A_{0,1} + B_{1,3,5}) c_0 $ \\     
  \vspace{2 mm}
  $c_7 = (A_{0,1,2,3,4,5,6}+ B_1 A_{2,3,4,5,6}+ B_2 A_{0,3,4,5,6}+ B_3 A_{0,1,4,5,6} + B_4 A_{0,1,2,5,6}+B_5 A_{0,1,2,3,6} $ \\
  \vspace{2 mm}
  \hspace{0.8 cm} $+ B_6 A_{0,1,2,3,4}+ B_{1,4} A_{2,5,6} + B_{1,3} A_{4,5,6}+ B_{2,4} A_{0,5,6}+ B_{1,5} A_{2,3,6} +B_{2,5} A_{0,3,6} $\\
  \vspace{2 mm}
  \hspace{0.8 cm} $+ B_{3,5} A_{0,1,6} + B_{1,6} A_{2,3,4}+ B_{2,6} A_{0,3,4} + B_{3,6} A_{0,1,4} + B_{4,6} A_{0,1,2}$\\
  \vspace{2 mm}                 
  \hspace{0.8 cm} $+ B_{1,3,6} A_4 + B_{1,3,5} A_6 + B_{1,4,6} A_2 + B_{2,4,6} A_0 ) c_0 $ \\
  \vspace{2 mm}
  $c_8 = (A_{0,1,2,3,4,5,6,7}+B_1 A_{2,3,4,5,6,7} +B_2 A_{0,3,4,5,6,7} + B_3 A_{0,1,4,5,6,7} + B_4 A_{0,1,2,5,6,7}+ B_5 A_{0,1,2,3,6,7} $ \\
  \vspace{2 mm}
  \hspace{0.8 cm} $+B_6 A_{0,1,2,3,4,7}+ B_{7} A_{0,1,2,3,4,5} + B_{1,4} A_{2,5,6,7} + B_{2,4} A_{0,5,6,7}+ B_{1,5} A_{2,3,6,7}+ B_{1,3} A_{4,5,6,7}  $\\
  \vspace{2 mm}
  \hspace{0.8 cm} $+ B_{2,5} A_{0,3,6,7} + B_{3,5} A_{0,1,6,7} + B_{1,6} A_{2,3,4,7} + B_{2,6} A_{0,3,4,7} + B_{3,6} A_{0,1,4,7} + B_{4,6} A_{0,1,2,7} $\\
  \vspace{2 mm}                 
  \hspace{0.8 cm} $+ B_{1,7} A_{2,3,4,5} +B_{3,7} A_{0,1,4,5}+ B_{2,7} A_{0,3,4,5}+ B_{4,7} A_{0,1,2,5}+ B_{5,7} A_{0,1,2,3} $\\
   \vspace{2 mm}
   \hspace{0.8 cm} $+ B_{1,4,6} A_{2,7} + B_{2,4,6} A_{0,7} +B_{1,3,7} A_{4,5}+ B_{1,4,7} A_{2,5}+ B_{1,3,6} A_{4,7} +B_{2,4,7} A_{0,5}$ \\
   \vspace{2 mm}
   \hspace{0.8 cm} $+ B_{1,3,5} A_{6,7}+ B_{1,5,7} A_{2,3}+ B_{2,5,7} A_{0,3}+ B_{3,5,7} A_{0,1} + B_{1,3,5,7}) c_0$\\ 
\hspace{2 mm}\large{\vdots} \hspace{5cm}\large{\vdots}\\ 
\end{tabular}
\label{eq:7}
\end{equation}
In Ref.\cite{Chou2012} I construct the Frobenius solution for infinite series and polynomial which makes $B_n$ term terminated in three term recurrence relation of linear ordinary differential equation by letting $A_n$ in sequence $c_n$ is the leading term in the analytic function $y(x)$. In this chapter I construct the power series for infinite series and polynomial which makes $A_n$ term terminated in three term recurrence relation by letting $B_n$ in the sequence $c_n$ is the leading term in $y(x)$.

In (\ref{eq:7}) the number of individual sequence $c_n$ follows Fibonacci sequence: 1,1,2,3,5,8,13,21,34,55,$\cdots$.
The sequence $c_n$ consists of combinations $A_n$ and $B_n$ in (\ref{eq:7}). 
First observe the term inside parentheses of sequence $c_n$ which does not include any $B_n$'s: $c_n$ with every index ($c_0$, $c_1$, $c_2$,$\cdots$). 

(a) Zero term of $B_n$'s

\begin{equation}
\begin{tabular}{  l  }
  \vspace{2 mm}
  $c_0$ \\
  \vspace{2 mm}
  $c_1 = A_0 c_0  $ \\
  \vspace{2 mm}
  $c_2 = A_{0,1} c_0  $ \\
  \vspace{2 mm}
  $c_3 = A_{0,1,2}c_0 $ \\
  \vspace{2 mm}
  $c_4 = A_{0,1,2,3}c_0 $\\
  \hspace{2 mm}
  \large{\vdots}\hspace{1cm}\large{\vdots} \\ 
\end{tabular}
\label{eq:8}
\end{equation}

When a function $y(x)$, analytic at $x=0$, is expanded in a power series, we write
\begin{equation}
y(x)= \sum_{n=0}^{\infty } c_n x^{n+\lambda }= \sum_{m=0}^{\infty } y_m(x) = y_0(x)+ y_1(x)+y_2(x)+ \cdots \label{eq:9}
\end{equation}
where
\begin{equation}
y_m(x)= \sum_{l=0}^{\infty } c_l^m x^{l+\lambda }\label{eq:104}
\end{equation}
$\lambda $ is the indicial root. $y_m(x)$ is sub-power series that have sequence $c_n$ including $m$ term of $B_n$'s in (\ref{eq:7}). For example $y_0(x)$ has sequences $c_n$ including zero term of $B_n$'s in (\ref{eq:7}), $y_1(x)$ has sequences $c_n$ including one term of $B_n$'s in (\ref{eq:7}), $y_2(x)$ has sequences $c_n$ including two term of $B_n$'s in (\ref{eq:7}), etc. Substitute (\ref{eq:8}) in (\ref{eq:104}) putting $m = 0$. 
\begin{equation}
y_0(x) = c_0 \sum_{n=0}^{\infty} \left\{ \prod _{i_0=0}^{n-1}A_{i_0} \right\} x^{n+\lambda }\label{eq:10}
\end{equation}
Observe the terms inside parentheses of sequence $c_n$ which include one term of $B_n$'s in (\ref{eq:7}): $c_n$ with every index except $c_0$ and $c_1$ ($c_2$, $c_3$, $c_4$,$\cdots$). 

(b) One term of $B_n$'s

\begin{equation}
\begin{tabular}{  l  }
  \vspace{2 mm}
  $c_2= B_1 c_0$ \\
  \vspace{2 mm}
  $c_3 = \left\{ B_1 \cdot 1\cdot \Big( \frac{A_2}{1}\Big) + B_2 A_0\Big(\frac{A_2}{A_2}\Big) \right\} c_0  $ \\
  \vspace{2 mm}
  $c_4 = \left\{ B_1 \cdot 1\cdot \Big( \frac{A_{2,3}}{1}\Big) + B_2 A_0 \Big(\frac{A_{2,3}}{A_2}\Big) + B_3 A_{0,1} \Big(\frac{A_{2,3}}{A_{2,3}}\Big) \right\} c_0  $ \\
  \vspace{2 mm}
  $c_5 = \left\{ B_1 \cdot 1\cdot \Big( \frac{A_{2,3,4}}{1}\Big)+ B_2 A_0 \Big(\frac{A_{2,3,4}}{A_2}\Big) + B_3 A_{0,1}\Big(\frac{A_{2,3,4}}{A_{2,3}}\Big) +  B_4 A_{0,1,2}\Big(\frac{A_{2,3,4}}{A_{2,3,4}}\Big) \right\} c_0  $ \\
  \vspace{2 mm}
  $c_6 = \bigg\{ B_1 \cdot 1\cdot \Big( \frac{A_{2,3,4,5}}{1}\Big)+ B_2 A_0 \Big(\frac{A_{2,3,4,5}}{A_2}\Big) + B_3 A_{0,1} \Big(\frac{A_{2,3,4,5}}{A_{2,3}}\Big) +  B_4 A_{0,1,2}\Big(\frac{A_{2,3,4,5}}{A_{2,3,4}}\Big)  $\\
  \vspace{2 mm}
  \hspace{0.8 cm} $+  B_5 A_{0,1,2,3}\Big(\frac{A_{2,3,4,5}}{A_{2,3,4,5}}\Big) \bigg\} c_0 $\\
  $c_7 =\bigg\{ B_1 \cdot 1\cdot  \Big( \frac{A_{2,3,4,5,6}}{1}\Big)+ B_2 A_0 \Big(\frac{A_{2,3,4,5,6}}{A_2}\Big) + B_3 A_{0,1}\Big(\frac{A_{2,3,4,5,6}}{A_{2,3}}\Big) +  B_4 A_{0,1,2}\Big(\frac{A_{2,3,4,5,6}}{A_{2,3,4}}\Big) $\\
   \vspace{2 mm}
   \hspace{0.8 cm} $+  B_5 A_{0,1,2,3} \Big(\frac{A_{2,3,4,5,6}}{A_{2,3,4,5}}\Big) +  B_6 A_{0,1,2,3,4}\Big(\frac{A_{2,3,4,5,6}}{A_{2,3,4,5,6}}\Big) \bigg\} c_0 $\\
  \hspace{2 mm}
  \large{\vdots}\hspace{3cm}\large{\vdots} \\ 
\end{tabular}
\label{eq:11}
\end{equation}

(\ref{eq:11}) is simply
\begin{equation}
c_{n+2}= c_0 \sum_{i_0=0}^{n} \left\{ B_{i_0+1} \prod _{i_1=0}^{i_0-1}A_{i_1} \prod _{i_2=i_0}^{n-1}A_{i_2+2} \right\}  
\label{eq:12}
\end{equation}
Substitute (\ref{eq:12}) in (\ref{eq:104}) putting $m = 1$. 
\begin{equation}
y_1(x)= c_0 \sum_{n=0}^{\infty}\left\{ \sum_{i_0=0}^{n} \left\{ B_{i_0+1} \prod _{i_1=0}^{i_0-1}A_{i_1} \prod _{i_2=i_0}^{n-1}A_{i_2+2} \right\} \right\} x^{n+2+\lambda } \label{eq:13}
\end{equation}
Observe the terms inside parentheses of sequence $c_n$ which include two terms of $B_n$'s in (\ref{eq:7}): $c_n$ with every index except $c_0$--$c_3$ ($c_4$, $c_5$, $c_6$,$\cdots$). 

(c) Two terms of $B_n$'s
\begin{equation}
\begin{tabular}{  l  }
  \vspace{2 mm}
  $c_4= B_{1,3} c_0$ \\
  \vspace{2 mm}
  $c_5 = \Bigg\{ B_1 \cdot 1\cdot \left\{ B_3 \left( \frac{1}{1}\right) \left(\frac{A_4}{1}\right) + B_4 \left(\frac{A_2}{1}\right) \left( \frac{A_4}{A_4}\right)\right\} + B_2 A_0\left\{ B_4 \left( \frac{A_2}{A_2}\right) \left( \frac{A_4}{A_4}\right) \right\} \Bigg\} c_0  $ \\
  \vspace{2 mm}
  $c_6 = \Bigg\{ B_1 \cdot 1\cdot \left\{ B_3 \left( \frac{1}{1}\right) \left(\frac{A_{4,5}}{1}\right)+ B_4 \left(\frac{A_2}{1}\right) \left( \frac{A_{4,5}}{A_4}\right)+ B_5 \left(\frac{A_{2,3}}{1}\right) \left( \frac{A_{4,5}}{A_{4,5}}\right)\right\}$\\
  \vspace{2 mm}
  \hspace{0.8 cm} $+ B_2 A_0\left\{ B_4 \left( \frac{A_2}{A_2}\right) \left( \frac{A_{4,5}}{A_4}\right) + B_5 \left( \frac{A_{2,3}}{A_2}\right) \left( \frac{A_{4,5}}{A_{4,5}}\right) \right\} + B_3  A_{0,1} \left\{ B_5 \left( \frac{A_{2,3}}{A_{2,3}}\right) \left( \frac{A_{4,5}}{A_{4,5}}\right)\right\} \Bigg\} c_0 $ \\
\vspace{2 mm}
    $c_7 = \Bigg\{ B_1\cdot  1 \cdot \bigg\{ B_3 \left( \frac{1}{1}\right) \left(\frac{A_{4,5,6}}{1}\right) + B_4 \left(\frac{A_2}{1}\right) \left( \frac{A_{4,5,6}}{A_4}\right) + B_5 \left(\frac{A_{2,3}}{1}\right) \left( \frac{A_{4,5,6}}{A_{4,5}}\right) $\\
\vspace{2 mm}
  \hspace{0.8 cm} $+  B_6 \left(\frac{A_{2,3,4}}{1}\right) \left( \frac{A_{4,5,6}}{A_{4,5,6}}\right) \bigg\}$\\
\vspace{2 mm}
  \hspace{0.8 cm} $+ B_2 A_0\left\{ B_4 \left( \frac{A_2}{A_2}\right) \left( \frac{A_{4,5,6}}{A_4}\right)  + B_5 \left( \frac{A_{2,3}}{A_2}\right) \left( \frac{A_{4,5,6}}{A_{4,5}}\right) + B_6 \left( \frac{A_{2,3,4}}{A_2}\right) \left( \frac{A_{4,5,6}}{A_{4,5,6}}\right) \right\}$\\
\vspace{2 mm}
  \hspace{0.8 cm} $+ B_3 A_{0,1} \left\{ B_5 \left( \frac{A_{2,3}}{A_{2,3}}\right) \left( \frac{A_{4,5,6}}{A_{4,5}}\right) +  B_6 \left( \frac{A_{2,3,4}}{A_{2,3}}\right) \left( \frac{A_{4,5,6}}{A_{4,5,6}}\right) \right\}$\\
\vspace{2 mm}
  \hspace{0.8 cm} $+ B_4 A_{0,1,2}\left\{ B_6 \Big( \frac{A_{2,3,4}}{A_{2,3,4}}\Big) \Big( \frac{A_{4,5,6}}{A_{4,5,6}}\Big)  \right\} \Bigg\} c_0$\\
  \hspace{2 mm}
  \large{\vdots}\hspace{5cm}\large{\vdots} \\  
\end{tabular}
\label{eq:14}
\end{equation}
(\ref{eq:14}) is simply
\begin{eqnarray}
 c_{n+4} &=& c_0 \sum_{i_0=0}^{n} \left\{ B_{i_0+1}\sum_{i_1=i_0}^{n} \left\{ B_{i_1+3} \prod _{i_2=0}^{i_0-1}A_{i_2} \prod _{i_3=i_0}^{i_1-1}A_{i_3+2}\prod _{i_4=i_1}^{n-1}A_{i_4+4} \right\}\right\}  
\nonumber\\
 &=&  c_0 \sum_{i_0=0}^{n} \left\{ B_{i_0+1}\prod _{i_1=0}^{i_0-1}A_{i_1} \sum_{i_2=i_0}^{n} \left\{ B_{i_2+3}  \prod _{i_3=i_0}^{i_2-1}A_{i_3+2}\prod _{i_4=i_2}^{n-1}A_{i_4+4} \right\}\right\} \label{eq:15}
\end{eqnarray}
Substitute (\ref{eq:15}) in (\ref{eq:104}) putting $m = 2$. 
\begin{eqnarray}
 y_2(x) &=&  c_0 \sum_{n=0}^{\infty}\left\{ \sum_{i_0=0}^{n} \left\{ B_{i_0+1}\prod _{i_1=0}^{i_0-1}A_{i_1} \sum_{i_2=i_0}^{n} \left\{ B_{i_2+3}  \prod _{i_3=i_0}^{i_2-1}A_{i_3+2} \right.\right.\right.\nonumber\\
&&\times \left.\left.\left. \prod _{i_4=i_2}^{n-1}A_{i_4+4} \right\}\right\} \right\} x^{n+4+\lambda } \label{eq:16}\
\end{eqnarray}
Observe the terms inside parentheses of sequence $c_n$ which include three terms of $B_n$'s in (\ref{eq:7}): $c_n$ with every index except $c_0$--$c_5$ ($c_6$, $c_7$, $c_8$,$\cdots$). 

(d) Three terms of $B_n$'s
\begin{equation}
\begin{tabular}{  l  }
  \vspace{2 mm}
  $c_6= B_{1,3,5} \;c_0$ \\
  \vspace{2 mm}
  $c_7 = \Bigg\{ B_1 \Big\{ B_3 \cdot 1\cdot \Big[ B_5 \left( \frac{1}{1}\right) \left( \frac{1}{1}\right) \left(\frac{A_{6}}{1}\right)+ B_6 \left( \frac{1}{1}\right) \left(\frac{A_4}{1}\right) \left( \frac{A_{6}}{A_{6}}\right)\Big] $\\
 \vspace{2 mm}
  \hspace{0.8 cm} $+ B_4 \cdot 1 \cdot \Big[ B_6 \left( \frac{A_2}{1}\right) \left( \frac{A_{4}}{A_{4}}\right) \left( \frac{A_{6}}{A_{6}}\right) \Big] \Big\} + B_2 \Big\{ B_4 A_0 \Big[ B_6 \left( \frac{A_{2}}{A_{2}}\right)\left( \frac{A_{4}}{A_{4}}\right)\left( \frac{A_{6}}{A_{6}}\right)\Big] \Big\} \Bigg\} c_0  $ \\
  \vspace{2 mm}
  $c_8 = \Bigg\{ B_1 \Big\{ B_3 \cdot 1\cdot  \Big[ B_5 \left( \frac{1}{1}\right) \left( \frac{1}{1}\right) \left(\frac{A_{6,7}}{1}\right)+ B_6 \left( \frac{1}{1}\right) \left(\frac{A_4}{1}\right) \left( \frac{A_{6,7}}{A_{6}}\right)+ B_7 \left( \frac{1}{1}\right) \left(\frac{A_{4,5}}{1}\right) \left( \frac{A_{6,7}}{A_{{6,7}}}\right)\Big]$\\
 \vspace{2 mm}
  \hspace{0.8 cm} $+ B_4 \cdot 1 \cdot \Big[ B_6 \left( \frac{A_2}{1}\right) \left( \frac{A_{4}}{A_{4}}\right) \left( \frac{A_{6,7}}{A_{6}}\right) + B_7 \left( \frac{A_2}{1}\right) \left( \frac{A_{4,5}}{A_{4}}\right) \left( \frac{A_{6,7}}{A_{6,7}}\right)\Big]+ B_5 \cdot 1 \cdot \Big[ B_7 \left( \frac{A_{2,3}}{1}\right) \left( \frac{A_{4,5}}{A_{4,5}}\right) \left( \frac{A_{6,7}}{A_{6,7}}\right) \Big] \Big\} $\\
 \vspace{2 mm}
  \hspace{0.8 cm} $+ B_2 \Big\{ B_4  A_0\Big[ B_6 \left( \frac{A_{2}}{A_{2}}\right)\left( \frac{A_{4}}{A_{4}}\right)\left( \frac{A_{6,7}}{A_{6}}\right)+ B_7 \left( \frac{A_{2}}{A_{2}}\right)\left( \frac{A_{4,5}}{A_{4}}\right)\left( \frac{A_{6,7}}{A_{6,7}}\right)\Big] $\\\vspace{2 mm}
  \hspace{0.8 cm}$ + B_5 A_0\Big[ B_7 \left( \frac{A_{2,3}}{A_{2}}\right) \left( \frac{A_{4,5}}{A_{4,5}}\right) \left( \frac{A_{6,7}}{A_{6,7}}\right) \Big] \Big\}+ B_3 \Big\{ B_5 A_{0,1} \Big[ B_7 \left( \frac{A_{2,3}}{A_{2,3}}\right) \left( \frac{A_{4,5}}{A_{4,5}}\right) \left( \frac{A_{6,7}}{A_{6,7}}\right) \Big]\Big\} \Bigg\} c_0 $ \\
 \vspace{2 mm}
   $c_9 = \Bigg\{ B_1 \bigg\{ B_3 \cdot 1 \cdot \Big[ B_5 \left( \frac{1}{1}\right) \left( \frac{1}{1}\right) \left(\frac{A_{6,7,8}}{1}\right)+ B_6 \left( \frac{1}{1}\right) \left(\frac{A_4}{1}\right) \left( \frac{A_{6,7,8}}{A_{6}}\right)+ B_7 \left( \frac{1}{1}\right) \left(\frac{A_{4,5}}{1}\right) \left( \frac{A_{6,7,8}}{A_{{6,7}}}\right) $\\
\vspace{2 mm}
\hspace{0.8 cm} $+ B_8 \left( \frac{1}{1}\right) \left(\frac{A_{4,5,6}}{1}\right) \left( \frac{A_{6,7,8}}{A_{{6,7,8}}}\right)\Big]$\\
\vspace{2 mm}
\hspace{0.8 cm} $+ B_4 \cdot 1\cdot  \Big[ B_6 \left( \frac{A_2}{1}\right) \left( \frac{A_{4}}{A_{4}}\right) \left( \frac{A_{6,7,8}}{A_{6}}\right) + B_7 \left( \frac{A_2}{1}\right) \left( \frac{A_{4,5}}{A_{4}}\right) \left( \frac{A_{6,7,8}}{A_{6,7}}\right)+ B_8 \left( \frac{A_2}{1}\right) \left( \frac{A_{4,5,6}}{A_{4}}\right) \left( \frac{A_{6,7,8}}{A_{6,7,8}}\right) \Big]$\\
  \vspace{2 mm}
  \hspace{0.8 cm} $+ B_5 \cdot 1 \cdot \Big[ B_7 \left( \frac{A_{2,3}}{1}\right) \left( \frac{A_{4,5}}{A_{4,5}}\right) \left( \frac{A_{6,7,8}}{A_{6,7}}\right)+ B_8 \left( \frac{A_{2,3}}{1}\right) \left( \frac{A_{4,5,6}}{A_{4,5}}\right) \left( \frac{A_{6,7,8}}{A_{6,7,8}}\right)\Big] $\\
 \vspace{2 mm}
  \hspace{0.8 cm} $+ B_6 \cdot 1 \cdot \Big[ B_8 \left( \frac{A_{2,3,4}}{1}\right) \left( \frac{A_{4,5,6}}{A_{4,5,6}}\right) \left( \frac{A_{6,7,8}}{A_{6,7,8}}\right) \Big] \bigg\} $\\
 \vspace{2 mm}
  \hspace{0.8 cm} $+ B_2 \bigg\{ B_4 A_0 \Big[ B_6  \left( \frac{A_{2}}{A_{2}}\right)\left( \frac{A_{4}}{A_{4}}\right)\left( \frac{A_{6,7,8}}{A_{6}}\right)+ B_7 \left( \frac{A_{2}}{A_{2}}\right)\left( \frac{A_{4,5}}{A_{4}}\right)\left( \frac{A_{6,7,8}}{A_{6,7}}\right) + B_8 \left( \frac{A_{2}}{A_{2}}\right)\left( \frac{A_{4,5,6}}{A_{4}}\right)\left( \frac{A_{6,7,8}}{A_{6,7,8}}\right)\Big]$\\
  \vspace{2 mm}
  \hspace{0.8 cm} $+ B_5 A_0 \Big[ B_7 \left( \frac{A_{2,3}}{A_{2}}\right) \left( \frac{A_{4,5}}{A_{4,5}}\right) \left( \frac{A_{6,7,8}}{A_{6,7}}\right) + B_8 \left( \frac{A_{2,3}}{A_{2}}\right) \left( \frac{A_{4,5,6}}{A_{4,5}}\right) \left( \frac{A_{6,7,8}}{A_{6,7,8}}\right)\Big] $\\
 \vspace{2 mm}
  \hspace{0.8 cm} $+ B_6 A_0 \Big[ B_8 \left( \frac{A_{2,3,4}}{A_{2}}\right) \left( \frac{A_{4,5,6}}{A_{4,5,6}}\right) \left( \frac{A_{6,7,8}}{A_{6,7,8}}\right) \Big] \bigg\}$\\
  \vspace{2 mm}
  \hspace{0.8 cm} $+ B_3 \bigg\{ B_5  A_{0,1}\Big[ B_7 \left( \frac{A_{2,3}}{A_{2,3}}\right) \left( \frac{A_{4,5}}{A_{4,5}}\right) \left( \frac{A_{6,7,8}}{A_{6,7}}\right) + B_8 \left( \frac{A_{2,3}}{A_{2,3}}\right) \left( \frac{A_{4,5,6}}{A_{4,5}}\right) \left( \frac{A_{6,7,8}}{A_{6,7,8}}\right)\Big] $\\
\end{tabular}
\nonumber
\end{equation}
\begin{equation}
\begin{tabular}{  l  }  
\vspace{2 mm}
 \hspace{0.8 cm} $+ B_6 A_{0,1}\Big[ B_8 \left( \frac{A_{2,3,4}}{A_{2,3}}\right) \left( \frac{A_{4,5,6}}{A_{4,5,6}}\right) \left( \frac{A_{6,7,8}}{A_{6,7,8}}\right) \Big] \bigg\}+ B_4 \bigg\{ B_6 A_{0,1,2} \Big[ B_8 \left( \frac{A_{2,3,4}}{A_{2,3,4}}\right) \left( \frac{A_{4,5,6}}{A_{4,5,6}}\right) \left( \frac{A_{6,7,8}}{A_{6,7,8}}\right)\Big] \bigg\} \Bigg\} c_0 $ \\
  \hspace{2 mm}
  \large{\vdots}\hspace{6cm}\large{\vdots} \\  
\end{tabular}
\label{eq:17}
\end{equation}
(\ref{eq:17}) is simply
\begin{eqnarray}
 c_{n+6} &=& c_0 \sum_{i_0=0}^{n} \left\{ B_{i_0+1}\sum_{i_1=i_0}^{n} \left\{ B_{i_1+3}\sum_{i_2=i_1}^{n} \left\{ B_{i_2+5} \prod _{i_3=0}^{i_0-1}A_{i_3} \prod _{i_4=i_0}^{i_1-1}A_{i_4+2} \prod _{i_5=i_1}^{i_2-1}A_{i_5+4}\prod _{i_6=i_2}^{n-1}A_{i_6+6} \right\} \right\} \right\} \nonumber\\
 &=& c_0 \sum_{i_0=0}^{n} \left\{ B_{i_0+1}\prod _{i_1=0}^{i_0-1}A_{i_1} \sum_{i_2=i_0}^{n} \left\{ B_{i_2+3}\prod _{i_3=i_0}^{i_2-1}A_{i_3+2} \sum_{i_4=i_2}^{n} \left\{ B_{i_4+5} \prod _{i_5=i_2}^{i_4-1}A_{i_5+4} \prod _{i_6=i_4}^{n-1}A_{i_6+6} \right\} \right\} \right\} \hspace{1.5cm}\label{eq:18}
\end{eqnarray}
Substitute (\ref{eq:18}) in (\ref{eq:104}) putting $m = 3$. 
\begin{eqnarray}
 y_3(x) &=& c_0 \sum_{n=0}^{\infty}\left\{\sum_{i_0=0}^{n} \left\{ B_{i_0+1}\prod _{i_1=0}^{i_0-1}A_{i_1} \sum_{i_2=i_0}^{n} \left\{ B_{i_2+3}\prod _{i_3=i_0}^{i_2-1}A_{i_3+2} \right.\right.\right. \nonumber\\
&&\times \left.\left.\left. \sum_{i_4=i_2}^{n} \left\{ B_{i_4+5} \prod _{i_5=i_2}^{i_4-1}A_{i_5+4} \prod _{i_6=i_4}^{n-1}A_{i_6+6} \right\} \right\} \right\}\right\} x^{n+6+\lambda } \label{eq:19}
\end{eqnarray}
By repeating this process for all higher terms of $B_n$'s, we obtain every $y_m(x)$ terms where $m \geq 4$. Substitute (\ref{eq:10}), (\ref{eq:13}), (\ref{eq:16}), (\ref{eq:19}) and including all $y_m(x)$ terms where $m \geq 4$ into (\ref{eq:9}). 
\begin{theorem}
The general expression of $y(x)$ for infinite series is
\begin{eqnarray}
 y(x) &=&  y_0(x)+ y_1(x)+ y_2(x)+y_3(x)+\cdots \nonumber\\
&=& c_0 \left\{ \sum_{n=0}^{\infty} \left\{ \prod _{i_0=0}^{n-1}A_{i_0} \right\} x^{n+\lambda} + \sum_{n=0}^{\infty}\left\{ \sum_{i_0=0}^{n} \left\{ B_{i_0+1} \prod _{i_1=0}^{i_0-1}A_{i_1} \prod _{i_2=i_0}^{n-1}A_{i_2+2} \right\} \right\} x^{n+2+\lambda} \right. \nonumber\\
&&+ \sum_{N=2}^{\infty } \left\{\sum_{n=0}^{\infty } \left\{ \sum_{i_0=0}^{n} \left\{B_{i_0+1}\prod _{i_1=0}^{i_0-1} A_{i_1} 
\prod _{k=1}^{N-1} \left( \sum_{i_{2k}= i_{2(k-1)}}^{n} B_{i_{2k}+2k+1}\prod _{i_{2k+1}=i_{2(k-1)}}^{i_{2k}-1}A_{i_{2k+1}+2k}\right) \right.\right.\right. \nonumber\\
&&\times \left.\left.\left.\left. \prod _{i_{2N} = i_{2(N-1)}}^{n-1} A_{i_{2N}+ 2N} \right\}\right\}\right\} x^{n+2N+\lambda} \right\} 
\label{eq:20}
\end{eqnarray}
\end{theorem}
\section{\label{sec:level4}Polynomial which makes $A_n$ term terminated}
Now let's investigate the polynomial case of $y(x)$. Assume that $A_n$ is terminated at certain value of $n$. Then each $y_i(x)$ where $i=0,1,2,\cdots$ will be polynomial. Examples of these are Heun's equation, GCH function\cite{Ch2012}, Lame function, etc. First $A_{k}$ is terminated at certain value of $k$. I choose eigenvalue $\alpha _0$ which $A_{k}$ is terminated where $\alpha _0 =0,1,2,\cdots$. $A_{k+2}$ is terminated at certain value of $k$. I choose eigenvalue $\alpha _1$ which $A_{k+2}$ is terminated where $\alpha _1 =0,1,2,\cdots$. Also $A_{k+4}$ is terminated at certain value of $k$. I choose eigenvalue $\alpha _2$ which $A_{k+4}$ is terminated where $\alpha _2 =0,1,2,\cdots$. By repeating this process I obtain
\begin{equation}
A_{\alpha _i+ 2i}=0 \hspace{1cm} \mathrm{where}\;i,\alpha _i =0,1,2,\cdots
\label{eq:21}
\end{equation}
In general, the two term recurrence relation for polynomial has only one eigenvalue: for example, the Laguerre function, confluent hypergeometric function, Legendre function, etc. 
In ``Generalization of the three-term recurrence formula and its applications''\cite{Chou2012}, polynomial which makes $B_n$ term terminated has infinite eigenvalues which is $\beta _i$ where $i,\beta _i =0,1,2,\cdots$. 
In this chapter polynomial which makes $A_n$ term terminated also has infinite eigenvalues which is $\alpha _i$.

First observe the term in sequence $c_n$ which does not include any $B_n$'s in (\ref{eq:8}): $c_n$ with every index ($c_0$, $c_1$, $c_2$,$\cdots$). 

(a) As $\alpha _0$=0, then $A_0$=0 in (\ref{eq:8}).
\begin{equation}
\begin{tabular}{  l  }
  \vspace{2 mm}
  $c_0$ \\ \hspace{1.6cm}
\end{tabular}
\label{eq:22}
\end{equation}
(b) As $\alpha _0$=1, then $A_1$=0 in (\ref{eq:8}).
\begin{equation}
\begin{tabular}{  l  }
  \vspace{2 mm}
  $c_0$ \\
  \vspace{2 mm}
  $c_1 = A_0 c_0  $ \\
\end{tabular}
\label{eq:23}
\end{equation}
(c) As $\alpha _0$=2, then $A_2$=0 in (\ref{eq:8}).
\begin{equation}
\begin{tabular}{  l  }
  \vspace{2 mm}
  $c_0$ \\
  \vspace{2 mm}
  $c_1 = A_0 c_0  $ \\
  \vspace{2 mm}
  $c_2 = A_{0,1} c_0  $ \\
\end{tabular}
\label{eq:24}
\end{equation}
Substitute (\ref{eq:22}),(\ref{eq:23}) and (\ref{eq:24}) in (\ref{eq:104}) putting $m = 0$.
\begin{equation}
y_0(x)= c_0 \sum_{n=0}^{\alpha _0} \left\{ \prod _{i_1=0}^{n-1}A_{i_1} \right\} x^{n+\lambda }
\label{eq:25}
\end{equation} 
Observe the terms inside curly brackets of sequence $c_n$ which include one term of $B_n$'s in (\ref{eq:11}): $c_n$ with every index except $c_0$ and $c_1$ ($c_1$, $c_3$, $c_5$,$\cdots$). 

(a) As $\alpha _0$=0, then $A_0$=0 in (\ref{eq:11}).
\begin{equation}
\begin{tabular}{  l  }
  \vspace{2 mm}
  $c_2= B_1 c_0$ \\
  \vspace{2 mm}
  $c_3 = B_1 \cdot 1\cdot  A_2 c_0  $ \\
  \vspace{2 mm}
  $c_4 = B_1 \cdot 1\cdot  A_{2,3}c_0 $ \\
  \vspace{2 mm}
  $c_5 = B_1 \cdot 1\cdot  A_{2,3,4}c_0  $ \\
  \vspace{2 mm}
  $c_6 = B_1 \cdot 1 \cdot A_{2,3,4,5}c_0  $ \\
  \hspace{2 mm}
  \large{\vdots}\hspace{1cm}\large{\vdots} \\ 
\end{tabular}
\label{eq:26}
\end{equation}
As i=1 in (\ref{eq:21}),
\begin{equation}
A_{\alpha _1 +2}=0 \hspace{1cm} \mathrm{where}\; \alpha _1 =0,1,2,\cdots
\label{eq:27}
\end{equation}
Substitute (\ref{eq:26}) in (\ref{eq:104}) putting $m = 1$ by using (\ref{eq:27}).
\begin{equation}
y_1^0(x)= c_0 B_1 \sum_{n=0}^{\alpha _1} \left\{ \prod _{i_1=0}^{n-1}A_{i_1+2} \right\} x^{n+2+\lambda }  
\label{eq:28}
\end{equation}
In (\ref{eq:28}) $y_1^0(x)$ is sub-power series, having sequences $c_n$ including one term of $B_n$'s in (\ref{eq:7}) as $\alpha _0$=0, for the polynomial case in which makes $A_n$ term terminated.

(b) As $\alpha _0$=1, then $A_1$=0 in (\ref{eq:11}).
\begin{equation}
\begin{tabular}{  l  }
  \vspace{2 mm}
  $c_2= B_1 c_0$ \\
  \vspace{2 mm}
  $c_3 = \{B_1 \cdot 1\cdot  A_2 + B_2 A_0 \cdot 1 \} c_0  $ \\
  \vspace{2 mm}
  $c_4 = \{B_1 \cdot 1 \cdot A_{2,3} + B_2 A_0 A_3 \} c_0 $ \\
  \vspace{2 mm}
  $c_5 = \{B_1 \cdot 1\cdot  A_{2,3,4} + B_2 A_0 A_{3,4} \} c_0  $ \\
  \vspace{2 mm}
  $c_6 = \{B_1 \cdot 1\cdot  A_{2,3,4,5} + B_2 A_0 A_{3,4,5} \} c_0  $ \\
  \hspace{2 mm}
  \large{\vdots}\hspace{2cm}\large{\vdots} \\ 
\end{tabular}
\label{eq:29}
\end{equation}
The first term in curly brackets of sequence $c_n$ in (\ref{eq:29}) is same as (\ref{eq:26}). Then, its solution is equal to (\ref{eq:28}). Substitute (\ref{eq:27}) into the second term in curly brackets of sequence $c_n$ in (\ref{eq:29}). Its power series expansion including the first and second terms in curly brackets of sequence $c_n$ in (\ref{eq:29}), analytic at $x=0$, is
\begin{equation}
y_1^1(x)= c_0 \left\{ B_1 \sum_{n=0}^{\alpha _1} \left\{ \prod _{i_1=0}^{n-1}A_{i_1+2} \right\} + B_2 A_0 \sum_{n=1}^{\alpha _1} \left\{ \prod _{i_1=1}^{n-1}A_{i_1+2} \right\}\right\} x^{n+2+\lambda }  
\label{eq:30}
\end{equation}
In (\ref{eq:30}) $y_1^1(x)$ is sub-power series, having sequences $c_n$ including one term of $B_n$'s in (\ref{eq:7}), as $\alpha _0$=1 for the polynomial case in which makes $A_n$ term terminated.

(c) As $\alpha _0$=2, then $A_2$=0 in (\ref{eq:11}).
\begin{equation}
\begin{tabular}{  l  }
  \vspace{2 mm}
  $c_2= B_1 c_0$ \\
  \vspace{2 mm}
  $c_3 = \{B_1 \cdot 1\cdot  A_2 + B_2 A_0 \cdot 1 \} c_0  $ \\
  \vspace{2 mm}
  $c_4 = \{B_1 \cdot 1 \cdot A_{2,3} + B_2 A_0 A_3 + B_3 A_{0,1} \cdot 1 \} c_0 $ \\
  \vspace{2 mm}
  $c_5 = \{B_1 \cdot 1\cdot  A_{2,3,4} + B_2 A_0 A_{3,4} + B_3 A_{0,1} A_4 \} c_0  $ \\
  \vspace{2 mm}
  $c_6 = \{B_1 \cdot 1\cdot  A_{2,3,4,5} + B_2 A_0 A_{3,4,5} + B_3 A_{0,1} A_{4,5} \} c_0  $ \\
  \vspace{2 mm}
  $c_7 = \{B_1 \cdot 1\cdot  A_{2,3,4,5,6} + B_2 A_0 A_{3,4,5,6} + B_3 A_{0,1} B_{4,5,6} \} c_0  $ \\
  \hspace{2 mm}
  \large{\vdots}\hspace{4cm}\large{\vdots} \\ 
\end{tabular}
\label{eq:31}
\end{equation}
The first and second term in curly brackets of sequence $c_n$ in (\ref{eq:31}) is same as (\ref{eq:29}). Then its power series expansion is same as (\ref{eq:30}). Substitute (\ref{eq:27}) into the third term in curly brackets of sequence $c_n$ in (\ref{eq:31}). Its power series expansion including the first, second and third terms in curly brackets of sequence $c_n$ in (\ref{eq:31}), analytic at $x=0$, is
\begin{eqnarray}
y_1^2(x)&=& c_0 \left\{  B_1 \sum_{n=0}^{\alpha _1} \left\{ \prod _{i_1=0}^{n-1}A_{i_1+2} \right\} + B_2 A_0 \sum_{n=1}^{\alpha _1} \left\{ \prod _{i_1=1}^{n-1}A_{i_1+2} \right\}\right.\nonumber\\
&&+\left. B_3 A_{0,1} \sum_{n=2}^{\alpha _1} \left\{ \prod _{i_1=2}^{n-1}A_{i_1+2} \right\}\right\} x^{n+2+\lambda } 
 \label{eq:32}
\end{eqnarray}
In (\ref{eq:32}) $y_1^2(x)$ is sub-power series, having sequence $c_n$ including one term of $B_n$'s in (\ref{eq:7}) as $\alpha _0$=2, for the polynomial case in which makes $A_n$ term terminated.
By repeating this process for all $\alpha _0 =3,4,5,\cdots$, I obtain every $y_1^j(x)$ terms where $j \geq 3$. 
According to (\ref{eq:28}), (\ref{eq:30}), (\ref{eq:32}) and every $y_1^j(x)$ where $j \geq 3$, the general expression of $y_1(x)$ for all $\alpha _0$ is 
\begin{equation}
y_1(x)= c_0 \sum_{i_0=0}^{\alpha _0}\left\{ B_{i_0+1} \prod _{i_1=0}^{i_0-1}A_{i_1}  \sum_{i_2=i_0}^{\alpha _1} \left\{ \prod _{i_3=i_0}^{i_2-1}A_{i_3+2} \right\}\right\} x^{i_2+2+\lambda }  \;\;\mbox{where}\; \alpha _0 \leq \alpha _1
 \label{eq:33}
\end{equation}

Observe the terms of sequence $c_n$ which include two terms of $B_n$'s in (\ref{eq:14}): $c_n$ with every index except $c_0$--$c_3$ ($c_4$, $c_5$, $c_6$,$\cdots$). 

(a) As $\alpha _0$=0, then $A_0$=0 in (\ref{eq:14}).
\begin{equation}
\begin{tabular}{  l  }
  \vspace{2 mm}
  $c_4= B_{1,3} c_0$ \\
  \vspace{2 mm}
  $c_5 = B_1\{B_3 \cdot 1\cdot 1\cdot  A_4 + B_4 \cdot 1 \cdot A_2 \cdot 1 \} c_0  $ \\
  \vspace{2 mm}
  $c_6 = B_1 \{B_3 \cdot 1\cdot 1 \cdot A_{4,5} + B_4 \cdot 1\cdot  A_2 A_5 + B_5 \cdot 1 \cdot A_{2,3} \cdot 1\} c_0 $ \\
  \vspace{2 mm}
  $c_7 = B_1 \{B_3 \cdot 1\cdot 1 \cdot A_{4,5,6} + B_4\cdot  1\cdot  A_2 A_{5,6} + B_5 \cdot 1 \cdot A_{2,3} A_6 + B_6 \cdot 1\cdot  A_{2,3,4} \cdot 1 \} c_0  $ \\
  \vspace{2 mm}
  $c_8 = B_1 \{B_3 \cdot 1\cdot 1 \cdot A_{4,5,6,7} + B_4 \cdot 1\cdot  A_2 A_{5,6,7} + B_5 \cdot 1\cdot  A_{2,3} A_{6,7} + B_6 \cdot 1 \cdot A_{2,3,4} A_7  $\\
\vspace{2 mm}
 \hspace{1 cm}$+ B_7 \cdot 1\cdot  A_{2,3,4,5}\cdot 1\} c_0$ \\
  \large{\vdots}\hspace{5cm}\large{\vdots} \\ 
\end{tabular}
\label{eq:34}
\end{equation}
(i) As $\alpha _1$=0, then $A_2$=0 in (\ref{eq:34}).
\begin{equation}
\begin{tabular}{  l  }
  \vspace{2 mm}
  $c_4= B_{1,3} \;c_0$ \\
  \vspace{2 mm}
  $c_5 = B_1 B_3\cdot 1 \cdot 1 \cdot A_4 \;c_0  $ \\
  \vspace{2 mm}
  $c_6 = B_1 B_3\cdot 1 \cdot 1 \cdot A_{4,5} \;c_0 $ \\
  \vspace{2 mm}
  $c_7 =  B_1 B_3\cdot 1 \cdot 1 \cdot A_{4,5,6} \;c_0  $ \\
  \vspace{2 mm}
  $c_8 =  B_1 B_3\cdot 1 \cdot 1 \cdot A_{4,5,6,7} \;c_0  $ \\
  \hspace{2 mm}
  \large{\vdots}\hspace{2cm}\large{\vdots} \\ 
\end{tabular}
\label{eq:35}
\end{equation}
As i=2 in (\ref{eq:21}),
\begin{equation}
A_{\alpha _2 +4}=0 \hspace{1cm} \mathrm{where}\; \alpha _2 =0,1,2,\cdots
\label{eq:36}
\end{equation}
Substitute (\ref{eq:35}) in (\ref{eq:104}) putting $m = 2$ by using (\ref{eq:36}).
\begin{equation}
y_2^{0,0}(x)= c_0 B_1 B_3 \sum_{n=0}^{\alpha _2} \left\{ \prod _{i_1=0}^{n-1}A_{i_1+4} \right\} x^{n+4+\lambda }  
\label{eq:37}
\end{equation} 
In (\ref{eq:37}) $y_2^{0,0}(x)$ is sub-power series, having sequences $c_n$ including two term of $B_n$'s in (\ref{eq:7}) as $\alpha _0$=0 and $\alpha _1$=0, for the polynomial case in which makes $A_n$ term terminated.

(ii) As $\alpha _1$=1, then $A_3$=0 in (\ref{eq:34}).
\begin{equation}
\begin{tabular}{  l  }
  \vspace{2 mm}
  $c_4= B_{1,3} c_0$ \\
  \vspace{2 mm}
  $c_5 = B_1\{B_3\cdot 1 \cdot 1 \cdot A_4 + B_4 \cdot 1 \cdot A_2 \cdot 1 \} c_0  $ \\
  \vspace{2 mm}
  $c_6 = B_1\{B_3\cdot 1 \cdot 1 \cdot A_{4,5} + B_4 \cdot 1 \cdot A_2 \cdot A_5 \} c_0 $ \\
  \vspace{2 mm}
  $c_7 = B_1\{B_3\cdot 1 \cdot 1 \cdot A_{4,5,6} + B_4 \cdot 1 \cdot A_2 \cdot A_{5,6} \} c_0  $ \\
  \vspace{2 mm}
  $c_8 = B_1\{B_3\cdot 1 \cdot 1 \cdot A_{4,5,6,7} + B_4 \cdot 1 \cdot A_2 \cdot A_{5,6,7} \} c_0  $ \\
  \hspace{2 mm}
  \large{\vdots}\hspace{3cm}\large{\vdots} \\ 
\end{tabular}
\label{eq:38}
\end{equation}
The first term in curly brackets of sequence $c_n$ in (\ref{eq:38}) is same as (\ref{eq:35}). Then its power series expansion is equal to (\ref{eq:37}). Substitute (\ref{eq:36}) into the second term in curly brackets of sequence $c_n$ in (\ref{eq:38}). Its power series expansion including the first and second terms in curly brackets of sequence $c_n$, analytic at $x=0$, is
\begin{equation}
y_2^{0,1}(x)= c_0 B_1 \left\{ B_3 \sum_{n=0}^{\alpha _2} \left\{ \prod _{i_1=0}^{n-1}A_{i_1+4} \right\} + B_4 A_2 \sum_{n=1}^{\alpha _2} \left\{ \prod _{i_1=1}^{n-1}A_{i_1+4} \right\} \right\} x^{n+4+\lambda }
\label{eq:39}
\end{equation}
In (\ref{eq:39}) $y_2^{0,1}(x)$ is sub-power series, having sequences $c_n$ including two term of $B_n$'s in (\ref{eq:7}) as $\alpha _0$=0 and $\alpha _1$=1, for the polynomial case in which makes $A_n$ term terminated.

(iii) As $\alpha _1$=2, then $A_4$=0 in (\ref{eq:34}).
\begin{equation}
\begin{tabular}{  l  }
  \vspace{2 mm}
  $c_4= B_{1,3} c_0$ \\
  \vspace{2 mm}
  $c_5 = B_1\{B_3\cdot 1 \cdot 1 \cdot A_4 + B_4 \cdot 1 \cdot A_2 \cdot 1 \} c_0  $ \\
  \vspace{2 mm}
  $c_6 = B_1\{B_3\cdot 1 \cdot 1 \cdot A_{4,5} + B_4 \cdot 1 \cdot A_2 \cdot A_5 + B_5\cdot 1\cdot A_{2,3}\cdot 1\} c_0 $ \\
  \vspace{2 mm}
  $c_7 = B_1\{B_3\cdot 1 \cdot 1 \cdot A_{4,5,6} + B_4 \cdot 1 \cdot A_2 \cdot A_{5,6} + B_5\cdot 1\cdot A_{2,3}\cdot A_6\} c_0  $ \\
  \vspace{2 mm}
  $c_8 = B_1\{B_3\cdot 1 \cdot 1 \cdot A_{4,5,6,7} + B_4 \cdot 1 \cdot A_2 \cdot A_{5,6,7} + B_5\cdot 1\cdot A_{2,3}\cdot A_{6,7}\} c_0  $ \\
  \hspace{2 mm}
  \large{\vdots}\hspace{3cm}\large{\vdots} \\ 
\end{tabular}
\label{eq:40}
\end{equation}
The first and second term in curly brackets of sequence $c_n$ in (\ref{eq:40}) is same as (\ref{eq:38}). Then its power series expansion is same as (\ref{eq:39}). Substitute (\ref{eq:36}) into the third term in curly brackets of sequence $c_n$ in (\ref{eq:40}). Its power series expansion including the first, second and third terms in curly brackets of sequence $c_n$ in (\ref{eq:40}), analytic at $x=0$, is
\begin{eqnarray}
y_2^{0,2}(x) &=& c_0 B_1 \left\{ B_3 \sum_{n=0}^{\alpha _2} \left\{ \prod _{i_1=0}^{n-1}A_{i_1+4} \right\} + B_4 A_2 \sum_{n=1}^{\alpha _2} \left\{ \prod _{i_1=1}^{n-1}A_{i_1+4} \right\}\right.\nonumber\\
&&+\left. B_5 A_{2,3} \sum_{n=2}^{\alpha _2} \left\{ \prod _{i_1=2}^{n-1}A_{i_1+4} \right\}\right\} x^{n+4+\lambda }
\label{eq:41}
\end{eqnarray}
In (\ref{eq:41}) $y_2^{0,2}(x)$ is sub-power series, having sequences $c_n$ including two term of $B_n$'s in (\ref{eq:7}) as $\alpha _0$=0 and $\alpha _1$=2, for the polynomial case in which makes $A_n$ term terminated.
By repeating this process for all $\alpha _1 =3,4,5,\cdots$, we obtain every $y_2^{0,j}(x)$ terms where $j \geq 3$. 
According to (\ref{eq:37}), (\ref{eq:39}), (\ref{eq:41}) and every $y_2^{0,j}(x)$ where $j \geq 3$, the general expression of $y_2^0(x)$ for the case of $\alpha _0=0$ replacing the index $n$ by $i_0$ is 
\begin{equation}
y_2^0(x)= c_0 B_1 \sum_{i_0=0}^{\alpha _1}\left\{ B_{i_0+3} \prod _{i_1=0}^{i_0-1}A_{i_1+2}  \sum_{i_2=i_0}^{\alpha _2} \left\{ \prod _{i_3=i_0}^{i_2-1}A_{i_3+4} \right\}\right\} x^{i_2+4+\lambda }  
\label{eq:42}
\end{equation}
In (\ref{eq:42}) $y_2^0(x)$ is sub-power series, having sequences $c_n$ including two term of $B_n$'s in (\ref{eq:7}) as $\alpha _0$=0, for the polynomial case in which makes $A_n$ term terminated.
\vspace{2 mm}

(b) As $\alpha _0$=1, then $A_1$=0 in (\ref{eq:14}).
\begin{equation}
\begin{tabular}{  l  }
  \vspace{2 mm}
  $c_4= B_{1,3} c_0$ \\
  \vspace{2 mm}
  $c_5 = \Big\{ B_1\Big[ B_3\cdot 1 \cdot 1 \cdot A_4 + B_4 \cdot 1 \cdot A_2 \cdot 1\Big]+ B_2\Big[ B_4 A_0 \cdot 1\cdot 1 \Big] \Big\} c_0  $ \\
  \vspace{2 mm}
  $c_6 = \Big\{ B_1 \Big[ B_3 \cdot 1 \cdot 1 \cdot A_{4,5} + B_4 \cdot 1 \cdot A_2 A_5 + B_5 \cdot 1\cdot A_{2,3}\cdot 1 \Big]$\\
\vspace{2 mm}
  \hspace{0.8 cm} $+ B_2\Big[ B_4 A_0 \cdot 1\cdot A_5 + B_5 A_0 A_3 \cdot 1 \Big] \Big\} c_0$ \\
  \vspace{2 mm}
  $c_7 = \Big\{ B_1 \Big[  B_3 \cdot 1 \cdot 1 \cdot A_{4,5,6} + B_4 \cdot 1 \cdot A_2 A_{5,6} + B_5 \cdot 1\cdot A_{2,3}\cdot A_6 +B_6 \cdot 1\cdot A_{2,3,4} \cdot 1 \Big] $ \\
 \vspace{2 mm}
  \hspace{0.8 cm} $+B_2 \Big[ B_4 A_0 \cdot 1\cdot A_{5,6} + B_5 A_0 A_3 A_6 + B_6 A_0 A_{3,4} \cdot 1 \Big] \Big\} c_0 $\\
  \vspace{2 mm}
  $c_8 = \Big\{ B_1 \Big[  B_3 \cdot 1 \cdot 1 \cdot A_{4,5,6,7} + B_4 \cdot 1 \cdot A_2 A_{5,6,7} + B_5 \cdot 1\cdot A_{2,3}\cdot A_{6,7} $\\
\vspace{2 mm}
  \hspace{0.8 cm} $+B_6 \cdot 1\cdot A_{2,3,4} \cdot A_7 + B_7 \cdot 1\cdot  A_{2,3,4,5} \cdot 1\Big]$ \\
  \vspace{2 mm}
  \hspace{0.8 cm} $+B_2 \Big[B_4 A_0 \cdot 1\cdot A_{5,6,7} + B_5 A_0 A_3 A_{6,7} + B_6 A_0 A_{3,4} \cdot A_7 + B_7 A_0 A_{3,4,5}\cdot 1 \Big] \Big\} c_0 $\\ 
  \hspace{2 mm}
  \large{\vdots}\hspace{7cm}\large{\vdots} \\ 
\end{tabular}
\label{eq:43}
\end{equation}
The first square brackets including $B_1$ inside curly brackets in sequence $c_n$ in (\ref{eq:43}) is same as (\ref{eq:34}). Then its power series expansion is same as (\ref{eq:42}). Observe the second square brackets including $B_2$ inside curly brackets in sequence $c_n$ in (\ref{eq:43}).

(i) As $\alpha _1$=1, then $A_3$=0 in the second square brackets including $B_2$ inside curly brackets in sequence $c_n$ in (\ref{eq:43}).
\begin{equation}
\begin{tabular}{  l  }
  \vspace{2 mm}
  $c_4= B_{2} A_0 \Big\{ B_4\cdot 1\cdot 1 \Big\} c_0$ \\
  \vspace{2 mm}
  $c_5 =B_{2} A_0 \Big\{ B_4\cdot 1\cdot A_5 \Big\}c_0 $ \\
  \vspace{2 mm}
  $c_6 = B_{2} A_0 \Big\{ B_4\cdot 1\cdot A_{5,6} \Big\} c_0  $ \\
  \vspace{2 mm}
  $c_7 = B_{2} A_0 \Big\{ B_4\cdot 1\cdot A_{5,6,7} \Big\} c_0  $ \\
  \hspace{2 mm}
  \large{\vdots}\hspace{2cm}\large{\vdots} \\ 
\end{tabular}
\label{eq:44}
\end{equation}
Its power series expansion of (\ref{eq:44}) by using (\ref{eq:36}), analytic at $x=0$, is
\begin{equation}
y_2^{1,1}(x)= c_0 A_0 B_{2,4} \sum_{n=1}^{\alpha _2}\left\{ \prod _{i_1=1}^{n-1}A_{i_1+4} \right\} x^{n+4+\lambda } 
\label{eq:45}
\end{equation}
In (\ref{eq:45}) $y_2^{1,1}(x)$ is sub-power series, for the second square brackets inside curly brackets in sequence $c_n$ including two term of $B_n$'s in (\ref{eq:43}) as $\alpha _0$=1 and $\alpha _1$=1, for the polynomial case in which makes $A_n$ term terminated.

(ii) As $\alpha _1$=2, then $A_4$=0 in second square brackets inside curly brackets in sequence $c_n$ including $B_2$ in (\ref{eq:43}).
\begin{equation}
\begin{tabular}{  l  }
  \vspace{2 mm}
  $c_5= B_{2} A_0 \Big\{ B_4\cdot 1\cdot 1 \Big\} c_0$ \\
  \vspace{2 mm}
  $c_6 =B_{2} A_0 \Big\{ B_4\cdot 1\cdot A_5 + B_5 A_3 \cdot 1 \Big\}c_0 $ \\
  \vspace{2 mm}
  $c_7 = B_{2} A_0 \Big\{ B_4\cdot 1\cdot A_{5,6} + B_5 A_3 \cdot A_6 \Big\} c_0  $ \\
  \vspace{2 mm}
  $c_8 = B_{2} A_0 \Big\{ B_4\cdot 1\cdot A_{5,6,7} + B_5 A_3 \cdot A_{6,7} \Big\} c_0  $ \\
  \hspace{2 mm}
  \large{\vdots}\hspace{3cm}\large{\vdots} \\ 
\end{tabular}
\label{eq:46}
\end{equation}
The first term in curly brackets of sequence $c_n$ in (\ref{eq:46}) is same as (\ref{eq:44}). Then its power series expansion is same as (\ref{eq:45}). Substitute (\ref{eq:36}) into the second term in curly brackets of sequence $c_n$ in (\ref{eq:46}). Its power series expansion including the first and second terms in curly brackets of sequence $c_n$ in (\ref{eq:46}), analytic at $x=0$, is
\begin{equation}
y_2^{1,2}(x) = c_0 B_2 A_0 \left\{ B_4 \sum_{n=1}^{\alpha _2}\left\{ \prod _{i_1=1}^{n-1}A_{i_1+4} \right\} + B_5 A_3 \sum_{n=2}^{\alpha _2}\left\{ \prod _{i_1=2}^{n-1}A_{i_1+4} \right\} \right\} x^{n+4+\lambda } 
\label{eq:47}
\end{equation}
In (\ref{eq:47}) $y_2^{1,2}(x)$ is sub-power series, for the second square brackets inside curly brackets in sequence $c_n$ including two term of $B_n$'s in (\ref{eq:43}) as $\alpha _0$=1 and $\alpha _1$=2, for the polynomial case which makes $A_n$ term terminated.

By using similar process as I did before, the solution for $\alpha _1$=3 with $A_5$=0 for the second square brackets inside curly brackets in sequence $c_n$ including $B_2$ in (\ref{eq:43}) is
\begin{eqnarray}
y_2^{1,3}(x) &=& c_0 B_2 A_0 \left\{ B_4 \sum_{n=1}^{\alpha _2}\left\{ \prod _{i_1=1}^{n-1}A_{i_1+4} \right\}+ B_5 A_3 \sum_{n=2}^{\alpha _2}\left\{ \prod _{i_1=2}^{n-1}A_{i_1+4} \right\}\right.\nonumber\\
&&+\left. B_6 A_{3,4} \sum_{n=3}^{\alpha _2}\left\{ \prod _{i_1=3}^{n-1}A_{i_1+4} \right\} \right\} x^{n+4+\lambda } 
 \label{eq:48}
\end{eqnarray}
By repeating this process for all $\alpha _1 =4,5,6,\cdots$, we obtain every $y_2^{1,j}(x)$ terms where $j \geq 4$ for the second square brackets inside curly brackets in sequence $c_n$ including two term of $B_n$'s in (\ref{eq:43}). 
According to (\ref{eq:42}), (\ref{eq:45}), (\ref{eq:47}), (\ref{eq:48}) and every $y_2^{1,j}(x)$ where $j \geq 4$, the general expression of $y_2^1(x)$ for all $\alpha _0=1$ replacing the index n by $i_0$ is 
\begin{eqnarray}
y_2^1(x)&=& c_0 \left\{B_1 \sum_{i_0=0}^{\alpha _1}\left\{ B_{i_0+3} \prod _{i_1=0}^{i_0-1}A_{i_1+2}  \sum_{i_2=i_0}^{\alpha _2} \left\{ \prod _{i_3=i_0}^{i_2-1}A_{i_3+4} \right\}\right\} \right. \label{eq:49}\\
&&+\left. B_2 A_0 \sum_{i_0=1}^{\alpha _1}\left\{ B_{i_0+3} \prod _{i_1=1}^{i_0-1}A_{i_1+2}  \sum_{i_2=i_0}^{\alpha _2} \left\{ \prod _{i_3=i_0}^{i_2-1}A_{i_3+4} \right\}\right\} \right\} x^{i_2+4+\lambda } \nonumber
\end{eqnarray}
(c) As $\alpha _0$=2, then $A_2$=0 in (\ref{eq:14}).
\begin{equation}
\begin{tabular}{  l  }
  \vspace{2 mm}
  $c_4= B_{1,3} c_0$ \\
  \vspace{2 mm}
  $c_5 = \Big\{ B_1\Big[ B_3\cdot 1 \cdot 1 \cdot A_4 + B_4 \cdot 1 \cdot A_2 \cdot 1\Big]+ B_2\Big[ B_4 A_0 \cdot 1\cdot 1 \Big] \Big\}  c_0  $ \\
  \vspace{2 mm}
  $c_6 = \Big\{ B_1\cdot 1\Big[ B_3\cdot 1 \cdot A_{4,5}+ B_4 A_2 A_5 + B_5 A_{2,3}\cdot 1 \Big] + B_2 A_0 \Big[ B_4\cdot 1\cdot A_5 + B_5 A_3 \cdot 1 \Big] $ \\
  \vspace{2 mm}
  \hspace{0.8 cm} $+ B_3 A_{0,1} \Big[ B_5 \cdot 1\cdot 1 \Big] \Big\} c_0 $\\ 
  \vspace{2 mm}
  $c_7 = \Big\{ B_1\cdot 1 \Big[ B_3\cdot 1 \cdot A_{4,5,6}+ B_4 A_2 A_{5,6}+ B_5 A_{2,3} A_6 + B_6 A_{2,3,4}\cdot 1 \Big] $ \\
 \vspace{2 mm}
  \hspace{0.8 cm} $+B_2 A_0 \Big[ B_4\cdot 1\cdot A_{5,6} + B_5 A_3 A_6 + B_6 A_{3,4} \cdot 1 \Big]+ B_3 A_{0,1} \Big[ B_5 \cdot 1\cdot A_6+ B_6 A_4\cdot 1\Big] \Big\} c_0 $\\ 
  \vspace{2 mm}
  $c_8 = \Big\{ B_1 \cdot 1\Big[ B_3\cdot 1 \cdot A_{4,5,6,7} + B_4 A_2 A_{5,6,7} + B_5 A_{2,3}A_{6,7} +B_6 A_{2,3,4}A_{7} + B_7 A_{2,3,4,5}\cdot 1 \Big] $ \\
 \vspace{2 mm}
  \hspace{0.8 cm} $+B_2 A_0 \Big[ B_4\cdot 1\cdot A_{5,6,7} + B_5 A_3 A_{6,7} + B_6 A_{3,4} A_7  + B_7 A_{3,4,5}\cdot 1 \Big] $\\ 
 \vspace{2 mm}
  \hspace{0.8 cm} $+B_3 A_{0,1} \Big[ B_5\cdot 1\cdot A_{6,7}+ B_6 A_4 A_7 + B_7 A_{4,5}\cdot 1\Big] \Big\} c_0 $\\ 
\hspace{2 mm}
  \large{\vdots}\hspace{7cm}\large{\vdots} \\ 
\end{tabular}
 \label{eq:50}
\end{equation}
By repeating similar process from the above, the general expression of $y_2^2(x)$ for all $\alpha _0=2$ in (\ref{eq:50}) is 
\begin{eqnarray}
y_2^2(x)&=& c_0 \left\{B_1 \sum_{i_0=0}^{\alpha _1}\left\{ B_{i_0+3} \prod _{i_1=0}^{i_0-1}A_{i_1+2}  \sum_{i_2=i_0}^{\alpha _2} \left\{ \prod _{i_3=i_0}^{i_2-1}A_{i_3+4} \right\}\right\}\right. \nonumber\\
&&+ B_2 A_0 \sum_{i_0=1}^{\alpha _1}\left\{ B_{i_0+3} \prod _{i_1=1}^{i_0-1}A_{i_1+2}  \sum_{i_2=i_0}^{\alpha _2} \left\{ \prod _{i_3=i_0}^{i_2-1}A_{i_3+4} \right\}\right\} \label{eq:51}\\
&&+\left. B_3 A_{0,1} \sum_{i_0=2}^{\alpha _1}\left\{ B_{i_0+3} \prod _{i_1=2}^{i_0-1}A_{i_1+2}  \sum_{i_2=i_0}^{\alpha _2} \left\{ \prod _{i_3=i_0}^{i_2-1}A_{i_3+4} \right\}\right\} \right\} x^{i_2+4+\lambda }\nonumber 
\end{eqnarray}
Again by repeating this process for all $\alpha _0 =3,4,5,\cdots$, I obtain every $y_2^{j}(x)$ terms where $j \geq 3$.
Then I have general expression $y_2(x)$ for all $\alpha _0$ of two term of $B_n$'s according (\ref{eq:42}), (\ref{eq:49}), (\ref{eq:51}) and $y_2^{j}(x)$ terms where $j \geq 3$. 
\begin{eqnarray}
y_2(x)&=& c_0 \sum_{i_0=0}^{\alpha _0}\left\{ B_{i_0+1} \prod _{i_1=0}^{i_0-1}A_{i_1}  \sum_{i_2=i_0}^{\alpha _1} \left\{ B_{i_2+3} \prod _{i_3=i_0}^{i_2-1}A_{i_3+2} \sum_{i_4=i_2}^{\alpha _2}\left\{ \prod _{i_5=i_2}^{i_4-1}A_{i_5+4}\right\}  \right\}\right\} x^{i_4+4+\lambda }\nonumber\\
&&\;\;\mbox{where}\; \alpha _0 \leq \alpha _1 \leq \alpha _2  
\label{eq:52}
\end{eqnarray}
By using similar process for the previous cases of zero, one and two terms of $B_n$'s, the function $y_3(x)$ for the case of three term of $B_n$'s is
\begin{eqnarray}
y_3(x) &=& c_0 \sum_{i_0=0}^{\alpha _0}\left\{ B_{i_0+1} \prod _{i_1=0}^{i_0-1}A_{i_1}  \sum_{i_2=i_0}^{\alpha _1} \left\{ B_{i_2+3} \prod _{i_3=i_0}^{i_2-1}A_{i_3+2} \right.\right. \label{eq:53}\\
&&\times \left.\left.  \sum_{i_4=i_2}^{\alpha _2}\left\{ B_{i_4+5}\prod _{i_5=i_2}^{i_4-1}A_{i_5+4} \sum_{i_6=i_4}^{\alpha _3} \left\{\prod _{i_7=i_4}^{i_6-1}A_{i_7+6}  \right\} \right\}  \right\}\right\} x^{i_6+6+\lambda }\nonumber\\
&&\;\;\mbox{where}\; \alpha _0 \leq \alpha _1 \leq \alpha _2 \leq \alpha _3 \nonumber  
\end{eqnarray}
By repeating this process for all higher terms of $B_n$'s, I obtain every $y_m(x)$ terms where $m > 3$. Substitute (\ref{eq:25}), (\ref{eq:33}), (\ref{eq:52}), (\ref{eq:53}) and including all $y_m(x)$ terms where $m > 3$ into (\ref{eq:9}). 
\begin{theorem}
The general expression of $y(x)$ for the polynomial case which makes $A_n$ term terminated in the three term recurrence relation is
\begin{eqnarray}
y(x) &=& y_0(x)+ y_1(x)+ y_2(x)+y_3(x)+\cdots \nonumber\\
&=&  c_0 \left\{ \sum_{i_0=0}^{\alpha _0} \left( \prod _{i_1=0}^{i_0-1}A_{i_1} \right) x^{i_0+\lambda}
+ \sum_{i_0=0}^{\alpha _0}\left\{ B_{i_0+1} \prod _{i_1=0}^{i_0-1}A_{i_1}  \sum_{i_2=i_0}^{\alpha _1} \left( \prod _{i_3=i_0}^{i_2-1}A_{i_3+2} \right)\right\} x^{i_2+2+\lambda} \right. \nonumber\\
&& + \sum_{N=2}^{\infty } \left\{ \sum_{i_0=0}^{\alpha _0} \left\{B_{i_0+1}\prod _{i_1=0}^{i_0-1} A_{i_1} 
\prod _{k=1}^{N-1} \left( \sum_{i_{2k}= i_{2(k-1)}}^{\alpha _k} B_{i_{2k}+2k+1}\prod _{i_{2k+1}=i_{2(k-1)}}^{i_{2k}-1}A_{i_{2k+1}+2k}\right) \right.\right.\nonumber\\
&& \times \left.\left.\left. \sum_{i_{2N} = i_{2(N-1)}}^{\alpha _N} \left( \prod _{i_{2N+1}=i_{2(N-1)}}^{i_{2N}-1} A_{i_{2N+1}+2N} \right) \right\} \right\} x^{i_{2N}+2N+\lambda}\right\}  \label{eq:54}
\end{eqnarray}
In the above, $\alpha _i\leq \alpha _j$ only if $i\leq j$ where $i,j,\alpha _i, \alpha _j \in \mathbb{N}_{0}$ 
\end{theorem}
\begin{theorem}
For infinite series, replacing $\alpha _0$,$\alpha _1$,$\alpha _k$ and $\alpha _N$ by $\infty $ in (\ref{eq:54})
\begin{eqnarray}
y(x) &=& y_0(x)+ y_1(x)+ y_2(x)+y_3(x)+\cdots \nonumber\\
&=& c_0 \left\{ \sum_{i_0=0}^{\infty } \left( \prod _{i_1=0}^{i_0-1}A_{i_1} \right) x^{i_0+\lambda}
+ \sum_{i_0=0}^{\infty }\left\{ B_{i_0+1} \prod _{i_1=0}^{i_0-1}A_{i_1}  \sum_{i_2=i_0}^{\infty } \left( \prod _{i_3=i_0}^{i_2-1}A_{i_3+2} \right)\right\} x^{i_2+2+\lambda} \right. \nonumber\\
&& + \sum_{N=2}^{\infty } \left\{ \sum_{i_0=0}^{\infty } \left\{B_{i_0+1}\prod _{i_1=0}^{i_0-1} A_{i_1} 
\prod _{k=1}^{N-1} \left( \sum_{i_{2k}= i_{2(k-1)}}^{\infty } B_{i_{2k}+2k+1}\prod _{i_{2k+1}=i_{2(k-1)}}^{i_{2k}-1}A_{i_{2k+1}+2k}\right)\right.\right.\nonumber\\
&& \times \left.\left.\left. \sum_{i_{2N} = i_{2(N-1)}}^{\infty } \left( \prod _{i_{2N+1}=i_{2(N-1)}}^{i_{2N}-1} A_{i_{2N+1}+2N} \right) \right\} \right\} x^{i_{2N}+2N+\lambda}\right\}    
\label{eq:55}
\end{eqnarray}
\end{theorem}
(\ref{eq:55}) is equivalent to (\ref{eq:20}). (\ref{eq:55}) is the another general expression of $y(x)$ for the infinite series.

\section{\label{sec:level5}Summary}
In Ref.\cite{Chou2012} I showed how to obtain the general expression of the power series for infinite series and polynomial which makes $B_n$ term terminated including all higher terms of $A_n$'s using 3TRF in a linear ordinary differential equation having a  recursive relation between a 3-term. This was done by letting $A_n$ in sequence $c_n$ is the leading term in the analytic function $y(x)$: the sequence $c_n$ consists of combinations $A_n$ and $B_n$.

In this chapter, I show how to construct the Frobenius solution in closed forms for infinite series and polynomial which makes $A_n$ term terminated including all higher terms of $B_n$'s using R3TRF in a linear ordinary differential equation having a 3-term recurrence relation between successive coefficients. This is done by letting $B_n$ in sequence $c_n$ is the leading term in the analytic function $y(x)$.

The infinite series in this chapter are equivalent to the infinite series in Ref.\cite{Chou2012}. The former is the letting $B_n$ in sequence $c_n$ is the leading term in the analytic function $y(x)$. And the latter is the letting $A_n$ in sequence $c_n$ is the leading term in the analytic function $y(x)$.

In chapters 2--9 I apply R3TRF to the power series expansions, integral forms and generating functions for Heun, Confluent Heun, Mathieu, etc equations; local solutions of these equations for infinite series and polynomial which makes $A_n$ term terminated are constructed analytically. 
 

\addcontentsline{toc}{section}{Bibliography}
\bibliographystyle{model1a-num-names}
\bibliography{<your-bib-database>}

\chapter{Heun function using reversible three-term recurrence formula}
\chaptermark{Heun function using R3TRF}
In Ref.\cite{Chou2012H11,Chou2012H21}, I construct the power series expansions in closed forms of Heun equation and its integral forms for infinite series and polynomial which makes $B_n$ term terminated including all higher terms of $A_n$'s\footnote{`` higher terms of $A_n$'s'' means at least two terms of $A_n$'s.} by applying three term recurrence formula (3TRF).

In this chapter, I will apply reversible three term recurrence formula to the Frobenius solutions in closed forms of Heun equation and its integral forms for infinite series and polynomial which makes $A_n$ term terminated including all higher terms of $B_n$'s\footnote{`` higher terms of $B_n$'s'' means at least two terms of $B_n$'s.}. 

The next chapter out of all 9 chapters describes the generating function of Heun polynomial which makes $A_n$ term terminated including all higher terms of $B_n$'s by using R3TRF.

Nine examples of 192 local solutions of the Heun equation (Maier, 2007) are provided in the appendices 1 and 2.  For each example, I construct the power series expansions of Heun equation and its integral forms by applying R3TRF.
\section{Introduction}

The Heun function generalizes all well-known special functions such as Spheroidal Wave, Lame, Mathieu, and hypergeometric $_2F_1$, $_1F_1$ and $_0F_1$ functions. The Heun function, having three term recurrence relations in its power series, is the most outstanding special functions in among every analytic functions.  

Recently Heun function started to appear in theoretical modern physics. For example the Heun function comes out in the hydrogen-molecule ion\cite{Wils19281}, in the Schr$\ddot{\mbox{o}}$dinger equation with doubly anharmonic potential\cite{Ronv19951} (its solution is the confluent forms of Heun equation), in the Stark effect\cite{Epst19261}, in perturbations of the Kerr metric\cite{Teuk19731,Leav19851,Bati20061,Bati20071,Bati20101}, in crystalline materials\cite{Slavy20001}, in Collogero-Moser-Sutherland systems\cite{Take20031}, etc., just to mention a few.\cite{Suzu19981,Suzu19991,rau20041,Birk20071}  Traditionally, we have explained all physical phenomenons by only using two term recursion relation in the power series expansion until 19th century. However, since modern physics (quantum gravity, SUSY, general relativity, etc) come out of the world, we have at least three or four recurrence relations in power series expansions. Furthermore these type of problems can not be reduced to two term recurrence relations by changing independent variables and coefficients.\cite{Hortacsu:2011rr1}

Due to its complexity Heun function was neglected for almost 100 years\cite{Heun18891}. According to Whittaker's hypothesis, `The Heun function can not be described in form of contour integrals of elementary functions even if it is the simplest class of special functions.' Because Heun equation consists of three recurrence relations in its power series expansion. Heun equation could not be described in the form of a definite or contour integral of any elementary function for the last 100 years. The 3-term recursive relation between successive coefficients in its power series expansion creates mathematical difficulty to be described it into direct integral representation. Instead in Fredholm integral equations, Heun equation is obtained by integral equations; such integral relationships express one analytic solution in terms of another analytic solution. More precisely, in earlier literature the integral representations of Heun equation were constructed by using two types of relations: (1) Linear relations using Fredholm integral equations. \cite{Lamb19341,Erde19421} (2) Non-linear relation (Malurkar-type integral relations) including Fredholm integral equations using two variables. \cite{Slee1969a1,Slee1969b1,Arsc19641,Schm19791}

In previous papers I show the analytic solutions of Heun equation for infinite series and polynomial which makes $B_n$ term terminated including all higher terms of $A_n$'s by applying 3TRF \cite{Chou2012H11,Chou2012H21}; (1) the power series expansion of Heun equation, (2) its integral representation and its asymptotic behavior including the boundary condition for an independent variable $x$. Indeed, according to analyzing Heun function into its combined direct and contour integral forms resulting in a precise, we are able to observe how Heun function relate to other well-known special functions such as Mathieu function, Lame function, confluent forms of Heun function and etc. Because a $_1F_1$ or $_2F_1$ function recurs in each of sub-integral forms of all these analytic functions. Also, an orthogonal relation of Heun function can be obtained from its integral representation.

In this chapter, by applying R3TRF, I construct the power series expansion in closed forms of Heun equation for infinite series and polynomial which makes $A_n$ term terminated including all higher terms of $B_n$'s analytically and its integral form. Heun's  differential equation is a second-order linear ordinary differential equation of the form \cite{Heun18891}.
\begin{equation}
\frac{d^2{y}}{d{x}^2} + \left(\frac{\gamma }{x} +\frac{\delta }{x-1} + \frac{\epsilon }{x-a}\right) \frac{d{y}}{d{x}} +  \frac{\alpha \beta x-q}{x(x-1)(x-a)} y = 0 \label{eq:1001}
\end{equation}
With the condition $\epsilon = \alpha +\beta -\gamma -\delta +1$. The parameters play different roles: $a \ne 0 $ is the singularity parameter, $\alpha $, $\beta $, $\gamma $, $\delta $, $\epsilon $ are exponent parameters, $q$ is the accessory parameter which in many physical applications appears as a spectral parameter. Also, $\alpha $ and $\beta $ are identical to each other. The total number of free parameters is six. It has four regular singular points which are 0, 1, $a$ and $\infty $ with exponents $\{ 0, 1-\gamma \}$, $\{ 0, 1-\delta \}$, $\{ 0, 1-\epsilon \}$ and $\{ \alpha, \beta \}$. Assume that its solution is
\begin{equation}
y(x)= \sum_{n=0}^{\infty } c_n x^{n+\lambda } \hspace{1cm}\mbox{where}\;\lambda  =\mbox{indicial}\;\mbox{root}\label{eq:1002}
\end{equation}
Plug (\ref{eq:1002})  into (\ref{eq:1001}).
\begin{equation}
c_{n+1}=A_n \;c_n +B_n \;c_{n-1} \hspace{1cm};n\geq 1 \label{eq:1003}
\end{equation}
where,
\begin{subequations}
\begin{eqnarray}
A_n &=& \frac{(n+\lambda )(n-1+\gamma +\epsilon +\lambda + a(n-1+\gamma +\lambda +\delta ))+q}{a(n+1+\lambda )(n+\gamma +\lambda )}\nonumber\\
&=& \frac{(n+\lambda )(n+\alpha +\beta -\delta +\lambda +a(n+\delta +\gamma -1+\lambda ))+q}{a(n+1+\lambda )(n+\gamma +\lambda )}\nonumber\\
&=& \frac{(1+a)}{a}\frac{\left( n- \frac{-(\varphi +2(1+a)\lambda )-\sqrt{\varphi ^2-4(1+a)q}}{2(1+a)}\right) \left( n- \frac{-(\varphi +2(1+a)\lambda )+\sqrt{\varphi ^2-4(1+a)q}}{2(1+a)}\right)}{(n+\lambda +1)(n+\lambda +\gamma )} \hspace{2cm}\label{eq:1004a}
\end{eqnarray}
and
\begin{equation}
\varphi =\alpha +\beta -\delta +a(\delta +\gamma -1)\nonumber
\end{equation}
\vspace{2mm}
\begin{eqnarray}
B_n &=& -\frac{(n-1+\lambda )(n+\gamma +\delta +\epsilon -2+\lambda )+\alpha \beta }{a(n+1+\lambda )(n+\gamma +\lambda )}\nonumber\\
&=&- \frac{(n-1+\lambda +\alpha )(n-1+\lambda +\beta )}{a(n+1+\lambda )(n+\gamma +\lambda )} \label{eq:1004b}
\end{eqnarray}
\begin{equation}
c_1= A_0 \;c_0 \label{eq:1004c}
\end{equation}
\end{subequations}
We have two indicial roots which are $\lambda = 0$ and $ 1-\gamma $.

\section{Power series}
\subsection{Polynomial of type 2}
There are three types of polynomials in three term recurrence relation of a linear ordinary differential equation: (1) polynomial which makes $B_n$ term terminated: $A_n$ term is not terminated, (2) polynomial which makes $A_n$ term terminated: $B_n$ term is not terminated, (3) polynomial which makes $A_n$ and $B_n$ terms terminated at the same time. In general Heun polynomial is defined as type 3 polynomial where $A_n$ and $B_n$ terms terminated. Heun polynomial comes from Heun equation that has a fixed integer value of $\alpha $ or $\beta $, just as it has a fixed value of $q$. In three term recurrence relation, polynomial of type 3 I categorize as complete polynomial. In future papers I will derive type 3 Heun polynomial. In this part I construct the power series expansion and an integral forms for Heun polynomial of type 2: I treat $\alpha $, $\beta $, $\gamma $ and $\delta $ as free variables and the accessory parameter $q$ as a fixed value. In my next papers I will work on the generating functions for Heun polynomial of type 2.\footnote{If $A_n$ and $B_n$ terms are not terminated, it turns to be infinite series.} 

In chapter 1, the general expression of power series of $y(x)$ for polynomial of type 2 is
\begin{eqnarray}
 y(x)&=& \sum_{n=0}^{\infty } y_{n}(x)= y_0(x)+ y_1(x)+ y_2(x)+y_3(x)+\cdots \nonumber\\
&=&  c_0 \Bigg\{ \sum_{i_0=0}^{\alpha _0} \left( \prod _{i_1=0}^{i_0-1}A_{i_1} \right) x^{i_0+\lambda }
+ \sum_{i_0=0}^{\alpha _0}\left\{ B_{i_0+1} \prod _{i_1=0}^{i_0-1}A_{i_1}  \sum_{i_2=i_0}^{\alpha _1} \left( \prod _{i_3=i_0}^{i_2-1}A_{i_3+2} \right)\right\} x^{i_2+2+\lambda }  \nonumber\\
&& + \sum_{N=2}^{\infty } \Bigg\{ \sum_{i_0=0}^{\alpha _0} \Bigg\{B_{i_0+1}\prod _{i_1=0}^{i_0-1} A_{i_1} 
\prod _{k=1}^{N-1} \Bigg( \sum_{i_{2k}= i_{2(k-1)}}^{\alpha _k} B_{i_{2k}+2k+1}\prod _{i_{2k+1}=i_{2(k-1)}}^{i_{2k}-1}A_{i_{2k+1}+2k}\Bigg)\nonumber\\
&& \times  \sum_{i_{2N} = i_{2(N-1)}}^{\alpha _N} \Bigg( \prod _{i_{2N+1}=i_{2(N-1)}}^{i_{2N}-1} A_{i_{2N+1}+2N} \Bigg) \Bigg\} \Bigg\} x^{i_{2N}+2N+\lambda }\Bigg\}  \label{eq:1005}
\end{eqnarray}
In the above, $\alpha _i\leq \alpha _j$ only if $i\leq j$ where $i,j,\alpha _i, \alpha _j \in \mathbb{N}_{0}$.

For a polynomial, we need a condition, which is:
\begin{equation}
A_{\alpha _i+ 2i}=0 \hspace{1cm} \mathrm{where}\;i,\alpha _i =0,1,2,\cdots
\label{eq:1006}
\end{equation}
In the above, $ \alpha _i$ is an eigenvalue that makes $A_n$ term terminated at certain value of index $n$. (\ref{eq:1006}) makes each $y_i(x)$ where $i=0,1,2,\cdots$ as the polynomial in (\ref{eq:1005}).
Replace $\alpha _i$ by $q_i$ and put $n=q_i+ 2i$ in (\ref{eq:1004a}) with the condition $A_{q_i+ 2i}=0$. Then, we obtain eigenvalues $q$ such as 
\begin{equation}
q_i +2i= \frac{-(\varphi +2(1+a)\lambda )\pm \sqrt{\varphi ^2-4(1+a)q}}{2(1+a)}
\nonumber
\end{equation}
\paragraph{The case of $ \sqrt{\varphi ^2-4(1+a)q}= -\{\varphi + 2(1+a)(q_i+2i+\lambda )\}$}
In (\ref{eq:1004a}) replace $\sqrt{\varphi ^2-4(1+a)q}$ by ${ \displaystyle -\{\varphi + 2(1+a)(q_i+2i+\lambda )\}}$. In (\ref{eq:1005}) replace index $\alpha _i$ by $q_i$. Take the new (\ref{eq:1004a}) and (\ref{eq:1004b}) in new (\ref{eq:1005}).
After the replacement process, the general expression of power series of Heun's equation for polynomial of type 2 is given by
\begin{eqnarray}
  y(x)&=& \sum_{n=0}^{\infty } y_{n}(x)= y_0(x)+ y_1(x)+ y_2(x)+y_3(x)+\cdots \nonumber\\
&=& c_0 x^{\lambda } \left\{\sum_{i_0=0}^{q_0} \frac{(-q_0)_{i_0} \left(q_0+ \frac{\varphi +2(1+a)\lambda }{(1+a)}\right)_{i_0}}{(1+\lambda )_{i_0}(\gamma +\lambda )_{i_0}} \eta ^{i_0}\right.\nonumber\\
&&+ \left\{ \sum_{i_0=0}^{q_0}\frac{(i_0+ \lambda +\alpha ) (i_0+ \lambda +\beta )}{(i_0+ \lambda +2)(i_0+ \lambda +1+\gamma )}\frac{(-q_0)_{i_0} \left(q_0+ \frac{\varphi +2(1+a)\lambda }{(1+a)}\right)_{i_0}}{(1+\lambda )_{i_0}(\gamma +\lambda )_{i_0}} \right.\nonumber\\
&&\times  \left. \sum_{i_1=i_0}^{q_1} \frac{(-q_1)_{i_1}\left(q_1+4+ \frac{\varphi +2(1+a)\lambda }{(1+a)}\right)_{i_1}(3+\lambda )_{i_0}(2+\gamma +\lambda )_{i_0}}{(-q_1)_{i_0}\left(q_1+4+ \frac{\varphi +2(1+a)\lambda }{(1+a)}\right)_{i_0}(3+\lambda )_{i_1}(2+\gamma +\lambda )_{i_1}} \eta ^{i_1}\right\} z\nonumber\\
&&+ \sum_{n=2}^{\infty } \left\{ \sum_{i_0=0}^{q_0} \frac{(i_0+ \lambda +\alpha ) (i_0+ \lambda +\beta )}{(i_0+ \lambda +2)(i_0+ \lambda +1+\gamma )}\frac{(-q_0)_{i_0} \left(q_0+ \frac{\varphi +2(1+a)\lambda }{(1+a)}\right)_{i_0}}{(1+\lambda )_{i_0}(\gamma +\lambda )_{i_0}}\right.\nonumber\\
&&\times \prod _{k=1}^{n-1} \left\{ \sum_{i_k=i_{k-1}}^{q_k} \frac{(i_k+ 2k+\lambda +\alpha ) (i_k+ 2k+\lambda +\beta )}{(i_k+ 2(k+1)+\lambda )(i_k+ 2k+1+\gamma +\lambda )}\right. \label{eq:1007}\\
&&\times \left.\frac{(-q_k)_{i_k}\left(q_k+4k+ \frac{\varphi +2(1+a)\lambda }{(1+a)}\right)_{i_k}(2k+1+\lambda )_{i_{k-1}}(2k+\gamma +\lambda )_{i_{k-1}}}{(-q_k)_{i_{k-1}}\left(q_k+4k+ \frac{\varphi +2(1+a)\lambda }{(1+a)}\right)_{i_{k-1}}(2k+1+\lambda )_{i_k}(2k+\gamma +\lambda )_{i_k}}\right\} \nonumber\\
&&\times \left. \left.\sum_{i_n= i_{n-1}}^{q_n} \frac{(-q_n)_{i_n}\left(q_n+4n+ \frac{\varphi +2(1+a)\lambda }{(1+a)}\right)_{i_n}(2n+1+\lambda )_{i_{n-1}}(2n+\gamma +\lambda )_{i_{n-1}}}{(-q_n)_{i_{n-1}}\left(q_n+4n+ \frac{\varphi +2(1+a)\lambda }{(1+a)}\right)_{i_{n-1}}(2n+1+\lambda )_{i_n}(2n+\gamma +\lambda )_{i_n}} \eta ^{i_n} \right\} z^n \right\}\nonumber
\end{eqnarray}
where
\begin{equation}
\begin{cases} z = -\frac{1}{a}x^2 \cr
\eta = \frac{(1+a)}{a} x \cr
q= -(q_j+2j+\lambda )\{\varphi +(1+a)(q_j+2j+\lambda ) \} \;\;\mbox{as}\;j,q_j\in \mathbb{N}_{0} \cr
q_i\leq q_j \;\;\mbox{only}\;\mbox{if}\;i\leq j\;\;\mbox{where}\;i,j\in \mathbb{N}_{0} 
\end{cases}\nonumber 
\end{equation}
\paragraph{The case of $ \sqrt{\varphi ^2-4(1+a)q}= \varphi + 2(1+a)(q_i+2i+\lambda )$}
In (\ref{eq:1004a}) replace $\sqrt{\varphi ^2-4(1+a)q}$ by ${ \displaystyle \varphi + 2(1+a)(q_i+2i+\lambda )}$. In (\ref{eq:1005}) replace index $\alpha _i$ by $q_i$. Take the new (\ref{eq:1004a}) and (\ref{eq:1004b}) in new (\ref{eq:1005}).
After the replacement process, its solution is equivalent to (\ref{eq:1007}).

Put $c_0$= 1 as $\lambda $=0 for the first kind of independent solutions of Heun equation and $\displaystyle{ c_0= \left( \frac{1+a}{a}\right)^{1-\gamma }}$ as $\lambda = 1-\gamma $ for the second one in (\ref{eq:1007}). 
\begin{remark}
The power series expansion of Heun equation of the first kind for polynomial of type 2 about $x=0$ as $q= -(q_j+2j)\{ \alpha +\beta -\delta +a(\delta +\gamma -1)+(1+a)(q_j+2j)\}$ where $j,q_j \in \mathbb{N}_{0}$ is
\begin{eqnarray}
 &y(x)&= HF_{q_j}^R \bigg( q_j =\frac{-\varphi \pm \sqrt{\varphi ^2-4(1+a)q}}{2(1+a)}-2j, \varphi =\alpha +\beta -\delta +a(\delta +\gamma -1), \Omega _1=\frac{\varphi }{(1+a)}\nonumber\\
&&; \eta = \frac{(1+a)}{a} x ; z= -\frac{1}{a} x^2 \bigg) \nonumber\\
&=& \sum_{i_0=0}^{q_0} \frac{(-q_0)_{i_0} \left(q_0+ \Omega _1\right)_{i_0}}{(1)_{i_0}(\gamma )_{i_0}} \eta ^{i_0}\nonumber\\
&+& \left\{ \sum_{i_0=0}^{q_0}\frac{(i_0+\alpha ) (i_0+\beta )}{(i_0+2)(i_0+ 1+\gamma )}\frac{(-q_0)_{i_0} \left(q_0+ \Omega _1\right)_{i_0}}{(1)_{i_0}(\gamma)_{i_0}} 
  \sum_{i_1=i_0}^{q_1} \frac{(-q_1)_{i_1}\left(q_1+4+ \Omega _1\right)_{i_1}(3)_{i_0}(2+\gamma)_{i_0}}{(-q_1)_{i_0}\left(q_1+4+ \Omega _1\right)_{i_0}(3)_{i_1}(2+\gamma)_{i_1}} \eta ^{i_1}\right\} z\nonumber\\
&+& \sum_{n=2}^{\infty } \left\{ \sum_{i_0=0}^{q_0} \frac{(i_0+\alpha ) (i_0+\beta )}{(i_0+ 2)(i_0+ 1+\gamma )}\frac{(-q_0)_{i_0} \left(q_0+ \Omega _1\right)_{i_0}}{(1)_{i_0}(\gamma )_{i_0}}\right.\nonumber\\
&\times& \prod _{k=1}^{n-1} \left\{ \sum_{i_k=i_{k-1}}^{q_k} \frac{(i_k+ 2k+\alpha ) (i_k+ 2k+\beta )}{(i_k+ 2(k+1))(i_k+ 2k+1+\gamma )} 
 \frac{(-q_k)_{i_k}\left(q_k + 4k + \Omega _1\right)_{i_k}(2k+1)_{i_{k-1}}(2k+\gamma )_{i_{k-1}}}{(-q_k)_{i_{k-1}}\left(q_k + 4k + \Omega _1\right)_{i_{k-1}}(2k+1)_{i_k}(2k+\gamma )_{i_k}}\right\} \nonumber\\
&\times& \left.\sum_{i_n= i_{n-1}}^{q_n} \frac{(-q_n)_{i_n}\left(q_n + 4n + \Omega _1\right)_{i_n}(2n+1)_{i_{n-1}}(2n+\gamma )_{i_{n-1}}}{(-q_n)_{i_{n-1}}\left(q_n + 4n + \Omega _1\right)_{i_{n-1}}(2n+1 )_{i_n}(2n+\gamma  )_{i_n}} \eta ^{i_n} \right\} z^n \label{eq:10010}
\end{eqnarray}
\end{remark}
For the minimum value of Heun equation of the first kind for a polynomial of type 2 about $x=0$, put $q_0=q_1=q_2=\cdots=0$ in (\ref{eq:10010}).
\begin{eqnarray}
y(x)&=& HF_{0}^R \left( q=-2j(\varphi +2(1+a)j), \varphi =\alpha +\beta -\delta +a(\delta +\gamma -1), \Omega _1=\frac{\varphi }{(1+a)}\right.\nonumber\\
&&\left.; \eta = \frac{(1+a)}{a} x ; z= -\frac{1}{a} x^2 \right) = \; _2F_1\left( \frac{\alpha }{2},\frac{\beta }{2},\frac{\gamma }{2}+\frac{1}{2},z \right) \hspace{.5cm}\mbox{where}\;\;|z| < 1 \label{cho:1001}
\end{eqnarray} 
\begin{remark}
The power series expansion of Heun equation of the second kind for polynomial of type 2 about $x=0$ as $q= -(q_j+2j+1-\gamma )\{ \alpha +\beta+1-\gamma -(1-a) \delta +(1+a)(q_j+2j)\}$ where $j,q_j \in \mathbb{N}_{0}$ is
\begin{eqnarray}
y(x)&=& HS_{q_j}^R \Bigg( q_j =\frac{-\{\varphi +2(1+a)(1-\gamma )\} \pm \sqrt{\varphi ^2-4(1+a)q}}{2(1+a)}-2j\nonumber\\
&&,\varphi =\alpha +\beta -\delta +a(\delta +\gamma -1), \Omega _2= \frac{\varphi +2(1+a)(1-\gamma )}{(1+a)}; \eta = \frac{(1+a)}{a} x ; z= -\frac{1}{a} x^2 \Bigg) \nonumber\\
&=& \eta ^{1-\gamma } \left\{\sum_{i_0=0}^{q_0} \frac{(-q_0)_{i_0} \left(q_0+ \Omega _2\right)_{i_0}}{(2-\gamma )_{i_0}(1)_{i_0}} \eta ^{i_0}\right.\nonumber\\
&+& \left\{ \sum_{i_0=0}^{q_0}\frac{(i_0+ 1-\gamma +\alpha ) (i_0+ 1-\gamma +\beta )}{(i_0+ 3-\gamma)(i_0+ 2)}\frac{(-q_0)_{i_0} \left( q_0+ \Omega _2\right)_{i_0}}{(2-\gamma )_{i_0}(1)_{i_0}} \right.\nonumber\\
&\times&  \left.\sum_{i_1=i_0}^{q_1} \frac{(-q_1)_{i_1}\left(q_1+4+ \Omega _2\right)_{i_1}(4-\gamma  )_{i_0}(3)_{i_0}}{(-q_1)_{i_0}\left(q_1+4+ \Omega _2\right)_{i_0}(4-\gamma  )_{i_1}(3)_{i_1}} \eta ^{i_1}\right\} z\nonumber\\
&+& \sum_{n=2}^{\infty } \left\{ \sum_{i_0=0}^{q_0} \frac{(i_0+ 1-\gamma +\alpha ) (i_0+ 1-\gamma +\beta )}{(i_0+ 3-\gamma )(i_0+ 2)}\frac{(-q_0)_{i_0} \left(q_0+ \Omega _2\right)_{i_0}}{(2-\gamma )_{i_0}(1)_{i_0}}\right.\nonumber\\
&\times& \prod _{k=1}^{n-1} \left\{ \sum_{i_k=i_{k-1}}^{q_k} \frac{(i_k+ 2k+1-\gamma+\alpha ) (i_k+ 2k+1-\gamma+\beta )}{(i_k+ 2(k+1)+1-\gamma)(i_k+ 2(k+1))} \right.\nonumber\\
&\times& \left.\frac{(-q_k)_{i_k}\left(q_k+4k+ \Omega _2\right)_{i_k}(2(k+1)-\gamma  )_{i_{k-1}}(2k+1 )_{i_{k-1}}}{(-q_k)_{i_{k-1}}\left( q_k+4k+ \Omega _2\right)_{i_{k-1}}(2(k+1)-\gamma )_{i_k}(2k+1 )_{i_k}}\right\} \label{eq:10011}\\
&\times& \left.\left.\sum_{i_n= i_{n-1}}^{q_n} \frac{(-q_n)_{i_n}\left(q_n+4n+ \Omega _2\right)_{i_n}(2(n+1)-\gamma )_{i_{n-1}}(2n+1 )_{i_{n-1}}}{(-q_n)_{i_{n-1}}\left(q_n+4n+ \Omega _2\right)_{i_{n-1}}(2(n+1)-\gamma  )_{i_n}(2n+1 )_{i_n}} \eta ^{i_n} \right\} z^n \right\}\nonumber
\end{eqnarray}
\end{remark}
For the minimum value of Heun equation of the second kind for a polynomial of type 2 about $x=0$, put $q_0=q_1=q_2=\cdots=0$ in (\ref{eq:10011}).
\begin{eqnarray}
y(x)&=& HS_{0}^R \Bigg( q=-(1+2j-\gamma )(\varphi +(1+a)(1+2j-\gamma )), \varphi =\alpha +\beta -\delta +a(\delta +\gamma -1)\nonumber\\
&&,\Omega _2= \frac{\varphi +2(1+a)(1-\gamma )}{(1+a)}; \eta = \frac{(1+a)}{a} x ; z= -\frac{1}{a} x^2 \Bigg) \nonumber\\
&=& \eta ^{1-\gamma } \; _2F_1\left( \frac{\alpha }{2}-\frac{\gamma }{2}+\frac{1}{2},\frac{\beta }{2}-\frac{\gamma }{2}+\frac{1}{2},-\frac{\gamma }{2}+\frac{3}{2},z \right) \hspace{.5cm}\mbox{where}\;\;|z| < 1 \label{cho:1002}
\end{eqnarray}
In (\ref{cho:1001}) and (\ref{cho:1002}), a polynomial of type 2 requires $\left| z\right|< 1$ for the convergence of the radius.

In Ref.\cite{Chou2012H11,Chou2012H21} I treat $\alpha $ and/or $\beta $ as a fixed value and $\gamma , \delta $, $q$ as free variables to construct Heun polynomial of type 1: (1) if $\alpha = -2\alpha_j -j $ and/or $\beta  = -2\beta _j -j $ where $j, \alpha_j, \beta _j \in \mathbb{N}_{0}$, an analytic solution of Heun equation turns to be the first kind of independent solution of Heun polynomial which makes $B_n$ term terminated. (2) if  $\alpha = -2\alpha_j -j -1+\gamma $ and/or  $\beta  = -2\beta _j -j -1+\gamma $, an analytic solution of Heun equation turns to be the second kind of independent solution of Heun polynomial of type 1. 

In this chapter I treat $q$ as a fixed value and $\alpha, \beta, \gamma, \delta $ as free variables to construct Heun polynomial of type 2: (1) if $q= -(q_j+2j)\{ \alpha +\beta -\delta +a(\delta +\gamma -1)+(1+a)(q_j+2j)\}$ where $j,q_j \in \mathbb{N}_{0}$, an analytic solution of Heun equation turns to be the first kind of independent solution of Heun polynomial of type 2. (2) if $q= -(q_j+2j+1-\gamma )\{ \alpha +\beta+1-\gamma -(1-a) \delta +(1+a)(q_j+2j)\}$, an analytic solution of Heun equation turns to be the second kind of independent solution of Heun polynomial of type 2.
\subsection{Infinite series}
In chapter 1, the general expression of power series of $y(x)$ for infinite series is
\begin{eqnarray}
 y(x)&=& \sum_{n=0}^{\infty } y_{n}(x)= y_0(x)+ y_1(x)+ y_2(x)+y_3(x)+\cdots \nonumber\\
&=& c_0 \Bigg\{ \sum_{i_0=0}^{\infty } \left( \prod _{i_1=0}^{i_0-1}A_{i_1} \right) x^{i_0+\lambda }
+ \sum_{i_0=0}^{\infty }\left\{ B_{i_0+1} \prod _{i_1=0}^{i_0-1}A_{i_1}  \sum_{i_2=i_0}^{\infty } \left( \prod _{i_3=i_0}^{i_2-1}A_{i_3+2} \right)\right\} x^{i_2+2+\lambda }  \nonumber\\
&&+ \sum_{N=2}^{\infty } \Bigg\{ \sum_{i_0=0}^{\infty } \Bigg\{B_{i_0+1}\prod _{i_1=0}^{i_0-1} A_{i_1} 
\prod _{k=1}^{N-1} \Bigg( \sum_{i_{2k}= i_{2(k-1)}}^{\infty } B_{i_{2k}+2k+1}\prod _{i_{2k+1}=i_{2(k-1)}}^{i_{2k}-1}A_{i_{2k+1}+2k}\Bigg)\nonumber\\
&&\times  \sum_{i_{2N} = i_{2(N-1)}}^{\infty } \Bigg( \prod _{i_{2N+1}=i_{2(N-1)}}^{i_{2N}-1} A_{i_{2N+1}+2N} \Bigg) \Bigg\} \Bigg\} x^{i_{2N}+2N+\lambda }\Bigg\}   \label{eq:10019}
\end{eqnarray}
Substitute (\ref{eq:1004a})-(\ref{eq:1004c}) into (\ref{eq:10019}). 
The general expression of power series of Heun equation for infinite series about $x=0$ is
\begin{eqnarray}
 y(x)&=& \sum_{n=0}^{\infty } y_n(x)= y_0(x)+ y_1(x)+ y_2(x)+ y_3(x)+\cdots \nonumber\\
&=& c_0 x^{\lambda } \left\{\sum_{i_0=0}^{\infty } \frac{\left(\Delta_0^{-}\right)_{i_0} \left(\Delta_0^{+}\right)_{i_0}}{(1+\lambda )_{i_0}(\gamma +\lambda )_{i_0}} \eta ^{i_0}\right.\nonumber\\
&+& \left\{ \sum_{i_0=0}^{\infty }\frac{(i_0+ \lambda +\alpha ) (i_0+ \lambda +\beta )}{(i_0+ \lambda +2)(i_0+ \lambda +1+\gamma )}\frac{\left(\Delta_0^{-}\right)_{i_0} \left(\Delta_0^{+}\right)_{i_0}}{(1+\lambda )_{i_0}(\gamma +\lambda )_{i_0}} \sum_{i_1=i_0}^{\infty } \frac{\left(\Delta_1^{-}\right)_{i_1} \left(\Delta_1^{+}\right)_{i_1}(3+\lambda )_{i_0}(2+\gamma +\lambda )_{i_0}}{\left(\Delta_1^{-}\right)_{i_0}  \left(\Delta_1^{+}\right)_{i_0}(3+\lambda )_{i_1}(2+\gamma +\lambda )_{i_1}}\eta ^{i_1}\right\} z\nonumber\\
&+& \sum_{n=2}^{\infty } \left\{ \sum_{i_0=0}^{\infty } \frac{(i_0+ \lambda +\alpha ) (i_0+ \lambda +\beta )}{(i_0+ \lambda +2)(i_0+ \lambda +1+\gamma )}\frac{\left(\Delta_0^{-}\right)_{i_0} \left(\Delta_0^{+}\right)_{i_0}}{(1+\lambda )_{i_0}(\gamma +\lambda )_{i_0}}\right.\nonumber\\
&\times& \prod _{k=1}^{n-1} \left\{ \sum_{i_k=i_{k-1}}^{\infty } \frac{(i_k+ 2k+\lambda +\alpha ) (i_k+ 2k+\lambda +\beta )}{(i_k+ 2(k+1)+\lambda )(i_k+ 2k+1+\gamma +\lambda )}\right.\nonumber\\
&\times& \left. \frac{ \left(\Delta_k^{-}\right)_{i_k} \left(\Delta_k^{+} \right)_{i_k}(2k+1+\lambda )_{i_{k-1}}(2k+\gamma +\lambda )_{i_{k-1}}}{\left(\Delta_k^{-}\right)_{i_{k-1}} \left(\Delta_k^{+} \right)_{i_{k-1}}(2k+1+\lambda )_{i_k}(2k+\gamma +\lambda )_{i_k}}\right\}\nonumber\\
&\times& \left.\left.\sum_{i_n= i_{n-1}}^{\infty } \frac{\left(\Delta_n^{-}\right)_{i_n}\left( \Delta_n^{+} \right)_{i_n}(2n+1+\lambda )_{i_{n-1}}(2n+\gamma +\lambda )_{i_{n-1}}}{\left(\Delta_n^{-}\right)_{i_{n-1}}\left(\Delta_n^{+} \right)_{i_{n-1}}3(2n+1+\lambda )_{i_n}(2n+\gamma +\lambda )_{i_n}} \eta ^{i_n} \right\} z^n \right\} \label{eq:10020}
\end{eqnarray}
where
\Large
\begin{equation}
\begin{cases} 
\Delta_0^{\pm}= \frac{\{\varphi +2(1+a)\lambda \} \pm\sqrt{\varphi ^2-4(1+a)q}}{2(1+a)} \cr
\Delta_1^{\pm}=  \frac{\{\varphi +2(1+a)(\lambda+2 ) \} \pm\sqrt{\varphi ^2-4(1+a)q}}{2(1+a)} \cr
\Delta_k^{\pm}=  \frac{\{\varphi +2(1+a)(\lambda+2k )\} \pm\sqrt{\varphi ^2-4(1+a)q}}{2(1+a)} \cr
\Delta_n^{\pm}=   \frac{\{\varphi +2(1+a)(\lambda+2n )\} \pm\sqrt{\varphi ^2-4(1+a)q}}{2(1+a)}
\end{cases}\nonumber 
\end{equation}
\normalsize
Put $c_0$= 1 as $\lambda $=0 for the first kind of independent solutions of Heun equation and $\displaystyle{ c_0= \left( \frac{1+a}{a}\right)^{1-\gamma }}$ as $\lambda = 1-\gamma $ for the second one in (\ref{eq:10020}).
\begin{remark}
The power series expansion of Heun equation of the first kind for infinite series about $x=0$ using R3TRF is
\footnotesize
\begin{eqnarray}
 y(x)&=& HF^R \left( \varphi =\alpha +\beta -\delta +a(\delta +\gamma -1); \eta = \frac{(1+a)}{a} x ; z= -\frac{1}{a} x^2 \right) \nonumber\\
&=& \sum_{i_0=0}^{\infty } \frac{\left(\Delta_0^{-}\right)_{i_0} \left(\Delta_0^{+}\right)_{i_0}}{(1)_{i_0}(\gamma )_{i_0}} \eta ^{i_0}\nonumber\\
&+& \left\{ \sum_{i_0=0}^{\infty }\frac{(i_0+ \alpha ) (i_0+ \beta )}{(i_0+ 2)(i_0+ 1+\gamma )}\frac{\left(\Delta_0^{-}\right)_{i_0} \left(\Delta_0^{+}\right)_{i_0}}{(1)_{i_0}(\gamma)_{i_0}} \sum_{i_1=i_0}^{\infty } \frac{\left( \Delta_1^{-}\right)_{i_1} \left( \Delta_1^{+}\right)_{i_1}(3)_{i_0}(2+\gamma )_{i_0}}{\left(\Delta_1^{-}\right)_{i_0} \left(\Delta_1^{+}\right)_{i_0}(3)_{i_1}(2+\gamma)_{i_1}} \eta ^{i_1}\right\} z\nonumber\\
&+& \sum_{n=2}^{\infty } \left\{ \sum_{i_0=0}^{\infty } \frac{(i_0 +\alpha ) (i_0 +\beta )}{(i_0 +2)(i_0 +1+\gamma )}\frac{\left(\Delta_0^{-}\right)_{i_0} \left(\Delta_0^{+}\right)_{i_0}}{(1)_{i_0}(\gamma )_{i_0}}\right.\nonumber\\
&\times& \prod _{k=1}^{n-1} \left\{ \sum_{i_k=i_{k-1}}^{\infty } \frac{(i_k+ 2k+\alpha ) (i_k+ 2k+\beta )}{(i_k+ 2(k+1) )(i_k+ 2k+1+\gamma )}\frac{\left(\Delta_k^{-}\right)_{i_k} \left( \Delta_k^{+}\right)_{i_k}(2k+1 )_{i_{k-1}}(2k+\gamma )_{i_{k-1}}}{ \left( \Delta_k^{-}\right)_{i_{k-1}}\left( \Delta_k^{+}\right)_{i_{k-1}}(2k+1)_{i_k}(2k+\gamma )_{i_k}}\right\}\nonumber \\
&\times& \left. \sum_{i_n= i_{n-1}}^{\infty } \frac{\left( \Delta_n^{-}\right)_{i_n} \left(\Delta_n^{+}\right)_{i_n}(2n+1 )_{i_{n-1}}(2n+\gamma )_{i_{n-1}}}{\left(\Delta_n^{-}\right)_{i_{n-1}}\left( \Delta_n^{+}\right)_{i_{n-1}}(2n+1)_{i_n}(2n+\gamma )_{i_n} } \eta ^{i_n} \right\} z^n  \label{eq:10021}
\end{eqnarray}
\normalsize
where
\Large
\begin{equation}
\begin{cases} 
\Delta_0^{\pm}= \frac{\varphi \pm\sqrt{\varphi ^2-4(1+a)q}}{2(1+a)} \cr
\Delta_1^{\pm}=  \frac{\{\varphi +4(1+a) \} \pm\sqrt{\varphi ^2-4(1+a)q}}{2(1+a)} \cr
\Delta_k^{\pm}=  \frac{\{\varphi +4(1+a)k \} \pm\sqrt{\varphi ^2-4(1+a)q}}{2(1+a)} \cr
\Delta_n^{\pm}=   \frac{\{\varphi +4(1+a)n\} \pm\sqrt{\varphi ^2-4(1+a)q}}{2(1+a)}
\end{cases}\nonumber 
\end{equation}
\normalsize
\end{remark}
\begin{remark}
The power series expansion of Heun equation of the second kind for infinite series about $x=0$ using R3TRF is
\footnotesize
\begin{eqnarray}
y(x)&=& HS^R \left( \varphi =\alpha +\beta -\delta +a(\delta +\gamma -1); \eta = \frac{(1+a)}{a} x ; z= -\frac{1}{a} x^2 \right) \nonumber\\
&=& \eta ^{1-\gamma } \left\{\sum_{i_0=0}^{\infty } \frac{\left(\Delta_0^{-}\right)_{i_0} \left(\Delta_0^{+}\right)_{i_0}}{(2-\gamma )_{i_0}(1)_{i_0}} \eta ^{i_0}\right.\nonumber\\
&+& \left\{ \sum_{i_0=0}^{\infty }\frac{(i_0+ 1-\gamma +\alpha ) (i_0+ 1-\gamma +\beta )}{(i_0+ 3-\gamma )(i_0+ 2)}\frac{\left(\Delta_0^{-}\right)_{i_0} \left(\Delta_0^{+}\right)_{i_0} }{(2-\gamma )_{i_0}(1)_{i_0}} \sum_{i_1=i_0}^{\infty } \frac{\left(\Delta_1^{-}\right)_{i_1} \left(\Delta_1^{+} \right)_{i_1}(4-\gamma  )_{i_0}(3)_{i_0}}{\left(\Delta_1^{-}\right)_{i_0} \left(\Delta_1^{+}\right)_{i_0}(4-\gamma  )_{i_1}(3)_{i_1}} \eta ^{i_1}\right\} z\nonumber\\
&+& \sum_{n=2}^{\infty } \left\{ \sum_{i_0=0}^{\infty } \frac{(i_0+ 1-\gamma +\alpha ) (i_0+ 1-\gamma +\beta )}{(i_0+ 3-\gamma )(i_0+ 2)}\frac{\left(\Delta_0^{-}\right)_{i_0} \left(\Delta_0^{+}\right)_{i_0} }{(2-\gamma )_{i_0}(1)_{i_0}}\right.\nonumber\\
&\times& \prod _{k=1}^{n-1} \left\{ \sum_{i_k=i_{k-1}}^{\infty } \frac{(i_k+ 2k+1-\gamma +\alpha ) (i_k+ 2k+1-\gamma +\beta )}{(i_k+ 2k+3-\gamma  )(i_k+ 2(k+1))}\frac{\left(\Delta_k^{-}\right)_{i_k}\left( \Delta_k^{+} \right)_{i_k} (2(k+1)-\gamma )_{i_{k-1}}(2k+1)_{i_{k-1}}}{\left( \Delta_k^{-}\right)_{i_{k-1}}\left( \Delta_k^{+}\right)_{i_{k-1}}(2(k+1)-\gamma )_{i_k}(2k+1)_{i_k}} \right\}\nonumber\\
&\times&  \left.\left.\sum_{i_n= i_{n-1}}^{\infty }\frac{\left( \Delta_n^{-}\right)_{i_n} \left( \Delta_n^{+} \right)_{i_n}(2(n+1)-\gamma )_{i_{n-1}}(2n+1 )_{i_{n-1}}}{\left( \Delta_n^{-}\right)_{i_{n-1}} \left( \Delta_n^{+}\right)_{i_{n-1}}(2(n+1)-\gamma )_{i_n}(2n+1)_{i_n}}\eta ^{i_n} \right\} z^n \right\} \label{eq:10022} 
\end{eqnarray}
\normalsize
where
\Large
\begin{equation}
\begin{cases} 
\Delta_0^{\pm}= \frac{\{\varphi +2(1+a)(1-\gamma ) \} \pm\sqrt{\varphi ^2-4(1+a)q}}{2(1+a)} \cr
\Delta_1^{\pm}=  \frac{\{\varphi +2(1+a)(3-\gamma  ) \} \pm\sqrt{\varphi ^2-4(1+a)q}}{2(1+a)} \cr
\Delta_k^{\pm}=  \frac{\{\varphi +2(1+a)(1+2k-\gamma  )\} \pm\sqrt{\varphi ^2-4(1+a)q}}{2(1+a)} \cr
\Delta_n^{\pm}=   \frac{\{\varphi +2(1+a)(1+2n-\gamma )\} \pm\sqrt{\varphi ^2-4(1+a)q}}{2(1+a)}
\end{cases}\nonumber 
\end{equation}
\normalsize
\end{remark}
It's required that $\gamma \ne 0,-1,-2,\cdots$ for the first kind of independent solution of Heun function for all cases. Because if it does not, its solution will be divergent. And it's required that $\gamma \ne 2,3,4,\cdots$ for the second kind of independent solution of Heun function for all cases.

The infinite series in this chapter are equivalent to the infinite series in Ref.\cite{Chou2012H11}. In this chapter $B_n$ is the leading term in sequence $c_n$ in the analytic function $y(x)$. In Ref.\cite{Chou2012H11} $A_n$ is the leading term in sequence $c_n$ in the analytic function $y(x)$.
\section{Integral Formalism}
\subsection{Polynomial of type 2}

Now let's investigate the integral formalism for the polynomial of type 2 at certain eigenvalue. There is a generalized hypergeometric function which is
\begin{eqnarray}
I_l &=& \sum_{i_l= i_{l-1}}^{q_l} \frac{(-q_l)_{i_l}\left( q_l+4l+ \frac{\varphi +2(1+a)\lambda }{(1+a)}\right)_{i_l}(2l+1+\lambda )_{i_{l-1}}(2l+\gamma +\lambda )_{i_{l-1}}}{(-q_l)_{i_{l-1}}\left( q_l+4l+ \frac{\varphi +2(1+a)\lambda }{(1+a)}\right)_{i_{l-1}}(2l+1+\lambda )_{i_l}(2l+\gamma +\lambda )_{i_l}} \eta^{i_l}\nonumber\\
&=& \eta ^{i_{l-1}} 
\sum_{j=0}^{\infty } \frac{B(i_{l-1}+2l+\lambda ,j+1) B(i_{l-1}+2l-1+\gamma +\lambda ,j+1)}{(i_{l-1}+2l+\lambda )^{-1}(i_{l-1}+2l-1+\gamma +\lambda )^{-1}} \nonumber\\
&&\times \frac{(i_{l-1}-q_l)_j \left( i_{l-1}+q_l+4l+\frac{\varphi +2(1+a)\lambda }{(1+a)}\right)_j}{(1)_j \;j!} \eta ^j\label{eq:10030}
\end{eqnarray}
By using integral form of beta function,
\begin{subequations}
\begin{equation}
B\left(i_{l-1}+2l+\lambda ,j+1\right)= \int_{0}^{1} dt_l\;t_l^{i_{l-1}+2l-1+\lambda } (1-t_l)^j \label{eq:10031a}
\end{equation}
\begin{equation}
B\left(i_{l-1}+2l-1+\gamma +\lambda ,j+1\right)= \int_{0}^{1} du_l\;u_l^{i_{l-1}+2(l-1)+\gamma +\lambda } (1-u_l)^j\label{eq:10031b}
\end{equation}
\end{subequations}
Substitute (\ref{eq:10031a}) and (\ref{eq:10031b}) into (\ref{eq:10030}). And divide $(i_{l-1}+2l+\lambda )(i_{l-1}+2l-1+\gamma +\lambda )$ into $I_l$.
\begin{eqnarray}
K_l&=& \frac{1}{(i_{l-1}+2l+\lambda )(i_{l-1}+2l-1+\gamma +\lambda )}\nonumber\\
&&\times \sum_{i_l= i_{l-1}}^{q_l} \frac{(-q_l)_{i_l}\left( q_l+4l+ \frac{\varphi +2(1+a)\lambda }{(1+a)}\right)_{i_l}(2l+1+\lambda )_{i_{l-1}}(2l+\gamma +\lambda )_{i_{l-1}}}{(-q_l)_{i_{l-1}}\left( q_l+4l+ \frac{\varphi +2(1+a)\lambda }{(1+a)}\right)_{i_{l-1}}(2l+1+\lambda )_{i_l}(2l+\gamma +\lambda )_{i_l}} \eta^{i_l}\nonumber\\
&=&  \int_{0}^{1} dt_l\;t_l^{2l-1+\lambda } \int_{0}^{1} du_l\;u_l^{2(l-1)+\gamma +\lambda } (\eta t_l u_l)^{i_{l-1}}\nonumber\\
&&\times \sum_{j=0}^{\infty } \frac{(i_{l-1}-q_l)_j \left( i_{l-1}+q_l+4l+\frac{\varphi +2(1+a)\lambda }{(1+a)}\right)_j}{(1)_j \;j!} [\eta (1-t_l)(1-u_l)]^j \nonumber 
\end{eqnarray}
The integral form of Gauss hypergeometric function is
\begin{eqnarray}
_2F_1 \left( \alpha ,\beta ; \gamma ; z \right) &=& \sum_{n=0}^{\infty } \frac{(\alpha )_n (\beta )_n}{(\gamma )_n (n!)} z^n \nonumber\\
&=& -\frac{1}{2\pi i} \frac{\Gamma(1-\alpha ) \Gamma(\gamma )}{\Gamma (\gamma -\alpha )} \oint dv_l\;(-v_l)^{\alpha -1} (1-v_l)^{\gamma -\alpha -1} (1-zv_l)^{-\beta }\hspace{1cm}\label{eq:10033}\\
&& \mbox{where} \;\mbox{Re}(\gamma -\alpha )>0 \nonumber
\end{eqnarray}
Replace $\alpha $, $\beta $, $\gamma $ and $z$ by $i_{l-1}-q_l$, $ { \displaystyle i_{l-1}+q_l+4l+\frac{\varphi +2(1+a)\lambda }{(1+a)}}$, 1 and $\eta (1-t_l)(1-u_l)$ in (\ref{eq:10033}).
\begin{eqnarray}
&& \sum_{j=0}^{\infty } \frac{\left(i_{l-1}-q_l)_j (i_{l-1}+q_l+4l+\frac{\varphi +2(1+a)\lambda }{(1+a)}\right)_j}{(1)_j \;j!} [\eta (1-t_l)(1-u_l)]^j \nonumber\\
&=& \frac{1}{2\pi i} \oint dv_l\;\frac{1}{v_l} \left(\frac{v_l-1}{v_l}\frac{1}{1-\eta (1-t_l)(1-u_l)v_l}\right)^{q_l} (1-\eta (1-t_l)(1-u_l)v_l)^{-\left( 4l+ \frac{\varphi +2(1+a)\lambda }{(1+a)}\right)}\nonumber\\
&&\times \left(\frac{v_l}{v_l-1} \frac{1}{1-\eta (1-t_l)(1-u_l)v_l}\right)^{i_{l-1}} \label{eq:10034}
\end{eqnarray}
Substitute (\ref{eq:10034}) into $K_l$.
\begin{eqnarray}
K_l&=& \frac{1}{(i_{l-1}+2l+\lambda )(i_{l-1}+2l-1+\gamma +\lambda )}\nonumber\\
&&\times \sum_{i_l= i_{l-1}}^{q_l} \frac{(-q_l)_{i_l}\left( q_l+4l+ \frac{\varphi +2(1+a)\lambda }{(1+a)}\right)_{i_l}(2l+1+\lambda )_{i_{l-1}}(2l+\gamma +\lambda )_{i_{l-1}}}{(-q_l)_{i_{l-1}}\left( q_l+4l+ \frac{\varphi +2(1+a)\lambda }{(1+a)}\right)_{i_{l-1}}(2l+1+\lambda )_{i_l}(2l+\gamma +\lambda )_{i_l}} \eta^{i_l}\nonumber\\
&=&  \int_{0}^{1} dt_l\;t_l^{2l-1+\lambda } \int_{0}^{1} du_l\;u_l^{2(l-1)+\gamma +\lambda } 
\frac{1}{2\pi i} \oint dv_l\;\frac{1}{v_l} \left(\frac{v_l-1}{v_l}\frac{1}{1-\eta (1-t_l)(1-u_l)v_l}\right)^{q_l} \nonumber\\
&&\times (1-\eta (1-t_l)(1-u_l)v_l)^{-\left( 4l+ \frac{\varphi +2(1+a)\lambda }{(1+a)}\right)} \left(\frac{v_l}{v_l-1} \frac{\eta t_l u_l}{1-\eta (1-t_l)(1-u_l)v_l}\right)^{i_{l-1}}\label{eq:10035} 
\end{eqnarray}
Substitute (\ref{eq:10035}) into (\ref{eq:1007}) where $l=1,2,3,\cdots$; apply $K_1$ into the second summation of sub-power series $y_1(x)$, apply $K_2$ into the third summation and $K_1$ into the second summation of sub-power series $y_2(x)$, apply $K_3$ into the forth summation, $K_2$ into the third summation and $K_1$ into the second summation of sub-power series $y_3(x)$, etc.\footnote{$y_1(x)$ means the sub-power series in (\ref{eq:1007}) contains one term of $B_n's$, $y_2(x)$ means the sub-power series in (\ref{eq:1007}) contains two terms of $B_n's$, $y_3(x)$ means the sub-power series in (\ref{eq:1007}) contains three terms of $B_n's$, etc.}
\begin{theorem}
The general representation in the form of integral of Heun polynomial of type 2 is given by
\begin{eqnarray}
 y(x)&=& \sum_{n=0}^{\infty } y_{n}(x)= y_0(x)+ y_1(x)+ y_2(x)+y_3(x)+\cdots \nonumber\\
&=& c_0 x^{\lambda } \Bigg\{ \sum_{i_0=0}^{q_0}\frac{(-q_0)_{i_0}\left(q_0+\frac{\varphi +2(1+a)\lambda }{(1+a)}\right)_{i_0}}{(1+\lambda )_{i_0}(\gamma +\lambda )_{i_0}}  \eta ^{i_0}\nonumber\\
&&+ \sum_{n=1}^{\infty } \Bigg\{\prod _{k=0}^{n-1} \Bigg\{ \int_{0}^{1} dt_{n-k}\;t_{n-k}^{2(n-k)-1+\lambda } \int_{0}^{1} du_{n-k}\;u_{n-k}^{2(n-k-1)+\gamma +\lambda } \nonumber\\
&&\times  \frac{1}{2\pi i}  \oint dv_{n-k} \frac{1}{v_{n-k}} \left( \frac{v_{n-k}-1}{v_{n-k}} \frac{1}{1-\overleftrightarrow {w}_{n-k+1,n}(1-t_{n-k})(1-u_{n-k})v_{n-k}}\right)^{q_{n-k}} \nonumber\\
&&\times \left( 1- \overleftrightarrow {w}_{n-k+1,n}(1-t_{n-k})(1-u_{n-k})v_{n-k}\right)^{-\left(4(n-k)+\frac{\varphi +2(1+a)\lambda }{(1+a)}\right)}\nonumber\\
&&\times \overleftrightarrow {w}_{n-k,n}^{-(2(n-k-1)+\alpha +\lambda )}\left(  \overleftrightarrow {w}_{n-k,n} \partial _{ \overleftrightarrow {w}_{n-k,n}}\right) \overleftrightarrow {w}_{n-k,n}^{\alpha -\beta} \left(  \overleftrightarrow {w}_{n-k,n} \partial _{ \overleftrightarrow {w}_{n-k,n}}\right) \overleftrightarrow {w}_{n-k,n}^{2(n-k-1)+\beta +\lambda } \Bigg\}\nonumber\\
&&\times \sum_{i_0=0}^{q_0}\frac{(-q_0)_{i_0}\left(q_0+\frac{\varphi +2(1+a)\lambda }{(1+a)}\right)_{i_0}}{(1+\lambda )_{i_0}(\gamma +\lambda )_{i_0}} \overleftrightarrow {w}_{1,n}^{i_0}\Bigg\} z^n \Bigg\} \label{eq:10039}
\end{eqnarray}
where
\begin{equation}\overleftrightarrow {w}_{i,j}=
\begin{cases} \displaystyle {\frac{v_i}{(v_i-1)}\; \frac{\overleftrightarrow w_{i+1,j} t_i u_i}{1- \overleftrightarrow w_{i+1,j} v_i (1-t_i)(1-u_i)}} \;\;\mbox{where}\; i\leq j\cr
\eta \;\;\mbox{only}\;\mbox{if}\; i>j
\end{cases}\nonumber 
\end{equation}
In the above, the first sub-integral form contains one term of $B_n's$, the second one contains two terms of $B_n$'s, the third one contains three terms of $B_n$'s, etc.
\end{theorem}
\begin{proof}
According to (\ref{eq:1007}), 
\begin{equation}
 y(x)= \sum_{n=0}^{\infty }y_n(x) = y_0(x)+ y_1(x)+ y_2(x)+y_3(x)+\cdots \label{eq:10041}
\end{equation}
In the above, sub-power series $y_0(x) $, $y_1(x)$, $y_2(x)$ and $y_3(x)$ of Heun polynomial which makes $A_n$ term terminated  about $x=0$ are
\begin{subequations}
\begin{equation}
 y_0(x)= c_0 x^{\lambda } \sum_{i_0=0}^{q_0} \frac{(-q_0)_{i_0} \left(q_0+ \frac{\varphi +2(1+a)\lambda }{(1+a)}\right)_{i_0}}{(1+\lambda )_{i_0}(\gamma +\lambda )_{i_0}} \eta ^{i_0} \label{eq:10042a}
\end{equation}
\begin{eqnarray}
 &y_1(x)&= c_0 x^{\lambda } \Bigg\{\sum_{i_0=0}^{q_0}\frac{(i_0+ \lambda +\alpha ) (i_0+ \lambda +\beta )}{(i_0+ \lambda +2)(i_0+ \lambda +1+\gamma )}\frac{(-q_0)_{i_0} \left(q_0+ \frac{\varphi +2(1+a)\lambda }{(1+a)}\right)_{i_0}}{(1+\lambda )_{i_0}(\gamma +\lambda )_{i_0}} \nonumber\\
&\times& \sum_{i_1=i_0}^{q_1} \frac{(-q_1)_{i_1}\left(q_1+4+ \frac{\varphi +2(1+a)\lambda }{(1+a)}\right)_{i_1}(3+\lambda )_{i_0}(2+\gamma +\lambda )_{i_0}}{(-q_1)_{i_0}\left(q_1+4+ \frac{\varphi +2(1+a)\lambda }{(1+a)}\right)_{i_0}(3+\lambda )_{i_1}(2+\gamma +\lambda )_{i_1}} \eta ^{i_1} \Bigg\}z  \label{eq:10042b}
\end{eqnarray}
\begin{eqnarray}
 &y_2(x)&= c_0 x^{\lambda } \Bigg\{\sum_{i_0=0}^{q_0}\frac{(i_0+ \lambda +\alpha ) (i_0+ \lambda +\beta )}{(i_0+ \lambda +2)(i_0+ \lambda +1+\gamma )}\frac{(-q_0)_{i_0} \left(q_0+ \frac{\varphi +2(1+a)\lambda }{(1+a)}\right)_{i_0}}{(1+\lambda )_{i_0}(\gamma +\lambda )_{i_0}} \nonumber\\
&\times&  \sum_{i_1=i_0}^{q_1} \frac{(i_1+2+ \lambda +\alpha ) (i_1+2+ \lambda +\beta )}{(i_1+ \lambda +4)(i_1+ \lambda +3+\gamma )} \frac{(-q_1)_{i_1}\left(q_1+4+ \frac{\varphi +2(1+a)\lambda }{(1+a)}\right)_{i_1}(3+\lambda )_{i_0}(2+\gamma +\lambda )_{i_0}}{(-q_1)_{i_0}\left(q_1+4+ \frac{\varphi +2(1+a)\lambda }{(1+a)}\right)_{i_0}(3+\lambda )_{i_1}(2+\gamma +\lambda )_{i_1}} \nonumber\\
&\times& \sum_{i_2=i_1}^{q_2} \frac{(-q_2)_{i_2}\left(q_2+8+ \frac{\varphi +2(1+a)\lambda }{(1+a)}\right)_{i_2}(5+\lambda )_{i_1}(4+\gamma +\lambda )_{i_1}}{(-q_2)_{i_1}\left(q_2+8+ \frac{\varphi +2(1+a)\lambda }{(1+a)}\right)_{i_1}(5+\lambda )_{i_2}(4+\gamma +\lambda )_{i_2}} \eta ^{i_2} \Bigg\} z^2  \label{eq:10042c}
\end{eqnarray}
\begin{eqnarray}
&y_3(x)&=  c_0 x^{\lambda } \Bigg\{\sum_{i_0=0}^{q_0}\frac{(i_0+ \lambda +\alpha ) (i_0+ \lambda +\beta )}{(i_0+ \lambda +2)(i_0+ \lambda +1+\gamma )}\frac{(-q_0)_{i_0} \left(q_0+ \frac{\varphi +2(1+a)\lambda }{(1+a)}\right)_{i_0}}{(1+\lambda )_{i_0}(\gamma +\lambda )_{i_0}} \nonumber\\
&\times&  \sum_{i_1=i_0}^{q_1} \frac{(i_1+2+ \lambda +\alpha ) (i_1+2+ \lambda +\beta )}{(i_1+ \lambda +4)(i_1+ \lambda +3+\gamma )} \frac{(-q_1)_{i_1}\left(q_1+4+ \frac{\varphi +2(1+a)\lambda }{(1+a)}\right)_{i_1}(3+\lambda )_{i_0}(2+\gamma +\lambda )_{i_0}}{(-q_1)_{i_0}\left(q_1+4+ \frac{\varphi +2(1+a)\lambda }{(1+a)}\right)_{i_0}(3+\lambda )_{i_1}(2+\gamma +\lambda )_{i_1}} \nonumber\\
&\times& \sum_{i_2=i_1}^{q_2} \frac{(i_2+4+ \lambda +\alpha ) (i_2+4+ \lambda +\beta )}{(i_2+ \lambda +6)(i_2+ \lambda +5+\gamma )}  \frac{(-q_2)_{i_2}\left(q_2+8+ \frac{\varphi +2(1+a)\lambda }{(1+a)}\right)_{i_2}(5+\lambda )_{i_1}(4+\gamma +\lambda )_{i_1}}{(-q_2)_{i_1}\left(q_2+8+ \frac{\varphi +2(1+a)\lambda }{(1+a)}\right)_{i_1}(5+\lambda )_{i_2}(4+\gamma +\lambda )_{i_2}} \nonumber\\
&\times& \sum_{i_3=i_2}^{q_3} \frac{(-q_3)_{i_3}\left(q_2+12+ \frac{\varphi +2(1+a)\lambda }{(1+a)}\right)_{i_3}(7+\lambda )_{i_2}(6+\gamma +\lambda )_{i_2}}{(-q_3)_{i_2}\left(q_2+12+ \frac{\varphi +2(1+a)\lambda }{(1+a)}\right)_{i_2}(7+\lambda )_{i_3}(6+\gamma +\lambda )_{i_3}}\eta ^{i_3} \Bigg\} z^3  \label{eq:10042d} 
\end{eqnarray}
\end{subequations}
Put $l=1$ in (\ref{eq:10035}). Take the new (\ref{eq:10035}) into (\ref{eq:10042b}).
\begin{eqnarray}
y_1(x) &=& \int_{0}^{1} dt_1\;t_1^{1+\lambda } \int_{0}^{1} du_1\;u_1^{\gamma +\lambda } \frac{1}{2\pi i} \oint dv_1 \;\frac{1}{v_1} 
\left( \frac{v_1-1}{v_1} \frac{1}{1-\eta (1-t_1)(1-u_1)v_1}\right)^{q_1}  \nonumber\\
&&\times (1-\eta (1-t_1)(1-u_1)v_1)^{-\left( 4+\frac{\varphi +2(1+a)\lambda }{(1+a)}\right)} \nonumber\\
&&\times \overleftrightarrow {w}_{1,1}^{-(\alpha +\lambda )} \left(\overleftrightarrow {w}_{1,1} \partial_{\overleftrightarrow {w}_{1,1}} \right) \overleftrightarrow {w}_{1,1}^{\alpha -\beta } \left(\overleftrightarrow {w}_{1,1} \partial_{\overleftrightarrow {w}_{1,1}} \right)\overleftrightarrow {w}_{1,1}^{\beta +\lambda } \nonumber\\
&&\times \left\{ c_0 x^{\lambda } \sum_{i_0=0}^{q_0} \frac{(-q_0)_{i_0} \left(q_0+ \frac{\varphi +2(1+a)\lambda }{(1+a)}\right)_{i_0}}{(1+\lambda )_{i_0}(\gamma +\lambda )_{i_0}} \overleftrightarrow {w}_{1,1} ^{i_0}\right\}z \label{eq:10043}
\end{eqnarray}
where
\begin{equation}
\overleftrightarrow {w}_{1,1} = \frac{v_1}{v_1-1} \frac{\eta t_1 u_1}{1-\eta (1-t_1)(1-u_1)v_1}\nonumber
\end{equation}
Put $l=2$ in (\ref{eq:10035}). Take the new (\ref{eq:10035}) into (\ref{eq:10042c}).
\begin{eqnarray}
y_2(x) &=& c_0 x^{\lambda } \int_{0}^{1} dt_2\;t_2^{3+\lambda } \int_{0}^{1} du_2\;u_2^{2+\gamma +\lambda } \frac{1}{2\pi i} \oint dv_2 \;\frac{1}{v_2} 
\left( \frac{v_2-1}{v_2} \frac{1}{1-\eta (1-t_2)(1-u_2)v_2}\right)^{q_2}  \nonumber\\
&&\times (1-\eta (1-t_2)(1-u_2)v_2)^{-\left( 8+\frac{\varphi +2(1+a)\lambda }{(1+a)}\right)} \nonumber\\
&&\times \overleftrightarrow {w}_{2,2}^{-(2+\alpha +\lambda) } \left(\overleftrightarrow {w}_{2,2} \partial_{\overleftrightarrow {w}_{2,2}} \right) \overleftrightarrow {w}_{2,2}^{\alpha -\beta } \left(\overleftrightarrow {w}_{2,2} \partial_{\overleftrightarrow {w}_{2,2}} \right)\overleftrightarrow {w}_{2,2}^{2+\beta +\lambda } \nonumber\\
&&\times \Bigg\{\sum_{i_0=0}^{q_0}\frac{(i_0+ \lambda +\alpha ) (i_0+ \lambda +\beta )}{(i_0+ \lambda +2)(i_0+ \lambda +1+\gamma )}\frac{(-q_0)_{i_0} \left(q_0+ \frac{\varphi +2(1+a)\lambda }{(1+a)}\right)_{i_0}}{(1+\lambda )_{i_0}(\gamma +\lambda )_{i_0}} \nonumber\\
&&\times  \sum_{i_1=i_0}^{q_1} \frac{(-q_1)_{i_1}\left(q_1+4+ \frac{\varphi +2(1+a)\lambda }{(1+a)}\right)_{i_1}(3+\lambda )_{i_0}(2+\gamma +\lambda )_{i_0}}{(-q_1)_{i_0}\left(q_1+4+ \frac{\varphi +2(1+a)\lambda }{(1+a)}\right)_{i_0}(3+\lambda )_{i_1}(2+\gamma +\lambda )_{i_1}} \overleftrightarrow {w}_{2,2}^{i_1} \Bigg\} z^2 \label{eq:10044}
\end{eqnarray}
where
\begin{equation}
\overleftrightarrow {w}_{2,2} = \frac{v_2}{v_2-1} \frac{\eta t_2 u_2}{1-\eta (1-t_2)(1-u_2)v_2}\nonumber
\end{equation}
Put $l=1$ and $\eta = \overleftrightarrow {w}_{2,2}$ in (\ref{eq:35}). Take the new (\ref{eq:35}) into (\ref{eq:44}).
\begin{eqnarray}
y_2(x) &=& \int_{0}^{1} dt_2\;t_2^{3+\lambda } \int_{0}^{1} du_2\;u_2^{2+\gamma +\lambda } \frac{1}{2\pi i} \oint dv_2 \;\frac{1}{v_2} 
\left( \frac{v_2-1}{v_2} \frac{1}{1-\eta (1-t_2)(1-u_2)v_2}\right)^{q_2}  \nonumber\\
&&\times (1-\eta (1-t_2)(1-u_2)v_2)^{-\left( 8+\frac{\varphi +2(1+a)\lambda }{(1+a)}\right)} \nonumber\\
&&\times \overleftrightarrow {w}_{2,2}^{-(2+\alpha +\lambda) } \left(\overleftrightarrow {w}_{2,2} \partial_{\overleftrightarrow {w}_{2,2}} \right) \overleftrightarrow {w}_{2,2}^{\alpha -\beta } \left(\overleftrightarrow {w}_{2,2} \partial_{\overleftrightarrow {w}_{2,2}} \right)\overleftrightarrow {w}_{2,2}^{2+\beta +\lambda } \nonumber\\
&&\times \int_{0}^{1} dt_1\;t_1^{1+\lambda } \int_{0}^{1} du_1\;u_1^{\gamma +\lambda } \frac{1}{2\pi i} \oint dv_1 \;\frac{1}{v_1} 
\left( \frac{v_1-1}{v_1} \frac{1}{1- \overleftrightarrow {w}_{2,2} (1-t_1)(1-u_1)v_1}\right)^{q_1}  \nonumber\\
&&\times (1- \overleftrightarrow {w}_{2,2} (1-t_1)(1-u_1)v_1)^{-\left( 4+\frac{\varphi +2(1+a)\lambda }{(1+a)}\right)} \nonumber\\
&&\times \overleftrightarrow {w}_{1,2}^{-(\alpha +\lambda)} \left(\overleftrightarrow {w}_{1,2} \partial_{\overleftrightarrow {w}_{1,2}} \right) \overleftrightarrow {w}_{1,2}^{\alpha -\beta } \left(\overleftrightarrow {w}_{1,2} \partial_{\overleftrightarrow {w}_{1,2}} \right)\overleftrightarrow {w}_{1,2}^{\beta +\lambda } \nonumber\\
&&\times \left\{ c_0 x^{\lambda } \sum_{i_0=0}^{q_0} \frac{(-q_0)_{i_0} \left(q_0+ \frac{\varphi +2(1+a)\lambda }{(1+a)}\right)_{i_0}}{(1+\lambda )_{i_0}(\gamma +\lambda )_{i_0}} \overleftrightarrow {w}_{1,2} ^{i_0}\right\} z^2 \label{eq:10045}
\end{eqnarray}
where
\begin{equation}
\overleftrightarrow {w}_{1,2} = \frac{v_1}{v_1-1} \frac{\overleftrightarrow {w}_{2,2} t_1 u_1}{1- \overleftrightarrow {w}_{2,2}(1-t_1)(1-u_1)v_1}\nonumber
\end{equation}
By using similar process for the previous cases of integral forms of $y_1(x)$ and $y_2(x)$, the integral form of sub-power series expansion of $y_3(x)$ is
\begin{eqnarray}
y_3(x) &=& \int_{0}^{1} dt_3\;t_3^{5+\lambda } \int_{0}^{1} du_3\;u_3^{4+\gamma +\lambda } \frac{1}{2\pi i} \oint dv_3 \;\frac{1}{v_3} 
\left( \frac{v_3-1}{v_3} \frac{1}{1-\eta (1-t_3)(1-u_3)v_3}\right)^{q_3}  \nonumber\\
&&\times (1-\eta (1-t_3)(1-u_3)v_3)^{-\left( 12+\frac{\varphi +2(1+a)\lambda }{(1+a)}\right)} \nonumber\\
&&\times \overleftrightarrow {w}_{3,3}^{-(4+\alpha +\lambda) } \left(\overleftrightarrow {w}_{3,3} \partial_{\overleftrightarrow {w}_{3,3}} \right) \overleftrightarrow {w}_{3,3}^{\alpha -\beta } \left(\overleftrightarrow {w}_{3,3} \partial_{\overleftrightarrow {w}_{3,3}} \right)\overleftrightarrow {w}_{3,3}^{4+\beta +\lambda } \nonumber\\
&&\times \int_{0}^{1} dt_2\;t_2^{3+\lambda } \int_{0}^{1} du_2\;u_2^{2+\gamma +\lambda } \frac{1}{2\pi i} \oint dv_2 \;\frac{1}{v_2} 
\left( \frac{v_2-1}{v_2} \frac{1}{1- \overleftrightarrow {w}_{3,3} (1-t_2)(1-u_2)v_2}\right)^{q_2}  \nonumber\\
&&\times (1- \overleftrightarrow {w}_{3,3} (1-t_2)(1-u_2)v_2)^{-\left( 8+\frac{\varphi +2(1+a)\lambda }{(1+a)}\right)} \nonumber\\
&&\times \overleftrightarrow {w}_{2,3}^{-(2+\alpha +\lambda)} \left(\overleftrightarrow {w}_{2,3} \partial_{\overleftrightarrow {w}_{2,3}} \right) \overleftrightarrow {w}_{2,3}^{\alpha -\beta } \left(\overleftrightarrow {w}_{2,3} \partial_{\overleftrightarrow {w}_{2,3}} \right)\overleftrightarrow {w}_{2,3}^{2+\beta +\lambda } \nonumber
\end{eqnarray}
\begin{eqnarray}
&&\times \int_{0}^{1} dt_1\;t_1^{1+\lambda } \int_{0}^{1} du_1\;u_1^{\gamma +\lambda } \frac{1}{2\pi i} \oint dv_1 \;\frac{1}{v_1} 
\left( \frac{v_1-1}{v_1} \frac{1}{1- \overleftrightarrow {w}_{2,3} (1-t_1)(1-u_1)v_1}\right)^{q_1}  \nonumber\\
&&\times (1- \overleftrightarrow {w}_{2,3} (1-t_1)(1-u_1)v_1)^{-\left( 4+\frac{\varphi +2(1+a)\lambda }{(1+a)}\right)}\nonumber\\
&&\times \overleftrightarrow {w}_{1,3}^{-(\alpha +\lambda)} \left(\overleftrightarrow {w}_{1,3} \partial_{\overleftrightarrow {w}_{1,3}} \right) \overleftrightarrow {w}_{1,3}^{\alpha -\beta } \left(\overleftrightarrow {w}_{1,3} \partial_{\overleftrightarrow {w}_{1,3}} \right)\overleftrightarrow {w}_{1,3}^{\beta +\lambda } \nonumber\\
&&\times \left\{ c_0 x^{\lambda } \sum_{i_0=0}^{q_0} \frac{(-q_0)_{i_0} \left(q_0+ \frac{\varphi +2(1+a)\lambda }{(1+a)}\right)_{i_0}}{(1+\lambda )_{i_0}(\gamma +\lambda )_{i_0}} \overleftrightarrow {w}_{1,3} ^{i_0}\right\} z^3 \label{eq:10046}
\end{eqnarray}
where
\begin{equation}
\begin{cases} \overleftrightarrow {w}_{3,3} = \frac{v_3}{v_3-1} \frac{\eta  t_3 u_3}{1- \eta (1-t_3)(1-u_3)v_3} \cr
\overleftrightarrow {w}_{2,3} = \frac{v_2}{v_2-1} \frac{\overleftrightarrow {w}_{3,3} t_2 u_2}{1- \overleftrightarrow {w}_{3,3}(1-t_2)(1-u_2)v_2} \cr
\overleftrightarrow {w}_{1,3} = \frac{v_1}{v_1-1} \frac{\overleftrightarrow {w}_{2,3} t_1 u_1}{1- \overleftrightarrow {w}_{2,3}(1-t_1)(1-u_1)v_1}
\end{cases}
\nonumber
\end{equation}
By repeating this process for all higher terms of integral forms of sub-summation $y_m(x)$ terms where $m \geq 4$, I obtain every integral forms of $y_m(x)$ terms. 
Since we substitute (\ref{eq:10042a}), (\ref{eq:10043}), (\ref{eq:10045}), (\ref{eq:10046}) and including all integral forms of $y_m(x)$ terms where $m \geq 4$ into (\ref{eq:10041}), We obtain (\ref{eq:10039}).
\qed
\end{proof}
Put $c_0$= 1 as $\lambda $=0 for the first kind of independent solutions of Heun equation and $\displaystyle{ c_0= \left( \frac{1+a}{a}\right)^{1-\gamma }}$ as $\lambda = 1-\gamma $ for the second one in (\ref{eq:10039}).
\begin{remark}
The integral representation of Heun equation of the first kind for polynomial of type 2 about $x=0$ as $q= -(q_j+2j)\{ \alpha +\beta -\delta +a(\delta +\gamma -1)+(1+a)(q_j+2j)\}$ where $j,q_j \in \mathbb{N}_{0}$ is
\begin{eqnarray}
 y(x)&=& HF_{q_j}^R \left( q_j =\frac{-\varphi \pm \sqrt{\varphi ^2-4(1+a)q}}{2(1+a)}-2j, \varphi =\alpha +\beta -\delta +a(\delta +\gamma -1) \right.\nonumber\\
&&, \left. \Omega _1=\frac{\varphi }{(1+a)}; \eta = \frac{(1+a)}{a} x ; z= -\frac{1}{a} x^2 \right) \nonumber\\
&=& _2F_1 \left(-q_0, q_0+\Omega _1; \gamma; \eta \right) \nonumber\\
&&+ \sum_{n=1}^{\infty } \Bigg\{\prod _{k=0}^{n-1} \Bigg\{ \int_{0}^{1} dt_{n-k}\;t_{n-k}^{2(n-k)-1} \int_{0}^{1} du_{n-k}\;u_{n-k}^{2(n-k-1)+\gamma } \nonumber\\
&&\times  \frac{1}{2\pi i}  \oint dv_{n-k} \frac{1}{v_{n-k}} \left( \frac{v_{n-k}-1}{v_{n-k}} \frac{1}{1-\overleftrightarrow {w}_{n-k+1,n}(1-t_{n-k})(1-u_{n-k})v_{n-k}}\right)^{q_{n-k}} \nonumber\\
&&\times \left( 1- \overleftrightarrow {w}_{n-k+1,n}(1-t_{n-k})(1-u_{n-k})v_{n-k}\right)^{-\left(4(n-k)+\Omega _1\right)}\nonumber\\
&&\times \overleftrightarrow {w}_{n-k,n}^{-(2(n-k-1)+\alpha )}\left( \overleftrightarrow {w}_{n-k,n} \partial _{ \overleftrightarrow {w}_{n-k,n}}\right) \overleftrightarrow {w}_{n-k,n}^{\alpha -\beta} \left( \overleftrightarrow {w}_{n-k,n} \partial _{ \overleftrightarrow {w}_{n-k,n}}\right) \overleftrightarrow {w}_{n-k,n}^{2(n-k-1)+\beta } \Bigg\}\nonumber\\
&&\times _2F_1 \left(-q_0, q_0+\Omega _1; \gamma; \overleftrightarrow {w}_{1,n} \right) \Bigg\} z^n \label{eq:10047}
\end{eqnarray}
\end{remark}
\begin{remark}
The integral representation of Heun equation of the second kind for polynomial of type 2 about $x=0$ as $q= -(q_j+2j+1-\gamma )\{ \alpha +\beta+1-\gamma -(1-a) \delta +(1+a)(q_j+2j)\}$ where $j,q_j \in \mathbb{N}_{0}$ is
\begin{eqnarray}
y(x)&=& HS_{q_j}^R \Bigg( q_j =\frac{-\{\varphi +2(1+a)(1-\gamma )\} \pm \sqrt{\varphi ^2-4(1+a)q}}{2(1+a)}-2j\nonumber\\
&&,\varphi =\alpha +\beta -\delta +a(\delta +\gamma -1), \Omega _2= \frac{\varphi +2(1+a)(1-\gamma )}{(1+a)}; \eta = \frac{(1+a)}{a} x; z= -\frac{1}{a} x^2 \Bigg) \nonumber\\
&=& \eta^{1-\gamma } \Bigg\{\; _2F_1 \left(-q_0, q_0+\Omega _2; 2-\gamma; \eta \right)+ \sum_{n=1}^{\infty } \Bigg\{\prod _{k=0}^{n-1} \Bigg\{ \int_{0}^{1} dt_{n-k}\;t_{n-k}^{2(n-k)-\gamma } \int_{0}^{1} du_{n-k}\;u_{n-k}^{2(n-k)-1} \nonumber\\
&&\times  \frac{1}{2\pi i}  \oint dv_{n-k} \frac{1}{v_{n-k}} \left( \frac{v_{n-k}-1}{v_{n-k}} \frac{1}{1-\overleftrightarrow {w}_{n-k+1,n}(1-t_{n-k})(1-u_{n-k})v_{n-k}}\right)^{q_{n-k}} \nonumber\\
&&\times \left( 1- \overleftrightarrow {w}_{n-k+1,n}(1-t_{n-k})(1-u_{n-k})v_{n-k}\right)^{-\left(4(n-k)+\Omega _2\right)}\nonumber\\
&&\times \overleftrightarrow {w}_{n-k,n}^{-(2(n-k)-1+\alpha -\gamma)}\left(  \overleftrightarrow {w}_{n-k,n} \partial _{ \overleftrightarrow {w}_{n-k,n}}\right) \overleftrightarrow {w}_{n-k,n}^{\alpha -\beta} \left(  \overleftrightarrow {w}_{n-k,n} \partial _{ \overleftrightarrow {w}_{n-k,n}}\right) \overleftrightarrow {w}_{n-k,n}^{2(n-k)-1+\beta -\gamma} \Bigg\}\nonumber\\
&&\times _2F_1 \left(-q_0, q_0+\Omega _2; 2-\gamma; \overleftrightarrow {w}_{1,n}\right)  \Bigg\} z^n \Bigg\}\label{eq:10048}
\end{eqnarray}
\end{remark}
\subsection{Infinite series}
Let's consider the integral representation of Heun equation about $x=0$ for infinite series by applying R3TRF.
There is a generalized hypergeometric function which is written by
\begin{eqnarray}
M_l &=& \sum_{i_l= i_{l-1}}^{\infty } \frac{\left(\Delta_l^{-}\right)_{i_l}\left(\Delta_l^{+}\right)_{i_l}(2l+1+\lambda )_{i_{l-1}}(2l+\gamma +\lambda )_{i_{l-1}}}{\left(\Delta_l^{-} \right)_{i_{l-1}}\left(\Delta_l^{+} \right)_{i_{l-1}}(2l+1+\lambda )_{i_l}(2l+\gamma +\lambda )_{i_l}} \eta^{i_l}\label{er:10030}\\
&=& \eta ^{i_{l-1}} 
\sum_{j=0}^{\infty } \frac{B(i_{l-1}+2l+\lambda ,j+1) B(i_{l-1}+2l-1+\gamma +\lambda ,j+1)\left(\Delta_l^{-} +i_{l-1}\right)_j \left(\Delta_l^{+} +i_{l-1}\right)_j}{(i_{l-1}+2l+\lambda )^{-1}(i_{l-1}+2l-1+\gamma +\lambda )^{-1}(1)_j \;j!} \;\eta ^j \nonumber
\end{eqnarray}
where
\begin{equation} 
\Delta_l^{\pm} = \frac{ \varphi +2(1+a)(\lambda +2l ) \pm\sqrt{\varphi ^2-4(1+a)q}}{2(1+a)}
\nonumber 
\end{equation}
Substitute (\ref{eq:10031a}) and (\ref{eq:10031b}) into (\ref{er:10030}). Divide $(i_{l-1}+2l+\lambda ) (i_{l-1}+2l-1+\gamma +\lambda ) $ into the new (\ref{er:10030}).
\begin{eqnarray}
V_l&=& \frac{1}{(i_{l-1}+2l+\lambda )(i_{l-1}+2l-1+\gamma +\lambda )} \sum_{i_l= i_{l-1}}^{\infty } \frac{\left(\Delta_l^{-}\right)_{i_l}\left(\Delta_l^{+}\right)_{i_l}(2l+1+\lambda )_{i_{l-1}}(2l+\gamma +\lambda )_{i_{l-1}}}{\left(\Delta_l^{-} \right)_{i_{l-1}}\left(\Delta_l^{+} \right)_{i_{l-1}}(2l+1+\lambda )_{i_l}(2l+\gamma +\lambda )_{i_l}} \eta^{i_l}\nonumber\\
&=&  \int_{0}^{1} dt_l\;t_l^{2l-1+\lambda } \int_{0}^{1} du_l\;u_l^{2(l-1)+\gamma +\lambda } (\eta t_l u_l)^{i_{l-1}}  \sum_{j=0}^{\infty } \frac{ \left(\Delta_l^{-} +i_{l-1}\right)_j \left(\Delta_l^{+} +i_{l-1}\right)_j \left(\eta (1-t_l)(1-u_l)\right)^j}{(i_{l-1}+2l+\lambda )^{-1}(i_{l-1}+2l-1+\gamma +\lambda )^{-1}(1)_j \;j!}  \nonumber 
\end{eqnarray}
The hypergeometric function is defined by
\begin{eqnarray}
_2F_1 \left( \alpha ,\beta ; \gamma ; z \right) &=& \sum_{n=0}^{\infty } \frac{(\alpha )_n (\beta )_n}{(\gamma )_n (n!)} z^n \nonumber\\
&=&  \frac{1}{2\pi i} \frac{\Gamma( 1+\alpha  -\gamma )}{\Gamma (\alpha )} \int_0^{(1+)} dv_l\; (-1)^{\gamma }(-v_l)^{\alpha -1} (1-v_l )^{\gamma -\alpha -1} (1-zv_l)^{-\beta }\hspace{1.5cm}\label{er:10031}\\
&& \mbox{where} \;\gamma -\alpha  \ne 1,2,3,\cdots, \;\mbox{Re}(\alpha )>0 \nonumber
\end{eqnarray}
Replace $\alpha $, $\beta $, $\gamma $ and $z$ by $\Delta_l^{-} +i_{l-1}$, $\Delta_l^{+} +i_{l-1}$, 1 and $\eta(1-t_l)(1-u_l)$ in (\ref{er:10031}). Take the new (\ref{er:10031}) into $V_l$.
\begin{eqnarray}
V_l&=& \frac{1}{(i_{l-1}+2l+\lambda )(i_{l-1}+2l-1+\gamma +\lambda )} \sum_{i_l= i_{l-1}}^{\infty } \frac{\left(\Delta_l^{-}\right)_{i_l}\left(\Delta_l^{+}\right)_{i_l}(2l+1+\lambda )_{i_{l-1}}(2l+\gamma +\lambda )_{i_{l-1}}}{\left(\Delta_l^{-} \right)_{i_{l-1}}\left(\Delta_l^{+} \right)_{i_{l-1}}(2l+1+\lambda )_{i_l}(2l+\gamma +\lambda )_{i_l}} \eta^{i_l}\nonumber\\
&=&  \int_{0}^{1} dt_l\;t_l^{2l-1+\lambda } \int_{0}^{1} du_l\;u_l^{2(l-1)+\gamma +\lambda } 
\frac{1}{2\pi i} \oint dv_l\;\frac{1}{v_l} \left( \frac{v_l-1}{v_l} \right)^{-\Delta_l^{-}} \nonumber\\
&&\times (1-\eta (1-t_l)(1-u_l)v_l)^{-\Delta_l^{+}} \left(\frac{v_l}{v_l-1} \frac{\eta t_l u_l}{1-\eta (1-t_l)(1-u_l)v_l}\right)^{i_{l-1}}\label{er:10032} 
\end{eqnarray}
Substitute (\ref{er:10032}) into (\ref{eq:10020}) where $l=1,2,3,\cdots$; apply $V_1$ into the second summation of sub-power series $y_1(x)$, apply $V_2$ into the third summation and $V_1$ into the second summation of sub-power series $y_2(x)$, apply $V_3$ into the forth summation, $V_2$ into the third summation and $V_1$ into the second summation of sub-power series $y_3(x)$, etc.\footnote{$y_1(x)$ means the sub-power series in (\ref{eq:10020}) contains one term of $B_n's$, $y_2(x)$ means the sub-power series in (\ref{eq:10020}) contains two terms of $B_n's$, $y_3(x)$ means the sub-power series in (\ref{eq:10020}) contains three terms of $B_n's$, etc.}
\begin{theorem}
The general representation in the form of integral of Heun equation for infinite series about $x=0$ using R3TRF is given by 
\begin{eqnarray}
 y(x)&=& \sum_{n=0}^{\infty } y_{n}(x)= y_0(x)+ y_1(x)+ y_2(x)+y_3(x)+\cdots \nonumber\\
&=& c_0 x^{\lambda } \left\{ \sum_{i_0=0}^{\infty }\frac{\left( \Delta_0^{-}\right)_{i_0}\left( \Delta_0^{+}\right)_{i_0}}{(1+\lambda )_{i_0}(\gamma +\lambda )_{i_0}}  \eta ^{i_0}\right.\nonumber\\
&+& \sum_{n=1}^{\infty } \left\{\prod _{k=0}^{n-1} \Bigg\{ \int_{0}^{1} dt_{n-k}\;t_{n-k}^{2(n-k)-1+\lambda } \int_{0}^{1} du_{n-k}\;u_{n-k}^{2(n-k-1)+\gamma +\lambda } \right.\nonumber\\
&\times& \frac{1}{2\pi i}  \oint dv_{n-k} \frac{1}{v_{n-k}}\left( \frac{v_{n-k}-1}{v_{n-k}}\right)^{ -\Delta_{n-k}^{-}}  \left( 1- \overleftrightarrow {w}_{n-k+1,n}(1-t_{n-k})(1-u_{n-k})v_{n-k}\right)^{ -\Delta_{n-k}^{+}}\nonumber\\
&\times& \overleftrightarrow {w}_{n-k,n}^{-(2(n-k-1)+\alpha +\lambda )}\left(  \overleftrightarrow {w}_{n-k,n} \partial _{ \overleftrightarrow {w}_{n-k,n}}\right) \overleftrightarrow {w}_{n-k,n}^{\alpha -\beta} \left(  \overleftrightarrow {w}_{n-k,n} \partial _{ \overleftrightarrow {w}_{n-k,n}}\right) \overleftrightarrow {w}_{n-k,n}^{2(n-k-1)+\beta +\lambda } \Bigg\}\nonumber\\
&\times& \left.\left.\sum_{i_0=0}^{\infty }\frac{\left( \Delta_0^{-}\right)_{i_0}\left( \Delta_0^{+}\right)_{i_0}}{(1+\lambda )_{i_0}(\gamma +\lambda )_{i_0}} \overleftrightarrow {w}_{1,n}^{i_0}\right\} z^n \right\} \label{eq:10049}
\end{eqnarray}
where
\begin{equation}
\begin{cases}
\Delta_0^{\pm}=  \frac{ \varphi +2(1+a)\lambda \pm\sqrt{\varphi ^2-4(1+a)q}}{2(1+a)} \cr
\Delta_{n-k}^{\pm}=  \frac{ \varphi +2(1+a)(\lambda+2(n-k) ) \pm\sqrt{\varphi ^2-4(1+a)q}}{2(1+a)}
\end{cases}\nonumber 
\end{equation}
In the above, the first sub-integral form contains one term of $B_n's$, the second one contains two terms of $B_n$'s, the third one contains three terms of $B_n$'s, etc.
\end{theorem}
\begin{proof}
In (\ref{eq:10020}) sub-power series $y_0(x) $, $y_1(x)$, $y_2(x)$ and $y_3(x)$ of Heun equation for infinite series about $x=0$ using R3TRF are given by
\begin{subequations}
\begin{equation}
 y_0(x)= c_0 x^{\lambda } \sum_{i_0=0}^{\infty } \frac{\left( \Delta_{0}^{-}\right)_{i_0} \left(\Delta_{0}^{+}\right)_{i_0}}{(1+\lambda )_{i_0}(\gamma +\lambda )_{i_0}} \eta ^{i_0} \label{er:10033a}
\end{equation}
\begin{eqnarray}
 y_1(x) &=& c_0 x^{\lambda } \left\{\sum_{i_0=0}^{\infty }\frac{(i_0+ \lambda +\alpha ) (i_0+ \lambda +\beta )}{(i_0+ \lambda +2)(i_0+ \lambda +1+\gamma )}\frac{\left( \Delta_{0}^{-}\right)_{i_0} \left(\Delta_{0}^{+}\right)_{i_0}}{(1+\lambda )_{i_0}(\gamma +\lambda )_{i_0}} \right. \nonumber\\
&&\times \left. \sum_{i_1=i_0}^{\infty } \frac{\left( \Delta_{1}^{-}\right)_{i_1}\left( \Delta_{1}^{+}\right)_{i_1}(3+\lambda )_{i_0}(2+\gamma +\lambda )_{i_0}}{\left(\Delta_{1}^{-}\right)_{i_0}\left(\Delta_{1}^{+}\right)_{i_0}(3+\lambda )_{i_1}(2+\gamma +\lambda )_{i_1}} \eta ^{i_1} \right\} z  \label{er:10033b}
\end{eqnarray}
\begin{eqnarray}
 y_2(x) &=& c_0 x^{\lambda } \left\{\sum_{i_0=0}^{\infty }\frac{(i_0+ \lambda +\alpha ) (i_0+ \lambda +\beta )}{(i_0+ \lambda +2)(i_0+ \lambda +1+\gamma )}\frac{\left(\Delta_{0}^{-}\right)_{i_0} \left(\Delta_{0}^{+}\right)_{i_0}}{(1+\lambda )_{i_0}(\gamma +\lambda )_{i_0}} \right. \nonumber\\
&&\times  \sum_{i_1=i_0}^{\infty } \frac{(i_1+2+ \lambda +\alpha ) (i_1+2+ \lambda +\beta )}{(i_1+ \lambda +4)(i_1+ \lambda +3+\gamma )}  \frac{\left(\Delta_{1}^{-}\right)_{i_1}\left(\Delta_{1}^{+}\right)_{i_1}(3+\lambda )_{i_0}(2+\gamma +\lambda )_{i_0}}{\left(\Delta_{1}^{-}\right)_{i_0}\left(\Delta_{1}^{+}\right)_{i_0}(3+\lambda )_{i_1}(2+\gamma +\lambda )_{i_1}} \nonumber\\
&&\times \left. \sum_{i_2=i_1}^{\infty } \frac{\left(\Delta_{2}^{-}\right)_{i_2}\left(\Delta_{2}^{+}\right)_{i_2}(5+\lambda )_{i_1}(4+\gamma +\lambda )_{i_1}}{\left(\Delta_{2}^{-}\right)_{i_1}\left(\Delta_{2}^{+}\right)_{i_1}(5+\lambda )_{i_2}(4+\gamma +\lambda )_{i_2}} \eta ^{i_2} \right\} z^2  \label{er:10033c}
\end{eqnarray}
\begin{eqnarray}
 y_3(x) &=&  c_0 x^{\lambda } \left\{\sum_{i_0=0}^{\infty }\frac{(i_0+ \lambda +\alpha ) (i_0+ \lambda +\beta )}{(i_0+ \lambda +2)(i_0+ \lambda +1+\gamma )}\frac{\left(\Delta_{0}^{-}\right)_{i_0} \left(\Delta_{0}^{+}\right)_{i_0}}{(1+\lambda )_{i_0}(\gamma +\lambda )_{i_0}} \right. \nonumber\\
&&\times  \sum_{i_1=i_0}^{\infty } \frac{(i_1+2+ \lambda +\alpha ) (i_1+2+ \lambda +\beta )}{(i_1+ \lambda +4)(i_1+ \lambda +3+\gamma )}  \frac{\left(\Delta_{1}^{-}\right)_{i_1}\left(\Delta_{1}^{+}\right)_{i_1}(3+\lambda )_{i_0}(2+\gamma +\lambda )_{i_0}}{\left(\Delta_{1}^{-}\right)_{i_0}\left(\Delta_{1}^{+}\right)_{i_0}(3+\lambda )_{i_1}(2+\gamma +\lambda )_{i_1}} \nonumber\\
&&\times \sum_{i_2=i_1}^{\infty } \frac{(i_2+4+ \lambda +\alpha ) (i_2+4+ \lambda +\beta )}{(i_2+ \lambda +6)(i_2+ \lambda +5+\gamma )}  \frac{\left(\Delta_{2}^{-}\right)_{i_2}\left(\Delta_{2}^{+}\right)_{i_2}(5+\lambda )_{i_1}(4+\gamma +\lambda )_{i_1}}{\left(\Delta_{2}^{-}\right)_{i_1}\left(\Delta_{2}^{+}\right)_{i_1}(5+\lambda )_{i_2}(4+\gamma +\lambda )_{i_2}} \nonumber\\
&&\times \left. \sum_{i_3=i_2}^{\infty } \frac{\left(\Delta_{3}^{-}\right)_{i_3}\left(\Delta_{3}^{+}\right)_{i_3}(7+\lambda )_{i_2}(6+\gamma +\lambda )_{i_2}}{\left(\Delta_{3}^{-}\right)_{i_2}\left(\Delta_{3}^{+}\right)_{i_2}(7+\lambda )_{i_3}(6+\gamma +\lambda )_{i_3}}\eta ^{i_3} \right\} z^3   \label{er:10033d} 
\end{eqnarray}
\end{subequations}
where
\begin{equation}
\begin{cases}
\Delta_{0}^{\pm}=  \frac{ \varphi +2(1+a) \lambda \pm\sqrt{\varphi ^2-4(1+a)q}}{2(1+a)} \cr
\Delta_{1}^{\pm}=  \frac{ \varphi +2(1+a)(\lambda +2) \pm\sqrt{\varphi ^2-4(1+a)q}}{2(1+a)} \cr
\Delta_{2}^{\pm}=  \frac{ \varphi +2(1+a)(\lambda +4) \pm\sqrt{\varphi ^2-4(1+a)q}}{2(1+a)} \cr
\Delta_{3}^{\pm}=  \frac{ \varphi +2(1+a)(\lambda +6) \pm\sqrt{\varphi ^2-4(1+a)q}}{2(1+a)} 
\end{cases}\nonumber 
\end{equation}
Put $l=1$ in (\ref{er:10032}). Take the new (\ref{er:10032}) into (\ref{er:10033b}).
\begin{eqnarray}
y_1(x) &=& \int_{0}^{1} dt_1\;t_1^{1+\lambda } \int_{0}^{1} du_1\;u_1^{\gamma +\lambda } \frac{1}{2\pi i} \oint dv_1 \;\frac{1}{v_1} 
\left( \frac{v_1-1}{v_1} \right)^{ -\Delta_1^{-}} \nonumber\\
&&\times (1-\eta (1-t_1)(1-u_1)v_1)^{-\Delta_1^{+}} \overleftrightarrow {w}_{1,1}^{-(\alpha +\lambda )} \left(\overleftrightarrow {w}_{1,1} \partial_{\overleftrightarrow {w}_{1,1}} \right) \overleftrightarrow {w}_{1,1}^{\alpha -\beta } \left(\overleftrightarrow {w}_{1,1} \partial_{\overleftrightarrow {w}_{1,1}} \right)\overleftrightarrow {w}_{1,1}^{\beta +\lambda }  \nonumber\\
&&\times \left\{ c_0 x^{\lambda }\sum_{i_0=0}^{\infty } \frac{\left(\Delta_1^{-}\right)_{i_0} \left(\Delta_1^{+}\right)_{i_0}}{(1+\lambda )_{i_0}(\gamma +\lambda )_{i_0}} \overleftrightarrow {w}_{1,1} ^{i_0}\right\}z \label{er:10034}
\end{eqnarray}
where
\begin{equation}
\overleftrightarrow {w}_{1,1} = \frac{v_1}{v_1-1} \frac{\eta t_1 u_1}{1-\eta (1-t_1)(1-u_1)v_1}
\nonumber 
\end{equation}
Put $l=2$ in (\ref{er:10032}). Take the new (\ref{er:10032}) into (\ref{er:10033c}).
\begin{eqnarray}
y_2(x) &=& c_0 x^{\lambda } \int_{0}^{1} dt_2\;t_2^{3+\lambda } \int_{0}^{1} du_2\;u_2^{2+\gamma +\lambda } \frac{1}{2\pi i} \oint dv_2 \;\frac{1}{v_2} 
\left( \frac{v_2-1}{v_2} \right)^{-\Delta_2^{-}} \nonumber\\
&&\times (1-\eta (1-t_2)(1-u_2)v_2)^{-\Delta_2^{+}}
 \overleftrightarrow {w}_{2,2}^{-(2+\alpha +\lambda) } \left(\overleftrightarrow {w}_{2,2} \partial_{\overleftrightarrow {w}_{2,2}} \right) \overleftrightarrow {w}_{2,2}^{\alpha -\beta } \left(\overleftrightarrow {w}_{2,2} \partial_{\overleftrightarrow {w}_{2,2}} \right)\overleftrightarrow {w}_{2,2}^{2+\beta +\lambda } \nonumber\\
&&\times \left\{\sum_{i_0=0}^{\infty }\frac{(i_0+ \lambda +\alpha ) (i_0+ \lambda +\beta )}{(i_0+ \lambda +2)(i_0+ \lambda +1+\gamma )}\frac{\left(\Delta_0^{-}\right)_{i_0} \left(\Delta_0^{+}\right)_{i_0}}{(1+\lambda )_{i_0}(\gamma +\lambda )_{i_0}} \right. \nonumber\\
&&\times \left. \sum_{i_1=i_0}^{\infty } \frac{\left(\Delta_1^{-}\right)_{i_1}\left(\Delta_1^{+}\right)_{i_1}(3+\lambda )_{i_0}(2+\gamma +\lambda )_{i_0}}{\left(\Delta_1^{-}\right)_{i_0}\left(\Delta_1^{+}\right)_{i_0}(3+\lambda )_{i_1}(2+\gamma +\lambda )_{i_1}} \overleftrightarrow {w}_{2,2}^{i_1} \right\} z^2 \hspace{1.5cm}\label{er:10035}
\end{eqnarray}
where
\begin{equation}
\overleftrightarrow {w}_{2,2} = \frac{v_2}{v_2-1} \frac{\eta t_2 u_2}{1-\eta (1-t_2)(1-u_2)v_2}
\nonumber 
\end{equation}
Put $l=1$ and $\eta = \overleftrightarrow {w}_{2,2}$ in (\ref{er:10032}). Take the new (\ref{er:10032}) into (\ref{er:10035}).
\begin{eqnarray}
y_2(x) &=& \int_{0}^{1} dt_2\;t_2^{3+\lambda } \int_{0}^{1} du_2\;u_2^{2+\gamma +\lambda } \frac{1}{2\pi i} \oint dv_2 \;\frac{1}{v_2} 
\left( \frac{v_2-1}{v_2} \right)^{-\Delta_2^{-}} \nonumber\\
&&\times (1-\eta (1-t_2)(1-u_2)v_2)^{-\Delta_2^{+}}
 \overleftrightarrow {w}_{2,2}^{-(2+\alpha +\lambda) } \left(\overleftrightarrow {w}_{2,2} \partial_{\overleftrightarrow {w}_{2,2}} \right) \overleftrightarrow {w}_{2,2}^{\alpha -\beta } \left(\overleftrightarrow {w}_{2,2} \partial_{\overleftrightarrow {w}_{2,2}} \right)\overleftrightarrow {w}_{2,2}^{2+\beta +\lambda } \nonumber\\
&&\times \int_{0}^{1} dt_1\;t_1^{1+\lambda } \int_{0}^{1} du_1\;u_1^{\gamma +\lambda } \frac{1}{2\pi i} \oint dv_1 \;\frac{1}{v_1} 
 \left( \frac{v_1-1}{v_1} \right)^{ -\Delta_1^{-}} \nonumber\\
&&\times (1-\overleftrightarrow {w}_{2,2} (1-t_1)(1-u_1)v_1)^{ -\Delta_1^{+}} \overleftrightarrow {w}_{1,2}^{-(\alpha +\lambda)} \left(\overleftrightarrow {w}_{1,2} \partial_{\overleftrightarrow {w}_{1,2}} \right) \overleftrightarrow {w}_{1,2}^{\alpha -\beta } \left(\overleftrightarrow {w}_{1,2} \partial_{\overleftrightarrow {w}_{1,2}} \right)\overleftrightarrow {w}_{1,2}^{\beta +\lambda } \nonumber\\
&&\times \left\{ c_0 x^{\lambda }\sum_{i_0=0}^{\infty } \frac{\left( \Delta_0^{-}\right)_{i_0} \left(\Delta_0^{+}\right)_{i_0}}{(1+\lambda )_{i_0}(\gamma +\lambda )_{i_0}} \overleftrightarrow {w}_{1,2} ^{i_0}\right\} z^2 \label{er:10036}
\end{eqnarray}
where
\begin{equation}
\overleftrightarrow {w}_{1,2} = \frac{v_1}{v_1-1} \frac{\overleftrightarrow {w}_{2,2} t_1 u_1}{1- \overleftrightarrow {w}_{2,2}(1-t_1)(1-u_1)v_1}
\nonumber 
\end{equation}
By using similar process for the previous cases of integral forms of $y_1(x)$ and $y_2(x)$, the integral form of sub-power series expansion of $y_3(x)$ is
\begin{eqnarray}
y_3(x) &=& \int_{0}^{1} dt_3\;t_3^{5+\lambda } \int_{0}^{1} du_3\;u_3^{4+\gamma +\lambda } \frac{1}{2\pi i} \oint dv_3 \;\frac{1}{v_3} 
 \left( \frac{v_3-1}{v_3} \right)^{ -\Delta_3^{-}} \nonumber\\
&&\times (1-\eta (1-t_3)(1-u_3)v_3)^{ -\Delta_3^{+}} \overleftrightarrow {w}_{3,3}^{-(4+\alpha +\lambda) } \left(\overleftrightarrow {w}_{3,3} \partial_{\overleftrightarrow {w}_{3,3}} \right) \overleftrightarrow {w}_{3,3}^{\alpha -\beta } \left(\overleftrightarrow {w}_{3,3} \partial_{\overleftrightarrow {w}_{3,3}} \right)\overleftrightarrow {w}_{3,3}^{4+\beta +\lambda }  \nonumber\\
&&\times \int_{0}^{1} dt_2\;t_2^{3+\lambda } \int_{0}^{1} du_2\;u_2^{2+\gamma +\lambda } \frac{1}{2\pi i} \oint dv_2 \;\frac{1}{v_2} 
 \left(  \frac{v_2-1}{v_2} \right)^{ -\Delta_2^{-}} \nonumber\\
&&\times (1-\overleftrightarrow {w}_{3,3} (1-t_2)(1-u_2)v_2)^{ -\Delta_2^{+}} \overleftrightarrow {w}_{2,3}^{-(2+\alpha +\lambda)} \left(\overleftrightarrow {w}_{2,3} \partial_{\overleftrightarrow {w}_{2,3}} \right) \overleftrightarrow {w}_{2,3}^{\alpha -\beta } \left(\overleftrightarrow {w}_{2,3} \partial_{\overleftrightarrow {w}_{2,3}} \right)\overleftrightarrow {w}_{2,3}^{2+\beta +\lambda } \nonumber\\
&&\times \int_{0}^{1} dt_1\;t_1^{1+\lambda } \int_{0}^{1} du_1\;u_1^{\gamma +\lambda } \frac{1}{2\pi i} \oint dv_1 \;\frac{1}{v_1} 
\left(  \frac{v_1-1}{v_1}\right)^{ -\Delta_1^{-}} \nonumber\\
&&\times (1-\overleftrightarrow {w}_{2,3} (1-t_1)(1-u_1)v_1)^{ -\Delta_1^{+}} \overleftrightarrow {w}_{1,3}^{-(\alpha +\lambda)} \left(\overleftrightarrow {w}_{1,3} \partial_{\overleftrightarrow {w}_{1,3}} \right) \overleftrightarrow {w}_{1,3}^{\alpha -\beta } \left(\overleftrightarrow {w}_{1,3} \partial_{\overleftrightarrow {w}_{1,3}} \right)\overleftrightarrow {w}_{1,3}^{\beta +\lambda } \nonumber\\
&&\times \left\{ c_0 x^{\lambda }\sum_{i_0=0}^{\infty } \frac{\left( \Delta_0^{-}\right)_{i_0} \left( \Delta_0^{+}\right)_{i_0}}{(1+\lambda )_{i_0}(\gamma +\lambda )_{i_0}} \overleftrightarrow {w}_{1,3} ^{i_0}\right\} z^3 \label{er:10037}
\end{eqnarray}
where
\begin{equation}
\begin{cases} \overleftrightarrow {w}_{3,3} = \frac{v_3}{v_3-1} \frac{\eta  t_3 u_3}{1- \eta (1-t_3)(1-u_3)v_3} \cr
\overleftrightarrow {w}_{2,3} = \frac{v_2}{v_2-1} \frac{\overleftrightarrow {w}_{3,3} t_2 u_2}{1- \overleftrightarrow {w}_{3,3}(1-t_2)(1-u_2)v_2} \cr
\overleftrightarrow {w}_{1,3} = \frac{v_1}{v_1-1} \frac{\overleftrightarrow {w}_{2,3} t_1 u_1}{1- \overleftrightarrow {w}_{2,3}(1-t_1)(1-u_1)v_1}
\end{cases}
\nonumber
\end{equation}
By repeating this process for all higher terms of integral forms of sub-summation $y_m(x)$ terms where $m \geq 4$, I obtain every integral forms of $y_m(x)$ terms. 
Since we substitute (\ref{er:10033a}), (\ref{er:10034}), (\ref{er:10036}), (\ref{er:10037}) and including all integral forms of $y_m(x)$ terms where $m \geq 4$ into (\ref{eq:10020}), we obtain (\ref{eq:10049}).\footnote{Or replace the finite summation with an interval $[0, q_0]$ by infinite summation with an interval  $[0,\infty ]$ in (\ref{eq:10039}).Replace $q_0$ and $q_{n-k}$ by $\frac{-(\varphi +2(1+a)\lambda )+\sqrt{\varphi ^2-4(1+a)q}}{2(1+a)}$ and $\frac{-(\varphi +2(1+a)(\lambda +2(n-k)))+\sqrt{\varphi ^2-4(1+a)q}}{2(1+a)}$ into the new (\ref{eq:10039}). Its solution is also equivalent to (\ref{eq:10049}).}
\end{proof}
Put $c_0$= 1 as $\lambda $=0 for the first kind of independent solutions of Heun equation and $\displaystyle{ c_0= \left( \frac{1+a}{a}\right)^{1-\gamma }}$ as $\lambda = 1-\gamma $ for the second one in (\ref{eq:10049}). 
\begin{remark}
The integral representation of Heun equation of the first kind for infinite series about $x=0$ using R3TRF is
\begin{eqnarray}
 y(x)&=& HF^R \left( \varphi =\alpha +\beta -\delta +a(\delta +\gamma -1); \eta = \frac{(1+a)}{a} x ; z= -\frac{1}{a} x^2 \right) \nonumber\\
&=& \; _2F_1 \left( \Delta_0^{-}, \Delta_0^{+}; \gamma; \eta \right) \nonumber\\
&&+ \sum_{n=1}^{\infty } \Bigg\{\prod _{k=0}^{n-1} \Bigg\{ \int_{0}^{1} dt_{n-k}\;t_{n-k}^{2(n-k)-1} \int_{0}^{1} du_{n-k}\;u_{n-k}^{2(n-k-1)+\gamma}   \nonumber\\
&&\times \frac{1}{2\pi i}  \oint dv_{n-k} \frac{1}{v_{n-k}} \left( \frac{v_{n-k}-1}{v_{n-k}}\right)^{ -\Delta_{n-k}^{-}} \left( 1- \overleftrightarrow {w}_{n-k+1,n}(1-t_{n-k})(1-u_{n-k})v_{n-k}\right)^{ -\Delta_{n-k}^{+}}\nonumber\\
&&\times \overleftrightarrow {w}_{n-k,n}^{-(2(n-k-1)+\alpha )}\left(  \overleftrightarrow {w}_{n-k,n} \partial _{ \overleftrightarrow {w}_{n-k,n}}\right) \overleftrightarrow {w}_{n-k,n}^{\alpha -\beta} \left(  \overleftrightarrow {w}_{n-k,n} \partial _{ \overleftrightarrow {w}_{n-k,n}}\right) \overleftrightarrow {w}_{n-k,n}^{2(n-k-1)+\beta } \Bigg\}\nonumber\\
&&\times \; _2F_1 \left( \Delta_0^{-}, \Delta_0^{+}; \gamma; \overleftrightarrow {w}_{1,n} \right)\Bigg\} z^n  \label{eq:10050}
\end{eqnarray}
where
\begin{equation}
\begin{cases} 
\Delta_0^{\pm}= \frac{ \varphi  \pm\sqrt{\varphi ^2-4(1+a)q}}{2(1+a)} \cr
\Delta_{n-k}^{\pm}= \frac{ \varphi +4(1+a)(n-k) \pm\sqrt{\varphi ^2-4(1+a)q}}{2(1+a)} 
\end{cases}\nonumber 
\end{equation}
\end{remark}
\begin{remark}
The integral representation of Heun equation of the second kind for infinite series about $x=0$ using R3TRF is
\begin{eqnarray}
y(x)&=& HS^R \left( \varphi =\alpha +\beta -\delta +a(\delta +\gamma -1); \eta = \frac{(1+a)}{a} x ; z= -\frac{1}{a} x^2 \right) \nonumber\\
&=& \eta ^{1-\gamma }  \Bigg\{ \; _2F_1 \left( \Delta_0^{-}, \Delta_0^{+}; 2-\gamma; \eta \right) + \sum_{n=1}^{\infty } \Bigg\{\prod _{k=0}^{n-1} \Bigg\{ \int_{0}^{1} dt_{n-k}\;t_{n-k}^{2(n-k)-\gamma } \int_{0}^{1} du_{n-k}\;u_{n-k}^{2(n-k)-1}\nonumber\\
&&\times \frac{1}{2\pi i}  \oint dv_{n-k} \frac{1}{v_{n-k}} \left( \frac{v_{n-k}-1}{v_{n-k}}\right)^{ -\Delta_{n-k}^{-}} \left( 1- \overleftrightarrow {w}_{n-k+1,n}(1-t_{n-k})(1-u_{n-k})v_{n-k}\right)^{ -\Delta_{n-k}^{+}}\nonumber\\
&&\times \overleftrightarrow {w}_{n-k,n}^{-(2(n-k)-1+\alpha -\gamma )}\left(  \overleftrightarrow {w}_{n-k,n} \partial _{ \overleftrightarrow {w}_{n-k,n}}\right) \overleftrightarrow {w}_{n-k,n}^{\alpha -\beta} \left(  \overleftrightarrow {w}_{n-k,n} \partial _{ \overleftrightarrow {w}_{n-k,n}}\right) \overleftrightarrow {w}_{n-k,n}^{2(n-k)-1+\beta 
-\gamma } \Bigg\}\nonumber\\
&&\times \; _2F_1 \left( \Delta_0^{-}, \Delta_0^{+}; 2-\gamma; \overleftrightarrow {w}_{1,n} \right) \Bigg\} z^n \Bigg\} \label{eq:10051}
\end{eqnarray}
where
\begin{equation}
\begin{cases} 
\Delta_0^{\pm}=  \frac{ \varphi +2(1+a)( 1-\gamma ) \pm\sqrt{\varphi ^2-4(1+a)q}}{2(1+a)} \cr
\Delta_{n-k}^{\pm}=  \frac{ \varphi +2(1+a)( 2(n-k)+1-\gamma ) \pm\sqrt{\varphi ^2-4(1+a)q}}{2(1+a)}
\end{cases}\nonumber 
\end{equation}
\end{remark}

\section{Summary}

In mathematical definition, Heun polynomial is the polynomial of type 3: for $\alpha $ or $\beta $ = $-N-\lambda $ where $N \in \mathbb{N}_{0}$ and $\lambda = 0$ or $1-\gamma $ there will be a polynomial equation of degree $N+1$ for the determination $q$ about $x=0$. 

In Ref.\cite{Chou2012H11,Chou2012H21} I show how to obtain the power series expansion in closed forms and its integral forms of Heun functions (infinite series and polynomial of type 1) including all higher terms of $A_n$'s by applying 3TRF. This was done by letting $A_n$ in sequence $c_n$ is the leading term in the analytic function $y(x)$: the sequence $c_n$ consists of combinations $A_n$ and $B_n$. For  polynomial of type 1, I treat $\gamma $, $\delta $ and $q$ as free variables and fixed values of $\alpha $ and/or $\beta $.

In this chapter I show how to construct the Frobenius solutions in closed forms of Heun equation and its integral forms for infinite series and polynomial of type 2 including all higher terms of $B_n$'s by applying R3TRF. This is done by letting $B_n$ in sequence $c_n$ is the leading term in the analytic function $y(x)$. For polynomial of type 2, I treat $\alpha $, $\beta $, $\gamma $ and $\delta $ as free variables and a fixed value of $q$. 

The power series expansion of Heun equation and its integral form for infinite series about $x=0$ in this chapter are equivalent to infinite series of Heun equation in Ref.\cite{Chou2012H11,Chou2012H21}. In this chapter $B_n$ is the leading term in sequence $c_n$ in the analytic function $y(x)$. In Ref.\cite{Chou2012H11,Chou2012H21} $A_n$ is the leading term in sequence $c_n$ in the analytic function $y(x)$.

In Ref.\cite{Chou2012H11,Chou2012H21} as we see the power series expansions of Heun equation for all cases of infinite series and polynomial, the denominators and numerators in all $B_n$ terms of Heun equation about $x=0$ arise with Pochhammer symbol. Again in this chapter the denominators and numerators in all $A_n$ terms of Heun function about $x=0$ arise with Pochhammer symbol. Since we construct the power series expansions with Pochhammer symbols in numerators and denominators, we are able to describe integral forms of Heun equation in simple way. As we observe representations in closed form integrals of Heun equation by applying either 3TRF or R3TRF, a $_2F_1$ function recurs in each of sub-integral forms of Heun function.
We can transforom Heun function into any special functions with two recursive coefficients because of $_2F_1$ functions in integral forms of Heun function. After we replace $_2F_1$ functions in integral forms of Heun function to other special functions, we can rebuild the power series expansion of Heun function in a backward.   

\begin{appendices}
\section{Power series expansion of 192 Heun functions}\label{App:AppendixA}
In this chapter, the fundamental power series expansion in closed forms and its integral forms of Heun function about $x=0$ is constructed by using R3TRF.
A machine-generated list of 192 (isomorphic to the Coxeter group of the Coxeter diagram $D_4$) local solutions of the Heun equation was obtained by Robert S. Maier(2007) \cite{Maie20071}. In appendix of Ref.\cite{Chou2012H21}, applying 3TRF, I analyzed the power series expansions and its integral forms of Heun function (infinite series and polynomial of type 1) of nine out of the 192 local solutions of Heun function in Table 2 \cite{Maie20071}.  

In these appendices 1 and 2, replacing coefficients in the general expression of power series expansion and integral representation of Heun function for polynomial and infinite series, I derive the analytic solutions for the power series expansions and its integrals (infinite series and polynomial of type 2) of the previous nine examples of the 192 local solutions of Heun function in Ref.\cite{Chou2012H21}.\footnote{I treat $\alpha $, $\beta $, $\gamma $ and $\delta $ as free variables and a fixed value of $q$ to construct the polynomial of type 2 and its integral formalism for all nine examples of the 192 local solutions of Heun function.}
\addtocontents{toc}{\protect\setcounter{tocdepth}{1}}
\subsection{ ${\displaystyle (1-x)^{1-\delta } Hl(a, q - (\delta  - 1)\gamma a; \alpha - \delta  + 1, \beta - \delta + 1, \gamma ,2 - \delta ; x)}$ }
\subsubsection{Polynomial of type 2}
Replacing coefficients $q$, $\alpha$, $\beta$, $\delta$, $c_0$ and $\lambda $ by $q - (\delta - 1)\gamma a $, $\alpha - \delta  + 1 $, $\beta - \delta + 1$, $2 - \delta$, 1 and zero into (\ref{eq:1007}). Multiply $(1-x)^{1-\delta }$ and the new (\ref{eq:1007}) together.
 \begin{eqnarray}
& &(1-x)^{1-\delta } y(x)\nonumber\\
&=& (1-x)^{1-\delta } Hl\left(a, q - (\delta  - 1)\gamma a; \alpha - \delta + 1, \beta - \delta + 1, \gamma ,2 - \delta ; x\right)\nonumber\\
&=& (1-x)^{1-\delta } \Bigg\{\sum_{i_0=0}^{q_0} \frac{(-q_0)_{i_0} \left(q_0+ \Omega\right)_{i_0}}{(1)_{i_0}(\gamma)_{i_0}} \eta ^{i_0}\nonumber\\
&+& \left\{ \sum_{i_0=0}^{q_0}\frac{(i_0+1+\alpha -\delta) (i_0+ 1+\beta -\delta)}{(i_0+ 2)(i_0+1+\gamma )}\frac{(-q_0)_{i_0} \left(q_0+ \Omega\right)_{i_0}}{(1)_{i_0}(\gamma )_{i_0}}\right.\nonumber\\
&\times&  \left. \sum_{i_1=i_0}^{q_1} \frac{(-q_1)_{i_1}\left(q_1+4+ \Omega\right)_{i_1}(3)_{i_0}(2+\gamma )_{i_0}}{(-q_1)_{i_0}\left(q_1+4+ \Omega\right)_{i_0}(3)_{i_1}(2+\gamma)_{i_1}} \eta ^{i_1}\right\} z\nonumber\\
&+& \sum_{n=2}^{\infty } \left\{ \sum_{i_0=0}^{q_0} \frac{(i_0+1+\alpha-\delta) (i_0+1+\beta-\delta)}{(i_0+2)(i_0+1+\gamma )}\frac{(-q_0)_{i_0} \left(q_0+ \Omega\right)_{i_0}}{(1)_{i_0}(\gamma)_{i_0}}\right.\nonumber\\
&\times& \prod _{k=1}^{n-1} \Bigg\{ \sum_{i_k=i_{k-1}}^{q_k} \frac{(i_k+ 2k+1+\alpha-\delta) (i_k+ 2k+1+\beta -\delta )}{(i_k+ 2(k+1))(i_k+ 2k+1+\gamma)} \nonumber\\
&\times& \frac{(-q_k)_{i_k}\left(q_k+4k+ \Omega\right)_{i_k}(2k+1)_{i_{k-1}}(2k+\gamma)_{i_{k-1}}}{(-q_k)_{i_{k-1}}\left(q_k+4k+\Omega \right)_{i_{k-1}}(2k+1)_{i_k}(2k+\gamma)_{i_k}}\Bigg\} \nonumber\\
&\times& \left.\left.\sum_{i_n= i_{n-1}}^{q_n} \frac{(-q_n)_{i_n}\left(q_n+4n+ \Omega\right)_{i_n}(2n+1)_{i_{n-1}}(2n+\gamma)_{i_{n-1}}}{(-q_n)_{i_{n-1}}\left(q_n+4n+ \Omega\right)_{i_{n-1}}(2n+1)_{i_n}(2n+\gamma)_{i_n}} \eta ^{i_n} \right\} z^n \right\} \label{eq:10052}
\end{eqnarray}
  where
 \begin{equation}
\begin{cases} z = -\frac{1}{a}x^2 \cr
\eta = \frac{(1+a)}{a} x \cr
\varphi = \alpha +\beta -\delta +a(\gamma -\delta +1) \cr
\Omega = \frac{\varphi }{(1+a)}\cr
q = (\delta - 1)\gamma a-(q_j+2j)\{\varphi +(1+a)(q_j+2j) \} \;\;\mbox{as}\;j,q_j\in \mathbb{N}_{0} \cr
q_i\leq q_j \;\;\mbox{only}\;\mbox{if}\;i\leq j\;\;\mbox{where}\;i,j\in \mathbb{N}_{0} 
\end{cases}\nonumber 
\end{equation}
 \subsubsection{Infinite series}
Replacing coefficients $q$, $\alpha$, $\beta$, $\delta$, $c_0$ and $\lambda $ by $q - (\delta - 1)\gamma a $, $\alpha - \delta  + 1 $, $\beta - \delta + 1$, $2-\delta$, 1 and zero into (\ref{eq:10020}). Multiply $(1-x)^{1-\delta }$ and the new (\ref{eq:10020}) together.
 \begin{eqnarray}
 & &(1-x)^{1-\delta } y(x)\nonumber\\
&=& (1-x)^{1-\delta } Hl\left(a, q - (\delta  - 1)\gamma a; \alpha - \delta + 1, \beta - \delta + 1, \gamma ,2 - \delta ; x\right)\nonumber\\
&=& (1-x)^{1-\delta } \left\{\sum_{i_0=0}^{\infty } \frac{\left(\Delta_0^{-}\right)_{i_0} \left(\Delta_0^{+}\right)_{i_0}}{(1)_{i_0}(\gamma)_{i_0}} \eta ^{i_0}\right.\nonumber\\
&+& \left\{ \sum_{i_0=0}^{\infty }\frac{(i_0 +1+\alpha - \delta) (i_0 +1+\beta - \delta)}{(i_0 +2)(i_0 +1+\gamma )}\frac{\left(\Delta_0^{-}\right)_{i_0} \left(\Delta_0^{+}\right)_{i_0}}{(1)_{i_0}(\gamma)_{i_0}}\sum_{i_1=i_0}^{\infty } \frac{\left(\Delta_1^{-}\right)_{i_1} \left(\Delta_1^{+}\right)_{i_1}(3)_{i_0}(2+\gamma)_{i_0}}{\left(\Delta_1^{-}\right)_{i_0} \left(\Delta_1^{+}\right)_{i_0}(3)_{i_1}(2+\gamma)_{i_1}} \eta ^{i_1}\right\} z\nonumber\\
&+& \sum_{n=2}^{\infty } \left\{ \sum_{i_0=0}^{\infty } \frac{(i_0 +1+\alpha - \delta) (i_0 +1+\beta - \delta)}{(i_0+2)(i_0 +1+\gamma )}\frac{\left(\Delta_0^{-}\right)_{i_0} \left(\Delta_0^{+}\right)_{i_0}}{(1)_{i_0}(\gamma)_{i_0}}\right.\nonumber\\
&\times& \prod _{k=1}^{n-1} \left\{ \sum_{i_k=i_{k-1}}^{\infty } \frac{(i_k+ 2k+1+\alpha -\delta) (i_k+ 2k+1+\beta - \delta)}{(i_k+ 2(k+1))(i_k+ 2k+1+\gamma)}\frac{\left(\Delta_k^{-}\right)_{i_k} \left(\Delta_k^{+}\right)_{i_k}(2k+1)_{i_{k-1}}(2k+\gamma)_{i_{k-1}}}{\left(\Delta_k^{-}\right)_{i_{k-1}} \left(\Delta_k^{+}\right)_{i_{k-1}}(2k+1)_{i_k}(2k+\gamma)_{i_k}}\right\} \nonumber\\
&\times& \left.\left.\sum_{i_n= i_{n-1}}^{\infty } \frac{\left(\Delta_n^{-}\right)_{i_n} \left(\Delta_n^{+}\right)_{i_n}(2n+1)_{i_{n-1}}(2n+\gamma)_{i_{n-1}}}{\left(\Delta_n^{-}\right)_{i_{n-1}} \left( \Delta_n^{+}\right)_{i_{n-1}}(2n+1)_{i_n}(2n+\gamma)_{i_n}} \eta ^{i_n} \right\} z^n \right\} \label{eq:10053}
\end{eqnarray}
where
 \begin{equation}
\begin{cases} 
\Delta_0^{\pm}= \frac{\varphi \pm\sqrt{\varphi ^2-4(1+a)q}}{2(1+a)} \cr
\Delta_1^{\pm}=  \frac{\{\varphi +4(1+a)\}\pm\sqrt{\varphi ^2-4(1+a)q}}{2(1+a)} \cr
\Delta_k^{\pm}=  \frac{\{\varphi +4(1+a)k\}\pm\sqrt{\varphi ^2-4(1+a)q}}{2(1+a)} \cr
\Delta_n^{\pm}=  \frac{\{\varphi +4(1+a)n \}\pm\sqrt{\varphi ^2-4(1+a)q}}{2(1+a)}
\end{cases}\nonumber 
\end{equation} 
\subsection{ \footnotesize ${\displaystyle x^{1-\gamma } (1-x)^{1-\delta }Hl(a, q-(\gamma +\delta -2)a -(\gamma -1)(\alpha +\beta -\gamma -\delta +1), \alpha - \gamma -\delta +2, \beta - \gamma -\delta +2, 2-\gamma, 2 - \delta ; x)}$ \normalsize}
\subsubsection{Polynomial of type 2}
Replacing coefficients $q$, $\alpha$, $\beta$, $\gamma $, $\delta$, $c_0$ and $\lambda $ by $q-(\gamma +\delta -2)a-(\gamma -1)(\alpha +\beta -\gamma -\delta +1)$, $\alpha - \gamma -\delta +2$, $\beta - \gamma -\delta +2, 2-\gamma$, $2 - \delta$,1 and zero into (\ref{eq:1007}). Multiply $x^{1-\gamma } (1-x)^{1-\delta }$ and the new (\ref{eq:1007}) together.
\begin{eqnarray}
 & &x^{1-\gamma } (1-x)^{1-\delta } y(x)\nonumber\\
&=& x^{1-\gamma } (1-x)^{1-\delta } Hl(a, q-(\gamma +\delta -2)a-(\gamma -1)(\alpha +\beta -\gamma -\delta +1); \alpha - \gamma -\delta +2 \nonumber\\
&&, \beta - \gamma -\delta +2, 2-\gamma, 2 - \delta ; x)\nonumber\\
&=& x^{1-\gamma } (1-x)^{1-\delta } \left\{\sum_{i_0=0}^{q_0} \frac{(-q_0)_{i_0} \left(q_0+ \Omega\right)_{i_0}}{(1)_{i_0}(2-\gamma)_{i_0}} \eta ^{i_0}\right.\nonumber\\
&&+ \left\{ \sum_{i_0=0}^{q_0}\frac{(i_0+2+\alpha-\gamma -\delta) (i_0+2+\beta -\gamma -\delta)}{(i_0+2)(i_0+3-\gamma )}\frac{(-q_0)_{i_0} \left(q_0+ \Omega\right)_{i_0}}{(1)_{i_0}(2-\gamma)_{i_0}} \right.\nonumber\\
&&\times \left. \sum_{i_1=i_0}^{q_1} \frac{(-q_1)_{i_1}\left(q_1+4+ \Omega\right)_{i_1}(3)_{i_0}(4-\gamma)_{i_0}}{(-q_1)_{i_0}\left(q_1+4+ \Omega\right)_{i_0}(3)_{i_1}(4-\gamma)_{i_1}} \eta ^{i_1}\right\} z\nonumber\\
&&+ \sum_{n=2}^{\infty } \left\{ \sum_{i_0=0}^{q_0} \frac{(i_0+2+\alpha-\gamma -\delta) (i_0+2+\beta -\gamma -\delta)}{(i_0+2)(i_0+3-\gamma )}\frac{(-q_0)_{i_0} \left(q_0+ \Omega\right)_{i_0}}{(1)_{i_0}(2-\gamma)_{i_0}}\right.\nonumber\\
&&\times \prod _{k=1}^{n-1} \left\{ \sum_{i_k=i_{k-1}}^{q_k} \frac{(i_k+ 2k+2+\alpha-\gamma -\delta) (i_k+ 2k+2+\beta -\gamma -\delta)}{(i_k+ 2(k+1))(i_k+ 2k+3-\gamma)}\right. \nonumber\\
&&\times \left.\frac{(-q_k)_{i_k}\left(q_k+4k+ \Omega\right)_{i_k}(2k+1)_{i_{k-1}}(2k+2-\gamma)_{i_{k-1}}}{(-q_k)_{i_{k-1}}\left(q_k+4k+ \Omega\right)_{i_{k-1}}(2k+1)_{i_k}(2k+2-\gamma)_{i_k}}\right\} \nonumber\\
&&\times \left.\left.\sum_{i_n= i_{n-1}}^{q_n} \frac{(-q_n)_{i_n}\left(q_n+4n+ \Omega\right)_{i_n}(2n+1)_{i_{n-1}}(2n+2-\gamma)_{i_{n-1}}}{(-q_n)_{i_{n-1}}\left(q_n+4n+ \Omega\right)_{i_{n-1}}(2n+1)_{i_n}(2n+2-\gamma)_{i_n}} \eta ^{i_n} \right\} z^n \right\}\label{eq:10054}
\end{eqnarray}
where
  \begin{equation}
\begin{cases} z = -\frac{1}{a}x^2 \cr
\eta = \frac{(1+a)}{a} x \cr
\varphi = \alpha +\beta -2\gamma -\delta +2 +a(3-\gamma -\delta ) \cr
\Omega = \frac{\varphi }{(1+a)}\cr
q = (\gamma +\delta -2)a+(\gamma -1)(\alpha +\beta -\gamma -\delta +1)\cr 
\hspace{0.6cm}-(q_j+2j)\{\varphi +(1+a)(q_j+2j) \} \;\;\mbox{as}\;j,q_j\in \mathbb{N}_{0} \cr
q_i\leq q_j \;\;\mbox{only}\;\mbox{if}\;i\leq j\;\;\mbox{where}\;i,j\in \mathbb{N}_{0} 
\end{cases}\nonumber 
\end{equation}
\subsubsection{Infinite series}
Replacing coefficients $q$, $\alpha$, $\beta$, $\gamma $, $\delta$, $c_0$ and $\lambda $ by $q-(\gamma +\delta -2)a-(\gamma -1)(\alpha +\beta -\gamma -\delta +1)$, $\alpha - \gamma -\delta +2$, $\beta - \gamma -\delta +2, 2-\gamma$, $2 - \delta$,1 and zero into (\ref{eq:10020}). Multiply $x^{1-\gamma } (1-x)^{1-\delta }$ and the new (\ref{eq:10020}) together.
 \begin{eqnarray}
 & &x^{1-\gamma } (1-x)^{1-\delta } y(x)\nonumber\\
&=& x^{1-\gamma } (1-x)^{1-\delta } Hl(a, q-(\gamma +\delta -2)a-(\gamma -1)(\alpha +\beta -\gamma -\delta +1); \alpha - \gamma -\delta +2\nonumber\\
&&, \beta - \gamma -\delta +2, 2-\gamma, 2 - \delta ; x)\nonumber\\
&=& x^{1-\gamma } (1-x)^{1-\delta } \left\{\sum_{i_0=0}^{\infty } \frac{\left(\Delta_0^{-}\right)_{i_0} \left(\Delta_0^{+}\right)_{i_0}}{(1)_{i_0}(2-\gamma)_{i_0}} \eta ^{i_0}\right.\nonumber\\
&+& \left\{ \sum_{i_0=0}^{\infty }\frac{(i_0+2+\alpha -\gamma -\delta ) (i_0+2+\beta -\gamma -\delta)}{(i_0+2)(i_0+3-\gamma )}\frac{\left(\Delta_0^{-}\right)_{i_0} \left(\Delta_0^{+}\right)_{i_0}}{(1)_{i_0}(2-\gamma )_{i_0}}  \sum_{i_1=i_0}^{\infty }\frac{\left(\Delta_1^{-}\right)_{i_1} \left(\Delta_1^{+}\right)_{i_1}(3)_{i_0}(4-\gamma)_{i_0}}{\left(\Delta_1^{-}\right)_{i_0} \left(\Delta_1^{+}\right)_{i_0}(3)_{i_1}(4-\gamma)_{i_1}} \eta ^{i_1}\right\} z\nonumber\\
&+& \sum_{n=2}^{\infty } \left\{ \sum_{i_0=0}^{\infty } \frac{(i_0+2+\alpha -\gamma -\delta ) (i_0+2+\beta -\gamma -\delta)}{(i_0+2)(i_0+3-\gamma )}\frac{\left(\Delta_0^{-}\right)_{i_0} \left(\Delta_0^{+}\right)_{i_0}}{(1)_{i_0}(2-\gamma )_{i_0}}\right.\nonumber\\
&\times& \prod _{k=1}^{n-1} \left\{ \sum_{i_k=i_{k-1}}^{\infty } \frac{(i_k+ 2k+2+\alpha -\gamma -\delta ) (i_k+ 2k+2+\beta -\gamma -\delta)}{(i_k+ 2(k+1) )(i_k+ 2k+3-\gamma )}\right. \nonumber\\
&\times& \left.\frac{\left(\Delta_k^{-}\right)_{i_k}\left(\Delta_k^{+}\right)_{i_k}(2k+1)_{i_{k-1}}(2k+2-\gamma)_{i_{k-1}}}{\left(\Delta_k^{-}\right)_{i_{k-1}} \left(\Delta_k^{+}\right)_{i_{k-1}}(2k+1)_{i_k}(2k+2-\gamma)_{i_k}}\right\}\nonumber\\
&\times& \left.\left. \sum_{i_n= i_{n-1}}^{\infty } \frac{\left(\Delta_n^{-}\right)_{i_n} \left(\Delta_n^{+}\right)_{i_n}(2n+1)_{i_{n-1}}(2n+2-\gamma)_{i_{n-1}}}{\left(\Delta_n^{-}\right)_{i_{n-1}} \left(\Delta_n^{+}\right)_{i_{n-1}}(2n+1)_{i_n}(2n+2-\gamma )_{i_n}} \eta ^{i_n} \right\} z^n \right\}\label{eq:10055}
\end{eqnarray}
  where
  \begin{equation}
\begin{cases} 
\Delta_0^{\pm}= \frac{\varphi \pm\sqrt{\varphi ^2-4(1+a)q}}{2(1+a)} \cr
\Delta_1^{\pm}= \frac{\{\varphi +4(1+a)\}\pm\sqrt{\varphi ^2-4(1+a)q}}{2(1+a)} \cr
\Delta_k^{\pm}=  \frac{\{\varphi +4(1+a)k\}\pm\sqrt{\varphi ^2-4(1+a)q}}{2(1+a)} \cr
\Delta_n^{\pm}=  \frac{\{\varphi +4(1+a)n \}\pm\sqrt{\varphi ^2-4(1+a)q}}{2(1+a)}
\end{cases}\nonumber 
\end{equation} 
\subsection{ ${\displaystyle  Hl(1-a,-q+\alpha \beta; \alpha,\beta, \delta, \gamma; 1-x)}$} 
\subsubsection{Polynomial of type 2}
Replacing coefficients $a$, $q$, $\gamma $, $\delta$, $x$, $c_0$ and $\lambda $ by $1-a$, $-q +\alpha \beta $, $\delta $, $\gamma $, $1-x$, 1 and zero into (\ref{eq:1007}).
 \begin{eqnarray}
y(\xi ) &=& Hl(1-a,-q+\alpha \beta; \alpha,\beta, \delta, \gamma; 1-x)\nonumber\\
&=& \sum_{i_0=0}^{q_0} \frac{(-q_0)_{i_0} \left(q_0+ \Omega\right)_{i_0}}{(1)_{i_0}(\delta )_{i_0}} \eta ^{i_0}\nonumber\\
&+& \left\{ \sum_{i_0=0}^{q_0}\frac{(i_0+\alpha) (i_0+\beta)}{(i_0+2)(i_0+1+\delta )}\frac{(-q_0)_{i_0} \left(q_0+ \Omega\right)_{i_0}}{(1)_{i_0}(\delta )_{i_0}}   \sum_{i_1=i_0}^{q_1} \frac{(-q_1)_{i_1}\left(q_1+4+\Omega\right)_{i_1}(3)_{i_0}(2+\delta )_{i_0}}{(-q_1)_{i_0}\left(q_1+4+ \Omega\right)_{i_0}(3)_{i_1}(2+\delta )_{i_1}} \eta ^{i_1}\right\} z\nonumber\\
&+& \sum_{n=2}^{\infty } \left\{ \sum_{i_0=0}^{q_0} \frac{(i_0+\alpha ) (i_0+\beta )}{(i_0+2)(i_0+1+\delta )}\frac{(-q_0)_{i_0} \left(q_0+\Omega\right)_{i_0}}{(1)_{i_0}(\delta )_{i_0}}\right.\nonumber\\
&\times& \prod _{k=1}^{n-1} \left\{ \sum_{i_k=i_{k-1}}^{q_k} \frac{(i_k+ 2k+\alpha ) (i_k+ 2k+\beta )}{(i_k+ 2(k+1))(i_k+ 2k+1+\delta )}  \frac{(-q_k)_{i_k}\left(q_k+4k+ \Omega\right)_{i_k}(2k+1)_{i_{k-1}}(2k+\delta )_{i_{k-1}}}{(-q_k)_{i_{k-1}}\left(q_k+4k+ \Omega\right)_{i_{k-1}}(2k+1)_{i_k}(2k+\delta )_{i_k}}\right\} \nonumber\\
&\times& \left.\sum_{i_n= i_{n-1}}^{q_n} \frac{(-q_n)_{i_n}\left(q_n+4n+\Omega\right)_{i_n}(2n+1)_{i_{n-1}}(2n+\delta )_{i_{n-1}}}{(-q_n)_{i_{n-1}}\left(q_n+4n+ \Omega\right)_{i_{n-1}}(2n+1)_{i_n}(2n+\delta )_{i_n}} \eta ^{i_n} \right\} z^n   \label{eq:10056}
\end{eqnarray}
 where
 \begin{equation}
\begin{cases} \xi =1-x \cr
z = \frac{-1}{1-a}\xi^2 \cr
\eta = \frac{2-a}{1-a}\xi \cr
\varphi = \alpha +\beta -\delta +(1-a)(\delta +\gamma -1) \cr
\Omega = \frac{\varphi }{(2-a)}\cr
q = \alpha \beta +(q_j+2j)\{\varphi +(2-a)(q_j+2j) \} \;\;\mbox{as}\;j,q_j\in \mathbb{N}_{0} \cr
q_i\leq q_j \;\;\mbox{only}\;\mbox{if}\;i\leq j\;\;\mbox{where}\;i,j\in \mathbb{N}_{0}
\end{cases}\nonumber 
\end{equation}
\subsubsection{Infinite series}
Replacing coefficients $a$, $q$, $\gamma $, $\delta$, $x$, $c_0$ and $\lambda $ by $1-a$, $-q +\alpha \beta $, $\delta $, $\gamma $, $1-x$,1 and zero into (\ref{eq:10020}).
 \begin{eqnarray}
y(\xi ) &=& Hl(1-a,-q+\alpha \beta; \alpha, \beta, \delta, \gamma; 1-x)\nonumber\\
&=& \sum_{i_0=0}^{\infty } \frac{\left(\Delta_0^{-}\right)_{i_0} \left(\Delta_0^{+}\right)_{i_0}}{(1)_{i_0}(\delta )_{i_0}} \eta ^{i_0}\nonumber\\
&+& \left\{ \sum_{i_0=0}^{\infty }\frac{(i_0 +\alpha ) (i_0+\beta )}{(i_0+2)(i_0+1+\delta)}\frac{\left(\Delta_0^{-}\right)_{i_0} \left(\Delta_0^{+}\right)_{i_0}}{(1)_{i_0}(\delta)_{i_0}} \sum_{i_1=i_0}^{\infty } \frac{\left(\Delta_1^{-}\right)_{i_1} \left(\Delta_1^{+} \right)_{i_1}(3)_{i_0}(2+\delta)_{i_0}}{\left(\Delta_1^{-}\right)_{i_0} \left(\Delta_1^{+}\right)_{i_0}(3)_{i_1}(2+\delta)_{i_1}} \eta ^{i_1} \right\} z \nonumber\\
&+& \sum_{n=2}^{\infty } \left\{ \sum_{i_0=0}^{\infty } \frac{(i_0+\alpha ) (i_0+\beta )}{(i_0+2)(i_0+1+\delta)}\frac{\left(\Delta_0^{-}\right)_{i_0} \left(\Delta_0^{+}\right)_{i_0}}{(1)_{i_0}(\delta)_{i_0}}\right.\nonumber\\
&\times& \prod _{k=1}^{n-1} \left\{ \sum_{i_k=i_{k-1}}^{\infty } \frac{(i_k+ 2k+\alpha ) (i_k+ 2k+\beta )}{(i_k+ 2(k+1))(i_k+ 2k+1+\delta)}\frac{\left(\Delta_k^{-}\right)_{i_k}\left(\Delta_k^{+}\right)_{i_k}(2k+1)_{i_{k-1}}(2k+\delta)_{i_{k-1}}}{ \left(\Delta_k^{-}\right)_{i_{k-1}}\left(\Delta_k^{+}\right)_{i_{k-1}}(2k+1)_{i_k}(2k+\delta)_{i_k}}\right\} \nonumber\\
&\times& \left. \sum_{i_n= i_{n-1}}^{\infty } \frac{ \left(\Delta_n^{-}\right)_{i_n} \left(\Delta_n^{+}\right)_{i_n}(2n+1)_{i_{n-1}}(2n+\delta)_{i_{n-1}}}{\left(\Delta_n^{-} \right)_{i_{n-1}} \left(\Delta_n^{+}\right)_{i_{n-1}}(2n+1)_{i_n}(2n+\delta)_{i_n}} \eta ^{i_n} \right\} z^n  \label{eq:10057}
\end{eqnarray}
 where
\begin{equation}
\begin{cases} 
\Delta_0^{\pm}= \frac{\varphi \pm\sqrt{\varphi ^2-4(2-a)q}}{2(2-a)}\cr
\Delta_1^{\pm}=  \frac{\{\varphi +4(2-a)\}\pm\sqrt{\varphi ^2-4(2-a)q}}{2(2-a)}\cr
\Delta_k^{\pm}=  \frac{\{\varphi +4(2-a)k\}\pm\sqrt{\varphi ^2-4(2-a)q}}{2(2-a)}\cr
\Delta_n^{\pm}=   \frac{\{\varphi +4(2-a)n\}\pm\sqrt{\varphi ^2-4(2-a)q}}{2(2-a)}
\end{cases}\nonumber 
\end{equation}
  \subsection{ \footnotesize ${\displaystyle (1-x)^{1-\delta } Hl(1-a,-q+(\delta -1)\gamma a+(\alpha -\delta +1)(\beta -\delta +1); \alpha-\delta +1,\beta-\delta +1, 2-\delta, \gamma; 1-x)}$ \normalsize}
\subsubsection{Polynomial of type 2}
Replacing coefficients $a$, $q$, $\alpha $, $\beta $, $\gamma $, $\delta$, $x$, $c_0$ and $\lambda $ by $1-a$, $-q+(\delta -1)\gamma a+(\alpha -\delta +1)(\beta -\delta +1)$, $\alpha-\delta +1 $, $\beta-\delta +1 $, $2-\delta$, $\gamma $, $1-x$, 1 and zero into (\ref{eq:1007}). Multiply $(1-x)^{1-\delta }$ and the new (\ref{eq:1007}) together.
\begin{eqnarray}
& &(1-x)^{1-\delta } y(\xi)\nonumber\\
&=& (1-x)^{1-\delta } Hl(1-a,-q+(\delta -1)\gamma a+(\alpha -\delta +1)(\beta -\delta +1); \alpha-\delta +1\nonumber\\
&&,\beta-\delta +1, 2-\delta, \gamma; 1-x) \nonumber\\
&=& (1-x)^{1-\delta } \left\{\sum_{i_0=0}^{q_0} \frac{(-q_0)_{i_0} \left(q_0+ \Omega\right)_{i_0}}{(1)_{i_0}(2-\delta )_{i_0}} \eta ^{i_0}\right.\nonumber\\
&&+ \left\{ \sum_{i_0=0}^{q_0}\frac{(i_0 +1+\alpha-\delta) (i_0 +1+\beta -\delta)}{(i_0+2)(i_0+3-\delta  )}\frac{(-q_0)_{i_0} \left(q_0+ \Omega\right)_{i_0}}{(1)_{i_0}(2-\delta )_{i_0}} \right.\nonumber\\
&&\times \left. \sum_{i_1=i_0}^{q_1} \frac{(-q_1)_{i_1}\left(q_1+4+ \Omega\right)_{i_1}(3)_{i_0}(4-\delta  )_{i_0}}{(-q_1)_{i_0}\left(q_1+4+ \Omega\right)_{i_0}(3)_{i_1}(4-\delta )_{i_1}} \eta ^{i_1}\right\} z\nonumber\\
&&+ \sum_{n=2}^{\infty } \left\{ \sum_{i_0=0}^{q_0} \frac{(i_0+1+\alpha-\delta) (i_0+1+\beta -\delta)}{(i_0+2)(i_0+3-\delta )}\frac{(-q_0)_{i_0} \left(q_0+ \Omega\right)_{i_0}}{(1)_{i_0}(2-\delta )_{i_0}}\right.\nonumber\\
&&\times \prod _{k=1}^{n-1} \left\{ \sum_{i_k=i_{k-1}}^{q_k} \frac{(i_k+ 2k+1+\alpha-\delta) (i_k+ 2k+1+\beta -\delta)}{(i_k+ 2(k+1))(i_k+ 2k+3-\delta )}\right. \nonumber\\
&&\times \left.\frac{(-q_k)_{i_k}\left( q_k+4k+ \Omega\right)_{i_k}(2k+1)_{i_{k-1}}(2k+2-\delta )_{i_{k-1}}}{(-q_k)_{i_{k-1}}\left(q_k+4k+ \Omega\right)_{i_{k-1}}(2k+1)_{i_k}(2k+2-\delta )_{i_k}}\right\} \nonumber\\
&&\times \left. \left.\sum_{i_n= i_{n-1}}^{q_n} \frac{(-q_n)_{i_n}\left(q_n+4n+ \Omega\right)_{i_n}(2n+1)_{i_{n-1}}(2n+2-\delta )_{i_{n-1}}}{(-q_n)_{i_{n-1}}\left(q_n+4n+ \Omega\right)_{i_{n-1}}(2n+1)_{i_n}(2n+2-\delta )_{i_n}} \eta ^{i_n} \right\} z^n \right\}\label{eq:10058}
\end{eqnarray}
where
 \begin{equation}
\begin{cases} \xi =1-x \cr
z = \frac{-1}{1-a}\xi^2 \cr
\eta = \frac{2-a}{1-a}\xi \cr
\varphi = \alpha +\beta -\gamma -2\delta +2+(1-a)(\gamma -\delta +1) \cr
\Omega = \frac{\varphi }{(2-a)}\cr
q =(\delta -1)\gamma a+(\alpha -\delta +1)(\beta -\delta +1)\cr
\hspace{0.6cm}+(q_j+2j)\{\varphi +(2-a)(q_j+2j ) \} \;\;\mbox{as}\;j,q_j\in \mathbb{N}_{0} \cr
q_i\leq q_j \;\;\mbox{only}\;\mbox{if}\;i\leq j\;\;\mbox{where}\;i,j\in \mathbb{N}_{0}
\end{cases}\nonumber 
\end{equation}  
\subsubsection{Infinite series}
Replacing coefficients $a$, $q$, $\alpha $, $\beta $, $\gamma $, $\delta$, $x$, $c_0$ and $\lambda $ by $1-a$, $-q+(\delta -1)\gamma a+(\alpha -\delta +1)(\beta -\delta +1)$, $\alpha-\delta +1 $, $\beta-\delta +1 $, $2-\delta$, $\gamma $, $1-x$, 1 and zero into (\ref{eq:10020}).  Multiply $(1-x)^{1-\delta }$ and the new (\ref{eq:10020}) together.
 \begin{eqnarray}
& &(1-x)^{1-\delta } y(\xi)\nonumber\\
&=& (1-x)^{1-\delta } Hl(1-a,-q+(\delta -1)\gamma a+(\alpha -\delta +1)(\beta -\delta +1); \alpha -\delta +1,\beta-\delta +1, 2-\delta, \gamma; 1-x) \nonumber\\
&=& (1-x)^{1-\delta } \left\{\sum_{i_0=0}^{\infty } \frac{\left(\Delta_0^{-}\right)_{i_0} \left(\Delta_0^{+}\right)_{i_0}}{(1)_{i_0}(2-\delta )_{i_0}} \eta ^{i_0}\right.\nonumber\\
&+& \left\{ \sum_{i_0=0}^{\infty }\frac{(i_0+1+\alpha-\delta ) (i_0+1+\beta -\delta )}{(i_0+2)(i_0+3-\delta)}\frac{\left(\Delta_0^{-}\right)_{i_0} \left(\Delta_0^{+}\right)_{i_0}}{(1)_{i_0}(2-\delta)_{i_0}}\sum_{i_1=i_0}^{\infty } \frac{\left(\Delta_1^{-}\right)_{i_1} \left(\Delta_1^{+}\right)_{i_1}(3)_{i_0}(4-\delta)_{i_0}}{\left(\Delta_1^{-}\right)_{i_0} \left(\Delta_1^{+}\right)_{i_0}(3)_{i_1}(4-\delta)_{i_1}} \eta ^{i_1}\right\} z\nonumber\\
&+& \sum_{n=2}^{\infty } \left\{ \sum_{i_0=0}^{\infty } \frac{(i_0+1+\alpha-\delta ) (i_0+1+\beta -\delta  )}{(i_0+2)(i_0+3-\delta )}\frac{\left(\Delta_0^{-}\right)_{i_0} \left(\Delta_0^{+}\right)_{i_0}}{(1)_{i_0}(2-\delta)_{i_0}}\right.\nonumber\\
&\times& \prod _{k=1}^{n-1} \left\{ \sum_{i_k=i_{k-1}}^{\infty } \frac{(i_k+ 2k+1+\alpha-\delta ) (i_k+ 2k+1+\beta -\delta )}{(i_k+ 2(k+1))(i_k+ 2k+3-\delta)}\frac{\left(\Delta_k^{-}\right)_{i_k} \left(\Delta_k^{+} \right)_{i_k}(2k+1)_{i_{k-1}}(2k+2-\delta)_{i_{k-1}}}{\left(\Delta_k^{-}\right)_{i_{k-1}} \left(\Delta_k^{+} \right)_{i_{k-1}}(2k+1)_{i_k}(2k+2-\delta)_{i_k}}\right\}\nonumber \\
&\times& \sum_{i_n= i_{n-1}}^{\infty } \left.\left.\frac{\left(\Delta_n^{-}\right)_{i_n} \left(\Delta_n^{+}\right)_{i_n}(2n+1)_{i_{n-1}}(2n+2-\delta)_{i_{n-1}}}{\left(\Delta_n^{-}\right)_{i_{n-1}} \left(\Delta_n^{+}\right)_{i_{n-1}}(2n+1)_{i_n}(2n+2-\delta)_{i_n}} \eta ^{i_n} \right\} z^n \right\} \label{eq:10059}
\end{eqnarray}
 where
\begin{equation}
\begin{cases} 
\Delta_0^{\pm}= \frac{\varphi \pm\sqrt{\varphi ^2-4(2-a)q}}{2(2-a)}\cr
\Delta_1^{\pm}=  \frac{\{\varphi +4(2-a) \}\pm\sqrt{\varphi ^2-4(2-a)q}}{2(2-a)}\cr
\Delta_k^{\pm}=   \frac{\{\varphi +4(2-a)k\}\pm\sqrt{\varphi ^2-4(2-a)q}}{2(2-a)}\cr
\Delta_n^{\pm}= \frac{\{\varphi +4(2-a)n \}\pm\sqrt{\varphi ^2-4(1+a)q}}{2(1+a)}
\end{cases}\nonumber 
\end{equation} 
\subsection{ \footnotesize ${\displaystyle x^{-\alpha } Hl\left(\frac{1}{a},\frac{q+\alpha [(\alpha -\gamma -\delta +1)a-\beta +\delta ]}{a}; \alpha , \alpha -\gamma +1, \alpha -\beta +1,\delta ;\frac{1}{x}\right)}$\normalsize}
\subsubsection{Polynomial of type 2}
Replacing coefficients $a$, $q$, $\beta $, $\gamma $, $x$, $c_0$ and $\lambda $ by $\frac{1}{a}$, $\frac{q+\alpha [(\alpha -\gamma -\delta +1)a-\beta +\delta ]}{a}$, $\alpha-\gamma +1 $, $\alpha -\beta +1 $, $\frac{1}{x}$, 1 and zero into (\ref{eq:1007}). Multiply $x^{-\alpha }$ and the new (\ref{eq:1007}) together.
\begin{eqnarray}
& &x^{-\alpha } y(\xi)\nonumber\\
&=& x^{-\alpha }  Hl\left(\frac{1}{a},\frac{q+\alpha [(\alpha -\gamma -\delta +1)a-\beta +\delta ]}{a}; \alpha , \alpha -\gamma +1, \alpha -\beta +1,\delta ;\frac{1}{x}\right) \nonumber\\
&=& x^{-\alpha } \left\{\sum_{i_0=0}^{q_0} \frac{(-q_0)_{i_0} \left(q_0+ \Omega\right)_{i_0}}{(1)_{i_0}(1+\alpha -\beta)_{i_0}} \eta ^{i_0}\right.\nonumber\\
&+& \left\{ \sum_{i_0=0}^{q_0}\frac{(i_0+\alpha ) (i_0+1+\alpha -\gamma )}{(i_0+2)(i_0+2+\alpha -\beta)}\frac{(-q_0)_{i_0} \left(q_0+ \Omega\right)_{i_0}}{(1)_{i_0}(1+\alpha -\beta)_{i_0}} \right. \nonumber\\
&\times& \left. \sum_{i_1=i_0}^{q_1} \frac{(-q_1)_{i_1}\left(q_1+4+ \Omega\right)_{i_1}(3)_{i_0}(3+\alpha -\beta )_{i_0}}{(-q_1)_{i_0}\left(q_1+4+ \Omega\right)_{i_0}(3)_{i_1}(3+\alpha -\beta )_{i_1}} \eta ^{i_1}\right\} z\nonumber\\
&+& \sum_{n=2}^{\infty } \left\{ \sum_{i_0=0}^{q_0} \frac{(i_0+\alpha ) (i_0+1+\alpha -\gamma )}{(i_0+2)(i_0+2+\alpha -\beta)}\frac{(-q_0)_{i_0} \left(q_0+ \Omega\right)_{i_0}}{(1)_{i_0}(1+\alpha -\beta)_{i_0}}\right.\nonumber\\
&\times& \prod _{k=1}^{n-1} \Bigg\{ \sum_{i_k=i_{k-1}}^{q_k} \frac{(i_k+ 2k+\alpha ) (i_k+ 2k+1+\alpha -\gamma)}{(i_k+ 2(k+1))(i_k+ 2k+2+\alpha -\beta )}\nonumber\\
&\times&  \frac{(-q_k)_{i_k}\left(q_k+4k+ \Omega\right)_{i_k}(2k+1)_{i_{k-1}}(2k+1+\alpha -\beta )_{i_{k-1}}}{(-q_k)_{i_{k-1}}\left(q_k+4k+ \Omega\right)_{i_{k-1}}(2k+1)_{i_k}(2k+1+\alpha -\beta )_{i_k}}\Bigg\} \label{eq:10060}\\
&\times& \left. \left.\sum_{i_n= i_{n-1}}^{q_n} \frac{(-q_n)_{i_n}\left(q_n+4n+ \Omega\right)_{i_n}(2n+1)_{i_{n-1}}(2n+1+\alpha -\beta )_{i_{n-1}}}{(-q_n)_{i_{n-1}}\left(q_n+4n+ \Omega\right)_{i_{n-1}}(2n+1)_{i_n}(2n+1+\alpha -\beta)_{i_n}} \eta ^{i_n} \right\} z^n \right\}\nonumber
\end{eqnarray}
 where
\begin{equation}
\begin{cases} \xi =\frac{1}{x} \cr
z = -a \xi^2 \cr
\eta = (1+a)\xi \cr
\varphi = 2\alpha -\gamma -\delta +1+\frac{1}{a}(\alpha -\beta +\delta ) \cr
\Omega = \frac{a\varphi }{(1+a)}\cr
q =-\alpha [(\alpha -\gamma -\delta +1)a-\beta +\delta]\cr
\hspace{0.6cm}-a(q_j+2j )\{\varphi +(1+1/a)(q_j+2j ) \} \;\;\mbox{as}\;j,q_j\in \mathbb{N}_{0} \cr
q_i\leq q_j \;\;\mbox{only}\;\mbox{if}\;i\leq j\;\;\mbox{where}\;i,j\in \mathbb{N}_{0}
\end{cases}\nonumber 
\end{equation} 
\subsubsection{Infinite series}
Replacing coefficients $a$, $q$, $\beta $, $\gamma $, $x$, $c_0$ and $\lambda $ by $\frac{1}{a}$, $\frac{q+\alpha [(\alpha -\gamma -\delta +1)a-\beta +\delta ]}{a}$, $\alpha-\gamma +1 $, $\alpha -\beta +1 $, $\frac{1}{x}$, 1 and zero into (\ref{eq:10020}). Multiply $x^{-\alpha }$ and the new (\ref{eq:10020}) together.
\begin{eqnarray}
& &x^{-\alpha } y(\xi)\nonumber\\
&=& x^{-\alpha } Hl\left(\frac{1}{a},\frac{q+\alpha [(\alpha -\gamma -\delta +1)a-\beta +\delta ]}{a}; \alpha , \alpha -\gamma +1, \alpha -\beta +1,\delta ;\frac{1}{x}\right) \nonumber\\
&=& x^{-\alpha } \left\{\sum_{i_0=0}^{\infty } \frac{\left(\Delta_0^{-}\right)_{i_0} \left(\Delta_0^{+}\right)_{i_0}}{(1)_{i_0}(1+\alpha -\beta )_{i_0}} \eta ^{i_0}\right.\nonumber\\
&+& \left\{ \sum_{i_0=0}^{\infty }\frac{(i_0+\alpha ) (i_0+1+\alpha -\gamma )}{(i_0+2)(i_0+2+\alpha -\beta  )}\frac{\left(\Delta_0^{-}\right)_{i_0} \left(\Delta_0^{+}\right)_{i_0}}{(1)_{i_0}(1+\alpha -\beta )_{i_0}}\sum_{i_1=i_0}^{\infty } \frac{\left(\Delta_1^{-}\right)_{i_1} \left(\Delta_1^{+}\right)_{i_1}(3)_{i_0}(3+\alpha -\beta )_{i_0}}{\left(\Delta_1^{-}\right)_{i_0} \left(\Delta_1^{+}\right)_{i_0}(3)_{i_1}(3+\alpha -\beta )_{i_1}} \eta ^{i_1}\right\} z\nonumber\\
&+& \sum_{n=2}^{\infty } \left\{ \sum_{i_0=0}^{\infty } \frac{(i_0+\alpha ) (i_0+1+\alpha -\gamma )}{(i_0+2)(i_0+2+\alpha -\beta )}\frac{\left(\Delta_0^{-}\right)_{i_0} \left(\Delta_0^{+}\right)_{i_0}}{(1)_{i_0}(1+\alpha -\beta )_{i_0}}\right.\nonumber\\
&\times& \prod _{k=1}^{n-1} \left\{ \sum_{i_k=i_{k-1}}^{\infty } \frac{(i_k+ 2k+\alpha ) (i_k+ 2k+1+\alpha -\gamma )}{(i_k+ 2(k+1))(i_k+ 2k+2+\alpha -\beta )} \frac{\left(\Delta_k^{-}\right)_{i_k} \left(\Delta_k^{+} \right)_{i_k}(2k+1)_{i_{k-1}}(2k+1+\alpha -\beta )_{i_{k-1}}}{\left( \Delta_k^{-}\right)_{i_{k-1}} \left(\Delta_k^{+}\right)_{i_{k-1}}(2k+1)_{i_k}(2k+1+\alpha -\beta )_{i_k}}\right\} \nonumber\\
&\times&  \left.\left.\sum_{i_n= i_{n-1}}^{\infty } \frac{\left(\Delta_n^{-}\right)_{i_n} \left(\Delta_n^{+}\right)_{i_n}(2n+1)_{i_{n-1}}(2n+1+\alpha -\beta )_{i_{n-1}}}{\left(\Delta_n^{-}\right)_{i_{n-1}} \left(\Delta_n^{+}\right)_{i_{n-1}}(2n+1)_{i_n}(2n+1+\alpha -\beta )_{i_n}} \eta ^{i_n} \right\} z^n \right\} \label{eq:10061}
\end{eqnarray}
where
\begin{equation}
\begin{cases} 
\Delta_0^{\pm}= \frac{\varphi \pm\sqrt{\varphi ^2-4(1+1/a)q}}{2(1+1/a)}\cr
\Delta_1^{\pm}=   \frac{\{\varphi +4(1+1/a)\}\pm\sqrt{\varphi ^2-4(1+1/a)q}}{2(1+1/a)}\cr
\Delta_k^{\pm}=  \frac{\{\varphi +4(1+1/a)k\}\pm\sqrt{\varphi ^2-4(1+1/a)q}}{2(1+1/a)}\cr
\Delta_n^{\pm}=  \frac{\{\varphi +4(1+1/a)n\}\pm\sqrt{\varphi ^2-4(1+1/a)q}}{2(1+1/a)}
\end{cases}\nonumber 
\end{equation}
\subsection{ ${\displaystyle \left(1-\frac{x}{a} \right)^{-\beta } Hl\left(1-a, -q+\gamma \beta; -\alpha +\gamma +\delta, \beta, \gamma, \delta; \frac{(1-a)x}{x-a} \right)}$}
\subsubsection{Polynomial of type 2}
Replacing coefficients $a$, $q$, $\alpha $, $x$, $c_0$ and $\lambda $ by $1-a$, $-q+\gamma \beta $, $-\alpha+\gamma +\delta $, $\frac{(1-a)x}{x-a}$, 1 and zero into (\ref{eq:1007}). Multiply $\left(1-\frac{x}{a} \right)^{-\beta }$ and the new (\ref{eq:1007}) together.
\begin{eqnarray}
 && \left(1-\frac{x}{a} \right)^{-\beta } y(\xi ) \nonumber\\
 &=& \left(1-\frac{x}{a} \right)^{-\beta } Hl\left(1-a, -q+\gamma \beta; -\alpha +\gamma +\delta, \beta, \gamma, \delta; \frac{(1-a)x}{x-a} \right) \nonumber\\
&=& \left(1-\frac{x}{a} \right)^{-\beta } \left\{\sum_{i_0=0}^{q_0} \frac{(-q_0)_{i_0} \left(q_0+\Omega\right)_{i_0}}{(1)_{i_0}(\gamma)_{i_0}} \eta ^{i_0}\right.\nonumber\\
&+& \left\{ \sum_{i_0=0}^{q_0}\frac{(i_0-\alpha +\gamma +\delta  ) (i_0 +\beta )}{(i_0 +2)(i_0 +1+\gamma )}\frac{(-q_0)_{i_0} \left(q_0+ \Omega\right)_{i_0}}{(1)_{i_0}(\gamma)_{i_0}} \right. \left. \sum_{i_1=i_0}^{q_1} \frac{(-q_1)_{i_1}\left(q_1+4+\Omega\right)_{i_1}(3)_{i_0}(2+\gamma)_{i_0}}{(-q_1)_{i_0}\left(q_1+4+\Omega\right)_{i_0}(3 )_{i_1}(2+\gamma)_{i_1}} \eta ^{i_1}\right\} z\nonumber\\
&+& \sum_{n=2}^{\infty } \left\{ \sum_{i_0=0}^{q_0} \frac{(i_0-\alpha +\gamma +\delta ) (i_0+\beta )}{(i_0+2)(i_0+1+\gamma )}\frac{(-q_0)_{i_0} \left(q_0+\Omega\right)_{i_0}}{(1)_{i_0}(\gamma)_{i_0}}\right.\nonumber\\
&\times& \prod _{k=1}^{n-1} \left\{ \sum_{i_k=i_{k-1}}^{q_k} \frac{(i_k+ 2k-\alpha +\gamma +\delta  ) (i_k+ 2k+\beta )}{(i_k+ 2(k+1))(i_k+ 2k+1+\gamma )}\right.  \left.\frac{(-q_k)_{i_k}\left(q_k+4k+ \Omega\right)_{i_k}(2k+1 )_{i_{k-1}}(2k+\gamma )_{i_{k-1}}}{(-q_k)_{i_{k-1}}\left( q_k+4k+ \Omega\right)_{i_{k-1}}(2k+1)_{i_k}(2k+\gamma )_{i_k}}\right\} \nonumber\\
&\times& \left. \left.\sum_{i_n= i_{n-1}}^{q_n} \frac{(-q_n)_{i_n}\left(q_n+4n+ \Omega\right)_{i_n}(2n+1)_{i_{n-1}}(2n+\gamma )_{i_{n-1}}}{(-q_n)_{i_{n-1}}\left(q_n+4n+ \Omega\right)_{i_{n-1}}(2n+1 )_{i_n}(2n+\gamma )_{i_n}} \eta ^{i_n} \right\} z^n \right\}\label{eq:100100}
\end{eqnarray}
 where
 \begin{equation}
\begin{cases} \xi = \frac{(1-a)x}{x-a} \cr
z = -\frac{1}{1-a}\xi^2 \cr
\eta = \frac{2-a}{1-a} \xi \cr
\varphi = -\alpha +\beta +\gamma +(1-a)(\gamma +\delta -1) \cr
\Omega = \frac{a\varphi }{(2-a)}\cr
q=\gamma \beta +(q_j+2j)\{\varphi +(2-a)(q_j+2j) \} \;\;\mbox{as}\;j,q_j\in \mathbb{N}_{0} \cr
q_i\leq q_j \;\;\mbox{only}\;\mbox{if}\;i\leq j\;\;\mbox{where}\;i,j\in \mathbb{N}_{0} 
\end{cases}\nonumber 
\end{equation} 
\subsubsection{Infinite series}
Replacing coefficients $a$, $q$, $\alpha $, $x$, $c_0$ and $\lambda $ by $1-a$, $-q+\gamma \beta $, $-\alpha+\gamma +\delta $, $\frac{(1-a)x}{x-a}$, 1 and zero into (\ref{eq:10020}). Multiply $\left(1-\frac{x}{a} \right)^{-\beta }$ and the new (\ref{eq:10020}) together.
 \begin{eqnarray}
 && \left(1-\frac{x}{a} \right)^{-\beta } y(\xi ) \nonumber\\
 &=& \left(1-\frac{x}{a} \right)^{-\beta } Hl\left(1-a, -q+\gamma \beta; -\alpha +\gamma +\delta, \beta, \gamma, \delta; \frac{(1-a)x}{x-a} \right) \nonumber\\
&=& \left(1-\frac{x}{a} \right)^{-\beta } \left\{\sum_{i_0=0}^{\infty } \frac{\left(\Delta_0^{-}\right)_{i_0} \left(\Delta_0^{+}\right)_{i_0}}{(1)_{i_0}(\gamma )_{i_0}} \eta ^{i_0}\right.\nonumber\\
&+& \left\{ \sum_{i_0=0}^{\infty }\frac{(i_0 -\alpha+\gamma +\delta) (i_0 +\beta )}{(i_0 +2)(i_0 +1+\gamma )}\frac{\left(\Delta_0^{-}\right)_{i_0} \left(\Delta_0^{+}\right)_{i_0}}{(1)_{i_0}(\gamma )_{i_0}}\sum_{i_1=i_0}^{\infty } \frac{\left(\Delta_1^{-}\right)_{i_1} \left(\Delta_1^{+}\right)_{i_1}(3 )_{i_0}(2+\gamma )_{i_0}}{\left(\Delta_1^{-}\right)_{i_0} \left(\Delta_1^{+}\right)_{i_0}(3)_{i_1}(2+\gamma )_{i_1}} \eta ^{i_1}\right\} z\nonumber\\
&+& \sum_{n=2}^{\infty } \left\{ \sum_{i_0=0}^{\infty } \frac{(i_0-\alpha+\gamma +\delta ) (i_0 +\beta )}{(i_0 +2)(i_0+1+\gamma )}\frac{\left(\Delta_0^{-}\right)_{i_0} \left(\Delta_0^{+}\right)_{i_0}}{(1)_{i_0}(\gamma)_{i_0}}\right.\nonumber\\
&\times& \prod _{k=1}^{n-1} \left\{ \sum_{i_k=i_{k-1}}^{\infty } \frac{(i_k+ 2k -\alpha+\gamma +\delta ) (i_k+ 2k +\beta )}{(i_k+ 2(k+1) )(i_k+ 2k+1+\gamma )} \frac{\left(\Delta_k^{-}\right)_{i_k} \left(\Delta_k^{+}\right)_{i_k}(2k+1 )_{i_{k-1}}(2k+\gamma )_{i_{k-1}}}{\left(\Delta_k^{-}\right)_{i_{k-1}} \left(\Delta_k^{+}\right)_{i_{k-1}}(2k+1)_{i_k}(2k+\gamma)_{i_k}}\right\} \nonumber\\
&\times& \left.\left.\sum_{i_n= i_{n-1}}^{\infty } \frac{\left(\Delta_n^{-}\right)_{i_n}  \left(\Delta_n^{+}\right)_{i_n}(2n+1)_{i_{n-1}}(2n+\gamma )_{i_{n-1}}}{\left(\Delta_n^{-}\right)_{i_{n-1}} \left(\Delta_n^{+}\right)_{i_{n-1}}(2n+1)_{i_n}(2n+\gamma )_{i_n}} \eta ^{i_n} \right\} z^n \right\} \label{eq:100101}
\end{eqnarray}
 where
 \begin{equation}
\begin{cases} 
\Delta_0^{\pm}= \frac{\varphi \pm\sqrt{\varphi ^2-4(2-a)q}}{2(2-a)}\cr
\Delta_1^{\pm}=  \frac{\{\varphi +4(2-a) \}\pm\sqrt{\varphi ^2-4(2-a)q}}{2(2-a)}\cr
\Delta_k^{\pm}=   \frac{\{\varphi +4(2-a)k\}\pm\sqrt{\varphi ^2-4(2-a)q}}{2(2-a)}\cr
\Delta_n^{\pm}=   \frac{\{\varphi +4(2-a)n \}\pm\sqrt{\varphi ^2-4(2-a)q}}{2(2-a)}
\end{cases}\nonumber 
\end{equation}
\subsection{ \footnotesize ${\displaystyle (1-x)^{1-\delta }\left(1-\frac{x}{a} \right)^{-\beta+\delta -1} Hl\left(1-a, -q+\gamma [(\delta -1)a+\beta -\delta +1]; -\alpha +\gamma +1, \beta -\delta+1, \gamma, 2-\delta; \frac{(1-a)x}{x-a} \right)}$ \normalsize}
\subsubsection{Polynomial of type 2}
Replacing coefficients $a$, $q$, $\alpha $, $\beta $, $\delta $, $x$, $c_0$ and $\lambda $ by $1-a$, $-q+\gamma [(\delta -1)a+\beta -\delta +1]$, $-\alpha +\gamma +1$, $\beta -\delta+1$, $2-\delta $, $\frac{(1-a)x}{x-a}$, 1 and zero into (\ref{eq:1007}). Multiply $(1-x)^{1-\delta }\left(1-\frac{x}{a} \right)^{-\beta+\delta -1}$ and the new (\ref{eq:1007}) together.
 \begin{eqnarray}
 && (1-x)^{1-\delta }\left(1-\frac{x}{a} \right)^{-\beta+\delta -1} y(\xi ) \nonumber\\
 &=& (1-x)^{1-\delta }\left(1-\frac{x}{a} \right)^{-\beta+\delta -1} Hl\bigg(1-a, -q+\gamma [(\delta -1)a+\beta -\delta +1]; -\alpha +\gamma +1, \beta -\delta+1, \gamma\nonumber\\
&&, 2-\delta; \frac{(1-a)x}{x-a} \bigg) \nonumber\\
&=& (1-x)^{1-\delta }\left(1-\frac{x}{a} \right)^{-\beta+\delta -1} \left\{\sum_{i_0=0}^{q_0} \frac{(-q_0)_{i_0} \left(q_0+ \Omega\right)_{i_0}}{(1)_{i_0}(\gamma )_{i_0}} \eta ^{i_0}\right.\nonumber\\
&+& \left\{ \sum_{i_0=0}^{q_0}\frac{(i_0-\alpha +\gamma +1) (i_0+\beta -\delta +1)}{(i_0 +2)(i_0 +1+\gamma )}\frac{(-q_0)_{i_0} \left(q_0+ \Omega\right)_{i_0}}{(1)_{i_0}(\gamma )_{i_0}} \right. \nonumber\\
&\times& \left. \sum_{i_1=i_0}^{q_1} \frac{(-q_1)_{i_1}\left(q_1+4+ \Omega\right)_{i_1}(3)_{i_0}(2+\gamma )_{i_0}}{(-q_1)_{i_0}\left(q_1+4+ \Omega\right)_{i_0}(3)_{i_1}(2+\gamma)_{i_1}} \eta ^{i_1}\right\} z\nonumber\\
&+& \sum_{n=2}^{\infty } \left\{ \sum_{i_0=0}^{q_0} \frac{(i_0-\alpha +\gamma +1) (i_0+\beta -\delta +1)}{(i_0 +2)(i_0 +1+\gamma )}\frac{(-q_0)_{i_0} \left(q_0+ \Omega\right)_{i_0}}{(1)_{i_0}(\gamma )_{i_0}}\right.\nonumber\\
&\times& \prod _{k=1}^{n-1} \Bigg\{ \sum_{i_k=i_{k-1}}^{q_k} \frac{(i_k+ 2k -\alpha +\gamma +1) (i_k+ 2k +\beta -\delta +1)}{(i_k+ 2(k+1) )(i_k+ 2k+1+\gamma )}\nonumber\\
&\times& \frac{(-q_k)_{i_k}\left(q_k+4k+ \Omega\right)_{i_k}(2k+1 )_{i_{k-1}}(2k+\gamma)_{i_{k-1}}}{(-q_k)_{i_{k-1}}\left(q_k+4k+ \Omega\right)_{i_{k-1}}(2k+1)_{i_k}(2k+\gamma )_{i_k}}\Bigg\} \nonumber\\
&\times& \left. \left.\sum_{i_n= i_{n-1}}^{q_n} \frac{(-q_n)_{i_n}\left(q_n+4n+ \Omega\right)_{i_n}(2n+1 )_{i_{n-1}}(2n+\gamma )_{i_{n-1}}}{(-q_n)_{i_{n-1}}\left(q_n+4n+ \Omega\right)_{i_{n-1}}(2n+1)_{i_n}(2n+\gamma)_{i_n}} \eta ^{i_n} \right\} z^n \right\} \label{eq:100102}
\end{eqnarray}
 where
 \begin{equation}
\begin{cases} 
\xi= \frac{(1-a)x}{x-a} \cr
z = -\frac{1}{1-a}\xi^2 \cr
\eta = \frac{(2-a)}{1-a} \xi \cr
\varphi = -\alpha +\beta +\gamma +(1-a)(\gamma -\delta +1) \cr
\Omega = \frac{\varphi }{(2-a)}\cr
q= \gamma [(\delta -1)a+\beta -\delta +1]+(q_j+2j)\{\varphi +(2-a)(q_j+2j) \} \;\;\mbox{as}\;j,q_j\in \mathbb{N}_{0} \cr
q_i\leq q_j \;\;\mbox{only}\;\mbox{if}\;i\leq j\;\;\mbox{where}\;i,j\in \mathbb{N}_{0} 
\end{cases}\nonumber 
\end{equation} 
\subsubsection{Infinite series}
Replacing coefficients $a$, $q$, $\alpha $, $\beta $, $\delta $, $x$, $c_0$ and $\lambda $ by $1-a$, $-q+\gamma [(\delta -1)a+\beta -\delta +1]$, $-\alpha +\gamma +1$, $\beta -\delta+1$, $2-\delta $, $\frac{(1-a)x}{x-a}$, 1 and zero into (\ref{eq:10020}). Multiply $(1-x)^{1-\delta }\left(1-\frac{x}{a} \right)^{-\beta+\delta -1}$ and the new (\ref{eq:10020}) together.
 \begin{eqnarray}
 && (1-x)^{1-\delta }\left(1-\frac{x}{a} \right)^{-\beta+\delta -1} y(\xi ) \nonumber\\
 &=& (1-x)^{1-\delta }\left(1-\frac{x}{a} \right)^{-\beta+\delta -1} Hl\bigg(1-a, -q+\gamma [(\delta -1)a+\beta -\delta +1]; -\alpha +\gamma +1, \beta -\delta+1\nonumber\\
&&, \gamma, 2-\delta; \frac{(1-a)x}{x-a} \bigg) \nonumber\\
&=& (1-x)^{1-\delta }\left(1-\frac{x}{a} \right)^{-\beta+\delta -1} \left\{\sum_{i_0=0}^{\infty } \frac{\left(\Delta_0^{-}\right)_{i_0} \left(\Delta_0^{+}\right)_{i_0}}{(1)_{i_0}(\gamma)_{i_0}} \eta ^{i_0}\right.\nonumber\\
&+& \left\{ \sum_{i_0=0}^{\infty }\frac{(i_0 -\alpha +\gamma +1) (i_0 +\beta-\delta +1 )}{(i_0+2)(i_0 +1+\gamma )}\frac{\left(\Delta_0^{-}\right)_{i_0} \left(\Delta_0^{+}\right)_{i_0}}{(1)_{i_0}(\gamma )_{i_0}}\sum_{i_1=i_0}^{\infty } \frac{\left(\Delta_1^{-}\right)_{i_1} \left(\Delta_1^{+} \right)_{i_1}(3 )_{i_0}(2+\gamma )_{i_0}}{\left(\Delta_1^{-}\right)_{i_0} \left(\Delta_1^{+}\right)_{i_0}(3)_{i_1}(2+\gamma)_{i_1}} \eta ^{i_1}\right\} z\nonumber\\
&+& \sum_{n=2}^{\infty } \left\{ \sum_{i_0=0}^{\infty } \frac{(i_0 -\alpha +\gamma +1) (i_0+\beta-\delta +1)}{(i_0+2)(i_0 +1+\gamma )}\frac{\left(\Delta_0^{-}\right)_{i_0} \left(\Delta_0^{+}\right)_{i_0}}{(1)_{i_0}(\gamma )_{i_0}}\right.\nonumber\\
&\times& \prod _{k=1}^{n-1} \left\{ \sum_{i_k=i_{k-1}}^{\infty } \frac{(i_k+ 2k -\alpha +\gamma +1) (i_k+ 2k+\beta-\delta +1)}{(i_k+ 2(k+1) )(i_k+ 2k+1+\gamma )} \frac{\left(\Delta_k^{-}\right)_{i_k} \left(\Delta_k^{+}\right)_{i_k}(2k+1)_{i_{k-1}}(2k+\gamma )_{i_{k-1}}}{\left(\Delta_k^{-}\right)_{i_{k-1}} \left(\Delta_k^{+}\right)_{i_{k-1}}(2k+1)_{i_k}(2k+\gamma )_{i_k}}\right\} \nonumber\\
&\times& \left.\left.\sum_{i_n= i_{n-1}}^{\infty } \frac{\left(\Delta_n^{-}\right)_{i_n} \left(\Delta_n^{+}\right)_{i_n}(2n+1)_{i_{n-1}}(2n+\gamma )_{i_{n-1}}}{ \left(\Delta_n^{-}\right)_{i_{n-1}} \left(\Delta_n^{+}\right)_{i_{n-1}}(2n+1)_{i_n}(2n+\gamma)_{i_n}} \eta ^{i_n} \right\} z^n \right\} \label{eq:100103}
\end{eqnarray}
where
  \begin{equation}
\begin{cases} 
\Delta_0^{\pm}= \frac{\varphi \pm\sqrt{\varphi ^2-4(2-a)q}}{2(2-a)}\cr
\Delta_1^{\pm}= \frac{\{\varphi +4(2-a) \}\pm\sqrt{\varphi ^2-4(2-a)q}}{2(2-a)}\cr
\Delta_k^{\pm}= \frac{\{\varphi +4(2-a)k\}\pm\sqrt{\varphi ^2-4(2-a)q}}{2(2-a)}\cr
\Delta_n^{\pm}=  \frac{\{\varphi +4(2-a)n\}\pm\sqrt{\varphi ^2-4(2-a)q}}{2(2-a)}
\end{cases}\nonumber 
\end{equation} 
\subsection{ \footnotesize ${\displaystyle x^{-\alpha } Hl\left(\frac{a-1}{a}, \frac{-q+\alpha (\delta a+\beta -\delta )}{a}; \alpha, \alpha -\gamma +1, \delta , \alpha -\beta +1; \frac{x-1}{x} \right)}$\normalsize}
\subsubsection{Polynomial of type 2}
Replacing coefficients $a$, $q$, $\beta $, $\gamma $, $\delta $, $x$, $c_0$ and $\lambda $ by $\frac{a-1}{a}$, $\frac{-q+\alpha (\delta a+\beta -\delta )}{a}$, $\alpha -\gamma +1$, $\delta $, $\alpha -\beta +1$, $\frac{x-1}{x}$, 1 and zero into (\ref{eq:1007}). Multiply $x^{-\alpha }$ and the new (\ref{eq:1007}) together.
  \begin{eqnarray}
&& x^{-\alpha } y(\xi ) \nonumber\\
 &=& x^{-\alpha } Hl\left(\frac{a-1}{a}, \frac{-q+\alpha (\delta a+\beta -\delta )}{a}; \alpha, \alpha -\gamma +1, \delta , \alpha -\beta +1; \frac{x-1}{x} \right) \nonumber\\
&=& x^{-\alpha } \left\{\sum_{i_0=0}^{q_0} \frac{(-q_0)_{i_0} \left(q_0+ \Omega\right)_{i_0}}{(1 )_{i_0}(\delta )_{i_0}} \eta ^{i_0}\right.\nonumber\\
&+& \left\{ \sum_{i_0=0}^{q_0}\frac{(i_0 +\alpha ) (i_0 +\alpha -\gamma +1)}{(i_0 +2)(i_0 +1+\delta )}\frac{(-q_0)_{i_0} \left(q_0+\Omega\right)_{i_0}}{(1)_{i_0}(\delta )_{i_0}} \sum_{i_1=i_0}^{q_1} \frac{(-q_1)_{i_1}\left(q_1+4+ \Omega\right)_{i_1}(3)_{i_0}(2+\delta )_{i_0}}{(-q_1)_{i_0}\left(q_1+4+\Omega\right)_{i_0}(3)_{i_1}(2+\delta )_{i_1}} \eta ^{i_1}\right\} z\nonumber\\
&+& \sum_{n=2}^{\infty } \left\{ \sum_{i_0=0}^{q_0} \frac{(i_0 +\alpha ) (i_0 +\alpha -\gamma +1)}{(i_0 +2)(i_0 +1+\delta )}\frac{(-q_0)_{i_0} \left(q_0+\Omega\right)_{i_0}}{(1)_{i_0}(\delta )_{i_0}}\right.\nonumber\\
&\times& \prod _{k=1}^{n-1} \left\{ \sum_{i_k=i_{k-1}}^{q_k} \frac{(i_k+ 2k +\alpha ) (i_k+ 2k +\alpha -\gamma +1)}{(i_k+ 2(k+1))(i_k+ 2k+1+\delta )}\frac{(-q_k)_{i_k}\left(q_k + 4k +\Omega\right)_{i_k}(2k+1)_{i_{k-1}}(2k+\delta )_{i_{k-1}}}{(-q_k)_{i_{k-1}}\left(q_k +4k +\Omega\right)_{i_{k-1}}(2k+1 )_{i_k}(2k+\delta )_{i_k}}\right\} \nonumber\\
&\times& \left. \left.\sum_{i_n= i_{n-1}}^{q_n} \frac{(-q_n)_{i_n}\left(q_n+4n+\Omega\right)_{i_n}(2n+1)_{i_{n-1}}(2n+\delta )_{i_{n-1}}}{(-q_n)_{i_{n-1}}\left(q_n+4n+\Omega\right)_{i_{n-1}}(2n+1)_{i_n}(2n+\delta )_{i_n}} \eta ^{i_n} \right\} z^n \right\}\label{eq:100104} 
\end{eqnarray}
 where
 \begin{equation}
\begin{cases} 
\xi= \frac{x-1}{x} \cr
z = \frac{-a}{a-1}\xi^2 \cr
\eta = \frac{2a-1}{a-1} \xi \cr
\varphi = \alpha +\beta -\gamma +\frac{a-1}{a}(\alpha -\beta +\delta ) \cr
\Omega = \frac{a\varphi }{(2a-1)}\cr
q= \alpha (\delta a+\beta -\delta )+(q_j+2j)\{a\varphi +(2a-1)(q_j+2j) \} \;\;\mbox{as}\;j,q_j\in \mathbb{N}_{0} \cr
q_i\leq q_j \;\;\mbox{only}\;\mbox{if}\;i\leq j\;\;\mbox{where}\;i,j\in \mathbb{N}_{0} 
\end{cases}\nonumber 
\end{equation}
\subsubsection{Infinite series}
Replacing coefficients $a$, $q$, $\beta $, $\gamma $, $\delta $, $x$, $c_0$ and $\lambda $ by $\frac{a-1}{a}$, $\frac{-q+\alpha (\delta a+\beta -\delta )}{a}$, $\alpha -\gamma +1$, $\delta $, $\alpha -\beta +1$, $\frac{x-1}{x}$, 1 and zero into (\ref{eq:10020}). Multiply $x^{-\alpha }$ and the new (\ref{eq:10020}) together.
 \begin{eqnarray}
&& x^{-\alpha } y(\xi ) \nonumber\\
 &=& x^{-\alpha } Hl\left(\frac{a-1}{a}, \frac{-q+\alpha (\delta a+\beta -\delta )}{a}; \alpha, \alpha -\gamma +1, \delta , \alpha -\beta +1; \frac{x-1}{x} \right) \nonumber\\
&=& x^{-\alpha } \left\{\sum_{i_0=0}^{\infty } \frac{\left(\Delta_0^{-}\right)_{i_0} \left(\Delta_0^{+}\right)_{i_0}}{(1)_{i_0}(\delta )_{i_0}} \eta ^{i_0}\right.\nonumber\\
&+& \left\{ \sum_{i_0=0}^{\infty }\frac{(i_0 +\alpha ) (i_0 +\alpha -\gamma +1)}{(i_0 +2)(i_0 +1+\delta )}\frac{\left(\Delta_0^{-}\right)_{i_0} \left(\Delta_0^{+}\right)_{i_0}}{(1)_{i_0}(\delta )_{i_0}}\sum_{i_1=i_0}^{\infty } \frac{\left(\Delta_1^{-}\right)_{i_1} \left(\Delta_1^{+}\right)_{i_1}(3)_{i_0}(2+\delta )_{i_0}}{\left(\Delta_1^{-}\right)_{i_0} \left(\Delta_1^{+}\right)_{i_0}(3)_{i_1}(2+\delta )_{i_1}} \eta ^{i_1}\right\} z\nonumber\\
&+& \sum_{n=2}^{\infty } \left\{ \sum_{i_0=0}^{\infty } \frac{(i_0 +\alpha ) (i_0 +\alpha -\gamma +1)}{(i_0 +2)(i_0 +1+\delta )}\frac{\left(\Delta_0^{-}\right)_{i_0} \left(\Delta_0^{+}\right)_{i_0}}{(1)_{i_0}(\delta )_{i_0}}\right.\nonumber\\
&\times& \prod _{k=1}^{n-1} \left\{ \sum_{i_k=i_{k-1}}^{\infty } \frac{(i_k+ 2k +\alpha ) (i_k+ 2k +\alpha -\gamma +1)}{(i_k+ 2(k+1))(i_k+ 2k+1+\delta )} \frac{\left(\Delta_k^{-}\right)_{i_k} \left(\Delta_k^{+}\right)_{i_k}(2k+1 )_{i_{k-1}}(2k+\delta )_{i_{k-1}}}{\left(\Delta_k^{-}\right)_{i_{k-1}} \left(\Delta_k^{+}\right)_{i_{k-1}}(2k+1 )_{i_k}(2k+\delta )_{i_k}}\right\} \nonumber\\
&\times& \left.\left.\sum_{i_n= i_{n-1}}^{\infty } \frac{\left(\Delta_n^{-}\right)_{i_n}\left(\Delta_n^{+}\right)_{i_n}(2n+1)_{i_{n-1}}(2n+\delta )_{i_{n-1}}}{\left(\Delta_n^{-}\right)_{i_{n-1}}\left(\Delta_n^{+}\right)_{i_{n-1}}(2n+1)_{i_n}(2n+\delta )_{i_n}} \eta ^{i_n} \right\} z^n \right\} \label{eq:100105}
\end{eqnarray}
 where
 \begin{equation}
\begin{cases} 
\Delta_0^{\pm}= \frac{a\varphi \pm a\sqrt{\varphi ^2-4\frac{(2a-1)}{a}q}}{2(2a-1)}\cr
\Delta_1^{\pm}=  \frac{a\left\{\varphi +4\frac{(2a-1)}{a}\right\}\pm a\sqrt{\varphi ^2-4\frac{(2a-1)}{a}q}}{2(2a-1)}\cr
\Delta_k^{\pm}= \frac{a\left\{\varphi +4\frac{(2a-1)}{a}k\right\}\pm a\sqrt{\varphi ^2-4\frac{(2a-1)}{a}q}}{2(2a-1)}\cr
\Delta_n^{\pm}=   \frac{a\left\{\varphi +4\frac{(2a-1)}{a}n\right\}\pm a\sqrt{\varphi ^2-4\frac{(2a-1)}{a}q}}{2(2a-1)}
\end{cases}\nonumber 
\end{equation}
\subsection{ ${\displaystyle \left(\frac{x-a}{1-a} \right)^{-\alpha } Hl\left(a, q-(\beta -\delta )\alpha ; \alpha , -\beta+\gamma +\delta , \delta , \gamma; \frac{a(x-1)}{x-a} \right)}$}
\subsubsection{Polynomial of type 2}
Replacing coefficients $q$, $\beta $, $\gamma $, $\delta $, $x$, $c_0$ and $\lambda $ by $q-(\beta -\delta )\alpha $, $-\beta+\gamma +\delta $, $\delta $,  $\gamma $, $\frac{a(x-1)}{x-a}$, 1 and zero into (\ref{eq:1007}). Multiply $\left(\frac{x-a}{1-a} \right)^{-\alpha }$ and the new (\ref{eq:1007}) together.
 \begin{eqnarray}
&& \left(\frac{x-a}{1-a} \right)^{-\alpha } y(\xi ) \nonumber\\
 &=& \left(\frac{x-a}{1-a} \right)^{-\alpha } Hl\left(a, q-(\beta -\delta )\alpha ; \alpha , -\beta+\gamma +\delta , \delta , \gamma; \frac{a(x-1)}{x-a} \right) \nonumber\\
&=& \left(\frac{x-a}{1-a} \right)^{-\alpha } \left\{\sum_{i_0=0}^{q_0} \frac{(-q_0)_{i_0} \left(q_0+ \Omega \right)_{i_0}}{(1)_{i_0}(\delta )_{i_0}} \eta ^{i_0}\right.\nonumber\\
&+& \left\{ \sum_{i_0=0}^{q_0}\frac{(i_0 +\alpha ) (i_0 -\beta +\gamma +\delta )}{(i_0 +2)(i_0 +1+\delta )}\frac{(-q_0)_{i_0} \left(q_0+ \Omega \right)_{i_0}}{(1)_{i_0}(\delta )_{i_0}} \right. \left. \sum_{i_1=i_0}^{q_1} \frac{(-q_1)_{i_1}\left(q_1+4+\Omega \right)_{i_1}(3)_{i_0}(2+\delta  )_{i_0}}{(-q_1)_{i_0}\left(q_1+4+\Omega \right)_{i_0}(3)_{i_1}(2+\delta )_{i_1}} \eta ^{i_1}\right\} z\nonumber\\
&\times& \prod _{k=1}^{n-1} \left\{ \sum_{i_k=i_{k-1}}^{q_k} \frac{(i_k+ 2k +\alpha ) (i_k+ 2k -\beta +\gamma +\delta )}{(i_k+ 2(k+1))(i_k+ 2k+1+\delta )}\right.  \left.\frac{(-q_k)_{i_k}\left(q_k+4k+ \Omega \right)_{i_k}(2k+1)_{i_{k-1}}(2k+\delta )_{i_{k-1}}}{(-q_k)_{i_{k-1}}\left(q_k+4k+ \Omega \right)_{i_{k-1}}(2k+1)_{i_k}(2k+\delta )_{i_k}}\right\} \nonumber\\
&\times& \left. \left.\sum_{i_n= i_{n-1}}^{q_n} \frac{(-q_n)_{i_n}\left(q_n+4n+ \Omega \right)_{i_n}(2n+1 )_{i_{n-1}}(2n+\delta )_{i_{n-1}}}{(-q_n)_{i_{n-1}}\left(q_n+4n+\Omega \right)_{i_{n-1}}(2n+1)_{i_n}(2n+\delta )_{i_n}} \eta ^{i_n} \right\} z^n \right\}\label{eq:100106}
\end{eqnarray}
 where
 \begin{equation}
\begin{cases} 
\xi= \frac{a(x-1)}{x-a} \cr
z = -\frac{1}{a}\xi^2 \cr
\eta = \frac{(1+a)}{a} \xi \cr
\varphi = \alpha -\beta +\delta +a(\gamma +\delta -1) \cr
\Omega = \frac{\varphi}{(1+a)}\cr
q= (\beta -\delta )\alpha -(q_j+2j)\{\varphi +(1+a)(q_j+2j ) \} \;\;\mbox{as}\;j,q_j\in \mathbb{N}_{0} \cr
q_i\leq q_j \;\;\mbox{only}\;\mbox{if}\;i\leq j\;\;\mbox{where}\;i,j\in \mathbb{N}_{0} 
\end{cases}\nonumber 
\end{equation}
\subsubsection{Infinite series}
Replacing coefficients $q$, $\beta $, $\gamma $, $\delta $, $x$, $c_0$ and $\lambda $ by $q-(\beta -\delta )\alpha $, $-\beta+\gamma +\delta $, $\delta $,  $\gamma $, $\frac{a(x-1)}{x-a}$, 1 and zero into (\ref{eq:10020}). Multiply $\left(\frac{x-a}{1-a} \right)^{-\alpha }$ and the new (\ref{eq:10020}) together.
 \begin{eqnarray}
&& \left(\frac{x-a}{1-a} \right)^{-\alpha } y(\xi ) \nonumber\\
 &=& \left(\frac{x-a}{1-a} \right)^{-\alpha } Hl\left(a, q-(\beta -\delta )\alpha ; \alpha , -\beta+\gamma +\delta , \delta , \gamma; \frac{a(x-1)}{x-a} \right) \nonumber\\
&=& \left(\frac{x-a}{1-a} \right)^{-\alpha } \left\{\sum_{i_0=0}^{\infty } \frac{\left(\Delta _0^-\right)_{i_0} \left(\Delta _0^+\right)_{i_0}}{(1)_{i_0}(\delta )_{i_0}} \eta ^{i_0}\right.\nonumber\\
&+& \left\{ \sum_{i_0=0}^{\infty }\frac{(i_0 +\alpha ) (i_0 -\beta +\gamma +\delta )}{(i_0 +2)(i_0 +1+\delta )}\frac{\left(\Delta _0^-\right)_{i_0} \left(\Delta _0^+\right)_{i_0}}{(1)_{i_0}(\delta )_{i_0}}\right.   \left.\sum_{i_1=i_0}^{\infty } \frac{\left(\Delta _1^-\right)_{i_1} \left(\Delta _1^+ \right)_{i_1}(3 )_{i_0}(2+\delta )_{i_0}}{\left(\Delta _1^-\right)_{i_0} \left(\Delta _1^+\right)_{i_0}(3 )_{i_1}(2+\delta )_{i_1}} \eta ^{i_1}\right\} z\nonumber\\
&+& \sum_{n=2}^{\infty } \left\{ \sum_{i_0=0}^{\infty } \frac{(i_0 +\alpha ) (i_0 -\beta +\gamma +\delta )}{(i_0 +2)(i_0 +1+\delta )}\frac{\left(\Delta _0^-\right)_{i_0} \left(\Delta _0^+\right)_{i_0}}{(1)_{i_0}(\delta )_{i_0}}\right.\nonumber\\
&\times& \prod _{k=1}^{n-1} \left\{ \sum_{i_k=i_{k-1}}^{\infty } \frac{(i_k+ 2k +\alpha ) (i_k+ 2k-\beta +\gamma +\delta )}{(i_k+ 2(k+1))(i_k+ 2k+1+\delta )} \right.\left.\frac{\left(\Delta _k^-\right)_{i_k} \left(\Delta _k^+ \right)_{i_k}(2k+1 )_{i_{k-1}}(2k+\delta )_{i_{k-1}}}{\left(\Delta _k^-\right)_{i_{k-1}} \left(\Delta _k^+ \right)_{i_{k-1}}(2k+1 )_{i_k}(2k+\delta )_{i_k}}\right\} \nonumber\\
&\times& \left.\left.\sum_{i_n= i_{n-1}}^{\infty } \frac{\left(\Delta _n^- \right)_{i_n}\left(\Delta _n^+ \right)_{i_n}(2n+1 )_{i_{n-1}}(2n+\delta )_{i_{n-1}}}{\left(\Delta _n^-\right)_{i_{n-1}} \left(\Delta _n^+ \right)_{i_{n-1}}(2n+1 )_{i_n}(2n+\delta )_{i_n}} \eta ^{i_n} \right\} z^n \right\} \label{eq:100107}
\end{eqnarray}
where
 \begin{equation}
\begin{cases} 
\Delta _0^{\pm}= \frac{\varphi \pm \sqrt{\varphi ^2-4(1+a)q}}{2(1+a)} \cr
\Delta _1^{\pm}=  \frac{\{\varphi +4(1+a) \} \pm \sqrt{\varphi ^2-4(1+a)q}}{2(1+a)} \cr
\Delta _k^{\pm}= \frac{\{\varphi +4(1+a)k \} \pm \sqrt{\varphi ^2-4(1+a)q}}{2(1+a)} \cr
\Delta _n^{\pm}=  \frac{\{\varphi +4(1+a)n \} \pm \sqrt{\varphi ^2-4(1+a)q}}{2(1+a)}
\end{cases}\nonumber 
\end{equation}
\section{Integral formalism of 192 Heun functions}\label{App:AppendixB}
\subsection{  ${\displaystyle (1-x)^{1-\delta } Hl(a, q - (\delta  - 1)\gamma a; \alpha - \delta  + 1, \beta - \delta + 1, \gamma ,2 - \delta ; x)}$ }
\subsubsection{Polynomial of type 2}
Replacing coefficients $q$, $\alpha$, $\beta$, $\delta$, $c_0$ and $\lambda $ by $q - (\delta - 1)\gamma a $, $\alpha - \delta  + 1 $, $\beta - \delta + 1$, $2 - \delta$, 1 and zero into (\ref{eq:10039}). Multiply $(1-x)^{1-\delta }$ and the new (\ref{eq:10039}) together.
\begin{eqnarray}
& &(1-x)^{1-\delta } y(x)\nonumber\\
&=& (1-x)^{1-\delta } Hl\left(a, q - (\delta  - 1)\gamma a; \alpha - \delta + 1, \beta - \delta + 1, \gamma ,2 - \delta ; x\right)\nonumber\\
&=& (1-x)^{1-\delta }  \Bigg\{ \;_2F_1 \left(-q_0, q_0+\Omega ; \gamma ; \eta\right) + \sum_{n=1}^{\infty } \Bigg\{\prod _{k=0}^{n-1} \Bigg\{ \int_{0}^{1} dt_{n-k}\;t_{n-k}^{2(n-k)-1} \int_{0}^{1} du_{n-k}\;u_{n-k}^{2(n-k-1)+\gamma } \nonumber\\
&\times&  \frac{1}{2\pi i}  \oint dv_{n-k} \frac{1}{v_{n-k}} \left( \frac{v_{n-k}-1}{v_{n-k}} \frac{1}{1-\overleftrightarrow {w}_{n-k+1,n}(1-t_{n-k})(1-u_{n-k})v_{n-k}}\right)^{q_{n-k}} \nonumber\\
&\times& \left( 1- \overleftrightarrow {w}_{n-k+1,n}(1-t_{n-k})(1-u_{n-k})v_{n-k}\right)^{-\left(4(n-k)+\Omega \right)}\nonumber\\
&\times&  \overleftrightarrow {w}_{n-k,n}^{-(2(n-k)-1+\alpha -\delta )}\left(  \overleftrightarrow {w}_{n-k,n} \partial _{ \overleftrightarrow {w}_{n-k,n}}\right) \overleftrightarrow {w}_{n-k,n}^{\alpha -\beta} \left(  \overleftrightarrow {w}_{n-k,n} \partial _{ \overleftrightarrow {w}_{n-k,n}}\right) \overleftrightarrow {w}_{n-k,n}^{2(n-k)-1+\beta -\delta } \Bigg\}\nonumber\\
&\times& \;_2F_1 \left(-q_0, q_0+\Omega ; \gamma ; \overleftrightarrow {w}_{1,n} \right) \Bigg\} z^n \Bigg\} \label{eq:10062}
\end{eqnarray}
where
 \begin{equation}
\begin{cases} z = -\frac{1}{a}x^2 \cr
\eta = \frac{(1+a)}{a} x \cr
\varphi = \alpha +\beta -\delta +a(\gamma -\delta +1) \cr
\Omega = \frac{\varphi }{(1+a)}\cr
q = (\delta - 1)\gamma a-(q_j+2j )\{\varphi +(1+a)(q_j+2j) \} \;\;\mbox{as}\;j,q_j\in \mathbb{N}_{0} \cr
q_i\leq q_j \;\;\mbox{only}\;\mbox{if}\;i\leq j\;\;\mbox{where}\;i,j\in \mathbb{N}_{0} 
\end{cases}\nonumber 
\end{equation}
 \subsubsection{Infinite series}
Replacing coefficients $q$, $\alpha$, $\beta$, $\delta$, $c_0$ and $\lambda $ by $q - (\delta - 1)\gamma a $, $\alpha - \delta  + 1 $, $\beta - \delta + 1$, $2 - \delta$, 1 and zero into (\ref{eq:10049}). Multiply $(1-x)^{1-\delta }$ and the new (\ref{eq:10049}) together.
\begin{eqnarray}
 & &(1-x)^{1-\delta } y(x)\nonumber\\
&=& (1-x)^{1-\delta } Hl\left(a, q - (\delta  - 1)\gamma a; \alpha - \delta + 1, \beta - \delta + 1, \gamma ,2 - \delta ; x\right)\nonumber\\
&=& (1-x)^{1-\delta } \Bigg\{ \;_2F_1 \left( -\Pi_0 ^{+}, -\Pi_0 ^{-}; \gamma ; \eta\right) \nonumber\\
&+& \sum_{n=1}^{\infty } \left\{\prod _{k=0}^{n-1} \Bigg\{ \int_{0}^{1} dt_{n-k}\;t_{n-k}^{2(n-k)-1} \int_{0}^{1} du_{n-k}\;u_{n-k}^{2(n-k-1)+\gamma} \right.\nonumber\\
&\times& \frac{1}{2\pi i}  \oint dv_{n-k} \frac{1}{v_{n-k}} \left( \frac{v_{n-k}-1}{v_{n-k}}\right)^{\Pi_{n-k} ^{+}} \left( 1- \overleftrightarrow {w}_{n-k+1,n}(1-t_{n-k})(1-u_{n-k})v_{n-k}\right)^{\Pi_{n-k} ^{-}}\nonumber\\
&\times& \overleftrightarrow {w}_{n-k,n}^{-(2(n-k)-1+\alpha -\delta )}\left(  \overleftrightarrow {w}_{n-k,n} \partial _{ \overleftrightarrow {w}_{n-k,n}}\right) \overleftrightarrow {w}_{n-k,n}^{\alpha -\beta} \left(  \overleftrightarrow {w}_{n-k,n} \partial _{ \overleftrightarrow {w}_{n-k,n}}\right) \overleftrightarrow {w}_{n-k,n}^{2(n-k)-1+\beta-\delta } \Bigg\}\nonumber\\
&\times&  \;_2F_1 \left( -\Pi_0 ^{+}, -\Pi_0 ^{-}; \gamma ; \overleftrightarrow {w}_{1,n}\right)\Bigg\} z^n \Bigg\} \label{eq:10063}
\end{eqnarray}
where
 \begin{equation}
\begin{cases} 
\Pi_0 ^{\pm}= \frac{-\varphi \pm\sqrt{\varphi ^2-4(1+a)(q - (\delta  - 1)\gamma a)}}{2(1+a)} \cr
\Pi_{n-k} ^{\pm}= \frac{-(\varphi +4(1+a)(n-k))\pm\sqrt{\varphi ^2-4(1+a)(q - (\delta  - 1)\gamma a)}}{2(1+a)} 
\end{cases}\nonumber 
\end{equation}
\subsection{  \footnotesize ${\displaystyle x^{1-\gamma } (1-x)^{1-\delta } Hl(a, q-(\gamma +\delta -2)a-(\gamma -1)(\alpha +\beta -\gamma -\delta +1); \alpha - \gamma -\delta +2, \beta - \gamma -\delta +2, 2-\gamma, 2 - \delta ; x)}$ \normalsize}
\subsubsection{Polynomial of type 2}
Replacing coefficients $q$, $\alpha$, $\beta$, $\gamma $, $\delta$, $c_0$ and $\lambda $ by $q-(\gamma +\delta -2)a-(\gamma -1)(\alpha +\beta -\gamma -\delta +1)$, $\alpha - \gamma -\delta +2$, $\beta - \gamma -\delta +2, 2-\gamma$, $2 - \delta$,1 and zero into (\ref{eq:10039}). Multiply $x^{1-\gamma } (1-x)^{1-\delta }$ and the new (\ref{eq:10039}) together.
\begin{eqnarray}
& &x^{1-\gamma } (1-x)^{1-\delta } y(x)\nonumber\\
&=& x^{1-\gamma } (1-x)^{1-\delta }  Hl(a, q-(\gamma +\delta -2)a-(\gamma -1)(\alpha +\beta -\gamma -\delta +1); \alpha - \gamma -\delta +2\nonumber\\
&&, \beta - \gamma -\delta +2, 2-\gamma, 2 - \delta ; x)\nonumber\\
&=& x^{1-\gamma } (1-x)^{1-\delta } \Bigg\{ \;_2F_1 \left(-q_0, q_0+\Omega; 2-\gamma ; \eta\right)\nonumber\\
&+& \sum_{n=1}^{\infty } \Bigg\{\prod _{k=0}^{n-1} \Bigg\{ \int_{0}^{1} dt_{n-k}\;t_{n-k}^{2(n-k)-1} \int_{0}^{1} du_{n-k}\;u_{n-k}^{2(n-k)-\gamma} \nonumber\\
&\times&  \frac{1}{2\pi i}  \oint dv_{n-k} \frac{1}{v_{n-k}} \left( \frac{v_{n-k}-1}{v_{n-k}} \frac{1}{1-\overleftrightarrow {w}_{n-k+1,n}(1-t_{n-k})(1-u_{n-k})v_{n-k}}\right)^{q_{n-k}} \nonumber\\
&\times& \left( 1- \overleftrightarrow {w}_{n-k+1,n}(1-t_{n-k})(1-u_{n-k})v_{n-k}\right)^{-\left(4(n-k)+\Omega\right)}\nonumber\\
&\times& \overleftrightarrow {w}_{n-k,n}^{-(2(n-k)+\alpha-\gamma -\delta)}\left(  \overleftrightarrow {w}_{n-k,n} \partial _{ \overleftrightarrow {w}_{n-k,n}}\right) \overleftrightarrow {w}_{n-k,n}^{\alpha -\beta} \left(  \overleftrightarrow {w}_{n-k,n} \partial _{ \overleftrightarrow {w}_{n-k,n}}\right) \overleftrightarrow {w}_{n-k,n}^{2(n-k)+\beta-\gamma -\delta  } \Bigg\}\nonumber\\
&\times& \;_2F_1 \left(-q_0, q_0+\Omega; 2-\gamma ; \overleftrightarrow {w}_{1,n} \right)\Bigg\} z^n \Bigg\} \label{eq:10064}
\end{eqnarray}
where
\begin{equation}
\begin{cases} z = -\frac{1}{a}x^2 \cr
\eta = \frac{(1+a)}{a} x \cr
\varphi = \alpha +\beta -2\gamma -\delta +2 +a(3-\gamma -\delta ) \cr
\Omega = \frac{\varphi}{(1+a)}\cr
q = (\gamma +\delta -2)a+(\gamma -1)(\alpha +\beta -\gamma -\delta +1)\cr
\hspace{0.6cm}-(q_j+2j)\{\varphi +(1+a)(q_j+2j) \} \;\;\mbox{as}\;j,q_j\in \mathbb{N}_{0} \cr
q_i\leq q_j \;\;\mbox{only}\;\mbox{if}\;i\leq j\;\;\mbox{where}\;i,j\in \mathbb{N}_{0} 
\end{cases}\nonumber 
\end{equation}
\subsubsection{Infinite series}
Replacing coefficients $q$, $\alpha$, $\beta$, $\gamma $, $\delta$, $c_0$ and $\lambda $ by $q-(\gamma +\delta -2)a-(\gamma -1)(\alpha +\beta -\gamma -\delta +1)$, $\alpha - \gamma -\delta +2$, $\beta - \gamma -\delta +2, 2-\gamma$, $2 - \delta$,1 and zero into (\ref{eq:10049}). Multiply $x^{1-\gamma } (1-x)^{1-\delta }$ and the new (\ref{eq:10049}) together.
\begin{eqnarray}
& &x^{1-\gamma } (1-x)^{1-\delta } y(x)\nonumber\\
&=& x^{1-\gamma } (1-x)^{1-\delta } Hl(a, q-(\gamma +\delta -2)a-(\gamma -1)(\alpha +\beta -\gamma -\delta +1); \alpha - \gamma -\delta +2\nonumber\\
&&, \beta - \gamma -\delta +2, 2-\gamma, 2 - \delta ; x)\nonumber\\
&=& x^{1-\gamma } (1-x)^{1-\delta } \Bigg\{ \;_2F_1 \left(-\Pi_{0} ^{+}, -\Pi_{0} ^{-}; 2-\gamma ; \eta \right)\nonumber\\
&+& \sum_{n=1}^{\infty } \left\{\prod _{k=0}^{n-1} \Bigg\{ \int_{0}^{1} dt_{n-k}\;t_{n-k}^{2(n-k)-1} \int_{0}^{1} du_{n-k}\;u_{n-k}^{2(n-k)-\gamma}\right. \nonumber\\
&\times& \left. \frac{1}{2\pi i} \oint dv_{n-k} \frac{1}{v_{n-k}} \left( \frac{v_{n-k}-1}{v_{n-k}}\right)^{\Pi_{n-k} ^{+}} \right.\left( 1- \overleftrightarrow {w}_{n-k+1,n}(1-t_{n-k})(1-u_{n-k})v_{n-k}\right)^{\Pi_{n-k} ^{-}}\nonumber\\
&\times& \overleftrightarrow {w}_{n-k,n}^{-(2(n-k)+\alpha-\gamma -\delta  )}\left(  \overleftrightarrow {w}_{n-k,n} \partial _{ \overleftrightarrow {w}_{n-k,n}}\right) \overleftrightarrow {w}_{n-k,n}^{\alpha -\beta} \left(  \overleftrightarrow {w}_{n-k,n} \partial _{ \overleftrightarrow {w}_{n-k,n}}\right) \overleftrightarrow {w}_{n-k,n}^{2(n-k)+\beta-\gamma -\delta  } \Bigg\}\nonumber\\
&\times& \;_2F_1 \left(-\Pi_{0} ^{+}, -\Pi_{0} ^{-}; 2-\gamma ; \overleftrightarrow {w}_{1,n} \right)\Bigg\} z^n \Bigg\} \label{eq:10065}
\end{eqnarray}
where
 \begin{equation}
\begin{cases} 
\Pi_{0} ^{\pm}= \frac{-\varphi \pm\sqrt{\varphi ^2-4(1+a)(q-(\gamma +\delta -2)a-(\gamma -1)(\alpha +\beta -\gamma -\delta +1))}}{2(1+a)} \cr
\Pi_{n-k} ^{\pm}= \frac{-(\varphi +4(1+a)(n-k))\pm\sqrt{\varphi ^2-4(1+a)(q-(\gamma +\delta -2)a-(\gamma -1)(\alpha +\beta -\gamma -\delta +1))}}{2(1+a)} 
\end{cases}\nonumber 
\end{equation} 
\subsection{ ${\displaystyle  Hl(1-a,-q+\alpha \beta; \alpha,\beta, \delta, \gamma; 1-x)}$} 
\subsubsection{Polynomial of type 2}
Replacing coefficients $a$, $q$, $\gamma $, $\delta$, $x$, $c_0$ and $\lambda $ by $1-a$, $-q +\alpha \beta $, $\delta $, $\gamma $, $1-x$, 1 and zero into (\ref{eq:10039}).
\begin{eqnarray}
y(\xi ) &=& Hl(1-a,-q+\alpha \beta; \alpha,\beta, \delta, \gamma; 1-x)\nonumber\\
&=&  \;_2F_1 \left(-q_0, q_0+\Omega ; \delta; \eta \right) + \sum_{n=1}^{\infty } \Bigg\{\prod _{k=0}^{n-1} \Bigg\{ \int_{0}^{1} dt_{n-k}\;t_{n-k}^{2(n-k)-1} \int_{0}^{1} du_{n-k}\;u_{n-k}^{2(n-k-1)+\delta } \nonumber\\
&\times&  \frac{1}{2\pi i}  \oint dv_{n-k} \frac{1}{v_{n-k}} \left( \frac{v_{n-k}-1}{v_{n-k}} \frac{1}{1-\overleftrightarrow {w}_{n-k+1,n}(1-t_{n-k})(1-u_{n-k})v_{n-k}}\right)^{q_{n-k}} \nonumber\\
&\times& \left( 1- \overleftrightarrow {w}_{n-k+1,n}(1-t_{n-k})(1-u_{n-k})v_{n-k}\right)^{-\left(4(n-k)+\Omega \right)}\nonumber\\
&\times& \overleftrightarrow {w}_{n-k,n}^{-(2(n-k-1)+\alpha)}\left(  \overleftrightarrow {w}_{n-k,n} \partial _{ \overleftrightarrow {w}_{n-k,n}}\right) \overleftrightarrow {w}_{n-k,n}^{\alpha -\beta} \left(  \overleftrightarrow {w}_{n-k,n} \partial _{ \overleftrightarrow {w}_{n-k,n}}\right) \overleftrightarrow {w}_{n-k,n}^{2(n-k-1)+\beta} \Bigg\}\nonumber\\
&\times& \;_2F_1 \left(-q_0, q_0+\Omega ; \delta; \overleftrightarrow {w}_{1,n} \right)\Bigg\} z^n \label{eq:10066}
\end{eqnarray}
where
 \begin{equation}
\begin{cases} \xi =1-x \cr
z = \frac{-1}{1-a}\xi^2 \cr
\eta = \frac{2-a}{1-a}\xi \cr
\varphi = \alpha +\beta -\delta +(1-a)(\delta +\gamma -1) \cr
\Omega = \frac{\varphi }{(2-a)}\cr
q = \alpha \beta +(q_j+2j)\{\varphi +(2-a)(q_j+2j ) \} \;\;\mbox{as}\;j,q_j\in \mathbb{N}_{0} \cr
q_i\leq q_j \;\;\mbox{only}\;\mbox{if}\;i\leq j\;\;\mbox{where}\;i,j\in \mathbb{N}_{0}
\end{cases}\nonumber 
\end{equation}
\subsubsection{Infinite series}
Replacing coefficients $a$, $q$, $\gamma $, $\delta$, $x$, $c_0$ and $\lambda $ by $1-a$, $-q +\alpha \beta $, $\delta $, $\gamma $, $1-x$,1 and zero into (\ref{eq:10049}).
\begin{eqnarray}
y(\xi ) &=& Hl(1-a,-q+\alpha \beta; \alpha, \beta, \delta, \gamma; 1-x)\nonumber\\
&=&  \;_2F_1 \left(-\Pi_{0} ^{+}, -\Pi_{0} ^{-}; \delta; \eta \right) \nonumber\\
&+& \sum_{n=1}^{\infty } \left\{\prod _{k=0}^{n-1} \Bigg\{ \int_{0}^{1} dt_{n-k}\;t_{n-k}^{2(n-k)-1} \int_{0}^{1} du_{n-k}\;u_{n-k}^{2(n-k-1)+\delta }\right.\nonumber\\
&\times& \frac{1}{2\pi i}  \oint dv_{n-k} \frac{1}{v_{n-k}} \left( \frac{v_{n-k}-1}{v_{n-k}}\right)^{\Pi_{n-k} ^{+}}  \left( 1- \overleftrightarrow {w}_{n-k+1,n}(1-t_{n-k})(1-u_{n-k})v_{n-k}\right)^{\Pi_{n-k} ^{-}} \nonumber\\
&\times& \overleftrightarrow {w}_{n-k,n}^{-(2(n-k-1)+\alpha)}\left(  \overleftrightarrow {w}_{n-k,n} \partial _{ \overleftrightarrow {w}_{n-k,n}}\right) \overleftrightarrow {w}_{n-k,n}^{\alpha -\beta} \left(  \overleftrightarrow {w}_{n-k,n} \partial _{ \overleftrightarrow {w}_{n-k,n}}\right) \overleftrightarrow {w}_{n-k,n}^{2(n-k-1)+\beta} \Bigg\}\nonumber\\
&\times& \;_2F_1 \left(-\Pi_{0} ^{+}, -\Pi_{0} ^{-}; \delta; \overleftrightarrow {w}_{1,n} \right) \Bigg\} z^n  \label{eq:10067}
\end{eqnarray}
where
 \begin{equation}
\begin{cases} 
\Pi_{0} ^{\pm}= \frac{-\varphi \pm\sqrt{\varphi ^2-4(2-a)(-q +\alpha \beta)}}{2(2-a)} \cr
\Pi_{n-k} ^{\pm}= \frac{-(\varphi +4(2-a)(n-k))\pm\sqrt{\varphi ^2-4(2-a)(-q +\alpha \beta)}}{2(2-a)} 
\end{cases}\nonumber 
\end{equation} 
\subsection{ \footnotesize ${\displaystyle (1-x)^{1-\delta } Hl(1-a,-q+(\delta -1)\gamma a+(\alpha -\delta +1)(\beta -\delta +1); \alpha-\delta +1,\beta-\delta +1, 2-\delta, \gamma; 1-x)}$ \normalsize}
\subsubsection{Polynomial of type 2}
Replacing coefficients $a$, $q$, $\alpha $, $\beta $, $\gamma $, $\delta$, $x$, $c_0$ and $\lambda $ by $1-a$, $-q+(\delta -1)\gamma a+(\alpha -\delta +1)(\beta -\delta +1)$, $\alpha-\delta +1 $, $\beta-\delta +1 $, $2-\delta$, $\gamma $, $1-x$, 1 and zero into (\ref{eq:10039}). Multiply $(1-x)^{1-\delta }$ and the new (\ref{eq:10039}) together.
\begin{eqnarray}
& &(1-x)^{1-\delta } y(\xi)\nonumber\\
&=& (1-x)^{1-\delta } Hl(1-a,-q+(\delta -1)\gamma a+(\alpha -\delta +1)(\beta -\delta +1); \alpha-\delta +1,\beta-\delta +1, 2-\delta, \gamma; 1-x) \nonumber\\
&=& (1-x)^{1-\delta } \Bigg\{ \;_2F_1 \left(-q_0, q_0+\Omega; 2-\delta ; \eta\right) \nonumber\\
&+& \sum_{n=1}^{\infty } \Bigg\{\prod _{k=0}^{n-1} \Bigg\{ \int_{0}^{1} dt_{n-k}\;t_{n-k}^{2(n-k)-1} \int_{0}^{1} du_{n-k}\;u_{n-k}^{2(n-k)-\delta} \nonumber\\
&\times&  \frac{1}{2\pi i}  \oint dv_{n-k} \frac{1}{v_{n-k}} \left( \frac{v_{n-k}-1}{v_{n-k}} \frac{1}{1-\overleftrightarrow {w}_{n-k+1,n}(1-t_{n-k})(1-u_{n-k})v_{n-k}}\right)^{q_{n-k}} \nonumber\\
&\times& \left( 1- \overleftrightarrow {w}_{n-k+1,n}(1-t_{n-k})(1-u_{n-k})v_{n-k}\right)^{-\left(4(n-k)+\Omega\right)}\nonumber\\
&\times& \overleftrightarrow {w}_{n-k,n}^{-(2(n-k)-1+\alpha-\delta )}\left(  \overleftrightarrow {w}_{n-k,n} \partial _{ \overleftrightarrow {w}_{n-k,n}}\right) \overleftrightarrow {w}_{n-k,n}^{\alpha -\beta} \left(  \overleftrightarrow {w}_{n-k,n} \partial _{ \overleftrightarrow {w}_{n-k,n}}\right) \overleftrightarrow {w}_{n-k,n}^{2(n-k)-1+\beta-\delta } \Bigg\}\nonumber\\
&\times& \;_2F_1 \left(-q_0, q_0+\Omega; 2-\delta ;  \overleftrightarrow {w}_{1,n}\right)\Bigg\} z^n \Bigg\} \label{eq:10068}
\end{eqnarray}
where
  \begin{equation}
\begin{cases} \xi =1-x \cr
z = \frac{-1}{1-a}\xi^2 \cr
\eta = \frac{2-a}{1-a}\xi \cr
\varphi = \alpha +\beta -\gamma -2\delta +2+(1-a)(\gamma -\delta +1) \cr
\Omega = \frac{\varphi }{(2-a)}\cr
q =(\delta -1)\gamma a+(\alpha -\delta +1)(\beta -\delta +1)\cr
\hspace{0.6cm}+(q_j+2j )\{\varphi +(2-a)(q_j+2j ) \} \;\;\mbox{as}\;j,q_j\in \mathbb{N}_{0} \cr
q_i\leq q_j \;\;\mbox{only}\;\mbox{if}\;i\leq j\;\;\mbox{where}\;i,j\in \mathbb{N}_{0}
\end{cases}\nonumber 
\end{equation} 
\subsubsection{Infinite series}
Replacing coefficients $a$, $q$, $\alpha $, $\beta $, $\gamma $, $\delta$, $x$, $c_0$ and $\lambda $ by $1-a$, $ -q+(\delta -1)\gamma a+(\alpha -\delta +1)(\beta -\delta +1)$, $\alpha-\delta +1 $, $\beta-\delta +1 $, $2-\delta$, $\gamma $, $1-x$, 1 and zero into (\ref{eq:10049}).  Multiply $(1-x)^{1-\delta }$ and the new (\ref{eq:10049}) together.
\begin{eqnarray}
& &(1-x)^{1-\delta } y(\xi)\nonumber\\
&=& (1-x)^{1-\delta } Hl(1-a, -q+(\delta -1)\gamma a+(\alpha -\delta +1)(\beta -\delta +1); \alpha -\delta +1,\beta-\delta +1 \nonumber\\
&&, 2-\delta, \gamma; 1-x) \nonumber\\
&=& (1-x)^{1-\delta } \Bigg\{ \;_2F_1 \left(-\Pi_{0} ^{+}, -\Pi_{0} ^{-}; 2-\delta ; \eta\right) \nonumber\\
&+& \sum_{n=1}^{\infty } \Bigg\{\prod _{k=0}^{n-1} \Bigg\{ \int_{0}^{1} dt_{n-k}\;t_{n-k}^{2(n-k)-1} \int_{0}^{1} du_{n-k}\;u_{n-k}^{2(n-k)-\delta } \nonumber\\
&\times& \frac{1}{2\pi i}  \oint dv_{n-k} \frac{1}{v_{n-k}} \left( \frac{v_{n-k}-1}{v_{n-k}}\right)^{\Pi_{n-k} ^{+}} \left( 1- \overleftrightarrow {w}_{n-k+1,n}(1-t_{n-k})(1-u_{n-k})v_{n-k}\right)^{\Pi_{n-k} ^{-}}\nonumber\\
&\times& \overleftrightarrow {w}_{n-k,n}^{-(2(n-k)-1+\alpha-\delta  )}\left(  \overleftrightarrow {w}_{n-k,n} \partial _{ \overleftrightarrow {w}_{n-k,n}}\right) \overleftrightarrow {w}_{n-k,n}^{\alpha -\beta} \left(  \overleftrightarrow {w}_{n-k,n} \partial _{ \overleftrightarrow {w}_{n-k,n}}\right) \overleftrightarrow {w}_{n-k,n}^{2(n-k)-1+\beta-\delta } \Bigg\}\nonumber\\
&\times& \;_2F_1 \left(-\Pi_{0} ^{+}, -\Pi_{0} ^{-}; 2-\delta ; \overleftrightarrow {w}_{1,n}\right)\Bigg\} z^n \Bigg\} \label{eq:10069}
\end{eqnarray}
where
 \begin{equation}
\begin{cases} 
\Pi_{0} ^{\pm}= \frac{-(\varphi +4(2-a)(n-k))\pm\sqrt{\varphi ^2-4(2-a)(-q+(\delta -1)\gamma a+(\alpha -\delta +1)(\beta -\delta +1))}}{2(2-a)} \cr
\Pi_{n-k} ^{\pm}= \frac{-(\varphi +4(2-a)(n-k))\pm\sqrt{\varphi ^2-4(2-a)(-q+(\delta -1)\gamma a+(\alpha -\delta +1)(\beta -\delta +1))}}{2(2-a)} 
\end{cases}\nonumber 
\end{equation}
\subsection{ \footnotesize ${\displaystyle x^{-\alpha } Hl\left(\frac{1}{a},\frac{q+\alpha [(\alpha -\gamma -\delta +1)a-\beta +\delta ]}{a}; \alpha , \alpha -\gamma +1, \alpha -\beta +1,\delta ;\frac{1}{x}\right)}$\normalsize}
\subsubsection{Polynomial of type 2}
Replacing coefficients $a$, $q$, $\beta $, $\gamma $, $x$, $c_0$ and $\lambda $ by $\frac{1}{a}$, $\frac{q+\alpha [(\alpha -\gamma -\delta +1)a-\beta +\delta ]}{a}$, $\alpha-\gamma +1 $, $\alpha -\beta +1 $, $\frac{1}{x}$, 1 and zero into (\ref{eq:10039}). Multiply $x^{-\alpha }$ and the new (\ref{eq:10039}) together.
\begin{eqnarray}
& &x^{-\alpha } y(\xi)\nonumber\\
&=& x^{-\alpha }  Hl\left(\frac{1}{a},\frac{q+\alpha [(\alpha -\gamma -\delta +1)a-\beta +\delta ]}{a}; \alpha , \alpha -\gamma +1, \alpha -\beta +1,\delta ;\frac{1}{x}\right) \nonumber\\
&=& x^{-\alpha } \Bigg\{ \;_2F_1 \left(-q_0, q_0+\Omega; 1+\alpha -\beta ; \eta\right)\nonumber\\
&+& \sum_{n=1}^{\infty } \Bigg\{\prod _{k=0}^{n-1} \Bigg\{ \int_{0}^{1} dt_{n-k}\;t_{n-k}^{2(n-k)-1} \int_{0}^{1} du_{n-k}\;u_{n-k}^{2(n-k)-1+\alpha -\beta } \nonumber\\
&\times&  \frac{1}{2\pi i}  \oint dv_{n-k} \frac{1}{v_{n-k}} \left( \frac{v_{n-k}-1}{v_{n-k}} \frac{1}{1-\overleftrightarrow {w}_{n-k+1,n}(1-t_{n-k})(1-u_{n-k})v_{n-k}}\right)^{q_{n-k}} \nonumber\\
&\times& \left( 1- \overleftrightarrow {w}_{n-k+1,n}(1-t_{n-k})(1-u_{n-k})v_{n-k}\right)^{-\left(4(n-k)+\Omega\right)}\nonumber\\
&\times& \overleftrightarrow {w}_{n-k,n}^{-(2(n-k-1)+\alpha )}\left(  \overleftrightarrow {w}_{n-k,n} \partial _{ \overleftrightarrow {w}_{n-k,n}}\right) \overleftrightarrow {w}_{n-k,n}^{\gamma -1} \left(  \overleftrightarrow {w}_{n-k,n} \partial _{ \overleftrightarrow {w}_{n-k,n}}\right) \overleftrightarrow {w}_{n-k,n}^{2(n-k)-1+\alpha -\gamma } \Bigg\}\nonumber\\
&\times& \;_2F_1 \left(-q_0, q_0+\Omega; 1+\alpha -\beta ;\overleftrightarrow {w}_{1,n}\right)\Bigg\} z^n \Bigg\} \label{eq:10070}
\end{eqnarray}
where
 \begin{equation}
\begin{cases} \xi= \frac{1}{x} \cr
z = -a \xi^2 \cr
\eta = (1+a)\xi \cr
\varphi = 2\alpha -\gamma -\delta +1+\frac{1}{a}(\alpha -\beta +\delta ) \cr
\Omega = \frac{a\varphi }{(1+a)}\cr
q =-\alpha [(\alpha -\gamma -\delta +1)a-\beta +\delta]\cr
\hspace{0.6cm}-a(q_j+2j )\{\varphi +(1+1/a)(q_j+2j ) \} \;\;\mbox{as}\;j,q_j\in \mathbb{N}_{0} \cr
q_i\leq q_j \;\;\mbox{only}\;\mbox{if}\;i\leq j\;\;\mbox{where}\;i,j\in \mathbb{N}_{0}
\end{cases}\nonumber 
\end{equation}
\subsubsection{Infinite series}
Replacing coefficients $a$, $q$, $\beta $, $\gamma $, $x$, $c_0$ and $\lambda $ by $\frac{1}{a}$, $\frac{q+\alpha [(\alpha -\gamma -\delta +1)a-\beta +\delta ]}{a}$, $\alpha-\gamma +1 $, $\alpha -\beta +1 $, $\frac{1}{x}$, 1 and zero into (\ref{eq:10049}). Multiply $x^{-\alpha }$ and the new (\ref{eq:10049}) together.
\begin{eqnarray}
& &x^{-\alpha } y(\xi)\nonumber\\
&=& x^{-\alpha } Hl\left(\frac{1}{a}, \frac{q+\alpha [(\alpha -\gamma -\delta +1)a-\beta +\delta ]}{a}; \alpha , \alpha -\gamma +1, \alpha -\beta +1,\delta ;\frac{1}{x}\right) \nonumber\\
&=& x^{-\alpha } \Bigg\{ \;_2F_1 \left(-\Pi_{0} ^{+}, -\Pi_{0} ^{-}; 1+\alpha -\beta ; \eta\right) \nonumber\\
&+& \sum_{n=1}^{\infty } \left\{\prod _{k=0}^{n-1} \Bigg\{ \int_{0}^{1} dt_{n-k}\;t_{n-k}^{2(n-k)-1} \int_{0}^{1} du_{n-k}\;u_{n-k}^{2(n-k)-1+\alpha -\beta} \right.\nonumber\\
&\times& \frac{1}{2\pi i} \oint dv_{n-k} \frac{1}{v_{n-k}} \left( \frac{v_{n-k}-1}{v_{n-k}}\right)^{\Pi_{n-k} ^{+}} \left( 1- \overleftrightarrow {w}_{n-k+1,n}(1-t_{n-k})(1-u_{n-k})v_{n-k}\right)^{\Pi_{n-k} ^{-}}\nonumber\\
&\times& \overleftrightarrow {w}_{n-k,n}^{-(2(n-k-1)+\alpha )}\left(  \overleftrightarrow {w}_{n-k,n} \partial _{ \overleftrightarrow {w}_{n-k,n}}\right) \overleftrightarrow {w}_{n-k,n}^{\gamma -1} \left(  \overleftrightarrow {w}_{n-k,n} \partial _{ \overleftrightarrow {w}_{n-k,n}}\right) \overleftrightarrow {w}_{n-k,n}^{2(n-k)-1+\alpha -\gamma} \Bigg\}\nonumber\\
&\times& \;_2F_1 \left(-\Pi_{0} ^{+}, -\Pi_{0} ^{-}; 1+\alpha -\beta ; \overleftrightarrow {w}_{1,n}\right) \Bigg\} z^n \Bigg\}  \label{eq:10071}
\end{eqnarray}
where
 \begin{equation}
\begin{cases} 
\Pi_{0} ^{\pm}= \frac{-\varphi \pm\sqrt{\varphi ^2-4\left( 1+\frac{1}{a}\right)\left(\frac{q+\alpha [(\alpha -\gamma -\delta +1)a-\beta +\delta ]}{a}\right)}}{2\left( 1+\frac{1}{a}\right)} \cr
\Pi_{n-k} ^{\pm}= \frac{-\left(\varphi +4\left( 1+\frac{1}{a}\right) (n-k)\right)\pm\sqrt{\varphi ^2-4\left( 1+\frac{1}{a}\right) \left(\frac{q+\alpha [(\alpha -\gamma -\delta +1)a-\beta +\delta ]}{a}\right)}}{2\left(1+\frac{1}{a}\right)} 
\end{cases}\nonumber 
\end{equation}
\subsection{ ${\displaystyle \left(1-\frac{x}{a} \right)^{-\beta } Hl\left(1-a, -q+\gamma \beta; -\alpha +\gamma +\delta, \beta, \gamma, \delta; \frac{(1-a)x}{x-a} \right)}$}
\subsubsection{Polynomial of type 2}
Replacing coefficients $a$, $q$, $\alpha $, $x$, $c_0$ and $\lambda $ by $1-a$, $-q+\gamma \beta $, $-\alpha+\gamma +\delta $, $\frac{(1-a)x}{x-a}$, 1 and zero into (\ref{eq:10039}). Multiply $\left(1-\frac{x}{a} \right)^{-\beta }$ and the new (\ref{eq:10039}) together.
\begin{eqnarray}
 && \left(1-\frac{x}{a} \right)^{-\beta } y(\xi ) \nonumber\\
 &=& \left(1-\frac{x}{a} \right)^{-\beta } Hl\left(1-a, -q+\gamma \beta; -\alpha +\gamma +\delta, \beta, \gamma, \delta; \frac{(1-a)x}{x-a} \right) \nonumber\\
&=& \left(1-\frac{x}{a} \right)^{-\beta } \Bigg\{ \;_2F_1 \left(-q_0, q_0+\Omega; \gamma; \eta \right) \nonumber\\
&+& \sum_{n=1}^{\infty } \Bigg\{\prod _{k=0}^{n-1} \Bigg\{ \int_{0}^{1} dt_{n-k}\;t_{n-k}^{2(n-k)-1} \int_{0}^{1} du_{n-k}\;u_{n-k}^{2(n-k-1)+\gamma } \nonumber\\
&\times&  \frac{1}{2\pi i}  \oint dv_{n-k} \frac{1}{v_{n-k}} \left( \frac{v_{n-k}-1}{v_{n-k}} \frac{1}{1-\overleftrightarrow {w}_{n-k+1,n}(1-t_{n-k})(1-u_{n-k})v_{n-k}}\right)^{q_{n-k}} \nonumber\\
&\times& \left( 1- \overleftrightarrow {w}_{n-k+1,n}(1-t_{n-k})(1-u_{n-k})v_{n-k}\right)^{-\left(4(n-k)+\Omega\right)}\nonumber\\
&\times& \overleftrightarrow {w}_{n-k,n}^{-(2(n-k-1)-\alpha+\gamma +\delta )}\left(  \overleftrightarrow {w}_{n-k,n} \partial _{ \overleftrightarrow {w}_{n-k,n}}\right) \overleftrightarrow {w}_{n-k,n}^{-\alpha -\beta+\gamma +\delta } \left(  \overleftrightarrow {w}_{n-k,n} \partial _{ \overleftrightarrow {w}_{n-k,n}}\right) \overleftrightarrow {w}_{n-k,n}^{2(n-k-1)+\beta } \Bigg\}\nonumber\\
&\times&  \;_2F_1 \left(-q_0, q_0+\Omega; \gamma; \overleftrightarrow {w}_{1,n} \right) \Bigg\} z^n \Bigg\} \label{eq:100108}
\end{eqnarray}
where
 \begin{equation}
\begin{cases} \xi = \frac{(1-a)x}{x-a} \cr
z = -\frac{1}{1-a}\xi^2 \cr
\eta = \frac{2-a}{1-a} \xi \cr
\varphi = -\alpha +\beta +\gamma +(1-a)(\gamma +\delta -1) \cr
\Omega = \frac{\varphi }{(2-a)}\cr
q=\gamma \beta +(q_j+2j)\{\varphi +(2-a)(q_j+2j) \} \;\;\mbox{as}\;j,q_j\in \mathbb{N}_{0} \cr
q_i\leq q_j \;\;\mbox{only}\;\mbox{if}\;i\leq j\;\;\mbox{where}\;i,j\in \mathbb{N}_{0} 
\end{cases}\nonumber 
\end{equation}
\subsubsection{Infinite series}
Replacing coefficients $a$, $q$, $\alpha $, $x$, $c_0$ and $\lambda $ by $1-a$, $-q+\gamma \beta $, $-\alpha+\gamma +\delta $, $\frac{(1-a)x}{x-a}$, 1 and zero into (\ref{eq:10049}). Multiply $\left(1-\frac{x}{a} \right)^{-\beta }$ and the new (\ref{eq:10049}) together.
\begin{eqnarray}
 && \left(1-\frac{x}{a} \right)^{-\beta } y(\xi ) \nonumber\\
 &=& \left(1-\frac{x}{a} \right)^{-\beta } Hl\left(1-a, -q+\gamma \beta; -\alpha +\gamma +\delta, \beta, \gamma, \delta; \frac{(1-a)x}{x-a} \right) \nonumber\\
&=& \left(1-\frac{x}{a} \right)^{-\beta } \Bigg\{  \;_2F_1 \left(-\Pi_{0} ^{+}, -\Pi_{0} ^{-}; \gamma; \eta \right) \nonumber\\
&+& \sum_{n=1}^{\infty } \left\{\prod _{k=0}^{n-1} \Bigg\{ \int_{0}^{1} dt_{n-k}\;t_{n-k}^{2(n-k)-1} \int_{0}^{1} du_{n-k}\;u_{n-k}^{2(n-k-1)+\gamma } \right.\nonumber\\
&\times& \frac{1}{2\pi i}  \oint dv_{n-k} \frac{1}{v_{n-k}} \left( \frac{v_{n-k}-1}{v_{n-k}}\right)^{\Pi_{n-k} ^{+}} \left( 1- \overleftrightarrow {w}_{n-k+1,n}(1-t_{n-k})(1-u_{n-k})v_{n-k}\right)^{\Pi_{n-k} ^{-}}\nonumber\\
&\times& \overleftrightarrow {w}_{n-k,n}^{-(2(n-k-1)-\alpha +\gamma +\delta )}\left(  \overleftrightarrow {w}_{n-k,n} \partial _{ \overleftrightarrow {w}_{n-k,n}}\right) \overleftrightarrow {w}_{n-k,n}^{-\alpha -\beta+\gamma +\delta } \left(  \overleftrightarrow {w}_{n-k,n} \partial _{ \overleftrightarrow {w}_{n-k,n}}\right) \overleftrightarrow {w}_{n-k,n}^{2(n-k-1)+\beta } \Bigg\}\nonumber\\
&\times& \;_2F_1 \left( -\Pi_{0} ^{+}, -\Pi_{0} ^{-}; \gamma; \overleftrightarrow {w}_{1,n}\right) \Bigg\} z^n \Bigg\} \label{eq:100109}
\end{eqnarray}
where
 \begin{equation}
\begin{cases} 
\Pi_{0} ^{\pm}= \frac{-\varphi \pm\sqrt{\varphi ^2-4(2-a)(-q+\gamma \beta)}}{2(2-a)} \cr
\Pi_{n-k} ^{\pm}= \frac{-(\varphi +4(2-a)(n-k))\pm\sqrt{\varphi ^2-4(2-a)(-q+\gamma \beta)}}{2(2-a)} 
\end{cases}\nonumber 
\end{equation} 
\subsection{ \footnotesize ${\displaystyle (1-x)^{1-\delta }\left(1-\frac{x}{a} \right)^{-\beta+\delta -1} Hl\left(1-a, -q+\gamma [(\delta -1)a+\beta -\delta +1]; -\alpha +\gamma +1, \beta -\delta+1, \gamma, 2-\delta; \frac{(1-a)x}{x-a} \right)}$ \normalsize}
\subsubsection{Polynomial of type 2}
Replacing coefficients $a$, $q$, $\alpha $, $\beta $, $\delta $, $x$, $c_0$ and $\lambda $ by $1-a$, $-q+\gamma [(\delta -1)a+\beta -\delta +1]$, $-\alpha +\gamma +1$, $\beta -\delta+1$, $2-\delta $, $\frac{(1-a)x}{x-a}$, 1 and zero into (\ref{eq:10039}). Multiply $(1-x)^{1-\delta }\left(1-\frac{x}{a} \right)^{-\beta+\delta -1}$ and the new (\ref{eq:10039}) together.
\begin{eqnarray}
 && (1-x)^{1-\delta }\left(1-\frac{x}{a} \right)^{-\beta+\delta -1} y(\xi ) \nonumber\\
 &=& (1-x)^{1-\delta }\left(1-\frac{x}{a} \right)^{-\beta+\delta -1} Hl\Big(1-a, -q+\gamma [(\delta -1)a+\beta -\delta +1]; -\alpha +\gamma +1, \beta -\delta+1, \gamma \nonumber\\
&&, 2-\delta; (1-a)x(x-a)^{-1}\Big) \nonumber\\
&=& (1-x)^{1-\delta }\left(1-\frac{x}{a} \right)^{-\beta+\delta -1} \Bigg\{ \;_2F_1 \left(-q_0, q_0+\Omega; \gamma; \eta \right) \nonumber\\
&+& \sum_{n=1}^{\infty } \Bigg\{\prod _{k=0}^{n-1} \Bigg\{ \int_{0}^{1} dt_{n-k}\;t_{n-k}^{2(n-k)-1} \int_{0}^{1} du_{n-k}\;u_{n-k}^{2(n-k-1)+\gamma} \nonumber\\
&\times&  \frac{1}{2\pi i}  \oint dv_{n-k} \frac{1}{v_{n-k}} \left( \frac{v_{n-k}-1}{v_{n-k}} \frac{1}{1-\overleftrightarrow {w}_{n-k+1,n}(1-t_{n-k})(1-u_{n-k})v_{n-k}}\right)^{q_{n-k}} \nonumber\\
&\times& \left( 1- \overleftrightarrow {w}_{n-k+1,n}(1-t_{n-k})(1-u_{n-k})v_{n-k}\right)^{-\left(4(n-k)+\Omega\right)}\nonumber\\
&\times& \overleftrightarrow {w}_{n-k,n}^{-(2(n-k)-1-\alpha +\gamma )}\left(  \overleftrightarrow {w}_{n-k,n} \partial _{ \overleftrightarrow {w}_{n-k,n}}\right) \overleftrightarrow {w}_{n-k,n}^{-\alpha -\beta+\gamma +\delta } \left(  \overleftrightarrow {w}_{n-k,n} \partial _{ \overleftrightarrow {w}_{n-k,n}}\right) \overleftrightarrow {w}_{n-k,n}^{2(n-k)-1+\beta -\delta }  \Bigg\}\nonumber\\
&\times& \;_2F_1 \left(-q_0, q_0+\Omega; \gamma; \overleftrightarrow {w}_{1,n} \right)\Bigg\} z^n \Bigg\} \label{eq:100110}
\end{eqnarray}
where
 \begin{equation}
\begin{cases} 
\xi= \frac{(1-a)x}{x-a} \cr
z = -\frac{1}{1-a}\xi^2 \cr
\eta = \frac{(2-a)}{1-a} \xi \cr
\varphi = -\alpha +\beta +\gamma +(1-a)(\gamma -\delta +1) \cr
\Omega = \frac{\varphi }{(2-a)}\cr
q= \gamma [(\delta -1)a+\beta -\delta +1]+(q_j+2j)\{\varphi +(2-a)(q_j+2j) \} \;\;\mbox{as}\;j,q_j\in \mathbb{N}_{0} \cr
q_i\leq q_j \;\;\mbox{only}\;\mbox{if}\;i\leq j\;\;\mbox{where}\;i,j\in \mathbb{N}_{0} 
\end{cases}\nonumber 
\end{equation}
\subsubsection{Infinite series}
Replacing coefficients $a$, $q$, $\alpha $, $\beta $, $\delta $, $x$, $c_0$ and $\lambda $ by $1-a$, $q^{\ast }= -q+\gamma [(\delta -1)a+\beta -\delta +1]$, $-\alpha +\gamma +1$, $\beta -\delta+1$, $2-\delta $, $\frac{(1-a)x}{x-a}$, 1 and zero into (\ref{eq:10049}). Multiply $(1-x)^{1-\delta }\left(1-\frac{x}{a} \right)^{-\beta+\delta -1}$ and the new (\ref{eq:10049}) together.
\begin{eqnarray}
 && (1-x)^{1-\delta }\left(1-\frac{x}{a} \right)^{-\beta+\delta -1} y(\xi ) \nonumber\\
 &=& (1-x)^{1-\delta }\left(1-\frac{x}{a} \right)^{-\beta+\delta -1} Hl\Big(1-a, q^{\ast }= -q+\gamma [(\delta -1)a+\beta -\delta +1]; -\alpha +\gamma +1 \nonumber\\
&&, \beta -\delta+1, \gamma, 2-\delta; (1-a)x(x-a)^{-1} \Big) \nonumber\\
&=& (1-x)^{1-\delta }\left(1-\frac{x}{a} \right)^{-\beta+\delta -1} \Bigg\{  \;_2F_1 \left(-\Pi_0 ^{+}, -\Pi_0 ^{-}; \gamma; \eta \right) \nonumber\\
&+& \sum_{n=1}^{\infty } \left\{\prod _{k=0}^{n-1} \Bigg\{ \int_{0}^{1} dt_{n-k}\;t_{n-k}^{2(n-k)-1} \int_{0}^{1} du_{n-k}\;u_{n-k}^{2(n-k-1)+\gamma} \right. \nonumber\\
&\times& \frac{1}{2\pi i}  \oint dv_{n-k} \frac{1}{v_{n-k}} \left( \frac{v_{n-k}-1}{v_{n-k}}\right)^{\Pi_{n-k} ^{+}}  \left( 1- \overleftrightarrow {w}_{n-k+1,n}(1-t_{n-k})(1-u_{n-k})v_{n-k}\right)^{\Pi_{n-k} ^{-}}\nonumber\\
&\times&  \overleftrightarrow {w}_{n-k,n}^{-(2(n-k)-1-\alpha +\gamma )}\left(  \overleftrightarrow {w}_{n-k,n} \partial _{ \overleftrightarrow {w}_{n-k,n}}\right) \overleftrightarrow {w}_{n-k,n}^{-\alpha -\beta+\gamma +\delta } \left(  \overleftrightarrow {w}_{n-k,n} \partial _{ \overleftrightarrow {w}_{n-k,n}}\right) \overleftrightarrow {w}_{n-k,n}^{2(n-k)-1+\beta-\delta } \Bigg\}\nonumber\\
&\times& \;_2F_1 \left(-\Pi_0 ^{+}, -\Pi_0 ^{-}; \gamma; \overleftrightarrow {w}_{1,n} \right) \Bigg\} z^n \Bigg\} \label{eq:100111}
\end{eqnarray}
where
 \begin{equation}
\begin{cases} 
\Pi_0 ^{\pm}= \frac{-\varphi \pm\sqrt{\varphi ^2-4(2-a)(-q+\gamma [(\delta -1)a+\beta -\delta +1])}}{2(2-a)} \cr
\Pi_{n-k} ^{\pm}= \frac{-(\varphi +4(2-a)(n-k))\pm\sqrt{\varphi ^2-4(2-a)(-q+\gamma [(\delta -1)a+\beta -\delta +1])}}{2(2-a)} 
\end{cases}\nonumber 
\end{equation} 
\subsection{ \footnotesize ${\displaystyle x^{-\alpha } Hl\left(\frac{a-1}{a}, \frac{-q+\alpha (\delta a+\beta -\delta )}{a}; \alpha, \alpha -\gamma +1, \delta , \alpha -\beta +1; \frac{x-1}{x} \right)}$\normalsize}
\subsubsection{Polynomial of type 2}
Replacing coefficients $a$, $q$, $\beta $, $\gamma $, $\delta $, $x$, $c_0$ and $\lambda $ by $\frac{a-1}{a}$, $\frac{-q+\alpha (\delta a+\beta -\delta )}{a}$, $\alpha -\gamma +1$, $\delta $, $\alpha -\beta +1$, $\frac{x-1}{x}$, 1 and zero into (\ref{eq:10039}). Multiply $x^{-\alpha }$ and the new (\ref{eq:10039}) together.
\begin{eqnarray}
&& x^{-\alpha } y(\xi ) \nonumber\\
 &=& x^{-\alpha } Hl\left(\frac{a-1}{a}, \frac{-q+\alpha (\delta a+\beta -\delta )}{a}; \alpha, \alpha -\gamma +1, \delta , \alpha -\beta +1; \frac{x-1}{x} \right) \nonumber\\
&=& x^{-\alpha } \Bigg\{\;_2F_1 \left(-q_0, q_0+\Omega; \delta; \eta \right) + \sum_{n=1}^{\infty } \Bigg\{\prod _{k=0}^{n-1} \Bigg\{ \int_{0}^{1} dt_{n-k}\;t_{n-k}^{2(n-k)-1} \int_{0}^{1} du_{n-k}\;u_{n-k}^{2(n-k-1)+\delta } \nonumber\\
&\times&  \frac{1}{2\pi i}  \oint dv_{n-k} \frac{1}{v_{n-k}} \left( \frac{v_{n-k}-1}{v_{n-k}} \frac{1}{1-\overleftrightarrow {w}_{n-k+1,n}(1-t_{n-k})(1-u_{n-k})v_{n-k}}\right)^{q_{n-k}} \nonumber\\
&\times& \left( 1- \overleftrightarrow {w}_{n-k+1,n}(1-t_{n-k})(1-u_{n-k})v_{n-k}\right)^{-\left(4(n-k)+\Omega\right)}\nonumber\\
&\times&  \overleftrightarrow {w}_{n-k,n}^{-(2(n-k-1)+\alpha )}\left(  \overleftrightarrow {w}_{n-k,n} \partial _{ \overleftrightarrow {w}_{n-k,n}}\right) \overleftrightarrow {w}_{n-k,n}^{\gamma -1} \left(  \overleftrightarrow {w}_{n-k,n} \partial _{ \overleftrightarrow {w}_{n-k,n}}\right) \overleftrightarrow {w}_{n-k,n}^{2(n-k)-1+\alpha -\gamma } \Bigg\}\nonumber\\
&\times& \;_2F_1 \left( -q_0, q_0+\Omega; \delta; \overleftrightarrow {w}_{1,n}\right) \Bigg\} z^n \Bigg\} \label{eq:100112}
\end{eqnarray}
where
 \begin{equation}
\begin{cases} 
\xi= \frac{x-1}{x} \cr
z = \frac{-a}{a-1}\xi^2 \cr
\eta = \frac{2a-1}{a-1} \xi \cr
\varphi = \alpha +\beta -\gamma +\frac{a-1}{a}(\alpha -\beta +\delta ) \cr
\Omega = \frac{a\varphi }{(2a-1)}\cr
q= \alpha (\delta a+\beta -\delta )+(q_j+2j)\{a\varphi +(2a-1)(q_j+2j) \} \;\;\mbox{as}\;j,q_j\in \mathbb{N}_{0} \cr
q_i\leq q_j \;\;\mbox{only}\;\mbox{if}\;i\leq j\;\;\mbox{where}\;i,j\in \mathbb{N}_{0} 
\end{cases}\nonumber 
\end{equation}
 \subsubsection{Infinite series}
Replacing coefficients $a$, $q$, $\beta $, $\gamma $, $\delta $, $x$, $c_0$ and $\lambda $ by $\frac{a-1}{a}$, $ \frac{-q+\alpha (\delta a+\beta -\delta )}{a}$, $\alpha -\gamma +1$, $\delta $, $\alpha -\beta +1$, $\frac{x-1}{x}$, 1 and zero into (\ref{eq:10049}). Multiply $x^{-\alpha }$ and the new (\ref{eq:10049}) together.
\begin{eqnarray}
&& x^{-\alpha } y(\xi ) \nonumber\\
 &=& x^{-\alpha } Hl\left(\frac{a-1}{a}, \frac{-q+\alpha (\delta a+\beta -\delta )}{a}; \alpha, \alpha -\gamma +1, \delta , \alpha -\beta +1; \frac{x-1}{x} \right) \nonumber\\
&=& x^{-\alpha } \Bigg\{ \;_2F_1 \left(-\Pi_0 ^{+}, -\Pi_0 ^{-}; \delta; \eta \right) \nonumber\\
&+& \sum_{n=1}^{\infty } \left\{\prod _{k=0}^{n-1} \Bigg\{ \int_{0}^{1} dt_{n-k}\;t_{n-k}^{2(n-k)-1} \int_{0}^{1} du_{n-k}\;u_{n-k}^{2(n-k-1)+\delta } \right. \nonumber\\
&\times& \frac{1}{2\pi i}  \oint dv_{n-k} \frac{1}{v_{n-k}} \left( \frac{v_{n-k}-1}{v_{n-k}}\right)^{\Pi_{n-k} ^{+}} \left( 1- \overleftrightarrow {w}_{n-k+1,n}(1-t_{n-k})(1-u_{n-k})v_{n-k}\right)^{\Pi_{n-k} ^{-}}\nonumber\\
&\times& \overleftrightarrow {w}_{n-k,n}^{-(2(n-k-1)+\alpha )}\left(  \overleftrightarrow {w}_{n-k,n} \partial _{ \overleftrightarrow {w}_{n-k,n}}\right) \overleftrightarrow {w}_{n-k,n}^{\gamma -1} \left(  \overleftrightarrow {w}_{n-k,n} \partial _{ \overleftrightarrow {w}_{n-k,n}}\right) \overleftrightarrow {w}_{n-k,n}^{2(n-k)-1+\alpha -\gamma } \Bigg\}\nonumber\\
&\times& \;_2F_1 \left(-\Pi_0 ^{+}, -\Pi_0 ^{-}; \delta; \overleftrightarrow {w}_{1,n} \right) \Bigg\} z^n \Bigg\} \label{eq:100113}
\end{eqnarray}
where
 \begin{equation}
\begin{cases} 
\Pi_0 ^{\pm}= \frac{-a\varphi \;\pm \;a\sqrt{\varphi ^2-4\frac{(2a-1)}{a}\left(\frac{-q+\alpha (\delta a+\beta -\delta )}{a}\right)}}{2(2a-1)} \cr
\Pi_{n-k} ^{\pm}= \frac{-a\left\{\varphi +4\frac{(2a-1)}{a}(n-k)\right\}\;\pm \;a\sqrt{\varphi ^2-4\frac{(2a-1)}{a}\left(\frac{-q+\alpha (\delta a+\beta -\delta )}{a}\right)}}{2(2a-1)} 
\end{cases}\nonumber 
\end{equation}
\subsection{ ${\displaystyle \left(\frac{x-a}{1-a} \right)^{-\alpha } Hl\left(a, q-(\beta -\delta )\alpha ; \alpha , -\beta+\gamma +\delta , \delta , \gamma; \frac{a(x-1)}{x-a} \right)}$}
\subsubsection{Polynomial of type 2}
Replacing coefficients $q$, $\beta $, $\gamma $, $\delta $, $x$, $c_0$ and $\lambda $ by $q-(\beta -\delta )\alpha $, $-\beta+\gamma +\delta $, $\delta $,  $\gamma $, $\frac{a(x-1)}{x-a}$, 1 and zero into (\ref{eq:10039}). Multiply $\left(\frac{x-a}{1-a} \right)^{-\alpha }$ and the new (\ref{eq:10039}) together.
\begin{eqnarray}
&& \left(\frac{x-a}{1-a} \right)^{-\alpha } y(\xi ) \nonumber\\
 &=& \left(\frac{x-a}{1-a} \right)^{-\alpha } Hl\left(a, q-(\beta -\delta )\alpha ; \alpha , -\beta+\gamma +\delta , \delta , \gamma; \frac{a(x-1)}{x-a} \right) \nonumber\\
&=& \left(\frac{x-a}{1-a} \right)^{-\alpha } \Bigg\{ \;_2F_1 \left(-q_0, q_0+\Omega; \delta; \eta\right)\nonumber\\
&+& \sum_{n=1}^{\infty } \Bigg\{\prod _{k=0}^{n-1} \Bigg\{ \int_{0}^{1} dt_{n-k}\;t_{n-k}^{2(n-k)-1} \int_{0}^{1} du_{n-k}\;u_{n-k}^{2(n-k-1)+\delta } \nonumber\\
&\times&  \frac{1}{2\pi i}  \oint dv_{n-k} \frac{1}{v_{n-k}} \left( \frac{v_{n-k}-1}{v_{n-k}} \frac{1}{1-\overleftrightarrow {w}_{n-k+1,n}(1-t_{n-k})(1-u_{n-k})v_{n-k}}\right)^{q_{n-k}} \nonumber\\
&\times& \left( 1- \overleftrightarrow {w}_{n-k+1,n}(1-t_{n-k})(1-u_{n-k})v_{n-k}\right)^{-\left(4(n-k)+\Omega \right)}\nonumber\\
&\times& \overleftrightarrow {w}_{n-k,n}^{-(2(n-k-1)+\alpha )}\left(  \overleftrightarrow {w}_{n-k,n} \partial _{ \overleftrightarrow {w}_{n-k,n}}\right) \overleftrightarrow {w}_{n-k,n}^{\alpha+\beta-\gamma -\delta } \left(  \overleftrightarrow {w}_{n-k,n} \partial _{ \overleftrightarrow {w}_{n-k,n}}\right) \overleftrightarrow {w}_{n-k,n}^{2(n-k-1)-\beta +\gamma +\delta } \Bigg\}\nonumber\\
&\times& \;_2F_1 \left(-q_0, q_0+\Omega; \delta; \overleftrightarrow {w}_{1,n}\right)\Bigg\} z^n \Bigg\} \label{eq:100114}
\end{eqnarray}
where
 \begin{equation}
\begin{cases} 
\xi= \frac{a(x-1)}{x-a} \cr
z = -\frac{1}{a}\xi^2 \cr
\eta = \frac{(1+a)}{a} \xi \cr
\varphi = \alpha -\beta +\delta +a(\gamma +\delta -1) \cr
\Omega = \frac{\varphi}{(1+a)}\cr
q= (\beta -\delta )\alpha -(q_j+2j)\{\varphi +(1+a)(q_j+2j ) \} \;\;\mbox{as}\;j,q_j\in \mathbb{N}_{0} \cr
q_i\leq q_j \;\;\mbox{only}\;\mbox{if}\;i\leq j\;\;\mbox{where}\;i,j\in \mathbb{N}_{0} 
\end{cases}\nonumber 
\end{equation}
 \subsubsection{Infinite series}
Replacing coefficients $q$, $\beta $, $\gamma $, $\delta $, $x$, $c_0$ and $\lambda $ by $q-(\beta -\delta )\alpha $, $-\beta+\gamma +\delta $, $\delta $,  $\gamma $, $\frac{a(x-1)}{x-a}$, 1 and zero into (\ref{eq:10049}). Multiply $\left(\frac{x-a}{1-a} \right)^{-\alpha }$ and the new (\ref{eq:10049}) together.
\begin{eqnarray}
&& \left(\frac{x-a}{1-a} \right)^{-\alpha } y(\xi ) \nonumber\\
 &=& \left(\frac{x-a}{1-a} \right)^{-\alpha } Hl\left(a, q-(\beta -\delta )\alpha ; \alpha , -\beta+\gamma +\delta , \delta , \gamma; \frac{a(x-1)}{x-a} \right) \nonumber\\
&=& \left(\frac{x-a}{1-a} \right)^{-\alpha } \Bigg\{ \;_2F_1 \left(-\Pi_0 ^{+}, -\Pi_0 ^{-}; \delta; \eta\right) \nonumber\\
&+& \sum_{n=1}^{\infty } \Bigg\{\prod _{k=0}^{n-1} \Bigg\{ \int_{0}^{1} dt_{n-k}\;t_{n-k}^{2(n-k)-1 } \int_{0}^{1} du_{n-k}\;u_{n-k}^{2(n-k-1)+\delta} \nonumber\\
&\times& \frac{1}{2\pi i}  \oint dv_{n-k} \frac{1}{v_{n-k}} \left( \frac{v_{n-k}-1}{v_{n-k}}\right)^{\Pi_{n-k} ^{+}} \left( 1- \overleftrightarrow {w}_{n-k+1,n}(1-t_{n-k})(1-u_{n-k})v_{n-k}\right)^{\Pi_{n-k} ^{-}}\nonumber\\
&\times&  \overleftrightarrow {w}_{n-k,n}^{-(2(n-k-1)+\alpha )}\left(  \overleftrightarrow {w}_{n-k,n} \partial _{ \overleftrightarrow {w}_{n-k,n}}\right) \overleftrightarrow {w}_{n-k,n}^{\alpha +\beta -\gamma -\delta } \left(  \overleftrightarrow {w}_{n-k,n} \partial _{ \overleftrightarrow {w}_{n-k,n}}\right) \overleftrightarrow {w}_{n-k,n}^{2(n-k-1)-\beta +\gamma +\delta } \Bigg\}\nonumber\\
&\times& \;_2F_1 \left(-\Pi_0 ^{+}, -\Pi_0 ^{-}; \delta; \overleftrightarrow {w}_{1,n}\right) \Bigg\} z^n \Bigg\}  \label{eq:100115}
\end{eqnarray}
where  
\begin{equation}
\begin{cases} 
\Pi_0 ^{\pm}= \frac{-\varphi \pm\sqrt{\varphi ^2-4(1+a)(q-(\beta -\delta )\alpha)}}{2(1+a)} \cr
\Pi_{n-k} ^{\pm}= \frac{-(\varphi +4(1+a)(n-k))\pm\sqrt{\varphi ^2-4(1+a)(q-(\beta -\delta )\alpha)}}{2(1+a)} 
\end{cases}\nonumber 
\end{equation}
  
\end{appendices} 

\addcontentsline{toc}{section}{Bibliography}
\bibliographystyle{model1a-num-names}
\bibliography{<your-bib-database>}
 
\chapter{The generating function for Heun polynomial using reversible three-term recurrence formula}
\chaptermark{Generating function for Heun polynomial}

In this chapter, by applying the generating function for Jacobi polynomial using hypergeometric functions into the integral representation of Heun polynomial of type 2\footnote{polynomial of type 2 is a polynomial which makes $A_n$ term terminated in three term recursion relation of the power series in a linear differential equation.}, I consider the generating function of Heun polynomial including all higher terms of $B_n$'s\footnote{`` higher terms of $B_n$'s'' means at least two terms of $B_n$'s.}.

Nine examples of 192 local solutions of the Heun equation (Maier, 2007) are provided in the appendices.  For each example, I construct the generating function for Heun polynomial of type 2 including all higher terms of $B_n$'s.
\section{Introduction}

In 1837, Gabriel Lam\'{e} introduced second ordinary differential equation which has four regular singular points in the method of separation of variables applied to the Laplace equation in elliptic coordinates\cite{Lame18372}. Various authors has called this equation as `Lam\'{e} equation' or `ellipsoidal harmonic equation'\cite{Erde1955}.
 
In 1889, K. Heun worked on the second ordinary differential equation which has four regular singular points. Later, the Heun function became a general function of all well-known special functions: Mathieu, Lame and Coulomb spheroidal functions.

Lam\'{e} ordinary differential equation represented either in the algebraic form or in Weierstrass's form is the special case of Heun's differential equation: Lam\'{e} equation is derived from Heun equation by changing all coefficients and an independent variable.
The solutions of Lam\'{e} and Heun equations are a form of power series that can be expressed as three term recurrence. In contrast, most of well-known special functions consist of two term recursion relation (Hypergeometric, Bessel, Legendre, Kummer functions, etc). Lam\'{e} and Heun functions arise from deriving the Laplace equation in general Jacobi ellipsoidal or conical coordinates.\cite{Kalnin20052} Also Heun polynomial arises in the separable coordinate systems on the n-sphere.\cite{Kaln19912}  

According to E. G. Kalnins and W. Miller Jr.(1990 \cite{Kalnin20052}) , ``Lam\'{e} \cite{Lame18372} and Heun functions have received relatively little attention, since they are rather intractable. Unfortunately the beautiful identities appearing have received little notice, probably because the methods of proof seemed obscure.'' 

Because Lam\'{e} and Heun equations have three term recursion relation in power series representation: a 3-term recursive relation between successive coefficients in its power series creates the complicated mathematical calculations. In general, most problems of physical and mechanical system in nature are nonlinear. For the simplified purposes and better approximation of physical systems, mathematically we have used linear (ordinary) differential equations to linearize these nonlinear system. Traditionally, we have explained these linear systems by only using two term recurrence relation of the power series in linear ordinary differential equations until 19th century. However, since modern physics (quantum gravity, SUSY, general relativity, etc) comes out the world, it seems to require at least three or four term recurrence relations in power series expansions. Furthermore these type of problems can not be reduced to two term recurrence relations.\cite{Hortacsu:2011rr2}

In 1939 Svartholm showed how Heun equation can be obtained in the form of series of degenerate hypergeometric functions. \cite{Svar1939} A few years later, His solution was developed by Erd\'{e}lyi in three papers (1942, 1944),  Sleeman (1969), Schmidt and Wolf (1979), Kalnins and Miller (1991).\cite{Kaln19912,Erde19421,Erde19422,Erde1944,Slee19691,Slee19692,Schm1979} 

In Ref.\cite{Chou2012H12,Chou2012H22}, by applying 3TRF, I show power series expansion in closed forms of Heun equation and its representation of the form of integrals for infinite series and polynomial of type 1.\footnote{polynomial of type 1 is a polynomial which makes $B_n$ term terminated in three term recursion relation of the power series in a linear differential equation.} I show that a $_2F_1$ function recurs in each of sub-integral forms of Heun function: the first sub-integral form contains zero term of $A_n's$, the second one contains one term of $A_n$'s, the third one contains two terms of $A_n$'s, etc. 

In chapter 2, by applying R3TRF, I show power series expansion in closed forms of Heun function and its representation of the form of integrals for infinite series and polynomial of type 2. I also show that a $_2F_1$ function recurs in each of sub-integral forms of Heun function: the first sub-integral form contains zero term of $B_n's$, the second one contains one term of $B_n$'s, the third one contains two terms of $B_n$'s, etc. 

In mathematics, a generating function $F(x,t)$ is a formal (it need not converge) power series in one indeterminate $t^n$, whose coefficients encode information about a sequence of numbers $f_n(x)$ that is indexed by the natural numbers (in quantum mechanically, the index $n$ corresponds to the eigenvalues denoted by natural numbers). In the case of Legendre equation, $f_n(x)$ corresponds to the Legendre polynomials. Since the expansion of $F(x,t)$ has generated the set $f_n(x)$ in $t$, $F(x,t)$ can be described as an analytic function, combined coefficients $x$ and $t$, in closed form without index $n$. In general, a generating function in closed form means that it is expressed in a direct form without summation signs. For instance, the generating function for Legendre polynomials is given by
\begin{eqnarray}
\sum_{n=0}^{\infty } P_n(x)\;t^n &=& \frac{1}{\sqrt{1-2xt+t^2}} \label{L1}\\
&=& 1+xt+ \frac{1}{2}(3x^2-1)t^2 + \frac{1}{2}(5x^3-3x)t^3 + \frac{1}{8}(35x^4-30x^2+3)t^4\nonumber\\
&& +\frac{1}{8}(63x^5-70x^3+15x)t^5+ +\frac{1}{16}(231x^6-315x^4+105x^2-5)t^6\cdots \hspace{2cm}\label{L2}
\end{eqnarray}
The power series representation (\ref{L2}) is not in closed form while an analytic function (\ref{L1}) is.     

In Ref.\cite{Chou2012h}, by apply a generating function for Jacobi polynomial using hypergeometric functions into the general representation in the form of integral of Lame polynomial of type 1 in Weierstrass's form, I construct the generating function for Lame polynomial analytically. As we see generating functions for Lame polynomial of type 1 of the first and second kinds, it is described as an analytic function in closed forms with removed index $\alpha _j$ 
 where $j,\alpha _j \in \mathbb{N}_{0}$. $\alpha _j$ is the natural numbers (eigenvalues) which makes $B_n$ term terminated at specific values of index $n$ in sequence $c_n$: the sequence $c_n$ consists of combinations $A_n$ and $B_n$ terms in three term recurrence relation of the power series representation of Lame equation. 

In this chapter, by using same technique for the generating function for Lame polynomial of type 1, I consider the generating function for Heun polynomial of type 2 including all higher terms of $B_n$'s by applying a generating function for Jacobi polynomial using hypergeometric functions into the general integral form of Heun polynomial of type 2. Since the generating function for Heun polynomial is derived, we might be possible to construct orthogonal relations of Heun polynomial.
In the physical point of view we might be possible to obtain the normalized constant for the wave function in modern physics, recursion relation and its expectation value of any physical quantities from the generating function for Heun polynomial. For the case of hydrogen-like atoms, the normalized wave function is derived from the generating function for associated Laguerre polynomials. And the expectation value of physical quantities such as position and momentum is constructed by applying the recursive relation of associated Laguerre polynomials.

In future papers I will derive Heun polynomial of type 3 and its generating function that has fixed integer values of $\alpha $ and/or $\beta $, just as it has a fixed value of $q$. In this paper I construct the generating function for Heun polynomial of type 2: I treat $\alpha $, $\beta $, $\gamma $ and $\delta $ as free variables and the accessory parameter $q$ as a fixed value.

Heun's equation is a second-order linear ordinary differential equation of the form \cite{Heun18892}.
\begin{equation}
\frac{d^2{y}}{d{x}^2} + \left(\frac{\gamma }{x} +\frac{\delta }{x-1} + \frac{\epsilon }{x-a}\right) \frac{d{y}}{d{x}} +  \frac{\alpha \beta x-q}{x(x-1)(x-a)} y = 0 \label{eq:2001}
\end{equation}
With the condition $\epsilon = \alpha +\beta -\gamma -\delta +1$. The parameters play different roles: $a \ne 0 $ is the singularity parameter, $\alpha $, $\beta $, $\gamma $, $\delta $, $\epsilon $ are exponent parameters, $q$ is the accessory parameter which in many physical applications appears as a spectral parameter. Also, $\alpha $ and $\beta $ are identical to each other. The total number of free parameters is six. It has four regular singular points which are 0, 1, a and $\infty $ with exponents $\{ 0, 1-\gamma \}$, $\{ 0, 1-\delta \}$, $\{ 0, 1-\epsilon \}$ and $\{ \alpha, \beta \}$.

\section[Generating function for Heun polynomial of type 2]{Generating function for Heun polynomial of type 2
  \sectionmark{Generating function for Heun polynomial of type 2}}
  \sectionmark{Generating function for Heun polynomial of type 2}
Let's investigate the generating function for Heun polynomial of type 2 of the first and second kinds. 
\begin{lemma}
The generating function for Jacobi polynomial using hypergeometric functions is defined by
\begin{eqnarray}
&&\sum_{q_0=0}^{\infty }\frac{(\gamma )_{q_0}}{q_0!} w^{q_0} \;_2F_1(-q_0, q_0+A; \gamma; x) \nonumber\\
&&= 2^{A -1}\frac{\left(1-w+\sqrt{w^2-2(1-2x)w+1}\right)^{1-\gamma } \left(1+w+\sqrt{w^2-2(1-2x)w+1}\right)^{\gamma -A}}{\sqrt{w^2-2(1-2x)w+1}} \nonumber\\
&&\;\;\mbox{where}\;|w|<1 \label{eq:2003}
\end{eqnarray}
\end{lemma}
\begin{proof}
Jacobi polynomial $P_n^{(\alpha, \beta )}(x)$ can be written in terms of hypergeometric function using
\begin{equation}
_2F_1 (-n, n+\alpha +\beta +1; \alpha +1;x) = \frac{n!}{(\alpha +1)_n}P_n^{(\alpha, \beta )}(1-2x) \label{eq:20020}
\end{equation}
And
\begin{equation}
P_n^{(\alpha, \beta )}(x) = \frac{\Gamma (n+\alpha +1)}{n! \Gamma (n+\alpha +\beta +1)} \sum_{m=0}^{n} \binom{n}{m} \frac{\Gamma (n+m+\alpha +\beta +1)}{\Gamma (m+\alpha +1)}\left( \frac{x-1}{2}\right)^m \label{eq:20021}
\end{equation}
The generating function for Jacobi polynomials is given by 
\begin{equation}
\sum_{n=0}^{\infty } P_n^{(\alpha, \beta )}(x) w^n = 2^{\alpha +\beta }\frac{\left(1-w+\sqrt{w^2-2xw+1}\right)^{-\alpha  } \left(1+w+\sqrt{w^2-2xw+1}\right)^{-\beta }}{\sqrt{w^2-2xw+1}} \label{eq:20022}
\end{equation}
Replace n, $\alpha $ and $\beta $ by $q_0$, $\gamma -1$ and $A-\gamma $ in (\ref{eq:20020}), and acting the summation operator ${\displaystyle \sum_{q_0=0}^{\infty }\frac{(\gamma )_{q_0}}{q_0!} w^{q_0} }$ on the new (\ref{eq:20020})
\begin{equation}
\sum_{q_0=0}^{\infty }\frac{(\gamma )_{q_0}}{q_0!} w^{q_0}\; _2F_1 (-q_0, q_0+ A; \gamma ;x) = \sum_{q_0=0}^{\infty } P_n^{(\gamma -1, A-\gamma  )}(1-2x) w^{q_0}\label{eq:20023}
\end{equation}
Replace $\alpha $, $\beta $ and $x$ by $\gamma -1$, $A-\gamma $ and $1-2x$ in (\ref{eq:20022}). As we take the new (\ref{eq:20022}) into (\ref{eq:20023}), we obtain (\ref{eq:2003}). 
\qed 
\end{proof}
\begin{definition}
I define that
\begin{equation}
\begin{cases}
\displaystyle { s_{a,b}} = \begin{cases} \displaystyle {  s_a\cdot s_{a+1}\cdot s_{a+2}\cdots s_{b-2}\cdot s_{b-1}\cdot s_b}\;\;\mbox{if}\;a>b \cr
s_a \;\;\mbox{if}\;a=b\end{cases}
\cr
\cr
\displaystyle { \widetilde{w}_{i,j}}  = 
\begin{cases} \displaystyle { \frac{ \widetilde{w}_{i+1,j}\; t_i u_i \left\{ 1+ (s_i+2\widetilde{w}_{i+1,j}(1-t_i)(1-u_i))s_i\right\}}{2(1-\widetilde{w}_{i+1,j}(1-t_i)(1-u_i))^2 s_i}} \cr
\displaystyle {-\frac{\widetilde{w}_{i+1,j}\; t_i u_i (1+s_i)\sqrt{s_i^2-2(1-2\widetilde{w}_{i+1,j}(1-t_i)(1-u_i))s_i+1}}{2(1-\widetilde{w}_{i+1,j}(1-t_i)(1-u_i))^2 s_i}} \;\;\mbox{where}\;i<j \cr
\cr
\displaystyle { \frac{\eta t_i u_i \left\{ 1+ (s_{i,\infty }+2\eta(1-t_i)(1-u_i))s_{i,\infty }\right\}}{2(1-\eta (1-t_i)(1-u_i))^2 s_{i,\infty }}} \cr
\displaystyle {-\frac{\eta t_i u_i(1+s_{i,\infty })\sqrt{s_{i,\infty }^2-2(1-2\eta (1-t_i)(1-u_i))s_{i,\infty }+1}}{2(1-\eta (1-t_i)(1-u_i))^2 s_{i,\infty }}} \;\;\mbox{where}\;i=j 
\end{cases}
\end{cases}\label{eq:2004}
\end{equation}
where
\begin{equation}
a,b,i,j\in \mathbb{N}_{0} \nonumber
\end{equation}
\end{definition}
And we have
\begin{equation}
\sum_{q_i = q_j}^{\infty } r_i^{q_i} = \frac{r_i^{q_j}}{(1-r_i)}\label{eq:2005}
\end{equation}
\begin{theorem}
In chapter 2, the representation in the form of integral of Heun polynomial of type 2 is given by\footnote{$y_1(x)$ means the sub-integral form in (\ref{eq:2006}) contains one term of $B_n's$, $y_2(x)$ means the sub-integral form in (\ref{eq:2006}) contains two terms of $B_n's$, $y_3(x)$ means the sub-integral form in (\ref{eq:2006}) contains three terms of $B_n's$, etc.}
\begin{eqnarray}
 y(x)&=& \sum_{n=0}^{\infty } y_{n}(x)= y_0(x)+ y_1(x)+ y_2(x)+y_3(x)+\cdots \nonumber\\
&=& c_0 x^{\lambda } \Bigg\{ \sum_{i_0=0}^{q_0}\frac{(-q_0)_{i_0}\left(q_0+\frac{\varphi +2(1+a)\lambda }{(1+a)}\right)_{i_0}}{(1+\lambda )_{i_0}(\gamma +\lambda )_{i_0}}  \eta ^{i_0}\nonumber\\
&&+ \sum_{n=1}^{\infty } \Bigg\{\prod _{k=0}^{n-1} \Bigg\{ \int_{0}^{1} dt_{n-k}\;t_{n-k}^{2(n-k)-1+\lambda } \int_{0}^{1} du_{n-k}\;u_{n-k}^{2(n-k-1)+\gamma +\lambda } \nonumber\\
&&\times  \frac{1}{2\pi i}  \oint dv_{n-k} \frac{1}{v_{n-k}} \left( \frac{v_{n-k}-1}{v_{n-k}} \frac{1}{1-\overleftrightarrow {w}_{n-k+1,n}(1-t_{n-k})(1-u_{n-k})v_{n-k}}\right)^{q_{n-k}} \nonumber\\
&&\times \left( 1- \overleftrightarrow {w}_{n-k+1,n}(1-t_{n-k})(1-u_{n-k})v_{n-k}\right)^{-\left(4(n-k)+\frac{\varphi +2(1+a)\lambda }{(1+a)}\right)}\nonumber\\
&&\times \overleftrightarrow {w}_{n-k,n}^{-(2(n-k-1)+\alpha +\lambda )}\left(  \overleftrightarrow {w}_{n-k,n} \partial _{ \overleftrightarrow {w}_{n-k,n}}\right) \overleftrightarrow {w}_{n-k,n}^{\alpha -\beta} \left(  \overleftrightarrow {w}_{n-k,n} \partial _{ \overleftrightarrow {w}_{n-k,n}}\right) \overleftrightarrow {w}_{n-k,n}^{2(n-k-1)+\beta +\lambda } \Bigg\}\nonumber\\
&&\times \sum_{i_0=0}^{q_0}\frac{(-q_0)_{i_0}\left(q_0+\frac{\varphi +2(1+a)\lambda }{(1+a)}\right)_{i_0}}{(1+\lambda )_{i_0}(\gamma +\lambda )_{i_0}} \overleftrightarrow {w}_{1,n}^{i_0}\Bigg\} z^n \Bigg\} \label{eq:2006}
\end{eqnarray}
where
\begin{equation}
\begin{cases} z = -\frac{1}{a}x^2 \cr
\eta = \frac{(1+a)}{a} x \cr
\varphi =\alpha +\beta -\delta +a(\delta +\gamma -1) \cr
q= -(q_j+2j+\lambda )\{\varphi +(1+a)(q_j+2j+\lambda ) \} \;\;\mbox{as}\;j,q_j\in \mathbb{N}_{0} \cr
q_i\leq q_j \;\;\mbox{only}\;\mbox{if}\;i\leq j\;\;\mbox{where}\;i,j\in \mathbb{N}_{0} 
\end{cases}\nonumber 
\end{equation}
and
\begin{equation}\overleftrightarrow {w}_{i,j}=
\begin{cases} \displaystyle {\frac{v_i}{(v_i-1)}\; \frac{\overleftrightarrow w_{i+1,j} t_i u_i}{1- \overleftrightarrow w_{i+1,j} v_i (1-t_i)(1-u_i)}} \;\;\mbox{where}\; i\leq j\cr
\eta \;\;\mbox{only}\;\mbox{if}\; i>j
\end{cases}\nonumber 
\end{equation}
In the above, the first sub-integral form contains one term of $B_n's$, the second one contains two terms of $B_n$'s, the third one contains three terms of $B_n$'s, etc.
\end{theorem}
Acting the summation operator $\displaystyle{ \sum_{q_0 =0}^{\infty } \frac{(\gamma')_{q_0}}{q_0!} s_0^{q_0} \prod _{n=1}^{\infty } \left\{ \sum_{ q_n = q_{n-1}}^{\infty } s_n^{q_n }\right\}}$ on (\ref{eq:2006}) where $|s_i|<1$ as $i=0,1,2,\cdots$ by using (\ref{eq:2004}) and (\ref{eq:2005}), 
\begin{theorem}  
The general expression of the generating function for Heun polynomial of type 2 is given by
\footnotesize
\begin{eqnarray}
&&\sum_{q_0 =0}^{\infty } \frac{(\gamma')_{q_0}}{q_0!} s_0^{q_0} \prod _{n=1}^{\infty } \left\{ \sum_{ q_n = q_{n-1}}^{\infty } s_n^{q_n }\right\} y(x) \nonumber\\
&&= \prod_{l=1}^{\infty } \frac{1}{(1-s_{l,\infty })} \mathbf{\Upsilon}(\lambda; s_{0,\infty } ;\eta) \nonumber\\
&&+ \Bigg\{ \prod_{l=2}^{\infty } \frac{1}{(1-s_{l,\infty })} \int_{0}^{1} dt_1\;t_1^{1+\lambda} \int_{0}^{1} du_1\;u_1^{\gamma +\lambda} \left(s_{1,\infty }^2-2(1-2\eta (1-t_1)(1-u_1))s_{1,\infty }+1\right)^{-\frac{1}{2}}\nonumber\\
&&\times \left( \frac{(1+s_{1,\infty })+\sqrt{s_{1,\infty }^2-2(1-2\eta (1-t_1)(1-u_1))s_{1,\infty }+1}}{2}\right)^{-\left(3+ \frac{\varphi +2(1+a)\lambda }{(1+a)}\right)}\nonumber\\
&&\times \widetilde{w}_{1,1}^{-(\alpha +\lambda )}\left(  \widetilde{w}_{1,1} \partial _{ \widetilde{w}_{1,1}}\right) \widetilde{w}_{1,1}^{\alpha -\beta } \left(  \widetilde{w}_{1,1} \partial _{ \widetilde{w}_{1,1}}\right) \widetilde{w}_{1,1}^{\beta +\lambda} \mathbf{\Upsilon}(\lambda ; s_0;\widetilde{w}_{1,1}) \Bigg\} z\nonumber\\
&&+ \sum_{n=2}^{\infty } \left\{ \prod_{l=n+1}^{\infty } \frac{1}{(1-s_{l,\infty })} \int_{0}^{1} dt_n\;t_n^{2n-1+\lambda } \int_{0}^{1} du_n\;u_n^{2(n-1)+\gamma +\lambda} \left( s_{n,\infty }^2-2(1-2\eta (1-t_n)(1-u_n))s_{n,\infty }+1\right)^{-\frac{1}{2}}\right.\nonumber\\
&&\times  \left( \frac{(1+s_{n,\infty })+\sqrt{s_{n,\infty }^2-2(1-2\eta (1-t_n)(1-u_n))s_{n,\infty }+1}}{2}\right)^{-\left(4n-1+\frac{\varphi +2(1+a)\lambda }{(1+a)} \right)}\nonumber\\
&&\times \widetilde{w}_{n,n}^{-(2(n-1)+\alpha +\lambda )}\left(  \widetilde{w}_{n,n} \partial _{ \widetilde{w}_{n,n}}\right) \widetilde{w}_{n,n}^{\alpha -\beta } \left(  \widetilde{w}_{n,n} \partial _{ \widetilde{w}_{n,n}}\right) \widetilde{w}_{n,n}^{2(n-1)+\beta +\lambda} \nonumber\\
&&\times \prod_{k=1}^{n-1} \Bigg\{ \int_{0}^{1} dt_{n-k}\;t_{n-k}^{2(n-k)-1+\lambda } \int_{0}^{1} du_{n-k} \;u_{n-k}^{2(n-k-1)+\gamma +\lambda }\left( s_{n-k}^2-2(1-2\widetilde{w}_{n-k+1,n} (1-t_{n-k})(1-u_{n-k}))s_{n-k}+1 \right)^{-\frac{1}{2}}\nonumber\\
&&\times \left( \frac{(1+s_{n-k})+\sqrt{s_{n-k}^2-2(1-2\widetilde{w}_{n-k+1,n} (1-t_{n-k})(1-u_{n-k}))s_{n-k}+1}}{2}\right)^{-\left(4(n-k)-1+\frac{\varphi +2(1+a)\lambda }{(1+a)} \right)} \nonumber\\
&&\times \left. \widetilde{w}_{n-k,n}^{-(2(n-k-1)+\alpha +\lambda )}\left(  \widetilde{w}_{n-k,n} \partial _{ \widetilde{w}_{n-k,n}}\right) \widetilde{w}_{n-k,n}^{\alpha -\beta } \left(  \widetilde{w}_{n-k,n} \partial _{ \widetilde{w}_{n-k,n}}\right) \widetilde{w}_{n-k,n}^{2(n-k-1)+\beta +\lambda} \Bigg\} \mathbf{\Upsilon}(\lambda ; s_0;\widetilde{w}_{1,n}) \right\} z^n\label{eq:2008}
\end{eqnarray}
\normalsize
where
\begin{equation}
\begin{cases} 
{ \displaystyle \mathbf{\Upsilon}(\lambda; s_{0,\infty } ;\eta)= \sum_{q_0 =0}^{\infty } \frac{(\gamma')_{q_0}}{q_0!} s_{0,\infty }^{q_0} \left( c_0 x^{\lambda } \sum_{i_0=0}^{q_0} \frac{(-q_0)_{i_0} \left(q_0+\frac{\varphi +2(1+a)\lambda }{(1+a)} \right)_{i_0}}{(1+\lambda )_{i_0}(\gamma +\lambda )_{i_0}} \eta ^{i_0} \right) }\cr
{ \displaystyle \mathbf{\Upsilon}(\lambda ; s_0;\widetilde{w}_{1,1}) = \sum_{q_0 =0}^{\infty } \frac{(\gamma')_{q_0}}{q_0!} s_0^{q_0}\left(c_0 x^{\lambda} \sum_{i_0=0}^{q_0} \frac{(-q_0)_{i_0} \left(q_0+\frac{\varphi +2(1+a)\lambda }{(1+a)} \right)_{i_0}}{(1+\lambda )_{i_0}(\gamma +\lambda )_{i_0}} \widetilde{w}_{1,1} ^{i_0} \right) }\cr
{ \displaystyle \mathbf{\Upsilon}(\lambda; s_0 ;\widetilde{w}_{1,n}) = \sum_{q_0 =0}^{\infty } \frac{(\gamma')_{q_0}}{q_0!} s_0^{q_0}\left(c_0 x^{\lambda} \sum_{i_0=0}^{q_0} \frac{(-q_0)_{i_0} \left(q_0+\frac{\varphi +2(1+a)\lambda }{(1+a)} \right)_{i_0}}{(1+\lambda )_{i_0}(\gamma +\lambda )_{i_0}} \widetilde{w}_{1,n} ^{i_0} \right)}
\end{cases}\nonumber 
\end{equation}
\end{theorem}
\begin{proof} 
Acting the summation operator $\displaystyle{ \sum_{q_0 =0}^{\infty } \frac{(\gamma')_{q_0}}{q_0!} s_0^{q_0} \prod _{n=1}^{\infty } \left\{ \sum_{ q_n = q_{n-1}}^{\infty } s_n^{q_n }\right\}}$ on the integral form of Heun polynomial of type 2 $y(x)$,
\begin{eqnarray}
&&\sum_{\alpha _0 =0}^{\infty } \frac{(\gamma')_{q_0}}{q_0!} s_0^{q_0} \prod _{n=1}^{\infty } \left\{ \sum_{ q_n = q_{n-1}}^{\infty } s_n^{q_n }\right\} y(x) \label{eq:20080}\\
&&= \sum_{q_0 =0}^{\infty } \frac{(\gamma')_{q_0}}{q_0!} s_0^{q_0} \prod _{n=1}^{\infty } \left\{ \sum_{ q_n = q_{n-1}}^{\infty } s_n^{q_n }\right\} \left\{ y_0(x)+y_1(x)+y_2(x)+y_3(x)+\cdots\right\} \nonumber
\end{eqnarray}
According to (\ref{eq:2006}), integral forms of sub-summation $y_0(x) $, $y_1(x)$, $y_2(x)$ and $y_3(x)$ are
\begin{subequations}
\begin{equation}
 y_0(x)= c_0 x^{\lambda } \sum_{i_0=0}^{q_0} \frac{(-q_0)_{i_0} \left(q_0+\frac{\varphi +2(1+a)\lambda }{(1+a)} \right)_{i_0}}{(1+\lambda )_{i_0}(\gamma +\lambda )_{i_0}} \eta ^{i_0}\label{eq:20024a}
\end{equation}
\begin{eqnarray}
 y_1(x) &=& \int_{0}^{1} dt_1\;t_1^{1+\lambda} \int_{0}^{1} du_1\;u_1^{\gamma +\lambda}
 \frac{1}{2\pi i} \oint dv_1 \;\frac{1}{v_1} \left( \frac{v_1-1}{v_1} \frac{1}{1-\eta (1-t_1)(1-u_1)v_1}\right)^{q_1} \nonumber\\
&\times& (1-\eta (1-t_1)(1-u_1)v_1)^{-\left(4+\frac{\varphi +2(1+a)\lambda }{(1+a)} \right)}\nonumber\\
&\times&  \overleftrightarrow {w}_{1,1}^{-(\alpha +\lambda )}\left(  \overleftrightarrow {w}_{1,1} \partial _{ \overleftrightarrow {w}_{1,1}}\right) \overleftrightarrow {w}_{1,1}^{\alpha -\beta } \left(  \overleftrightarrow {w}_{1,1} \partial _{ \overleftrightarrow {w}_{1,1}}\right) \overleftrightarrow {w}_{1,1}^{\beta +\lambda} \nonumber\\
&\times& \left( c_0 x^{\lambda } \sum_{i_0=0}^{q_0} \frac{(-q_0)_{i_0} \left(q_0+\frac{\varphi +2(1+a)\lambda }{(1+a)} \right)_{i_0}}{(1+\lambda )_{i_0}(\gamma +\lambda )_{i_0}} \overleftrightarrow {w}_{1,1} ^{i_0} \right) z\label{eq:20024b}
\end{eqnarray}
\begin{eqnarray}
 y_2(x) &=& \int_{0}^{1} dt_2\;t_2^{3+\lambda} \int_{0}^{1} du_2\;u_2^{2+\gamma +\lambda }
 \frac{1}{2\pi i} \oint dv_2 \;\frac{1}{v_2} \left( \frac{v_2-1}{v_2} \frac{1}{1-\eta (1-t_2)(1-u_2)v_2}\right)^{q_2} \nonumber\\
&\times& (1-\eta (1-t_2)(1-u_2)v_2)^{-\left(8+\frac{\varphi +2(1+a)\lambda }{(1+a)} \right)}\nonumber\\
&\times& \overleftrightarrow {w}_{2,2}^{-(2+\alpha +\lambda )}\left(  \overleftrightarrow {w}_{2,2} \partial _{ \overleftrightarrow {w}_{2,2}}\right) \overleftrightarrow {w}_{2,2}^{\alpha -\beta } \left(  \overleftrightarrow {w}_{2,2} \partial _{ \overleftrightarrow {w}_{2,2}}\right) \overleftrightarrow {w}_{2,2}^{2+\beta +\lambda} \nonumber\\
&\times& \int_{0}^{1} dt_1\;t_1^{1+\lambda} \int_{0}^{1} du_1\;u_1^{\gamma +\lambda}
 \frac{1}{2\pi i} \oint dv_1 \;\frac{1}{v_1} \left( \frac{v_1-1}{v_1} \frac{1}{1-\overleftrightarrow {w}_{2,2} (1-t_1)(1-u_1)v_1}\right)^{q_1} \nonumber\\
&\times& (1-\overleftrightarrow {w}_{2,2} (1-t_1)(1-u_1)v_1)^{-\left(4+\frac{\varphi +2(1+a)\lambda }{(1+a)} \right)}\nonumber\\
&\times& \overleftrightarrow {w}_{1,1}^{-(\alpha +\lambda )}\left(  \overleftrightarrow {w}_{1,1} \partial _{ \overleftrightarrow {w}_{1,1}}\right) \overleftrightarrow {w}_{1,1}^{\alpha -\beta } \left(  \overleftrightarrow {w}_{1,1} \partial _{ \overleftrightarrow {w}_{1,1}}\right) \overleftrightarrow {w}_{1,1}^{\beta +\lambda} \nonumber\\
&\times& \left( c_0 x^{\lambda } \sum_{i_0=0}^{q_0} \frac{(-q_0)_{i_0} \left(q_0+\frac{\varphi +2(1+a)\lambda }{(1+a)} \right)_{i_0}}{(1+\lambda )_{i_0}(\gamma +\lambda )_{i_0}} \overleftrightarrow {w}_{1,2} ^{i_0} \right)z^2 \label{eq:20024c}
\end{eqnarray}
\begin{eqnarray}
 y_3(x) &=& \int_{0}^{1} dt_3\;t_3^{5+\lambda} \int_{0}^{1} du_3\;u_3^{4+\gamma +\lambda }
 \frac{1}{2\pi i} \oint dv_3 \;\frac{1}{v_3} \left( \frac{v_3-1}{v_3} \frac{1}{1-\eta (1-t_3)(1-u_3)v_3}\right)^{q_3} \nonumber\\
&\times& (1-\eta (1-t_3)(1-u_3)v_3)^{-\left(12+\frac{\varphi +2(1+a)\lambda }{(1+a)} \right)} \nonumber\\
&\times& \overleftrightarrow {w}_{3,3}^{-(4+\alpha +\lambda )}\left(  \overleftrightarrow {w}_{3,3} \partial _{ \overleftrightarrow {w}_{3,3}}\right) \overleftrightarrow {w}_{3,3}^{\alpha -\beta } \left(  \overleftrightarrow {w}_{3,3} \partial _{ \overleftrightarrow {w}_{3,3}}\right) \overleftrightarrow {w}_{3,3}^{4+\beta +\lambda} \nonumber\\
&\times& \int_{0}^{1} dt_2\;t_2^{3+\lambda} \int_{0}^{1} du_2\;u_2^{2+\gamma +\lambda }
 \frac{1}{2\pi i} \oint dv_2 \;\frac{1}{v_2} \left( \frac{v_2-1}{v_2} \frac{1}{1-\overleftrightarrow {w}_{3,3} (1-t_2)(1-u_2)v_2}\right)^{q_2} \nonumber\\
&\times& (1-\overleftrightarrow {w}_{3,3}(1-t_2)(1-u_2)v_2)^{-\left(8+\frac{\varphi +2(1+a)\lambda }{(1+a)} \right)} \nonumber\\
&\times&  \overleftrightarrow {w}_{2,3}^{-(2+\alpha +\lambda )}\left(  \overleftrightarrow {w}_{2,3} \partial _{ \overleftrightarrow {w}_{2,3}}\right) \overleftrightarrow {w}_{2,3}^{\alpha -\beta } \left(  \overleftrightarrow {w}_{2,3} \partial _{ \overleftrightarrow {w}_{2,3}}\right) \overleftrightarrow {w}_{2,3}^{2+\beta +\lambda} \nonumber\\
&\times&  \int_{0}^{1} dt_1\;t_1^{1+\lambda} \int_{0}^{1} du_1\;u_1^{\gamma +\lambda}
 \frac{1}{2\pi i} \oint dv_1 \;\frac{1}{v_1} \left( \frac{v_1-1}{v_1} \frac{1}{1-\overleftrightarrow {w}_{2,3} (1-t_1)(1-u_1)v_1}\right)^{q_1} \nonumber\\
&\times& (1-\overleftrightarrow {w}_{2,3} (1-t_1)(1-u_1)v_1)^{-\left(4+\frac{\varphi +2(1+a)\lambda }{(1+a)} \right)} \nonumber\\
&\times& \overleftrightarrow {w}_{1,3}^{-(\alpha +\lambda )}\left(  \overleftrightarrow {w}_{1,3} \partial _{ \overleftrightarrow {w}_{1,3}}\right) \overleftrightarrow {w}_{1,3}^{\alpha -\beta } \left(  \overleftrightarrow {w}_{1,3} \partial _{ \overleftrightarrow {w}_{1,3}}\right) \overleftrightarrow {w}_{1,3}^{\beta +\lambda} \nonumber\\
&\times& \left( c_0 x^{\lambda } \sum_{i_0=0}^{q_0} \frac{(-q_0)_{i_0} \left(q_0+\frac{\varphi +2(1+a)\lambda }{(1+a)} \right)_{i_0}}{(1+\lambda )_{i_0}(\gamma +\lambda )_{i_0}} \overleftrightarrow {w}_{1,3} ^{i_0} \right) z^3 \label{eq:20024d}
\end{eqnarray}
\end{subequations}
Acting the summation operator $\displaystyle{ \sum_{q_0 =0}^{\infty } \frac{(\gamma')_{q_0}}{q_0!} s_0^{q_0} \prod _{n=1}^{\infty } \left\{ \sum_{ q_n = q_{n-1}}^{\infty } s_n^{q_n }\right\}}$ on (\ref{eq:20024a}),
\begin{eqnarray}
&&\sum_{q_0 =0}^{\infty } \frac{(\gamma')_{q_0}}{q_0!} s_0^{q_0} \prod _{n=1}^{\infty } \left\{ \sum_{ q_n = q_{n-1}}^{\infty } s_n^{q_n }\right\} y_0(x) \nonumber\\
&&= \prod_{l=1}^{\infty } \frac{1}{(1-s_{l,\infty })} \sum_{q_0 =0}^{\infty } \frac{(\gamma')_{q_0}}{q_0!} s_{0,\infty }^{q_0} \left( c_0 x^{\lambda } \sum_{i_0=0}^{q_0}  \frac{(-q_0)_{i_0} \left(q_0+\frac{\varphi +2(1+a)\lambda }{(1+a)} \right)_{i_0}}{(1+\lambda )_{i_0}(\gamma +\lambda )_{i_0}} \eta ^{i_0} \right) \hspace{2cm}\label{eq:20025}
\end{eqnarray}
Acting the summation operator $\displaystyle{ \sum_{q_0 =0}^{\infty } \frac{(\gamma')_{q_0}}{q_0!} s_0^{q_0} \prod _{n=1}^{\infty } \left\{ \sum_{ q_n = q_{n-1}}^{\infty } s_n^{q_n }\right\}}$ on (\ref{eq:20024b}),
\begin{eqnarray}
&&\sum_{q_0 =0}^{\infty } \frac{(\gamma')_{q_0}}{q_0!} s_0^{q_0} \prod _{n=1}^{\infty } \left\{ \sum_{ q_n = q_{n-1}}^{\infty } s_n^{q_n }\right\} y_1(x) \nonumber\\
&&= \prod_{l=2}^{\infty } \frac{1}{(1-s_{l,\infty })} \int_{0}^{1} dt_1\;t_1^{1+\lambda } \int_{0}^{1} du_1\;u_1^{\gamma +\lambda}
 \frac{1}{2\pi i} \oint dv_1 \;\frac{1}{v_1} (1-\eta (1-t_1)(1-u_1)v_1)^{-\left(4+\frac{\varphi +2(1+a)\lambda }{(1+a)}\right)} \nonumber\\
&&\times \sum_{q_1 =q_0}^{\infty }\left( \frac{v_1-1}{v_1} \frac{s_{1,\infty }}{1-\eta (1-t_1)(1-u_1)v_1}\right)^{q_1}  \overleftrightarrow {w}_{1,1}^{-(\alpha +\lambda )}\left(  \overleftrightarrow {w}_{1,1} \partial _{ \overleftrightarrow {w}_{1,1}}\right) \overleftrightarrow {w}_{1,1}^{\alpha -\beta } \left(  \overleftrightarrow {w}_{1,1} \partial _{ \overleftrightarrow {w}_{1,1}}\right)\overleftrightarrow {w}_{1,1}^{\beta +\lambda} \nonumber\\
&&\times  \sum_{q_0 =0}^{\infty } \frac{(\gamma' )_{q_0}}{q_0!}s_0^{q_0}\left( c_0 x^{\lambda } \sum_{i_0=0}^{q_0} \frac{(-q_0)_{i_0} \left(q_0+\frac{\varphi +2(1+a)\lambda }{(1+a)} \right)_{i_0}}{(1+\lambda )_{i_0}(\gamma +\lambda )_{i_0}} \overleftrightarrow {w}_{1,1} ^{i_0} \right) z \label{eq:20026}
\end{eqnarray}
Replace $q_i$, $q_j$ and $r_i$ by $q_1$, $q_0$ and ${ \displaystyle \frac{v_1-1}{v_1} \frac{s_{1,\infty }}{1-\eta (1-t_1)(1-u_1)v_1}}$ in (\ref{eq:2005}). Take the new (\ref{eq:2005}) into (\ref{eq:20026}).
\begin{eqnarray}
&&\sum_{q_0 =0}^{\infty } \frac{(\gamma')_{q_0}}{q_0!} s_0^{q_0} \prod _{n=1}^{\infty } \left\{ \sum_{ q_n = q_{n-1}}^{\infty } s_n^{q_n }\right\} y_1(x) \nonumber\\
&&= \prod_{l=2}^{\infty } \frac{1}{(1-s_{l,\infty })} \int_{0}^{1} dt_1\;t_1^{1+\lambda } \int_{0}^{1} du_1\;u_1^{\gamma +\lambda}
 \frac{1}{2\pi i} \oint dv_1 \;\frac{(1-\eta (1-t_1)(1-u_1)v_1)^{-\left(3+\frac{\varphi +2(1+a)\lambda }{(1+a)} \right)} }{-\eta (1-t_1)(1-u_1)v_1^2+ (1-s_{1,\infty })v_1+s_{1,\infty } } \nonumber\\
&&\times  \overleftrightarrow {w}_{1,1}^{-(\alpha +\lambda )}\left(  \overleftrightarrow {w}_{1,1} \partial _{ \overleftrightarrow {w}_{1,1}}\right) \overleftrightarrow {w}_{1,1}^{\alpha -\beta } \left(  \overleftrightarrow {w}_{1,1} \partial _{ \overleftrightarrow {w}_{1,1}}\right)\overleftrightarrow {w}_{1,1}^{\beta +\lambda} \label{eq:20027}\\
&&\times  \sum_{q_0 =0}^{\infty } \frac{(\gamma' )_{q_0}}{q_0!} \left( \frac{v_1-1}{v_1} \frac{s_{0,\infty }}{1-\eta (1-t_1)(1-u_1)v_1}\right)^{q_0}\left( c_0 x^{\lambda } \sum_{i_0=0}^{q_0} \frac{(-q_0)_{i_0} \left(q_0+\frac{\varphi +2(1+a)\lambda }{(1+a)} \right)_{i_0}}{(1+\lambda )_{i_0}(\gamma +\lambda )_{i_0}} \overleftrightarrow {w}_{1,1} ^{i_0} \right) z \nonumber
\end{eqnarray}
By using Cauchy's integral formula, the contour integrand has poles at
\footnotesize ${\displaystyle v_1= \frac{1-s_{1,\infty }-\sqrt{(1-s_{1,\infty })^2+4\eta (1-t_1)(1-u_1)s_{1,\infty }}}{2\eta (1-t_1)(1-u_1)}  \;\;\mbox{or}\;\frac{1-s_{1,\infty }+\sqrt{(1-s_{1,\infty })^2+4\eta (1-t_1)(1-u_1)s_{1,\infty }}}{2\eta (1-t_1)(1-u_1)} }$\normalsize and ${ \displaystyle \frac{1-s_{1,\infty }-\sqrt{(1-s_{1,\infty })^2+4\eta (1-t_1)(1-u_1)s_{1,\infty }}}{2\eta (1-t_1)(1-u_1)}}$ is only inside the unit circle. As we compute the residue there in (\ref{eq:20027}) we obtain
\begin{eqnarray}
&&\sum_{q_0 =0}^{\infty } \frac{(\gamma')_{q_0}}{q_0!} s_0^{q_0} \prod _{n=1}^{\infty } \left\{ \sum_{ q_n = q_{n-1}}^{\infty } s_n^{q_n }\right\} y_1(x) \nonumber\\
&&= \prod_{l=2}^{\infty } \frac{1}{(1-s_{l,\infty })} \int_{0}^{1} dt_1\;t_1^{1+\lambda} \int_{0}^{1} du_1\;u_1^{\gamma +\lambda}
 \left( s_{1,\infty }^2-2(1-2\eta (1-t_1)(1-u_1))s_{1,\infty }+1 \right)^{-\frac{1}{2}} \nonumber\\
&&\times \left(\frac{1+s_{1,\infty }+\sqrt{s_{1,\infty }^2-2(1-2\eta (1-t_1)(1-u_1))s_{1,\infty }+1}}{2}\right)^{-\left(3+\frac{\varphi +2(1+a)\lambda }{(1+a)} \right)}\nonumber\\
&&\times  \widetilde{w}_{1,1}^{-(\alpha +\lambda )}\left(  \widetilde{w}_{1,1} \partial _{ \widetilde{w}_{1,1}}\right) \widetilde{w}_{1,1}^{\alpha -\beta } \left(  \widetilde{w}_{1,1} \partial _{ \widetilde{w}_{1,1}}\right)\widetilde{w}_{1,1}^{\beta +\lambda} \nonumber\\
&&\times  \sum_{q_0 =0}^{\infty } \frac{(\gamma' )_{q_0}}{q_0!} s_0^{q_0}\left( c_0 x^{\lambda } \sum_{i_0=0}^{q_0} \frac{(-q_0)_{i_0} \left(q_0+\frac{\varphi +2(1+a)\lambda }{(1+a)} \right)_{i_0}}{(1+\lambda )_{i_0}(\gamma +\lambda )_{i_0}} \widetilde{w}_{1,1} ^{i_0} \right) z \label{eq:20028}
\end{eqnarray}
where
\begin{eqnarray}
\widetilde{w}_{1,1} &=& \frac{v_1}{(v_1-1)}\; \frac{\eta t_1 u_1}{1- \eta v_1 (1-t_1)(1-u_1)}\Bigg|_{v_1=\frac{1-s_{1,\infty }-\sqrt{(1-s_{1,\infty })^2+4\eta (1-t_1)(1-u_1)s_{1,\infty }}}{2\eta (1-t_1)(1-u_1)}}\nonumber\\
&=& \frac{\eta t_1 u_1 \left\{ 1+ (s_{1,\infty }+2\eta(1-t_1)(1-u_1) )s_{1,\infty }\right\}}{2(1-\eta (1-t_1)(1-u_1))^2 s_{1,\infty }}\nonumber\\
&&-\frac{\eta t_1 u_1(1+s_{1,\infty })\sqrt{s_{1,\infty }^2-2(1-2\eta (1-t_1)(1-u_1))s_{1,\infty }+1}}{2(1-\eta (1-t_1)(1-u_1))^2 s_{1,\infty }}\nonumber
\end{eqnarray}
Acting the summation operator $\displaystyle{ \sum_{q_0 =0}^{\infty } \frac{(\gamma')_{q_0}}{q_0!} s_0^{q_0} \prod _{n=1}^{\infty } \left\{ \sum_{ q_n = q_{n-1}}^{\infty } s_n^{q_n }\right\}}$ on (\ref{eq:20024c}),
\begin{eqnarray}
&&\sum_{q_0 =0}^{\infty } \frac{(\gamma')_{q_0}}{q_0!} s_0^{q_0} \prod _{n=1}^{\infty } \left\{ \sum_{ q_n = q_{n-1}}^{\infty } s_n^{q_n }\right\} y_2(x) \nonumber\\
&&= \prod_{l=3}^{\infty } \frac{1}{(1-s_{l,\infty })} \int_{0}^{1} dt_2\;t_2^{3+\lambda} \int_{0}^{1} du_2\;u_2^{2+\gamma +\lambda }
 \frac{1}{2\pi i} \oint dv_2 \;\frac{1}{v_2} (1-\eta (1-t_2)(1-u_2)v_2)^{-\left(8+\frac{\varphi +2(1+a)\lambda }{(1+a)} \right)} \nonumber\\
&&\times \sum_{q_2 =q_1}^{\infty }\left( \frac{v_2-1}{v_2} \frac{s_{2,\infty }}{1-\eta (1-t_2)(1-u_2)v_2}\right)^{q_2} \overleftrightarrow {w}_{2,2}^{-(2+\alpha +\lambda )}\left(  \overleftrightarrow {w}_{2,2} \partial _{ \overleftrightarrow {w}_{2,2}}\right) \overleftrightarrow {w}_{2,2}^{\alpha -\beta } \left(  \overleftrightarrow {w}_{2,2} \partial _{ \overleftrightarrow {w}_{2,2}}\right)\overleftrightarrow {w}_{2,2}^{2+\beta +\lambda}  \nonumber\\
&&\times \int_{0}^{1} dt_1\;t_1^{1+\lambda} \int_{0}^{1} du_1\;u_1^{\gamma +\lambda}
 \frac{1}{2\pi i} \oint dv_1 \;\frac{1}{v_1} (1-\overleftrightarrow {w}_{2,2} (1-t_1)(1-u_1)v_1)^{-\left( 4+\frac{\varphi +2(1+a)\lambda }{(1+a)}\right)} \nonumber\\
&&\times \sum_{q_1 =q_0}^{\infty }\left( \frac{v_1-1}{v_1} \frac{s_1}{1-\overleftrightarrow {w}_{2,2}(1-t_1)(1-u_1)v_1}\right)^{q_1} \overleftrightarrow {w}_{1,2}^{-(\alpha +\lambda )}\left(  \overleftrightarrow {w}_{1,2} \partial _{ \overleftrightarrow {w}_{1,2}}\right) \overleftrightarrow {w}_{1,2}^{\alpha -\beta } \left(  \overleftrightarrow {w}_{1,2} \partial _{ \overleftrightarrow {w}_{1,2}}\right)\overleftrightarrow {w}_{1,2}^{\beta +\lambda} \nonumber\\
&&\times  \sum_{q_0 =0}^{\infty } \frac{(\gamma' )_{q_0}}{q_0!} s_0^{q_0}\left( c_0 x^{\lambda } \sum_{i_0=0}^{q_0} \frac{(-q_0)_{i_0} \left(q_0+\frac{\varphi +2(1+a)\lambda }{(1+a)} \right)_{i_0}}{(1+\lambda )_{i_0}(\gamma +\lambda )_{i_0}} \overleftrightarrow {w}_{1,2} ^{i_0} \right) z^2 \label{eq:20029}
\end{eqnarray}
Replace $q_i$, $q_j$ and $r_i$ by $q_2$, $q_1$ and ${ \displaystyle \frac{v_2-1}{v_2} \frac{s_{2,\infty }}{1-\eta (1-t_2)(1-u_2)v_2}}$ in (\ref{eq:2005}). Take the new (\ref{eq:2005}) into (\ref{eq:20029}).
\begin{eqnarray}
&&\sum_{q_0 =0}^{\infty } \frac{(\gamma')_{q_0}}{q_0!} s_0^{q_0} \prod _{n=1}^{\infty } \left\{ \sum_{ q_n = q_{n-1}}^{\infty } s_n^{q_n }\right\} y_2(x) \nonumber\\
&&= \prod_{l=3}^{\infty } \frac{1}{(1-s_{l,\infty })} \int_{0}^{1} dt_2\;t_2^{3+\lambda } \int_{0}^{1} du_2\;u_2^{2+\gamma +\lambda }
 \frac{1}{2\pi i} \oint dv_2 \;\frac{\left(1-\eta (1-t_2)(1-u_2)v_2\right)^{-\left(7+\frac{\varphi +2(1+a)\lambda }{(1+a)} \right)}}{-\eta (1-t_2)(1-u_2)v_2^2+ (1-s_{2,\infty })v_2+s_{2,\infty } } \nonumber\\
&&\times \overleftrightarrow {w}_{2,2}^{-(2+\alpha +\lambda )}\left(  \overleftrightarrow {w}_{2,2} \partial _{ \overleftrightarrow {w}_{2,2}}\right) \overleftrightarrow {w}_{2,2}^{\alpha -\beta } \left(  \overleftrightarrow {w}_{2,2} \partial _{ \overleftrightarrow {w}_{2,2}}\right)\overleftrightarrow {w}_{2,2}^{2+\beta +\lambda} \nonumber\\
&&\times \int_{0}^{1} dt_1\;t_1^{1+\lambda} \int_{0}^{1} du_1\;u_1^{\gamma +\lambda}
 \frac{1}{2\pi i} \oint dv_1 \;\frac{1}{v_1} \left(1-\overleftrightarrow {w}_{2,2} (1-t_1)(1-u_1)v_1\right)^{-\left( 4+\frac{\varphi +2(1+a)\lambda }{(1+a)}\right)} \nonumber\\
&&\times \sum_{q_1 =q_0}^{\infty }\left( \frac{v_2-1}{v_2} \frac{s_{1,\infty }}{1-\eta (1-t_2)(1-u_2)v_2} \frac{v_1-1}{v_1}\frac{1}{1-\overleftrightarrow {w}_{2,2}(1-t_1)(1-u_1)v_1}\right)^{q_1} \nonumber\\
&&\times \overleftrightarrow {w}_{1,2}^{-(\alpha +\lambda )}\left(  \overleftrightarrow {w}_{1,2} \partial _{ \overleftrightarrow {w}_{1,2}}\right) \overleftrightarrow {w}_{1,2}^{\alpha -\beta } \left(  \overleftrightarrow {w}_{1,2} \partial _{ \overleftrightarrow {w}_{1,2}}\right)\overleftrightarrow {w}_{1,2}^{\beta +\lambda}  \nonumber\\
&&\times  \sum_{q_0 =0}^{\infty } \frac{(\gamma' )_{q_0}}{q_0!} s_0^{q_0}\left( c_0 x^{\lambda } \sum_{i_0=0}^{q_0} \frac{(-q_0)_{i_0} \left(q_0+\frac{\varphi +2(1+a)\lambda }{(1+a)} \right)_{i_0}}{(1+\lambda )_{i_0}(\gamma +\lambda )_{i_0}} \overleftrightarrow {w}_{1,2} ^{i_0} \right) z^2 \label{eq:20030}
\end{eqnarray}
By using Cauchy's integral formula, the contour integrand has poles at
\footnotesize ${\displaystyle
 v_2= \frac{1-s_{2,\infty }-\sqrt{(1-s_{2,\infty })^2+4\eta (1-t_2)(1-u_2)s_{2,\infty }}}{2\eta (1-t_2)(1-u_2)}  \;\;\mbox{or}\;\frac{1-s_{2,\infty }+\sqrt{(1-s_{2,\infty })^2+4\eta (1-t_2)(1-u_2)s_{2,\infty }}}{2\eta (1-t_2)(1-u_2)} }$ \normalsize
and ${ \displaystyle\frac{1-s_{2,\infty }-\sqrt{(1-s_{2,\infty })^2+4\eta (1-t_2)(1-u_2)s_{2,\infty }}}{2\eta (1-t_2)(1-u_2)}}$ is only inside the unit circle. As we compute the residue there in (\ref{eq:20030}) we obtain
\begin{eqnarray}
&&\sum_{q_0 =0}^{\infty } \frac{(\gamma')_{q_0}}{q_0!} s_0^{q_0} \prod _{n=1}^{\infty } \left\{ \sum_{ q_n = q_{n-1}}^{\infty } s_n^{q_n }\right\} y_2(x) \nonumber\\
&&= \prod_{l=3}^{\infty } \frac{1}{(1-s_{l,\infty })} \int_{0}^{1} dt_2\;t_2^{3+\lambda} \int_{0}^{1} du_2\;u_2^{2+\gamma +\lambda }
 \left( s_{2,\infty }^2-2(1-2\eta (1-t_2)(1-u_2))s_{2,\infty }+1\right)^{-\frac{1}{2}}\nonumber\\
&&\times \left(\frac{1+s_{2,\infty }+\sqrt{s_{2,\infty }^2-2(1-2\eta (1-t_2)(1-u_2))s_{2,\infty }+1}}{2}\right)^{-\left(7+\frac{\varphi +2(1+a)\lambda }{(1+a)} \right)} \nonumber\\
&&\times  \widetilde{w}_{2,2}^{-(2+\alpha +\lambda )}\left( \widetilde{w}_{2,2} \partial _{ \widetilde{w}_{2,2}}\right) \widetilde{w}_{2,2}^{\alpha -\beta } \left( \widetilde{w}_{2,2} \partial _{ \widetilde{w}_{2,2}}\right)\widetilde{w}_{2,2}^{2+\beta +\lambda}  \nonumber\\
&&\times \int_{0}^{1} dt_1\;t_1^{1+\lambda} \int_{0}^{1} du_1\;u_1^{\gamma +\lambda}
 \frac{1}{2\pi i} \oint dv_1 \;\frac{1}{v_1} \left(1-\widetilde{w}_{2,2} (1-t_1)(1-u_1)v_1\right)^{-\left(4+\frac{\varphi +2(1+a)\lambda }{(1+a)} \right)} \nonumber\\
&&\times \sum_{q_1 =q_0}^{\infty }\left( \frac{v_1-1}{v_1}\frac{s_1}{1-\widetilde{w}_{2,2}(1-t_1)(1-u_1)v_1}\right)^{q_1}
 \ddot{w}_{1,2}^{-(\alpha +\lambda )}\left( \ddot{w}_{1,2} \partial _{ \ddot{w}_{1,2}}\right) \ddot{w}_{1,2}^{\alpha -\beta } \left( \ddot{w}_{1,2} \partial _{ \ddot{w}_{1,2}}\right)\ddot{w}_{1,2}^{\beta +\lambda}  \nonumber\\
&&\times  \sum_{q_0 =0}^{\infty } \frac{(\gamma' )_{q_0}}{q_0!} s_0^{q_0}\left( c_0 x^{\lambda } \sum_{i_0=0}^{q_0} \frac{(-q_0)_{i_0} \left(q_0+\frac{\varphi +2(1+a)\lambda }{(1+a)} \right)_{i_0}}{(1+\lambda )_{i_0}(\gamma +\lambda )_{i_0}} \ddot{w}_{1,2} ^{i_0} \right) z^2 \label{eq:20031}
\end{eqnarray}
where
\begin{eqnarray}
\widetilde{w}_{2,2} &=& \frac{v_2}{(v_2-1)}\; \frac{\eta t_2 u_2}{1- \eta v_2 (1-t_2)(1-u_2)}\Bigg|_{v_2=\frac{1-s_{2,\infty }-\sqrt{(1-s_{2,\infty })^2+4\eta (1-t_2)(1-u_2)s_{2,\infty }}}{2\eta (1-t_2)(1-u_2)}}\nonumber\\
&=& \frac{\eta t_2 u_2 \left\{ 1+ (s_{2,\infty }+2\eta(1-t_2)(1-u_2) )s_{2,\infty }\right\}}{2(1-\eta (1-t_2)(1-u_2))^2 s_{2,\infty }}\nonumber\\
&&- \frac{\eta t_2 u_2  (1+s_{2,\infty })\sqrt{s_{2,\infty }^2-2(1-2\eta (1-t_2)(1-u_2))s_{2,\infty }+1}}{2(1-\eta (1-t_2)(1-u_2))^2 s_{2,\infty }}\nonumber
\end{eqnarray}
and
\begin{equation}
\ddot{w}_{1,2} = \frac{v_1}{(v_1-1)}\; \frac{\widetilde{w}_{2,2} t_1 u_1}{1- \widetilde{w}_{2,2}v_1 (1-t_1)(1-u_1)}\nonumber
\end{equation}
Replace $q_i$, $q_j$ and $r_i$ by $q_1$, $q_0$ and ${ \displaystyle \frac{v_1-1}{v_1}\frac{s_1}{1-\widetilde{w}_{2,2}(1-t_1)(1-u_1)v_1}}$ in (\ref{eq:2005}). Take the new (\ref{eq:2005}) into (\ref{eq:20031}).
\begin{eqnarray}
&&\sum_{q_0 =0}^{\infty } \frac{(\gamma')_{q_0}}{q_0!} s_0^{q_0} \prod _{n=1}^{\infty } \left\{ \sum_{ q_n = q_{n-1}}^{\infty } s_n^{q_n }\right\} y_2(x) \nonumber\\
&&= \prod_{l=3}^{\infty } \frac{1}{(1-s_{l,\infty })} \int_{0}^{1} dt_2\;t_2^{3+\lambda} \int_{0}^{1} du_2\;u_2^{2+\gamma +\lambda }
 \left( s_{2,\infty }^2-2(1-2\eta (1-t_2)(1-u_2))s_{2,\infty }+1\right)^{-\frac{1}{2}}\nonumber\\
&&\times \left(\frac{1+s_{2,\infty }+\sqrt{s_{2,\infty }^2-2(1-2\eta (1-t_2)(1-u_2))s_{2,\infty }+1}}{2}\right)^{-\left(7+\frac{\varphi +2(1+a)\lambda }{(1+a)} \right)} \nonumber\\
&&\times  \widetilde{w}_{2,2}^{-(2+\alpha +\lambda )}\left( \widetilde{w}_{2,2} \partial _{ \widetilde{w}_{2,2}}\right) \widetilde{w}_{2,2}^{\alpha -\beta } \left( \widetilde{w}_{2,2} \partial _{ \widetilde{w}_{2,2}}\right)\widetilde{w}_{2,2}^{2+\beta +\lambda} \nonumber\\
&&\times \int_{0}^{1} dt_1\;t_1^{1+\lambda} \int_{0}^{1} du_1\;u_1^{\gamma +\lambda}
 \frac{1}{2\pi i} \oint dv_1 \;\frac{\left(1-\widetilde{w}_{2,2} (1-t_1)(1-u_1)v_1\right)^{-\left(3+\frac{\varphi +2(1+a)\lambda }{(1+a)}\right)}}{-\widetilde{w}_{2,2} (1-t_1)(1-u_1)v_1^2+(1-s_1)v_1+s_1} \nonumber\\
&&\times \ddot{w}_{1,2}^{-(\alpha +\lambda )}\left( \ddot{w}_{1,2} \partial _{ \ddot{w}_{1,2}}\right) \ddot{w}_{1,2}^{\alpha -\beta } \left( \ddot{w}_{1,2} \partial _{ \ddot{w}_{1,2}}\right)\ddot{w}_{1,2}^{\beta +\lambda} \label{eq:20032}\\
&&\times  \sum_{q_0 =0}^{\infty } \frac{(\gamma')_{q_0}}{q_0!} \left( \frac{v_1-1}{v_1}\frac{s_{0,1}}{1-\widetilde{w}_{2,2}(1-t_1)(1-u_1)v_1}\right)^{q_0} \left( c_0 x^{\lambda } \sum_{i_0=0}^{q_0} \frac{(-q_0)_{i_0} \left(q_0+\frac{\varphi +2(1+a)\lambda }{(1+a)} \right)_{i_0}}{(1+\lambda )_{i_0}(\gamma +\lambda )_{i_0}} \ddot{w}_{1,2} ^{i_0} \right) z^2 \nonumber
\end{eqnarray}
By using Cauchy's integral formula, the contour integrand has poles at \footnotesize ${\displaystyle
 v_1= \frac{1-s_1-\sqrt{(1-s_1)^2+4\widetilde{w}_{2,2} (1-t_1)(1-u_1)s_1}}{2\widetilde{w}_{2,2} (1-t_1)(1-u_1)}  \;\;\mbox{or}\;\frac{1-s_1+\sqrt{(1-s_1)^2+4\widetilde{w}_{2,2} (1-t_1)(1-u_1)s_1}}{2\widetilde{w}_{2,2} (1-t_1)(1-u_1)}}$ \normalsize
and ${ \displaystyle \frac{1-s_1-\sqrt{(1-s_1)^2+4\widetilde{w}_{2,2} (1-t_1)(1-u_1)s_1}}{2\widetilde{w}_{2,2} (1-t_1)(1-u_1)}}$ is only inside the unit circle. As we compute the residue there in (\ref{eq:20032}) we obtain
\begin{eqnarray}
&&\sum_{q_0 =0}^{\infty } \frac{(\gamma')_{q_0}}{q_0!} s_0^{q_0} \prod _{n=1}^{\infty } \left\{ \sum_{ q_n = q_{n-1}}^{\infty } s_n^{q_n }\right\} y_2(x) \nonumber\\
&&= \prod_{l=3}^{\infty } \frac{1}{(1-s_{l,\infty })} \int_{0}^{1} dt_2\;t_2^{3+\lambda} \int_{0}^{1} du_2\;u_2^{2+\gamma +\lambda }
\left( s_{2,\infty }^2-2(1-2\eta (1-t_2)(1-u_2))s_{2,\infty }+1\right)^{-\frac{1}{2}}\nonumber\\
&&\times \left(\frac{1+s_{2,\infty }+\sqrt{s_{2,\infty }^2-2(1-2\eta (1-t_2)(1-u_2))s_{2,\infty }+1}}{2}\right)^{-\left(7+\frac{\varphi +2(1+a)\lambda }{(1+a)} \right)} \nonumber\\
&&\times \widetilde{w}_{2,2}^{-(2+\alpha +\lambda )}\left( \widetilde{w}_{2,2} \partial _{ \widetilde{w}_{2,2}}\right) \widetilde{w}_{2,2}^{\alpha -\beta } \left( \widetilde{w}_{2,2} \partial _{ \widetilde{w}_{2,2}}\right)\widetilde{w}_{2,2}^{2+\beta +\lambda} \nonumber\\
&&\times \int_{0}^{1} dt_1\;t_1^{1+\lambda } \int_{0}^{1} du_1\;u_1^{\gamma +\lambda}
 \left( s_1^2-2(1-2\widetilde{w}_{2,2} (1-t_1)(1-u_1))s_1+1\right)^{-\frac{1}{2}}\nonumber\\
&&\times \left(\frac{1+s_1+\sqrt{s_1^2-2(1-2\widetilde{w}_{2,2}(1-t_1)(1-u_1))s_1+1}}{2}\right)^{-\left(3+\frac{\varphi +2(1+a)\lambda }{(1+a)}\right)} \nonumber\\
&&\times \widetilde{w}_{1,2}^{-(\alpha +\lambda )}\left( \widetilde{w}_{1,2} \partial _{ \widetilde{w}_{1,2}}\right) \widetilde{w}_{1,2}^{\alpha -\beta } \left( \widetilde{w}_{1,2} \partial _{ \widetilde{w}_{1,2}}\right)\widetilde{w}_{1,2}^{\beta +\lambda} \nonumber\\
&&\times \sum_{q_0 =0}^{\infty } \frac{(\gamma' )_{q_0}}{q_0!} s_0^{q_0}\left( c_0 x^{\lambda } \sum_{i_0=0}^{q_0} \frac{(-q_0)_{i_0} \left(q_0+\frac{\varphi +2(1+a)\lambda }{(1+a)} \right)_{i_0}}{(1+\lambda )_{i_0}(\gamma +\lambda )_{i_0}} \widetilde{w}_{1,2} ^{i_0} \right) z^2  \label{eq:20033}
\end{eqnarray}
where
\begin{eqnarray}
\widetilde{w}_{1,2} &=& \frac{v_1}{(v_1-1)}\; \frac{\widetilde{w}_{2,2} t_1 u_1}{1- \widetilde{w}_{2,2} v_1 (1-t_1)(1-u_1)}\Bigg|_{v_1=\frac{1-s_1-\sqrt{(1-s_1)^2+4\widetilde{w}_{2,2} (1-t_1)(1-u_1)s_1}}{2\widetilde{w}_{2,2} (1-t_1)(1-u_1)}}\nonumber\\
&=& \frac{\widetilde{w}_{2,2} t_1 u_1 \left\{ 1+ (s_1+2\widetilde{w}_{2,2}(1-t_1)(1-u_1) )s_1\right\}}{2(1-\widetilde{w}_{2,2}(1-t_1)(1-u_1))^2 s_1}\nonumber\\
&&-\frac{\widetilde{w}_{2,2} t_1 u_1  (1+s_1)\sqrt{s_1^2-2(1-2\widetilde{w}_{2,2} (1-t_1)(1-u_1))s_1+1} }{2(1-\widetilde{w}_{2,2}(1-t_1)(1-u_1))^2 s_1}\nonumber
\end{eqnarray}
Acting the summation operator $\displaystyle{ \sum_{q_0 =0}^{\infty } \frac{(\gamma')_{q_0}}{q_0!} s_0^{q_0} \prod _{n=1}^{\infty } \left\{ \sum_{ q_n = q_{n-1}}^{\infty } s_n^{q_n }\right\}}$ on (\ref{eq:20024d}),
\begin{eqnarray}
&&\sum_{q_0 =0}^{\infty } \frac{(\gamma')_{q_0}}{q_0!} s_0^{q_0} \prod _{n=1}^{\infty } \left\{ \sum_{ q_n = q_{n-1}}^{\infty } s_n^{q_n }\right\} y_3(x) \nonumber\\
&&= \prod_{l=4}^{\infty } \frac{1}{(1-s_{l,\infty })} \int_{0}^{1} dt_3\;t_3^{5+\lambda} \int_{0}^{1} du_3\;u_3^{4+\gamma +\lambda}\left( s_{3,\infty }^2-2(1-2\eta (1-t_3)(1-u_3))s_{3,\infty }+1\right)^{-\frac{1}{2}}\nonumber\\
&&\times \left(\frac{1+s_{3,\infty }+\sqrt{s_{3,\infty }^2-2(1-2\eta (1-t_3)(1-u_3))s_{3,\infty }+1}}{2}\right)^{-\left(11+\frac{\varphi +2(1+a)\lambda }{(1+a)} \right)} \nonumber\\
&&\times \widetilde{w}_{3,3}^{-(4+\alpha +\lambda )}\left( \widetilde{w}_{3,3} \partial _{ \widetilde{w}_{3,3}}\right) \widetilde{w}_{3,3}^{\alpha -\beta } \left( \widetilde{w}_{3,3} \partial _{ \widetilde{w}_{3,3}}\right)\widetilde{w}_{3,3}^{4+\beta +\lambda}  \nonumber\\
&&\times \int_{0}^{1} dt_2\;t_2^{3+\lambda} \int_{0}^{1} du_2\;u_2^{2+\gamma +\lambda }\left(s_2^2-2(1-2\widetilde{w}_{3,3}(1-t_2)(1-u_2))s_2+1\right)^{-\frac{1}{2}}\nonumber\\
&&\times \left(\frac{1+s_2+\sqrt{s_2^2-2(1-2\widetilde{w}_{3,3}(1-t_2)(1-u_2))s_2+1}}{2}\right)^{-\left(7+\frac{\varphi +2(1+a)\lambda }{(1+a)} \right)} \nonumber\\
&&\times  \widetilde{w}_{2,3}^{-(2+\alpha +\lambda )}\left( \widetilde{w}_{2,3} \partial _{ \widetilde{w}_{2,3}}\right) \widetilde{w}_{2,3}^{\alpha -\beta } \left( \widetilde{w}_{2,3} \partial _{ \widetilde{w}_{2,3}}\right)\widetilde{w}_{2,3}^{2+\beta +\lambda}  \nonumber\\
&&\times \int_{0}^{1} dt_1\;t_1^{1+\lambda } \int_{0}^{1} du_1\;u_1^{\gamma +\lambda}\left( s_1^2-2(1-2\widetilde{w}_{2,3} (1-t_1)(1-u_1))s_1+1\right)^{-\frac{1}{2}}\nonumber\\
&&\times \left(\frac{1+s_1+\sqrt{s_1^2-2(1-2\widetilde{w}_{2,3}(1-t_1)(1-u_1))s_1+1}}{2}\right)^{-\left(3+\frac{\varphi +2(1+a)\lambda }{(1+a)}\right)} \nonumber\\
&&\times  \widetilde{w}_{1,3}^{-(\alpha +\lambda )}\left( \widetilde{w}_{1,3} \partial _{ \widetilde{w}_{1,3}}\right) \widetilde{w}_{1,3}^{\alpha -\beta } \left( \widetilde{w}_{1,3} \partial _{ \widetilde{w}_{1,3}}\right)\widetilde{w}_{1,3}^{\beta +\lambda} \nonumber\\
&&\times \sum_{q_0 =0}^{\infty } \frac{(\gamma' )_{q_0}}{q_0!} s_0^{q_0}\left( c_0 x^{\lambda } \sum_{i_0=0}^{q_0} \frac{(-q_0)_{i_0} \left(q_0+\frac{\varphi +2(1+a)\lambda }{(1+a)} \right)_{i_0}}{(1+\lambda )_{i_0}(\gamma +\lambda )_{i_0}} \widetilde{w}_{1,3} ^{i_0} \right) z^3  \label{eq:20034}
\end{eqnarray}

\vspace{1cm}
where
\begin{eqnarray}
\widetilde{w}_{3,3} &=& \frac{v_3}{(v_3-1)}\; \frac{\eta  t_3 u_3}{1- \eta (1-t_3)(1-u_3)v_3}\Bigg|_{v_3=\frac{1-s_{3,\infty }-\sqrt{(1-s_{3,\infty })^2+4\eta (1-t_3)(1-u_3)s_{3,\infty }}}{2\eta (1-t_3)(1-u_3)}}\nonumber\\
&=& \frac{\eta t_3 u_3 \left\{ 1+ (s_{3,\infty }+2\eta (1-t_3)(1-u_3) )s_{3,\infty } \right\}}{2(1-\eta (1-t_3)(1-u_3))^2 s_{3,\infty }}\nonumber\\
&&- \frac{\eta t_3 u_3 (1+s_{3,\infty })\sqrt{s_{3,\infty }^2-2(1-2\eta  (1-t_3)(1-u_3))s_{3,\infty }+1}}{2(1-\eta (1-t_3)(1-u_3))^2 s_{3,\infty }}\nonumber
\end{eqnarray}
\begin{eqnarray}
\widetilde{w}_{2,3} &=& \frac{v_2}{(v_2-1)}\; \frac{\widetilde{w}_{3,3} t_2 u_2}{1- \widetilde{w}_{3,3} (1-t_2)(1-u_2)v_2 }\Bigg|_{v_2=\frac{1-s_2-\sqrt{(1-s_2)^2+4\widetilde{w}_{3,3} (1-t_2)(1-u_2)s_2}}{2\widetilde{w}_{3,3} (1-t_2)(1-u_2)}}\nonumber\\
&=& \frac{\widetilde{w}_{3,3} t_2 u_2 \left\{ 1+ (s_2+2\widetilde{w}_{3,3}(1-t_2)(1-u_2) )s_2 \right\}}{2(1-\widetilde{w}_{3,3}(1-t_2)(1-u_2))^2 s_2}\nonumber\\
&&- \frac{\widetilde{w}_{3,3} t_2 u_2 (1+s_2)\sqrt{s_2^2-2(1-2\widetilde{w}_{3,3} (1-t_2)(1-u_2))s_2+1}}{2(1-\widetilde{w}_{3,3}(1-t_2)(1-u_2))^2 s_2}\nonumber
\end{eqnarray}
\begin{eqnarray}
\widetilde{w}_{1,3} &=& \frac{v_1}{(v_1-1)}\; \frac{\widetilde{w}_{2,3} t_1 u_1}{1- \widetilde{w}_{2,3} (1-t_1)(1-u_1)v_1 }\Bigg|_{v_1=\frac{1-s_1-\sqrt{(1-s_1)^2+4\widetilde{w}_{2,3} (1-t_1)(1-u_1)s_1}}{2\widetilde{w}_{2,3} (1-t_1)(1-u_1)}}\nonumber\\
&=& \frac{\widetilde{w}_{2,3} t_1 u_1 \left\{ 1+ (s_1+2\widetilde{w}_{2,3}(1-t_1)(1-u_1) )s_1 \right\}}{2(1-\widetilde{w}_{2,3}(1-t_1)(1-u_1))^2 s_1}\nonumber\\
&&- \frac{\widetilde{w}_{2,3} t_1 u_1 (1+s_1)\sqrt{s_1^2-2(1-2\widetilde{w}_{2,3} (1-t_1)(1-u_1))s_1+1} }{2(1-\widetilde{w}_{2,3}(1-t_1)(1-u_1))^2 s_1}\nonumber
\end{eqnarray}
By repeating this process for all higher terms of integral forms of sub-summation $y_m(x)$ terms where $m > 3$, I obtain every  $\displaystyle{ \sum_{q_0 =0}^{\infty } \frac{(\gamma')_{q_0}}{q_0!} s_0^{q_0} \prod _{n=1}^{\infty } \left\{ \sum_{ q_n = q_{n-1}}^{\infty } s_n^{q_n }\right\}}  y_m(x)$ terms. 
Substitute (\ref{eq:20025}), (\ref{eq:20028}), (\ref{eq:20033}), (\ref{eq:20034}) and including all $\displaystyle{ \sum_{q_0 =0}^{\infty } \frac{(\gamma')_{q_0}}{q_0!} s_0^{q_0} \prod _{n=1}^{\infty } \left\{ \sum_{ q_n = q_{n-1}}^{\infty } s_n^{q_n }\right\}}  y_m(x)$ terms where $m > 3$ into (\ref{eq:20080}). 
\qed
\end{proof}
\begin{remark}
The generating function for Heun polynomial of type 2 of the first kind as $q =-(q_j+2j)\{\alpha +\beta -\delta +a(\delta +\gamma -1)+(1+a)(q_j+2j)\} $ where $j,q_j \in \mathbb{N}_{0}$ is
\begin{eqnarray}
&&\sum_{q_0 =0}^{\infty } \frac{(\gamma)_{q_0}}{q_0!} s_0^{q_0} \prod _{n=1}^{\infty } \left\{ \sum_{ q_n = q_{n-1}}^{\infty } s_n^{q_n }\right\} HF_{q_j}^R \Bigg( q_j =\frac{-\varphi \pm \sqrt{\varphi ^2-4(1+a)q}}{2(1+a)}-2j\nonumber\\
&&, \varphi =\alpha +\beta -\delta +a(\delta +\gamma -1), \Omega _1=\frac{\varphi }{(1+a)}; \eta = \frac{(1+a)}{a} x ; z= -\frac{1}{a} x^2 \Bigg) \nonumber\\
&&=2^{\frac{\varphi }{(1+a)}-1}\Bigg\{ \prod_{l=1}^{\infty } \frac{1}{(1-s_{l,\infty })}  \mathbf{A}\left( s_{0,\infty } ;\eta\right) + \Bigg\{  \prod_{l=2}^{\infty } \frac{1}{(1-s_{l,\infty })} \int_{0}^{1} dt_1\;t_1 \int_{0}^{1} du_1\;u_1^{\gamma} \overleftrightarrow {\mathbf{\Gamma}}_1 \left(s_{1,\infty };t_1,u_1,\eta\right)\nonumber\\
&&\times  \widetilde{w}_{1,1}^{-\alpha}\left( \widetilde{w}_{1,1} \partial _{ \widetilde{w}_{1,1}}\right) \widetilde{w}_{1,1}^{\alpha -\beta } \left( \widetilde{w}_{1,1} \partial _{ \widetilde{w}_{1,1}}\right)\widetilde{w}_{1,1}^{\beta } \mathbf{A}\left( s_{0} ;\widetilde{w}_{1,1}\right)\Bigg\} z \nonumber\\
&&+ \sum_{n=2}^{\infty } \Bigg\{ \prod_{l=n+1}^{\infty } \frac{1}{(1-s_{l,\infty })} \int_{0}^{1} dt_n\;t_n^{2n-1} \int_{0}^{1} du_n\;u_n^{2(n-1)+\gamma } \overleftrightarrow {\mathbf{\Gamma}}_n \left(s_{n,\infty };t_n,u_n,\eta \right)\nonumber\\
&&\times \widetilde{w}_{n,n}^{-(2(n-1)+\alpha)}\left( \widetilde{w}_{n,n} \partial _{ \widetilde{w}_{n,n}}\right) \widetilde{w}_{n,n}^{\alpha -\beta } \left( \widetilde{w}_{n,n} \partial _{ \widetilde{w}_{n,n}}\right)\widetilde{w}_{n,n}^{2(n-1)+\beta } \label{eq:2009}\\
&&\times \prod_{k=1}^{n-1} \Bigg\{ \int_{0}^{1} dt_{n-k}\;t_{n-k}^{2(n-k)-1} \int_{0}^{1} du_{n-k} \;u_{n-k}^{2(n-k-1)+\gamma}\overleftrightarrow {\mathbf{\Gamma}}_{n-k} \left(s_{n-k};t_{n-k},u_{n-k},\widetilde{w}_{n-k+1,n} \right)\nonumber\\
&&\times \widetilde{w}_{n-k,n}^{-(2(n-k-1)+\alpha)}\left( \widetilde{w}_{n-k,n} \partial _{ \widetilde{w}_{n-k,n}}\right) \widetilde{w}_{n-k,n}^{\alpha -\beta } \left( \widetilde{w}_{n-k,n} \partial _{ \widetilde{w}_{n-k,n}}\right)\widetilde{w}_{n-k,n}^{2(n-k-1)+\beta }  \Bigg\} \mathbf{A} \left( s_{0} ;\widetilde{w}_{1,n}\right) \Bigg\} z^n \Bigg\} \nonumber   
\end{eqnarray}
where
\begin{equation}
\begin{cases} 
{ \displaystyle \overleftrightarrow {\mathbf{\Gamma}}_1 \left(s_{1,\infty };t_1,u_1,\eta\right)= \frac{\left( \frac{1+s_{1,\infty }+\sqrt{s_{1,\infty }^2-2(1-2\eta (1-t_1)(1-u_1))s_{1,\infty }+1}}{2}\right)^{-\left(3+\Omega _1\right)}}{\sqrt{s_{1,\infty }^2-2(1-2\eta (1-t_1)(1-u_1))s_{1,\infty }+1}}}\cr
{ \displaystyle  \overleftrightarrow {\mathbf{\Gamma}}_n \left(s_{n,\infty };t_n,u_n,\eta \right) =\frac{\left( \frac{1+s_{n,\infty }+\sqrt{s_{n,\infty }^2-2(1-2\eta (1-t_n)(1-u_n))s_{n,\infty }+1}}{2}\right)^{-\left(4n-1+\Omega _1\right)}}{\sqrt{ s_{n,\infty }^2-2(1-2\eta (1-t_n)(1-u_n))s_{n,\infty }+1}}}\cr
{ \displaystyle \overleftrightarrow {\mathbf{\Gamma}}_{n-k} \left(s_{n-k};t_{n-k},u_{n-k},\widetilde{w}_{n-k+1,n} \right)}\cr
{ \displaystyle = \frac{ \left( \frac{1+s_{n-k}+\sqrt{s_{n-k}^2-2(1-2\widetilde{w}_{n-k+1,n} (1-t_{n-k})(1-u_{n-k}))s_{n-k}+1}}{2}\right)^{-\left(4(n-k)-1+\Omega _1\right)}}{\sqrt{ s_{n-k}^2-2(1-2\widetilde{w}_{n-k+1,n} (1-t_{n-k})(1-u_{n-k}))s_{n-k}+1}}}
\end{cases}\nonumber 
\end{equation}
and
\begin{equation}
\begin{cases} 
{ \displaystyle \mathbf{A} \left( s_{0,\infty } ;\eta\right)= \frac{\left(1- s_{0,\infty }+\sqrt{s_{0,\infty }^2-2(1-2\eta )s_{0,\infty }+1}\right)^{1-\gamma } }{\left(1+s_{0,\infty }+\sqrt{s_{0,\infty }^2-2(1-2\eta )s_{0,\infty }+1}\right)^{\Omega _1 -\gamma }\sqrt{s_{0,\infty }^2-2(1-2\eta )s_{0,\infty }+1}}}\cr
{ \displaystyle  \mathbf{A} \left( s_{0} ;\widetilde{w}_{1,1}\right) = \frac{\left(1- s_0+\sqrt{s_0^2-2(1-2\widetilde{w}_{1,1})s_0+1}\right)^{1-\gamma}}{ \left(1+s_0+\sqrt{s_0^2-2(1-2\widetilde{w}_{1,1} )s_0+1}\right)^{\Omega _1 -\gamma}\sqrt{s_0^2-2(1-2\widetilde{w}_{1,1})s_0+1}}} \cr
{ \displaystyle \mathbf{A} \left( s_{0} ;\widetilde{w}_{1,n}\right) = \frac{\left(1- s_0+\sqrt{s_0^2-2(1-2\widetilde{w}_{1,n})s_0+1}\right)^{1-\gamma} }{\left(1+s_0+\sqrt{s_0^2-2(1-2\widetilde{w}_{1,n} )s_0+1}\right)^{\Omega _1 -\gamma}\sqrt{s_0^2-2(1-2\widetilde{w}_{1,n})s_0+1}}}
\end{cases}\nonumber 
\end{equation}
\end{remark}
\begin{proof}
Replace A, $w$ and $x$  by $\displaystyle{\frac{\varphi }{(1+a)}}$, $s_{0,\infty }$ and $\eta $ in (\ref{eq:2003}). 
\begin{eqnarray}
&&\sum_{q_0=0}^{\infty }\frac{(\gamma )_{q_0}}{q_0!} s_{0,\infty }^{q_0} \;_2F_1\left(-q_0, q_0+\frac{\varphi }{(1+a)}; \gamma ; \eta \right) \label{eq:20035}\\
&&= 2^{\frac{\varphi }{(1+a)}-1}\frac{\left(1- s_{0,\infty }+\sqrt{s_{0,\infty }^2-2(1-2\eta )s_{0,\infty }+1}\right)^{1-\gamma } }{\left(1+s_{0,\infty }+\sqrt{s_{0,\infty }^2-2(1-2\eta )s_{0,\infty }+1}\right)^{\frac{\varphi }{(1+a)}-\gamma }\sqrt{s_{0,\infty }^2-2(1-2\eta )s_{0,\infty }+1}} \nonumber
\end{eqnarray} 
Replace A, $w$ and $x$  by $\displaystyle{\frac{\varphi }{(1+a)}}$, $s_0$ and $\widetilde{w}_{1,1}$ in (\ref{eq:2003}). 
\begin{eqnarray}
&&\sum_{q_0=0}^{\infty }\frac{(\gamma )_{q_0}}{q_0!} s_0^{q_0} \;_2F_1\left(-q_0, q_0+\frac{\varphi }{(1+a)}; \gamma; \widetilde{w}_{1,1} \right) \label{eq:20036}\\
&&= 2^{\frac{\varphi }{(1+a)}-1}\frac{\left(1- s_0+\sqrt{s_0^2-2(1-2\widetilde{w}_{1,1})s_0+1}\right)^{1-\gamma} }{\left(1+s_0+\sqrt{s_0^2-2(1-2\widetilde{w}_{1,1} )s_0+1}\right)^{\frac{\varphi }{(1+a)}-\gamma }\sqrt{s_0^2-2(1-2\widetilde{w}_{1,1})s_0+1}} \nonumber
\end{eqnarray} 
Replace A, $w$ and $x$  by $\displaystyle{\frac{\varphi }{(1+a)}}$, $s_0$ and $\widetilde{w}_{1,n}$  in (\ref{eq:2003}). 
\begin{eqnarray}
&&\sum_{q_0=0}^{\infty }\frac{(\gamma )_{q_0}}{q_0!} s_0^{q_0} \;_2F_1\left(-q_0, q_0+\frac{\varphi }{(1+a)}; \gamma; \widetilde{w}_{1,n} \right) \label{eq:20037}\\
&&= 2^{\frac{\varphi }{(1+a)}-1}\frac{\left(1- s_0+\sqrt{s_0^2-2(1-2\widetilde{w}_{1,n})s_0+1}\right)^{1-\gamma} }{\left(1+s_0+\sqrt{s_0^2-2(1-2\widetilde{w}_{1,n} )s_0+1}\right)^{\frac{\varphi }{(1+a)}-\gamma }\sqrt{s_0^2-2(1-2\widetilde{w}_{1,n})s_0+1}} \nonumber
\end{eqnarray} 
Put $c_0$= 1, $\lambda $=0 and $\gamma' = \gamma $ in (\ref{eq:2008}). Substitute (\ref{eq:20035}), (\ref{eq:20036}) and (\ref{eq:20037}) into the new (\ref{eq:2008}).
\qed
\end{proof}
\begin{remark}
The generating function for Heun polynomial of type 2 of the second kind as $q =-(q_j+2j+1-\gamma )\{\alpha +\beta +1-\gamma -(1-a)\delta +(1+a)(q_j+2j)\} $ where $j, q_j \in \mathbb{N}_{0}$ is
\begin{eqnarray}
&&\sum_{q_0 =0}^{\infty } \frac{(2-\gamma )_{q_0}}{q_0!} s_0^{q_0} \prod _{n=1}^{\infty } \left\{ \sum_{ q_n = q_{n-1}}^{\infty } s_n^{q_n }\right\} HS_{q_j}^R \Bigg( q_j =\frac{-\{\varphi +2(1+a)(1-\gamma )\} \pm \sqrt{\varphi ^2-4(1+a)q}}{2(1+a)}-2j\nonumber\\
&&, \varphi =\alpha +\beta -\delta +a(\delta +\gamma -1), \Omega _2= \frac{\varphi +2(1+a)(1-\gamma )}{(1+a)}; \eta = \frac{(1+a)}{a} x ; z= -\frac{1}{a} x^2 \Bigg) \nonumber\\
&&= 2^{\frac{\varphi +2(1+a)(1/2-\gamma )}{(1+a)}}\eta^{1-\gamma }\left\{ \prod_{l=1}^{\infty } \frac{1}{(1-s_{l,\infty })} \mathbf{B}\left( s_{0,\infty } ;\eta\right) \right. \nonumber\\
&&+ \Bigg\{ \prod_{l=2}^{\infty } \frac{1}{(1-s_{l,\infty })} \int_{0}^{1} dt_1\;t_1^{2-\gamma } \int_{0}^{1} du_1\;u_1 \overleftrightarrow {\mathbf{\Psi}}_1 \left(s_{1,\infty };t_1,u_1,\eta\right)\nonumber\\
&& \times \widetilde{w}_{1,1}^{-(\alpha-\gamma +1)}\left( \widetilde{w}_{1,1} \partial _{ \widetilde{w}_{1,1}}\right) \widetilde{w}_{1,1}^{\alpha -\beta } \left( \widetilde{w}_{1,1} \partial _{ \widetilde{w}_{1,1}}\right)\widetilde{w}_{1,1}^{\beta -\gamma +1} \mathbf{B}\left( s_{0} ;\widetilde{w}_{1,1}\right) \Bigg\} z \nonumber\\
&&+ \sum_{n=2}^{\infty } \Bigg\{ \prod_{l=n+1}^{\infty } \frac{1}{(1-s_{l,\infty })} \int_{0}^{1} dt_n\;t_n^{2n-\gamma } \int_{0}^{1} du_n\;u_n^{2n-1} \overleftrightarrow {\mathbf{\Psi}}_n \left(s_{n,\infty };t_n,u_n,\eta \right)\nonumber\\
&&\times \widetilde{w}_{n,n}^{-(2n-1+\alpha-\gamma )}\left( \widetilde{w}_{n,n} \partial _{ \widetilde{w}_{n,n}}\right) \widetilde{w}_{n,n}^{\alpha -\beta } \left( \widetilde{w}_{n,n} \partial _{ \widetilde{w}_{n,n}}\right)\widetilde{w}_{n,n}^{2n-1+\beta -\gamma } \nonumber\\
&&\times \prod_{k=1}^{n-1} \Bigg\{  \int_{0}^{1} dt_{n-k}\;t_{n-k}^{2(n-k)-\gamma} \int_{0}^{1} du_{n-k} \;u_{n-k}^{2(n-k)-1} \overleftrightarrow {\mathbf{\Psi}}_{n-k} \left(s_{n-k};t_{n-k},u_{n-k},\widetilde{w}_{n-k+1,n} \right) \nonumber\\
&&\times \widetilde{w}_{n-k,n}^{-(2(n-k)-1+\alpha-\gamma )}\left( \widetilde{w}_{n-k,n} \partial _{ \widetilde{w}_{n-k,n}}\right) \widetilde{w}_{n-k,n}^{\alpha -\beta } \left( \widetilde{w}_{n-k,n} \partial _{ \widetilde{w}_{n-k,n}}\right)\widetilde{w}_{n-k,n}^{2(n-k)-1+\beta -\gamma } \Bigg\}\nonumber\\
&&\times \left. \mathbf{B}\left( s_{0} ;\widetilde{w}_{1,n}\right)\Bigg\} z^n  \right\}
 \label{eq:20010}
\end{eqnarray}
where
\begin{equation}
\begin{cases} 
{ \displaystyle \overleftrightarrow {\mathbf{\Psi}}_1 \left(s_{1,\infty };t_1,u_1,\eta\right)= \frac{\left( \frac{1+s_{1,\infty }+\sqrt{s_{1,\infty }^2-2(1-2\eta (1-t_1)(1-u_1))s_{1,\infty }+1}}{2}\right)^{-\left( 3+\Omega _2\right)}}{\sqrt{s_{1,\infty }^2-2(1-2\eta (1-t_1)(1-u_1))s_{1,\infty }+1}} }\cr
{ \displaystyle  \overleftrightarrow {\mathbf{\Psi}}_n \left(s_{n,\infty };t_n,u_n,\eta \right) = \frac{\left( \frac{1+s_{n,\infty }+\sqrt{s_{n,\infty }^2-2(1-2\eta (1-t_n)(1-u_n))s_{n,\infty }+1}}{2}\right)^{-\left(4n-1+\Omega _2\right)}}{\sqrt{s_{n,\infty }^2-2(1-2\eta (1-t_n)(1-u_n))s_{n,\infty }+1}}}\cr
{ \displaystyle \overleftrightarrow {\mathbf{\Psi}}_{n-k} \left(s_{n-k};t_{n-k},u_{n-k},\widetilde{w}_{n-k+1,n} \right)}\cr
{ \displaystyle = \frac{\left( \frac{(1+s_{n-k})+\sqrt{s_{n-k}^2-2(1-2\widetilde{w}_{n-k+1,n} (1-t_{n-k})(1-u_{n-k}))s_{n-k}+1}}{2}\right)^{-\left( 4(n-k)-1+\Omega _2\right)}}{\sqrt{s_{n-k}^2-2(1-2\widetilde{w}_{n-k+1,n} (1-t_{n-k})(1-u_{n-k}))s_{n-k}+1}}}
\end{cases}\nonumber 
\end{equation}
and
\begin{equation}
\begin{cases} 
{ \displaystyle \mathbf{B} \left( s_{0,\infty } ;\eta\right)= \frac{ \left(1+s_{0,\infty }+\sqrt{s_{0,\infty }^2-2(1-2\eta )s_{0,\infty }+1}\right)^{2-\gamma -\Omega _2}}{\left(1- s_{0,\infty }+\sqrt{s_{0,\infty }^2-2(1-2\eta )s_{0,\infty }+1}\right)^{1-\gamma}\sqrt{s_{0,\infty }^2-2(1-2\eta )s_{0,\infty }+1}}}\cr
{ \displaystyle  \mathbf{B} \left( s_{0} ;\widetilde{w}_{1,1}\right) = \frac{\left(1+s_0+\sqrt{s_0^2-2(1-2\widetilde{w}_{1,1} )s_0+1}\right)^{2-\gamma -\Omega _2}}{\left(1- s_0+\sqrt{s_0^2-2(1-2\widetilde{w}_{1,1})s_0+1}\right)^{1-\gamma}\sqrt{s_0^2-2(1-2\widetilde{w}_{1,1})s_0+1}}} \cr
{ \displaystyle \mathbf{B} \left( s_{0} ;\widetilde{w}_{1,n}\right) = \frac{\left(1+s_0+\sqrt{s_0^2-2(1-2\widetilde{w}_{1,n} )s_0+1}\right)^{2-\gamma -\Omega _2}}{\left(1- s_0+\sqrt{s_0^2-2(1-2\widetilde{w}_{1,n})s_0+1}\right)^{1-\gamma}\sqrt{s_0^2-2(1-2\widetilde{w}_{1,n})s_0+1}}}
\end{cases}\nonumber 
\end{equation}
\end{remark}
\begin{proof}
Replace $\gamma $, A, $w$ and $x$  by $2-\gamma $, $\displaystyle {\frac{\varphi +2(1+a)(1-\gamma )}{(1+a)}}$, $s_{0,\infty }$ and $\eta $ in (\ref{eq:2003}). 
\begin{eqnarray}
&&\sum_{q_0=0}^{\infty }\frac{(2-\gamma)_{q_0}}{q_0!} s_{0,\infty }^{q_0} \;_2F_1\left(-q_0, q_0+\frac{\varphi +2(1+a)(1-\gamma )}{(1+a)}; 2-\gamma; \eta \right) \label{eq:20038}\\
&&=2^{\frac{\varphi +2(1+a)(1/2-\gamma )}{(1+a)}}\frac{ \left(1+s_{0,\infty }+\sqrt{s_{0,\infty }^2-2(1-2\eta )s_{0,\infty }+1}\right)^{-\frac{\varphi -(1+a)\gamma }{(1+a)}}}{\left(1- s_{0,\infty }+\sqrt{s_{0,\infty }^2-2(1-2\eta )s_{0,\infty }+1}\right)^{1-\gamma}\sqrt{s_{0,\infty }^2-2(1-2\eta )s_{0,\infty }+1}} \nonumber
\end{eqnarray} 
Replace $\gamma $, A, $w$ and $x$  by $2-\gamma $, $\displaystyle {\frac{\varphi +2(1+a)(1-\gamma )}{(1+a)}}$, $s_0$ and $\widetilde{w}_{1,1}$ in (\ref{eq:2003}). 
\begin{eqnarray}
&&\sum_{q_0=0}^{\infty }\frac{(2-\gamma)_{q_0}}{q_0!} s_0^{q_0} \;_2F_1\left(-q_0, q_0+\frac{\varphi +2(1+a)(1-\gamma )}{(1+a)}; 2-\gamma; \widetilde{w}_{1,1} \right) \label{eq:20039}\\
&&=2^{\frac{\varphi +2(1+a)(1/2-\gamma )}{(1+a)}}\frac{ \left(1+s_0+\sqrt{s_0^2-2(1-2\widetilde{w}_{1,1} )s_0+1}\right)^{-\frac{\varphi -(1+a)\gamma }{(1+a)}}}{\left(1- s_0+\sqrt{s_0^2-2(1-2\widetilde{w}_{1,1})s_0+1}\right)^{1-\gamma }\sqrt{s_0^2-2(1-2\widetilde{w}_{1,1})s_0+1}} \nonumber
\end{eqnarray} 
Replace $\gamma $, A, $w$ and $x$  by $2-\gamma $, $\displaystyle { \frac{\varphi +2(1+a)(1-\gamma )}{(1+a)}}$, $s_0$ and $\widetilde{w}_{1,n}$ in (\ref{eq:2003}). 
\begin{eqnarray}
&&\sum_{q_0=0}^{\infty }\frac{(2-\gamma)_{q_0}}{q_0!} s_0^{q_0} \;_2F_1\left(-q_0, q_0+\frac{\varphi +2(1+a)(1-\gamma )}{(1+a)}; 2-\gamma; \widetilde{w}_{1,n} \right) \label{eq:20040}\\
&&= 2^{\frac{\varphi +2(1+a)(1/2-\gamma )}{(1+a)}}\frac{ \left(1+s_0+\sqrt{s_0^2-2(1-2\widetilde{w}_{1,n} )s_0+1}\right)^{-\frac{\varphi -(1+a)\gamma }{(1+a)}}}{\left(1- s_0+\sqrt{s_0^2-2(1-2\widetilde{w}_{1,n})s_0+1}\right)^{1-\gamma }\sqrt{s_0^2-2(1-2\widetilde{w}_{1,n})s_0+1}} \nonumber
\end{eqnarray} 
Put $\displaystyle {c_0= \left( \frac{1+a}{a}\right)^{1-\gamma }}$, $\lambda =1-\gamma $ and $\gamma' = 2-\gamma $ in (\ref{eq:2008}). Substitute (\ref{eq:20038}), (\ref{eq:20039}) and (\ref{eq:20040}) into the new (\ref{eq:2008}).
\qed
\end{proof}

\section{Summary} 

As we see all solutions of power series expansions in Heun function by using 3TRF, denominators and numerators in all $B_n$ terms arise with Pochhammer symbol.\cite{Chou2012H12} Also the Frobenius solutions in Heun funtion by using R3TRF, denominators and numerators in all $A_n$ terms arise with Pochhammer symbol: the meaning of this is that the analytic solutions of Heun's ordinary differential equations with three recursive coefficients can be described as Hypergoemetric functions in a strict mathematical way. 

We can express representations in closed form integrals of Heun function since we have power series expansions with Pochhammer symbols in numerators and denominators: (1) look at Remark 1 through Remark 6 in Ref.\cite{Chou2012H22}.  (2) look at Remark 2.3.2 through Remark 2.3.5 in chapter 2. We can transform any special functions into all other well-known special functions such as Legendre, Bessel, Laguerre, Kummer functions and etc because a $_2F_1$ function recurs in each of sub-integral forms of Heun function. After we replace $_2F_1$ functions in integral forms of Heun function to other special functions, we can rebuild the power series expansion of Heun function and its linear ordinary differential equation in a backward: understanding the connection between other special functions is important in the mathematical and physical points of views as we all know.

Since integral representation of Heun function is derived from power series expansion in closed forms, I construct generating function for Heun polynomial of type 2 by applying the generating function for Jacobi polynomial using hypergeometric functions into the general integral representation of Heun polynomial. The generating function is really helpful in order to derive orthogonal relations, recursion relations and expectation values of any physical quantities as we all recognize; i.e. the normalized wave function of hydrogen-like atoms and expectation values of its physical quantities such as position and momentum.

\begin{appendices}
\section[Generating functions for 192 Heun polynomials of type 2]{Generating functions for 192 Heun polynomials of type 2
  \sectionmark{Generating functions for 192 Heun polynomials of type 2}}
  \sectionmark{Generating functions for 192 Heun polynomials of type 2}
In appendix of chapter 2, by applying R3TRF, I construct the fundamental power series expansions in closed forms and integral forms of Heun function (infinite series and polynomial of type 2) of nine out of the 192 local solution of Heun function in Table 2, obtained by Robert S. Maier(2007) \cite{Maie20072}.  

In this appendix, replacing coefficients in the general expression of generating function for Heun polynomial of type 2, I derive generating functions for Heun polynomial of type 2 for the previous nine examples of the 192 local solution of Heun equation \cite{Maie20072}.\footnote{I treat $\alpha $, $\beta $, $\gamma $ and $\delta $ as free variables and a fixed value of $q$ to construct the generating functions for Heun polynomial of type 2 for all nine examples of the 192 local solution of Heun function.}
\addtocontents{toc}{\protect\setcounter{tocdepth}{1}}
\subsection{ ${\displaystyle (1-x)^{1-\delta } Hl(a, q - (\delta  - 1)\gamma a; \alpha - \delta  + 1, \beta - \delta + 1, \gamma ,2 - \delta ; x)}$ }
Replacing coefficients $q$, $\alpha$, $\beta$ and $\delta$ by $q - (\delta - 1)\gamma a $, $\alpha - \delta  + 1 $, $\beta - \delta + 1$ and $2 - \delta$ into (\ref{eq:2009}). Multiply $(1-x)^{1-\delta }$ and the new (\ref{eq:2009}) together.
\begin{eqnarray}
&&\sum_{q_0 =0}^{\infty } \frac{(\gamma)_{q_0}}{q_0!} s_0^{q_0} \prod _{n=1}^{\infty } \left\{ \sum_{ q_n = q_{n-1}}^{\infty } s_n^{q_n }\right\} (1-x)^{1-\delta } Hl(a, q - (\delta  - 1)\gamma a; \alpha - \delta  + 1, \beta - \delta + 1, \gamma ,2 - \delta ; x) \nonumber\\
&&=\frac{2^{\frac{\varphi }{(1+a)}-1}}{(1-x)^{\delta -1}}\left\{ \prod_{l=1}^{\infty } \frac{1}{(1-s_{l,\infty })} \mathbf{A} \left( s_{0,\infty } ;\eta\right)\right.  \nonumber\\
&&+ \Bigg\{ \prod_{l=2}^{\infty } \frac{1}{(1-s_{l,\infty })} \int_{0}^{1} dt_1\;t_1 \int_{0}^{1} du_1\;u_1^{\gamma} \overleftrightarrow {\mathbf{\Gamma}}_1 \left(s_{1,\infty };t_1,u_1,\eta\right)\nonumber\\
&&\times \widetilde{w}_{1,1}^{-\alpha+\delta -1}\left( \widetilde{w}_{1,1} \partial _{ \widetilde{w}_{1,1}}\right) \widetilde{w}_{1,1}^{\alpha -\beta } \left( \widetilde{w}_{1,1} \partial _{ \widetilde{w}_{1,1}}\right)\widetilde{w}_{1,1}^{\beta -\delta +1}  \mathbf{A} \left( s_{0} ;\widetilde{w}_{1,1}\right) \Bigg\}z \nonumber\\
&&+ \sum_{n=2}^{\infty } \left\{ \prod_{l=n+1}^{\infty } \frac{1}{(1-s_{l,\infty })} \int_{0}^{1} dt_n\;t_n^{2n-1} \int_{0}^{1} du_n\;u_n^{2(n-1)+\gamma }\right.  \overleftrightarrow {\mathbf{\Gamma}}_n \left(s_{n,\infty };t_n,u_n,\eta \right)\nonumber\\
&&\times \widetilde{w}_{n,n}^{-(2n-1+\alpha-\delta )}\left( \widetilde{w}_{n,n} \partial _{ \widetilde{w}_{n,n}}\right) \widetilde{w}_{n,n}^{\alpha -\beta } \left( \widetilde{w}_{n,n} \partial _{ \widetilde{w}_{n,n}}\right)\widetilde{w}_{n,n}^{2n-1+\beta -\delta } \nonumber\\
&&\times \prod_{k=1}^{n-1} \Bigg\{  \int_{0}^{1} dt_{n-k}\;t_{n-k}^{2(n-k)-1} \int_{0}^{1} du_{n-k} \;u_{n-k}^{2(n-k-1)+\gamma}  \overleftrightarrow {\mathbf{\Gamma}}_{n-k} \left(s_{n-k};t_{n-k},u_{n-k},\widetilde{w}_{n-k+1,n} \right)  \nonumber\\
&&\times \widetilde{w}_{n-k,n}^{-(2(n-k)-1+\alpha-\delta )}\left( \widetilde{w}_{n-k,n} \partial _{ \widetilde{w}_{n-k,n}}\right) \widetilde{w}_{n-k,n}^{\alpha -\beta } \left( \widetilde{w}_{n-k,n} \partial _{ \widetilde{w}_{n-k,n}}\right)\widetilde{w}_{n-k,n}^{2(n-k)-1+\beta-\delta}  \Bigg\}\nonumber\\
&&\times \left.\mathbf{A} \left( s_{0} ;\widetilde{w}_{1,n}\right) \Bigg\} z^n \right\}   \label{eq:200200}
\end{eqnarray}
where
\begin{equation}
\begin{cases} z = -\frac{1}{a}x^2 \cr
\eta = \frac{(1+a)}{a} x \cr
\varphi = \alpha +\beta -\delta +a(\gamma -\delta +1) \cr
\Omega_1 = \frac{\varphi }{(1+a)} \cr
q = (\delta - 1)\gamma a-(q_j+2j)\{\varphi +(1+a)(q_j+2j) \} \;\;\mbox{as}\;j,q_j\in \mathbb{N}_{0} \
\end{cases}\nonumber 
\end{equation}
\begin{equation}
\begin{cases} 
{ \displaystyle \overleftrightarrow {\mathbf{\Gamma}}_1 \left(s_{1,\infty };t_1,u_1,\eta\right)= \frac{\left( \frac{1+s_{1,\infty }+\sqrt{s_{1,\infty }^2-2(1-2\eta (1-t_1)(1-u_1))s_{1,\infty }+1}}{2}\right)^{-\left(3+\Omega_1\right)}}{\sqrt{s_{1,\infty }^2-2(1-2\eta (1-t_1)(1-u_1))s_{1,\infty }+1}}}\cr
{ \displaystyle  \overleftrightarrow {\mathbf{\Gamma}}_n \left(s_{n,\infty };t_n,u_n,\eta \right) =\frac{ \left( \frac{1+s_{n,\infty }+\sqrt{s_{n,\infty }^2-2(1-2\eta (1-t_n)(1-u_n))s_{n,\infty }+1}}{2}\right)^{-\left(4n-1+\Omega_1\right)} }{\sqrt{s_{n,\infty }^2-2(1-2\eta (1-t_n)(1-u_n))s_{n,\infty }+1}}}\cr
{ \displaystyle \overleftrightarrow {\mathbf{\Gamma}}_{n-k} \left(s_{n-k};t_{n-k},u_{n-k},\widetilde{w}_{n-k+1,n} \right)}\cr
{ \displaystyle = \frac{ \left( \frac{1+s_{n-k}+\sqrt{s_{n-k}^2-2(1-2\widetilde{w}_{n-k+1,n} (1-t_{n-k})(1-u_{n-k}))s_{n-k}+1}}{2}\right)^{-\left(4(n-k)-1+\Omega_1\right)} }{\sqrt{s_{n-k}^2-2(1-2\widetilde{w}_{n-k+1,n} (1-t_{n-k})(1-u_{n-k}))s_{n-k}+1}} }
\end{cases}\nonumber 
\end{equation}
  and
 \begin{equation}
\begin{cases} 
{ \displaystyle \mathbf{A} \left( s_{0,\infty } ;\eta\right)= \frac{\left(1+s_{0,\infty }+\sqrt{s_{0,\infty }^2-2(1-2\eta )s_{0,\infty }+1}\right)^{\gamma-\Omega_1}}{\left(1- s_{0,\infty }+\sqrt{s_{0,\infty }^2-2(1-2\eta )s_{0,\infty }+1}\right)^{\gamma -1}  \sqrt{s_{0,\infty }^2-2(1-2\eta )s_{0,\infty }+1}}}\cr
{ \displaystyle  \mathbf{A} \left( s_{0} ;\widetilde{w}_{1,1}\right) = \frac{ \left(1+s_0+\sqrt{s_0^2-2(1-2\widetilde{w}_{1,1} )s_0+1}\right)^{\gamma -\Omega_1}}{\left(1- s_0+\sqrt{s_0^2-2(1-2\widetilde{w}_{1,1})s_0+1}\right)^{\gamma -1} \sqrt{s_0^2-2(1-2\widetilde{w}_{1,1})s_0+1}}} \cr
{ \displaystyle \mathbf{A} \left( s_{0} ;\widetilde{w}_{1,n}\right) = \frac{ \left(1+s_0+\sqrt{s_0^2-2(1-2\widetilde{w}_{1,n} )s_0+1}\right)^{\gamma -\Omega_1}}{\left(1- s_0+\sqrt{s_0^2-2(1-2\widetilde{w}_{1,n})s_0+1}\right)^{\gamma -1} \sqrt{s_0^2-2(1-2\widetilde{w}_{1,n})s_0+1}}}
\end{cases}\nonumber 
\end{equation}
 \subsection{\footnotesize ${\displaystyle x^{1-\gamma } (1-x)^{1-\delta } Hl(a, q-(\gamma +\delta -2)a-(\gamma -1)(\alpha +\beta -\gamma -\delta +1); \alpha - \gamma -\delta +2, \beta - \gamma -\delta +2, 2-\gamma, 2 - \delta ; x)}$ \normalsize}
Replacing coefficients $q$, $\alpha$, $\beta$, $\gamma $ and $\delta$ by $q-(\gamma +\delta -2)a-(\gamma -1)(\alpha +\beta -\gamma -\delta +1)$, $\alpha - \gamma -\delta +2$, $\beta - \gamma -\delta +2, 2-\gamma$ and $2 - \delta$ into (\ref{eq:2009}). Multiply $x^{1-\gamma } (1-x)^{1-\delta }$ and the new (\ref{eq:2009}) together.
\begin{eqnarray}
&& \sum_{q_0 =0}^{\infty } \frac{(\gamma)_{q_0}}{q_0!} s_0^{q_0} \prod _{n=1}^{\infty } \left\{ \sum_{ q_n = q_{n-1}}^{\infty } s_n^{q_n }\right\} x^{1-\gamma } (1-x)^{1-\delta } Hl(a, q-(\gamma +\delta -2)a\nonumber\\
&&\hspace{4mm}-(\gamma -1)(\alpha +\beta -\gamma -\delta +1); \alpha - \gamma -\delta +2, \beta - \gamma -\delta +2, 2-\gamma, 2 - \delta ; x)  \nonumber\\
&& =\frac{2^{\frac{\varphi }{(1+a)}-1}}{x^{\gamma -1} (1-x)^{\delta -1}}\left\{ \prod_{l=1}^{\infty } \frac{1}{(1-s_{l,\infty })} \mathbf{A} \left( s_{0,\infty } ;\eta\right)\right.  \nonumber\\
&&+ \Bigg\{\prod_{l=2}^{\infty } \frac{1}{(1-s_{l,\infty })} \int_{0}^{1} dt_1\;t_1 \int_{0}^{1} du_1\;u_1^{2-\gamma} \overleftrightarrow {\mathbf{\Gamma}}_1 \left(s_{1,\infty };t_1,u_1,\eta\right)\nonumber\\
&&\times \widetilde{w}_{1,1}^{-\alpha+\gamma +\delta -2}\left( \widetilde{w}_{1,1} \partial _{ \widetilde{w}_{1,1}}\right) \widetilde{w}_{1,1}^{\alpha -\beta } \left( \widetilde{w}_{1,1} \partial _{ \widetilde{w}_{1,1}}\right)\widetilde{w}_{1,1}^{\beta -\gamma -\delta +2}  \mathbf{A} \left( s_{0} ;\widetilde{w}_{1,1}\right) \Bigg\}z \nonumber\\
&&+ \sum_{n=2}^{\infty } \left\{ \prod_{l=n+1}^{\infty } \frac{1}{(1-s_{l,\infty })} \int_{0}^{1} dt_n\;t_n^{2n-1} \int_{0}^{1} du_n\;u_n^{2n-\gamma } \overleftrightarrow {\mathbf{\Gamma}}_n \left(s_{n,\infty };t_n,u_n,\eta \right)\right.\nonumber\\
&&\times  \widetilde{w}_{n,n}^{-(2n+\alpha-\gamma -\delta)}\left( \widetilde{w}_{n,n} \partial _{ \widetilde{w}_{n,n}}\right) \widetilde{w}_{n,n}^{\alpha -\beta } \left( \widetilde{w}_{n,n} \partial _{ \widetilde{w}_{n,n}}\right)\widetilde{w}_{n,n}^{2n+\beta-\gamma -\delta }  \nonumber\\
&&\times \prod_{k=1}^{n-1} \Bigg\{  \int_{0}^{1} dt_{n-k}\;t_{n-k}^{2(n-k)-1} \int_{0}^{1} du_{n-k} \;u_{n-k}^{2(n-k)-\gamma}  \overleftrightarrow {\mathbf{\Gamma}}_{n-k} \left(s_{n-k};t_{n-k},u_{n-k},\widetilde{w}_{n-k+1,n} \right)  \nonumber\\
&&\times \widetilde{w}_{n-k,n}^{-(2(n-k)+\alpha-\gamma -\delta )}\left( \widetilde{w}_{n-k,n} \partial _{ \widetilde{w}_{n-k,n}}\right) \widetilde{w}_{n-k,n}^{\alpha -\beta } \left( \widetilde{w}_{n-k,n} \partial _{ \widetilde{w}_{n-k,n}}\right)\widetilde{w}_{n-k,n}^{2(n-k)+\beta -\gamma -\delta } \Bigg\}\nonumber\\
&&\times  \mathbf{A} \left( s_{0} ;\widetilde{w}_{1,n}\right) \Bigg\} z^n \Bigg\}  
 \label{eq:200201}
\end{eqnarray}
where
\begin{equation}
\begin{cases} z = -\frac{1}{a}x^2 \cr
\eta = \frac{(1+a)}{a} x \cr
\varphi = \alpha +\beta -2\gamma -\delta +2 +a(3-\gamma -\delta ) \cr
\Omega_1 = \frac{\varphi }{(1+a)}\cr
q = (\gamma +\delta -2)a+(\gamma -1)(\alpha +\beta -\gamma -\delta +1)-(q_j+2j)\{\varphi +(1+a)(q_j+2j) \} \;\;\mbox{as}\;j,q_j\in \mathbb{N}_{0} 
\end{cases}\nonumber 
\end{equation}
\begin{equation}
\begin{cases} 
{ \displaystyle \overleftrightarrow {\mathbf{\Gamma}}_1 \left(s_{1,\infty };t_1,u_1,\eta\right)= \frac{\left( \frac{1+s_{1,\infty }+\sqrt{s_{1,\infty }^2-2(1-2\eta (1-t_1)(1-u_1))s_{1,\infty }+1}}{2}\right)^{-(3+\Omega_1 )}}{\sqrt{s_{1,\infty }^2-2(1-2\eta (1-t_1)(1-u_1))s_{1,\infty }+1}}}\cr
{ \displaystyle  \overleftrightarrow {\mathbf{\Gamma}}_n \left(s_{n,\infty };t_n,u_n,\eta \right) = \frac{ \left( \frac{1+s_{n,\infty }+\sqrt{s_{n,\infty }^2-2(1-2\eta (1-t_n)(1-u_n))s_{n,\infty }+1}}{2}\right)^{-\left(4n-1+\Omega_1 \right)} }{\sqrt{s_{n,\infty }^2-2(1-2\eta (1-t_n)(1-u_n))s_{n,\infty }+1}}}\cr
{ \displaystyle \overleftrightarrow {\mathbf{\Gamma}}_{n-k} \left(s_{n-k};t_{n-k},u_{n-k},\widetilde{w}_{n-k+1,n} \right)}\cr
{ \displaystyle = \frac{ \left( \frac{1+s_{n-k}+\sqrt{s_{n-k}^2-2(1-2\widetilde{w}_{n-k+1,n} (1-t_{n-k})(1-u_{n-k}))s_{n-k}+1}}{2}\right)^{-\left(4(n-k)-1+\Omega_1\right)} }{\sqrt{s_{n-k}^2-2(1-2\widetilde{w}_{n-k+1,n} (1-t_{n-k})(1-u_{n-k}))s_{n-k}+1}} }
\end{cases}\nonumber 
\end{equation}
  and
 \begin{equation}
\begin{cases} 
{ \displaystyle \mathbf{A} \left( s_{0,\infty } ;\eta\right)= \frac{ \left(1+s_{0,\infty }+\sqrt{s_{0,\infty }^2-2(1-2\eta )s_{0,\infty }+1}\right)^{2-\gamma -\Omega_1}}{\left(1- s_{0,\infty }+\sqrt{s_{0,\infty }^2-2(1-2\eta )s_{0,\infty }+1}\right)^{1-\gamma }\sqrt{s_{0,\infty }^2-2(1-2\eta )s_{0,\infty }+1}}}\cr
{ \displaystyle  \mathbf{A} \left( s_{0} ;\widetilde{w}_{1,1}\right) = \frac{ \left(1+s_0+\sqrt{s_0^2-2(1-2\widetilde{w}_{1,1} )s_0+1}\right)^{2-\gamma -\Omega_1}}{\left(1- s_0+\sqrt{s_0^2-2(1-2\widetilde{w}_{1,1})s_0+1}\right)^{1-\gamma }\sqrt{s_0^2-2(1-2\widetilde{w}_{1,1})s_0+1}}} \cr
{ \displaystyle \mathbf{A} \left( s_{0} ;\widetilde{w}_{1,n}\right) = \frac{ \left(1+s_0+\sqrt{s_0^2-2(1-2\widetilde{w}_{1,n} )s_0+1}\right)^{2-\gamma -\Omega_1}}{\left(1- s_0+\sqrt{s_0^2-2(1-2\widetilde{w}_{1,n})s_0+1}\right)^{1-\gamma }\sqrt{s_0^2-2(1-2\widetilde{w}_{1,n})s_0+1}}}
\end{cases}\nonumber 
\end{equation}
 \subsection{ ${\displaystyle  Hl(1-a,-q+\alpha \beta; \alpha,\beta, \delta, \gamma; 1-x)}$} 
Replacing coefficients $a$, $q$, $\gamma $, $\delta$ and $x$ by $1-a$, $-q +\alpha \beta $, $\delta $, $\gamma $ and $1-x$ into (\ref{eq:2009}).
\begin{eqnarray}
&&\sum_{q_0 =0}^{\infty } \frac{(\gamma)_{q_0}}{q_0!} s_0^{q_0} \prod _{n=1}^{\infty } \left\{ \sum_{ q_n = q_{n-1}}^{\infty } s_n^{q_n }\right\} Hl(1-a,-q+\alpha \beta; \alpha,\beta, \delta, \gamma; 1-x) \nonumber\\
&&=2^{\frac{\varphi }{(2-a)}-1}\left\{ \prod_{l=1}^{\infty } \frac{1}{(1-s_{l,\infty })} \mathbf{A} \left( s_{0,\infty } ;\eta\right)\right. \nonumber\\
&&+ \Bigg\{\prod_{l=2}^{\infty } \frac{1}{(1-s_{l,\infty })} \int_{0}^{1} dt_1\;t_1 \int_{0}^{1} du_1\;u_1^{\delta } \overleftrightarrow {\mathbf{\Gamma}}_1 \left(s_{1,\infty };t_1,u_1,\eta\right) \nonumber\\
&&\times  \widetilde{w}_{1,1}^{-\alpha}\left( \widetilde{w}_{1,1} \partial _{ \widetilde{w}_{1,1}}\right) \widetilde{w}_{1,1}^{\alpha -\beta } \left( \widetilde{w}_{1,1} \partial _{ \widetilde{w}_{1,1}}\right)\widetilde{w}_{1,1}^{\beta }  \mathbf{A} \left( s_{0} ;\widetilde{w}_{1,1}\right) \Bigg\}z \nonumber\\
&&+ \sum_{n=2}^{\infty } \left\{ \prod_{l=n+1}^{\infty } \frac{1}{(1-s_{l,\infty })} \int_{0}^{1} dt_n\;t_n^{2n-1} \int_{0}^{1} du_n\;u_n^{2(n-1)+\delta  } \right. \overleftrightarrow {\mathbf{\Gamma}}_n \left(s_{n,\infty };t_n,u_n,\eta \right)\nonumber\\
&&\times  \widetilde{w}_{n,n}^{-(2(n-1)+\alpha)}\left( \widetilde{w}_{n,n} \partial _{ \widetilde{w}_{n,n}}\right) \widetilde{w}_{n,n}^{\alpha -\beta } \left( \widetilde{w}_{n,n} \partial _{ \widetilde{w}_{n,n}}\right)\widetilde{w}_{n,n}^{2(n-1)+\beta } \nonumber\\
&&\times \prod_{k=1}^{n-1} \Bigg\{  \int_{0}^{1} dt_{n-k}\;t_{n-k}^{2(n-k)-1} \int_{0}^{1} du_{n-k} \;u_{n-k}^{2(n-k-1)+\delta }  \overleftrightarrow {\mathbf{\Gamma}}_{n-k} \left(s_{n-k};t_{n-k},u_{n-k},\widetilde{w}_{n-k+1,n} \right)  \nonumber\\
&&\times \widetilde{w}_{n-k,n}^{-(2(n-k-1)+\alpha)}\left( \widetilde{w}_{n-k,n} \partial _{ \widetilde{w}_{n-k,n}}\right) \widetilde{w}_{n-k,n}^{\alpha -\beta } \left( \widetilde{w}_{n-k,n} \partial _{ \widetilde{w}_{n-k,n}}\right)\widetilde{w}_{n-k,n}^{2(n-k-1)+\beta }  \Bigg\}\nonumber\\
&&\times \left.\mathbf{A} \left( s_{0} ;\widetilde{w}_{1,n}\right) \Bigg\} z^n \right\}  
 \label{eq:200202}
\end{eqnarray}
where
\begin{equation}
\begin{cases} \xi =1-x \cr
z = \frac{-1}{1-a}\xi^2 \cr
\eta = \frac{2-a}{1-a}\xi \cr
\varphi = \alpha +\beta -\delta +(1-a)(\delta +\gamma -1) \cr
\Omega _1= \frac{\varphi }{(2-a)}\cr
q = \alpha \beta +(q_j+2j)\{\varphi +(2-a)(q_j+2j) \} \;\;\mbox{as}\;j,q_j\in \mathbb{N}_{0} 
\end{cases}\nonumber 
\end{equation}
\begin{equation}
\begin{cases} 
{ \displaystyle \overleftrightarrow {\mathbf{\Gamma}}_1 \left(s_{1,\infty };t_1,u_1,\eta\right)= \frac{\left( \frac{1+s_{1,\infty }+\sqrt{s_{1,\infty }^2-2(1-2\eta (1-t_1)(1-u_1))s_{1,\infty }+1}}{2}\right)^{-\left(3+\Omega _1\right)}}{\sqrt{s_{1,\infty }^2-2(1-2\eta (1-t_1)(1-u_1))s_{1,\infty }+1}}}\cr
{ \displaystyle  \overleftrightarrow {\mathbf{\Gamma}}_n \left(s_{n,\infty };t_n,u_n,\eta \right) = \frac{ \left( \frac{1+s_{n,\infty }+\sqrt{s_{n,\infty }^2-2(1-2\eta (1-t_n)(1-u_n))s_{n,\infty }+1}}{2}\right)^{-\left(4n-1+\Omega _1\right)} }{\sqrt{s_{n,\infty }^2-2(1-2\eta (1-t_n)(1-u_n))s_{n,\infty }+1}}}\cr
{ \displaystyle \overleftrightarrow {\mathbf{\Gamma}}_{n-k} \left(s_{n-k};t_{n-k},u_{n-k},\widetilde{w}_{n-k+1,n} \right)}\cr
{ \displaystyle =  \frac{ \left( \frac{1+s_{n-k}+\sqrt{s_{n-k}^2-2(1-2\widetilde{w}_{n-k+1,n} (1-t_{n-k})(1-u_{n-k}))s_{n-k}+1}}{2}\right)^{-\left(4(n-k)-1+\Omega _1\right)} }{\sqrt{s_{n-k}^2-2(1-2\widetilde{w}_{n-k+1,n} (1-t_{n-k})(1-u_{n-k}))s_{n-k}+1}}}
\end{cases}\nonumber 
\end{equation}
 and
  \begin{equation}
\begin{cases} 
{ \displaystyle \mathbf{A} \left( s_{0,\infty } ;\eta\right)= \frac{\left(1+s_{0,\infty }+\sqrt{s_{0,\infty }^2-2(1-2\eta )s_{0,\infty }+1}\right)^{\delta  -\Omega _1}}{\left(1- s_{0,\infty }+\sqrt{s_{0,\infty }^2-2(1-2\eta )s_{0,\infty }+1}\right)^{\delta -1} \sqrt{s_{0,\infty }^2-2(1-2\eta )s_{0,\infty }+1}}}\cr
{ \displaystyle  \mathbf{A} \left( s_{0} ;\widetilde{w}_{1,1}\right) = \frac{ \left(1+s_0+\sqrt{s_0^2-2(1-2\widetilde{w}_{1,1} )s_0+1}\right)^{\delta -\Omega _1}}{\left(1- s_0+\sqrt{s_0^2-2(1-2\widetilde{w}_{1,1})s_0+1}\right)^{\delta -1}\sqrt{s_0^2-2(1-2\widetilde{w}_{1,1})s_0+1}}} \cr
{ \displaystyle \mathbf{A} \left( s_{0} ;\widetilde{w}_{1,n}\right) = \frac{\left(1+s_0+\sqrt{s_0^2-2(1-2\widetilde{w}_{1,n} )s_0+1}\right)^{\delta -\Omega _1}}{\left(1- s_0+\sqrt{s_0^2-2(1-2\widetilde{w}_{1,n})s_0+1}\right)^{\delta -1} \sqrt{s_0^2-2(1-2\widetilde{w}_{1,n})s_0+1}} }
\end{cases}\nonumber 
\end{equation}
 \subsection{ \footnotesize ${\displaystyle (1-x)^{1-\delta } Hl(1-a,-q+(\delta -1)\gamma a+(\alpha -\delta +1)(\beta -\delta +1); \alpha-\delta +1,\beta-\delta +1, 2-\delta, \gamma; 1-x)}$ \normalsize}
Replacing coefficients $a$, $q$, $\alpha $, $\beta $, $\gamma $, $\delta$ and $x$ by $1-a$, $-q+(\delta -1)\gamma a+(\alpha -\delta +1)(\beta -\delta +1)$, $\alpha-\delta +1 $, $\beta-\delta +1 $, $2-\delta$, $\gamma $ and $1-x$ into (\ref{eq:2009}). Multiply $(1-x)^{1-\delta }$ and the new (\ref{eq:2009}) together.
\begin{eqnarray}
&&\sum_{q_0 =0}^{\infty } \frac{(\gamma)_{q_0}}{q_0!} s_0^{q_0} \prod _{n=1}^{\infty } \left\{ \sum_{ q_n = q_{n-1}}^{\infty } s_n^{q_n }\right\} (1-x)^{1-\delta } \nonumber\\
&&\times Hl(1-a,-q+(\delta -1)\gamma a+(\alpha -\delta +1)(\beta -\delta +1); \alpha-\delta +1,\beta-\delta +1, 2-\delta, \gamma; 1-x) \nonumber\\ 
&&=\frac{2^{\frac{\varphi }{(2-a)}-1}}{(1-x)^{\delta -1}}\left\{ \prod_{l=1}^{\infty } \frac{1}{(1-s_{l,\infty })}\right. \mathbf{A} \left( s_{0,\infty } ;\eta\right) \nonumber\\
&&+ \Bigg\{\prod_{l=2}^{\infty } \frac{1}{(1-s_{l,\infty })} \int_{0}^{1} dt_1\;t_1 \int_{0}^{1} du_1\;u_1^{2-\delta } \overleftrightarrow {\mathbf{\Gamma}}_1 \left(s_{1,\infty };t_1,u_1,\eta\right)\nonumber\\
&&\times \widetilde{w}_{1,1}^{-\alpha+\delta -1}\left( \widetilde{w}_{1,1} \partial _{ \widetilde{w}_{1,1}}\right) \widetilde{w}_{1,1}^{\alpha -\beta } \left( \widetilde{w}_{1,1} \partial _{ \widetilde{w}_{1,1}}\right)\widetilde{w}_{1,1}^{\beta-\delta +1} \mathbf{A} \left( s_{0} ;\widetilde{w}_{1,1}\right)\Bigg\} z \nonumber\\
&&+ \sum_{n=2}^{\infty } \Bigg\{ \prod_{l=n+1}^{\infty } \frac{1}{(1-s_{l,\infty })} \int_{0}^{1} dt_n\;t_n^{2n-1} \int_{0}^{1} du_n\;u_n^{2n-\delta } \overleftrightarrow {\mathbf{\Gamma}}_n \left(s_{n,\infty };t_n,u_n,\eta \right) \nonumber\\
&&\times  \widetilde{w}_{n,n}^{-(2n-1+\alpha-\delta )}\left( \widetilde{w}_{n,n} \partial _{ \widetilde{w}_{n,n}}\right) \widetilde{w}_{n,n}^{\alpha -\beta } \left( \widetilde{w}_{n,n} \partial _{ \widetilde{w}_{n,n}}\right)\widetilde{w}_{n,n}^{2n-1+\beta-\delta }  \nonumber\\
&&\times \prod_{k=1}^{n-1} \Bigg\{  \int_{0}^{1} dt_{n-k}\;t_{n-k}^{2(n-k)-1} \int_{0}^{1} du_{n-k} \;u_{n-k}^{2(n-k)-\delta }  \overleftrightarrow {\mathbf{\Gamma}}_{n-k} \left(s_{n-k};t_{n-k},u_{n-k},\widetilde{w}_{n-k+1,n} \right)  \nonumber\\
&&\times \widetilde{w}_{n-k,n}^{-(2(n-k)-1+\alpha-\delta )}\left( \widetilde{w}_{n-k,n} \partial _{ \widetilde{w}_{n-k,n}}\right) \widetilde{w}_{n-k,n}^{\alpha -\beta } \left( \widetilde{w}_{n-k,n} \partial _{ \widetilde{w}_{n-k,n}}\right)\widetilde{w}_{n-k,n}^{2(n-k)-1+\beta-\delta  } \Bigg\}\nonumber\\
&&\times \left.\mathbf{A} \left( s_{0} ;\widetilde{w}_{1,n}\right) \Bigg\} z^n \right\}  
 \label{eq:200203}
\end{eqnarray}
where
\begin{equation}
\begin{cases} \xi =1-x \cr
z = \frac{-1}{1-a}\xi^2 \cr
\eta = \frac{2-a}{1-a}\xi \cr
\varphi = \alpha +\beta -\gamma -2\delta +2+(1-a)(\gamma -\delta +1) \cr
\Omega _1= \frac{\varphi }{(2-a)}\cr
q =(\delta -1)\gamma a+(\alpha -\delta +1)(\beta -\delta +1)+(q_j+2j)\{\varphi +(2-a)(q_j+2j ) \} \;\;\mbox{as}\;j,q_j\in \mathbb{N}_{0}
\end{cases}\nonumber 
\end{equation}
\begin{equation}
\begin{cases} 
{ \displaystyle \overleftrightarrow {\mathbf{\Gamma}}_1 \left(s_{1,\infty };t_1,u_1,\eta\right)= \frac{\left( \frac{1+s_{1,\infty }+\sqrt{s_{1,\infty }^2-2(1-2\eta (1-t_1)(1-u_1))s_{1,\infty }+1}}{2}\right)^{-\left(3+\Omega _1\right)}}{\sqrt{s_{1,\infty }^2-2(1-2\eta (1-t_1)(1-u_1))s_{1,\infty }+1}}}\cr
{ \displaystyle \overleftrightarrow {\mathbf{\Gamma}}_n \left(s_{n,\infty };t_n,u_n,\eta \right) = \frac{ \left( \frac{1+s_{n,\infty }+\sqrt{s_{n,\infty }^2-2(1-2\eta (1-t_n)(1-u_n))s_{n,\infty }+1}}{2}\right)^{-\left(4n-1+\Omega _1\right)} }{\sqrt{s_{n,\infty }^2-2(1-2\eta (1-t_n)(1-u_n))s_{n,\infty }+1}}}\cr
{ \displaystyle \overleftrightarrow {\mathbf{\Gamma}}_{n-k} \left(s_{n-k};t_{n-k},u_{n-k},\widetilde{w}_{n-k+1,n} \right)}\cr
{ \displaystyle = \frac{ \left( \frac{1+s_{n-k}+\sqrt{s_{n-k}^2-2(1-2\widetilde{w}_{n-k+1,n} (1-t_{n-k})(1-u_{n-k}))s_{n-k}+1}}{2}\right)^{-\left(4(n-k)-1+\Omega _1\right)} }{\sqrt{s_{n-k}^2-2(1-2\widetilde{w}_{n-k+1,n} (1-t_{n-k})(1-u_{n-k}))s_{n-k}+1}} }
\end{cases}\nonumber 
\end{equation}
  and
 \begin{equation}
\begin{cases} 
{ \displaystyle \mathbf{A} \left( s_{0,\infty } ;\eta\right)= \frac{\left(1+s_{0,\infty }+\sqrt{s_{0,\infty }^2-2(1-2\eta )s_{0,\infty }+1}\right)^{2-\delta -\Omega _1}}{\left(1- s_{0,\infty }+\sqrt{s_{0,\infty }^2-2(1-2\eta )s_{0,\infty }+1}\right)^{1-\delta} \sqrt{s_{0,\infty }^2-2(1-2\eta )s_{0,\infty }+1}}}\cr
{ \displaystyle \mathbf{A} \left( s_{0} ;\widetilde{w}_{1,1}\right) = \frac{\left(1+s_0+\sqrt{s_0^2-2(1-2\widetilde{w}_{1,1} )s_0+1}\right)^{2-\delta -\Omega _1}}{\left(1- s_0+\sqrt{s_0^2-2(1-2\widetilde{w}_{1,1})s_0+1}\right)^{1-\delta} \sqrt{s_0^2-2(1-2\widetilde{w}_{1,1})s_0+1}}} \cr
{ \displaystyle \mathbf{A} \left( s_{0} ;\widetilde{w}_{1,n}\right) = \frac{\left(1+s_0+\sqrt{s_0^2-2(1-2\widetilde{w}_{1,n} )s_0+1}\right)^{2-\delta -\Omega _1}}{\left(1- s_0+\sqrt{s_0^2-2(1-2\widetilde{w}_{1,n})s_0+1}\right)^{1-\delta} \sqrt{s_0^2-2(1-2\widetilde{w}_{1,n})s_0+1}} }
\end{cases}\nonumber 
\end{equation}
 \subsection{ \footnotesize ${\displaystyle x^{-\alpha } Hl\left(\frac{1}{a},\frac{q+\alpha [(\alpha -\gamma -\delta +1)a-\beta +\delta ]}{a}; \alpha , \alpha -\gamma +1, \alpha -\beta +1,\delta ;\frac{1}{x}\right)}$\normalsize}
Replacing coefficients $a$, $q$, $\beta $, $\gamma $ and $x$ by $\frac{1}{a}$, $\frac{q+\alpha [(\alpha -\gamma -\delta +1)a-\beta +\delta ]}{a}$, $\alpha-\gamma +1 $, $\alpha -\beta +1 $ and $\frac{1}{x}$ into (\ref{eq:2009}). Multiply $x^{-\alpha }$ and the new (\ref{eq:2009}) together.
\begin{eqnarray}
&&\sum_{q_0 =0}^{\infty } \frac{(\gamma)_{q_0}}{q_0!} s_0^{q_0} \prod _{n=1}^{\infty } \left\{ \sum_{ q_n = q_{n-1}}^{\infty } s_n^{q_n }\right\}x^{-\alpha } Hl\Big(\frac{1}{a},\frac{q+\alpha [(\alpha -\gamma -\delta +1)a-\beta +\delta ]}{a}; \alpha , \alpha -\gamma +1 \nonumber\\
&&, \alpha -\beta +1,\delta ;\frac{1}{x}\Big) \nonumber\\
&&=\frac{2^{\frac{a\varphi }{(1+a)}-1}}{x^{\alpha }}\left\{ \prod_{l=1}^{\infty } \frac{1}{(1-s_{l,\infty })} \mathbf{A} \left( s_{0,\infty } ;\eta\right)\right. \nonumber\\
&&+ \Bigg\{\prod_{l=2}^{\infty } \frac{1}{(1-s_{l,\infty })} \int_{0}^{1} dt_1\;t_1 \int_{0}^{1} du_1\;u_1^{\alpha -\beta +1} \overleftrightarrow {\mathbf{\Gamma}}_1 \left(s_{1,\infty };t_1,u_1,\eta\right)\nonumber\\
&&\times \widetilde{w}_{1,1}^{-\alpha}\left( \widetilde{w}_{1,1} \partial _{ \widetilde{w}_{1,1}}\right) \widetilde{w}_{1,1}^{\gamma -1} \left( \widetilde{w}_{1,1} \partial _{ \widetilde{w}_{1,1}}\right)\widetilde{w}_{1,1}^{\alpha -\gamma +1}  \mathbf{A} \left( s_{0} ;\widetilde{w}_{1,1}\right) \Bigg\}z \nonumber\\
&&+ \sum_{n=2}^{\infty } \Bigg\{ \prod_{l=n+1}^{\infty } \frac{1}{(1-s_{l,\infty })} \int_{0}^{1} dt_n\;t_n^{2n-1} \int_{0}^{1} du_n\;u_n^{2n-1+\alpha -\beta } \overleftrightarrow {\mathbf{\Gamma}}_n \left(s_{n,\infty };t_n,u_n,\eta \right) \nonumber\\
&&\times \widetilde{w}_{n,n}^{-(2(n-1)+\alpha)}\left( \widetilde{w}_{n,n} \partial _{ \widetilde{w}_{n,n}}\right) \widetilde{w}_{n,n}^{\gamma -1} \left( \widetilde{w}_{n,n} \partial _{ \widetilde{w}_{n,n}}\right)\widetilde{w}_{n,n}^{2n-1+\alpha -\gamma } \nonumber\\
&&\times \prod_{k=1}^{n-1} \Bigg\{  \int_{0}^{1} dt_{n-k}\;t_{n-k}^{2(n-k)-1} \int_{0}^{1} du_{n-k} \;u_{n-k}^{2(n-k)-1+\alpha -\beta } \overleftrightarrow {\mathbf{\Gamma}}_{n-k} \left(s_{n-k};t_{n-k},u_{n-k},\widetilde{w}_{n-k+1,n} \right)  \nonumber\\
&&\times \widetilde{w}_{n-k,n}^{-(2(n-k-1)+\alpha)}\left( \widetilde{w}_{n-k,n} \partial _{ \widetilde{w}_{n-k,n}}\right) \widetilde{w}_{n-k,n}^{\gamma -1} \left( \widetilde{w}_{n-k,n} \partial _{ \widetilde{w}_{n-k,n}}\right)\widetilde{w}_{n-k,n}^{2(n-k)-1+\alpha -\gamma } \Bigg\}\nonumber\\
&&\times \left. \mathbf{A} \left( s_{0} ;\widetilde{w}_{1,n}\right) \Bigg\} z^n \right\}  \label{eq:200204} 
\end{eqnarray}
where
\begin{equation}
\begin{cases} \xi =\frac{1}{x} \cr
z = -a \xi^2 \cr
\eta = (1+a)\xi \cr
\varphi = 2\alpha -\gamma -\delta +1+\frac{1}{a}(\alpha -\beta +\delta ) \cr
\Omega _1= \frac{a\varphi }{(1+a)}\cr
q =-\alpha [(\alpha -\gamma -\delta +1)a-\beta +\delta]-a(q_j+2j )\{\varphi +(1+1/a)(q_j+2j ) \} \;\;\mbox{as}\;j,q_j\in \mathbb{N}_{0} 
\end{cases}\nonumber 
\end{equation}
\begin{equation}
\begin{cases} 
{ \displaystyle \overleftrightarrow {\mathbf{\Gamma}}_1 \left(s_{1,\infty };t_1,u_1,\eta\right)= \frac{\left( \frac{1+s_{1,\infty }+\sqrt{s_{1,\infty }^2-2(1-2\eta (1-t_1)(1-u_1))s_{1,\infty }+1}}{2}\right)^{-\left(3+\Omega _1\right)}}{\sqrt{s_{1,\infty }^2-2(1-2\eta (1-t_1)(1-u_1))s_{1,\infty }+1}}}\cr
{ \displaystyle \overleftrightarrow {\mathbf{\Gamma}}_n \left(s_{n,\infty };t_n,u_n,\eta \right) = \frac{ \left( \frac{1+s_{n,\infty }+\sqrt{s_{n,\infty }^2-2(1-2\eta (1-t_n)(1-u_n))s_{n,\infty }+1}}{2}\right)^{-\left(4n-1+\Omega _1\right)} }{\sqrt{s_{n,\infty }^2-2(1-2\eta (1-t_n)(1-u_n))s_{n,\infty }+1}}}\cr
{ \displaystyle \overleftrightarrow {\mathbf{\Gamma}}_{n-k} \left(s_{n-k};t_{n-k},u_{n-k},\widetilde{w}_{n-k+1,n} \right)}\cr
{ \displaystyle =  \frac{ \left( \frac{1+s_{n-k}+\sqrt{s_{n-k}^2-2(1-2\widetilde{w}_{n-k+1,n} (1-t_{n-k})(1-u_{n-k}))s_{n-k}+1}}{2}\right)^{-\left(4(n-k)-1+\Omega _1\right)} }{\sqrt{s_{n-k}^2-2(1-2\widetilde{w}_{n-k+1,n} (1-t_{n-k})(1-u_{n-k}))s_{n-k}+1}}}
\end{cases}\nonumber 
\end{equation}
 and
 \begin{equation}
\begin{cases} 
{ \displaystyle \mathbf{A} \left( s_{0,\infty } ;\eta\right)= \frac{ \left(1+s_{0,\infty }+\sqrt{s_{0,\infty }^2-2(1-2\eta )s_{0,\infty }+1}\right)^{\alpha -\beta +1 -\Omega _1}}{\left(1- s_{0,\infty }+\sqrt{s_{0,\infty }^2-2(1-2\eta )s_{0,\infty }+1}\right)^{\alpha-\beta }\sqrt{s_{0,\infty }^2-2(1-2\eta )s_{0,\infty }+1}}}\cr
{ \displaystyle \mathbf{A} \left( s_{0} ;\widetilde{w}_{1,1}\right) = \frac{\left(1+s_0+\sqrt{s_0^2-2(1-2\widetilde{w}_{1,1} )s_0+1}\right)^{\alpha -\beta +1 -\Omega _1}}{\left(1- s_0+\sqrt{s_0^2-2(1-2\widetilde{w}_{1,1})s_0+1}\right)^{\alpha-\beta } \sqrt{s_0^2-2(1-2\widetilde{w}_{1,1})s_0+1}}} \cr
{ \displaystyle \mathbf{A} \left( s_{0} ;\widetilde{w}_{1,n}\right) = \frac{ \left(1+s_0+\sqrt{s_0^2-2(1-2\widetilde{w}_{1,n} )s_0+1}\right)^{\alpha -\beta +1 -\Omega _1}}{\left(1- s_0+\sqrt{s_0^2-2(1-2\widetilde{w}_{1,n})s_0+1}\right)^{\alpha-\beta }\sqrt{s_0^2-2(1-2\widetilde{w}_{1,n})s_0+1}}}
\end{cases}\nonumber 
\end{equation}
 \subsection{ ${\displaystyle \left(1-\frac{x}{a} \right)^{-\beta } Hl\left(1-a, -q+\gamma \beta; -\alpha +\gamma +\delta, \beta, \gamma, \delta; \frac{(1-a)x}{x-a} \right)}$}
Replacing coefficients $a$, $q$, $\alpha $ and $x$ by $1-a$, $-q+\gamma \beta $, $-\alpha+\gamma +\delta $ and $\frac{(1-a)x}{x-a}$ into (\ref{eq:2009}). Multiply $\left(1-\frac{x}{a} \right)^{-\beta }$ and the new (\ref{eq:2009}) together.
\begin{eqnarray}
&&\sum_{q_0 =0}^{\infty } \frac{(\gamma)_{q_0}}{q_0!} s_0^{q_0} \prod _{n=1}^{\infty } \left\{ \sum_{ q_n = q_{n-1}}^{\infty } s_n^{q_n }\right\} \left(1-\frac{x}{a} \right)^{-\beta } Hl\left(1-a, -q+\gamma \beta; -\alpha +\gamma +\delta, \beta, \gamma, \delta; \frac{(1-a)x}{x-a} \right) \nonumber\\
&&=\frac{2^{\frac{\varphi }{(2-a)}-1}}{\left(1-\frac{x}{a} \right)^{\beta }}\left\{ \prod_{l=1}^{\infty } \frac{1}{(1-s_{l,\infty })} \mathbf{A} \left( s_{0,\infty } ;\eta\right)\right. \nonumber\\
&&+ \Bigg\{\prod_{l=2}^{\infty } \frac{1}{(1-s_{l,\infty })} \int_{0}^{1} dt_1\;t_1 \int_{0}^{1} du_1\;u_1^{\gamma} \overleftrightarrow {\mathbf{\Gamma}}_1 \left(s_{1,\infty };t_1,u_1,\eta\right)\nonumber\\
&&\times \widetilde{w}_{1,1}^{\alpha-\gamma -\delta }\left( \widetilde{w}_{1,1} \partial _{ \widetilde{w}_{1,1}}\right) \widetilde{w}_{1,1}^{-\alpha -\beta +\gamma +\delta } \left( \widetilde{w}_{1,1} \partial _{ \widetilde{w}_{1,1}}\right)\widetilde{w}_{1,1}^{\beta }  \mathbf{A} \left( s_{0} ;\widetilde{w}_{1,1}\right) \Bigg\}z \nonumber\\
&&+ \sum_{n=2}^{\infty } \bigg\{ \prod_{l=n+1}^{\infty } \frac{1}{(1-s_{l,\infty })} \int_{0}^{1} dt_n\;t_n^{2n-1} \int_{0}^{1} du_n\;u_n^{2(n-1)+\gamma }  \overleftrightarrow {\mathbf{\Gamma}}_n \left(s_{n,\infty };t_n,u_n,\eta \right)\nonumber\\
&&\times \widetilde{w}_{n,n}^{-(2(n-1)-\alpha+\gamma +\delta )}\left( \widetilde{w}_{n,n} \partial _{ \widetilde{w}_{n,n}}\right) \widetilde{w}_{n,n}^{-\alpha -\beta +\gamma +\delta } \left( \widetilde{w}_{n,n} \partial _{ \widetilde{w}_{n,n}}\right)\widetilde{w}_{n,n}^{2(n-1)+\beta } \nonumber\\
&&\times \prod_{k=1}^{n-1} \Bigg\{  \int_{0}^{1} dt_{n-k}\;t_{n-k}^{2(n-k)-1} \int_{0}^{1} du_{n-k} \;u_{n-k}^{2(n-k-1)+\gamma} \overleftrightarrow {\mathbf{\Gamma}}_{n-k} \left(s_{n-k};t_{n-k},u_{n-k},\widetilde{w}_{n-k+1,n} \right)  \nonumber\\
&&\times \widetilde{w}_{n-k,n}^{-(2(n-k-1)-\alpha+\gamma +\delta )}\left( \widetilde{w}_{n-k,n} \partial _{ \widetilde{w}_{n-k,n}}\right) \widetilde{w}_{n-k,n}^{-\alpha -\beta +\gamma +\delta } \left( \widetilde{w}_{n-k,n} \partial _{ \widetilde{w}_{n-k,n}}\right)\widetilde{w}_{n-k,n}^{2(n-k-1)+\beta } \Bigg\} \nonumber\\
&&\times \left. \mathbf{A} \left( s_{0} ;\widetilde{w}_{1,n}\right) \Bigg\} z^n \right\} \label{eq:200205}
\end{eqnarray}
where
\begin{equation}
\begin{cases} \xi = \frac{(1-a)x}{x-a} \cr
z = -\frac{1}{1-a}\xi^2 \cr
\eta = \frac{2-a}{1-a} \xi \cr
\varphi = -\alpha +\beta +\gamma +(1-a)(\gamma +\delta -1) \cr
\Omega _1= \frac{\varphi }{(2-a)}\cr
q=\gamma \beta +(q_j+2j)\{\varphi +(2-a)(q_j+2j) \} \;\;\mbox{as}\;j,q_j\in \mathbb{N}_{0} 
\end{cases}\nonumber 
\end{equation}
\begin{equation}
\begin{cases} 
{ \displaystyle \overleftrightarrow {\mathbf{\Gamma}}_1 \left(s_{1,\infty };t_1,u_1,\eta\right)= \frac{\left( \frac{1+s_{1,\infty }+\sqrt{s_{1,\infty }^2-2(1-2\eta (1-t_1)(1-u_1))s_{1,\infty }+1}}{2}\right)^{-\left(3+\Omega _1\right)}}{\sqrt{s_{1,\infty }^2-2(1-2\eta (1-t_1)(1-u_1))s_{1,\infty }+1}}}\cr
{ \displaystyle \overleftrightarrow {\mathbf{\Gamma}}_n \left(s_{n,\infty };t_n,u_n,\eta \right) = \frac{ \left( \frac{1+s_{n,\infty }+\sqrt{s_{n,\infty }^2-2(1-2\eta (1-t_n)(1-u_n))s_{n,\infty }+1}}{2}\right)^{-\left(4n-1+\Omega _1\right)} }{\sqrt{s_{n,\infty }^2-2(1-2\eta (1-t_n)(1-u_n))s_{n,\infty }+1}}}\cr
{ \displaystyle \overleftrightarrow {\mathbf{\Gamma}}_{n-k} \left(s_{n-k};t_{n-k},u_{n-k},\widetilde{w}_{n-k+1,n} \right)}\cr
{ \displaystyle = \frac{ \left( \frac{1+s_{n-k}+\sqrt{s_{n-k}^2-2(1-2\widetilde{w}_{n-k+1,n} (1-t_{n-k})(1-u_{n-k}))s_{n-k}+1}}{2}\right)^{-\left(4(n-k)-1+\Omega _1\right)} }{\sqrt{s_{n-k}^2-2(1-2\widetilde{w}_{n-k+1,n} (1-t_{n-k})(1-u_{n-k}))s_{n-k}+1}} }
\end{cases}\nonumber 
\end{equation}
 and
  \begin{equation}
\begin{cases} 
{ \displaystyle \mathbf{A} \left( s_{0,\infty } ;\eta\right)= \frac{ \left(1+s_{0,\infty }+\sqrt{s_{0,\infty }^2-2(1-2\eta )s_{0,\infty }+1}\right)^{\gamma -\Omega _1}}{\left(1- s_{0,\infty }+\sqrt{s_{0,\infty }^2-2(1-2\eta )s_{0,\infty }+1}\right)^{\gamma -1}\sqrt{s_{0,\infty }^2-2(1-2\eta )s_{0,\infty }+1}}}\cr
{ \displaystyle \mathbf{A} \left( s_{0} ;\widetilde{w}_{1,1}\right) = \frac{ \left(1+s_0+\sqrt{s_0^2-2(1-2\widetilde{w}_{1,1} )s_0+1}\right)^{\gamma -\Omega _1}}{\left(1- s_0+\sqrt{s_0^2-2(1-2\widetilde{w}_{1,1})s_0+1}\right)^{ \gamma -1}\sqrt{s_0^2-2(1-2\widetilde{w}_{1,1})s_0+1}}} \cr
{ \displaystyle \mathbf{A} \left( s_{0} ;\widetilde{w}_{1,n}\right) = \frac{\left(1+s_0+\sqrt{s_0^2-2(1-2\widetilde{w}_{1,n} )s_0+1}\right)^{\gamma -\Omega _1}}{\left(1- s_0+\sqrt{s_0^2-2(1-2\widetilde{w}_{1,n})s_0+1}\right)^{ \gamma -1} \sqrt{s_0^2-2(1-2\widetilde{w}_{1,n})s_0+1}}}
\end{cases}\nonumber 
\end{equation}
 \subsection{ \footnotesize ${\displaystyle (1-x)^{1-\delta }\left(1-\frac{x}{a} \right)^{-\beta+\delta -1} Hl\left(1-a, -q+\gamma [(\delta -1)a+\beta -\delta +1]; -\alpha +\gamma +1, \beta -\delta+1, \gamma, 2-\delta; \frac{(1-a)x}{x-a} \right)}$\normalsize}
Replacing coefficients $a$, $q$, $\alpha $, $\beta $, $\delta $ and $x$ by $1-a$, $-q+\gamma [(\delta -1)a+\beta -\delta +1]$, $-\alpha +\gamma +1$, $\beta -\delta+1$, $2-\delta $ and $\frac{(1-a)x}{x-a}$ into (\ref{eq:2009}). Multiply $(1-x)^{1-\delta }\left(1-\frac{x}{a} \right)^{-\beta+\delta -1}$ and the new (\ref{eq:2009}) together.
\begin{eqnarray}
&&\sum_{q_0 =0}^{\infty } \frac{(\gamma)_{q_0}}{q_0!} s_0^{q_0} \prod _{n=1}^{\infty } \left\{ \sum_{ q_n = q_{n-1}}^{\infty } s_n^{q_n }\right\} (1-x)^{1-\delta }\left(1-\frac{x}{a} \right)^{-\beta+\delta -1}\nonumber\\
&&\times Hl\left(1-a, -q+\gamma [(\delta -1)a+\beta -\delta +1]; -\alpha +\gamma +1, \beta -\delta+1, \gamma, 2-\delta; \frac{(1-a)x}{x-a} \right) \nonumber\\
&&=\frac{2^{\frac{\varphi }{(2-a)}-1}\left(1-\frac{x}{a} \right)^{-\beta+\delta -1}}{(1-x)^{\delta -1}}\left\{ \prod_{l=1}^{\infty } \frac{1}{(1-s_{l,\infty })}\right. \mathbf{A} \left( s_{0,\infty } ;\eta\right)  \nonumber\\
&&+ \Bigg\{\prod_{l=2}^{\infty } \frac{1}{(1-s_{l,\infty })} \int_{0}^{1} dt_1\;t_1 \int_{0}^{1} du_1\;u_1^{\gamma} \overleftrightarrow {\mathbf{\Gamma}}_1 \left(s_{1,\infty };t_1,u_1,\eta\right)\nonumber\\
&&\times \widetilde{w}_{1,1}^{\alpha-\gamma -1}\left( \widetilde{w}_{1,1} \partial _{ \widetilde{w}_{1,1}}\right) \widetilde{w}_{1,1}^{-\alpha -\beta+\gamma +\delta  } \left( \widetilde{w}_{1,1} \partial _{ \widetilde{w}_{1,1}}\right)\widetilde{w}_{1,1}^{\beta -\delta +1}  \mathbf{A} \left( s_{0} ;\widetilde{w}_{1,1}\right) \Bigg\}z \nonumber\\
&&+ \sum_{n=2}^{\infty } \Bigg\{ \prod_{l=n+1}^{\infty } \frac{1}{(1-s_{l,\infty })} \int_{0}^{1} dt_n\;t_n^{2n-1} \int_{0}^{1} du_n\;u_n^{2(n-1)+\gamma } \overleftrightarrow {\mathbf{\Gamma}}_n \left(s_{n,\infty };t_n,u_n,\eta \right) \nonumber\\
&&\times \widetilde{w}_{n,n}^{-(2n-1-\alpha+\gamma )}\left( \widetilde{w}_{n,n} \partial _{ \widetilde{w}_{n,n}}\right) \widetilde{w}_{n,n}^{-\alpha -\beta+\gamma +\delta } \left( \widetilde{w}_{n,n} \partial _{ \widetilde{w}_{n,n}}\right)\widetilde{w}_{n,n}^{2n-1+\beta -\delta } \nonumber\\
&&\times \prod_{k=1}^{n-1} \Bigg\{  \int_{0}^{1} dt_{n-k}\;t_{n-k}^{2(n-k)-1} \int_{0}^{1} du_{n-k} \;u_{n-k}^{2(n-k-1)+\gamma}  \overleftrightarrow {\mathbf{\Gamma}}_{n-k} \left(s_{n-k};t_{n-k},u_{n-k},\widetilde{w}_{n-k+1,n} \right) \nonumber\\
&&\times \widetilde{w}_{n-k,n}^{-(2(n-k)-1-\alpha+\gamma )}\left( \widetilde{w}_{n-k,n} \partial _{ \widetilde{w}_{n-k,n}}\right) \widetilde{w}_{n-k,n}^{-\alpha -\beta +\gamma +\delta } \left( \widetilde{w}_{n-k,n} \partial _{ \widetilde{w}_{n-k,n}}\right)\widetilde{w}_{n-k,n}^{2(n-k)-1+\beta-\delta } \Bigg\} \nonumber\\
&&\times \left. \mathbf{A} \left( s_{0} ;\widetilde{w}_{1,n}\right) \Bigg\} z^n \right\}  
 \label{eq:200206}
\end{eqnarray}
where
\begin{equation}
\begin{cases} 
\xi= \frac{(1-a)x}{x-a} \cr
z = -\frac{1}{1-a}\xi^2 \cr
\eta = \frac{(2-a)}{1-a} \xi \cr
\varphi = -\alpha +\beta +\gamma +(1-a)(\gamma -\delta +1) \cr
\Omega _1= \frac{\varphi }{(2-a)}\cr
q= \gamma [(\delta -1)a+\beta -\delta +1]+(q_j+2j)\{\varphi +(2-a)(q_j+2j) \} \;\;\mbox{as}\;j,q_j\in \mathbb{N}_{0} 
\end{cases}\nonumber 
\end{equation}
\begin{equation}
\begin{cases} 
{ \displaystyle \overleftrightarrow {\mathbf{\Gamma}}_1 \left(s_{1,\infty };t_1,u_1,\eta\right)= \frac{\left( \frac{1+s_{1,\infty }+\sqrt{s_{1,\infty }^2-2(1-2\eta (1-t_1)(1-u_1))s_{1,\infty }+1}}{2}\right)^{-\left(3+\Omega _1\right)}}{\sqrt{s_{1,\infty }^2-2(1-2\eta (1-t_1)(1-u_1))s_{1,\infty }+1}}}\cr
{ \displaystyle \overleftrightarrow {\mathbf{\Gamma}}_n \left(s_{n,\infty };t_n,u_n,\eta \right) = \frac{ \left( \frac{1+s_{n,\infty }+\sqrt{s_{n,\infty }^2-2(1-2\eta (1-t_n)(1-u_n))s_{n,\infty }+1}}{2}\right)^{-\left(4n-1+\Omega _1\right)} }{\sqrt{s_{n,\infty }^2-2(1-2\eta (1-t_n)(1-u_n))s_{n,\infty }+1}}}\cr
{ \displaystyle \overleftrightarrow {\mathbf{\Gamma}}_{n-k} \left(s_{n-k};t_{n-k},u_{n-k},\widetilde{w}_{n-k+1,n} \right)}\cr
{ \displaystyle = \frac{ \left( \frac{1+s_{n-k}+\sqrt{s_{n-k}^2-2(1-2\widetilde{w}_{n-k+1,n} (1-t_{n-k})(1-u_{n-k}))s_{n-k}+1}}{2}\right)^{-\left(4(n-k)-1+\Omega _1\right)} }{\sqrt{s_{n-k}^2-2(1-2\widetilde{w}_{n-k+1,n} (1-t_{n-k})(1-u_{n-k}))s_{n-k}+1}} }
\end{cases}\nonumber 
\end{equation}
 and
 \begin{equation}
\begin{cases} 
{ \displaystyle \mathbf{A} \left( s_{0,\infty } ;\eta\right) =  \frac{\left(1+s_{0,\infty }+\sqrt{s_{0,\infty }^2-2(1-2\eta )s_{0,\infty }+1}\right)^{\gamma -\Omega _1}}{\left(1- s_{0,\infty }+\sqrt{s_{0,\infty }^2-2(1-2\eta )s_{0,\infty }+1}\right)^{ \gamma -1} \sqrt{s_{0,\infty }^2-2(1-2\eta )s_{0,\infty }+1}}}\cr
{ \displaystyle \mathbf{A} \left( s_{0} ;\widetilde{w}_{1,1}\right) = \frac{\left(1+s_0+\sqrt{s_0^2-2(1-2\widetilde{w}_{1,1} )s_0+1}\right)^{\gamma -\Omega _1}}{\left(1- s_0+\sqrt{s_0^2-2(1-2\widetilde{w}_{1,1})s_0+1}\right)^ {\gamma -1} \sqrt{s_0^2-2(1-2\widetilde{w}_{1,1})s_0+1}}} \cr
{ \displaystyle \mathbf{A} \left( s_{0} ;\widetilde{w}_{1,n}\right) = \frac{\left(1+s_0+\sqrt{s_0^2-2(1-2\widetilde{w}_{1,n} )s_0+1}\right)^{\gamma -\Omega _1}}{\left(1- s_0+\sqrt{s_0^2-2(1-2\widetilde{w}_{1,n})s_0+1}\right)^{ \gamma -1} \sqrt{s_0^2-2(1-2\widetilde{w}_{1,n})s_0+1}}}
\end{cases}\nonumber 
\end{equation}
 \subsection{ \footnotesize ${\displaystyle x^{-\alpha } Hl\left(\frac{a-1}{a}, \frac{-q+\alpha (\delta a+\beta -\delta )}{a}; \alpha, \alpha -\gamma +1, \delta , \alpha -\beta +1; \frac{x-1}{x} \right)}$\normalsize}
Replacing coefficients $a$, $q$, $\beta $, $\gamma $, $\delta $ and $x$ by $\frac{a-1}{a}$, $\frac{-q+\alpha (\delta a+\beta -\delta )}{a}$, $\alpha -\gamma +1$, $\delta $, $\alpha -\beta +1$ and $\frac{x-1}{x}$ into (\ref{eq:2009}). Multiply $x^{-\alpha }$ and the new (\ref{eq:2009}) together.
\begin{eqnarray}
&&\sum_{q_0 =0}^{\infty } \frac{(\gamma)_{q_0}}{q_0!} s_0^{q_0} \prod _{n=1}^{\infty } \left\{ \sum_{ q_n = q_{n-1}}^{\infty } s_n^{q_n }\right\}x^{-\alpha } Hl\Big(\frac{a-1}{a}, \frac{-q+\alpha (\delta a+\beta -\delta )}{a}; \alpha, \alpha -\gamma +1, \delta \nonumber\\
&&\hspace{3mm}, \alpha -\beta +1; \frac{x-1}{x} \Big) \nonumber\\
&&=\frac{2^{\frac{a \varphi }{(2a-1)}-1}}{x^{\alpha }}\left\{ \prod_{l=1}^{\infty } \frac{1}{(1-s_{l,\infty })} \right. \mathbf{A} \left( s_{0,\infty } ;\eta\right)  \nonumber\\
&&+ \Bigg\{\prod_{l=2}^{\infty } \frac{1}{(1-s_{l,\infty })} \int_{0}^{1} dt_1\;t_1 \int_{0}^{1} du_1\;u_1^{\delta } \overleftrightarrow {\mathbf{\Gamma}}_1 \left(s_{1,\infty };t_1,u_1,\eta\right)\nonumber\\
&&\times \widetilde{w}_{1,1}^{-\alpha}\left( \widetilde{w}_{1,1} \partial _{ \widetilde{w}_{1,1}}\right) \widetilde{w}_{1,1}^{\gamma -1} \left( \widetilde{w}_{1,1} \partial _{ \widetilde{w}_{1,1}}\right)\widetilde{w}_{1,1}^{\alpha -\gamma +1}  \mathbf{A} \left( s_{0} ;\widetilde{w}_{1,1}\right) \Bigg\}z \nonumber\\
&&+ \sum_{n=2}^{\infty } \Bigg\{ \prod_{l=n+1}^{\infty } \frac{1}{(1-s_{l,\infty })} \int_{0}^{1} dt_n\;t_n^{2n-1} \int_{0}^{1} du_n\;u_n^{2(n-1)+\delta  } \overleftrightarrow {\mathbf{\Gamma}}_n \left(s_{n,\infty };t_n,u_n,\eta \right)\nonumber\\
&&\times \widetilde{w}_{n,n}^{-(2(n-1)+\alpha)}\left( \widetilde{w}_{n,n} \partial _{ \widetilde{w}_{n,n}}\right) \widetilde{w}_{n,n}^{\gamma -1} \left( \widetilde{w}_{n,n} \partial _{ \widetilde{w}_{n,n}}\right)\widetilde{w}_{n,n}^{2n-1+\alpha -\gamma } \nonumber\\
&&\times \prod_{k=1}^{n-1} \Bigg\{  \int_{0}^{1} dt_{n-k}\;t_{n-k}^{2(n-k)-1} \int_{0}^{1} du_{n-k} \;u_{n-k}^{2(n-k-1)+\delta }  \overleftrightarrow {\mathbf{\Gamma}}_{n-k} \left(s_{n-k};t_{n-k},u_{n-k},\widetilde{w}_{n-k+1,n} \right)  \nonumber\\
&&\times  \widetilde{w}_{n-k,n}^{-(2(n-k-1)+\alpha)}\left( \widetilde{w}_{n-k,n} \partial _{ \widetilde{w}_{n-k,n}}\right) \widetilde{w}_{n-k,n}^{\gamma -1} \left( \widetilde{w}_{n-k,n} \partial _{ \widetilde{w}_{n-k,n}}\right)\widetilde{w}_{n-k,n}^{2(n-k)-1+\alpha -\gamma } \Bigg\} \nonumber\\
&&\times \left.\mathbf{A} \left( s_{0} ;\widetilde{w}_{1,n}\right) \Bigg\} z^n \right\}  
\label{eq:200207}
\end{eqnarray}
where
\begin{equation}
\begin{cases} 
\xi= \frac{x-1}{x} \cr
z = \frac{-a}{a-1}\xi^2 \cr
\eta = \frac{2a-1}{a-1} \xi \cr
\varphi = \alpha +\beta -\gamma +\frac{a-1}{a}(\alpha -\beta +\delta ) \cr
\Omega _1= \frac{a \varphi }{(2a-1)}\cr
q= \alpha (\delta a+\beta -\delta )+(q_j+2j)\{a\varphi +(2a-1)(q_j+2j) \} \;\;\mbox{as}\;j,q_j\in \mathbb{N}_{0} 
\end{cases}\nonumber 
\end{equation}
\begin{equation}
\begin{cases} 
{ \displaystyle \overleftrightarrow {\mathbf{\Gamma}}_1 \left(s_{1,\infty };t_1,u_1,\eta\right)= \frac{\left( \frac{1+s_{1,\infty }+\sqrt{s_{1,\infty }^2-2(1-2\eta (1-t_1)(1-u_1))s_{1,\infty }+1}}{2}\right)^{-\left(3+\Omega _1\right)}}{\sqrt{s_{1,\infty }^2-2(1-2\eta (1-t_1)(1-u_1))s_{1,\infty }+1}}}\cr
{ \displaystyle \overleftrightarrow {\mathbf{\Gamma}}_n \left(s_{n,\infty };t_n,u_n,\eta \right) = \frac{ \left( \frac{1+s_{n,\infty }+\sqrt{s_{n,\infty }^2-2(1-2\eta (1-t_n)(1-u_n))s_{n,\infty }+1}}{2}\right)^{-\left(4n-1+\Omega _1\right)} }{\sqrt{s_{n,\infty }^2-2(1-2\eta (1-t_n)(1-u_n))s_{n,\infty }+1}}}\cr
{ \displaystyle \overleftrightarrow {\mathbf{\Gamma}}_{n-k} \left(s_{n-k};t_{n-k},u_{n-k},\widetilde{w}_{n-k+1,n} \right)}\cr
{ \displaystyle = \frac{ \left( \frac{1+s_{n-k}+\sqrt{s_{n-k}^2-2(1-2\widetilde{w}_{n-k+1,n} (1-t_{n-k})(1-u_{n-k}))s_{n-k}+1}}{2}\right)^{-\left(4(n-k)-1+\Omega _1\right)} }{\sqrt{s_{n-k}^2-2(1-2\widetilde{w}_{n-k+1,n} (1-t_{n-k})(1-u_{n-k}))s_{n-k}+1}} }
\end{cases}\nonumber 
\end{equation}
  and
 \begin{equation}
\begin{cases} 
{ \displaystyle \mathbf{A} \left( s_{0,\infty } ;\eta\right) = \frac{\left(1+s_{0,\infty }+\sqrt{s_{0,\infty }^2-2(1-2\eta )s_{0,\infty }+1}\right)^{\delta  -\Omega _1}}{\left(1- s_{0,\infty }+\sqrt{s_{0,\infty }^2-2(1-2\eta )s_{0,\infty }+1}\right)^{ \delta -1}\sqrt{s_{0,\infty }^2-2(1-2\eta )s_{0,\infty }+1}}}\cr
{ \displaystyle \mathbf{A} \left( s_{0} ;\widetilde{w}_{1,1}\right) = \frac{\left(1+s_0+\sqrt{s_0^2-2(1-2\widetilde{w}_{1,1} )s_0+1}\right)^{\delta -\Omega _1}}{\left(1- s_0+\sqrt{s_0^2-2(1-2\widetilde{w}_{1,1})s_0+1}\right)^{ \delta -1} \sqrt{s_0^2-2(1-2\widetilde{w}_{1,1})s_0+1}}} \cr
{ \displaystyle \mathbf{A} \left( s_{0} ;\widetilde{w}_{1,n}\right) = \frac{\left(1+s_0+\sqrt{s_0^2-2(1-2\widetilde{w}_{1,n} )s_0+1}\right)^{\delta -\Omega _1}}{\left(1- s_0+\sqrt{s_0^2-2(1-2\widetilde{w}_{1,n})s_0+1}\right)^{ \delta -1} \sqrt{s_0^2-2(1-2\widetilde{w}_{1,n})s_0+1}}}
\end{cases}\nonumber 
\end{equation}
 
\subsection{ ${\displaystyle \left(\frac{x-a}{1-a} \right)^{-\alpha } Hl\left(a, q-(\beta -\delta )\alpha ; \alpha , -\beta+\gamma +\delta , \delta , \gamma; \frac{a(x-1)}{x-a} \right)}$}
Replacing coefficients $q$, $\beta $, $\gamma $, $\delta $ and $x$ by $q-(\beta -\delta )\alpha $, $-\beta+\gamma +\delta $, $\delta $,  $\gamma $ and $\frac{a(x-1)}{x-a}$ into (\ref{eq:2009}). Multiply $\left(\frac{x-a}{1-a} \right)^{-\alpha }$ and the new (\ref{eq:2009}) together.
\begin{eqnarray}
&&\sum_{q_0 =0}^{\infty } \frac{(\gamma)_{q_0}}{q_0!} s_0^{q_0} \prod _{n=1}^{\infty } \left\{ \sum_{ q_n = q_{n-1}}^{\infty } s_n^{q_n }\right\} \left(\frac{x-a}{1-a} \right)^{-\alpha } Hl\left(a, q-(\beta -\delta )\alpha ; \alpha , -\beta+\gamma +\delta , \delta , \gamma; \frac{a(x-1)}{x-a} \right) \nonumber\\
&&=\frac{2^{\frac{\varphi }{(1+a)}-1}}{\left(\frac{x-a}{1-a} \right)^{\alpha }}\left\{ \prod_{l=1}^{\infty } \frac{1}{(1-s_{l,\infty })} \right. \mathbf{A} \left( s_{0,\infty } ;\eta\right) \nonumber\\
&&+ \Bigg\{\prod_{l=2}^{\infty } \frac{1}{(1-s_{l,\infty })} \int_{0}^{1} dt_1\;t_1 \int_{0}^{1} du_1\;u_1^{\delta } \overleftrightarrow {\mathbf{\Gamma}}_1 \left(s_{1,\infty };t_1,u_1,\eta\right)\nonumber\\
&&\times \widetilde{w}_{1,1}^{-\alpha}\left( \widetilde{w}_{1,1} \partial _{ \widetilde{w}_{1,1}}\right) \widetilde{w}_{1,1}^{\alpha +\beta -\gamma -\delta } \left( \widetilde{w}_{1,1} \partial _{ \widetilde{w}_{1,1}}\right)\widetilde{w}_{1,1}^{-\beta +\gamma +\delta }  \mathbf{A} \left( s_{0} ;\widetilde{w}_{1,1}\right) \Bigg\}z \nonumber\\
&&+ \sum_{n=2}^{\infty } \Bigg\{ \prod_{l=n+1}^{\infty } \frac{1}{(1-s_{l,\infty })} \int_{0}^{1} dt_n\;t_n^{2n-1} \int_{0}^{1} du_n\;u_n^{2(n-1)+\delta  } \overleftrightarrow {\mathbf{\Gamma}}_n \left(s_{n,\infty };t_n,u_n,\eta \right)\nonumber\\
&&\times \widetilde{w}_{n,n}^{-(2(n-1)+\alpha)}\left( \widetilde{w}_{n,n} \partial _{ \widetilde{w}_{n,n}}\right) \widetilde{w}_{n,n}^{\alpha +\beta -\gamma -\delta } \left( \widetilde{w}_{n,n} \partial _{ \widetilde{w}_{n,n}}\right)\widetilde{w}_{n,n}^{2(n-1)-\beta +\gamma +\delta } \nonumber\\
&&\times \prod_{k=1}^{n-1} \Bigg\{  \int_{0}^{1} dt_{n-k}\;t_{n-k}^{2(n-k)-1} \int_{0}^{1} du_{n-k} \;u_{n-k}^{2(n-k-1)+\delta } \overleftrightarrow {\mathbf{\Gamma}}_{n-k} \left(s_{n-k};t_{n-k},u_{n-k},\widetilde{w}_{n-k+1,n} \right)  \nonumber\\
&&\times \widetilde{w}_{n-k,n}^{-(2(n-k-1)+\alpha)}\left( \widetilde{w}_{n-k,n} \partial _{ \widetilde{w}_{n-k,n}}\right) \widetilde{w}_{n-k,n}^{\alpha+\beta -\gamma -\delta } \left( \widetilde{w}_{n-k,n} \partial _{ \widetilde{w}_{n-k,n}}\right)\widetilde{w}_{n-k,n}^{2(n-k-1)-\beta +\gamma +\delta } \Bigg\} \nonumber\\
&&\times \left.\mathbf{A} \left( s_{0} ;\widetilde{w}_{1,n}\right) \Bigg\} z^n \right\}  
\label{eq:200208}
\end{eqnarray}
where
\begin{equation}
\begin{cases} 
\xi= \frac{a(x-1)}{x-a} \cr
z = -\frac{1}{a}\xi^2 \cr
\eta = \frac{(1+a)}{a} \xi \cr
\varphi = \alpha -\beta +\delta +a(\gamma +\delta -1) \cr
\Omega _1= \frac{\varphi }{(1+a)}\cr
q= (\beta -\delta )\alpha -(q_j+2j)\{\varphi +(1+a)(q_j+2j ) \} \;\;\mbox{as}\;j,q_j\in \mathbb{N}_{0} 
\end{cases}\nonumber 
\end{equation}
\begin{equation}
\begin{cases} 
{ \displaystyle \overleftrightarrow {\mathbf{\Gamma}}_1 \left(s_{1,\infty };t_1,u_1,\eta\right)= \frac{\left( \frac{1+s_{1,\infty }+\sqrt{s_{1,\infty }^2-2(1-2\eta (1-t_1)(1-u_1))s_{1,\infty }+1}}{2}\right)^{-\left(3+\Omega _1\right)}}{\sqrt{s_{1,\infty }^2-2(1-2\eta (1-t_1)(1-u_1))s_{1,\infty }+1}}}\cr
{ \displaystyle \overleftrightarrow {\mathbf{\Gamma}}_n \left(s_{n,\infty };t_n,u_n,\eta \right) = \frac{ \left( \frac{1+s_{n,\infty }+\sqrt{s_{n,\infty }^2-2(1-2\eta (1-t_n)(1-u_n))s_{n,\infty }+1}}{2}\right)^{-\left(4n-1+\Omega _1\right)} }{\sqrt{s_{n,\infty }^2-2(1-2\eta (1-t_n)(1-u_n))s_{n,\infty }+1}}}\cr
{ \displaystyle \overleftrightarrow {\mathbf{\Gamma}}_{n-k} \left(s_{n-k};t_{n-k},u_{n-k},\widetilde{w}_{n-k+1,n} \right)}\cr
{ \displaystyle = \frac{ \left( \frac{1+s_{n-k}+\sqrt{s_{n-k}^2-2(1-2\widetilde{w}_{n-k+1,n} (1-t_{n-k})(1-u_{n-k}))s_{n-k}+1}}{2}\right)^{-\left(4(n-k)-1+\Omega _1\right)} }{\sqrt{s_{n-k}^2-2(1-2\widetilde{w}_{n-k+1,n} (1-t_{n-k})(1-u_{n-k}))s_{n-k}+1}}  }
\end{cases}\nonumber 
\end{equation}
 and
 \begin{equation}
\begin{cases} 
{ \displaystyle \mathbf{A} \left( s_{0,\infty } ;\eta\right) = \frac{\left(1+s_{0,\infty }+\sqrt{s_{0,\infty }^2-2(1-2\eta )s_{0,\infty }+1}\right)^{\delta  -\Omega _1}}{\left(1- s_{0,\infty }+\sqrt{s_{0,\infty }^2-2(1-2\eta )s_{0,\infty }+1}\right)^{ \delta -1} \sqrt{s_{0,\infty }^2-2(1-2\eta )s_{0,\infty }+1}}}\cr
{ \displaystyle \mathbf{A} \left( s_{0} ;\widetilde{w}_{1,1}\right) = \frac{\left(1+s_0+\sqrt{s_0^2-2(1-2\widetilde{w}_{1,1} )s_0+1}\right)^{\delta -\Omega _1}}{\left(1- s_0+\sqrt{s_0^2-2(1-2\widetilde{w}_{1,1})s_0+1}\right)^{ \delta -1} \sqrt{s_0^2-2(1-2\widetilde{w}_{1,1})s_0+1}}} \cr
{ \displaystyle \mathbf{A} \left( s_{0} ;\widetilde{w}_{1,n}\right) = \frac{\left(1+s_0+\sqrt{s_0^2-2(1-2\widetilde{w}_{1,n} )s_0+1}\right)^{\delta -\Omega _1}}{\left(1- s_0+\sqrt{s_0^2-2(1-2\widetilde{w}_{1,n})s_0+1}\right)^{ \delta -1} \sqrt{s_0^2-2(1-2\widetilde{w}_{1,n})s_0+1}}}
\end{cases}\nonumber 
\end{equation}
 \end{appendices}

\addcontentsline{toc}{section}{Bibliography}
\bibliographystyle{model1a-num-names}
\bibliography{<your-bib-database>}
\bibliographystyle{model1a-num-names}
\bibliography{<your-bib-database>}
 
\chapter{Confluent Heun function using three-term recurrence formula}
\chaptermark{Confluent Heun function using 3TRF} 

The Confluent Heun equation (CHE) is one of 4 confluent forms of Heun's differential equation. The Confluent Heun function (CHF) is applicable to areas such as general relativity, theory of rotating/non-rotating black hole, the gauge theories on thick brane words, hydrogen molecule ion in Stark effect and etc.

In this chapter I will apply three term recurrence formula (3TRF) \cite{Chou2012} to the power series expansion in closed forms and its integral forms of the CHE for infinite series and polynomial which makes $B_n$ term terminated including all higher terms of $A_n$'s.\footnote{`` higher terms of $A_n$'s'' means at least two terms of $A_n$'s.} Indeed, I construct the generating function for Confluent Heun polynomial (CHP) which makes $B_n$ term terminated. 

\addtocontents{toc}{\protect\setcounter{tocdepth}{3}}
\section{Introduction}
\subsection[Non-symmetrical canonical form of the CHE]{Non-symmetrical canonical form of Confluent Heun Equation}
In Ref.\cite{Heun1889,Ronv1995}, Heun's equation is a second-order linear ordinary differential equation of the form
\begin{equation}
\frac{d^2{y}}{d{x}^2} + \left(\frac{\gamma }{x} +\frac{\delta }{x-1} + \frac{\epsilon }{x-a}\right) \frac{d{y}}{d{x}} +  \frac{\alpha \beta x-q}{x(x-1)(x-a)} y = 0 \label{eq:4001}
\end{equation}
With the condition $\epsilon = \alpha +\beta -\gamma -\delta +1$. The parameters play different roles: $a \ne 0 $ is the singularity parameter, $\alpha $, $\beta $, $\gamma $, $\delta $, $\epsilon $ are exponent parameters, q is the accessory parameter. Also, $\alpha $ and $\beta $ are identical to each other. The total number of free parameters is six. It has four regular singular points which are 0, 1, $a$ and $\infty $ with exponents $\{ 0, 1-\gamma \}$, $\{ 0, 1-\delta \}$, $\{ 0, 1-\epsilon \}$ and $\{ \alpha, \beta \}$.

Heun equation has the four kinds of confluent forms: (1) Confluent Heun (two regular and one irregular singularities), (2) Doubly Confluent Heun (two irregular singularities), (3) Biconfluent Heun (one regular and one irregular singularities)\footnote{Biconfluent Heun equation is derived from the grand confluent hypergeometric equation by changing all coefficients $\mu =1$ and $ \varepsilon \omega = -q$.\cite{Chou2012i,Chou2012j}}, (4) Triconfluent Heun equations (one irregular singularity).
We can derive these four confluent forms from Heun equation by combining two or more regular singularities to take form an irregular singularity. Its process, converting Heun equation to other confluent forms, is similar to deriving of confluent hypergeometric equation from the hypergeometric equation. For the CHE, divide (\ref{eq:4001}) by $a$ to get 
\begin{equation}
x(x-1)\left(\frac{x}{a}-1\right)\frac{d^2{y}}{d{x}^2}+\left( \gamma (x-1)\left(\frac{x}{a}-1\right)+ \delta x\left(\frac{x}{a}-1\right)+\frac{\epsilon }{a}x(x-1)\right)\frac{d{y}}{d{x}}+ \left( \alpha \frac{\beta }{a}x-\frac{q}{a}\right)y =0 \nonumber
\end{equation}
Let $a, \beta, \epsilon, q \rightarrow \infty $ in such a way that
\begin{equation}
\frac{\beta }{a} \rightarrow \frac{\beta }{a}\rightarrow -\beta , \frac{q}{a} \rightarrow -q \nonumber
\end{equation}
the new $\beta $ and $q$ being constants, we find the non-symmetrical canonical form of the Confluent Heun Equation (CHE).\cite{Ronv1995,Deca1978,Decar1978}  
\begin{equation}
\frac{d^2{y}}{d{x}^2} + \left(\beta  +\frac{\gamma }{x} + \frac{\delta }{x-1}\right) \frac{d{y}}{d{x}} +  \frac{\alpha \beta x-q}{x(x-1)} y = 0 \label{eq:4002}
\end{equation}
(\ref{eq:4002}) has three singular points: two regular singular points which are 0 and 1 with exponents $\{0, 1-\gamma\}$ and $\{0, 1-\delta \}$, and one irregular singular point which is $\infty$ with an exponent $\alpha$. Its solution is denoted as $H_{c}(\alpha,\beta,\gamma,\delta,q;x)$.\footnote{Several authors denote as coefficients $4p$ and $\sigma $ instead of $\beta $ and $q$. And they define the solution of the CHE as $H_{c}^{(a)}(p,\alpha,\gamma,\delta,\sigma;x)$.}  
\subsection{Generalized spheroidal equation (GSE)}
The B\^{o}cher symmetrical form of the CHE in applications of mathematical physics is written as \cite{Wils1928,Wils1928a}
\begin{equation}
\frac{d}{d{z}}(z^2-1)\frac{d}{d{z}} y(z) + \left( -p^2(z^2-1)+2p\beta z-\lambda -\frac{m^2+s^2+2msz}{(z^2-1)}\right) y(z) = 0 \label{eq:4003}
\end{equation}
It is called the generalized spheroidal equation (GSE). The GSE has three singular points: two regular singular points which are $\pm 1$ and one irregular singular point which is $\infty$.
Putting ${\displaystyle y(z)= (1-z^2)^{-\frac{1}{2}}v(z)}$ in (\ref{eq:4002}), we obtain the normal symmetrical form of the CHE:
\begin{equation}
(z^2-1)\frac{d^2{v}}{d{z}^2} + \left( -p^2(z^2-1)+2p\beta z-\lambda -\frac{m^2+s^2-1+2msz}{(z^2-1)}\right)v(z) = 0 \label{eq:4004}
\end{equation}
It is called a normal form of the GSE.
By putting ${ \displaystyle y(z)= (z-1)^{\pm \frac{1}{2}(m+s)} (z+1)^{\pm \frac{1}{2}(m-s)} e^{\pm pz}w(z)}$ into (\ref{eq:4003}), The solution of the GSE can be described in the form of non-symmetrical canonical form of the CHE. The 16 possible solutions are as follows.
\begin{table}
\begin{center}
\begin{tabular}{cc}
\hline
$ $\\
$1.$ & $y(z) = (z-1)^{\frac{1}{2}(m+s)} (z+1)^{\frac{1}{2}(m-s)} e^{pz} H_{c}\big( \beta +m+1, 4p, m-s+1, m+s+1$\\
& $, 2p(\beta +1)-m(m+1)-\lambda; \frac{1+z}{2}\big)$\\ 
$2.$ & $y(z) = (z-1)^{\frac{1}{2}(m+s)} (z+1)^{\frac{1}{2}(m-s)} e^{-pz} H_{c}\big( -\beta +m+1, -4p, m-s+1, m+s+1$\\ 
& $, 2p(\beta -1)-m(m+1)-\lambda; \frac{1+z}{2}\big)$\\
$3.$ & $y(z) = (z-1)^{\frac{1}{2}(m+s)} (z+1)^{\frac{1}{2}(m-s)} e^{pz} H_{c}\big( \beta +m+1, -4p, m+s+1, m-s+1$\\
& $, -2p(\beta +2m+1)-m(m+1)-\lambda; \frac{1-z}{2}\big)$\\ 
$4.$ & $y(z) = (z-1)^{\frac{1}{2}(m+s)} (z+1)^{\frac{1}{2}(m-s)} e^{-pz} H_{c}\big( -\beta +m+1, 4p, m+s+1, m-s+1$\\ 
& $, -2p(\beta -2m-1)-m(m+1)-\lambda; \frac{1-z}{2}\big)$\\
$5.$ & $y(z) = (z-1)^{-\frac{1}{2}(m+s)} (z+1)^{-\frac{1}{2}(m-s)} e^{pz} H_{c}\big( \beta -m+1, 4p, -m+s+1, -m-s+1$\\
& $, 2p(\beta +1)+m(-m+1)-\lambda; \frac{1+z}{2}\big)$\\ 
$6.$ & $y(z) = (z-1)^{-\frac{1}{2}(m+s)} (z+1)^{-\frac{1}{2}(m-s)} e^{-pz} H_{c}\big( -\beta -m+1, -4p, -m+s+1, -m-s+1$\\ 
& $, 2p(\beta -1)+m(-m+1)-\lambda; \frac{1+z}{2}\big)$\\
$7.$ & $y(z) = (z-1)^{-\frac{1}{2}(m+s)} (z+1)^{-\frac{1}{2}(m-s)} e^{pz} H_{c}\big( \beta -m+1, -4p, -m-s+1, -m+s+1$\\ 
& $, -2p(\beta -2m +1)+m(-m+1)-\lambda; \frac{1-z}{2}\big)$\\
$8.$ & $y(z) = (z-1)^{-\frac{1}{2}(m+s)} (z+1)^{-\frac{1}{2}(m-s)} e^{-pz} H_{c}\big( -\beta -m+1, 4p, -m-s+1, -m+s+1$\\ 
& $, -2p(\beta -2m -1)+m(-m+1)-\lambda; \frac{1-z}{2}\big)$\\
$9.$ & $y(z) = (z-1)^{\frac{1}{2}(m+s)} (z+1)^{-\frac{1}{2}(m-s)} e^{pz} H_{c}\big( \beta +s+1, 4p, -m+s+1, m+s+1$\\ 
& $, 2p(\beta +1)-s(s+1)-\lambda; \frac{1+z}{2}\big)$\\
$10.$ & $y(z) = (z-1)^{\frac{1}{2}(m+s)} (z+1)^{-\frac{1}{2}(m-s)} e^{-pz} H_{c}\big( -\beta +s+1, -4p, -m+s+1, m+s+1$\\ 
& $, 2p(\beta -1)-s(s+1)-\lambda; \frac{1+z}{2}\big)$\\
$11.$ & $y(z) = (z-1)^{\frac{1}{2}(m+s)} (z+1)^{-\frac{1}{2}(m-s)} e^{pz} H_{c}\big( \beta +s+1, -4p, m+s+1, -m+s+1$\\ 
& $, -2p(\beta +2s+1)-s(s+1)-\lambda; \frac{1-z}{2}\big)$\\
$12.$ & $y(z) = (z-1)^{\frac{1}{2}(m+s)} (z+1)^{-\frac{1}{2}(m-s)} e^{-pz} H_{c}\big( -\beta +s+1, 4p, m+s+1, -m+s+1$\\ 
& $, -2p(\beta -2s-1)-s(s+1)-\lambda; \frac{1-z}{2}\big)$\\
$13.$ & $y(z) = (z-1)^{-\frac{1}{2}(m+s)} (z+1)^{\frac{1}{2}(m-s)} e^{pz} H_{c}\big( \beta -s+1, 4p, m-s+1, -m-s+1$\\ 
& $, 2p(\beta +1)+s(-s+1)-\lambda; \frac{1+z}{2}\big)$\\
$14.$ & $y(z) = (z-1)^{-\frac{1}{2}(m+s)} (z+1)^{\frac{1}{2}(m-s)} e^{-pz} H_{c}\big( -\beta -s+1, -4p, m-s+1, -m-s+1$\\ 
& $, 2p(\beta -1)+s(-s+1)-\lambda; \frac{1+z}{2}\big)$\\
$15.$ & $y(z) = (z-1)^{-\frac{1}{2}(m+s)} (z+1)^{\frac{1}{2}(m-s)} e^{pz} H_{c}\big( \beta -s+1, -4p, -m-s+1, m-s+1$\\ 
& $, -2p(\beta -2s+1)+s(-s+1)-\lambda; \frac{1-z}{2}\big)$\\
$16.$ & $y(z) = (z-1)^{-\frac{1}{2}(m+s)} (z+1)^{\frac{1}{2}(m-s)} e^{-pz} H_{c}\big( -\beta -s+1, 4p, -m-s+1, m-s+1$\\ 
& $, -2p(\beta +2s+1)+s(-s+1)-\lambda; \frac{1-z}{2}\big)$\\
$ $\\
\hline
\end{tabular}
\end{center}
\caption{The 16 possible solutions of the GSE using non-symmetrical canonical form of the CHE}
\end{table}   
\newpage
\subsection{Generalized spheroidal equation in the Leaver version}
The generalized spheroidal wave equation (GSWE), one of standard forms of the CHE, in the Leaver version reads \cite{Leav1986}
\begin{equation}
z(z-z_0)\frac{d^2{U}}{d{z}^2} + (B_1+B_2 z)\frac{d{U}}{d{z}}+ \left( B_3 -2\eta \omega (z-z_0)+\omega ^2z(z-z_0)\right) U(z) = 0 \hspace{.5cm}\mbox{where}\; \omega \ne0 \label{eq:4005}
\end{equation}
where $B_1, B_2, B_3, \eta $ and $\omega $ are constant. (\ref{eq:4005}) has three singular points: two regular singular points which are 0 and $z_0$ with exponents $\left\{ 0, 1+\frac{B_1}{z_0}\right\}$ and $\left\{ 0, 1-B_2-\frac{B_1}{z_0}\right\}$, and one irregular singular point which is $\infty $.
To permit the non-symmetrical canonical form of the CHE limit, the first step is given by the substitutions
\begin{equation}
z=z_0(1-x),\hspace{1cm} U(z)= \exp\left( iwz_0(1-x)\right) y(x) \nonumber
\end{equation}
which convert (\ref{eq:4005}) into
\begin{equation}
\frac{d^2{y}}{d{x}^2} + \left( -2 i\omega z_0  +\frac{B_2+\frac{B_1}{z_0}}{x} + \frac{\frac{-B_1}{z_0}}{x-1} \right) \frac{d{y}}{d{x}} +   \frac{-2 i\omega z_0\left( i\eta +\frac{B_2}{2}\right)x- \left( -B_2-i\omega z_0\left( B_2+\frac{B_1}{z_0}\right) \right) }{x(x-1)} y(x) = 0 \label{eq:4006}
\end{equation}
As we compare all coefficients in (\ref{eq:4006}) with (\ref{eq:4002}), we find that
\begin{equation}
\begin{cases} \alpha \rightarrow i\eta +\frac{B_2}{2}  \cr
\beta \rightarrow -2 i\omega z_0  \cr
\gamma \rightarrow B_2 +\frac{B_1}{z_0}  \cr
\delta \rightarrow -\frac{B_1}{z_0}  \cr
q \rightarrow -B_2 - i z_0 \left( B_2+ \frac{B_1}{2}\right)  
\end{cases}\nonumber 
\end{equation}
Then the solution $U(z)$ for the GSWE in the Leaver version is described by the non-symmetrical canonical form of the CHE.
\begin{eqnarray}
&&U(B_1, B_2, B_3; z_0, \omega, \eta; z)\label{eq:4007}\\ 
&&= \exp\left( iwz_0(1-x)\right)H_c\left( i\eta +\frac{B_2}{2},-2 i\omega z_0,B_2 +\frac{B_1}{z_0},-\frac{B_1}{z_0},-B_2 - i z_0 \left( B_2+ \frac{B_1}{2}\right); z_0(1-x)\right) \nonumber
\end{eqnarray} 
\section[CHF with a regular singular point at zero]{Confluent Heun function with a regular singular point at zero} 
In this chapter, by applying 3TRF \cite{Chou2012}, I'll construct the power series expansion in closed forms, the integral form and the generating function for non-symmetrical canonical form of the CHE for infinite series and polynomial which makes $B_n$ term terminated including all higher terms of $A_n$'s: I prefer (\ref{eq:4002}) as an analytic solution of the CHE rather than the GSE and its Leaver version because of more convenient mathematical calculations. 

Assume that its solution is
\begin{equation}
y(x)= \sum_{n=0}^{\infty } c_n x^{n+\lambda } \hspace{1cm}\mbox{where}\; \lambda =\mbox{indicial}\;\mbox{root}\label{eq:4008}
\end{equation}
Plug (\ref{eq:4008})  into (\ref{eq:4002}). Then we get a three-term recurrence relation for the coefficients $c_n$:
\begin{equation}
c_{n+1}=A_n \;c_n +B_n \;c_{n-1} \hspace{1cm};n\geq 1 \label{eq:4009}
\end{equation}
where,
\begin{subequations}
\begin{equation}
A_n =\frac{(n+\lambda )(n+\lambda -\beta +\gamma +\delta -1)-q}{(n+\lambda +1)(n+\lambda +\gamma )}\label{eq:40010a}
\end{equation}
\begin{equation}
B_n = \frac{\beta (n+\lambda +\alpha -1)}{(n+\lambda +1)(n+\lambda +\gamma )} \label{eq:40010b}
\end{equation}
\begin{equation}
c_1= A_0 \;c_0 \label{eq:40010c}
\end{equation}
\end{subequations}
We have two indicial roots which are $\lambda _1= 0$ and $\lambda _2= 1-\gamma $.

Now let's test for convergence of infinite series of the analytic function $y(x)$. As $n\gg 1$ (for sufficiently large), (\ref{eq:40010a}) and (\ref{eq:40010b}) are
\begin{subequations}
\begin{equation}
\lim_{n\gg 1} A_n = 1  
 \label{eq:40011a}
\end{equation}
\begin{equation}
  \lim_{n\gg 1} B_n = \frac{\beta }{n}
 \label{eq:40011b}
\end{equation}
\end{subequations}
(\ref{eq:40011b}) is negligible because of $n\gg 1$. Substitute (\ref{eq:40011a}) into (\ref{eq:4009}) with $B_n =0$. For $n=0,1,2,\cdots$, it give
\begin{equation}
\begin{tabular}{ l }
  \vspace{2 mm}
  $c_0$ \\
  \vspace{2 mm}
  $c_1 = c_0 $ \\
  \vspace{2 mm}
  $c_2 = c_0 $ \\
  \vspace{2 mm}
  $c_3 = c_0 $ \\
  \vspace{2 mm}                       
 \; \vdots \hspace{5mm} \vdots \\
\end{tabular}\label{eq:40012}
\end{equation}
When a function $y(x)$, analytic at $x=0$, is expanded in a power series putting $c_0=1$ by using (\ref{eq:40012}), we write
\begin{equation}
\lim_{n\gg 1}y(x) = \frac{1}{1-x}\hspace{1cm}\mbox{where}\;|x|<1 \label{eq:40013}
\end{equation}
For polynomial which makes $A_n$ terminated, (\ref{eq:40011b}) is only available for the asymptotic behavior of the minimum $y(x)$.\footnote{(\ref{eq:40011a}) is negligible for the minimum $y(x)$ because $A_n$ term will be terminated at the specific eigenvalues. } Substitute (\ref{eq:40011b}) into (\ref{eq:4009}) with $A_n =0$. For $n=0,1,2,\cdots$, it gives
\begin{equation}
\begin{tabular}{  l  l }
  \vspace{2 mm}
   $c_0$ &\hspace{1cm}  $c_1$ \\
  \vspace{2 mm}
   $c_2 = \frac{\Gamma \left(\frac{1}{2}\right)}{\Gamma \left(\frac{3}{2}\right)} \frac{\beta }{2}c_0 $  &\hspace{1cm}  $c_3 = \frac{1}{1!} \frac{\beta }{2} c_1$  \\
  \vspace{2 mm}
  $c_4 = \frac{\Gamma \left(\frac{1}{2}\right)}{\Gamma \left(\frac{5}{2}\right)} \left(\frac{\beta }{2}\right)^2 c_0 $ &\hspace{1cm}  $c_5 = \frac{1}{2!} \left(\frac{\beta }{2}\right)^2 c_1$\\
  \vspace{2 mm}
  $c_6 = \frac{\Gamma \left(\frac{1}{2}\right)}{\Gamma \left(\frac{7}{2}\right)} \left(\frac{\beta }{2}\right)^3 c_0 $ &\hspace{1cm}  $c_7 = \frac{1}{3!} \left(\frac{\beta }{2}\right)^3 c_1$\\
 \hspace{2 mm} \large{\vdots} & \hspace{1.5 cm}\large{\vdots} \\
\end{tabular}
\label{eq:40069}
\end{equation}
Put (\ref{eq:40069}) into (\ref{eq:4008}) taking $c_0=1$ and $\lambda = 0$ for simplicity by letting $c_1=A_0 c_0 =0$.
\begin{equation}
\mbox{min}\left( \lim_{n\gg 1}y(z) \right) = 1+  e^{\eta} \sqrt{\pi \eta } \; \mbox{Erf}\left(\sqrt{\eta}\right) \hspace{1cm}\mbox{where}\;\eta =\frac{1}{2}\beta x^2 \label{eq:40070}
\end{equation}
On the above $\mbox{Erf(y)} $ is an error function which is
\begin{equation}
\mbox{Erf(y)} = \frac{2}{\sqrt{\pi }} \int_{0}^{y} dt\; e^{-t^2}\nonumber
\end{equation} 
\subsection{Power series}
\subsubsection{Polynomial of type 1}
As I mentioned in chapter 2, there are three types of polynomials in three term recurrence relation of a linear ordinary differential equation: (1) type 1 is the polynomial which makes $B_n$ term terminated: $A_n$ term is not terminated, (2) type 2 is the polynomial which makes $A_n$ term terminated: $B_n$ term is not terminated, (3) type 3 is the polynomial which makes $A_n$ and $B_n$ terms terminated at the same time.\footnote{If $A_n$ and $B_n$ terms are not terminated, it turns to be infinite series.} In general the CHP is defined as type 3 polynomial where $A_n$ and $B_n$ terms terminated. The CHP comes from the CHE that has a fixed integer value of $\alpha $, just as it has a fixed value of $q$. In three term recurrence relation, polynomial of type 3 I categorize as complete polynomial. In future papers I will derive type 3 CHP. In this chapter I construct the power series expansion, its integral forms and the generating function for the CHP of type 1: I treat $\beta $, $\gamma $, $\delta $ and $q$ as free variables and $\alpha $ as a fixed value. 

In Ref.\cite{Chou2012}, the general expression of power series of $y(x)$ for polynomial of type 1 is given by
\begin{eqnarray}
 y(x)&=& \sum_{n=0}^{\infty } y_n(x)= y_0(x)+ y_1(x)+ y_2(x)+ y_3(x)+\cdots \nonumber\\
&=&  c_0 \Bigg\{ \sum_{i_0=0}^{\beta _0} \left( \prod _{i_1=0}^{i_0-1}B_{2i_1+1} \right) x^{2i_0+\lambda } 
+ \sum_{i_0=0}^{\beta _0}\left\{ A_{2i_0} \prod _{i_1=0}^{i_0-1}B_{2i_1+1}  \sum_{i_2=i_0}^{\beta _1} \left( \prod _{i_3=i_0}^{i_2-1}B_{2i_3+2} \right)\right\} x^{2i_2+1+\lambda }\nonumber\\
 && + \sum_{N=2}^{\infty } \Bigg\{ \sum_{i_0=0}^{\beta _0} \Bigg\{A_{2i_0}\prod _{i_1=0}^{i_0-1} B_{2i_1+1} \prod _{k=1}^{N-1} \Bigg( \sum_{i_{2k}= i_{2(k-1)}}^{\beta _k} A_{2i_{2k}+k}\prod _{i_{2k+1}=i_{2(k-1)}}^{i_{2k}-1}B_{2i_{2k+1}+(k+1)}\Bigg)\nonumber\\
 &&\times  \sum_{i_{2N} = i_{2(N-1)}}^{\beta _N} \Bigg( \prod _{i_{2N+1}=i_{2(N-1)}}^{i_{2N}-1} B_{2i_{2N+1}+(N+1)} \Bigg) \Bigg\} \Bigg\} x^{2i_{2N}+N+\lambda }\Bigg\}
  \label{eq:40014}
\end{eqnarray}
In the above, $\beta _i\leq \beta _j$ only if $i\leq j$ where $i,j,\beta _i, \beta _j \in \mathbb{N}_{0}$.

For a polynomial, we need a condition, which is:
\begin{equation}
 B_{2\beta _i + (i+1)}=0 \hspace{1cm} \mathrm{where}\; i= 0,1,2,\cdots, \beta _i=0,1,2,\cdots
 \label{eq:40015}
\end{equation}
In the above, $ \beta _i$ is an eigenvalue that makes $B_n$ term terminated at certain value of index $n$. (\ref{eq:40015}) makes each $y_i(x)$ where $i=0,1,2,\cdots$ as the polynomial in (\ref{eq:40014}).
Replace $\beta _i$ by $\alpha _i$ and put $n=2\alpha _i + (i+1)$ in (\ref{eq:40010b}) with the condition $B_{2\alpha  _i + (i+1)}=0$. Then, we obtain eigenvalues $\alpha = -2 \alpha _j- j-\lambda $. 

In (\ref{eq:40010b}) replace $\alpha $ by $-2 \alpha _i-i -\lambda $. In (\ref{eq:40014}) replace index $\beta _i$ by $\alpha _i$. Take (\ref{eq:40010a}) and the new (\ref{eq:40010b}) in new (\ref{eq:40014}).
After the replacement process, the general expression of power series of the CHE for polynomial of type 1 is given by
\begin{eqnarray}
 y(x)&=& \sum_{n=0}^{\infty } y_n(x)= y_0(x)+ y_1(x)+ y_2(x)+ y_3(x)+\cdots \nonumber\\
&=&  c_0 x^{\lambda } \left\{\sum_{i_0=0}^{\alpha _0} \frac{(-\alpha _0)_{i_0}}{(1+\frac{\lambda }{2})_{i_0}(\frac{1}{2}+ \frac{\gamma}{2} +\frac{\lambda }{2})_{i_0}} \eta ^{i_0}\right.\nonumber\\
&&+ \Bigg\{ \sum_{i_0=0}^{\alpha _0}\frac{(i_0+ \frac{\lambda }{2})\left( i_0+ \frac{1}{2}(-\beta +\gamma+\delta -1+\lambda )\right)- \frac{q}{4}}{(i_0+ \frac{1}{2}+ \frac{\lambda }{2})(i_0 + \frac{\gamma }{2}+ \frac{\lambda }{2})} \nonumber\\
&&\times \frac{(-\alpha _0)_{i_0}}{(1+\frac{\lambda }{2})_{i_0}(\frac{1}{2}+ \frac{\gamma}{2} +\frac{\lambda }{2})_{i_0}} \sum_{i_1=i_0}^{\alpha _1} \frac{(-\alpha _1)_{i_1}(\frac{3}{2}+\frac{\lambda }{2})_{i_0}(1+\frac{\gamma }{2}+ \frac{\lambda }{2})_{i_0}}{(-\alpha _1)_{i_0}(\frac{3}{2}+\frac{\lambda }{2})_{i_1}(1+ \frac{\gamma}{2} +\frac{\lambda }{2})_{i_1}} \eta ^{i_1} \Bigg\} x \nonumber\\
&&+ \sum_{n=2}^{\infty } \Bigg\{ \sum_{i_0=0}^{\alpha _0} \frac{(i_0+ \frac{\lambda }{2}) \left( i_0+ \frac{1}{2}(-\beta +\gamma+\delta -1+\lambda )\right)- \frac{q}{4}}{(i_0+ \frac{1}{2}+ \frac{\lambda }{2})(i_0 + \frac{\gamma }{2}+ \frac{\lambda }{2})}
   \frac{(-\alpha _0)_{i_0}}{(1+\frac{\lambda }{2})_{i_0}(\frac{1}{2}+ \frac{\gamma}{2} +\frac{\lambda }{2})_{i_0}}\nonumber\\
&&\times \prod _{k=1}^{n-1} \Bigg\{ \sum_{i_k=i_{k-1}}^{\alpha _k} \frac{(i_k+\frac{k}{2}+ \frac{\lambda }{2}) \left( i_k+ \frac{1}{2}(-\beta +\gamma +\delta +\lambda +k-1)\right)- \frac{q}{4}}{(i_k+ \frac{k}{2}+\frac{1}{2}+\frac{\lambda }{2})(i_k +\frac{k}{2}+\frac{\gamma }{2}+\frac{\lambda }{2})} \nonumber\\
&&\times \frac{(-\alpha _k)_{i_k}(1+ \frac{k}{2}+\frac{\lambda }{2})_{i_{k-1}}(\frac{1}{2}+\frac{k}{2}+\frac{\gamma }{2}+ \frac{\lambda }{2})_{i_{k-1}}}{(-\alpha _k)_{i_{k-1}}(1+\frac{k}{2}+\frac{\lambda }{2})_{i_k}(\frac{1}{2}+ \frac{k}{2}+ \frac{\gamma}{2} +\frac{\lambda }{2})_{i_k}}\Bigg\} \nonumber\\
&&\times \left.\sum_{i_n= i_{n-1}}^{\alpha _n} \frac{(-\alpha _n)_{i_n}(1+ \frac{n}{2}+\frac{\lambda }{2})_{i_{n-1}}(\frac{1}{2}+\frac{n}{2}+\frac{\gamma }{2}+ \frac{\lambda }{2})_{i_{n-1}}}{(-\alpha _n)_{i_{n-1}}(1+\frac{n}{2}+\frac{\lambda }{2})_{i_n}(\frac{1}{2}+ \frac{n}{2}+ \frac{\gamma}{2} +\frac{\lambda }{2})_{i_n}} \eta ^{i_n} \Bigg\} x^n \right\}\label{eq:40016}
\end{eqnarray}
where
\begin{equation}
\begin{cases} \eta = \frac{1}{2} \beta x^2 \cr
 \alpha = -2 \alpha _j- j-\lambda  \cr
\alpha _i\leq \alpha _j \;\;\mbox{only}\;\mbox{if}\;i\leq j\;\;\mbox{where}\;i,j,\alpha _i,\alpha _j\in \mathbb{N}_{0} 
\end{cases}\nonumber
\end{equation}
Put $c_0$= 1 as $\lambda =0$  for the first kind of independent solutions of the CHE and $\lambda =1-\gamma$ for the second one in (\ref{eq:40016}). 
\begin{remark}
The power series expansion of the CHE of the first kind for polynomial of type 1 about $x=0$ as $\alpha = -(2 \alpha _j+j) $ where $j,\alpha _j \in \mathbb{N}_{0}$ is
\begin{eqnarray}
 y(x)&=& H_c^{(a)}F_{\alpha _j}\left(\alpha = -(2 \alpha _j+j), \beta, \gamma, \delta, q; \eta =\frac{1}{2}\beta x^2  \right)\nonumber\\
&=&  \sum_{i_0=0}^{\alpha _0} \frac{(-\alpha _0)_{i_0}}{(1 )_{i_0}(\frac{1}{2}+ \frac{\gamma}{2} )_{i_0}} \eta ^{i_0} \nonumber\\
&+& \Bigg\{ \sum_{i_0=0}^{\alpha _0}\frac{ i_0 \left( i_0+ \frac{1}{2}(-\beta +\gamma+\delta -1 )\right)- \frac{q}{4}}{(i_0+ \frac{1}{2} )(i_0 + \frac{\gamma }{2} )}  \frac{(-\alpha _0)_{i_0}}{(1 )_{i_0}(\frac{1}{2}+ \frac{\gamma}{2} )_{i_0}} \sum_{i_1=i_0}^{\alpha _1} \frac{(-\alpha _1)_{i_1}(\frac{3}{2} )_{i_0}(1+\frac{\gamma }{2} )_{i_0}}{(-\alpha _1)_{i_0}(\frac{3}{2} )_{i_1}(1+ \frac{\gamma}{2} )_{i_1}} \eta ^{i_1} \Bigg\} x \nonumber\\
&+& \sum_{n=2}^{\infty } \Bigg\{ \sum_{i_0=0}^{\alpha _0} \frac{ i_0 \left(i_0+ \frac{1}{2}(-\beta +\gamma+\delta -1 )\right)- \frac{q}{4}}{(i_0+ \frac{1}{2} )(i_0 + \frac{\gamma }{2} )}
   \frac{(-\alpha _0)_{i_0}}{(1 )_{i_0}(\frac{1}{2}+ \frac{\gamma}{2} )_{i_0}}\nonumber\\
&\times& \prod _{k=1}^{n-1} \Bigg\{ \sum_{i_k=i_{k-1}}^{\alpha _k} \frac{(i_k+\frac{k}{2} ) \left(i_k+ \frac{1}{2}(-\beta +\gamma +\delta +k-1)\right)- \frac{q}{4}}{(i_k+ \frac{k}{2}+\frac{1}{2} )(i_k +\frac{k}{2}+\frac{\gamma }{2} )}  \frac{(-\alpha _k)_{i_k}(\frac{k}{2} +1)_{i_{k-1}}(\frac{k}{2}+\frac{1}{2}+\frac{\gamma }{2} )_{i_{k-1}}}{(-\alpha _k)_{i_{k-1}}(\frac{k}{2}+1 )_{i_k}(\frac{k}{2}+\frac{1}{2}+  \frac{\gamma}{2} )_{i_k}}\Bigg\} \nonumber\\
&\times&  \sum_{i_n= i_{n-1}}^{\alpha _n} \frac{(-\alpha _n)_{i_n}(\frac{n}{2}+1)_{i_{n-1}}(\frac{n}{2}+\frac{1}{2}+\frac{\gamma }{2} )_{i_{n-1}}}{(-\alpha _n)_{i_{n-1}}(\frac{n}{2}+1)_{i_n}(\frac{n}{2}+\frac{1}{2}+ \frac{\gamma}{2} )_{i_n}} \eta ^{i_n} \Bigg\} x^n  \label{eq:40017}
\end{eqnarray}
\end{remark}
For the minimum value of the CHE of the first kind for a polynomial of type 1 about $x=0$, put $\alpha _0=\alpha _1=\alpha _2=\cdots=0$ in (\ref{eq:40017}).
\begin{eqnarray}
y(x)&=& H_c^{(a)}F_{0}\left(\alpha = -j, \beta, \gamma, \delta, q; \eta =\frac{1}{2}\beta x^2  \right)\nonumber\\ 
&=& \; _2F_1\left( \frac{1}{2}\left( \varphi +\sqrt{\varphi ^2+4q}\right), \frac{1}{2}\left( \varphi -\sqrt{\varphi ^2+4q}\right), \gamma ,x\right) \label{choun:4001}
\end{eqnarray} 
where $\varphi =-\beta +\gamma +\delta -1$ and $|x|< 1$.
\begin{remark}
The power series expansion of the CHE of the second kind for polynomial of type 1 about $x=0$ as $\alpha = -(2 \alpha _j+j +1-\gamma ) $ where $j,\alpha _j \in \mathbb{N}_{0}$ is
\begin{eqnarray}
 y(x)&=& H_c^{(a)}S_{\alpha _j}\left(\alpha = -(2 \alpha _j+j+1-\gamma ), \beta, \gamma, \delta, q; \eta =\frac{1}{2}\beta x^2  \right)\nonumber\\
&=&  x^{1-\gamma } \left\{\sum_{i_0=0}^{\alpha _0} \frac{(-\alpha _0)_{i_0}}{( \frac{3 }{2}-\frac{\gamma }{2})_{i_0}(1)_{i_0}} \eta ^{i_0}\right.\nonumber\\
&+& \Bigg\{ \sum_{i_0=0}^{\alpha _0}\frac{\left(i_0+ \frac{1 }{2}-\frac{\gamma }{2}\right) \left( i_0+ \frac{1}{2}(-\beta +\delta )\right)- \frac{q}{4}}{(i_0+1- \frac{\gamma }{2})(i_0 + \frac{1}{2} )}  \frac{(-\alpha _0)_{i_0}}{( \frac{3 }{2}-\frac{\gamma }{2})_{i_0}(1)_{i_0}} \sum_{i_1=i_0}^{\alpha _1} \frac{(-\alpha _1)_{i_1}(2-\frac{\gamma }{2})_{i_0}( \frac{3}{2} )_{i_0}}{(-\alpha _1)_{i_0}(2-\frac{\gamma }{2})_{i_1}( \frac{3}{2} )_{i_1}} \eta ^{i_1} \Bigg\} x \nonumber\\
&+& \sum_{n=2}^{\infty } \Bigg\{ \sum_{i_0=0}^{\alpha _0} \frac{(i_0+ \frac{1}{2}- \frac{\gamma }{2}) \left(i_0+ \frac{1}{2}(-\beta +\delta )\right)- \frac{q}{4}}{(i_0+ 1- \frac{\gamma }{2})(i_0 + \frac{1}{2} )}
   \frac{(-\alpha _0)_{i_0}}{(\frac{3}{2}-\frac{\gamma }{2})_{i_0}(1)_{i_0}}\nonumber\\
&\times& \prod _{k=1}^{n-1} \Bigg\{ \sum_{i_k=i_{k-1}}^{\alpha _k} \frac{(i_k+\frac{k}{2}+ \frac{1}{2}-\frac{\gamma }{2}) \left(i_k+ \frac{1}{2}(-\beta +\delta +k )\right)- \frac{q}{4}}{(i_k+ \frac{k}{2}+1-\frac{\gamma }{2} )(i_k +\frac{k}{2}+\frac{1}{2} )}  \frac{(-\alpha _k)_{i_k}( \frac{k}{2}+\frac{3}{2}-\frac{\gamma }{2})_{i_{k-1}}( \frac{k}{2}+1)_{i_{k-1}}}{(-\alpha _k)_{i_{k-1}}(\frac{k}{2}+\frac{3}{2}-\frac{\gamma }{2})_{i_k}( \frac{k}{2}+1)_{i_k}}\Bigg\} \nonumber\\
&\times& \left.\sum_{i_n= i_{n-1}}^{\alpha _n} \frac{(-\alpha _n)_{i_n}( \frac{n}{2}+\frac{3}{2}-\frac{\gamma }{2})_{i_{n-1}}( \frac{n}{2} + 1)_{i_{n-1}}}{(-\alpha _n)_{i_{n-1}}( \frac{n}{2}+\frac{3}{2}-\frac{\gamma }{2})_{i_n}( \frac{n}{2}+1)_{i_n}} \eta ^{i_n} \Bigg\} x^n \right\}\label{eq:40018}
\end{eqnarray}
\end{remark}
For the minimum value of the CHE of the second kind for a polynomial of type 1 about $x=0$, put $\alpha _0=\alpha _1=\alpha _2=\cdots=0$ in (\ref{eq:40018}).
\begin{eqnarray}
y(x)&=& H_c^{(a)}S_{0}\left(\alpha = -(j+1-\gamma ), \beta, \gamma, \delta, q; \eta =\frac{1}{2}\beta x^2  \right)\nonumber\\ 
&=& x^{1-\gamma }\; _2F_1\left( \frac{1}{2}\left( \varphi +2(1-\gamma )+\sqrt{\varphi ^2+4q}\right), \frac{1}{2}\left( \varphi +2(1-\gamma )-\sqrt{\varphi ^2+4q}\right), 2-\gamma ,x\right) \hspace{1.5cm}\label{choun:4002}
\end{eqnarray} 
where $\varphi =-\beta +\gamma +\delta -1$ and $|x|< 1$.
In (\ref{choun:4001}) and (\ref{choun:4002}), a polynomial of type 1 requires $\left| x\right|< 1$ for the convergence of the radius.
\subsubsection{Infinite series}
In Ref.\cite{Chou2012}, the general expression of power series of $y(x)$ for infinite series is defined by
\begin{eqnarray}
y(x)  &=& \sum_{n=0}^{\infty } y_{n}(x)= y_0(x)+ y_1(x)+ y_2(x)+ y_3(x)+\cdots \nonumber\\
&=& c_0 \Bigg\{ \sum_{i_0=0}^{\infty } \left( \prod _{i_1=0}^{i_0-1}B_{2i_1+1} \right) x^{2i_0+\lambda } 
+ \sum_{i_0=0}^{\infty }\left\{ A_{2i_0} \prod _{i_1=0}^{i_0-1}B_{2i_1+1}  \sum_{i_2=i_0}^{\infty } \left( \prod _{i_3=i_0}^{i_2-1}B_{2i_3+2} \right)\right\} x^{2i_2+1+\lambda }  \nonumber\\
&& + \sum_{N=2}^{\infty } \Bigg\{ \sum_{i_0=0}^{\infty } \Bigg\{A_{2i_0}\prod _{i_1=0}^{i_0-1} B_{2i_1+1} 
 \prod _{k=1}^{N-1} \Bigg( \sum_{i_{2k}= i_{2(k-1)}}^{\infty } A_{2i_{2k}+k}\prod _{i_{2k+1}=i_{2(k-1)}}^{i_{2k}-1}B_{2i_{2k+1}+(k+1)}\Bigg)\nonumber\\
&& \times  \sum_{i_{2N} = i_{2(N-1)}}^{\infty } \Bigg( \prod _{i_{2N+1}=i_{2(N-1)}}^{i_{2N}-1} B_{2i_{2N+1}+(N+1)} \Bigg) \Bigg\} \Bigg\} x^{2i_{2N}+N+\lambda }\Bigg\} 
\label{eq:40019}
\end{eqnarray}
Substitute (\ref{eq:40010a})--(\ref{eq:40010c}) into (\ref{eq:40019}). 
The general expression of power series of the CHE for infinite series about $x=0$ is
\begin{eqnarray}
 y(x)&=& \sum_{n=0}^{\infty } y_n(x)= y_0(x)+ y_1(x)+ y_2(x)+ y_3(x)+\cdots \nonumber\\
&=& c_0 x^{\lambda } \left\{\sum_{i_0=0}^{\infty } \frac{(\frac{\alpha }{2}+\frac{\lambda }{2})_{i_0} }{(1+\frac{\lambda }{2})_{i_0}(\frac{1}{2}+ \frac{\gamma}{2} +\frac{\lambda }{2})_{i_0}} \eta ^{i_0}\right.\nonumber\\
&&+ \Bigg\{ \sum_{i_0=0}^{\infty }\frac{(i_0+ \frac{\lambda }{2}) \left( i_0+ \frac{1}{2}(-\beta +\gamma +\delta -1+\lambda )\right)-\frac{q}{4}}{(i_0+ \frac{1}{2}+ \frac{\lambda }{2})(i_0 + \frac{\gamma }{2}+ \frac{\lambda }{2})}\frac{(\frac{\alpha }{2}+\frac{\lambda }{2})_{i_0}}{(1+\frac{\lambda }{2})_{i_0}(\frac{1}{2}+ \frac{\gamma}{2} +\frac{\lambda }{2})_{i_0}}  \nonumber\\
&&\times  \sum_{i_1=i_0}^{\infty } \frac{(\frac{1}{2}+\frac{\alpha }{2}+ \frac{\lambda }{2})_{i_1} (\frac{3}{2}+\frac{\lambda }{2})_{i_0}(1+\frac{\gamma }{2}+ \frac{\lambda }{2})_{i_0}}{(\frac{1}{2}+\frac{\alpha }{2}+ \frac{\lambda }{2})_{i_0} (\frac{3}{2}+\frac{\lambda }{2})_{i_1}(1+ \frac{\gamma}{2} +\frac{\lambda }{2})_{i_1}} \eta ^{i_1} \Bigg\} x\nonumber\\
&&+ \sum_{n=2}^{\infty } \Bigg\{ \sum_{i_0=0}^{\infty } \frac{(i_0+ \frac{\lambda }{2}) \left( i_0+ \frac{1}{2}(-\beta +\gamma +\delta -1+\lambda )\right)-\frac{q}{4}}{(i_0+ \frac{1}{2}+ \frac{\lambda }{2})(i_0 + \frac{\gamma }{2}+ \frac{\lambda }{2})}
 \frac{(\frac{\alpha }{2}+\frac{\lambda }{2})_{i_0}}{(1+\frac{\lambda }{2})_{i_0}(\frac{1}{2}+ \frac{\gamma}{2} +\frac{\lambda }{2})_{i_0}}\nonumber\\
&&\times \prod _{k=1}^{n-1} \left\{ \sum_{i_k=i_{k-1}}^{\infty } \frac{(i_k+\frac{k}{2}+ \frac{\lambda }{2}) \left( i_k+ \frac{1}{2}(-\beta +\gamma +\delta -1+k+\lambda )\right)-\frac{q}{4}}{(i_k+ \frac{k}{2}+\frac{1}{2}+\frac{\lambda }{2} )(i_k +\frac{k}{2}+\frac{\gamma }{2}+\frac{\lambda }{2})} \right. \nonumber\\
&&\times \left. \frac{(\frac{k}{2}+\frac{\alpha }{2}+ \frac{\lambda }{2})_{i_k} ( \frac{k}{2}+1+\frac{\lambda }{2})_{i_{k-1}}(\frac{k}{2}+\frac{1}{2}+\frac{\gamma }{2}+ \frac{\lambda }{2})_{i_{k-1}}}{(\frac{k}{2}+\frac{\alpha }{2}+ \frac{\lambda }{2})_{i_{k-1}} ( \frac{k}{2}+1+\frac{\lambda }{2})_{i_k}( \frac{k}{2}+ \frac{1}{2}+\frac{\gamma}{2} +\frac{\lambda }{2})_{i_k}}\right\} \nonumber\\
&&\times \left.\sum_{i_n= i_{n-1}}^{\infty } \frac{(\frac{n}{2}+\frac{\alpha }{2}+ \frac{\lambda }{2})_{i_n} ( \frac{n}{2}+1+\frac{\lambda }{2})_{i_{n-1}}(\frac{n}{2}+\frac{1}{2}+\frac{\gamma }{2}+ \frac{\lambda }{2})_{i_{n-1}}}{(\frac{n}{2}+\frac{\alpha }{2}+ \frac{\lambda }{2})_{i_{n-1}} (\frac{n}{2}+1+\frac{\lambda }{2})_{i_n}(\frac{n}{2}+ \frac{1}{2}+ \frac{\gamma}{2} +\frac{\lambda }{2})_{i_n}} \eta ^{i_n} \Bigg\} x^n \right\}\label{eq:40020}
\end{eqnarray}
Put $c_0$= 1 as $\lambda =0$ for the first kind of independent solutions of the CHE and $\lambda = 1-\gamma $ in (\ref{eq:40020}). 
\begin{remark}
The power series expansion of the CHE of the first kind for infinite series about $x=0$ using 3TRF is
\begin{eqnarray}
 y(x)&=& H_c^{(a)}F\left(\alpha, \beta, \gamma, \delta, q; \eta =\frac{1}{2}\beta x^2  \right) \nonumber\\
&=& \sum_{i_0=0}^{\infty } \frac{(\frac{\alpha }{2} )_{i_0} }{(1 )_{i_0}(\frac{1}{2}+ \frac{\gamma}{2} )_{i_0}} \eta ^{i_0} \nonumber\\
&&+ \Bigg\{ \sum_{i_0=0}^{\infty }\frac{ i_0  \left( i_0+ \frac{1}{2}(-\beta +\gamma +\delta -1 )\right)-\frac{q}{4}}{(i_0+ \frac{1}{2} )(i_0 + \frac{\gamma }{2} )}\frac{(\frac{\alpha }{2} )_{i_0}}{(1 )_{i_0}(\frac{1}{2}+ \frac{\gamma}{2} )_{i_0}}  \sum_{i_1=i_0}^{\infty } \frac{(\frac{1}{2}+\frac{\alpha }{2} )_{i_1} (\frac{3}{2} )_{i_0}(1+\frac{\gamma }{2} )_{i_0}}{(\frac{1}{2}+\frac{\alpha }{2} )_{i_0} (\frac{3}{2} )_{i_1}(1+ \frac{\gamma}{2} )_{i_1}} \eta ^{i_1} \Bigg\} x\nonumber\\
&&+ \sum_{n=2}^{\infty } \Bigg\{ \sum_{i_0=0}^{\infty } \frac{ i_0 \left( i_0+ \frac{1}{2}(-\beta +\gamma +\delta -1 )\right)-\frac{q}{4}}{(i_0+ \frac{1}{2} )(i_0 + \frac{\gamma }{2} )}
 \frac{(\frac{\alpha }{2} )_{i_0}}{(1 )_{i_0}(\frac{1}{2}+ \frac{\gamma}{2} )_{i_0}}\nonumber\\
&&\times \prod _{k=1}^{n-1} \left\{ \sum_{i_k=i_{k-1}}^{\infty } \frac{(i_k+\frac{k}{2} ) \left( i_k+ \frac{1}{2}(-\beta +\gamma +\delta -1+k )\right)-\frac{q}{4}}{(i_k+ \frac{k}{2}+\frac{1}{2} )(i_k +\frac{k}{2}+\frac{\gamma }{2} )} \right. \nonumber\\
&&\times \left. \frac{(\frac{k}{2}+\frac{\alpha }{2} )_{i_k} ( \frac{k}{2}+1 )_{i_{k-1}}(\frac{k}{2}+\frac{1}{2}+\frac{\gamma }{2} )_{i_{k-1}}}{(\frac{k}{2}+\frac{\alpha }{2} )_{i_{k-1}} ( \frac{k}{2}+1 )_{i_k}( \frac{k}{2}+ \frac{1}{2}+\frac{\gamma}{2} )_{i_k}}\right\} \nonumber\\
&&\times  \sum_{i_n= i_{n-1}}^{\infty } \frac{(\frac{n}{2}+\frac{\alpha }{2} )_{i_n} ( \frac{n}{2}+1 )_{i_{n-1}}(\frac{n}{2}+\frac{1}{2}+\frac{\gamma }{2} )_{i_{n-1}}}{(\frac{n}{2}+\frac{\alpha }{2} )_{i_{n-1}} (\frac{n}{2}+1 )_{i_n}(\frac{n}{2}+ \frac{1}{2}+ \frac{\gamma}{2} )_{i_n}} \eta ^{i_n} \Bigg\} x^n  \label{eq:40021}
\end{eqnarray}
\end{remark}
\begin{remark}
The power series expansion of the CHE of the second kind for infinite series for infinite series about $x=0$ using 3TRF is
\begin{eqnarray}
 y(x)&=&  H_c^{(a)}S\left(\alpha, \beta, \gamma, \delta, q; \eta =\frac{1}{2}\beta x^2  \right)\nonumber\\
&=& x^{1-\gamma  } \left\{\sum_{i_0=0}^{\infty } \frac{(\frac{\alpha }{2}+\frac{1}{2}-\frac{\gamma }{2})_{i_0} }{(\frac{3}{2}-\frac{\gamma  }{2})_{i_0}(1)_{i_0}} \eta ^{i_0}\right.\nonumber\\
&&+ \Bigg\{ \sum_{i_0=0}^{\infty }\frac{(i_0+ \frac{1}{2}-\frac{\gamma }{2}) \left( i_0+ \frac{1}{2}(-\beta +\delta )\right)-\frac{q}{4}}{(i_0+ 1-\frac{\gamma }{2})(i_0 + \frac{1}{2} )}\frac{(\frac{\alpha }{2}+\frac{1}{2}-\frac{\gamma }{2})_{i_0}}{( \frac{3}{2}-\frac{\gamma }{2})_{i_0}(1)_{i_0}}\nonumber\\
&&\times \sum_{i_1=i_0}^{\infty } \frac{(1+\frac{\alpha }{2}-\frac{\gamma }{2})_{i_1} (2-\frac{\gamma }{2})_{i_0}( \frac{3}{2})_{i_0}}{(1+\frac{\alpha }{2}-\frac{\gamma }{2})_{i_0} (2-\frac{\gamma }{2})_{i_1}(\frac{3}{2})_{i_1}} \eta ^{i_1} \Bigg\} x\nonumber\\
&&+ \sum_{n=2}^{\infty } \Bigg\{ \sum_{i_0=0}^{\infty } \frac{(i_0+ \frac{1}{2}-\frac{\gamma }{2}) \left( i_0+ \frac{1}{2}(-\beta +\delta )\right)-\frac{q}{4}}{(i_0+1-\frac{\gamma }{2})(i_0 + \frac{1}{2} )}
 \frac{(\frac{\alpha }{2}+\frac{1}{2}-\frac{\gamma }{2})_{i_0}}{( \frac{3}{2}-\frac{\gamma }{2})_{i_0}(1)_{i_0}}\nonumber\\
&&\times \prod _{k=1}^{n-1} \left\{ \sum_{i_k=i_{k-1}}^{\infty } \frac{(i_k+\frac{k}{2}+ \frac{1}{2}-\frac{\gamma }{2}) \left( i_k+ \frac{1}{2}(-\beta +\delta +k )\right)-\frac{q}{4}}{(i_k+ \frac{k}{2}+1-\frac{\gamma }{2} )(i_k +\frac{k}{2}+\frac{1}{2} )}\right. \nonumber\\
&&\times \left. \frac{(\frac{k}{2}+\frac{\alpha }{2}+ \frac{1}{2}-\frac{\gamma }{2})_{i_k} ( \frac{k}{2}+ \frac{3}{2}-\frac{\gamma }{2})_{i_{k-1}}(\frac{k}{2}+1)_{i_{k-1}}}{(\frac{k}{2}+\frac{\alpha }{2}+ \frac{1}{2}-\frac{\gamma }{2})_{i_{k-1}} ( \frac{k}{2}+ \frac{3}{2}-\frac{\gamma }{2})_{i_k}( \frac{k}{2}+ 1)_{i_k}}\right\} \nonumber\\
&&\times \left.\sum_{i_n= i_{n-1}}^{\infty } \frac{(\frac{n}{2}+\frac{\alpha }{2}+ \frac{1}{2}-\frac{\gamma }{2})_{i_n} ( \frac{n}{2}+ \frac{3}{2}-\frac{\gamma }{2})_{i_{n-1}}(\frac{n}{2}+1)_{i_{n-1}}}{(\frac{n}{2}+\frac{\alpha }{2}+ \frac{1}{2}-\frac{\gamma }{2})_{i_{n-1}} (\frac{n}{2} +\frac{3}{2}-\frac{\gamma }{2})_{i_n}(\frac{n}{2}+1)_{i_n}} \eta ^{i_n} \Bigg\} x^n \right\}\label{eq:40022}
\end{eqnarray}
\end{remark}
It is required that $\gamma \ne 0,-1,-2,\cdots$ for the first kind of independent solution of the CHF for all cases. Because if it does not, its solution will be divergent. And it is required that $\gamma \ne 2,3,4,\cdots$ for the second kind of independent solution of the CHF for all cases.
\subsection{Integral formalism}
The CHE could not be constructed in the definite or contour integral form of any well-known simple function because of a 3-term recursive relation between successive coefficients in its power series. The three term recurrence relation in power series of the CHE creates mathematical difficulty to be analyzed it into direct or contour integral formalism.

Instead in earlier literature the integral relations of the CHE were constructed by using Fredholm integral equations; such integral relationships express one analytic solution in terms of another analytic solution such as a confluent hypergeometric function with a branch point at zero. \cite{Ronv1995} There are many other forms of integral relations in the CHE. \cite{Lamb1934,Erde1942,Slee1969,Abra1976,Schm1979,Kaza1986}

Now I consider the combined definite and contour integral forms of the CHF by using 3TRF.
 Expressing the CHF in integral forms resulting in a precise and simplified transformation of the CHF to other well-known special functions such as confluent hypergeometric function.  
\subsubsection{Polynomial of type 1}
Let's investigate the integral representation for polynomial of type 1. There is a generalized hypergeometric function which is
\begin{eqnarray}
I_l &=& \sum_{i_l= i_{l-1}}^{\alpha _l} \frac{(-\alpha _l)_{i_l} ( 1+\frac{l}{2}+ \frac{\lambda}{2} )_{i_{l-1}}(\frac{1}{2}+\frac{\gamma }{2}+\frac{l}{2}+\frac{\lambda }{2})_{i_{l-1}}}{(-\alpha _l)_{i_{l-1}} (1+\frac{l}{2}+ \frac{\lambda}{2} )_{i_l} (\frac{1}{2}+\frac{\gamma }{2}+\frac{l}{2}+\frac{\lambda }{2})_l} \eta^{i_l}\nonumber\\
&=& \eta ^{i_{l-1}} 
\sum_{j=0}^{\infty } \frac{B(i_{l-1}+\frac{l}{2}+\frac{\lambda}{2},j+1) B(i_{l-1}+\frac{l}{2}-\frac{1}{2}+\frac{\gamma }{2}+\frac{\lambda }{2},j+1)(i_{l-1}-\alpha _l)_j}{(i_{l-1}+\frac{l}{2}+\frac{\lambda}{2})^{-1}(i_{l-1}+\frac{l}{2}-\frac{1}{2}+\frac{\gamma}{2} +\frac{\lambda}{2} )^{-1}(1)_j \;j!} \eta ^j \hspace{2cm}\label{eq:40023}
\end{eqnarray}
By using integral form of beta function,
\begin{subequations}
\begin{equation}
B\left(i_{l-1}+\frac{l}{2}+\frac{\lambda}{2},j+1\right)= \int_{0}^{1} dt_l\;t_l^{i_{l-1}+\frac{l}{2}-1+\frac{\lambda }{2} } (1-t_l)^j \label{eq:40024a}
\end{equation}
\begin{equation}
B\left(i_{l-1}+\frac{l}{2}-\frac{1}{2}+\frac{\gamma }{2}+\frac{\lambda }{2},j+1\right)= \int_{0}^{1} du_l\;u_l^{i_{l-1}+\frac{l}{2}-\frac{3}{2}+\frac{\gamma }{2} +\frac{\lambda }{2}} (1-u_l)^j\label{eq:40024b}
\end{equation}
\end{subequations}
Substitute (\ref{eq:40024a}) and (\ref{eq:40024b}) into (\ref{eq:40023}). And divide $(i_{l-1}+\frac{l}{2}+\frac{\lambda}{2})(i_{l-1}+\frac{l}{2}-\frac{1}{2}+\frac{\gamma}{2} +\frac{\lambda}{2} )$ into $I_l$.
\begin{eqnarray}
K_l&=& \frac{1}{(i_{l-1}+\frac{l}{2}+\frac{\lambda}{2})(i_{l-1}+\frac{l}{2}-\frac{1}{2}+\frac{\gamma}{2} +\frac{\lambda}{2} )}
\sum_{i_l= i_{l-1}}^{\alpha _l} \frac{(-\alpha _l)_{i_l} ( 1+\frac{l}{2}+ \frac{\lambda}{2} )_{i_{l-1}}(\frac{1}{2}+\frac{\gamma }{2}+\frac{l}{2}+\frac{\lambda }{2})_{i_{l-1}}}{(-\alpha _l)_{i_{l-1}} (1+\frac{l}{2}+ \frac{\lambda}{2} )_{i_l} (\frac{1}{2}+\frac{\gamma }{2}+\frac{l}{2}+\frac{\lambda }{2})i_l} \eta^{i_l}\nonumber\\
&=& \int_{0}^{1} dt_l\;t_l^{\frac{l}{2}-1+\frac{\lambda }{2}} \int_{0}^{1} du_l\;u_l^{\frac{l}{2}-\frac{3}{2}+\frac{\gamma }{2}+\frac{\lambda }{2}} (\eta t_l u_l)^{i_{l-1}}
 \sum_{j=0}^{\infty } \frac{(i_{l-1}-\alpha _l)_j }{(1)_j \;j!} [\eta (1-t_l)(1-u_l)]^j \nonumber 
\end{eqnarray}
Confluent hypergeometric polynomial of the first kind is defined by
\begin{equation}
F_{\alpha _0}(\gamma ;z) = \frac{\Gamma (\alpha _0+\gamma )}{\Gamma (\gamma )}\sum_{j=0}^{\alpha _0}\frac{(-\alpha _0)_j}{(\gamma )_j j!}z^j= \frac{\Gamma ( \alpha _0+1) }{2\pi i}\oint dv_l\;\frac{\exp\left(-\frac{z v_l}{(1-v_l)} \right)}{v_l^{\alpha _0+1} (1-v_l)^{\gamma }}\label{eq:40025}
\end{equation}
replaced $\alpha _0 $, $\gamma $ and z by $ \alpha _l-i_{l-1}$, 1 and $\eta (1-t_l)(1-u_l)$ in (\ref{eq:40025})
\begin{equation}
  \sum_{j=0}^{\infty }\frac{(i_{l-1}-\alpha _l)_j }{(1)_j \;j!} [\eta (1-t_l)(1-u_l)]^j = \frac{1}{2\pi i}\oint dv_l\;\frac{\exp\left(-\frac{ v_l}{(1-v_l)}\eta (1-t_l)(1-u_l) \right)}{v_l^{\alpha _l-i_{l-1}+1} (1-v_l)}\label{eq:40026}
\end{equation}
Substitute (\ref{eq:40026}) into $K_l$.
\begin{eqnarray}
K_l&=& \frac{1}{(i_{l-1}+\frac{l}{2}+\frac{\lambda}{2})(i_{l-1}+\frac{l}{2}-\frac{1}{2}+\frac{\gamma}{2} +\frac{\lambda}{2} )}
\sum_{i_l= i_{l-1}}^{\alpha _l} \frac{(-\alpha _l)_{i_l} ( 1+\frac{l}{2}+ \frac{\lambda}{2} )_{i_{l-1}}(\frac{1}{2}+\frac{\gamma }{2}+\frac{l}{2}+\frac{\lambda }{2})_{i_{l-1}}}{(-\alpha _l)_{i_{l-1}} (1+\frac{l}{2}+ \frac{\lambda}{2} )_{i_l} (\frac{1}{2}+\frac{\gamma }{2}+\frac{l}{2}+\frac{\lambda }{2})i_l} \eta^{i_l}\nonumber\\
&=& \int_{0}^{1} dt_l\;t_l^{\frac{l}{2}-1+\frac{\lambda }{2}} \int_{0}^{1} du_l\;u_l^{\frac{l}{2}-\frac{3}{2}+\frac{\gamma }{2}+\frac{\lambda }{2}} \frac{1}{2\pi i}\oint dv_l\;\frac{\exp\left(-\frac{ v_l}{(1-v_l)}\eta (1-t_l)(1-u_l) \right)}{v_l^{\alpha _l +1} (1-v_l)}(\eta t_l u_l v_l)^{i_{l-1}}\hspace{1.5cm}\label{eq:40027}
\end{eqnarray}
Substitute (\ref{eq:40027}) into (\ref{eq:40016}) where $l=1,2,3,\cdots$; apply $K_1$ into the second summation of sub-power series $y_1(x)$, apply $K_2$ into the third summation and $K_1$ into the second summation of sub-power series $y_2(x)$, apply $K_3$ into the forth summation, $K_2$ into the third summation and $K_1$ into the second summation of sub-power series $y_3(x)$, etc.\footnote{$y_1(x)$ means the sub-power series in (\ref{eq:40016}) contains one term of $A_n's$, $y_2(x)$ means the sub-power series in (\ref{eq:40016}) contains two terms of $A_n's$, $y_3(x)$ means the sub-power series in (\ref{eq:40016}) contains three terms of $A_n's$, etc.}
\begin{theorem}
The general representation in the form of integral of the CHP of type 1 about $x=0$ is given by
\begin{eqnarray}
 y(x)&=& \sum_{n=0}^{\infty } y_n(x)= y_0(x)+ y_1(x)+ y_2(x)+ y_3(x)+\cdots \nonumber\\
&=&  c_0 x^{\lambda } \left\{ \sum_{i_0=0}^{\alpha _0}\frac{(-\alpha _0)_{i_0} }{(1+\frac{\lambda }{2})_{i_0}(\frac{1}{2}+ \frac{\gamma }{2} +\frac{\lambda }{2})_{i_0}} \eta^{i_0}\right. \nonumber\\
&&+ \left.\sum_{n=1}^{\infty } \Bigg\{\prod _{k=0}^{n-1} \Bigg\{ \int_{0}^{1} dt_{n-k}\;t_{n-k}^{\frac{1}{2}(n-k-2+\lambda )} \int_{0}^{1} du_{n-k}\;u_{n-k}^{\frac{1}{2}(n-k-3+\gamma +\lambda )}\right. \nonumber\\
&&\times  \frac{1}{2\pi i}  \oint dv_{n-k} \frac{\exp\left( -\frac{v_{n-k}}{(1-v_{n-k})}w_{n-k+1,n}(1-t_{n-k})(1-u_{n-k})\right) }{v_{n-k}^{\alpha_{n-k}+1} (1-v_{n-k})} \nonumber\\
&&\times \Bigg( w_{n-k,n}^{-\frac{1}{2}(n-k-1+\lambda )}\left( w_{n-k,n} \partial _{ w_{n-k,n}}\right) w_{n-k,n}^{\frac{1}{2}(\beta -\gamma -\delta +1)} \left( w_{n-k,n} \partial _{ w_{n-k,n}}\right) w_{n-k,n}^{\frac{1}{2}(-\beta +\gamma +\delta +\lambda +n-k-2)} -\frac{q}{4}\Bigg) \Bigg\}\nonumber\\
&&\times\left.\sum_{i_0=0}^{\alpha _0}\frac{(-\alpha _0)_{i_0} }{(1+\frac{\lambda }{2})_{i_0}(\frac{1}{2}+ \frac{\gamma }{2} +\frac{\lambda }{2})_{i_0}} w_{1,n}^{i_0}\Bigg\} x^n \right\} \label{eq:40028}
\end{eqnarray}
where
\begin{equation} w_{a,b}=
\begin{cases} \displaystyle {\eta \prod_{l=a}^{b} t_l u_l v_l\;\;\mbox{where}\; a\leq b}\cr
\eta \;\;\mbox{only}\;\mbox{if}\; a>b
\end{cases}\nonumber 
\end{equation}
In the above, the first sub-integral form contains one term of $A_n's$, the second one contains two terms of $A_n$'s, the third one contains three terms of $A_n$'s, etc.
\end{theorem}
\begin{proof} 
According to (\ref{eq:40016}), 
\begin{equation}
 y(x)= \sum_{n=0}^{\infty }y_n(x) = y_0(x)+ y_1(x)+ y_2(x)+y_3(x)+\cdots \label{eq:40029}
\end{equation}
In the above, sub-power series $y_0(x) $, $y_1(x)$, $y_2(x)$ and $y_3(x)$ of the CHP of type 1 using 3TRF about $x=0$ are
\begin{subequations}
\begin{equation}
 y_0(x)= c_0 x^{\lambda } \sum_{i_0=0}^{\alpha _0} \frac{(-\alpha _0)_{i_0}}{(1+\frac{\lambda}{2} )_{i_0}(\frac{1}{2}+\frac{\gamma }{2}+\frac{\lambda }{2} )_{i_0}} \eta ^{i_0} \label{eq:40030a}
\end{equation}
\begin{eqnarray}
 y_1(x) &=& c_0 x^{\lambda } \Bigg\{\sum_{i_0=0}^{\alpha _0} \frac{(i_0+\frac{\lambda }{2})(i_0+\frac{1}{2}(-\beta +\gamma +\delta -1+\lambda ))-\frac{q}{4}}{(i_0+\frac{1}{2}+\frac{\lambda }{2})(i_0+\frac{\gamma }{2}+\frac{\lambda }{2})} \frac{(-\alpha _0)_{i_0}}{(1+\frac{\lambda}{2} )_{i_0}(\frac{1}{2}+\frac{\gamma }{2}+\frac{\lambda }{2} )_{i_0}}\nonumber\\
&&\times  \sum_{i_1=i_0}^{\alpha _1}  \frac{(-\alpha _1)_{i_1} (\frac{3}{2}+\frac{\lambda}{2} )_{i_0}(1+\frac{\gamma }{2} +\frac{\lambda }{2})_{i_0}}{(-\alpha _1)_{i_0} (\frac{3}{2}+\frac{\lambda}{2})_{i_1}(1+\frac{\gamma }{2} +\frac{\lambda }{2})_{i_1}} \eta ^{i_1} \Bigg\}x  \label{eq:40030b}
\end{eqnarray}
\begin{eqnarray}
 y_2(x) &=& c_0 x^{\lambda } \Bigg\{\sum_{i_0=0}^{\alpha _0} \frac{(i_0+\frac{\lambda }{2})(i_0+\frac{1}{2}(-\beta +\gamma +\delta -1+\lambda ))-\frac{q}{4}}{(i_0+\frac{1}{2}+\frac{\lambda }{2})(i_0+\frac{\gamma }{2}+\frac{\lambda }{2})} \frac{(-\alpha _0)_{i_0}}{(1+\frac{\lambda}{2} )_{i_0}(\frac{1}{2}+\frac{\gamma }{2}+\frac{\lambda }{2} )_{i_0}} \nonumber\\
&&\times  \sum_{i_1=i_0}^{\alpha _1} \frac{(i_1+\frac{1}{2}+\frac{\lambda }{2})(i_1+\frac{1}{2}(-\beta +\gamma +\delta +\lambda ))-\frac{q}{4} }{ (i_1+1+\frac{\lambda }{2})(i_1+\frac{1}{2}+\frac{\gamma }{2}+\frac{\lambda }{2})} \frac{(-\alpha _1)_{i_1} (\frac{3}{2}+\frac{\lambda}{2} )_{i_0}(1+\frac{\gamma }{2} +\frac{\lambda }{2})_{i_0}}{(-\alpha _1)_{i_0} (\frac{3}{2}+\frac{\lambda}{2})_{i_1}(1+\frac{\gamma }{2} +\frac{\lambda }{2})_{i_1}}   \nonumber\\
&&\times \sum_{i_2=i_1}^{\alpha _2} \frac{(-\alpha _2)_{i_2} (2+\frac{\lambda}{2} )_{i_1}(\frac{3}{2}+\frac{\gamma }{2} +\frac{\lambda }{2})_{i_1}}{(-\alpha _2)_{i_1} (2+\frac{\lambda}{2})_{i_2}(\frac{3}{2}+\frac{\gamma }{2} +\frac{\lambda }{2})_{i_2}} \eta ^{i_2} \Bigg\} x^2  \label{eq:40030c}
\end{eqnarray}
\begin{eqnarray}
 y_3(x) &=&  c_0 x^{\lambda } \Bigg\{\sum_{i_0=0}^{\alpha _0} \frac{(i_0+\frac{\lambda }{2})(i_0+\frac{1}{2}(-\beta +\gamma +\delta -1+\lambda ))-\frac{q}{4}}{(i_0+\frac{1}{2}+\frac{\lambda }{2})(i_0+\frac{\gamma }{2}+\frac{\lambda }{2})} \frac{(-\alpha _0)_{i_0}}{(1+\frac{\lambda}{2} )_{i_0}(\frac{1}{2}+\frac{\gamma }{2}+\frac{\lambda }{2} )_{i_0}} \nonumber\\
&&\times  \sum_{i_1=i_0}^{\alpha _1} \frac{(i_1+\frac{1}{2}+\frac{\lambda }{2})(i_1+\frac{1}{2}(-\beta +\gamma +\delta +\lambda ))-\frac{q}{4} }{ (i_1+1+\frac{\lambda }{2})(i_1+\frac{1}{2}+\frac{\gamma }{2}+\frac{\lambda }{2})} \frac{(-\alpha _1)_{i_1} (\frac{3}{2}+\frac{\lambda}{2} )_{i_0}(1+\frac{\gamma }{2} +\frac{\lambda }{2})_{i_0}}{(-\alpha _1)_{i_0} (\frac{3}{2}+\frac{\lambda}{2})_{i_1}(1+\frac{\gamma }{2} +\frac{\lambda }{2})_{i_1}} \nonumber\\
&&\times \sum_{i_2=i_1}^{\alpha _2} \frac{(i_2+1+\frac{\lambda }{2})(i_1+\frac{1}{2}(-\beta +\gamma +\delta +1+\lambda ))-\frac{q}{4} }{ (i_2+\frac{3}{2}+\frac{\lambda }{2})(i_2+1+\frac{\gamma }{2}+\frac{\lambda }{2})} \frac{(-\alpha _2)_{i_2} (2+\frac{\lambda}{2} )_{i_1}(\frac{3}{2}+\frac{\gamma }{2} +\frac{\lambda }{2})_{i_1}}{(-\alpha _2)_{i_1} (2+\frac{\lambda}{2})_{i_2}(\frac{3}{2}+\frac{\gamma }{2} +\frac{\lambda }{2})_{i_2}} \nonumber\\
&&\times \sum_{i_3=i_2}^{\alpha _3} \frac{(-\alpha _3)_{i_3} (\frac{5}{2}+\frac{\lambda}{2} )_{i_2}(2+\frac{\gamma }{2} +\frac{\lambda }{2})_{i_2}}{(-\alpha _3)_{i_2} (\frac{5}{2}+\frac{\lambda}{2})_{i_3}(2+\frac{\gamma }{2} +\frac{\lambda }{2})_{i_3}} \eta ^{i_3} \Bigg\} x^3  \label{eq:40030d} 
\end{eqnarray}
\end{subequations}
Put $l=1$ in (\ref{eq:40027}). Take the new (\ref{eq:40027}) into (\ref{eq:40030b}).
\begin{eqnarray}
y_1(x) &=& \int_{0}^{1} dt_1\;t_1^{\frac{1}{2}(-1+\lambda )} \int_{0}^{1} du_1\;u_1^{\frac{1}{2}(-2+\gamma +\lambda )} \frac{1}{2\pi i} \oint dv_1 \; \frac{\exp\left(-\frac{v_1}{(1-v_1)}(1-t_1)(1-u_1)\eta \right)}{v_1^{\alpha _1+1}(1-v_1)}\nonumber\\
&&\times \left( w_{1,1}^{-\frac{\lambda }{2}} \left( w_{1,1} \partial_{w_{1,1}}\right) w_{1,1}^{\frac{1}{2}(\beta -\gamma -\delta +1)} \left( w_{1,1} \partial_{w_{1,1}}\right) w_{1,1}^{\frac{1}{2}(-\beta +\gamma +\delta -1+\lambda )}-\frac{q}{4}\right) \nonumber\\
&&\times \Bigg\{ c_0 x^{\lambda }  \sum_{i_0=0}^{\alpha _0} \frac{(-\alpha _0)_{i_0}}{(1+\frac{\lambda}{2} )_{i_0}(\frac{1}{2}+\frac{\gamma }{2}+\frac{\lambda }{2} )_{i_0}} w_{1,1}^{i_0} \Bigg\} x\label{eq:40031}
\end{eqnarray}
where
\begin{equation}
w_{1,1} = \eta \prod_{l=1}^{1} t_l u_l v_l\nonumber
\end{equation}
Put $l=2$ in (\ref{eq:40027}). Take the new (\ref{eq:40027}) into (\ref{eq:40030c}).
\begin{eqnarray}
y_2(x) &=& c_0 x^{\lambda } \int_{0}^{1} dt_2\;t_2^{\frac{\lambda }{2}} \int_{0}^{1} du_2\;u_2^{\frac{1}{2}(-1+\gamma +\lambda )} \frac{1}{2\pi i} \oint dv_2 \; \frac{\exp\left(-\frac{v_2}{(1-v_2)}(1-t_2)(1-u_2)\eta \right)}{v_2^{\alpha _2+1}(1-v_2)} \nonumber\\
&&\times \left( w_{2,2}^{-\frac{1}{2}(1+\lambda )} \left( w_{2,2} \partial_{w_{2,2}}\right) w_{2,2}^{\frac{1}{2}(\beta -\gamma -\delta +1)} \left( w_{2,2} \partial_{w_{2,2}}\right) w_{2,2}^{\frac{1}{2}(-\beta +\gamma +\delta +\lambda )}-\frac{q}{4}\right) \nonumber\\
&&\times \Bigg\{\sum_{i_0=0}^{\alpha _0} \frac{(i_0+\frac{\lambda }{2})(i_0+\frac{1}{2}(-\beta +\gamma +\delta -1+\lambda ))-\frac{q}{4}}{(i_0+\frac{1}{2}+\frac{\lambda }{2})(i_0+\frac{\gamma }{2}+\frac{\lambda }{2})} \frac{(-\alpha _0)_{i_0}}{(1+\frac{\lambda}{2} )_{i_0}(\frac{1}{2}+\frac{\gamma }{2}+\frac{\lambda }{2} )_{i_0}} \nonumber\\
&&\times  \sum_{i_1=i_0}^{\alpha _1} \frac{(-\alpha _1)_{i_1} (\frac{3}{2}+\frac{\lambda}{2} )_{i_0}(1+\frac{\gamma }{2} +\frac{\lambda }{2})_{i_0}}{(-\alpha _1)_{i_0} (\frac{3}{2}+\frac{\lambda}{2})_{i_1}(1+\frac{\gamma }{2} +\frac{\lambda }{2})_{i_1}} w_{2,2}^{i_1} \Bigg\} x^2 \label{eq:40032}
\end{eqnarray}
where
\begin{equation}
w_{2,2} = \eta \prod_{l=2}^{2} t_l u_l v_l\nonumber
\end{equation}
Put $l=1$ and $\eta = w_{2,2}$ in (\ref{eq:40027}). Take the new (\ref{eq:40027}) into (\ref{eq:40032}).
\begin{eqnarray}
y_2(x) &=& \int_{0}^{1} dt_2\;t_2^{\frac{\lambda }{2}} \int_{0}^{1} du_2\;u_2^{\frac{1}{2}(-1+\gamma +\lambda )} \frac{1}{2\pi i} \oint dv_2 \; \frac{\exp\left(-\frac{v_2}{(1-v_2)}(1-t_2)(1-u_2)\eta \right)}{v_2^{\alpha _2+1}(1-v_2)} \nonumber\\
&&\times \left( w_{2,2}^{-\frac{1}{2}(1+\lambda )} \left( w_{2,2} \partial_{w_{2,2}}\right) w_{2,2}^{\frac{1}{2}(\beta -\gamma -\delta +1)} \left( w_{2,2} \partial_{w_{2,2}}\right) w_{2,2}^{\frac{1}{2}(-\beta +\gamma +\delta +\lambda )}-\frac{q}{4}\right) \nonumber\\ 
&&\times \int_{0}^{1} dt_1\;t_1^{\frac{1}{2}(-1+\lambda )} \int_{0}^{1} du_1\;u_1^{\frac{1}{2}(-2+\gamma +\lambda )} \frac{1}{2\pi i} \oint dv_1 \; \frac{\exp\left(-\frac{v_1}{(1-v_1)}(1-t_1)(1-u_1)w_{2,2} \right)}{v_1^{\alpha _1+1}(1-v_1)} \nonumber\\
&&\times \left( w_{1,2}^{-\frac{\lambda }{2}} \left( w_{1,2} \partial_{w_{1,2}}\right) w_{1,2}^{\frac{1}{2}(\beta -\gamma -\delta +1)} \left( w_{1,2} \partial_{w_{1,2}}\right) w_{1,2}^{\frac{1}{2}(-\beta +\gamma +\delta -1+\lambda )}-\frac{q}{4}\right) \nonumber\\
&&\times \left\{ c_0 x^{\lambda } \sum_{i_0=0}^{\alpha _0} \frac{(-\alpha _0)_{i_0}}{(1+\frac{\lambda}{2} )_{i_0}(\frac{1}{2}+\frac{\gamma }{2}+\frac{\lambda }{2} )_{i_0}} w_{1,2}^{i_0} \right\} x^2 \label{eq:40033}
\end{eqnarray}
where
\begin{equation}
w_{1,2} = \eta \prod_{l=1}^{2} t_l u_l v_l \nonumber
\end{equation}
By using similar process for the previous cases of integral forms of $y_1(x)$ and $y_2(x)$, the integral form of sub-power series expansion of $y_3(x)$ is
\begin{eqnarray}
y_3(x) &=& \int_{0}^{1} dt_3\;t_3^{\frac{1}{2}(1+\lambda )} \int_{0}^{1} du_3\;u_3^{\frac{1}{2}(\gamma +\lambda )} \frac{1}{2\pi i} \oint dv_3 \; \frac{\exp\left(-\frac{v_3}{(1-v_3)}(1-t_3)(1-u_3)\eta \right)}{v_3^{\alpha _3+1}(1-v_3)} \nonumber\\
&&\times \left( w_{3,3}^{-\frac{1}{2}(2+\lambda )} \left( w_{3,3} \partial_{w_{3,3}}\right) w_{3,3}^{\frac{1}{2}(\beta -\gamma -\delta +1)} \left( w_{3,3} \partial_{w_{3,3}}\right) w_{3,3}^{\frac{1}{2}(-\beta +\gamma +\delta +1+\lambda )}-\frac{q}{4}\right) \nonumber\\ 
&&\times \int_{0}^{1} dt_2\;t_2^{\frac{\lambda }{2}} \int_{0}^{1} du_2\;u_2^{\frac{1}{2}(-1+\gamma +\lambda )} \frac{1}{2\pi i} \oint dv_2 \; \frac{\exp\left(-\frac{v_2}{(1-v_2)}(1-t_2)(1-u_2)w_{3,3} \right)}{v_2^{\alpha _2+1}(1-v_2)} \nonumber\\
&&\times \left( w_{2,3}^{-\frac{1}{2}(1+\lambda )} \left( w_{2,3} \partial_{w_{2,3}}\right) w_{2,3}^{\frac{1}{2}(\beta -\gamma -\delta +1)} \left( w_{2,3} \partial_{w_{2,3}}\right) w_{2,3}^{\frac{1}{2}(-\beta +\gamma +\delta +\lambda )}-\frac{q}{4}\right) \nonumber\\
&&\times \int_{0}^{1} dt_1\;t_1^{\frac{1}{2}(-1+\lambda )} \int_{0}^{1} du_1\;u_1^{\frac{1}{2}(-2+\gamma +\lambda )} \frac{1}{2\pi i} \oint dv_1 \; \frac{\exp\left(-\frac{v_1}{(1-v_1)}(1-t_1)(1-u_1)w_{2,3} \right)}{v_1^{\alpha _1+1}(1-v_1)} \nonumber\\
&&\times \left( w_{1,3}^{-\frac{\lambda }{2}} \left( w_{1,3} \partial_{w_{1,3}}\right) w_{1,3}^{\frac{1}{2}(\beta -\gamma -\delta +1)} \left( w_{1,3} \partial_{w_{1,3}}\right) w_{1,3}^{\frac{1}{2}(-\beta +\gamma +\delta -1+\lambda )}-\frac{q}{4}\right) \nonumber\\
&&\times \left\{ c_0 x^{\lambda } \sum_{i_0=0}^{\alpha _0} \frac{(-\alpha _0)_{i_0}}{(1+\frac{\lambda}{2} )_{i_0}(\frac{1}{2}+\frac{\gamma }{2}+\frac{\lambda }{2} )_{i_0}} w_{1,3}^{i_0} \right\} x^3 \label{eq:40034}
\end{eqnarray}
where
\begin{equation}
 w_{3,3} = \eta \prod_{l=3}^{3} t_l u_l v_l \hspace{1cm}
w_{2,3} = \eta \prod_{l=2}^{3} t_l u_l v_l \hspace{1cm}
w_{1,3} = \eta \prod_{l=1}^{3} t_l u_l v_l \nonumber
\end{equation}
By repeating this process for all higher terms of integral forms of sub-summation $y_m(x)$ terms where $m \geq 4$, we obtain every integral forms of $y_m(x)$ terms. 
Substitute (\ref{eq:40030a}), (\ref{eq:40031}), (\ref{eq:40033}), (\ref{eq:40034}) and including all integral forms of $y_m(x)$ terms where $m \geq 4$ into (\ref{eq:40029}). 
\qed
\end{proof}
Put $c_0$= 1 as $\lambda =0$  for the first kind of independent solutions of the CHE and $\lambda = 1-\gamma $ in (\ref{eq:40028}).  
\begin{remark}
The integral representation of the CHE of the first kind for polynomial of type 1 about $x=0$ as $\alpha = -(2 \alpha _j+j) $ where $j,\alpha _j \in \mathbb{N}_{0}$ is
\begin{eqnarray}
 y(x)&=& H_c^{(a)}F_{\alpha _j}\left(\alpha = -(2 \alpha _j+j), \beta, \gamma, \delta, q; \eta =\frac{1}{2}\beta x^2  \right)\nonumber\\
&=& _1F_1 \left(-\alpha _0; \frac{1}{2}+ \frac{\gamma }{2}; \eta \right) + \sum_{n=1}^{\infty } \Bigg\{\prod _{k=0}^{n-1} \Bigg\{ \int_{0}^{1} dt_{n-k}\;t_{n-k}^{\frac{1}{2}(n-k-2)} \int_{0}^{1} du_{n-k}\;u_{n-k}^{\frac{1}{2}(n-k-3+\gamma )} \nonumber\\
&&\times  \frac{1}{2\pi i}  \oint dv_{n-k} \frac{\exp\left( -\frac{v_{n-k}}{(1-v_{n-k})}w_{n-k+1,n}(1-t_{n-k})(1-u_{n-k})\right) }{v_{n-k}^{\alpha_{n-k}+1} (1-v_{n-k})} \nonumber\\
&&\times \Bigg( w_{n-k,n}^{-\frac{1}{2}(n-k-1)}\left( w_{n-k,n} \partial _{ w_{n-k,n}}\right) w_{n-k,n}^{\frac{1}{2}(\beta -\gamma -\delta +1)} \left( w_{n-k,n} \partial _{ w_{n-k,n}}\right) w_{n-k,n}^{\frac{1}{2}(-\beta +\gamma +\delta +n-k-2)} -\frac{q}{4}\Bigg) \Bigg\}\nonumber\\
&&\times _1F_1 \left(-\alpha _0; \frac{1}{2}+ \frac{\gamma }{2}; w_{1,n} \right) \Bigg\} x^n  \label{eq:40035}
\end{eqnarray}
\end{remark}
\begin{remark}
The integral representation of the CHE of the second kind for polynomial of type 1 about $x=0$ as $\alpha = -(2 \alpha _j+j +1-\gamma ) $ where $j,\alpha _j \in \mathbb{N}_{0}$ is
\begin{eqnarray}
 y(x)&=& H_c^{(a)}S_{\alpha _j}\left(\alpha = -(2 \alpha _j+j+1-\gamma ), \beta, \gamma, \delta, q; \eta =\frac{1}{2}\beta x^2  \right)\nonumber\\
&=&  x^{1-\gamma } \left\{ _1F_1 \left(-\alpha _0; \frac{3}{2}-\frac{\gamma }{2}; \eta \right) 
+ \sum_{n=1}^{\infty } \Bigg\{\prod _{k=0}^{n-1} \Bigg\{ \int_{0}^{1} dt_{n-k}\;t_{n-k}^{\frac{1}{2}(n-k-1-\gamma )} \int_{0}^{1} du_{n-k}\;u_{n-k}^{\frac{1}{2}(n-k-2)}\right. \nonumber\\
&&\times  \frac{1}{2\pi i}  \oint dv_{n-k} \frac{\exp\left( -\frac{v_{n-k}}{(1-v_{n-k})}w_{n-k+1,n}(1-t_{n-k})(1-u_{n-k})\right) }{v_{n-k}^{\alpha_{n-k}+1} (1-v_{n-k})} \nonumber\\
&&\times \Bigg( w_{n-k,n}^{-\frac{1}{2}(n-k-\gamma )}\left( w_{n-k,n} \partial _{ w_{n-k,n}}\right) w_{n-k,n}^{\frac{1}{2}(\beta -\gamma -\delta +1)} \left( w_{n-k,n} \partial _{ w_{n-k,n}}\right) w_{n-k,n}^{\frac{1}{2}(-\beta + \delta + n-k-1)} -\frac{q}{4}\Bigg) \Bigg\}\nonumber\\
&&\times\left. _1F_1 \left(-\alpha _0; \frac{3}{2}-\frac{\gamma }{2}; w_{1,n} \right) \Bigg\} x^n \right\}\label{eq:40036}
\end{eqnarray}
\end{remark}
In the above, $_1F_1 \left(a; b; z \right)$ is a Kummer function of the first kind which is defined by
\begin{eqnarray}
_1F_1 \left(a; b; z \right) &=& M(a,b,z) = \sum_{n=0}^{\infty } \frac{(a)_n}{(b)_n n!} z^n = e^z M(b-a,b,-z) \nonumber\\
&=& -\frac{1}{2\pi i}\frac{\Gamma \left( 1-a\right)\Gamma \left(b\right)}{\Gamma \left(b-a\right) } \oint  dv_j\; e^{z v_j}(-v_j)^{a-1} \left( 1- v_j \right)^{b-a-1} \nonumber\\
&=&  \frac{\Gamma \left(a\right)}{2\pi i} \oint  dv_j\; e^{v_j}v_j^{-b} \left( 1-\frac{z}{v_j}\right)^{-a} \nonumber\\
&=& \frac{1}{2\pi i} \frac{\Gamma \left( 1-a\right)\Gamma \left(b\right)}{\Gamma \left(b-a\right) } \oint  dv_j\; e^{-\frac{z \;v_j}{1-v_j}}v_j^{a-1} \left( 1- v_j \right)^{-b} \label{er:40024}
\end{eqnarray}
\subsubsection{Infinite series}
Let's consider the integral representation of the CHE about $x=0$ for infinite series by applying 3TRF.
There is a generalized hypergeometric function which is written by
\begin{eqnarray}
M_l &=& \sum_{i_l= i_{l-1}}^{\infty } \frac{\left( \frac{\alpha }{2}+\frac{l}{2}+\frac{\lambda }{2} \right)_{i_l} ( 1+\frac{l}{2}+ \frac{\lambda}{2} )_{i_{l-1}}(\frac{1}{2}+\frac{\gamma }{2}+\frac{l}{2}+\frac{\lambda }{2})_{i_{l-1}}}{\left( \frac{\alpha }{2}+\frac{l}{2}+\frac{\lambda }{2} \right)_{i_{l-1}} (1+\frac{l}{2}+ \frac{\lambda}{2} )_{i_l} (\frac{1}{2}+\frac{\gamma }{2}+\frac{l}{2}+\frac{\lambda }{2})_l} \eta^{i_l}\nonumber\\
&=& \eta ^{i_{l-1}} 
\sum_{j=0}^{\infty } \frac{B(i_{l-1}+\frac{l}{2}+\frac{\lambda}{2},j+1) B(i_{l-1}+\frac{l}{2}-\frac{1}{2}+\frac{\gamma }{2}+\frac{\lambda }{2},j+1)\left( \frac{\alpha }{2}+\frac{l}{2}+\frac{\lambda }{2}+i_{l-1} \right)_j}{(i_{l-1}+\frac{l}{2}+\frac{\lambda}{2})^{-1}(i_{l-1}+\frac{l}{2}-\frac{1}{2}+\frac{\gamma}{2} +\frac{\lambda}{2} )^{-1}(1)_j \;j!} \eta ^j \hspace{1.5cm}\label{er:40023}
\end{eqnarray}
Substitute (\ref{eq:40024a}) and (\ref{eq:40024b}) into (\ref{er:40023}). And divide $(i_{l-1}+\frac{l}{2}+\frac{\lambda}{2})(i_{l-1}+\frac{l}{2}-\frac{1}{2}+\frac{\gamma}{2} +\frac{\lambda}{2} )$ into the new (\ref{er:40023}).
\begin{eqnarray}
V_l&=& \frac{1}{(i_{l-1}+\frac{l}{2}+\frac{\lambda}{2})(i_{l-1}+\frac{l}{2}-\frac{1}{2}+\frac{\gamma}{2} +\frac{\lambda}{2} )}
\sum_{i_l= i_{l-1}}^{\infty } \frac{\left( \frac{\alpha }{2}+\frac{l}{2}+\frac{\lambda }{2} \right)_{i_l} ( 1+\frac{l}{2}+ \frac{\lambda}{2} )_{i_{l-1}}(\frac{1}{2}+\frac{\gamma }{2}+\frac{l}{2}+\frac{\lambda }{2})_{i_{l-1}}}{\left( \frac{\alpha }{2}+\frac{l}{2}+\frac{\lambda }{2} \right)_{i_{l-1}} (1+\frac{l}{2}+ \frac{\lambda}{2} )_{i_l} (\frac{1}{2}+\frac{\gamma }{2}+\frac{l}{2}+\frac{\lambda }{2})_l} \eta^{i_l}\nonumber\\
&=& \int_{0}^{1} dt_l\;t_l^{\frac{l}{2}-1+\frac{\lambda }{2}} \int_{0}^{1} du_l\;u_l^{\frac{l}{2}-\frac{3}{2}+\frac{\gamma }{2}+\frac{\lambda }{2}} (\eta t_l u_l)^{i_{l-1}}
 \sum_{j=0}^{\infty } \frac{\left( \frac{\alpha }{2}+\frac{l}{2}+\frac{\lambda }{2}+i_{l-1} \right)_j }{(1)_j \;j!} [\eta (1-t_l)(1-u_l)]^j \nonumber 
\end{eqnarray}
Replace $a$, $b$ and $z$ by $ \frac{\alpha }{2}+\frac{l}{2}+\frac{\lambda }{2}+i_{l-1}$, 1 and $\eta(1-t_j)(1-u_j)$ in (\ref{er:40024}).
Take the new (\ref{er:40024}) into $V_l$
\begin{eqnarray}
V_l&=& \frac{1}{(i_{l-1}+\frac{l}{2}+\frac{\lambda}{2})(i_{l-1}+\frac{l}{2}-\frac{1}{2}+\frac{\gamma}{2} +\frac{\lambda}{2} )}
\sum_{i_l= i_{l-1}}^{\infty } \frac{\left( \frac{\alpha }{2}+\frac{l}{2}+\frac{\lambda }{2} \right)_{i_l} ( 1+\frac{l}{2}+ \frac{\lambda}{2} )_{i_{l-1}}(\frac{1}{2}+\frac{\gamma }{2}+\frac{l}{2}+\frac{\lambda }{2})_{i_{l-1}}}{\left( \frac{\alpha }{2}+\frac{l}{2}+\frac{\lambda }{2} \right)_{i_{l-1}} (1+\frac{l}{2}+ \frac{\lambda}{2} )_{i_l} (\frac{1}{2}+\frac{\gamma }{2}+\frac{l}{2}+\frac{\lambda }{2})_l} \eta^{i_l}\nonumber\\
&=& \int_{0}^{1} dt_l\;t_l^{\frac{l}{2}-1+\frac{\lambda }{2}} \int_{0}^{1} du_l\;u_l^{\frac{l}{2}-\frac{3}{2}+\frac{\gamma }{2}+\frac{\lambda }{2}} \frac{1}{2\pi i}\oint dv_l\;\frac{\exp\left(-\frac{ v_l}{(1-v_l)}\eta (1-t_l)(1-u_l) \right)}{v_l^{-\frac{1}{2}(\alpha +l+\lambda ) +1} (1-v_l)}(\eta t_l u_l v_l)^{i_{l-1}}\hspace{1cm} \label{er:40025}
\end{eqnarray}
Substitute (\ref{er:40025}) into (\ref{eq:40020}) where $l=1,2,3,\cdots$; apply $V_1$ into the second summation of sub-power series $y_1(x)$, apply $V_2$ into the third summation and $V_1$ into the second summation of sub-power series $y_2(x)$, apply $V_3$ into the forth summation, $V_2$ into the third summation and $V_1$ into the second summation of sub-power series $y_3(x)$, etc.\footnote{$y_1(x)$ means the sub-power series in (\ref{eq:40020}) contains one term of $A_n's$, $y_2(x)$ means the sub-power series in (\ref{eq:40020}) contains two terms of $A_n's$, $y_3(x)$ means the sub-power series in (\ref{eq:40020}) contains three terms of $A_n's$, etc.}
\begin{theorem}
The general expression of an integral form of the CHE for infinite series about $x=0$ using 3TRF is given by 
\begin{eqnarray}
 y(x) &=& \sum_{n=0}^{\infty } y_n(x)= y_0(x)+ y_1(x)+ y_2(x)+ y_3(x)+\cdots \nonumber\\
&=&  c_0 x^{\lambda } \left\{ \sum_{i_0=0}^{\infty }\frac{(\frac{\alpha }{2}+\frac{\lambda }{2})_{i_0} }{(1+\frac{\lambda }{2})_{i_0}(\frac{1}{2}+ \frac{\gamma }{2} +\frac{\lambda }{2})_{i_0}} \eta^{i_0}\right.\nonumber\\
&&+ \left.\sum_{n=1}^{\infty } \Bigg\{\prod _{k=0}^{n-1} \Bigg\{ \int_{0}^{1} dt_{n-k}\;t_{n-k}^{\frac{1}{2}(n-k-2+\lambda )} \int_{0}^{1} du_{n-k}\;u_{n-k}^{\frac{1}{2}(n-k-3+\gamma +\lambda )}\right. \nonumber\\
&&\times  \frac{1}{2\pi i}  \oint dv_{n-k} \frac{\exp\left( -\frac{v_{n-k}}{(1-v_{n-k})}w_{n-k+1,n}(1-t_{n-k})(1-u_{n-k})\right) }{v_{n-k}^{-\frac{1}{2}(\alpha +n-k+\lambda )+1} (1-v_{n-k})} \nonumber\\
&&\times \Bigg( w_{n-k,n}^{-\frac{1}{2}(n-k-1+\lambda )}\left( w_{n-k,n} \partial _{ w_{n-k,n}}\right) w_{n-k,n}^{\frac{1}{2}(\beta -\gamma -\delta +1)} \left( w_{n-k,n} \partial _{ w_{n-k,n}}\right) w_{n-k,n}^{\frac{1}{2}(-\beta +\gamma +\delta +\lambda +n-k-2)} -\frac{q}{4}\Bigg) \Bigg\}\nonumber\\
&&\times\left.\sum_{i_0=0}^{\infty }\frac{(\frac{\alpha }{2}+\frac{\lambda }{2})_{i_0} }{(1+\frac{\lambda }{2})_{i_0}(\frac{1}{2}+ \frac{\gamma }{2} +\frac{\lambda }{2})_{i_0}} w_{1,n}^{i_0}\Bigg\} x^n \right\} \label{eq:40037}
\end{eqnarray}
In the above, the first sub-integral form contains one term of $A_n's$, the second one contains two terms of $A_n$'s, the third one contains three terms of $A_n$'s, etc.
\end{theorem}
\begin{proof}
In (\ref{eq:40020}) sub-power series $y_0(x) $, $y_1(x)$, $y_2(x)$ and $y_3(x)$ of the CHE for infinite series using 3TRF about $x=0$ are
\begin{subequations}
\begin{equation}
 y_0(x)= c_0 x^{\lambda } \sum_{i_0=0}^{\infty } \frac{\left( \frac{\alpha }{2} +\frac{\lambda }{2} \right)_{i_0}}{(1+\frac{\lambda}{2} )_{i_0}(\frac{1}{2}+\frac{\gamma }{2}+\frac{\lambda }{2} )_{i_0}} \eta ^{i_0} \label{er:40026a}
\end{equation}
\begin{eqnarray}
 y_1(x) &=& c_0 x^{\lambda } \left\{\sum_{i_0=0}^{\infty} \frac{(i_0+\frac{\lambda }{2})(i_0+\frac{1}{2}(-\beta +\gamma +\delta -1+\lambda ))-\frac{q}{4}}{(i_0+\frac{1}{2}+\frac{\lambda }{2})(i_0+\frac{\gamma }{2}+\frac{\lambda }{2})} \frac{\left( \frac{\alpha }{2} +\frac{\lambda }{2} \right)_{i_0}}{(1+\frac{\lambda}{2} )_{i_0}(\frac{1}{2}+\frac{\gamma }{2}+\frac{\lambda }{2} )_{i_0}}\right.\nonumber\\
&&\times  \left. \sum_{i_1=i_0}^{\infty}  \frac{\left( \frac{\alpha }{2}+\frac{1}{2}+\frac{\lambda }{2} \right)_{i_1} (\frac{3}{2}+\frac{\lambda}{2} )_{i_0}(1+\frac{\gamma }{2} +\frac{\lambda }{2})_{i_0}}{\left( \frac{\alpha }{2}+\frac{1}{2}+\frac{\lambda }{2} \right)_{i_0} (\frac{3}{2}+\frac{\lambda}{2})_{i_1}(1+\frac{\gamma }{2} +\frac{\lambda }{2})_{i_1}} \eta ^{i_1} \right\}x  \label{er:40026b}
\end{eqnarray}
\begin{eqnarray}
 y_2(x) &=& c_0 x^{\lambda } \left\{\sum_{i_0=0}^{\infty} \frac{(i_0+\frac{\lambda }{2})(i_0+\frac{1}{2}(-\beta +\gamma +\delta -1+\lambda ))-\frac{q}{4}}{(i_0+\frac{1}{2}+\frac{\lambda }{2})(i_0+\frac{\gamma }{2}+\frac{\lambda }{2})} \frac{\left( \frac{\alpha }{2} +\frac{\lambda }{2} \right)_{i_0}}{(1+\frac{\lambda}{2} )_{i_0}(\frac{1}{2}+\frac{\gamma }{2}+\frac{\lambda }{2} )_{i_0}} \right.\nonumber\\
&&\times  \sum_{i_1=i_0}^{\infty} \frac{(i_1+\frac{1}{2}+\frac{\lambda }{2})(i_1+\frac{1}{2}(-\beta +\gamma +\delta +\lambda ))-\frac{q}{4} }{ (i_1+1+\frac{\lambda }{2})(i_1+\frac{1}{2}+\frac{\gamma }{2}+\frac{\lambda }{2})} \frac{\left( \frac{\alpha }{2}+\frac{1}{2}+\frac{\lambda }{2} \right)_{i_1} (\frac{3}{2}+\frac{\lambda}{2} )_{i_0}(1+\frac{\gamma }{2} +\frac{\lambda }{2})_{i_0}}{\left( \frac{\alpha }{2}+\frac{1}{2}+\frac{\lambda }{2} \right)_{i_0} (\frac{3}{2}+\frac{\lambda}{2})_{i_1}(1+\frac{\gamma }{2} +\frac{\lambda }{2})_{i_1}}   \nonumber\\
&&\times \left.\sum_{i_2=i_1}^{\infty} \frac{\left( \frac{\alpha }{2}+1+\frac{\lambda }{2} \right)_{i_2} (2+\frac{\lambda}{2} )_{i_1}(\frac{3}{2}+\frac{\gamma }{2} +\frac{\lambda }{2})_{i_1}}{\left( \frac{\alpha }{2}+1+\frac{\lambda }{2} \right)_{i_1} (2+\frac{\lambda}{2})_{i_2}(\frac{3}{2}+\frac{\gamma }{2} +\frac{\lambda }{2})_{i_2}} \eta ^{i_2} \right\} x^2  \label{er:40026c}
\end{eqnarray}
\begin{eqnarray}
 y_3(x) &=&  c_0 x^{\lambda } \left\{\sum_{i_0=0}^{\infty} \frac{(i_0+\frac{\lambda }{2})(i_0+\frac{1}{2}(-\beta +\gamma +\delta -1+\lambda ))-\frac{q}{4}}{(i_0+\frac{1}{2}+\frac{\lambda }{2})(i_0+\frac{\gamma }{2}+\frac{\lambda }{2})} \frac{\left( \frac{\alpha }{2} +\frac{\lambda }{2} \right)_{i_0}}{(1+\frac{\lambda}{2} )_{i_0}(\frac{1}{2}+\frac{\gamma }{2}+\frac{\lambda }{2} )_{i_0}} \right.\nonumber\\
&&\times  \sum_{i_1=i_0}^{\infty} \frac{(i_1+\frac{1}{2}+\frac{\lambda }{2})(i_1+\frac{1}{2}(-\beta +\gamma +\delta +\lambda ))-\frac{q}{4} }{ (i_1+1+\frac{\lambda }{2})(i_1+\frac{1}{2}+\frac{\gamma }{2}+\frac{\lambda }{2})} \frac{\left( \frac{\alpha }{2}+\frac{1}{2}+\frac{\lambda }{2} \right)_{i_1} (\frac{3}{2}+\frac{\lambda}{2} )_{i_0}(1+\frac{\gamma }{2} +\frac{\lambda }{2})_{i_0}}{\left( \frac{\alpha }{2}+\frac{1}{2}+\frac{\lambda }{2} \right)_{i_0} (\frac{3}{2}+\frac{\lambda}{2})_{i_1}(1+\frac{\gamma }{2} +\frac{\lambda }{2})_{i_1}} \nonumber\\
&&\times \sum_{i_2=i_1}^{\infty} \frac{(i_2+1+\frac{\lambda }{2})(i_1+\frac{1}{2}(-\beta +\gamma +\delta +1+\lambda ))-\frac{q}{4} }{ (i_2+\frac{3}{2}+\frac{\lambda }{2})(i_2+1+\frac{\gamma }{2}+\frac{\lambda }{2})} \frac{\left( \frac{\alpha }{2}+1+\frac{\lambda }{2} \right)_{i_2} (2+\frac{\lambda}{2} )_{i_1}(\frac{3}{2}+\frac{\gamma }{2} +\frac{\lambda }{2})_{i_1}}{\left( \frac{\alpha }{2}+1+\frac{\lambda }{2} \right)_{i_1} (2+\frac{\lambda}{2})_{i_2}(\frac{3}{2}+\frac{\gamma }{2} +\frac{\lambda }{2})_{i_2}} \nonumber\\
&&\times \left.\sum_{i_3=i_2}^{\infty} \frac{\left( \frac{\alpha }{2}+\frac{3}{2}+\frac{\lambda }{2} \right)_{i_3} (\frac{5}{2}+\frac{\lambda}{2} )_{i_2}(2+\frac{\gamma }{2} +\frac{\lambda }{2})_{i_2}}{\left( \frac{\alpha }{2}+\frac{3}{2}+\frac{\lambda }{2} \right)_{i_2} (\frac{5}{2}+\frac{\lambda}{2})_{i_3}(2+\frac{\gamma }{2} +\frac{\lambda }{2})_{i_3}} \eta ^{i_3} \right\} x^3  \label{er:40026d} 
\end{eqnarray}
\end{subequations}
Put $l=1$ in (\ref{er:40025}). Take the new (\ref{er:40025}) into (\ref{er:40026b}).
\begin{eqnarray}
y_1(x) &=& \int_{0}^{1} dt_1\;t_1^{\frac{1}{2}(-1+\lambda )} \int_{0}^{1} du_1\;u_1^{\frac{1}{2}(-2+\gamma +\lambda )} \frac{1}{2\pi i} \oint dv_1 \; \frac{\exp\left(-\frac{v_1}{(1-v_1)}(1-t_1)(1-u_1)\eta \right)}{v_1^{-\frac{1}{2}(\alpha +1+\lambda )+1}(1-v_1)}\nonumber\\
&&\times \left( w_{1,1}^{-\frac{\lambda }{2}} \left( w_{1,1} \partial_{w_{1,1}}\right) w_{1,1}^{\frac{1}{2}(\beta -\gamma -\delta +1)} \left( w_{1,1} \partial_{w_{1,1}}\right) w_{1,1}^{\frac{1}{2}(-\beta +\gamma +\delta -1+\lambda )}-\frac{q}{4}\right) \nonumber\\
&&\times \Bigg\{ c_0 x^{\lambda }  \sum_{i_0=0}^{\infty } \frac{\left( \frac{\alpha }{2} +\frac{\lambda }{2} \right)_{i_0}}{(1+\frac{\lambda}{2} )_{i_0}(\frac{1}{2}+\frac{\gamma }{2}+\frac{\lambda }{2} )_{i_0}} w_{1,1}^{i_0} \Bigg\} x\label{er:40027}
\end{eqnarray}
where
\begin{equation}
w_{1,1} = \eta \prod_{l=1}^{1} t_l u_l v_l\nonumber
\end{equation}
Put $l=2$ in (\ref{er:40025}). Take the new (\ref{er:40025}) into (\ref{er:40026c}).
\begin{eqnarray}
y_2(x) &=& c_0 x^{\lambda } \int_{0}^{1} dt_2\;t_2^{\frac{\lambda }{2}} \int_{0}^{1} du_2\;u_2^{\frac{1}{2}(-1+\gamma +\lambda )} \frac{1}{2\pi i} \oint dv_2 \; \frac{\exp\left(-\frac{v_2}{(1-v_2)}(1-t_2)(1-u_2)\eta \right)}{v_2^{-\frac{1}{2}(\alpha +2+\lambda )+1}(1-v_2)} \nonumber\\
&&\times \left( w_{2,2}^{-\frac{1}{2}(1+\lambda )} \left( w_{2,2} \partial_{w_{2,2}}\right) w_{2,2}^{\frac{1}{2}(\beta -\gamma -\delta +1)} \left( w_{2,2} \partial_{w_{2,2}}\right) w_{2,2}^{\frac{1}{2}(-\beta +\gamma +\delta +\lambda )}-\frac{q}{4}\right) \nonumber\\
&&\times \Bigg\{\sum_{i_0=0}^{\infty } \frac{(i_0+\frac{\lambda }{2})(i_0+\frac{1}{2}(-\beta +\gamma +\delta -1+\lambda ))-\frac{q}{4}}{(i_0+\frac{1}{2}+\frac{\lambda }{2})(i_0+\frac{\gamma }{2}+\frac{\lambda }{2})} \frac{\left( \frac{\alpha }{2} +\frac{\lambda }{2} \right)_{i_0}}{(1+\frac{\lambda}{2} )_{i_0}(\frac{1}{2}+\frac{\gamma }{2}+\frac{\lambda }{2} )_{i_0}} \nonumber\\
&&\times  \sum_{i_1=i_0}^{\infty } \frac{\left( \frac{\alpha }{2} +\frac{1}{2}+\frac{\lambda }{2} \right)_{i_1} (\frac{3}{2}+\frac{\lambda}{2} )_{i_0}(1+\frac{\gamma }{2} +\frac{\lambda }{2})_{i_0}}{\left( \frac{\alpha }{2} +\frac{1}{2}+\frac{\lambda }{2} \right)_{i_0} (\frac{3}{2}+\frac{\lambda}{2})_{i_1}(1+\frac{\gamma }{2} +\frac{\lambda }{2})_{i_1}} w_{2,2}^{i_1} \Bigg\} x^2 \label{er:40028}
\end{eqnarray}
where
\begin{equation}
w_{2,2} = \eta \prod_{l=2}^{2} t_l u_l v_l\nonumber
\end{equation}
Put $l=1$ and $\eta = w_{2,2}$ in (\ref{er:40025}). Take the new (\ref{er:40025}) into (\ref{er:40028}).
\begin{eqnarray}
y_2(x) &=& \int_{0}^{1} dt_2\;t_2^{\frac{\lambda }{2}} \int_{0}^{1} du_2\;u_2^{\frac{1}{2}(-1+\gamma +\lambda )} \frac{1}{2\pi i} \oint dv_2 \; \frac{\exp\left(-\frac{v_2}{(1-v_2)}(1-t_2)(1-u_2)\eta \right)}{v_2^{-\frac{1}{2}(\alpha +2+\lambda )+1}(1-v_2)} \nonumber\\
&&\times \left( w_{2,2}^{-\frac{1}{2}(1+\lambda )} \left( w_{2,2} \partial_{w_{2,2}}\right) w_{2,2}^{\frac{1}{2}(\beta -\gamma -\delta +1)} \left( w_{2,2} \partial_{w_{2,2}}\right) w_{2,2}^{\frac{1}{2}(-\beta +\gamma +\delta +\lambda )}-\frac{q}{4}\right) \nonumber\\
&&\times \int_{0}^{1} dt_1\;t_1^{\frac{1}{2}(-1+\lambda )} \int_{0}^{1} du_1\;u_1^{\frac{1}{2}(-2+\gamma +\lambda )} \frac{1}{2\pi i} \oint dv_1 \; \frac{\exp\left(-\frac{v_1}{(1-v_1)}(1-t_1)(1-u_1)w_{2,2} \right)}{v_1^{-\frac{1}{2}(\alpha +1+\lambda )+1}(1-v_1)} \nonumber\\
&&\times \left( w_{1,2}^{-\frac{\lambda }{2}} \left( w_{1,2} \partial_{w_{1,2}}\right) w_{1,2}^{\frac{1}{2}(\beta -\gamma -\delta +1)} \left( w_{1,2} \partial_{w_{1,2}}\right) w_{1,2}^{\frac{1}{2}(-\beta +\gamma +\delta -1+\lambda )}-\frac{q}{4}\right) \nonumber\\
&&\times \left\{ c_0 x^{\lambda } \sum_{i_0=0}^{\infty } \frac{\left( \frac{\alpha }{2} +\frac{\lambda }{2} \right)_{i_0}}{(1+\frac{\lambda}{2} )_{i_0}(\frac{1}{2}+\frac{\gamma }{2}+\frac{\lambda }{2} )_{i_0}} w_{1,2}^{i_0} \right\} x^2 \label{er:40029}
\end{eqnarray}
where
\begin{equation}
w_{1,2} = \eta \prod_{l=1}^{2} t_l u_l v_l \nonumber
\end{equation}
By using similar process for the previous cases of integral forms of $y_1(x)$ and $y_2(x)$, the integral form of sub-power series expansion of $y_3(x)$ is
\begin{eqnarray}
y_3(x) &=& \int_{0}^{1} dt_3\;t_3^{\frac{1}{2}(1+\lambda )} \int_{0}^{1} du_3\;u_3^{\frac{1}{2}(\gamma +\lambda )} \frac{1}{2\pi i} \oint dv_3 \; \frac{\exp\left(-\frac{v_3}{(1-v_3)}(1-t_3)(1-u_3)\eta \right)}{v_3^{-\frac{1}{2}(\alpha +3+\lambda )+1}(1-v_3)} \nonumber\\
&&\times \left( w_{3,3}^{-\frac{1}{2}(2+\lambda )} \left( w_{3,3} \partial_{w_{3,3}}\right) w_{3,3}^{\frac{1}{2}(\beta -\gamma -\delta +1)} \left( w_{3,3} \partial_{w_{3,3}}\right) w_{3,3}^{\frac{1}{2}(-\beta +\gamma +\delta +1+\lambda )}-\frac{q}{4}\right) \nonumber\\ 
&&\times \int_{0}^{1} dt_2\;t_2^{\frac{\lambda }{2}} \int_{0}^{1} du_2\;u_2^{\frac{1}{2}(-1+\gamma +\lambda )} \frac{1}{2\pi i} \oint dv_2 \; \frac{\exp\left(-\frac{v_2}{(1-v_2)}(1-t_2)(1-u_2)w_{3,3} \right)}{v_2^{-\frac{1}{2}(\alpha +2+\lambda )+1}(1-v_2)} \nonumber\\
&&\times \left( w_{2,3}^{-\frac{1}{2}(1+\lambda )} \left( w_{2,3} \partial_{w_{2,3}}\right) w_{2,3}^{\frac{1}{2}(\beta -\gamma -\delta +1)} \left( w_{2,3} \partial_{w_{2,3}}\right) w_{2,3}^{\frac{1}{2}(-\beta +\gamma +\delta +\lambda )}-\frac{q}{4}\right) \nonumber\\
&&\times \int_{0}^{1} dt_1\;t_1^{\frac{1}{2}(-1+\lambda )} \int_{0}^{1} du_1\;u_1^{\frac{1}{2}(-2+\gamma +\lambda )} \frac{1}{2\pi i} \oint dv_1 \; \frac{\exp\left(-\frac{v_1}{(1-v_1)}(1-t_1)(1-u_1)w_{2,3} \right)}{v_1^{-\frac{1}{2}(\alpha +1+\lambda )+1}(1-v_1)} \nonumber\\
&&\times \left( w_{1,3}^{-\frac{\lambda }{2}} \left( w_{1,3} \partial_{w_{1,3}}\right) w_{1,3}^{\frac{1}{2}(\beta -\gamma -\delta +1)} \left( w_{1,3} \partial_{w_{1,3}}\right) w_{1,3}^{\frac{1}{2}(-\beta +\gamma +\delta -1+\lambda )}-\frac{q}{4}\right) \nonumber\\
&&\times \left\{ c_0 x^{\lambda } \sum_{i_0=0}^{\infty } \frac{\left( \frac{\alpha }{2} +\frac{\lambda }{2} \right)_{i_0}}{(1+\frac{\lambda}{2} )_{i_0}(\frac{1}{2}+\frac{\gamma }{2}+\frac{\lambda }{2} )_{i_0}} w_{1,3}^{i_0} \right\} x^3 \label{er:40030}
\end{eqnarray}
where
\begin{equation}
 w_{3,3} = \eta \prod_{l=3}^{3} t_l u_l v_l \hspace{1cm}
w_{2,3} = \eta \prod_{l=2}^{3} t_l u_l v_l \hspace{1cm}
w_{1,3} = \eta \prod_{l=1}^{3} t_l u_l v_l \nonumber
\end{equation}
By repeating this process for all higher terms of integral forms of sub-summation $y_m(x)$ terms where $m \geq 4$, we obtain every integral forms of $y_m(x)$ terms. 
Since we substitute (\ref{er:40026a}), (\ref{er:40027}), (\ref{er:40029}), (\ref{er:40030}) and including all integral forms of $y_m(x)$ terms where $m \geq 4$ into (\ref{eq:40020}), we obtain (\ref{eq:40037}).\footnote{Or replace the finite summation with an interval $[0, \alpha _0]$ by infinite summation with an interval  $[0,\infty ]$ in (\ref{eq:40028}). Replace $\alpha _0$ and $\alpha _{n-k}$ by $-\frac{1}{2}(\alpha +\lambda ) $ and $-\frac{1}{2}(\alpha +n-k+\lambda ) $ into the new (\ref{eq:40028}). Its solution is also equivalent to (\ref{eq:40037})} 
\qed
\end{proof}
Put $c_0$= 1 as $\lambda =0$  for the first kind of independent solutions of the CHE and $\lambda = 1-\gamma $ in (\ref{eq:40037}).  
\begin{remark}
The integral representation of the CHE of the first kind for infinite series about $x=0$ using 3TRF is
\begin{eqnarray}
 y(x)&=& H_c^{(a)}F\left(\alpha, \beta, \gamma, \delta, q; \eta =\frac{1}{2}\beta x^2  \right) \nonumber\\
&=& _1F_1 \left(\frac{\alpha }{2}; \frac{1}{2}+ \frac{\gamma }{2}; \eta \right) 
+ \sum_{n=1}^{\infty } \Bigg\{\prod _{k=0}^{n-1} \Bigg\{ \int_{0}^{1} dt_{n-k}\;t_{n-k}^{\frac{1}{2}(n-k-2)} \int_{0}^{1} du_{n-k}\;u_{n-k}^{\frac{1}{2}(n-k-3+\gamma )} \nonumber\\
&&\times  \frac{1}{2\pi i}  \oint dv_{n-k} \frac{\exp\left( -\frac{v_{n-k}}{(1-v_{n-k})}w_{n-k+1,n}(1-t_{n-k})(1-u_{n-k})\right) }{v_{n-k}^{-\frac{1}{2}(\alpha +n-k-2)} (1-v_{n-k})} \nonumber\\
&&\times \Bigg( w_{n-k,n}^{-\frac{1}{2}(n-k-1)}\left( w_{n-k,n} \partial _{ w_{n-k,n}}\right) w_{n-k,n}^{\frac{1}{2}(\beta -\gamma -\delta +1)} \left( w_{n-k,n} \partial _{ w_{n-k,n}}\right) w_{n-k,n}^{\frac{1}{2}(-\beta +\gamma +\delta +n-k-2)} -\frac{q}{4}\Bigg) \Bigg\}\nonumber\\
&&\times _1F_1 \left(\frac{\alpha }{2}; \frac{1}{2}+ \frac{\gamma }{2};  w_{1,n} \right) \Bigg\} x^n   \label{eq:40038}
\end{eqnarray}
\end{remark}
\begin{remark}
The integral representation of the CHE of the second kind for infinite series about $x=0$ using 3TRF is
\begin{eqnarray}
y(x)&=& H_c^{(a)}S\left(\alpha, \beta, \gamma, \delta, q; \eta =\frac{1}{2}\beta x^2  \right) \nonumber\\
&=& x^{1-\gamma } \Bigg\{ \;_1F_1 \left(\frac{\alpha }{2}+\frac{1}{2}-\frac{\gamma }{2}; \frac{3}{2}-\frac{\gamma }{2}; \eta \right) \nonumber\\ 
&&+ \sum_{n=1}^{\infty } \Bigg\{\prod _{k=0}^{n-1} \Bigg\{ \int_{0}^{1} dt_{n-k}\;t_{n-k}^{\frac{1}{2}(n-k-1-\gamma )} \int_{0}^{1} du_{n-k}\;u_{n-k}^{\frac{1}{2}(n-k-2)}  \nonumber\\
&&\times  \frac{1}{2\pi i}  \oint dv_{n-k} \frac{\exp\left( -\frac{v_{n-k}}{(1-v_{n-k})}w_{n-k+1,n}(1-t_{n-k})(1-u_{n-k})\right) }{v_{n-k}^{-\frac{1}{2}(\alpha +n-k-1-\gamma )} (1-v_{n-k})} \nonumber\\
&&\times \Bigg( w_{n-k,n}^{-\frac{1}{2}(n-k-\gamma )}\left( w_{n-k,n} \partial _{ w_{n-k,n}}\right) w_{n-k,n}^{\frac{1}{2}(\beta -\gamma -\delta +1)} \left( w_{n-k,n} \partial _{ w_{n-k,n}}\right) w_{n-k,n}^{\frac{1}{2}(-\beta +\delta +n-k-1)} -\frac{q}{4}\Bigg) \Bigg\}\nonumber\\
&&\times _1F_1 \left(\frac{\alpha }{2}+\frac{1}{2}-\frac{\gamma }{2}; \frac{3}{2}-\frac{\gamma }{2}; w_{1,n} \right) \Bigg\} x^n \Bigg\} \label{eq:40039}
\end{eqnarray}
\end{remark}
\subsection[Generating function for the CHP of type 1]{Generating function for the CHP of type 1}
I consider the generating function for the CHP of type 1. Since the generating function for the CHP is derived, we might be possible to construct orthogonal relations of the CHP.
In the physical point of view we might be possible to obtain the normalized constant for the wave function in modern physics     
(especially black hole problems) recursion relation and its expectation value of any physical quantities from the generating function for the CHP.\footnote{In this chapter I consider the generating functions for type 1 CHP. In the next chapter and future paper, I will construct the generating functions for type 2 and 3 CHP's. For the type 1 polynomial, I treat $\beta $, $\gamma $, $\delta $ and $q$ as free variables and $\alpha $ as a fixed value. For the type 2 polynomial, I treat $\alpha $, $\beta $, $\gamma $ and $\delta $  as free variables and  $q$ as a fixed value. For the type 3 polynomial, I treat $\beta $, $\gamma $ and $\delta $  as free variables and  $\alpha $ and $q$ as fixed values.} For the case of hydrogen-like atoms, the normalized wave function is derived from the generating function for associated Laguerre polynomial. And the expectation value of physical quantities such as position and momentum is constructed by applying the recursive relation of associated Laguerre polynomial.

Let's investigate the generating function for the first and second kinds of the CHP about $x=0$ as $B_n's$ term terminated at certain eigenvalue. 
\begin{definition}
I define that
\begin{equation}
\begin{cases}
\displaystyle { s_{a,b}} = \begin{cases} \displaystyle {  s_a\cdot s_{a+1}\cdot s_{a+2}\cdots s_{b-2}\cdot s_{b-1}\cdot s_b}\;\;\mbox{where}\;a>b \cr
s_a \;\;\mbox{only}\;\mbox{if}\;a=b\end{cases}
\cr
\cr
\displaystyle {  \widetilde{w}_{a,b}= \eta s_{a,\infty }\prod_{l=a}^{b}t_l u_l}
\end{cases}\label{eq:40040}
\end{equation}
where
\begin{equation}
a,b\in \mathbb{N}_{0} \nonumber
\end{equation}
\end{definition}
And we have
\begin{equation}
\sum_{\alpha _i = \alpha _j}^{\infty } s_i^{\alpha _i} = \frac{s_i^{\alpha _j}}{(1-s_i)}\label{eq:40041}
\end{equation}
Acting the summation operator $\displaystyle{ \sum_{\alpha _0 =0}^{\infty } \frac{s_0^{\alpha _0}}{\alpha _0!} \frac{\Gamma (\alpha _0+\gamma')}{\Gamma(\gamma') }  \prod _{n=1}^{\infty } \left\{ \sum_{ \alpha _n = \alpha _{n-1}}^{\infty } s_n^{\alpha _n }\right\}}$ on (\ref{eq:40028}) where $|s_i|<1$ as $i=0,1,2,\cdots$ by using (\ref{eq:40040}) and (\ref{eq:40041}),
\begin{theorem} 
The general expression of the generating function for the CHP of type 1 about $x=0$ is given by
\begin{eqnarray}
&&\sum_{\alpha _0 =0}^{\infty } \frac{s_0^{\alpha _0}}{\alpha _0!} \frac{\Gamma (\alpha _0+\gamma')}{\Gamma(\gamma') }  \prod _{n=1}^{\infty } \left\{ \sum_{ \alpha _n = \alpha _{n-1}}^{\infty } s_n^{\alpha _n }\right\} y(x) \nonumber\\
&&= \prod_{l=1}^{\infty } \frac{1}{(1-s_{l,\infty })} \mathbf{\Upsilon}(\gamma', \lambda; s_{0,\infty } ;\eta) \nonumber\\
&&+ \Bigg\{ \prod_{l=1}^{\infty } \frac{1}{(1-s_{l,\infty })} \int_{0}^{1} dt_1\;t_1^{\frac{1}{2}(-1+\lambda )} \int_{0}^{1} du_1\;u_1^{\frac{1}{2}(-2+\gamma +\lambda )} \exp\left( -\frac{s_{1,\infty }}{(1-s_{1,\infty })}\eta (1-t_1)(1-u_1)\right) \nonumber\\
&&\times \Bigg(  \widetilde{w}_{1,1}^{-\frac{\lambda }{2}}\left(  \widetilde{w}_{1,1} \partial _{ \widetilde{w}_{1,1}}\right) \widetilde{w}_{1,1}^{\frac{1}{2}(\beta -\gamma -\delta +1)} \left(  \widetilde{w}_{1,1} \partial _{ \widetilde{w}_{1,1}}\right) \widetilde{w}_{1,1}^{\frac{1}{2}(-\beta +\gamma +\delta -1+\lambda )}-\frac{q}{4}\Bigg)\; \mathbf{\Upsilon}(\gamma', \lambda ; s_0;\widetilde{w}_{1,1})\Bigg\} x\nonumber\\
&&+ \sum_{n=2}^{\infty } \Bigg\{ \prod_{l=n}^{\infty } \frac{1}{(1-s_{l,\infty })} \int_{0}^{1} dt_n\;t_n^{\frac{1}{2}(n-2+\lambda )} \int_{0}^{1} du_n\;u_n^{\frac{1}{2}(n-3+\gamma +\lambda )} \exp\left( -\frac{s_{n,\infty }}{(1-s_{n,\infty })}\eta (1-t_n)(1-u_n)\right) \nonumber\\
&&\times \Bigg(  \widetilde{w}_{n,n}^{-\frac{1}{2}(n-1+\lambda )}\left(  \widetilde{w}_{n,n} \partial _{ \widetilde{w}_{n,n}}\right) \widetilde{w}_{n,n}^{\frac{1}{2}(\beta -\gamma -\delta +1)} \left(  \widetilde{w}_{n,n} \partial _{ \widetilde{w}_{n,n}}\right) \widetilde{w}_{n,n}^{\frac{1}{2}(-\beta +\gamma +\delta +n-2+\lambda )}-\frac{q}{4}\Bigg) \nonumber\\
&&\times \prod_{k=1}^{n-1} \Bigg\{ \int_{0}^{1} dt_{n-k}\;t_{n-k}^{\frac{1}{2}(n-k-2+\lambda )} \int_{0}^{1} du_{n-k} \;u_{n-k}^{\frac{1}{2}(n-k-3+\gamma +\lambda ) } \nonumber\\
&&\times \frac{ \exp\left( -\frac{s_{n-k}}{(1-s_{n-k})}\widetilde{w}_{n-k+1,n} (1-t_{n-k})(1-u_{n-k})\right)}{(1-s_{n-k})}\nonumber\\
&&\times \left. \Bigg( \widetilde{w}_{n-k,n}^{-\frac{1}{2}(n-k-1+\lambda )}\left(  \widetilde{w}_{n-k,n} \partial _{ \widetilde{w}_{n-k,n}}\right) \widetilde{w}_{n-k,n}^{\frac{1}{2}(\beta -\gamma -\delta +1)} \left(  \widetilde{w}_{n-k,n} \partial _{ \widetilde{w}_{n-k,n}}\right) \widetilde{w}_{n-k,n}^{\frac{1}{2}(-\beta +\gamma +\delta +n-k-2+\lambda )}-\frac{q}{4}\Bigg) \Bigg\} \right.\nonumber\\
&&\times \mathbf{\Upsilon}(\gamma', \lambda ; s_0;\widetilde{w}_{1,n}) \Bigg\} x^n \label{eq:40042}
\end{eqnarray}
\normalsize
where
\begin{equation}
\begin{cases} 
{ \displaystyle \mathbf{\Upsilon}(\gamma',\lambda; s_{0,\infty } ;\eta)=  \sum_{\alpha _0 =0}^{\infty } \frac{s_{0,\infty }^{\alpha _0}}{\alpha _0!} \frac{\Gamma (\alpha _0+\gamma')}{\Gamma(\gamma') }\left( c_0 x^{\lambda } \sum_{i_0=0}^{\alpha _0} \frac{(-\alpha _0)_{i_0} }{(1+\frac{\lambda }{2} )_{i_0}(\frac{1}{2}+\frac{\gamma }{2}+\frac{\lambda }{2})_{i_0}} \eta ^{i_0} \right) }\cr
{ \displaystyle \mathbf{\Upsilon}(\gamma',\lambda ; s_0;\widetilde{w}_{1,1}) = \sum_{\alpha _0 =0}^{\infty } \frac{s_0^{\alpha _0}}{\alpha _0!} \frac{\Gamma (\alpha _0+\gamma')}{\Gamma(\gamma') }\left( c_0 x^{\lambda } \sum_{i_0=0}^{\alpha _0} \frac{(-\alpha _0)_{i_0} }{(1+\frac{\lambda }{2} )_{i_0}(\frac{1}{2}+\frac{\gamma }{2}+\frac{\lambda }{2})_{i_0}} \widetilde{w}_{1,1} ^{i_0} \right) }\cr
{ \displaystyle \mathbf{\Upsilon}(\gamma',\lambda; s_0 ;\widetilde{w}_{1,n}) = \sum_{\alpha _0 =0}^{\infty } \frac{s_0^{\alpha _0}}{\alpha _0!} \frac{\Gamma (\alpha _0+\gamma')}{\Gamma(\gamma') }\left( c_0 x^{\lambda } \sum_{i_0=0}^{\alpha _0} \frac{(-\alpha _0)_{i_0} }{(1+\frac{\lambda }{2} )_{i_0}(\frac{1}{2}+\frac{\gamma }{2}+\frac{\lambda }{2})_{i_0}} \widetilde{w}_{1,n} ^{i_0} \right)}
\end{cases}\nonumber 
\end{equation}
\end{theorem}
\begin{proof} 
Acting the summation operator $\displaystyle{ \sum_{\alpha _0 =0}^{\infty } \frac{s_0^{\alpha _0}}{\alpha _0!} \frac{\Gamma (\alpha _0+\gamma')}{\Gamma(\gamma') }  \prod _{n=1}^{\infty } \left\{ \sum_{ \alpha _n = \alpha _{n-1}}^{\infty } s_n^{\alpha _n }\right\}}$ on the form of integral of the CHP of type 1 $y(x)$,
\begin{eqnarray}
&&\sum_{\alpha _0 =0}^{\infty } \frac{s_0^{\alpha _0}}{\alpha _0!} \frac{\Gamma (\alpha _0+\gamma')}{\Gamma(\gamma') }  \prod _{n=1}^{\infty } \left\{ \sum_{ \alpha _n = \alpha _{n-1}}^{\infty } s_n^{\alpha _n }\right\} y(x) \label{eq:40043}\\
&&= \sum_{\alpha _0 =0}^{\infty } \frac{s_0^{\alpha _0}}{\alpha _0!} \frac{\Gamma (\alpha _0+\gamma')}{\Gamma(\gamma') }  \prod _{n=1}^{\infty } \left\{ \sum_{ \alpha _n = \alpha _{n-1}}^{\infty } s_n^{\alpha _n }\right\} \Big\{ y_0(x)+y_1(x)+y_2(x)+y_3(x)+\cdots\Big\} \nonumber
\end{eqnarray}
Acting the summation operator $\displaystyle{ \sum_{\alpha _0 =0}^{\infty } \frac{s_0^{\alpha _0}}{\alpha _0!} \frac{\Gamma (\alpha _0+\gamma')}{\Gamma(\gamma') }  \prod _{n=1}^{\infty } \left\{ \sum_{ \alpha _n = \alpha _{n-1}}^{\infty } s_n^{\alpha _n }\right\}}$ on (\ref{eq:40030a}) by applying (\ref{eq:40041}),
\begin{eqnarray}
&&\sum_{\alpha _0 =0}^{\infty } \frac{s_0^{\alpha _0}}{\alpha _0!} \frac{\Gamma (\alpha _0+\gamma')}{\Gamma(\gamma') }  \prod _{n=1}^{\infty } \left\{ \sum_{ \alpha _n = \alpha _{n-1}}^{\infty } s_n^{\alpha _n }\right\} y_0(x) \nonumber\\
&&= \prod_{l=1}^{\infty } \frac{1}{(1-s_{l,\infty })} \sum_{\alpha _0 =0}^{\infty } \frac{s_{0,\infty }^{\alpha _0}}{\alpha _0!}  \left( c_0 x^{\lambda } \frac{\Gamma (\alpha _0+\gamma')}{\Gamma(\gamma') } \sum_{i_0=0}^{\alpha _0}  \frac{(-\alpha _0)_{i_0} }{(1+\frac{\lambda}{2} )_{i_0}(\frac{1}{2}+\frac{\gamma }{2}+\frac{\lambda }{2} )_{i_0}} \eta ^{i_0} \right) \hspace{2cm}\label{eq:40044}
\end{eqnarray}
Acting the summation operator $\displaystyle{ \sum_{\alpha _0 =0}^{\infty } \frac{s_0^{\alpha _0}}{\alpha _0!} \frac{\Gamma (\alpha _0+\gamma')}{\Gamma(\gamma') }  \prod _{n=1}^{\infty } \left\{ \sum_{ \alpha _n = \alpha _{n-1}}^{\infty } s_n^{\alpha _n }\right\}}$ on (\ref{eq:40031}),
\begin{eqnarray}
&&\sum_{\alpha _0 =0}^{\infty } \frac{s_0^{\alpha _0}}{\alpha _0!} \frac{\Gamma (\alpha _0+\gamma')}{\Gamma(\gamma') }  \prod _{n=1}^{\infty } \left\{ \sum_{ \alpha _n = \alpha _{n-1}}^{\infty } s_n^{\alpha _n }\right\} y_1(x) \nonumber\\
&&= \prod_{l=2}^{\infty } \frac{1}{(1-s_{l,\infty })} \int_{0}^{1} dt_1\;t_1^{\frac{1}{2}(-1+\lambda )} \int_{0}^{1} du_1\;u_1^{\frac{1}{2}(-2+\gamma +\lambda )}
 \frac{1}{2\pi i} \oint dv_1 \;\frac{\exp\left( -\frac{v_1}{(1-v_1)}\eta (1-t_1)(1-u_1)\right)}{v_1(1-v_1)} \nonumber\\
&&\times \sum_{\alpha _1 =\alpha _0}^{\infty }\left( \frac{s_{1,\infty }}{v_1}\right)^{\alpha _1} \Bigg( w_{1,1}^{-\frac{\lambda }{2}}\left(  w_{1,1} \partial _{ w_{1,1}}\right) w_{1,1}^{\frac{1}{2}(\beta -\gamma -\delta +1)} \left(  w_{1,1} \partial _{ w_{1,1}}\right) w_{1,1}^{\frac{1}{2}(-\beta +\gamma +\delta -1+\lambda )}-\frac{q}{4} \Bigg) \nonumber\\
&&\times  \sum_{\alpha _0 =0}^{\infty } \frac{s_0^{\alpha _0} }{\alpha _0!} \frac{\Gamma (\alpha _0+\gamma')}{\Gamma(\gamma') }  \left( c_0 x^{\lambda } \sum_{i_0=0}^{\alpha _0} \frac{(-\alpha _0)_{i_0} }{(1+\frac{\lambda}{2} )_{i_0}(\frac{1}{2}+\frac{\gamma }{2}+\frac{\lambda }{2} )_{i_0}} w_{1,1} ^{i_0} \right) x \label{eq:40045}
\end{eqnarray}
Replace $\alpha _i$, $\alpha _j$ and $s_i$ by $\alpha _1$, $\alpha _0$ and ${ \displaystyle \frac{s_{1,\infty }}{v_1}}$ in (\ref{eq:40041}). Take the new (\ref{eq:40041}) into (\ref{eq:40045}).
\begin{eqnarray}
&&\sum_{\alpha _0 =0}^{\infty } \frac{s_0^{\alpha _0}}{\alpha _0!} \frac{\Gamma (\alpha _0+\gamma')}{\Gamma(\gamma') }  \prod _{n=1}^{\infty } \left\{ \sum_{ \alpha _n = \alpha _{n-1}}^{\infty } s_n^{\alpha _n }\right\} y_1(x) \nonumber\\
&&= \prod_{l=2}^{\infty } \frac{1}{(1-s_{l,\infty })} \int_{0}^{1} dt_1\;t_1^{\frac{1}{2}(-1+\lambda )} \int_{0}^{1} du_1\;u_1^{\frac{1}{2}(-2+\gamma +\lambda )}
 \frac{1}{2\pi i} \oint dv_1 \;\frac{\exp\left( -\frac{v_1}{(1-v_1)}\eta (1-t_1)(1-u_1)\right)}{(1-v_1)(v_1-s_{1,\infty })} \nonumber\\
&&\times \Bigg( w_{1,1}^{-\frac{\lambda }{2}}\left(  w_{1,1} \partial _{ w_{1,1}}\right) w_{1,1}^{\frac{1}{2}(\beta -\gamma -\delta +1)} \left(  w_{1,1} \partial _{ w_{1,1}}\right) w_{1,1}^{\frac{1}{2}(-\beta +\gamma +\delta -1+\lambda )}-\frac{q}{4} \Bigg) \nonumber\\
&&\times  \sum_{\alpha _0 =0}^{\infty } \frac{1}{\alpha _0!}\left( \frac{s_{0,\infty }}{v_1}\right)^{\alpha _0} \frac{\Gamma (\alpha _0+\gamma')}{\Gamma(\gamma') }  \left( c_0 x^{\lambda } \sum_{i_0=0}^{\alpha _0} \frac{(-\alpha _0)_{i_0} }{(1+\frac{\lambda}{2} )_{i_0}(\frac{1}{2}+\frac{\gamma }{2}+\frac{\lambda }{2} )_{i_0}} w_{1,1} ^{i_0} \right) x \label{eq:40046}
\end{eqnarray}
By using Cauchy's integral formula, the contour integrand has poles at
 $ v_1= 1 $  or $ s_{1,\infty }$ and $ s_{1,\infty }$ is only inside the unit circle. As we compute the residue there in (\ref{eq:40046}) we obtain
\begin{eqnarray}
&&\sum_{\alpha _0 =0}^{\infty } \frac{s_0^{\alpha _0}}{\alpha _0!} \frac{\Gamma (\alpha _0+\gamma')}{\Gamma(\gamma') }  \prod _{n=1}^{\infty } \left\{ \sum_{ \alpha _n = \alpha _{n-1}}^{\infty } s_n^{\alpha _n }\right\} y_1(x) \nonumber\\
&&= \prod_{l=1}^{\infty } \frac{1}{(1-s_{l,\infty })} \int_{0}^{1} dt_1\;t_1^{\frac{1}{2}(-1+\lambda )} \int_{0}^{1} du_1\;u_1^{\frac{1}{2}(-2+\gamma +\lambda )}
 \; \exp\left( -\frac{s_{1,\infty }}{(1-s_{1,\infty })}\eta (1-t_1)(1-u_1)\right)  \nonumber\\
&&\times \Bigg( \widetilde{w}_{1,1} ^{-\frac{\lambda }{2}}\left(  \widetilde{w}_{1,1}  \partial _{\widetilde{w}_{1,1} }\right) \widetilde{w}_{1,1} ^{\frac{1}{2}(\beta -\gamma -\delta +1)} \left(\widetilde{w}_{1,1}  \partial _{ \widetilde{w}_{1,1} }\right) \widetilde{w}_{1,1} ^{\frac{1}{2}(-\beta +\gamma +\delta -1+\lambda )}-\frac{q}{4} \Bigg) \nonumber\\
&&\times  \sum_{\alpha _0 =0}^{\infty } \frac{s_0^{\alpha _0}}{\alpha _0!} \frac{\Gamma (\alpha _0+\gamma')}{\Gamma(\gamma') }  \left( c_0 x^{\lambda } \sum_{i_0=0}^{\alpha _0} \frac{(-\alpha _0)_{i_0} }{(1+\frac{\lambda}{2} )_{i_0}(\frac{1}{2}+\frac{\gamma }{2}+\frac{\lambda }{2} )_{i_0}} \widetilde{w}_{1,1}^{i_0} \right) x \label{eq:40047}\\
&&\mbox{where}\;\widetilde{w}_{1,1}= \eta s_{1,\infty }\prod_{l=1}^{1}t_l u_l \nonumber
\end{eqnarray} 
Acting the summation operator $\displaystyle{  \sum_{\alpha _0 =0}^{\infty } \frac{s_0^{\alpha _0}}{\alpha _0!} \frac{\Gamma (\alpha _0+\gamma')}{\Gamma(\gamma') }  \prod _{n=1}^{\infty } \left\{ \sum_{ \alpha _n = \alpha _{n-1}}^{\infty } s_n^{\alpha _n }\right\}}$ on (\ref{eq:40033}),
\begin{eqnarray}
&&  \sum_{\alpha _0 =0}^{\infty } \frac{s_0^{\alpha _0}}{\alpha _0!} \frac{\Gamma (\alpha _0+\gamma')}{\Gamma(\gamma') }  \prod _{n=1}^{\infty } \left\{ \sum_{ \alpha _n = \alpha _{n-1}}^{\infty } s_n^{\alpha _n }\right\} y_2(x) \nonumber\\
&&= \prod_{l=3}^{\infty } \frac{1}{(1-s_{l,\infty })} \int_{0}^{1} dt_2\;t_2^{\frac{\lambda}{2}} \int_{0}^{1} du_2\;u_2^{\frac{1}{2}(-1+\gamma +\lambda )}
 \frac{1}{2\pi i} \oint dv_2 \;\frac{\exp\left( -\frac{v_2}{(1-v_2)}\eta (1-t_2)(1-u_2)\right)}{v_2(1-v_2)}  \nonumber\\
&&\times \sum_{\alpha _2 =\alpha _1}^{\infty }\left(\frac{s_{2,\infty }}{v_2} \right)^{\alpha _2} \Bigg( w_{2,2}^{-\frac{1}{2}(1+\lambda )}\left(  w_{2,2} \partial _{ w_{2,2}}\right) w_{2,2}^{\frac{1}{2}(\beta -\gamma -\delta +1)} \left(  w_{2,2} \partial _{ w_{2,2}}\right) w_{2,2}^{\frac{1}{2}(-\beta +\gamma +\delta +\lambda )}-\frac{q}{4}\Bigg) \nonumber\\
&&\times \int_{0}^{1} dt_1\;t_1^{\frac{1}{2}(-1+\lambda )} \int_{0}^{1} du_1\;u_1^{\frac{1}{2}(-2+\gamma +\lambda )}
 \frac{1}{2\pi i} \oint dv_1 \;\frac{\exp\left( -\frac{v_1}{(1-v_1)} w_{2,2} (1-t_1)(1-u_1)\right)}{v_1(1-v_1)} \nonumber\\
&&\times \sum_{\alpha _1 =\alpha _0}^{\infty }\left(\frac{s_1}{v_1}\right)^{\alpha _1} \Bigg(  w_{1,2}^{-\frac{\lambda }{2}}\left(  w_{1,2} \partial _{ w_{1,2}}\right) w_{1,2}^{\frac{1}{2}(\beta -\gamma -\delta +1)} \left(  w_{1,2} \partial _{ w_{1,2}}\right) w_{1,2}^{\frac{1}{2}(-\beta +\gamma +\delta -1+\lambda )}-\frac{q}{4}\Bigg) \nonumber\\
&&\times  \sum_{\alpha _0 =0}^{\infty } \frac{s_0^{\alpha _0}}{\alpha _0!} \frac{\Gamma (\alpha _0+\gamma')}{\Gamma(\gamma') } \left( c_0 x^{\lambda } \sum_{i_0=0}^{\alpha _0} \frac{(-\alpha _0)_{i_0} }{(1+\frac{\lambda }{2} )_{i_0}(\frac{1}{2}+\frac{\gamma }{2}+\frac{\lambda }{2} )_{i_0}}  w_{1,2} ^{i_0} \right) x^2 \label{eq:40048}
\end{eqnarray}
Replace $\alpha _i$, $\alpha _j$ and $s_i$ by $\alpha _2$, $\alpha _1$ and ${ \displaystyle \frac{s_{2,\infty }}{v_2}}$ in (\ref{eq:40041}). Take the new (\ref{eq:40041}) into (\ref{eq:40048}).
\begin{eqnarray}
&&  \sum_{\alpha _0 =0}^{\infty } \frac{s_0^{\alpha _0}}{\alpha _0!} \frac{\Gamma (\alpha _0+\gamma')}{\Gamma(\gamma') }  \prod _{n=1}^{\infty } \left\{ \sum_{ \alpha _n = \alpha _{n-1}}^{\infty } s_n^{\alpha _n }\right\} y_2(x) \nonumber\\
&&= \prod_{l=3}^{\infty } \frac{1}{(1-s_{l,\infty })} \int_{0}^{1} dt_2\;t_2^{\frac{\lambda}{2}} \int_{0}^{1} du_2\;u_2^{\frac{1}{2}(-1+\gamma +\lambda )}
 \frac{1}{2\pi i} \oint dv_2 \;\frac{\exp\left( -\frac{v_2}{(1-v_2)}\eta (1-t_2)(1-u_2)\right)}{ (1-v_2)(v_2-s_{2,\infty })}  \nonumber\\
&&\times \Bigg( w_{2,2}^{-\frac{1}{2}(1+\lambda )}\left(  w_{2,2} \partial _{ w_{2,2}}\right) w_{2,2}^{\frac{1}{2}(\beta -\gamma -\delta +1)} \left(  w_{2,2} \partial _{ w_{2,2}}\right) w_{2,2}^{\frac{1}{2}(-\beta +\gamma +\delta +\lambda )}-\frac{q}{4}\Bigg) \nonumber\\
&&\times \int_{0}^{1} dt_1\;t_1^{\frac{1}{2}(-1+\lambda )} \int_{0}^{1} du_1\;u_1^{\frac{1}{2}(-2+\gamma +\lambda )}
 \frac{1}{2\pi i} \oint dv_1 \;\frac{\exp\left( -\frac{v_1}{(1-v_1)} w_{2,2} (1-t_1)(1-u_1)\right)}{v_1(1-v_1)} \nonumber\\
&&\times \sum_{\alpha _1 =\alpha _0}^{\infty }\left(\frac{s_{1,\infty }}{v_1 v_2}\right)^{\alpha _1} \Bigg(  w_{1,2}^{-\frac{\lambda }{2}}\left(  w_{1,2} \partial _{ w_{1,2}}\right) w_{1,2}^{\frac{1}{2}(\beta -\gamma -\delta +1)} \left(  w_{1,2} \partial _{ w_{1,2}}\right) w_{1,2}^{\frac{1}{2}(-\beta +\gamma +\delta -1+\lambda )}-\frac{q}{4}\Bigg) \nonumber\\
&&\times  \sum_{\alpha _0 =0}^{\infty } \frac{s_0^{\alpha _0}}{\alpha _0!} \frac{\Gamma (\alpha _0+\gamma')}{\Gamma(\gamma') } \left( c_0 x^{\lambda } \sum_{i_0=0}^{\alpha _0} \frac{(-\alpha _0)_{i_0} }{(1+\frac{\lambda }{2} )_{i_0}(\frac{1}{2}+\frac{\gamma }{2}+\frac{\lambda }{2} )_{i_0}}  w_{1,2} ^{i_0} \right) x^2 \label{eq:40049}
\end{eqnarray}
By using Cauchy's integral formula, the contour integrand has poles at
 $ v_2= 1 $  or $ s_{2,\infty }$ and $ s_{2,\infty }$ is only inside the unit circle. As we compute the residue there in (\ref{eq:40049}) we obtain
\begin{eqnarray}
&&  \sum_{\alpha _0 =0}^{\infty } \frac{s_0^{\alpha _0}}{\alpha _0!} \frac{\Gamma (\alpha _0+\gamma')}{\Gamma(\gamma') }  \prod _{n=1}^{\infty } \left\{ \sum_{ \alpha _n = \alpha _{n-1}}^{\infty } s_n^{\alpha _n }\right\} y_2(x) \nonumber\\
&&= \prod_{l=2}^{\infty } \frac{1}{(1-s_{l,\infty })} \int_{0}^{1} dt_2\;t_2^{\frac{\lambda}{2}} \int_{0}^{1} du_2\;u_2^{\frac{1}{2}(-1+\gamma +\lambda )}
   \exp\left( -\frac{s_{2,\infty }}{(1-s_{2,\infty })}\eta (1-t_2)(1-u_2)\right)   \nonumber\\
&&\times \Bigg( \widetilde{w}_{2,2}^{-\frac{1}{2}(1+\lambda )}\left( \widetilde{w}_{2,2} \partial _{\widetilde{w}_{2,2}}\right) \widetilde{w}_{2,2}^{\frac{1}{2}(\beta -\gamma -\delta +1)} \left(  \widetilde{w}_{2,2} \partial _{ \widetilde{w}_{2,2}}\right) \widetilde{w}_{2,2}^{\frac{1}{2}(-\beta +\gamma +\delta +\lambda )}-\frac{q}{4}\Bigg) \nonumber\\
&&\times \int_{0}^{1} dt_1\;t_1^{\frac{1}{2}(-1+\lambda )} \int_{0}^{1} du_1\;u_1^{\frac{1}{2}(-2+\gamma +\lambda )}
 \frac{1}{2\pi i} \oint dv_1 \;\frac{\exp\left( -\frac{v_1}{(1-v_1)} \widetilde{w}_{2,2} (1-t_1)(1-u_1)\right)}{v_1(1-v_1)} \nonumber\\
&&\times \sum_{\alpha _1 =\alpha _0}^{\infty }\left(\frac{s_1}{v_1 }\right)^{\alpha _1}  \Bigg(  \ddot{w}_{1,2}^{-\frac{\lambda }{2}}\left( \ddot{w}_{1,2} \partial _{ \ddot{w}_{1,2}}\right) \ddot{w}_{1,2}^{\frac{1}{2}(\beta -\gamma -\delta +1)} \left( \ddot{w}_{1,2} \partial _{ \ddot{w}_{1,2}}\right) \ddot{w}_{1,2}^{\frac{1}{2}(-\beta +\gamma +\delta -1+\lambda )}-\frac{q}{4}\Bigg) \nonumber\\
&&\times  \sum_{\alpha _0 =0}^{\infty } \frac{s_0^{\alpha _0}}{\alpha _0!} \frac{\Gamma (\alpha _0+\gamma')}{\Gamma(\gamma') } \left( c_0 x^{\lambda } \sum_{i_0=0}^{\alpha _0} \frac{(-\alpha _0)_{i_0} }{(1+\frac{\lambda }{2} )_{i_0}(\frac{1}{2}+\frac{\gamma }{2}+\frac{\lambda }{2} )_{i_0}} \ddot{w}_{1,2} ^{i_0} \right) x^2 \label{eq:40050}\\
&&\mbox{where}\;\widetilde{w}_{2,2}= \eta s_{2,\infty }\prod_{l=2}^{2}t_l u_l \;\;\mbox{and}\; \ddot{w}_{1,2}= \eta v_1  s_{2,\infty }\prod_{l=1}^{2}t_l u_l \nonumber
\end{eqnarray}
Replace $\alpha _i$, $\alpha _j$ and $s_i$ by $\alpha _1$, $\alpha _0$ and ${ \displaystyle \frac{s_1}{v_1 }}$ in (\ref{eq:40041}). Take the new (\ref{eq:40041}) into (\ref{eq:40050}). 
\begin{eqnarray}
&&  \sum_{\alpha _0 =0}^{\infty } \frac{s_0^{\alpha _0}}{\alpha _0!} \frac{\Gamma (\alpha _0+\gamma')}{\Gamma(\gamma') }  \prod _{n=1}^{\infty } \left\{ \sum_{ \alpha _n = \alpha _{n-1}}^{\infty } s_n^{\alpha _n }\right\} y_2(x) \nonumber\\
&&= \prod_{l=2}^{\infty } \frac{1}{(1-s_{l,\infty })} \int_{0}^{1} dt_2\;t_2^{\frac{\lambda}{2}} \int_{0}^{1} du_2\;u_2^{\frac{1}{2}(-1+\gamma +\lambda )}
   \exp\left( -\frac{s_{2,\infty }}{(1-s_{2,\infty })}\eta (1-t_2)(1-u_2)\right)   \nonumber\\
&&\times \Bigg( \widetilde{w}_{2,2}^{-\frac{1}{2}(1+\lambda )}\left( \widetilde{w}_{2,2} \partial _{\widetilde{w}_{2,2}}\right) \widetilde{w}_{2,2}^{\frac{1}{2}(\beta -\gamma -\delta +1)} \left(  \widetilde{w}_{2,2} \partial _{ \widetilde{w}_{2,2}}\right) \widetilde{w}_{2,2}^{\frac{1}{2}(-\beta +\gamma +\delta +\lambda )}-\frac{q}{4}\Bigg) \nonumber\\
&&\times \int_{0}^{1} dt_1\;t_1^{\frac{1}{2}(-1+\lambda )} \int_{0}^{1} du_1\;u_1^{\frac{1}{2}(-2+\gamma +\lambda )}
 \frac{1}{2\pi i} \oint dv_1 \;\frac{\exp\left( -\frac{v_1}{(1-v_1)} \widetilde{w}_{2,2} (1-t_1)(1-u_1)\right)}{ (1-v_1)(v_1-s_1)} \nonumber\\
&&\times  \Bigg( \ddot{w}_{1,2}^{-\frac{\lambda }{2}}\left( \ddot{w}_{1,2} \partial _{ \ddot{w}_{1,2}}\right) \ddot{w}_{1,2}^{\frac{1}{2}(\beta -\gamma -\delta +1)} \left( \ddot{w}_{1,2} \partial _{ \ddot{w}_{1,2}}\right) \ddot{w}_{1,2}^{\frac{1}{2}(-\beta +\gamma +\delta -1+\lambda )}-\frac{q}{4}\Bigg) \nonumber\\
&&\times  \sum_{\alpha _0 =0}^{\infty } \frac{1}{\alpha _0!} \left( \frac{s_{0,1}}{v_1}\right)^{\alpha _0} \frac{\Gamma (\alpha _0+\gamma')}{\Gamma(\gamma') } \left( c_0 x^{\lambda } \sum_{i_0=0}^{\alpha _0} \frac{(-\alpha _0)_{i_0} }{(1+\frac{\lambda }{2} )_{i_0}(\frac{1}{2}+\frac{\gamma }{2}+\frac{\lambda }{2} )_{i_0}} \ddot{w}_{1,2} ^{i_0} \right) x^2 \label{eq:40051}
\end{eqnarray} 
By using Cauchy's integral formula, the contour integrand has poles at
 $ v_1= 1 $  or $ s_1$ and $ s_1$ is only inside the unit circle. As we compute the residue there in (\ref{eq:40051}) we obtain 
\begin{eqnarray}
&&  \sum_{\alpha _0 =0}^{\infty } \frac{s_0^{\alpha _0}}{\alpha _0!} \frac{\Gamma (\alpha _0+\gamma')}{\Gamma(\gamma') }  \prod _{n=1}^{\infty } \left\{ \sum_{ \alpha _n = \alpha _{n-1}}^{\infty } s_n^{\alpha _n }\right\} y_2(x) \nonumber\\
&&= \prod_{l=2}^{\infty } \frac{1}{(1-s_{l,\infty })} \int_{0}^{1} dt_2\;t_2^{\frac{\lambda}{2}} \int_{0}^{1} du_2\;u_2^{\frac{1}{2}(-1+\gamma +\lambda )}
   \exp\left( -\frac{s_{2,\infty }}{(1-s_{2,\infty })}\eta (1-t_2)(1-u_2)\right)   \nonumber\\
&&\times \Bigg( \widetilde{w}_{2,2}^{-\frac{1}{2}(1+\lambda )}\left( \widetilde{w}_{2,2} \partial _{\widetilde{w}_{2,2}}\right) \widetilde{w}_{2,2}^{\frac{1}{2}(\beta -\gamma -\delta +1)} \left(  \widetilde{w}_{2,2} \partial _{ \widetilde{w}_{2,2}}\right) \widetilde{w}_{2,2}^{\frac{1}{2}(-\beta +\gamma +\delta +\lambda )}-\frac{q}{4}\Bigg) \nonumber\\
&&\times \int_{0}^{1} dt_1\;t_1^{\frac{1}{2}(-1+\lambda )} \int_{0}^{1} du_1\;u_1^{\frac{1}{2}(-2+\gamma +\lambda )}
  \frac{\exp\left( -\frac{s_1}{(1-s_1)} \widetilde{w}_{2,2} (1-t_1)(1-u_1)\right)}{ (1-s_1) } \nonumber\\
&&\times  \Bigg( \widetilde{w}_{1,2}^{-\frac{\lambda }{2}}\left( \widetilde{w}_{1,2} \partial _{ \widetilde{w}_{1,2}}\right) \widetilde{w}_{1,2}^{\frac{1}{2}(\beta -\gamma -\delta +1)} \left( \widetilde{w}_{1,2} \partial _{ \widetilde{w}_{1,2}}\right) \widetilde{w}_{1,2}^{\frac{1}{2}(-\beta +\gamma +\delta -1+\lambda )}-\frac{q}{4}\Bigg) \nonumber\\
&&\times  \sum_{\alpha _0 =0}^{\infty } \frac{s_0^{\alpha _0}}{\alpha _0!} \frac{\Gamma (\alpha _0+\gamma')}{\Gamma(\gamma') } \left( c_0 x^{\lambda } \sum_{i_0=0}^{\alpha _0} \frac{(-\alpha _0)_{i_0} }{(1+\frac{\lambda }{2} )_{i_0}(\frac{1}{2}+\frac{\gamma }{2}+\frac{\lambda }{2} )_{i_0}} \widetilde{w}_{1,2} ^{i_0} \right) x^2 \label{eq:40052}\\
&&\mbox{where}\;  \widetilde{w}_{1,2}= \eta s_{1,\infty }\prod_{l=1}^{2}t_l u_l \nonumber
\end{eqnarray} 
Acting the summation operator $\displaystyle{ \sum_{\alpha _0 =0}^{\infty } \frac{s_0^{\alpha _0}}{\alpha _0!} \frac{\Gamma (\alpha _0+\gamma')}{\Gamma(\gamma') }  \prod _{n=1}^{\infty } \left\{ \sum_{ \alpha _n = \alpha _{n-1}}^{\infty } s_n^{\alpha _n }\right\}}$ on (\ref{eq:40034}),
\begin{eqnarray}
&&  \sum_{\alpha _0 =0}^{\infty } \frac{s_0^{\alpha _0}}{\alpha _0!} \frac{\Gamma (\alpha _0+\gamma')}{\Gamma(\gamma') }  \prod _{n=1}^{\infty } \left\{ \sum_{ \alpha _n = \alpha _{n-1}}^{\infty } s_n^{\alpha _n }\right\} y_3(x) \nonumber\\
&&= \prod_{l=3}^{\infty } \frac{1}{(1-s_{l,\infty })} \int_{0}^{1} dt_3\;t_3^{\frac{1}{2}(1+\lambda )} \int_{0}^{1} du_3\;u_3^{\frac{1}{2}( \gamma +\lambda )}
   \exp\left( -\frac{s_{3,\infty }}{(1-s_{3,\infty })}\eta (1-t_3)(1-u_3)\right)   \nonumber\\
&&\times \Bigg( \widetilde{w}_{3,3}^{-\frac{1}{2}(2+\lambda )}\left( \widetilde{w}_{3,3} \partial _{\widetilde{w}_{3,3}}\right) \widetilde{w}_{3,3}^{\frac{1}{2}(\beta -\gamma -\delta +1)} \left(  \widetilde{w}_{3,3} \partial _{ \widetilde{w}_{3,3}}\right) \widetilde{w}_{3,3}^{\frac{1}{2}(-\beta +\gamma +\delta +1+\lambda )}-\frac{q}{4}\Bigg) \nonumber\\
&&\times \int_{0}^{1} dt_2\;t_2^{\frac{\lambda }{2}} \int_{0}^{1} du_2\;u_2^{\frac{1}{2}(-1+\gamma +\lambda )}
  \frac{\exp\left( -\frac{s_2}{(1-s_2)} \widetilde{w}_{3,3} (1-t_2)(1-u_2)\right)}{ (1-s_2) } \nonumber\\
&&\times  \Bigg( \widetilde{w}_{2,3}^{-\frac{1}{2}(-1+\lambda )}\left( \widetilde{w}_{2,3} \partial _{ \widetilde{w}_{2,3}}\right) \widetilde{w}_{2,3}^{\frac{1}{2}(\beta -\gamma -\delta +1)} \left( \widetilde{w}_{2,3} \partial _{ \widetilde{w}_{2,3}}\right) \widetilde{w}_{2,3}^{\frac{1}{2}(-\beta +\gamma +\delta +\lambda )}-\frac{q}{4}\Bigg) \nonumber\\
&&\times \int_{0}^{1} dt_1\;t_1^{\frac{1}{2}(-1+\lambda )} \int_{0}^{1} du_1\;u_1^{\frac{1}{2}(-2+\gamma +\lambda )}
  \frac{\exp\left( -\frac{s_1}{(1-s_1)} \widetilde{w}_{2,3} (1-t_1)(1-u_1)\right)}{ (1-s_1) } \nonumber\\
&&\times  \Bigg( \widetilde{w}_{1,3}^{-\frac{\lambda }{2}}\left( \widetilde{w}_{1,3} \partial _{ \widetilde{w}_{1,3}}\right) \widetilde{w}_{1,3}^{\frac{1}{2}(\beta -\gamma -\delta +1)} \left( \widetilde{w}_{1,3} \partial _{ \widetilde{w}_{1,3}}\right) \widetilde{w}_{1,3}^{\frac{1}{2}(-\beta +\gamma +\delta -1+\lambda )}-\frac{q}{4}\Bigg) \nonumber\\
&&\times  \sum_{\alpha _0 =0}^{\infty } \frac{s_0^{\alpha _0}}{\alpha _0!} \frac{\Gamma (\alpha _0+\gamma')}{\Gamma(\gamma') } \left( c_0 x^{\lambda } \sum_{i_0=0}^{\alpha _0} \frac{(-\alpha _0)_{i_0} }{(1+\frac{\lambda }{2} )_{i_0}(\frac{1}{2}+\frac{\gamma }{2}+\frac{\lambda }{2} )_{i_0}} \widetilde{w}_{1,3} ^{i_0} \right) x^3 \label{eq:40053}\\
&&\mbox{where}\;\widetilde{w}_{3,3}= \eta s_{3,\infty }\prod_{l=3}^{3}t_l u_l , \;\;\widetilde{w}_{2,3}= \eta s_{2,\infty }\prod_{l=2}^{3}t_l u_l \;\;\mbox{and}\;\widetilde{w}_{1,3}= \eta s_{1,\infty }\prod_{l=1}^{3}t_l u_l \nonumber
\end{eqnarray} 
By repeating this process for all higher terms of integral forms of sub-summation $y_m(x)$ terms where $m > 3$, I obtain every  $\displaystyle{ \sum_{\alpha _0 =0}^{\infty } \frac{s_0^{\alpha _0}}{\alpha _0!} \frac{\Gamma (\alpha _0+\gamma')}{\Gamma(\gamma') }  \prod _{n=1}^{\infty } \left\{ \sum_{ \alpha _n = \alpha _{n-1}}^{\infty } s_n^{\alpha _n }\right\}}  y_m(x)$ terms. 
Substitute (\ref{eq:40044}), (\ref{eq:40047}), (\ref{eq:40052}), (\ref{eq:40053}) and including all $\displaystyle{ \sum_{q_0 =0}^{\infty } \frac{(\gamma')_{q_0}}{(q_0)!} s_0^{q_0} \prod _{n=1}^{\infty } \left\{ \sum_{ q_n = q_{n-1}}^{\infty } s_n^{q_n }\right\}}  y_m(x)$ terms where $m > 3$ into (\ref{eq:40043}).
\qed
\end{proof}
\begin{remark}
The generating function for the CHP of type 1 of the first kind about $x=0$ as $\alpha = -(2 \alpha _j+j) $ where $j,\alpha _j \in \mathbb{N}_{0}$ is
\begin{eqnarray}
&&\sum_{\alpha _0 =0}^{\infty } \frac{s_0^{\alpha _0}}{\alpha _0!} \frac{\Gamma (\alpha _0+\frac{1}{2}+\frac{\gamma }{2})}{\Gamma(\frac{1}{2}+\frac{\gamma }{2}) }  \prod _{n=1}^{\infty } \left\{ \sum_{ \alpha _n = \alpha _{n-1}}^{\infty } s_n^{\alpha _n }\right\} H_c^{(a)}F_{\alpha _j}\left(\alpha = -(2 \alpha _j+j), \beta, \gamma, \delta, q; \eta =\frac{1}{2}\beta x^2  \right)  \nonumber\\
&&= \prod_{l=1}^{\infty } \frac{1}{(1-s_{l,\infty })} \mathbf{A} \left( s_{0,\infty } ;\eta\right)\nonumber\\
&&+\Bigg\{ \prod_{l=1}^{\infty } \frac{1}{(1-s_{l,\infty })} \int_{0}^{1} dt_1\;t_1^{-\frac{1}{2} } \int_{0}^{1} du_1\;u_1^{\frac{1}{2}(-2+\gamma )} \overleftrightarrow {\mathbf{\Gamma}}_1 \left(s_{1,\infty };t_1,u_1,\eta\right) \nonumber\\
&&\times \Bigg( \left(  \widetilde{w}_{1,1} \partial _{ \widetilde{w}_{1,1}}\right) \widetilde{w}_{1,1}^{\frac{1}{2}(\beta -\gamma -\delta +1)} \left(  \widetilde{w}_{1,1} \partial _{ \widetilde{w}_{1,1}}\right) \widetilde{w}_{1,1}^{\frac{1}{2}(-\beta +\gamma +\delta -1 )}-\frac{q}{4}\Bigg)\; \mathbf{A} \left( s_{0} ;\widetilde{w}_{1,1}\right)\Bigg\}x\nonumber\\
&&+ \sum_{n=2}^{\infty } \Bigg\{ \prod_{l=n}^{\infty } \frac{1}{(1-s_{l,\infty })} \int_{0}^{1} dt_n\;t_n^{\frac{1}{2}(n-2)} \int_{0}^{1} du_n\;u_n^{\frac{1}{2}(n-3+\gamma )} \overleftrightarrow {\mathbf{\Gamma}}_n \left(s_{n,\infty };t_n,u_n,\eta \right)  \nonumber\\
&&\times \Bigg( \widetilde{w}_{n,n}^{-\frac{1}{2}(n-1 )}\left(  \widetilde{w}_{n,n} \partial _{ \widetilde{w}_{n,n}}\right) \widetilde{w}_{n,n}^{\frac{1}{2}(\beta -\gamma -\delta +1)} \left(  \widetilde{w}_{n,n} \partial _{ \widetilde{w}_{n,n}}\right) \widetilde{w}_{n,n}^{\frac{1}{2}(-\beta +\gamma +\delta +n-2 )}-\frac{q}{4}\Bigg)\nonumber\\
&&\times \prod_{k=1}^{n-1} \Bigg\{ \int_{0}^{1} dt_{n-k}\;t_{n-k}^{\frac{1}{2}(n-k-2 )} \int_{0}^{1} du_{n-k} \;u_{n-k}^{\frac{1}{2}(n-k-3+\gamma ) } \overleftrightarrow {\mathbf{\Gamma}}_{n-k} \left(s_{n-k};t_{n-k},u_{n-k},\widetilde{w}_{n-k+1,n} \right)\nonumber\\
&&\times \Bigg( \widetilde{w}_{n-k,n}^{-\frac{1}{2}(n-k-1 )}\left(  \widetilde{w}_{n-k,n} \partial _{ \widetilde{w}_{n-k,n}}\right) \widetilde{w}_{n-k,n}^{\frac{1}{2}(\beta -\gamma -\delta +1)} \left(  \widetilde{w}_{n-k,n} \partial _{ \widetilde{w}_{n-k,n}}\right) \widetilde{w}_{n-k,n}^{\frac{1}{2}(-\beta +\gamma +\delta +n-k-2 )}-\frac{q}{4}\Bigg) \Bigg\}\nonumber\\
&&\times \mathbf{A} \left( s_{0} ;\widetilde{w}_{1,n}\right) \Bigg\} x^n  \label{eq:40054}
\end{eqnarray}
where
\begin{equation}
\begin{cases} 
{ \displaystyle \overleftrightarrow {\mathbf{\Gamma}}_1 \left(s_{1,\infty };t_1,u_1,\eta\right)= \exp\left( -\frac{s_{1,\infty }}{(1-s_{1,\infty })}\eta (1-t_1)(1-u_1)\right) }\cr
{ \displaystyle  \overleftrightarrow {\mathbf{\Gamma}}_n \left(s_{n,\infty };t_n,u_n,\eta \right) = \exp\left( -\frac{s_{n,\infty }}{(1-s_{n,\infty })}\eta (1-t_n)(1-u_n)\right)}\cr
{ \displaystyle \overleftrightarrow {\mathbf{\Gamma}}_{n-k} \left(s_{n-k};t_{n-k},u_{n-k},\widetilde{w}_{n-k+1,n} \right) = \frac{ \exp\left( -\frac{s_{n-k}}{(1-s_{n-k})}\widetilde{w}_{n-k+1,n} (1-t_{n-k})(1-u_{n-k})\right)}{(1-s_{n-k})} }
\end{cases}\nonumber 
\end{equation}
and
\begin{equation}
\begin{cases} 
{ \displaystyle \mathbf{A} \left( s_{0,\infty } ;\eta\right)=  (1-s_{0,\infty })^{-\frac{1}{2}(1+\gamma )} \exp\left( -\frac{\eta s_{0,\infty }}{(1-s_{0,\infty })}\right)  }\cr
{ \displaystyle  \mathbf{A} \left( s_{0} ;\widetilde{w}_{1,1}\right) = (1-s_0 )^{-\frac{1}{2}(1+\gamma )} \exp\left( -\frac{\widetilde{w}_{1,1} s_0}{(1-s_0)}\right) } \cr
{ \displaystyle \mathbf{A} \left( s_{0} ;\widetilde{w}_{1,n}\right) = (1-s_0 )^{-\frac{1}{2}(1+\gamma )} \exp\left( -\frac{\widetilde{w}_{1,n} s_0}{(1-s_0)}\right)  }
\end{cases}\nonumber 
\end{equation}
\end{remark}
\begin{proof}
The generating function for Confluent hypergeometric (Kummer's) polynomial of the first kind is given by
\begin{equation}
\sum_{\alpha _0=0}^{\infty }\frac{t^{\alpha _0}}{\alpha _0!}\frac{\Gamma (\alpha _0+\gamma')}{\Gamma (\gamma')}\; _1F_1\left(-\alpha _0; \gamma'; z \right)= \frac{\exp\left( -\frac{z t}{(1-t)}\right)}{(1-t)^{\gamma'}} \label{eq:40055}
\end{equation}

Replace $t$, $\gamma'$ and $z$  by $\displaystyle{s_{0,\infty }}$, $\displaystyle{\frac{1}{2}+\frac{\gamma }{2}}$ and $\eta $ in (\ref{eq:40055}). 
\begin{equation}
\sum_{\alpha _0=0}^{\infty }\frac{s_{0,\infty }^{\alpha _0}}{\alpha _0!} \frac{\Gamma (\alpha _0+\frac{1}{2}+\frac{\gamma }{2})}{\Gamma (\frac{1}{2}+\frac{\gamma }{2})}\; _1F_1\left(-\alpha _0; \frac{1}{2}+\frac{\gamma }{2}; \eta  \right)= (1-s_{0,\infty })^{-\frac{1}{2}(1+\gamma )} \exp\left( -\frac{\eta  s_{0,\infty }}{(1-s_{0,\infty })}\right)   \label{eq:40056}
\end{equation}
Replace $t$, $\gamma'$ and $z$  by $\displaystyle{s_0}$, $\displaystyle{\frac{1}{2}+\frac{\gamma }{2}}$ and $\widetilde{w}_{1,1} $ in (\ref{eq:40055}). 
\begin{equation}
\sum_{\alpha _0=0}^{\infty }\frac{s_0^{\alpha _0}}{\alpha _0!} \frac{\Gamma (\alpha _0+\frac{1}{2}+\frac{\gamma }{2})}{\Gamma (\frac{1}{2}+\frac{\gamma }{2})}\; _1F_1\left(-\alpha _0; \frac{1}{2}+\frac{\gamma }{2}; \widetilde{w}_{1,1} \right)= (1-s_0)^{-\frac{1}{2}(1+\gamma )} \exp\left( -\frac{\widetilde{w}_{1,1}  s_0}{(1-s_0)}\right)   \label{eq:40057}
\end{equation}
Replace $t$, $\gamma'$ and $z$  by $\displaystyle{s_0}$, $\displaystyle{\frac{1}{2}+\frac{\gamma }{2}}$ and $\widetilde{w}_{1,n} $ in (\ref{eq:40055}). 
\begin{equation}
\sum_{\alpha _0=0}^{\infty }\frac{s_0^{\alpha _0}}{\alpha _0!} \frac{\Gamma (\alpha _0+\frac{1}{2}+\frac{\gamma }{2})}{\Gamma (\frac{1}{2}+\frac{\gamma }{2})}\; _1F_1\left(-\alpha _0; \frac{1}{2}+\frac{\gamma }{2}; \widetilde{w}_{1,n} \right)= (1-s_0)^{-\frac{1}{2}(1+\gamma )} \exp\left( -\frac{\widetilde{w}_{1,n}  s_0}{(1-s_0)}\right)   \label{eq:40058}
\end{equation}
Put $c_0$= 1, $\lambda $=0 and $\displaystyle{\gamma' =\frac{1}{2}+\frac{\gamma }{2}}$ in (\ref{eq:40042}). Substitute (\ref{eq:40056}), (\ref{eq:40057}) and (\ref{eq:40058}) into the new (\ref{eq:40042}).\qed
\end{proof}
\begin{remark}
The generating function for the CHP of type 1 of the second kind about $x=0$ as     $\alpha = -(2 \alpha _j+j +1-\gamma ) $ where $j,\alpha _j \in \mathbb{N}_{0}$ is
\begin{eqnarray}
&&\sum_{\alpha _0 =0}^{\infty } \frac{s_0^{\alpha _0}}{\alpha _0!} \frac{\Gamma (\alpha _0+\frac{3}{2}-\frac{\gamma }{2})}{\Gamma(\frac{3}{2}-\frac{\gamma }{2}) }  \prod _{n=1}^{\infty } \left\{ \sum_{ \alpha _n = \alpha _{n-1}}^{\infty } s_n^{\alpha _n }\right\} H_c^{(a)}S_{\alpha _j}\left(\alpha = -(2 \alpha _j+j+1-\gamma ), \beta, \gamma, \delta, q; \eta =\frac{1}{2}\beta x^2  \right) \nonumber\\
&&= x^{1-\gamma } \Bigg\{ \prod_{l=1}^{\infty } \frac{1}{(1-s_{l,\infty })} \mathbf{B} \left( s_{0,\infty } ;\eta\right)  \nonumber\\
&&+ \Bigg\{ \prod_{l=1}^{\infty } \frac{1}{(1-s_{l,\infty })} \int_{0}^{1} dt_1\;t_1^{-\frac{\gamma }{2}} \int_{0}^{1} du_1\;u_1^{-\frac{1}{2}} \overleftrightarrow {\mathbf{\Gamma}}_1 \left(s_{1,\infty };t_1,u_1,\eta\right) \nonumber\\
&&\times \Bigg(  \widetilde{w}_{1,1}^{-\frac{1}{2}(1-\gamma )}\left(  \widetilde{w}_{1,1} \partial _{ \widetilde{w}_{1,1}}\right) \widetilde{w}_{1,1}^{\frac{1}{2}(\beta -\gamma -\delta +1)} \left(  \widetilde{w}_{1,1} \partial _{ \widetilde{w}_{1,1}}\right) \widetilde{w}_{1,1}^{\frac{1}{2}(-\beta +\delta )}-\frac{q}{4}\Bigg)\; \mathbf{B} \left( s_{0} ;\widetilde{w}_{1,1}\right)\Bigg\}x\nonumber\\
&&+ \sum_{n=2}^{\infty } \Bigg\{ \prod_{l=n}^{\infty } \frac{1}{(1-s_{l,\infty })} \int_{0}^{1} dt_n\;t_n^{\frac{1}{2}(n-1-\gamma  )} \int_{0}^{1} du_n\;u_n^{\frac{1}{2}(n-2)} \overleftrightarrow {\mathbf{\Gamma}}_n \left(s_{n,\infty };t_n,u_n,\eta \right)  \nonumber\\
&&\times \Bigg( \widetilde{w}_{n,n}^{-\frac{1}{2}(n-\gamma )}\left(  \widetilde{w}_{n,n} \partial _{ \widetilde{w}_{n,n}}\right) \widetilde{w}_{n,n}^{\frac{1}{2}(\beta -\gamma -\delta +1)} \left(  \widetilde{w}_{n,n} \partial _{ \widetilde{w}_{n,n}}\right) \widetilde{w}_{n,n}^{\frac{1}{2}(-\beta +\delta +n-1)}-\frac{q}{4}\Bigg)\nonumber\\
&&\times \prod_{k=1}^{n-1} \Bigg\{ \int_{0}^{1} dt_{n-k}\;t_{n-k}^{\frac{1}{2}(n-k-1-\gamma )} \int_{0}^{1} du_{n-k} \;u_{n-k}^{\frac{1}{2}(n-k-2) } \overleftrightarrow {\mathbf{\Gamma}}_{n-k} \left(s_{n-k};t_{n-k},u_{n-k},\widetilde{w}_{n-k+1,n} \right)\nonumber\\
&&\times  \Bigg( \widetilde{w}_{n-k,n}^{-\frac{1}{2}(n-k-\gamma )}\left(  \widetilde{w}_{n-k,n} \partial _{ \widetilde{w}_{n-k,n}}\right) \widetilde{w}_{n-k,n}^{\frac{1}{2}(\beta -\gamma -\delta +1)} \left(  \widetilde{w}_{n-k,n} \partial _{ \widetilde{w}_{n-k,n}}\right) \widetilde{w}_{n-k,n}^{\frac{1}{2}(-\beta +\delta +n-k-1)}-\frac{q}{4}\Bigg) \Bigg\}\nonumber\\
&&\times \mathbf{B} \left( s_{0} ;\widetilde{w}_{1,n}\right) \Bigg\} x^n \Bigg\}  \label{eq:40059}
\end{eqnarray}
where
\begin{equation}
\begin{cases} 
{ \displaystyle \overleftrightarrow {\mathbf{\Gamma}}_1 \left(s_{1,\infty };t_1,u_1,\eta\right)= \exp\left( -\frac{s_{1,\infty }}{(1-s_{1,\infty })}\eta (1-t_1)(1-u_1)\right) }\cr
{ \displaystyle  \overleftrightarrow {\mathbf{\Gamma}}_n \left(s_{n,\infty };t_n,u_n,\eta \right) = \exp\left( -\frac{s_{n,\infty }}{(1-s_{n,\infty })}\eta (1-t_n)(1-u_n)\right)}\cr
{ \displaystyle \overleftrightarrow {\mathbf{\Gamma}}_{n-k} \left(s_{n-k};t_{n-k},u_{n-k},\widetilde{w}_{n-k+1,n} \right) = \frac{ \exp\left( -\frac{s_{n-k}}{(1-s_{n-k})}\widetilde{w}_{n-k+1,n} (1-t_{n-k})(1-u_{n-k})\right)}{(1-s_{n-k})} }
\end{cases}\nonumber 
\end{equation}
and
\begin{equation}
\begin{cases} 
{ \displaystyle \mathbf{B} \left( s_{0,\infty } ;\eta\right)= (1-s_{0,\infty })^{-\frac{1}{2}(3-\gamma )} \exp\left( -\frac{\eta s_{0,\infty }}{(1-s_{0,\infty })}\right) }\cr
{ \displaystyle  \mathbf{B} \left( s_{0} ;\widetilde{w}_{1,1}\right) = (1-s_0 )^{-\frac{1}{2}(3-\gamma )} \exp\left( -\frac{\widetilde{w}_{1,1} s_0}{(1-s_0)}\right)  } \cr
{ \displaystyle \mathbf{B} \left( s_{0} ;\widetilde{w}_{1,n}\right) = (1-s_0 )^{-\frac{1}{2}(3-\gamma )} \exp\left( -\frac{\widetilde{w}_{1,n} s_0}{(1-s_0)}\right)  }
\end{cases}\nonumber 
\end{equation}
\end{remark}
\begin{proof}
Replace $t$, $\gamma'$ and $z$  by $\displaystyle{s_{0,\infty }}$, $\displaystyle{\frac{3}{2}-\frac{\gamma }{2}}$ and $\eta $ in (\ref{eq:40055}). 
\begin{equation}
\sum_{\alpha _0=0}^{\infty }\frac{s_{0,\infty }^{\alpha _0}}{\alpha _0!} \frac{\Gamma (\alpha _0+\frac{3}{2}-\frac{\gamma }{2})}{\Gamma (\frac{3}{2}-\frac{\gamma }{2})}\; _1F_1\left(-\alpha _0; \frac{3}{2}-\frac{\gamma }{2}; \eta  \right)= (1-s_{0,\infty })^{-\frac{1}{2}(3-\gamma )} \exp\left( -\frac{\eta  s_{0,\infty }}{(1-s_{0,\infty })}\right)   \label{eq:40060}
\end{equation}
Replace $t$, $\gamma'$ and $z$  by $\displaystyle{s_0}$, $\displaystyle{\frac{3}{2}-\frac{\gamma }{2}}$ and $\widetilde{w}_{1,1} $ in (\ref{eq:40055}). 
\begin{equation}
\sum_{\alpha _0=0}^{\infty }\frac{s_0^{\alpha _0}}{\alpha _0!} \frac{\Gamma (\alpha _0+\frac{3}{2}-\frac{\gamma }{2})}{\Gamma (\frac{3}{2}-\frac{\gamma }{2})}\; _1F_1\left(-\alpha _0; \frac{3}{2}-\frac{\gamma }{2}; \widetilde{w}_{1,1} \right)= (1-s_0)^{-\frac{1}{2}(3-\gamma )} \exp\left( -\frac{\widetilde{w}_{1,1}  s_0}{(1-s_0)}\right)   \label{eq:40061}
\end{equation}
Replace $t$, $\gamma'$ and $z$  by $\displaystyle{s_0}$, $\displaystyle{\frac{3}{2}-\frac{\gamma }{2}}$ and $\widetilde{w}_{1,n} $ in (\ref{eq:40055}). 
\begin{equation}
\sum_{\alpha _0=0}^{\infty }\frac{s_0^{\alpha _0}}{\alpha _0!} \frac{\Gamma (\alpha _0+\frac{3}{2}-\frac{\gamma }{2})}{\Gamma (\frac{3}{2}-\frac{\gamma }{2})}\; _1F_1\left(-\alpha _0; \frac{3}{2}-\frac{\gamma }{2}; \widetilde{w}_{1,n} \right)= (1-s_0)^{-\frac{1}{2}(3-\gamma )} \exp\left( -\frac{\widetilde{w}_{1,n}  s_0}{(1-s_0)}\right)   \label{eq:40062}
\end{equation}
Put $c_0$= 1, $\lambda =1-\gamma $  and $\displaystyle{\gamma' =\frac{3}{2}-\frac{\gamma }{2}}$ in (\ref{eq:40042}). Substitute (\ref{eq:40060}), (\ref{eq:40061}) and (\ref{eq:40062}) into the new (\ref{eq:40042}).\qed
\end{proof}
\section{Summary}
In this chapter I apply the 3TRF to (1) the power series expansion in closed forms of the CHE about regular singularity $x=0$, (2) its integral representation for infinite series and polynomial of type 1 and (3) generating functions of the CHP of type 1. 

 As we see the power series expansions of the CHE for infinite series and polynomial, denominators and numerators in every $B_n$ terms arise with Pochhammer symbol: the meaning of this is that the analytic solutions of the CHF can be described as hypergoemetric functions in a strict mathematical way. I show how to construct representations in closed form integrals resulting in a precise since we have power series expansions with Pochhammer symbols in numerators and denominators. We can transform the CHF into all other well-known special functions with two recursive coefficients because a $_1F_1$ function recurs in each of sub-integral forms of the CHF.

I show how to analyze generating functions for the CHP of type 1 from the general representation of an integral form of the CHP. We can derive orthogonal relations, recursion relations and expectation values of physical quantities from the generating functions for the CHP: the processes in order to obtain orthogonal and recursion relations of the CHP are similar to the case of a normalized wave function for the hydrogen-like atoms.  

The second independent solution of the CHF is only valid for non-integer values of $1-\gamma $ in the domain $|x|<1$. Its analytic solution is defined by its behavior at infinity in the domain $|x|>1$. The local solution of the CHF about regular singularity $x=0$ is called angular solution of the CHE, denoted as $Hc^{(a)}(\alpha ,\beta ,\gamma ,\delta ,q; x)$. And the local solution of the CHE about the singular point at infinity is called radial solution of the CHE, denoted as $Hc^{(r)}(\alpha ,\beta ,\gamma ,\delta ,q; x)$. \cite{Ronv1995} Various local solutions of the CHE can be constructed using the combination of transformed independent parameters which resembles Kummer's 24 local solutions of the hypergeometric equation. All possible local solutions of the CHE (Regge-Wheeler and Teukolsky equations) were constructed by Fiziev.\cite{Fizi2009,Fizi2010}
In general the power series expansion for the CHE about the singular point at infinity is not convergent by only asymptotic. \cite{Ronv1995,Fizi2010a} Putting $z=\frac{1}{x}$ into (\ref{eq:4002}),
\begin{equation}
z^2 \frac{d^2{y}}{d{z}^2} + z\left((2-\gamma ) -\frac{\beta }{z} + \frac{\delta }{z-1}\right) \frac{d{y}}{d{z}} +  \frac{q z-\alpha \beta }{z(z-1)} y = 0 \label{eq:40063}
\end{equation}
Assume that its solution is
\begin{equation}
y(z)= \sum_{n=0}^{\infty } c_n z^{n+\lambda } \hspace{1cm}\mbox{where}\; \lambda =\mbox{indicial}\;\mbox{root}\label{eq:40064}
\end{equation}
Plug (\ref{eq:40064})  into (\ref{eq:40063}). We again obtain a three-term recurrence relation for the coefficients $c_n$:
\begin{equation}
c_{n+1}=A_n \;c_n +B_n \;c_{n-1} \hspace{1cm};n\geq 1 \label{eq:40065}
\end{equation}
where,
\begin{subequations}
\begin{equation}
A_n =\frac{(n+\lambda )(n+\lambda +\beta -\gamma -\delta +1)-q}{\beta (n+1-\alpha +\lambda )}\label{eq:40066a}
\end{equation}
\begin{equation}
B_n = \frac{-(n+\lambda -1)(n+\lambda -\gamma )}{\beta (n+\lambda +1-\alpha )} \label{eq:40066b}
\end{equation}
\begin{equation}
c_1= A_0 \;c_0 \label{eq:40066c}
\end{equation}
\end{subequations}
We only have an indicial root which is $\lambda =\alpha $.
Let's test for the asymptotic behavior of the function $y(z)$. As $n\gg 1$ (for sufficiently large),
\begin{subequations}
\begin{equation}
A= \lim_{n\gg 1}A_n =\frac{n}{\beta }\label{eq:40067a}
\end{equation}
\begin{equation}
B= \lim_{n\gg 1}B_n = -\frac{n}{\beta} \label{eq:40067b}
\end{equation}
\end{subequations} 
As we see (\ref{eq:40067a}) and (\ref{eq:40067b}), the function $y(z)$ for infinite series is divergent. Also, for polynomial of type 1, (\ref{eq:40067b}) is negligible for the minimum $y(z)$ because $B_n$ term will be terminated at the specific eigenvalues. (\ref{eq:40067a}) is only survived. And its recurrence relation is
\begin{equation}
c_{n+1}= \frac{n}{\beta } \;c_n \label{eq:40068}
\end{equation}
The local solution for (\ref{eq:40068}) is divergent. By similar reason, polynomial of type 2 also will be divergent.
Therefore, the local solution of the CHE about the singular point at infinity is only valid for the type 3 polynomial where $A_n$ and $B_n$ terms terminated at the same time. Also we need a condition $\alpha >0$ in order to make a function $y(z)$ convergent. If it does not, the function $y(z)$ will be divergent as $x$ goes to infinity. For the black hole problems, $x$ is correspond to the radius.\cite{Fizi2009,Fizi2010,Fizi2010a} In the quantum mechanical points of views any physical quantities should be disappeared as the distance goes to the infinity as we all know.  In the next series I will construct the Frobenius solution in closed forms of the CHE about the singular point at infinity, its integral representation and its generating functions for the type 3 polynomial. I also derive the mathematical formula of the spectral parameter $q$ for this type polynomial.
\addcontentsline{toc}{section}{Bibliography}
\bibliographystyle{model1a-num-names}
\bibliography{<your-bib-database>}
\bibliographystyle{model1a-num-names}
\bibliography{<your-bib-database>}

\chapter{Confluent Heun function using reversible three-term recurrence formula}
\chaptermark{Confluent Heun function using  R3TRF} 


In this chapter I apply reversible three term recurrence formula to (1) the power series expansions in closed forms, (2) its integral forms of the Confluent Heun function (CHF) for infinite series and polynomial of type 2\footnote{polynomial of type 2 is a polynomial which makes $A_n$ term terminated in three term recursion relation of the power series in a linear differential equation.} including all higher terms of $B_n$'s\footnote{`` higher terms of $B_n$'s'' means at least two terms of $B_n$'s.} and (3) the generating function for Confluent Heun polynomial (CHP) of type 2.
 
\section{Non-symmetrical canonical form of Confluent Heun Equation}
In general, there are three types of the Confluent Heun equation (CHE): (1)Generalized spheroidal equation (GSE), (2) the generalized spheroidal wave equation (GSWE) in the Leaver version, (3) Non-symmetrical canonical form of the CHE.

In 1928, A.H. Wilson worked on the wave equation for the ion of the hydrogen molecule $H_2^{+}$ with considering the protons as at rest: he assume that the mass of the electron is negligible compared with that of the protons. Neglecting the mass of electrons in his wave equation, he obtained the second ordinary differential equation in elliptic coordinates for the approximative wave equation. And he introduced on the solution of the more general form of the Helmholtz equation in prolate spheroidal coordinates which is written by\cite{Wils1928,Wils1928a}
\begin{equation}
\frac{d}{d{\xi }}\left( (1-\xi^2 )\frac{d X}{d{\xi}} \right) + \left( \lambda ^2 \xi^2 -2p\lambda \xi -\frac{n_3^2}{1-\xi ^2}+\mu ^{'}\right) X = 0 \nonumber
\end{equation}
When $p=0$ this equation becomes the equation giving the solution of $\nabla ^2 X- \lambda ^2 X =0$ in prolate spheroidal coordinates. Therefore, this equation is called the generalized spheroidal equation (GSE). The GSE has three singular points: two regular singular points which are $\pm 1$ and one irregular singular point which is $\infty$.

In 1986, E.W. Leaver introduced the different version of the GSE which is given by \cite{Leav1986}
\begin{equation}
z(z-z_0)\frac{d^2{U}}{d{z}^2} + (B_1+B_2 z)\frac{d{U}}{d{z}}+ \left( B_3 -2\eta \omega (z-z_0)+\omega ^2z(z-z_0)\right) U(z) = 0 \hspace{.5cm}\mbox{where}\; \omega \ne0 \nonumber
\end{equation}
This equation has three singular points: two regular singular points which are 0 and $z_0$ with exponents $\left\{ 0, 1+\frac{B_1}{z_0}\right\}$ and $\left\{ 0, 1-B_2-\frac{B_1}{z_0}\right\}$, and one irregular singular point which is $\infty $. It arises in the scattering problem of non-relativistic electrons from polar molecule and ions, in wave equation on the background geometry of the Schwarzschild black hole (Teukolsky's equation) and etc.

In Ref.\cite{Heun1889,Ronv1995}, Heun equation is a second-order linear ordinary differential equation of the form
\begin{equation}
\frac{d^2{y}}{d{x}^2} + \left(\frac{\gamma }{x} +\frac{\delta }{x-1} + \frac{\epsilon }{x-a}\right) \frac{d{y}}{d{x}} +  \frac{\alpha \beta x-q}{x(x-1)(x-a)} y = 0 \label{eq:3001}
\end{equation}
It has four regular singular points which are 0, 1, a and $\infty $ with exponents $\{ 0, 1-\gamma \}$, $\{ 0, 1-\delta \}$, $\{ 0, 1-\epsilon \}$ and $\{ \alpha, \beta \}$.

Heun equation has the four kinds of confluent forms: (1) Confluent Heun (two regular and one irregular singularities), (2) Doubly Confluent Heun (two irregular singularities), (3) Biconfluent Heun (one regular and one irregular singularities),\footnote{Biconfluent Heun equation is derived from the grand confluent hypergeometric equation by changing all coefficients $\mu =1$ and $ \varepsilon \omega = -q$.\cite{Chou2012i,Chou2012j}} (4) Triconfluent Heun equations (one irregular singularity).
We can derive these four confluent forms from Heun equation by combining two or more regular singularities to take form an irregular singularity. Its process, converting Heun equation to other confluent forms, is similar to deriving of confluent hypergeometric equation from the hypergeometric equation. For the non-symmetrical canonical form of the Confluent Heun Equation (CHE).\cite{Ronv1995,Deca1978,Decar1978}  
\begin{equation}
\frac{d^2{y}}{d{x}^2} + \left(\beta  +\frac{\gamma }{x} + \frac{\delta }{x-1}\right) \frac{d{y}}{d{x}} +  \frac{\alpha \beta x-q}{x(x-1)} y = 0 \label{eq:3002}
\end{equation}
(\ref{eq:3002}) has three singular points: two regular singular points which are 0 and 1 with exponents $\{0, 1-\gamma\}$ and $\{0, 1-\delta \}$, and one irregular singular point which is $\infty$ with an exponent $\alpha$. Its solution is denoted as $H_{c}(\alpha,\beta,\gamma,\delta,q;x)$.\footnote{Several authors denote as coefficients $4p$ and $\sigma $ instead of $\beta $ and $q$. And they define the solution of the CHE as $H_{c}^{(a)}(p,\alpha,\gamma,\delta,\sigma;x)$.} We can convert the GSE and its Leaver version to the non-symmetrical canonical form of the CHE by changing coefficients, adding the new variables and combining regular singularities. 

\section[CHF with a regular singular point at 0]{Confluent Heun function with a regular singular point at 0} 
In chapter 4 I showed the analytic solutions of the CHE for infinite series and polynomial of type 1\footnote{polynomial of type 1 is a polynomial which makes $B_n$ term terminated in three term recursion relation of the power series in a linear differential equation.} including all higher terms of $A_n$'s by applying 3TRF \cite{Chou2012}; (1) the power series expansion of the CHE, (2) the integral representation of the CHF, (3) its generating function for the Confluent Heun polynomial (CHP) of type 1. Expressing the CHF in integral forms resulting in a precise and simplified transformation of the CHF to a confluent hypergeometric function. Also, the orthogonal relations of the  CHF can be obtained from the integral forms.

In this chapter, by applying R3TRF, I construct the power series expansion in closed forms of non-symmetrical canonical form of the CHE (for infinite series and polynomial of type 2 including all higher terms of $B_n$'s), its integral forms and the generating function for the CHP of type 2 analytically. 

Assume that its solution is
\begin{equation}
y(x)= \sum_{n=0}^{\infty } c_n x^{n+\lambda } \hspace{1cm}\mbox{where}\; \lambda =\mbox{indicial}\;\mbox{root} \label{eq:3008}
\end{equation}
Above $\lambda$ is indicial root.  Plug (\ref{eq:3008})  into (\ref{eq:3002}). Then we get a three-term recurrence relation for the coefficients $c_n$:
\begin{equation}
c_{n+1}=A_n \;c_n +B_n \;c_{n-1} \hspace{1cm};n\geq 1 \label{eq:3009}
\end{equation}
where,
\begin{subequations}
\begin{eqnarray}
A_n &=&\frac{(n+\lambda )(n+\lambda -\beta +\gamma +\delta -1)-q}{(n+\lambda +1)(n+\lambda +\gamma )}\nonumber\\
&=& \frac{\left( n- \frac{-(\varphi +2\lambda )-\sqrt{\varphi ^2+ 4q}}{2}\right) \left( n- \frac{-(\varphi +2\lambda )+\sqrt{\varphi ^2+ 4q}}{2}\right)}{(n+\lambda +1)(n+\lambda +\gamma )} 
\label{eq:30010a}
\end{eqnarray}
and
\begin{equation}
\varphi = -\beta + \gamma +\delta -1 \nonumber
\end{equation}
\vspace{2mm}
\begin{equation}
B_n = \frac{\beta (n+\lambda +\alpha -1)}{(n+\lambda +1)(n+\lambda +\gamma )} \label{eq:30010b}
\end{equation}
\begin{equation}
c_1= A_0 \;c_0 \label{eq:30010c}
\end{equation}
\end{subequations}
We have two indicial roots which are $\lambda _1= 0$ and $\lambda _2= 1-\gamma $.
\subsection{Power series}
\subsubsection{Polynomial of type 2}
In the previous chapter, by applying 3TRF, I construct the power series expansion, its integral forms and the generating function for the CHP of type 1: I treat $\beta $, $\gamma $, $\delta $ and $q$ as free variables and $\alpha $ as a fixed value
In this chapter, by applying R3TRF, I will show how to derive the power series expansion, the integral representation and the generating function for the CHP of type 2: I treat $\alpha $, $\beta $, $\gamma $ and $\delta $ as free variables and $q$ as a fixed value.  

In chapter 1 the general expression of power series of $y(x)$ for polynomial of type 2 is defined by
\begin{eqnarray}
y(x) &=& \sum_{n=0}^{\infty } y_n(x)= y_0(x)+ y_1(x)+ y_2(x)+ y_3(x)+\cdots \nonumber\\
&=&  c_0 \Bigg\{ \sum_{i_0=0}^{\alpha _0} \left( \prod _{i_1=0}^{i_0-1}A_{i_1} \right) x^{i_0+\lambda }
+ \sum_{i_0=0}^{\alpha _0}\left\{ B_{i_0+1} \prod _{i_1=0}^{i_0-1}A_{i_1}  \sum_{i_2=i_0}^{\alpha _1} \left( \prod _{i_3=i_0}^{i_2-1}A_{i_3+2} \right)\right\} x^{i_2+2+\lambda }  \nonumber\\
&& + \sum_{N=2}^{\infty } \Bigg\{ \sum_{i_0=0}^{\alpha _0} \Bigg\{B_{i_0+1}\prod _{i_1=0}^{i_0-1} A_{i_1} 
\prod _{k=1}^{N-1} \Bigg( \sum_{i_{2k}= i_{2(k-1)}}^{\alpha _k} B_{i_{2k}+2k+1}\prod _{i_{2k+1}=i_{2(k-1)}}^{i_{2k}-1}A_{i_{2k+1}+2k}\Bigg)\nonumber\\
&& \times  \sum_{i_{2N} = i_{2(N-1)}}^{\alpha _N} \Bigg( \prod _{i_{2N+1}=i_{2(N-1)}}^{i_{2N}-1} A_{i_{2N+1}+2N} \Bigg) \Bigg\} \Bigg\} x^{i_{2N}+2N+\lambda }\Bigg\}  \label{eq:30014}
\end{eqnarray}
In the above, $\alpha _i\leq \alpha _j$ only if $i\leq j$ where $i,j,\alpha _i, \alpha _j \in \mathbb{N}_{0}$.

For a polynomial, we need a condition which is:
\begin{equation}
A_{\alpha _i+ 2i}=0 \hspace{1cm} \mathrm{where}\;i,\alpha _i =0,1,2,\cdots
\label{eq:30015}
\end{equation}
In the above, $ \alpha _i$ is an eigenvalue that makes $A_n$ term terminated at certain value of n. (\ref{eq:30015}) makes each $y_i(x)$ where $i=0,1,2,\cdots$ as the polynomial in (\ref{eq:30014}).
Replace $\alpha _i$ by $q_i$ and put $n=q_i+ 2i$ in (\ref{eq:30010a}) with the condition $A_{q_i+ 2i}=0$. Then, we obtain eigenvalues $q$ such as $\sqrt{\varphi ^2+ 4q}= \pm\{\varphi +2\lambda + 2(q_i+2i)\} $.

\paragraph{The case of $ \sqrt{\varphi ^2+ 4q}= -\{\varphi +2\lambda + 2(q_i+2i)\}$}
In (\ref{eq:30010a}) replace $ \sqrt{\varphi ^2+ 4q}$ by ${ \displaystyle -\{\varphi +2\lambda + 2(q_i+2i)\}}$. In (\ref{eq:30014}) replace index $\alpha _i$ by $q_i$. Take the new (\ref{eq:30010a}) and (\ref{eq:30010b}) in new (\ref{eq:30014}).
After the replacement process, the general expression of power series of the CHE for polynomial of type 2 is given by
\begin{eqnarray}
 y(x)&=& \sum_{n=0}^{\infty } y_n(x)= y_0(x)+ y_1(x)+ y_2(x)+ y_3(x)+\cdots \nonumber\\
&=&  c_0 x^{\lambda } \left\{\sum_{i_0=0}^{q_0} \frac{(-q_0)_{i_0} \left( q_0+ \varphi +2\lambda \right)_{i_0}}{(1+\lambda )_{i_0}(\gamma +\lambda )_{i_0}} x^{i_0}\right.\nonumber\\
&&+ \left\{ \sum_{i_0=0}^{q_0}\frac{(i_0+\alpha + \lambda )}{ (i_0+ 2+\lambda )(i_0+ 1+\gamma + \lambda )}\frac{(-q_0)_{i_0} \left( q_0+ \varphi +2\lambda \right)_{i_0}}{(1+\lambda )_{i_0}(\gamma +\lambda )_{i_0}} \right.\nonumber\\
&&\times  \left. \sum_{i_1=i_0}^{q_1} \frac{(-q_1)_{i_1}\left(q_1+\varphi +4+2\lambda \right)_{i_1}(3+\lambda )_{i_0}(2+\gamma +\lambda )_{i_0}}{(-q_1)_{i_0}\left(q_1+\varphi +4+2\lambda \right)_{i_0}(3+\lambda )_{i_1}(2+\gamma +\lambda )_{i_1}} x^{i_1}\right\} z\nonumber\\
&&+ \sum_{n=2}^{\infty } \left\{ \sum_{i_0=0}^{q_0}\frac{(i_0+\alpha + \lambda )}{ (i_0+ 2+\lambda )(i_0+ 1+\gamma + \lambda )}\frac{(-q_0)_{i_0} \left( q_0+ \varphi +2\lambda \right)_{i_0}}{(1+\lambda )_{i_0}(\gamma +\lambda )_{i_0}}\right.\nonumber\\
&&\times \prod _{k=1}^{n-1} \left\{ \sum_{i_k=i_{k-1}}^{q_k} \frac{(i_k+ 2k+\alpha +\lambda )}{(i_k+ 2(k+1)+\lambda )(i_k+ 2k+1+\gamma +\lambda )}\right. \label{eq:30016}\\
&&\times \left.\frac{(-q_k)_{i_k}\left( q_k+4k +\varphi +2\lambda \right)_{i_k}(2k+1+\lambda )_{i_{k-1}}(2k+\gamma +\lambda )_{i_{k-1}}}{(-q_k)_{i_{k-1}}\left( q_k+4k +\varphi +2\lambda \right)_{i_{k-1}}(2k+1+\lambda )_{i_k}(2k+\gamma +\lambda )_{i_k}}\right\} \nonumber\\
&&\times \left. \left.\sum_{i_n= i_{n-1}}^{q_n} \frac{(-q_n)_{i_n}\left( q_n+4n +\varphi +2\lambda \right)_{i_n}(2n+1+\lambda )_{i_{n-1}}(2n+\gamma +\lambda )_{i_{n-1}}}{(-q_n)_{i_{n-1}}\left( q_n+4n +\varphi +2\lambda \right)_{i_{n-1}}(2n+1+\lambda )_{i_n}(2n+\gamma +\lambda )_{i_n}} x^{i_n} \right\} z^n \right\} \nonumber
\end{eqnarray}
where
\begin{equation}
\begin{cases} z = \beta x^2 \cr
\varphi = -\beta +\gamma +\delta -1 \cr
q= (q_j+2j+\lambda )(q_j+2j+\lambda +\varphi ) \;\;\mbox{as}\;j,q_j\in \mathbb{N}_{0} \cr
q_i\leq q_j \;\;\mbox{only}\;\mbox{if}\;i\leq j\;\;\mbox{where}\;i,j\in \mathbb{N}_{0} 
\end{cases}\nonumber 
\end{equation}
\paragraph{The case of  $ \sqrt{\varphi ^2+ 4q}= \varphi +2\lambda + 2(q_i+2i)$}
In (\ref{eq:30010a}) replace $ \sqrt{\varphi ^2+ 4q}$ by ${ \displaystyle \varphi +2\lambda + 2(q_i+2i)}$. In (\ref{eq:30014}) replace index $\alpha _i$ by $q_i$. Take the new (\ref{eq:30010a}) and (\ref{eq:30010b}) in new (\ref{eq:30014}).
After the replacement process, its solution is equivalent to (\ref{eq:30016}).

Put $c_0$= 1 as $\lambda =0$  for the first kind of independent solutions of the CHE and $\lambda =\frac{1}{2}$ for the second one in (\ref{eq:30016}).  
\begin{remark}
The power series expansion of the CHE of the first kind for polynomial of type 2 about $x=0$ as $q= (q_j+2j)(q_j+2j-\beta +\gamma +\delta -1) $ where $j,q_j \in \mathbb{N}_{0}$ is
\begin{eqnarray}
 y(x)&=& H_c^{(a)}F_{q_j}^R\left(\alpha, \beta, \gamma, \delta, q= (q_j+2j)(q_j+2j+\varphi ); \varphi = -\beta +\gamma +\delta -1; z = \beta x^2  \right)\nonumber\\
&=& \sum_{i_0=0}^{q_0} \frac{(-q_0)_{i_0} \left( q_0+ \varphi\right)_{i_0}}{(1)_{i_0}(\gamma)_{i_0}} x^{i_0}\nonumber\\
&&+ \left\{ \sum_{i_0=0}^{q_0}\frac{(i_0+\alpha )}{ (i_0+ 2)(i_0+ 1+\gamma )}\frac{(-q_0)_{i_0} \left( q_0+ \varphi \right)_{i_0}}{(1)_{i_0}(\gamma  )_{i_0}} \right. \nonumber\\
&&\times \left. \sum_{i_1=i_0}^{q_1} \frac{(-q_1)_{i_1}\left(q_1+\varphi +4 \right)_{i_1}(3)_{i_0}(2+\gamma)_{i_0}}{(-q_1)_{i_0}\left(q_1+\varphi +4 \right)_{i_0}(3+\lambda )_{i_1}(2+\gamma)_{i_1}} x^{i_1}\right\} z\nonumber\\
&&+ \sum_{n=2}^{\infty } \left\{ \sum_{i_0=0}^{q_0}\frac{(i_0+\alpha)}{ (i_0+ 2)(i_0+ 1+\gamma)}\frac{(-q_0)_{i_0} \left( q_0+ \varphi \right)_{i_0}}{(1)_{i_0}(\gamma)_{i_0}}\right.\nonumber\\
&&\times \prod _{k=1}^{n-1} \left\{ \sum_{i_k=i_{k-1}}^{q_k} \frac{(i_k+ 2k+\alpha)}{(i_k+ 2(k+1))(i_k+ 2k+1+\gamma)}\right.\nonumber\\
&&\times \left.\frac{(-q_k)_{i_k}\left( q_k+4k +\varphi \right)_{i_k}(2k+1)_{i_{k-1}}(2k+\gamma )_{i_{k-1}}}{(-q_k)_{i_{k-1}}\left( q_k+4k +\varphi \right)_{i_{k-1}}(2k+1)_{i_k}(2k+\gamma)_{i_k}}\right\} \nonumber\\
&&\times \left.\sum_{i_n= i_{n-1}}^{q_n} \frac{(-q_n)_{i_n}\left( q_n+4n +\varphi \right)_{i_n}(2n+1)_{i_{n-1}}(2n+\gamma)_{i_{n-1}}}{(-q_n)_{i_{n-1}}\left( q_n+4n +\varphi \right)_{i_{n-1}}(2n+1)_{i_n}(2n+\gamma)_{i_n}} x^{i_n} \right\} z^n \label{eq:30017}
\end{eqnarray}
\end{remark}
For the minimum value of the CHE of the first kind for polynomial of type 2 about $x=0$, put $q_0=q_1=q_2=\cdots=0$ in (\ref{eq:30017}).
\begin{eqnarray}
y(x)&=& H_c^{(a)}F_{0}^R\left(\alpha, \beta, \gamma, \delta, q= 2j(2j+\varphi ); \varphi = -\beta +\gamma +\delta -1; z = \beta x^2  \right)\nonumber\\
&=& \; _1F_1\left( \frac{\alpha }{2}, \frac{\gamma }{2}+\frac{1}{2}, \frac{1}{2}\beta x^2 \right) \;\;\mbox{where}\;-\infty < x< \infty \nonumber 
\end{eqnarray}
On the above,  $_1F_1( a,b,x)= \sum_{n=0}^{\infty }\frac{(a)_n}{(b)_n}\frac{x^n}{n!}$.
\begin{remark}
The power series expansion of the CHE of the first kind for polynomial of type 2 about $x=0$ as $q= (q_j+2j+1-\gamma )(q_j+2j -\beta +\delta ) $ where $j,q_j \in \mathbb{N}_{0}$ is
\begin{eqnarray}
 y(x)&=& H_c^{(a)}S_{\alpha _j}^R\left(\alpha, \beta, \gamma, \delta, q= (q_j+2j+1-\gamma )(q_j+2j+1-\gamma +\varphi ); \varphi = -\beta +\gamma +\delta -1; z= \beta x^2  \right)\nonumber\\
&=&  x^{1-\gamma } \left\{\sum_{i_0=0}^{q_0} \frac{(-q_0)_{i_0} \left( q_0+ \varphi +2(1-\gamma ) \right)_{i_0}}{(2-\gamma )_{i_0}(1)_{i_0}} x^{i_0}\right.\nonumber\\
&&+ \left\{ \sum_{i_0=0}^{q_0}\frac{(i_0+\alpha +1-\gamma )}{ (i_0+ 3-\gamma )(i_0+2)}\frac{(-q_0)_{i_0} \left( q_0+ \varphi +2(1-\gamma ) \right)_{i_0}}{(2-\gamma )_{i_0}(1)_{i_0}} \right.\nonumber\\
&&\times \left.\sum_{i_1=i_0}^{q_1} \frac{(-q_1)_{i_1}\left(q_1+\varphi +2(3-\gamma ) \right)_{i_1}(4-\gamma )_{i_0}(3)_{i_0}}{(-q_1)_{i_0}\left(q_1+\varphi +2(3-\gamma ) \right)_{i_0}(4-\gamma )_{i_1}(3)_{i_1}} x^{i_1}\right\} z\nonumber\\
&&+ \sum_{n=2}^{\infty } \left\{ \sum_{i_0=0}^{q_0}\frac{(i_0+\alpha +1-\gamma )}{ (i_0+ 3-\gamma )(i_0+2)}\frac{(-q_0)_{i_0} \left( q_0+ \varphi +2(1-\gamma ) \right)_{i_0}}{(2-\gamma )_{i_0}(1)_{i_0}}\right.\nonumber\\
&&\times \prod _{k=1}^{n-1} \left\{ \sum_{i_k=i_{k-1}}^{q_k} \frac{(i_k+ 2k+\alpha +1-\gamma )}{(i_k+ 2k+3-\gamma )(i_k+ 2k+2 )}\right. \nonumber\\
&&\times \left. \frac{(-q_k)_{i_k}\left( q_k +\varphi +2(2k+1-\gamma ) \right)_{i_k}(2k+2-\gamma )_{i_{k-1}}(2k+1)_{i_{k-1}}}{(-q_k)_{i_{k-1}}\left( q_k +\varphi +2(2k+1-\gamma ) \right)_{i_{k-1}}(2k+2-\gamma )_{i_k}(2k+1)_{i_k}}\right\} \nonumber\\
&&\times \left. \left.\sum_{i_n= i_{n-1}}^{q_n} \frac{(-q_n)_{i_n}\left( q_n +\varphi +2(2n+1-\gamma ) \right)_{i_n}(2n+2-\gamma  )_{i_{n-1}}(2n+1)_{i_{n-1}}}{(-q_n)_{i_{n-1}}\left( q_n +\varphi +2(2n+1-\gamma ) \right)_{i_{n-1}}(2n+2-\gamma)_{i_n}(2n+1)_{i_n}} x^{i_n} \right\} z^n \right\}\label{eq:30018}
\end{eqnarray}
\end{remark}
For the minimum value of the CHE of the second kind for polynomial of type 2 about $x=0$, put $q_0=q_1=q_2=\cdots=0$ in 
(\ref{eq:30018}). 
\begin{eqnarray}
y(x)&=& H_c^{(a)}S_{0}^R\left(\alpha, \beta, \gamma, \delta, q= (2j+1-\gamma )(2j+1-\gamma +\varphi ); \varphi = -\beta +\gamma +\delta -1; z= \beta x^2  \right)\nonumber\\
&=& x^{1-\gamma }\; _1F_1\left( \frac{\alpha }{2}-\frac{\gamma }{2}+\frac{1}{2}, -\frac{\gamma }{2}+\frac{3}{2}, \frac{1}{2}\beta x^2 \right) \;\;\mbox{where}\;-\infty < x< \infty \nonumber 
\end{eqnarray}
In chapter 4 I treat $\alpha $ as a fixed value and $\beta$, $\gamma , \delta $, $q$ as free variables to construct the CHP of type 1: (1) if $\alpha = -(2 \alpha _j+j) $ where $j, \alpha_j \in \mathbb{N}_{0}$, an analytic solution of the CHE turns to be the first kind of independent solution of the CHP of type 1. (2) if  $\alpha = -(2 \alpha _j+j +1-\gamma ) $, an analytic solution of the CHE turns to be the second kind of independent solution of the CHP of type 1. 

In this chapter I treat $q$ as a fixed value and $\alpha, \beta, \gamma, \delta $ as free variables to construct the CHP of type 2: (1) if $q= (q_j+2j)(q_j+2j-\beta +\gamma +\delta -1) $ where $j,q_j \in \mathbb{N}_{0}$, an analytic solution of the CHE turns to be the first kind of independent solution of the CHP of type 2. (2) if $q= (q_j+2j+1-\gamma )(q_j+2j -\beta +\delta ) $, an analytic solution of the CHE turns to be the second kind of independent solution of the CHP of type 2.
\subsubsection{Infinite series}
In chapter 1 the general expression of power series of $y(x)$ for infinite series is defined by
\begin{eqnarray}
y(x) &=& \sum_{n=0}^{\infty } y_n(x)= y_0(x)+ y_1(x)+ y_2(x)+ y_3(x)+\cdots \nonumber\\
&=& c_0 \Bigg\{ \sum_{i_0=0}^{\infty } \left( \prod _{i_1=0}^{i_0-1}A_{i_1} \right) x^{i_0+\lambda }
+ \sum_{i_0=0}^{\infty }\left\{ B_{i_0+1} \prod _{i_1=0}^{i_0-1}A_{i_1}  \sum_{i_2=i_0}^{\infty } \left( \prod _{i_3=i_0}^{i_2-1}A_{i_3+2} \right)\right\} x^{i_2+2+\lambda }  \nonumber\\
&& + \sum_{N=2}^{\infty } \Bigg\{ \sum_{i_0=0}^{\infty } \Bigg\{B_{i_0+1}\prod _{i_1=0}^{i_0-1} A_{i_1} 
\prod _{k=1}^{N-1} \Bigg( \sum_{i_{2k}= i_{2(k-1)}}^{\infty } B_{i_{2k}+2k+1}\prod _{i_{2k+1}=i_{2(k-1)}}^{i_{2k}-1}A_{i_{2k+1}+2k}\Bigg)\nonumber\\
&& \times  \sum_{i_{2N} = i_{2(N-1)}}^{\infty } \Bigg( \prod _{i_{2N+1}=i_{2(N-1)}}^{i_{2N}-1} A_{i_{2N+1}+2N} \Bigg) \Bigg\} \Bigg\} x^{i_{2N}+2N+\lambda }\Bigg\}   \label{eq:30019}
\end{eqnarray}
Substitute (\ref{eq:30010a})--(\ref{eq:30010c}) into (\ref{eq:30019}). 
The general expression of power series of the CHE for infinite series about $x=0$ is given by
\begin{eqnarray}
 y(x)&=&\sum_{n=0}^{\infty } y_n(x)= y_0(x)+ y_1(x)+ y_2(x)+ y_3(x)+\cdots \nonumber\\
&=& c_0 x^{\lambda } \left\{\sum_{i_0=0}^{\infty } \frac{\left(\Delta_0^{-}\right)_{i_0} \left(\Delta_0^{+}\right)_{i_0}}{(1+\lambda )_{i_0}(\gamma +\lambda )_{i_0}} x^{i_0}\right.\nonumber\\
&&+ \left\{ \sum_{i_0=0}^{\infty }\frac{(i_0+ \lambda +\alpha ) }{(i_0+ \lambda +2)(i_0+ \lambda +1+\gamma )}\frac{\left(\Delta_0^{-}\right)_{i_0} \left(\Delta_0^{+}\right)_{i_0}}{(1+\lambda )_{i_0}(\gamma +\lambda )_{i_0}} \right. \nonumber\\
&&\times \left. \sum_{i_1=i_0}^{\infty } \frac{\left(\Delta_1^{-}\right)_{i_1} \left(\Delta_1^{+}\right)_{i_1}(3+\lambda )_{i_0}(2+\gamma +\lambda )_{i_0}}{\left(\Delta_1^{-}\right)_{i_0}  \left(\Delta_1^{+}\right)_{i_0}(3+\lambda )_{i_1}(2+\gamma +\lambda )_{i_1}}x^{i_1}\right\} z\nonumber\\
&&+ \sum_{n=2}^{\infty } \left\{ \sum_{i_0=0}^{\infty } \frac{(i_0+ \lambda +\alpha ) }{(i_0+ \lambda +2)(i_0+ \lambda +1+\gamma )}\frac{\left(\Delta_0^{-}\right)_{i_0} \left(\Delta_0^{+}\right)_{i_0}}{(1+\lambda )_{i_0}(\gamma +\lambda )_{i_0}}\right.\nonumber\\
&&\times \prod _{k=1}^{n-1} \left\{ \sum_{i_k=i_{k-1}}^{\infty } \frac{(i_k+ 2k+\lambda +\alpha ) }{(i_k+ 2(k+1)+\lambda )(i_k+ 2k+1+\gamma +\lambda )}\right.\nonumber\\
&&\times \left. \frac{ \left(\Delta_k^{-}\right)_{i_k} \left(\Delta_k^{+} \right)_{i_k}(2k+1+\lambda )_{i_{k-1}}(2k+\gamma +\lambda )_{i_{k-1}}}{\left(\Delta_k^{-}\right)_{i_{k-1}} \left(\Delta_k^{+} \right)_{i_{k-1}}(2k+1+\lambda )_{i_k}(2k+\gamma +\lambda )_{i_k}}\right\}\nonumber\\
&&\times \left.\left.\sum_{i_n= i_{n-1}}^{\infty } \frac{\left(\Delta_n^{-}\right)_{i_n}\left( \Delta_n^{+} \right)_{i_n}(2n+1+\lambda )_{i_{n-1}}(2n+\gamma +\lambda )_{i_{n-1}}}{\left(\Delta_n^{-}\right)_{i_{n-1}}\left(\Delta_n^{+} \right)_{i_{n-1}} (2n+1+\lambda )_{i_n}(2n+\gamma +\lambda )_{i_n}} x^{i_n} \right\} z^n \right\} \hspace{1cm}\label{eq:30020}
\end{eqnarray}
where
\begin{equation}
\begin{cases} 
\Delta_0^{\pm}= \frac{ \varphi +2 \lambda  \pm\sqrt{\varphi ^2+4q}}{2} \cr
\Delta_1^{\pm}=  \frac{ \varphi +4+ 2 \lambda  \pm\sqrt{\varphi ^2+4q}}{2} \cr
\Delta_k^{\pm}=   \frac{ \varphi + 2 \lambda +4k \pm\sqrt{\varphi ^2+4q}}{2} \cr
\Delta_n^{\pm}=   \frac{ \varphi + 2 \lambda +4n \pm\sqrt{\varphi ^2+4q}}{2}
\end{cases}\nonumber 
\end{equation}
Put $c_0$= 1 as $\lambda =0$  for the first kind of independent solutions of the CHE and $\lambda =1-\gamma$ for the second one in (\ref{eq:30020}). 
\begin{remark}
The power series expansion of the CHE of the first kind for infinite series about $x=0$ using R3TRF is
\begin{eqnarray}
 y(x)&=& H_c^{(a)}F^R\left(\alpha, \beta, \gamma, \delta, q; \varphi = -\beta +\gamma +\delta -1; z= \beta x^2 \right) \nonumber\\
&=&  \left\{\sum_{i_0=0}^{\infty } \frac{\left(\Delta_0^{-}\right)_{i_0} \left(\Delta_0^{+}\right)_{i_0}}{(1)_{i_0}(\gamma )_{i_0}} x^{i_0}\right.\nonumber\\
&+& \left\{ \sum_{i_0=0}^{\infty }\frac{(i_0 +\alpha ) }{(i_0 +2)(i_0 +1+\gamma )}\frac{\left(\Delta_0^{-}\right)_{i_0} \left(\Delta_0^{+}\right)_{i_0}}{(1)_{i_0}(\gamma )_{i_0}} \sum_{i_1=i_0}^{\infty } \frac{\left(\Delta_1^{-}\right)_{i_1} \left(\Delta_1^{+}\right)_{i_1}(3)_{i_0}(2+\gamma)_{i_0}}{\left(\Delta_1^{-}\right)_{i_0}  \left(\Delta_1^{+}\right)_{i_0}(3)_{i_1}(2+\gamma)_{i_1}}x^{i_1}\right\} z\nonumber\\
&+& \sum_{n=2}^{\infty } \left\{ \sum_{i_0=0}^{\infty } \frac{(i_0 +\alpha ) }{(i_0 +2)(i_0 +1+\gamma )}\frac{\left(\Delta_0^{-}\right)_{i_0} \left(\Delta_0^{+}\right)_{i_0}}{(1)_{i_0}(\gamma)_{i_0}}\right.\nonumber\\
&\times& \prod _{k=1}^{n-1} \left\{ \sum_{i_k=i_{k-1}}^{\infty } \frac{(i_k+ 2k+\alpha ) }{(i_k+ 2(k+1))(i_k+ 2k+1+\gamma )} \frac{ \left(\Delta_k^{-}\right)_{i_k} \left(\Delta_k^{+} \right)_{i_k}(2k+1)_{i_{k-1}}(2k+\gamma )_{i_{k-1}}}{\left(\Delta_k^{-}\right)_{i_{k-1}} \left(\Delta_k^{+} \right)_{i_{k-1}}(2k+1)_{i_k}(2k+\gamma )_{i_k}}\right\}\nonumber\\
&\times& \left.\sum_{i_n= i_{n-1}}^{\infty } \frac{\left(\Delta_n^{-}\right)_{i_n}\left( \Delta_n^{+} \right)_{i_n}(2n+1)_{i_{n-1}}(2n+\gamma)_{i_{n-1}}}{\left(\Delta_n^{-}\right)_{i_{n-1}}\left(\Delta_n^{+} \right)_{i_{n-1}} (2n+1)_{i_n}(2n+\gamma)_{i_n}} x^{i_n} \right\} z^n  \label{eq:30021}
\end{eqnarray}
where
\begin{equation}
\begin{cases} 
\Delta_0^{\pm}= \frac{ \varphi \pm\sqrt{\varphi ^2+4q}}{2} \cr
\Delta_1^{\pm}=  \frac{ \varphi +4 \pm\sqrt{\varphi ^2+4q}}{2} \cr
\Delta_k^{\pm}=   \frac{ \varphi +4k \pm\sqrt{\varphi ^2+4q}}{2} \cr
\Delta_n^{\pm}=   \frac{ \varphi +4n \pm\sqrt{\varphi ^2+4q}}{2}
\end{cases}\nonumber 
\end{equation}
\end{remark}
\begin{remark}
The power series expansion of the CHE of the second kind for infinite series about $x=0$ using R3TRF is
\begin{eqnarray}
 y(x)&=&  H_c^{(a)}S^R\left(\alpha, \beta, \gamma, \delta, q; \varphi = -\beta +\gamma +\delta -1; z= \beta x^2   \right)\nonumber\\
&=& x^{1-\gamma  }\left\{\sum_{i_0=0}^{\infty } \frac{\left(\Delta_0^{-}\right)_{i_0} \left(\Delta_0^{+}\right)_{i_0}}{(2-\gamma)_{i_0}(1)_{i_0}} x^{i_0}\right.\nonumber\\
&+& \left\{ \sum_{i_0=0}^{\infty }\frac{(i_0+1-\gamma +\alpha ) }{(i_0+3-\gamma )(i_0+2)}\frac{\left(\Delta_0^{-}\right)_{i_0} \left(\Delta_0^{+}\right)_{i_0}}{(2-\gamma )_{i_0}(1)_{i_0}} \sum_{i_1=i_0}^{\infty } \frac{\left(\Delta_1^{-}\right)_{i_1} \left(\Delta_1^{+}\right)_{i_1}(4-\gamma )_{i_0}(3)_{i_0}}{\left(\Delta_1^{-}\right)_{i_0}  \left(\Delta_1^{+}\right)_{i_0}(4-\gamma )_{i_1}(3)_{i_1}}x^{i_1}\right\} z\nonumber\\
&+& \sum_{n=2}^{\infty } \left\{ \sum_{i_0=0}^{\infty } \frac{(i_0+1-\gamma +\alpha ) }{(i_0+3-\gamma )(i_0+2)}\frac{\left(\Delta_0^{-}\right)_{i_0} \left(\Delta_0^{+}\right)_{i_0}}{(2-\gamma )_{i_0}(1)_{i_0}}\right.\nonumber\\
&\times& \prod _{k=1}^{n-1} \left\{ \sum_{i_k=i_{k-1}}^{\infty } \frac{(i_k+ 2k+1-\gamma +\alpha ) }{(i_k+ 2k+3-\gamma  )(i_k+ 2k+2)} \frac{ \left(\Delta_k^{-}\right)_{i_k} \left(\Delta_k^{+} \right)_{i_k}(2k+2-\gamma )_{i_{k-1}}(2k+1)_{i_{k-1}}}{\left(\Delta_k^{-}\right)_{i_{k-1}} \left(\Delta_k^{+} \right)_{i_{k-1}}(2k+2-\gamma )_{i_k}(2k+1)_{i_k}}\right\}\nonumber\\
&\times& \left.\left.\sum_{i_n= i_{n-1}}^{\infty } \frac{\left(\Delta_n^{-}\right)_{i_n}\left( \Delta_n^{+} \right)_{i_n}(2n+2-\gamma )_{i_{n-1}}(2n+1)_{i_{n-1}}}{\left(\Delta_n^{-}\right)_{i_{n-1}}\left(\Delta_n^{+} \right)_{i_{n-1}} (2n+2-\gamma )_{i_n}(2n+1)_{i_n}} x^{i_n} \right\} z^n \right\} \label{eq:30022}
\end{eqnarray}
where
\begin{equation}
\begin{cases} 
\Delta_0^{\pm}= \frac{ \varphi +2(1-\gamma )  \pm\sqrt{\varphi ^2+4q}}{2} \cr
\Delta_1^{\pm}=  \frac{ \varphi + 2(3-\gamma )  \pm\sqrt{\varphi ^2+4q}}{2} \cr
\Delta_k^{\pm}=   \frac{ \varphi + 2(2k+1-\gamma ) \pm\sqrt{\varphi ^2+4q}}{2} \cr
\Delta_n^{\pm}=   \frac{ \varphi + 2(2n+1-\gamma ) \pm\sqrt{\varphi ^2+4q}}{2}
\end{cases}\nonumber 
\end{equation}
\end{remark}
It is required that $\gamma \ne 0,-1,-2,\cdots$ for the first kind of independent solution of the CHE for all cases. Because if it does not, its solution will be divergent. And it is required that $\gamma \ne 2,3,4,\cdots$ for the second kind of independent solution of the CHE for all cases.

The infinite series in this chapter are equivalent to the infinite series in chapter 4. In this chapter $B_n$ is the leading term in sequence $c_n$ of the analytic function $y(x)$. In chapter 4 $A_n$ is the leading term in sequence $c_n$ of the analytic function $y(x)$.

For the special case, as $|\beta| \ll 1$ in (\ref{eq:30021}) and (\ref{eq:30022}), we have
\begin{subequations}
\begin{eqnarray}
&&\lim_{|\beta| \ll 1}  H_c^{(a)}F^R\left(\alpha, \beta, \gamma, \delta, q; \varphi = -\beta +\gamma +\delta -1; z= \beta x^2 \right) \nonumber\\
&&\approx \; _2F_1 \left(\frac{ \gamma +\delta -1 -\sqrt{(\gamma +\delta -1)^2+4q}}{2},  \frac{ \gamma +\delta -1 +\sqrt{(\gamma +\delta -1)^2+4q}}{2}; \gamma ; x \right)\hspace{2cm}\label{eq:30023a}
\end{eqnarray}
And,
\begin{eqnarray}
&&\lim_{|\beta| \ll 1} H_c^{(a)}S^R\left(\alpha, \beta, \gamma, \delta, q; \varphi = -\beta +\gamma +\delta -1; z= \beta x^2 \right)\label{eq:30023b}\\
&&\approx  x^{1-\gamma } \;_2F_1 \left(\frac{ -\gamma +\delta +1 -\sqrt{(\gamma +\delta -1)^2+4q}}{2}, \frac{ -\gamma +\delta +1 +\sqrt{(\gamma +\delta -1)^2+4q}}{2}; 2-\gamma; x \right) \nonumber
\end{eqnarray}
\end{subequations}
(\ref{eq:30023a}) and (\ref{eq:30023b}) are Gauss Hypergeometric function.
\subsection{Integral formalism}
In earlier literature the integral relations of the CHE were constructed by using Fredholm integral equations; such integral relationships express one analytic solution in terms of another analytic solution such as a confluent hypergeometric function with a branch point at zero. \cite{Ronv1995} There are many other forms of integral relations in the CHE. \cite{Lamb1934,Erde1942,Slee1969,Abra1976,Schm1979,Kaza1986}

In chapter 4 I derive the combined definite and contour integral forms of the CHF by applying 3TRF: a $_1F_1$ function recurs in each of sub-integral forms of the CHF.
Now I consider the combined definite and contour integral representations of the CHF by using R3TRF.
\subsubsection{Polynomial of type 2}

There is a generalized hypergeometric function which is given by
\begin{eqnarray}
I_l &=& \sum_{i_l= i_{l-1}}^{q_l} \frac{(-q_l)_{i_l}\left( q_l+4l+ \varphi +2\lambda \right)_{i_l}(2l+1+\lambda )_{i_{l-1}}(2l+\gamma +\lambda )_{i_{l-1}}}{(-q_l)_{i_{l-1}}\left( q_l+4l+ \varphi +2\lambda \right)_{i_{l-1}}(2l+1+\lambda )_{i_l}(2l+\gamma +\lambda )_{i_l}} x^{i_l}\nonumber\\
&=& x^{i_{l-1}} \sum_{j=0}^{\infty } \frac{B(i_{l-1}+2l+\lambda ,j+1) B(i_{l-1}+2l-1+\gamma +\lambda ,j+1)}{(i_{l-1}+2l+\lambda )^{-1}(i_{l-1}+2l-1+\gamma +\lambda )^{-1}} \nonumber\\
&&\times \frac{(i_{l-1}-q_l)_j \left( i_{l-1}+q_l+4l+ \varphi +2\lambda \right)_j}{(1)_j \;j!} x^j \label{eq:30024}
\end{eqnarray}
By using integral form of beta function,
\begin{subequations}
\begin{equation}
B\left(i_{l-1}+2l+\lambda ,j+1\right)= \int_{0}^{1} dt_l\;t_l^{i_{l-1}+2l-1+\lambda } (1-t_l)^j \label{eq:30025a}
\end{equation}
\begin{equation}
B\left(i_{l-1}+2l-1+\gamma +\lambda ,j+1\right)= \int_{0}^{1} du_l\;u_l^{i_{l-1}+2(l-1)+\gamma +\lambda } (1-u_l)^j\label{eq:30025b}
\end{equation}
\end{subequations}
Substitute (\ref{eq:30025a}) and (\ref{eq:30025b}) into (\ref{eq:30024}). And divide $(i_{l-1}+2l+\lambda )(i_{l-1}+2l-1+\gamma +\lambda )$ into $I_l$.
\begin{eqnarray}
K_l&=& \frac{1}{(i_{l-1}+2l+\lambda )(i_{l-1}+2l-1+\gamma +\lambda )}\nonumber\\
&&\times \sum_{i_l= i_{l-1}}^{q_l} \frac{(-q_l)_{i_l}\left( q_l+4l+ \varphi +2\lambda \right)_{i_l}(2l+1+\lambda )_{i_{l-1}}(2l+\gamma +\lambda )_{i_{l-1}}}{(-q_l)_{i_{l-1}}\left( q_l+4l+ \varphi +2\lambda \right)_{i_{l-1}}(2l+1+\lambda )_{i_l}(2l+\gamma +\lambda )_{i_l}} x^{i_l}\nonumber\\
&=&  \int_{0}^{1} dt_l\;t_l^{2l-1+\lambda } \int_{0}^{1} du_l\;u_l^{2(l-1)+\gamma +\lambda } (x t_l u_l)^{i_{l-1}}\nonumber\\
&&\times \sum_{j=0}^{\infty } \frac{(i_{l-1}-q_l)_j \left( i_{l-1}+q_l+4l+\varphi +2\lambda \right)_j}{(1)_j \;j!} 
[x(1-t_l)(1-u_l)]^j \nonumber 
\end{eqnarray}
The integral form of Gauss hypergeometric function is defined by
\begin{eqnarray}
_2F_1 \left( \alpha ,\beta ; \gamma ; z \right) &=& \sum_{j=0}^{\infty } \frac{(\alpha )_j (\beta )_j}{(\gamma )_j (j!)} z^j \nonumber\\
&=& -\frac{1}{2\pi i} \frac{\Gamma(1-\alpha ) \Gamma(\gamma )}{\Gamma (\gamma -\alpha )} \oint dv_l\;(-v_l)^{\alpha -1} (1-v_l)^{\gamma -\alpha -1} (1-zv_l)^{-\beta }\hspace{1cm}\label{eq:30026}\\
&& \mbox{where} \;\mbox{Re}(\gamma -\alpha )>0 \nonumber
\end{eqnarray}
replaced $\alpha $, $\beta $, $\gamma $ and z by $i_{l-1}-q_l$, $ { \displaystyle i_{l-1}+q_l+4l+ \varphi +2\lambda }$, 1 and $x (1-t_l)(1-u_l)$ in (\ref{eq:30026})
\begin{eqnarray}
&& \sum_{j=0}^{\infty } \frac{\left(i_{l-1}-q_l)_j (i_{l-1}+q_l+4l+ \varphi +2\lambda \right)_j}{(1)_j \;j!} [x(1-t_l)(1-u_l)]^j \nonumber\\
&=& \frac{1}{2\pi i} \oint dv_l\;\frac{1}{v_l} \left(\frac{v_l-1}{v_l}\frac{1}{1-x(1-t_l)(1-u_l)v_l}\right)^{q_l} (1-x(1-t_l)(1-u_l)v_l)^{-\left( 4l+ \varphi +2\lambda \right)}\nonumber\\
&&\times \left(\frac{v_l}{v_l-1} \frac{1}{1-x(1-t_l)(1-u_l)v_l}\right)^{i_{l-1}} \label{eq:30027}
\end{eqnarray}
Substitute (\ref{eq:30027}) into $K_l$.
\begin{eqnarray}
K_l&=& \frac{1}{(i_{l-1}+2l+\lambda )(i_{l-1}+2l-1+\gamma +\lambda )}\nonumber\\
&&\times \sum_{i_l= i_{l-1}}^{q_l} \frac{(-q_l)_{i_l}\left( q_l+4l+ \varphi +2\lambda \right)_{i_l}(2l+1+\lambda )_{i_{l-1}}(2l+\gamma +\lambda )_{i_{l-1}}}{(-q_l)_{i_{l-1}}\left( q_l+4l+ \varphi +2\lambda \right)_{i_{l-1}}(2l+1+\lambda )_{i_l}(2l+\gamma +\lambda )_{i_l}} x^{i_l}\nonumber\\
&=&  \int_{0}^{1} dt_l\;t_l^{2l-1+\lambda } \int_{0}^{1} du_l\;u_l^{2(l-1)+\gamma +\lambda } 
\frac{1}{2\pi i} \oint dv_l\;\frac{1}{v_l} \left(\frac{v_l-1}{v_l}\frac{1}{1-x(1-t_l)(1-u_l)v_l}\right)^{q_l} \nonumber\\
&&\times (1-x(1-t_l)(1-u_l)v_l)^{-\left( 4l+ \varphi +2\lambda \right)} \left(\frac{v_l}{v_l-1} \frac{xt_l u_l}{1-x(1-t_l)(1-u_l)v_l}\right)^{i_{l-1}}\label{eq:30028} 
\end{eqnarray}
Substitute (\ref{eq:30028}) into (\ref{eq:30016}) where $l=1,2,3,\cdots$; apply $K_1$ into the second summation of sub-power series $y_1(x)$, apply $K_2$ into the third summation and $K_1$ into the second summation of sub-power series $y_2(x)$, apply $K_3$ into the forth summation, $K_2$ into the third summation and $K_1$ into the second summation of sub-power series $y_3(x)$, etc.\footnote{$y_1(x)$ means the sub-power series in (\ref{eq:30016}) contains one term of $B_n's$, $y_2(x)$ means the sub-power series in (\ref{eq:30016}) contains two terms of $B_n's$, $y_3(x)$ means the sub-power series in (\ref{eq:30016}) contains three terms of $B_n's$, etc.}
\begin{theorem}
The general representation in the form of integral of the CHP of type 2 is given by
\begin{eqnarray}
 y(x)&=& \sum_{n=0}^{\infty } y_{n}(x)= y_0(x)+ y_1(x)+ y_2(x)+y_3(x)+\cdots \nonumber\\
&=& c_0 x^{\lambda } \Bigg\{ \sum_{i_0=0}^{q_0}\frac{(-q_0)_{i_0}\left(q_0+ \varphi +2\lambda \right)_{i_0}}{(1+\lambda )_{i_0}(\gamma +\lambda )_{i_0}} x^{i_0}\nonumber\\
&&+ \sum_{n=1}^{\infty } \Bigg\{\prod _{k=0}^{n-1} \Bigg\{ \int_{0}^{1} dt_{n-k}\;t_{n-k}^{2(n-k)-1+\lambda } \int_{0}^{1} du_{n-k}\;u_{n-k}^{2(n-k-1)+\gamma +\lambda } \nonumber\\
&&\times  \frac{1}{2\pi i}  \oint dv_{n-k} \frac{1}{v_{n-k}} \left( \frac{v_{n-k}-1}{v_{n-k}} \frac{1}{1-\overleftrightarrow {w}_{n-k+1,n}(1-t_{n-k})(1-u_{n-k})v_{n-k}}\right)^{q_{n-k}} \nonumber\\
&&\times \left( 1- \overleftrightarrow {w}_{n-k+1,n}(1-t_{n-k})(1-u_{n-k})v_{n-k}\right)^{-\left(4(n-k)+ \varphi +2\lambda \right)}\nonumber\\
&&\times  \overleftrightarrow {w}_{n-k,n}^{-(2(n-k-1)+\alpha +\lambda )}\left(  \overleftrightarrow {w}_{n-k,n} \partial _{ \overleftrightarrow {w}_{n-k,n}}\right) \overleftrightarrow {w}_{n-k,n}^{2(n-k-1)+\alpha +\lambda } \Bigg\} \nonumber\\
&&\times \sum_{i_0=0}^{q_0}\frac{(-q_0)_{i_0}\left(q_0+ \varphi +2\lambda \right)_{i_0}}{(1+\lambda )_{i_0}(\gamma +\lambda )_{i_0}} \overleftrightarrow {w}_{1,n}^{i_0}\Bigg\} z^n \Bigg\} \label{eq:30029}
\end{eqnarray}
where
\begin{equation}\overleftrightarrow {w}_{i,j}=
\begin{cases} \displaystyle {\frac{v_i}{(v_i-1)}\; \frac{\overleftrightarrow w_{i+1,j} t_i u_i}{1- \overleftrightarrow w_{i+1,j} v_i (1-t_i)(1-u_i)}} \;\;\mbox{where}\; i\leq j\cr
x \;\;\mbox{only}\;\mbox{if}\; i>j
\end{cases}\nonumber 
\end{equation}
In the above, the first sub-integral form contains one term of $B_n's$, the second one contains two terms of $B_n$'s, the third one contains three terms of $B_n$'s, etc.
\end{theorem}
\begin{proof} 
According to (\ref{eq:30016}), 
\begin{equation}
 y(x)= \sum_{n=0}^{\infty }y_n(x) = y_0(x)+ y_1(x)+ y_2(x)+y_3(x)+\cdots \label{eq:30030}
\end{equation}
In the above, the power series expansion of sub-summation $y_0(x) $, $y_1(x)$, $y_2(x)$ and $y_3(x)$ of the type 2 CHP using R3TRF about $x=0$ are
\begin{subequations}
\begin{equation}
 y_0(x)= c_0 x^{\lambda } \sum_{i_0=0}^{q_0} \frac{(-q_0)_{i_0} \left(q_0+ \varphi +2\lambda \right)_{i_0}}{(1+\lambda )_{i_0}(\gamma +\lambda )_{i_0}} x^{i_0} \label{eq:30031a}
\end{equation}
\begin{eqnarray}
 y_1(x) &=& c_0 x^{\lambda } \Bigg\{\sum_{i_0=0}^{q_0}\frac{(i_0+ \lambda +\alpha ) }{(i_0+ \lambda +2)(i_0+ \lambda +1+\gamma )}\frac{(-q_0)_{i_0} \left(q_0+ \varphi +2\lambda \right)_{i_0}}{(1+\lambda )_{i_0}(\gamma +\lambda )_{i_0}} \nonumber\\
&&\times  \sum_{i_1=i_0}^{q_1} \frac{(-q_1)_{i_1}\left(q_1+4+ \varphi +2\lambda \right)_{i_1}(3+\lambda )_{i_0}(2+\gamma +\lambda )_{i_0}}{(-q_1)_{i_0}\left(q_1+4+ \varphi +2\lambda \right)_{i_0}(3+\lambda )_{i_1}(2+\gamma +\lambda )_{i_1}} x^{i_1} \Bigg\} z  \hspace{1cm}\label{eq:30031b}
\end{eqnarray}
\begin{eqnarray}
 y_2(x) &=& c_0 x^{\lambda } \Bigg\{\sum_{i_0=0}^{q_0}\frac{(i_0+ \lambda +\alpha ) }{(i_0+ \lambda +2)(i_0+ \lambda +1+\gamma )}\frac{(-q_0)_{i_0} \left(q_0+ \varphi +2\lambda \right)_{i_0}}{(1+\lambda )_{i_0}(\gamma +\lambda )_{i_0}} \nonumber\\
&&\times  \sum_{i_1=i_0}^{q_1} \frac{(i_1+2+ \lambda +\alpha ) }{(i_1+ \lambda +4)(i_1+ \lambda +3+\gamma )} \frac{(-q_1)_{i_1}\left(q_1+4+  \varphi +2\lambda \right)_{i_1}(3+\lambda )_{i_0}(2+\gamma +\lambda )_{i_0}}{(-q_1)_{i_0}\left(q_1+4+ \varphi +2\lambda \right)_{i_0}(3+\lambda )_{i_1}(2+\gamma +\lambda )_{i_1}} \nonumber\\
&&\times \sum_{i_2=i_1}^{q_2} \frac{(-q_2)_{i_2}\left(q_2+8+ \varphi +2\lambda \right)_{i_2}(5+\lambda )_{i_1}(4+\gamma +\lambda )_{i_1}}{(-q_2)_{i_1}\left(q_2+8+ \varphi +2\lambda \right)_{i_1}(5+\lambda )_{i_2}(4+\gamma +\lambda )_{i_2}} x^{i_2} \Bigg\} z^2  \label{eq:30031c}
\end{eqnarray}
\begin{eqnarray}
 y_3(x) &=&  c_0 x^{\lambda } \Bigg\{\sum_{i_0=0}^{q_0}\frac{(i_0+ \lambda +\alpha ) }{(i_0+ \lambda +2)(i_0+ \lambda +1+\gamma )}\frac{(-q_0)_{i_0} \left(q_0+ \varphi +2\lambda \right)_{i_0}}{(1+\lambda )_{i_0}(\gamma +\lambda )_{i_0}} \nonumber\\
&&\times  \sum_{i_1=i_0}^{q_1} \frac{(i_1+2+ \lambda +\alpha ) }{(i_1+ \lambda +4)(i_1+ \lambda +3+\gamma )} \frac{(-q_1)_{i_1}\left(q_1+4+  \varphi +2\lambda \right)_{i_1}(3+\lambda )_{i_0}(2+\gamma +\lambda )_{i_0}}{(-q_1)_{i_0}\left(q_1+4+ \varphi +2\lambda \right)_{i_0}(3+\lambda )_{i_1}(2+\gamma +\lambda )_{i_1}} \nonumber\\
&&\times \sum_{i_2=i_1}^{q_2} \frac{(i_2+4+ \lambda +\alpha ) }{(i_2+ \lambda +6)(i_2+ \lambda +5+\gamma )}  \frac{(-q_2)_{i_2}\left(q_2+8+  \varphi +2\lambda \right)_{i_2}(5+\lambda )_{i_1}(4+\gamma +\lambda )_{i_1}}{(-q_2)_{i_1}\left(q_2+8+ \varphi +2\lambda \right)_{i_1}(5+\lambda )_{i_2}(4+\gamma +\lambda )_{i_2}} \nonumber\\
&&\times \sum_{i_3=i_2}^{q_3} \frac{(-q_3)_{i_3}\left(q_2+12+ \varphi +2\lambda \right)_{i_3}(7+\lambda )_{i_2}(6+\gamma +\lambda )_{i_2}}{(-q_3)_{i_2}\left(q_2+12+ \varphi +2\lambda \right)_{i_2}(7+\lambda )_{i_3}(6+\gamma +\lambda )_{i_3}}x^{i_3} \Bigg\} z^3  \label{eq:30031d} 
\end{eqnarray}
\end{subequations}
Put $l=1$ in (\ref{eq:30028}). Take the new (\ref{eq:30028}) into (\ref{eq:30031b}).
\begin{eqnarray}
y_1(x) &=& \int_{0}^{1} dt_1\;t_1^{1+\lambda } \int_{0}^{1} du_1\;u_1^{\gamma +\lambda } \frac{1}{2\pi i} \oint dv_1 \;\frac{1}{v_1} 
\left( \frac{v_1-1}{v_1} \frac{1}{1-x(1-t_1)(1-u_1)v_1}\right)^{q_1}  \nonumber\\
&&\times (1-x(1-t_1)(1-u_1)v_1)^{-\left( 4+ \varphi +2\lambda \right)} \overleftrightarrow {w}_{1,1}^{-(\alpha +\lambda )} \left(\overleftrightarrow {w}_{1,1} \partial_{\overleftrightarrow {w}_{1,1}} \right) \overleftrightarrow {w}_{1,1}^{\alpha +\lambda } \nonumber\\
&&\times \left\{ c_0 x^{\lambda } \sum_{i_0=0}^{q_0} \frac{(-q_0)_{i_0} \left(q_0+ \varphi +2 \lambda \right)_{i_0}}{(1+\lambda )_{i_0}(\gamma +\lambda )_{i_0}} \overleftrightarrow {w}_{1,1} ^{i_0}\right\}z \label{eq:30032}
\end{eqnarray}
where
\begin{equation}
\overleftrightarrow {w}_{1,1} = \frac{v_1}{v_1-1} \frac{x t_1 u_1}{1-x(1-t_1)(1-u_1)v_1}\nonumber
\end{equation}
Put $l=2$ in (\ref{eq:30028}). Take the new (\ref{eq:30028}) into (\ref{eq:30031c}).
\begin{eqnarray}
y_2(x) &=& c_0 x^{\lambda } \int_{0}^{1} dt_2\;t_2^{3+\lambda } \int_{0}^{1} du_2\;u_2^{2+\gamma +\lambda } \frac{1}{2\pi i} \oint dv_2 \;\frac{1}{v_2} 
\left( \frac{v_2-1}{v_2} \frac{1}{1-x(1-t_2)(1-u_2)v_2}\right)^{q_2}  \nonumber\\
&&\times (1-x(1-t_2)(1-u_2)v_2)^{-\left( 8+ \varphi +2 \lambda \right)} 
 \overleftrightarrow {w}_{2,2}^{-(2+\alpha +\lambda) } \left(\overleftrightarrow {w}_{2,2} \partial_{\overleftrightarrow {w}_{2,2}} \right) \overleftrightarrow {w}_{2,2}^{2+\alpha +\lambda } \nonumber\\
&&\times \Bigg\{\sum_{i_0=0}^{q_0}\frac{(i_0+ \lambda +\alpha ) }{(i_0+ \lambda +2)(i_0+ \lambda +1+\gamma )}\frac{(-q_0)_{i_0} \left(q_0+  \varphi +2 \lambda \right)_{i_0}}{(1+\lambda )_{i_0}(\gamma +\lambda )_{i_0}} \nonumber\\
&&\times  \sum_{i_1=i_0}^{q_1} \frac{(-q_1)_{i_1}\left(q_1 + 4+ \varphi +2 \lambda \right)_{i_1}(3+\lambda )_{i_0}(2+\gamma +\lambda )_{i_0}}{(-q_1)_{i_0}\left(q_1 + 4+ \varphi +2 \lambda \right)_{i_0}(3+\lambda )_{i_1}(2+\gamma +\lambda )_{i_1}} \overleftrightarrow {w}_{2,2}^{i_1} \Bigg\} z^2 \label{eq:30033}
\end{eqnarray}
where
\begin{equation}
\overleftrightarrow {w}_{2,2} = \frac{v_2}{v_2-1} \frac{x t_2 u_2}{1-x(1-t_2)(1-u_2)v_2}\nonumber
\end{equation}
Put $l=1$ and $\eta = \overleftrightarrow {w}_{2,2}$ in (\ref{eq:30028}). Take the new (\ref{eq:30028}) into (\ref{eq:30033}).
\begin{eqnarray}
y_2(x) &=& \int_{0}^{1} dt_2\;t_2^{3+\lambda } \int_{0}^{1} du_2\;u_2^{2+\gamma +\lambda } \frac{1}{2\pi i} \oint dv_2 \;\frac{1}{v_2} 
\left( \frac{v_2-1}{v_2} \frac{1}{1-x(1-t_2)(1-u_2)v_2}\right)^{q_2}  \nonumber\\
&&\times (1-x(1-t_2)(1-u_2)v_2)^{-\left( 8+ \varphi +2 \lambda \right)} \overleftrightarrow {w}_{2,2}^{-(2+\alpha +\lambda) } \left(\overleftrightarrow {w}_{2,2} \partial_{\overleftrightarrow {w}_{2,2}} \right) \overleftrightarrow {w}_{2,2}^{2+\alpha +\lambda }  \nonumber\\
&&\times \int_{0}^{1} dt_1\;t_1^{1+\lambda } \int_{0}^{1} du_1\;u_1^{\gamma +\lambda } \frac{1}{2\pi i} \oint dv_1 \;\frac{1}{v_1} 
\left( \frac{v_1-1}{v_1} \frac{1}{1- \overleftrightarrow {w}_{2,2} (1-t_1)(1-u_1)v_1}\right)^{q_1}  \nonumber\\
&&\times (1- \overleftrightarrow {w}_{2,2} (1-t_1)(1-u_1)v_1)^{-\left( 4+ \varphi +2 \lambda \right)}  \overleftrightarrow {w}_{1,2}^{-(\alpha +\lambda)} \left(\overleftrightarrow {w}_{1,2} \partial_{\overleftrightarrow {w}_{1,2}} \right) \overleftrightarrow {w}_{1,2}^{\alpha +\lambda } \nonumber\\
&&\times \left\{ c_0 x^{\lambda } \sum_{i_0=0}^{q_0} \frac{(-q_0)_{i_0} \left(q_0+ \varphi +2 \lambda \right)_{i_0}}{(1+\lambda )_{i_0}(\gamma +\lambda )_{i_0}} \overleftrightarrow {w}_{1,2} ^{i_0}\right\} z^2 \label{eq:30034}
\end{eqnarray}
where
\begin{equation}
\overleftrightarrow {w}_{1,2} = \frac{v_1}{v_1-1} \frac{\overleftrightarrow {w}_{2,2} t_1 u_1}{1- \overleftrightarrow {w}_{2,2}(1-t_1)(1-u_1)v_1}\nonumber
\end{equation}
By using similar process for the previous cases of integral forms of $y_1(x)$ and $y_2(x)$, the integral form of sub-power series expansion of $y_3(x)$ is
\begin{eqnarray}
y_3(x) &=& \int_{0}^{1} dt_3\;t_3^{5+\lambda } \int_{0}^{1} du_3\;u_3^{4+\gamma +\lambda } \frac{1}{2\pi i} \oint dv_3 \;\frac{1}{v_3} 
\left( \frac{v_3-1}{v_3} \frac{1}{1-x(1-t_3)(1-u_3)v_3}\right)^{q_3}  \nonumber\\
&&\times (1-x(1-t_3)(1-u_3)v_3)^{-\left( 12+ \varphi +2 \lambda \right)} \overleftrightarrow {w}_{3,3}^{-(4+\alpha +\lambda) } \left(\overleftrightarrow {w}_{3,3} \partial_{\overleftrightarrow {w}_{3,3}} \right) \overleftrightarrow {w}_{3,3}^{4+\alpha +\lambda }  \nonumber\\
&&\times \int_{0}^{1} dt_2\;t_2^{3+\lambda } \int_{0}^{1} du_2\;u_2^{2+\gamma +\lambda } \frac{1}{2\pi i} \oint dv_2 \;\frac{1}{v_2} 
\left( \frac{v_2-1}{v_2} \frac{1}{1- \overleftrightarrow {w}_{3,3} (1-t_2)(1-u_2)v_2}\right)^{q_2}  \nonumber\\
&&\times (1- \overleftrightarrow {w}_{3,3} (1-t_2)(1-u_2)v_2)^{-\left( 8+ \varphi +2 \lambda \right)} \overleftrightarrow {w}_{2,3}^{-(2+\alpha +\lambda)} \left(\overleftrightarrow {w}_{2,3} \partial_{\overleftrightarrow {w}_{2,3}} \right) \overleftrightarrow {w}_{2,3}^{2+\alpha +\lambda } \nonumber\\
&&\times \int_{0}^{1} dt_1\;t_1^{1+\lambda } \int_{0}^{1} du_1\;u_1^{\gamma +\lambda } \frac{1}{2\pi i} \oint dv_1 \;\frac{1}{v_1} 
\left( \frac{v_1-1}{v_1} \frac{1}{1- \overleftrightarrow {w}_{2,3} (1-t_1)(1-u_1)v_1}\right)^{q_1}  \nonumber\\
&&\times (1- \overleftrightarrow {w}_{2,3} (1-t_1)(1-u_1)v_1)^{-\left( 4+ \varphi +2 \lambda \right)} \overleftrightarrow {w}_{1,3}^{-(\alpha +\lambda)} \left(\overleftrightarrow {w}_{1,3} \partial_{\overleftrightarrow {w}_{1,3}} \right) \overleftrightarrow {w}_{1,3}^{\alpha +\lambda } \nonumber\\
&&\times \left\{ c_0 x^{\lambda } \sum_{i_0=0}^{q_0} \frac{(-q_0)_{i_0} \left(q_0+ \varphi +2 \lambda \right)_{i_0}}{(1+\lambda )_{i_0}(\gamma +\lambda )_{i_0}} \overleftrightarrow {w}_{1,3} ^{i_0}\right\} z^3 \label{eq:30035}
\end{eqnarray}
where
\begin{equation}
\begin{cases} \overleftrightarrow {w}_{3,3} = \frac{v_3}{v_3-1} \frac{ x t_3 u_3}{1- x(1-t_3)(1-u_3)v_3} \cr
\overleftrightarrow {w}_{2,3} = \frac{v_2}{v_2-1} \frac{\overleftrightarrow {w}_{3,3} t_2 u_2}{1- \overleftrightarrow {w}_{3,3}(1-t_2)(1-u_2)v_2} \cr
\overleftrightarrow {w}_{1,3} = \frac{v_1}{v_1-1} \frac{\overleftrightarrow {w}_{2,3} t_1 u_1}{1- \overleftrightarrow {w}_{2,3}(1-t_1)(1-u_1)v_1}
\end{cases}
\nonumber
\end{equation}
By repeating this process for all higher terms of integral forms of sub-summation $y_m(x)$ terms where $m \geq 4$, we obtain every integral forms of $y_m(x)$ terms. 
Substitute (\ref{eq:30031a}), (\ref{eq:30032}), (\ref{eq:30034}), (\ref{eq:30035}) and including all integral forms of $y_m(x)$ terms where $m \geq 4$ into (\ref{eq:30030}). 
\end{proof}
Put $c_0$= 1 as $\lambda =0$  for the first kind of independent solutions of the CHE and $\lambda = 1-\gamma $  for the second one in (\ref{eq:30029}).
\begin{remark}
The integral representation of the CHE of the first kind for polynomial of type 2 about $x=0$ as $q= (q_j+2j)(-\beta +\gamma +\delta +q_j+2j-1)$ where $j,q_j \in \mathbb{N}_{0}$ is
\begin{eqnarray}
  y(x)&=& H_c^{(a)}F_{q_j}^R\left(\alpha, \beta, \gamma, \delta, q= (q_j+2j)(q_j+2j+\varphi ); \varphi = -\beta +\gamma +\delta -1; z = \beta x^2  \right)\nonumber\\
&=& _2F_1 \left( -q_0, q_0+\varphi ; \gamma; x \right) \nonumber\\
&&+ \sum_{n=1}^{\infty } \Bigg\{\prod _{k=0}^{n-1} \Bigg\{ \int_{0}^{1} dt_{n-k}\;t_{n-k}^{2(n-k)-1 } \int_{0}^{1} du_{n-k}\;u_{n-k}^{2(n-k-1)+\gamma } \nonumber\\
&&\times  \frac{1}{2\pi i}  \oint dv_{n-k} \frac{1}{v_{n-k}} \left( \frac{v_{n-k}-1}{v_{n-k}} \frac{1}{1-\overleftrightarrow {w}_{n-k+1,n}(1-t_{n-k})(1-u_{n-k})v_{n-k}}\right)^{q_{n-k}} \nonumber\\
&&\times \left( 1- \overleftrightarrow {w}_{n-k+1,n}(1-t_{n-k})(1-u_{n-k})v_{n-k}\right)^{-\left(4(n-k)+ \varphi \right)}\nonumber\\
&&\times \overleftrightarrow {w}_{n-k,n}^{-(2(n-k-1)+\alpha )}\left(  \overleftrightarrow {w}_{n-k,n} \partial _{ \overleftrightarrow {w}_{n-k,n}}\right) \overleftrightarrow {w}_{n-k,n}^{2(n-k-1)+\alpha } \Bigg\}\nonumber\\
&&\times\; _2F_1 \left( -q_0, q_0+\varphi ; \gamma; \overleftrightarrow {w}_{1,n} \right) \Bigg\} z^n \label{eq:30036}
\end{eqnarray}
\end{remark} 
\begin{remark}
The integral representation of the CHE of the first kind for polynomial of type 2 about $x=0$ as $q= (q_j+2j+1-\gamma )(-\beta +\delta +q_j+2j)$ where $j,q_j \in \mathbb{N}_{0}$ is
\begin{eqnarray}
 y(x)&=& H_c^{(a)}S_{\alpha _j}^R\left(\alpha, \beta, \gamma, \delta, q= (q_j+2j+1-\gamma )(q_j+2j+1-\gamma +\varphi ); \varphi = -\beta +\gamma +\delta -1; z= \beta x^2  \right)\nonumber\\
&=& x^{1-\gamma } \Bigg\{\; _2F_1 \left( -q_0, q_0+\varphi +2(1-\gamma ); 2-\gamma; x \right) \nonumber\\
&&+ \sum_{n=1}^{\infty } \Bigg\{\prod _{k=0}^{n-1} \Bigg\{ \int_{0}^{1} dt_{n-k}\;t_{n-k}^{2(n-k)-\gamma } \int_{0}^{1} du_{n-k}\;u_{n-k}^{2(n-k)-1} \nonumber\\
&&\times  \frac{1}{2\pi i}  \oint dv_{n-k} \frac{1}{v_{n-k}} \left( \frac{v_{n-k}-1}{v_{n-k}} \frac{1}{1-\overleftrightarrow {w}_{n-k+1,n}(1-t_{n-k})(1-u_{n-k})v_{n-k}}\right)^{q_{n-k}} \nonumber\\
&&\times \left( 1- \overleftrightarrow {w}_{n-k+1,n}(1-t_{n-k})(1-u_{n-k})v_{n-k}\right)^{-\left(4(n-k)+ \varphi +2(1-\gamma ) \right)} \nonumber\\
&&\times \overleftrightarrow {w}_{n-k,n}^{-(2(n-k)-1+\alpha -\gamma )}\left(  \overleftrightarrow {w}_{n-k,n} \partial _{ \overleftrightarrow {w}_{n-k,n}}\right) \overleftrightarrow {w}_{n-k,n}^{2(n-k)-1+\alpha -\gamma } \Bigg\}\nonumber\\
&&\times\; _2F_1 \left( -q_0, q_0+\varphi +2(1-\gamma ); 2-\gamma; \overleftrightarrow {w}_{1,n} \right) \Bigg\} z^n \Bigg\} \label{eq:30037}
\end{eqnarray}
\end{remark}
\subsubsection{Infinite series}
Let's consider the integral representation of the CHE about $x=0$ for infinite series by applying R3TRF.
There is a generalized hypergeometric function which is written by
\begin{eqnarray}
M_l &=& \sum_{i_l= i_{l-1}}^{\infty } \frac{\left( \Delta_{l}^{-}\right)_{i_l}\left(\Delta_{l}^{+} \right)_{i_l}(2l+1+\lambda )_{i_{l-1}}(2l+\gamma +\lambda )_{i_{l-1}}}{\left( \Delta_{l}^{-}\right)_{i_{l-1}}\left( \Delta_{l}^{+}\right)_{i_{l-1}}(2l+1+\lambda )_{i_l}(2l+\gamma +\lambda )_{i_l}} x^{i_l}\label{er:30024}\\
&=& x^{i_{l-1}} 
\sum_{j=0}^{\infty } \frac{B(i_{l-1}+2l+\lambda ,j+1) B(i_{l-1}+2l-1+\gamma +\lambda ,j+1) \left( \Delta_{l}^{-}+ i_{l-1} \right)_j \left( \Delta_{l}^{+} + i_{l-1} \right)_j}{(i_{l-1}+2l+\lambda )^{-1}(i_{l-1}+2l-1+\gamma +\lambda )^{-1}(1)_j \;j!} x^j\nonumber
\end{eqnarray}
where
\begin{equation} 
\Delta_{l}^{\pm}= \frac{ \varphi + 2 \lambda +4l \pm\sqrt{\varphi ^2+4q}}{2} \nonumber 
\end{equation}
Substitute (\ref{eq:30025a}) and (\ref{eq:30025b}) into (\ref{er:30024}). And divide $(i_{l-1}+2l+\lambda )(i_{l-1}+2l-1+\gamma +\lambda )$ into the new (\ref{er:30024}).
\begin{eqnarray}
V_l&=& \frac{1}{(i_{l-1}+2l+\lambda )(i_{l-1}+2l-1+\gamma +\lambda )} \sum_{i_l= i_{l-1}}^{\infty } \frac{\left( \Delta_{l}^{-}\right)_{i_l}\left(\Delta_{l}^{+} \right)_{i_l}(2l+1+\lambda )_{i_{l-1}}(2l+\gamma +\lambda )_{i_{l-1}}}{\left( \Delta_{l}^{-}\right)_{i_{l-1}}\left( \Delta_{l}^{+}\right)_{i_{l-1}}(2l+1+\lambda )_{i_l}(2l+\gamma +\lambda )_{i_l}} x^{i_l}\nonumber\\
&=&  \int_{0}^{1} dt_l\;t_l^{2l-1+\lambda } \int_{0}^{1} du_l\;u_l^{2(l-1)+\gamma +\lambda } (x t_l u_l)^{i_{l-1}} \sum_{j=0}^{\infty } \frac{\left( \Delta_{l}^{-}+ i_{l-1} \right)_j \left( \Delta_{l}^{+} + i_{l-1} \right)_j}{(1)_j \;j!} 
(x(1-t_l)(1-u_l))^j \nonumber 
\end{eqnarray}
The hypergeometric function is defined by
\begin{eqnarray}
_2F_1 \left( \alpha ,\beta ; \gamma ; z \right) &=& \sum_{n=0}^{\infty } \frac{(\alpha )_n (\beta )_n}{(\gamma )_n (n!)} z^n \nonumber\\
&=&  \frac{1}{2\pi i} \frac{\Gamma( 1+\alpha  -\gamma )}{\Gamma (\alpha )} \int_0^{(1+)} dv_l\; (-1)^{\gamma }(-v_l)^{\alpha -1} (1-v_l )^{\gamma -\alpha -1} (1-zv_l)^{-\beta }\hspace{1.5cm}\label{er:30025}\\
&& \mbox{where} \;\gamma -\alpha  \ne 1,2,3,\cdots, \;\mbox{Re}(\alpha )>0 \nonumber
\end{eqnarray}
Replace $\alpha $, $\beta $, $\gamma $ and $z$ by $\Delta_{l}^{-}+ i_{l-1}$, $ { \displaystyle  \Delta_{l}^{+} + i_{l-1} }$, 1 and $x (1-t_l)(1-u_l)$ in (\ref{er:30025}). Take the new (\ref{er:30025}) into $V_l$.
\begin{eqnarray}
V_l&=&  \frac{1}{(i_{l-1}+2l+\lambda )(i_{l-1}+2l-1+\gamma +\lambda )} \sum_{i_l= i_{l-1}}^{\infty } \frac{\left( \Delta_{l}^{-}\right)_{i_l}\left(\Delta_{l}^{+} \right)_{i_l}(2l+1+\lambda )_{i_{l-1}}(2l+\gamma +\lambda )_{i_{l-1}}}{\left( \Delta_{l}^{-}\right)_{i_{l-1}}\left( \Delta_{l}^{+}\right)_{i_{l-1}}(2l+1+\lambda )_{i_l}(2l+\gamma +\lambda )_{i_l}} x^{i_l} \nonumber\\
&=&  \int_{0}^{1} dt_l\;t_l^{2l-1+\lambda } \int_{0}^{1} du_l\;u_l^{2(l-1)+\gamma +\lambda } 
\frac{1}{2\pi i} \oint dv_l\;\frac{1}{v_l} \left(\frac{v_l-1}{v_l} \right)^{-\Delta_{l}^{-}} \nonumber\\
&&\times (1-x(1-t_l)(1-u_l)v_l)^{-\Delta_{l}^{+}} \left(\frac{v_l}{v_l-1} \frac{xt_l u_l}{1-x(1-t_l)(1-u_l)v_l}\right)^{i_{l-1}}\label{er:30026} 
\end{eqnarray}
Substitute (\ref{er:30026}) into (\ref{eq:30020}) where $l=1,2,3,\cdots$; apply $V_1$ into the second summation of sub-power series $y_1(x)$, apply $V_2$ into the third summation and $V_1$ into the second summation of sub-power series $y_2(x)$, apply $V_3$ into the forth summation, $V_2$ into the third summation and $V_1$ into the second summation of sub-power series $y_3(x)$, etc.\footnote{$y_1(x)$ means the sub-power series in (\ref{eq:30020}) contains one term of $B_n's$, $y_2(x)$ means the sub-power series in (\ref{eq:30020}) contains two terms of $B_n's$, $y_3(x)$ means the sub-power series in (\ref{eq:30020}) contains three terms of $B_n's$, etc.}
\begin{theorem}
The general representation in the form of integral of the CHE for infinite series about $x=0$ using R3TRF is given by
\begin{eqnarray}
 y(x) &=& \sum_{n=0}^{\infty } y_n(x)= y_0(x)+ y_1(x)+ y_2(x)+ y_3(x)+\cdots \nonumber\\
&=&  c_0 x^{\lambda } \left\{ \sum_{i_0=0}^{\infty }\frac{\left( \Delta_0^{-}\right)_{i_0}\left( \Delta_0^{+}\right)_{i_0}}{(1+\lambda )_{i_0}(\gamma +\lambda )_{i_0}}  x^{i_0}\right. \nonumber\\
&&+ \sum_{n=1}^{\infty } \left\{\prod _{k=0}^{n-1} \Bigg\{ \int_{0}^{1} dt_{n-k}\;t_{n-k}^{2(n-k)-1+\lambda } \int_{0}^{1} du_{n-k}\;u_{n-k}^{2(n-k-1)+\gamma +\lambda } \right.\nonumber\\
&\times& \frac{1}{2\pi i}  \oint dv_{n-k} \frac{1}{v_{n-k}}\left( \frac{v_{n-k}-1}{v_{n-k}}\right)^{ -\Delta_{n-k}^{-}}  \left( 1- \overleftrightarrow {w}_{n-k+1,n}(1-t_{n-k})(1-u_{n-k})v_{n-k}\right)^{ -\Delta_{n-k}^{+}}\nonumber\\
&\times& \overleftrightarrow {w}_{n-k,n}^{-(2(n-k-1)+\alpha +\lambda )}\left(  \overleftrightarrow {w}_{n-k,n} \partial _{ \overleftrightarrow {w}_{n-k,n}}\right) \overleftrightarrow {w}_{n-k,n}^{2(n-k-1)+\alpha +\lambda } \Bigg\}\nonumber\\
&\times& \left.\left.\sum_{i_0=0}^{\infty }\frac{\left( \Delta_0^{-}\right)_{i_0}\left( \Delta_0^{+}\right)_{i_0}}{(1+\lambda )_{i_0}(\gamma +\lambda )_{i_0}} \overleftrightarrow {w}_{1,n}^{i_0}\right\} z^n \right\} \label{eq:30038}
\end{eqnarray}
where
\begin{equation}
\begin{cases} 
\Delta_0^{\pm}= \frac{ \varphi +2 \lambda  \pm\sqrt{\varphi ^2+4q}}{2} \cr
\Delta_{n-k}^{\pm}=   \frac{ \varphi + 2 \lambda +4(n-k) \pm\sqrt{\varphi ^2+4q}}{2}
\end{cases}\nonumber 
\end{equation}
In the above, the first sub-integral form contains one term of $B_n's$, the second one contains two terms of $B_n$'s, the third one contains three terms of $B_n$'s, etc.
\end{theorem}
\begin{proof} 
In (\ref{eq:30020}) sub-power series $y_0(x) $, $y_1(x)$, $y_2(x)$ and $y_3(x)$ of the CHE for infinite series about $x=0$ using R3TRF are given by
\begin{subequations}
\begin{equation}
 y_0(x)= c_0 x^{\lambda } \sum_{i_0=0}^{\infty } \frac{\left( \Delta_{0}^{-}\right)_{i_0} \left( \Delta_{0}^{+}\right)_{i_0}}{(1+\lambda )_{i_0}(\gamma +\lambda )_{i_0}} x^{i_0} \label{er:30027a}
\end{equation}
\begin{eqnarray}
 y_1(x) &=& c_0 x^{\lambda } \left\{\sum_{i_0=0}^{\infty}\frac{(i_0+ \lambda +\alpha ) }{(i_0+ \lambda +2)(i_0+ \lambda +1+\gamma )}\frac{\left( \Delta_{0}^{-}\right)_{i_0} \left( \Delta_{0}^{+}\right)_{i_0}}{(1+\lambda )_{i_0}(\gamma +\lambda )_{i_0}} \right. \nonumber\\
&&\times \left. \sum_{i_1=i_0}^{\infty} \frac{\left( \Delta_{1}^{-}\right)_{i_1}\left( \Delta_{1}^{+}\right)_{i_1}(3+\lambda )_{i_0}(2+\gamma +\lambda )_{i_0}}{\left( \Delta_{1}^{-}\right)_{i_0}\left( \Delta_{1}^{+}\right)_{i_0}(3+\lambda )_{i_1}(2+\gamma +\lambda )_{i_1}} x^{i_1} \right\} z  \label{er:30027b}
\end{eqnarray}
\begin{eqnarray}
 y_2(x) &=& c_0 x^{\lambda } \left\{\sum_{i_0=0}^{\infty}\frac{(i_0+ \lambda +\alpha ) }{(i_0+ \lambda +2)(i_0+ \lambda +1+\gamma )}\frac{\left( \Delta_{0}^{-}\right)_{i_0} \left( \Delta_{0}^{+}\right)_{i_0}}{(1+\lambda )_{i_0}(\gamma +\lambda )_{i_0}} \right. \nonumber\\
&&\times  \sum_{i_1=i_0}^{\infty} \frac{(i_1+2+ \lambda +\alpha ) }{(i_1+ \lambda +4)(i_1+ \lambda +3+\gamma )} \frac{\left( \Delta_{1}^{-}\right)_{i_1}\left( \Delta_{1}^{+}\right)_{i_1}(3+\lambda )_{i_0}(2+\gamma +\lambda )_{i_0}}{\left( \Delta_{1}^{-}\right)_{i_0}\left( \Delta_{1}^{+}\right)_{i_0}(3+\lambda )_{i_1}(2+\gamma +\lambda )_{i_1}} \nonumber\\
&&\times \left.\sum_{i_2=i_1}^{\infty} \frac{\left( \Delta_{2}^{-}\right)_{i_2}\left( \Delta_{2}^{+}\right)_{i_2}(5+\lambda )_{i_1}(4+\gamma +\lambda )_{i_1}}{\left( \Delta_{2}^{-}\right)_{i_1}\left( \Delta_{2}^{+}\right)_{i_1}(5+\lambda )_{i_2}(4+\gamma +\lambda )_{i_2}} x^{i_2} \right\} z^2  \label{er:30027c}
\end{eqnarray}
\begin{eqnarray}
 y_3(x) &=&  c_0 x^{\lambda } \left\{\sum_{i_0=0}^{\infty}\frac{(i_0+ \lambda +\alpha ) }{(i_0+ \lambda +2)(i_0+ \lambda +1+\gamma )}\frac{\left( \Delta_{0}^{-}\right)_{i_0} \left( \Delta_{0}^{+}\right)_{i_0}}{(1+\lambda )_{i_0}(\gamma +\lambda )_{i_0}} \right.\nonumber\\
&&\times  \sum_{i_1=i_0}^{\infty} \frac{(i_1+2+ \lambda +\alpha ) }{(i_1+ \lambda +4)(i_1+ \lambda +3+\gamma )} \frac{\left( \Delta_{1}^{-}\right)_{i_1}\left( \Delta_{1}^{+}\right)_{i_1}(3+\lambda )_{i_0}(2+\gamma +\lambda )_{i_0}}{\left( \Delta_{1}^{-}\right)_{i_0}\left( \Delta_{1}^{+}\right)_{i_0}(3+\lambda )_{i_1}(2+\gamma +\lambda )_{i_1}} \nonumber\\
&&\times \sum_{i_2=i_1}^{\infty} \frac{(i_2+4+ \lambda +\alpha ) }{(i_2+ \lambda +6)(i_2+ \lambda +5+\gamma )}  \frac{\left( \Delta_{2}^{-}\right)_{i_2}\left( \Delta_{2}^{+}\right)_{i_2}(5+\lambda )_{i_1}(4+\gamma +\lambda )_{i_1}}{\left( \Delta_{2}^{-}\right)_{i_1}\left( \Delta_{2}^{+}\right)_{i_1}(5+\lambda )_{i_2}(4+\gamma +\lambda )_{i_2}} \nonumber\\
&&\times \left.\sum_{i_3=i_2}^{\infty} \frac{\left( \Delta_{3}^{-}\right)_{i_3}\left( \Delta_{3}^{+}\right)_{i_3}(7+\lambda )_{i_2}(6+\gamma +\lambda )_{i_2}}{\left( \Delta_{3}^{-}\right)_{i_2}\left( \Delta_{3}^{+}\right)_{i_2}(7+\lambda )_{i_3}(6+\gamma +\lambda )_{i_3}}x^{i_3} \right\} z^3  \label{er:30027d} 
\end{eqnarray}
\end{subequations}
where
\begin{equation}
\begin{cases}
\Delta_{0}^{\pm}=    \frac{ \varphi + 2 \lambda \pm\sqrt{\varphi ^2+4q}}{2} \cr
\Delta_{1}^{\pm}=    \frac{ \varphi + 2 \lambda +4 \pm\sqrt{\varphi ^2+4q}}{2} \cr
\Delta_{2}^{\pm}=    \frac{ \varphi + 2 \lambda +8 \pm\sqrt{\varphi ^2+4q}}{2} \cr
\Delta_{3}^{\pm}=    \frac{ \varphi + 2 \lambda +12 \pm\sqrt{\varphi ^2+4q}}{2} 
\end{cases}\nonumber 
\end{equation}
Put $l=1$ in (\ref{er:30026}). Take the new (\ref{er:30026}) into (\ref{er:30027b}).
\begin{eqnarray}
y_1(x) &=& \int_{0}^{1} dt_1\;t_1^{1+\lambda } \int_{0}^{1} du_1\;u_1^{\gamma +\lambda } \frac{1}{2\pi i} \oint dv_1 \;\frac{1}{v_1} 
\left( \frac{v_1-1}{v_1}  \right)^{-\Delta_{1}^{-}}  \nonumber\\
&&\times (1-x(1-t_1)(1-u_1)v_1)^{-\Delta_{1}^{+}} \overleftrightarrow {w}_{1,1}^{-(\alpha +\lambda )} \left(\overleftrightarrow {w}_{1,1} \partial_{\overleftrightarrow {w}_{1,1}} \right) \overleftrightarrow {w}_{1,1}^{\alpha +\lambda } \nonumber\\
&&\times \left\{ c_0 x^{\lambda } \sum_{i_0=0}^{\infty } \frac{\left( \Delta_{0}^{-}\right)_{i_0} \left( \Delta_{0}^{+}\right)_{i_0}}{(1+\lambda )_{i_0}(\gamma +\lambda )_{i_0}} \overleftrightarrow {w}_{1,1} ^{i_0}\right\}z \label{er:30028}
\end{eqnarray}
where
\begin{equation}
\overleftrightarrow {w}_{1,1} = \frac{v_1}{v_1-1} \frac{x t_1 u_1}{1-x(1-t_1)(1-u_1)v_1}\nonumber
\end{equation}
Put $l=2$ in (\ref{er:30026}). Take the new (\ref{er:30026}) into (\ref{er:30027c}).
\begin{eqnarray}
y_2(x) &=& c_0 x^{\lambda } \int_{0}^{1} dt_2\;t_2^{3+\lambda } \int_{0}^{1} du_2\;u_2^{2+\gamma +\lambda } \frac{1}{2\pi i} \oint dv_2 \;\frac{1}{v_2} \left( \frac{v_2-1}{v_2} \right)^{-\Delta_{2}^{-}}  \nonumber\\
&&\times (1-x(1-t_2)(1-u_2)v_2)^{-\Delta_{2}^{+}} 
 \overleftrightarrow {w}_{2,2}^{-(2+\alpha +\lambda) } \left(\overleftrightarrow {w}_{2,2} \partial_{\overleftrightarrow {w}_{2,2}} \right) \overleftrightarrow {w}_{2,2}^{2+\alpha +\lambda } \nonumber\\
&&\times \left\{\sum_{i_0=0}^{\infty}\frac{(i_0+ \lambda +\alpha ) }{(i_0+ \lambda +2)(i_0+ \lambda +1+\gamma )}\frac{\left( \Delta_{0}^{-}\right)_{i_0} \left( \Delta_{0}^{+}\right)_{i_0}}{(1+\lambda )_{i_0}(\gamma +\lambda )_{i_0}} \right. \nonumber\\
&&\times \left. \sum_{i_1=i_0}^{\infty} \frac{\left( \Delta_{1}^{-}\right)_{i_1}\left( \Delta_{1}^{+}\right)_{i_1}(3+\lambda )_{i_0}(2+\gamma +\lambda )_{i_0}}{\left( \Delta_{1}^{-}\right)_{i_0}\left( \Delta_{1}^{+}\right)_{i_0}(3+\lambda )_{i_1}(2+\gamma +\lambda )_{i_1}} \overleftrightarrow {w}_{2,2}^{i_1} \right\} z^2 \label{er:30029}
\end{eqnarray}
where
\begin{equation}
\overleftrightarrow {w}_{2,2} = \frac{v_2}{v_2-1} \frac{x t_2 u_2}{1-x(1-t_2)(1-u_2)v_2}\nonumber
\end{equation}
Put $l=1$ and $\eta = \overleftrightarrow {w}_{2,2}$ in (\ref{er:30026}). Take the new (\ref{er:30026}) into (\ref{er:30029}).
\begin{eqnarray}
y_2(x) &=& \int_{0}^{1} dt_2\;t_2^{3+\lambda } \int_{0}^{1} du_2\;u_2^{2+\gamma +\lambda } \frac{1}{2\pi i} \oint dv_2 \;\frac{1}{v_2} 
\left( \frac{v_2-1}{v_2} \right)^{-\Delta_{2}^{-}}  \nonumber\\
&&\times (1-x(1-t_2)(1-u_2)v_2)^{-\Delta_{2}^{+}} \overleftrightarrow {w}_{2,2}^{-(2+\alpha +\lambda) } \left(\overleftrightarrow {w}_{2,2} \partial_{\overleftrightarrow {w}_{2,2}} \right) \overleftrightarrow {w}_{2,2}^{2+\alpha +\lambda }  \nonumber\\
&&\times \int_{0}^{1} dt_1\;t_1^{1+\lambda } \int_{0}^{1} du_1\;u_1^{\gamma +\lambda } \frac{1}{2\pi i} \oint dv_1 \;\frac{1}{v_1} 
\left( \frac{v_1-1}{v_1}  \right)^{-\Delta_{1}^{-}}  \nonumber\\
&&\times (1- \overleftrightarrow {w}_{2,2} (1-t_1)(1-u_1)v_1)^{-\Delta_{1}^{+}}  \overleftrightarrow {w}_{1,2}^{-(\alpha +\lambda)} \left(\overleftrightarrow {w}_{1,2} \partial_{\overleftrightarrow {w}_{1,2}} \right) \overleftrightarrow {w}_{1,2}^{\alpha +\lambda } \nonumber\\
&&\times \left\{ \sum_{i_0=0}^{\infty } \frac{\left( \Delta_{0}^{-}\right)_{i_0} \left( \Delta_{0}^{+}\right)_{i_0}}{(1+\lambda )_{i_0}(\gamma +\lambda )_{i_0}} \overleftrightarrow {w}_{1,2} ^{i_0}\right\} z^2 \label{er:30030}
\end{eqnarray}
where
\begin{equation}
\overleftrightarrow {w}_{1,2} = \frac{v_1}{v_1-1} \frac{\overleftrightarrow {w}_{2,2} t_1 u_1}{1- \overleftrightarrow {w}_{2,2}(1-t_1)(1-u_1)v_1}\nonumber
\end{equation}
By using similar process for the previous cases of integral forms of $y_1(x)$ and $y_2(x)$, the integral form of sub-power series expansion of $y_3(x)$ is
\begin{eqnarray}
y_3(x) &=& \int_{0}^{1} dt_3\;t_3^{5+\lambda } \int_{0}^{1} du_3\;u_3^{4+\gamma +\lambda } \frac{1}{2\pi i} \oint dv_3 \;\frac{1}{v_3} 
\left( \frac{v_3-1}{v_3} \right)^{-\Delta_{3}^{-}}  \nonumber\\
&&\times (1-x(1-t_3)(1-u_3)v_3)^{-\Delta_{3}^{+}} \overleftrightarrow {w}_{3,3}^{-(4+\alpha +\lambda) } \left(\overleftrightarrow {w}_{3,3} \partial_{\overleftrightarrow {w}_{3,3}} \right) \overleftrightarrow {w}_{3,3}^{4+\alpha +\lambda }  \nonumber\\
&&\times \int_{0}^{1} dt_2\;t_2^{3+\lambda } \int_{0}^{1} du_2\;u_2^{2+\gamma +\lambda } \frac{1}{2\pi i} \oint dv_2 \;\frac{1}{v_2} 
\left( \frac{v_2-1}{v_2} \right)^{-\Delta_{2}^{-}}  \nonumber\\
&&\times (1- \overleftrightarrow {w}_{3,3} (1-t_2)(1-u_2)v_2)^{-\Delta_{2}^{+}} \overleftrightarrow {w}_{2,3}^{-(2+\alpha +\lambda)} \left(\overleftrightarrow {w}_{2,3} \partial_{\overleftrightarrow {w}_{2,3}} \right) \overleftrightarrow {w}_{2,3}^{2+\alpha +\lambda } \nonumber\\
&&\times \int_{0}^{1} dt_1\;t_1^{1+\lambda } \int_{0}^{1} du_1\;u_1^{\gamma +\lambda } \frac{1}{2\pi i} \oint dv_1 \;\frac{1}{v_1} 
\left( \frac{v_1-1}{v_1} \right)^{-\Delta_{1}^{-}}  \nonumber\\
&&\times (1- \overleftrightarrow {w}_{2,3} (1-t_1)(1-u_1)v_1)^{-\Delta_{1}^{+}} \overleftrightarrow {w}_{1,3}^{-(\alpha +\lambda)} \left(\overleftrightarrow {w}_{1,3} \partial_{\overleftrightarrow {w}_{1,3}} \right) \overleftrightarrow {w}_{1,3}^{\alpha +\lambda } \nonumber\\
&&\times \left\{ c_0 x^{\lambda } \sum_{i_0=0}^{\infty } \frac{\left( \Delta_{0}^{-}\right)_{i_0} \left( \Delta_{0}^{+}\right)_{i_0}}{(1+\lambda )_{i_0}(\gamma +\lambda )_{i_0}} \overleftrightarrow {w}_{1,3} ^{i_0}\right\} z^3 \label{er:30031}
\end{eqnarray}
where
\begin{equation}
\begin{cases} \overleftrightarrow {w}_{3,3} = \frac{v_3}{v_3-1} \frac{ x t_3 u_3}{1- x(1-t_3)(1-u_3)v_3} \cr
\overleftrightarrow {w}_{2,3} = \frac{v_2}{v_2-1} \frac{\overleftrightarrow {w}_{3,3} t_2 u_2}{1- \overleftrightarrow {w}_{3,3}(1-t_2)(1-u_2)v_2} \cr
\overleftrightarrow {w}_{1,3} = \frac{v_1}{v_1-1} \frac{\overleftrightarrow {w}_{2,3} t_1 u_1}{1- \overleftrightarrow {w}_{2,3}(1-t_1)(1-u_1)v_1}
\end{cases}
\nonumber
\end{equation}
By repeating this process for all higher terms of integral forms of sub-summation $y_m(x)$ terms where $m \geq 4$, we obtain every integral forms of $y_m(x)$ terms. 
Since we substitute (\ref{er:30027a}), (\ref{er:30028}), (\ref{er:30030}), (\ref{er:30031}) and including all integral forms of $y_m(x)$ terms where $m \geq 4$ into (\ref{eq:30020}), we obtain (\ref{eq:30038}). \footnote{Or replace the finite summation with an interval $[0, q_0]$ by infinite summation with an interval  $[0,\infty ]$ in (\ref{eq:30029}). Replace $q_0$ and $q_{n-k}$ by $\frac{-(\varphi +2\lambda )+\sqrt{\varphi ^2+4q}}{2}$ and $\frac{-(\varphi +2(\lambda +2(n-k)))+\sqrt{\varphi ^2+4q}}{2}$ into the new (\ref{eq:30029}). Its solution is also equivalent to (\ref{eq:30038})}
\end{proof}
Put $c_0$= 1 as $\lambda =0$  for the first kind of independent solutions of the CHE and $\lambda = 1-\gamma $  for the second one in (\ref{eq:30038}). 
\begin{remark}
The integral representation of the CHE of the first kind for infinite series about $x=0$ using R3TRF is
\begin{eqnarray}
y(x)&=& H_c^{(a)}F^R\left(\alpha, \beta, \gamma, \delta, q; \varphi = -\beta +\gamma +\delta -1; z= \beta x^2 \right) \nonumber\\
&=& \; _2F_1 \left( \Delta_0^{-}, \Delta_0^{+}; \gamma; x \right) + \sum_{n=1}^{\infty } \Bigg\{\prod _{k=0}^{n-1} \Bigg\{ \int_{0}^{1} dt_{n-k}\;t_{n-k}^{2(n-k)-1} \int_{0}^{1} du_{n-k}\;u_{n-k}^{2(n-k-1)+\gamma } \nonumber\\
&\times& \frac{1}{2\pi i}  \oint dv_{n-k} \frac{1}{v_{n-k}}\left( \frac{v_{n-k}-1}{v_{n-k}}\right)^{ -\Delta_{n-k}^{-}}  \left( 1- \overleftrightarrow {w}_{n-k+1,n}(1-t_{n-k})(1-u_{n-k})v_{n-k}\right)^{-\Delta_{n-k}^{+}}\nonumber\\
&\times& \overleftrightarrow {w}_{n-k,n}^{-(2(n-k-1)+\alpha )}\left(  \overleftrightarrow {w}_{n-k,n} \partial _{ \overleftrightarrow {w}_{n-k,n}}\right) \overleftrightarrow {w}_{n-k,n}^{2(n-k-1)+\alpha } \Bigg\}\; _2F_1 \left( \Delta_0^{-}, \Delta_0^{+}; \gamma; \overleftrightarrow {w}_{1,n} \right) \Bigg\} z^n \hspace{2cm}\label{eq:30039}
\end{eqnarray}
where
\begin{equation}
\begin{cases} 
\Delta_0^{\pm} = \frac{ \varphi \pm\sqrt{\varphi ^2+4q}}{2} \cr
\Delta_{n-k}^{\pm} = \frac{ \varphi +4(n-k) \pm\sqrt{\varphi ^2+4q}}{2}
\end{cases}\nonumber 
\end{equation}
\end{remark}
\begin{remark}
The integral representation of the CHE of the first kind for infinite series about $x=0$ using R3TRF is
\begin{eqnarray}
 y(x)&=&  H_c^{(a)}S^R\left(\alpha, \beta, \gamma, \delta, q; \varphi = -\beta +\gamma +\delta -1; z= \beta x^2   \right)\nonumber\\
&=& x^{1-\gamma }  \Bigg\{ \; _2F_1 \left( \Delta_0^{-}, \Delta_0^{+}; 2-\gamma; x \right)  + \sum_{n=1}^{\infty } \Bigg\{\prod _{k=0}^{n-1} \Bigg\{ \int_{0}^{1} dt_{n-k}\;t_{n-k}^{2(n-k)-\gamma } \int_{0}^{1} du_{n-k}\;u_{n-k}^{2(n-k)-1} \nonumber\\
&\times& \frac{1}{2\pi i}  \oint dv_{n-k} \frac{1}{v_{n-k}}\left( \frac{v_{n-k}-1}{v_{n-k}}\right)^{-\Delta_{n-k}^{-}}  \left( 1- \overleftrightarrow {w}_{n-k+1,n}(1-t_{n-k})(1-u_{n-k})v_{n-k}\right)^{-\Delta_{n-k}^{+}}\nonumber\\
&\times& \overleftrightarrow {w}_{n-k,n}^{-(2(n-k)-1 +\alpha -\gamma )}\left(  \overleftrightarrow {w}_{n-k,n} \partial _{ \overleftrightarrow {w}_{n-k,n}}\right) \overleftrightarrow {w}_{n-k,n}^{2(n-k)-1 +\alpha -\gamma  } \Bigg\} \nonumber\\
&\times& \; _2F_1 \left( \Delta_0^{-}, \Delta_0^{+}; 2-\gamma; \overleftrightarrow {w}_{1,n} \right) \Bigg\} z^n \Bigg\}  \label{eq:30040}
\end{eqnarray}
where
\begin{equation}
\begin{cases} 
\Delta_0^{\pm} = \frac{ \varphi +2 (1-\gamma )  \pm\sqrt{\varphi ^2+4q}}{2} \cr
\Delta_{n-k}^{\pm} = \frac{ \varphi + 2(1-\gamma ) +4(n-k) \pm\sqrt{\varphi ^2+4q}}{2}
\end{cases}\nonumber 
\end{equation}
\end{remark}
\subsection[Generating function for the CHP of type 2]{Generating function for the CHP of type 2}
In chapter 4 I construct the generating function for the CHP of type 1 by applying 3TRF. Now I consider the generating function for the CHP of type 2. Since the generating function for the CHP is derived, we might be possible to construct orthogonal relations of the CHP.\footnote{Polynomial of type 3 is a polynomial which makes $A_n$ and $B_n$ terms terminated at the same time in three term recursion relation of the power series in a linear differential equation. For the type 1 polynomial, I treat $\beta $, $\gamma $, $\delta $ and $q$ as free variables and $\alpha $ as a fixed value. For the type 2 polynomial, I treat $\alpha $, $\beta $, $\gamma $ and $\delta $  as free variables and  $q$ as a fixed value. For the type 3 polynomial, I treat $\beta $, $\gamma $ and $\delta $  as free variables and  $\alpha $ and $q$ as fixed values. In the next series, I will construct the generating functions for the CHP of type 3.} 
\begin{lemma}
The generating function for Jacobi polynomial using hypergeometric functions is given by
\begin{eqnarray}
&&\sum_{q _0=0}^{\infty }\frac{(\gamma )_{q_0}}{q_0!} w^{q_0} \;_2F_1(-q_0, q_0+A; \gamma; x) \label{eq:30041}\\
&&= 2^{A -1}\frac{\left(1-w+\sqrt{w^2-2(1-2x)w+1}\right)^{1-\gamma } \left(1+w+\sqrt{w^2-2(1-2x)w+1}\right)^{\gamma -A}}{\sqrt{w^2-2(1-2x)w+1}} \nonumber\\
&& \hspace{.5cm} \mbox{where}\;|w|<1 \nonumber
\end{eqnarray}
\end{lemma}
\begin{proof}
The proof of this lemma is given in Lemma 3.2.1
\end{proof}

\begin{definition}
I define that
\begin{equation}
\begin{cases}
\displaystyle { s_{a,b}} = \begin{cases} \displaystyle {  s_a\cdot s_{a+1}\cdot s_{a+2}\cdots s_{b-2}\cdot s_{b-1}\cdot s_b}\;\;\mbox{if}\;a>b \cr
s_a \;\;\mbox{if}\;a=b\end{cases}
\cr
\cr
\displaystyle { \widetilde{w}_{i,j}}  = 
\begin{cases} \displaystyle { \frac{ \widetilde{w}_{i+1,j}\; t_i u_i \left\{ 1+ (s_i+2\widetilde{w}_{i+1,j}(1-t_i)(1-u_i))s_i\right\}}{2(1-\widetilde{w}_{i+1,j}(1-t_i)(1-u_i))^2 s_i}} \cr
\displaystyle {-\frac{\widetilde{w}_{i+1,j}\; t_i u_i (1+s_i)\sqrt{s_i^2-2(1-2\widetilde{w}_{i+1,j}(1-t_i)(1-u_i))s_i+1}}{2(1-\widetilde{w}_{i+1,j}(1-t_i)(1-u_i))^2 s_i}} \;\;\mbox{where}\;i<j \cr
\cr
\displaystyle { \frac{x t_i u_i \left\{ 1+ (s_{i,\infty }+2x(1-t_i)(1-u_i))s_{i,\infty }\right\}}{2(1-x(1-t_i)(1-u_i))^2 s_{i,\infty }}} \cr
\displaystyle {-\frac{x t_i u_i(1+s_{i,\infty })\sqrt{s_{i,\infty }^2-2(1-2x (1-t_i)(1-u_i))s_{i,\infty }+1}}{2(1-x (1-t_i)(1-u_i))^2 s_{i,\infty }}} \;\;\mbox{where}\;i=j 
\end{cases}
\end{cases}\label{eq:30046}
\end{equation}
where
\begin{equation}
a,b,i,j\in \mathbb{N}_{0} \nonumber
\end{equation}
\end{definition}
And we have
\begin{equation}
\sum_{q_i = q_j}^{\infty } r_i^{q_i} = \frac{r_i^{q_j}}{(1-r_i)}\label{eq:30047}
\end{equation}
Acting the summation operator $\displaystyle{ \sum_{q_0 =0}^{\infty } \frac{(\gamma')_{q_0}}{q_0!} s_0^{q_0} \prod _{n=1}^{\infty } \left\{ \sum_{ q_n = q_{n-1}}^{\infty } s_n^{q_n }\right\}}$ on (\ref{eq:30029}) where $|s_i|<1$ as $i=0,1,2,\cdots$ by using (\ref{eq:30046}) and (\ref{eq:30047}),
\begin{theorem} 
The general expression of the generating function for the CHP of type 2 is given by
\begin{eqnarray}
&&\sum_{q_0 =0}^{\infty } \frac{(\gamma')_{q_0}}{q_0!} s_0^{q_0} \prod _{n=1}^{\infty } \left\{ \sum_{ q_n = q_{n-1}}^{\infty } s_n^{q_n }\right\} y(x) \nonumber\\
&&= \prod_{l=1}^{\infty } \frac{1}{(1-s_{l,\infty })} \mathbf{\Upsilon}(\lambda; s_{0,\infty } ;x) \nonumber\\
&&+\Bigg\{ \prod_{l=2}^{\infty } \frac{1}{(1-s_{l,\infty })} \int_{0}^{1} dt_1\;t_1^{1+\lambda} \int_{0}^{1} du_1\;u_1^{\gamma +\lambda} \left(s_{1,\infty }^2-2(1-2 x(1-t_1)(1-u_1))s_{1,\infty }+1\right)^{-\frac{1}{2}}\nonumber\\
&&\times \left( \frac{ 1+s_{1,\infty } +\sqrt{s_{1,\infty }^2-2(1-2x(1-t_1)(1-u_1))s_{1,\infty }+1}}{2}\right)^{-\left(3+ \varphi +2 \lambda  \right)}\nonumber\\
&&\times \widetilde{w}_{1,1}^{-(\alpha +\lambda )}\left(  \widetilde{w}_{1,1} \partial _{ \widetilde{w}_{1,1}}\right) \widetilde{w}_{1,1}^{\alpha +\lambda} \mathbf{\Upsilon}(\lambda ; s_0;\widetilde{w}_{1,1}) \Bigg\}z\nonumber\\
&&+ \sum_{n=2}^{\infty } \left\{ \prod_{l=n+1}^{\infty } \frac{1}{(1-s_{l,\infty })} \int_{0}^{1} dt_n\;t_n^{2n-1+\lambda } \int_{0}^{1} du_n\;u_n^{2(n-1)+\gamma +\lambda} \right.\nonumber\\
&&\times \left( s_{n,\infty }^2-2(1-2x(1-t_n)(1-u_n))s_{n,\infty }+1\right)^{-\frac{1}{2}} \nonumber\\
&&\times  \left( \frac{ 1+s_{n,\infty } +\sqrt{s_{n,\infty }^2-2(1-2x(1-t_n)(1-u_n))s_{n,\infty }+1}}{2}\right)^{-\left( 4n-1+ \varphi +2 \lambda \right)} \nonumber\\
&&\times \widetilde{w}_{n,n}^{-(2(n-1)+\alpha +\lambda )}\left(  \widetilde{w}_{n,n} \partial _{ \widetilde{w}_{n,n}}\right)  \widetilde{w}_{n,n}^{2(n-1)+\alpha +\lambda} \nonumber\\
&&\times \prod_{k=1}^{n-1} \Bigg\{ \int_{0}^{1} dt_{n-k}\;t_{n-k}^{2(n-k)-1+\lambda } \int_{0}^{1} du_{n-k} \;u_{n-k}^{2(n-k-1)+\gamma +\lambda }\nonumber\\
&&\times \left( s_{n-k}^2-2(1-2\widetilde{w}_{n+1-k,n} (1-t_{n-k})(1-u_{n-k}))s_{n-k}+1 \right)^{-\frac{1}{2}}\nonumber\\
&&\times \left( \frac{ 1+s_{n-k} +\sqrt{s_{n-k}^2-2(1-2\widetilde{w}_{n+1-k,n} (1-t_{n-k})(1-u_{n-k}))s_{n-k}+1}}{2}\right)^{-\left( 4(n-k)-1+ \varphi +2 \lambda \right)} \nonumber\\
&&\times \left. \widetilde{w}_{n-k,n}^{-(2(n-k-1)+\alpha +\lambda )}\left(  \widetilde{w}_{n-k,n} \partial _{ \widetilde{w}_{n-k,n}}\right) \widetilde{w}_{n-k,n}^{2(n-k-1)+\alpha +\lambda} \Bigg\} \mathbf{\Upsilon}(\lambda ; s_0;\widetilde{w}_{1,n}) \right\} z^n\label{eq:30048}
\end{eqnarray}
where
\begin{equation}
\begin{cases} 
{ \displaystyle \mathbf{\Upsilon}(\lambda; s_{0,\infty } ;x)= \sum_{q_0 =0}^{\infty } \frac{(\gamma')_{q_0}}{q_0!} s_{0,\infty }^{q_0} \left( c_0 x^{\lambda } \sum_{i_0=0}^{q_0} \frac{(-q_0)_{i_0} \left(q_0+ \varphi +2 \lambda \right)_{i_0}}{(1+\lambda )_{i_0}(\gamma +\lambda )_{i_0}} x^{i_0} \right) }\cr
{ \displaystyle \mathbf{\Upsilon}(\lambda ; s_0;\widetilde{w}_{1,1}) = \sum_{q_0 =0}^{\infty } \frac{(\gamma')_{q_0}}{q_0!} s_0^{q_0}\left(c_0 x^{\lambda} \sum_{i_0=0}^{q_0} \frac{(-q_0)_{i_0} \left(q_0+ \varphi +2 \lambda \right)_{i_0}}{(1+\lambda )_{i_0}(\gamma +\lambda )_{i_0}} \widetilde{w}_{1,1} ^{i_0} \right) }\cr
{ \displaystyle \mathbf{\Upsilon}(\lambda; s_0 ;\widetilde{w}_{1,n}) = \sum_{q_0 =0}^{\infty } \frac{(\gamma')_{q_0}}{q_0!} s_0^{q_0}\left(c_0 x^{\lambda} \sum_{i_0=0}^{q_0} \frac{(-q_0)_{i_0} \left(q_0+ \varphi +2 \lambda \right)_{i_0}}{(1+\lambda )_{i_0}(\gamma +\lambda )_{i_0}} \widetilde{w}_{1,n} ^{i_0} \right)}
\end{cases}\nonumber 
\end{equation}
\end{theorem}
\begin{proof} 
Acting the summation operator $\displaystyle{ \sum_{q_0 =0}^{\infty } \frac{(\gamma')_{q_0}}{q_0!} s_0^{q_0} \prod _{n=1}^{\infty } \left\{ \sum_{ q_n = q_{n-1}}^{\infty } s_n^{q_n }\right\}}$ on the form of integral of the type 2 CHP  $y(x)$,
\begin{eqnarray}
&&\sum_{\alpha _0 =0}^{\infty } \frac{(\gamma')_{q_0}}{q_0!} s_0^{q_0} \prod _{n=1}^{\infty } \left\{ \sum_{ q_n = q_{n-1}}^{\infty } s_n^{q_n }\right\} y(x) \label{eq:30049}\\
&&= \sum_{q_0 =0}^{\infty } \frac{(\gamma')_{q_0}}{q_0!} s_0^{q_0} \prod _{n=1}^{\infty } \left\{ \sum_{ q_n = q_{n-1}}^{\infty } s_n^{q_n }\right\} \Big\{ y_0(x)+y_1(x)+y_2(x)+y_3(x)+\cdots\Big\} \nonumber
\end{eqnarray}
Acting the summation operator $\displaystyle{ \sum_{q_0 =0}^{\infty } \frac{(\gamma')_{q_0}}{q_0!} s_0^{q_0} \prod _{n=1}^{\infty } \left\{ \sum_{ q_n = q_{n-1}}^{\infty } s_n^{q_n }\right\}}$ on (\ref{eq:30031a}),
\begin{eqnarray}
&&\sum_{q_0 =0}^{\infty } \frac{(\gamma')_{q_0}}{q_0!} s_0^{q_0} \prod _{n=1}^{\infty } \left\{ \sum_{ q_n = q_{n-1}}^{\infty } s_n^{q_n }\right\} y_0(x) \nonumber\\
&&= \prod_{l=1}^{\infty } \frac{1}{(1-s_{l,\infty })} \sum_{q_0 =0}^{\infty } \frac{(\gamma')_{q_0}}{q_0!} s_{0,\infty }^{q_0} \left( c_0 x^{\lambda } \sum_{i_0=0}^{q_0}  \frac{(-q_0)_{i_0} \left(q_0+ \varphi +2 \lambda \right)_{i_0}}{(1+\lambda )_{i_0}(\gamma +\lambda )_{i_0}} x^{i_0} \right) \hspace{2cm}\label{eq:30051}
\end{eqnarray}
Acting the summation operator $\displaystyle{ \sum_{q_0 =0}^{\infty } \frac{(\gamma')_{q_0}}{q_0!} s_0^{q_0} \prod _{n=1}^{\infty } \left\{ \sum_{ q_n = q_{n-1}}^{\infty } s_n^{q_n }\right\}}$ on (\ref{eq:30032}),
\begin{eqnarray}
&&\sum_{q_0 =0}^{\infty } \frac{(\gamma')_{q_0}}{q_0!} s_0^{q_0} \prod _{n=1}^{\infty } \left\{ \sum_{ q_n = q_{n-1}}^{\infty } s_n^{q_n }\right\} y_1(x) \nonumber\\
&&= \prod_{l=2}^{\infty } \frac{1}{(1-s_{l,\infty })} \int_{0}^{1} dt_1\;t_1^{1+\lambda } \int_{0}^{1} du_1\;u_1^{\gamma +\lambda}
 \frac{1}{2\pi i} \oint dv_1 \;\frac{1}{v_1} (1-x(1-t_1)(1-u_1)v_1)^{-\left(4+ \varphi +2 \lambda \right)} \nonumber\\
&&\times \sum_{q_1 =q_0}^{\infty }\left( \frac{v_1-1}{v_1} \frac{s_{1,\infty }}{1-x(1-t_1)(1-u_1)v_1}\right)^{q_1}  \overleftrightarrow {w}_{1,1}^{-(\alpha +\lambda )}\left(  \overleftrightarrow {w}_{1,1} \partial _{ \overleftrightarrow {w}_{1,1}}\right) \overleftrightarrow {w}_{1,1}^{\alpha +\lambda}  \nonumber\\
&&\times  \sum_{q_0 =0}^{\infty } \frac{(\gamma' )_{q_0}}{q_0!}s_0^{q_0}\left( c_0 x^{\lambda } \sum_{i_0=0}^{q_0} \frac{(-q_0)_{i_0} \left(q_0+ \varphi +2 \lambda \right)_{i_0}}{(1+\lambda )_{i_0}(\gamma +\lambda )_{i_0}} \overleftrightarrow {w}_{1,1} ^{i_0} \right) z \label{eq:30052}
\end{eqnarray}
Replace $q_i$, $q_j$ and $r_i$ by $q_1$, $q_0$ and ${ \displaystyle \frac{v_1-1}{v_1} \frac{s_{1,\infty }}{1-x(1-t_1)(1-u_1)v_1}}$ in (\ref{eq:30047}). Take the new (\ref{eq:30047}) into (\ref{eq:30052}).
\begin{eqnarray}
&&\sum_{q_0 =0}^{\infty } \frac{(\gamma')_{q_0}}{q_0!} s_0^{q_0} \prod _{n=1}^{\infty } \left\{ \sum_{ q_n = q_{n-1}}^{\infty } s_n^{q_n }\right\} y_1(x) \nonumber\\
&&= \prod_{l=2}^{\infty } \frac{1}{(1-s_{l,\infty })} \int_{0}^{1} dt_1\;t_1^{1+\lambda } \int_{0}^{1} du_1\;u_1^{\gamma +\lambda}
 \frac{1}{2\pi i} \oint dv_1 \;\frac{(1-x(1-t_1)(1-u_1)v_1)^{-\left(3+ \varphi +2 \lambda \right)} }{-x(1-t_1)(1-u_1)v_1^2+ (1-s_{1,\infty })v_1+s_{1,\infty } } \nonumber\\
&&\times \overleftrightarrow {w}_{1,1}^{-(\alpha +\lambda )}\left(  \overleftrightarrow {w}_{1,1} \partial _{ \overleftrightarrow {w}_{1,1}}\right)  \overleftrightarrow {w}_{1,1}^{\alpha +\lambda}  \sum_{q_0 =0}^{\infty } \frac{(\gamma' )_{q_0}}{q_0!} \left( \frac{v_1-1}{v_1} \frac{s_{0,\infty }}{1-x(1-t_1)(1-u_1)v_1}\right)^{q_0} \label{eq:30053}\\
&&\times \left( c_0 x^{\lambda } \sum_{i_0=0}^{q_0} \frac{(-q_0)_{i_0} \left(q_0+ \varphi +2 \lambda \right)_{i_0}}{(1+\lambda )_{i_0}(\gamma +\lambda )_{i_0}} \overleftrightarrow {w}_{1,1} ^{i_0} \right) z \nonumber
\end{eqnarray}
By using Cauchy's integral formula, the contour integrand has poles at\\
 ${\displaystyle v_1= \frac{1-s_{1,\infty }-\sqrt{(1-s_{1,\infty })^2+4x(1-t_1)(1-u_1)s_{1,\infty }}}{2x(1-t_1)(1-u_1)} }$\\  or ${\displaystyle \frac{1-s_{1,\infty }+\sqrt{(1-s_{1,\infty })^2+4x(1-t_1)(1-u_1)s_{1,\infty }}}{2x(1-t_1)(1-u_1)} }$ and ${ \displaystyle \frac{1-s_{1,\infty }-\sqrt{(1-s_{1,\infty })^2+4x(1-t_1)(1-u_1)s_{1,\infty }}}{2x(1-t_1)(1-u_1)}}$ is only inside the unit circle. As we compute the residue there in (\ref{eq:30053}) we obtain
\begin{eqnarray}
&&\sum_{q_0 =0}^{\infty } \frac{(\gamma')_{q_0}}{q_0!} s_0^{q_0} \prod _{n=1}^{\infty } \left\{ \sum_{ q_n = q_{n-1}}^{\infty } s_n^{q_n }\right\} y_1(x) \nonumber\\
&&= \prod_{l=2}^{\infty } \frac{1}{(1-s_{l,\infty })} \int_{0}^{1} dt_1\;t_1^{1+\lambda} \int_{0}^{1} du_1\;u_1^{\gamma +\lambda}
 \left( s_{1,\infty }^2-2(1-2x(1-t_1)(1-u_1))s_{1,\infty }+1 \right)^{-\frac{1}{2}} \nonumber\\
&&\times \left(\frac{1+s_{1,\infty }+\sqrt{s_{1,\infty }^2-2(1-2x(1-t_1)(1-u_1))s_{1,\infty }+1}}{2}\right)^{-\left(3+ \varphi +2 \lambda  \right)} \label{eq:30054}\\
&&\times \widetilde{w}_{1,1}^{-(\alpha +\lambda )}\left(  \widetilde{w}_{1,1} \partial _{ \widetilde{w}_{1,1}}\right) \widetilde{w}_{1,1}^{\alpha +\lambda} \sum_{q_0 =0}^{\infty } \frac{(\gamma' )_{q_0}}{q_0!} s_0^{q_0}\left( c_0 x^{\lambda } \sum_{i_0=0}^{q_0} \frac{(-q_0)_{i_0} \left(q_0+ \varphi +2 \lambda \right)_{i_0}}{(1+\lambda )_{i_0}(\gamma +\lambda )_{i_0}} \widetilde{w}_{1,1} ^{i_0} \right) z \nonumber 
\end{eqnarray}
where
\begin{eqnarray}
\widetilde{w}_{1,1} &=& \frac{v_1}{(v_1-1)}\; \frac{x t_1 u_1}{1- x v_1 (1-t_1)(1-u_1)}\Bigg|_{\Large v_1=\frac{1-s_{1,\infty }-\sqrt{(1-s_{1,\infty })^2+4x(1-t_1)(1-u_1)s_{1,\infty }}}{2x(1-t_1)(1-u_1)}\normalsize}\nonumber\\
&=& \frac{x t_1 u_1 \left\{ 1+ (s_{1,\infty }+2x(1-t_1)(1-u_1) )s_{1,\infty }\right\}}{2(1-x(1-t_1)(1-u_1))^2 s_{1,\infty }}\nonumber\\
&&-\frac{xt_1 u_1(1+s_{1,\infty })\sqrt{s_{1,\infty }^2-2(1-2x(1-t_1)(1-u_1))s_{1,\infty }+1}}{2(1-x(1-t_1)(1-u_1))^2 s_{1,\infty }}\nonumber
\end{eqnarray}
Acting the summation operator $\displaystyle{ \sum_{q_0 =0}^{\infty } \frac{(\gamma')_{q_0}}{q_0!} s_0^{q_0} \prod _{n=1}^{\infty } \left\{ \sum_{ q_n = q_{n-1}}^{\infty } s_n^{q_n }\right\}}$ on (\ref{eq:30034}),
\begin{eqnarray}
&&\sum_{q_0 =0}^{\infty } \frac{(\gamma')_{q_0}}{q_0!} s_0^{q_0} \prod _{n=1}^{\infty } \left\{ \sum_{ q_n = q_{n-1}}^{\infty } s_n^{q_n }\right\} y_2(x) \nonumber\\
&&= \prod_{l=3}^{\infty } \frac{1}{(1-s_{l,\infty })} \int_{0}^{1} dt_2\;t_2^{3+\lambda} \int_{0}^{1} du_2\;u_2^{2+\gamma +\lambda }
 \frac{1}{2\pi i} \oint dv_2 \;\frac{1}{v_2} (1-x(1-t_2)(1-u_2)v_2)^{-\left(8+ \varphi +2 \lambda \right)} \nonumber\\
&&\times \sum_{q_2 =q_1}^{\infty }\left( \frac{v_2-1}{v_2} \frac{s_{2,\infty }}{1-x(1-t_2)(1-u_2)v_2}\right)^{q_2} \overleftrightarrow {w}_{2,2}^{-(2+\alpha +\lambda )}\left(  \overleftrightarrow {w}_{2,2} \partial _{ \overleftrightarrow {w}_{2,2}}\right) \overleftrightarrow {w}_{2,2}^{2+\alpha +\lambda} \nonumber\\
&&\times \int_{0}^{1} dt_1\;t_1^{1+\lambda} \int_{0}^{1} du_1\;u_1^{\gamma +\lambda}
 \frac{1}{2\pi i} \oint dv_1 \;\frac{1}{v_1} (1-\overleftrightarrow {w}_{2,2} (1-t_1)(1-u_1)v_1)^{-\left( 4+ \varphi +2 \lambda  \right)} \nonumber\\
&&\times \sum_{q_1 =q_0}^{\infty }\left( \frac{v_1-1}{v_1} \frac{s_1}{1-\overleftrightarrow {w}_{2,2}(1-t_1)(1-u_1)v_1}\right)^{q_1} \overleftrightarrow {w}_{1,2}^{-(\alpha +\lambda )}\left(  \overleftrightarrow {w}_{1,2} \partial _{ \overleftrightarrow {w}_{1,2}}\right)  \overleftrightarrow {w}_{1,2}^{\alpha +\lambda} \nonumber\\
&&\times  \sum_{q_0 =0}^{\infty } \frac{(\gamma' )_{q_0}}{q_0!} s_0^{q_0}\left( c_0 x^{\lambda } \sum_{i_0=0}^{q_0} \frac{(-q_0)_{i_0} \left(q_0+ \varphi +2 \lambda \right)_{i_0}}{(1+\lambda )_{i_0}(\gamma +\lambda )_{i_0}} \overleftrightarrow {w}_{1,2} ^{i_0} \right) z^2 \label{eq:30055}
\end{eqnarray}
Replace $q_i$, $q_j$ and $r_i$ by $q_2$, $q_1$ and ${ \displaystyle \frac{v_2-1}{v_2} \frac{s_{2,\infty }}{1-x(1-t_2)(1-u_2)v_2}}$ in (\ref{eq:30047}). Take the new (\ref{eq:30047}) into (\ref{eq:30055}).
\begin{eqnarray}
&&\sum_{q_0 =0}^{\infty } \frac{(\gamma')_{q_0}}{q_0!} s_0^{q_0} \prod _{n=1}^{\infty } \left\{ \sum_{ q_n = q_{n-1}}^{\infty } s_n^{q_n }\right\} y_2(x) \nonumber\\
&&= \prod_{l=3}^{\infty } \frac{1}{(1-s_{l,\infty })} \int_{0}^{1} dt_2\;t_2^{3+\lambda } \int_{0}^{1} du_2\;u_2^{2+\gamma +\lambda }
 \frac{1}{2\pi i} \oint dv_2 \;\frac{\left(1-x (1-t_2)(1-u_2)v_2\right)^{-\left(7+ \varphi +2 \lambda \right)}}{-x (1-t_2)(1-u_2)v_2^2+ (1-s_{2,\infty })v_2+s_{2,\infty } } \nonumber\\
&&\times \overleftrightarrow {w}_{2,2}^{-(2+\alpha +\lambda )}\left(  \overleftrightarrow {w}_{2,2} \partial _{ \overleftrightarrow {w}_{2,2}}\right) \overleftrightarrow {w}_{2,2}^{2+\alpha +\lambda} \nonumber\\
&&\times \int_{0}^{1} dt_1\;t_1^{1+\lambda} \int_{0}^{1} du_1\;u_1^{\gamma +\lambda}
 \frac{1}{2\pi i} \oint dv_1 \;\frac{1}{v_1} \left( 1-\overleftrightarrow {w}_{2,2} (1-t_1)(1-u_1)v_1\right)^{-\left( 4+ \varphi +2 \lambda  \right)} \nonumber\\
&&\times \sum_{q_1 =q_0}^{\infty }\left( \frac{v_2-1}{v_2} \frac{s_{1,\infty }}{1-x(1-t_2)(1-u_2)v_2} \frac{v_1-1}{v_1}\frac{1}{1-\overleftrightarrow {w}_{2,2}(1-t_1)(1-u_1)v_1}\right)^{q_1} \nonumber\\
&&\times \overleftrightarrow {w}_{1,2}^{-(\alpha +\lambda )}\left(  \overleftrightarrow {w}_{1,2} \partial _{ \overleftrightarrow {w}_{1,2}}\right)  \overleftrightarrow {w}_{1,2}^{\alpha +\lambda} \nonumber\\
&&\times  \sum_{q_0 =0}^{\infty } \frac{(\gamma' )_{q_0}}{q_0!} s_0^{q_0}\left( c_0 x^{\lambda } \sum_{i_0=0}^{q_0} \frac{(-q_0)_{i_0} \left(q_0 + \varphi +2 \lambda \right)_{i_0}}{(1+\lambda )_{i_0}(\gamma +\lambda )_{i_0}} \overleftrightarrow {w}_{1,2} ^{i_0} \right) z^2 \label{eq:30056}
\end{eqnarray}
By using Cauchy's integral formula, the contour integrand has poles at\\
 ${\displaystyle
 v_2= \frac{1-s_{2,\infty }-\sqrt{(1-s_{2,\infty })^2+4x(1-t_2)(1-u_2)s_{2,\infty }}}{2x(1-t_2)(1-u_2)}}$ \\or${\displaystyle\frac{1-s_{2,\infty }+\sqrt{(1-s_{2,\infty })^2+4x(1-t_2)(1-u_2)s_{2,\infty }}}{2x(1-t_2)(1-u_2)} }$ 
and ${ \displaystyle\frac{1-s_{2,\infty }-\sqrt{(1-s_{2,\infty })^2+4x(1-t_2)(1-u_2)s_{2,\infty }}}{2x(1-t_2)(1-u_2)}}$ is only inside the unit circle. As we compute the residue there in (\ref{eq:30056}) we obtain
\begin{eqnarray}
&&\sum_{q_0 =0}^{\infty } \frac{(\gamma')_{q_0}}{q_0!} s_0^{q_0} \prod _{n=1}^{\infty } \left\{ \sum_{ q_n = q_{n-1}}^{\infty } s_n^{q_n }\right\} y_2(x) \nonumber\\
&&= \prod_{l=3}^{\infty } \frac{1}{(1-s_{l,\infty })} \int_{0}^{1} dt_2\;t_2^{3+\lambda} \int_{0}^{1} du_2\;u_2^{2+\gamma +\lambda }
 \left( s_{2,\infty }^2-2(1-2x(1-t_2)(1-u_2))s_{2,\infty }+1\right)^{-\frac{1}{2}}\nonumber\\
&&\times \left(\frac{1+s_{2,\infty }+\sqrt{s_{2,\infty }^2-2(1-2x(1-t_2)(1-u_2))s_{2,\infty }+1}}{2}\right)^{-\left(7+ \varphi +2 \lambda  \right)}\nonumber\\
&&\times  \widetilde{w}_{2,2}^{-(2+\alpha +\lambda )}\left( \widetilde{w}_{2,2} \partial _{ \widetilde{w}_{2,2}}\right) \widetilde{w}_{2,2}^{2+\alpha +\lambda} \nonumber\\
&&\times \int_{0}^{1} dt_1\;t_1^{1+\lambda} \int_{0}^{1} du_1\;u_1^{\gamma +\lambda}
 \frac{1}{2\pi i} \oint dv_1 \;\frac{1}{v_1} \left( 1-\widetilde{w}_{2,2} (1-t_1)(1-u_1)v_1\right)^{-\left( 4+ \varphi +2 \lambda \right)} \nonumber\\
&&\times \sum_{q_1 =q_0}^{\infty }\left( \frac{v_1-1}{v_1}\frac{s_1}{1-\widetilde{w}_{2,2}(1-t_1)(1-u_1)v_1}\right)^{q_1}
 \ddot{w}_{1,2}^{-(\alpha +\lambda )}\left( \ddot{w}_{1,2} \partial _{ \ddot{w}_{1,2}}\right) \ddot{w}_{1,2}^{\alpha +\lambda} \nonumber\\
&&\times  \sum_{q_0 =0}^{\infty } \frac{(\gamma' )_{q_0}}{q_0!} s_0^{q_0}\left( c_0 x^{\lambda } \sum_{i_0=0}^{q_0} \frac{(-q_0)_{i_0} \left(q_0+ \varphi +2 \lambda \right)_{i_0}}{(1+\lambda )_{i_0}(\gamma +\lambda )_{i_0}} \ddot{w}_{1,2} ^{i_0} \right) z^2 \label{eq:30057}
\end{eqnarray}
where
\begin{eqnarray}
\widetilde{w}_{2,2} &=& \frac{v_2}{(v_2-1)}\; \frac{x t_2 u_2}{1- x v_2 (1-t_2)(1-u_2)}\Bigg|_{\Large v_2=\frac{1-s_{2,\infty }-\sqrt{(1-s_{2,\infty })^2+4x (1-t_2)(1-u_2)s_{2,\infty }}}{2x(1-t_2)(1-u_2)}\normalsize}\nonumber\\
&=& \frac{x t_2 u_2 \left\{ 1+ (s_{2,\infty }+2x(1-t_2)(1-u_2) )s_{2,\infty }\right\}}{2(1-x(1-t_2)(1-u_2))^2 s_{2,\infty }}\nonumber\\
&&- \frac{x t_2 u_2 (1+s_{2,\infty })\sqrt{s_{2,\infty }^2-2(1-2x (1-t_2)(1-u_2))s_{2,\infty }+1}}{2(1-x(1-t_2)(1-u_2))^2 s_{2,\infty }} \nonumber
\end{eqnarray}
and
\begin{equation}
\ddot{w}_{1,2} = \frac{v_1}{(v_1-1)}\; \frac{\widetilde{w}_{2,2} t_1 u_1}{1- \widetilde{w}_{2,2}v_1 (1-t_1)(1-u_1)}\nonumber
\end{equation}
Replace $q_i$, $q_j$ and $r_i$ by $q_1$, $q_0$ and ${ \displaystyle \frac{v_1-1}{v_1}\frac{s_1}{1-\widetilde{w}_{2,2}(1-t_1)(1-u_1)v_1}}$ in (\ref{eq:30047}). Take the new (\ref{eq:30047}) into (\ref{eq:30057}).
\begin{eqnarray}
&&\sum_{q_0 =0}^{\infty } \frac{(\gamma')_{q_0}}{q_0!} s_0^{q_0} \prod _{n=1}^{\infty } \left\{ \sum_{ q_n = q_{n-1}}^{\infty } s_n^{q_n }\right\} y_2(x) \nonumber\\
&&= \prod_{l=3}^{\infty } \frac{1}{(1-s_{l,\infty })} \int_{0}^{1} dt_2\;t_2^{3+\lambda} \int_{0}^{1} du_2\;u_2^{2+\gamma +\lambda }
 \left( s_{2,\infty }^2-2(1-2x(1-t_2)(1-u_2))s_{2,\infty }+1\right)^{-\frac{1}{2}}\nonumber\\
&&\times \left(\frac{1+s_{2,\infty }+\sqrt{s_{2,\infty }^2-2(1-2x(1-t_2)(1-u_2))s_{2,\infty }+1}}{2}\right)^{-\left(7+ \varphi +2 \lambda  \right)} \nonumber\\
&&\times \widetilde{w}_{2,2}^{-(2+\alpha +\lambda )}\left( \widetilde{w}_{2,2} \partial _{ \widetilde{w}_{2,2}}\right) \widetilde{w}_{2,2}^{2+\alpha +\lambda} \nonumber\\
&&\times \int_{0}^{1} dt_1\;t_1^{1+\lambda} \int_{0}^{1} du_1\;u_1^{\gamma +\lambda}
 \frac{1}{2\pi i} \oint dv_1 \;\frac{\left(1-\widetilde{w}_{2,2} (1-t_1)(1-u_1)v_1\right)^{-\left(3+ \varphi +2 \lambda \right)}}{-\widetilde{w}_{2,2} (1-t_1)(1-u_1)v_1^2+(1-s_1)v_1+s_1} \nonumber\\
&&\times  \ddot{w}_{1,2}^{-(\alpha +\lambda )}\left( \ddot{w}_{1,2} \partial _{ \ddot{w}_{1,2}}\right) \ddot{w}_{1,2}^{\alpha +\lambda}  \sum_{q_0 =0}^{\infty } \frac{(\gamma')_{q_0}}{q_0!} \left( \frac{v_1-1}{v_1}\frac{s_{0,1}}{1-\widetilde{w}_{2,2}(1-t_1)(1-u_1)v_1}\right)^{q_0} \nonumber\\
&&\times \left( c_0 x^{\lambda } \sum_{i_0=0}^{q_0} \frac{(-q_0)_{i_0} \left(q_0+ \varphi +2 \lambda  \right)_{i_0}}{(1+\lambda )_{i_0}(\gamma +\lambda )_{i_0}} \ddot{w}_{1,2} ^{i_0} \right) z^2  \label{eq:30058}
\end{eqnarray}
By using Cauchy's integral formula, the contour integrand has poles at\\ ${\displaystyle
 v_1= \frac{1-s_1-\sqrt{(1-s_1)^2+4\widetilde{w}_{2,2} (1-t_1)(1-u_1)s_1}}{2\widetilde{w}_{2,2} (1-t_1)(1-u_1)}}$\\ or ${\displaystyle\frac{1-s_1+\sqrt{(1-s_1)^2+4\widetilde{w}_{2,2} (1-t_1)(1-u_1)s_1}}{2\widetilde{w}_{2,2} (1-t_1)(1-u_1)}}$ 
and ${ \displaystyle \frac{1-s_1-\sqrt{(1-s_1)^2+4\widetilde{w}_{2,2} (1-t_1)(1-u_1)s_1}}{2\widetilde{w}_{2,2} (1-t_1)(1-u_1)}}$ is only inside the unit circle. As we compute the residue there in (\ref{eq:30058}) we obtain
\begin{eqnarray}
&&\sum_{q_0 =0}^{\infty } \frac{(\gamma')_{q_0}}{q_0!} s_0^{q_0} \prod _{n=1}^{\infty } \left\{ \sum_{ q_n = q_{n-1}}^{\infty } s_n^{q_n }\right\} y_2(x) \nonumber\\
&&= \prod_{l=3}^{\infty } \frac{1}{(1-s_{l,\infty })} \int_{0}^{1} dt_2\;t_2^{3+\lambda} \int_{0}^{1} du_2\;u_2^{2+\gamma +\lambda }
\left( s_{2,\infty }^2-2(1-2x(1-t_2)(1-u_2))s_{2,\infty }+1\right)^{-\frac{1}{2}}\nonumber\\
&&\times \left(\frac{1+s_{2,\infty }+\sqrt{s_{2,\infty }^2-2(1-2x(1-t_2)(1-u_2))s_{2,\infty }+1}}{2}\right)^{-\left(7+ \varphi +2 \lambda  \right)}\nonumber\\
&&\times \widetilde{w}_{2,2}^{-(2+\alpha +\lambda )}\left( \widetilde{w}_{2,2} \partial _{ \widetilde{w}_{2,2}}\right) \widetilde{w}_{2,2}^{2+\alpha +\lambda} \nonumber\\
&&\times \int_{0}^{1} dt_1\;t_1^{1+\lambda } \int_{0}^{1} du_1\;u_1^{\gamma +\lambda}
 \left( s_1^2-2(1-2\widetilde{w}_{2,2} (1-t_1)(1-u_1))s_1+1\right)^{-\frac{1}{2}}\nonumber\\
&&\times \left(\frac{1+s_1+\sqrt{s_1^2-2(1-2\widetilde{w}_{2,2}(1-t_1)(1-u_1))s_1+1}}{2}\right)^{-\left(3+ \varphi +2 \lambda \right)} \nonumber\\
&&\times \widetilde{w}_{1,2}^{-(\alpha +\lambda )}\left( \widetilde{w}_{1,2} \partial _{ \widetilde{w}_{1,2}}\right) \widetilde{w}_{1,2}^{\alpha +\lambda}  \nonumber\\
&&\times \sum_{q_0 =0}^{\infty } \frac{(\gamma' )_{q_0}}{q_0!} s_0^{q_0}\left( c_0 x^{\lambda } \sum_{i_0=0}^{q_0} \frac{(-q_0)_{i_0} \left(q_0+ \varphi +2 \lambda \right)_{i_0}}{(1+\lambda )_{i_0}(\gamma +\lambda )_{i_0}} \widetilde{w}_{1,2} ^{i_0} \right) z^2  \label{eq:30059}
\end{eqnarray}
where
\begin{eqnarray}
\widetilde{w}_{1,2} &=& \frac{v_1}{(v_1-1)}\; \frac{\widetilde{w}_{2,2} t_1 u_1}{1- \widetilde{w}_{2,2} v_1 (1-t_1)(1-u_1)}\Bigg|_{\Large v_1=\frac{1-s_1-\sqrt{(1-s_1)^2+4\widetilde{w}_{2,2} (1-t_1)(1-u_1)s_1}}{2\widetilde{w}_{2,2} (1-t_1)(1-u_1)}\normalsize}\nonumber\\
&=& \frac{\widetilde{w}_{2,2} t_1 u_1 \left\{ 1+ (s_1+2\widetilde{w}_{2,2}(1-t_1)(1-u_1) )s_1\right\}}{2(1-\widetilde{w}_{2,2}(1-t_1)(1-u_1))^2 s_1}\nonumber\\
&&- \frac{\widetilde{w}_{2,2} t_1 u_1 (1+s_1)\sqrt{s_1^2-2(1-2\widetilde{w}_{2,2} (1-t_1)(1-u_1))s_1+1}}{2(1-\widetilde{w}_{2,2}(1-t_1)(1-u_1))^2 s_1} \nonumber
\end{eqnarray}
Acting the summation operator $\displaystyle{ \sum_{q_0 =0}^{\infty } \frac{(\gamma')_{q_0}}{q_0!} s_0^{q_0} \prod _{n=1}^{\infty } \left\{ \sum_{ q_n = q_{n-1}}^{\infty } s_n^{q_n }\right\}}$ on (\ref{eq:30035}),
\begin{eqnarray}
&&\sum_{q_0 =0}^{\infty } \frac{(\gamma')_{q_0}}{q_0!} s_0^{q_0} \prod _{n=1}^{\infty } \left\{ \sum_{ q_n = q_{n-1}}^{\infty } s_n^{q_n }\right\} y_3(x) \nonumber\\
&&= \prod_{l=4}^{\infty } \frac{1}{(1-s_{l,\infty })} \int_{0}^{1} dt_3\;t_3^{5+\lambda} \int_{0}^{1} du_3\;u_3^{4+\gamma +\lambda}\left( s_{3,\infty }^2-2(1-2x(1-t_3)(1-u_3))s_{3,\infty }+1\right)^{-\frac{1}{2}}\nonumber\\
&&\times \left(\frac{1+s_{3,\infty }+\sqrt{s_{3,\infty }^2-2(1-2x(1-t_3)(1-u_3))s_{3,\infty }+1}}{2}\right)^{-\left(11+ \varphi +2 \lambda \right)} \nonumber\\
&&\times \widetilde{w}_{3,3}^{-(4+\alpha +\lambda )}\left( \widetilde{w}_{3,3} \partial _{ \widetilde{w}_{3,3}}\right) \widetilde{w}_{3,3}^{4+\alpha +\lambda}  \nonumber\\
&&\times \int_{0}^{1} dt_2\;t_2^{3+\lambda} \int_{0}^{1} du_2\;u_2^{2+\gamma +\lambda }\left(s_2^2-2(1-2\widetilde{w}_{3,3}(1-t_2)(1-u_2))s_2+1\right)^{-\frac{1}{2}}\nonumber\\
&&\times \left(\frac{1+s_2+\sqrt{s_2^2-2(1-2\widetilde{w}_{3,3}(1-t_2)(1-u_2))s_2+1}}{2}\right)^{-\left(7+ \varphi +2 \lambda \right)} \widetilde{w}_{2,3}^{-(2+\alpha +\lambda )}\left( \widetilde{w}_{2,3} \partial _{ \widetilde{w}_{2,3}}\right) \widetilde{w}_{2,3}^{2+\alpha +\lambda} \nonumber\\
&&\times \int_{0}^{1} dt_1\;t_1^{1+\lambda } \int_{0}^{1} du_1\;u_1^{\gamma +\lambda}\left( s_1^2-2(1-2\widetilde{w}_{2,3} (1-t_1)(1-u_1))s_1+1\right)^{-\frac{1}{2}}\nonumber\\
&&\times \left(\frac{1+s_1+\sqrt{s_1^2-2(1-2\widetilde{w}_{2,3}(1-t_1)(1-u_1))s_1+1}}{2}\right)^{-\left(3+ \varphi +2 \lambda \right)}   \widetilde{w}_{1,3}^{-(\alpha +\lambda )}\left( \widetilde{w}_{1,3} \partial _{ \widetilde{w}_{1,3}}\right) \widetilde{w}_{1,3}^{\alpha +\lambda} \nonumber\\
&&\times \sum_{q_0 =0}^{\infty } \frac{(\gamma' )_{q_0}}{q_0!} s_0^{q_0}\left( c_0 x^{\lambda } \sum_{i_0=0}^{q_0} \frac{(-q_0)_{i_0} \left(q_0+ \varphi +2 \lambda \right)_{i_0}}{(1+\lambda )_{i_0}(\gamma +\lambda )_{i_0}} \widetilde{w}_{1,3} ^{i_0} \right) z^3   \label{eq:30060}
\end{eqnarray}

\vspace{1cm}
where
\begin{eqnarray}
\widetilde{w}_{3,3} &=& \frac{v_3}{(v_3-1)}\; \frac{x t_3 u_3}{1- x(1-t_3)(1-u_3)v_3}\Bigg|_{\Large v_3=\frac{1-s_{3,\infty }-\sqrt{(1-s_{3,\infty })^2+4x(1-t_3)(1-u_3)s_{3,\infty }}}{2x(1-t_3)(1-u_3)}\normalsize}\nonumber\\
&=& \frac{x t_3 u_3 \left\{ 1+ (s_{3,\infty }+2x(1-t_3)(1-u_3) )s_{3,\infty } \right\}}{2(1-x(1-t_3)(1-u_3))^2 s_{3,\infty }}\nonumber\\
&&- \frac{x t_3 u_3 (1+s_{3,\infty })\sqrt{s_{3,\infty }^2-2(1-2x (1-t_3)(1-u_3))s_{3,\infty }+1}}{2(1-x(1-t_3)(1-u_3))^2 s_{3,\infty }} \nonumber
\end{eqnarray}
\begin{eqnarray}
\widetilde{w}_{2,3} &=& \frac{v_2}{(v_2-1)}\; \frac{\widetilde{w}_{3,3} t_2 u_2}{1- \widetilde{w}_{3,3} (1-t_2)(1-u_2)v_2 }\Bigg|_{\Large v_2=\frac{1-s_2-\sqrt{(1-s_2)^2+4\widetilde{w}_{3,3} (1-t_2)(1-u_2)s_2}}{2\widetilde{w}_{3,3} (1-t_2)(1-u_2)}\normalsize}\nonumber\\
&=& \frac{\widetilde{w}_{3,3} t_2 u_2 \left\{ 1+ (s_2+2\widetilde{w}_{3,3}(1-t_2)(1-u_2) )s_2\right\}}{2(1-\widetilde{w}_{3,3}(1-t_2)(1-u_2))^2 s_2}\nonumber\\
&&- \frac{\widetilde{w}_{3,3} t_2 u_2 (1+s_2)\sqrt{s_2^2-2(1-2\widetilde{w}_{3,3} (1-t_2)(1-u_2))s_2+1}}{2(1-\widetilde{w}_{3,3}(1-t_2)(1-u_2))^2 s_2} \nonumber
\end{eqnarray}
\begin{eqnarray}
\widetilde{w}_{1,3} &=& \frac{v_1}{(v_1-1)}\; \frac{\widetilde{w}_{2,3} t_1 u_1}{1- \widetilde{w}_{2,3} (1-t_1)(1-u_1)v_1 }\Bigg|_{\Large v_1=\frac{1-s_1-\sqrt{(1-s_1)^2+4\widetilde{w}_{2,3} (1-t_1)(1-u_1)s_1}}{2\widetilde{w}_{2,3} (1-t_1)(1-u_1)}\normalsize}\nonumber\\
&=& \frac{\widetilde{w}_{2,3} t_1 u_1 \left\{ 1+ (s_1+2\widetilde{w}_{2,3}(1-t_1)(1-u_1) )s_1\right\}}{2(1-\widetilde{w}_{2,3}(1-t_1)(1-u_1))^2 s_1}\nonumber\\
&&- \frac{\widetilde{w}_{2,3} t_1 u_1 (1+s_1)\sqrt{s_1^2-2(1-2\widetilde{w}_{2,3} (1-t_1)(1-u_1))s_1+1}}{2(1-\widetilde{w}_{2,3}(1-t_1)(1-u_1))^2 s_1} \nonumber
\end{eqnarray}
By repeating this process for all higher terms of integral forms of sub-summation $y_m(x)$ terms where $m \geq 4$, I obtain every  $\displaystyle{ \sum_{q_0 =0}^{\infty } \frac{(\gamma')_{q_0}}{q_0!} s_0^{q_0} \prod _{n=1}^{\infty } \left\{ \sum_{ q_n = q_{n-1}}^{\infty } s_n^{q_n }\right\}}  y_m(x)$ terms. 
Substitute (\ref{eq:30051}), (\ref{eq:30054}), (\ref{eq:30059}), (\ref{eq:30060}) and including all $\displaystyle{ \sum_{q_0 =0}^{\infty } \frac{(\gamma')_{q_0}}{q_0!} s_0^{q_0} \prod _{n=1}^{\infty } \left\{ \sum_{ q_n = q_{n-1}}^{\infty } s_n^{q_n }\right\}}  y_m(x)$ terms where $m > 3$ into (\ref{eq:30049}). 
\qed
\end{proof}
\begin{remark}
The generating function for the CHP of type 2 of the first kind about $x=0$ as $q = (q_j+2j)(-\beta +\gamma +\delta -1 + q_j+2j ) $ where $j,q_j \in \mathbb{N}_{0}$ is
\begin{eqnarray}
&&\sum_{q_0 =0}^{\infty } \frac{(\gamma)_{q_0}}{q_0!} s_0^{q_0} \prod _{n=1}^{\infty } \left\{ \sum_{ q_n = q_{n-1}}^{\infty } s_n^{q_n }\right\} H_c^{(a)}F_{q_j}^R\left(\alpha, \beta, \gamma, \delta, q= (q_j+2j)(q_j+2j+\varphi )\right.\nonumber\\
&&; \left.\varphi = -\beta +\gamma +\delta -1; z = \beta x^2  \right) \nonumber\\
&&=2^{\varphi -1}\Bigg\{ \prod_{l=1}^{\infty } \frac{1}{(1-s_{l,\infty })}  \mathbf{A}\left( s_{0,\infty } ;x\right) + \Bigg\{ \prod_{l=2}^{\infty } \frac{1}{(1-s_{l,\infty })} \int_{0}^{1} dt_1\;t_1 \int_{0}^{1} du_1\;u_1^{\gamma} \overleftrightarrow {\mathbf{\Gamma}}_1 \left(s_{1,\infty };t_1,u_1,x\right)\nonumber\\
&&\times \widetilde{w}_{1,1}^{-\alpha}\left( \widetilde{w}_{1,1} \partial _{ \widetilde{w}_{1,1}}\right) \widetilde{w}_{1,1}^{\alpha } \mathbf{A}\left( s_{0} ;\widetilde{w}_{1,1}\right) \Bigg\}z \nonumber\\
&&+ \sum_{n=2}^{\infty } \Bigg\{ \prod_{l=n+1}^{\infty } \frac{1}{(1-s_{l,\infty })} \int_{0}^{1} dt_n\;t_n^{2n-1} \int_{0}^{1} du_n\;u_n^{2(n-1)+\gamma } \overleftrightarrow {\mathbf{\Gamma}}_n \left(s_{n,\infty };t_n,u_n,x \right) \nonumber\\
&&\times \widetilde{w}_{n,n}^{-(2(n-1)+\alpha)}\left( \widetilde{w}_{n,n} \partial _{ \widetilde{w}_{n,n}}\right) \widetilde{w}_{n,n}^{2(n-1)+\alpha }  \nonumber\\
&&\times \prod_{k=1}^{n-1} \Bigg\{ \int_{0}^{1} dt_{n-k}\;t_{n-k}^{2(n-k)-1} \int_{0}^{1} du_{n-k} \;u_{n-k}^{2(n-k-1)+\gamma}\overleftrightarrow {\mathbf{\Gamma}}_{n-k} \left(s_{n-k};t_{n-k},u_{n-k},\widetilde{w}_{n-k+1,n} \right)\nonumber\\
&&\times \widetilde{w}_{n-k,n}^{-(2(n-k-1)+\alpha )}\left( \widetilde{w}_{n-k,n} \partial _{ \widetilde{w}_{n-k,n}}\right)  \widetilde{w}_{n-k,n}^{2(n-k-1)+\alpha } \Bigg\} \mathbf{A} \left( s_{0} ;\widetilde{w}_{1,n}\right) \Bigg\} z^n \Bigg\}   \label{eq:30061} 
\end{eqnarray}
where
\begin{equation}
\begin{cases} 
{ \displaystyle \overleftrightarrow {\mathbf{\Gamma}}_1 \left(s_{1,\infty };t_1,u_1,x\right)= \frac{\left( \frac{1+s_{1,\infty }+\sqrt{s_{1,\infty }^2-2(1-2x(1-t_1)(1-u_1))s_{1,\infty }+1}}{2}\right)^{-\left(3+\varphi \right)}}{\sqrt{s_{1,\infty }^2-2(1-2x (1-t_1)(1-u_1))s_{1,\infty }+1}}}\cr
{ \displaystyle  \overleftrightarrow {\mathbf{\Gamma}}_n \left(s_{n,\infty };t_n,u_n,x\right) =\frac{\left( \frac{1+s_{n,\infty }+\sqrt{s_{n,\infty }^2-2(1-2x(1-t_n)(1-u_n))s_{n,\infty }+1}}{2}\right)^{-\left( 4n-1+\varphi \right)}}{\sqrt{ s_{n,\infty }^2-2(1-2x(1-t_n)(1-u_n))s_{n,\infty }+1}}}\cr
{ \displaystyle \overleftrightarrow {\mathbf{\Gamma}}_{n-k} \left(s_{n-k};t_{n-k},u_{n-k},\widetilde{w}_{n-k+1,n} \right)}\cr
{ \displaystyle = \frac{ \left( \frac{1+s_{n-k}+\sqrt{s_{n-k}^2-2(1-2\widetilde{w}_{n-k+1,n} (1-t_{n-k})(1-u_{n-k}))s_{n-k}+1}}{2}\right)^{-\left( 4(n-k)-1 +\varphi \right)}}{\sqrt{ s_{n-k}^2-2(1-2\widetilde{w}_{n-k+1,n} (1-t_{n-k})(1-u_{n-k}))s_{n-k}+1}}}
\end{cases}\nonumber 
\end{equation}
and
\begin{equation}
\begin{cases} 
{ \displaystyle \mathbf{A} \left( s_{0,\infty } ;x\right)= \frac{\left(1+s_{0,\infty }+\sqrt{s_{0,\infty }^2-2(1-2x )s_{0,\infty }+1}\right)^{\gamma -\varphi}}{\left(1- s_{0,\infty }+\sqrt{s_{0,\infty }^2-2(1-2x)s_{0,\infty }+1}\right)^{\gamma -1} \sqrt{s_{0,\infty }^2-2(1-2x)s_{0,\infty }+1}}}\cr
{ \displaystyle  \mathbf{A} \left( s_{0} ;\widetilde{w}_{1,1}\right) = \frac{\left(1+s_0+\sqrt{s_0^2-2(1-2\widetilde{w}_{1,1} )s_0+1}\right)^{\gamma -\varphi}}{\left(1- s_0+\sqrt{s_0^2-2(1-2\widetilde{w}_{1,1})s_0+1}\right)^{\gamma -1} \sqrt{s_0^2-2(1-2\widetilde{w}_{1,1})s_0+1}}} \cr
{ \displaystyle \mathbf{A} \left( s_{0} ;\widetilde{w}_{1,n}\right) = \frac{\left(1+s_0+\sqrt{s_0^2-2(1-2\widetilde{w}_{1,n} )s_0+1}\right)^{\gamma -\varphi}}{\left(1- s_0+\sqrt{s_0^2-2(1-2\widetilde{w}_{1,n})s_0+1}\right)^{\gamma -1} \sqrt{s_0^2-2(1-2\widetilde{w}_{1,n})s_0+1}}}
\end{cases}\nonumber 
\end{equation}
\end{remark}
\begin{proof}
Replace A and $w$ by $\varphi $ and $s_{0,\infty } $ in (\ref{eq:30041}). 
\begin{eqnarray}
&&\sum_{q_0=0}^{\infty }\frac{(\gamma )_{q_0}}{q_0!} s_{0,\infty }^{q_0} \;_2F_1\left(-q_0, q_0+\varphi; \gamma ; x\right) \label{eq:30062}\\
&&= 2^{\varphi -1}\frac{\left(1+s_{0,\infty }+\sqrt{s_{0,\infty }^2-2(1-2x)s_{0,\infty }+1}\right)^{\gamma - \varphi }}{\left(1- s_{0,\infty }+\sqrt{s_{0,\infty }^2-2(1-2x)s_{0,\infty }+1}\right)^{ \gamma -1} \sqrt{s_{0,\infty }^2-2(1-2x)s_{0,\infty }+1}} \nonumber
\end{eqnarray} 
Replace A, $w$ and $x$  by $ \varphi $, $s_0$ and $\widetilde{w}_{1,1}$ in (\ref{eq:30041}). 
\begin{eqnarray}
&&\sum_{q_0=0}^{\infty }\frac{(\gamma )_{q_0}}{q_0!} s_0^{q_0} \;_2F_1\left(-q_0, q_0+ \varphi; \gamma; \widetilde{w}_{1,1} \right) \label{eq:30063}\\
&&= 2^{\varphi -1}\frac{\left(1+s_0+\sqrt{s_0^2-2(1-2\widetilde{w}_{1,1} )s_0+1}\right)^{\gamma - \varphi }}{\left(1- s_0+\sqrt{s_0^2-2(1-2\widetilde{w}_{1,1})s_0+1}\right)^{ \gamma -1} \sqrt{s_0^2-2(1-2\widetilde{w}_{1,1})s_0+1}} \nonumber
\end{eqnarray} 
Replace A, $w$ and $x$  by $ \varphi $, $s_0$ and $\widetilde{w}_{1,n}$  in (\ref{eq:30041}). 
\begin{eqnarray}
&&\sum_{q_0=0}^{\infty }\frac{(\gamma )_{q_0}}{q_0!} s_0^{q_0} \;_2F_1\left(-q_0, q_0+ \varphi; \gamma; \widetilde{w}_{1,n} \right) \label{eq:30064}\\
&&= 2^{\varphi -1}\frac{ \left(1+s_0+\sqrt{s_0^2-2(1-2\widetilde{w}_{1,n} )s_0+1}\right)^{\gamma - \varphi}}{\left(1- s_0+\sqrt{s_0^2-2(1-2\widetilde{w}_{1,n})s_0+1}\right)^{ \gamma -1} \sqrt{s_0^2-2(1-2\widetilde{w}_{1,n})s_0+1}} \nonumber
\end{eqnarray} 
Put $c_0$= 1, $\lambda $=0 and $\gamma' =\gamma $ in (\ref{eq:30048}). Substitute (\ref{eq:30062}), (\ref{eq:30063}) and (\ref{eq:30064}) into the new (\ref{eq:30048}).\qed
\end{proof}
\begin{remark}
The generating function for the CHP of type 2 of the second kind about $x=0$ as $q =(q_j+2j+1-\gamma )(-\beta +\delta + q_j+2j ) $ where $j,q_j \in \mathbb{N}_{0}$ is
\begin{eqnarray}
&&\sum_{q_0 =0}^{\infty } \frac{(2-\gamma )_{q_0}}{q_0!} s_0^{q_0} \prod _{n=1}^{\infty } \left\{ \sum_{ q_n = q_{n-1}}^{\infty } s_n^{q_n }\right\}  H_c^{(a)}S_{\alpha _j}^R\left(\alpha, \beta, \gamma, \delta \right.\nonumber\\
&&, q= (q_j+2j+1-\gamma )(q_j+2j+1-\gamma +\varphi );\left. \varphi = -\beta +\gamma +\delta -1; z= \beta x^2  \right) \nonumber\\
&&= \frac{2^{ \varphi +1-2\gamma }}{x^{\gamma -1}}\left\{ \prod_{l=1}^{\infty } \frac{1}{(1-s_{l,\infty })} \mathbf{B}\left( s_{0,\infty } ;x\right) \right. \nonumber\\
&&+ \Bigg\{ \prod_{l=2}^{\infty } \frac{1}{(1-s_{l,\infty })} \int_{0}^{1} dt_1\;t_1^{2-\gamma } \int_{0}^{1} du_1\;u_1 \overleftrightarrow {\mathbf{\Psi}}_1 \left(s_{1,\infty };t_1,u_1,x\right)\nonumber\\
&& \times \widetilde{w}_{1,1}^{-(\alpha-\gamma +1)}\left( \widetilde{w}_{1,1} \partial _{ \widetilde{w}_{1,1}}\right) \widetilde{w}_{1,1}^{\alpha -\gamma +1} \mathbf{B}\left( s_{0} ;\widetilde{w}_{1,1}\right) \Bigg\}z \nonumber\\
&&+ \sum_{n=2}^{\infty } \Bigg\{ \prod_{l=n+1}^{\infty } \frac{1}{(1-s_{l,\infty })} \int_{0}^{1} dt_n\;t_n^{2n-\gamma } \int_{0}^{1} du_n\;u_n^{2n-1} \overleftrightarrow {\mathbf{\Psi}}_n \left(s_{n,\infty };t_n,u_n,x \right) \nonumber\\
&&\times \widetilde{w}_{n,n}^{-(2n-1+\alpha-\gamma )}\left( \widetilde{w}_{n,n} \partial _{ \widetilde{w}_{n,n}}\right) \widetilde{w}_{n,n}^{2n-1+\alpha -\gamma } \nonumber\\
&&\times \prod_{k=1}^{n-1} \Bigg\{ \int_{0}^{1} dt_{n-k}\;t_{n-k}^{2(n-k)-\gamma} \int_{0}^{1} du_{n-k} \;u_{n-k}^{2(n-k)-1} \overleftrightarrow {\mathbf{\Psi}}_{n-k} \left( s_{n-k};t_{n-k},u_{n-k},\widetilde{w}_{n-k+1,n} \right) \nonumber\\
&&\times \widetilde{w}_{n-k,n}^{-(2(n-k)-1+\alpha-\gamma )}\left( \widetilde{w}_{n-k,n} \partial _{ \widetilde{w}_{n-k,n}}\right)  \widetilde{w}_{n-k,n}^{2(n-k)-1+\alpha -\gamma } \Bigg\} \left. \mathbf{B}\left( s_{0} ;\widetilde{w}_{1,n}\right)\Bigg\} z^n  \right\}\hspace{1cm}
 \label{eq:30065}
\end{eqnarray}
where
\begin{equation}
\begin{cases} 
{ \displaystyle \overleftrightarrow {\mathbf{\Psi}}_1 \left(s_{1,\infty };t_1,u_1,x\right)= \frac{\left( \frac{1+s_{1,\infty }+\sqrt{s_{1,\infty }^2-2(1-2x(1-t_1)(1-u_1))s_{1,\infty }+1}}{2}\right)^{-\left(5-2\gamma +\varphi \right)}}{\sqrt{s_{1,\infty }^2-2(1-2x(1-t_1)(1-u_1))s_{1,\infty }+1}} }\cr
{ \displaystyle  \overleftrightarrow {\mathbf{\Psi}}_n \left(s_{n,\infty };t_n,u_n,x \right) = \frac{\left( \frac{1+s_{n,\infty }+\sqrt{s_{n,\infty }^2-2(1-2x(1-t_n)(1-u_n))s_{n,\infty }+1}}{2}\right)^{-\left( 4n+1-2\gamma+ \varphi \right)}}{\sqrt{s_{n,\infty }^2-2(1-2x(1-t_n)(1-u_n))s_{n,\infty }+1}}}\cr
{ \displaystyle \overleftrightarrow {\mathbf{\Psi}}_{n-k} \left(s_{n-k};t_{n-k},u_{n-k},\widetilde{w}_{n-k+1,n} \right) }\cr
{ \displaystyle = \frac{\left( \frac{ 1+s_{n-k} +\sqrt{s_{n-k}^2-2(1-2\widetilde{w}_{n-k+1,n} (1-t_{n-k})(1-u_{n-k}))s_{n-k}+1}}{2}\right)^{-\left( 4(n-k)+1-2\gamma +\varphi \right)}}{\sqrt{s_{n-k}^2-2(1-2\widetilde{w}_{n-k+1,n} (1-t_{n-k})(1-u_{n-k}))s_{n-k}+1}}}
\end{cases}\nonumber 
\end{equation}
and
\begin{equation}
\begin{cases} 
{ \displaystyle \mathbf{B} \left( s_{0,\infty } ;x\right)= \frac{\left(1+s_{0,\infty }+\sqrt{s_{0,\infty }^2-2(1-2x )s_{0,\infty }+1}\right)^{-2+3\gamma -\varphi }}{\left(1- s_{0,\infty }+\sqrt{s_{0,\infty }^2-2(1-2x)s_{0,\infty }+1}\right)^{1-\gamma } \sqrt{s_{0,\infty }^2-2(1-2x)s_{0,\infty }+1}}}\cr
{ \displaystyle  \mathbf{B} \left( s_{0} ;\widetilde{w}_{1,1}\right) = \frac{\left( 1+s_0+\sqrt{s_0^2-2(1-2\widetilde{w}_{1,1} )s_0+1}\right)^{-2+3\gamma -\varphi }}{\left(1- s_0+\sqrt{s_0^2-2(1-2\widetilde{w}_{1,1})s_0+1}\right)^{1-\gamma} \sqrt{s_0^2-2(1-2\widetilde{w}_{1,1})s_0+1}}} \cr
{ \displaystyle \mathbf{B} \left( s_{0} ;\widetilde{w}_{1,n}\right) = \frac{\left(1+s_0+\sqrt{s_0^2-2(1-2\widetilde{w}_{1,n} )s_0+1}\right)^{-2+3\gamma -\varphi }}{\left(1- s_0+\sqrt{s_0^2-2(1-2\widetilde{w}_{1,n})s_0+1}\right)^{1-\gamma } \sqrt{s_0^2-2(1-2\widetilde{w}_{1,n})s_0+1}}}
\end{cases}\nonumber 
\end{equation}
\end{remark}
\begin{proof}
Replace $\gamma $, A and $w$ by $2-\gamma $, $ \varphi +2(1-\gamma )$ and $s_{0,\infty }$ in (\ref{eq:30041}). 
\begin{eqnarray}
&&\sum_{q_0=0}^{\infty }\frac{(2-\gamma)_{q_0}}{q_0!} s_{0,\infty }^{q_0} \;_2F_1\left(-q_0, q_0+ \varphi +2(1-\gamma ); 2-\gamma; x\right) \label{eq:30066}\\
&&=2^{ \varphi +1-2\gamma } \frac{ \left(1+s_{0,\infty }+\sqrt{s_{0,\infty }^2-2(1-2x)s_{0,\infty }+1}\right)^{-2+3\gamma -\varphi }}{\left(1- s_{0,\infty }+\sqrt{s_{0,\infty }^2-2(1-2x)s_{0,\infty }+1}\right)^{1-\gamma }\sqrt{s_{0,\infty }^2-2(1-2x )s_{0,\infty }+1}} \nonumber
\end{eqnarray} 
Replace $\gamma $, A, $w$ and $x$  by $2-\gamma $, $\varphi +2(1-\gamma )$, $s_0$ and $\widetilde{w}_{1,1}$ in (\ref{eq:30041}). 
\begin{eqnarray}
&&\sum_{q_0=0}^{\infty }\frac{(2-\gamma)_{q_0}}{q_0!} s_0^{q_0} \;_2F_1\left(-q_0, q_0+ \varphi +2(1-\gamma ); 2-\gamma; \widetilde{w}_{1,1} \right) \label{eq:30067}\\
&&=2^{ \varphi +1-2\gamma } \frac{\left(1+s_0+\sqrt{s_0^2-2(1-2\widetilde{w}_{1,1} )s_0+1}\right)^{-2+3\gamma -\varphi }}{\left(1- s_0+\sqrt{s_0^2-2(1-2\widetilde{w}_{1,1})s_0+1}\right)^{1-\gamma } \sqrt{s_0^2-2(1-2\widetilde{w}_{1,1})s_0+1}} \nonumber
\end{eqnarray} 
Replace $\gamma $, A, $w$ and $x$  by $2-\gamma $, $\varphi +2(1-\gamma )$, $s_0$ and $\widetilde{w}_{1,n}$ in (\ref{eq:30041}). 
\begin{eqnarray}
&&\sum_{q_0=0}^{\infty }\frac{(2-\gamma)_{q_0}}{q_0!} s_0^{q_0} \;_2F_1\left(-q_0, q_0+ \varphi +2(1-\gamma ); 2-\gamma; \widetilde{w}_{1,n} \right) \label{eq:30068}\\
&&= 2^{ \varphi +1-2\gamma }\frac{\left(1+s_0+\sqrt{s_0^2-2(1-2\widetilde{w}_{1,n} )s_0+1}\right)^{-2+3\gamma -\varphi }}{\left(1- s_0+\sqrt{s_0^2-2(1-2\widetilde{w}_{1,n})s_0+1}\right)^{1-\gamma } \sqrt{s_0^2-2(1-2\widetilde{w}_{1,n})s_0+1}} \nonumber
\end{eqnarray} 
Put $ c_0= 1 $, $\lambda =1-\gamma $ and $\gamma' =2-\gamma $ in (\ref{eq:30048}). Substitute (\ref{eq:30066}), (\ref{eq:30067}) and (\ref{eq:30068}) into the new (\ref{eq:30048}).\qed
\end{proof}
\section{Summary}
The second independent solutions of the CHF by using 3TRF and R3TRF are only available for non-integer value of $\gamma $ in the domain $|x|<1$. If $\gamma $ is an integer, one of two solutions will be dependent. Its analytic solution is defined by its behavior at infinity in the domain $|x|>1$.
All possible local solutions of the CHE (Regge-Wheeler and Teukolsky equations) were constructed by Fiziev.\cite{Fizi2009,Fizi2010}
The power series representation of the CHE about the singular point at infinity is not convergent by only asymptotic \cite{Ronv1995,Fizi2010a}: Its solution will be constructed in the future paper for (1) the power series expansion in closed forms of the CHF about regular singularity $x=\infty $, (2) its integral representation and (3) the generating functions for type 3 CHP.

In chapter 4 I show how to obtain the power series expansion in closed forms and its integral forms of the CHE for infinite series and polynomial of type 1 including all higher terms of $A_n$'s by applying 3TRF. This was done by letting $A_n$ in sequence $c_n$ is the leading term in the analytic function $y(x)$: the sequence $c_n$ consists of combinations $A_n$ and $B_n$. For polynomial of type 1, I treat $\beta $, $\gamma $, $\delta $ and $q$ as free variables and a fixed value of $\alpha $.

In this chapter I show how to construct the power series expansion in closed forms and its integral forms of the CHE for  infinite series and polynomial of type 2 including all higher terms of $B_n$'s by applying R3TRF. This is done by letting $B_n$ in sequence $c_n$ is the leading term in the analytic function $y(x)$. For polynomial of type 2, I treat $\alpha $, $\beta $, $\gamma $ and $\delta $ as free variables and a fixed value of $q$. 

The power series expansion and integral forms of the CHE for infinite series about $x=0$ in this chapter are equivalent to infinite series of the CHE in chapter 4. In this chapter $B_n$ is the leading term in sequence $c_n$ in the analytic function $y(x)$. In chapter 4 $A_n$ is the leading term in sequence $c_n$ in the analytic function $y(x)$.
 
In chapter 4 as we see the power series expansions of the CHE about $x=0$ for infinite series and polynomial of type 1, the denominators and numerators in all $B_n$ terms of $y(x)$ arise with Pochhammer symbol. Again in this chapter the denominators and numerators in all $A_n$ terms of the CHF about $x=0$ arise with Pochhammer symbol. Since we construct the power series expansion with Pochhammer symbols in numerators and denominators, we are able to describe integral forms of the CHF in the mathematical rigor. As we observe representations in closed form integrals of the CHF by applying 3TRF, a $_1F_1$ function recurs in each of sub-integral forms of the CHF. Also $_2F_1$ function recurs in each of sub-integral forms of the CHF as we observe integral form of the CHF by applying R3TRF.

We can transform the CHF into any special functions having two recursive coefficients because of $_1F_1$ and $_2F_1$ functions in each of sub-integral forms of the CHF. After we replace $_1F_1$ and $_2F_1$ functions in integral forms of the CHF to other special functions, we are able to rebuild the power series expansion of the CHF in a backward.   

In chapter 4 I show how to analyze the generating function for the CHP of type 1 from its general integral representation. And in this chapter I construct the generating function for the CHP of type 2 from its integral form. We are able to derive orthogonal relations, recursion relations and expectation values of physical quantities from these two generating functions compared with the case of a normalized wave function for the hydrogen-like atoms.\footnote{The normalized constant for this wave function is obtained by applying the generating function for associated Laguerre polynomials into the Helmholtz equation in the spherical coordinates. It only has one eigenvalue. And the CHP of either type 1 or 2 has infinite eigenvalues because of a 3-term recursive relation between successive coefficients in its power series.}   

Mathematical structure of the generating function for the CHP of type 2 closely resembles the generating function for the Heun polynomial of type 2 in chapter 3. If $a\rightarrow 0$, $\eta \rightarrow x $, $\alpha =\beta $, $z= -\frac{1}{a}x^2 \rightarrow \beta x^2$, $\varphi = \alpha +\beta -\delta +a(\delta +\gamma -1)\rightarrow -\beta +\gamma +\delta -1$, $q= -(q_j+2j+\lambda )\{\varphi +(1+a)(q_j+2j+\lambda ) \} \rightarrow (q_j+2j+\lambda )(\varphi + q_j+2j+\lambda ) $, $\left(  \widetilde{w}_{1,1} \partial _{ \widetilde{w}_{1,1}}\right)^2 \rightarrow \widetilde{w}_{1,1} \partial _{ \widetilde{w}_{1,1}}$, $\left(  \widetilde{w}_{n,n} \partial _{ \widetilde{w}_{n,n}}\right)^2 \rightarrow \widetilde{w}_{n,n} \partial _{ \widetilde{w}_{n,n}}$ and  $\left(  \widetilde{w}_{n-k,n} \partial _{ \widetilde{w}_{n-k,n}}\right)^2 \rightarrow \widetilde{w}_{n-k,n} \partial _{ \widetilde{w}_{n-k,n}}$ in the general expression of the generating function for the Heun polynomial of type 2, its solution is equivalent to the generating function for the CHP of type 2; compare theorem 3.2.4 in chapter 3 with theorem 5.2.12 in this chapter. Because I derive the both generating functions for the CHP and Heun polynomial of type 2 by applying the generating function for the Jacobi polynomial using $_2F_1$ functions into integral forms of the CHP and Heun polynomial: $_2F_1$ function recurs in each of sub-integral forms of them. 
By similar reason, mathematical structure of the generating function for the CHP of type 1 closely resembles the generating function for the GCH polynomial of type 1 \cite{Chou2012j}. The generating function for confluent hypergeometric polynomial is applied into the integral representations of the CHP and GCH polynomial: $_1F_1$ function recurs in each of sub-integral forms of them.
\addcontentsline{toc}{section}{Bibliography}
\bibliographystyle{model1a-num-names}
\bibliography{<your-bib-database>}
\bibliographystyle{model1a-num-names}
\bibliography{<your-bib-database>}

\chapter{Grand Confluent Hypergeometric function using reversible three-term recurrence formula}
\chaptermark{GCH function using  R3TRF} 

Biconfluent Heun (BCH) equation, a confluent form of Heun equation\cite{Heun1889,Ronv1995}, is the special case of Grand Confluent Hypergeometric (GCH) equation\footnote{For the canonical form of the BCH equation\cite{Ronv1995}, replace $\mu $, $\varepsilon $, $\nu $, $\Omega $ and $\omega $ by $-2$, $-\beta  $, $ 1+\alpha $, $\gamma -\alpha -2 $ and $ 1/2 (\delta /\beta +1+\alpha )$ in (\ref{eq:5001}). For DLFM version \cite{NIST} or in Ref.\cite{Slavy2000}, replace $\mu $ and $\omega $ by 1 and $-q/\varepsilon $ in (\ref{eq:5001}).}: this has a regular singularity at $x=0$, and an irregular singularity at $\infty$ of rank 2.  

In this chapter I will apply reversible three term recurrence formula (R3TRF) in chapter 1 to (1) the power series expansion in closed forms, (2) its integral forms of the GCH equation for infinite series and polynomial which makes $A_n$ term terminated  including all higher terms of $B_n$'s\footnote{`` higher terms of $B_n$'s'' means at least two terms of $B_n$'s.}  and (3) the generating function for the GCH polynomial which makes $A_n$ term terminated. 

\section{Introduction}
\begin{equation}
x \frac{d^2{y}}{d{x}^2} + \left( \mu x^2 + \varepsilon x + \nu  \right) \frac{d{y}}{d{x}} + \left( \Omega x + \varepsilon \omega \right) y = 0
\label{eq:5001}
\end{equation}
(\ref{eq:5001}) is  Grand Confluent Hypergeometric (GCH) differential equation where $\mu$, $\varepsilon$, $\nu $, $\Omega$ and $\omega$ are real or complex parameters \cite{Chou2012a,Chou2012i}.
GCH ordinary differential equation is of Fuchsian types with two singular points: one regular singular point which is zero with exponents $\{ 0,1-\nu \}$, and one irregular singular point which is infinity with an exponent $\frac{\Omega }{\mu }$. In contrast, Heun equation of Fuchsian types has the four regular singularities. Heun equation has the four kind of confluent forms: (1) Confluent Heun (two regular and one irregular singularities), (2) Doubly confluent Heun (two irregular singularities), (3) Biconfluent Heun (one regular and one irregular singularities), (4) Triconfluent Heun equations (one irregular singularity). 

BCH equation is derived from the GCH equation by changing all coefficients.\cite{Ronv1995} The GCH (or BCH) equation is applicable into the modern physics. \cite{Slav1996,Ralk2002,Kand2005,Hortacsu:2011rr,Arri1991} The BCH equation appears in the radial Schr$\ddot{\mbox{o}}$dinger equation with those potentials such as the rotating harmonic oscillator \cite{Mass1983}, the doubly anharmonic oscillator \cite{Chau1983,Chau1984,Leau1986}, a three-dimensional anharmonic oscillator \cite{Fles1980,Fles1982,Leau1986}, Coulomb potential with a linear confining potential \cite{Leau1986,Pons1988} and other kinds of potentials \cite{Leau1990,Lemi1969}.  

The fundamental solution of the BCH equation for infinite series and the BCH spectral polynomial about $x=0$ in the canonical form was obtained by applying the power series expansion.\cite{Deca1978,Bato1977,Haut1969,Urwi1975} For the case of the irregular singular point $x=\infty $, the three term recurrence of the power series in the BCH equation was derived.\cite{Maro1967,Meix1933} And they left the analytic solution of the BCH equation as solutions of recurrences because of a 3-term recursive relation between successive coefficients in its power series expansion of the BCH equation.\footnote{For the special case, the explicit solutions of the BCH equation in the canonical form was constructed when one of coefficients $\beta =0$.\cite{Exto1989}} In comparison with the two term recursion relation of the power series in the linear differential equation, the analytic solution in closed forms on the three term recurrence relation of the power series is unknown currently because of its complex mathematical calculation.  

As I know, there are no examples for the analytic solution of the BCH equation about $x=0$ and $x=\infty $ in the form of a definite or contour integral containing the well-known special functions which consists of two term recursion relation in its power series of the linear differential equation such as a $_2F_1$ or $_1F_1$ function.
In place of describing the integral representation of the BCH equation involving only simple functions, especially for confluent hypergeometric functions, the BCH equation is obtained by Fredholm-type integral equations; such integral relationships express one analytic solution in terms of another analytic solution.\cite{Batola1982,Batola.a.1982,Maroni1979,Maroni.a.1979,Maroni.b.1984,Chiang2013,Belm2004} 

In Ref.\cite{Chou2012i,Chou2012j} I construct analytic solutions of GCH equation about the regular singular point at $x=0$ by applying three term recurrence formula (3TRF). \cite{Chou2012b}; (1) the power series expansion in closed forms, (2) its integral forms of the GCH equation (for infinite series and polynomial which makes $B_n$ term terminated including all higher terms of $A_n$'s\footnote{`` higher terms of $A_n$'s'' means at least two terms of $A_n$'s.}) and (3) the generating function for GCH polynomial  which makes $B_n$ term terminated. 
And I show four examples of the analytic wave functions and its eigenvalues in the radial Schr$\ddot{\mbox{o}}$dinger equation with certain potentials: (1) Schr$\ddot{\mbox{o}}$dinger equation with the rotating harmonic oscillator and a class of confinement potentials, (2) The spin free Hamiltonian involving only scalar potential for the $q-\bar{q}$ system, (3) The radial Schr$\ddot{\mbox{o}}$dinger equation with Confinement potentials, (4) Two interacting electrons in a uniform magnetic field and a parabolic potential. The Frobenius solutions in closed forms and its combined definite and contour integrals of these four quantum mechanical wave functions are derived analytically.

In this chapter, by applying R3TRF in chapter 1, I construct the power series expansion in closed forms of GCH equation about the regular singular point at $x=0$ for infinite series and polynomial which makes $A_n$ term terminated. The integral forms of GCH equation and its generating function for the GCH polynomial which makes $A_n$ term terminated are derived analytically. Also the Frobenius solution of GCH equation about irregular singular point at $x=\infty $ by applying 3TRF \cite{Chou2012b} is constructed analytically including its integral representation and the generating function for the GCH polynomial polynomial which makes $B_n$ term terminated.

\section{GCH equation about regular singular point at zero}

Assume that its solution is
\begin{equation}
y(x)= \sum_{n=0}^{\infty } c_n x^{n+\lambda }  \label{eq:5002}
\end{equation}
where $\lambda$ is an indicial root.  Plug (\ref{eq:5002})  into (\ref{eq:5001}). We obtain a three-term recurrence relation for the coefficients $c_n$:
\begin{equation}
c_{n+1}=A_n \;c_n +B_n \;c_{n-1} \hspace{1cm};n\geq 1 \label{eq:5003}
\end{equation}
where,
\begin{subequations}
\begin{equation}
A_n = \frac{-\varepsilon (n+\omega +\lambda ) }{ (n+1+\lambda )(n+\nu +\lambda)}\nonumber
\label{eq:5004a}
\end{equation}
\begin{equation}
B_n = -\frac{\Omega +\mu (n-1+\lambda ) }{ (n+1+\lambda )(n+\nu +\lambda)} \label{eq:5004b}
\end{equation}
\begin{equation}
c_1= A_0 \;c_0 \label{eq:5004c}
\end{equation}
\end{subequations} 
We have two indicial roots which are $\lambda = 0$ and $  1-\nu $.

\subsection{Power series}
\subsubsection{Polynomial of type 2}
There are two types of power series expansion using the two term recurrence relation in a linear ordinary differential equation which are a polynomial and an infinite series. In contrast there are three types of polynomials and an infinite series in three term recurrence relation of linear ordinary differential equation: (1) polynomial which makes $B_n$ term terminated: $A_n$ term is not terminated, (2) polynomial which makes $A_n$ term terminated: $B_n$ term is not terminated, (3) polynomial which makes $A_n$ and $B_n$ terms terminated at the same time.\footnote{If $A_n$ and $B_n$ terms are not terminated, it turns to be infinite series.} In general the GCH polynomial is defined as type 3 polynomial where $A_n$ and $B_n$ terms terminated. For the type 3 GCH  polynomial about $x=0$, it has a fixed integer value of $\Omega $, just as it has a fixed value of $\omega $. In three term recurrence relation, polynomial of type 3 I categorize as complete polynomial. In the previous paper \cite{Chou2012i,Chou2012j} I construct the type 1 GCH polynomial about $x=0$. In future papers I will derive type 3 polynomial about $x=0$. In this chapter I construct the power series expansion, its integral forms and the generating function for the GCH polynomial of type 2 about $x=0$: I treat $\mu $, $\varepsilon  $, $\nu $ and $\Omega $ as free variables and $\omega $ as a fixed value.  

In chapter 1 the general expression of power series of $y(x)$ for polynomial of type 2 is defined by
\begin{eqnarray}
y(x) &=&  \sum_{n=0}^{\infty } y_{n}(x) = y_0(x)+ y_1(x)+ y_2(x)+y_3(x)+\cdots \nonumber\\
&=&  c_0 \Bigg\{ \sum_{i_0=0}^{\alpha _0} \left( \prod _{i_1=0}^{i_0-1}A_{i_1} \right) x^{i_0+\lambda }
+ \sum_{i_0=0}^{\alpha _0}\left\{ B_{i_0+1} \prod _{i_1=0}^{i_0-1}A_{i_1}  \sum_{i_2=i_0}^{\alpha _1} \left( \prod _{i_3=i_0}^{i_2-1}A_{i_3+2} \right)\right\} x^{i_2+2+\lambda }  \nonumber\\
&& + \sum_{N=2}^{\infty } \Bigg\{ \sum_{i_0=0}^{\alpha _0} \Bigg\{B_{i_0+1}\prod _{i_1=0}^{i_0-1} A_{i_1} 
\prod _{k=1}^{N-1} \Bigg( \sum_{i_{2k}= i_{2(k-1)}}^{\alpha _k} B_{i_{2k}+2k+1}\prod _{i_{2k+1}=i_{2(k-1)}}^{i_{2k}-1}A_{i_{2k+1}+2k}\Bigg)\nonumber\\
&& \times  \sum_{i_{2N} = i_{2(N-1)}}^{\alpha _N} \Bigg( \prod _{i_{2N+1}=i_{2(N-1)}}^{i_{2N}-1} A_{i_{2N+1}+2N} \Bigg) \Bigg\} \Bigg\} x^{i_{2N}+2N+\lambda }\Bigg\}  \label{eq:5005}
\end{eqnarray}
In the above, $\alpha _i\leq \alpha _j$ only if $i\leq j$ where $i,j,\alpha _i, \alpha _j \in \mathbb{N}_{0}$.

For a polynomial, we need a condition which is:
\begin{equation}
A_{\alpha _i+ 2i}=0 \hspace{1cm} \mathrm{where}\;i,\alpha _i =0,1,2,\cdots
\label{eq:5006}
\end{equation}
In the above, $ \alpha _i$ is an eigenvalue that makes $A_n$ term terminated at certain value of index $n$. (\ref{eq:5006}) makes each $y_i(x)$ where $i=0,1,2,\cdots$ as the polynomial in (\ref{eq:5005}).
In (\ref{eq:5004a})-(\ref{eq:5004c}) replace $\omega $ by ${ \displaystyle  -(\omega_i +2i+\lambda )}$. In (\ref{eq:5006}) replace  index $\alpha _i$ by $\omega _i$.  Take the new (\ref{eq:5004a})-(\ref{eq:5004c}), (\ref{eq:5006}) and put them in (\ref{eq:5005}).
After the replacement process, the general expression of power series of the GCH equation for polynomial of type 2 is given by
\begin{eqnarray}
 y(x)&=& \sum_{n=0}^{\infty } y_{n}(x) = y_0(x)+ y_1(x)+ y_2(x)+y_3(x)+\cdots \nonumber\\ 
&=& c_0 x^{\lambda } \left\{\sum_{i_0=0}^{\omega _0} \frac{(-\omega _0)_{i_0} }{(1+\lambda )_{i_0}(\nu +\lambda )_{i_0}} \eta ^{i_0}\right.\nonumber\\
&&+ \left\{ \sum_{i_0=0}^{\omega _0}\frac{(i_0+ \Omega / \mu + \lambda )}{ (i_0+ 2+\lambda )(i_0+ 1+\nu + \lambda )}\frac{(-\omega _0)_{i_0} }{(1+\lambda )_{i_0}(\nu +\lambda )_{i_0}} \right.   \left. \sum_{i_1=i_0}^{\omega _1} \frac{(-\omega _1)_{i_1} (3+\lambda )_{i_0}(2+\nu +\lambda )_{i_0}}{(-\omega _1)_{i_0} (3+\lambda )_{i_1}(2+\nu +\lambda )_{i_1}} \eta ^{i_1}\right\} \rho \nonumber\\
&&+ \sum_{n=2}^{\infty } \left\{ \sum_{i_0=0}^{\omega _0}\frac{(i_0+\Omega / \mu + \lambda )}{ (i_0+ 2+\lambda )(i_0+ 1+\nu + \lambda )}\frac{(-\omega_0)_{i_0}  }{(1+\lambda )_{i_0}(\nu +\lambda )_{i_0}}\right.\nonumber\\
&&\times \prod _{k=1}^{n-1} \left\{ \sum_{i_k=i_{k-1}}^{\omega_k} \frac{(i_k+ 2k +\Omega / \mu +\lambda )}{(i_k+ 2 k+2+\lambda )(i_k+ 2k+1+\nu +\lambda )}\right.  \left.\frac{(-\omega _k)_{i_k} (2k+1+\lambda )_{i_{k-1}}(2k+\nu +\lambda )_{i_{k-1}}}{(-\omega _k)_{i_{k-1}}(2k+1+\lambda )_{i_k}(2k+\nu +\lambda )_{i_k}}\right\} \nonumber\\
&&\times \left. \left.\sum_{i_n= i_{n-1}}^{\omega _n} \frac{(-\omega_n)_{i_n} (2n+1+\lambda )_{i_{n-1}}(2n+\nu +\lambda )_{i_{n-1}}}{(-\omega _n)_{i_{n-1}} (2n+1+\lambda )_{i_n}(2n+\nu +\lambda )_{i_n}} \eta ^{i_n} \right\} \rho ^n \right\}\label{eq:5007}
\end{eqnarray}
where
\begin{equation}
\begin{cases} \eta  = -\varepsilon  x \cr
\rho = -\mu x^2 \cr
\omega = -(\omega_j +2 j+\lambda ) \;\;\mbox{as}\;j,\omega _j\in \mathbb{N}_{0} \cr
\omega _i\leq \omega _j \;\;\mbox{only}\;\mbox{if}\;i\leq j\;\;\mbox{where}\;i,j\in \mathbb{N}_{0} 
\end{cases}\nonumber 
\end{equation}

Put $c_0$= 1 as $\lambda =0$  for the first kind of independent solutions of the GCH equation and $\lambda =1-\nu$ for the second one in (\ref{eq:5007}).  
\begin{remark}
The power series expansion of the GCH equation of the first kind for polynomial of type 2 about $x=0$ as $\omega = -(\omega_j +2 j) $ where $j,\omega _j \in \mathbb{N}_{0}$ is
\begin{eqnarray}
 y(x)&=& QW_{\omega _j}^R\left(\mu ,\varepsilon ,\nu ,\Omega ,\omega =-(\omega_j +2 j); \rho =-\mu x^2; \eta = -\varepsilon x \right)\nonumber\\
&=& \sum_{i_0=0}^{\omega _0} \frac{(-\omega _0)_{i_0} }{(1 )_{i_0}(\nu )_{i_0}} \eta ^{i_0} + \left\{ \sum_{i_0=0}^{\omega _0}\frac{(i_0+ \Omega / \mu  )}{ (i_0+ 2 )(i_0+ 1+\nu )}\frac{(-\omega _0)_{i_0} }{(1 )_{i_0}(\nu )_{i_0}} \right.   \left. \sum_{i_1=i_0}^{\omega _1} \frac{(-\omega _1)_{i_1} (3 )_{i_0}(2+\nu )_{i_0}}{(-\omega _1)_{i_0} (3 )_{i_1}(2+\nu )_{i_1}} \eta ^{i_1}\right\} \rho \nonumber\\
&&+ \sum_{n=2}^{\infty } \left\{ \sum_{i_0=0}^{\omega _0}\frac{(i_0+\Omega / \mu )}{ (i_0+ 2 )(i_0+ 1+\nu )}\frac{(-\omega_0)_{i_0}  }{(1 )_{i_0}(\nu  )_{i_0}}\right.\nonumber\\
&&\times \prod _{k=1}^{n-1} \left\{ \sum_{i_k=i_{k-1}}^{\omega_k} \frac{(i_k+ 2k +\Omega / \mu )}{(i_k+ 2 k+2 )(i_k+ 2k+1+\nu )}\right.  \left.\frac{(-\omega _k)_{i_k} (2k+1 )_{i_{k-1}}(2k+\nu )_{i_{k-1}}}{(-\omega _k)_{i_{k-1}}(2k+1 )_{i_k}(2k+\nu )_{i_k}}\right\} \nonumber\\
&&\times  \left.\sum_{i_n= i_{n-1}}^{\omega _n} \frac{(-\omega_n)_{i_n} (2n+1 )_{i_{n-1}}(2n+\nu )_{i_{n-1}}}{(-\omega _n)_{i_{n-1}} (2n+1 )_{i_n}(2n+\nu )_{i_n}} \eta ^{i_n} \right\} \rho ^n  \label{eq:5008}
\end{eqnarray}
\end{remark}
For the minimum value of the GCH equation of the first kind for polynomial of type 2 about $x=0$, put $\omega _0=\omega _1=\omega _2=\cdots=0$ in (\ref{eq:5008}).
\begin{eqnarray}
y(x)&=& QW_{0}^R\left(\mu ,\varepsilon ,\nu ,\Omega ,\omega =-2 j; \rho =-\mu x^2; \eta = -\varepsilon x \right)\nonumber\\
&=& \; _1F_1\left( \frac{\Omega }{2\mu }, \frac{\nu }{2}+\frac{1}{2}, -\frac{1}{2}\mu x^2 \right) \;\;\mbox{where}\;-\infty < x< \infty \nonumber 
\end{eqnarray}
On the above,  $_1F_1( a,b,x)= \sum_{n=0}^{\infty }\frac{(a)_n}{(b)_n}\frac{x^n}{n!}$.
\begin{remark}
The power series expansion of the GCH equation of the second kind for polynomial of type 2 about $x=0$ as $ \omega = -(\omega_j +2 j+1-\nu ) $ where $j,\omega _j \in \mathbb{N}_{0}$ is
\begin{eqnarray}
 y(x)&=& RW_{\omega _j}^R\left(\mu ,\varepsilon ,\nu ,\Omega ,\omega =-(\omega_j +2 j+1-\nu); \rho =-\mu x^2; \eta = -\varepsilon x \right)\nonumber\\
&=&  x^{1-\nu } \left\{\sum_{i_0=0}^{\omega _0} \frac{(-\omega _0)_{i_0} }{(2-\nu  )_{i_0}(1)_{i_0}} \eta ^{i_0}\right.\nonumber\\
&&+ \left\{ \sum_{i_0=0}^{\omega _0}\frac{(i_0+ 1+\Omega / \mu -\nu )}{ (i_0+ 3-\nu )(i_0+ 2)}\frac{(-\omega _0)_{i_0} }{(2-\nu )_{i_0}(1)_{i_0}} \right.   \left. \sum_{i_1=i_0}^{\omega _1} \frac{(-\omega _1)_{i_1} (4-\nu )_{i_0}(3)_{i_0}}{(-\omega _1)_{i_0} (4-\nu  )_{i_1}(3)_{i_1}} \eta ^{i_1}\right\} \rho \nonumber\\
&&+ \sum_{n=2}^{\infty } \left\{ \sum_{i_0=0}^{\omega _0}\frac{(i_0+1+\Omega / \mu -\nu )}{ (i_0+ 3-\nu )(i_0+ 2)}\frac{(-\omega_0)_{i_0}  }{(2-\nu  )_{i_0}(1)_{i_0}}\right.\nonumber\\
&&\times \prod _{k=1}^{n-1} \left\{ \sum_{i_k=i_{k-1}}^{\omega_k} \frac{(i_k+ 2k +1+\Omega / \mu -\nu )}{(i_k+ 2 k+3-\nu )(i_k+ 2k+2)}\right.  \left.\frac{(-\omega _k)_{i_k} (2k+2-\nu )_{i_{k-1}}(2k+1)_{i_{k-1}}}{(-\omega _k)_{i_{k-1}}(2k+2-\nu  )_{i_k}(2k+1)_{i_k}}\right\} \nonumber\\
&&\times \left. \left.\sum_{i_n= i_{n-1}}^{\omega _n} \frac{(-\omega_n)_{i_n} (2n+2-\nu )_{i_{n-1}}(2n+1)_{i_{n-1}}}{(-\omega _n)_{i_{n-1}} (2n+2-\nu )_{i_n}(2n+1)_{i_n}} \eta ^{i_n} \right\} \rho ^n \right\}\label{eq:5009}
\end{eqnarray}
\end{remark}
For the minimum value of the GCH equation of the second kind for polynomial of type 2 about $x=0$, put $\omega _0=\omega _1=\omega _2=\cdots=0$ in (\ref{eq:5009}). 
\begin{eqnarray}
y(x)&=& RW_{0}^R\left(\mu ,\varepsilon ,\nu ,\Omega ,\omega =-(2 j+1-\nu); \rho =-\mu x^2; \eta = -\varepsilon x \right)\nonumber\\
&=& x^{1-\nu }\; _1F_1\left( \frac{\Omega }{2\mu }-\frac{\nu }{2}+\frac{1}{2}, -\frac{\nu }{2}+\frac{3}{2}, -\frac{1}{2}\mu x^2 \right) \;\;\mbox{where}\;-\infty < x< \infty \nonumber 
\end{eqnarray}
In Ref.\cite{Chou2012i,Chou2012j} I treat $\Omega $ as a fixed value and $\mu$, $\varepsilon$, $\nu $, $\omega$ as free variables to construct the GCH polynomial of type 1: (1) if $\Omega = -\mu (2 \beta _j+j) $ where $j, \beta _j \in \mathbb{N}_{0}$, an analytic solution of the GCH equation turns to be the first kind of independent solution of the GCH polynomial of type 1. (2) if  $\Omega = -\mu (2 \psi _j+j+1-\nu ) $ where $j, \psi _j \in \mathbb{N}_{0}$, an analytic solution of the GCH equation turns to be the second kind of independent solution of the GCH polynomial of type 1. 

In this chapter I treat $\omega $ as a fixed value and $\mu$, $\varepsilon$, $\nu $, $\Omega$ as free variables to construct the GCH polynomial of type 2: (1) if $\omega =-(\omega_j +2 j)$ where $j,\omega _j \in \mathbb{N}_{0}$, an analytic solution of the GCH equation turns to be the first kind of independent solution of the GCH polynomial of type 2. (2) if $ \omega = -(\omega_j +2 j+1-\nu ) $, an analytic solution of the GCH equation turns to be the second kind of independent solution of the GCH polynomial of type 2.
\subsubsection{Infinite series}
In chapter 1 the general expression of power series of $y(x)$ for infinite series is defined by
\begin{eqnarray}
y(x) &=& \sum_{n=0}^{\infty } y_{n}(x) = y_0(x)+ y_1(x)+ y_2(x)+y_3(x)+\cdots \nonumber\\
&=& c_0 \Bigg\{ \sum_{i_0=0}^{\infty } \left( \prod _{i_1=0}^{i_0-1}A_{i_1} \right) x^{i_0+\lambda }
+ \sum_{i_0=0}^{\infty }\left\{ B_{i_0+1} \prod _{i_1=0}^{i_0-1}A_{i_1}  \sum_{i_2=i_0}^{\infty } \left( \prod _{i_3=i_0}^{i_2-1}A_{i_3+2} \right)\right\} x^{i_2+2+\lambda }  \nonumber\\
&& + \sum_{N=2}^{\infty } \Bigg\{ \sum_{i_0=0}^{\infty } \Bigg\{B_{i_0+1}\prod _{i_1=0}^{i_0-1} A_{i_1} 
\prod _{k=1}^{N-1} \Bigg( \sum_{i_{2k}= i_{2(k-1)}}^{\infty } B_{i_{2k}+2k+1}\prod _{i_{2k+1}=i_{2(k-1)}}^{i_{2k}-1}A_{i_{2k+1}+2k}\Bigg)\nonumber\\
&& \times  \sum_{i_{2N} = i_{2(N-1)}}^{\infty } \Bigg( \prod _{i_{2N+1}=i_{2(N-1)}}^{i_{2N}-1} A_{i_{2N+1}+2N} \Bigg) \Bigg\} \Bigg\} x^{i_{2N}+2N+\lambda }\Bigg\}   \label{eq:50010}
\end{eqnarray}
Substitute (\ref{eq:5004a})--(\ref{eq:5004c}) into (\ref{eq:50010}). 
The general expression of power series of the GCH equation for infinite series about $x=0$ is given by
\begin{eqnarray}
 y(x)&=& \sum_{n=0}^{\infty } y_n(x)= y_0(x)+ y_1(x)+ y_2(x)+ y_3(x)+\cdots \nonumber\\
&=& c_0 x^{\lambda } \left\{\sum_{i_0=0}^{\infty } \frac{ (\omega +\lambda )_{i_0}  }{(1+\lambda )_{i_0}(\nu +\lambda )_{i_0}} \eta ^{i_0}\right.\nonumber\\
&+& \left\{ \sum_{i_0=0}^{\infty }\frac{(i_0+ \Omega /\mu +\lambda ) }{(i_0+ 2+\lambda )(i_0+1+\nu + \lambda )}\frac{ ( \omega +\lambda  )_{i_0}  }{(1+\lambda )_{i_0}(\nu +\lambda )_{i_0}} \right.\nonumber\\
&&\times \left. \sum_{i_1=i_0}^{\infty } \frac{ ( \omega +2+\lambda )_{i_1}  (3+\lambda )_{i_0}(2+\nu +\lambda )_{i_0}}{ ( \omega +2+\lambda )_{i_0} (3+\lambda )_{i_1}(2+\nu +\lambda )_{i_1}}\eta ^{i_1}\right\} \rho \nonumber\\
&+& \sum_{n=2}^{\infty } \left\{ \sum_{i_0=0}^{\infty } \frac{(i_0+ \Omega /\mu +\lambda ) }{(i_0+ 2+\lambda )(i_0+ 1+\nu +\lambda )}\frac{ (\omega +\lambda  )_{i_0}  }{(1+\lambda )_{i_0}(\nu +\lambda )_{i_0}}\right.\nonumber\\
&&\times \prod _{k=1}^{n-1} \left\{ \sum_{i_k=i_{k-1}}^{\infty } \frac{(i_k+ 2k+ \Omega /\mu +\lambda) }{(i_k+ 2 k+2+\lambda )(i_k+ 2k+1+\nu +\lambda )} \right. \nonumber\\
&&\times \left. \frac{ (\omega +2k+\lambda  )_{i_k}  (2k+1+\lambda )_{i_{k-1}}(2k+\nu +\lambda )_{i_{k-1}}}{(\omega +2k+\lambda  )_{i_{k-1}} (2k+1+\lambda )_{i_{k-1}}(2k+\nu +\lambda )_{i_k}}\right\}\nonumber\\
&&\times \left.\left.\sum_{i_n= i_{n-1}}^{\infty }  \frac{ (\omega +2n+\lambda  )_{i_n}  (2n+1+\lambda )_{i_{n-1}}(2n+\nu +\lambda )_{i_{n-1}}}{(\omega +2n+\lambda  )_{i_{n-1}} (2n+1+\lambda )_{i_{n-1}}(2n+\nu +\lambda )_{i_n}} \eta ^{i_n} \right\} \rho ^n \right\} \hspace{1cm}\label{eq:50011}
\end{eqnarray}
Put $c_0$= 1 as $\lambda =0$  for the first kind of independent solutions of the GCH equation and $\lambda =1-\nu$ for the second one in (\ref{eq:50011}).  
\begin{remark}
The power series expansion of the GCH equation of the first kind for infinite series about $x=0$ using R3TRF is 
\begin{eqnarray}
 y(x)&=& QW^R\left(\mu ,\varepsilon ,\nu ,\Omega ,\omega; \rho =-\mu x^2; \eta = -\varepsilon x \right) \nonumber\\
&=&  \sum_{i_0=0}^{\infty } \frac{ (\omega )_{i_0}}{(1 )_{i_0}(\nu )_{i_0}} \eta ^{i_0} + \left\{ \sum_{i_0=0}^{\infty }\frac{(i_0+ \Omega /\mu ) }{(i_0+ 2 )(i_0+1+\nu )}\frac{ ( \omega )_{i_0}  }{(1 )_{i_0}(\nu )_{i_0}} \sum_{i_1=i_0}^{\infty } \frac{ ( \omega +2 )_{i_1}  (3 )_{i_0}(2+\nu )_{i_0}}{ (\omega +2 )_{i_0} (3 )_{i_1}(2+\nu )_{i_1}}\eta ^{i_1}\right\} \rho \nonumber\\
&+& \sum_{n=2}^{\infty } \left\{ \sum_{i_0=0}^{\infty } \frac{(i_0+ \Omega /\mu ) }{(i_0+ 2 )(i_0+ 1+\nu )}\frac{ (\omega )_{i_0}  }{(1)_{i_0}(\nu )_{i_0}}\right.\nonumber\\
&\times& \prod _{k=1}^{n-1} \left\{ \sum_{i_k=i_{k-1}}^{\infty } \frac{(i_k+ 2k+ \Omega /\mu ) }{(i_k+ 2 k+2 )(i_k+ 2k+1+\nu )} \frac{ (\omega +2k )_{i_k}  (2k+1 )_{i_{k-1}}(2k+\nu )_{i_{k-1}}}{(\omega +2k )_{i_{k-1}} (2k+1 )_{i_{k-1}}(2k+\nu )_{i_k}}\right\}\nonumber\\
&\times&  \left.\sum_{i_n= i_{n-1}}^{\infty }  \frac{ (\omega +2n )_{i_n}  (2n+1 )_{i_{n-1}}(2n+\nu )_{i_{n-1}}}{(\omega +2n)_{i_{n-1}} (2n+1 )_{i_{n-1}}(2n+\nu )_{i_n}} \eta ^{i_n} \right\} \rho ^n  \label{eq:50012}
\end{eqnarray}
\end{remark}
\begin{remark}
The power series expansion of the GCH equation of the second kind for infinite series about $x=0$ using R3TRF is
\begin{eqnarray}
 y(x)&=&  RW^R\left(\mu ,\varepsilon ,\nu ,\Omega ,\omega; \rho =-\mu x^2; \eta = -\varepsilon x \right)\nonumber\\
&=& x^{1-\nu } \left\{\sum_{i_0=0}^{\infty } \frac{ (\omega +1-\nu  )_{i_0}  }{(2-\nu)_{i_0}(1)_{i_0}} \eta ^{i_0}\right.\nonumber\\
&+& \left\{ \sum_{i_0=0}^{\infty }\frac{(i_0+ 1+\Omega /\mu -\nu ) }{(i_0+ 3-\nu )(i_0+2)}\frac{ ( \omega +1-\nu )_{i_0}  }{(2-\nu )_{i_0}(1)_{i_0}} \sum_{i_1=i_0}^{\infty } \frac{ ( \omega +3-\nu)_{i_1}  (4-\nu)_{i_0}(3)_{i_0}}{ ( \omega +3-\nu )_{i_0} (4-\nu )_{i_1}(3)_{i_1}}\eta ^{i_1}\right\} \rho \nonumber\\
&+& \sum_{n=2}^{\infty } \left\{ \sum_{i_0=0}^{\infty } \frac{(i_0+ 1+\Omega /\mu -\nu ) }{(i_0+ 3-\nu )(i_0+ 2)}\frac{ (\omega +1-\nu )_{i_0}  }{(2-\nu )_{i_0}(1)_{i_0}}\right.\nonumber\\
&\times& \prod _{k=1}^{n-1} \left\{ \sum_{i_k=i_{k-1}}^{\infty } \frac{(i_k+ 2k+1+ \Omega /\mu -\nu ) }{(i_k+ 2 k+3-\nu )(i_k+ 2k+2)} \frac{ (\omega +2k+1-\nu )_{i_k}  (2k+2-\nu )_{i_{k-1}}(2k+1)_{i_{k-1}}}{(\omega +2k+1-\nu )_{i_{k-1}} (2k+2-\nu )_{i_{k-1}}(2k+1)_{i_k}}\right\}\nonumber\\
&\times& \left.\left.\sum_{i_n= i_{n-1}}^{\infty }  \frac{ (\omega +2n+1-\nu )_{i_n}  (2n+2-\nu )_{i_{n-1}}(2n+1)_{i_{n-1}}}{(\omega +2n+1-\nu)_{i_{n-1}} (2n+2-\nu)_{i_{n-1}}(2n+1)_{i_n}} \eta ^{i_n} \right\} \rho ^n \right\} \label{eq:50013}
\end{eqnarray}
\end{remark}
It is required that $\nu \ne 0,-1,-2,\cdots$ for the first kind of independent solutions of the GCH equation for all cases. Because if it does not, its solution will be divergent. And it is required that $\nu \ne 2,3,4,\cdots$ for the second kind of independent solutions of the GCH equation for all cases.

The infinite series in this chapter are equivalent to the infinite series in Ref.\cite{Chou2012i,Chou2012j}. In this chapter $B_n$ is the leading term in sequence $c_n$ of the analytic function $y(x)$. In Ref.\cite{Chou2012i,Chou2012j} $A_n$ is the leading term in sequence $c_n$ of the analytic function $y(x)$.\footnote{As ${\displaystyle \frac{\Gamma (1/2+\nu/2-\Omega /(2\mu))}{\Gamma (1/2+\nu/2)}}$  is multiplied on 
(\ref{eq:50012}), the new (\ref{eq:50012}) is equivalent to the first kind solution of the GCH equation for the infinite series using 3TRF.\cite{Chou2012i} Again, As ${\displaystyle \left( -\mu /2\right)^{1/2(1-\nu)} \frac{\Gamma (1-\Omega /(2\mu ))}{\Gamma (3/2-\nu/2)} }$  is multiplied on (\ref{eq:50013}), the new (\ref{eq:50013}) corresponds to the second kind solution of the GCH equation for the infinite series using 3TRF.\cite{Chou2012i} }
%
%
\subsection{Integral formalism}
\subsubsection{Polynomial of type 2}
There is a generalized hypergeometric function which is
\begin{eqnarray}
I_l &=& \sum_{i_l= i_{l-1}}^{\omega _l} \frac{(-\omega _l)_{i_l}(2l+1+ \lambda )_{i_{l-1}}(2l+\nu + \lambda)_{i_{l-1}}}{(-\omega _l)_{i_{l-1}}(2l+1+ \lambda )_{i_l}(2l+\nu + \lambda)_{i_l}} \eta ^{i_l}\nonumber\\
&=& \eta^{i_{l-1}} 
\sum_{j=0}^{\infty } \frac{B(i_{l-1}+2l+ \lambda ,j+1) B(i_{l-1}+2l-1+\nu +\lambda , j+1)(i_{l-1}-\omega _l)_j}{(i_{l-1}+2l+ \lambda )^{-1}(i_{l-1}+2l-1+\nu +\lambda )^{-1}(1)_j \;j!} \eta^j \hspace{1.5cm}
\label{eq:50015}
\end{eqnarray}
By using integral form of beta function,
\begin{subequations}
\begin{equation}
B\left(i_{l-1}+2l+ \lambda , j+1\right)= \int_{0}^{1} dt_l\;t_l^{i_{l-1}+2l-1+ \lambda } (1-t_l)^j
\label{eq:50016a}
\end{equation}
\begin{equation}
 B\left(i_{l-1}+2l-1+\nu + \lambda , j+1\right)= \int_{0}^{1} du_l\;u_l^{i_{l-1}+2l -2+\nu + \lambda } (1-u_l)^j
 \label{eq:50016b}
\end{equation}
\end{subequations}
Substitute (\ref{eq:50016a}) and (\ref{eq:50016b}) into (\ref{eq:50015}), and divide $(i_{l-1}+2l+ \lambda )(i_{l-1} +2l-1+\nu + \lambda )$ into $I_l$.
\begin{eqnarray}
&& \frac{1}{(i_{l-1}+2l+ \lambda )(i_{l-1} +2l-1+\nu + \lambda )}
 \sum_{i_l= i_{l-1}}^{\omega _l} \frac{(-\omega _l)_{i_l}(2l+1+ \lambda )_{i_{l-1}}(2l+\nu + \lambda)_{i_{l-1}}}{(-\omega _l)_{i_{l-1}}(2l+1+ \lambda )_{i_l}(2l+\nu + \lambda)_{i_l}} \eta ^{i_l}\nonumber\\
&=&  \int_{0}^{1} dt_l\;t_l^{2l-1+ \lambda } \int_{0}^{1} du_l\;u_l^{2l-2+\nu+ \lambda } (\eta t_l u_l)^{i_{l-1}}
 \sum_{j=0}^{\infty } \frac{(i_{l-1}- \omega _l)_j}{(1)_j \;j!} [\eta (1-t_l)(1-u_l)]^j
 \label{eq:50017}
\end{eqnarray}
The integral form of confluent hypergeometric function of the first kind is given by
\begin{equation}
 \sum_{j=0}^{\infty }\frac{(-\alpha _0)_j}{(\gamma )_j j!}z^j= \frac{\Gamma ( \alpha _0+1) \Gamma (\gamma )}{2\pi i\; \Gamma (\alpha _0+\gamma )}\oint dv_l\;\frac{\exp\left(-\frac{z v_l}{(1-v_l)} \right)}{v_l^{\alpha _0+1} (1-v_l)^{\gamma }} \label{eq:50018}
\end{equation}
replace $\alpha_0$, $\gamma $ and z by $\omega _l-i_{l-1}$, 1 and $\eta (1-t_l)(1-u_l)$ in (\ref{eq:50018}).
\begin{equation}
\sum_{j=0}^{\infty } \frac{(i_{l-1}- \omega _l)_j}{(1)_j \;j!} [\eta (1-t_l)(1-u_l)]^j = \frac{1}{2\pi i}  \oint dv_l \frac{\exp\left(-\frac{v_l}{(1-v_l)}\eta (1-t_l)(1-u_l)\right)}{v_l^{\omega _l+1-i_{l-1}}(1-v_l)}  \label{eq:50019}
\end{equation}
Substitute (\ref{eq:50019}) into (\ref{eq:50017}).
\begin{eqnarray}
K_l &=& \frac{1}{(i_{l-1}+2l+ \lambda )(i_{l-1} +2l-1+\nu + \lambda )}
 \sum_{i_l= i_{l-1}}^{\omega _l} \frac{(-\omega _l)_{i_l}(2l+1+ \lambda )_{i_{l-1}}(2l+\nu + \lambda)_{i_{l-1}}}{(-\omega _l)_{i_{l-1}}(2l+1+ \lambda )_{i_l}(2l+\nu + \lambda)_{i_l}} \eta ^{i_l} \nonumber\\
&=&  \int_{0}^{1} dt_l\;t_l^{2l-1+ \lambda } \int_{0}^{1} du_l\;u_l^{2l-2+\nu+ \lambda }
\frac{1}{2\pi i}  \oint dv_l \frac{\exp\left(-\frac{v_l}{(1-v_l)}\eta (1-t_l)(1-u_l)\right)}{v_l^{\omega _l+1 }(1-v_l)} (\eta t_l u_l v_l)^{i_{l-1}}
 \hspace{1.5cm}\label{eq:50020}
\end{eqnarray}
Substitute (\ref{eq:50020}) into (\ref{eq:5007}) where $l=1,2,3,\cdots$; apply $K_1$ into the second summation of sub-power series $y_1(x)$, apply $K_2$ into the third summation and $K_1$ into the second summation of sub-power series $y_2(x)$, apply $K_3$ into the forth summation, $K_2$ into the third summation and $K_1$ into the second summation of sub-power series $y_3(x)$, etc.\footnote{$y_1(x)$ means the sub-power series in (\ref{eq:5007}) contains one term of $B_n's$, $y_2(x)$ means the sub-power series in (\ref{eq:5007}) contains two terms of $B_n's$, $y_3(x)$ means the sub-power series in (\ref{eq:5007}) contains three terms of $B_n's$, etc.}
\begin{theorem}
The general representation in the form of integral of the GCH polynomial of type 2 is given by
\begin{eqnarray}
 y(x)&=& \sum_{n=0}^{\infty } y_n(x)= y_0(x)+ y_1(x)+ y_2(x)+ y_3(x)+\cdots \nonumber\\
&=& c_0 x^{\lambda } \left\{ \sum_{i_0=0}^{\omega _0 }\frac{(-\omega _0)_{i_0}}{(1+ \lambda )_{i_0}(\nu + \lambda )_{i_0}}  \eta^{i_0} \right.\nonumber\\
&&+ \sum_{n=1}^{\infty } \left\{\prod _{k=0}^{n-1} \Bigg\{ \int_{0}^{1} dt_{n-k}\;t_{n-k}^{2(n-k)-1+ \lambda } \int_{0}^{1} du_{n-k}\;u_{n-k}^{ 2(n-k-1)+\nu+\lambda } \right.\nonumber\\
&&\times \frac{1}{2\pi i}  \oint dv_{n-k} \frac{\exp\left(-\frac{v_{n-k}}{(1-v_{n-k})}w_{n-k+1,n}(1-t_{n-k})(1-u_{n-k})\right)}{v_{n-k}^{\omega_{n-k}+1}(1-v_{n-k})} \nonumber\\
&&\times  w_{n-k,n}^{-\left(\Omega /\mu +2(n-k-1)+\lambda \right)} \left(  w_{n-k,n} \partial _{w_{n-k,n}} \right) w_{n-k,n}^{ \Omega /\mu +2(n-k-1)+\lambda }\Bigg\} \nonumber\\
&&\times \left.\left. \sum_{i_0=0}^{\omega _0}\frac{(-\omega _0)_{i_0}}{(1+ \lambda )_{i_0}(\nu + \lambda )_{i_0}}  w_{1,n}^{i_0}\right\} \rho ^n \right\}
\label{eq:50021}
\end{eqnarray}
where
\begin{equation}w_{a,b}=
\begin{cases} \displaystyle {\eta \prod _{l=a}^{b} t_l u_l v_l }\;\;\mbox{where}\; a\leq b\cr
\eta \;\;\mbox{only}\;\mbox{if}\; a>b
\end{cases}
\nonumber
\end{equation}
In the above, the first sub-integral form contains one term of $B_n's$, the second one contains two terms of $B_n$'s, the third one contains three terms of $B_n$'s, etc.
\end{theorem}
\begin{proof} 
In (\ref{eq:5007}) the power series expansions of sub-summation $y_0(x) $, $y_1(x)$, $y_2(x)$ and $y_3(x)$ of the GCH polynomial of type 2 are
\begin{equation}
 y(x)= \sum_{n=0}^{\infty } y_{n}(x) = y_0(x)+ y_1(x)+ y_2(x)+y_3(x)+\cdots \label{eq:50022}
\end{equation}
where
\begin{subequations}
\begin{equation}
 y_0(x)= c_0 x^{\lambda } \sum_{i_0=0}^{\omega _0 }\frac{(-\omega _0)_{i_0}}{(1+ \lambda )_{i_0}(\nu + \lambda )_{i_0}}  \eta^{i_0} \label{eq:50023a}
\end{equation}
\begin{eqnarray}
 y_1(x)&=&  c_0 x^{\lambda } \Bigg\{ \sum_{i_0=0}^{\omega _0}\frac{(i_0+\Omega /\mu +\lambda )}{(i_0+ 2 +\lambda )(i_0 +1+\nu + \lambda )} \frac{(-\omega _0)_{i_0}}{(1+ \lambda )_{i_0}(\nu + \lambda )_{i_0}}\nonumber\\
&&\times  \sum_{i_1=i_0}^{\omega_1} \frac{(-\omega_1)_{i_1}(3+ \lambda )_{i_0}(2+\nu + \lambda )_{i_0}}{(-\omega_1)_{i_0}(3+ \lambda )_{i_1}(2+\nu + \lambda)_{i_1}} \eta^{i_1} \Bigg\}\rho \label{eq:50023b}
\end{eqnarray}
\begin{eqnarray}
 y_2(x) &=& c_0 x^{\lambda }\Bigg\{ \sum_{i_0=0}^{\omega _0}\frac{(i_0+\Omega /\mu +\lambda )}{(i_0+ 2 +\lambda )(i_0 +1+\nu + \lambda )} \frac{(-\omega _0)_{i_0}}{(1+ \lambda )_{i_0}(\nu + \lambda )_{i_0}} \nonumber\\
&&\times  \sum_{i_1=i_0}^{\omega_1} \frac{(i_1+ 2 +\Omega /\mu +\lambda )}{(i_1+4+ \lambda )(i_1+3 +\nu + \lambda )}  \frac{(-\omega_1)_{i_1}(3+ \lambda )_{i_0}(2+\nu + \lambda )_{i_0}}{(-\omega_1)_{i_0}(3+ \lambda )_{i_1}(2+\nu + \lambda)_{i_1}}  \nonumber\\
&&\times \sum_{i_2=i_1}^{\omega _2} \frac{(-\omega_2)_{i_2}(5+ \lambda )_{i_1}(4+\nu +  \lambda )_{i_1}}{(-\omega _2)_{i_1}(5+ \lambda )_{i_2}(4+\nu + \lambda )_{i_2}} \eta^{i_2} \Bigg\} \rho ^2 
\label{eq:50023c}
\end{eqnarray}
\begin{eqnarray}
 y_3(x)&=&  c_0 x^{\lambda } \Bigg\{ \sum_{i_0=0}^{\omega _0}\frac{(i_0+\Omega /\mu +\lambda )}{(i_0+ 2 +\lambda )(i_0 +1+\nu + \lambda )} \frac{(-\omega _0)_{i_0}}{(1+ \lambda )_{i_0}(\nu + \lambda )_{i_0}}\nonumber\\
&&\times  \sum_{i_1=i_0}^{\omega_1} \frac{(i_1+ 2 +\Omega /\mu +\lambda )}{(i_1+4+ \lambda )(i_1+3 +\nu + \lambda )} \frac{(-\omega_1)_{i_1}(3+ \lambda )_{i_0}(2+\nu + \lambda )_{i_0}}{(-\omega_1)_{i_0}(3+ \lambda )_{i_1}(2+\nu + \lambda)_{i_1}} \nonumber\\
&&\times \sum_{i_2=i_1}^{\omega _2} \frac{(i_2+4+\Omega /\mu +\lambda )}{(i_2+6+ \lambda )(i_2+5+ \nu + \lambda )} \frac{(-\omega_2)_{i_2}(5+ \lambda )_{i_1}(4+\nu +  \lambda )_{i_1}}{(-\omega _2)_{i_1}(5+ \lambda )_{i_2}(4+\nu + \lambda )_{i_2}}\nonumber\\
&&\times \sum_{i_3=i_2}^{\omega _3} \frac{(-\omega _3)_{i_3}(7+ \lambda )_{i_2}(6 +\nu + \lambda )_{i_2}}{(-\omega _3)_{i_2}(7+ \lambda )_{i_3}(6 +\nu + \lambda )_{i_3}} \eta^{i_3} \Bigg\} \rho ^3 
\label{eq:50023d}
\end{eqnarray}
\end{subequations}
Put $l=1$ in (\ref{eq:50020}). Take the new (\ref{eq:50020}) into (\ref{eq:50023b}).
\begin{eqnarray}
 y_1(x)&=& c_0 x^{\lambda }  \int_{0}^{1} dt_1\;t_1^{1+ \lambda } \int_{0}^{1} du_1\;u_1^{\nu + \lambda } 
\frac{1}{2\pi i}  \oint dv_1 \frac{\exp\left(-\frac{v_1}{(1-v_1)}\eta (1-t_1)(1-u_1)\right)}{v_1^{\omega _1+1}(1-v_1)} \nonumber\\
&&\times \left\{ \sum_{i_0=0}^{\omega _0} \left(i_0+\Omega /\mu + \lambda \right) \frac{(-\omega _0)_{i_0}}{(1+ \lambda )_{i_0}(\nu + \lambda  )_{i_0}} (\eta t_1 u_1 v_1)^{i_0} \right\} \rho \nonumber\\
&=& c_0 x^{\lambda }  \int_{0}^{1} dt_1\;t_1^{1+ \lambda } \int_{0}^{1} du_1\;u_1^{\nu + \lambda } 
\frac{1}{2\pi i}  \oint dv_1 \frac{\exp\left(-\frac{v_1}{(1-v_1)}\eta (1-t_1)(1-u_1)\right)}{v_1^{\omega _1+1}(1-v_1)}  \nonumber\\
&&\times  w_{1,1}^{-(\Omega /\mu +\lambda )} \left( w_{1,1}\partial_{w_{1,1}} \right) w_{1,1}^{ \Omega /\mu +\lambda } \left\{ \sum_{i_0=0}^{\omega _0} \frac{(-\omega _0)_{i_0}}{(1+ \lambda )_{i_0}(\nu + \lambda )_{i_0}} w_{1,1}^{i_0} \right\} \rho \label{eq:50024}\\
&& \mathrm{where}\hspace{.5cm} w_{1,1}=\eta \prod _{l=1}^{1} t_l u_l v_l \nonumber
\end{eqnarray}
Put $l=2$ in (\ref{eq:50020}). Take the new (\ref{eq:50020}) into (\ref{eq:50023c}). 
\begin{eqnarray}
 y_2(x) &=& c_0 x^{\lambda } \int_{0}^{1} dt_2\;t_2^{3+ \lambda } \int_{0}^{1} du_2\;u_2^{2+\nu + \lambda } 
\frac{1}{2\pi i}  \oint dv_2 \frac{\exp\left(-\frac{v_2}{(1-v_2)}\eta (1-t_2)(1-u_2)\right)}{v_2^{\omega _2+1}(1-v_2)}\nonumber\\
&&\times   w_{2,2}^{-(\Omega /\mu +2+\lambda )} \left( w_{2,2}\partial_{w_{2,2}} \right) w_{2,2}^{ \Omega /\mu +2 +\lambda } \nonumber\\
&&\times  \left\{ \sum_{i_0=0}^{\omega _0}\frac{(i_0+\Omega /\mu +\lambda )}{(i_0+ 2 +\lambda )(i_0 +1+\nu + \lambda )} \frac{(-\omega _0)_{i_0}}{(1+ \lambda )_{i_0}(\nu + \lambda )_{i_0}}\right.\nonumber\\
&&\times  \left.\sum_{i_1=i_0}^{\omega_1} \frac{(-\omega_1)_{i_1}(3+ \lambda )_{i_0}(2+\nu + \lambda )_{i_0}}{(-\omega_1)_{i_0}(3+ \lambda )_{i_1}(2+\nu + \lambda)_{i_1}} w_{2,2}^{i_1} \right\} \rho ^2 \label{eq:50025}\\
&& \mathrm{where}\hspace{.5cm} w_{2,2}= \eta \prod _{l=2}^{2} t_l u_l v_l \nonumber
\end{eqnarray}
Put $l=1$ and $\eta = w_{2,2}$ in (\ref{eq:50020}). Take the new (\ref{eq:50020}) into (\ref{eq:50025}).
\begin{eqnarray}
 y_2(x)&=& c_0 x^{\lambda } \int_{0}^{1} dt_2\;t_2^{3+ \lambda } \int_{0}^{1} du_2\;u_2^{2+\nu + \lambda } 
\frac{1}{2\pi i}  \oint dv_2 \frac{\exp\left(-\frac{v_2}{(1-v_2)}\eta (1-t_2)(1-u_2)\right)}{v_2^{\omega _2+1}(1-v_2)}  \nonumber\\
&&\times w_{2,2}^{-(\Omega /\mu +2+\lambda )} \left( w_{2,2}\partial_{w_{2,2}} \right) w_{2,2}^{ \Omega /\mu +2 +\lambda } \nonumber\\
&&\times  \int_{0}^{1} dt_1\;t_1^{1+ \lambda } \int_{0}^{1} du_1\;u_1^{\nu + \lambda } 
\frac{1}{2\pi i}  \oint dv_1 \frac{\exp\left(-\frac{v_1}{(1-v_1)}w_{2,2}(1-t_1)(1-u_1)\right)}{v_1^{\omega _1+1}(1-v_1)} \nonumber\\
&&\times w_{1,2}^{-(\Omega /\mu +\lambda )} \left( w_{1,2}\partial_{w_{1,2}} \right) w_{2,2}^{ \Omega /\mu +\lambda } \Bigg\{ \sum_{i_0=0}^{\omega _0} \frac{(-\omega _0)_{i_0}}{(1+ \lambda )_{i_0}(\nu + \lambda )_{i_0}} w_{1,2}^{i_0} \Bigg\} \rho^2 \label{eq:50026}\\
&& \mathrm{where}\hspace{.5cm} w_{1,2}=\eta \prod _{l=1}^{2} t_l u_l v_l \nonumber
\end{eqnarray}
By using similar process for the previous cases of integral forms of $y_1(x)$ and $y_2(x)$, the integral form of sub-power series expansion of $y_3(x)$ is
\begin{eqnarray}
 y_3(x)&=& c_0 x^{\lambda } \int_{0}^{1} dt_3\;t_3^{5+ \lambda } \int_{0}^{1} du_3\;u_3^{4+\nu + \lambda } \frac{1}{2\pi i}  \oint dv_3 \frac{\exp\left(-\frac{v_3}{(1-v_3)}\eta (1-t_3)(1-u_3)\right)}{v_3^{\omega _3+1}(1-v_3)}\nonumber\\
&&\times w_{3,3}^{-(\Omega /\mu +4+\lambda )} \left( w_{3,3}\partial_{w_{3,3}} \right) w_{3,3}^{ \Omega /\mu +4 +\lambda } \nonumber\\
&&\times \int_{0}^{1} dt_2\;t_2^{3+ \lambda } \int_{0}^{1} du_2\;u_2^{2+\nu + \lambda } 
\frac{1}{2\pi i}  \oint dv_2 \frac{\exp\left(-\frac{v_2}{(1-v_2)}w_{3,3} (1-t_2)(1-u_2)\right)}{v_2^{\omega _2+1}(1-v_2)} \nonumber\\
&&\times w_{2,3}^{-(\Omega /\mu +2+\lambda )} \left( w_{2,3}\partial_{w_{2,3}} \right) w_{2,3}^{ \Omega /\mu +2 +\lambda } \nonumber\\
&&\times  \int_{0}^{1} dt_1\;t_1^{1+ \lambda } \int_{0}^{1} du_1\;u_1^{\nu + \lambda } 
\frac{1}{2\pi i}  \oint dv_1 \frac{\exp\left(-\frac{v_1}{(1-v_1)}w_{2,3}(1-t_1)(1-u_1)\right)}{v_1^{\omega _1+1}(1-v_1)} \nonumber\\
&&\times w_{1,3}^{-(\Omega /\mu +\lambda )} \left( w_{1,3}\partial_{w_{1,3}} \right) w_{1,3}^{ \Omega /\mu +\lambda } \Bigg\{ \sum_{i_0=0}^{\omega _0} \frac{(-\omega _0)_{i_0}}{(1+ \lambda )_{i_0}(\nu + \lambda )_{i_0}} w_{1,3}^{i_0}\Bigg\} \rho ^3 \label{eq:50027}
\end{eqnarray}
where
\begin{equation}
\begin{cases} \displaystyle{w_{3,3} =\eta \prod _{l=3}^{3} t_l u_l v_l} \cr
\displaystyle{w_{2,3} = \eta \prod _{l=2}^{3} t_l u_l v_l} \cr
\displaystyle{w_{1,3}= \eta \prod _{l=1}^{3} t_l u_l v_l}
\end{cases}
\nonumber
\end{equation}
By repeating this process for all higher terms of integral forms of sub-summation $y_m(x)$ terms where $m \geq 4$, we obtain every integral forms of $y_m(x)$ terms. 
Since we substitute (\ref{eq:50023a}), (\ref{eq:50024}), (\ref{eq:50026}), (\ref{eq:50027}) and including all integral forms of $y_m(x)$ terms where $m \geq 4$ into (\ref{eq:50022}), we obtain (\ref{eq:50021}).
\qed
\end{proof} 
Put $c_0$= 1 as $\lambda =0$  for the first kind of independent solutions of the GCH equation and $\lambda = 1-\nu $  for the second one in (\ref{eq:50021}). 
\begin{remark}
The integral representation of the GCH equation of the first kind for polynomial of type 2 about $x=0$ as $\omega =-(\omega_j +2 j)$ where $j,\omega_j = 0,1,2,\cdots$ is
\begin{eqnarray}
 y(x)&=&  QW_{\omega _j}^R\left(\mu ,\varepsilon ,\nu ,\Omega ,\omega =-(\omega_j +2 j); \rho =-\mu x^2; \eta = -\varepsilon x \right) \nonumber\\
&=& _1F_1 \left(-\omega _0; \nu ; \eta \right) + \sum_{n=1}^{\infty } \Bigg\{\prod _{k=0}^{n-1} \Bigg\{ \int_{0}^{1} dt_{n-k}\;t_{n-k}^{2(n-k)-1 } \int_{0}^{1} du_{n-k}\;u_{n-k}^{ 2(n-k-1)+\nu }  \nonumber\\
&&\times \frac{1}{2\pi i}  \oint dv_{n-k} \frac{\exp\left(-\frac{v_{n-k}}{(1-v_{n-k})}w_{n-k+1,n}(1-t_{n-k})(1-u_{n-k})\right)}{v_{n-k}^{\omega_{n-k}+1}(1-v_{n-k})} \nonumber\\
&&\times  w_{n-k,n}^{-\left(\Omega /\mu +2(n-k-1) \right)} \left(  w_{n-k,n} \partial _{w_{n-k,n}} \right) w_{n-k,n}^{ \Omega /\mu +2(n-k-1) }\Bigg\}\;  _1F_1 \left(-\omega _0; \nu ; w_{1,n} \right) \Bigg\} \rho ^n  \hspace{2cm}
\label{eq:50028}
\end{eqnarray}
\end{remark}
\begin{remark}
The integral representation of the GCH equation of the second kind for polynomial of type 2 about $x=0$ as $\omega =-(\omega_j +2 j+1-\nu)$ where $j,\omega_j = 0,1,2,\cdots$ is
\begin{eqnarray}
 y(x)&=&  RW_{\omega _j}^R\left(\mu ,\varepsilon ,\nu ,\Omega ,\omega =-(\omega_j +2 j+1-\nu); \rho =-\mu x^2; \eta = -\varepsilon x \right) \nonumber\\
&=& x^{1-\nu} \Bigg\{\;  _1F_1 \left(-\omega _0; 2-\nu ; \eta \right) + \sum_{n=1}^{\infty } \Bigg\{\prod _{k=0}^{n-1} \Bigg\{ \int_{0}^{1} dt_{n-k}\;t_{n-k}^{2(n-k)-\nu } \int_{0}^{1} du_{n-k}\;u_{n-k}^{ 2(n-k)-1 } \nonumber\\
&&\times \frac{1}{2\pi i}  \oint dv_{n-k} \frac{\exp\left(-\frac{v_{n-k}}{(1-v_{n-k})}w_{n-k+1,n}(1-t_{n-k})(1-u_{n-k})\right)}{v_{n-k}^{\omega_{n-k}+1}(1-v_{n-k})} \nonumber\\
&&\times  w_{n-k,n}^{-\left(\Omega /\mu +2(n-k)-1 -\nu \right)} \left(  w_{n-k,n} \partial _{w_{n-k,n}} \right) w_{n-k,n}^{ \Omega /\mu +2(n-k)-1-\nu }\Bigg\} \nonumber\\
&&\times\; _1F_1 \left(-\omega _0; 2-\nu ; w_{1,n} \right) \Bigg\} \rho ^n \Bigg\} \label{eq:50029}
\end{eqnarray}
\end{remark}
In the above, $_1F_1 \left(a; b; z \right)$ is a Kummer function of the first kind which is defined by
\begin{eqnarray}
_1F_1 \left(a; b; z \right) &=& M(a,b,z) = \sum_{n=0}^{\infty } \frac{(a)_n}{(b)_n n!} z^n = e^z M(b-a,b,-z) \nonumber\\
&=& -\frac{1}{2\pi i}\frac{\Gamma \left( 1-a\right)\Gamma \left(b\right)}{\Gamma \left(b-a\right) } \oint  dv_j\; e^{z v_j}(-v_j)^{a-1} \left( 1- v_j \right)^{b-a-1} \nonumber\\
&=&  \frac{\Gamma \left(a\right)}{2\pi i} \oint  dv_j\; e^{v_j}v_j^{-b} \left( 1-\frac{z}{v_j}\right)^{-a} \nonumber\\
&=& \frac{1}{2\pi i} \frac{\Gamma \left( 1-a\right)\Gamma \left(b\right)}{\Gamma \left(b-a\right) } \oint  dv_j\; e^{-\frac{z \;v_j}{1-v_j}}v_j^{a-1} \left( 1- v_j \right)^{-b} \label{er:50024}
\end{eqnarray}
\subsubsection{Infinite series}
Let's consider the integral representation of the GCH equation about $x=0$ for infinite series by applying R3TRF.
There is a generalized hypergeometric function which is written by
\begin{eqnarray}
M_l &=& \sum_{i_l= i_{l-1}}^{\infty } \frac{\left( \omega +2l+\lambda \right)_{i_l}(2l+1+ \lambda )_{i_{l-1}}(2l+\nu + \lambda)_{i_{l-1}}}{\left( \omega +2l+\lambda \right)_{i_{l-1}}(2l+1+ \lambda )_{i_l}(2l+\nu + \lambda)_{i_l}} \eta ^{i_l}\nonumber\\
&=& \eta^{i_{l-1}} 
\sum_{j=0}^{\infty } \frac{B(i_{l-1}+2l+ \lambda ,j+1) B(i_{l-1}+2l-1+\nu +\lambda , j+1)(\omega +2l+\lambda+ i_{l-1} )_j}{(i_{l-1}+2l+ \lambda )^{-1}(i_{l-1}+2l-1+\nu +\lambda )^{-1}(1)_j \;j!} \eta^j \hspace{1.5cm} \label{er:50025}
\end{eqnarray}
Substitute (\ref{eq:50016a}) and (\ref{eq:50016b}) into (\ref{er:50025}), and divide $(i_{l-1}+2l+ \lambda )(i_{l-1} +2l-1+\nu + \lambda )$ into the new (\ref{er:50025}).
\begin{eqnarray}
&& \frac{1}{(i_{l-1}+2l+ \lambda )(i_{l-1} +2l-1+\nu + \lambda )}
 \sum_{i_l= i_{l-1}}^{\infty } \frac{\left( \omega +2l+\lambda \right)_{i_l}(2l+1+ \lambda )_{i_{l-1}}(2l+\nu + \lambda)_{i_{l-1}}}{\left( \omega +2l+\lambda \right)_{i_{l-1}}(2l+1+ \lambda )_{i_l}(2l+\nu + \lambda)_{i_l}} \eta ^{i_l}\nonumber\\
&=&  \int_{0}^{1} dt_l\;t_l^{2l-1+ \lambda } \int_{0}^{1} du_l\;u_l^{2l-2+\nu+ \lambda } (\eta t_l u_l)^{i_{l-1}}
 \sum_{j=0}^{\infty } \frac{(\omega +2l+\lambda+ i_{l-1} )_j}{(1)_j \;j!} (\eta (1-t_l)(1-u_l))^j\hspace{1.5cm}
 \label{er:50026}
\end{eqnarray}
Replace $a$, $b$ and $z$ by $ \omega +2l+\lambda+ i_{l-1}$, 1 and $\eta(1-t_j)(1-u_j)$ in (\ref{er:50024}). Take the new (\ref{er:50024}) into (\ref{er:50026}).
\begin{eqnarray}
V_l &=& \frac{1}{(i_{l-1}+2l+ \lambda )(i_{l-1} +2l-1+\nu + \lambda )}
 \sum_{i_l= i_{l-1}}^{\infty } \frac{\left( \omega +2l+\lambda \right)_{i_l}(2l+1+ \lambda )_{i_{l-1}}(2l+\nu + \lambda)_{i_{l-1}}}{\left( \omega +2l+\lambda \right)_{i_{l-1}}(2l+1+ \lambda )_{i_l}(2l+\nu + \lambda)_{i_l}} \eta ^{i_l} \nonumber\\
&=&  \int_{0}^{1} dt_l\;t_l^{2l-1+ \lambda } \int_{0}^{1} du_l\;u_l^{2l-2+\nu+ \lambda }
\frac{1}{2\pi i}  \oint dv_l \frac{\exp\left(-\frac{v_l}{(1-v_l)}\eta (1-t_l)(1-u_l)\right)}{v_l^{-(\omega +2l-1+\lambda )}(1-v_l)} (\eta t_l u_l v_l)^{i_{l-1}}
 \hspace{1.5cm}\label{er:50027}
\end{eqnarray}
Substitute (\ref{er:50027}) into (\ref{eq:50011}) where $l=1,2,3,\cdots$; apply $V_1$ into the second summation of sub-power series $y_1(x)$, apply $V_2$ into the third summation and $V_1$ into the second summation of sub-power series $y_2(x)$, apply $V_3$ into the forth summation, $V_2$ into the third summation and $V_1$ into the second summation of sub-power series $y_3(x)$, etc.\footnote{$y_1(x)$ means the sub-power series in (\ref{eq:50011}) contains one term of $B_n's$, $y_2(x)$ means the sub-power series in (\ref{eq:50011}) contains two terms of $B_n's$, $y_3(x)$ means the sub-power series in (\ref{eq:50011}) contains three terms of $B_n's$, etc.}
\begin{theorem}
The general representation in the form of integral of the GCH equation for infinite series about $x=0$ using R3TRF is given by
\begin{eqnarray}
 y(x) &=& \sum_{n=0}^{\infty } y_n(x)= y_0(x)+ y_1(x)+ y_2(x)+ y_3(x)+ \cdots \nonumber\\
&=& c_0 x^{\lambda } \left\{ \sum_{i_0=0}^{\infty }\frac{(\omega +\lambda )_{i_0}}{(1+ \lambda )_{i_0}(\nu + \lambda )_{i_0}}  \eta^{i_0} \right. \nonumber\\
&&+ \sum_{n=1}^{\infty } \left\{\prod _{k=0}^{n-1} \Bigg\{ \int_{0}^{1} dt_{n-k}\;t_{n-k}^{2(n-k)-1+ \lambda } \int_{0}^{1} du_{n-k}\;u_{n-k}^{ 2(n-k-1)+\nu+\lambda } \right.\nonumber\\
&&\times \frac{1}{2\pi i}  \oint dv_{n-k} \frac{\exp\left(-\frac{v_{n-k}}{(1-v_{n-k})}w_{n-k+1,n}(1-t_{n-k})(1-u_{n-k})\right)}{v_{n-k}^{-(\omega +2(n-k)-1+\lambda )}(1-v_{n-k})} \nonumber\\
&&\times  w_{n-k,n}^{-\left(\Omega /\mu +2(n-k-1)+\lambda \right)} \left(  w_{n-k,n} \partial _{w_{n-k,n}} \right) w_{n-k,n}^{ \Omega /\mu +2(n-k-1)+\lambda }\Bigg\}  \nonumber\\
&&\times \left.\left. \sum_{i_0=0}^{\infty}\frac{(\omega +\lambda)_{i_0}}{(1+ \lambda )_{i_0}(\nu + \lambda )_{i_0}}  w_{1,n}^{i_0}\right\} \rho ^n \right\} \label{eq:50030}
\end{eqnarray}
In the above, the first sub-integral form contains one term of $B_n's$, the second one contains two terms of $B_n$'s, the third one contains three terms of $B_n$'s, etc.
\end{theorem}
\begin{proof} 
In (\ref{eq:50011}) sub-power series $y_0(x) $, $y_1(x)$, $y_2(x)$ and $y_3(x)$ of the CHE for infinite series using 3TRF about $x=0$ are given by
\begin{subequations}
\begin{equation}
 y_0(x)= c_0 x^{\lambda } \sum_{i_0=0}^{\infty }\frac{\left( \omega +\lambda \right)_{i_0}}{(1+ \lambda )_{i_0}(\nu + \lambda )_{i_0}}  \eta^{i_0} \label{er:50028a}
\end{equation}
\begin{eqnarray}
 y_1(x)&=&  c_0 x^{\lambda } \left\{ \sum_{i_0=0}^{\infty}\frac{(i_0+\Omega /\mu +\lambda )}{(i_0+ 2 +\lambda )(i_0 +1+\nu + \lambda )} \frac{\left( \omega +\lambda \right)_{i_0}}{(1+ \lambda )_{i_0}(\nu + \lambda )_{i_0}}\right.\nonumber\\
&&\times \left. \sum_{i_1=i_0}^{\infty} \frac{\left( \omega +2+\lambda \right)_{i_1}(3+ \lambda )_{i_0}(2+\nu + \lambda )_{i_0}}{\left( \omega +2+\lambda \right)_{i_0}(3+ \lambda )_{i_1}(2+\nu + \lambda)_{i_1}} \eta^{i_1} \right\}\rho \label{er:50028b}
\end{eqnarray}
\begin{eqnarray}
 y_2(x) &=& c_0 x^{\lambda }\left\{ \sum_{i_0=0}^{\infty}\frac{(i_0+\Omega /\mu +\lambda )}{(i_0+ 2 +\lambda )(i_0 +1+\nu + \lambda )} \frac{\left( \omega +\lambda \right)_{i_0}}{(1+ \lambda )_{i_0}(\nu + \lambda )_{i_0}} \right.\nonumber\\
&&\times  \sum_{i_1=i_0}^{\infty} \frac{(i_1+ 2 +\Omega /\mu +\lambda )}{(i_1+4+ \lambda )(i_1+3 +\nu + \lambda )}  \frac{\left( \omega +2+\lambda \right)_{i_1}(3+ \lambda )_{i_0}(2+\nu + \lambda )_{i_0}}{\left( \omega +2+\lambda \right)_{i_0}(3+ \lambda )_{i_1}(2+\nu + \lambda)_{i_1}}  \nonumber\\
&&\times \left.\sum_{i_2=i_1}^{\infty} \frac{\left( \omega +4+\lambda \right)_{i_2}(5+ \lambda )_{i_1}(4+\nu +  \lambda )_{i_1}}{\left( \omega +4+\lambda \right)_{i_1}(5+ \lambda )_{i_2}(4+\nu + \lambda )_{i_2}} \eta^{i_2} \right\} \rho ^2 
\label{er:50028c}
\end{eqnarray}
\begin{eqnarray}
 y_3(x)&=&  c_0 x^{\lambda } \left\{ \sum_{i_0=0}^{\infty}\frac{(i_0+\Omega /\mu +\lambda )}{(i_0+ 2 +\lambda )(i_0 +1+\nu + \lambda )} \frac{\left( \omega +\lambda \right)_{i_0}}{(1+ \lambda )_{i_0}(\nu + \lambda )_{i_0}}\right.\nonumber\\
&&\times  \sum_{i_1=i_0}^{\infty} \frac{(i_1+ 2 +\Omega /\mu +\lambda )}{(i_1+4+ \lambda )(i_1+3 +\nu + \lambda )} \frac{\left( \omega +2+\lambda \right)_{i_1}(3+ \lambda )_{i_0}(2+\nu + \lambda )_{i_0}}{\left( \omega +2+\lambda \right)_{i_0}(3+ \lambda )_{i_1}(2+\nu + \lambda)_{i_1}} \nonumber\\
&&\times \sum_{i_2=i_1}^{\infty } \frac{(i_2+4+\Omega /\mu +\lambda )}{(i_2+6+ \lambda )(i_2+5+ \nu + \lambda )} \frac{\left( \omega +4+\lambda \right)_{i_2}(5+ \lambda )_{i_1}(4+\nu +  \lambda )_{i_1}}{\left( \omega +4+\lambda \right)_{i_1}(5+ \lambda )_{i_2}(4+\nu + \lambda )_{i_2}}\nonumber\\
&&\times \left.\sum_{i_3=i_2}^{\infty} \frac{\left( \omega +6+\lambda \right)_{i_3}(7+ \lambda )_{i_2}(6 +\nu + \lambda )_{i_2}}{\left( \omega +6+\lambda \right)_{i_2}(7+ \lambda )_{i_3}(6 +\nu + \lambda )_{i_3}} \eta^{i_3} \right\} \rho ^3 
\label{er:50028d}
\end{eqnarray}
\end{subequations}
Put $l=1$ in (\ref{er:50027}). Take the new (\ref{er:50027}) into (\ref{er:50028b}).
\begin{eqnarray}
 y_1(x)&=& c_0 x^{\lambda }  \int_{0}^{1} dt_1\;t_1^{1+ \lambda } \int_{0}^{1} du_1\;u_1^{\nu + \lambda } 
\frac{1}{2\pi i}  \oint dv_1 \frac{\exp\left(-\frac{v_1}{(1-v_1)}\eta (1-t_1)(1-u_1)\right)}{v_1^{-(\omega +1+\lambda )}(1-v_1)} \nonumber\\
&&\times \left\{ \sum_{i_0=0}^{\infty } \left(i_0+\Omega /\mu + \lambda \right) \frac{( \omega +\lambda )_{i_0}}{(1+ \lambda )_{i_0}(\nu + \lambda  )_{i_0}} (\eta t_1 u_1 v_1)^{i_0} \right\} \rho \nonumber\\
&=& c_0 x^{\lambda }  \int_{0}^{1} dt_1\;t_1^{1+ \lambda } \int_{0}^{1} du_1\;u_1^{\nu + \lambda } 
\frac{1}{2\pi i}  \oint dv_1 \frac{\exp\left(-\frac{v_1}{(1-v_1)}\eta (1-t_1)(1-u_1)\right)}{v_1^{-( \omega +1+\lambda )}(1-v_1)}  \nonumber\\
&&\times  w_{1,1}^{-(\Omega /\mu +\lambda )} \left( w_{1,1}\partial_{w_{1,1}} \right) w_{1,1}^{ \Omega /\mu +\lambda } \left\{ \sum_{i_0=0}^{\infty } \frac{( \omega +\lambda )_{i_0}}{(1+ \lambda )_{i_0}(\nu + \lambda )_{i_0}} w_{1,1}^{i_0} \right\} \rho \label{er:50029}\\
&& \mathrm{where}\hspace{.5cm} w_{1,1}=\eta \prod _{l=1}^{1} t_l u_l v_l \nonumber
\end{eqnarray}
Put $l=2$ in (\ref{er:50027}). Take the new (\ref{er:50027}) into (\ref{er:50028c}). 
\begin{eqnarray}
 y_2(x) &=& c_0 x^{\lambda } \int_{0}^{1} dt_2\;t_2^{3+ \lambda } \int_{0}^{1} du_2\;u_2^{2+\nu + \lambda } 
\frac{1}{2\pi i}  \oint dv_2 \frac{\exp\left(-\frac{v_2}{(1-v_2)}\eta (1-t_2)(1-u_2)\right)}{v_2^{-(\omega +3+\lambda )}(1-v_2)}\nonumber\\
&&\times   w_{2,2}^{-(\Omega /\mu +2+\lambda )} \left( w_{2,2}\partial_{w_{2,2}} \right) w_{2,2}^{ \Omega /\mu +2 +\lambda } \nonumber\\
&&\times  \left\{ \sum_{i_0=0}^{\infty }\frac{(i_0+\Omega /\mu +\lambda )}{(i_0+ 2 +\lambda )(i_0 +1+\nu + \lambda )} \frac{ (\omega +\lambda )_{i_0}}{(1+ \lambda )_{i_0}(\nu + \lambda )_{i_0}}\right.\nonumber\\
&&\times  \left.\sum_{i_1=i_0}^{\infty } \frac{(\omega +2+\lambda )_{i_1}(3+ \lambda )_{i_0}(2+\nu + \lambda )_{i_0}}{(\omega +2+\lambda )_{i_0}(3+ \lambda )_{i_1}(2+\nu + \lambda)_{i_1}} w_{2,2}^{i_1} \right\} \rho ^2 \label{er:50030}\\
&& \mathrm{where}\hspace{.5cm} w_{2,2}= \eta \prod _{l=2}^{2} t_l u_l v_l \nonumber
\end{eqnarray}
Put $l=1$ and $\eta = w_{2,2}$ in (\ref{er:50027}). Take the new (\ref{er:50027}) into (\ref{er:50030}).
\begin{eqnarray}
 y_2(x)&=& c_0 x^{\lambda } \int_{0}^{1} dt_2\;t_2^{3+ \lambda } \int_{0}^{1} du_2\;u_2^{2+\nu + \lambda } 
\frac{1}{2\pi i}  \oint dv_2 \frac{\exp\left(-\frac{v_2}{(1-v_2)}\eta (1-t_2)(1-u_2)\right)}{v_2^{-(\omega +3+\lambda )}(1-v_2)}\nonumber\\
&&\times  w_{2,2}^{-(\Omega /\mu +2+\lambda )} \left( w_{2,2}\partial_{w_{2,2}} \right) w_{2,2}^{ \Omega /\mu +2 +\lambda } \nonumber\\
&&\times  \int_{0}^{1} dt_1\;t_1^{1+ \lambda } \int_{0}^{1} du_1\;u_1^{\nu + \lambda } 
\frac{1}{2\pi i}  \oint dv_1 \frac{\exp\left(-\frac{v_1}{(1-v_1)}w_{2,2}(1-t_1)(1-u_1)\right)}{v_1^{-(\omega +1+\lambda )}(1-v_1)}  \nonumber\\
&&\times w_{1,2}^{-(\Omega /\mu +\lambda )} \left( w_{1,2}\partial_{w_{1,2}} \right) w_{2,2}^{ \Omega /\mu +\lambda } \Bigg\{ \sum_{i_0=0}^{\infty } \frac{(\omega +\lambda )_{i_0}}{(1+ \lambda )_{i_0}(\nu + \lambda )_{i_0}} w_{1,2}^{i_0} \Bigg\} \rho^2 \label{er:50031}\\
&& \mathrm{where}\hspace{.5cm} w_{1,2}=\eta \prod _{l=1}^{2} t_l u_l v_l \nonumber
\end{eqnarray}
By using similar process for the previous cases of integral forms of $y_1(x)$ and $y_2(x)$, the integral form of sub-power series expansion of $y_3(x)$ is
\begin{eqnarray}
 y_3(x)&=& c_0 x^{\lambda } \int_{0}^{1} dt_3\;t_3^{5+ \lambda } \int_{0}^{1} du_3\;u_3^{4+\nu + \lambda } \frac{1}{2\pi i}  \oint dv_3 \frac{\exp\left(-\frac{v_3}{(1-v_3)}\eta (1-t_3)(1-u_3)\right)}{v_3^{-(\omega +5+\lambda )}(1-v_3)}\nonumber\\
&&\times w_{3,3}^{-(\Omega /\mu +4+\lambda )} \left( w_{3,3}\partial_{w_{3,3}} \right) w_{3,3}^{ \Omega /\mu +4 +\lambda } \nonumber\\
&&\times \int_{0}^{1} dt_2\;t_2^{3+ \lambda } \int_{0}^{1} du_2\;u_2^{2+\nu + \lambda } 
\frac{1}{2\pi i}  \oint dv_2 \frac{\exp\left(-\frac{v_2}{(1-v_2)}w_{3,3} (1-t_2)(1-u_2)\right)}{v_2^{-(\omega +3+\lambda )}(1-v_2)}\nonumber\\
&&\times  w_{2,3}^{-(\Omega /\mu +2+\lambda )} \left( w_{2,3}\partial_{w_{2,3}} \right) w_{2,3}^{ \Omega /\mu +2 +\lambda } \nonumber\\
&&\times  \int_{0}^{1} dt_1\;t_1^{1+ \lambda } \int_{0}^{1} du_1\;u_1^{\nu + \lambda } 
\frac{1}{2\pi i}  \oint dv_1 \frac{\exp\left(-\frac{v_1}{(1-v_1)}w_{2,3}(1-t_1)(1-u_1)\right)}{v_1^{-(\omega +1+\lambda )}(1-v_1)} \nonumber\\
&&\times  w_{1,3}^{-(\Omega /\mu +\lambda )} \left( w_{1,3}\partial_{w_{1,3}} \right) w_{1,3}^{ \Omega /\mu +\lambda }   \Bigg\{ \sum_{i_0=0}^{\infty } \frac{(\omega +\lambda )_{i_0}}{(1+ \lambda )_{i_0}(\nu + \lambda )_{i_0}} w_{1,3}^{i_0}\Bigg\} \rho ^3 \label{er:50032}
\end{eqnarray}
where
\begin{equation}
\begin{cases} \displaystyle{w_{3,3} =\eta \prod _{l=3}^{3} t_l u_l v_l} \cr
\displaystyle{w_{2,3} = \eta \prod _{l=2}^{3} t_l u_l v_l} \cr
\displaystyle{w_{1,3}= \eta \prod _{l=1}^{3} t_l u_l v_l}
\end{cases}
\nonumber
\end{equation}
By repeating this process for all higher terms of integral forms of sub-summation $y_m(x)$ terms where $m \geq 4$, we obtain every integral forms of $y_m(x)$ terms. 
Since we substitute (\ref{er:50028a}), (\ref{er:50029}), (\ref{er:50031}), (\ref{er:50032}) and including all integral forms of $y_m(x)$ terms where $m \geq 4$ into (\ref{eq:50011}), we obtain (\ref{eq:50030}).\footnote{Or replace the finite summation with an interval $[0,\omega _0]$ by infinite summation with an interval  $[0,\infty ]$ in (\ref{eq:50021}). Replace $ \omega _0$ and $ \omega _{n-j}$ by $-(\omega +\lambda )$ and $-(\omega +2(n-k)+\lambda )$ into the new (\ref{eq:50021}). Its solution is also equivalent to (\ref{eq:50030})}
\qed
\end{proof} 
Put $c_0$= 1 as $\lambda =0$  for the first kind of independent solutions of the GCH equation and $\lambda = 1-\nu $  for the second one in (\ref{eq:50030}).
\begin{remark}
The integral representation of the GCH equation of the first kind for infinite series about $x=0$ using R3TRF is
\begin{eqnarray}
 y(x) &=&  QW^R\left(\mu ,\varepsilon ,\nu ,\Omega ,\omega; \rho =-\mu x^2; \eta = -\varepsilon x \right) \nonumber\\
&=& _1F_1 \left( \omega ; \nu ; \eta \right) + \sum_{n=1}^{\infty } \Bigg\{\prod _{k=0}^{n-1} \Bigg\{ \int_{0}^{1} dt_{n-k}\;t_{n-k}^{2(n-k)-1 } \int_{0}^{1} du_{n-k}\;u_{n-k}^{ 2(n-k-1)+\nu }  \nonumber\\
&&\times \frac{1}{2\pi i}  \oint dv_{n-k} \frac{\exp\left(-\frac{v_{n-k}}{(1-v_{n-k})}w_{n-k+1,n}(1-t_{n-k})(1-u_{n-k})\right)}{v_{n-k}^{-(\omega +2(n-k)-1 )}(1-v_{n-k})} \nonumber\\
&&\times  w_{n-k,n}^{-\left(\Omega /\mu +2(n-k-1) \right)} \left(  w_{n-k,n} \partial _{w_{n-k,n}} \right) w_{n-k,n}^{ \Omega /\mu +2(n-k-1) }\Bigg\}\;  _1F_1 \left( \omega ; \nu ; w_{1,n} \right) \Bigg\} \rho ^n \hspace{1.5cm}
\label{eq:50031}
\end{eqnarray}
\end{remark}
\begin{remark}
The integral representation of the GCH equation of the second kind for infinite series about $x=0$ using R3TRF is
\begin{eqnarray}
 y(x)&=& RW^R\left(\mu ,\varepsilon ,\nu ,\Omega ,\omega; \rho =-\mu x^2; \eta = -\varepsilon x \right) \nonumber\\
&=& x^{1-\nu} \Bigg\{\; _1F_1 \left( \omega +1-\nu ; 2-\nu ; \eta \right) + \sum_{n=1}^{\infty } \Bigg\{\prod _{k=0}^{n-1} \Bigg\{ \int_{0}^{1} dt_{n-k}\;t_{n-k}^{2(n-k)-\nu } \int_{0}^{1} du_{n-k}\;u_{n-k}^{ 2(n-k)-1 }  \nonumber\\
&&\times \frac{1}{2\pi i}  \oint dv_{n-k} \frac{\exp\left(-\frac{v_{n-k}}{(1-v_{n-k})}w_{n-k+1,n}(1-t_{n-k})(1-u_{n-k})\right)}{v_{n-k}^{-(\omega +2(n-k)-\nu )}(1-v_{n-k})} \nonumber\\
&&\times  w_{n-k,n}^{-\left(\Omega /\mu +2(n-k)-1 -\nu \right)} \left(  w_{n-k,n} \partial _{w_{n-k,n}} \right) w_{n-k,n}^{ \Omega /\mu +2(n-k)-1-\nu }\Bigg\}\nonumber\\
&&\times\; _1F_1 \left( \omega +1-\nu ; 2-\nu ; w_{1,n} \right) \Bigg\} \rho ^n \Bigg\} \hspace{1cm} \label{eq:50032}
\end{eqnarray}
\end{remark}
 As ${\displaystyle \frac{\Gamma (1/2+\nu/2-\Omega /(2\mu))}{\Gamma (1/2+\nu/2)}}$  is multiplied on 
(\ref{eq:50031}), the new (\ref{eq:50031}) is equivalent to the integral form of the first kind solution of the GCH equation for the infinite series using 3TRF.\cite{Chou2012i} Again, As ${\displaystyle \left( -\mu /2\right)^{1/2(1-\nu)} \frac{\Gamma (1-\Omega /(2\mu ))}{\Gamma (3/2-\nu/2)} }$  is multiplied on (\ref{eq:50032}), the new (\ref{eq:50032}) corresponds to the integral representation of the second kind solution of the GCH equation for the infinite series using 3TRF.\cite{Chou2012i} 
\subsection{Generating function for the GCH polynomial of type 2}
Now let's investigate the generating function for the type 2 GCH polynomials of the first and second kind about $x=0$. 
\begin{definition}
I define that
\begin{equation}
\begin{cases}
\displaystyle { s_{a,b}} = \begin{cases} \displaystyle {  s_a\cdot s_{a+1}\cdot s_{a+2}\cdots s_{b-2}\cdot s_{b-1}\cdot s_b}\;\;\mbox{where}\;a>b \cr
s_a \;\;\mbox{only}\;\mbox{if}\;a=b\end{cases}
\cr
\cr
\displaystyle {  \widetilde{w}_{a,b}= \eta s_{a,\infty }\prod_{l=a}^{b}t_l u_l}
\end{cases}\label{eq:50033}
\end{equation}
where
\begin{equation}
a,b\in \mathbb{N}_{0} \nonumber
\end{equation}
\end{definition}
And we have
\begin{equation}
\sum_{\omega _i = \omega _j}^{\infty } s_i^{\omega _i} = \frac{s_i^{\omega _j}}{(1-s_i)}\label{eq:50034}
\end{equation}
Acting the summation operator $\displaystyle{ \sum_{\omega_0 =0}^{\infty } \frac{s_0^{\omega _0}}{ \omega _0 !} \frac{\Gamma (\omega _0+\gamma')}{\Gamma(\gamma') }  \prod _{n=1}^{\infty } \left\{ \sum_{ \omega _n = \omega _{n-1}}^{\infty } s_n^{\omega _n }\right\}}$ on (\ref{eq:50021}) where $|s_i|<1$ as $i=0,1,2,\cdots$ by using (\ref{eq:50033}) and (\ref{eq:50034}),
\begin{theorem}
The general expression of the generating function for the GCH polynomial of type 2 about $x=0$ is given by
\begin{eqnarray}
&&\sum_{\omega _0 =0}^{\infty } \frac{s_0^{\omega _0}}{ \omega _0 !} \frac{\Gamma (\omega _0+\gamma')}{\Gamma(\gamma') }  \prod _{n=1}^{\infty } \left\{ \sum_{ \omega _n = \omega _{n-1}}^{\infty } s_n^{\omega _n }\right\} y(x) \nonumber\\
&&= \prod_{k=1}^{\infty } \frac{1}{(1-s_{k,\infty })} \mathbf{\Upsilon}(\lambda; s_{0,\infty } ;\eta )  \nonumber\\
&&+ \Bigg\{ \prod_{k=1}^{\infty } \frac{1}{(1-s_{k,\infty })} \int_{0}^{1} dt_1\;t_1^{1+ \lambda } \int_{0}^{1} du_1\;u_1^{\nu + \lambda } \exp\left(-\frac{s_{1,\infty }}{(1-s_{1,\infty })}\eta (1-t_1)(1-u_1)\right)\nonumber\\
&&\times \widetilde{w}_{1,1}^{-\left(\Omega /\mu +\lambda \right)} \left(  \widetilde{w}_{1,1} \partial _{\widetilde{w}_{1,1}} \right) \widetilde{w}_{1,1}^{ \Omega /\mu  +\lambda }\; \mathbf{\Upsilon}(\lambda ; s_0; \widetilde{w}_{1,1}) \Bigg\} \rho  \nonumber\\
&&+ \sum_{n=2}^{\infty } \Bigg\{ \prod_{k=n}^{\infty } \frac{1}{(1-s_{k,\infty })} \int_{0}^{1} dt_n\;t_n^{2n-1+ \lambda } \int_{0}^{1} du_n \;u_n^{2(n-1)+\nu + \lambda } \exp\left(-\frac{s_{n,\infty }}{(1-s_{n,\infty })}\eta (1-t_n)(1-u_n)\right)\nonumber\\
&&\times \widetilde{w}_{n,n}^{-\left(\Omega /\mu +2(n-1)+ \lambda \right)} \left(  \widetilde{w}_{n,n} \partial _{\widetilde{w}_{n,n}} \right) \widetilde{w}_{n,n}^{ \Omega /\mu +2(n-1)+ \lambda} \nonumber\\
&&\times  \prod_{j=1}^{n-1} \Bigg\{ \int_{0}^{1} dt_{n-j}\;t_{n-j}^{2(n-j)-1+ \lambda } \int_{0}^{1} du_{n-j} \;u_{n-j}^{2(n-j-1)+\nu + \lambda }\frac{\exp\left(-\frac{s_{n-j}}{(1-s_{n-j})}\widetilde{w}_{n-j+1,n}(1-t_{n-j})(1-u_{n-j})\right)}{(1-s_{n-j})} \nonumber\\
&&\times \widetilde{w}_{n-j,n}^{-\left(\Omega /\mu +2(n-j-1)+ \lambda \right)} \left( \widetilde{w}_{n-j,n} \partial _{\widetilde{w}_{n-j,n}} \right) \widetilde{w}_{n-j,n}^{ \Omega /\mu +2(n-j-1)+ \lambda} \Bigg\}
 \mathbf{\Upsilon}(\lambda; s_0 ;\widetilde{w}_{1,n}) \Bigg\} \rho ^n  \label{eq:50035}
\end{eqnarray}
where
\begin{equation}
\begin{cases} 
{ \displaystyle \mathbf{\Upsilon}(\lambda; s_{0,\infty } ;\eta)= \sum_{\omega _0=0}^{\infty } \frac{s_{0,\infty }^{\omega _0}}{\omega_0!} \frac{\Gamma (\omega _0+\gamma')}{\Gamma(\gamma') } \left\{  c_0 x^{\lambda }  \sum_{i_0=0}^{\omega_0} \frac{(-\omega_0)_{i_0}}{(1+ \lambda )_{i_0}(\nu + \lambda )_{i_0}} \eta^{i_0} \right\} }
\cr
{ \displaystyle \mathbf{\Upsilon}(\lambda ; s_0;\widetilde{w}_{1,1}) =  \sum_{\omega_0=0}^{\infty }\frac{s_0^{\omega_0}}{\omega _0!}\frac{\Gamma (\omega _0+\gamma')}{\Gamma(\gamma') } \left\{  c_0 x^{\lambda } \sum_{i_0=0}^{\omega_0} \frac{(-\omega_0)_{i_0}}{(1+ \lambda )_{i_0}(\nu + \lambda )_{i_0}} \widetilde{w}_{1,1}^{i_0}\right\} } \cr
{ \displaystyle \mathbf{\Upsilon}(\lambda; s_0 ;\widetilde{w}_{1,n}) = \sum_{\omega_0 =0}^{\infty }\frac{s_0^{\omega_0}}{\omega _0!} \frac{\Gamma (\omega _0+\gamma')}{\Gamma(\gamma') } \left\{  c_0 x^{\lambda } \sum_{i_0=0}^{\omega _0}\frac{(-\omega_0)_{i_0}}{(1+ \lambda )_{i_0}(\nu + \lambda )_{i_0}} \widetilde{w}_{1,n}^{i_0}\right\}}
\end{cases}\nonumber 
\end{equation}
\end{theorem}
\begin{proof} 
Acting the summation operator  $\displaystyle{ \sum_{\omega_0 =0}^{\infty } \frac{s_0^{\omega _0}}{ \omega _0 !} \frac{\Gamma (\omega _0+\gamma')}{\Gamma(\gamma') }  \prod _{n=1}^{\infty } \left\{ \sum_{ \omega _n = \omega _{n-1}}^{\infty } s_n^{\omega _n }\right\}}$ on the form of integral of the type 2 GCH polynomial $y(x)$,
\begin{eqnarray}
&& \sum_{\omega_0 =0}^{\infty } \frac{s_0^{\omega _0}}{ \omega _0 !} \frac{\Gamma (\omega _0+\gamma')}{\Gamma(\gamma') }  \prod _{n=1}^{\infty } \left\{ \sum_{ \omega _n = \omega _{n-1}}^{\infty } s_n^{\omega _n }\right\} y(x) \nonumber\\
&&=  \sum_{\omega_0 =0}^{\infty } \frac{s_0^{\omega _0}}{ \omega _0 !} \frac{\Gamma (\omega _0+\gamma')}{\Gamma(\gamma') }  \prod _{n=1}^{\infty } \left\{ \sum_{ \omega _n = \omega _{n-1}}^{\infty } s_n^{\omega _n }\right\} \Big( y_0(x)+y_1(x)+y_2(x)+y_3(x)+ \cdots \Big)\hspace{1.5cm} \label{eq:50036}
\end{eqnarray}
Acting the summation operator $\displaystyle{ \sum_{\omega_0 =0}^{\infty } \frac{s_0^{\omega _0}}{ \omega _0 !} \frac{\Gamma (\omega _0+\gamma')}{\Gamma(\gamma') }  \prod _{n=1}^{\infty } \left\{ \sum_{ \omega _n = \omega _{n-1}}^{\infty } s_n^{\omega _n }\right\}}$ on (\ref{eq:50023a}) by using (\ref{eq:50033}) and (\ref{eq:50034}),
\begin{eqnarray}
&&\sum_{\omega_0 =0}^{\infty } \frac{s_0^{\omega _0}}{ \omega _0 !} \frac{\Gamma (\omega _0+\gamma')}{\Gamma(\gamma') }  \prod _{n=1}^{\infty } \left\{ \sum_{ \omega _n = \omega _{n-1}}^{\infty } s_n^{\omega _n }\right\} y_0(x)\nonumber\\
&& = \prod_{k=1}^{\infty } \frac{1}{(1-s_{k,\infty })} \sum_{\omega_0 =0}^{\infty } \frac{s_{0,\infty }^{\omega_0}}{\omega_0!}\frac{\Gamma (\omega _0+\gamma')}{\Gamma(\gamma') } \Bigg\{ c_0 x^{\lambda } \sum_{i_0=0}^{\omega _0} \frac{(-\omega _0)_{i_0}}{(1+ \lambda )_{i_0}(\nu + \lambda )_{i_0}}\eta^{i_0} \Bigg\} \hspace{.5cm}\label{eq:50037}
\end{eqnarray}
Acting the summation operator $\displaystyle{ \sum_{\omega_0 =0}^{\infty } \frac{s_0^{\omega _0}}{ \omega _0 !} \frac{\Gamma (\omega _0+\gamma')}{\Gamma(\gamma') }  \prod _{n=1}^{\infty } \left\{ \sum_{ \omega _n = \omega _{n-1}}^{\infty } s_n^{\omega _n }\right\}}$ on (\ref{eq:50024}) by using (\ref{eq:50033}) and (\ref{eq:50034}),
\begin{eqnarray}
&&\sum_{\omega_0 =0}^{\infty } \frac{s_0^{\omega _0}}{ \omega _0 !} \frac{\Gamma (\omega _0+\gamma')}{\Gamma(\gamma') }  \prod _{n=1}^{\infty } \left\{ \sum_{ \omega _n = \omega _{n-1}}^{\infty } s_n^{\omega _n }\right\} y_1(x) \nonumber\\
&&=  \prod_{k=2}^{\infty } \frac{1}{(1-s_{k,\infty })} \int_{0}^{1} dt_1\;t_1^{1+ \lambda } \int_{0}^{1} du_1\;u_1^{\nu + \lambda } 
  \frac{1}{2\pi i}  \oint dv_1 \frac{\exp\left(-\frac{v_1}{(1-v_1)}\eta (1-t_1)(1-u_1)\right)}{v_1(1-v_1)}\nonumber\\
&&\times \sum_{\omega _1=\omega _0}^{\infty }\left( \frac{s_{1,\infty }}{v_1}\right)^{\omega _1}  w_{1,1}^{-\left(\Omega /\mu + \lambda \right)} \left(  w_{1,1} \partial _{w_{1,1}} \right) w_{1,1}^{ \Omega /\mu + \lambda}\nonumber\\
&&\times \sum_{\omega  _0=0}^{\infty }\frac{s_0^{\omega _0}}{\omega _0!} \frac{\Gamma (\omega _0+\gamma')}{\Gamma(\gamma') } \Bigg\{ c_0 x^{\lambda} \sum_{i_0=0}^{\omega _0} \frac{(-\omega _0)_{i_0}}{(1+ \lambda )_{i_0}(\nu + \lambda )_{i_0}} w_{1,1}^{i_0}\Bigg\} \rho  \label{eq:50038}
\end{eqnarray}
Replace $\omega_i$, $\omega _j$ and $s_i$ by $\omega_1$, $\omega _0$ and ${ \displaystyle \frac{s_{1,\infty }}{v_1}}$ in (\ref{eq:50034}). Take the new (\ref{eq:50034}) into (\ref{eq:50038}).
\begin{eqnarray}
&&\sum_{\omega_0 =0}^{\infty } \frac{s_0^{\omega _0}}{ \omega _0 !} \frac{\Gamma (\omega _0+\gamma')}{\Gamma(\gamma') }  \prod _{n=1}^{\infty } \left\{ \sum_{ \omega _n = \omega _{n-1}}^{\infty } s_n^{\omega _n }\right\} y_1(x) \nonumber\\
&&=  \prod_{k=2}^{\infty } \frac{1}{(1-s_{k,\infty })} \int_{0}^{1} dt_1\;t_1^{1+ \lambda } \int_{0}^{1} du_1\;u_1^{\nu + \lambda } 
  \frac{1}{2\pi i}  \oint dv_1 \frac{\exp\left(-\frac{v_1}{(1-v_1)}\eta (1-t_1)(1-u_1)\right)}{(1-v_1)( v_1-s_{1,\infty })}   \nonumber\\
&&\times  w_{1,1}^{-\left(\Omega /\mu + \lambda \right)} \left(  w_{1,1} \partial _{w_{1,1}} \right) w_{1,1}^{ \Omega /\mu + \lambda}  \nonumber\\
&&\times \sum_{\omega  _0=0}^{\infty }\frac{1}{\omega _0!} \left( \frac{s_{0,\infty }}{v_1}\right)^{\omega _0} \frac{\Gamma (\omega _0+\gamma')}{\Gamma(\gamma') } \Bigg\{ c_0 x^{\lambda} \sum_{i_0=0}^{\omega _0} \frac{(-\omega _0)_{i_0}}{(1+ \lambda )_{i_0}(\nu + \lambda )_{i_0}} w_{1,1}^{i_0}\Bigg\} \rho  \label{eq:50039}
\end{eqnarray}
By using Cauchy's integral formula, the contour integrand has poles at $v_1 =1$ or $s_{1,\infty}$,
and $s_{1,\infty}$ is only inside the unit circle. As we compute the residue there in (\ref{eq:50039}) we obtain
\begin{eqnarray}
&&\sum_{\omega_0 =0}^{\infty } \frac{s_0^{\omega _0}}{ \omega _0 !} \frac{\Gamma (\omega _0+\gamma')}{\Gamma(\gamma') }  \prod _{n=1}^{\infty } \left\{ \sum_{ \omega _n = \omega _{n-1}}^{\infty } s_n^{\omega _n }\right\} y_1(x) \nonumber\\
&&=  \prod_{k=1}^{\infty } \frac{1}{(1-s_{k,\infty })} \int_{0}^{1} dt_1\;t_1^{1+ \lambda } \int_{0}^{1} du_1\;u_1^{\nu + \lambda } 
    \exp\left(-\frac{s_{1,\infty}}{(1-s_{1,\infty})}\eta (1-t_1)(1-u_1)\right) \nonumber\\
&&\times  \widetilde{w}_{1,1}^{-\left(\Omega /\mu + \lambda \right)} \left(\widetilde{w}_{1,1} \partial _{\widetilde{w}_{1,1}} \right) \widetilde{w}_{1,1}^{ \Omega /\mu + \lambda} \sum_{\omega  _0=0}^{\infty }\frac{s_0^{\omega _0}}{\omega _0!} \frac{\Gamma (\omega _0+\gamma')}{\Gamma(\gamma') } \Bigg\{ c_0 x^{\lambda} \sum_{i_0=0}^{\omega _0} \frac{(-\omega _0)_{i_0}}{(1+ \lambda )_{i_0}(\nu + \lambda )_{i_0}} \widetilde{w}_{1,1}^{i_0}\Bigg\} \rho \hspace{1.5cm} \label{eq:50040}
\end{eqnarray}
where
\begin{eqnarray}
\widetilde{w}_{1,1} &=& \eta s_{1,\infty } \prod _{l=1}^{1} t_l u_l\nonumber
\end{eqnarray}
Acting the summation operator $\displaystyle{ \sum_{\omega_0 =0}^{\infty } \frac{s_0^{\omega _0}}{ \omega _0 !} \frac{\Gamma (\omega _0+\gamma')}{\Gamma(\gamma') }  \prod _{n=1}^{\infty } \left\{ \sum_{ \omega _n = \omega _{n-1}}^{\infty } s_n^{\omega _n }\right\}}$ on (\ref{eq:50026}) by using (\ref{eq:50033}) and (\ref{eq:50034}),
\begin{eqnarray}
&&\sum_{\omega_0 =0}^{\infty } \frac{s_0^{\omega _0}}{ \omega _0 !} \frac{\Gamma (\omega _0+\gamma')}{\Gamma(\gamma') }  \prod _{n=1}^{\infty } \left\{ \sum_{ \omega _n = \omega _{n-1}}^{\infty } s_n^{\omega _n }\right\} y_2(x) \nonumber\\
&&=  \prod_{k=3}^{\infty } \frac{1}{(1-s_{k,\infty })} \int_{0}^{1} dt_2\;t_2^{3+ \lambda } \int_{0}^{1} du_2\;u_2^{2+\nu + \lambda } \frac{1}{2\pi i}  \oint dv_2 \frac{\exp\left(-\frac{v_2}{(1-v_2)}\eta (1-t_2)(1-u_2)\right)}{v_2(1-v_2)}\nonumber\\
&&\times \sum_{\omega _2=\omega  _1}^{\infty }\left( \frac{s_{2,\infty }}{v_2}\right)^{\omega _2} w_{2,2}^{-\left(\Omega /\mu +2+ \lambda \right)} \left(  w_{2,2} \partial _{w_{2,2}} \right) w_{2,2}^{ \Omega /\mu +2+ \lambda}\nonumber\\
&&\times  \int_{0}^{1} dt_1\;t_1^{1+ \lambda } \int_{0}^{1} du_1\;u_1^{\nu + \lambda } 
\frac{1}{2\pi i}  \oint dv_1 \frac{\exp\left(-\frac{v_1}{(1-v_1)}w_{2,2}(1-t_1)(1-u_1)\right)}{v_1 (1-v_1)} \nonumber\\
&&\times \sum_{\omega  _1=\omega _0}^{\infty } \left( \frac{s_1}{v_1}\right)^{\omega _1} w_{1,2}^{-\left(\Omega /\mu  + \lambda \right)} \left(  w_{1,2} \partial _{w_{1,2}} \right) w_{1,2}^{ \Omega /\mu + \lambda}\nonumber\\
&&\times \sum_{\omega _0=0}^{\infty }\frac{s_0^{\omega _0}}{\omega _0!}\frac{\Gamma (\omega _0+\gamma')}{\Gamma(\gamma') }\Bigg\{c_0 x^{\lambda} \sum_{i_0=0}^{\omega_0} \frac{(-\omega_0)_{i_0}}{(1+ \lambda )_{i_0}(\nu + \lambda )_{i_0}}w_{1,2}^{i_0} \Bigg\} \rho ^2  \label{eq:50041}
\end{eqnarray}
Replace $\omega_i$, $\omega_j$ and $s_i$ by $\omega_2$, $\omega_1$ and ${ \displaystyle \frac{s_{2,\infty }}{v_2}}$ in (\ref{eq:50034}). Take the new (\ref{eq:50034}) into (\ref{eq:50041}).
\begin{eqnarray}
&&\sum_{\omega_0 =0}^{\infty } \frac{s_0^{\omega _0}}{ \omega _0 !} \frac{\Gamma (\omega _0+\gamma')}{\Gamma(\gamma') }  \prod _{n=1}^{\infty } \left\{ \sum_{ \omega _n = \omega _{n-1}}^{\infty } s_n^{\omega _n }\right\} y_2(x) \nonumber\\
&&=  \prod_{k=3}^{\infty } \frac{1}{(1-s_{k,\infty })} \int_{0}^{1} dt_2\;t_2^{3+ \lambda } \int_{0}^{1} du_2\;u_2^{2+\nu + \lambda } \frac{1}{2\pi i}  \oint dv_2 \frac{\exp\left(-\frac{v_2}{(1-v_2)}\eta (1-t_2)(1-u_2)\right)}{ (1-v_2)( v_2-s_{2,\infty })}\nonumber\\
&&\times w_{2,2}^{-\left(\Omega /\mu +2+ \lambda \right)} \left(  w_{2,2} \partial _{w_{2,2}} \right) w_{2,2}^{ \Omega /\mu +2+ \lambda}\nonumber\\
&&\times  \int_{0}^{1} dt_1\;t_1^{1+ \lambda } \int_{0}^{1} du_1\;u_1^{\nu + \lambda } 
\frac{1}{2\pi i}  \oint dv_1 \frac{\exp\left(-\frac{v_1}{(1-v_1)}w_{2,2}(1-t_1)(1-u_1)\right)}{v_1 (1-v_1)} \nonumber\\
&&\times \sum_{\omega  _1=\omega _0}^{\infty } \left( \frac{s_{1,\infty }}{v_1 v_2}\right)^{\omega _1} w_{1,2}^{-\left(\Omega /\mu  + \lambda \right)} \left(  w_{1,2} \partial _{w_{1,2}} \right) w_{1,2}^{ \Omega /\mu + \lambda}\nonumber\\
&&\times \sum_{\omega _0=0}^{\infty }\frac{s_0^{\omega _0}}{\omega _0!}\frac{\Gamma (\omega _0+\gamma')}{\Gamma(\gamma') }\Bigg\{c_0 x^{\lambda} \sum_{i_0=0}^{\omega_0} \frac{(-\omega_0)_{i_0}}{(1+ \lambda )_{i_0}(\nu + \lambda )_{i_0}}w_{1,2}^{i_0} \Bigg\} \rho ^2  \label{eq:50042}
\end{eqnarray}
By using Cauchy's integral formula, the contour integrand has poles at $v_2 =1$ or $s_{2,\infty}$,
and $s_{2,\infty}$ is only inside the unit circle. As we compute the residue there in (\ref{eq:50042}) we obtain
\begin{eqnarray}
&&\sum_{\omega_0 =0}^{\infty } \frac{s_0^{\omega _0}}{ \omega _0 !} \frac{\Gamma (\omega _0+\gamma')}{\Gamma(\gamma') }  \prod _{n=1}^{\infty } \left\{ \sum_{ \omega _n = \omega _{n-1}}^{\infty } s_n^{\omega _n }\right\} y_2(x) \nonumber\\
&&=   \prod_{k=2}^{\infty } \frac{1}{(1-s_{k,\infty })} \int_{0}^{1} dt_2\;t_2^{3+ \lambda } \int_{0}^{1} du_2\;u_2^{2+ \nu + \lambda } \exp\left(-\frac{s_{2,\infty }}{(1-s_{2,\infty })}\eta (1-t_2)(1-u_2)\right)\nonumber\\
&&\times   \widetilde{w}_{2,2}^{-\left(\Omega /\mu +2+ \lambda \right)} \left(\widetilde{w}_{2,2} \partial _{\widetilde{w}_{2,2}} \right) \widetilde{w}_{2,2}^{ \Omega /\mu +2+ \lambda}\nonumber\\
&&\times  \int_{0}^{1} dt_1\;t_1^{ 1 + \lambda } \int_{0}^{1} du_1\;u_1^{\nu + \lambda } 
\frac{1}{2\pi i}  \oint dv_1 \frac{\exp\left(-\frac{v_1}{(1-v_1)}\widetilde{w}_{2,2} (1-t_1)(1-u_1)\right)}{v_1(1-v_1)} \sum_{\omega _1=\omega _0}^{\infty } \left( \frac{s_1}{v_1}\right)^{\omega_1}\nonumber \\
&&\times \ddot{w}_{1,2}^{-\left(\Omega /\mu + \lambda \right)} \left(\ddot{w}_{1,2} \partial _{\ddot{w}_{1,2}} \right) \ddot{w}_{1,2}^{ \Omega /\mu + \lambda} \nonumber\\
&&\times \sum_{\omega _0=0}^{\infty }\frac{s_0^{\omega_0}}{\omega _0!} \frac{\Gamma (\omega _0+\gamma')}{\Gamma(\gamma') }\Bigg\{ c_0 x^{\lambda}\sum_{i_0=0}^{\omega_0} \frac{(-\omega_0)_{i_0}}{(1+ \lambda )_{i_0}(\nu + \lambda )_{i_0}}\ddot{w}_{1,2}^{i_0}\Bigg\}\rho ^2  \label{eq:50043}
\end{eqnarray}
where
\begin{eqnarray}
\widetilde{w}_{2,2} &=& \eta s_{2,\infty } \prod _{l=2}^{2} t_l u_l\hspace{2cm}\ddot{w}_{1,2} = \eta s_{2,\infty } v_1\prod _{l=1}^{2} t_l u_l\nonumber
\end{eqnarray}
Replace $\omega _i$, $\omega _j$ and $s_i$ by $\omega _1$, $\omega _0$ and ${ \displaystyle \frac{s_1}{v_1}}$ in (\ref{eq:50034}). Take the new (\ref{eq:50034}) into (\ref{eq:50043}).
\begin{eqnarray}
&&\sum_{\omega_0 =0}^{\infty } \frac{s_0^{\omega _0}}{ \omega _0 !} \frac{\Gamma (\omega _0+\gamma')}{\Gamma(\gamma') }  \prod _{n=1}^{\infty } \left\{ \sum_{ \omega _n = \omega _{n-1}}^{\infty } s_n^{\omega _n }\right\} y_2(x) \nonumber\\
&&=   \prod_{k=2}^{\infty } \frac{1}{(1-s_{k,\infty })} \int_{0}^{1} dt_2\;t_2^{3+ \lambda } \int_{0}^{1} du_2\;u_2^{2+ \nu + \lambda } \exp\left(-\frac{s_{2,\infty }}{(1-s_{2,\infty })}\eta (1-t_2)(1-u_2)\right)\nonumber\\
&&\times   \widetilde{w}_{2,2}^{-\left(\Omega /\mu +2+ \lambda \right)} \left(\widetilde{w}_{2,2} \partial _{\widetilde{w}_{2,2}} \right) \widetilde{w}_{2,2}^{ \Omega /\mu +2+ \lambda}\nonumber\\
&&\times  \int_{0}^{1} dt_1\;t_1^{ 1 + \lambda } \int_{0}^{1} du_1\;u_1^{\nu + \lambda } 
\frac{1}{2\pi i}  \oint dv_1 \frac{\exp\left(-\frac{v_1}{(1-v_1)}\widetilde{w}_{2,2} (1-t_1)(1-u_1)\right)}{ (1-v_1)( v_1 -s_1)} \nonumber \\
&&\times \ddot{w}_{1,2}^{-\left(\Omega /\mu + \lambda \right)} \left(\ddot{w}_{1,2} \partial _{\ddot{w}_{1,2}} \right) \ddot{w}_{1,2}^{ \Omega /\mu + \lambda}\nonumber\\
&&\times  \sum_{\omega _0=0}^{\infty }\frac{1}{\omega _0!} \left( \frac{s_{0,1}}{v_1}\right)^{\omega_0} \frac{\Gamma (\omega _0+\gamma')}{\Gamma(\gamma') }\Bigg\{ c_0 x^{\lambda}\sum_{i_0=0}^{\omega_0} \frac{(-\omega_0)_{i_0}}{(1+ \lambda )_{i_0}(\nu + \lambda )_{i_0}}\ddot{w}_{1,2}^{i_0}\Bigg\}\rho ^2 \label{eq:50044}
\end{eqnarray}
By using Cauchy's integral formula, the contour integrand has poles at $v_1 =1$ or $s_1$,
and $s_1$ is only inside the unit circle. As we compute the residue there in (\ref{eq:50044}) we obtain
\begin{eqnarray}
&&\sum_{\omega_0 =0}^{\infty } \frac{s_0^{\omega _0}}{ \omega _0 !} \frac{\Gamma (\omega _0+\gamma')}{\Gamma(\gamma') }  \prod _{n=1}^{\infty } \left\{ \sum_{ \omega _n = \omega _{n-1}}^{\infty } s_n^{\omega _n }\right\} y_2(x) \nonumber\\
&&=   \prod_{k=2}^{\infty } \frac{1}{(1-s_{k,\infty })} \int_{0}^{1} dt_2\;t_2^{3+ \lambda } \int_{0}^{1} du_2\;u_2^{2+ \nu + \lambda } \exp\left(-\frac{s_{2,\infty }}{(1-s_{2,\infty })}\eta (1-t_2)(1-u_2)\right)\nonumber\\
&&\times   \widetilde{w}_{2,2}^{-\left(\Omega /\mu +2+ \lambda \right)} \left(\widetilde{w}_{2,2} \partial _{\widetilde{w}_{2,2}} \right) \widetilde{w}_{2,2}^{ \Omega /\mu +2+ \lambda}\nonumber\\
&&\times  \int_{0}^{1} dt_1\;t_1^{ 1 + \lambda } \int_{0}^{1} du_1\;u_1^{\nu + \lambda } 
 \frac{\exp\left(-\frac{s_1}{(1-s_1)}\widetilde{w}_{2,2}(1-t_1)(1-u_1)\right)}{(1-s_1)} \widetilde{w}_{1,2}^{-\left(\Omega /\mu + \lambda \right)} \left(\widetilde{w}_{1,2} \partial _{\widetilde{w}_{1,2}} \right) \widetilde{w}_{1,2}^{ \Omega /\mu + \lambda} \nonumber\\
&&\times  \sum_{\omega _0=0}^{\infty }\frac{s_0^{\omega _0}}{\omega _0!}  \frac{\Gamma (\omega _0+\gamma')}{\Gamma(\gamma') }\Bigg\{ c_0 x^{\lambda}\sum_{i_0=0}^{\omega _0} \frac{(-\omega _0)_{i_0}}{(1+ \lambda )_{i_0}(\nu + \lambda )_{i_0}} \widetilde{w}_{1,2}^{i_0}\Bigg\}\rho ^2 \label{eq:50045}
\end{eqnarray}
where
\begin{eqnarray}
\widetilde{w}_{1,2} &=& \eta s_{1,\infty } \prod _{l=1}^{2} t_l u_l\nonumber
\end{eqnarray}
Acting the summation operator $\displaystyle{ \sum_{\omega_0 =0}^{\infty } \frac{s_0^{\omega _0}}{ \omega _0 !} \frac{\Gamma (\omega _0+\gamma')}{\Gamma(\gamma') }  \prod _{n=1}^{\infty } \left\{ \sum_{ \omega _n = \omega _{n-1}}^{\infty } s_n^{\omega _n }\right\}}$ on (\ref{eq:50045}) by using (\ref{eq:50033}) and (\ref{eq:50034}),
\begin{eqnarray}
&&\sum_{\omega_0 =0}^{\infty } \frac{s_0^{\omega _0}}{ \omega _0 !} \frac{\Gamma (\omega _0+\gamma')}{\Gamma(\gamma') }  \prod _{n=1}^{\infty } \left\{ \sum_{ \omega _n = \omega _{n-1}}^{\infty } s_n^{\omega _n }\right\} y_3(x) \nonumber\\
&&= \prod_{k=3}^{\infty } \frac{1}{(1-s_{k,\infty })} \int_{0}^{1} dt_3\;t_3^{5+ \lambda } \int_{0}^{1} du_3\;u_3^{4+ \nu + \lambda } \exp\left(-\frac{s_{3,\infty }}{(1-s_{3,\infty })}\eta (1-t_3)(1-u_3)\right)\nonumber\\
&&\times   \widetilde{w}_{3,3}^{-\left(\Omega /\mu +4+ \lambda \right)} \left(\widetilde{w}_{3,3} \partial _{\widetilde{w}_{3,3}} \right) \widetilde{w}_{3,3}^{ \Omega /\mu +4+ \lambda}\nonumber\\
&&\times  \int_{0}^{1} dt_2\;t_2^{3 + \lambda } \int_{0}^{1} du_2\;u_2^{2+ \nu + \lambda } 
 \frac{\exp\left(-\frac{s_2}{(1-s_2)}\widetilde{w}_{3,3}(1-t_2)(1-u_2)\right)}{(1-s_2)} \widetilde{w}_{2,3}^{-\left(\Omega /\mu +2+ \lambda \right)} \left(\widetilde{w}_{2,3} \partial _{\widetilde{w}_{2,3}} \right) \widetilde{w}_{2,3}^{ \Omega /\mu +2+ \lambda} \nonumber\\
 &&\times  \int_{0}^{1} dt_1\;t_1^{1 + \lambda } \int_{0}^{1} du_1\;u_1^{ \nu + \lambda } 
 \frac{\exp\left(-\frac{s_1}{(1-s_1)}\widetilde{w}_{2,3}(1-t_1)(1-u_1)\right)}{(1-s_1)} \widetilde{w}_{1,3}^{-\left(\Omega /\mu + \lambda \right)} \left(\widetilde{w}_{1,3} \partial _{\widetilde{w}_{1,3}} \right) \widetilde{w}_{1,3}^{ \Omega /\mu + \lambda} \nonumber\\
&&\times  \sum_{\omega _0=0}^{\infty }\frac{s_0^{\omega _0}}{\omega _0!}  \frac{\Gamma (\omega _0+\gamma')}{\Gamma(\gamma') }\Bigg\{ c_0 x^{\lambda}\sum_{i_0=0}^{\omega _0} \frac{(-\omega _0)_{i_0}}{(1+ \lambda )_{i_0}(\nu + \lambda )_{i_0}} \widetilde{w}_{1,3}^{i_0}\Bigg\}\rho ^3 \label{eq:50046}
\end{eqnarray}
where
\begin{equation}
\widetilde{w}_{3,3} = \eta s_{3,\infty } \prod _{l=3}^{3} t_l u_l \hspace{1cm} \widetilde{w}_{2,3} = \eta s_{2,\infty } \prod _{l=2}^{3} t_l u_l\hspace{1cm} \widetilde{w}_{1,3} = \eta s_{1,\infty } \prod _{l=1}^{3} t_l u_l\nonumber
\end{equation}
By repeating this process for all higher terms of integral forms of sub-summation $y_m(x)$ terms where $m > 3$, we obtain every  $\displaystyle{ \sum_{\omega_0 =0}^{\infty } \frac{s_0^{\omega _0}}{ \omega _0 !} \frac{\Gamma (\omega _0+\gamma')}{\Gamma(\gamma') }  \prod _{n=1}^{\infty } \left\{ \sum_{ \omega _n = \omega _{n-1}}^{\infty } s_n^{\omega _n }\right\}  y_m(x)}$ terms. 
Since we substitute (\ref{eq:50037}), (\ref{eq:50040}), (\ref{eq:50045}), (\ref{eq:50046}) and including all $\displaystyle{ \sum_{\omega_0 =0}^{\infty } \frac{s_0^{\omega _0}}{ \omega _0 !} \frac{\Gamma (\omega _0+\gamma')}{\Gamma(\gamma') }  \prod _{n=1}^{\infty } \left\{ \sum_{ \omega _n = \omega _{n-1}}^{\infty } s_n^{\omega _n }\right\}  y_m(x)}$ terms where $m > 3$ into (\ref{eq:50036}), we obtain (\ref{eq:50035})
\qed
\end{proof}
\begin{remark}
The generating function for the GCH polynomial of type 2 of the first kind about $x=0$ as $\omega =-(\omega_j +2 j)$ where $j,\omega _j= 0,1,2,\cdots$ is
\begin{eqnarray}
&&\sum_{\omega_0 =0}^{\infty } \frac{s_0^{\omega _0}}{ \omega _0 !} \frac{\Gamma (\omega _0+\nu)}{\Gamma(\nu) }  \prod _{n=1}^{\infty } \left\{ \sum_{ \omega _n = \omega _{n-1}}^{\infty } s_n^{\omega _n }\right\} QW_{\omega _j}^R\left(\mu ,\varepsilon ,\nu ,\Omega ,\omega =-(\omega_j +2 j); \rho =-\mu x^2; \eta = -\varepsilon x \right)\nonumber\\
&&=  \prod_{k=1}^{\infty } \frac{1}{(1-s_{k,\infty })} \displaystyle \mathbf{A} \left( s_{0,\infty } ;\eta \right)\nonumber\\
&&+ \Bigg\{ \prod_{k=1}^{\infty } \frac{1}{(1-s_{k,\infty })} \int_{0}^{1} dt_1\;t_1 \int_{0}^{1} du_1\;u_1^{\nu } \overleftrightarrow {\mathbf{\Gamma}}_1 \left(s_{1,\infty };t_1,u_1,\eta \right) \widetilde{w}_{1,1}^{- \Omega /\mu } \left(\widetilde{w}_{1,1} \partial _{\widetilde{w}_{1,1}} \right) \widetilde{w}_{1,1}^{ \Omega /\mu } \mathbf{A} \left(s_{0} ; \widetilde{w}_{1,1}\right) \Bigg\} \rho \nonumber\\
&&+ \sum_{n=2}^{\infty } \Bigg\{ \prod_{k=n}^{\infty } \frac{1}{(1-s_{k,\infty })} \int_{0}^{1} dt_n\;t_n^{2n-1} \int_{0}^{1} du_n \;u_n^{2(n-1) +\nu}   \overleftrightarrow {\mathbf{\Gamma}}_n \left(s_{n,\infty };t_n,u_n,\eta\right) \nonumber\\
&&\times \widetilde{w}_{n,n}^{-\left( \Omega /\mu +2(n-1)\right) } \left(\widetilde{w}_{n,n} \partial _{\widetilde{w}_{n,n}} \right) \widetilde{w}_{n,n}^{ \Omega /\mu +2(n-1)}\nonumber\\
&&\times  \prod_{j=1}^{n-1} \Bigg\{ \int_{0}^{1} dt_{n-j}\;t_{n-j}^{2(n-j)-1} \int_{0}^{1} du_{n-j} \;u_{n-j}^{2(n-j-1)+\nu}  \overleftrightarrow {\mathbf{\Gamma}}_{n-j} \left(s_{n-j};t_{n-j},u_{n-j},\widetilde{w}_{n-j+1,n} \right) \nonumber\\
&&\times \widetilde{w}_{n-j,n}^{-\left( \Omega /\mu +2(n-j-1)\right) } \left(\widetilde{w}_{n-j,n} \partial _{\widetilde{w}_{n-j,n}} \right) \widetilde{w}_{n-j,n}^{ \Omega /\mu +2(n-j-1)} \Bigg\}
  \mathbf{A} \left(s_{0} ;\widetilde{w}_{1,n} \right)  \Bigg\} \rho ^n 
 \label{eq:50047}
\end{eqnarray}
where
\begin{equation}
\begin{cases} 
{ \displaystyle \overleftrightarrow {\mathbf{\Gamma}}_1 \left(s_{1,\infty };t_1,u_1,\eta \right)=  \exp \left(-\frac{s_{1,\infty }}{(1-s_{1,\infty })}\eta (1-t_1)(1-u_1)\right) }\cr
{ \displaystyle  \overleftrightarrow {\mathbf{\Gamma}}_n \left(s_{n,\infty };t_n,u_n,\eta \right) = \exp\left( -\frac{s_{n,\infty }}{(1-s_{n,\infty })}\eta (1-t_n)(1-u_n)\right)}\cr
{ \displaystyle \overleftrightarrow {\mathbf{\Gamma}}_{n-j} \left(s_{n-j};t_{n-j},u_{n-j},\widetilde{w}_{n-j+1,n} \right) = \frac{\exp\left(-\frac{s_{n-j}}{(1-s_{n-j})}\widetilde{w}_{n-j+1,n} (1-t_{n-j})(1-u_{n-j})\right)}{(1-s_{n-j})}}
\end{cases}\nonumber 
\end{equation}
and
\begin{equation}
\begin{cases} 
{ \displaystyle \mathbf{A} \left( s_{0,\infty } ;\eta \right)= (1-s_{0,\infty })^{-\nu }\exp\left(-\frac{\eta s_{0,\infty }}{(1-s_{0,\infty })}\right)}\cr
{ \displaystyle  \mathbf{A} \left(s_{0} ;\widetilde{w}_{1,1} \right) = (1-s_0)^{-\nu }\exp\left(-\frac{\widetilde{w}_{1,1} s_0}{(1-s_0)}\right)} \cr
{ \displaystyle \mathbf{A} \left(s_{0} ;\widetilde{w}_{1,n} \right) = (1-s_0)^{-\nu }\exp\left(-\frac{\widetilde{w}_{1,n} s_0}{(1-s_0)}\right)}
\end{cases}\nonumber 
\end{equation}
\end{remark}
\begin{proof}
The generating function for confluent Hypergeometric polynomial of the first kind is given by
\begin{equation}
\sum_{\omega _0=0}^{\infty } \frac{t^{\omega _0}}{\omega _0!} \frac{\Gamma (\omega _0+\gamma )}{\Gamma(\gamma ) }\; _1F_1 \left(-\omega _0; \gamma ; z \right)  = (1-t)^{-\gamma } \exp\left(-\frac{z t}{(1-t)}\right) 
 \label{eq:50048}
\end{equation}
Replace $t$, $\gamma$ and $z$  by $s_{0,\infty }$ , $\nu$ and $\eta$ in (\ref{eq:50048}). 
\begin{equation}
\sum_{\omega _0=0}^{\infty } \frac{s_{0,\infty }^{\omega _0}}{\omega _0!} \frac{\Gamma (\omega _0+\nu )}{\Gamma(\nu ) }\; _1F_1 \left(-\omega _0; \nu ; \eta \right)  = (1-s_{0,\infty })^{-\nu } \exp\left(-\frac{\eta s_{0,\infty }}{(1-s_{0,\infty })}\right)  \label{eq:50049}
\end{equation} 
Replace $t$, $\gamma$ and $z$  by $s_0$ , $\nu$ and $\widetilde{w}_{1,1}$ in (\ref{eq:50048}). 
\begin{equation}
\sum_{\omega _0=0}^{\infty } \frac{s_0^{\omega _0}}{\omega _0!} \frac{\Gamma (\omega _0+\nu )}{\Gamma(\nu ) }\; _1F_1 \left(-\omega _0; \nu ; \widetilde{w}_{1,1} \right)  = (1-s_0)^{-\nu } \exp\left(-\frac{\widetilde{w}_{1,1} s_0}{(1-s_0)}\right)  \label{eq:50050}
\end{equation} 
Replace $t$, $\gamma$ and $z$  by $s_0$ , $\nu$ and $\widetilde{w}_{1,n}$ in (\ref{eq:50048}). 
\begin{equation}
\sum_{\omega _0=0}^{\infty } \frac{s_0^{\omega _0}}{\omega _0!} \frac{\Gamma (\omega _0+\nu )}{\Gamma(\nu ) }\; _1F_1 \left(-\omega _0; \nu ; \widetilde{w}_{1,n} \right)  = (1-s_0)^{-\nu } \exp\left(-\frac{\widetilde{w}_{1,n} s_0}{(1-s_0)}\right)  \label{eq:50051}
\end{equation} 
Put $c_0$= 1, $\lambda $=0 and $\displaystyle{\gamma' =\nu}$ in (\ref{eq:50035}). Substitute (\ref{eq:50049}), (\ref{eq:50050}) and (\ref{eq:50051}) into the new (\ref{eq:50035}).
\qed
\end{proof}
\begin{remark}
The generating function for the GCH polynomial of type 2 of the second kind about $x=0$ as $\omega =-(\omega_j +2 j+1-\nu)$ where $j,\omega _j= 0,1,2,\cdots$ is
\begin{eqnarray}
&&\sum_{\omega_0 =0}^{\infty } \frac{s_0^{\omega _0}}{ \omega _0 !} \frac{\Gamma (\omega _0+2-\nu)}{\Gamma( 2-\nu ) }  \prod _{n=1}^{\infty } \left\{ \sum_{ \omega _n = \omega _{n-1}}^{\infty } s_n^{\omega _n }\right\} RW_{\omega _j}^R\left(\mu ,\varepsilon ,\nu ,\Omega ,\omega =-(\omega_j +2 j+1-\nu)\right. \nonumber\\
&&;\left. \rho =-\mu x^2; \eta = -\varepsilon x \right) \nonumber\\
&&= x^{1-\nu } \Bigg\{ \prod_{k=1}^{\infty } \frac{1}{(1-s_{k,\infty })} \mathbf{B} \left( s_{0,\infty } ;\eta\right)  \nonumber\\
 &&+ \Bigg\{ \prod_{k=1}^{\infty } \frac{1}{(1-s_{k,\infty })} \int_{0}^{1} dt_1\;t_1^{2-\nu } \int_{0}^{1} du_1\;u_1 \overleftrightarrow {\mathbf{\Gamma}}_1 \left(s_{1,\infty };t_1,u_1,\eta\right)\nonumber\\
&&\times  \widetilde{w}_{1,1}^{- (\Omega /\mu +1-\nu ) } \left(\widetilde{w}_{1,1} \partial _{\widetilde{w}_{1,1}} \right) \widetilde{w}_{1,1}^{ \Omega /\mu +1-\nu} \mathbf{B} \left( s_{0} ;\widetilde{w}_{1,1} \right)\Bigg\} \rho \nonumber\\
&&+ \sum_{n=2}^{\infty } \Bigg\{ \prod_{k=n}^{\infty } \frac{1}{(1-s_{k,\infty })} \int_{0}^{1} dt_n\;t_n^{2n-\nu } \int_{0}^{1} du_n \;u_n^{2n-1} \overleftrightarrow {\mathbf{\Gamma}}_n \left(s_{n,\infty };t_n,u_n,\eta \right)\nonumber\\
&&\times \widetilde{w}_{n,n}^{- (\Omega /\mu +2n-1-\nu ) } \left(\widetilde{w}_{n,n} \partial _{\widetilde{w}_{n,n}} \right) \widetilde{w}_{n,n}^{ \Omega /\mu +2n-1-\nu} \nonumber\\
&&\times  \prod_{j=1}^{n-1} \Bigg\{ \int_{0}^{1} dt_{n-j}\;t_{n-j}^{2(n-j)-\nu } \int_{0}^{1} du_{n-j} \;u_{n-j}^{2(n-j)-1} \overleftrightarrow {\mathbf{\Gamma}}_{n-j} \left(s_{n-j};t_{n-j},u_{n-j},\widetilde{w}_{n-j+1,n} \right)\nonumber\\
&&\times \widetilde{w}_{n-j,n}^{- (\Omega /\mu +2(n-j)-1-\nu ) } \left(\widetilde{w}_{n-j,n} \partial _{\widetilde{w}_{n-j,n}} \right) \widetilde{w}_{n-j,n}^{ \Omega /\mu +2(n-j)-1-\nu}\Bigg\} \mathbf{B} \left( s_{0} ;\widetilde{w}_{1,n} \right) \Bigg\} \rho ^n \Bigg\}
  \label{eq:50052}
\end{eqnarray}
where
\begin{equation}
\begin{cases} 
{ \displaystyle \overleftrightarrow {\mathbf{\Gamma}}_1 \left(s_{1,\infty };t_1,u_1,\eta \right)=  \exp \left(-\frac{s_{1,\infty }}{(1-s_{1,\infty })}\eta (1-t_1)(1-u_1)\right) }\cr
{ \displaystyle  \overleftrightarrow {\mathbf{\Gamma}}_n \left(s_{n,\infty };t_n,u_n,\eta \right) = \exp\left( -\frac{s_{n,\infty }}{(1-s_{n,\infty })}\eta (1-t_n)(1-u_n)\right)}\cr
{ \displaystyle \overleftrightarrow {\mathbf{\Gamma}}_{n-j} \left(s_{n-j};t_{n-j},u_{n-j},\widetilde{w}_{n-j+1,n}\right) = \frac{\exp\left(-\frac{s_{n-j}}{(1-s_{n-j})}\widetilde{w}_{n-j+1,n}(1-t_{n-j})(1-u_{n-j})\right)}{(1-s_{n-j})}}
\end{cases}\nonumber 
\end{equation}
and
\begin{equation}
\begin{cases} 
{ \displaystyle \mathbf{B} \left( s_{0,\infty } ;\eta \right)= (1-s_{0,\infty })^{ \nu -2}\exp\left(-\frac{\eta s_{0,\infty }}{(1-s_{0,\infty })}\right)}\cr
{ \displaystyle  \mathbf{B} \left( s_{0} ;\widetilde{w}_{1,1} \right) = (1-s_0)^{ \nu -2} \exp\left(-\frac{\widetilde{w}_{1,1} s_0}{(1-s_0)}\right)   } \cr
{ \displaystyle \mathbf{B} \left( s_{0} ;\widetilde{w}_{1,n}\right) = (1-s_0)^{ \nu -2}   \exp\left(-\frac{\widetilde{w}_{1,n} s_0}{(1-s_0)} \right) }
\end{cases}\nonumber 
\end{equation}
\end{remark}
\begin{proof}
Replace $t$, $\gamma$ and $z$  by $s_{0,\infty }$ , $2-\nu$ and $\eta$ in (\ref{eq:50048}). 
\begin{equation}
\sum_{\omega _0=0}^{\infty } \frac{s_{0,\infty }^{\omega _0}}{\omega _0!} \frac{\Gamma (\omega _0+2-\nu )}{\Gamma(2-\nu ) }\; _1F_1 \left(-\omega _0; 2-\nu ; \eta \right)  = (1-s_{0,\infty })^{\nu -2} \exp\left(-\frac{\eta s_{0,\infty }}{(1-s_{0,\infty })}\right)  \label{eq:50053}
\end{equation} 
Replace $t$, $\gamma$ and $z$  by $s_0$ , $2-\nu$ and $\widetilde{w}_{1,1}$ in (\ref{eq:50048}). 
\begin{equation}
\sum_{\omega _0=0}^{\infty } \frac{s_0^{\omega _0}}{\omega _0!} \frac{\Gamma (\omega _0+2-\nu )}{\Gamma( 2-\nu ) }\; _1F_1 \left(-\omega _0; 2-\nu ; \widetilde{w}_{1,1} \right)  = (1-s_0)^{\nu -2} \exp\left(-\frac{\widetilde{w}_{1,1} s_0}{(1-s_0)}\right)  \label{eq:50054}
\end{equation} 
Replace $t$, $\gamma$ and $z$  by $s_0$ , $2-\nu$ and $\widetilde{w}_{1,n}$ in (\ref{eq:50048}). 
\begin{equation}
\sum_{\omega _0=0}^{\infty } \frac{s_0^{\omega _0}}{\omega _0!} \frac{\Gamma (\omega _0+2-\nu )}{\Gamma( 2-\nu ) }\; _1F_1 \left(-\omega _0; 2-\nu ; \widetilde{w}_{1,n} \right)  = (1-s_0)^{\nu -2 } \exp\left(-\frac{\widetilde{w}_{1,n} s_0}{(1-s_0)}\right)  \label{eq:50055}
\end{equation} 
Put $c_0$= 1, $\lambda =1-\nu$  and $\displaystyle{\gamma' =2-\nu}$ in (\ref{eq:50035}). Substitute (\ref{eq:50053}), (\ref{eq:50054}) and (\ref{eq:50055}) into the new (\ref{eq:50035}).
\qed
\end{proof}

\section{GCH equation about irregular singular point at infinity}
Let $z=\frac{1}{x}$ in (\ref{eq:5001}) in order to obtain the analytic solution of the GCH equation about $x=\infty $.
\begin{equation}
z^4 \frac{d^2{y}}{d{z}^2} + \left( (2-\nu ) z^3 - \varepsilon z^2 - \mu z\right) \frac{d{y}}{d{z}} + \left( \Omega + \varepsilon \omega z\right) y = 0
\label{eq:50056}
\end{equation}
Assume that its solution is
\begin{equation}
y(z)= \sum_{n=0}^{\infty } c_n z^{n+\lambda }  \label{eq:50057}
\end{equation}
where $\lambda$ is indicial root.  Plug (\ref{eq:50057})  into (\ref{eq:50056}). We get a three-term recurrence relation for the coefficients $c_n$:
\begin{equation}
c_{n+1}=A_n \;c_n +B_n \;c_{n-1} \hspace{1cm};n\geq 1 \label{eq:50058}
\end{equation}
where,
\begin{subequations}
\begin{equation}
A_n = -\frac{\varepsilon}{\mu } \frac{(n-\omega +\lambda )}{(n+1- \Omega /\mu +\lambda )} 
\label{eq:50059a}
\end{equation}
\begin{equation}
B_n = \frac{1}{\mu} \frac{(n-1 +\lambda )(n-\nu +\lambda )}{(n+1- \Omega /\mu +\lambda )} \label{eq:50059b}
\end{equation}
\begin{equation}
c_1= A_0 \;c_0 \label{eq:50059c}
\end{equation}
\end{subequations} 
We have an indicial roots which is $\lambda = \Omega /\mu $.

Now let's test for convergence of the analytic function $y(z)$. As $n\gg 1$ (for sufficiently large $n$), (\ref{eq:50059a}) and (\ref{eq:50059b}) are
\begin{subequations}
\begin{equation}
 \lim_{n\gg 1} A_n = -\frac{\varepsilon}{\mu } 
 \label{eq:50060a}
\end{equation}
\begin{equation}
 \lim_{n\gg 1} B_n = \frac{n}{\mu}
 \label{eq:50060b}
\end{equation}
\end{subequations}
There are no analytic solutions for type 2 polynomial and infinite series. Because the $y(z)$ will be divergent as $n\gg 1$ in (\ref{eq:50060b}). So there are only two types of the analytic solutions of the GCH equation about $x=\infty $ which are type 1  and type 3 polynomials.
In future papers I will derive type 3 polynomial about $x=\infty $: I treat $\mu $, $\varepsilon  $ and $\Omega $ as free variables and $\nu $ and $\omega $ as  fixed values. In this chapter I construct the power series expansion, its integral forms and the generating function for the GCH polynomial of type 1 about $x=\infty $: I treat $\mu $, $\varepsilon  $, $\omega $ and $\Omega $ as free variables and $\nu $ as a fixed value.  

For type 1 polynomial, (\ref{eq:50060a}) is only available for the asymptotic behavior of the minimum $y(z)$.\footnote{(\ref{eq:50060b}) is negligible for the minimum $y(z)$ because $B_n$ term will be terminated at the specific eigenvalues.} Substitute (\ref{eq:50060a}) into (\ref{eq:50058}) with $B_n =0$. For $n=0,1,2,\cdots$, it give
\begin{equation}
\begin{tabular}{ l }
  \vspace{2 mm}
  $c_0$ \\
  \vspace{2 mm}
  $c_1 = -\frac{\varepsilon}{\mu } c_0 $ \\
  \vspace{2 mm}
  $c_2 = \left(-\frac{\varepsilon}{\mu }\right)^{2} c_0 $ \\
  \vspace{2 mm}
  $c_3 = \left(-\frac{\varepsilon}{\mu }\right)^{3} c_0 $ \\
  \vspace{2 mm}  
  $c_4 = \left(-\frac{\varepsilon}{\mu }\right)^{4} c_0 $ \\
  \vspace{2 mm}                        
 \; \vdots \hspace{5mm} \vdots \\
\end{tabular}\label{eq:50061}
\end{equation}
Put (\ref{eq:50061}) into (\ref{eq:50057}) taking $c_0=1$ for simplicity and $\lambda = \Omega /\mu$.
\begin{equation}
\lim_{n\gg 1}y(z) > z^{\Omega /\mu} \sum_{n=0}^{\infty }\left(-\frac{\varepsilon}{\mu }\right)^{n} z^n = \frac{z^{\Omega /\mu}}{1+ \frac{\varepsilon}{\mu }z}\hspace{1cm}\mbox{where}\;\left| \frac{\varepsilon}{\mu }z \right| <1 \label{eq:50062}
\end{equation}
Since $|\varepsilon |,|\mu |\gg 1$ for a polynomial of type 1, (\ref{eq:50060b}) is negligible relative to (\ref{eq:50060a}). The approximative $y(z)$ as $|\varepsilon |,|\mu |\gg 1$ is
\begin{equation}
\lim_{\substack{n\gg 1\\ |\varepsilon |,|\mu | \gg 1}} y(z)  = \frac{z^{\Omega /\mu}}{1+ \frac{\varepsilon}{\mu }z} \label{eq:50063}
\end{equation}

\subsection{Power series for polynomial of type 1}
In Ref.\cite{Chou2012b}, the general expression of power series of $y(x)$ for polynomial of type 1 is given by
\begin{eqnarray}
 y(x)&=&  \sum_{n=0}^{\infty } y_{n}(x) = y_0(x)+ y_1(x)+ y_2(x)+y_3(x)+\cdots \nonumber\\
&=& c_0 \Bigg\{ \sum_{i_0=0}^{\beta _0} \left( \prod _{i_1=0}^{i_0-1}B_{2i_1+1} \right) x^{2i_0+\lambda } 
+ \sum_{i_0=0}^{\beta _0}\left\{ A_{2i_0} \prod _{i_1=0}^{i_0-1}B_{2i_1+1}  \sum_{i_2=i_0}^{\beta _1} \left( \prod _{i_3=i_0}^{i_2-1}B_{2i_3+2} \right)\right\} x^{2i_2+1+\lambda }\nonumber\\
 && + \sum_{N=2}^{\infty } \Bigg\{ \sum_{i_0=0}^{\beta _0} \Bigg\{A_{2i_0}\prod _{i_1=0}^{i_0-1} B_{2i_1+1} \prod _{k=1}^{N-1} \Bigg( \sum_{i_{2k}= i_{2(k-1)}}^{\beta _k} A_{2i_{2k}+k}\prod _{i_{2k+1}=i_{2(k-1)}}^{i_{2k}-1}B_{2i_{2k+1}+(k+1)}\Bigg)\nonumber\\
 &&\times  \sum_{i_{2N} = i_{2(N-1)}}^{\beta _N} \Bigg( \prod _{i_{2N+1}=i_{2(N-1)}}^{i_{2N}-1} B_{2i_{2N+1}+(N+1)} \Bigg) \Bigg\} \Bigg\} x^{2i_{2N}+N+\lambda }\Bigg\}
  \label{eq:50064}
\end{eqnarray}
In the above, $\beta _i\leq \beta _j$ only if $i\leq j$ where $i,j,\beta _i, \beta _j \in \mathbb{N}_{0}$.

For a polynomial, we need a condition, which is:
\begin{equation}
 B_{2\beta _i + (i+1)}=0 \hspace{1cm} \mathrm{where}\; i= 0,1,2,\cdots, \beta _i=0,1,2,\cdots
 \label{eq:50065}
\end{equation}
In the above, $ \beta _i$ is an eigenvalue that makes $B_n$ term terminated at certain value of index $n$. (\ref{eq:50065}) makes each $y_i(x)$ where $i=0,1,2,\cdots$ as the polynomial in (\ref{eq:50064}).

In (\ref{eq:50059b}) replace $\nu $ by $2 \nu _i+i +1+ \lambda $. In (\ref{eq:50064}) and (\ref{eq:50065}) replace an independent variable $x$ and an index $\beta _i$ by $z$ and $\nu _i$. Take the new (\ref{eq:50059a})--(\ref{eq:50059c}), (\ref{eq:50065}) and put them in (\ref{eq:50064}).
After the replacement process, the general expression of power series of the GCH equation for polynomial of type 1 about $x=\infty $ is given by
\begin{eqnarray}
 y(z)&=& \sum_{n=0}^{\infty } y_n(z)= y_0(z)+ y_1(z)+ y_2(z)+y_3(z)+\cdots \nonumber\\
&=& c_0 z^{\lambda } \left\{\sum_{i_0=0}^{\nu _0} \frac{(-\nu _0)_{i_0}( \frac{\lambda }{2})_{i_0}}{(1- \frac{\Omega }{2\mu } +\frac{\lambda }{2})_{i_0}} \eta ^{i_0}\right.\nonumber\\
&&+ \left\{ \sum_{i_0=0}^{\nu _0}\frac{(i_0 -\frac{\omega }{2}+ \frac{\lambda }{2})}{(i_0+ \frac{1}{2}-\frac{\Omega }{2\mu}+ \frac{\lambda }{2})}  \frac{(-\nu _0)_{i_0}( \frac{\lambda }{2})_{i_0}}{(1- \frac{\Omega }{2\mu } +\frac{\lambda }{2})_{i_0}} \sum_{i_1=i_0}^{\nu _1} \frac{(-\nu _1)_{i_1}(\frac{1}{2}+\frac{\lambda }{2})_{i_1}(\frac{3}{2}-\frac{\Omega}{2\mu}+ \frac{\lambda }{2})_{i_0}}{(-\nu _1)_{i_0}(\frac{1}{2}+\frac{\lambda }{2})_{i_0}(\frac{3}{2}-\frac{\Omega}{2\mu}+ \frac{\lambda }{2})_{i_1}} \eta ^{i_1} \right\} \xi  \nonumber\\
&&+ \sum_{n=2}^{\infty } \left\{ \sum_{i_0=0}^{\nu _0}\frac{(i_0 -\frac{\omega }{2}+ \frac{\lambda }{2})}{(i_0+ \frac{1}{2}-\frac{\Omega }{2\mu}+ \frac{\lambda }{2})}  \frac{(-\nu _0)_{i_0}( \frac{\lambda }{2})_{i_0}}{(1- \frac{\Omega }{2\mu } +\frac{\lambda }{2})_{i_0}}\right.\nonumber\\
&&\times \prod _{k=1}^{n-1} \left\{ \sum_{i_k=i_{k-1}}^{\nu _k} \frac{(i_k+\frac{k}{2} -\frac{\omega }{2}+ \frac{\lambda }{2}) }{(i_k+ \frac{k}{2}+\frac{1}{2}-\frac{\Omega }{2\mu}+\frac{\lambda }{2})} \frac{(-\nu _k)_{i_k}( \frac{k}{2}+\frac{\lambda }{2})_{i_k}( \frac{k}{2}+1-\frac{\Omega }{2\mu}+ \frac{\lambda }{2})_{i_{k-1}}}{(-\nu _k)_{i_{k-1}}( \frac{k}{2}+\frac{\lambda }{2})_{i_{k-1}}( \frac{k}{2}+1-\frac{\Omega }{2\mu}+ \frac{\lambda }{2})_{i_k}}\right\} \nonumber\\
&&\times \left.\left.\sum_{i_n= i_{n-1}}^{\nu _n} \frac{(-\nu _n)_{i_n}( \frac{n}{2}+\frac{\lambda }{2})_{i_n}( \frac{n}{2}+1-\frac{\Omega }{2\mu}+ \frac{\lambda }{2})_{i_{n-1}}}{(-\nu _n)_{i_{n-1}}( \frac{n}{2}+\frac{\lambda }{2})_{i_{n-1}}( \frac{n}{2}+1-\frac{\Omega }{2\mu}+ \frac{\lambda }{2})_{i_n}} \eta ^{i_n} \right\} \xi ^n \right\}\label{eq:50066}
\end{eqnarray}
where
\begin{equation}
\begin{cases} \eta = \frac{2}{\mu } z^2 \cr
\xi = -\frac{\varepsilon }{\mu } z \cr
 \nu = 2 \nu _j+j +1+ \lambda  \cr
 z= \frac{1}{x} \cr
\nu _i\leq \nu _j \;\;\mbox{only}\;\mbox{if}\;i\leq j\;\;\mbox{where}\;i,j,\nu _i,\nu _j\in \mathbb{N}_{0}\cdots
\end{cases}\nonumber
\end{equation}
Put $c_0$= 1 as $\lambda = \Omega /\mu$ in (\ref{eq:50066}). 
\begin{remark}
The power series expansion of the GCH equation of the first kind for polynomial of type 1 about $x=\infty $ as $\nu = 2 \nu _j+j +1+ \Omega /\mu $ where $j,\nu _j \in \mathbb{N}_{0}$ is
\begin{eqnarray}
 y(z)&=& Q^{(i)}W_{\nu _j}\left( \mu ,\varepsilon, \Omega, \omega ,\nu = 2 \nu _j+j +1+ \frac{\Omega}{\mu}; z=\frac{1}{x}; \xi = -\frac{\varepsilon }{\mu}z; \eta = \frac{2}{\mu } z^2 \right)\nonumber\\
&=& z^{\frac{\Omega}{\mu}} \left\{\sum_{i_0=0}^{\nu _0} \frac{(-\nu _0)_{i_0}( \frac{\Omega }{2\mu})_{i_0}}{(1)_{i_0}} \eta ^{i_0}\right.\nonumber\\
&&+ \left\{ \sum_{i_0=0}^{\nu _0}\frac{(i_0 -\frac{\omega }{2}+ \frac{\Omega }{2\mu })}{(i_0+ \frac{1}{2})}  \frac{(-\nu _0)_{i_0}( \frac{\Omega }{2\mu})_{i_0}}{(1)_{i_0}} \sum_{i_1=i_0}^{\nu _1} \frac{(-\nu _1)_{i_1}(\frac{1}{2}+\frac{\Omega }{2\mu})_{i_1}(\frac{3}{2})_{i_0}}{(-\nu _1)_{i_0}(\frac{1}{2}+\frac{\Omega }{2\mu})_{i_0}(\frac{3}{2})_{i_1}} \eta ^{i_1} \right\} \xi  \nonumber\\
&&+ \sum_{n=2}^{\infty } \left\{ \sum_{i_0=0}^{\nu _0}\frac{(i_0 -\frac{\omega }{2}+ \frac{\Omega }{2\mu })}{(i_0+ \frac{1}{2})}  \frac{(-\nu _0)_{i_0}( \frac{\Omega }{2\mu})_{i_0}}{(1)_{i_0}}\right.\nonumber\\
&&\times \prod _{k=1}^{n-1} \left\{ \sum_{i_k=i_{k-1}}^{\nu _k} \frac{(i_k+\frac{k}{2} -\frac{\omega }{2}+ \frac{\Omega }{2\mu}) }{(i_k+ \frac{k}{2}+\frac{1}{2})} \frac{(-\nu _k)_{i_k}( \frac{k}{2}+\frac{\Omega }{2\mu})_{i_k}( \frac{k}{2}+1)_{i_{k-1}}}{(-\nu _k)_{i_{k-1}}( \frac{k}{2}+\frac{\Omega }{2\mu})_{i_{k-1}}( \frac{k}{2}+1)_{i_k}}\right\} \nonumber\\
&&\times \left.\left.\sum_{i_n= i_{n-1}}^{\nu _n} \frac{(-\nu _n)_{i_n}( \frac{n}{2}+\frac{\Omega }{2\mu})_{i_n}( \frac{n}{2}+1)_{i_{n-1}}}{(-\nu _n)_{i_{n-1}}( \frac{n}{2}+\frac{\Omega }{2\mu})_{i_{n-1}}( \frac{n}{2}+1)_{i_n}} \eta ^{i_n} \right\} \xi ^n \right\}\label{eq:50067}
\end{eqnarray}
\end{remark}
For the minimum value of the GCH equation of the first kind for polynomial of type 1 about $x=\infty $, put $\nu_0=\nu _1=\nu_2=\cdots=0$ in (\ref{eq:50067}).
\begin{eqnarray}
y(z)&=& Q^{(i)}W_{0}\left( \mu ,\varepsilon, \Omega, \omega ,\nu = j +1+ \frac{\Omega}{\mu}; z=\frac{1}{x}; \xi = -\frac{\varepsilon }{\mu}z; \eta = \frac{2}{\mu } z^2 \right)\nonumber\\
&=& z^{\frac{\Omega}{\mu}}\left( 1+\frac{\varepsilon }{\mu }z\right)^{-\left( \frac{\Omega}{\mu}-\omega \right)} \hspace{1cm}\mbox{where}\;\left| \frac{\varepsilon}{\mu }z \right| <1\nonumber 
\end{eqnarray}
\subsection{Integral formalism for polynomial of type 1}
There is a generalized hypergeometric function which is given by
\begin{eqnarray}
L_l &=& \sum_{i_l= i_{l-1}}^{\nu _l} \frac{(-\nu _l)_{i_l}(\frac{l}{2}+\frac{\lambda }{2})_{i_l}( \frac{l}{2}+1-\frac{\Omega }{2\mu}+ \frac{\lambda }{2})_{i_{l-1}}}{(-\nu _l)_{i_{l-1}}(\frac{l}{2}+\frac{\lambda }{2})_{i_{l-1}}( \frac{l}{2}+1-\frac{\Omega }{2\mu}+ \frac{\lambda }{2})_{i_l}} \eta ^{i_l}\nonumber\\
&=& \eta^{i_{l-1}} 
\sum_{j=0}^{\infty } \frac{B\left(i_{l-1}+\frac{l}{2}-\frac{\Omega }{2\mu}+\frac{\lambda }{2},j+1\right) (i_{l-1}-\nu _l)_j (i_{l-1}+\frac{l}{2}+\frac{\lambda }{2})_j}{(i_{l-1}+\frac{l}{2} -\frac{\Omega }{2\mu}+\frac{\lambda }{2} )^{-1}(1)_j} \eta^j 
\label{eq:50068}
\end{eqnarray}
By using integral form of beta function,
\begin{equation}
B\left(i_{l-1}+\frac{l}{2}-\frac{\Omega }{2\mu}+\frac{\lambda }{2},j+1\right) = \int_{0}^{1} dt_l\;t_l^{i_{l-1}+\frac{l}{2}-1-\frac{\Omega }{2\mu}+ \frac{\lambda }{2} } (1-t_l)^j
\label{eq:50069}
\end{equation}
Substitute (\ref{eq:50069}) into (\ref{eq:50068}), and divide $(i_{l-1}+\frac{l}{2} -\frac{\Omega }{2\mu}+\frac{\lambda }{2} )$ into $L_l$.
\begin{eqnarray}
G_l &=& \frac{1}{(i_{l-1}+\frac{l}{2} -\frac{\Omega }{2\mu}+\frac{\lambda }{2} )}
 \sum_{i_l= i_{l-1}}^{\nu _l} \frac{(-\nu _l)_{i_l}(\frac{l}{2}+\frac{\lambda }{2})_{i_l}( \frac{l}{2}+1-\frac{\Omega }{2\mu}+ \frac{\lambda }{2})_{i_{l-1}}}{(-\nu _l)_{i_{l-1}}(\frac{l}{2}+\frac{\lambda }{2})_{i_{l-1}}( \frac{l}{2}+1-\frac{\Omega }{2\mu}+ \frac{\lambda }{2})_{i_l}} \eta ^{i_l} \nonumber\\
&=&  \int_{0}^{1} dt_l\;t_l^{\frac{l}{2}-1-\frac{\Omega }{2\mu}+ \frac{\lambda }{2}} \left(\eta t_l \right)^{i_{l-1}} \sum_{j=0}^{\infty } \frac{(i_{l-1}-\nu _l)_j (i_{l-1}+\frac{l}{2}+\frac{\lambda }{2})_j}{(1)_j} \left(\eta (1-t_l) \right)^j \hspace{1cm} \label{eq:50070}
\end{eqnarray}
Tricomi's function is defined by
\begin{equation}
U(a,b,z)=  \frac{\Gamma (1-b)}{\Gamma (a-b+1)} M(a,b,z) +\frac{\Gamma (b-1)}{\Gamma (a)}z^{1-b} M(a-b+1,2-b,z)
\label{eq:50072}
\end{equation}
The contour integral form of (\ref{eq:50072}) is given by\cite{NIST1}
\begin{equation}
U(a,b,z)=  e^{-a\pi i}\frac{\Gamma (1-a)}{2\pi i} \int_{\infty }^{(0+)} dp_l\; e^{-z p_l} p_l^{a-1} (1+p_l)^{b-a-1}\;\mbox{where}\;a\ne 1,2,3, \cdots, \;\;\; |\mbox{ph} z|< \frac{1}{2}\pi
\label{eq:50073}
\end{equation}
Also (\ref{eq:50072}) is written by \cite{NIST1}
\begin{equation}
U(a,b,z)=  z^{-a}  \sum_{j= 0}^{\infty } \frac{(a)_j (a-b+1)_j}{(1)_j} (-z^{-1})^j = z^{-a}\; _2F_0 (a,a-b+1;-;-z^{-1})
\label{eq:50074}
\end{equation}
Replace $a$, $b$ and $z$ by $i_{l-1} -\nu_l$, $-\nu_l +1- \frac{l}{2}- \frac{\lambda }{2}$ and $\frac{-1}{\eta (1-t_l)}$ into (\ref{eq:50074}).
\begin{equation}
\sum_{j= 0}^{\infty } \frac{(i_{l-1} -\nu_l )_j (i_{l-1}+ \frac{l}{2}+ \frac{\lambda }{2} )_j}{(1)_j} (\eta (1-t_l))^j =
\left( \frac{-1}{\eta (1-t_l)}\right)^{i_{l-1} -\nu_l} U\left( i_{l-1} -\nu_l, -\nu_l +1- \frac{l}{2}- \frac{\lambda }{2},\frac{-1}{\eta (1-t_l)} \right)
\label{eq:50075}
\end{equation}
Replace $a$, $b$ and $z$ by $i_{l-1} -\nu_l$, $-\nu_l +1- \frac{l}{2}- \frac{\lambda }{2}$ and $\frac{-1}{\eta (1-t_l)}$ into (\ref{eq:50073}).
Take the new (\ref{eq:50073}) into (\ref{eq:50075})
\begin{eqnarray}
&& \sum_{j= 0}^{\infty } \frac{(i_{l-1} -\nu_l )_j (i_{l-1}+ \frac{l}{2}+ \frac{\lambda }{2} )_j}{(1)_j} (\eta (1-t_l))^j \nonumber\\ 
&&= \frac{\Gamma (\nu_l -i_{l-1}+1)}{2\pi i} \int_{\infty }^{(0+)} dp_l\; \exp\left(\frac{ p_l}{\eta (1-t_l)}\right) p_l^{-1} (1+p_l)^{-\frac{1}{2}(l+\lambda )} \nonumber\\
&&\times \left( \frac{\eta (1-t_l)}{p_l}\right)^{\nu_l} \left( \frac{p_l}{\eta (1-t_l)(1+p_l)}\right)^{i_{l-1}}
 \label{eq:50076}
\end{eqnarray}
The definition of the Gamma function $\Gamma (z)$ is written by
\begin{equation}
\Gamma (z) = \int_{0}^{\infty} du_l\; e^{-u_l} u_l^{z-1} \;\mbox{where}\; Re(z)>0
\label{eq:50077}
\end{equation}
Put $z=\nu_l -i_{l-1}+1$ in (\ref{eq:50077}).
\begin{equation}
\Gamma (\nu_l -i_{l-1}+1) = \int_{0}^{\infty} du_l\; e^{-u_l} u_l^{\nu_l -i_{l-1}} \label{eq:50078}
\end{equation}
Put (\ref{eq:50078}) in (\ref{eq:50076}). Take the new (\ref{eq:50076}) into (\ref{eq:50070})
\begin{eqnarray}
G_l &=& \frac{1}{(i_{l-1}+\frac{l}{2} -\frac{\Omega }{2\mu}+\frac{\lambda }{2} )}
 \sum_{i_l= i_{l-1}}^{\nu _l} \frac{(-\nu _l)_{i_l}(\frac{l}{2}+\frac{\lambda }{2})_{i_l}( \frac{l}{2}+1-\frac{\Omega }{2\mu}+ \frac{\lambda }{2})_{i_{l-1}}}{(-\nu _l)_{i_{l-1}}(\frac{l}{2}+\frac{\lambda }{2})_{i_{l-1}}( \frac{l}{2}+1-\frac{\Omega }{2\mu}+ \frac{\lambda }{2})_{i_l}} \eta ^{i_l} \nonumber\\
&=&  \int_{0}^{1} dt_l\;t_l^{\frac{1}{2}\left( l-2-\frac{\Omega }{\mu}+ \lambda \right)} \int_{0}^{\infty} du_l\; e^{-u_l}
\frac{1}{2\pi i} \int_{\infty }^{(0+)} dp_l\; \exp\left(\frac{ p_l}{\eta (1-t_l)}\right) p_l^{-1} (1+p_l)^{-\frac{1}{2}(l+\lambda )} \nonumber\\
&&\times \left( \frac{\eta u_l (1-t_l)}{p_l}\right)^{\nu_l} \left( \frac{t_l p_l}{u_l (1-t_l)(1+p_l)}\right)^{i_{l-1}}
 \label{eq:50079}
\end{eqnarray}
Substitute (\ref{eq:50079}) into (\ref{eq:50066}) where $l=1,2,3,\cdots$; apply $G_1$ into the second summation of sub-power series $y_1(z)$, apply $G_2$ into the third summation and $G_1$ into the second summation of sub-power series $y_2(z)$, apply $G_3$ into the forth summation, $G_2$ into the third summation and $G_1$ into the second summation of sub-power series $y_3(z)$, etc.\footnote{$y_1(z)$ means the sub-power series in (\ref{eq:50066}) contains one term of $A_n's$, $y_2(z)$ means the sub-power series in (\ref{eq:50066}) contains two terms of $A_n's$, $y_3(z)$ means the sub-power series in (\ref{eq:50066}) contains three terms of $A_n's$, etc.}
\begin{theorem}
The general representation in the form of integral of the GCH polynomial of type 1 about $x=\infty $ is given by 
\begin{eqnarray}
 y(z)&=& \sum_{n=0}^{\infty } y_n(z)= y_0(z)+ y_1(z)+ y_2(z)+y_3(z)+\cdots \nonumber\\
&=& c_0 z^{\lambda } \left\{ \sum_{i_0=0}^{\nu _0 }\frac{(-\nu _0)_{i_0}\left(\frac{\lambda }{2}\right)_{i_0}}{\left(1-\frac{\Omega }{2\mu }+\frac{\lambda }{2}\right)_{i_0}}  \eta^{i_0} 
+ \sum_{n=1}^{\infty } \left\{\prod _{k=0}^{n-1} \Bigg\{ \int_{0}^{1} dt_{n-k}\;t_{n-k}^{\frac{1}{2}\left( n-k-2 -\frac{\Omega }{\mu }+ \lambda \right)} \int_{0}^{\infty } du_{n-k}\;e^{-u_{n-k}}  \right.\right.\nonumber\\
&&\times \frac{1}{2\pi i} \int_{\infty }^{(0+)} dp_{n-k}\; p_{n-k}^{-1} (1+p_{n-k})^{-\frac{1}{2}(n-k+\lambda )} \exp \left( \frac{p_{n-k}}{w_{n-k+1,n}(1-t_{n-k})}\right)\nonumber\\
&&\times  \left( \frac{w_{n-k+1,n} u_{n-k} (1-t_{n-k})}{p_{n-k}}\right)^{\nu_{n-k}}   \nonumber\\
&&\times  w_{n-k,n}^{-\frac{1}{2} (n-k-1-\omega +\lambda )} \left(  w_{n-k,n} \partial _{w_{n-k,n}} \right) w_{n-k,n}^{\frac{1}{2} (n-k-1-\omega +\lambda )}\Bigg\}  \left.\left. \sum_{i_0=0}^{\nu _0}\frac{(-\nu _0)_{i_0}\left(\frac{\lambda }{2}\right)_{i_0}}{\left(1-\frac{\Omega }{2\mu }+\frac{\lambda }{2}\right)_{i_0}}  w_{1,n}^{i_0}\right\} \xi ^n \right\}\hspace{1.5cm}
\label{eq:50080}
\end{eqnarray}
where
\begin{equation}w_{i,j}=
\begin{cases} \displaystyle {\frac{t_i p_i}{u_i (1-t_i)(1+p_i)} }\;\;\mbox{where}\; i\leq j\cr
\eta \;\;\mbox{only}\;\mbox{if}\; i>j
\end{cases}
\nonumber
\end{equation}
In the above, the first sub-integral form contains one term of $A_n's$, the second one contains two terms of $A_n$'s, the third one contains three terms of $A_n$'s, etc.
\end{theorem}

\begin{proof} 
In (\ref{eq:50066}) the power series expansions of sub-summation $y_0(z) $, $y_1(z)$, $y_2(z)$ and $y_3(z)$ of the GCH polynomial of type 1 about $x=\infty $ are
\begin{equation}
 y(z)= \sum_{n=0}^{\infty } y_n(z) = y_0(z)+ y_1(z)+ y_2(z)+y_3(z)+\cdots  \label{eq:50081}
\end{equation}
where
\begin{subequations}
\begin{equation}
 y_0(z)= \sum_{i_0=0}^{\nu _0 }\frac{(-\nu _0)_{i_0}\left(\frac{\lambda }{2}\right)_{i_0}}{\left(1-\frac{\Omega }{2\mu }+\frac{\lambda }{2}\right)_{i_0}}  \eta^{i_0} \label{eq:50082a}
\end{equation}
\begin{eqnarray}
 y_1(z)&=&  c_0 z^{\lambda } \left\{ \sum_{i_0=0}^{\nu _0}\frac{\left(i_0- \frac{\omega }{2} +\frac{\lambda }{2} \right)}{\left(i_0 + \frac{1}{2} -\frac{\Omega }{2\mu } +\frac{\lambda }{2}\right)} \frac{(-\nu _0)_{i_0}\left(\frac{\lambda }{2}\right)_{i_0}}{\left(1-\frac{\Omega }{2\mu }+\frac{\lambda }{2}\right)_{i_0}} \right. \nonumber\\
&&\times \left. \sum_{i_1=i_0}^{\nu_1} \frac{(-\nu_1)_{i_1}\left(\frac{1}{2}+\frac{\lambda }{2} \right)_{i_1} \left( \frac{3}{2}-\frac{\Omega }{2\mu }+\frac{\lambda }{2} \right)_{i_0}}{(-\nu_1)_{i_0}\left(\frac{1}{2}+\frac{\lambda }{2} \right)_{i_0} \left( \frac{3}{2}-\frac{\Omega }{2\mu }+\frac{\lambda }{2} \right)_{i_1}} \eta^{i_1} \right\}\xi  \label{eq:50082b}
\end{eqnarray}
\begin{eqnarray}
 y_2(z) &=& c_0 z^{\lambda }\left\{\sum_{i_0=0}^{\nu _0}\frac{\left(i_0- \frac{\omega }{2} +\frac{\lambda }{2} \right)}{\left(i_0 + \frac{1}{2} -\frac{\Omega }{2\mu } +\frac{\lambda }{2}\right)} \frac{(-\nu _0)_{i_0}\left(\frac{\lambda }{2}\right)_{i_0}}{\left(1-\frac{\Omega }{2\mu }+\frac{\lambda }{2}\right)_{i_0}} \right.\nonumber\\
&&\times  \sum_{i_1=i_0}^{\nu_1} \frac{\left( i_1 +\frac{1}{2}-\frac{\omega }{2}+\frac{\lambda }{2}\right)}{\left( i_1 +1-\frac{\Omega }{2\mu} +\frac{\lambda }{2} \right)} \frac{(-\nu_1)_{i_1}\left(\frac{1}{2}+\frac{\lambda }{2} \right)_{i_1} \left( \frac{3}{2}-\frac{\Omega }{2\mu }+\frac{\lambda }{2} \right)_{i_0}}{(-\nu_1)_{i_0}\left(\frac{1}{2}+\frac{\lambda }{2} \right)_{i_0} \left( \frac{3}{2}-\frac{\Omega }{2\mu }+\frac{\lambda }{2} \right)_{i_1}}  \nonumber\\
&&\times \left. \sum_{i_2=i_1}^{\nu _2} \frac{(-\nu_2)_{i_2}\left(1+\frac{\lambda }{2} \right)_{i_2} \left( 2-\frac{\Omega }{2\mu }+\frac{\lambda }{2} \right)_{i_1}}{(-\nu_2)_{i_1}\left(1+\frac{\lambda }{2} \right)_{i_1} \left( 2-\frac{\Omega }{2\mu }+\frac{\lambda }{2} \right)_{i_2}}  \eta^{i_2} \right\} \xi ^2 
\label{eq:50082c}
\end{eqnarray}
\begin{eqnarray}
 y_3(z)&=&  c_0 z^{\lambda } \left\{ \sum_{i_0=0}^{\nu _0}\frac{\left(i_0- \frac{\omega }{2} +\frac{\lambda }{2} \right)}{\left(i_0 + \frac{1}{2} -\frac{\Omega }{2\mu } +\frac{\lambda }{2}\right)} \frac{(-\nu _0)_{i_0}\left(\frac{\lambda }{2}\right)_{i_0}}{\left(1-\frac{\Omega }{2\mu }+\frac{\lambda }{2}\right)_{i_0}}\right. \nonumber\\
&&\times  \sum_{i_1=i_0}^{\nu_1} \frac{\left( i_1 +\frac{1}{2}-\frac{\omega }{2}+\frac{\lambda }{2}\right)}{\left( i_1 +1-\frac{\Omega }{2\mu} +\frac{\lambda }{2} \right)} \frac{(-\nu_1)_{i_1}\left(\frac{1}{2}+\frac{\lambda }{2} \right)_{i_1} \left( \frac{3}{2}-\frac{\Omega }{2\mu }+\frac{\lambda }{2} \right)_{i_0}}{(-\nu_1)_{i_0}\left(\frac{1}{2}+\frac{\lambda }{2} \right)_{i_0} \left( \frac{3}{2}-\frac{\Omega }{2\mu }+\frac{\lambda }{2} \right)_{i_1}}  \nonumber\\
&&\times \sum_{i_2=i_1}^{\nu _2} \frac{\left( i_2 +1-\frac{\omega }{2}+\frac{\lambda }{2}\right)}{\left( i_2 +\frac{3}{2}-\frac{\Omega }{2\mu} +\frac{\lambda }{2} \right)} \frac{(-\nu_2)_{i_2}\left(1+\frac{\lambda }{2} \right)_{i_2} \left( 2-\frac{\Omega }{2\mu }+\frac{\lambda }{2} \right)_{i_1}}{(-\nu_2)_{i_1}\left(1+\frac{\lambda }{2} \right)_{i_1} \left( 2-\frac{\Omega }{2\mu }+\frac{\lambda }{2} \right)_{i_2}}\nonumber\\
&&\times \left. \sum_{i_3=i_2}^{\nu _3} \frac{(-\nu_3)_{i_3}\left(\frac{3}{2}+\frac{\lambda }{2} \right)_{i_3} \left( \frac{5}{2}-\frac{\Omega }{2\mu }+\frac{\lambda }{2} \right)_{i_2}}{(-\nu_3)_{i_2}\left(\frac{3}{2}+\frac{\lambda }{2} \right)_{i_2} \left( \frac{5}{2}-\frac{\Omega }{2\mu }+\frac{\lambda }{2} \right)_{i_3}} \eta^{i_3} \right\} \xi ^3 
\label{eq:50082d}
\end{eqnarray}
\end{subequations}
Put $l=1$ in (\ref{eq:50079}). Take the new (\ref{eq:50079}) into (\ref{eq:50082b}).
\begin{eqnarray}
 y_1(z)&=& c_0 z^{\lambda }  \int_{0}^{1} dt_1\;t_1^{\frac{1}{2}\left( -1-\frac{\Omega }{\mu }+\lambda \right)} \int_{0}^{\infty } du_1\;e^{-u_1} 
\frac{1}{2\pi i}  \int_{\infty }^{(0+)} dp_1\; \exp\left(\frac{ p_1}{\eta (1-t_1)}\right) p_1^{-1} (1+p_1)^{-\frac{1}{2}(1+\lambda )} \nonumber\\
&&\times \left( \frac{\eta u_1 (1-t_1)}{p_1}\right)^{\nu_1} \left\{ \sum_{i_0=0}^{\nu _0} \left( i_0-\frac{\omega }{2}+\frac{\lambda }{2}\right) \frac{(-\nu _0)_{i_0}\left(\frac{\lambda }{2}\right)_{i_0}}{\left(1-\frac{\Omega }{2\mu }+\frac{\lambda }{2}\right)_{i_0}} \left(\frac{t_1 p_1}{u_1 (1-t_1)(1+p_1)}\right)^{i_0} \right\} \xi \nonumber\\
&=&c_0 z^{\lambda }  \int_{0}^{1} dt_1\;t_1^{\frac{1}{2}\left( -1-\frac{\Omega }{\mu }+\lambda \right)} \int_{0}^{\infty } du_1\;e^{-u_1} 
\frac{1}{2\pi i}  \int_{\infty }^{(0+)} dp_1\; \exp\left(\frac{ p_1}{\eta (1-t_1)}\right) p_1^{-1} (1+p_1)^{-\frac{1}{2}(1+\lambda )} \nonumber\\
&&\times \left( \frac{\eta u_1 (1-t_1)}{p_1}\right)^{\nu_1} w_{1,1}^{-\frac{1}{2}\left(-\omega +\lambda \right)} \left( w_{1,1}\partial_{w_{1,1}} \right) w_{1,1}^{\frac{1}{2}\left(-\omega +\lambda \right) } \nonumber\\
&&\times \left\{ \sum_{i_0=0}^{\nu _0} \frac{(-\nu _0)_{i_0}\left(\frac{\lambda }{2}\right)_{i_0}}{\left(1-\frac{\Omega }{2\mu }+\frac{\lambda }{2}\right)_{i_0}} w_{1,1}^{i_0} \right\} \xi \label{eq:50083}\\
&& \mathrm{where}\hspace{.5cm} w_{1,1}= \frac{t_1 p_1}{u_1 (1-t_1)(1+p_1)}  \nonumber
\end{eqnarray}
Put $l=2$ in (\ref{eq:50079}). Take the new (\ref{eq:50079}) into (\ref{eq:50082c}). 
\begin{eqnarray}
 y_2(z) &=& c_0 z^{\lambda } \int_{0}^{1} dt_2\;t_2^{\frac{1}{2}\left(-\frac{\Omega }{\mu }+\lambda \right)} \int_{0}^{\infty } du_2\;e^{-u_2} 
\frac{1}{2\pi i}  \int_{\infty }^{(0+)} dp_2\; \exp\left(\frac{ p_2}{\eta (1-t_2)}\right) p_2^{-1} (1+p_2)^{-\frac{1}{2}(2+\lambda )} \nonumber\\
&&\times \left( \frac{\eta u_2 (1-t_2)}{p_2}\right)^{\nu_2}  w_{2,2}^{-\frac{1}{2}\left(1-\omega +\lambda \right)} \left( w_{2,2}\partial_{w_{2,2}} \right) w_{2,2}^{\frac{1}{2}\left(1-\omega +\lambda \right) } \label{eq:50084}\\
&&\times \left\{  \sum_{i_0=0}^{\nu _0}\frac{\left(i_0- \frac{\omega }{2} +\frac{\lambda }{2} \right)}{\left(i_0 + \frac{1}{2} -\frac{\Omega }{2\mu } +\frac{\lambda }{2}\right)} \frac{(-\nu _0)_{i_0}\left(\frac{\lambda }{2}\right)_{i_0}}{\left(1-\frac{\Omega }{2\mu }+\frac{\lambda }{2}\right)_{i_0}} \right.  \left. \sum_{i_1=i_0}^{\nu_1} \frac{(-\nu_1)_{i_1}\left(\frac{1}{2}+\frac{\lambda }{2} \right)_{i_1} \left( \frac{3}{2}-\frac{\Omega }{2\mu }+\frac{\lambda }{2} \right)_{i_0}}{(-\nu_1)_{i_0}\left(\frac{1}{2}+\frac{\lambda }{2} \right)_{i_0} \left( \frac{3}{2}-\frac{\Omega }{2\mu }+\frac{\lambda }{2} \right)_{i_1}} w_{2,2}^{i_1} \right\} \xi ^2  \nonumber\\
&& \mathrm{where}\hspace{.5cm} w_{2,2}=  \frac{t_2 p_2}{u_2 (1-t_2)(1+p_2)} \nonumber
\end{eqnarray}
Put $l=1$ and $\eta = w_{2,2}$ in (\ref{eq:50079}). Take the new (\ref{eq:50079}) into (\ref{eq:50084}).
\begin{eqnarray}
 y_2(z)&=& c_0 z^{\lambda } \int_{0}^{1} dt_2\;t_2^{\frac{1}{2}\left(-\frac{\Omega }{\mu }+\lambda \right)} \int_{0}^{\infty } du_2\;e^{-u_2} 
\frac{1}{2\pi i}  \int_{\infty }^{(0+)} dp_2\; \exp\left(\frac{ p_2}{\eta (1-t_2)}\right) p_2^{-1} (1+p_2)^{-\frac{1}{2}(2+\lambda )} \nonumber\\
&&\times \left( \frac{\eta u_2 (1-t_2)}{p_2}\right)^{\nu_2}  w_{2,2}^{-\frac{1}{2}\left(1-\omega +\lambda \right)} \left( w_{2,2}\partial_{w_{2,2}} \right) w_{2,2}^{\frac{1}{2}\left(1-\omega +\lambda \right) } \nonumber\\
&&\times  \int_{0}^{1} dt_1\;t_1^{\frac{1}{2}\left( -1-\frac{\Omega }{\mu }+\lambda \right)} \int_{0}^{\infty } du_1\;e^{-u_1} 
\frac{1}{2\pi i}  \int_{\infty }^{(0+)} dp_1\; \exp\left(\frac{ p_1}{w_{2,2} (1-t_1)}\right) p_1^{-1} (1+p_1)^{-\frac{1}{2}(1+\lambda )} \nonumber\\
&&\times \left( \frac{w_{2,2} u_1 (1-t_1)}{p_1}\right)^{\nu_1}  w_{1,2}^{-\frac{1}{2}\left( -\omega +\lambda \right)} \left( w_{1,2}\partial_{w_{1,2}} \right) w_{1,2}^{\frac{1}{2}\left( -\omega +\lambda \right) }  \nonumber\\
&&\times \left\{ \sum_{i_0=0}^{\nu _0} \frac{(-\nu _0)_{i_0}\left(\frac{\lambda }{2}\right)_{i_0}}{\left(1-\frac{\Omega }{2\mu }+\frac{\lambda }{2}\right)_{i_0}} w_{1,2}^{i_0} \right\} \xi^2 \label{eq:50085}\\
&& \mathrm{where}\hspace{.5cm} w_{1,2}= \frac{t_1 p_1}{u_1 (1-t_1)(1+p_1)} \nonumber
\end{eqnarray}
By using similar process for the previous cases of integral forms of $y_1(z)$ and $y_2(z)$, the integral form of sub-power series expansion of $y_3(z)$ is
\begin{eqnarray}
 y_3(z)&=& c_0 z^{\lambda } \int_{0}^{1} dt_3\;t_3^{\frac{1}{2}\left( 1-\frac{\Omega }{\mu }+\lambda \right)} \int_{0}^{\infty } du_3\;e^{-u_3} 
\frac{1}{2\pi i}  \int_{\infty }^{(0+)} dp_3\; \exp\left(\frac{ p_3}{\eta (1-t_3)}\right) p_3^{-1} (1+p_3)^{-\frac{1}{2}(3+\lambda )} \nonumber\\
&&\times \left( \frac{\eta u_3 (1-t_3)}{p_3}\right)^{\nu_3}  w_{3,3}^{-\frac{1}{2}\left(2-\omega +\lambda \right)} \left( w_{3,3}\partial_{w_{3,3}} \right) w_{3,3}^{\frac{1}{2}\left( 2-\omega +\lambda \right) } \nonumber\\
&&\times \int_{0}^{1} dt_2\;t_2^{\frac{1}{2}\left(-\frac{\Omega }{\mu }+\lambda \right)} \int_{0}^{\infty } du_2\;e^{-u_2} 
\frac{1}{2\pi i}  \int_{\infty }^{(0+)} dp_2\; \exp\left(\frac{ p_2}{w_{3,3} (1-t_2)}\right) p_2^{-1} (1+p_2)^{-\frac{1}{2}(2+\lambda )} \nonumber\\
&&\times \left( \frac{w_{3,3} u_2 (1-t_2)}{p_2}\right)^{\nu_2}  w_{2,3}^{-\frac{1}{2}\left(1-\omega +\lambda \right)} \left( w_{2,3}\partial_{w_{2,3}} \right) w_{2,3}^{\frac{1}{2}\left(1-\omega +\lambda \right) } \nonumber\\
&&\times  \int_{0}^{1} dt_1\;t_1^{\frac{1}{2}\left( -1-\frac{\Omega }{\mu }+\lambda \right)} \int_{0}^{\infty } du_1\;e^{-u_1} 
\frac{1}{2\pi i}  \int_{\infty }^{(0+)} dp_1\; \exp\left(\frac{ p_1}{w_{2,3} (1-t_1)}\right) p_1^{-1} (1+p_1)^{-\frac{1}{2}(1+\lambda )} \nonumber\\
&&\times \left( \frac{w_{2,3} u_1 (1-t_1)}{p_1}\right)^{\nu_1}  w_{1,3}^{-\frac{1}{2}\left( -\omega +\lambda \right)} \left( w_{1,3}\partial_{w_{1,3}} \right) w_{1,3}^{\frac{1}{2}\left( -\omega +\lambda \right) } \nonumber\\
&&\times  \left\{ \sum_{i_0=0}^{\nu _0} \frac{(-\nu _0)_{i_0}\left(\frac{\lambda }{2}\right)_{i_0}}{\left(1-\frac{\Omega }{2\mu }+\frac{\lambda }{2}\right)_{i_0}} w_{1,3}^{i_0} \right\} \xi^3 \label{eq:50086} 
\end{eqnarray}
where
\begin{equation}
\begin{cases} \displaystyle{w_{3,3} =\frac{t_3 p_3}{u_3 (1-t_3)(1+p_3)}} \cr
\displaystyle{w_{2,3} = \frac{t_2 p_2}{u_2 (1-t_2)(1+p_2)}} \cr
\displaystyle{w_{1,3}= \frac{t_1 p_1}{u_1 (1-t_1)(1+p_1)}}
\end{cases}
\nonumber
\end{equation}
By repeating this process for all higher terms of integral forms of sub-summation $y_m(z)$ terms where $m \geq 4$, we obtain every integral forms of $y_m(z)$ terms. 
Since we substitute (\ref{eq:50082a}), (\ref{eq:50083}), (\ref{eq:50085}), (\ref{eq:50086}) and including all integral forms of $y_m(z)$ terms where $m \geq 4$ into (\ref{eq:50081}), we obtain (\ref{eq:50080}).
\qed
\end{proof} 
\begin{remark}
The integral representation of the GCH equation of the first kind for polynomial of type 1 about $x=\infty $ as $\nu = 2 \nu _j+j +1+ \Omega /\mu $ where $j,\nu _j \in \mathbb{N}_{0}$ is
\begin{eqnarray}
y(z)&=& Q^{(i)}W_{\nu _j}\left( \mu ,\varepsilon , \Omega, \omega ,\nu = 2 \nu _j+j +1+ \frac{\Omega}{\mu}; z=\frac{1}{x}; \xi = -\frac{\varepsilon }{\mu}z; \eta = \frac{2}{\mu } z^2 \right)\nonumber\\
&=& z^{\frac{\Omega }{\mu }} \Bigg\{ \left( -\eta\right)^{\nu_0} U\left( -\nu_0, -\nu_0 +1-\frac{\Omega}{\mu}, -\eta^{-1} \right)\nonumber\\ 
&&+ \sum_{n=1}^{\infty } \Bigg\{\prod _{k=0}^{n-1} \Bigg\{ \int_{0}^{1} dt_{n-k}\;t_{n-k}^{\frac{1}{2}\left( n-k-2 \right)} \int_{0}^{\infty } du_{n-k}\;e^{-u_{n-k}}   \nonumber\\
&&\times \frac{1}{2\pi i} \int_{\infty }^{(0+)} dp_{n-k}\; p_{n-k}^{-1} (1+p_{n-k})^{-\frac{1}{2}\left(n-k+\frac{\Omega }{\mu }\right)} \exp \left( \frac{p_{n-k}}{w_{n-k+1,n}(1-t_{n-k})}\right)  \nonumber\\
&&\times \left( \frac{w_{n-k+1,n} u_{n-k} (1-t_{n-k})}{p_{n-k}}\right)^{\nu_{n-k}}   w_{n-k,n}^{-\frac{1}{2} \left(n-k-1-\omega +\frac{\Omega }{\mu}\right)} \left(  w_{n-k,n} \partial _{w_{n-k,n}} \right) w_{n-k,n}^{\frac{1}{2} \left(n-k-1-\omega +\frac{\Omega }{\mu}\right)}\Bigg\} \nonumber\\
&&\times   \left( -w_{1,n}\right)^{\nu_0} U\left( -\nu_0, -\nu_0 +1-\frac{\Omega}{\mu}, -w_{1,n}^{-1} \right) \Bigg\} \xi ^n \Bigg\}
\label{eq:50087}
\end{eqnarray}
\end{remark}
\begin{proof}
Replace $a$, $b$ and $z$ by $ -\nu_0$, $-\nu_0 +1 - \frac{\Omega }{2\mu }$ and $ -\eta^{-1} $ into (\ref{eq:50074}).
\begin{equation}
\sum_{j=0}^{\nu_0} \frac{(-\nu _0)_{i_0}\left(\frac{\lambda }{2}\right)_{i_0}}{\left(1-\frac{\Omega }{2\mu }+\frac{\lambda }{2}\right)_{i_0}} \eta^j  = \left( -\eta\right)^{\nu_0} U\left( -\nu_0, -\nu_0 +1-\frac{\Omega}{\mu}, -\eta^{-1} \right)  \label{eq:50088}
\end{equation} 
Replace $a$, $b$ and $z$ by $ -\nu_0$, $-\nu_0 +1 - \frac{\Omega }{2\mu }$ and $ -w_{1,n}^{-1} $ into (\ref{eq:50074}).
\begin{equation}
\sum_{j=0}^{\nu_0} \frac{(-\nu _0)_{i_0}\left(\frac{\lambda }{2}\right)_{i_0}}{\left(1-\frac{\Omega }{2\mu }+\frac{\lambda }{2}\right)_{i_0}} w_{1,n}^j  = \left( -w_{1,n}\right)^{\nu_0} U\left( -\nu_0, -\nu_0 +1-\frac{\Omega}{\mu}, -w_{1,n}^{-1} \right)  \label{eq:50089}
\end{equation} 
Put $c_0$= 1 and $\lambda =\Omega /\mu $ in (\ref{eq:50080}). Substitute (\ref{eq:50088}) and (\ref{eq:50089}) into the new (\ref{eq:50080}).
\qed
\end{proof}
\subsection{Generating function of the GCH polynomial of type 1}
Let's investigate the generating function for the GCH polynomial of type 1 of the first kind about $x=\infty $. 
\begin{definition}
I define that
\begin{equation}
\begin{cases}
\displaystyle { s_{a,b}} = \begin{cases} \displaystyle {  s_a\cdot s_{a+1}\cdot s_{a+2}\cdots s_{b-2}\cdot s_{b-1}\cdot s_b}\;\;\mbox{where}\;a>b \cr
s_a \;\;\mbox{only}\;\mbox{if}\;a=b\end{cases}
\cr
\cr
\displaystyle {  \widetilde{w}_{i,j}= \begin{cases} \displaystyle { \frac{s_i  t_i \widetilde{w}_{i+1,j}}{1+s_i u_i (1-t_i) \widetilde{w}_{i+1,j} }}\;\;\mbox{where}\;i<j \cr
\displaystyle { \frac{s_{i,\infty } t_i\eta }{1+s_{i,\infty } u_i (1-t_i) \eta }}  \;\;\mbox{only}\;\mbox{if}\;i=j\end{cases}}
\end{cases}\label{eq:50090}
\end{equation}
where
\begin{equation}
a,b,i,j\in \mathbb{N}_{0} \nonumber
\end{equation}
\end{definition}
And we have
\begin{equation}
\sum_{\nu _i = \nu _j}^{\infty } s_i^{\nu _i} = \frac{s_i^{\nu _j}}{(1-s_i)}\label{eq:50091}
\end{equation}
Acting the summation operator $\displaystyle{ \sum_{\nu_0 =0}^{\infty } \frac{s_0^{\nu _0}}{ \nu _0 !}  \prod _{n=1}^{\infty } \left\{ \sum_{ \nu _n = \nu _{n-1}}^{\infty } s_n^{\nu _n }\right\}}$ on (\ref{eq:50080}) where $|s_i|<1$ as $i=0,1,2,\cdots$ by using (\ref{eq:50090}) and (\ref{eq:50091}),
\begin{theorem}
The general expression of the generating function for the GCH polynomial of type 1 about $x=\infty $ is given by
\begin{eqnarray}
&& \sum_{\nu_0 =0}^{\infty } \frac{s_0^{\nu _0}}{ \nu _0 !}  \prod _{n=1}^{\infty } \left\{ \sum_{ \nu _n = \nu _{n-1}}^{\infty } s_n^{\nu _n }\right\} y(z) \nonumber\\
&&= \prod_{k=1}^{\infty } \frac{1}{(1-s_{k,\infty })} \mathbf{\Upsilon}(\lambda; s_{0,\infty } ;\eta )  \nonumber\\
&&+ \Bigg\{ \prod_{k=2}^{\infty } \frac{1}{(1-s_{k,\infty })} \int_{0}^{1} dt_1\;t_1^{\frac{1}{2}\left( -1-\frac{\Omega }{\mu}+\lambda \right)} \int_{0}^{\infty } du_1\; \exp\left( -(1-s_{1,\infty })u_1 \right) \left( 1+ s_{1,\infty } u_1 (1-t_1)\eta \right)^{-\frac{1}{2}\left(1+\lambda \right)}\nonumber\\
&&\times \widetilde{w}_{1,1}^{-\frac{1}{2}\left(-\omega +\lambda \right)} \left(  \widetilde{w}_{1,1} \partial _{\widetilde{w}_{1,1}} \right) \widetilde{w}_{1,1}^{ \frac{1}{2}\left(-\omega +\lambda \right)}\; \mathbf{\Upsilon}(\lambda ; s_0; \widetilde{w}_{1,1}) \Bigg\} \xi  \nonumber\\
&&+ \sum_{n=2}^{\infty } \Bigg\{ \prod_{k=n+1}^{\infty } \frac{1}{(1-s_{k,\infty })}  \int_{0}^{1} dt_n\;t_n^{\frac{1}{2}\left( n-2-\frac{\Omega }{\mu}+\lambda \right)} \int_{0}^{\infty } du_n\; \exp\left( -(1-s_{n,\infty })u_n \right) \nonumber\\
&&\times \left( 1+ s_{n,\infty }  u_n (1-t_n)\eta\right)^{-\frac{1}{2}\left(n +\lambda \right)} \widetilde{w}_{n,n}^{-\frac{1}{2}\left( n-1-\omega +\lambda \right)} \left(  \widetilde{w}_{n,n} \partial _{\widetilde{w}_{n,n}} \right) \widetilde{w}_{n,n}^{ \frac{1}{2}\left( n-1-\omega +\lambda \right)} \nonumber\\
&&\times  \prod_{j=1}^{n-1} \Bigg\{ \int_{0}^{1} dt_{n-j}\;t_{n-j}^{\frac{1}{2}\left( n-j-2-\frac{\Omega }{\mu}+\lambda \right)} \int_{0}^{\infty } du_{n-j}\; \exp\left( -(1-s_{n-j} )u_{n-j} \right)\nonumber\\
&&\times \left( 1+ s_{n-j} u_{n-j} (1-t_{n-j})\widetilde{w}_{n-j+1,n} \right)^{-\frac{1}{2}\left({n-j} +\lambda \right)}\nonumber\\
&&\times \widetilde{w}_{n-j,n}^{-\frac{1}{2}\left( n-j-1-\omega +\lambda \right)} \left(  \widetilde{w}_{n-j,n} \partial _{\widetilde{w}_{n-j,n}} \right) \widetilde{w}_{n-j,n}^{ \frac{1}{2}\left( n-j-1-\omega +\lambda \right)} \Bigg\}
 \mathbf{\Upsilon}(\lambda; s_0 ;\widetilde{w}_{1,n}) \Bigg\} \xi ^n  \label{eq:50092}
\end{eqnarray}
where
\begin{equation}
\begin{cases} 
{ \displaystyle \mathbf{\Upsilon}(\lambda; s_{0,\infty } ;\eta)= \sum_{\nu _0=0}^{\infty } \frac{s_{0,\infty }^{\nu _0}}{\nu_0!}  \left\{  c_0 z^{\lambda }  \sum_{i_0=0}^{\nu _0} \frac{(-\nu _0)_{i_0}\left(\frac{\lambda }{2}\right)_{i_0}}{\left(1-\frac{\Omega }{2\mu }+\frac{\lambda }{2}\right)_{i_0}} \eta^{i_0} \right\} }
\cr
{ \displaystyle \mathbf{\Upsilon}(\lambda ; s_0;\widetilde{w}_{1,1}) =  \sum_{\nu _0=0}^{\infty } \frac{s_0^{\nu _0}}{\nu_0!} \left\{  c_0 z^{\lambda } \sum_{i_0=0}^{\nu _0} \frac{(-\nu _0)_{i_0}\left(\frac{\lambda }{2}\right)_{i_0}}{\left(1-\frac{\Omega }{2\mu }+\frac{\lambda }{2}\right)_{i_0}} \widetilde{w}_{1,1}^{i_0}\right\} } \cr
{ \displaystyle \mathbf{\Upsilon}(\lambda; s_0 ;\widetilde{w}_{1,n}) = \sum_{\nu _0=0}^{\infty } \frac{s_0^{\nu _0}}{\nu_0!} \left\{  c_0 z^{\lambda }  \sum_{i_0=0}^{\nu _0} \frac{(-\nu _0)_{i_0}\left(\frac{\lambda }{2}\right)_{i_0}}{\left(1-\frac{\Omega }{2\mu }+\frac{\lambda }{2}\right)_{i_0}} \widetilde{w}_{1,n}^{i_0}\right\}}
\end{cases}\nonumber 
\end{equation}
\end{theorem}
\begin{proof} 
Acting the summation operator  $\displaystyle{ \sum_{\nu_0 =0}^{\infty } \frac{s_0^{\nu _0}}{ \nu _0 !}  \prod _{n=1}^{\infty } \left\{ \sum_{ \nu _n = \nu _{n-1}}^{\infty } s_n^{\nu _n }\right\}}$ on the form of integral of the type 1 GCH polynomial $y(z)$,
\begin{eqnarray}
&& \sum_{\nu_0 =0}^{\infty } \frac{s_0^{\nu _0}}{ \nu _0 !}  \prod _{n=1}^{\infty } \left\{ \sum_{ \nu _n = \nu _{n-1}}^{\infty } s_n^{\nu _n }\right\} y(z) \nonumber\\
&&=  \sum_{\nu_0 =0}^{\infty } \frac{s_0^{\nu _0}}{ \nu _0 !}  \prod _{n=1}^{\infty } \left\{ \sum_{ \nu _n = \nu _{n-1}}^{\infty } s_n^{\nu _n }\right\} \Big( y_0(z)+y_1(z)+y_2(z)+y_3(z)+ \cdots \Big) \label{eq:50093}
\end{eqnarray}
Acting the summation operator $\displaystyle{ \sum_{\nu_0 =0}^{\infty } \frac{s_0^{\nu _0}}{ \nu _0 !}  \prod _{n=1}^{\infty } \left\{ \sum_{ \nu _n = \nu _{n-1}}^{\infty } s_n^{\nu _n }\right\}}$ on  (\ref{eq:50082a}) by using (\ref{eq:50090}) and (\ref{eq:50091}),
\begin{eqnarray}
&&\sum_{\nu_0 =0}^{\infty } \frac{s_0^{\nu _0}}{ \nu _0 !}  \prod _{n=1}^{\infty } \left\{ \sum_{ \nu _n = \nu _{n-1}}^{\infty } s_n^{\nu _n }\right\} y_0(z)\nonumber\\
&& = \prod_{k=1}^{\infty } \frac{1}{(1-s_{k,\infty })} \sum_{\nu_0 =0}^{\infty } \frac{s_{0,\infty }^{\nu_0}}{\nu_0!} \left\{ c_0 z^{\lambda } \sum_{i_0=0}^{\nu _0} \frac{(-\nu _0)_{i_0}\left(\frac{\lambda }{2}\right)_{i_0}}{\left(1-\frac{\Omega }{2\mu }+\frac{\lambda }{2}\right)_{i_0}} \eta^{i_0} \right\} \label{eq:50094}
\end{eqnarray}
Acting the summation operator $\displaystyle{ \sum_{\nu_0 =0}^{\infty } \frac{s_0^{\nu _0}}{ \nu _0 !}  \prod _{n=1}^{\infty } \left\{ \sum_{ \nu _n = \nu _{n-1}}^{\infty } s_n^{\nu _n }\right\}}$ on (\ref{eq:50083}) by using (\ref{eq:50090}) and (\ref{eq:50091}),
\begin{eqnarray}
&&\sum_{\nu_0 =0}^{\infty } \frac{s_0^{\nu _0}}{ \nu _0 !}  \prod _{n=1}^{\infty } \left\{ \sum_{ \nu _n = \nu _{n-1}}^{\infty } s_n^{\nu _n }\right\} y_1(z) \nonumber\\
&&=  \prod_{k=2}^{\infty } \frac{1}{(1-s_{k,\infty })}  \int_{0}^{1} dt_1\;t_1^{\frac{1}{2}\left( -1-\frac{\Omega }{\mu }+\lambda \right)} \int_{0}^{\infty } du_1\;e^{-u_1} 
\frac{1}{2\pi i}  \int_{\infty }^{(0+)} dp_1\; \exp\left(\frac{ p_1}{\eta (1-t_1)}\right) \nonumber\\
&&\times p_1^{-1} (1+p_1)^{-\frac{1}{2}(1+\lambda )} \sum_{ \nu _1 = \nu _0}^{\infty } \left( \frac{s_{1,\infty }\eta u_1 (1-t_1)}{p_1}\right)^{\nu_1} w_{1,1}^{-\frac{1}{2}\left(-\omega +\lambda \right)} \left( w_{1,1}\partial_{w_{1,1}} \right) w_{1,1}^{\frac{1}{2}\left(-\omega +\lambda \right) } \nonumber\\
&&\times \sum_{\nu_0 =0}^{\infty } \frac{s_0^{\nu _0}}{ \nu _0 !} \left\{c_0 z^{\lambda } \sum_{i_0=0}^{\nu _0} \frac{(-\nu _0)_{i_0}\left(\frac{\lambda }{2}\right)_{i_0}}{\left(1-\frac{\Omega }{2\mu }+\frac{\lambda }{2}\right)_{i_0}} w_{1,1}^{i_0} \right\} \xi  \label{eq:50095}
\end{eqnarray}
Replace $\nu_i$, $\nu _j$ and $s_i$ by $\nu_1$, $\nu _0$ and ${ \displaystyle \frac{s_{1,\infty }\eta u_1 (1-t_1)}{p_1}}$ in (\ref{eq:50091}). Take the new (\ref{eq:50091}) into (\ref{eq:50095}).
\begin{eqnarray}
&&\sum_{\nu_0 =0}^{\infty } \frac{s_0^{\nu _0}}{ \nu _0 !}  \prod _{n=1}^{\infty } \left\{ \sum_{ \nu _n = \nu _{n-1}}^{\infty } s_n^{\nu _n }\right\} y_1(z) \nonumber\\
&&=  \prod_{k=2}^{\infty } \frac{1}{(1-s_{k,\infty })}  \int_{0}^{1} dt_1\;t_1^{\frac{1}{2}\left( -1-\frac{\Omega }{\mu }+\lambda \right)} \int_{0}^{\infty } du_1\;e^{-u_1} 
\frac{1}{2\pi i}  \int_{\infty }^{(0+)} dp_1\; \exp\left(\frac{ p_1}{\eta (1-t_1)}\right) \nonumber\\
&&\times \frac{(1+p_1)^{-\frac{1}{2}(1+\lambda )}}{p_1-s_{1,\infty } \eta u_1 (1-t_1)}  w_{1,1}^{-\frac{1}{2}\left(-\omega +\lambda \right)} \left( w_{1,1}\partial_{w_{1,1}} \right) w_{1,1}^{\frac{1}{2}\left(-\omega +\lambda \right) } \nonumber\\
&&\times  \sum_{\nu_0 =0}^{\infty } \left( \frac{s_{0,\infty }\eta u_1 (1-t_1)}{p_1}\right)^{\nu_0}  \frac{1}{ \nu _0 !} \left\{c_0 z^{\lambda } \sum_{i_0=0}^{\nu _0} \frac{(-\nu _0)_{i_0}\left(\frac{\lambda }{2}\right)_{i_0}}{\left(1-\frac{\Omega }{2\mu }+\frac{\lambda }{2}\right)_{i_0}} w_{1,1}^{i_0} \right\} \xi  \label{eq:50096}
\end{eqnarray}
By using Cauchy's integral formula, the contour integrand has poles at $p_1 = s_{1,\infty }\eta u_1 (1-t_1)$,
and $s_{1,\infty }\eta u_1 (1-t_1)$ is inside the unit circle. As we compute the residue there in (\ref{eq:50096}) we obtain
\begin{eqnarray}
&&\sum_{\nu_0 =0}^{\infty } \frac{s_0^{\nu _0}}{ \nu _0 !}  \prod _{n=1}^{\infty } \left\{ \sum_{ \nu _n = \nu _{n-1}}^{\infty } s_n^{\nu _n }\right\} y_1(z) \nonumber\\
&&=  \prod_{k=2}^{\infty } \frac{1}{(1-s_{k,\infty })}  \int_{0}^{1} dt_1\;t_1^{\frac{1}{2}\left( -1-\frac{\Omega }{\mu }+\lambda \right)} \int_{0}^{\infty } du_1\;\exp\left( -(1-s_{1,\infty })u_1\right) (1+s_{1,\infty } u_1 (1-t_1)\eta)^{-\frac{1}{2}(1+\lambda )}  \nonumber\\
&&\times  \widetilde{w}_{1,1}^{-\frac{1}{2}\left(-\omega +\lambda \right)} \left( \widetilde{w}_{1,1}\partial_{\widetilde{w}_{1,1}} \right) \widetilde{w}_{1,1}^{\frac{1}{2}\left(-\omega +\lambda \right) }  \sum_{\nu_0 =0}^{\infty }  \frac{s_0^{\nu _0}}{ \nu _0 !} \left\{c_0 z^{\lambda } \sum_{i_0=0}^{\nu _0} \frac{(-\nu _0)_{i_0}\left(\frac{\lambda }{2}\right)_{i_0}}{\left(1-\frac{\Omega }{2\mu }+\frac{\lambda }{2}\right)_{i_0}} \widetilde{w}_{1,1}^{i_0} \right\} \xi  \label{eq:50097}
\end{eqnarray}
where
\begin{equation}
\widetilde{w}_{1,1} =  \frac{t_1 p_1}{u_1 (1-t_1)(1+p_1)}\Bigg|_{p_1 = s_{1,\infty }\eta u_1 (1-t_1)} = \frac{s_{1,\infty }t_1 \eta}{1+s_{1,\infty }u_1 (1-t_1)\eta} \nonumber
\end{equation}
Acting the summation operator $\displaystyle{ \sum_{\nu_0 =0}^{\infty } \frac{s_0^{\nu _0}}{ \nu _0 !}  \prod _{n=1}^{\infty } \left\{ \sum_{ \nu _n = \nu _{n-1}}^{\infty } s_n^{\nu _n }\right\}}$ on (\ref{eq:50085}) by using (\ref{eq:50090}) and (\ref{eq:50091}),
\begin{eqnarray}
&&\sum_{\nu_0 =0}^{\infty } \frac{s_0^{\nu _0}}{ \nu _0 !}  \prod _{n=1}^{\infty } \left\{ \sum_{ \nu _n = \nu _{n-1}}^{\infty } s_n^{\nu _n }\right\} y_2(z) \nonumber\\
&&=  \prod_{k=3}^{\infty } \frac{1}{(1-s_{k,\infty })}  \int_{0}^{1} dt_2\;t_2^{\frac{1}{2}\left( -\frac{\Omega }{\mu }+\lambda \right)} \int_{0}^{\infty } du_2\;e^{-u_2} 
\frac{1}{2\pi i}  \int_{\infty }^{(0+)} dp_2\; \exp\left(\frac{ p_2}{\eta (1-t_2)}\right) \nonumber\\
&&\times p_2^{-1} (1+p_2)^{-\frac{1}{2}(2+\lambda )}  \sum_{ \nu _2 = \nu _1}^{\infty } \left( \frac{s_{2,\infty }\eta u_2 (1-t_2)}{p_2}\right)^{\nu_2} w_{2,2}^{-\frac{1}{2}\left( 1-\omega +\lambda \right)} \left( w_{2,2}\partial_{w_{2,2}} \right) w_{2,2}^{\frac{1}{2}\left( 1-\omega +\lambda \right) } \nonumber\\
&&\times  \int_{0}^{1} dt_1\;t_1^{\frac{1}{2}\left( -1-\frac{\Omega }{\mu }+\lambda \right)} \int_{0}^{\infty } du_1\;e^{-u_1} 
\frac{1}{2\pi i}  \int_{\infty }^{(0+)} dp_1\; \exp\left(\frac{ p_1}{w_{2,2} (1-t_1)}\right) p_1^{-1} (1+p_1)^{-\frac{1}{2}(1+\lambda )} \nonumber\\
&&\times \sum_{ \nu _1 = \nu _0}^{\infty } \left( \frac{s_1 w_{2,2} u_1 (1-t_1)}{p_1}\right)^{\nu_1} w_{1,2}^{-\frac{1}{2}\left(-\omega +\lambda \right)} \left( w_{1,2}\partial_{w_{1,2}} \right) w_{1,2}^{\frac{1}{2}\left(-\omega +\lambda \right) } \nonumber\\
&&\times \sum_{\nu_0 =0}^{\infty } \frac{s_0^{\nu _0}}{ \nu _0 !} \left\{c_0 z^{\lambda } \sum_{i_0=0}^{\nu _0} \frac{(-\nu _0)_{i_0}\left(\frac{\lambda }{2}\right)_{i_0}}{\left(1-\frac{\Omega }{2\mu }+\frac{\lambda }{2}\right)_{i_0}} w_{1,2}^{i_0} \right\} \xi^2  \label{eq:50098}
\end{eqnarray}
Replace $\nu_i$, $\nu_j$ and $s_i$ by $\nu_2$, $\nu_1$ and ${ \displaystyle \frac{s_{2,\infty }\eta u_2 (1-t_2)}{p_2}}$ in (\ref{eq:50091}). Take the new (\ref{eq:50091}) into (\ref{eq:50098}).
\begin{eqnarray}
&&\sum_{\nu_0 =0}^{\infty } \frac{s_0^{\nu _0}}{ \nu _0 !}  \prod _{n=1}^{\infty } \left\{ \sum_{ \nu _n = \nu _{n-1}}^{\infty } s_n^{\nu _n }\right\} y_2(z) \nonumber\\
&&=  \prod_{k=3}^{\infty } \frac{1}{(1-s_{k,\infty })}  \int_{0}^{1} dt_2\;t_2^{\frac{1}{2}\left( -\frac{\Omega }{\mu }+\lambda \right)} \int_{0}^{\infty } du_2\;e^{-u_2} 
\frac{1}{2\pi i}  \int_{\infty }^{(0+)} dp_2\; \exp\left(\frac{ p_2}{\eta (1-t_2)}\right) \nonumber\\
&&\times \frac{(1+p_2)^{-\frac{1}{2}(2+\lambda )}}{p_2-s_{2,\infty } \eta u_2 (1-t_2)}\; w_{2,2}^{-\frac{1}{2}\left( 1-\omega +\lambda \right)} \left( w_{2,2}\partial_{w_{2,2}} \right) w_{2,2}^{\frac{1}{2}\left( 1-\omega +\lambda \right) } \nonumber\\
&&\times \int_{0}^{1} dt_1\;t_1^{\frac{1}{2}\left( -1-\frac{\Omega }{\mu }+\lambda \right)} \int_{0}^{\infty } du_1\;e^{-u_1} 
\frac{1}{2\pi i}  \int_{\infty }^{(0+)} dp_1\; \exp\left(\frac{ p_1}{w_{2,2} (1-t_1)}\right) p_1^{-1} (1+p_1)^{-\frac{1}{2}(1+\lambda )} \nonumber\\ \nonumber\\
&&\times \sum_{ \nu _1 = \nu _0}^{\infty } \left( \frac{s_{1,\infty } \eta u_2(1-t_2) w_{2,2} u_1 (1-t_1)}{p_1 p_2}\right)^{\nu_1} w_{1,2}^{-\frac{1}{2}\left( -\omega +\lambda \right)} \left( w_{1,2}\partial_{w_{1,2}} \right) w_{1,2}^{\frac{1}{2}\left( -\omega +\lambda \right) } \nonumber\\
&&\times  \sum_{\nu_0 =0}^{\infty } \frac{s_0^{\nu_0}}{ \nu _0 !} \left\{c_0 z^{\lambda } \sum_{i_0=0}^{\nu _0} \frac{(-\nu _0)_{i_0}\left(\frac{\lambda }{2}\right)_{i_0}}{\left(1-\frac{\Omega }{2\mu }+\frac{\lambda }{2}\right)_{i_0}} w_{1,2}^{i_0} \right\} \xi^2 \label{eq:50099}
\end{eqnarray}
By using Cauchy's integral formula, the contour integrand has poles at $p_2 = s_{2,\infty }\eta u_2 (1-t_2)$,
and $s_{2,\infty }\eta u_2 (1-t_2)$ is inside the unit circle. As we compute the residue there in (\ref{eq:50099}) we obtain
\begin{eqnarray}
&&\sum_{\nu_0 =0}^{\infty } \frac{s_0^{\nu _0}}{ \nu _0 !}  \prod _{n=1}^{\infty } \left\{ \sum_{ \nu _n = \nu _{n-1}}^{\infty } s_n^{\nu _n }\right\} y_2(z) \nonumber\\
&&=   \prod_{k=3}^{\infty } \frac{1}{(1-s_{k,\infty })}  \int_{0}^{1} dt_2\;t_2^{\frac{1}{2}\left( -\frac{\Omega }{\mu }+\lambda \right)} \int_{0}^{\infty } du_2\;\exp\left( -(1-s_{2,\infty })u_2\right) (1+s_{2,\infty } u_2 (1-t_2)\eta)^{-\frac{1}{2}(2+\lambda )}  \nonumber\\
&&\times  \widetilde{w}_{2,2}^{-\frac{1}{2}\left( 1-\omega +\lambda \right)} \left( \widetilde{w}_{2,2}\partial_{\widetilde{w}_{2,2}} \right) \widetilde{w}_{2,2}^{\frac{1}{2}\left( 1-\omega +\lambda \right) } \nonumber\\
&&\times \int_{0}^{1} dt_1\;t_1^{\frac{1}{2}\left( -1-\frac{\Omega }{\mu }+\lambda \right)} \int_{0}^{\infty } du_1\;e^{-u_1} \frac{1}{2\pi i}  \int_{\infty }^{(0+)} dp_1\; \exp\left(\frac{ p_1}{\widetilde{w}_{2,2} (1-t_1)}\right) p_1^{-1} (1+p_1)^{-\frac{1}{2}(1+\lambda )} \nonumber\\
&&\times \sum_{ \nu _1 = \nu _0}^{\infty } \left( \frac{s_1 \widetilde{w}_{2,2} u_1 (1-t_1)}{p_1}\right)^{\nu_1} w_{1,2}^{-\frac{1}{2}\left( -\omega +\lambda \right)} \left( w_{1,2}\partial_{w_{1,2}} \right) w_{1,2}^{\frac{1}{2}\left( -\omega +\lambda \right) } \nonumber\\
&&\times \sum_{\nu_0 =0}^{\infty }  \frac{s_0^{\nu _0}}{ \nu _0 !} \left\{c_0 z^{\lambda } \sum_{i_0=0}^{\nu _0} \frac{(-\nu _0)_{i_0}\left(\frac{\lambda }{2}\right)_{i_0}}{\left(1-\frac{\Omega }{2\mu }+\frac{\lambda }{2}\right)_{i_0}} w_{1,2}^{i_0} \right\} \xi^2 \label{eq:500100}
\end{eqnarray}
where
\begin{equation}
\widetilde{w}_{2,2}  =   \frac{t_2 p_2}{u_2 (1-t_2)(1+p_2)}\Bigg|_{p_2 = s_{2,\infty }\eta u_2 (1-t_2)} = \frac{s_{2,\infty }t_2 \eta}{1+s_{2,\infty }u_2 (1-t_2)\eta} \nonumber
\end{equation}
Replace $\nu _i$, $\nu _j$ and $s_i$ by $\nu _1$, $\nu _0$ and ${ \displaystyle \frac{s_1 \widetilde{w}_{2,2} u_1 (1-t_1)}{p_1}}$ in (\ref{eq:50091}). Take the new (\ref{eq:50091}) into (\ref{eq:500100}).
\begin{eqnarray}
&&\sum_{\nu_0 =0}^{\infty } \frac{s_0^{\nu _0}}{ \nu _0 !}  \prod _{n=1}^{\infty } \left\{ \sum_{ \nu _n = \nu _{n-1}}^{\infty } s_n^{\nu _n }\right\} y_2(z) \nonumber\\
&&=   \prod_{k=3}^{\infty } \frac{1}{(1-s_{k,\infty })}  \int_{0}^{1} dt_2\;t_2^{\frac{1}{2}\left( -\frac{\Omega }{\mu }+\lambda \right)} \int_{0}^{\infty } du_2\;\exp\left( -(1-s_{2,\infty })u_2\right) (1+s_{2,\infty } u_2 (1-t_2)\eta)^{-\frac{1}{2}(2+\lambda )}  \nonumber\\
&&\times  \widetilde{w}_{2,2}^{-\frac{1}{2}\left( 1-\omega +\lambda \right)} \left( \widetilde{w}_{2,2}\partial_{\widetilde{w}_{2,2}} \right) \widetilde{w}_{2,2}^{\frac{1}{2}\left( 1-\omega +\lambda \right) } \nonumber\\
&&\times \int_{0}^{1} dt_1\;t_1^{\frac{1}{2}\left( -1-\frac{\Omega }{\mu }+\lambda \right)} \int_{0}^{\infty } du_1\;e^{-u_1} 
\frac{1}{2\pi i}  \int_{\infty }^{(0+)} dp_1\; \exp\left(\frac{ p_1}{\widetilde{w}_{2,2} (1-t_1)}\right)\frac{(1+p_1)^{-\frac{1}{2}(1+\lambda )}}{p_1-s_1 \widetilde{w}_{2,2} u_1 (1-t_1)} \nonumber\\
&&\times  w_{1,2}^{-\frac{1}{2}\left(-\omega +\lambda \right)} \left( w_{1,2}\partial_{w_{1,2}} \right) w_{1,2}^{\frac{1}{2}\left(-\omega +\lambda \right) } \nonumber\\
&&\times  \sum_{\nu_0 =0}^{\infty } \left( \frac{s_{0,1}\widetilde{w}_{2,2} u_1 (1-t_1)}{p_1}\right)^{\nu_0}  \frac{1}{ \nu _0 !} \left\{c_0 z^{\lambda } \sum_{i_0=0}^{\nu _0} \frac{(-\nu _0)_{i_0}\left(\frac{\lambda }{2}\right)_{i_0}}{\left(1-\frac{\Omega }{2\mu }+\frac{\lambda }{2}\right)_{i_0}} w_{1,2}^{i_0} \right\} \xi^2 \label{eq:500101}
\end{eqnarray}
By using Cauchy's integral formula, the contour integrand has poles at $p_1 = s_1\widetilde{w}_{2,2} u_1 (1-t_1)$,
and $s_1\widetilde{w}_{2,2} u_1 (1-t_1)$ is inside the unit circle. As we compute the residue there in (\ref{eq:500101}) we obtain
\begin{eqnarray}
&&\sum_{\nu_0 =0}^{\infty } \frac{s_0^{\nu _0}}{ \nu _0 !}  \prod _{n=1}^{\infty } \left\{ \sum_{ \nu _n = \nu _{n-1}}^{\infty } s_n^{\nu _n }\right\} y_2(z) \nonumber\\
&&=   \prod_{k=3}^{\infty } \frac{1}{(1-s_{k,\infty })}  \int_{0}^{1} dt_2\;t_2^{\frac{1}{2}\left( -\frac{\Omega }{\mu }+\lambda \right)} \int_{0}^{\infty } du_2\;\exp\left( -(1-s_{2,\infty })u_2\right) (1+s_{2,\infty } u_2 (1-t_2)\eta)^{-\frac{1}{2}(2+\lambda )}  \nonumber\\
&&\times  \widetilde{w}_{2,2}^{-\frac{1}{2}\left( 1-\omega +\lambda \right)} \left( \widetilde{w}_{2,2}\partial_{\widetilde{w}_{2,2}} \right) \widetilde{w}_{2,2}^{\frac{1}{2}\left( 1-\omega +\lambda \right) } \nonumber\\
&&\times \int_{0}^{1} dt_1\;t_1^{\frac{1}{2}\left( -1-\frac{\Omega }{\mu }+\lambda \right)} \int_{0}^{\infty } du_1\;\exp\left( -(1-s_1)u_1\right) (1+s_1 u_1 (1-t_1)\widetilde{w}_{2,2})^{-\frac{1}{2}(1+\lambda )}  \nonumber\\
&&\times  \widetilde{w}_{1,2}^{-\frac{1}{2}\left(-\omega +\lambda \right)} \left( \widetilde{w}_{1,2}\partial_{\widetilde{w}_{1,2}} \right) \widetilde{w}_{1,2}^{\frac{1}{2}\left(-\omega +\lambda \right) }  \sum_{\nu_0 =0}^{\infty }  \frac{s_0^{\nu _0}}{ \nu _0 !} \left\{c_0 z^{\lambda } \sum_{i_0=0}^{\nu _0} \frac{(-\nu _0)_{i_0}\left(\frac{\lambda }{2}\right)_{i_0}}{\left(1-\frac{\Omega }{2\mu }+\frac{\lambda }{2}\right)_{i_0}} \widetilde{w}_{1,2}^{i_0} \right\} \xi^2 \label{eq:500102}
\end{eqnarray}
where
\begin{equation}
\widetilde{w}_{1,2}  =   \frac{t_1 p_1}{u_1 (1-t_1)(1+p_1)}\Bigg|_{p_1 = s_1\widetilde{w}_{2,2} u_1 (1-t_1)} = \frac{s_1 t_1 \widetilde{w}_{2,2}}{1+s_1 u_1 (1-t_1)\widetilde{w}_{2,2}} \nonumber
\end{equation}
Acting the summation operator $\displaystyle{ \sum_{\nu_0 =0}^{\infty } \frac{s_0^{\nu _0}}{ \nu _0 !}  \prod _{n=1}^{\infty } \left\{ \sum_{ \nu _n = \nu _{n-1}}^{\infty } s_n^{\nu _n }\right\}}$ on (\ref{eq:50086}) by using (\ref{eq:50090}) and (\ref{eq:50091}),
\begin{eqnarray}
&&\sum_{\nu_0 =0}^{\infty } \frac{s_0^{\nu _0}}{ \nu _0 !}  \prod _{n=1}^{\infty } \left\{ \sum_{ \nu _n = \nu _{n-1}}^{\infty } s_n^{\nu _n }\right\} y_3(z) \nonumber\\
&&= \prod_{k=4}^{\infty } \frac{1}{(1-s_{k,\infty })}\int_{0}^{1} dt_3\;t_3^{\frac{1}{2}\left( 1-\frac{\Omega }{\mu }+\lambda \right)} \int_{0}^{\infty } du_3\;\exp\left( -(1-s_{3,\infty })u_3\right) (1+s_{3,\infty } u_3 (1-t_3)\eta)^{-\frac{1}{2}(3+\lambda )}  \nonumber\\
&&\times  \widetilde{w}_{3,3}^{-\frac{1}{2}\left( 2-\omega +\lambda \right)} \left( \widetilde{w}_{3,3}\partial_{\widetilde{w}_{3,3}} \right) \widetilde{w}_{3,3}^{\frac{1}{2}\left( 2-\omega +\lambda \right) } \nonumber\\
&&\times \int_{0}^{1} dt_2\;t_2^{\frac{1}{2}\left( -\frac{\Omega }{\mu }+\lambda \right)} \int_{0}^{\infty } du_2\;\exp\left( -(1-s_2)u_2\right) (1+s_2 u_2 (1-t_2)\widetilde{w}_{3,3})^{-\frac{1}{2}(2+\lambda )}  \nonumber\\
&&\times  \widetilde{w}_{2,3}^{-\frac{1}{2}\left( 1-\omega +\lambda \right)} \left( \widetilde{w}_{2,3}\partial_{\widetilde{w}_{2,3}} \right) \widetilde{w}_{2,3}^{\frac{1}{2}\left( 1-\omega +\lambda \right) } \nonumber\\
&&\times \int_{0}^{1} dt_1\;t_1^{\frac{1}{2}\left( -1-\frac{\Omega }{\mu }+\lambda \right)} \int_{0}^{\infty } du_1\;\exp\left( -(1-s_1)u_1\right) (1+s_1 u_1 (1-t_1)\widetilde{w}_{2,3})^{-\frac{1}{2}(1+\lambda )}  \nonumber\\
&&\times  \widetilde{w}_{1,3}^{-\frac{1}{2}\left( -\omega +\lambda \right)} \left( \widetilde{w}_{1,3}\partial_{\widetilde{w}_{1,3}} \right) \widetilde{w}_{1,3}^{\frac{1}{2}\left( -\omega +\lambda \right) } \sum_{\nu_0 =0}^{\infty }  \frac{s_0^{\nu _0}}{ \nu _0 !} \left\{c_0 z^{\lambda } \sum_{i_0=0}^{\nu _0} \frac{(-\nu _0)_{i_0}\left(\frac{\lambda }{2}\right)_{i_0}}{\left(1-\frac{\Omega }{2\mu }+\frac{\lambda }{2}\right)_{i_0}} \widetilde{w}_{1,3}^{i_0} \right\} \xi^3 \label{eq:500103}
\end{eqnarray}
where
\begin{equation}
\begin{cases} \displaystyle{\widetilde{w}_{3,3} = \frac{t_3 p_3}{u_3 (1-t_3)(1+p_3)}\Bigg|_{p_3 = s_{3,\infty }\eta u_3 (1-t_3)} = \frac{s_{3,\infty }t_3 \eta}{1+s_{3,\infty }u_3 (1-t_3)\eta}}  \cr
\displaystyle{ \widetilde{w}_{2,3}  =   \frac{t_2 p_2}{u_2 (1-t_2)(1+p_2)}\Bigg|_{p_2 = s_2\widetilde{w}_{3,3} u_2 (1-t_2)} = \frac{s_2 t_2 \widetilde{w}_{3,3}}{1+s_2 u_2 (1-t_2)\widetilde{w}_{3,3}}} \cr
\displaystyle{ \widetilde{w}_{1,3}  =   \frac{t_1 p_1}{u_1 (1-t_1)(1+p_1)}\Bigg|_{p_1 = s_1\widetilde{w}_{2,3} u_1 (1-t_1)} = \frac{s_1 t_1 \widetilde{w}_{2,3}}{1+s_1 u_1 (1-t_1)\widetilde{w}_{2,3}}}
\end{cases}\nonumber 
\end{equation}
By repeating this process for all higher terms of integral forms of sub-summation $y_m(z)$ terms where $m > 3$, we obtain every  $\displaystyle{ \sum_{\nu_0 =0}^{\infty } \frac{s_0^{\nu _0}}{ \nu _0 !}  \prod _{n=1}^{\infty } \left\{ \sum_{ \nu _n = \nu _{n-1}}^{\infty } s_n^{\nu _n }\right\} y_m(z)}$ terms. 
Since we substitute (\ref{eq:50094}), (\ref{eq:50097}), (\ref{eq:500102}), (\ref{eq:500103}) and including all $\displaystyle{ \sum_{\nu_0 =0}^{\infty } \frac{s_0^{\nu _0}}{ \nu _0 !}  \prod _{n=1}^{\infty } \left\{ \sum_{ \nu _n = \nu _{n-1}}^{\infty } s_n^{\nu _n }\right\} y_m(z)}$ terms where $m > 3$ into (\ref{eq:50093}), we obtain (\ref{eq:50092})
\qed
\end{proof}
\begin{remark}
The generating function for the GCH polynomial of type 1 of the first kind about $x=\infty $ as $\nu = 2 \nu _j+j +1+ \Omega /\mu $ where $j,\nu _j \in \mathbb{N}_{0}$ is
\begin{eqnarray}
&&\sum_{\nu_0 =0}^{\infty } \frac{s_0^{\nu _0}}{ \nu _0 !}  \prod _{n=1}^{\infty } \left\{ \sum_{ \nu _n = \nu _{n-1}}^{\infty } s_n^{\nu _n }\right\} Q^{(i)}W_{\nu _j}\left( \mu ,\varepsilon, \Omega, \omega, \nu = 2 \nu _j+j +1+ \frac{\Omega}{\mu}; z=\frac{1}{x}; \xi = -\frac{\varepsilon }{\mu}z; \eta = \frac{2}{\mu } z^2 \right) \nonumber\\
&&= z^{\frac{\Omega }{\mu}} \Bigg\{ \prod_{k=1}^{\infty } \frac{1}{(1-s_{k,\infty })} \mathbf{A} \left( s_{0,\infty } ;\eta \right) \nonumber\\
&&+ \Bigg\{ \prod_{k=2}^{\infty } \frac{1}{(1-s_{k,\infty })} \int_{0}^{1} dt_1\;t_1^{-\frac{1}{2}} \int_{0}^{\infty } du_1\; \overleftrightarrow {\mathbf{\Gamma}}_1 \left(s_{1,\infty };t_1,u_1,\eta \right) \nonumber\\
&&\times \widetilde{w}_{1,1}^{-\frac{1}{2}\left(-\omega +\frac{\Omega }{\mu} \right)} \left(  \widetilde{w}_{1,1} \partial _{\widetilde{w}_{1,1}} \right) \widetilde{w}_{1,1}^{ \frac{1}{2}\left(-\omega +\frac{\Omega }{\mu} \right)}\; \mathbf{A} \left(s_{0} ;\widetilde{w}_{1,1} \right) \Bigg\} \xi  \nonumber\\
&&+ \sum_{n=2}^{\infty } \Bigg\{ \prod_{k=n+1}^{\infty } \frac{1}{(1-s_{k,\infty })}  \int_{0}^{1} dt_n\;t_n^{\frac{1}{2}\left( n-2\right)} \int_{0}^{\infty } du_n\; \overleftrightarrow {\mathbf{\Gamma}}_n \left(s_{n,\infty };t_n,u_n,\eta \right)\nonumber\\
&&\times \widetilde{w}_{n,n}^{-\frac{1}{2}\left( n-1-\omega +\frac{\Omega}{\mu} \right)} \left(  \widetilde{w}_{n,n} \partial _{\widetilde{w}_{n,n}} \right) \widetilde{w}_{n,n}^{ \frac{1}{2}\left( n-1-\omega +\frac{\Omega}{\mu} \right)} \nonumber\\
&&\times  \prod_{j=1}^{n-1} \Bigg\{ \int_{0}^{1} dt_{n-j}\;t_{n-j}^{\frac{1}{2}\left( n-j-2 \right)} \int_{0}^{\infty } du_{n-j}\; \overleftrightarrow {\mathbf{\Gamma}}_{n-j} \left(s_{n-j};t_{n-j},u_{n-j},\widetilde{w}_{n-j+1,n} \right) \nonumber\\
&&\times \widetilde{w}_{n-j,n}^{-\frac{1}{2}\left( n-j-1-\omega +\frac{\Omega}{\mu} \right)} \left(  \widetilde{w}_{n-j,n} \partial _{\widetilde{w}_{n-j,n}} \right) \widetilde{w}_{n-j,n}^{ \frac{1}{2}\left( n-j-1-\omega +\frac{\Omega}{\mu}\right)} \Bigg\}\mathbf{A} \left(s_{0} ;\widetilde{w}_{1,n} \right) \Bigg\} \xi ^n \Bigg\}  \label{eq:500104}
\end{eqnarray}
where
\begin{equation}
\begin{cases} 
{ \displaystyle \overleftrightarrow {\mathbf{\Gamma}}_1 \left(s_{1,\infty };t_1,u_1,\eta \right)=  \exp\left( -(1-s_{1,\infty })u_1 \right) \left( 1+ s_{1,\infty } u_1 (1-t_1)\eta \right)^{-\frac{1}{2}\left(1+\frac{\Omega }{\mu} \right)} }\cr
{ \displaystyle  \overleftrightarrow {\mathbf{\Gamma}}_n \left(s_{n,\infty };t_n,u_n,\eta \right) = \exp\left( -(1-s_{n,\infty })u_n \right) \left( 1+ s_{n,\infty }  u_n (1-t_n)\eta\right)^{-\frac{1}{2}\left(n +\frac{\Omega}{\mu} \right)} }\cr
{ \displaystyle \overleftrightarrow {\mathbf{\Gamma}}_{n-j} \left(s_{n-j};t_{n-j},u_{n-j},\widetilde{w}_{n-j+1,n} \right)  } \cr
{ \displaystyle = \exp\left( -(1-s_{n-j} )u_{n-j} \right) \left( 1+ s_{n-j} u_{n-j} (1-t_{n-j})\widetilde{w}_{n-j+1,n} \right)^{-\frac{1}{2}\left({n-j} +\frac{\Omega}{\mu} \right)}}
\end{cases}\nonumber 
\end{equation}
and
\begin{equation}
\begin{cases} 
{ \displaystyle \mathbf{A} \left( s_{0,\infty } ;\eta \right)= \exp\left( s_{0,\infty }\right)  (1+s_{0,\infty } \eta)^{-\frac{\Omega }{2\mu}}}\cr
{ \displaystyle  \mathbf{A} \left(s_{0} ;\widetilde{w}_{1,1} \right) = \exp\left( s_0 \right) (1+s_0 \widetilde{w}_{1,1})^{-\frac{\Omega }{2\mu}} } \cr
{ \displaystyle \mathbf{A} \left(s_{0} ;\widetilde{w}_{1,n} \right) = \exp\left( s_0 \right) (1+s_0 \widetilde{w}_{1,n})^{-\frac{\Omega }{2\mu}} }
\end{cases}\nonumber 
\end{equation}
\end{remark}
\begin{proof}
Replace $a$, $b$, $j$  and $z$ by $-\nu_0$, $-\nu_0+1-a$, $i_0$ and $-z^{-1}$ in (\ref{eq:50074}). Acting the summation operator $\displaystyle{ \sum_{\nu_0 =0}^{\infty } \frac{s_0^{\nu _0}}{ \nu _0 !}}$ on the new (\ref{eq:50074}) 
\begin{equation}
\sum_{\nu_0 =0}^{\infty } \frac{s_0^{\nu _0}}{ \nu _0 !} \sum_{i_0=0}^{\nu_0} \frac{(-\nu_0)_{i_0} (a)_{i_0}}{(1)_{i_0}} z^{i_0} = \sum_{\nu_0 =0}^{\infty } \frac{(-s_0 z)^{\nu _0}}{ \nu _0 !} U\left(-\nu_0, -\nu_0 +1-a,-z^{-1}\right)   
\label{eq:500105}
\end{equation}
Replace $a$, $b$, $p_l$ and $z$ by $-\nu_0$, $-\nu_0+1-a$, $p$ and $-z^{-1}$ in (\ref{eq:50073}).
\begin{equation}
U\left(-\nu_0, -\nu_0 +1-a,-z^{-1}\right) =  e^{\nu_0\pi i}\frac{\nu_0 !}{2\pi i} \int_{\infty }^{(0+)} dp\; e^{\frac{p}{z}} p^{-\nu_0-1} (1+p)^{-a} 
\label{eq:500106}
\end{equation}
Put (\ref{eq:500106}) in (\ref{eq:500105}).
\begin{eqnarray}
\sum_{\nu_0 =0}^{\infty } \frac{s_0^{\nu _0}}{ \nu _0 !} \sum_{i_0=0}^{\nu_0}\frac{(-\nu_0)_{i_0} (a)_{i_0}}{(1)_{i_0}} z^{i_0} &=& \frac{1}{2\pi i} \int_{\infty }^{(0+)} dp\; e^{\frac{p}{z}} p^{-1} (1+p)^{-a} \sum_{\nu_0 =0}^{\infty } \left(\frac{s_0 z}{p}\right)^{\nu_0} \nonumber\\
&=& \frac{1}{2\pi i} \int_{\infty }^{(0+)} dp\; e^{\frac{p}{z}} \frac{(1+p)^{-a}}{(p-s_0 z)} \nonumber\\
&=& \exp(s_0) (1+s_0 z)^{-a}
\label{eq:500107}
\end{eqnarray}
Replace $s_0$, $a$ and $z$  by $s_{0,\infty }$ , $\frac{\Omega }{2\mu}$ and $\eta$ in (\ref{eq:500107}).
\begin{equation}
\sum_{\nu_0 =0}^{\infty } \frac{s_{0,\infty }^{\nu _0}}{ \nu _0 !} \sum_{i_0=0}^{\nu_0} \frac{(-\nu_0)_{i_0} \left(\frac{\Omega }{2\mu}\right)_{i_0}}{(1)_{i_0}} \eta^{i_0} =  \exp(s_{0,\infty }) (1+s_{0,\infty } \eta)^{-\frac{\Omega }{2\mu}}
\label{eq:500108}
\end{equation}
Replace $a$ and $z$  by $\frac{\Omega }{2\mu}$ and $\widetilde{w}_{1,1}$ in (\ref{eq:500107}).
\begin{equation}
\sum_{\nu_0 =0}^{\infty } \frac{s_0^{\nu _0}}{ \nu _0 !} \sum_{i_0=0}^{\nu_0} \frac{(-\nu_0)_{i_0} \left(\frac{\Omega }{2\mu}\right)_{i_0}}{(1)_{i_0}} \widetilde{w}_{1,1}^{i_0} = \exp(s_0) (1+s_0 \widetilde{w}_{1,1})^{-\frac{\Omega }{2\mu}}
\label{eq:500109}
\end{equation}
Replace $a$ and $z$  by $\frac{\Omega }{2\mu}$ and $\widetilde{w}_{1,n}$ in (\ref{eq:500107}).
\begin{equation}
\sum_{\nu_0 =0}^{\infty } \frac{s_0^{\nu _0}}{ \nu _0 !} \sum_{i_0=0}^{\nu_0}\frac{(-\nu_0)_{i_0} \left(\frac{\Omega }{2\mu}\right)_{i_0}}{(1)_{i_0}} \widetilde{w}_{1,n}^{i_0} = \exp(s_0) (1+s_0 \widetilde{w}_{1,n})^{-\frac{\Omega }{2\mu}}
\label{eq:500110}
\end{equation}
Put $c_0$= 1 and $\lambda = \frac{\Omega }{\mu}$ in (\ref{eq:50092}). Substitute (\ref{eq:500108}), (\ref{eq:500109}) and (\ref{eq:500110}) into the new (\ref{eq:50092}).
\qed
\end{proof}
\section{Summary}
The canonical form of the biconfluent Heun equation is defined by \cite{Ronv1995,Maro1967}
\begin{equation}
x \frac{d^2{y}}{d{x}^2} + \left( 1+\alpha -\beta x-2 x^2 \right) \frac{d{y}}{d{x}} + \left( (\gamma -\alpha -2)x- \frac{1}{2}[\delta +(1+\alpha )\beta ]\right) y = 0
\label{eq:500111}
\end{equation}
in which $(\alpha ,\beta ,\gamma ,\delta ) \in \mathbb{C}^{4}$. This equation has two singular points which are a regular singularity at $x=0$ and an irregular singularity at $\infty$. This equation is derived from the GCH equation \cite{Chou2012a} by replacing all coefficients $\mu $, $\varepsilon $, $\nu $, $\Omega $ and $\omega $ by $-2$, $-\beta  $, $ 1+\alpha $, $\gamma -\alpha -2 $ and $ \frac{1}{2} \left( \frac{\delta}{\beta} +1+\alpha \right)$ in (\ref{eq:5001}).

In my previous two papers \cite{Chou2012i,Chou2012j} I show how to derive the power series expansion in closed forms, its integral forms (each sub-integral $y_m(x)$ where $m=0,1,2,\cdots$ is composed of $2m$ terms of definite integrals and $m$ terms of contour integrals) and the generating function for the GCH equation about $x=0$ by applying 3TRF. \cite{Chou2012b}     

In this chapter I show how to construct the power series expansion in closed forms and its integral forms of the GCH equation about $x=0$ for infinite series and polynomial of type 2 by applying R3TRF. This is done by letting $B_n$ in sequence $c_n$ is the leading term in the analytic function $y(x)$. For polynomial of type 2, I treat $\omega $ as a fixed value and $\mu$, $\varepsilon$, $\nu $, $\Omega$ as free variables. 

The power series expansion and the integral representation of the GCH equation for infinite series about $x=0$ in this chapter are equivalent to infinite series of the GCH equation in Ref.\cite{Chou2012i,Chou2012j}. In this chapter $B_n$ is the leading term in sequence $c_n$ in the analytic function $y(x)$. In Ref.\cite{Chou2012i,Chou2012j} $A_n$ is the leading term in sequence $c_n$ in the analytic function $y(x)$.
 
In Ref.\cite{Chou2012i,Chou2012j}, as we see the power series expansions of the GCH equation for infinite series and polynomial, the denominators and numerators in all $B_n$ terms of each sub-power series expansion $y_m(x)$ where $m=0,1,2,\cdots$ arise with Pochhammer symbol. And in this chapter the denominators and numerators in all $A_n$ terms of each sub-power series expansion $y_m(x)$ also arise with Pochhammer symbol. Since we construct the power series expansions with Pochhammer symbols in numerators and denominators, we are able to describe the integral representation of the GCH equation analytically. As we observe representations in closed form integrals of the GCH equation about $x=0$ by applying either 3TRF or R3TRF, a $_1F_1$ function (the Kummer function of the first kind) recurs in each of sub-integral forms of the GCH equation. It means that we are able to transform the GCH (or BCH) function about $x=0$ into any well-known special functions having two recursive coefficients in the power series of its ordinary differential equation because $_1F_1$ function arises in each of sub-integral forms on the GCH equation. After we replace $_1F_1$ functions in its integral forms to other special functions, we can rebuild the Frobenius solution of the GCH equation about $x=0$ in a backward.
 
In Ref.\cite{Chou2012j} and this chapter, I show how to derive the generating function for the type 1 and type 2 GCH polynomials from its analytic integral representation. And in this chapter I construct the generating function for the type 2 GCH polynomial from its integral representation. We are able to derive orthogonal relations, recursion relations and expectation values of the physical quantities from these two generating functions: the processes in order to obtain orthogonal and recursion relations of the GCH polynomial are similar as the case of a normalized wave function for the hydrogen-like atoms.\footnote{For instance, in the quantum mechanical aspects, if the eigenenergy is contained in $B_n$ term in a 3-term recursive relation between successive coefficients of the power series expansion, we have to apply the type 1 GCH polynomial. If the eigenenergy is included in $A_n$ term in a 3-term recursive relation, we should apply the type 2 GCH polynomial. If the first eigenenergy (mathematically, it is denoted by a spectral parameter) is included in $A_n$ term and the second one is involved in $B_n$ terms, we must apply the type 3 GCH polynomial. In my future paper I will discuss about the type 3 GCH polynomial.}  

In section 6.3 I construct the Frobenius solution of the GCH equation about $x=\infty $ for the type 1 polynomial by applying 3TRF analytically \cite{Chou2012b}. Also its integral representation and the generating function for the GCH equation are derived analytically. There are no such solutions for infinite series and the type 2 polynomial because the $B_n$ term is divergent in (\ref{eq:50060b}) since the index $n\gg 1$.\footnote{There are only two types of the analytic solution of the GCH equation about $x=\infty $ which are type 1 and type 3 polynomials. Its type 3 polynomial will be constructed in the future paper.}  
In comparison with integral forms of the GCH polynomials of type 1 and 2 about $x=0$, a Tricomi's function (Kummer's function of the second kind) recurs in each of sub-integral forms of the GCH polynomial of type 1 about $x=\infty $. 

\addcontentsline{toc}{section}{Bibliography}
\bibliographystyle{model1a-num-names}
\bibliography{<your-bib-database>}
\bibliographystyle{model1a-num-names}
\bibliography{<your-bib-database>}

\chapter{Mathieu function using reversible three-term recurrence formula}
\chaptermark{Mathieu function using R3TRF} 

In Ref.\cite{1Chou2012e}, by applying three term recurrence formula (3TRF)\cite{1chou2012b}, I construct the power series expansion in closed forms and its integral forms of Mathieu equation for infinite series including all higher terms of $A_n$'s.\footnote{`` higher terms of $A_n$'s'' means at least two terms of $A_n$'s.}  

In this chapter I will apply reversible three term recurrence formula (R3TRF) in chapter 1 to (1) the power series expansion in closed forms, (2) its integral forms of Mathieu equation (for infinite series and polynomial which makes $A_n$ term terminated including all higher terms of $B_n$'s\footnote{`` higher terms of $B_n$'s'' means at least two terms of $B_n$'s.})  and (3) the generating function for Mathieu polynomial which makes $A_n$ term terminated.

\section{Introduction}

Mathieu ordinary differential equation is of Fuchsian types with the two regular and one irregular singularities. In contrast, Heun equation of Fuchsian types has the four regular singularities. Heun equation has the four kind of confluent forms: (1) Confluent Heun (two regular and one irregular singularities), (2) Doubly confluent Heun (two irregular singularities), (3) Biconfluent Heun (one regular and one irregular singularities), (4) Triconfluent Heun equations (one irregular singularity).  For DLFM version \cite{NIST}, Mathieu equation in algebraic forms is also derived from the Confluent Heun equation by changing all coefficients $\delta =\gamma =\frac{1}{2}$, $\epsilon =0$, $\alpha =q$ and $q=\frac{\lambda +2q}{4}$.

The simplest example of three therm recursion relations in the power series of a linear ordinary differential equation is the Mathieu equation, introduced by Mathieu (1868)\cite{Math1868}, while investigating the vibration of an elliptical drumhead: it is derived from the Helmholtz equation in elliptic cylinder coordinates by using the method of separation of variables. Mathieu equation, known for the elliptic cylinder equation appears in diverse areas such as astronomy and physical problems involving Schr$\ddot{\mbox{o}}$dinger equation for a periodic potentials \cite{Conn1984}, the parametric Resonance in the reheating process of universe \cite{Zlat1998}, and wave equations in general relativity\cite{Birk2007}, etc. Mathieu function has been used in various areas in modern physics and mathematics.\cite{McLa1947,Guti2003,Daym1955,Troe1973,Alha1995,Shen1981,Bhat1988,Ragh1991,Aliev1999,Hort2007,Newm1962} 

Unfortunately, even though Mathieu equation has been observed in various areas mentioned above, the Mathieu function could not be described in the form of a definite or contour integral of any elementary functions analytically. Besides the analytic solutions of its power series expansion in closed forms have not be constructed yet for almost 150 years because of a 3-term recursive relation between successive coefficients in the Frobenius solution. Three recursive relation in its power series creates the mathematical complex for the calculation. Instead, the analytic solutions for any linear ordinary differential equation having two term recurrence relation in its power series expansion can be derived.
Mathieu function has only been described in numerical approximations (Whittaker 1914\cite{Whit1914}, Frenkel and Portugal 2001\cite{Fren2001}). Sips 1949\cite{Sips1949}, Frenkel and Portugal 2001\cite{Fren2001} argued that it is not possible to represent analytically the Mathieu function in a simple and handy way. All these authors treat solutions of Mathieu equation as the periodic-function with period $2\pi$. It is not possible for arbitrary parameters $\lambda $ and $q$ in (\ref{eq:6001}). They leave the fundamental solution of Mathieu equation by applying Floquet's theorem as solutions of recurrences because of a 3-term recursive relation in its power series expansion. In place of using the Bloch's theorem I will derive the analytic solution of Mathieu equation for any arbitrary $\lambda $ and $q$ by applying R3TRF in this chapter.  

In Ref.\cite{1Chou2012e} I construct analytic solutions of Mathieu equation about the regular singular point at $x=0$ by applying 3TRF\cite{1chou2012b}; for the power series expansion in closed forms and its integral representation for infinite series including all higher terms of $A_n$'s. 

In this chapter, by applying R3TRF in chapter 1, I construct the power series expansion in closed forms of Mathieu equation about the regular singular point at $x=0$ (for infinite series and polynomial which makes $A_n$ term terminated including all higher terms of $B_n$'s) analytically. The integral representation of Mathieu equation and the generating function for Mathieu polynomial which makes $A_n$ term terminated are derived in a mathematical rigor. Also Frobenius solution of the Mathieu equation about regular singular point at $x=1$ and irregular singular one at $x=\infty $ by applying R3TRF is constructed including its integral forms and the generating function of it.

Mathieu equation is 
\begin{equation}
 \frac{d^2{y}}{d{z}^2} + \left( \lambda - 2q\;\mathrm{\cos}2z \right) y = 0 \label{eq:6001}
\end{equation}
where $\lambda $ and $q$ are parameters. This is an equation with periodic-function coefficient. Mathieu equation also can be described in algebraic forms putting $x=\mathrm{\cos}^2z$:
\begin{equation}
 4x (1-x ) \frac{d^2{y}}{d{x}^2} + 2( 1-2x) \frac{d{y}}{d{x}} + ( \lambda + 2 q - 4 q x ) y = 0\label{eq:6002}
\end{equation}
This equation has two regular singularities: $x=0$ and $x=1$; the other singularity $x=\infty $ is irregular. Assume that its solution is
\begin{equation}
y(x)= \sum_{n=0}^{\infty } c_n x^{n+\nu }\label{eq:6003}
\end{equation}
where $\nu$ is an indicial root. Plug (\ref{eq:6003}) into (\ref{eq:6002}).
\begin{equation}
c_{n+1}=A_n \;c_n +B_n \;c_{n-1} \hspace{1cm};n\geq 1\label{eq:6004}
\end{equation}
where,
\begin{subequations}
\begin{equation}
A_n = \frac{4(n+\nu )^2-(\lambda +2q)}{2(n+1+\nu )(2(n+\nu )+1)} = \frac{(n+\nu-\varphi )(n+\nu +\varphi )}{ (n+1+\nu )\left( n+\nu +\frac{1}{2}\right)}\label{eq:6005a}
\end{equation}
and 
\begin{equation}
\varphi  = \frac{\sqrt{\lambda + 2q}}{2}\nonumber
\end{equation}
\begin{equation}
B_n = \frac{q}{ (n+1+\nu )\left( n+\nu +\frac{1}{2}\right)}\label{eq:6005b}
\end{equation}
\begin{equation}
c_1= A_0 \;c_0 \label{eq:6005c}
\end{equation}
\end{subequations}
We have two indicial roots which are $\nu  = 0$ and $ \frac{1}{2} $. As we see (\ref{eq:6005b}), there are no way to build the analytic solution of Mathieu equation for polynomial which makes $B_n$ term terminated at certain value of $n$. Because the numerator of (\ref{eq:6005b}) is just consist of constant $q$ parameter.\footnote{Whenever index $n$ increase in (\ref{eq:6005b}), $B_n$ term never be terminated with a fixed constant parameter $q$.}

Now let's test for convergence of the analytic function $y(x)$. As $n\gg 1$ (for sufficiently large), (\ref{eq:6005a}) and (\ref{eq:6005b}) are
\begin{subequations}
\begin{equation}
 \lim_{n\gg 1} A_n = 1 
 \label{eq:6006a}
\end{equation}
\begin{equation}
 \lim_{n\gg 1} B_n = \frac{q}{n^2}
 \label{eq:6006b}
\end{equation}
\end{subequations}

As $n\gg 1$, (\ref{eq:6006b}) is negligible. Put (\ref{eq:6006a}) with $B_n =0$ into (\ref{eq:6004}). 
\begin{equation}
\begin{tabular}{ l }
  \vspace{2 mm}
  $c_0$ \\
  \vspace{2 mm}
  $c_1 = c_0 $ \\
  \vspace{2 mm}
  $c_2 = c_0 $ \\
  \vspace{2 mm}
  $c_3 = c_0 $ \\
  \vspace{2 mm}                       
 \; \vdots \hspace{5mm} \vdots \\
\end{tabular}\label{eq:6007}
\end{equation}
When a function $y(x)$, analytic at $x=0$, is expanded in a power series putting $c_0=1$ by using (\ref{eq:6007}), we write
\begin{equation}
\lim_{n\gg 1}y(x) \approx \sum_{n=0}^{\infty } x^n = \frac{1}{1-x} \hspace{1cm}\mbox{where}\;0\leq x=\mathrm{\cos}^2z < 1
\label{eq:6008}
\end{equation}
For being convergent of $y(x)$ in (\ref{eq:6008}), an independent variable $x=\cos^2z$ should be less than 1.

The function $y(x)$ should be convergent as $x=1$ in terms of the application of the mathematical physics. As $A_n$ term is terminated as certain value of index $n$, the function $y(x)$ will be convergent even if $x=1$. (\ref{eq:6008}) is not applicable for the asymptotic behavior of $y(x)$ any more from the quantum mechanical point of view. Put (\ref{eq:6006b}) into (\ref{eq:6004}) with $A_n=0$.\footnote{Actually $A_n $ can not be zero for $n\gg 1$. However, I only interesting in the minimum value of $y(x)$ for asymptotic behavior.} For $n=0,1,2,\cdots$, it give
\begin{equation}
\begin{tabular}{  l  l }
  \vspace{2 mm}
   $c_2 = \frac{q}{1^2 } c_0 $  &\hspace{1cm}  $c_3 = \frac{q}{2^2 } c_1$  \\
  \vspace{2 mm}
  $c_4 = \frac{q}{1^2\cdot 3^2} c_0 $ &\hspace{1cm}  $c_5 =  \frac{q}{2^2\cdot 4^2 } c_1 $\\
  \vspace{2 mm}
  $c_6 =  \frac{q}{1^2\cdot 3^2\cdot 5^2 } c_0 $ &\hspace{1cm}  $c_7 = \frac{q}{2^2\cdot 4^2\cdot 6^2 } c_1 $\\
 \hspace{2 mm} \large{\vdots} & \hspace{1.5 cm}\large{\vdots} \\
 \vspace{2 mm}
  $c_{2n} = \left( \frac{ \Gamma \left( \frac{1}{2}\right)}{\Gamma \left( n+ \frac{1}{2}\right) } \right)^2 \left( \frac{ q}{4 } \right)^n c_0 $ &\hspace{1cm}  $c_{2n+1} = \frac{1}{(n!)^2} \left( \frac{ q}{4 } \right)^n c_1 $\\
\end{tabular}
\label{eq:6009}
\end{equation}
Put (\ref{eq:6009}) into the power series expansion where $\displaystyle {\sum_{n=0}^{\infty }c_n x^n}$ putting $c_0=1$ for simplicity and letting $c_1=A_o c_0 =0$ for the minimum value of $y(x)$. 
\begin{equation}
\mbox{min}\lim_{n\gg 1}y(x)=\; _1F_2 \left(1; \frac{1}{2},\frac{1}{2};  \frac{q}{4}x^2 \right) \label{eq:60010}
\end{equation}
On the above $_1F_2 (a;b,c;z) $  is a generalized hypergeometric function which is given by
\begin{equation}
_1F_2 (a;b,c;z) = \sum_{n=0}^{\infty } \frac{(a)_n}{\left(b\right)_n \left(c\right)_n} \frac{z^n}{n!} \nonumber
\end{equation} 
\section{Mathieu equation about regular singular point at zero}
\subsection{Power series}
\subsubsection{Polynomial of type 2}
There are two types of power series expansions using the two term recurrence relation in a linear ordinary differential equation which are a polynomial and an infinite series. In contrast there are three types of polynomials and an infinite series in three term recurrence relation of it: (1) polynomial which makes $B_n$ term terminated: $A_n$ term is not terminated, (2) polynomial which makes $A_n$ term terminated: $B_n$ term is not terminated, (3) polynomial which makes $A_n$ and $B_n$ terms terminated at the same time.\footnote{If $A_n$ and $B_n$ terms are not terminated, it turns to be infinite series.} In three term recurrence relation, polynomial of type 3 I categorize as complete polynomial. The solutions for the type 1 and 3 polynomials of Mathieu equation in the algebraic form does not exist because of the numerator in $B_n$ term. There are only two type solutions of  Mathieu equation which are type 2 polynomial and infinite series. For type 2 Mathieu polynomial about $x=0$, I treat a parameter $q$ as a free variable and $\lambda $ as a fixed value.

In chapter 1 the general expression of power series of $y(x)$ for polynomial of type 2 is defined by
\begin{eqnarray}
y(x) &=& \sum_{n=0}^{\infty } y_{n}(x) = y_0(x)+ y_1(x)+ y_2(x)+y_3(x)+\cdots \nonumber\\
&=&  c_0 \Bigg\{ \sum_{i_0=0}^{\alpha _0} \left( \prod _{i_1=0}^{i_0-1}A_{i_1} \right) x^{i_0+\lambda }
+ \sum_{i_0=0}^{\alpha _0}\left\{ B_{i_0+1} \prod _{i_1=0}^{i_0-1}A_{i_1}  \sum_{i_2=i_0}^{\alpha _1} \left( \prod _{i_3=i_0}^{i_2-1}A_{i_3+2} \right)\right\} x^{i_2+2+\lambda }  \nonumber\\
&& + \sum_{N=2}^{\infty } \Bigg\{ \sum_{i_0=0}^{\alpha _0} \Bigg\{B_{i_0+1}\prod _{i_1=0}^{i_0-1} A_{i_1} 
\prod _{k=1}^{N-1} \Bigg( \sum_{i_{2k}= i_{2(k-1)}}^{\alpha _k} B_{i_{2k}+2k+1}\prod _{i_{2k+1}=i_{2(k-1)}}^{i_{2k}-1}A_{i_{2k+1}+2k}\Bigg)\nonumber\\
&& \times  \sum_{i_{2N} = i_{2(N-1)}}^{\alpha _N} \Bigg( \prod _{i_{2N+1}=i_{2(N-1)}}^{i_{2N}-1} A_{i_{2N+1}+2N} \Bigg) \Bigg\} \Bigg\} x^{i_{2N}+2N+\lambda }\Bigg\}  \label{eq:60011}
\end{eqnarray}
In the above, $\alpha _i\leq \alpha _j$ only if $i\leq j$ where $i,j,\alpha _i, \alpha _j \in \mathbb{N}_{0}$.

For a polynomial, we need a condition which is:
\begin{equation}
A_{\alpha _i+ 2i}=0 \hspace{1cm} \mathrm{where}\;i,\alpha _i =0,1,2,\cdots
\label{eq:60012}
\end{equation}
In the above, $ \alpha _i$ is an eigenvalue that makes $A_n$ term terminated at certain value of index $n$. (\ref{eq:60012}) makes each $y_i(x)$ where $i=0,1,2,\cdots$ as the polynomial in (\ref{eq:60011}).

\paragraph{The case of  $ \varphi =\frac{1}{2}\sqrt{\lambda + 2q}= \lambda _i +2i +\nu$}
In (\ref{eq:6005a})-(\ref{eq:6005c}) replace $\varphi =\frac{1}{2}\sqrt{\lambda + 2q}$ by ${ \displaystyle \lambda _i +2i +\nu}$. In (\ref{eq:60012}) replace index $\alpha _i$ by $\lambda _i$. Take the new (\ref{eq:6005a})-(\ref{eq:6005c}), (\ref{eq:60012}) and put them in (\ref{eq:60011}). After the replacement process, the general expression of power series of Mathieu equation for polynomial of type 2 is given by
\begin{eqnarray}
 y(x)&=&  \sum_{n=0}^{\infty } y_{n}(x) = y_0(x)+ y_1(x)+ y_2(x)+y_3(x)+\cdots \nonumber\\ 
&=& c_0 x^{\nu } \left\{\sum_{i_0=0}^{\lambda _0} \frac{(-\lambda_0)_{i_0} \left( \lambda _0+ 2\nu \right)_{i_0}}{(1+\nu )_{i_0}\left( \frac{1}{2}+ \nu \right)_{i_0}} x^{i_0}\right.\nonumber\\
&&+ \left\{ \sum_{i_0=0}^{\lambda_0}\frac{1}{ (i_0+ 2+\nu )\left( i_0+\frac{3}{2} + \nu \right)}\frac{(-\lambda _0)_{i_0} \left( \lambda_0+  2\nu \right)_{i_0}}{(1+\nu )_{i_0}\left(\frac{1}{2}+\nu \right)_{i_0}} \right.\nonumber\\
&&\times \left. \sum_{i_1=i_0}^{\lambda _1} \frac{(-\lambda_1)_{i_1}\left(\lambda_1 + 4+2\nu \right)_{i_1}(3+\nu )_{i_0}\left(\frac{5}{2} +\nu \right)_{i_0}}{(-\lambda_1)_{i_0}\left(\lambda_1 + 4+2\nu \right)_{i_0}(3+\nu )_{i_1}\left(\frac{5}{2} +\nu \right)_{i_1}} x^{i_1}\right\}\eta \nonumber\\
&&+ \sum_{n=2}^{\infty } \left\{ \sum_{i_0=0}^{\lambda_0}\frac{1}{ (i_0+ 2+\nu )\left( i_0+\frac{3}{2} + \nu \right)}\frac{(-\lambda _0)_{i_0} \left( \lambda_0+  2\nu \right)_{i_0}}{(1+\nu )_{i_0}\left(\frac{1}{2}+\nu \right)_{i_0}}\right.\nonumber\\
&&\times \prod _{k=1}^{n-1} \left\{ \sum_{i_k=i_{k-1}}^{\lambda _k} \frac{1}{(i_k+ 2k+2+\nu )\left( i_k+ 2k+\frac{3}{2} +\nu \right)}\right.\nonumber\\
&&\times \left.\frac{(-\lambda _k)_{i_k}\left( \lambda_k +4k +2\nu \right)_{i_k}(2k+1+\nu )_{i_{k-1}}\left( 2k+\frac{1}{2} +\nu \right)_{i_{k-1}}}{(-\lambda _k)_{i_{k-1}}\left( \lambda_k +4k +2\nu \right)_{i_{k-1}}(2k+1+\nu )_{i_k}\left( 2k+\frac{1}{2} +\nu \right)_{i_k}}\right\} \nonumber\\
&&\times \left. \left.\sum_{i_n= i_{n-1}}^{\lambda _n} \frac{(-\lambda _n)_{i_n}\left( \lambda_n +4n +2\nu \right)_{i_n}(2n+1+\nu )_{i_{n-1}}\left( 2n+\frac{1}{2} +\nu \right)_{i_{n-1}}}{(-\lambda _n)_{i_{n-1}}\left( \lambda_n +4n +2\nu \right)_{i_{n-1}}(2n+1+\nu )_{i_n}\left( 2n+\frac{1}{2} +\nu \right)_{i_n}} x^{i_n} \right\} \eta ^n \right\} \hspace{2cm}\label{eq:60013}
\end{eqnarray}
where
\begin{equation}
\begin{cases} x= \mathrm{\cos}^2z \cr
\eta  = q x^2 \cr
\lambda = 2^2\left(\lambda_j +2j+\nu \right)^2- 2 q \;\;\mbox{as}\;j,\lambda _j\in \mathbb{N}_{0} \cr
\lambda_i\leq \lambda_j \;\;\mbox{only}\;\mbox{if}\;i\leq j\;\;\mbox{where}\;i,j\in \mathbb{N}_{0} 
\end{cases}\nonumber 
\end{equation}
\paragraph{The case of  $ \varphi =\frac{1}{2}\sqrt{\lambda + 2q}= -(\lambda _i +2i +\nu)$}
In (\ref{eq:6005a})-(\ref{eq:6005c}) replace $\frac{1}{2}\sqrt{\lambda + 2q}$ by $-(\lambda _i +2i +\nu)$. In (\ref{eq:60012}) replace index $\alpha _i$ by $\lambda _i$. Take the new (\ref{eq:6005a})-(\ref{eq:6005c}), (\ref{eq:60012}) and put them in (\ref{eq:60011}). After the replacement process, its solution is equivalent to (\ref{eq:60013}).

Put $c_0$= 1 as $\nu =0$ for the first kind of independent solutions of Mathieu equation and $\nu = 1/2 $ for the second one in (\ref{eq:60013}).  
\begin{remark}
The power series expansion of Mathieu equation of the first kind for polynomial of type 2 about $x=0$ as $\lambda = 2^2(\lambda_j +2j )^2- 2 q$ where $j,\lambda_j \in \mathbb{N}_{0}$ is
\begin{eqnarray}
 y(x)&=& MF_{\lambda _j}^R\left( q,\lambda = 2^2(\lambda_j +2j )^2- 2 q; \eta  = q x^2, x= \mathrm{\cos}^2z \right)\nonumber\\
&=& \sum_{i_0=0}^{\lambda _0} \frac{(-\lambda_0)_{i_0} \left( \lambda _0 \right)_{i_0}}{(1)_{i_0}\left( \frac{1}{2}\right)_{i_0}} x^{i_0} \nonumber\\
&+& \left\{ \sum_{i_0=0}^{\lambda_0}\frac{1}{ (i_0+ 2)\left( i_0+\frac{3}{2} \right)}\frac{(-\lambda _0)_{i_0} \left( \lambda_0 \right)_{i_0}}{(1)_{i_0}\left(\frac{1}{2}\right)_{i_0}} \right. \left. \sum_{i_1=i_0}^{\lambda _1} \frac{(-\lambda_1)_{i_1}\left(\lambda_1 + 4 \right)_{i_1}(3)_{i_0}\left(\frac{5}{2}\right)_{i_0}}{(-\lambda_1)_{i_0}\left(\lambda_1 + 4\right)_{i_0}(3)_{i_1}\left(\frac{5}{2}\right)_{i_1}} x^{i_1}\right\}\eta \nonumber\\
&+& \sum_{n=2}^{\infty } \left\{ \sum_{i_0=0}^{\lambda_0}\frac{1}{ (i_0+ 2)\left( i_0+\frac{3}{2}\right)}\frac{(-\lambda _0)_{i_0} \left( \lambda_0 \right)_{i_0}}{(1)_{i_0}\left(\frac{1}{2} \right)_{i_0}}\right.\nonumber\\
&\times& \prod _{k=1}^{n-1} \left\{ \sum_{i_k=i_{k-1}}^{\lambda _k} \frac{1}{(i_k+ 2k+2 )\left( i_k+ 2k+\frac{3}{2}\right)}\right.   \left.\frac{(-\lambda _k)_{i_k}\left( \lambda_k +4k \right)_{i_k}(2k+1)_{i_{k-1}}\left( 2k+\frac{1}{2} \right)_{i_{k-1}}}{(-\lambda _k)_{i_{k-1}}\left( \lambda_k +4k \right)_{i_{k-1}}(2k+1)_{i_k}\left( 2k+\frac{1}{2} \right)_{i_k}}\right\} \nonumber\\
&\times&  \left.\sum_{i_n= i_{n-1}}^{\lambda _n} \frac{(-\lambda _n)_{i_n}\left( \lambda_n +4n \right)_{i_n}(2n+1 )_{i_{n-1}}\left( 2n+\frac{1}{2} \right)_{i_{n-1}}}{(-\lambda _n)_{i_{n-1}}\left( \lambda_n +4n \right)_{i_{n-1}}(2n+1)_{i_n}\left( 2n+\frac{1}{2}\right)_{i_n}} x^{i_n} \right\} \eta ^n  \label{eq:60014}
\end{eqnarray}
\end{remark}
For the minimum value of Mathieu equation of the first kind for polynomial of type 2 about $x=0$, put $\lambda _0=\lambda _1=\lambda _2=\cdots=0$ in (\ref{eq:60014}).
\begin{eqnarray}
y(x)&=& MF_{0}^R\left( q,\lambda = 2(8j^2- q); \eta  = q x^2, x= \mathrm{\cos}^2z \right)\nonumber\\
&=&  \sum_{n=0}^{\infty } \frac{1}{\left( \frac{3}{4}\right)_n}\frac{\left( \frac{q}{4}x^2\right)^n}{n!}
= \Gamma \left( 3/4\right) \left( \frac{i}{2}\sqrt{q}x\right)^{\frac{1}{4}}J_{-1/4}\left( i\sqrt{q}x\right)\label{jjj:6001}
\end{eqnarray} 
On the above, $J_{\alpha }(x)$ is Bessel functions of the first kind.
\begin{remark}
The power series expansion of Mathieu equation of the second kind for polynomial of type 2 about $x=0$ as $\lambda = 2^2\left(\lambda_j +2j+1/2 \right)^2- 2 q $ where $j,\lambda_j \in \mathbb{N}_{0}$ is
\begin{eqnarray}
 y(x)&=& MS_{\lambda _j}^R\left( q,\lambda =2^2\left(\lambda_j +2j+1/2 \right)^2- 2 q; \eta  = q x^2, x= \mathrm{\cos}^2z \right)\nonumber\\
&=& x^{\frac{1}{2}} \left\{\sum_{i_0=0}^{\lambda _0} \frac{(-\lambda_0)_{i_0} \left( \lambda _0+ 1 \right)_{i_0}}{\left(\frac{3}{2} \right)_{i_0}\left( 1\right)_{i_0}} x^{i_0}\right.\nonumber\\
&+& \left\{ \sum_{i_0=0}^{\lambda_0}\frac{1}{ \left(i_0+ \frac{5}{2} \right)\left( i_0 +2 \right)}\frac{(-\lambda _0)_{i_0} \left( \lambda_0+ 1 \right)_{i_0}}{\left(\frac{3}{2} \right)_{i_0}\left(1 \right)_{i_0}} \right. \left. \sum_{i_1=i_0}^{\lambda _1} \frac{(-\lambda_1)_{i_1}\left(\lambda_1 + 5 \right)_{i_1}\left(\frac{7}{2} \right)_{i_0}\left( 3\right)_{i_0}}{(-\lambda_1)_{i_0}\left(\lambda_1 +5 \right)_{i_0}\left(\frac{7}{2} \right)_{i_1}\left( 3 \right)_{i_1}} x^{i_1}\right\}\eta \nonumber\\
&+& \sum_{n=2}^{\infty } \left\{ \sum_{i_0=0}^{\lambda_0}\frac{1}{ \left(i_0+ \frac{5}{2} \right)\left( i_0+2 \right)}\frac{(-\lambda _0)_{i_0} \left( \lambda_0+  1 \right)_{i_0}}{\left(\frac{3}{2} \right)_{i_0}\left( 1 \right)_{i_0}}\right.\nonumber\\
&\times& \prod _{k=1}^{n-1} \left\{ \sum_{i_k=i_{k-1}}^{\lambda _k} \frac{1}{\left( i_k+ 2k+\frac{5}{2} \right)\left( i_k+ 2k+2 \right)}\right.   \left.\frac{(-\lambda _k)_{i_k}\left( \lambda_k +4k +1\right)_{i_k}\left( 2k+\frac{3}{2} \right)_{i_{k-1}}\left( 2k+1 \right)_{i_{k-1}}}{(-\lambda _k)_{i_{k-1}}\left( \lambda_k +4k +1 \right)_{i_{k-1}}\left( 2k+\frac{3}{2} \right)_{i_k}\left( 2k+1 \right)_{i_k}}\right\} \nonumber\\
&\times& \left. \left.\sum_{i_n= i_{n-1}}^{\lambda _n} \frac{(-\lambda _n)_{i_n}\left( \lambda_n +4n +1 \right)_{i_n}\left(2n+\frac{3}{2} \right)_{i_{n-1}}\left( 2n+1 \right)_{i_{n-1}}}{(-\lambda _n)_{i_{n-1}}\left( \lambda_n +4n +1 \right)_{i_{n-1}}\left( 2n+\frac{3}{2} \right)_{i_n}\left( 2n+1 \right)_{i_n}} x^{i_n} \right\} \eta ^n \right\} \label{eq:60015}
\end{eqnarray}
\end{remark}
For the minimum value of Mathieu equation of the first kind for polynomial of type 2 about $x=0$, put $\lambda _0=\lambda _1=\lambda _2=\cdots=0$ in (\ref{eq:60015}).
\begin{eqnarray}
y(x)&=& MS_{0}^R\left( q,\lambda =4\left( 2j+1/2 \right)^2- 2 q; \eta  = q x^2, x= \mathrm{\cos}^2z \right)\nonumber\\
&=&  x^{\frac{1}{2}}\sum_{n=0}^{\infty } \frac{1}{\left( \frac{5}{4}\right)_n}\frac{\left( \frac{q}{4}x^2\right)^n}{n!}
= \Gamma \left( 5/4\right) x^{\frac{1}{2}} \left( \frac{i}{2}\sqrt{q}x\right)^{-\frac{1}{4}}J_{1/4}\left( i\sqrt{q}x\right)\label{jjj:6002}
\end{eqnarray}
\subsubsection{Infinite series}
In chapter 1 the general expression of power series of $y(x)$ for infinite series is
\begin{eqnarray}
y(x) &=&  \sum_{n=0}^{\infty } y_{n}(x) = y_0(x)+ y_1(x)+ y_2(x)+y_3(x)+\cdots \nonumber\\
&=& c_0 \Bigg\{ \sum_{i_0=0}^{\infty } \left( \prod _{i_1=0}^{i_0-1}A_{i_1} \right) x^{i_0+\lambda }
+ \sum_{i_0=0}^{\infty }\left\{ B_{i_0+1} \prod _{i_1=0}^{i_0-1}A_{i_1}  \sum_{i_2=i_0}^{\infty } \left( \prod _{i_3=i_0}^{i_2-1}A_{i_3+2} \right)\right\} x^{i_2+2+\lambda }  \nonumber\\
&& + \sum_{N=2}^{\infty } \Bigg\{ \sum_{i_0=0}^{\infty } \Bigg\{B_{i_0+1}\prod _{i_1=0}^{i_0-1} A_{i_1} 
\prod _{k=1}^{N-1} \Bigg( \sum_{i_{2k}= i_{2(k-1)}}^{\infty } B_{i_{2k}+2k+1}\prod _{i_{2k+1}=i_{2(k-1)}}^{i_{2k}-1}A_{i_{2k+1}+2k}\Bigg)\nonumber\\
&& \times  \sum_{i_{2N} = i_{2(N-1)}}^{\infty } \Bigg( \prod _{i_{2N+1}=i_{2(N-1)}}^{i_{2N}-1} A_{i_{2N+1}+2N} \Bigg) \Bigg\} \Bigg\} x^{i_{2N}+2N+\lambda }\Bigg\}   \label{eq:60016}
\end{eqnarray}
Substitute (\ref{eq:6005a})--(\ref{eq:6005c}) into (\ref{eq:60016}). 
The general expression of power series of Mathieu equation for infinite series about $x=0$ is given by
\begin{eqnarray}
 y(x)&=&\sum_{n=0}^{\infty } y_n(x)= y_0(x)+ y_1(x)+ y_2(x)+ y_3(x)+\cdots \nonumber\\
&=& c_0 x^{\nu } \left\{\sum_{i_0=0}^{\infty } \frac{\left(\nu - \varphi \right)_{i_0} \left(\nu + \varphi \right)_{i_0}}{\left(1+\nu \right)_{i_0}\left(\frac{1}{2} +\nu \right)_{i_0}} x^{i_0}\right.\nonumber\\
&+& \left\{ \sum_{i_0=0}^{\infty }\frac{1}{\left( i_0+ 2+\nu\right) \left( i_0+ \frac{3}{2} +\nu \right)}\frac{\left(\nu-\varphi \right)_{i_0} \left( \nu +\varphi \right)_{i_0}}{(1+\nu )_{i_0} \left( \frac{1}{2} +\nu \right)_{i_0}}\right.\nonumber\\
&\times& \left. \sum_{i_1=i_0}^{\infty } \frac{\left( \nu +2-\varphi \right)_{i_1} \left(\nu +2+\varphi \right)_{i_1}\left( 3+\nu \right)_{i_0}\left( \frac{5}{2} +\nu \right)_{i_0}}{\left( \nu +2-\varphi \right)_{i_0} \left(\nu +2+\varphi \right)_{i_0}\left( 3+\nu \right)_{i_1}\left( \frac{5}{2} +\nu \right)_{i_1}}x^{i_1}\right\} \eta \nonumber\\
&+& \sum_{n=2}^{\infty } \left\{ \sum_{i_0=0}^{\infty }\frac{1}{\left( i_0+ 2+\nu\right) \left( i_0+ \frac{3}{2} +\nu \right)}\frac{\left(\nu-\varphi \right)_{i_0} \left( \nu +\varphi \right)_{i_0}}{(1+\nu )_{i_0} \left( \frac{1}{2} +\nu \right)_{i_0}}\right.\nonumber\\
&\times& \prod _{k=1}^{n-1} \left\{ \sum_{i_k=i_{k-1}}^{\infty } \frac{1}{\left( i_k+ 2k+2+\nu \right) \left( i_k+ 2k+\frac{3}{2} +\nu \right)}\right.\nonumber\\
&\times& \left. \frac{ \left(\nu +2k -\varphi \right)_{i_k} \left(\nu +2k +\varphi \right)_{i_k} \left( 2k+1+\nu \right)_{i_{k-1}}\left( 2k+\frac{1}{2} +\nu \right)_{i_{k-1}}}{ \left(\nu +2k -\varphi \right)_{i_{k-1}} \left(\nu +2k +\varphi \right)_{i_{k-1}} \left( 2k+1+\nu \right)_{i_k}\left( 2k+\frac{1}{2} +\nu \right)_{i_k}}\right\} \nonumber\\
&\times& \left.\left.\sum_{i_n= i_{n-1}}^{\infty }\frac{ \left(\nu +2 n -\varphi \right)_{i_n} \left(\nu +2n +\varphi \right)_{i_n} \left( 2n+1+\nu \right)_{i_{n-1}}\left( 2n+\frac{1}{2} +\nu \right)_{i_{n-1}}}{ \left(\nu +2n -\varphi \right)_{i_{n-1}} \left(\nu +2n +\varphi \right)_{i_{n-1}} \left( 2n+1+\nu \right)_{i_n}\left( 2n+\frac{1}{2} +\nu \right)_{i_n}} x^{i_n} \right\} \eta ^n \right\} \hspace{2cm}\label{eq:60017}
\end{eqnarray}
where
\begin{equation}
\begin{cases} x= \mathrm{\cos}^2z \cr
\eta  = q x^2 \cr
\varphi  = \frac{\sqrt{\lambda + 2q}}{2} 
\end{cases}\nonumber 
\end{equation}
Put $c_0$= 1 as $\nu =0$ for the first kind of independent solutions of Mathieu equation and $\nu = 1/2 $ for the second one in (\ref{eq:60017}).  
\begin{remark}
The power series expansion of Mathieu equation of the first kind for infinite series about $x=0$ using R3TRF is
\begin{eqnarray}
 y(x)&=& MF^R\left( q,\lambda, \varphi  = \frac{1}{2} \sqrt{\lambda + 2q} ; \eta  = q x^2, x= \mathrm{\cos}^2z \right) \nonumber\\
&=& \sum_{i_0=0}^{\infty } \frac{\left( - \varphi \right)_{i_0} \left( \varphi \right)_{i_0}}{\left(1 \right)_{i_0}\left(\frac{1}{2} \right)_{i_0}} x^{i_0} \nonumber\\
&+& \left\{ \sum_{i_0=0}^{\infty }\frac{1}{\left( i_0+ 2 \right) \left( i_0+ \frac{3}{2} \right)}\frac{\left( -\varphi \right)_{i_0} \left(  \varphi \right)_{i_0}}{(1 )_{i_0} \left( \frac{1}{2} \right)_{i_0}} \sum_{i_1=i_0}^{\infty } \frac{\left( 2-\varphi \right)_{i_1} \left( 2+\varphi \right)_{i_1}\left( 3 \right)_{i_0}\left( \frac{5}{2} \right)_{i_0}}{\left( 2-\varphi \right)_{i_0} \left( 2+\varphi \right)_{i_0}\left( 3 \right)_{i_1}\left( \frac{5}{2} \right)_{i_1}}x^{i_1}\right\} \eta \nonumber\\
&+& \sum_{n=2}^{\infty } \left\{ \sum_{i_0=0}^{\infty }\frac{1}{\left( i_0+ 2 \right) \left( i_0+ \frac{3}{2} \right)}\frac{\left( -\varphi \right)_{i_0} \left( \varphi \right)_{i_0}}{(1 )_{i_0} \left( \frac{1}{2} \right)_{i_0}}\right.\nonumber\\
&\times& \prod _{k=1}^{n-1} \left\{ \sum_{i_k=i_{k-1}}^{\infty } \frac{1}{\left( i_k+ 2k+2 \right) \left( i_k+ 2k+\frac{3}{2} \right)} \frac{ \left( 2k -\varphi \right)_{i_k} \left( 2k +\varphi \right)_{i_k} \left( 2k+1 \right)_{i_{k-1}}\left( 2k+\frac{1}{2} \right)_{i_{k-1}}}{ \left( 2k -\varphi \right)_{i_{k-1}} \left( 2k +\varphi \right)_{i_{k-1}} \left( 2k+1 \right)_{i_k}\left( 2k+\frac{1}{2} \right)_{i_k}}\right\} \nonumber\\
&\times&  \left.\sum_{i_n= i_{n-1}}^{\infty }\frac{ \left( 2 n -\varphi \right)_{i_n} \left( 2n +\varphi \right)_{i_n} \left( 2n+1 \right)_{i_{n-1}}\left( 2n+\frac{1}{2} \right)_{i_{n-1}}}{ \left( 2n -\varphi \right)_{i_{n-1}} \left( 2n +\varphi \right)_{i_{n-1}} \left( 2n+1  \right)_{i_n}\left( 2n+\frac{1}{2} \right)_{i_n}} x^{i_n} \right\} \eta ^n  \label{eq:60018}
\end{eqnarray}
\end{remark}
\begin{remark}
The power series expansion of Mathieu equation of the second kind for infinite series about $x=0$ using R3TRF is
\begin{eqnarray}
y(x)&=& MS^R\left( q,\lambda, \varphi  = \frac{1}{2} \sqrt{\lambda + 2q} ; \eta  = q x^2, x= \mathrm{\cos}^2z \right) \nonumber\\ \nonumber\\
&=& x^{\frac{1}{2}} \left\{\sum_{i_0=0}^{\infty } \frac{\left(\frac{1}{2} - \varphi \right)_{i_0} \left(\frac{1}{2} + \varphi \right)_{i_0}}{\left(\frac{3}{2} \right)_{i_0}\left(1\right)_{i_0}} x^{i_0}\right.\nonumber\\
&+& \left\{ \sum_{i_0=0}^{\infty }\frac{1}{\left( i_0+ \frac{5}{2}\right) \left( i_0+ 2 \right)}\frac{\left(\frac{1}{2}-\varphi \right)_{i_0} \left( \frac{1}{2} +\varphi \right)_{i_0}}{(\frac{3}{2} )_{i_0} \left( 1 \right)_{i_0}} \sum_{i_1=i_0}^{\infty } \frac{\left( \frac{5}{2}-\varphi \right)_{i_1} \left(\frac{5}{2}+\varphi \right)_{i_1}\left( \frac{7}{2} \right)_{i_0}\left( 3 \right)_{i_0}}{\left( \frac{5}{2}-\varphi \right)_{i_0} \left(\frac{5}{2}+\varphi \right)_{i_0}\left( \frac{7}{2} \right)_{i_1}\left( 3 \right)_{i_1}}x^{i_1}\right\} \eta \nonumber\\
&+& \sum_{n=2}^{\infty } \left\{ \sum_{i_0=0}^{\infty }\frac{1}{\left( i_0+ \frac{5}{2}\right) \left( i_0+ 2 \right)}\frac{\left(\frac{1}{2}-\varphi \right)_{i_0} \left( \frac{1}{2} +\varphi \right)_{i_0}}{(\frac{3}{2} )_{i_0} \left( 1 \right)_{i_0}}\right.\nonumber\\
&\times& \prod _{k=1}^{n-1} \left\{ \sum_{i_k=i_{k-1}}^{\infty } \frac{1}{\left( i_k+ 2k+\frac{5}{2} \right) \left( i_k+ 2k+2 \right)} \frac{ \left( 2k+\frac{1}{2} -\varphi \right)_{i_k} \left( 2k+\frac{1}{2} +\varphi \right)_{i_k} \left( 2k+\frac{3}{2} \right)_{i_{k-1}}\left( 2k+1 \right)_{i_{k-1}}}{ \left( 2k+\frac{1}{2} -\varphi \right)_{i_{k-1}} \left( 2k+ \frac{1}{2} +\varphi \right)_{i_{k-1}} \left( 2k+\frac{3}{2} \right)_{i_k}\left( 2k+1 \right)_{i_k}}\right\} \nonumber\\
&\times& \left.\left.\sum_{i_n= i_{n-1}}^{\infty }\frac{ \left( 2n+\frac{1}{2} -\varphi \right)_{i_n} \left( 2n+\frac{1}{2} +\varphi \right)_{i_n} \left( 2n+\frac{3}{2} \right)_{i_{n-1}}\left( 2n+1 \right)_{i_{n-1}}}{ \left( 2n+\frac{1}{2} -\varphi \right)_{i_{n-1}} \left( 2n+\frac{1}{2} +\varphi \right)_{i_{n-1}} \left( 2n+\frac{3}{2} \right)_{i_n}\left( 2n+1 \right)_{i_n}} x^{i_n} \right\} \eta ^n \right\} \label{eq:60019}
\end{eqnarray}
\end{remark}
Two power series expansions of Mathieu equation for infinite series in this chapter and Ref.\cite{1Chou2012e} are equivalent to each other. In this chapter, $B_n$ is the leading term in sequence $c_n$ of the analytic function $y(x)$. In Ref.\cite{1Chou2012e}, $A_n$ is the leading term in sequence $c_n$ of the analytic function $y(x)$.

For the special case, as $|q| \ll 1$ in (\ref{eq:60018}) and (\ref{eq:60019}), we have
\begin{subequations}
\begin{eqnarray}
&&\lim_{|q| \ll 1} MF^R\left( q,\lambda, \varphi  = \frac{1}{2} \sqrt{\lambda + 2q} ; \eta  = q x^2, x= \mathrm{\cos}^2z \right) \nonumber\\
&&\approx  \sum_{i_0=0}^{\infty } \frac{\left( - \frac{\sqrt{\lambda }}{2} \right)_{i_0} \left( \frac{\sqrt{\lambda }}{2} \right)_{i_0}}{\left(1 \right)_{i_0}\left(\frac{1}{2} \right)_{i_0}} x^{i_0} \nonumber\\
&&+  \eta \sum_{i_0=0}^{\infty }\frac{1}{\left( i_0+ 2 \right) \left( i_0+ \frac{3}{2} \right)}\frac{\left( -\frac{\sqrt{\lambda }}{2} \right)_{i_0} \left(  \frac{\sqrt{\lambda }}{2} \right)_{i_0}}{(1 )_{i_0} \left( \frac{1}{2} \right)_{i_0}} \sum_{i_1=i_0}^{\infty } \frac{\left( 2-\frac{\sqrt{\lambda }}{2} \right)_{i_1} \left( 2+\frac{\sqrt{\lambda }}{2} \right)_{i_1}\left( 3 \right)_{i_0}\left( \frac{5}{2} \right)_{i_0}}{\left( 2-\frac{\sqrt{\lambda }}{2}\right)_{i_0} \left( 2+\frac{\sqrt{\lambda }}{2} \right)_{i_0}\left( 3 \right)_{i_1}\left( \frac{5}{2} \right)_{i_1}}x^{i_1} \nonumber\\
&&> \mathrm{\cos}\left( \sqrt{\lambda } \mathrm{\sin}^{-1}\left(\sqrt{x}\right)\right) \label{eq:60020a}
\end{eqnarray}
And,
\begin{eqnarray}
&&\lim_{|q| \ll 1}  MS^R\left( q,\lambda, \varphi  = \frac{1}{2} \sqrt{\lambda + 2q} ; \eta  = q x^2, x= \mathrm{\cos}^2z \right) \nonumber\\
&&\approx  x^{\frac{1}{2}} \left\{\sum_{i_0=0}^{\infty } \frac{\left(\frac{1}{2} - \frac{\sqrt{\lambda }}{2} \right)_{i_0} \left(\frac{1}{2} + \frac{\sqrt{\lambda }}{2} \right)_{i_0}}{\left(\frac{3}{2} \right)_{i_0}\left(1\right)_{i_0}} x^{i_0}\right.\nonumber\\
&&+ \left. \eta \sum_{i_0=0}^{\infty }\frac{1}{\left( i_0+ \frac{5}{2}\right) \left( i_0+ 2 \right)}\frac{\left(\frac{1}{2}-\frac{\sqrt{\lambda }}{2} \right)_{i_0} \left( \frac{1}{2} +\frac{\sqrt{\lambda }}{2} \right)_{i_0}}{(\frac{3}{2} )_{i_0} \left( 1 \right)_{i_0}} \sum_{i_1=i_0}^{\infty } \frac{\left( \frac{5}{2}-\frac{\sqrt{\lambda }}{2} \right)_{i_1} \left(\frac{5}{2}+\frac{\sqrt{\lambda }}{2} \right)_{i_1}\left( \frac{7}{2} \right)_{i_0}\left( 3 \right)_{i_0}}{\left( \frac{5}{2}-\frac{\sqrt{\lambda }}{2} \right)_{i_0} \left(\frac{5}{2}+\frac{\sqrt{\lambda }}{2} \right)_{i_0}\left( \frac{7}{2} \right)_{i_1}\left( 3 \right)_{i_1}}x^{i_1}\right\} \nonumber\\
&&> \sqrt{\frac{x}{1-x}}  \mathrm{\cos}\left( \sqrt{\lambda } \mathrm{\sin}^{-1}\left(\sqrt{x}\right)\right)  \label{eq:60020b}
\end{eqnarray}
\end{subequations}
\subsection{Integral formalism}
The Mathieu function could not be constructed in a definite or contour integral form of any well-known simple functions because of a 3-term recursive relation in its power series expansion. The three term recurrence relation in power series of Mathieu equation creates mathematical difficulty to be analyzed it into a direct or contour integrals.

In place of describing the Mathieu function into the integral representation of any simple functions which have two term recursion in its power series expansion of a linear ordinary differential equation, in earlier literature the integral equations of the Mathieu function were constructed by using trigonometric functions kernels, Bessel-function kernels and etc; such integral relationships express one analytic solution in terms of another analytic solution such as the Fourier series of the periodic Mathieu functions; one of $\mbox{ce}$, $\mbox{se}$, $\mbox{me}$, $\mbox{Ce}$ and $\mbox{Se}$, $\mbox{Me}$, $\mbox{Ge}$ and $\mbox{Ge}$ functions.\cite{Prud1990,Erde1955,Grad2000,Volk1983} There are many other forms of integral relations in Mathieu equation. \cite{Meix1980,Sips1970,Arsc1964,Wang1989}

In Ref.\cite{1Chou2012e}, I show integral representation (each sub-integral $y_m(x)$ where $m=0,1,2,\cdots$ is composed of $2m$ terms of definite integrals and $m$ terms of contour integrals) of the Mathieu function using 3TRF: a Modified Bessel function recurs in each of sub-integral forms of the Mathieu function.

Now I consider integral forms of the Mathieu function by using R3TRF.
\subsubsection{Polynomial of type 2}
Now I consider an integral form of the Mathieu polynomial of type 2 by using R3TRF. There is a generalized hypergeometric function which is written by
\begin{eqnarray}
I_l &=& \sum_{i_l= i_{l-1}}^{\lambda _l} \frac{(-\lambda _l)_{i_l}\left( \lambda _l+4l +2\nu \right)_{i_l}(2l+1+\nu )_{i_{l-1}}\left( 2l+\frac{1}{2} +\nu \right)_{i_{l-1}}}{(-\lambda _l)_{i_{l-1}}\left( \lambda _l+4l +2\nu \right)_{i_{l-1}}(2l+1+\nu )_{i_l}(2l+\frac{1}{2} +\nu )_{i_l}} x^{i_l}\label{eq:60021}\\
&=& x^{i_{l-1}} 
\sum_{j=0}^{\infty } \frac{B\left(i_{l-1}+2l-\frac{1}{2}+\nu ,j+1\right) B\left( i_{l-1}+2l +\nu ,j+1\right) (i_{l-1}-\lambda _l)_j \left( i_{l-1}+\lambda _l+4l +2\nu \right)_j}{(i_{l-1}+2l+\nu )^{-1}\left( i_{l-1}+2l-\frac{1}{2}+\nu \right)^{-1}(1)_j \;j!} x^j\nonumber
\end{eqnarray}
By using integral form of beta function,
\begin{subequations}
\begin{equation}
B\left( i_{l-1}+2l +\nu ,j+1\right) = \int_{0}^{1} dt_l\;t_l^{i_{l-1}+2l -1 +\nu } (1-t_l)^j \label{eq:60022a}
\end{equation}
\begin{equation}
B\left(i_{l-1}+2l-\frac{1}{2}+\nu ,j+1\right) = \int_{0}^{1} du_l\;u_l^{i_{l-1}+2l -\frac{3}{2} +\nu } (1-u_l)^j\label{eq:60022b}
\end{equation}
\end{subequations}
Substitute (\ref{eq:60022a}) and (\ref{eq:60022b}) into (\ref{eq:60021}). And divide $(i_{l-1}+2l+\nu ) \left( i_{l-1}+2l-\frac{1}{2}+\nu \right) $ into the new (\ref{eq:60021}).
\begin{eqnarray}
K_l&=& \frac{1}{(i_{l-1}+2l+\nu ) \left( i_{l-1}+2l-\frac{1}{2}+\nu \right)}\nonumber\\
&&\times \sum_{i_l= i_{l-1}}^{\lambda _l} \frac{(-\lambda _l)_{i_l}\left( \lambda _l+4l +2\nu \right)_{i_l}(2l+1+\nu )_{i_{l-1}}\left( 2l+\frac{1}{2} +\nu \right)_{i_{l-1}}}{(-\lambda _l)_{i_{l-1}}\left( \lambda _l+4l +2\nu \right)_{i_{l-1}}(2l+1+\nu )_{i_l}(2l+\frac{1}{2} +\nu )_{i_l}} x^{i_l}\nonumber\\
&=&  \int_{0}^{1} dt_l\;t_l^{2l-1+\nu } \int_{0}^{1} du_l\;u_l^{2l -\frac{3}{2} +\nu } (x t_l u_l)^{i_{l-1}}\nonumber\\
&&\times  \sum_{j=0}^{\infty } \frac{(i_{l-1}-\lambda _l)_j \left( i_{l-1}+\lambda _l+4l+ 2\nu \right)_j}{(1)_j \;j!} 
[x(1-t_l)(1-u_l)]^j \nonumber 
\end{eqnarray}
The contour integral form of Gauss hypergeometric function is given by
\begin{eqnarray}
_2F_1 \left( \alpha ,\beta ; \gamma ; z \right) &=& \sum_{n=0}^{\infty } \frac{(\alpha )_n (\beta )_n}{(\gamma )_n (n!)} z^n \nonumber\\
&=& -\frac{1}{2\pi i} \frac{\Gamma(1-\alpha ) \Gamma(\gamma )}{\Gamma (\gamma -\alpha )} \oint dp_l\;(-p_l)^{\alpha -1} (1-p_l)^{\gamma -\alpha -1} (1-z p_l)^{-\beta }\hspace{1cm}\label{eq:60023}\\
&& \mbox{where} \;\mbox{Re}(\gamma -\alpha )>0 \nonumber
\end{eqnarray}
replaced $\alpha $, $\beta $, $\gamma $ and z by $i_{l-1}-\lambda _l$, $ { \displaystyle i_{l-1}+\lambda _l+4l +2\nu }$, 1 and $x (1-t_l)(1-u_l)$ in (\ref{eq:60023})
\begin{eqnarray}
&& \sum_{j=0}^{\infty } \frac{\left(i_{l-1}-\lambda _l)_j (i_{l-1}+\lambda _l+4l +2\nu \right)_j}{(1)_j \;j!} [x(1-t_l)(1-u_l)]^j \nonumber\\
&=& \frac{1}{2\pi i} \oint dp_l\;\frac{1}{p_l} (1-x(1-t_l)(1-u_l)p_l)^{-\left( 4l +2\nu \right)} \nonumber\\
&&\times \left(\frac{p_l-1}{p_l}\frac{1}{1-x(1-t_l)(1-u_l)p_l}\right)^{\lambda _l} \left(\frac{p_l}{p_l-1} \frac{1}{1-x(1-t_l)(1-u_l)p_l}\right)^{i_{l-1}} \label{eq:60024}
\end{eqnarray}
Substitute (\ref{eq:60024}) into $K_l$.
\begin{eqnarray}
K_l&=& \frac{1}{(i_{l-1}+2l+\nu ) \left( i_{l-1}+2l-\frac{1}{2}+\nu \right)}\nonumber\\
&&\times \sum_{i_l= i_{l-1}}^{\lambda _l} \frac{(-\lambda _l)_{i_l}\left( \lambda _l+4l +2\nu \right)_{i_l}(2l+1+\nu )_{i_{l-1}}\left( 2l+\frac{1}{2} +\nu \right)_{i_{l-1}}}{(-\lambda _l)_{i_{l-1}}\left( \lambda _l+4l +2\nu \right)_{i_{l-1}}(2l+1+\nu )_{i_l}(2l+\frac{1}{2} +\nu )_{i_l}} x^{i_l}\nonumber\\
&=&  \int_{0}^{1} dt_l\;t_l^{2l-1+\nu } \int_{0}^{1} du_l\;u_l^{2l-\frac{3}{2} +\nu } 
\frac{1}{2\pi i} \oint dp_l\;\frac{1}{p_l} (1-x(1-t_l)(1-u_l)p_l)^{-\left( 4l+ 2\nu \right)}  \nonumber\\
&&\times \left(\frac{p_l-1}{p_l}\frac{1}{1-x(1-t_l)(1-u_l)p_l}\right)^{\lambda _l} \left(\frac{p_l}{p_l-1} \frac{xt_l u_l}{1-x(1-t_l)(1-u_l)p_l}\right)^{i_{l-1}}\label{eq:60025} 
\end{eqnarray}
Substitute (\ref{eq:60025}) into (\ref{eq:60013}) where $l=1,2,3,\cdots$; apply $K_1$ into the second summation of sub-power series $y_1(x)$, apply $K_2$ into the third summation and $K_1$ into the second summation of sub-power series $y_2(x)$, apply $K_3$ into the forth summation, $K_2$ into the third summation and $K_1$ into the second summation of sub-power series $y_3(x)$, etc.\footnote{$y_1(x)$ means the sub-power series in (\ref{eq:60013}) contains one term of $B_n's$, $y_2(x)$ means the sub-power series in (\ref{eq:60013}) contains two terms of $B_n's$, $y_3(x)$ means the sub-power series in (\ref{eq:60013}) contains three terms of $B_n's$, etc.}
\begin{theorem}
The general representation in the form of integral of the Mathieu polynomial of type 2 is given by
\begin{eqnarray}
 y(x)&=& \sum_{n=0}^{\infty } y_{n}(x)= y_0(x)+ y_1(x)+ y_2(x)+y_3(x)+\cdots \nonumber\\
&=& c_0 x^{\nu } \left\{ \sum_{i_0=0}^{\lambda _0}\frac{(-\lambda_0)_{i_0}\left( \lambda_0 +2\nu \right)_{i_0}}{\left(1+\nu \right)_{i_0}\left(\nu +\frac{1}{2}\right)_{i_0}} x^{i_0}\right.\nonumber\\
&&+ \sum_{n=1}^{\infty } \left\{\prod _{k=0}^{n-1} \Bigg\{ \int_{0}^{1} dt_{n-k}\;t_{n-k}^{2(n-k)-1+\nu } \int_{0}^{1} du_{n-k}\;u_{n-k}^{2(n-k )-\frac{3}{2} +\nu }\right.  \nonumber\\
&&\times  \frac{1}{2\pi i}  \oint dp_{n-k} \frac{1}{p_{n-k}} \left( 1-  w_{n-k+1,n}(1-t_{n-k})(1-u_{n-k})p_{n-k}\right)^{-\left( 4(n-k) +2\nu \right)}  \nonumber\\
&&\times  \left( \frac{p_{n-k}-1}{p_{n-k}} \frac{1}{1- w_{n-k+1,n}(1-t_{n-k})(1-u_{n-k})p_{n-k}}\right)^{\lambda_{n-k}}  \Bigg\} \nonumber\\
&&\times  \left.\left. \sum_{i_0=0}^{\lambda _0}\frac{(-\lambda_0)_{i_0}\left( \lambda_0 +2\nu \right)_{i_0}}{\left(1+\nu \right)_{i_0}\left(\nu +\frac{1}{2}\right)_{i_0}} w_{1,n}^{i_0}\right\} \eta ^n \right\}  \label{eq:60026}
\end{eqnarray}
where
\begin{equation} w _{i,j}=
\begin{cases} \displaystyle {\frac{p_i}{(p_i-1)}\; \frac{ w_{i+1,j} t_i u_i}{1-  w_{i+1,j} p_i (1-t_i)(1-u_i)}} \;\;\mbox{where}\; i\leq j\cr
x \;\;\mbox{only}\;\mbox{if}\; i>j
\end{cases}\nonumber 
\end{equation}
In the above, the first sub-integral form contains one term of $B_n's$, the second one contains two terms of $B_n$'s, the third one contains three terms of $B_n$'s, etc.
\end{theorem}
\begin{proof} 
According to (\ref{eq:60013}), 
\begin{equation}
 y(x)=  \sum_{n=0}^{\infty } y_{n}(x) = y_0(x)+ y_1(x)+ y_2(x)+y_3(x)+\cdots \label{eq:60027}
\end{equation}
In the above, the power series expansions of sub-summation $y_0(x) $, $y_1(x)$, $y_2(x)$ and $y_3(x)$ of Mathieu equation using R3TRF about $x=0$ are
\begin{subequations}
\begin{equation}
 y_0(x)= c_0 x^{\nu } \sum_{i_0=0}^{\lambda _0} \frac{(-\lambda_0)_{i_0} \left( \lambda _0+ 2\nu \right)_{i_0}}{(1+\nu )_{i_0}\left( \frac{1}{2}+ \nu \right)_{i_0}} x^{i_0} \label{eq:60028a}
\end{equation}
\begin{eqnarray}
 y_1(x) &=& c_0 x^{\nu } \left\{ \sum_{i_0=0}^{\lambda_0}\frac{1}{ (i_0+ 2+\nu )\left( i_0+\frac{3}{2} + \nu \right)}\frac{(-\lambda _0)_{i_0} \left( \lambda_0+  2\nu \right)_{i_0}}{(1+\nu )_{i_0}\left(\frac{1}{2}+\nu \right)_{i_0}} \right.\nonumber\\
&&\times \left. \sum_{i_1=i_0}^{\lambda _1} \frac{(-\lambda_1)_{i_1}\left(\lambda_1 + 4+2\nu \right)_{i_1}(3+\nu )_{i_0}\left(\frac{5}{2} +\nu \right)_{i_0}}{(-\lambda_1)_{i_0}\left(\lambda_1 + 4+2\nu \right)_{i_0}(3+\nu )_{i_1}\left(\frac{5}{2} +\nu \right)_{i_1}} x^{i_1}\right\}\eta  \label{eq:60028b}
\end{eqnarray}
\begin{eqnarray}
 y_2(x) &=& c_0 x^{\nu } \left\{ \sum_{i_0=0}^{\lambda_0}\frac{1}{ (i_0+ 2+\nu )\left( i_0+\frac{3}{2} + \nu \right)}\frac{(-\lambda _0)_{i_0} \left( \lambda_0+  2\nu \right)_{i_0}}{(1+\nu )_{i_0}\left(\frac{1}{2}+\nu \right)_{i_0}} \right.\nonumber\\
&&\times \sum_{i_1=i_0}^{\lambda _1} \frac{1}{ (i_1+ 4+\nu )\left( i_1+\frac{7}{2} + \nu \right)} \frac{(-\lambda_1)_{i_1}\left(\lambda_1 + 4+2\nu \right)_{i_1}(3+\nu )_{i_0}\left(\frac{5}{2} +\nu \right)_{i_0}}{(-\lambda_1)_{i_0}\left(\lambda_1 + 4+2\nu \right)_{i_0}(3+\nu )_{i_1}\left(\frac{5}{2} +\nu \right)_{i_1}} \nonumber\\
&&\times \left. \sum_{i_2=i_1}^{\lambda _2} \frac{(-\lambda _2)_{i_2}\left( \lambda_2 +8  +2\nu \right)_{i_2}\left( 5+\nu \right)_{i_1}\left( \frac{9}{2}+\nu \right)_{i_1}}{(-\lambda _2)_{i_1}\left( \lambda_2 +8  +2\nu \right)_{i_1}\left( 5+\nu \right)_{i_2}\left( \frac{9}{2}+\nu \right)_{i_2}} x^{i_2} \right\} \eta ^2  \label{eq:60028c}
\end{eqnarray}
\begin{eqnarray}
 y_3(x) &=&  c_0 x^{\nu } \left\{ \sum_{i_0=0}^{\lambda_0}\frac{1}{ (i_0+ 2+\nu )\left( i_0+\frac{3}{2} + \nu \right)}\frac{(-\lambda _0)_{i_0} \left( \lambda_0+  2\nu \right)_{i_0}}{(1+\nu )_{i_0}\left(\frac{1}{2}+\nu \right)_{i_0}} \right.\nonumber\\
&&\times \sum_{i_1=i_0}^{\lambda _1} \frac{1}{ (i_1+ 4+\nu )\left( i_1+\frac{7}{2} + \nu \right)} \frac{(-\lambda_1)_{i_1}\left(\lambda_1 + 4+2\nu \right)_{i_1}(3+\nu )_{i_0}\left(\frac{5}{2} +\nu \right)_{i_0}}{(-\lambda_1)_{i_0}\left(\lambda_1 + 4+2\nu \right)_{i_0}(3+\nu )_{i_1}\left(\frac{5}{2} +\nu \right)_{i_1}} \nonumber\\
&&\times  \sum_{i_2=i_1}^{\lambda _2} \frac{1}{ (i_2+ 6+\nu )\left( i_2+\frac{11}{2} + \nu \right)} \frac{(-\lambda _2)_{i_2}\left( \lambda_2 +8  +2\nu \right)_{i_2}\left( 5+\nu \right)_{i_1}\left( \frac{9}{2}+\nu \right)_{i_1}}{(-\lambda _2)_{i_1}\left( \lambda_2 +8  +2\nu \right)_{i_1}\left( 5+\nu \right)_{i_2}\left( \frac{9}{2}+\nu \right)_{i_2}}  \nonumber\\
&&\times \sum_{i_3=i_2}^{\lambda _3} \frac{(-\lambda _3)_{i_3}\left( \lambda_3 +12  +2\nu \right)_{i_3}\left( 7+\nu \right)_{i_2}\left( \frac{13}{2}+\nu \right)_{i_2}}{(-\lambda _3)_{i_2}\left( \lambda_3 +12  +2\nu \right)_{i_2}\left( 7+\nu \right)_{i_3}\left( \frac{13}{2}+\nu \right)_{i_3}} x^{i_3} \Bigg\} \eta ^3  \label{eq:60028d} 
\end{eqnarray}
\end{subequations}
Put $l=1$ in (\ref{eq:60025}). Take the new (\ref{eq:60025}) into (\ref{eq:60028b}).
\begin{eqnarray}
y_1(x) &=& \int_{0}^{1} dt_1\;t_1^{1+\nu } \int_{0}^{1} du_1\;u_1^{\frac{1}{2} +\nu } \frac{1}{2\pi i} \oint dp_1 \;\frac{1}{p_1} 
 (1-x(1-t_1)(1-u_1)p_1)^{-\left( 4+ 2\nu \right)} \nonumber\\
&&\times \left( \frac{p_1-1}{p_1} \frac{1}{1-x(1-t_1)(1-u_1)p_1}\right)^{\lambda _1}  \left\{ c_0 x^{\nu } \sum_{i_0=0}^{\lambda _0} \frac{(-\lambda_0)_{i_0} \left( \lambda _0+ 2\nu \right)_{i_0}}{(1+\nu )_{i_0}\left( \frac{1}{2}+ \nu \right)_{i_0}} w_{1,1}^{i_0} \right\} \eta \hspace{2cm}\label{eq:60029}
\end{eqnarray}
where
\begin{equation}
w_{1,1} = \frac{p_1}{p_1-1} \frac{x t_1 u_1}{1-x(1-t_1)(1-u_1)p_1}\nonumber
\end{equation}
Put $l=2$ in (\ref{eq:60025}). Take the new (\ref{eq:60025}) into (\ref{eq:60028c}).
\begin{eqnarray}
y_2(x) &=& c_0 x^{\nu } \int_{0}^{1} dt_2\;t_2^{3+\nu } \int_{0}^{1} du_2\;u_2^{\frac{5}{2} +\nu } \frac{1}{2\pi i} \oint dp_2 \;\frac{1}{p_2} (1-x(1-t_2)(1-u_2)p_2)^{-\left( 8 +2 \nu \right)}
 \nonumber\\
&&\times \left( \frac{p_2-1}{p_2} \frac{1}{1-x(1-t_2)(1-u_2)p_2}\right)^{\lambda _2}  \nonumber\\
&&\times \left\{ \sum_{i_0=0}^{\lambda_0}\frac{1}{ (i_0+ 2+\nu )\left( i_0+\frac{3}{2} + \nu \right)}\frac{(-\lambda _0)_{i_0} \left( \lambda_0+  2\nu \right)_{i_0}}{(1+\nu )_{i_0}\left(\frac{1}{2}+\nu \right)_{i_0}} \right.\nonumber\\
&&\times \left. \sum_{i_1=i_0}^{\lambda _1} \frac{(-\lambda_1)_{i_1}\left(\lambda_1 + 4+2\nu \right)_{i_1}(3+\nu )_{i_0}\left(\frac{5}{2} +\nu \right)_{i_0}}{(-\lambda_1)_{i_0}\left(\lambda_1 + 4+2\nu \right)_{i_0}(3+\nu )_{i_1}\left(\frac{5}{2} +\nu \right)_{i_1}} w_{2,2}^{i_1}\right\} \eta^2 \label{eq:60030}
\end{eqnarray}
where
\begin{equation}
w_{2,2} = \frac{p_2}{p_2-1} \frac{x t_2 u_2}{1-x(1-t_2)(1-u_2)p_2}\nonumber
\end{equation}
Put $l=1$ and $\eta = w_{2,2}$ in (\ref{eq:60025}). Take the new (\ref{eq:60025}) into (\ref{eq:60030}).
\begin{eqnarray}
y_2(x) &=& \int_{0}^{1} dt_2\;t_2^{3+\nu } \int_{0}^{1} du_2\;u_2^{\frac{5}{2} +\nu } \frac{1}{2\pi i} \oint dp_2 \;\frac{1}{p_2} (1-x(1-t_2)(1-u_2)p_2)^{-\left( 8 +2 \nu \right)}
 \nonumber\\
&&\times \left( \frac{p_2-1}{p_2} \frac{1}{1-x(1-t_2)(1-u_2)p_2}\right)^{\lambda _2}  \nonumber\\
&&\times \int_{0}^{1} dt_1\;t_1^{1+\nu } \int_{0}^{1} du_1\;u_1^{\frac{1}{2} +\nu } \frac{1}{2\pi i} \oint dp_1 \;\frac{1}{p_1} (1- w_{2,2}(1-t_1)(1-u_1)p_1)^{-\left( 4 +2 \nu \right)} \nonumber\\
&&\times \left( \frac{p_1-1}{p_1} \frac{1}{1- w_{2,2}(1-t_1)(1-u_1)p_1}\right)^{\lambda _1} \left\{ c_0 x^{\nu } \sum_{i_0=0}^{\lambda _0} \frac{(-\lambda_0)_{i_0} \left( \lambda _0+ 2\nu \right)_{i_0}}{(1+\nu )_{i_0}\left( \frac{1}{2}+ \nu \right)_{i_0}} w_{1,2}^{i_0} \right\} \eta^2 \hspace{2cm}\label{eq:60031}
\end{eqnarray}
where
\begin{equation}
w_{1,2} = \frac{p_1}{p_1-1} \frac{ w_{2,2} t_1 u_1}{1- w_{2,2}(1-t_1)(1-u_1)p_1}\nonumber
\end{equation}
By using similar process for the previous cases of integral forms of $y_1(x)$ and $y_2(x)$, the integral form of sub-power series expansion of $y_3(x)$ is
\begin{eqnarray}
y_3(x) &=& \int_{0}^{1} dt_3\;t_3^{5+\nu } \int_{0}^{1} du_3\;u_3^{\frac{9}{2} +\nu } \frac{1}{2\pi i} \oint dp_3 \;\frac{1}{p_3} (1-x(1-t_3)(1-u_3)p_3)^{-\left( 12 +2 \nu \right)}
 \nonumber\\
&&\times \left( \frac{p_3-1}{p_3} \frac{1}{1-x(1-t_3)(1-u_3)p_3}\right)^{\lambda _3}  \nonumber\\
&&\times \int_{0}^{1} dt_2\;t_2^{3+\nu } \int_{0}^{1} du_2\;u_2^{\frac{5}{2} +\nu } \frac{1}{2\pi i} \oint dp_2 \;\frac{1}{p_2} (1- w_{3,3}(1-t_2)(1-u_2)p_2)^{-\left( 8 +2 \nu \right)} \nonumber\\
&&\times \left( \frac{p_2-1}{p_2} \frac{1}{1- w_{3,3}(1-t_2)(1-u_2)p_2}\right)^{\lambda _2} \nonumber\\
&&\times \int_{0}^{1} dt_1\;t_1^{1+\nu } \int_{0}^{1} du_1\;u_1^{\frac{1}{2} +\nu } \frac{1}{2\pi i} \oint dp_1 \;\frac{1}{p_1} (1- w_{2,3}(1-t_1)(1-u_1)p_1)^{-\left( 4 +2 \nu \right)} \nonumber\\
&&\times \left( \frac{p_1-1}{p_1} \frac{1}{1- w_{2,3}(1-t_1)(1-u_1)p_1}\right)^{\lambda _1}
\left\{ c_0 x^{\nu } \sum_{i_0=0}^{\lambda _0} \frac{(-\lambda_0)_{i_0} \left( \lambda _0+ 2\nu \right)_{i_0}}{(1+\nu )_{i_0}\left( \frac{1}{2}+ \nu \right)_{i_0}} w_{1,3}^{i_0} \right\} \eta^3 \hspace{2cm}\label{eq:60032}
\end{eqnarray}
where
\begin{equation}
\begin{cases} w_{3,3} = \frac{p_3}{p_3-1} \frac{ x t_3 u_3}{1- x(1-t_3)(1-u_3)p_3} \cr
w_{2,3} = \frac{p_2}{p_2-1} \frac{w_{3,3} t_2 u_2}{1- w_{3,3}(1-t_2)(1-u_2)p_2} \cr
w_{1,3} = \frac{p_1}{p_1-1} \frac{w_{2,3} t_1 u_1}{1- w_{2,3}(1-t_1)(1-u_1)p_1}
\end{cases}
\nonumber
\end{equation}
By repeating this process for all higher terms of integral forms of sub-summation $y_m(x)$ terms where $m \geq 4$, we obtain every integral forms of $y_m(x)$ terms. 
Substitute (\ref{eq:60028a}), (\ref{eq:60029}), (\ref{eq:60031}), (\ref{eq:60032}) and including all integral forms of $y_m(x)$ terms where $m \geq 4$ into (\ref{eq:60026}). 
\end{proof}
Put $c_0$= 1 as $\nu =0$ for the first kind of independent solutions of Mathieu equation and $\nu = 1/2 $ for the second one in (\ref{eq:60026}).
\begin{remark}
The integral representation of Mathieu equation of the first kind for polynomial of type 2 about $x=0$ as $\lambda = 2^2(\lambda_j +2j )^2- 2 q$ where $j,\lambda _j \in \mathbb{N}_{0}$ is
\begin{eqnarray}
  y(x)&=& MF_{\lambda _j}^R\left( q,\lambda = 2^2(\lambda_j +2j )^2- 2 q; \eta  = q x^2, x= \mathrm{\cos}^2z \right)\nonumber\\
&=&\; _2F_1 \left( -\lambda _0, \lambda_0 ; \frac{1}{2}; x \right)  + \sum_{n=1}^{\infty } \Bigg\{\prod _{k=0}^{n-1} \Bigg\{ \int_{0}^{1} dt_{n-k}\;t_{n-k}^{2(n-k)-1 } \int_{0}^{1} du_{n-k}\;u_{n-k}^{2(n-k )-\frac{3}{2} } \nonumber\\
&&\times  \frac{1}{2\pi i}  \oint dp_{n-k} \frac{1}{p_{n-k}} \left( 1-  w_{n-k+1,n}(1-t_{n-k})(1-u_{n-k})p_{n-k}\right)^{- 4(n-k) }  \nonumber\\
&&\times  \left( \frac{p_{n-k}-1}{p_{n-k}} \frac{1}{1- w_{n-k+1,n}(1-t_{n-k})(1-u_{n-k})p_{n-k}}\right)^{\lambda_{n-k}}  \Bigg\} \nonumber\\
&&\times  \left.  _2F_1 \left( -\lambda _0, \lambda_0 ; \frac{1}{2}; w_{1,n} \right) \right\} \eta ^n \hspace{1cm}\label{eq:60033}
\end{eqnarray}
\end{remark} 
\begin{remark}
The integral representation of Mathieu equation of the second kind for polynomial of type 2 about $x=0$ as $\lambda = 2^2\left(\lambda_j +2j+1/2 \right)^2- 2 q $ where $j,\lambda _j \in \mathbb{N}_{0}$ is
\begin{eqnarray}
 y(x)&=& MS_{\lambda _j}^R\left( q,\lambda =2^2\left(\lambda_j +2j+1/2 \right)^2- 2 q; \eta  = q x^2, x= \mathrm{\cos}^2z \right) \nonumber\\
&=& x^{\frac{1}{2}} \left\{\; _2F_1 \left( -\lambda _0, \lambda_0 +1; \frac{3}{2}; x \right) \right. + \sum_{n=1}^{\infty } \Bigg\{\prod _{k=0}^{n-1} \Bigg\{ \int_{0}^{1} dt_{n-k}\;t_{n-k}^{2(n-k)-\frac{1}{2} } \int_{0}^{1} du_{n-k}\;u_{n-k}^{2(n-k )-1 } \nonumber\\
&&\times  \frac{1}{2\pi i}  \oint dp_{n-k} \frac{1}{p_{n-k}} \left( 1-  w_{n-k+1,n}(1-t_{n-k})(1-u_{n-k})p_{n-k}\right)^{-\left( 4(n-k) +1 \right)}  \nonumber\\
&&\times  \left( \frac{p_{n-k}-1}{p_{n-k}} \frac{1}{1- w_{n-k+1,n}(1-t_{n-k})(1-u_{n-k})p_{n-k}}\right)^{\lambda_{n-k}}  \Bigg\} \nonumber\\
&&\times  \left.\left.  _2F_1 \left( -\lambda _0, \lambda_0 +1; \frac{3}{2}; w_{1,n} \right) \right\} \eta ^n \right\} \label{eq:60034}
\end{eqnarray}
\end{remark}
\subsubsection{Infinite series}
Let's consider the integral representation of Mathieu equation about $x=0$ for infinite series by applying R3TRF.
There is a generalized hypergeometric function which is written by
\begin{eqnarray}
M_l &=& \sum_{i_l= i_{l-1}}^{\infty } \frac{\left(\nu +2l-\varphi \right)_{i_l}\left(\nu +2l+\varphi \right)_{i_l}(2l+1+\nu )_{i_{l-1}}\left( 2l+\frac{1}{2} +\nu \right)_{i_{l-1}}}{\left(\nu +2l-\varphi \right)_{i_{l-1}}\left(\nu +2l+\varphi \right)_{i_{l-1}}(2l+1+\nu )_{i_l}(2l+\frac{1}{2} +\nu )_{i_l}} x^{i_l}\label{er:60021}\\
&=& x^{i_{l-1}} 
\sum_{j=0}^{\infty } \frac{B\left(i_{l-1}+2l-\frac{1}{2}+\nu ,j+1\right) B\left( i_{l-1}+2l +\nu ,j+1\right) ( \nu +2l-\varphi  +i_{l-1} )_j \left( \nu +2l+\varphi  + i_{l-1} \right)_j}{(i_{l-1}+2l+\nu )^{-1}\left( i_{l-1}+2l-\frac{1}{2}+\nu \right)^{-1}(1)_j \;j!} x^j\nonumber
\end{eqnarray}
Substitute (\ref{eq:60022a}) and (\ref{eq:60022b}) into (\ref{er:60021}). And divide $(i_{l-1}+2l+\nu ) \left( i_{l-1}+2l-\frac{1}{2}+\nu \right) $ into the new (\ref{er:60021}).
\begin{eqnarray}
V_l&=& \frac{1}{(i_{l-1}+2l+\nu ) \left( i_{l-1}+2l-\frac{1}{2}+\nu \right)} \nonumber\\
&&\times \sum_{i_l= i_{l-1}}^{\infty } \frac{\left(\nu +2l-\varphi \right)_{i_l}\left(\nu +2l+\varphi \right)_{i_l}(2l+1+\nu )_{i_{l-1}}\left( 2l+\frac{1}{2} +\nu \right)_{i_{l-1}}}{\left(\nu +2l-\varphi \right)_{i_{l-1}}\left(\nu +2l+\varphi \right)_{i_{l-1}}(2l+1+\nu )_{i_l}(2l+\frac{1}{2} +\nu )_{i_l}} x^{i_l}\nonumber\\
&=&  \int_{0}^{1} dt_l\;t_l^{2l-1+\nu } \int_{0}^{1} du_l\;u_l^{2l -\frac{3}{2} +\nu } (x t_l u_l)^{i_{l-1}} \nonumber\\
&&\times  \sum_{j=0}^{\infty } \frac{(\nu +2l-\varphi  +i_{l-1})_j \left( \nu +2l+\varphi + i_{l-1} \right)_j}{(1)_j \;j!} 
(x(1-t_l)(1-u_l))^j \nonumber 
\end{eqnarray}
The hypergeometric function is defined by
\begin{eqnarray}
_2F_1 \left( \alpha ,\beta ; \gamma ; z \right) &=& \sum_{n=0}^{\infty } \frac{(\alpha )_n (\beta )_n}{(\gamma )_n (n!)} z^n \nonumber\\
&=&  \frac{1}{2\pi i} \frac{\Gamma( 1+\alpha  -\gamma )}{\Gamma (\alpha )} \int_0^{(1+)} dv_l\; (-1)^{\gamma }(-v_l)^{\alpha -1} (1-v_l )^{\gamma -\alpha -1} (1-zv_l)^{-\beta }\hspace{1.5cm}\label{er:60022}\\
&& \mbox{where} \;\gamma -\alpha  \ne 1,2,3,\cdots, \;\mbox{Re}(\alpha )>0 \nonumber
\end{eqnarray}
Replace $\alpha $, $\beta $, $\gamma $ and $z$ by $\nu +2l-\varphi  +i_{l-1}$, $\nu +2l+\varphi + i_{l-1}$, 1 and $x(1-t_l)(1-u_l)$ in (\ref{er:60022}). Take the new (\ref{er:60022}) into $V_l$.
\begin{eqnarray}
V_l&=& \frac{1}{(i_{l-1}+2l+\nu ) \left( i_{l-1}+2l-\frac{1}{2}+\nu \right)} \nonumber\\
&&\times \sum_{i_l= i_{l-1}}^{\infty } \frac{\left(\nu +2l-\varphi \right)_{i_l}\left(\nu +2l+\varphi \right)_{i_l}(2l+1+\nu )_{i_{l-1}}\left( 2l+\frac{1}{2} +\nu \right)_{i_{l-1}}}{\left(\nu +2l-\varphi \right)_{i_{l-1}}\left(\nu +2l+\varphi \right)_{i_{l-1}}(2l+1+\nu )_{i_l}(2l+\frac{1}{2} +\nu )_{i_l}} x^{i_l}\nonumber\\
&=&  \int_{0}^{1} dt_l\;t_l^{2l-1+\nu } \int_{0}^{1} du_l\;u_l^{2l-\frac{3}{2} +\nu } 
\frac{1}{2\pi i} \oint dp_l\;\frac{1}{p_l} \left(\frac{p_l-1}{p_l} \right)^{-(\nu +2l)+\varphi }  \nonumber\\
&&\times (1-x(1-t_l)(1-u_l)p_l)^{-(\nu +2l)-\varphi} \left(\frac{p_l}{p_l-1} \frac{xt_l u_l}{1-x(1-t_l)(1-u_l)p_l}\right)^{i_{l-1}}
\label{er:60023} 
\end{eqnarray}
Substitute (\ref{er:60023}) into (\ref{eq:60017}) where $l=1,2,3,\cdots$; apply $V_1$ into the second summation of sub-power series $y_1(x)$, apply $V_2$ into the third summation and $V_1$ into the second summation of sub-power series $y_2(x)$, apply $V_3$ into the forth summation, $V_2$ into the third summation and $V_1$ into the second summation of sub-power series $y_3(x)$, etc.\footnote{$y_1(x)$ means the sub-power series in (\ref{eq:60017}) contains one term of $B_n's$, $y_2(x)$ means the sub-power series in (\ref{eq:60017}) contains two terms of $B_n's$, $y_3(x)$ means the sub-power series in (\ref{eq:60017}) contains three terms of $B_n's$, etc.}
\begin{theorem}
The general expression of an integral form of Mathieu equation for infinite series about $x=0$ using R3TRF is given by
\begin{eqnarray}
y(x) &=&  \sum_{n=0}^{\infty } y_{n}(x) = y_0(x)+ y_1(x)+ y_2(x)+y_3(x)+\cdots \nonumber\\
&=& c_0 x^{\nu } \left\{ \sum_{i_0=0}^{\infty }\frac{(\nu -\varphi )_{i_0}\left( \nu +\varphi \right)_{i_0}}{\left(1+\nu \right)_{i_0}\left(\nu +\frac{1}{2}\right)_{i_0}} x^{i_0} \right. \nonumber\\
&&+ \sum_{n=1}^{\infty } \Bigg\{\prod _{k=0}^{n-1} \Bigg\{ \int_{0}^{1} dt_{n-k}\;t_{n-k}^{2(n-k)-1+\nu } \int_{0}^{1} du_{n-k}\;u_{n-k}^{2(n-k )-\frac{3}{2} +\nu }  \nonumber\\
&&\times  \frac{1}{2\pi i}  \oint dp_{n-k} \frac{1}{p_{n-k}} \left( \frac{p_{n-k}-1}{p_{n-k}} \right)^{-(\nu + 2(n-k))+ \varphi } \nonumber\\
&&\times \left( 1- w_{n-k+1,n}(1-t_{n-k})(1-u_{n-k})p_{n-k}\right)^{-(\nu + 2(n-k))- \varphi} \Bigg\}\nonumber\\
&&\times \left. \sum_{i_0=0}^{\infty }\frac{(\nu -\varphi )_{i_0}\left( \nu +\varphi \right)_{i_0}}{\left(1+\nu \right)_{i_0}\left(\nu +\frac{1}{2}\right)_{i_0}} w_{1,n}^{i_0}\Bigg\} \eta ^n \right\}  \label{eq:60035}
\end{eqnarray}
where
\begin{equation}
\varphi  = \frac{1}{2}\sqrt{\lambda +2q} \nonumber
\end{equation}
In the above, the first sub-integral form contains one term of $B_n's$, the second one contains two terms of $B_n$'s, the third one contains three terms of $B_n$'s, etc.
\end{theorem}
\begin{proof} 
In (\ref{eq:60017}) sub-power series $y_0(x) $, $y_1(x)$, $y_2(x)$ and $y_3(x)$ of Mathieu equation for infinite series about $x=0$ using R3TRF are given by
\begin{subequations}
\begin{equation}
 y_0(x)= c_0 x^{\nu } \sum_{i_0=0}^{\infty } \frac{(\nu -\varphi )_{i_0} \left( \nu + \varphi \right)_{i_0}}{(1+\nu )_{i_0}\left( \frac{1}{2}+ \nu \right)_{i_0}} x^{i_0} \label{er:60024a}
\end{equation}
\begin{eqnarray}
 y_1(x) &=& c_0 x^{\nu } \left\{ \sum_{i_0=0}^{\infty }\frac{1}{ (i_0+ 2+\nu )\left( i_0+\frac{3}{2} + \nu \right)}\frac{(\nu -\varphi)_{i_0} \left( \nu + \varphi\right)_{i_0}}{(1+\nu )_{i_0}\left(\frac{1}{2}+\nu \right)_{i_0}} \right.\nonumber\\
&&\times \left. \sum_{i_1=i_0}^{\infty } \frac{(\nu +2 -\varphi)_{i_1}\left(\nu +2+\varphi \right)_{i_1}(3+\nu )_{i_0}\left(\frac{5}{2} +\nu \right)_{i_0}}{(\nu +2 -\varphi)_{i_0}\left(\nu +2+\varphi\right)_{i_0}(3+\nu )_{i_1}\left(\frac{5}{2} +\nu \right)_{i_1}} x^{i_1}\right\}\eta  \label{er:60024b}
\end{eqnarray}
\begin{eqnarray}
 y_2(x) &=& c_0 x^{\nu } \left\{ \sum_{i_0=0}^{\infty }\frac{1}{ (i_0+ 2+\nu )\left( i_0+\frac{3}{2} + \nu \right)}\frac{(\nu -\varphi)_{i_0} \left( \nu + \varphi\right)_{i_0}}{(1+\nu )_{i_0}\left(\frac{1}{2}+\nu \right)_{i_0}} \right.\nonumber\\
&&\times \sum_{i_1=i_0}^{\infty } \frac{1}{ (i_1+ 4+\nu )\left( i_1+\frac{7}{2} + \nu \right)} \frac{(\nu +2 -\varphi)_{i_1}\left(\nu +2+\varphi\right)_{i_1}(3+\nu )_{i_0}\left(\frac{5}{2} +\nu \right)_{i_0}}{(\nu +2 -\varphi)_{i_0}\left(\nu +2+\varphi\right)_{i_0}(3+\nu )_{i_1}\left(\frac{5}{2} +\nu \right)_{i_1}} \nonumber\\
&&\times \left. \sum_{i_2=i_1}^{\infty } \frac{(\nu +4-\varphi)_{i_2}\left( \nu +4+\varphi \right)_{i_2}\left( 5+\nu \right)_{i_1}\left( \frac{9}{2}+\nu \right)_{i_1}}{(\nu +4-\varphi)_{i_1}\left( \nu +4+\varphi \right)_{i_1}\left( 5+\nu \right)_{i_2}\left( \frac{9}{2}+\nu \right)_{i_2}} x^{i_2} \right\} \eta ^2  \label{er:60024c}
\end{eqnarray}
\begin{eqnarray}
 y_3(x) &=&  c_0 x^{\nu } \left\{ \sum_{i_0=0}^{\infty }\frac{1}{ (i_0+ 2+\nu )\left( i_0+\frac{3}{2} + \nu \right)}\frac{(\nu -\varphi)_{i_0} \left( \nu + \varphi\right)_{i_0}}{(1+\nu )_{i_0}\left(\frac{1}{2}+\nu \right)_{i_0}} \right.\nonumber\\
&&\times \sum_{i_1=i_0}^{\infty } \frac{1}{ (i_1+ 4+\nu )\left( i_1+\frac{7}{2} + \nu \right)} \frac{(\nu +2 -\varphi)_{i_1}\left(\nu +2+\varphi\right)_{i_1}(3+\nu )_{i_0}\left(\frac{5}{2} +\nu \right)_{i_0}}{(\nu +2 -\varphi)_{i_0}\left(\nu +2+\varphi\right)_{i_0}(3+\nu )_{i_1}\left(\frac{5}{2} +\nu \right)_{i_1}} \nonumber\\
&&\times  \sum_{i_2=i_1}^{\infty } \frac{1}{ (i_2+ 6+\nu )\left( i_2+\frac{11}{2} + \nu \right)} \frac{(\nu +4-\varphi)_{i_2}\left( \nu +4+\varphi \right)_{i_2}\left( 5+\nu \right)_{i_1}\left( \frac{9}{2}+\nu \right)_{i_1}}{(\nu +4-\varphi)_{i_1}\left( \nu +4+\varphi \right)_{i_1}\left( 5+\nu \right)_{i_2}\left( \frac{9}{2}+\nu \right)_{i_2}}  \nonumber\\
&&\times \left. \sum_{i_3=i_2}^{\infty } \frac{(\nu +6-\varphi)_{i_3}\left( \nu +6+\varphi\right)_{i_3}\left( 7+\nu \right)_{i_2}\left( \frac{13}{2}+\nu \right)_{i_2}}{(\nu +6-\varphi)_{i_2}\left( \nu +6+\varphi \right)_{i_2}\left( 7+\nu \right)_{i_3}\left( \frac{13}{2}+\nu \right)_{i_3}} x^{i_3} \right\} \eta ^3  \label{er:60024d} 
\end{eqnarray}
\end{subequations}
Put $l=1$ in (\ref{er:60023}). Take the new (\ref{er:60023}) into (\ref{er:60024b}).
\begin{eqnarray}
y_1(x) &=& \int_{0}^{1} dt_1\;t_1^{1+\nu } \int_{0}^{1} du_1\;u_1^{\frac{1}{2} +\nu } \frac{1}{2\pi i} \oint dp_1 \;\frac{1}{p_1}\left( \frac{p_1-1}{p_1} \right)^{-(\nu+2 )+\varphi }  \nonumber\\
&&\times (1-x(1-t_1)(1-u_1)p_1)^{-(\nu+2 )-\varphi} \left\{ c_0 x^{\nu } \sum_{i_0=0}^{\infty } \frac{(\nu -\varphi )_{i_0} \left( \nu + \varphi \right)_{i_0}}{(1+\nu )_{i_0}\left( \frac{1}{2}+ \nu \right)_{i_0}} w_{1,1}^{i_0} \right\} \eta \hspace{1.5cm}\label{er:60025}
\end{eqnarray}
where
\begin{equation}
w_{1,1} = \frac{p_1}{p_1-1} \frac{x t_1 u_1}{1-x(1-t_1)(1-u_1)p_1}\nonumber
\end{equation}
Put $l=2$ in (\ref{er:60023}). Take the new (\ref{er:60023}) into (\ref{er:60024c}).
\begin{eqnarray}
y_2(x) &=& c_0 x^{\nu } \int_{0}^{1} dt_2\;t_2^{3+\nu } \int_{0}^{1} du_2\;u_2^{\frac{5}{2} +\nu } \frac{1}{2\pi i} \oint dp_2 \;\frac{1}{p_2} \left( \frac{p_2-1}{p_2} \right)^{-(\nu +4)+\varphi} \nonumber\\
&&\times  (1-x(1-t_2)(1-u_2)p_2)^{-(\nu +4)-\varphi} \nonumber\\
&&\times \left\{ \sum_{i_0=0}^{\infty }\frac{1}{ (i_0+ 2+\nu )\left( i_0+\frac{3}{2} + \nu \right)}\frac{(\nu -\varphi)_{i_0} \left( \nu + \varphi\right)_{i_0}}{(1+\nu )_{i_0}\left(\frac{1}{2}+\nu \right)_{i_0}} \right.\nonumber\\
&&\times \left. \sum_{i_1=i_0}^{\infty } \frac{(\nu +2 -\varphi)_{i_1}\left(\nu +2+\varphi \right)_{i_1}(3+\nu )_{i_0}\left(\frac{5}{2} +\nu \right)_{i_0}}{(\nu +2 -\varphi)_{i_0}\left(\nu +2+\varphi\right)_{i_0}(3+\nu )_{i_1}\left(\frac{5}{2} +\nu \right)_{i_1}} w_{2,2}^{i_1}\right\} \eta^2 \label{er:60026}
\end{eqnarray}
where
\begin{equation}
w_{2,2} = \frac{p_2}{p_2-1} \frac{x t_2 u_2}{1-x(1-t_2)(1-u_2)p_2}\nonumber
\end{equation}
Put $l=1$ and $\eta = w_{2,2}$ in (\ref{er:60023}). Take the new (\ref{er:60023}) into (\ref{er:60026}).
\begin{eqnarray}
y_2(x) &=& \int_{0}^{1} dt_2\;t_2^{3+\nu } \int_{0}^{1} du_2\;u_2^{\frac{5}{2} +\nu } \frac{1}{2\pi i} \oint dp_2 \;\frac{1}{p_2} \left( \frac{p_2-1}{p_2} \right)^{-(\nu +4)+\varphi} \nonumber\\
&&\times  (1-x(1-t_2)(1-u_2)p_2)^{-(\nu +4)-\varphi} \nonumber\\
&&\times \int_{0}^{1} dt_1\;t_1^{1+\nu } \int_{0}^{1} du_1\;u_1^{\frac{1}{2} +\nu } \frac{1}{2\pi i} \oint dp_1 \;\frac{1}{p_1} \left( \frac{p_1-1}{p_1} \right)^{-(\nu +2)+\varphi}  \nonumber\\
&&\times (1- w_{2,2}(1-t_1)(1-u_1)p_1)^{-(\nu +2)-\varphi} \left\{ c_0 x^{\nu } \sum_{i_0=0}^{\infty } \frac{(\nu -\varphi )_{i_0} \left( \nu + \varphi \right)_{i_0}}{(1+\nu )_{i_0}\left( \frac{1}{2}+ \nu \right)_{i_0}} w_{1,2}^{i_0} \right\} \eta^2 \hspace{1.5cm}\label{er:60027}
\end{eqnarray}
where
\begin{equation}
w_{1,2} = \frac{p_1}{p_1-1} \frac{ w_{2,2} t_1 u_1}{1- w_{2,2}(1-t_1)(1-u_1)p_1}\nonumber
\end{equation}
By using similar process for the previous cases of integral forms of $y_1(x)$ and $y_2(x)$, the integral form of sub-power series expansion of $y_3(x)$ is
\begin{eqnarray}
y_3(x) &=& \int_{0}^{1} dt_3\;t_3^{5+\nu } \int_{0}^{1} du_3\;u_3^{\frac{9}{2} +\nu } \frac{1}{2\pi i} \oint dp_3 \;\frac{1}{p_3} \left( \frac{p_3-1}{p_3}  \right)^{-(\nu +6)+\varphi}  \nonumber\\
&&\times (1-x(1-t_3)(1-u_3)p_3)^{-(\nu +6)-\varphi} \nonumber\\
&&\times \int_{0}^{1} dt_2\;t_2^{3+\nu } \int_{0}^{1} du_2\;u_2^{\frac{5}{2} +\nu } \frac{1}{2\pi i} \oint dp_2 \;\frac{1}{p_2} \left( \frac{p_2-1}{p_2}  \right)^{-(\nu +4)+\varphi} \nonumber\\
&&\times (1- w_{3,3}(1-t_2)(1-u_2)p_2)^{-(\nu +4)-\varphi} \nonumber\\
&&\times \int_{0}^{1} dt_1\;t_1^{1+\nu } \int_{0}^{1} du_1\;u_1^{\frac{1}{2} +\nu } \frac{1}{2\pi i} \oint dp_1 \;\frac{1}{p_1} \left( \frac{p_1-1}{p_1} \right)^{-(\nu +2)+\varphi} \nonumber\\
&&\times (1- w_{2,3}(1-t_1)(1-u_1)p_1)^{-(\nu +2)-\varphi}
\left\{ c_0 x^{\nu } \sum_{i_0=0}^{\infty } \frac{(\nu -\varphi )_{i_0} \left( \nu + \varphi \right)_{i_0}}{(1+\nu )_{i_0}\left( \frac{1}{2}+ \nu \right)_{i_0}} w_{1,3}^{i_0} \right\} \eta^3 \hspace{1.5cm}\label{er:60028}
\end{eqnarray}
where
\begin{equation}
\begin{cases} w_{3,3} = \frac{p_3}{p_3-1} \frac{ x t_3 u_3}{1- x(1-t_3)(1-u_3)p_3} \cr
w_{2,3} = \frac{p_2}{p_2-1} \frac{w_{3,3} t_2 u_2}{1- w_{3,3}(1-t_2)(1-u_2)p_2} \cr
w_{1,3} = \frac{p_1}{p_1-1} \frac{w_{2,3} t_1 u_1}{1- w_{2,3}(1-t_1)(1-u_1)p_1}
\end{cases}
\nonumber
\end{equation}
By repeating this process for all higher terms of integral forms of sub-summation $y_m(x)$ terms where $m \geq 4$, we obtain every integral forms of $y_m(x)$ terms. 
Since we substitute (\ref{er:60024a}), (\ref{er:60025}), (\ref{er:60027}), (\ref{er:60028}) and including all integral forms of $y_m(x)$ terms where $m \geq 4$ into (\ref{eq:60017}),  we obtain (\ref{eq:60035}).\footnote{Or replace the finite summation with an interval $[0, \lambda _0]$ by infinite summation with an interval  $[0,\infty ]$ in (\ref{eq:60026}). Replace $\lambda _0$ and $\lambda _{n-k}$ by $\frac{1}{2}\sqrt{\lambda +2q}-\nu$ and $\frac{1}{2}\sqrt{\lambda +2q}-2(n-k)-\nu$ into the new (\ref{eq:60026}). Its solution is also equivalent to (\ref{eq:60035})}
\end{proof}
Put $c_0$= 1 as $\nu =0$ for the first kind of independent solutions of Mathieu equation and $\nu = 1/2 $ for the second one in (\ref{eq:60035}). 
\begin{remark}
The integral representation of Mathieu equation of the first kind for infinite series about $x=0$ using R3TRF is
\begin{eqnarray}
y(x)&=&  MF^R\left( q,\lambda, \varphi  = \frac{1}{2} \sqrt{\lambda + 2q} ; \eta  = q x^2, x= \mathrm{\cos}^2z \right) \nonumber\\
&=& \; _2F_1 \left( -\varphi, \varphi ; \frac{1}{2}; x \right)  + \sum_{n=1}^{\infty } \Bigg\{\prod _{k=0}^{n-1} \Bigg\{ \int_{0}^{1} dt_{n-k}\;t_{n-k}^{2(n-k)-1 } \int_{0}^{1} du_{n-k}\;u_{n-k}^{2(n-k )-\frac{3}{2} }  \nonumber\\
&&\times  \frac{1}{2\pi i}  \oint dp_{n-k} \frac{1}{p_{n-k}} \left( \frac{p_{n-k}-1}{p_{n-k}} \right)^{-2(n-k)+ \varphi } \nonumber\\
&&\times  \left( 1- w_{n-k+1,n}(1-t_{n-k})(1-u_{n-k})p_{n-k}\right)^{-2(n-k)-\varphi}   \Bigg\} \nonumber\\
&&\times  \; _2F_1 \left( -\varphi, \varphi ; \frac{1}{2}; w_{1,n} \right)  \Bigg\} \eta ^n  \label{eq:60036}
\end{eqnarray}
\end{remark}
\begin{remark}
The integral representation of Mathieu equation of the second kind for infinite series about $x=0$ using R3TRF is
\begin{eqnarray}
 y(x)&=& MS^R\left( q,\lambda, \varphi  = \frac{1}{2} \sqrt{\lambda + 2q}; \eta  = q x^2, x= \mathrm{\cos}^2z \right)\nonumber\\
&=& x^{\frac{1}{2}} \left\{\; _2F_1 \left( \frac{1}{2}-\varphi, \frac{1}{2}+\varphi ; \frac{3}{2}; x \right) \right. + \sum_{n=1}^{\infty } \Bigg\{\prod _{k=0}^{n-1} \Bigg\{ \int_{0}^{1} dt_{n-k}\;t_{n-k}^{2(n-k)-\frac{1}{2} } \int_{0}^{1} du_{n-k}\;u_{n-k}^{2(n-k )-1 }  \nonumber\\
&&\times  \frac{1}{2\pi i}  \oint dp_{n-k} \frac{1}{p_{n-k}} \left( \frac{p_{n-k}-1}{p_{n-k}} \right)^{-(2(n-k) +\frac{1}{2})+ \varphi}  \nonumber\\
&&\times \left( 1- w_{n-k+1,n}(1-t_{n-k})(1-u_{n-k})p_{n-k}\right)^{-(2(n-k) +\frac{1}{2})-\varphi}   \Bigg\} \nonumber\\
&&\times \left. \; _2F_1 \left( \frac{1}{2}-\varphi, \frac{1}{2}+\varphi ; \frac{3}{2}; w_{1,n} \right) \Bigg\} \eta ^n \right\}   \label{eq:60037}
\end{eqnarray}
\end{remark}
Two integral forms of Mathieu equation for infinite series in this chapter and Ref.\cite{1Chou2012e} are equivalent to each other. In this chapter, $_2F_1$ function recurs in each of sub-integral forms of the Mathieu function. In Ref.\cite{1Chou2012e}, the Modified Bessel function recurs in each of sub-integral forms of it.
\subsection{Generating function for the Mathieu polynomial of type 2}
I consider the generating function for the Mathieu polynomial of type 2. Since its generating function is derived, we might be possible to construct orthogonal relations of it.
\begin{lemma}
The generating function for Jacobi polynomial using hypergeometric functions is given by
\begin{eqnarray}
&&\sum_{\lambda _0=0}^{\infty }\frac{(\gamma )_{\lambda _0}}{ \lambda_0 !} w^{\lambda_0} \;_2F_1(-\lambda_0, \lambda_0+A; \gamma; x) \label{eq:60038}\\
&&= 2^{A -1}\frac{\left(1-w+\sqrt{w^2-2(1-2x)w+1}\right)^{1-\gamma } \left(1+w+\sqrt{w^2-2(1-2x)w+1}\right)^{\gamma -A}}{\sqrt{w^2-2(1-2x)w+1}} \nonumber\\
&& \hspace{.5cm} \mbox{where}\;|w|<1 \nonumber
\end{eqnarray}
\end{lemma}
\begin{proof}
Jacobi polynomial $P_n^{(\alpha, \beta )}(x)$ can be written in terms of hypergeometric function using
\begin{equation}
_2F_1 (-n, n+\alpha +\beta +1; \alpha +1;x) = \frac{n!}{(\alpha +1)_n}P_n^{(\alpha, \beta )}(1-2x) \label{eq:60039}
\end{equation}
And
\begin{equation}
P_n^{(\alpha, \beta )}(x) = \frac{\Gamma (n+\alpha +1)}{n! \Gamma (n+\alpha +\beta +1)} \sum_{m=0}^{n} \binom{n}{m} \frac{\Gamma (n+m+\alpha +\beta +1)}{\Gamma (m+\alpha +1)}\left( \frac{x-1}{2}\right)^m \label{eq:60040}
\end{equation}
The generating function for the Jacobi polynomials is witten by 
\begin{equation}
\sum_{n=0}^{\infty } P_n^{(\alpha, \beta )}(x) w^n = 2^{\alpha +\beta }\frac{\left(1-w+\sqrt{w^2-2xw+1}\right)^{-\alpha  } \left(1+w+\sqrt{w^2-2xw+1}\right)^{-\beta }}{\sqrt{w^2-2xw+1}} \label{eq:60041}
\end{equation}
Replace $n$, $\alpha $ and $\beta $ by $\lambda_0$, $\gamma -1$ and $A-\gamma $ in (\ref{eq:60039}), and acting the summation operator ${\displaystyle \sum_{\lambda_0=0}^{\infty }\frac{(\gamma )_{\lambda_0}}{ \lambda_0 !} w^{\lambda_0} }$ on the new (\ref{eq:60039})
\begin{equation}
\sum_{\lambda _0=0}^{\infty }\frac{(\gamma )_{\lambda_0}}{ \lambda_0 !} w^{\lambda_0}\; _2F_1 (-\lambda_0, \lambda_0+ A; \gamma ;x) = \sum_{\lambda_0=0}^{\infty } P_n^{(\gamma -1, A-\gamma  )}(1-2x) w^{\lambda_0}\label{eq:60042}
\end{equation}
Replace $\alpha $, $\beta $ and $x$ by $\gamma -1$, $A-\gamma $ and $1-2x$ in (\ref{eq:60041}). As we take the new (\ref{eq:60041}) into (\ref{eq:60042}), we obtain (\ref{eq:60038}).\qed
\end{proof}

\begin{definition}
I define that
\begin{equation}
\begin{cases}
\displaystyle { s_{a,b}} = \begin{cases} \displaystyle {  s_a\cdot s_{a+1}\cdot s_{a+2}\cdots s_{b-2}\cdot s_{b-1}\cdot s_b}\;\;\mbox{if}\;a>b \cr
s_a \;\;\mbox{if}\;a=b\end{cases}
\cr
\cr
\displaystyle { \widetilde{w}_{i,j}}  = 
\begin{cases} \displaystyle { \frac{ \widetilde{w}_{i+1,j}\; t_i u_i \left\{ 1+ (s_i+2\widetilde{w}_{i+1,j}(1-t_i)(1-u_i))s_i\right\}}{2(1-\widetilde{w}_{i+1,j}(1-t_i)(1-u_i))^2 s_i}} \cr
\displaystyle {-\frac{\widetilde{w}_{i+1,j}\; t_i u_i (1+s_i)\sqrt{s_i^2-2(1-2\widetilde{w}_{i+1,j}(1-t_i)(1-u_i))s_i+1}}{2(1-\widetilde{w}_{i+1,j}(1-t_i)(1-u_i))^2 s_i}} \;\;\mbox{where}\;i<j \cr
\cr
\displaystyle { \frac{x t_i u_i \left\{ 1+ (s_{i,\infty }+2x(1-t_i)(1-u_i))s_{i,\infty }\right\}}{2(1-x(1-t_i)(1-u_i))^2 s_{i,\infty }}} \cr
\displaystyle {-\frac{x t_i u_i(1+s_{i,\infty })\sqrt{s_{i,\infty }^2-2(1-2x (1-t_i)(1-u_i))s_{i,\infty }+1}}{2(1-x (1-t_i)(1-u_i))^2 s_{i,\infty }}} \;\;\mbox{where}\;i=j 
\end{cases}
\end{cases}\label{eq:60043}
\end{equation}
where
\begin{equation}
a,b,i,j\in \mathbb{N}_{0} \nonumber
\end{equation}
\end{definition}
And we have
\begin{equation}
\sum_{\lambda _i = \lambda_j}^{\infty } s_i^{\lambda_i} = \frac{s_i^{\lambda_j}}{(1-s_i)}\label{eq:60044}
\end{equation}
Acting the summation operator $\displaystyle{ \sum_{\lambda _0 =0}^{\infty } \frac{(\gamma')_{\lambda _0}}{ \lambda _0 !} s_0^{\lambda _0} \prod _{n=1}^{\infty } \left\{ \sum_{ \lambda _n = \lambda _{n-1}}^{\infty } s_n^{\lambda _n }\right\}}$ on (\ref{eq:60026}) where $|s_i|<1$ as $i=0,1,2,\cdots$ by using (\ref{eq:60043}) and (\ref{eq:60044}),
\begin{theorem} 
The general expression of the generating function for the Mathieu polynomial of type 2 is given by
\begin{eqnarray}
&&\sum_{\lambda _0 =0}^{\infty } \frac{(\gamma')_{\lambda _0}}{ \lambda _0 !} s_0^{\lambda _0} \prod _{n=1}^{\infty } \left\{ \sum_{ \lambda _n = \lambda _{n-1}}^{\infty } s_n^{\lambda _n }\right\} y(x) \nonumber\\
&&= \prod_{l=1}^{\infty } \frac{1}{(1-s_{l,\infty })} \mathbf{\Upsilon}(\nu; s_{0,\infty } ;x) \nonumber\\
&&+ \Bigg\{ \prod_{l=2}^{\infty } \frac{1}{(1-s_{l,\infty })} \int_{0}^{1} dt_1\;t_1^{1+\nu} \int_{0}^{1} du_1\;u_1^{\frac{1}{2} +\nu} \overleftrightarrow {\mathbf{\Xi}}_1 \left(\nu; s_{1,\infty };t_1,u_1,x\right)  \mathbf{\Upsilon}(\nu ; s_0;\widetilde{w}_{1,1})\Bigg\} \eta\nonumber\\
&&+ \sum_{n=2}^{\infty } \Bigg\{ \prod_{l=n+1}^{\infty } \frac{1}{(1-s_{l,\infty })} \int_{0}^{1} dt_n\;t_n^{2n-1+\nu } \int_{0}^{1} du_n\;u_n^{2n -\frac{3}{2}+\nu} \overleftrightarrow {\mathbf{\Xi}}_n \left(\nu; s_{n,\infty };t_n,u_n,x \right) \nonumber\\
&&\times \prod_{k=1}^{n-1} \left\{ \int_{0}^{1} dt_{n-k}\;t_{n-k}^{2(n-k)-1+\nu } \int_{0}^{1} du_{n-k} \;u_{n-k}^{2(n-k)-\frac{3}{2} +\nu } \overleftrightarrow {\mathbf{\Xi}}_{n-k} \left(\nu; s_{n-k};t_{n-k},u_{n-k},\widetilde{w}_{n-k+1,n} \right) \right\}\nonumber\\
&&\times   \mathbf{\Upsilon}(\nu ; s_0;\widetilde{w}_{1,n}) \Bigg\}\; \eta^n\label{eq:60045}
\end{eqnarray}
where
\begin{equation}
\begin{cases} 
{ \displaystyle \overleftrightarrow {\mathbf{\Xi}}_1 \left(\nu; s_{1,\infty };t_1,u_1,x\right) = 
\frac{\left( \frac{ 1+s_{1,\infty } +\sqrt{s_{1,\infty }^2-2(1-2x(1-t_1)(1-u_1))s_{1,\infty }+1}}{2}\right)^{-\left( 3+2\nu  \right)}}{\sqrt{s_{1,\infty }^2-2(1-2 x(1-t_1)(1-u_1))s_{1,\infty }+1}  }}\cr
{ \displaystyle \overleftrightarrow {\mathbf{\Xi}}_n \left(\nu; s_{n,\infty };t_n,u_n,x \right) = 
\frac{\left( \frac{ 1+s_{n,\infty } +\sqrt{s_{n,\infty }^2-2(1-2x(1-t_n)(1-u_n))s_{n,\infty }+1}}{2}\right)^{-\left( 4n-1+2\nu \right)}}{\sqrt{s_{n,\infty }^2-2(1-2x(1-t_n)(1-u_n))s_{n,\infty }+1}} }\cr
{ \displaystyle \overleftrightarrow {\mathbf{\Xi}}_{n-k} \left(\nu; s_{n-k};t_{n-k},u_{n-k},\widetilde{w}_{n-k+1,n} \right)
 } \cr
{ \displaystyle = \frac{\left( \frac{ 1+s_{n-k} +\sqrt{s_{n-k}^2-2(1-2\widetilde{w}_{n-k+1,n} (1-t_{n-k})(1-u_{n-k}))s_{n-k}+1}}{2}\right)^{-\left( 4(n-k)-1 +2\nu \right)}}{\sqrt{s_{n-k}^2-2(1-2\widetilde{w}_{n-k+1,n} (1-t_{n-k})(1-u_{n-k}))s_{n-k}+1}}}
\end{cases}\nonumber 
\end{equation}
and
\begin{equation}
\begin{cases} 
{ \displaystyle \mathbf{\Upsilon}(\nu; s_{0,\infty } ;x)= \sum_{\lambda _0 =0}^{\infty } \frac{(\gamma')_{\lambda_0}}{ \lambda_0 !} s_{0,\infty }^{\lambda_0} \left( c_0 x^{\nu } \sum_{i_0=0}^{\lambda_0} \frac{(-\lambda_0)_{i_0} \left( \lambda_0+  2 \nu \right)_{i_0}}{\left(1+\nu \right)_{i_0} \left(\frac{1}{2} +\nu \right)_{i_0}} x^{i_0} \right) }\cr
{ \displaystyle \mathbf{\Upsilon}(\nu ; s_0;\widetilde{w}_{1,1}) = \sum_{\lambda_0 =0}^{\infty } \frac{(\gamma')_{\lambda_0}}{ \lambda_0 !} s_0^{\lambda_0}\left( c_0 x^{\nu} \sum_{i_0=0}^{\lambda_0} \frac{(-\lambda_0)_{i_0} \left( \lambda_0 +2 \nu \right)_{i_0}}{\left( 1+\nu \right)_{i_0}\left(\frac{1}{2} +\nu \right)_{i_0}} \widetilde{w}_{1,1} ^{i_0} \right) }\cr
{ \displaystyle \mathbf{\Upsilon}(\nu; s_0 ;\widetilde{w}_{1,n}) = \sum_{\lambda_0 =0}^{\infty } \frac{(\gamma')_{\lambda_0}}{ \lambda_0 !} s_0^{\lambda_0}\left( c_0 x^{\nu} \sum_{i_0=0}^{\lambda_0} \frac{(-\lambda_0)_{i_0} \left( \lambda_0 +2 \nu \right)_{i_0}}{\left( 1+\nu \right)_{i_0}\left(\frac{1}{2} +\nu \right)_{i_0}} \widetilde{w}_{1,n} ^{i_0} \right)}
\end{cases}\nonumber 
\end{equation}
\end{theorem}
\begin{proof} 
Acting the summation operator $\displaystyle{ \sum_{\lambda _0 =0}^{\infty } \frac{(\gamma')_{\lambda _0}}{ \lambda _0 !} s_0^{\lambda _0} \prod _{n=1}^{\infty } \left\{ \sum_{ \lambda _n = \lambda _{n-1}}^{\infty } s_n^{\lambda _n }\right\}}$ on the form of integral of the Mathieu polynomial of type 2 $y(x)$,
\begin{eqnarray}
&&\sum_{\lambda _0 =0}^{\infty } \frac{(\gamma')_{\lambda _0}}{ \lambda _0 !} s_0^{\lambda _0} \prod _{n=1}^{\infty } \left\{ \sum_{ \lambda _n = \lambda _{n-1}}^{\infty } s_n^{\lambda _n }\right\} y(x) \label{eq:60046}\\
&&= \sum_{\lambda _0 =0}^{\infty } \frac{(\gamma')_{\lambda _0}}{ \lambda _0 !} s_0^{\lambda _0} \prod _{n=1}^{\infty } \left\{ \sum_{ \lambda _n = \lambda _{n-1}}^{\infty } s_n^{\lambda _n }\right\} \Big\{ y_0(x)+y_1(x)+y_2(x)+y_3(x)+\cdots\Big\} \nonumber
\end{eqnarray}
Acting the summation operator $\displaystyle{ \sum_{\lambda _0 =0}^{\infty } \frac{(\gamma')_{\lambda _0}}{ \lambda _0 !} s_0^{\lambda _0} \prod _{n=1}^{\infty } \left\{ \sum_{ \lambda _n = \lambda _{n-1}}^{\infty } s_n^{\lambda _n }\right\}}$ on (\ref{eq:60028a}),
\begin{eqnarray}
&&\sum_{\lambda _0 =0}^{\infty } \frac{(\gamma')_{\lambda _0}}{ \lambda _0 !} s_0^{\lambda _0} \prod _{n=1}^{\infty } \left\{ \sum_{ \lambda _n = \lambda _{n-1}}^{\infty } s_n^{\lambda _n }\right\} y_0(x) \nonumber\\
&&= \prod_{l=1}^{\infty } \frac{1}{(1-s_{l,\infty })} \sum_{\lambda_0 =0}^{\infty } \frac{(\gamma')_{\lambda_0}}{ \lambda_0 !} s_{0,\infty }^{\lambda_0} \left( c_0 x^{\nu } \sum_{i_0=0}^{\lambda _0} \frac{(-\lambda_0)_{i_0} \left( \lambda _0+ 2\nu \right)_{i_0}}{(1+\nu )_{i_0}\left( \frac{1}{2}+ \nu \right)_{i_0}} x^{i_0} \right) \hspace{2cm}\label{eq:60047}
\end{eqnarray}
Acting the summation operator $\displaystyle{ \sum_{\lambda _0 =0}^{\infty } \frac{(\gamma')_{\lambda _0}}{ \lambda _0 !} s_0^{\lambda _0} \prod _{n=1}^{\infty } \left\{ \sum_{ \lambda _n = \lambda _{n-1}}^{\infty } s_n^{\lambda _n }\right\}}$ on (\ref{eq:60029}),
\begin{eqnarray}
&&\sum_{\lambda _0 =0}^{\infty } \frac{(\gamma')_{\lambda _0}}{ \lambda _0 !} s_0^{\lambda _0} \prod _{n=1}^{\infty } \left\{ \sum_{ \lambda _n = \lambda _{n-1}}^{\infty } s_n^{\lambda _n }\right\} y_1(x) \nonumber\\
&&= \prod_{l=2}^{\infty } \frac{1}{(1-s_{l,\infty })} \int_{0}^{1} dt_1\;t_1^{1+\nu } \int_{0}^{1} du_1\;u_1^{\frac{1}{2} +\nu}
 \frac{1}{2\pi i} \oint dp_1 \;\frac{1}{p_1} (1-x(1-t_1)(1-u_1)p_1)^{-\left( 4 +2 \nu \right)} \nonumber\\
&&\times \sum_{\lambda _1 =\lambda_0}^{\infty }\left( \frac{p_1-1}{p_1} \frac{s_{1,\infty }}{1-x(1-t_1)(1-u_1)p_1}\right)^{\lambda _1}     \sum_{\lambda_0 =0}^{\infty } \frac{(\gamma')_{\lambda_0}}{ \lambda_0 !} s_0^{\lambda_0} \left( c_0 x^{\nu } \sum_{i_0=0}^{\lambda _0} \frac{(-\lambda_0)_{i_0} \left( \lambda _0+ 2\nu \right)_{i_0}}{(1+\nu )_{i_0}\left( \frac{1}{2}+ \nu \right)_{i_0}} w_{1,1}^{i_0} \right)\eta \hspace{1.5cm}  \label{eq:60048}
\end{eqnarray}
Replace $\lambda _i$, $\lambda_j$ and $s_i$ by $\lambda_1$, $\lambda_0$ and ${ \displaystyle \frac{p_1-1}{p_1} \frac{s_{1,\infty }}{1-x(1-t_1)(1-u_1)p_1}}$ in (\ref{eq:60044}). Take the new (\ref{eq:60044}) into (\ref{eq:60048}).
\begin{eqnarray}
&&\sum_{\lambda _0 =0}^{\infty } \frac{(\gamma')_{\lambda _0}}{ \lambda _0 !} s_0^{\lambda _0} \prod _{n=1}^{\infty } \left\{ \sum_{ \lambda _n = \lambda _{n-1}}^{\infty } s_n^{\lambda _n }\right\} y_1(x) \nonumber\\
&&= \prod_{l=2}^{\infty } \frac{1}{(1-s_{l,\infty })} \int_{0}^{1} dt_1\;t_1^{1+\nu } \int_{0}^{1} du_1\;u_1^{\frac{1}{2} +\nu}
 \frac{1}{2\pi i} \oint dp_1 \; \frac{(1-x(1-t_1)(1-u_1)p_1)^{-\left( 3 +2 \nu \right)}}{-x(1-t_1)(1-u_1)p_1^2+(1-s_{1,\infty })p_1+s_{1,\infty }} \nonumber\\
&&\times  \sum_{\lambda_0 =0}^{\infty } \left( \frac{p_1-1}{p_1} \frac{s_{0,\infty }}{1-x(1-t_1)(1-u_1)p_1}\right)^{\lambda _0}     \frac{(\gamma')_{\lambda_0}}{ \lambda_0 !} \left( c_0 x^{\nu } \sum_{i_0=0}^{\lambda _0} \frac{(-\lambda_0)_{i_0} \left( \lambda _0+ 2\nu \right)_{i_0}}{(1+\nu )_{i_0}\left( \frac{1}{2}+ \nu \right)_{i_0}} w_{1,1}^{i_0} \right)\eta\hspace{1cm}  \label{eq:60049}
\end{eqnarray}
By using Cauchy's integral formula, the contour integrand has poles at\\
 ${\displaystyle p_1= \frac{1-s_{1,\infty }-\sqrt{(1-s_{1,\infty })^2+4x(1-t_1)(1-u_1)s_{1,\infty }}}{2x(1-t_1)(1-u_1)} }$\\  or ${\displaystyle \frac{1-s_{1,\infty }+\sqrt{(1-s_{1,\infty })^2+4x(1-t_1)(1-u_1)s_{1,\infty }}}{2x(1-t_1)(1-u_1)} }$ and ${ \displaystyle \frac{1-s_{1,\infty }-\sqrt{(1-s_{1,\infty })^2+4x(1-t_1)(1-u_1)s_{1,\infty }}}{2x(1-t_1)(1-u_1)}}$ is only inside the unit circle. As we compute the residue there in (\ref{eq:60049}) we obtain
\begin{eqnarray}
&&\sum_{\lambda _0 =0}^{\infty } \frac{(\gamma')_{\lambda _0}}{ \lambda _0 !} s_0^{\lambda _0} \prod _{n=1}^{\infty } \left\{ \sum_{ \lambda _n = \lambda _{n-1}}^{\infty } s_n^{\lambda _n }\right\} y_1(x) \nonumber\\
&&= \prod_{l=2}^{\infty } \frac{1}{(1-s_{l,\infty })} \int_{0}^{1} dt_1\;t_1^{1+\nu} \int_{0}^{1} du_1\;u_1^{\frac{1}{2} +\nu}
 \left( s_{1,\infty }^2-2(1-2x(1-t_1)(1-u_1))s_{1,\infty }+1 \right)^{-\frac{1}{2}} \nonumber\\
&&\times \left(\frac{1+s_{1,\infty }+\sqrt{s_{1,\infty }^2-2(1-2x(1-t_1)(1-u_1))s_{1,\infty }+1}}{2}\right)^{-\left(3+ 2\nu \right)} \nonumber \\
&&\times \sum_{\lambda _0 =0}^{\infty } \frac{(\gamma' )_{\lambda_0}}{ \lambda_0 !} s_0^{\lambda_0}\left( c_0 x^{\nu } \sum_{i_0=0}^{\lambda_0} \frac{(-\lambda_0)_{i_0} \left(\lambda_0 +2\nu \right)_{i_0}}{(1+\nu )_{i_0}\left(\frac{1}{2} +\nu \right)_{i_0}} \widetilde{w}_{1,1} ^{i_0} \right) \eta \label{eq:60050}
\end{eqnarray}
where
\begin{eqnarray}
\widetilde{w}_{1,1} &=& \frac{p_1}{(p_1-1)}\; \frac{x t_1 u_1}{1- x (1-t_1)(1-u_1) p_1}\Bigg|_{\Large p_1=\frac{1-s_{1,\infty }-\sqrt{(1-s_{1,\infty })^2+4x(1-t_1)(1-u_1)s_{1,\infty }}}{2x(1-t_1)(1-u_1)}\normalsize}\nonumber\\
&=& \frac{x t_1 u_1 \left\{ 1+ (s_{1,\infty }+2x(1-t_1)(1-u_1) )s_{1,\infty }\right\}}{2(1-x(1-t_1)(1-u_1))^2 s_{1,\infty }}\nonumber\\
&&-\frac{xt_1 u_1(1+s_{1,\infty })\sqrt{s_{1,\infty }^2-2(1-2x(1-t_1)(1-u_1))s_{1,\infty }+1}}{2(1-x(1-t_1)(1-u_1))^2 s_{1,\infty }}\nonumber
\end{eqnarray}
Acting the summation operator $\displaystyle{ \sum_{\lambda _0 =0}^{\infty } \frac{(\gamma')_{\lambda _0}}{ \lambda _0 !} s_0^{\lambda _0} \prod _{n=1}^{\infty } \left\{ \sum_{ \lambda _n = \lambda _{n-1}}^{\infty } s_n^{\lambda _n }\right\}}$ on (\ref{eq:60031}),
\begin{eqnarray}
&&\sum_{\lambda _0 =0}^{\infty } \frac{(\gamma')_{\lambda _0}}{ \lambda _0 !} s_0^{\lambda _0} \prod _{n=1}^{\infty } \left\{ \sum_{ \lambda _n = \lambda _{n-1}}^{\infty } s_n^{\lambda _n }\right\} y_2(x) \nonumber\\
&&= \prod_{l=3}^{\infty } \frac{1}{(1-s_{l,\infty })} \int_{0}^{1} dt_2\;t_2^{3+\nu} \int_{0}^{1} du_2\;u_2^{\frac{5}{2} +\nu }
 \frac{1}{2\pi i} \oint dp_2 \;\frac{1}{p_2} (1-x(1-t_2)(1-u_2)p_2)^{-\left(8 +2 \nu \right)} \nonumber\\
&&\times \sum_{\lambda _2 =\lambda _1}^{\infty }\left( \frac{p_2-1}{p_2} \frac{s_{2,\infty }}{1-x(1-t_2)(1-u_2)p_2}\right)^{\lambda _2}  \nonumber\\
&&\times \int_{0}^{1} dt_1\;t_1^{1+\nu} \int_{0}^{1} du_1\;u_1^{\frac{1}{2} +\nu}
 \frac{1}{2\pi i} \oint dp_1 \;\frac{1}{p_1} (1-w_{2,2} (1-t_1)(1-u_1)p_1)^{-\left( 4 +2 \nu  \right)} \label{eq:60051}\\
&&\times \sum_{\lambda _1 =\lambda _0}^{\infty }\left( \frac{p_1-1}{p_1} \frac{s_1}{1-w_{2,2}(1-t_1)(1-u_1)p_1}\right)^{\lambda _1}     \sum_{\lambda _0 =0}^{\infty } \frac{(\gamma' )_{\lambda _0}}{ \lambda _0 !} s_0^{\lambda _0}\left( c_0 x^{\nu } \sum_{i_0=0}^{\lambda _0} \frac{(-\lambda _0)_{i_0} \left( \lambda _0 +2 \nu \right)_{i_0}}{(1+\nu )_{i_0}\left(\frac{1}{2} +\nu \right)_{i_0}} w_{1,2} ^{i_0} \right) \eta^2  \nonumber
\end{eqnarray}
Replace $\lambda _i$, $\lambda _j$ and $s_i$ by $\lambda _2$, $\lambda _1$ and ${ \displaystyle \frac{p_2-1}{p_2} \frac{s_{2,\infty }}{1-x(1-t_2)(1-u_2)p_2}}$ in (\ref{eq:60044}). Take the new (\ref{eq:60044}) into (\ref{eq:60051}).
\begin{eqnarray}
&&\sum_{\lambda _0 =0}^{\infty } \frac{(\gamma')_{\lambda _0}}{ \lambda _0 !} s_0^{\lambda _0} \prod _{n=1}^{\infty } \left\{ \sum_{ \lambda _n = \lambda _{n-1}}^{\infty } s_n^{\lambda _n }\right\} y_2(x) \nonumber\\
&&= \prod_{l=3}^{\infty } \frac{1}{(1-s_{l,\infty })} \int_{0}^{1} dt_2\;t_2^{3+\nu } \int_{0}^{1} du_2\;u_2^{\frac{5}{2} +\nu }
 \frac{1}{2\pi i} \oint dp_2 \;\frac{\left(1-x (1-t_2)(1-u_2)p_2\right)^{-\left(7 +2 \nu \right)}}{-x (1-t_2)(1-u_2)p_2^2+ (1-s_{2,\infty })p_2+s_{2,\infty } } \nonumber\\
&&\times \int_{0}^{1} dt_1\;t_1^{1+\nu} \int_{0}^{1} du_1\;u_1^{\frac{1}{2} +\nu}
 \frac{1}{2\pi i} \oint dp_1 \;\frac{1}{p_1} \left( 1-w_{2,2} (1-t_1)(1-u_1)p_1\right)^{-\left( 4+2\nu \right)} \nonumber\\
&&\times \sum_{\lambda _1 = \lambda_0}^{\infty }\left( \frac{p_2-1}{p_2} \frac{s_{1,\infty }}{1-x(1-t_2)(1-u_2)p_2} \frac{p_1-1}{p_1}\frac{1}{1-w_{2,2}(1-t_1)(1-u_1)p_1}\right)^{\lambda_1} \nonumber\\
&&\times  \sum_{\lambda_0 =0}^{\infty } \frac{(\gamma' )_{\lambda_0}}{ \lambda_0 !} s_0^{\lambda_0} \left( c_0 x^{\nu } \sum_{i_0=0}^{\lambda_0} \frac{(-\lambda_0)_{i_0} \left( \lambda_0 +2\nu \right)_{i_0}}{(1+\nu )_{i_0}\left( \frac{1}{2} +\nu \right)_{i_0}} w_{1,2} ^{i_0} \right) \eta^2 \label{eq:60052}
\end{eqnarray}
By using Cauchy's integral formula, the contour integrand has poles at\\
 ${\displaystyle
 p_2= \frac{1-s_{2,\infty }-\sqrt{(1-s_{2,\infty })^2+4x(1-t_2)(1-u_2)s_{2,\infty }}}{2x(1-t_2)(1-u_2)}}$ \\or${\displaystyle\frac{1-s_{2,\infty }+\sqrt{(1-s_{2,\infty })^2+4x(1-t_2)(1-u_2)s_{2,\infty }}}{2x(1-t_2)(1-u_2)} }$ 
and ${ \displaystyle\frac{1-s_{2,\infty }-\sqrt{(1-s_{2,\infty })^2+4x(1-t_2)(1-u_2)s_{2,\infty }}}{2x(1-t_2)(1-u_2)}}$ is only inside the unit circle. As we compute the residue there in (\ref{eq:60052}) we obtain
\begin{eqnarray}
&&\sum_{\lambda _0 =0}^{\infty } \frac{(\gamma')_{\lambda _0}}{ \lambda _0 !} s_0^{\lambda _0} \prod _{n=1}^{\infty } \left\{ \sum_{ \lambda _n = \lambda _{n-1}}^{\infty } s_n^{\lambda _n }\right\} y_2(x) \nonumber\\
&&= \prod_{l=3}^{\infty } \frac{1}{(1-s_{l,\infty })} \int_{0}^{1} dt_2\;t_2^{3+\nu} \int_{0}^{1} du_2\;u_2^{\frac{5}{2} +\nu }
 \left( s_{2,\infty }^2-2(1-2x(1-t_2)(1-u_2))s_{2,\infty }+1\right)^{-\frac{1}{2}}\nonumber\\
&&\times \left(\frac{1+s_{2,\infty }+\sqrt{s_{2,\infty }^2-2(1-2x(1-t_2)(1-u_2))s_{2,\infty }+1}}{2}\right)^{-\left(7 +2\nu \right)}   \nonumber\\
&&\times \int_{0}^{1} dt_1\;t_1^{1+\nu} \int_{0}^{1} du_1\;u_1^{\frac{1}{2} +\nu}
 \frac{1}{2\pi i} \oint dp_1 \;\frac{1}{p_1} \left( 1-\widetilde{w}_{2,2} (1-t_1)(1-u_1)p_1\right)^{-\left( 4 +2 \nu \right)} \label{eq:60053}\\
&&\times \sum_{\lambda _1 = \lambda_0}^{\infty }\left( \frac{p_1-1}{p_1}\frac{s_1}{1-\widetilde{w}_{2,2}(1-t_1)(1-u_1)p_1}\right)^{\lambda _1}
   \sum_{\lambda_0 =0}^{\infty } \frac{(\gamma' )_{\lambda_0}}{ \lambda_0 !} s_0^{\lambda_0}\left( c_0 x^{\nu } \sum_{i_0=0}^{\lambda_0} \frac{(-\lambda_0)_{i_0} \left( \lambda_0 +2 \nu \right)_{i_0}}{(1+\nu )_{i_0}\left(\frac{1}{2} +\nu \right)_{i_0}} \ddot{w}_{1,2} ^{i_0} \right) \eta^2 \nonumber
\end{eqnarray}
where
\begin{eqnarray}
\widetilde{w}_{2,2} &=& \frac{p_2}{(p_2-1)}\; \frac{x t_2 u_2}{1- x (1-t_2)(1-u_2) p_2}\Bigg|_{\Large p_2=\frac{1-s_{2,\infty }-\sqrt{(1-s_{2,\infty })^2+4x (1-t_2)(1-u_2)s_{2,\infty }}}{2x(1-t_2)(1-u_2)}\normalsize}\nonumber\\
&=& \frac{x t_2 u_2 \left\{1+ (s_{2,\infty }+2x(1-t_2)(1-u_2) )s_{2,\infty }\right\}  }{2(1-x(1-t_2)(1-u_2))^2 s_{2,\infty }} \nonumber\\
&&- \frac{x t_2 u_2 (1+s_{2,\infty })\sqrt{s_{2,\infty }^2-2(1-2x (1-t_2)(1-u_2))s_{2,\infty }+1}}{2(1-x(1-t_2)(1-u_2))^2 s_{2,\infty }}
\nonumber
\end{eqnarray}
and
\begin{equation}
\ddot{w}_{1,2} = \frac{p_1}{(p_1-1)}\; \frac{\widetilde{w}_{2,2} t_1 u_1}{1- \widetilde{w}_{2,2}(1-t_1)(1-u_1)p_1 }\nonumber
\end{equation}
Replace $\lambda _i$, $\lambda_j$ and $s_i$ by $\lambda_1$, $\lambda_0$ and ${ \displaystyle \frac{p_1-1}{p_1}\frac{s_1}{1-\widetilde{w}_{2,2}(1-t_1)(1-u_1)p_1}}$ in (\ref{eq:60044}). Take the new (\ref{eq:60044}) into (\ref{eq:60053}).
\begin{eqnarray}
&&\sum_{\lambda _0 =0}^{\infty } \frac{(\gamma')_{\lambda _0}}{ \lambda _0 !} s_0^{\lambda _0} \prod _{n=1}^{\infty } \left\{ \sum_{ \lambda _n = \lambda _{n-1}}^{\infty } s_n^{\lambda _n }\right\} y_2(x) \nonumber\\
&&= \prod_{l=3}^{\infty } \frac{1}{(1-s_{l,\infty })} \int_{0}^{1} dt_2\;t_2^{3+\nu} \int_{0}^{1} du_2\;u_2^{\frac{5}{2} +\nu }
 \left( s_{2,\infty }^2-2(1-2x(1-t_2)(1-u_2))s_{2,\infty }+1\right)^{-\frac{1}{2}}\nonumber\\
&&\times \left(\frac{1+s_{2,\infty }+\sqrt{s_{2,\infty }^2-2(1-2x(1-t_2)(1-u_2))s_{2,\infty }+1}}{2}\right)^{-\left( 7 +2\nu \right)}  \nonumber\\
&&\times \int_{0}^{1} dt_1\;t_1^{1+\nu} \int_{0}^{1} du_1\;u_1^{\frac{1}{2} +\nu}
 \frac{1}{2\pi i} \oint dp_1 \;\frac{\left(1-\widetilde{w}_{2,2} (1-t_1)(1-u_1)p_1\right)^{-\left(3 +2 \nu \right)}}{-\widetilde{w}_{2,2} (1-t_1)(1-u_1)p_1^2+(1-s_1)p_1+s_1} \nonumber\\
&&\times   \sum_{\lambda _0 =0}^{\infty } \frac{(\gamma')_{\lambda_0}}{ \lambda_0 !} \left( \frac{p_1-1}{p_1}\frac{s_{0,1}}{1-\widetilde{w}_{2,2}(1-t_1)(1-u_1)p_1}\right)^{\lambda_0} \left( c_0 x^{\nu } \sum_{i_0=0}^{\lambda_0} \frac{(-\lambda_0)_{i_0} \left( \lambda_0 +2 \nu  \right)_{i_0}}{(1+\nu )_{i_0}\left(\frac{1}{2} +\lambda \right)_{i_0}} \ddot{w}_{1,2} ^{i_0} \right) \eta^2  \hspace{1.5cm}\label{eq:60054}
\end{eqnarray}
By using Cauchy's integral formula, the contour integrand has poles at\\ ${\displaystyle
 p_1= \frac{1-s_1-\sqrt{(1-s_1)^2+4\widetilde{w}_{2,2} (1-t_1)(1-u_1)s_1}}{2\widetilde{w}_{2,2} (1-t_1)(1-u_1)}}$\\ or ${\displaystyle\frac{1-s_1+\sqrt{(1-s_1)^2+4\widetilde{w}_{2,2} (1-t_1)(1-u_1)s_1}}{2\widetilde{w}_{2,2} (1-t_1)(1-u_1)}}$ 
and ${ \displaystyle \frac{1-s_1-\sqrt{(1-s_1)^2+4\widetilde{w}_{2,2} (1-t_1)(1-u_1)s_1}}{2\widetilde{w}_{2,2} (1-t_1)(1-u_1)}}$ is only inside the unit circle. As we compute the residue there in (\ref{eq:60054}) we obtain
\begin{eqnarray}
&&\sum_{\lambda _0 =0}^{\infty } \frac{(\gamma')_{\lambda _0}}{ \lambda _0 !} s_0^{\lambda _0} \prod _{n=1}^{\infty } \left\{ \sum_{ \lambda _n = \lambda _{n-1}}^{\infty } s_n^{\lambda _n }\right\} y_2(x) \nonumber\\
&&= \prod_{l=3}^{\infty } \frac{1}{(1-s_{l,\infty })} \int_{0}^{1} dt_2\;t_2^{3+\nu} \int_{0}^{1} du_2\;u_2^{\frac{5}{2} +\nu }
\left( s_{2,\infty }^2-2(1-2x(1-t_2)(1-u_2))s_{2,\infty }+1\right)^{-\frac{1}{2}}\nonumber\\
&&\times \left(\frac{1+s_{2,\infty }+\sqrt{s_{2,\infty }^2-2(1-2x(1-t_2)(1-u_2))s_{2,\infty }+1}}{2}\right)^{-\left(7 +2 \nu  \right)} \nonumber\\
&&\times \int_{0}^{1} dt_1\;t_1^{1+\nu } \int_{0}^{1} du_1\;u_1^{\frac{1}{2} +\nu}
 \left( s_1^2-2(1-2\widetilde{w}_{2,2} (1-t_1)(1-u_1))s_1+1\right)^{-\frac{1}{2}}\nonumber\\
&&\times \left(\frac{1+s_1+\sqrt{s_1^2-2(1-2\widetilde{w}_{2,2}(1-t_1)(1-u_1))s_1+1}}{2}\right)^{-\left(3 +2 \nu \right)} \nonumber\\
&&\times \sum_{\lambda _0 =0}^{\infty } \frac{(\gamma' )_{\lambda_0}}{ \lambda_0 !} s_0^{\lambda_0}\left( c_0 x^{\nu } \sum_{i_0=0}^{\lambda_0} \frac{(-\lambda_0)_{i_0} \left( \lambda_0  +2 \nu \right)_{i_0}}{(1+\nu )_{i_0}\left(\frac{1}{2} +\nu \right)_{i_0}} \widetilde{w}_{1,2} ^{i_0} \right) \eta^2  \label{eq:60055}
\end{eqnarray}
where
\begin{eqnarray}
\widetilde{w}_{1,2} &=& \frac{p_1}{(p_1-1)}\; \frac{\widetilde{w}_{2,2} t_1 u_1}{1- \widetilde{w}_{2,2} (1-t_1)(1-u_1) p_1}\Bigg|_{\Large p_1=\frac{1-s_1-\sqrt{(1-s_1)^2+4\widetilde{w}_{2,2} (1-t_1)(1-u_1)s_1}}{2\widetilde{w}_{2,2} (1-t_1)(1-u_1)}\normalsize}\nonumber\\
&=& \frac{\widetilde{w}_{2,2} t_1 u_1 \left\{ 1+ (s_1+2\widetilde{w}_{2,2}(1-t_1)(1-u_1) )s_1-(1+s_1)\sqrt{s_1^2-2(1-2\widetilde{w}_{2,2} (1-t_1)(1-u_1))s_1+1}\right\}}{2(1-\widetilde{w}_{2,2}(1-t_1)(1-u_1))^2 s_1}\nonumber
\end{eqnarray}
Acting the summation operator $\displaystyle{ \sum_{\lambda _0 =0}^{\infty } \frac{(\gamma')_{\lambda _0}}{ \lambda _0 !} s_0^{\lambda _0} \prod _{n=1}^{\infty } \left\{ \sum_{ \lambda _n = \lambda _{n-1}}^{\infty } s_n^{\lambda _n }\right\}}$ on (\ref{eq:60032}),
\begin{eqnarray}
&&\sum_{\lambda _0 =0}^{\infty } \frac{(\gamma')_{\lambda _0}}{ \lambda _0 !} s_0^{\lambda _0} \prod _{n=1}^{\infty } \left\{ \sum_{ \lambda _n = \lambda _{n-1}}^{\infty } s_n^{\lambda _n }\right\} y_3(x) \nonumber\\
&&= \prod_{l=4}^{\infty } \frac{1}{(1-s_{l,\infty })} \int_{0}^{1} dt_3\;t_3^{5+\nu} \int_{0}^{1} du_3\;u_3^{\frac{9}{2} +\nu}\left( s_{3,\infty }^2-2(1-2x(1-t_3)(1-u_3))s_{3,\infty }+1\right)^{-\frac{1}{2}}\nonumber\\
&&\times \left(\frac{1+s_{3,\infty }+\sqrt{s_{3,\infty }^2-2(1-2x(1-t_3)(1-u_3))s_{3,\infty }+1}}{2}\right)^{-\left(11 +2\nu \right)}  \nonumber\\
&&\times \int_{0}^{1} dt_2\;t_2^{3+\nu} \int_{0}^{1} du_2\;u_2^{\frac{5}{2} +\nu }\left(s_2^2-2(1-2\widetilde{w}_{3,3}(1-t_2)(1-u_2))s_2+1\right)^{-\frac{1}{2}}\nonumber\\
&&\times \left(\frac{1+s_2+\sqrt{s_2^2-2(1-2\widetilde{w}_{3,3}(1-t_2)(1-u_2))s_2+1}}{2}\right)^{-\left( 7 +2\nu \right)}   \nonumber\\
&&\times \int_{0}^{1} dt_1\;t_1^{1+\nu } \int_{0}^{1} du_1\;u_1^{\frac{1}{2} +\nu}\left( s_1^2-2(1-2\widetilde{w}_{2,3} (1-t_1)(1-u_1))s_1+1\right)^{-\frac{1}{2}}\nonumber\\
&&\times \left(\frac{1+s_1+\sqrt{s_1^2-2(1-2\widetilde{w}_{2,3}(1-t_1)(1-u_1))s_1+1}}{2}\right)^{-\left(3 +2 \nu \right)}     \nonumber\\
&&\times \sum_{\lambda _0 =0}^{\infty } \frac{(\gamma' )_{\lambda_0}}{ \lambda_0 !} s_0^{\lambda_0}\left( c_0 x^{\nu } \sum_{i_0=0}^{\lambda_0} \frac{(-\lambda_0)_{i_0} \left( \lambda_0 +2 \nu \right)_{i_0}}{(1+\nu )_{i_0}\left( \frac{1}{2} +\nu \right)_{i_0}} \widetilde{w}_{1,3} ^{i_0} \right) \eta^3  \label{eq:60056}
\end{eqnarray}

\vspace{1cm}
where
\begin{eqnarray}
\widetilde{w}_{3,3} &=& \frac{p_3}{(p_3-1)}\; \frac{x t_3 u_3}{1- x(1-t_3)(1-u_3)p_3}\Bigg|_{\Large p_3=\frac{1-s_{3,\infty }-\sqrt{(1-s_{3,\infty })^2+4x(1-t_3)(1-u_3)s_{3,\infty }}}{2x(1-t_3)(1-u_3)}\normalsize}\nonumber\\
&=& \frac{x t_3 u_3 \left\{ 1+ (s_{3,\infty }+2x(1-t_3)(1-u_3) )s_{3,\infty }\right\}}{2(1-x(1-t_3)(1-u_3))^2 s_{3,\infty }} \nonumber\\
&&-\frac{x t_3 u_3 (1+s_{3,\infty })\sqrt{s_{3,\infty }^2-2(1-2x (1-t_3)(1-u_3))s_{3,\infty }+1}}{2(1-x(1-t_3)(1-u_3))^2 s_{3,\infty }} \nonumber
\end{eqnarray}
\begin{eqnarray}
\widetilde{w}_{2,3} &=& \frac{p_2}{(p_2-1)}\; \frac{\widetilde{w}_{3,3} t_2 u_2}{1- \widetilde{w}_{3,3} (1-t_2)(1-u_2)p_2 }\Bigg|_{\Large p_2=\frac{1-s_2-\sqrt{(1-s_2)^2+4\widetilde{w}_{3,3} (1-t_2)(1-u_2)s_2}}{2\widetilde{w}_{3,3} (1-t_2)(1-u_2)}\normalsize}\nonumber\\
&=& \frac{\widetilde{w}_{3,3} t_2 u_2 \left\{ 1+ (s_2+2\widetilde{w}_{3,3}(1-t_2)(1-u_2) )s_2-(1+s_2)\sqrt{s_2^2-2(1-2\widetilde{w}_{3,3} (1-t_2)(1-u_2))s_2+1}\right\}}{2(1-\widetilde{w}_{3,3}(1-t_2)(1-u_2))^2 s_2}\nonumber
\end{eqnarray}
\begin{eqnarray}
\widetilde{w}_{1,3} &=& \frac{p_1}{(p_1-1)}\; \frac{\widetilde{w}_{2,3} t_1 u_1}{1- \widetilde{w}_{2,3} (1-t_1)(1-u_1)p_1 }\Bigg|_{\Large p_1=\frac{1-s_1-\sqrt{(1-s_1)^2+4\widetilde{w}_{2,3} (1-t_1)(1-u_1)s_1}}{2\widetilde{w}_{2,3} (1-t_1)(1-u_1)}\normalsize}\nonumber\\
&=& \frac{\widetilde{w}_{2,3} t_1 u_1 \left\{ 1+ (s_1+2\widetilde{w}_{2,3}(1-t_1)(1-u_1) )s_1-(1+s_1)\sqrt{s_1^2-2(1-2\widetilde{w}_{2,3} (1-t_1)(1-u_1))s_1+1}\right\}}{2(1-\widetilde{w}_{2,3}(1-t_1)(1-u_1))^2 s_1}\nonumber
\end{eqnarray}
By repeating this process for all higher terms of integral forms of sub-summation $y_m(x)$ terms where $m \geq 4$, I obtain every  $\displaystyle{ \sum_{\lambda _0 =0}^{\infty } \frac{(\gamma')_{\lambda _0}}{ \lambda _0 !} s_0^{\lambda _0} \prod _{n=1}^{\infty } \left\{ \sum_{ \lambda _n = \lambda _{n-1}}^{\infty } s_n^{\lambda _n }\right\}}  y_m(x)$ terms. 
Substitute (\ref{eq:60047}), (\ref{eq:60050}), (\ref{eq:60055}), (\ref{eq:60056}) and including all $\displaystyle{ \sum_{\lambda _0 =0}^{\infty } \frac{(\gamma')_{\lambda _0}}{ \lambda _0 !} s_0^{\lambda _0} \prod _{n=1}^{\infty } \left\{ \sum_{ \lambda _n = \lambda _{n-1}}^{\infty } s_n^{\lambda _n }\right\}}  y_m(x)$ terms where $m > 3$ into (\ref{eq:60046}). 
\qed
\end{proof}
\begin{remark}
The generating function for the Mathieu polynomial of type 2 of the first kind about $x=0$ as $\lambda  = 2^2(\lambda _j+2j )^2-2q $ where $j,\lambda _j \in \mathbb{N}_{0}$ is
\begin{eqnarray}
&&\sum_{\lambda _0 =0}^{\infty } \frac{ (\frac{1}{2} )_{\lambda _0}}{ \lambda _0 !} s_0^{\lambda _0} \prod _{n=1}^{\infty } \left\{ \sum_{ \lambda _n = \lambda _{n-1}}^{\infty } s_n^{\lambda _n }\right\} MF_{\lambda _j}^R\left( q,\lambda = 2^2(\lambda_j +2j )^2- 2 q; \eta  = q x^2, x= \mathrm{\cos}^2z \right) \nonumber\\
&&= \frac{1}{2} \Bigg\{ \prod_{l=1}^{\infty } \frac{1}{(1-s_{l,\infty })}  \mathbf{A}\left( s_{0,\infty } ;x\right) \nonumber\\
&&+ \Bigg\{ \prod_{l=2}^{\infty } \frac{1}{(1-s_{l,\infty })} \int_{0}^{1} dt_1\;t_1 \int_{0}^{1} du_1\;u_1^{\frac{1}{2}} \overleftrightarrow {\mathbf{\Gamma}}_1 \left(s_{1,\infty };t_1,u_1,x\right) \mathbf{A}\left( s_{0} ;\widetilde{w}_{1,1}\right)\Bigg\} \eta \nonumber\\
&&+ \sum_{n=2}^{\infty } \Bigg\{ \prod_{l=n+1}^{\infty } \frac{1}{(1-s_{l,\infty })} \int_{0}^{1} dt_n\;t_n^{2n-1} \int_{0}^{1} du_n\;u_n^{ 2n-\frac{3}{2} } \overleftrightarrow {\mathbf{\Gamma}}_n \left(s_{n,\infty };t_n,u_n,x \right) \label{eq:60057}\\
&&\times \prod_{k=1}^{n-1} \Bigg\{ \int_{0}^{1} dt_{n-k}\;t_{n-k}^{2(n-k)-1} \int_{0}^{1} du_{n-k} \;u_{n-k}^{2(n-k)-\frac{3}{2}}\overleftrightarrow {\mathbf{\Gamma}}_{n-k} \left(s_{n-k};t_{n-k},u_{n-k},\widetilde{w}_{n-k+1,n} \right)\Bigg\}  \mathbf{A} \left( s_{0} ;\widetilde{w}_{1,n}\right) \Bigg\} \eta^n \Bigg\} \nonumber
\end{eqnarray}
where
\begin{equation}
\begin{cases} 
{ \displaystyle \overleftrightarrow {\mathbf{\Gamma}}_1 \left(s_{1,\infty };t_1,u_1,x\right)= \frac{\left( \frac{1+s_{1,\infty }+\sqrt{s_{1,\infty }^2-2(1-2x(1-t_1)(1-u_1))s_{1,\infty }+1}}{2}\right)^{-3}}{\sqrt{s_{1,\infty }^2-2(1-2x (1-t_1)(1-u_1))s_{1,\infty }+1}}}\cr
{ \displaystyle  \overleftrightarrow {\mathbf{\Gamma}}_n \left(s_{n,\infty };t_n,u_n,x\right) =\frac{\left( \frac{1+s_{n,\infty }+\sqrt{s_{n,\infty }^2-2(1-2x(1-t_n)(1-u_n))s_{n,\infty }+1}}{2}\right)^{-\left( 4n-1 \right)}}{\sqrt{ s_{n,\infty }^2-2(1-2x(1-t_n)(1-u_n))s_{n,\infty }+1}}}\cr
{ \displaystyle \overleftrightarrow {\mathbf{\Gamma}}_{n-k} \left(s_{n-k};t_{n-k},u_{n-k},\widetilde{w}_{n-k+1,n} \right)= \frac{ \left( \frac{1+s_{n-k}+\sqrt{s_{n-k}^2-2(1-2\widetilde{w}_{n-k+1,n} (1-t_{n-k})(1-u_{n-k}))s_{n-k}+1}}{2}\right)^{-\left( 4(n-k)-1 \right)}}{\sqrt{ s_{n-k}^2-2(1-2\widetilde{w}_{n-k+1,n} (1-t_{n-k})(1-u_{n-k}))s_{n-k}+1}}} 
\end{cases}\nonumber 
\end{equation}
and
\begin{equation}
\begin{cases} 
{ \displaystyle \mathbf{A} \left( s_{0,\infty } ;x\right)= \frac{\left(1- s_{0,\infty }+\sqrt{s_{0,\infty }^2-2(1-2x)s_{0,\infty }+1}\right)^{\frac{1}{2}} \left(1+s_{0,\infty }+\sqrt{s_{0,\infty }^2-2(1-2x )s_{0,\infty }+1}\right)^{\frac{1}{2}}}{\sqrt{s_{0,\infty }^2-2(1-2x)s_{0,\infty }+1}}}\cr
{ \displaystyle  \mathbf{A} \left( s_{0} ;\widetilde{w}_{1,1}\right) = \frac{\left(1- s_0+\sqrt{s_0^2-2(1-2\widetilde{w}_{1,1})s_0+1}\right)^{\frac{1}{2}} \left( 1+s_0+\sqrt{s_0^2-2(1-2\widetilde{w}_{1,1} )s_0+1}\right)^{\frac{1}{2}}}{\sqrt{s_0^2-2(1-2\widetilde{w}_{1,1})s_0+1}}} \cr
{ \displaystyle \mathbf{A} \left( s_{0} ;\widetilde{w}_{1,n}\right) = \frac{\left( 1- s_0+\sqrt{s_0^2-2(1-2\widetilde{w}_{1,n})s_0+1}\right)^{\frac{1}{2}} \left(1+s_0+\sqrt{s_0^2-2(1-2\widetilde{w}_{1,n} )s_0+1}\right)^{\frac{1}{2}}}{\sqrt{s_0^2-2(1-2\widetilde{w}_{1,n})s_0+1}}}
\end{cases}\nonumber 
\end{equation}
\end{remark}
\begin{proof}
Replace $A$, $w$ and $\gamma  $ by 0, $s_{0,\infty } $ and $\frac{1}{2}$ in (\ref{eq:60038}). 
\begin{eqnarray}
&&\sum_{\lambda _0=0}^{\infty }\frac{(\frac{1}{2})_{\lambda_0}}{\lambda_0 !} s_{0,\infty }^{\lambda_0} \;_2F_1\left( -\lambda_0, \lambda_0; \frac{1}{2} ; x\right) \label{eq:60058}\\
&&=  \frac{\left(1- s_{0,\infty }+\sqrt{s_{0,\infty }^2-2(1-2x)s_{0,\infty }+1}\right)^{\frac{1}{2}} \left(1+s_{0,\infty }+\sqrt{s_{0,\infty }^2-2(1-2x)s_{0,\infty }+1}\right)^{\frac{1}{2}}}{\sqrt{s_{0,\infty }^2-2(1-2x)s_{0,\infty }+1}} \nonumber
\end{eqnarray} 
Replace $A$, $w$, $\gamma  $, $x$ by 0, $s_0 $, $\frac{1}{2}$ and $\widetilde{w}_{1,1}$ in (\ref{eq:60038}). 
\begin{eqnarray}
&&\sum_{\lambda _0=0}^{\infty }\frac{(\frac{1}{2})_{\lambda_0}}{\lambda_0 !} s_0^{\lambda_0} \;_2F_1\left( -\lambda_0, \lambda_0; \frac{1}{2} ; \widetilde{w}_{1,1} \right) \label{eq:60059}\\
&&=  \frac{\left(1- s_0+\sqrt{s_0^2-2(1-2\widetilde{w}_{1,1})s_0+1}\right)^{\frac{1}{2}} \left( 1+ s_0+\sqrt{s_0^2-2(1-2\widetilde{w}_{1,1})s_0+1}\right)^{\frac{1}{2}}}{\sqrt{s_0^2-2(1-2\widetilde{w}_{1,1})s_0+1}} \nonumber
\end{eqnarray} 
Replace $A$, $w$, $\gamma  $, $x$ by 0, $s_0 $, $\frac{1}{2}$ and $\widetilde{w}_{1,n}$ in (\ref{eq:60038}). 
\begin{eqnarray}
&&\sum_{\lambda _0=0}^{\infty }\frac{(\frac{1}{2})_{\lambda_0}}{\lambda_0 !} s_0^{\lambda_0} \;_2F_1\left( -\lambda_0, \lambda_0; \frac{1}{2} ; \widetilde{w}_{1,n} \right) \label{eq:60060}\\
&&=  \frac{\left(1- s_0+\sqrt{s_0^2-2(1-2\widetilde{w}_{1,n})s_0+1}\right)^{\frac{1}{2}} \left( 1+ s_0+\sqrt{s_0^2-2(1-2\widetilde{w}_{1,n})s_0+1}\right)^{\frac{1}{2}}}{\sqrt{s_0^2-2(1-2\widetilde{w}_{1,n})s_0+1}} \nonumber
\end{eqnarray} 
Put $c_0$= 1, $\lambda =0$ and $\gamma' =\frac{1}{2} $ in (\ref{eq:60045}). Substitute (\ref{eq:60058}), (\ref{eq:60059}) and (\ref{eq:60060}) into the new (\ref{eq:60045}).\qed
\end{proof}
\begin{remark}
The generating function for the Mathieu polynomial of type 2 of the second kind about $x=0$  as $\lambda  = 2^2(\lambda _j+2j +1/2)^2-2q $ where $j,\lambda _j \in \mathbb{N}_{0}$ is
\begin{eqnarray}
&&\sum_{\lambda _0 =0}^{\infty } \frac{ (\frac{3}{2} )_{\lambda _0}}{ \lambda _0 !} s_0^{\lambda _0} \prod _{n=1}^{\infty } \left\{ \sum_{ \lambda _n = \lambda _{n-1}}^{\infty } s_n^{\lambda _n }\right\}  MS_{\lambda _j}^R\left( q,\lambda =2^2\left(\lambda_j +2j+1/2 \right)^2- 2 q; \eta  = q x^2, x= \mathrm{\cos}^2z \right)  \nonumber\\
&&= x^{\frac{1}{2}}\Bigg\{ \prod_{l=1}^{\infty } \frac{1}{(1-s_{l,\infty })} \mathbf{B}\left( s_{0,\infty } ;x\right)  + \Bigg\{\prod_{l=2}^{\infty } \frac{1}{(1-s_{l,\infty })} \int_{0}^{1} dt_1\;t_1^{\frac{3}{2} } \int_{0}^{1} du_1\;u_1 \overleftrightarrow {\mathbf{\Psi}}_1 \left(s_{1,\infty };t_1,u_1,x\right) \mathbf{B}\left( s_{0} ;\widetilde{w}_{1,1}\right) \Bigg\}\eta \nonumber\\
&&+ \sum_{n=2}^{\infty } \Bigg\{ \prod_{l=n+1}^{\infty } \frac{1}{(1-s_{l,\infty })} \int_{0}^{1} dt_n\;t_n^{2n-\frac{1}{2} } \int_{0}^{1} du_n\;u_n^{2n-1} \overleftrightarrow {\mathbf{\Psi}}_n \left(s_{n,\infty };t_n,u_n,x \right)  \label{eq:60061}\\
&&\times \prod_{k=1}^{n-1} \Bigg\{ \int_{0}^{1} dt_{n-k}\;t_{n-k}^{2(n-k)-\frac{1}{2}} \int_{0}^{1} du_{n-k} \;u_{n-k}^{2(n-k)-1} \overleftrightarrow {\mathbf{\Psi}}_{n-k} \left( s_{n-k};t_{n-k},u_{n-k},\widetilde{w}_{n-k+1,n} \right) \Bigg\} \left. \mathbf{B}\left( s_{0} ;\widetilde{w}_{1,n}\right)\Bigg\} \eta^n  \right\} \nonumber
\end{eqnarray}
where
\begin{equation}
\begin{cases} 
{ \displaystyle \overleftrightarrow {\mathbf{\Psi}}_1 \left(s_{1,\infty };t_1,u_1,x\right)= \frac{\left( \frac{1+s_{1,\infty }+\sqrt{s_{1,\infty }^2-2(1-2x(1-t_1)(1-u_1))s_{1,\infty }+1}}{2}\right)^{-4}}{\sqrt{s_{1,\infty }^2-2(1-2x(1-t_1)(1-u_1))s_{1,\infty }+1}} }\cr
{ \displaystyle  \overleftrightarrow {\mathbf{\Psi}}_n \left(s_{n,\infty };t_n,u_n,x \right) = \frac{\left( \frac{1+s_{n,\infty }+\sqrt{s_{n,\infty }^2-2(1-2x(1-t_n)(1-u_n))s_{n,\infty }+1}}{2}\right)^{- 4n }}{\sqrt{s_{n,\infty }^2-2(1-2x(1-t_n)(1-u_n))s_{n,\infty }+1}}}\cr
{ \displaystyle \overleftrightarrow {\mathbf{\Psi}}_{n-k} \left(s_{n-k};t_{n-k},u_{n-k},\widetilde{w}_{n-k+1,n} \right) = \frac{\left( \frac{(1+s_{n-k})+\sqrt{s_{n-k}^2-2(1-2\widetilde{w}_{n-k+1,n} (1-t_{n-k})(1-u_{n-k}))s_{n-k}+1}}{2}\right)^{-4(n-k)}}{\sqrt{s_{n-k}^2-2(1-2\widetilde{w}_{n-k+1,n} (1-t_{n-k})(1-u_{n-k}))s_{n-k}+1}}}
\end{cases}\nonumber 
\end{equation}
and
\begin{equation}
\begin{cases} 
{ \displaystyle \mathbf{B} \left( s_{0,\infty } ;x\right)= \frac{\left( 1- s_{0,\infty }+\sqrt{s_{0,\infty }^2-2(1-2x)s_{0,\infty }+1}\right)^{-\frac{1}{2}} \left( 1+s_{0,\infty }+\sqrt{s_{0,\infty }^2-2(1-2x )s_{0,\infty }+1}\right)^{\frac{1}{2}}}{\sqrt{s_{0,\infty }^2-2(1-2x)s_{0,\infty }+1}}}\cr
{ \displaystyle  \mathbf{B} \left( s_{0} ;\widetilde{w}_{1,1}\right) = \frac{\left( 1- s_0+\sqrt{s_0^2-2(1-2\widetilde{w}_{1,1})s_0+1}\right)^{-\frac{1}{2}} \left( 1+s_0+\sqrt{s_0^2-2(1-2\widetilde{w}_{1,1} )s_0+1}\right)^{\frac{1}{2}}}{\sqrt{s_0^2-2(1-2\widetilde{w}_{1,1})s_0+1}}} \cr
{ \displaystyle \mathbf{B} \left( s_{0} ;\widetilde{w}_{1,n}\right) = \frac{\left( 1- s_0+\sqrt{s_0^2-2(1-2\widetilde{w}_{1,n})s_0+1}\right)^{ -\frac{1}{2}} \left(1+s_0+\sqrt{s_0^2-2(1-2\widetilde{w}_{1,n} )s_0+1}\right)^{\frac{1}{2}}}{\sqrt{s_0^2-2(1-2\widetilde{w}_{1,n})s_0+1}}}
\end{cases}\nonumber 
\end{equation}
\end{remark}
\begin{proof}
Replace $A$, $w$ and $\gamma  $  by 1, $s_{0,\infty } $ and $\frac{3}{2}$ in (\ref{eq:60038}).  
\begin{eqnarray}
&&\sum_{\lambda _0=0}^{\infty }\frac{(\frac{3}{2})_{\lambda_0}}{\lambda_0 !} s_{0,\infty }^{\lambda_0} \;_2F_1\left( -\lambda_0, \lambda_0 +1; \frac{3}{2} ; x\right) \label{eq:60062}\\
&&=  \frac{\left(1- s_{0,\infty }+\sqrt{s_{0,\infty }^2-2(1-2x)s_{0,\infty }+1}\right)^{-\frac{1}{2}} \left(1+s_{0,\infty }+\sqrt{s_{0,\infty }^2-2(1-2x)s_{0,\infty }+1}\right)^{\frac{1}{2}}}{\sqrt{s_{0,\infty }^2-2(1-2x)s_{0,\infty }+1}} \nonumber
\end{eqnarray} 
Replace $A$, $w$, $\gamma  $, $x$ by 1, $s_0 $, $\frac{3}{2}$ and $\widetilde{w}_{1,1}$ in (\ref{eq:60038}). 
\begin{eqnarray}
&&\sum_{\lambda _0=0}^{\infty }\frac{(\frac{3}{2})_{\lambda_0}}{\lambda_0 !} s_0^{\lambda_0} \;_2F_1\left( -\lambda_0, \lambda_0 +1; \frac{3}{2} ; \widetilde{w}_{1,1} \right) \label{eq:60063}\\
&&=  \frac{\left(1- s_0+\sqrt{s_0^2-2(1-2\widetilde{w}_{1,1})s_0+1}\right)^{-\frac{1}{2}} \left( 1+ s_0+\sqrt{s_0^2-2(1-2\widetilde{w}_{1,1})s_0+1}\right)^{\frac{1}{2}}}{\sqrt{s_0^2-2(1-2\widetilde{w}_{1,1})s_0+1}} \nonumber
\end{eqnarray} 
Replace $A$, $w$, $\gamma  $, $x$ by 1, $s_0 $, $\frac{3}{2}$ and $\widetilde{w}_{1,n}$ in (\ref{eq:60038}). 
\begin{eqnarray}
&&\sum_{\lambda _0=0}^{\infty }\frac{(\frac{3}{2})_{\lambda_0}}{\lambda_0 !} s_0^{\lambda_0} \;_2F_1\left( -\lambda_0, \lambda_0 +1; \frac{3}{2} ; \widetilde{w}_{1,n} \right) \label{eq:60064}\\
&&=  \frac{\left(1- s_0+\sqrt{s_0^2-2(1-2\widetilde{w}_{1,n})s_0+1}\right)^{-\frac{1}{2}} \left( 1+ s_0+\sqrt{s_0^2-2(1-2\widetilde{w}_{1,n})s_0+1}\right)^{\frac{1}{2}}}{\sqrt{s_0^2-2(1-2\widetilde{w}_{1,n})s_0+1}} \nonumber
\end{eqnarray} 
Put $ c_0= 1 $, $\lambda =1/2 $ and $\gamma' =3/2 $ in (\ref{eq:60045}). Substitute (\ref{eq:60062}), (\ref{eq:60063}) and (\ref{eq:60064}) into the new (\ref{eq:60045}).\qed
\end{proof}

\section{Mathieu equation about regular singular point at one}
Let $\zeta =1-x$ in (\ref{eq:6002}) to obtain analytic solutions of the Mathieu equation about $x=1 $.
\begin{equation}
4\zeta (1-\zeta) \frac{d^2{y}}{d{\zeta}^2} +2(1-2\zeta) \frac{d{y}}{d{\zeta}} + \left( \lambda -2q +4q\zeta \right) y = 0
\label{eq:60065}
\end{equation}
As we compare (\ref{eq:60065}) with (\ref{eq:6002}), an independent variable and a coefficient on the above are correspondent to the following way.
\begin{equation}
\begin{split}
& x  \longrightarrow   \zeta= 1-x \\ & q \longrightarrow  -q 
\end{split}\label{eq:60066}   
\end{equation}
For the asymptotic behavior of the function $y(x)$ about $x=1$ for infinite series, replace $x$ by $\zeta= 1-x$ in (\ref{eq:6008}). 
\begin{equation}
\lim_{n\gg 1}y(\zeta) \approx  \frac{1}{x} \hspace{1cm}\mbox{where}\;  x=\mathrm{\cos}^2z  
\nonumber
\end{equation}
For polynomial of type 2, replace $x$ and $q$ by $\zeta= 1-x$ and $-q$ in (\ref{eq:60010}). 
\begin{equation}
\lim_{n\gg 1}y(\zeta) >  \; _1F_2 \left(1; \frac{1}{2},\frac{1}{2};  -\frac{q}{4}\zeta^2 \right) 
\nonumber
\end{equation}
As we put (\ref{eq:60066}) into analytic solutions of Mathieu equation about $x=0$ using 3TRF\cite{1Chou2012e} and R3TRF, we obtain  (1) the Frobenius solution in closed form of Mathieu equation about $x=1$, (2) its integral form, (3) the generating function for the Mathieu polynomial of type 2. The local solutions are as follow.
\subsection{Mathieu equation about $x=1$ using 3TRF} 
\subsubsection{Power series for infinite series}
\begin{remark}
The power series expansion of Mathieu equation of the first kind for infinite series about $x=1$ using 3TRF is
\begin{eqnarray}
 y(\zeta )&=& M^{(o)}F\left( q, \lambda; \zeta= 1-x, \eta = -\frac{1}{4}q\zeta^2, x=\mathrm{\cos}^2z \right) \nonumber\\
&=&  \sum_{i_0=0}^{\infty} \frac{1}{(1)_{i_0} (\frac{3}{4})_{i_0}} \eta^{i_0} +   \left\{ \sum_{i_0=0}^{\infty }\frac{i_0^2 - \frac{1}{4^2}(\lambda -2q)}{(i_0+\frac{1}{2})(i_0+\frac{1}{4})} \frac{1}{(1)_{i_0} (\frac{3}{4})_{i_0}}
 \sum_{i_1=i_0}^{\infty } \frac{(\frac{3}{2})_{i_0}(\frac{5}{4})_{i_0}}{(\frac{3}{2})_{i_1}(\frac{5}{4})_{i_1}} \eta^{i_1} \right\}\zeta\nonumber\\
&&+ \sum_{n=2}^{\infty } \left\{ \sum_{i_0=0}^{\infty }  \frac{i_0^2 - \frac{1}{4^2}(\lambda -2q)}{(i_0+\frac{1}{2})(i_0+\frac{1}{4})} \frac{1}{(1)_{i_0} (\frac{3}{4})_{i_0}}\right. \nonumber\\
&&\times \prod _{k=1}^{n-1} \left\{ \sum_{i_k=i_{k-1}}^{\infty }  \frac{(i_k+\frac{k}{2})^2 - \frac{1}{4^2}(\lambda -2q)}{(i_k+\frac{k}{2}+\frac{1}{2})(i_k +\frac{k}{2}+\frac{1}{4})} \frac{(1+\frac{k}{2})_{i_{k-1}}(\frac{3}{4}+\frac{k}{2})_{i_{k-1}}}{(1+\frac{k}{2})_{i_{k}}(\frac{3}{4}+\frac{k}{2})_{i_{k}}}\right\} \nonumber\\
&&\times \left. \sum_{i_n= i_{n-1}}^{\infty } \frac{(1+\frac{n}{2})_{i_{n-1}}(\frac{3}{4}+\frac{n}{2})_{i_{n-1}}}{(1+\frac{n}{2})_{i_{n}}(\frac{3}{4}+\frac{n}{2})_{i_{n}}} \eta ^{i_n}\right\}  \zeta^n  \label{eq:60067}
\end{eqnarray}
\end{remark}
\begin{remark}
The power series expansion of Mathieu equation of the second kind for infinite series about $x=1$ using 3TRF is
\begin{eqnarray}
 y(\zeta)&=& M^{(o)}S\left( q, \lambda; \zeta= 1-x, \eta = -\frac{1}{4}q\zeta^2, x=\mathrm{\cos}^2z \right) \nonumber\\
&=& \zeta^{\frac{1}{2}} \left\{\sum_{i_0=0}^{\infty} \frac{1}{(\frac{5}{4})_{i_0} (1)_{i_0}} \eta^{i_0} \right. +  \left\{ \sum_{i_0=0}^{\infty }\frac{(i_0+\frac{1}{4})^2 - \frac{1}{4^2}(\lambda -2q)}{(i_0+\frac{3}{4})(i_0+\frac{1}{2})} \frac{1}{(\frac{5}{4})_{i_0} (1)_{i_0}}
  \sum_{i_1=i_0}^{\infty } \frac{(\frac{7}{4})_{i_0}(\frac{3}{2})_{i_0}}{(\frac{7}{4})_{i_1}(\frac{3}{2})_{i_1}} \eta^{i_1} \right\}\zeta\nonumber\\
&&+ \sum_{n=2}^{\infty } \left\{ \sum_{i_0=0}^{\infty } \frac{(i_0+\frac{1}{4})^2 - \frac{1}{4^2}(\lambda -2q)}{(i_0+\frac{3}{4})(i_0+\frac{1}{2})} \frac{1}{(\frac{5}{4})_{i_0} (1)_{i_0}}\right.\nonumber\\
&&\times \prod _{k=1}^{n-1} \left\{ \sum_{i_k=i_{k-1}}^{\infty }  \frac{(i_k+\frac{1}{4}+\frac{k}{2})^2 - \frac{1}{4^2}(\lambda -2q)}{(i_k+\frac{3}{4}+\frac{k}{2})(i_k +\frac{1}{2}+\frac{k}{2})} \frac{(\frac{5}{4}+\frac{k}{2})_{i_{k-1}}(1+\frac{k}{2})_{i_{k-1}}}{(\frac{5}{4}+\frac{k}{2})_{i_{k}}(1+\frac{k}{2})_{i_{k}}}\right\} \nonumber\\
&&\times \left.\left. \sum_{i_n= i_{n-1}}^{\infty } \frac{(\frac{5}{4}+\frac{n}{2})_{i_{n-1}}(1+\frac{n}{2})_{i_{n-1}}}{( \frac{5}{4}+\frac{n}{2})_{i_{n}}(1+\frac{n}{2})_{i_{n}}} \eta ^{i_n}\right\} \zeta^n \right\}\label{eq:60068}
\end{eqnarray}
\end{remark}
\subsubsection{Integral formalism for infinite series}
\begin{remark}
The integral representation of Mathieu equation of the first kind for infinite series about $x=1$ using 3TRF is
\begin{eqnarray}
 y(\zeta)&=& M^{(o)}F\left( q, \lambda; \zeta= 1-x, \eta = -\frac{1}{4}q\zeta^2, x=\mathrm{\cos}^2z \right) \nonumber\\
 &=& \Gamma \left( 3/4 \right) \Bigg\{ \eta ^{\frac{1}{8}} I_{-\frac{1}{4}}\left(2\sqrt{\eta}\right)+ \sum_{n=1}^{\infty } \Bigg\{\prod _{k=0}^{n-1} \left\{ \int_{0}^{1} dt_{n-k}\;t_{n-k}^{-\frac{5}{4}+\frac{1}{2}(n-k)} \int_{0}^{1} du_{n-k}\;u_{n-k}^{-1+\frac{1}{2}(n-k)}\right. \nonumber\\
&&\times  I_0\left(2\sqrt{w_{n+1-k,n}(1-t_{n-k})(1-u_{n-k})}\right)\nonumber\\
&&\times  \left. \left( w_{n-k,n}^{-\frac{1}{2}(n-1-k)}\left(w_{n-k,n}\partial_{w_{n-k,n}}\right)^2 w_{n-k,n}^{\frac{1}{2}(n-1-k)}- \frac{1}{4^2}(\lambda -2q) \right)\right\}
  w_{1,n}^{\frac{1}{8}} I_{-\frac{1}{4}}\left(2\sqrt{w_{1,n}}\right) \Bigg\} \zeta^n \Bigg\} \hspace{1.5cm}\label{eq:60069}
\end{eqnarray}
\end{remark}
\begin{remark}
The integral representation of Mathieu equation of the second kind for infinite series about $x=1$ using 3TRF is
\begin{eqnarray}
 y(\zeta)&=& M^{(o)}S\left( q, \lambda; \zeta= 1-x, \eta = -\frac{1}{4}q\zeta^2, x=\mathrm{\cos}^2z \right) \nonumber\\
&=& \Gamma ( 5/4 ) \zeta^{\frac{1}{2}} \Bigg\{ \eta ^{-\frac{1}{8}}I_{\frac{1}{4}}\left(2\sqrt{\eta}\right)+ \sum_{n=1}^{\infty } \Bigg\{\prod _{k=0}^{n-1} \left\{ \int_{0}^{1} dt_{n-k}\;t_{n-k}^{-1+\frac{1}{2}(n-k)} \int_{0}^{1} du_{n-k}\;u_{n-k}^{-\frac{3}{4}+\frac{1}{2}(n-k)} \right.\nonumber\\
&&\times I_0\left(2\sqrt{w_{n+1-k,n}(1-t_{n-k})(1-u_{n-k})}\right)\nonumber\\
&&\times \left. \left( w_{n-k,n}^{-\frac{1}{2}(n-k-\frac{1}{2})}\left(w_{n-k,n}\partial_{w_{n-k,n}}\right)^2 w_{n-k,n}^{\frac{1}{2}(n-k-\frac{1}{2})}- \frac{1}{4^2}(\lambda -2q) \right)\right\} 
  w_{1,n}^{-\frac{1}{8}} I_{\frac{1}{4}}\left(2\sqrt{w_{1,n}}\right) \Bigg\}  \zeta^n \Bigg\}  \hspace{1.5cm}\label{eq:60070}
\end{eqnarray}
\end{remark}
In (\ref{eq:60069}) and (\ref{eq:60070}),
\begin{equation}w_{a,b}=
\begin{cases} \displaystyle {\eta \prod _{l=a}^{b} t_l u_l }\;\;\mbox{where}\; a\leq b\cr
\eta  \;\;\mbox{only}\;\mbox{if}\; a>b
\end{cases}
\nonumber
\end{equation}
\subsection{Mathieu equation about $x=1$ using R3TRF}
\subsubsection{Power series}
\paragraph{Polynomial of type 2}

\begin{remark}
The power series expansion of Mathieu equation of the first kind for polynomial of type 2 about $x=1$ as $\lambda = 2^2(\lambda_j +2j )^2+ 2 q$ where $j,\lambda_j \in \mathbb{N}_{0}$ is
\begin{eqnarray}
 y(\zeta )&=& M^{(o)}F_{\lambda _j}^R\left( q,\lambda = 2^2(\lambda_j +2j )^2+ 2 q;\zeta= 1-x, \eta  = -q \zeta^2, x= \mathrm{\cos}^2z \right)\nonumber\\
&=& \sum_{i_0=0}^{\lambda _0} \frac{(-\lambda_0)_{i_0} \left( \lambda _0 \right)_{i_0}}{(1)_{i_0}\left( \frac{1}{2}\right)_{i_0}} \zeta^{i_0} \nonumber\\
&&+ \left\{ \sum_{i_0=0}^{\lambda_0}\frac{1}{ (i_0+ 2)\left( i_0+\frac{3}{2} \right)}\frac{(-\lambda _0)_{i_0} \left( \lambda_0 \right)_{i_0}}{(1)_{i_0}\left(\frac{1}{2}\right)_{i_0}} \right. \left. \sum_{i_1=i_0}^{\lambda _1} \frac{(-\lambda_1)_{i_1}\left(\lambda_1 + 4 \right)_{i_1}(3)_{i_0}\left(\frac{5}{2}\right)_{i_0}}{(-\lambda_1)_{i_0}\left(\lambda_1 + 4\right)_{i_0}(3)_{i_1}\left(\frac{5}{2}\right)_{i_1}} \zeta^{i_1}\right\}\eta \nonumber\\
&&+ \sum_{n=2}^{\infty } \left\{ \sum_{i_0=0}^{\lambda_0}\frac{1}{ (i_0+ 2)\left( i_0+\frac{3}{2}\right)}\frac{(-\lambda _0)_{i_0} \left( \lambda_0 \right)_{i_0}}{(1)_{i_0}\left(\frac{1}{2} \right)_{i_0}}\right.\nonumber\\
&&\times \prod _{k=1}^{n-1} \left\{ \sum_{i_k=i_{k-1}}^{\lambda _k} \frac{1}{(i_k+ 2k+2 )\left( i_k+ 2k+\frac{3}{2}\right)}\right.   \left.\frac{(-\lambda _k)_{i_k}\left( \lambda_k +4k \right)_{i_k}(2k+1)_{i_{k-1}}\left( 2k+\frac{1}{2} \right)_{i_{k-1}}}{(-\lambda _k)_{i_{k-1}}\left( \lambda_k +4k \right)_{i_{k-1}}(2k+1)_{i_k}\left( 2k+\frac{1}{2} \right)_{i_k}}\right\} \nonumber\\
&&\times  \left.\sum_{i_n= i_{n-1}}^{\lambda _n} \frac{(-\lambda _n)_{i_n}\left( \lambda_n +4n \right)_{i_n}(2n+1 )_{i_{n-1}}\left( 2n+\frac{1}{2} \right)_{i_{n-1}}}{(-\lambda _n)_{i_{n-1}}\left( \lambda_n +4n \right)_{i_{n-1}}(2n+1)_{i_n}\left( 2n+\frac{1}{2}\right)_{i_n}} \zeta^{i_n} \right\} \eta ^n  \label{eq:60071}
\end{eqnarray}
\end{remark}
\begin{remark}
The power series expansion of Mathieu equation of the second kind for polynomial of type 2 about $x=1$ as $\lambda = 2^2\left(\lambda_j +2j+1/2 \right)^2+ 2 q $ where $j,\lambda_j \in \mathbb{N}_{0}$ is
\begin{eqnarray}
 y(\zeta)&=& M^{(o)}S_{\lambda _j}^R\left( q,\lambda =2^2\left(\lambda_j +2j+1/2 \right)^2 + 2 q; \zeta= 1-x, \eta =-q \zeta^2, x= \mathrm{\cos}^2z \right)\nonumber\\
&=& \zeta^{\frac{1}{2}} \left\{\sum_{i_0=0}^{\lambda _0} \frac{(-\lambda_0)_{i_0} \left( \lambda _0+ 1 \right)_{i_0}}{\left(\frac{3}{2} \right)_{i_0}\left( 1\right)_{i_0}} \zeta^{i_0}\right.\nonumber\\
&&+ \left\{ \sum_{i_0=0}^{\lambda_0}\frac{1}{ \left(i_0+ \frac{5}{2} \right)\left( i_0 +2 \right)}\frac{(-\lambda _0)_{i_0} \left( \lambda_0+ 1 \right)_{i_0}}{\left(\frac{3}{2} \right)_{i_0}\left(1 \right)_{i_0}} \right. \left. \sum_{i_1=i_0}^{\lambda _1} \frac{(-\lambda_1)_{i_1}\left(\lambda_1 + 5 \right)_{i_1}\left(\frac{7}{2} \right)_{i_0}\left( 3\right)_{i_0}}{(-\lambda_1)_{i_0}\left(\lambda_1 +5 \right)_{i_0}\left(\frac{7}{2} \right)_{i_1}\left( 3 \right)_{i_1}} \zeta^{i_1}\right\}\eta \nonumber\\
&&+ \sum_{n=2}^{\infty } \left\{ \sum_{i_0=0}^{\lambda_0}\frac{1}{ \left(i_0+ \frac{5}{2} \right)\left( i_0+2 \right)}\frac{(-\lambda _0)_{i_0} \left( \lambda_0+  1 \right)_{i_0}}{\left(\frac{3}{2} \right)_{i_0}\left( 1 \right)_{i_0}}\right.\nonumber\\
&&\times \prod _{k=1}^{n-1} \left\{ \sum_{i_k=i_{k-1}}^{\lambda _k} \frac{1}{\left( i_k+ 2k+\frac{5}{2} \right)\left( i_k+ 2k+2 \right)}\right.   \left.\frac{(-\lambda _k)_{i_k}\left( \lambda_k +4k +1\right)_{i_k}\left( 2k+\frac{3}{2} \right)_{i_{k-1}}\left( 2k+1 \right)_{i_{k-1}}}{(-\lambda _k)_{i_{k-1}}\left( \lambda_k +4k +1 \right)_{i_{k-1}}\left( 2k+\frac{3}{2} \right)_{i_k}\left( 2k+1 \right)_{i_k}}\right\} \nonumber\\
&&\times \left. \left.\sum_{i_n= i_{n-1}}^{\lambda _n} \frac{(-\lambda _n)_{i_n}\left( \lambda_n +4n +1 \right)_{i_n}\left(2n+\frac{3}{2} \right)_{i_{n-1}}\left( 2n+1 \right)_{i_{n-1}}}{(-\lambda _n)_{i_{n-1}}\left( \lambda_n +4n +1 \right)_{i_{n-1}}\left( 2n+\frac{3}{2} \right)_{i_n}\left( 2n+1 \right)_{i_n}} \zeta^{i_n} \right\} \eta ^n \right\} \label{eq:60072}
\end{eqnarray}
\end{remark}

\paragraph{Infinite series}

\begin{remark}
The power series expansion of Mathieu equation of the first kind for infinite series about $x=1$ using R3TRF is
\begin{eqnarray}
 y(\zeta)&=& M^{(o)}F^R\left( q,\lambda, \varphi  = \frac{1}{2} \sqrt{\lambda - 2q} ; \zeta =1-x, \eta  = -q \zeta^2, x= \mathrm{\cos}^2z \right) \nonumber\\
&=& \sum_{i_0=0}^{\infty } \frac{\left( - \varphi \right)_{i_0} \left( \varphi \right)_{i_0}}{\left(1 \right)_{i_0}\left(\frac{1}{2} \right)_{i_0}} \zeta^{i_0} + \left\{ \sum_{i_0=0}^{\infty }\frac{1}{\left( i_0+ 2 \right) \left( i_0+ \frac{3}{2} \right)}\frac{\left( -\varphi \right)_{i_0} \left(  \varphi \right)_{i_0}}{(1 )_{i_0} \left( \frac{1}{2} \right)_{i_0}} \sum_{i_1=i_0}^{\infty } \frac{\left( 2-\varphi \right)_{i_1} \left( 2+\varphi \right)_{i_1}\left( 3 \right)_{i_0}\left( \frac{5}{2} \right)_{i_0}}{\left( 2-\varphi \right)_{i_0} \left( 2+\varphi \right)_{i_0}\left( 3 \right)_{i_1}\left( \frac{5}{2} \right)_{i_1}}\zeta^{i_1}\right\} \eta \nonumber\\
&+& \sum_{n=2}^{\infty } \left\{ \sum_{i_0=0}^{\infty }\frac{1}{\left( i_0+ 2 \right) \left( i_0+ \frac{3}{2} \right)}\frac{\left( -\varphi \right)_{i_0} \left( \varphi \right)_{i_0}}{(1 )_{i_0} \left( \frac{1}{2} \right)_{i_0}}\right.\nonumber\\
&\times& \prod _{k=1}^{n-1} \left\{ \sum_{i_k=i_{k-1}}^{\infty } \frac{1}{\left( i_k+ 2k+2 \right) \left( i_k+ 2k+\frac{3}{2} \right)} \frac{ \left( 2k -\varphi \right)_{i_k} \left( 2k +\varphi \right)_{i_k} \left( 2k+1 \right)_{i_{k-1}}\left( 2k+\frac{1}{2} \right)_{i_{k-1}}}{ \left( 2k -\varphi \right)_{i_{k-1}} \left( 2k +\varphi \right)_{i_{k-1}} \left( 2k+1 \right)_{i_k}\left( 2k+\frac{1}{2} \right)_{i_k}}\right\} \nonumber\\
&\times&  \left.\sum_{i_n= i_{n-1}}^{\infty }\frac{ \left( 2 n -\varphi \right)_{i_n} \left( 2n +\varphi \right)_{i_n} \left( 2n+1 \right)_{i_{n-1}}\left( 2n+\frac{1}{2} \right)_{i_{n-1}}}{ \left( 2n -\varphi \right)_{i_{n-1}} \left( 2n +\varphi \right)_{i_{n-1}} \left( 2n+1  \right)_{i_n}\left( 2n+\frac{1}{2} \right)_{i_n}} \zeta^{i_n} \right\} \eta ^n  \label{eq:60073}
\end{eqnarray}
\end{remark}
\begin{remark}
The power series expansion of Mathieu equation of the second kind for infinite series about $x=1$ using R3TRF is
\begin{eqnarray}
y(\zeta)&=& M^{(o)}S^R\left( q,\lambda, \varphi  = \frac{1}{2} \sqrt{\lambda - 2q} ; \zeta =1-x, \eta  = -q \zeta^2, x= \mathrm{\cos}^2z  \right) \nonumber\\
&=& \zeta^{\frac{1}{2}} \left\{\sum_{i_0=0}^{\infty } \frac{\left(\frac{1}{2} - \varphi \right)_{i_0} \left(\frac{1}{2} + \varphi \right)_{i_0}}{\left(\frac{3}{2} \right)_{i_0}\left(1\right)_{i_0}} \zeta^{i_0}\right.\nonumber\\
&+& \left\{ \sum_{i_0=0}^{\infty }\frac{1}{\left( i_0+ \frac{5}{2}\right) \left( i_0+ 2 \right)}\frac{\left(\frac{1}{2}-\varphi \right)_{i_0} \left( \frac{1}{2} +\varphi \right)_{i_0}}{(\frac{3}{2} )_{i_0} \left( 1 \right)_{i_0}} \sum_{i_1=i_0}^{\infty } \frac{\left( \frac{5}{2}-\varphi \right)_{i_1} \left(\frac{5}{2}+\varphi \right)_{i_1}\left( \frac{7}{2} \right)_{i_0}\left( 3 \right)_{i_0}}{\left( \frac{5}{2}-\varphi \right)_{i_0} \left(\frac{5}{2}+\varphi \right)_{i_0}\left( \frac{7}{2} \right)_{i_1}\left( 3 \right)_{i_1}}\zeta^{i_1}\right\} \eta \nonumber\\
&+& \sum_{n=2}^{\infty } \left\{ \sum_{i_0=0}^{\infty }\frac{1}{\left( i_0+ \frac{5}{2}\right) \left( i_0+ 2 \right)}\frac{\left(\frac{1}{2}-\varphi \right)_{i_0} \left( \frac{1}{2} +\varphi \right)_{i_0}}{(\frac{3}{2} )_{i_0} \left( 1 \right)_{i_0}}\right.\nonumber\\
&\times& \prod _{k=1}^{n-1} \left\{ \sum_{i_k=i_{k-1}}^{\infty } \frac{1}{\left( i_k+ 2k+\frac{5}{2} \right) \left( i_k+ 2k+2 \right)} \frac{ \left( 2k+\frac{1}{2} -\varphi \right)_{i_k} \left( 2k+\frac{1}{2} +\varphi \right)_{i_k} \left( 2k+\frac{3}{2} \right)_{i_{k-1}}\left( 2k+1 \right)_{i_{k-1}}}{ \left( 2k+\frac{1}{2} -\varphi \right)_{i_{k-1}} \left( 2k+ \frac{1}{2} +\varphi \right)_{i_{k-1}} \left( 2k+\frac{3}{2} \right)_{i_k}\left( 2k+1 \right)_{i_k}}\right\} \nonumber\\
&\times& \left.\left.\sum_{i_n= i_{n-1}}^{\infty }\frac{ \left( 2n+\frac{1}{2} -\varphi \right)_{i_n} \left( 2n+\frac{1}{2} +\varphi \right)_{i_n} \left( 2n+\frac{3}{2} \right)_{i_{n-1}}\left( 2n+1 \right)_{i_{n-1}}}{ \left( 2n+\frac{1}{2} -\varphi \right)_{i_{n-1}} \left( 2n+\frac{1}{2} +\varphi \right)_{i_{n-1}} \left( 2n+\frac{3}{2} \right)_{i_n}\left( 2n+1 \right)_{i_n}} \zeta^{i_n} \right\} \eta ^n \right\} \label{eq:60074}
\end{eqnarray}
\end{remark}
\paragraph{The case of small values of $q$}
As $|q| \ll 1$ in (\ref{eq:60073}) and (\ref{eq:60074}), we have
\begin{subequations}
\begin{eqnarray}
&&\lim_{|q| \ll 1} M^{(o)}F^R\left( q,\lambda, \varphi  = \frac{1}{2} \sqrt{\lambda - 2q} ; \zeta =1-x, \eta  = -q \zeta^2, x= \mathrm{\cos}^2z \right) \nonumber\\
&&\approx  \sum_{i_0=0}^{\infty } \frac{\left( -\frac{\sqrt{\lambda}}{2}\right)_{i_0} \left( \frac{\sqrt{\lambda}}{2}\right)_{i_0}}{\left(1 \right)_{i_0}\left(\frac{1}{2} \right)_{i_0}} \zeta^{i_0} \nonumber\\
&&+  \eta \sum_{i_0=0}^{\infty }\frac{1}{\left( i_0+ 2 \right) \left( i_0+ \frac{3}{2} \right)}\frac{\left( -\frac{\sqrt{\lambda }}{2} \right)_{i_0} \left(  \frac{\sqrt{\lambda }}{2} \right)_{i_0}}{(1 )_{i_0} \left( \frac{1}{2} \right)_{i_0}} \sum_{i_1=i_0}^{\infty } \frac{\left( 2-\frac{\sqrt{\lambda }}{2} \right)_{i_1} \left( 2+\frac{\sqrt{\lambda }}{2} \right)_{i_1}\left( 3 \right)_{i_0}\left( \frac{5}{2} \right)_{i_0}}{\left( 2-\frac{\sqrt{\lambda }}{2}\right)_{i_0} \left( 2+\frac{\sqrt{\lambda }}{2} \right)_{i_0}\left( 3 \right)_{i_1}\left( \frac{5}{2} \right)_{i_1}}\zeta^{i_1} \nonumber\\ 
&&> \mathrm{\cos}\left( \sqrt{\lambda } \mathrm{\sin}^{-1}\left(\sqrt{\zeta}\right)\right) \label{eq:60075a}
\end{eqnarray}
And,
\begin{eqnarray}
&&\lim_{|q| \ll 1}  M^{(o)}S^R\left( q,\lambda, \varphi  = \frac{1}{2} \sqrt{\lambda - 2q} ; \zeta =1-x, \eta  = -q \zeta^2, x= \mathrm{\cos}^2z  \right) \nonumber\\
&&\approx   \zeta^{\frac{1}{2}} \left\{ \sum_{i_0=0}^{\infty } \frac{\left(\frac{1}{2} - \frac{\sqrt{\lambda }}{2} \right)_{i_0} \left(\frac{1}{2} + \frac{\sqrt{\lambda }}{2} \right)_{i_0}}{\left(\frac{3}{2} \right)_{i_0}\left(1\right)_{i_0}} \zeta^{i_0} \right. \nonumber\\
&&+ \left. \eta \sum_{i_0=0}^{\infty }\frac{1}{\left( i_0+ \frac{5}{2}\right) \left( i_0+ 2 \right)}\frac{\left(\frac{1}{2}-\frac{\sqrt{\lambda }}{2} \right)_{i_0} \left( \frac{1}{2} +\frac{\sqrt{\lambda }}{2} \right)_{i_0}}{(\frac{3}{2} )_{i_0} \left( 1 \right)_{i_0}} \sum_{i_1=i_0}^{\infty } \frac{\left( \frac{5}{2}-\frac{\sqrt{\lambda }}{2} \right)_{i_1} \left(\frac{5}{2}+\frac{\sqrt{\lambda }}{2} \right)_{i_1}\left( \frac{7}{2} \right)_{i_0}\left( 3 \right)_{i_0}}{\left( \frac{5}{2}-\frac{\sqrt{\lambda }}{2} \right)_{i_0} \left(\frac{5}{2}+\frac{\sqrt{\lambda }}{2} \right)_{i_0}\left( \frac{7}{2} \right)_{i_1}\left( 3 \right)_{i_1}}\zeta^{i_1}\right\} \nonumber\\
&&> \sqrt{\frac{\zeta}{1-\zeta}}  \mathrm{\cos}\left( \sqrt{\lambda } \mathrm{\sin}^{-1}\left(\sqrt{\zeta}\right)\right)  \label{eq:60075b}
\end{eqnarray}
\end{subequations}
\subsubsection{Integral formalism}
\paragraph{Polynomial of type 2}

\begin{remark}
The integral representation of Mathieu equation of the first kind for polynomial of type 2 about $x=1$ as $\lambda = 2^2(\lambda_j +2j )^2+ 2 q$ where $j,\lambda _j \in \mathbb{N}_{0}$ is
\begin{eqnarray}
  y(\zeta)&=& M^{(o)}F_{\lambda _j}^R\left( q,\lambda = 2^2(\lambda_j +2j )^2+ 2 q;\zeta= 1-x, \eta  = -q \zeta^2, x= \mathrm{\cos}^2z \right) \nonumber\\
&=&\; _2F_1 \left( -\lambda _0, \lambda_0 ; \frac{1}{2}; \zeta \right)  + \sum_{n=1}^{\infty } \Bigg\{\prod _{k=0}^{n-1} \Bigg\{ \int_{0}^{1} dt_{n-k}\;t_{n-k}^{2(n-k)-1 } \int_{0}^{1} du_{n-k}\;u_{n-k}^{2(n-k )-\frac{3}{2} } \nonumber\\
&&\times  \frac{1}{2\pi i}  \oint dp_{n-k} \frac{1}{p_{n-k}} \left( 1-  w_{n-k+1,n}(1-t_{n-k})(1-u_{n-k})p_{n-k}\right)^{- 4(n-k) }  \nonumber\\
&&\times  \left( \frac{p_{n-k}-1}{p_{n-k}} \frac{1}{1- w_{n-k+1,n}(1-t_{n-k})(1-u_{n-k})p_{n-k}}\right)^{\lambda_{n-k}}  \Bigg\} \nonumber\\
&&\times  \left.  _2F_1 \left( -\lambda _0, \lambda_0 ; \frac{1}{2}; w_{1,n} \right) \right\} \eta ^n \hspace{1cm}\label{eq:60076}
\end{eqnarray}
\end{remark} 
\begin{remark}
The integral representation of Mathieu equation of the second kind for polynomial of type 2 about $x=1$ as $\lambda = 2^2\left(\lambda_j +2j+1/2 \right)^2+ 2 q $ where $j,\lambda _j \in \mathbb{N}_{0}$ is
\begin{eqnarray}
 y(\zeta)&=& M^{(o)}S_{\lambda _j}^R\left(  q,\lambda =2^2\left(\lambda_j +2j+1/2 \right)^2 + 2 q; \zeta= 1-x, \eta   =-q \zeta^2, x= \mathrm{\cos}^2z \right) \nonumber\\
&=& \zeta^{\frac{1}{2}} \left\{\; _2F_1 \left( -\lambda _0, \lambda_0 +1; \frac{3}{2}; \zeta \right) \right. + \sum_{n=1}^{\infty } \Bigg\{\prod _{k=0}^{n-1} \Bigg\{ \int_{0}^{1} dt_{n-k}\;t_{n-k}^{2(n-k)-\frac{1}{2} } \int_{0}^{1} du_{n-k}\;u_{n-k}^{2(n-k )-1 } \nonumber\\
&&\times  \frac{1}{2\pi i}  \oint dp_{n-k} \frac{1}{p_{n-k}} \left( 1-  w_{n-k+1,n}(1-t_{n-k})(1-u_{n-k})p_{n-k}\right)^{-\left( 4(n-k) +1 \right)}  \nonumber\\
&&\times  \left( \frac{p_{n-k}-1}{p_{n-k}} \frac{1}{1- w_{n-k+1,n}(1-t_{n-k})(1-u_{n-k})p_{n-k}}\right)^{\lambda_{n-k}}  \Bigg\} \nonumber\\
&&\times  \left.\left.  _2F_1 \left( -\lambda _0, \lambda_0 +1; \frac{3}{2}; w_{1,n} \right) \right\} \eta ^n \right\} \label{eq:60077}
\end{eqnarray}
\end{remark}
\paragraph{Infinite series}

\begin{remark}
The integral representation of Mathieu equation of the first kind for infinite series about $x=1$ using R3TRF is
\begin{eqnarray}
y(\zeta)&=& M^{(o)}F^R\left( q,\lambda, \varphi  = \frac{1}{2} \sqrt{\lambda - 2q} ; \zeta =1-x, \eta  = -q \zeta^2, x= \mathrm{\cos}^2z \right)\nonumber\\
&=& \; _2F_1 \left( -\varphi, \varphi ; \frac{1}{2}; \zeta \right)  + \sum_{n=1}^{\infty } \Bigg\{\prod _{k=0}^{n-1} \Bigg\{ \int_{0}^{1} dt_{n-k}\;t_{n-k}^{2(n-k)-1 } \int_{0}^{1} du_{n-k}\;u_{n-k}^{2(n-k )-\frac{3}{2} }  \nonumber\\
&&\times  \frac{1}{2\pi i}  \oint dp_{n-k} \frac{1}{p_{n-k}} \left( 1-  w_{n-k+1,n}(1-t_{n-k})(1-u_{n-k})p_{n-k}\right)^{-4(n-k) }  \nonumber\\
&&\times  \left( \frac{p_{n-k}-1}{p_{n-k}} \frac{1}{1- w_{n-k+1,n}(1-t_{n-k})(1-u_{n-k})p_{n-k}}\right)^{-2(n-k)+ \varphi }  \Bigg\} \nonumber\\
&&\times  \; _2F_1 \left( -\varphi, \varphi ; \frac{1}{2}; w_{1,n} \right)  \Bigg\} \eta ^n  \label{eq:60078}
\end{eqnarray}
\end{remark}
\begin{remark}
The integral representation of Mathieu equation of the second kind for infinite series about $x=1$ using R3TRF is
\begin{eqnarray}
 y(\zeta)&=& M^{(o)}S^R\left( q,\lambda, \varphi  = \frac{1}{2} \sqrt{\lambda - 2q} ; \zeta =1-x, \eta  = -q \zeta^2, x= \mathrm{\cos}^2z  \right) \nonumber\\
&=& \zeta^{\frac{1}{2}} \left\{\; _2F_1 \left( \frac{1}{2}-\varphi, \frac{1}{2}+\varphi ; \frac{3}{2}; \zeta \right) \right. + \sum_{n=1}^{\infty } \Bigg\{\prod _{k=0}^{n-1} \Bigg\{ \int_{0}^{1} dt_{n-k}\;t_{n-k}^{2(n-k)-\frac{1}{2} } \int_{0}^{1} du_{n-k}\;u_{n-k}^{2(n-k )-1 }  \nonumber\\
&&\times  \frac{1}{2\pi i}  \oint dp_{n-k} \frac{1}{p_{n-k}} \left( 1-  w_{n-k+1,n}(1-t_{n-k})(1-u_{n-k})p_{n-k}\right)^{-\left( 4(n-k) +1 \right)}  \nonumber\\
&&\times  \left( \frac{p_{n-k}-1}{p_{n-k}} \frac{1}{1- w_{n-k+1,n}(1-t_{n-k})(1-u_{n-k})p_{n-k}}\right)^{-2(n-k) -\frac{1}{2}+ \varphi}  \Bigg\} \nonumber\\
&&\times \left. \; _2F_1 \left( \frac{1}{2}-\varphi, \frac{1}{2}+\varphi ; \frac{3}{2}; w_{1,n} \right) \Bigg\} \eta ^n \right\}  \label{eq:60079}
\end{eqnarray}
\end{remark}
In (\ref{eq:60076})--(\ref{eq:60079}),
\begin{equation} w _{i,j}=
\begin{cases} \displaystyle {\frac{p_i}{(p_i-1)}\; \frac{ w_{i+1,j} t_i u_i}{1-  w_{i+1,j} p_i (1-t_i)(1-u_i)}} \;\;\mbox{where}\; i\leq j\cr
\zeta \;\;\mbox{only}\;\mbox{if}\; i>j
\end{cases}\nonumber 
\end{equation}
\subsubsection{Generating function for the Mathieu polynomial of type 2}
\begin{remark}
The generating function for the Mathieu polynomial of type 2 of the first kind about $x=1$  as $\lambda  = 2^2(\lambda _j+2j )^2+2q $ where $j,\lambda _j \in \mathbb{N}_{0}$ is
\begin{eqnarray}
&&\sum_{\lambda _0 =0}^{\infty } \frac{ (\frac{1}{2} )_{\lambda _0}}{ \lambda _0 !} s_0^{\lambda _0} \prod _{n=1}^{\infty } \left\{ \sum_{ \lambda _n = \lambda _{n-1}}^{\infty } s_n^{\lambda _n }\right\} M^{(o)}F_{\lambda _j}^R\left( q,\lambda = 2^2(\lambda_j +2j )^2+ 2 q;\zeta= 1-x, \eta  = -q \zeta^2, x= \mathrm{\cos}^2z \right) \nonumber\\
&&= \frac{1}{2} \Bigg\{ \prod_{l=1}^{\infty } \frac{1}{(1-s_{l,\infty })}  \mathbf{A}\left( s_{0,\infty } ;\zeta\right) \nonumber\\
&&+ \Bigg\{ \prod_{l=2}^{\infty } \frac{1}{(1-s_{l,\infty })} \int_{0}^{1} dt_1\;t_1 \int_{0}^{1} du_1\;u_1^{\frac{1}{2}} \overleftrightarrow {\mathbf{\Gamma}}_1 \left(s_{1,\infty };t_1,u_1,\zeta\right) \mathbf{A}\left( s_{0} ;\widetilde{w}_{1,1}\right)\Bigg\} \eta \nonumber\\
&&+ \sum_{n=2}^{\infty } \Bigg\{ \prod_{l=n+1}^{\infty } \frac{1}{(1-s_{l,\infty })} \int_{0}^{1} dt_n\;t_n^{2n-1} \int_{0}^{1} du_n\;u_n^{ 2n-\frac{3}{2} } \overleftrightarrow {\mathbf{\Gamma}}_n \left(s_{n,\infty };t_n,u_n,\zeta \right)   \label{eq:60080}\\
&&\times \prod_{k=1}^{n-1} \Bigg\{ \int_{0}^{1} dt_{n-k}\;t_{n-k}^{2(n-k)-1} \int_{0}^{1} du_{n-k} \;u_{n-k}^{2(n-k)-\frac{3}{2}}\overleftrightarrow {\mathbf{\Gamma}}_{n-k} \left(s_{n-k};t_{n-k},u_{n-k},\widetilde{w}_{n-k+1,n} \right)\Bigg\}  \mathbf{A} \left( s_{0} ;\widetilde{w}_{1,n}\right) \Bigg\} \eta^n \Bigg\}  \nonumber 
\end{eqnarray}
where
\begin{equation}
\begin{cases} 
{ \displaystyle \overleftrightarrow {\mathbf{\Gamma}}_1 \left(s_{1,\infty };t_1,u_1,\zeta\right)= \frac{\left( \frac{1+s_{1,\infty }+\sqrt{s_{1,\infty }^2-2(1-2\zeta(1-t_1)(1-u_1))s_{1,\infty }+1}}{2}\right)^{-3}}{\sqrt{s_{1,\infty }^2-2(1-2\zeta (1-t_1)(1-u_1))s_{1,\infty }+1}}}\cr
{ \displaystyle  \overleftrightarrow {\mathbf{\Gamma}}_n \left(s_{n,\infty };t_n,u_n,\zeta\right) =\frac{\left( \frac{1+s_{n,\infty }+\sqrt{s_{n,\infty }^2-2(1-2\zeta(1-t_n)(1-u_n))s_{n,\infty }+1}}{2}\right)^{-\left( 4n-1 \right)}}{\sqrt{ s_{n,\infty }^2-2(1-2\zeta(1-t_n)(1-u_n))s_{n,\infty }+1}}}\cr
{ \displaystyle \overleftrightarrow {\mathbf{\Gamma}}_{n-k} \left(s_{n-k};t_{n-k},u_{n-k},\widetilde{w}_{n-k+1,n} \right) = \frac{ \left( \frac{1+s_{n-k}+\sqrt{s_{n-k}^2-2(1-2\widetilde{w}_{n-k+1,n} (1-t_{n-k})(1-u_{n-k}))s_{n-k}+1}}{2}\right)^{-\left( 4(n-k)-1 \right)}}{\sqrt{ s_{n-k}^2-2(1-2\widetilde{w}_{n-k+1,n} (1-t_{n-k})(1-u_{n-k}))s_{n-k}+1}}}
\end{cases}\nonumber 
\end{equation}
and
\begin{equation}
\begin{cases} 
{ \displaystyle \mathbf{A} \left( s_{0,\infty } ;\zeta\right)= \frac{\left(1- s_{0,\infty }+\sqrt{s_{0,\infty }^2-2(1-2\zeta)s_{0,\infty }+1}\right)^{\frac{1}{2}} \left(1+s_{0,\infty }+\sqrt{s_{0,\infty }^2-2(1-2\zeta )s_{0,\infty }+1}\right)^{\frac{1}{2}}}{\sqrt{s_{0,\infty }^2-2(1-2\zeta )s_{0,\infty }+1}}}\cr
{ \displaystyle  \mathbf{A} \left( s_{0} ;\widetilde{w}_{1,1}\right) = \frac{\left(1- s_0+\sqrt{s_0^2-2(1-2\widetilde{w}_{1,1})s_0+1}\right)^{\frac{1}{2}} \left( 1+s_0+\sqrt{s_0^2-2(1-2\widetilde{w}_{1,1} )s_0+1}\right)^{\frac{1}{2}}}{\sqrt{s_0^2-2(1-2\widetilde{w}_{1,1})s_0+1}}} \cr
{ \displaystyle \mathbf{A} \left( s_{0} ;\widetilde{w}_{1,n}\right) = \frac{\left( 1- s_0+\sqrt{s_0^2-2(1-2\widetilde{w}_{1,n})s_0+1}\right)^{\frac{1}{2}} \left(1+s_0+\sqrt{s_0^2-2(1-2\widetilde{w}_{1,n} )s_0+1}\right)^{\frac{1}{2}}}{\sqrt{s_0^2-2(1-2\widetilde{w}_{1,n})s_0+1}}}
\end{cases}\nonumber 
\end{equation}
\end{remark}

\begin{remark}
The generating function for the Mathieu polynomial of type 2 of the second kind about $x=1$ as $\lambda  = 2^2(\lambda _j+2j +1/2)^2+2q $ where $j,\lambda _j \in \mathbb{N}_{0}$ is
\begin{eqnarray}
&&\sum_{\lambda _0 =0}^{\infty } \frac{ (\frac{3}{2} )_{\lambda _0}}{ \lambda _0 !} s_0^{\lambda _0} \prod _{n=1}^{\infty } \left\{ \sum_{ \lambda _n = \lambda _{n-1}}^{\infty } s_n^{\lambda _n }\right\}  M^{(o)}S_{\lambda _j}^R\bigg(  q,\lambda =2^2\left(\lambda_j +2j+1/2 \right)^2 + 2 q; \zeta= 1-x\nonumber\\
&&,\eta =-q \zeta^2, x= \mathrm{\cos}^2z \bigg)  \nonumber\\
&&= \zeta^{\frac{1}{2}}\Bigg\{ \prod_{l=1}^{\infty } \frac{1}{(1-s_{l,\infty })} \mathbf{B}\left( s_{0,\infty } ;\zeta\right) \nonumber\\
&&+ \Bigg\{\prod_{l=2}^{\infty } \frac{1}{(1-s_{l,\infty })} \int_{0}^{1} dt_1\;t_1^{\frac{3}{2} } \int_{0}^{1} du_1\;u_1 \overleftrightarrow {\mathbf{\Psi}}_1 \left(s_{1,\infty };t_1,u_1,\zeta\right) \mathbf{B}\left( s_{0} ;\widetilde{w}_{1,1}\right) \Bigg\}\eta \nonumber\\
&&+ \sum_{n=2}^{\infty } \Bigg\{ \prod_{l=n+1}^{\infty } \frac{1}{(1-s_{l,\infty })} \int_{0}^{1} dt_n\;t_n^{2n-\frac{1}{2} } \int_{0}^{1} du_n\;u_n^{2n-1} \overleftrightarrow {\mathbf{\Psi}}_n \left(s_{n,\infty };t_n,u_n,\zeta \right) \nonumber\\
&&\times \prod_{k=1}^{n-1} \Bigg\{ \int_{0}^{1} dt_{n-k}\;t_{n-k}^{2(n-k)-\frac{1}{2}} \int_{0}^{1} du_{n-k} \;u_{n-k}^{2(n-k)-1} \overleftrightarrow {\mathbf{\Psi}}_{n-k} \left( s_{n-k};t_{n-k},u_{n-k},\widetilde{w}_{n-k+1,n} \right) \Bigg\} \nonumber\\
&&\times \left. \mathbf{B}\left( s_{0} ;\widetilde{w}_{1,n}\right)\Bigg\} \eta^n  \right\} \label{eq:60081}
\end{eqnarray}
where
\begin{equation}
\begin{cases} 
{ \displaystyle \overleftrightarrow {\mathbf{\Psi}}_1 \left(s_{1,\infty };t_1,u_1,\zeta\right)= \frac{\left( \frac{1+s_{1,\infty }+\sqrt{s_{1,\infty }^2-2(1-2\zeta(1-t_1)(1-u_1))s_{1,\infty }+1}}{2}\right)^{-4}}{\sqrt{s_{1,\infty }^2-2(1-2\zeta(1-t_1)(1-u_1))s_{1,\infty }+1}} }\cr
{ \displaystyle  \overleftrightarrow {\mathbf{\Psi}}_n \left(s_{n,\infty };t_n,u_n,\zeta \right) = \frac{\left( \frac{1+s_{n,\infty }+\sqrt{s_{n,\infty }^2-2(1-2\zeta(1-t_n)(1-u_n))s_{n,\infty }+1}}{2}\right)^{- 4n }}{\sqrt{s_{n,\infty }^2-2(1-2\zeta(1-t_n)(1-u_n))s_{n,\infty }+1}}}\cr
{ \displaystyle \overleftrightarrow {\mathbf{\Psi}}_{n-k} \left(s_{n-k};t_{n-k},u_{n-k},\widetilde{w}_{n-k+1,n} \right) = \frac{\left( \frac{(1+s_{n-k})+\sqrt{s_{n-k}^2-2(1-2\widetilde{w}_{n-k+1,n} (1-t_{n-k})(1-u_{n-k}))s_{n-k}+1}}{2}\right)^{-4(n-k)}}{\sqrt{s_{n-k}^2-2(1-2\widetilde{w}_{n-k+1,n} (1-t_{n-k})(1-u_{n-k}))s_{n-k}+1}}}
\end{cases}\nonumber 
\end{equation}
and
\begin{equation}
\begin{cases} 
{ \displaystyle \mathbf{B} \left( s_{0,\infty } ;\zeta\right)= \frac{\left( 1- s_{0,\infty }+\sqrt{s_{0,\infty }^2-2(1-2\zeta)s_{0,\infty }+1}\right)^{-\frac{1}{2}} \left( 1+s_{0,\infty }+\sqrt{s_{0,\infty }^2-2(1-2\zeta )s_{0,\infty }+1}\right)^{\frac{1}{2}}}{\sqrt{s_{0,\infty }^2-2(1-2\zeta)s_{0,\infty }+1}}}\cr
{ \displaystyle  \mathbf{B} \left( s_{0} ;\widetilde{w}_{1,1}\right) = \frac{\left( 1- s_0+\sqrt{s_0^2-2(1-2\widetilde{w}_{1,1})s_0+1}\right)^{-\frac{1}{2}} \left( 1+s_0+\sqrt{s_0^2-2(1-2\widetilde{w}_{1,1} )s_0+1}\right)^{\frac{1}{2}}}{\sqrt{s_0^2-2(1-2\widetilde{w}_{1,1})s_0+1}}} \cr
{ \displaystyle \mathbf{B} \left( s_{0} ;\widetilde{w}_{1,n}\right) = \frac{\left( 1- s_0+\sqrt{s_0^2-2(1-2\widetilde{w}_{1,n})s_0+1}\right)^{ -\frac{1}{2}} \left(1+s_0+\sqrt{s_0^2-2(1-2\widetilde{w}_{1,n} )s_0+1}\right)^{\frac{1}{2}}}{\sqrt{s_0^2-2(1-2\widetilde{w}_{1,n})s_0+1}}}
\end{cases}\nonumber 
\end{equation}
\end{remark}
In (\ref{eq:60080}) and (\ref{eq:60081}),
\begin{equation}
\begin{cases}
\displaystyle { s_{a,b}} = \begin{cases} \displaystyle {  s_a\cdot s_{a+1}\cdot s_{a+2}\cdots s_{b-2}\cdot s_{b-1}\cdot s_b}\;\;\mbox{if}\;a>b \cr
s_a \;\;\mbox{if}\;a=b\end{cases}
\cr
\cr
\displaystyle { \widetilde{w}_{i,j}}  = 
\begin{cases} \displaystyle { \frac{ \widetilde{w}_{i+1,j}\; t_i u_i \left\{ 1+ (s_i+2\widetilde{w}_{i+1,j}(1-t_i)(1-u_i))s_i\right\}}{2(1-\widetilde{w}_{i+1,j}(1-t_i)(1-u_i))^2 s_i}} \cr
\displaystyle {-\frac{\widetilde{w}_{i+1,j}\; t_i u_i (1+s_i)\sqrt{s_i^2-2(1-2\widetilde{w}_{i+1,j}(1-t_i)(1-u_i))s_i+1}}{2(1-\widetilde{w}_{i+1,j}(1-t_i)(1-u_i))^2 s_i}} \;\;\mbox{where}\;i<j \cr
\cr
\displaystyle { \frac{\zeta t_i u_i \left\{ 1+ (s_{i,\infty }+2\zeta(1-t_i)(1-u_i))s_{i,\infty }\right\}}{2(1-\zeta(1-t_i)(1-u_i))^2 s_{i,\infty }}} \cr
\displaystyle {-\frac{\zeta t_i u_i(1+s_{i,\infty })\sqrt{s_{i,\infty }^2-2(1-2\zeta (1-t_i)(1-u_i))s_{i,\infty }+1}}{2(1-\zeta (1-t_i)(1-u_i))^2 s_{i,\infty }}} \;\;\mbox{where}\;i=j 
\end{cases}
\end{cases}\nonumber
\end{equation}
where
\begin{equation}
a,b,i,j\in \mathbb{N}_{0} \nonumber
\end{equation}
\section{Mathieu equation about irregular singular point at infinity}
\subsection{Mathieu equation about $x=\infty $}
Let $\varsigma  =\frac{1}{x}$ in (\ref{eq:6002}) to obtain the analytic solution of Mathieu equation about $x=\infty $.
\begin{equation}
4\varsigma^3 (\varsigma -1) \frac{d^2{y}}{d{\varsigma}^2} +2\varsigma^2 (3\varsigma-2) \frac{d{y}}{d{\varsigma}} + \left( (\lambda +2q)\varsigma -4q \right) y = 0
\label{eq:60082}
\end{equation}
Assume that its solution is
\begin{equation}
y(\varsigma)= \sum_{n=0}^{\infty } c_n \varsigma^{n+\nu }\label{eq:60083}
\end{equation}
Plug (\ref{eq:60083}) into (\ref{eq:60082}).
\begin{equation}
c_{n+1}= B_n \;c_{n-1} \hspace{1cm};n\geq 1\label{eq:60084}
\end{equation}
where,
\begin{subequations}
\begin{equation}
B_n = \frac{\left( n+\frac{1}{2}\sqrt{\lambda }-1\right)\left( n+\frac{1}{2}\sqrt{\lambda } -\frac{1}{2}\right)}{n \left( n+\sqrt{\lambda }\right)}\label{eq:60085a}
\end{equation}
with
\begin{equation}
q=0 \hspace{2cm} \nu= \frac{1}{2}\sqrt{\lambda }\label{eq:60085b}
\end{equation}
\end{subequations}
The Frobenius solution of (\ref{eq:60085a}) with (\ref{eq:60085b}) is
\begin{equation}
y(\varsigma) = M^{(i,1)}F\left( q=0, \lambda; \varsigma = \frac{1}{x}, x= \mathrm{\cos}^2z \right)
 = \varsigma^{\frac{1}{2}\sqrt{\lambda }}\; _2F_1 \left( \frac{1}{2}\sqrt{\lambda },\frac{1}{2}\sqrt{\lambda }+\frac{1}{2}; \sqrt{\lambda }+1;  \varsigma \right) \label{eq:60086}
\end{equation}
The power series of (\ref{eq:60085a}) for polynomial does not exist. If $\sqrt{\lambda }$ is equal to $-2m$ or $-2m-1$ where $m=0,1,2,\cdots$, the function $y(\varsigma)$ does not convergent because Pochhammer symbol $(\sqrt{\lambda }+1)_n$ on the denominator will be zero at certain value of $m$. For infinite series, the series in $y(\varsigma)$ does not have circle $|\varsigma|= |\mathrm{\cos}^{-2}z|<1$ as its circle of convergence.   
   
\begin{lemma}
The hypergeometric functions is defined by
\begin{equation}
_2F_1(a,b;c;z)= \sum_{n=0}^{\infty } \frac{(a)_n (b)_n}{(c)_n} \frac{z^n}{n!} \hspace{1cm}\mbox{where}\;c\ne 0,-1,-2,\cdots \nonumber
\end{equation}
converges for $|z|<1$. On the boundary $z=1$ of the region of convergence, a sufficient condition for absolute convergence of the series is $Re(c-a-b)>0$. \cite{Rain1960} 
\end{lemma}
By using this lemma, the analytic solution of (\ref{eq:60086}) is given by
\begin{eqnarray}
y(\varsigma) &=& M^{(i,1)}F\left( q=0, \lambda; \varsigma = \frac{1}{x}=1, x= \mathrm{\cos}^2z \right)\nonumber\\
 &=&  \; _2F_1 \left( \frac{1}{2}\sqrt{\lambda },\frac{1}{2}\sqrt{\lambda }+\frac{1}{2}; \sqrt{\lambda }+1;  1 \right) = 2^{\sqrt{\lambda }} 
\hspace{.5cm}\mbox{where}\; z= 0,\pm \pi, \pm 2\pi, \cdots \hspace{2cm}  \label{eq:60087}
\end{eqnarray}
\subsection{Mathieu equation about $\zeta=\infty $ }
Let $\phi  =\frac{1}{\zeta}$ in (\ref{eq:60065}) to obtain the analytic solution of Mathieu equation about $\zeta=\infty $.
\begin{equation}
4\phi^3 (\phi -1) \frac{d^2{y}}{d{\phi}^2} +2\phi^2 (3\phi-2) \frac{d{y}}{d{\phi}} + \left( (\lambda +2q)\phi -4q \right) y = 0
\label{eq:60088}
\end{equation}
As $\varsigma$ is replaced by $\phi$ in (\ref{eq:60082}), its equation is correspondent to (\ref{eq:60088}). The series of Mathieu equation about $\zeta=\infty $ is given by
\begin{eqnarray}
y(\varsigma) &=& M^{(i,2)}F\left( q=0, \lambda; \phi = \frac{1}{\zeta}=1, \zeta= 1-x, x= \mathrm{\cos}^2z \right)\nonumber\\
 &=&  \; _2F_1 \left( \frac{1}{2}\sqrt{\lambda },\frac{1}{2}\sqrt{\lambda }+\frac{1}{2}; \sqrt{\lambda }+1;  1 \right) = 2^{\sqrt{\lambda }} 
\hspace{.5cm}\mbox{where}\; z= \pm \frac{1}{2}\pi, \pm \frac{3}{2}\pi, \pm \frac{5}{2}\pi, \cdots \hspace{2cm}  \label{eq:60089}
\end{eqnarray}
\section{Summary}

In Ref.\cite{1Chou2012e} I show how to obtain the power series expansion in closed forms and its integral form of Mathieu equation for infinite series including all higher terms of $A_n$'s by applying 3TRF. This was done by letting $A_n$ in sequence $c_n$ is the leading term in the analytic function $y(x)$: the sequence $c_n$ consists of combinations $A_n$ and $B_n$. The series for the type 1 and 3 polynomials does not exist by using either 3TRF or R3TRF. Because there are no ways to terminate $B_n$ term with a fixed parameter $q$ of numerator in $B_n$ whenever index $n$ increases in (\ref{eq:6005b}).

In this chapter I show how to construct the power series expansion in closed forms and its integral form of Mathieu equation for infinite series and polynomial of type 2 including all higher terms of $B_n$'s by applying R3TRF. This is done by letting $B_n$ in sequence $c_n$ is the leading term in the analytic function $y(x)$. For the type 2 polynomial, I treat $q$ as a free variable and a fixed value of $\lambda $. 

The power series expansion and its integral form of Mathieu equation for infinite series about $x=0$ in this chapter are equivalent to infinite series of Mathieu equation in Ref.\cite{1Chou2012e}. In this chapter $B_n$ is the leading term in sequence $c_n$ in the analytic function $y(x)$. In Ref.\cite{1Chou2012e} $A_n$ is the leading term in sequence $c_n$ in the analytic function $y(x)$.
 
In Ref.\cite{1Chou2012e} and this chapter, as we see the power series expansion of Mathieu equation about $x=0$ for either polynomial or infinite series, the denominators and numerators in all $A_n$ or $B_n$ terms of $y(x)$ arise with Pochhammer symbol. Since we construct the power series expansion with Pochhammer symbols in numerators and denominators, we are able to describe integral representation of Mathieu equation analytically. As we observe representations in closed form integrals of Mathieu equation by applying 3TRF, a modified Bessel function recurs in each of sub-integral forms of the $y(x)$. And $_2F_1$ function recurs in each of sub-integral forms of $y(x)$ as we observe  integral forms of Mathieu equation by applying R3TRF (each sub-integral $y_m(x)$ where $m=0,1,2,\cdots$ is composed of $2m$ terms of definite integrals and $m$ terms of contour integrals).
We are able to transform the Mathieu function into any well-known special functions having two recurrence relation in its power series of a linear differential equation because of modified Bessel and $_2F_1$ functions in each of sub-integral forms of the Mathieu function. After we replace modified Bessel and $_2F_1$ functions in its integral forms to other special functions (such as  Bessel function, Kummer function, hypergeometric function, Laguerre function and etc), we are able to rebuild the power series expansion of  Mathieu equation in a backward.   
 
Indeed, the generating function for the Mathieu polynomial of type 2 is constructed analytically from its integral representation by applying the generating function for the Jacobi polynomial using hypergeometric functions. Also it might be possible to obtain orthogonal relations, recursion relations and expectation values of physical quantities from the generating function for the Mathieu polynomial of type 2: the processes in order to obtain orthogonal and recursion relations of the Mathieu polynomial are similar as the case of a normalized wave function for the hydrogen-like atoms.  

Mathematical structure of the Frobenius solution, its integral representation and the generating function for Mathieu equation closely resembles the case of Confluent Heun equation using 3TRF and R3TRF. As $\beta \rightarrow 0$, $\delta =\gamma = \frac{1}{2}$, $\alpha = \frac{q}{0}$, $q\rightarrow \frac{\lambda +2q}{4}$ and $\lambda \rightarrow \nu$ in the power series expansion and its integral forms of the Confluent Heun equation using R3TRF in chapter 5, its analytic solution are correspondent to the Frobenius solution and its integral representation of Mathieu equation in this chapter; compare remarks 5.2.1--5.2.9 in chapter 5 with remarks 7.2.1--7.2.9 in this chapter. By similar process, the analytic solution of the Confluent Heun equation using 3TRF for infinite series by changing all coefficients in chapter 4 is equivalent to the power series and its integral form of Mathieu equation in Ref.\cite{1Chou2012e}; compare remarks 4.2.3, 4.2.4 and 4.2.8 and 4.2.9 in chapter 4 with remarks 1--4 in Ref.\cite{1Chou2012e}. 

As $\displaystyle {\widetilde{w}_{1,1}^{-(\alpha +\lambda )}\left(  \widetilde{w}_{1,1} \partial _{ \widetilde{w}_{1,1}}\right) \widetilde{w}_{1,1}^{\alpha +\lambda}, \widetilde{w}_{n,n}^{-(2(n-1)+\alpha +\lambda )}\left(  \widetilde{w}_{n,n} \partial _{ \widetilde{w}_{n,n}}\right)  \widetilde{w}_{n,n}^{2(n-1)+\alpha +\lambda},}$ 

 $\displaystyle { \widetilde{w}_{n-k,n}^{-(2(n-k-1)+\alpha +\lambda )}\left(  \widetilde{w}_{n-k,n} \partial _{ \widetilde{w}_{n-k,n}}\right) \widetilde{w}_{n-k,n}^{2(n-k-1)+\alpha +\lambda} \rightarrow 1}$ with $\beta \rightarrow 0$, $\delta =\gamma = \frac{1}{2}$, $\alpha = \frac{q}{0}$, $q\rightarrow \frac{\lambda +2q}{4}$ and $\lambda \rightarrow \nu$  in the general expression of the generating function for the type 2 Confluent Heun polynomial using R3TRF in the chapter 5, its solution is equivalent to the generating function for the type 2 Mathieu polynomial in this chapter; compare theorem 5.2.12 in chapter 5 with theorem 7.2.12 in this chapter.
Because, by applying the generating function for the Jacobi polynomial using $_2F_1$ functions into both integral forms of the Confluent Heun and Mathieu polynomials of type 2 ($_2F_1$ function recurs in each of sub-integral forms of them), I derive generating functions for these polynomials of type 2.   

\addcontentsline{toc}{section}{Bibliography}
\bibliographystyle{model1a-num-names}
\bibliography{<your-bib-database>}
\bibliographystyle{model1a-num-names}
\bibliography{<your-bib-database>}

\chapter{Lame function in the algebraic form using reversible three-term recurrence formula}
\chaptermark{Lame function in the algebraic form using R3TRF}
In Ref.\cite{zChou2012f} I apply three term recurrence formula (3TRF)\cite{zchou2012b} to the power series expansion in closed forms of Lame equation in the algebraic form (for infinite series and polynomial which makes $B_n$ term terminated  including all higher terms of $A_n$'s\footnote{`` higher terms of $A_n$'s'' means at least two terms of $A_n$'s.}) and its integral form.  I show how to transform power series expansion in Lame equation to its integral representation for infinite series and polynomial in mathematical rigour.

In this chapter I will apply reversible three term recurrence formula (R3TRF) in chapter 1 to (1) the Frobenius solution in closed forms of Lame equation, (2) its integral form (for infinite series and polynomial which makes $A_n$ term terminated including all higher terms of $B_n$'s\footnote{`` higher terms of $B_n$'s'' means at least two terms of $B_n$'s.}), (3) the generating function for the Lame polynomial which makes $A_n$ term terminated. 
               
\section{Introduction}

The sphere is a geometrical perfect shape, the set of points which are all equidistant from its center (a fixed point) in three-dimensional space. In contrast, an ellipsoid is a imperfect one, a surface whose plane sections are all ellipses or circles; the set of points are not same distance from the center of the ellipsoid any more. 
As we all recognize, the nature is nonlinear and imperfect geometrically. For the purpose of simplification, we usually linearize those system in order to take a step to the future with a good numerical approximation. Actually, many geometrical spherical objects (earth, sun, black hole, etc) are not perfectly sphere in nature. The shape of those objects are closely better interpreted by an ellipsoid because of their rotations by themselves. For example, the ellipsoidal harmonics are represented in calculations of gravitational potential\cite{zRomai2001}. However spherical harmonic is preferred over the more mathematically  complex ellipsoid harmonics (the coefficients in the power series expansion of Lame equation have a relation between a 3-term).

In 1837, Gabriel Lame introduced a second ordinary differential equation which has four regular singular points in the method of separation of variables applied to the Laplace equation in elliptic coordinates\cite{zLame1837}. Various authors has called this equation as `Lame equation' or `ellipsoidal harmonic equation'\cite{zErde1955}. Lame equation is applicable to diverse areas such as boundary value problems in ellipsoidal geometry, chaotic Hamiltonian systems, the theory of Bose-Einstein condensates, etc \cite{zBrac2001,zBron2001,zQian2003}.

Until the present day, there are no analytic solutions in closed forms of the Lame function\cite{zHobs1931,zWhit1952,zErde1955}. Using Frobenius method to obtain an analytic solution of Lame equation represented either in the algebraic form or in Weierstrass's form, a solution automatically comes out three term recurrence relation\cite{zHobs1931,zWhit1952}. They left the analytic solution of Lame equation as solutions of recurrences because of a 3-term recursive relation between successive coefficients in its power series expansion. In comparison with two term recursion relation of power series in any linear ordinary differential equations, an analytic solution in closed forms on three term recurrence relation of power series is unknown currently because of its complex mathematical calculation.  

In Ref.\cite{zChou2012f} I construct analytic solutions of Lame equation in the algebraic form about the regular singular point at $x=a$ by applying 3TRF. \cite{zchou2012b}: the power series expansion in closed forms of Lame equation and its integral form for infinite series and polynomial which makes $B_n$ term terminated including all higher terms of $A_n$'s.
 One interesting observation resulting from the calculations is the fact that a $_2F_1$ function recurs in each of sub-integral forms: the first sub-integral form contains zero term of $A_n's$, the second one contains one term of $A_n$'s, the third one contains two terms of $A_n$'s, etc. 

In this chapter, by applying R3TRF in chapter 1, I construct the power series expansion in closed forms of Lame equation in the algebraic form about a regular singular point at $x=a$ for infinite series and polynomial which makes $A_n$ term terminated including all higher terms of $B_n$'s. Indeed, an integral form of Lame equation and the generating function of Lame polynomial which makes $A_n$ term terminated will be derived analytically. 

Lame equation is a second-order linear ordinary differential equation of the algebraic form\cite{zLame1837,zStie1885,zZwil1997}
\begin{equation}
\frac{d^2{y}}{d{x}^2} + \frac{1}{2}\left(\frac{1}{x-a} +\frac{1}{x-b} + \frac{1}{x-c}\right) \frac{d{y}}{d{x}} +  \frac{-\alpha (\alpha +1) x+q}{4 (x-a)(x-b)(x-c)} y = 0\label{eq:7001}
\end{equation}
Lame equation has four regular singular points: $a$, $b$, $c$ and $\infty $; the exponents at the first three are all $0$ and $\frac{1}{2}$, and those at infinity are $-\frac{1}{2}\alpha $ and $\frac{1}{2}(\alpha +1)$.\cite{zMoon1961,zBoch1894} (\ref{eq:7001}) is a special case of Heun equation. Heun equation is a second-order linear ordinary differential equation of the form \cite{zHeun1889}.
\begin{equation}
\frac{d^2{y}}{d{x}^2} + \left(\frac{\gamma }{x} +\frac{\delta }{x-1} + \frac{\epsilon }{x-a}\right) \frac{d{y}}{d{x}} +  \frac{\alpha \beta x-q}{x(x-1)(x-a)} y = 0 \label{eq:7002}
\end{equation}
With the condition $\epsilon = \alpha +\beta -\gamma -\delta +1$. The parameters play different roles: $a \ne 0 $ is the singularity parameter, $\alpha $, $\beta $, $\gamma $, $\delta $, $\epsilon $ are exponent parameters, $q$ is the accessory parameter which in many physical applications appears as a spectral parameter. Also, $\alpha $ and $\beta $ are identical to each other. The total number of free parameters is six. It has four regular singular points which are 0, 1, $a$ and $\infty $ with exponents $\{ 0, 1-\gamma \}$, $\{ 0, 1-\delta \}$, $\{ 0, 1-\epsilon \}$ and $\{ \alpha, \beta \}$.

Let $z=x-a$ in (\ref{eq:7001}).
\begin{equation}
\frac{d^2{y}}{d{z}^2} + \frac{1}{2}\left(\frac{1}{z} +\frac{1}{z-(b-a)} + \frac{1}{z-(c-a)}\right) \frac{d{y}}{d{z}} +  \frac{-\frac{1}{4}\alpha (\alpha +1)z- \frac{1}{4}\left(- q+ \alpha (\alpha +1)a\right)}{ z(z-(b-a))(z-(c-a))} y = 0\label{eq:7003}
\end{equation}
As we compare (\ref{eq:7002}) with (\ref{eq:7003}), all coefficients on the above are correspondent to the following way.
\begin{equation}
\begin{split}
& \gamma ,\delta ,\epsilon  \longleftrightarrow   \frac{1}{2} \\ & 1\longleftrightarrow  b-a \\ & a\longleftrightarrow  c-a \\ & \alpha  \longleftrightarrow \frac{1}{2}(\alpha +1) \\
& \beta   \longleftrightarrow -\frac{1}{2}\alpha \\
& q \longleftrightarrow  \frac{1}{4}\left(- q+ \alpha (\alpha +1)a\right) \\ & x \longleftrightarrow z
\end{split}\label{eq:7004}   
\end{equation}
Assume that the solution of (\ref{eq:7003}) is
\begin{equation}
y(z)= \sum_{n=0}^{\infty } c_n z^{n+\lambda } \label{eq:7005}
\end{equation}
where $\lambda $ is an indicial root. Plug (\ref{eq:7005}) into (\ref{eq:7003}).
\begin{equation}
c_{n+1}=A_n \;c_n +B_n \;c_{n-1} \hspace{1cm};n\geq 1\label{eq:7006}
\end{equation}
where,
\begin{subequations}
\begin{eqnarray}
A_n &=& \frac{\frac{1}{4}(\alpha (\alpha +1)a-q)-(2a-b-c)(n+\lambda )^2}{(a-b)(a-c)(n+1+\lambda )(n+\frac{1}{2}+\lambda )}\nonumber\\
 &=& -\frac{(2a-b-c)}{(a-b)(a-c)} \frac{(n-\varphi +\lambda )(n+ \varphi +\lambda)}{(n+1+\lambda )(n+\frac{1}{2}+\lambda )} \label{eq:7007a}
\end{eqnarray}
and
\begin{equation}
\varphi = \sqrt{\frac{\alpha (\alpha +1)a-q}{4(2a-b-c)}} \nonumber
\end{equation}
\vspace{2mm}
\begin{eqnarray}
B_n &=& \frac{[\alpha -(1-2(n+\lambda ))][\alpha -2(n-1+\lambda )]}{2^2(a-b)(a-c)(n+1+\lambda )(n+\frac{1}{2}+\lambda )}\nonumber\\
 &=& \frac{-1}{(a-b)(a-c)}\frac{\left( n+\frac{\alpha }{2}-\frac{1}{2}+\lambda \right) \left( n-\frac{\alpha }{2}-1+\lambda\right)}{(n+1+\lambda )(n+\frac{1}{2}+\lambda )}\label{eq:7007b}
\end{eqnarray}
\begin{equation}
c_1= A_0 \;c_0\label{eq:7007c}
\end{equation}
\end{subequations}
We have two indicial roots which are $\lambda = 0$ and $ \frac{1}{2}$.
\section{Power series}
\subsection{Polynomial of type 2}
There are three types of polynomials in three-term recurrence relation of a linear ordinary differential equation: (1) polynomial which makes $B_n$ term terminated: $A_n$ term is not terminated, (2) polynomial which makes $A_n$ term terminated: $B_n$ term is not terminated, (3) polynomial which makes $A_n$ and $B_n$ terms terminated at the same time.\footnote{If $A_n$ and $B_n$ terms are not terminated, it turns to be infinite series.} In general Lame (spectral) polynomial represented either in the algebraic form or in Weierstrass's form is defined as type 3 polynomial where $A_n$ and $B_n$ terms terminated. Lame polynomial comes from a Lame equation that has a fixed integer value of $\alpha $, just as it has a fixed value of $q$. In three-term recurrence formula, polynomial of type 3 I categorize as complete polynomial. In future papers I will derive type 3 Lame polynomial. In Ref.\cite{zChou2012f} I derive the Frobenius solution and its integral form of the Lame polynomial of type 1: I treat $q$ as a free variable  and $\alpha $ as a fixed value. In this chapter I construct the power series expansion and integral representation for the Lame polynomial of type 2:  I treat $\alpha $ as a free variable and $q$ as a fixed value.   

In chapter 1 the general expression of power series of $y(x)$ for polynomial of type 2 is defined by
\begin{eqnarray}
y(x) &=&   \sum_{n=0}^{\infty } y_{n}(x) = y_0(x)+ y_1(x)+ y_2(x)+y_3(x)+\cdots \nonumber\\
&=&  c_0 \Bigg\{ \sum_{i_0=0}^{\alpha _0} \left( \prod _{i_1=0}^{i_0-1}A_{i_1} \right) x^{i_0+\lambda }
+ \sum_{i_0=0}^{\alpha _0}\left\{ B_{i_0+1} \prod _{i_1=0}^{i_0-1}A_{i_1}  \sum_{i_2=i_0}^{\alpha _1} \left( \prod _{i_3=i_0}^{i_2-1}A_{i_3+2} \right)\right\} x^{i_2+2+\lambda }  \nonumber\\
&& + \sum_{N=2}^{\infty } \Bigg\{ \sum_{i_0=0}^{\alpha _0} \Bigg\{B_{i_0+1}\prod _{i_1=0}^{i_0-1} A_{i_1} 
\prod _{k=1}^{N-1} \Bigg( \sum_{i_{2k}= i_{2(k-1)}}^{\alpha _k} B_{i_{2k}+2k+1}\prod _{i_{2k+1}=i_{2(k-1)}}^{i_{2k}-1}A_{i_{2k+1}+2k}\Bigg)\nonumber\\
&& \times  \sum_{i_{2N} = i_{2(N-1)}}^{\alpha _N} \Bigg( \prod _{i_{2N+1}=i_{2(N-1)}}^{i_{2N}-1} A_{i_{2N+1}+2N} \Bigg) \Bigg\} \Bigg\} x^{i_{2N}+2N+\lambda }\Bigg\}  \label{eq:7008}
\end{eqnarray}
In the above, $\alpha _i\leq \alpha _j$ only if $i\leq j$ where $i,j,\alpha _i, \alpha _j \in \mathbb{N}_{0}$.

For a polynomial, we need a condition which is:
\begin{equation}
A_{\alpha _i+ 2i}=0 \hspace{1cm} \mathrm{where}\;i,\alpha _i =0,1,2,\cdots
\label{eq:7009}
\end{equation}
In the above, $ \alpha _i$ is an eigenvalue that makes $A_n$ term terminated at certain value of index $n$. (\ref{eq:7009}) makes each $y_i(x)$ where $i=0,1,2,\cdots$ as the polynomial in (\ref{eq:7008}).
\subsubsection*{The case of $ \varphi = -( q_i +2i+\lambda )$ where $i,q_i = 0,1,2,\cdots$}

In (\ref{eq:7009}) replace  index $\alpha _i$ by $q_i$. In (\ref{eq:7007a})-(\ref{eq:7007c}) replace $ \varphi $ by ${ \displaystyle -( q_i +2i+\lambda ) }$.   Take the new (\ref{eq:7007a})-(\ref{eq:7007c}), (\ref{eq:7009}) and put them in (\ref{eq:7008}) with replacing variable $x$ by $z$.
After the replacement process, the general expression of power series of Lame equation in the algebraic form for polynomial of type 2 is given by
\begin{eqnarray}
 y(z)&=&  \sum_{n=0}^{\infty } y_{n}(z) = y_0(z)+ y_1(z)+ y_2(z)+y_3(z)+\cdots \nonumber\\ 
&=& c_0 z^{\lambda } \left\{\sum_{i_0=0}^{q_0} \frac{(-q_0)_{i_0} \left( q_0 +2\lambda \right)_{i_0}}{\left(1+\lambda \right)_{i_0}\left(\frac{1}{2} +\lambda \right)_{i_0}} \rho ^{i_0}\right.\nonumber\\
&&+ \left\{ \sum_{i_0=0}^{q_0}\frac{\left( i_0+\frac{\alpha}{2}+\frac{1}{2} + \lambda \right)\left( i_0-\frac{\alpha}{2} + \lambda \right)}{ \left(i_0+ 2+\lambda \right) \left( i_0+ \frac{3}{2}+ \lambda \right)}\frac{(-q_0)_{i_0} \left( q_0+ 2\lambda \right)_{i_0}}{\left(1+\lambda \right)_{i_0}\left(\frac{1}{2} +\lambda \right)_{i_0}} \right.\nonumber\\
&&\times  \left. \sum_{i_1=i_0}^{q_1} \frac{(-q_1)_{i_1}\left( q_1 +4+2\lambda \right)_{i_1}\left( 3+\lambda \right)_{i_0}\left( \frac{5}{2} +\lambda \right)_{i_0}}{(-q_1)_{i_0}\left( q_1 +4+2\lambda \right)_{i_0}\left( 3+\lambda \right)_{i_1}\left( \frac{5}{2} +\lambda \right)_{i_1}} \rho ^{i_1}\right\} \eta \nonumber\\
&&+ \sum_{n=2}^{\infty } \left\{ \sum_{i_0=0}^{q_0}\frac{\left( i_0+\frac{\alpha}{2}+\frac{1}{2} + \lambda \right)\left( i_0-\frac{\alpha}{2} + \lambda \right)}{ \left(i_0+ 2+\lambda \right) \left( i_0+ \frac{3}{2}+ \lambda \right)}\frac{(-q_0)_{i_0} \left( q_0+ 2\lambda \right)_{i_0}}{\left(1+\lambda \right)_{i_0}\left(\frac{1}{2} +\lambda \right)_{i_0}} \right.\nonumber\\
&&\times \prod _{k=1}^{n-1} \left\{ \sum_{i_k=i_{k-1}}^{q_k} \frac{\left( i_k+2k+\frac{\alpha}{2}+\frac{1}{2} + \lambda \right)\left( i_k+2k-\frac{\alpha}{2} + \lambda \right)}{ \left(i_k+ 2k+2+\lambda \right) \left( i_k+2k +\frac{3}{2}+ \lambda \right)}\right. \nonumber\\
&&\times \left.\frac{(-q_k)_{i_k}\left( q_k+4k + 2\lambda \right)_{i_k}\left(2k+1+\lambda \right)_{i_{k-1}}\left( 2k+\frac{1}{2} +\lambda \right)_{i_{k-1}}}{(-q_k)_{i_{k-1}}\left( q_k+4k +2\lambda \right)_{i_{k-1}}\left(2k+1+\lambda \right)_{i_k}\left( 2k+\frac{1}{2} +\lambda \right)_{i_k}}\right\} \nonumber\\
&&\times \left. \left.\sum_{i_n= i_{n-1}}^{q_n} \frac{(-q_n)_{i_n}\left( q_n+4n + 2\lambda \right)_{i_n}\left( 2n +1+\lambda \right)_{i_{n-1}}\left( 2n+\frac{1}{2} +\lambda \right)_{i_{n-1}}}{(-q_n)_{i_{n-1}}\left( q_n+4n +2\lambda \right)_{i_{n-1}}\left( 2n +1+\lambda \right)_{i_n}\left( 2n +\frac{1}{2} +\lambda \right)_{i_n}} \rho ^{i_n} \right\} \eta ^n \right\} \hspace{1.5cm}\label{eq:70010}
\end{eqnarray}
where
\begin{equation}
\begin{cases} z = x-a \cr
\rho = -\frac{2a-b-c}{(a-b)(a-c)} z \cr
\eta = -\frac{1}{(a-b)(a-c)} z^2 \cr
\varphi = -\beta +\gamma +\delta -1 \cr
q= \alpha (\alpha +1)a- 4(2a-b-c)(q_j+2j+\lambda )^2 \;\;\mbox{as}\;j,q_j\in \mathbb{N}_{0} \cr
q_i\leq q_j \;\;\mbox{only}\;\mbox{if}\;i\leq j\;\;\mbox{where}\;i,j\in \mathbb{N}_{0} 
\end{cases}\nonumber 
\end{equation}
\subsubsection*{The case of  $ \varphi = q_i +2i+\lambda $}

In (\ref{eq:7009}) replace  index $\alpha _i$ by $q_i$. In (\ref{eq:7007a})-(\ref{eq:7007c}) replace $ \varphi $ by ${ \displaystyle q_i +2i+\lambda }$. Take the new (\ref{eq:7007a})-(\ref{eq:7007c}), (\ref{eq:7009}) and put them in (\ref{eq:7008}) with replacing variable $x$ by $z$. After the replacement process, its solution is equivalent to (\ref{eq:70010}).

Put $c_0$= 1 as $\lambda =0$  for the first independent solution of Lame equation and $\lambda =\frac{1}{2}$ for the second one in (\ref{eq:70010}).
\begin{remark}
The power series expansion of Lame equation in the algebraic form of the first kind for polynomial of type 2 about $x=a$ as $q= \alpha (\alpha +1)a- 4(2a-b-c)(q_j+2j )^2 $ where $j,q_j \in \mathbb{N}_{0}$ is
\begin{eqnarray}
y(z)&=& LF_{q_j}^R\left( a, b, c, \alpha, q= \alpha (\alpha +1)a- 4(2a-b-c)(q_j+2j )^2; z= x-a \right. \nonumber\\
&&, \rho = -\frac{2a-b-c}{(a-b)(a-c)} z, \left. \eta = \frac{-z^2}{(a-b)(a-c)} \right) \nonumber\\
&=&  \sum_{i_0=0}^{q_0} \frac{(-q_0)_{i_0} \left( q_0  \right)_{i_0}}{\left(1 \right)_{i_0}\left(\frac{1}{2} \right)_{i_0}} \rho ^{i_0} \nonumber\\
&+& \left\{ \sum_{i_0=0}^{q_0}\frac{\left( i_0+\frac{\alpha}{2}+\frac{1}{2} \right)\left( i_0-\frac{\alpha}{2} \right)}{ \left(i_0+ 2 \right) \left( i_0+ \frac{3}{2} \right)}\frac{(-q_0)_{i_0} \left( q_0 \right)_{i_0}}{\left( 1 \right)_{i_0}\left(\frac{1}{2} \right)_{i_0}} \right. \left. \sum_{i_1=i_0}^{q_1} \frac{(-q_1)_{i_1}\left( q_1 +4 \right)_{i_1}\left( 3 \right)_{i_0}\left( \frac{5}{2} \right)_{i_0}}{(-q_1)_{i_0}\left( q_1 +4 \right)_{i_0}\left( 3 \right)_{i_1}\left( \frac{5}{2} \right)_{i_1}} \rho ^{i_1}\right\} \eta \nonumber\\
&+& \sum_{n=2}^{\infty } \left\{ \sum_{i_0=0}^{q_0}\frac{\left( i_0+\frac{\alpha}{2}+\frac{1}{2} \right)\left( i_0-\frac{\alpha}{2} \right)}{ \left(i_0+ 2 \right) \left( i_0+ \frac{3}{2} \right)}\frac{(-q_0)_{i_0} \left( q_0 \right)_{i_0}}{\left(1 \right)_{i_0}\left(\frac{1}{2} \right)_{i_0}} \right.\nonumber\\
&\times& \prod _{k=1}^{n-1} \left\{ \sum_{i_k=i_{k-1}}^{q_k} \frac{\left( i_k+2k+\frac{\alpha}{2}+\frac{1}{2} \right)\left( i_k+2k-\frac{\alpha}{2} \right)}{ \left(i_k+ 2k+2 \right) \left( i_k+2k +\frac{3}{2} \right)}\right.  \left.\frac{(-q_k)_{i_k}\left( q_k+4k \right)_{i_k}\left( 2k+1 \right)_{i_{k-1}}\left( 2k+\frac{1}{2} \right)_{i_{k-1}}}{(-q_k)_{i_{k-1}}\left( q_k+4k \right)_{i_{k-1}}\left( 2k+1 \right)_{i_k}\left( 2k+\frac{1}{2} \right)_{i_k}}\right\} \nonumber\\
&\times&  \left.\sum_{i_n= i_{n-1}}^{q_n} \frac{(-q_n)_{i_n}\left( q_n+4n \right)_{i_n}\left( 2n +1 \right)_{i_{n-1}}\left( 2n+\frac{1}{2} \right)_{i_{n-1}}}{(-q_n)_{i_{n-1}}\left( q_n+4n  \right)_{i_{n-1}}\left( 2n +1 \right)_{i_n}\left( 2n +\frac{1}{2} \right)_{i_n}} \rho ^{i_n} \right\} \eta ^n  \label{eq:70011}
\end{eqnarray}
\end{remark}
For the minimum value of Lame equation in the algebraic form of the first kind for a polynomial of type 2 about $x=a$, put $q_0=q_1=q_2=\cdots=0$ in (\ref{eq:70011}).
\begin{eqnarray}
y(z)&=& LF_{0}^R\left( a, b, c, \alpha, q= \alpha (\alpha +1)a- 16(2a-b-c)j^2; z= x-a, \rho = -\frac{2a-b-c}{(a-b)(a-c)} z \right. \nonumber\\
&&, \left. \eta = \frac{-z^2}{(a-b)(a-c)} \right) = \; _2F_1\left( -\frac{\alpha }{4},\frac{\alpha }{4}+\frac{1}{4},\frac{3}{4},\eta \right) \hspace{1cm}\mbox{where}\;\;|\eta |<1 \label{ysc:7001}
\end{eqnarray} 
For the special case, if $\eta =1$ in (\ref{ysc:7001}),
\begin{eqnarray}
y(z)&=& LF_{0}^R\Bigg( a, b, c, \alpha, q= \alpha (\alpha +1)a- 16(2a-b-c)j^2; z= \sqrt{-(a-b)(a-c)} \nonumber\\
&&, \rho = -\frac{2a-b-c}{\sqrt{-(a-b)(a-c)}}, \eta =1 \Bigg) =  \frac{\sqrt{\pi}\Gamma \left( \frac{3}{4}\right)}{\Gamma \left( \frac{2-\alpha }{4}\right)\Gamma \left( \frac{3+\alpha }{4}\right)}\label{ysc:7002}
\end{eqnarray}
\begin{remark}
The power series expansion of Lame equation in the algebraic form of the second kind for polynomial of type 2 about $x=a$ as $q= \alpha (\alpha +1)a- 4(2a-b-c)\left( q_j+2j+\frac{1}{2} \right)^2$ where $j,q_j \in \mathbb{N}_{0}$ is
\begin{eqnarray}
y(z)&=& LS_{q_j}^R\left( a, b, c, \alpha, q= \alpha (\alpha +1)a- 4(2a-b-c)\left( q_j+2j+\frac{1}{2} \right)^2; z= x-a \right. \nonumber\\
&&,\left. \rho = -\frac{2a-b-c}{(a-b)(a-c)} z, \eta = \frac{-z^2}{(a-b)(a-c)} \right) \nonumber\\
&=& z^{\frac{1}{2}} \left\{\sum_{i_0=0}^{q_0} \frac{(-q_0)_{i_0} \left( q_0 +1 \right)_{i_0}}{\left(\frac{3}{2} \right)_{i_0}\left(1 \right)_{i_0}} \rho ^{i_0}\right.\nonumber\\
&+& \left\{ \sum_{i_0=0}^{q_0}\frac{\left( i_0+\frac{\alpha}{2}+1 \right)\left( i_0-\frac{\alpha}{2} + \frac{1}{2} \right)}{ \left(i_0+ \frac{5}{2} \right) \left( i_0+ 2 \right)}\frac{(-q_0)_{i_0} \left( q_0+ 1 \right)_{i_0}}{\left(\frac{3}{2} \right)_{i_0}\left( 1\right)_{i_0}} \right.   \left. \sum_{i_1=i_0}^{q_1} \frac{(-q_1)_{i_1}\left( q_1 +5 \right)_{i_1}\left( \frac{7}{2} \right)_{i_0}\left( 3 \right)_{i_0}}{(-q_1)_{i_0}\left( q_1 +5 \right)_{i_0}\left( \frac{7}{2} \right)_{i_1}\left( 3 \right)_{i_1}} \rho ^{i_1}\right\} \eta \nonumber\\
&+& \sum_{n=2}^{\infty } \left\{ \sum_{i_0=0}^{q_0}\frac{\left( i_0+\frac{\alpha}{2}+1 \right)\left( i_0-\frac{\alpha}{2} +\frac{1}{2} \right)}{ \left(i_0+ \frac{5}{2} \right) \left( i_0+ 2 \right)}\frac{(-q_0)_{i_0} \left( q_0+ 1 \right)_{i_0}}{\left(\frac{3}{2} \right)_{i_0}\left(1\right)_{i_0}} \right.\nonumber\\
&\times& \prod _{k=1}^{n-1} \left\{ \sum_{i_k=i_{k-1}}^{q_k} \frac{\left( i_k+2k+\frac{\alpha}{2}+1\right)\left( i_k+2k-\frac{\alpha}{2} + \frac{1}{2} \right)}{ \left(i_k+ 2k+\frac{5}{2} \right) \left( i_k+2k +2 \right)}\right. \nonumber\\
&\times&  \left.\frac{(-q_k)_{i_k}\left( q_k+4k + 1 \right)_{i_k}\left( 2k+\frac{3}{2} \right)_{i_{k-1}}\left( 2k+1 \right)_{i_{k-1}}}{(-q_k)_{i_{k-1}}\left( q_k+4k +1 \right)_{i_{k-1}}\left( 2k+\frac{3}{2} \right)_{i_k}\left( 2k+1\right)_{i_k}}\right\} \nonumber\\
&\times& \left. \left.\sum_{i_n= i_{n-1}}^{q_n} \frac{(-q_n)_{i_n}\left( q_n+4n + 1 \right)_{i_n}\left( 2n +\frac{3}{2} \right)_{i_{n-1}}\left( 2n+1 \right)_{i_{n-1}}}{(-q_n)_{i_{n-1}}\left( q_n+4n +1 \right)_{i_{n-1}}\left( 2n +\frac{3}{2} \right)_{i_n}\left( 2n +1 \right)_{i_n}} \rho ^{i_n} \right\} \eta ^n \right\} \label{eq:70012}
\end{eqnarray}
\end{remark}
For the minimum value of Lame equation in the algebraic form of the second kind for a polynomial of type 2 about $x=a$, put $q_0=q_1=q_2=\cdots=0$ in (\ref{eq:70012}).
\begin{eqnarray}
y(z)&=& LS_{0}^R\left( a, b, c, \alpha, q= \alpha (\alpha +1)a- 4(2a-b-c)\left( 2j+\frac{1}{2} \right)^2; z= x-a, \rho = -\frac{2a-b-c}{(a-b)(a-c)} z \right. \nonumber\\
&&, \left. \eta = \frac{-z^2}{(a-b)(a-c)} \right) =z^{\frac{1}{2}}\; _2F_1\left( -\frac{\alpha }{4}+\frac{1}{4},\frac{\alpha }{4}+\frac{1}{2},\frac{5}{4},\eta \right) \hspace{1cm}\mbox{where}\;\;|\eta |<1 \label{ysc:7003}
\end{eqnarray} 
For the special case, if $\eta =1$ in (\ref{ysc:7003}),
\begin{eqnarray}
y(z)&=& LS_{0}^R\Bigg( a, b, c, \alpha, q= \alpha (\alpha +1)a- 4(2a-b-c)\left( 2j+\frac{1}{2} \right)^2; z= \sqrt{-(a-b)(a-c)} \nonumber\\
&&, \rho = -\frac{2a-b-c}{\sqrt{-(a-b)(a-c)}}, \eta =1 \Bigg) =  \frac{\sqrt{\pi}\Gamma \left( \frac{5}{4}\right)}{\Gamma \left( \frac{3-\alpha }{4}\right)\Gamma \left( \frac{4+\alpha }{4}\right)}\left( -(a-b)(a-c)\right)^{\frac{1}{4}}\label{ysc:7004}
\end{eqnarray}
In (\ref{ysc:7001}) and (\ref{ysc:7003}), a polynomial of type 2 requires $\left| \eta = \frac{-(x-a)^2}{(a-b)(a-c)}\right|<1$ for the convergence of the radius. For huge values $q_j$, Lame functions for a polynomial of type 2 will not convergent at $\eta =1$ any more. Its radius of convergence is $\left| \eta \right|< 1$. For more details about this issue, it is explained in chapter 3 and 4 \cite{zChoun2014}. 

In Ref.\cite{zChou2012f} I treat $\alpha $ as a fixed value and  $q$ as a free variable  to construct the Lame polynomial of type 1: (1) if $\alpha = 2(2\alpha_j +j)$ or  $-2(2\alpha_j +j)-1$ where $j, \alpha_j \in \mathbb{N}_{0}$, an analytic solution of Lame equation turns to be the first kind of independent solution of the Lame polynomial of type 1. (2) if  $\alpha = 2(2\alpha_j +j)+1$ or $-2(2\alpha_j +j+1)$, an analytic solution of Lame equation turns to be the second kind of independent solution of the Lame polynomial  of type 1. 

In this chapter I treat $q$ as a fixed value and $\alpha$ as a free variable  to construct the Lame polynomial of type 2: (1) if $q= \alpha (\alpha +1)a- 4(2a-b-c)(q_j+2j )^2 $ where $j,q_j \in \mathbb{N}_{0}$, an analytic solution of Lame equation turns to be the first kind of independent solution of the Lame polynomial of type 2. (2) if $q= \alpha (\alpha +1)a- 4(2a-b-c)\left( q_j+2j+\frac{1}{2} \right)^2$, an analytic solution of Lame equation turns to be the second kind of independent solution of the Lame polynomial of type 2.
\subsection{Infinite series}
In chapter 1 the general expression of power series of $y(x)$ for infinite series is defined by
\begin{eqnarray}
y(x) &=&  \sum_{n=0}^{\infty } y_{n}(x) = y_0(x)+ y_1(x)+ y_2(x)+y_3(x)+\cdots \nonumber\\
&=& c_0 \Bigg\{ \sum_{i_0=0}^{\infty } \left( \prod _{i_1=0}^{i_0-1}A_{i_1} \right) x^{i_0+\lambda }
+ \sum_{i_0=0}^{\infty }\left\{ B_{i_0+1} \prod _{i_1=0}^{i_0-1}A_{i_1}  \sum_{i_2=i_0}^{\infty } \left( \prod _{i_3=i_0}^{i_2-1}A_{i_3+2} \right)\right\} x^{i_2+2+\lambda }  \nonumber\\
&& + \sum_{N=2}^{\infty } \Bigg\{ \sum_{i_0=0}^{\infty } \Bigg\{B_{i_0+1}\prod _{i_1=0}^{i_0-1} A_{i_1} 
\prod _{k=1}^{N-1} \Bigg( \sum_{i_{2k}= i_{2(k-1)}}^{\infty } B_{i_{2k}+2k+1}\prod _{i_{2k+1}=i_{2(k-1)}}^{i_{2k}-1}A_{i_{2k+1}+2k}\Bigg)\nonumber\\
&& \times  \sum_{i_{2N} = i_{2(N-1)}}^{\infty } \Bigg( \prod _{i_{2N+1}=i_{2(N-1)}}^{i_{2N}-1} A_{i_{2N+1}+2N} \Bigg) \Bigg\} \Bigg\} x^{i_{2N}+2N+\lambda }\Bigg\}   \label{eq:70013}
\end{eqnarray}
Substitute (\ref{eq:7007a})--(\ref{eq:7007c}) into (\ref{eq:70013}) with replacing variable $x$ by $z$. 
The general expression of power series of Lame equation in the algebraic form for infinite series about $x=a$ is given by
\begin{eqnarray}
  y(z)&=&\sum_{n=0}^{\infty } y_n(z)= y_0(z)+ y_1(z)+ y_2(z)+ y_3(z)+\cdots \nonumber\\
&=& c_0 z^{\lambda } \left\{\sum_{i_0=0}^{\infty } \frac{\left( -\varphi +\lambda \right)_{i_0} \left(  \varphi+\lambda \right)_{i_0}}{\left( 1+\lambda \right)_{i_0}\left(\frac{1}{2} +\lambda \right)_{i_0}} \rho ^{i_0}\right.\nonumber\\
&&+ \left\{ \sum_{i_0=0}^{\infty }\frac{\left( i_0+ \frac{\alpha}{2}+\frac{1}{2}+\lambda \right) \left( i_0- \frac{\alpha}{2} +\lambda \right)}{\left( i_0+ 2+ \lambda \right)\left( i_0+ \frac{3}{2}+\lambda \right)}\frac{\left(-\varphi +\lambda \right)_{i_0} \left( \varphi+\lambda \right)_{i_0}}{\left( 1+\lambda \right)_{i_0}\left(\frac{1}{2} +\lambda \right)_{i_0}} \right.\nonumber\\
&&\times \left.\sum_{i_1=i_0}^{\infty } \frac{\left(-\varphi +2+\lambda \right)_{i_1} \left( \varphi +2+\lambda \right)_{i_1}\left( 3+\lambda \right)_{i_0}\left(\frac{5}{2} +\lambda \right)_{i_0}}{\left(-\varphi +2+\lambda \right)_{i_0} \left( \varphi +2+\lambda \right)_{i_0}\left( 3+\lambda \right)_{i_1}\left(\frac{5}{2} +\lambda \right)_{i_1}}\rho ^{i_1}\right\} \eta \nonumber\\
&&+ \sum_{n=2}^{\infty } \left\{ \sum_{i_0=0}^{\infty } \frac{\left( i_0+ \frac{\alpha}{2}+\frac{1}{2}+\lambda \right) \left( i_0- \frac{\alpha}{2} +\lambda \right)}{\left( i_0+ 2+ \lambda \right)\left( i_0+ \frac{3}{2}+\lambda \right)}\frac{\left(-\varphi +\lambda \right)_{i_0} \left( \varphi+\lambda \right)_{i_0}}{\left( 1+\lambda \right)_{i_0}\left(\frac{1}{2} +\lambda \right)_{i_0}}\right.\nonumber\\
&&\times \prod _{k=1}^{n-1} \left\{ \sum_{i_k=i_{k-1}}^{\infty } \frac{\left( i_k+ 2k+\frac{\alpha }{2}+\frac{1}{2}+\lambda \right) \left( i_k+ 2k-\frac{\alpha }{2} +\lambda \right)}{\left( i_k+ 2k+2+\lambda \right)\left( i_k+ 2k+\frac{3}{2}+\lambda \right)}\right.\nonumber\\
&&\times \left.\frac{ \left(-\varphi +2k+\lambda \right)_{i_k} \left( \varphi +2k+\lambda \right)_{i_k}\left( 2k+1+\lambda \right)_{i_{k-1}}\left( 2k+\frac{1}{2} +\lambda \right)_{i_{k-1}}}{\left(-\varphi +2k+\lambda \right)_{i_{k-1}} \left( \varphi +2k+\lambda \right)_{i_{k-1}}\left( 2k+1+\lambda \right)_{i_k}\left( 2k+\frac{1}{2} +\lambda \right)_{i_k}}\right\}\nonumber\\
&&\times \left.\left.\sum_{i_n= i_{n-1}}^{\infty } \frac{ \left( -\varphi +2n+\lambda \right)_{i_n} \left( \varphi +2n+\lambda \right)_{i_n}\left( 2n+1+\lambda \right)_{i_{n-1}}\left( 2n+\frac{1}{2} +\lambda \right)_{i_{n-1}}}{\left(-\varphi +2n+\lambda \right)_{i_{n-1}} \left( \varphi +2n+\lambda \right)_{i_{n-1}}\left( 2n+1+\lambda \right)_{i_n}\left( 2n+\frac{1}{2} +\lambda \right)_{i_n}} \rho ^{i_n} \right\} \eta ^n \right\} \hspace{1.5cm} \label{eq:70014}
\end{eqnarray}
Put $c_0$= 1 as $\lambda =0$ for the first kind of independent solutions of Lame equation and $\lambda =\frac{1}{2}$ for the second one in (\ref{eq:70014}).
\begin{remark}
The power series expansion of Lame equation in the algebraic form of the first kind for infinite series about $x=a$ using R3TRF is
\begin{eqnarray}
 y(z)&=& LF^R\left( a, b, c, \alpha, q, \varphi = \sqrt{\frac{\alpha (\alpha +1)a-q}{4(2a-b-c)}}; z= x-a, \rho = -\frac{2a-b-c}{(a-b)(a-c)} z, \eta = \frac{-z^2}{(a-b)(a-c)} \right) \nonumber\\
&=& \sum_{i_0=0}^{\infty } \frac{\left(-\varphi \right)_{i_0} \left( \varphi \right)_{i_0}}{\left( 1 \right)_{i_0}\left(\frac{1}{2} \right)_{i_0}} \rho ^{i_0} \nonumber\\
&&+ \left\{ \sum_{i_0=0}^{\infty }\frac{\left( i_0+ \frac{\alpha}{2}+\frac{1}{2} \right) \left( i_0- \frac{\alpha}{2} \right)}{\left( i_0+ 2 \right)\left( i_0+ \frac{3}{2} \right)}\frac{\left(-\varphi \right)_{i_0} \left( \varphi \right)_{i_0}}{\left( 1 \right)_{i_0}\left(\frac{1}{2} \right)_{i_0}} \right.  \left.\sum_{i_1=i_0}^{\infty } \frac{\left(-\varphi +2 \right)_{i_1} \left( \varphi +2 \right)_{i_1}\left( 3 \right)_{i_0}\left(\frac{5}{2} \right)_{i_0}}{\left(-\varphi +2 \right)_{i_0} \left( \varphi +2 \right)_{i_0}\left( 3 \right)_{i_1}\left(\frac{5}{2} \right)_{i_1}}\rho ^{i_1}\right\} \eta \nonumber\\
&&+ \sum_{n=2}^{\infty } \left\{ \sum_{i_0=0}^{\infty } \frac{\left( i_0+ \frac{\alpha}{2}+\frac{1}{2} \right) \left( i_0- \frac{\alpha}{2} \right)}{\left( i_0+ 2 \right)\left( i_0+ \frac{3}{2} \right)}\frac{\left(-\varphi \right)_{i_0} \left( \varphi \right)_{i_0}}{\left( 1 \right)_{i_0}\left(\frac{1}{2} \right)_{i_0}}\right.\nonumber\\
&&\times \prod _{k=1}^{n-1} \left\{ \sum_{i_k=i_{k-1}}^{\infty } \frac{\left( i_k+ 2k+\frac{\alpha }{2}+\frac{1}{2} \right) \left( i_k+ 2k-\frac{\alpha }{2} \right)}{\left( i_k+ 2k+2 \right)\left( i_k+ 2k+\frac{3}{2} \right)}\right. \left.\frac{ \left(-\varphi +2k \right)_{i_k} \left( \varphi +2k \right)_{i_k}\left( 2k+1 \right)_{i_{k-1}}\left( 2k+\frac{1}{2} \right)_{i_{k-1}}}{\left(-\varphi +2k \right)_{i_{k-1}} \left( \varphi +2k \right)_{i_{k-1}}\left( 2k+1 \right)_{i_k}\left( 2k+\frac{1}{2} \right)_{i_k}}\right\}\nonumber\\
&&\times \left.\sum_{i_n= i_{n-1}}^{\infty } \frac{ \left(-\varphi +2n \right)_{i_n} \left( \varphi +2n \right)_{i_n}\left( 2n+1 \right)_{i_{n-1}}\left( 2n+\frac{1}{2} \right)_{i_{n-1}}}{\left(-\varphi +2n \right)_{i_{n-1}} \left( \varphi +2n \right)_{i_{n-1}}\left( 2n+1 \right)_{i_n}\left( 2n+\frac{1}{2}  \right)_{i_n}} \rho ^{i_n} \right\} \eta ^n \label{eq:70015}
\end{eqnarray}
\end{remark}
\begin{remark}
The power series expansion of Lame equation in the algebraic form of the second kind for infinite series about $x=a$ using R3TRF is
\begin{eqnarray}
 y(z)&=& LS^R\left( a, b, c, \alpha, q, \varphi = \sqrt{\frac{\alpha (\alpha +1)a-q}{4(2a-b-c)}}; z= x-a, \rho = -\frac{2a-b-c}{(a-b)(a-c)} z, \eta = \frac{-z^2}{(a-b)(a-c)} \right) \nonumber\\
&=& z^{\frac{1}{2}} \left\{\sum_{i_0=0}^{\infty } \frac{\left(-\varphi +\frac{1}{2} \right)_{i_0} \left( \varphi +\frac{1}{2} \right)_{i_0}}{\left(\frac{3}{2} \right)_{i_0}\left(1 \right)_{i_0}} \rho ^{i_0}\right.\nonumber\\
&&+ \left\{ \sum_{i_0=0}^{\infty }\frac{\left( i_0+ \frac{\alpha}{2}+1 \right) \left( i_0- \frac{\alpha}{2} +\frac{1}{2} \right)}{\left( i_0 + \frac{5}{2}\right)\left( i_0+ 2 \right)}\frac{\left(-\varphi +\frac{1}{2} \right)_{i_0} \left( \varphi +\frac{1}{2} \right)_{i_0}}{\left( \frac{3}{2} \right)_{i_0}\left( 1 \right)_{i_0}} \right. \nonumber\\
&&\times \left.\sum_{i_1=i_0}^{\infty } \frac{\left(-\varphi  +\frac{5}{2} \right)_{i_1} \left( \varphi +\frac{5}{2} \right)_{i_1}\left( \frac{7}{2}\right)_{i_0}\left( 3 \right)_{i_0}}{\left(-\varphi +\frac{5}{2} \right)_{i_0} \left( \varphi +\frac{5}{2} \right)_{i_0}\left( \frac{7}{2} \right)_{i_1}\left( 3 \right)_{i_1}}\rho ^{i_1}\right\} \eta \nonumber\\
&&+ \sum_{n=2}^{\infty } \left\{ \sum_{i_0=0}^{\infty } \frac{\left( i_0+ \frac{\alpha}{2}+1 \right) \left( i_0- \frac{\alpha}{2} +\frac{1}{2} \right)}{\left( i_0 + \frac{5}{2} \right)\left( i_0+ 2 \right)}\frac{\left(-\varphi +\frac{1}{2} \right)_{i_0} \left( \varphi +\frac{1}{2} \right)_{i_0}}{\left(  \frac{3}{2} \right)_{i_0}\left( 1 \right)_{i_0}}\right.\nonumber\\
&&\times \prod _{k=1}^{n-1} \left\{ \sum_{i_k=i_{k-1}}^{\infty } \frac{\left( i_k+ 2k+\frac{\alpha }{2} +1 \right) \left( i_k+ 2k-\frac{\alpha }{2} +\frac{1}{2} \right)}{\left( i_k+ 2k +\frac{5}{2} \right)\left( i_k+ 2k+2 \right)}\right. \nonumber\\
&&\times \left.\frac{ \left(-\varphi +2k+\frac{1}{2} \right)_{i_k} \left( \varphi +2k+\frac{1}{2} \right)_{i_k}\left( 2k +\frac{3}{2}\right)_{i_{k-1}}\left( 2k+1 \right)_{i_{k-1}}}{\left(-\varphi +2k+\frac{1}{2} \right)_{i_{k-1}} \left( \varphi +2k+\frac{1}{2} \right)_{i_{k-1}}\left( 2k +\frac{3}{2}\right)_{i_k}\left( 2k+1 \right)_{i_k}}\right\}\nonumber\\
&&\times \left.\left.\sum_{i_n= i_{n-1}}^{\infty } \frac{ \left(-\varphi +2n+\frac{1}{2} \right)_{i_n} \left( \varphi +2n+\frac{1}{2} \right)_{i_n}\left( 2n+\frac{3}{2} \right)_{i_{n-1}}\left( 2n+1 \right)_{i_{n-1}}}{\left(-\varphi +2n+\frac{1}{2} \right)_{i_{n-1}} \left( \varphi +2n+\frac{1}{2} \right)_{i_{n-1}}\left( 2n+\frac{3}{2} \right)_{i_n}\left( 2n+1 \right)_{i_n}} \rho ^{i_n} \right\} \eta ^n \right\}\label{eq:70016}
\end{eqnarray}
\end{remark}
The infinite series in this chapter are equivalent to the infinite series in Ref.\cite{zChou2012f}. In this chapter $B_n$ is the leading term in sequence $c_n$ of the analytic function $y(z)$. In Ref.\cite{zChou2012f} $A_n$ is the leading term in sequence $c_n$ of the analytic function $y(z)$.
%
\section{Integral formalism}
Lame equation could not be built in a definite or contour integral form of any well-known simple functions such as Gauss hypergeometric, Kummer functions and etc because of a 3-term recursive relation in its power series. The three term recurrence relation in the power series expansion of Lame equation brings mathematical difficulty to be derived into a direct or contour integral.
 
In place of analyzing Lame equation into its integral representation of any simple functions, in earlier literature the integral equations of Lame equation were constructed by using simple kernels involving Legendre functions of the Jacobian elliptic function or ellipsoidal harmonic functions: such integral relationships express one analytic solution in terms of another analytic solution involving Jacobian elliptic functions. There are many other forms of integral relations in Lame equation.\cite{zErde1955,zArsc1964a,zShai1980,zSlee1968a,zVolk1982,zVolk1983,zVolk1984,zWang1989,zWhit1996}

In Ref.\cite{zChou2012f}, I show the integral form (each sub-integral $y_m(x)$ where $m=0,1,2,\cdots$ is composed of $2m$ terms of definite integrals and $m$ terms of contour integrals) of Lame equation using 3TRF: a $_2F_1$ function recurs in each of sub-integral forms of Lame equation.
\subsection{Polynomial of type 2}

Now I consider the combined definite and contour integral representation of Lame equation by using R3TRF.
There is a generalized hypergeometric function which is given by
\begin{eqnarray}
I_l &=& \sum_{i_l= i_{l-1}}^{q_l} \frac{(-q_l)_{i_l}\left( q_l+4l +2\lambda \right)_{i_l}\left(2l+1+\lambda \right)_{i_{l-1}}\left( 2l+\frac{1}{2} +\lambda \right)_{i_{l-1}}}{(-q_l)_{i_{l-1}}\left( q_l+4l +2\lambda \right)_{i_{l-1}}\left(2l+1+\lambda \right)_{i_l}\left( 2l+\frac{1}{2} +\lambda \right)_{i_l}} \rho ^{i_l}\label{eq:70018}\\
&=& \rho ^{i_{l-1}} 
\sum_{j=0}^{\infty } \frac{B(i_{l-1}+2l+\lambda ,j+1) B\left(i_{l-1}+2l- \frac{1}{2} +\lambda ,j+1\right)(i_{l-1}-q_l)_j \left( i_{l-1}+q_l+4l +2\lambda \right)_j}{(i_{l-1}+2l+\lambda )^{-1}\left( i_{l-1}+2l-\frac{1}{2} +\lambda \right)^{-1}(1)_j \;j!} \rho ^j\nonumber
\end{eqnarray}
By using integral form of beta function,
\begin{subequations}
\begin{equation}
B\left(i_{l-1}+2l+\lambda ,j+1\right)= \int_{0}^{1} dt_l\;t_l^{i_{l-1}+2l-1+\lambda } (1-t_l)^j \label{eq:70019a}
\end{equation}
\begin{equation}
B\left(i_{l-1}+2l-\frac{1}{2} +\lambda ,j+1\right)= \int_{0}^{1} du_l\;u_l^{i_{l-1}+2l-\frac{3}{2} +\lambda } (1-u_l)^j\label{eq:70019b}
\end{equation}
\end{subequations}
Substitute (\ref{eq:70019a}) and (\ref{eq:70019b}) into (\ref{eq:70018}). And divide $(i_{l-1}+2l+\lambda )\left( i_{l-1}+2l-\frac{1}{2} +\lambda \right)$ into $I_l$.
\begin{eqnarray}
K_l&=& \frac{1}{(i_{l-1}+2l+\lambda )\left( i_{l-1}+2l-\frac{1}{2} +\lambda \right)}\nonumber\\
&&\times \sum_{i_l= i_{l-1}}^{q_l} \frac{(-q_l)_{i_l}\left( q_l+4l +2\lambda \right)_{i_l}\left(2l+1+\lambda \right)_{i_{l-1}}\left( 2l+\frac{1}{2} +\lambda \right)_{i_{l-1}}}{(-q_l)_{i_{l-1}}\left( q_l+4l +2\lambda \right)_{i_{l-1}}\left(2l+1+\lambda \right)_{i_l}\left( 2l+\frac{1}{2} +\lambda \right)_{i_l}} \rho ^{i_l}\nonumber\\
&=&  \int_{0}^{1} dt_l\;t_l^{2l-1+\lambda } \int_{0}^{1} du_l\;u_l^{2l-\frac{3}{2}+\lambda } (\rho t_l u_l)^{i_{l-1}}   \sum_{j=0}^{\infty } \frac{(i_{l-1}-q_l)_j \left( i_{l-1}+q_l+4l +2\lambda \right)_j}{(1)_j \;j!} 
(\rho (1-t_l)(1-u_l))^j \nonumber 
\end{eqnarray}
The integral form of Gauss hypergeometric function is written by
\begin{eqnarray}
_2F_1 \left( \alpha ,\beta ; \gamma ; z \right) &=& \sum_{n=0}^{\infty } \frac{(\alpha )_n (\beta )_n}{(\gamma )_n (n!)} z^n \nonumber\\
&=& -\frac{1}{2\pi i} \frac{\Gamma(1-\alpha ) \Gamma(\gamma )}{\Gamma (\gamma -\alpha )} \oint dv_l\;(-v_l)^{\alpha -1} (1-v_l)^{\gamma -\alpha -1} (1-zv_l)^{-\beta }\hspace{1cm}\label{eq:70020}\\
&& \mbox{where} \;\mbox{Re}(\gamma -\alpha )>0 \nonumber
\end{eqnarray}
replaced $\alpha $, $\beta $, $\gamma $ and $z$ by $i_{l-1}-q_l$, $ { \displaystyle i_{l-1}+q_l+4l +2\lambda }$, 1 and $\rho  (1-t_l)(1-u_l)$ in (\ref{eq:70020})
\begin{eqnarray}
&& \sum_{j=0}^{\infty } \frac{\left(i_{l-1}-q_l)_j (i_{l-1}+q_l+4l +2\lambda \right)_j}{(1)_j \;j!} (\rho (1-t_l)(1-u_l))^j \nonumber\\
&=& \frac{1}{2\pi i} \oint dv_l\;\frac{1}{v_l} \left(\frac{v_l-1}{v_l}\frac{1}{1-\rho (1-t_l)(1-u_l)v_l}\right)^{q_l} (1-\rho (1-t_l)(1-u_l)v_l)^{-\left( 4l +2\lambda \right)}\nonumber\\
&&\times \left(\frac{v_l}{v_l-1} \frac{1}{1-\rho (1-t_l)(1-u_l)v_l}\right)^{i_{l-1}} \label{eq:70021}
\end{eqnarray}
Substitute (\ref{eq:70021}) into $K_l$.
\begin{eqnarray}
K_l&=& \frac{1}{(i_{l-1}+2l+\lambda )\left( i_{l-1}+2l-\frac{1}{2} +\lambda \right)}\nonumber\\
&&\times \sum_{i_l= i_{l-1}}^{q_l} \frac{(-q_l)_{i_l}\left( q_l+4l +2\lambda \right)_{i_l}\left(2l+1+\lambda \right)_{i_{l-1}}\left( 2l+\frac{1}{2} +\lambda \right)_{i_{l-1}}}{(-q_l)_{i_{l-1}}\left( q_l+4l +2\lambda \right)_{i_{l-1}}\left(2l+1+\lambda \right)_{i_l}\left( 2l+\frac{1}{2} +\lambda \right)_{i_l}} \rho ^{i_l}\nonumber\\
&=&  \int_{0}^{1} dt_l\;t_l^{2l-1+\lambda } \int_{0}^{1} du_l\;u_l^{2l-\frac{3}{2} +\lambda } 
\frac{1}{2\pi i} \oint dv_l\;\frac{1}{v_l} \left(\frac{v_l-1}{v_l}\frac{1}{1-\rho (1-t_l)(1-u_l)v_l}\right)^{q_l} \nonumber\\
&&\times (1-\rho (1-t_l)(1-u_l)v_l)^{-\left( 4l +2\lambda \right)} \left(\frac{v_l}{v_l-1} \frac{\rho t_l u_l}{1-\rho (1-t_l)(1-u_l)v_l}\right)^{i_{l-1}}\label{eq:70022} 
\end{eqnarray}
Substitute (\ref{eq:70022}) into (\ref{eq:70010}) where $l=1,2,3,\cdots$; apply $K_1$ into the second summation of sub-power series $y_1(z)$, apply $K_2$ into the third summation and $K_1$ into the second summation of sub-power series $y_2(z)$, apply $K_3$ into the forth summation, $K_2$ into the third summation and $K_1$ into the second summation of sub-power series $y_3(z)$, etc.\footnote{$y_1(z)$ means the sub-power series in (\ref{eq:70010}) contains one term of $B_n's$, $y_2(z)$ means the sub-power series in (\ref{eq:70010}) contains two terms of $B_n's$, $y_3(z)$ means the sub-power series in (\ref{eq:70010}) contains three terms of $B_n's$, etc.}
\begin{theorem}
The general representation in the form of integral of the Lame polynomial of type 2 in the algebraic form is given by
\begin{eqnarray}
 y(z)&=& \sum_{n=0}^{\infty } y_{n}(z)= y_0(z)+ y_1(z)+ y_2(z)+y_3(z)+\cdots \nonumber\\
&=& c_0 z^{\lambda } \left\{ \sum_{i_0=0}^{q_0}\frac{(-q_0)_{i_0}\left(q_0 +2\lambda \right)_{i_0}}{(1+\lambda )_{i_0}\left(\frac{1}{2} +\lambda \right)_{i_0}} \rho ^{i_0}\right.\nonumber\\
&&+ \sum_{n=1}^{\infty } \left\{\prod _{k=0}^{n-1} \Bigg\{ \int_{0}^{1} dt_{n-k}\;t_{n-k}^{2(n-k)-1+\lambda } \int_{0}^{1} du_{n-k}\;u_{n-k}^{2(n-k)- \frac{3}{2} +\lambda }\right. \nonumber\\
&&\times  \frac{1}{2\pi i}  \oint dv_{n-k} \frac{1}{v_{n-k}} \left( \frac{v_{n-k}-1}{v_{n-k}} \frac{1}{1-\overleftrightarrow {w}_{n-k+1,n}(1-t_{n-k})(1-u_{n-k})v_{n-k}}\right)^{q_{n-k}} \nonumber\\
&&\times \left( 1- \overleftrightarrow {w}_{n-k+1,n}(1-t_{n-k})(1-u_{n-k})v_{n-k}\right)^{-\left( 4(n-k) +2\lambda \right)}  \nonumber\\
&&\times \overleftrightarrow {w}_{n-k,n}^{-(2(n-k)- \frac{3}{2}+\frac{\alpha }{2} +\lambda )}\left(  \overleftrightarrow {w}_{n-k,n} \partial _{ \overleftrightarrow {w}_{n-k,n}}\right) \overleftrightarrow {w}_{n-k,n}^{\alpha +\frac{1}{2}}\left(  \overleftrightarrow {w}_{n-k,n} \partial _{ \overleftrightarrow {w}_{n-k,n}}\right) \overleftrightarrow {w}_{n-k,n}^{2(n-k)-2-\frac{\alpha}{2} +\lambda } \Bigg\} \nonumber\\
&&\times \left.\left.\sum_{i_0=0}^{q_0}\frac{(-q_0)_{i_0}\left(q_0 +2\lambda \right)_{i_0}}{(1+\lambda )_{i_0}\left(\frac{1}{2} +\lambda \right)_{i_0}} \overleftrightarrow {w}_{1,n}^{i_0}\right\} \eta ^n \right\}  \label{eq:70023}
\end{eqnarray}
where
\begin{equation}\overleftrightarrow {w}_{i,j}=
\begin{cases} \displaystyle {\frac{v_i}{(v_i-1)}\; \frac{\overleftrightarrow w_{i+1,j} t_i u_i}{1- \overleftrightarrow w_{i+1,j} (1-t_i)(1-u_i) v_i}} \;\;\mbox{where}\; i\leq j\cr
\rho  \;\;\mbox{only}\;\mbox{if}\; i>j
\end{cases}\nonumber 
\end{equation}
In the above, the first sub-integral form contains one term of $B_n's$, the second one contains two terms of $B_n$'s, the third one contains three terms of $B_n$'s, etc.
\end{theorem}
\begin{proof} 
According to (\ref{eq:70010}), 
\begin{equation}
 y(z)=  \sum_{n=0}^{\infty } y_{n}(z) = y_0(z)+ y_1(z)+ y_2(z)+y_3(z)+\cdots \label{eq:70024}
\end{equation}
In the above, sub-power series $y_0(z) $, $y_1(z)$, $y_2(z)$ and $y_3(z)$ of the Lame polynomial using R3TRF about $x=a$ are given by
\begin{subequations}
\begin{equation}
 y_0(z)= c_0 z^{\lambda } \sum_{i_0=0}^{q_0} \frac{(-q_0)_{i_0} \left(q_0 +2\lambda \right)_{i_0}}{\left( 1+\lambda \right)_{i_0}\left( \frac{1}{2} +\lambda \right)_{i_0}} \rho ^{i_0} \label{eq:70025a}
\end{equation}
\begin{eqnarray}
 y_1(z) &=& c_0 z^{\lambda } \left\{\sum_{i_0=0}^{q_0}\frac{\left( i_0 +\frac{\alpha}{2}+ \frac{1}{2}+\lambda  \right) \left( i_0-\frac{\alpha}{2}+\lambda  \right)}{\left(i_0 +2+ \lambda\right) \left( i_0+ \frac{3}{2}+\lambda  \right)}\frac{(-q_0)_{i_0} \left(q_0 +2\lambda \right)_{i_0}}{\left( 1+\lambda \right)_{i_0}\left( \frac{1}{2} +\lambda \right)_{i_0}}\right. \nonumber\\
&&\times \left. \sum_{i_1=i_0}^{q_1} \frac{(-q_1)_{i_1}\left( q_1+4 +2\lambda \right)_{i_1}\left( 3+\lambda \right)_{i_0}\left( \frac{5}{2} +\lambda \right)_{i_0}}{(-q_1)_{i_0}\left( q_1+4 +2\lambda \right)_{i_0}\left( 3+\lambda \right)_{i_1}\left( \frac{5}{2} +\lambda \right)_{i_1}} \rho ^{i_1} \right\}\eta  \label{eq:70025b}
\end{eqnarray}
\begin{eqnarray}
 y_2(z) &=& c_0 z^{\lambda } \left\{\sum_{i_0=0}^{q_0}\frac{\left( i_0+\frac{\alpha}{2}+ \frac{1}{2} +\lambda  \right) \left( i_0-\frac{\alpha}{2}+\lambda  \right)}{\left(i_0 +2+ \lambda\right) \left( i_0+ \frac{3}{2}+\lambda  \right)}\frac{(-q_0)_{i_0} \left(q_0 +2\lambda \right)_{i_0}}{\left( 1+\lambda \right)_{i_0}\left( \frac{1}{2} +\lambda \right)_{i_0}}\right. \nonumber\\
&&\times  \sum_{i_1=i_0}^{q_1} \frac{\left( i_1 +\frac{\alpha}{2}+ \frac{5}{2}+\lambda  \right) \left( i_1-\frac{\alpha}{2}+2+\lambda  \right)}{\left(i_1 +4+ \lambda\right) \left( i_1+ \frac{7}{2}+\lambda  \right)} \frac{(-q_1)_{i_1}\left( q_1+4 +2\lambda \right)_{i_1}\left( 3+\lambda \right)_{i_0}\left( \frac{5}{2} +\lambda \right)_{i_0}}{(-q_1)_{i_0}\left( q_1+4 +2\lambda \right)_{i_0}\left( 3+\lambda \right)_{i_1}\left( \frac{5}{2} +\lambda \right)_{i_1}} \nonumber\\
&&\times \left.\sum_{i_2=i_1}^{q_2} \frac{(-q_2)_{i_2}\left( q_2+8 +2\lambda \right)_{i_2}\left( 5+\lambda \right)_{i_1}\left(\frac{9}{2}+\lambda \right)_{i_1}}{(-q_2)_{i_1}\left( q_2+8 +2\lambda \right)_{i_1}\left( 5+\lambda \right)_{i_2}\left(\frac{9}{2}+\lambda \right)_{i_2}} \rho ^{i_2} \right\} \eta ^2  \label{eq:70025c}
\end{eqnarray}
\begin{eqnarray}
 y_3(z) &=&  c_0 z^{\lambda } \left\{\sum_{i_0=0}^{q_0}\frac{\left( i_0+\frac{\alpha}{2}+ \frac{1}{2} +\lambda  \right) \left( i_0-\frac{\alpha}{2}+\lambda  \right)}{\left(i_0 +2+ \lambda\right) \left( i_0+ \frac{3}{2}+\lambda  \right)}\frac{(-q_0)_{i_0} \left(q_0 +2\lambda \right)_{i_0}}{\left( 1+\lambda \right)_{i_0}\left( \frac{1}{2} +\lambda \right)_{i_0}} \right. \nonumber\\
&&\times  \sum_{i_1=i_0}^{q_1} \frac{\left( i_1 +\frac{\alpha}{2}+ \frac{5}{2}+\lambda  \right) \left( i_1-\frac{\alpha}{2}+2+\lambda  \right)}{\left(i_1 +4+ \lambda\right) \left( i_1+ \frac{7}{2}+\lambda  \right)} \frac{(-q_1)_{i_1}\left( q_1+4 +2\lambda \right)_{i_1}\left( 3+\lambda \right)_{i_0}\left( \frac{5}{2} +\lambda \right)_{i_0}}{(-q_1)_{i_0}\left( q_1+4 +2\lambda \right)_{i_0}\left( 3+\lambda \right)_{i_1}\left( \frac{5}{2} +\lambda \right)_{i_1}} \nonumber\\
&&\times \sum_{i_2=i_1}^{q_2} \frac{\left( i_2 +\frac{\alpha}{2}+ \frac{9}{2}+\lambda  \right) \left( i_2-\frac{\alpha}{2}+4+\lambda  \right)}{\left(i_2 +6+ \lambda\right) \left( i_2+ \frac{11}{2}+\lambda  \right)}  \frac{(-q_2)_{i_2}\left( q_2+8 +2\lambda \right)_{i_2}\left( 5+\lambda \right)_{i_1}\left(\frac{9}{2}+\lambda \right)_{i_1}}{(-q_2)_{i_1}\left( q_2+8 +2\lambda \right)_{i_1}\left( 5+\lambda \right)_{i_2}\left(\frac{9}{2}+\lambda \right)_{i_2}} \nonumber\\
&&\times \left.\sum_{i_3=i_2}^{q_3}\frac{(-q_3)_{i_3}\left( q_3+12 +2\lambda \right)_{i_3}\left( 7+\lambda \right)_{i_2}\left(\frac{13}{2}+\lambda \right)_{i_2}}{(-q_3)_{i_2}\left( q_3+12 +2\lambda \right)_{i_3}\left( 7+\lambda \right)_{i_3}\left(\frac{13}{2}+\lambda \right)_{i_3}}\rho ^{i_3} \right\} \eta ^3  \label{eq:70025d} 
\end{eqnarray}
\end{subequations}
Put $l=1$ in (\ref{eq:70022}). Take the new (\ref{eq:70022}) into (\ref{eq:70025b}).
\begin{eqnarray}
y_1(z) &=& \int_{0}^{1} dt_1\;t_1^{1+\lambda } \int_{0}^{1} du_1\;u_1^{\frac{1}{2} +\lambda } \frac{1}{2\pi i} \oint dv_1 \;\frac{1}{v_1} 
\left( \frac{v_1-1}{v_1} \frac{1}{1-\rho (1-t_1)(1-u_1)v_1}\right)^{q_1}  \nonumber\\
&&\times (1-\rho (1-t_1)(1-u_1)v_1)^{-\left( 4 +2\lambda \right)} \overleftrightarrow {w}_{1,1}^{-(\frac{1}{2}+\frac{\alpha }{2}+\lambda )} \left(\overleftrightarrow {w}_{1,1} \partial_{\overleftrightarrow {w}_{1,1}} \right) \overleftrightarrow {w}_{1,1}^{\alpha +\frac{1}{2} } \left(\overleftrightarrow {w}_{1,1} \partial_{\overleftrightarrow {w}_{1,1}} \right) \overleftrightarrow {w}_{1,1}^{- \frac{\alpha }{2}+\lambda }\nonumber\\
&&\times \left\{ c_0 z^{\lambda } \sum_{i_0=0}^{q_0} \frac{(-q_0)_{i_0} \left(q_0 +2\lambda \right)_{i_0}}{\left( 1+\lambda \right)_{i_0}\left( \frac{1}{2} +\lambda \right)_{i_0}} \overleftrightarrow {w}_{1,1} ^{i_0}\right\}\eta  \label{eq:70026}
\end{eqnarray}
where
\begin{equation}
\overleftrightarrow {w}_{1,1} = \frac{v_1}{v_1-1} \frac{\rho  t_1 u_1}{1-\rho (1-t_1)(1-u_1)v_1}\nonumber
\end{equation}
Put $l=2$ in (\ref{eq:70022}). Take the new (\ref{eq:70022}) into (\ref{eq:70025c}).
\begin{eqnarray}
y_2(z) &=& c_0 z^{\lambda } \int_{0}^{1} dt_2\;t_2^{3+\lambda } \int_{0}^{1} du_2\;u_2^{\frac{5}{2} +\lambda } \frac{1}{2\pi i} \oint dv_2 \;\frac{1}{v_2} 
\left( \frac{v_2-1}{v_2} \frac{1}{1-\rho (1-t_2)(1-u_2)v_2}\right)^{q_2}  \nonumber\\
&&\times (1-\rho (1-t_2)(1-u_2)v_2)^{-\left( 8 +2 \lambda \right)} 
 \overleftrightarrow {w}_{2,2}^{-(\frac{5}{2}+\frac{\alpha }{2}+\lambda) } \left(\overleftrightarrow {w}_{2,2} \partial_{\overleftrightarrow {w}_{2,2}} \right) \overleftrightarrow {w}_{2,2}^{ \alpha +\frac{1}{2}}\left(\overleftrightarrow {w}_{2,2} \partial_{\overleftrightarrow {w}_{2,2}} \right) \overleftrightarrow {w}_{2,2}^{ 2-\frac{\alpha}{2} +\lambda } \nonumber\\
&&\times \left\{\sum_{i_0=0}^{q_0}\frac{\left( i_0+\frac{\alpha}{2}+ \frac{1}{2} +\lambda  \right) \left( i_0-\frac{\alpha}{2}+\lambda  \right)}{\left(i_0 +2+ \lambda\right) \left( i_0+ \frac{3}{2}+\lambda  \right)}\frac{(-q_0)_{i_0} \left(q_0 +2\lambda \right)_{i_0}}{\left( 1+\lambda \right)_{i_0}\left( \frac{1}{2} +\lambda \right)_{i_0}} \right.\nonumber\\
&&\times \left. \sum_{i_1=i_0}^{q_1} \frac{(-q_1)_{i_1}\left( q_1+4 +2\lambda \right)_{i_1}\left( 3+\lambda \right)_{i_0}\left( \frac{5}{2} +\lambda \right)_{i_0}}{(-q_1)_{i_0}\left( q_1+4 +2\lambda \right)_{i_0}\left( 3+\lambda \right)_{i_1}\left( \frac{5}{2} +\lambda \right)_{i_1}} \overleftrightarrow {w}_{2,2}^{i_1} \right\} \eta ^2 \label{eq:70027}
\end{eqnarray}
where
\begin{equation}
\overleftrightarrow {w}_{2,2} = \frac{v_2}{v_2-1} \frac{\rho  t_2 u_2}{1-\rho (1-t_2)(1-u_2)v_2}\nonumber
\end{equation}
Put $l=1$ and $\rho  = \overleftrightarrow {w}_{2,2}$ in (\ref{eq:70022}). Take the new (\ref{eq:70022}) into (\ref{eq:70027}).
\begin{eqnarray}
y_2(z) &=& c_0 z^{\lambda } \int_{0}^{1} dt_2\;t_2^{3+\lambda } \int_{0}^{1} du_2\;u_2^{\frac{5}{2} +\lambda } \frac{1}{2\pi i} \oint dv_2 \;\frac{1}{v_2} 
\left( \frac{v_2-1}{v_2} \frac{1}{1-\rho (1-t_2)(1-u_2)v_2}\right)^{q_2}  \nonumber\\
&&\times (1-\rho (1-t_2)(1-u_2)v_2)^{-\left( 8 +2 \lambda \right)} 
 \overleftrightarrow {w}_{2,2}^{-(\frac{5}{2}+\frac{\alpha }{2}+\lambda) } \left(\overleftrightarrow {w}_{2,2} \partial_{\overleftrightarrow {w}_{2,2}} \right) \overleftrightarrow {w}_{2,2}^{ \alpha +\frac{1}{2}}\left(\overleftrightarrow {w}_{2,2} \partial_{\overleftrightarrow {w}_{2,2}} \right) \overleftrightarrow {w}_{2,2}^{ 2-\frac{\alpha}{2} +\lambda } \nonumber\\
&&\times \int_{0}^{1} dt_1\;t_1^{1+\lambda } \int_{0}^{1} du_1\;u_1^{\frac{1}{2} +\lambda } \frac{1}{2\pi i} \oint dv_1 \;\frac{1}{v_1} 
\left( \frac{v_1-1}{v_1} \frac{1}{1- \overleftrightarrow {w}_{2,2} (1-t_1)(1-u_1)v_1}\right)^{q_1}  \nonumber\\
&&\times (1- \overleftrightarrow {w}_{2,2} (1-t_1)(1-u_1)v_1)^{-\left( 4 +2 \lambda \right)} \overleftrightarrow {w}_{1,2}^{-(\frac{1}{2}+\frac{\alpha }{2}+\lambda) } \left(\overleftrightarrow {w}_{1,2} \partial_{\overleftrightarrow {w}_{1,2}} \right) \overleftrightarrow {w}_{1,2}^{ \alpha +\frac{1}{2}}\left(\overleftrightarrow {w}_{1,2} \partial_{\overleftrightarrow {w}_{1,2}} \right) \overleftrightarrow {w}_{1,2}^{-\frac{\alpha }{2}+\lambda } \nonumber\\
&&\times \left\{ c_0 z^{\lambda } \sum_{i_0=0}^{q_0} \frac{(-q_0)_{i_0} \left( q_0 +2 \lambda \right)_{i_0}}{(1+\lambda )_{i_0}\left(\frac{1}{2} +\lambda \right)_{i_0}} \overleftrightarrow {w}_{1,2} ^{i_0}\right\} \eta ^2 \label{eq:70028}
\end{eqnarray}
where
\begin{equation}
\overleftrightarrow {w}_{1,2} = \frac{v_1}{v_1-1} \frac{\overleftrightarrow {w}_{2,2} t_1 u_1}{1- \overleftrightarrow {w}_{2,2}(1-t_1)(1-u_1)v_1}\nonumber
\end{equation}
By using similar process for the previous cases of integral forms of $y_1(z)$ and $y_2(z)$, the integral form of sub-power series expansion of $y_3(z)$ is given by
\begin{eqnarray}
y_3(z) &=& \int_{0}^{1} dt_3\;t_3^{5+\lambda } \int_{0}^{1} du_3\;u_3^{\frac{9}{2} +\lambda } \frac{1}{2\pi i} \oint dv_3 \;\frac{1}{v_3} 
\left( \frac{v_3-1}{v_3} \frac{1}{1-\rho (1-t_3)(1-u_3)v_3}\right)^{q_3}  \nonumber\\
&\times& (1-\rho (1-t_3)(1-u_3)v_3)^{-\left( 12 +2 \lambda \right)} \overleftrightarrow {w}_{3,3}^{-(\frac{9}{2}+\frac{\alpha }{2}+\lambda) } \left(\overleftrightarrow {w}_{3,3} \partial_{\overleftrightarrow {w}_{3,3}} \right) \overleftrightarrow {w}_{3,3}^{ \alpha +\frac{1}{2}}\left(\overleftrightarrow {w}_{3,3} \partial_{\overleftrightarrow {w}_{3,3}} \right) \overleftrightarrow {w}_{3,3}^{ 4-\frac{\alpha}{2} +\lambda } \nonumber\\
&\times& \int_{0}^{1} dt_2\;t_2^{3+\lambda } \int_{0}^{1} du_2\;u_2^{\frac{5}{2} +\lambda } \frac{1}{2\pi i} \oint dv_2 \;\frac{1}{v_2} 
\left( \frac{v_2-1}{v_2} \frac{1}{1- \overleftrightarrow {w}_{3,3} (1-t_2)(1-u_2)v_2}\right)^{q_2}  \nonumber\\
&\times& (1- \overleftrightarrow {w}_{3,3} (1-t_2)(1-u_2)v_2)^{-\left( 8 +2 \lambda \right)} \overleftrightarrow {w}_{2,3}^{-(\frac{5}{2}+\frac{\alpha }{2}+\lambda) } \left(\overleftrightarrow {w}_{2,3} \partial_{\overleftrightarrow {w}_{2,3}} \right) \overleftrightarrow {w}_{2,3}^{ \alpha +\frac{1}{2}}\left(\overleftrightarrow {w}_{2,3} \partial_{\overleftrightarrow {w}_{2,3}} \right) \overleftrightarrow {w}_{2,3}^{ 2-\frac{\alpha}{2} +\lambda } \nonumber\\
&\times& \int_{0}^{1} dt_1\;t_1^{1+\lambda } \int_{0}^{1} du_1\;u_1^{\frac{1}{2} +\lambda } \frac{1}{2\pi i} \oint dv_1 \;\frac{1}{v_1} 
\left( \frac{v_1-1}{v_1} \frac{1}{1- \overleftrightarrow {w}_{2,3} (1-t_1)(1-u_1)v_1}\right)^{q_1}  \nonumber\\
&\times& (1- \overleftrightarrow {w}_{2,3} (1-t_1)(1-u_1)v_1)^{-\left( 4 +2 \lambda \right)} \overleftrightarrow {w}_{1,3}^{-(\frac{1}{2}+\frac{\alpha }{2}+\lambda) } \left(\overleftrightarrow {w}_{1,3} \partial_{\overleftrightarrow {w}_{1,3}} \right) \overleftrightarrow {w}_{1,3}^{ \alpha +\frac{1}{2}}\left(\overleftrightarrow {w}_{1,3} \partial_{\overleftrightarrow {w}_{1,3}} \right) \overleftrightarrow {w}_{1,3}^{ -\frac{\alpha}{2} +\lambda } \nonumber\\
&\times& \left\{ c_0 z^{\lambda } \sum_{i_0=0}^{q_0} \frac{(-q_0)_{i_0} \left( q_0 +2 \lambda \right)_{i_0}}{(1+\lambda )_{i_0}\left(\frac{1}{2} +\lambda \right)_{i_0}} \overleftrightarrow {w}_{1,3} ^{i_0}\right\} \eta ^3  \label{eq:70029}
\end{eqnarray}
where
\begin{equation}
\begin{cases} \overleftrightarrow {w}_{3,3} = \frac{v_3}{v_3-1} \frac{ \rho  t_3 u_3}{1- \rho (1-t_3)(1-u_3)v_3} \cr
\overleftrightarrow {w}_{2,3} = \frac{v_2}{v_2-1} \frac{\overleftrightarrow {w}_{3,3} t_2 u_2}{1- \overleftrightarrow {w}_{3,3}(1-t_2)(1-u_2)v_2} \cr
\overleftrightarrow {w}_{1,3} = \frac{v_1}{v_1-1} \frac{\overleftrightarrow {w}_{2,3} t_1 u_1}{1- \overleftrightarrow {w}_{2,3}(1-t_1)(1-u_1)v_1}
\end{cases}
\nonumber
\end{equation}
By repeating this process for all higher terms of integral forms of sub-summation $y_m(x)$ terms where $m \geq 4$, we obtain every integral forms of $y_m(x)$ terms. 
Since we substitute (\ref{eq:70025a}), (\ref{eq:70026}), (\ref{eq:70028}), (\ref{eq:70029}) and including all integral forms of $y_m(x)$ terms where $m \geq 4$ into (\ref{eq:70024}), we obtain (\ref{eq:70023}).
\end{proof}
Put $c_0$= 1 as $\lambda =0$ for the first kind of independent solutions of Lame equation and $\lambda =\frac{1}{2}$ for the second one in (\ref{eq:70023}).
\begin{remark}
The integral representation of Lame equation in the algebraic form of the first kind for polynomial of type 2 about $x=a$ as $q= \alpha (\alpha +1)a- 4(2a-b-c)(q_j+2j )^2 $ where $j,q_j \in \mathbb{N}_{0}$ is
\begin{eqnarray}
y(z)&=& LF_{q_j}^R\left( a, b, c, \alpha, q= \alpha (\alpha +1)a- 4(2a-b-c)(q_j+2j )^2; z= x-a, \rho = -\frac{2a-b-c}{(a-b)(a-c)} z \right. \nonumber\\
&&, \left. \eta = \frac{-z^2}{(a-b)(a-c)} \right) \nonumber\\
&=&\; _2F_1 \left( -q_0, q_0; \frac{1}{2}; \rho \right) + \sum_{n=1}^{\infty } \left\{\prod _{k=0}^{n-1} \Bigg\{ \int_{0}^{1} dt_{n-k}\;t_{n-k}^{2(n-k)-1 } \int_{0}^{1} du_{n-k}\;u_{n-k}^{2(n-k)- \frac{3}{2} }\right. \nonumber\\
&&\times  \frac{1}{2\pi i}  \oint dv_{n-k} \frac{1}{v_{n-k}} \left( \frac{v_{n-k}-1}{v_{n-k}} \frac{1}{1-\overleftrightarrow {w}_{n-k+1,n}(1-t_{n-k})(1-u_{n-k})v_{n-k}}\right)^{q_{n-k}} \nonumber\\
&&\times \left( 1- \overleftrightarrow {w}_{n-k+1,n}(1-t_{n-k})(1-u_{n-k})v_{n-k}\right)^{-4(n-k)}  \nonumber\\
&&\times \overleftrightarrow {w}_{n-k,n}^{-(2(n-k)-\frac{3}{2}+\frac{\alpha }{2} )}\left(  \overleftrightarrow {w}_{n-k,n} \partial _{ \overleftrightarrow {w}_{n-k,n}}\right) \overleftrightarrow {w}_{n-k,n}^{\alpha +\frac{1}{2}}\left(  \overleftrightarrow {w}_{n-k,n} \partial _{ \overleftrightarrow {w}_{n-k,n}}\right) \overleftrightarrow {w}_{n-k,n}^{2(n-k)-2-\frac{\alpha}{2}} \Bigg\} \nonumber\\
&&\times \left. \; _2F_1 \left( -q_0, q_0; \frac{1}{2}; \overleftrightarrow {w}_{1,n} \right) \right\} \eta ^n  \label{eq:70030}
\end{eqnarray}
\end{remark} 
\begin{remark}
The integral representation of Lame equation in the algebraic form of the second kind for polynomial of type 2 about $x=a$ as $q= \alpha (\alpha +1)a- 4(2a-b-c)\left( q_j+2j+\frac{1}{2} \right)^2$ where $j,q_j \in \mathbb{N}_{0}$ is
\begin{eqnarray}
y(z)&=& LS_{q_j}^R\left( a, b, c, \alpha, q= \alpha (\alpha +1)a- 4(2a-b-c)\left( q_j+2j+\frac{1}{2} \right)^2; z= x-a, \rho = -\frac{2a-b-c}{(a-b)(a-c)} z \right. \nonumber\\
&&, \left. \eta = \frac{-z^2}{(a-b)(a-c)} \right) \nonumber\\
&=& z^{\frac{1}{2} } \left\{ \; _2F_1 \left( -q_0, q_0 +1; \frac{3}{2}; \rho \right) \right. + \sum_{n=1}^{\infty } \left\{\prod _{k=0}^{n-1} \Bigg\{ \int_{0}^{1} dt_{n-k}\;t_{n-k}^{2(n-k)- \frac{1}{2} } \int_{0}^{1} du_{n-k}\;u_{n-k}^{2(n-k)-1 }\right. \nonumber\\
&&\times  \frac{1}{2\pi i}  \oint dv_{n-k} \frac{1}{v_{n-k}} \left( \frac{v_{n-k}-1}{v_{n-k}} \frac{1}{1-\overleftrightarrow {w}_{n-k+1,n}(1-t_{n-k})(1-u_{n-k})v_{n-k}}\right)^{q_{n-k}} \nonumber\\
&&\times \left( 1- \overleftrightarrow {w}_{n-k+1,n}(1-t_{n-k})(1-u_{n-k})v_{n-k}\right)^{-\left( 4(n-k) +1 \right)}  \nonumber\\
&&\times \overleftrightarrow {w}_{n-k,n}^{-(2(n-k)-1 +\frac{\alpha }{2} )}\left(  \overleftrightarrow {w}_{n-k,n} \partial _{ \overleftrightarrow {w}_{n-k,n}}\right) \overleftrightarrow {w}_{n-k,n}^{\alpha +\frac{1}{2}}\left(  \overleftrightarrow {w}_{n-k,n} \partial _{ \overleftrightarrow {w}_{n-k,n}}\right) \overleftrightarrow {w}_{n-k,n}^{2(n-k)- \frac{3}{2}-\frac{\alpha}{2}  } \Bigg\} \nonumber\\
&&\times \left.\left.\; _2F_1 \left( -q_0, q_0 +1; \frac{3}{2}; \overleftrightarrow {w}_{1,n} \right) \right\} \eta ^n \right\}  \label{eq:70031}
\end{eqnarray}
\end{remark}
\subsection{Infinite series}
Let's consider the integral representation of Lame equation about $x=0$ for infinite series by applying R3TRF.
There is a generalized hypergeometric function which is written by
\begin{eqnarray}
M_l &=& \sum_{i_l= i_{l-1}}^{\infty } \frac{(-\varphi +2l+\lambda )_{i_l}\left( \varphi +2l+\lambda \right)_{i_l}\left(2l+1+\lambda \right)_{i_{l-1}}\left( 2l+\frac{1}{2} +\lambda \right)_{i_{l-1}}}{(-\varphi +2l+\lambda)_{i_{l-1}}\left( \varphi +2l+\lambda\right)_{i_{l-1}}\left(2l+1+\lambda \right)_{i_l}\left( 2l+\frac{1}{2} +\lambda \right)_{i_l}} \rho ^{i_l} \nonumber\\
&=& \rho ^{i_{l-1}} 
\sum_{j=0}^{\infty } \frac{B(i_{l-1}+2l+\lambda ,j+1) B\left( i_{l-1}+2l- \frac{1}{2} +\lambda ,j+1\right)}{(i_{l-1}+2l+\lambda )^{-1}\left( i_{l-1}+2l-\frac{1}{2} +\lambda \right)^{-1}} \nonumber\\
&&\times \frac{( -\varphi +2l+\lambda+ i_{l-1} )_j \left( \varphi +2l+\lambda +i_{l-1} \right)_j}{(1)_j \;j!}\rho ^j \label{er:70018}
\end{eqnarray}
Substitute (\ref{eq:70019a}) and (\ref{eq:70019b}) into (\ref{er:70018}). And divide $(i_{l-1}+2l+\lambda )\left( i_{l-1}+2l-\frac{1}{2} +\lambda \right)$ into the new (\ref{er:70018}).
\begin{eqnarray}
V_l&=& \frac{1}{(i_{l-1}+2l+\lambda )\left( i_{l-1}+2l-\frac{1}{2} +\lambda \right)} \nonumber\\
&&\times \sum_{i_l= i_{l-1}}^{\infty } \frac{(-\varphi +2l+\lambda )_{i_l}\left( \varphi +2l+\lambda \right)_{i_l}\left(2l+1+\lambda \right)_{i_{l-1}}\left( 2l+\frac{1}{2} +\lambda \right)_{i_{l-1}}}{(-\varphi +2l+\lambda)_{i_{l-1}}\left( \varphi +2l+\lambda\right)_{i_{l-1}}\left(2l+1+\lambda \right)_{i_l}\left( 2l+\frac{1}{2} +\lambda \right)_{i_l}} \rho ^{i_l} \nonumber\\
&=&  \int_{0}^{1} dt_l\;t_l^{2l-1+\lambda } \int_{0}^{1} du_l\;u_l^{2l-\frac{3}{2}+\lambda } (\rho t_l u_l)^{i_{l-1}}  \nonumber\\
&&\times \sum_{j=0}^{\infty } \frac{(-\varphi +2l+\lambda+ i_{l-1} )_j \left( \varphi +2l+\lambda +i_{l-1} \right)_j}{(1)_j \;j!} 
(\rho (1-t_l)(1-u_l))^j \nonumber 
\end{eqnarray}
The hypergeometric function is defined by
\begin{eqnarray}
_2F_1 \left( \alpha ,\beta ; \gamma ; z \right) &=& \sum_{n=0}^{\infty } \frac{(\alpha )_n (\beta )_n}{(\gamma )_n (n!)} z^n \nonumber\\
&=&  \frac{1}{2\pi i} \frac{\Gamma( 1+\alpha  -\gamma )}{\Gamma (\alpha )} \int_0^{(1+)} dv_l\; (-1)^{\gamma }(-v_l)^{\alpha -1} (1-v_l )^{\gamma -\alpha -1} (1-zv_l)^{-\beta }\hspace{1.5cm}\label{er:70019}\\
&& \mbox{where} \;\gamma -\alpha  \ne 1,2,3,\cdots, \;\mbox{Re}(\alpha )>0 \nonumber
\end{eqnarray}
Replace $\alpha $, $\beta $, $\gamma $ and $z$ by $-\varphi +2l+\lambda+ i_{l-1}$, $ \varphi +2l+\lambda +i_{l-1}$, 1 and $\rho (1-t_l)(1-u_l)$ in (\ref{er:70019}). Take the new (\ref{er:70019}) into $V_l$.
\begin{eqnarray}
V_l&=& \frac{1}{(i_{l-1}+2l+\lambda )\left( i_{l-1}+2l-\frac{1}{2} +\lambda \right)} \nonumber\\
&&\times  \sum_{i_l= i_{l-1}}^{\infty } \frac{(-\varphi +2l+\lambda )_{i_l}\left( \varphi +2l+\lambda \right)_{i_l}\left(2l+1+\lambda \right)_{i_{l-1}}\left( 2l+\frac{1}{2} +\lambda \right)_{i_{l-1}}}{(-\varphi +2l+\lambda)_{i_{l-1}}\left( \varphi +2l+\lambda\right)_{i_{l-1}}\left(2l+1+\lambda \right)_{i_l}\left( 2l+\frac{1}{2} +\lambda \right)_{i_l}} \rho ^{i_l}\nonumber\\
&=&  \int_{0}^{1} dt_l\;t_l^{2l-1+\lambda } \int_{0}^{1} du_l\;u_l^{2l-\frac{3}{2} +\lambda } 
\frac{1}{2\pi i} \oint dv_l\;\frac{1}{v_l} \left(\frac{v_l-1}{v_l} \right)^{ \varphi -(2l+\lambda ) } \nonumber\\
&&\times (1-\rho (1-t_l)(1-u_l)v_l)^{-\varphi -(2l+\lambda )} \left(\frac{v_l}{v_l-1} \frac{\rho t_l u_l}{1-\rho (1-t_l)(1-u_l)v_l}\right)^{i_{l-1}}\label{er:70020} 
\end{eqnarray}
Substitute (\ref{er:70020}) into (\ref{eq:70014}) where $l=1,2,3,\cdots$; apply $V_1$ into the second summation of sub-power series $y_1(z)$, apply $V_2$ into the third summation and $V_1$ into the second summation of sub-power series $y_2(z)$, apply $V_3$ into the forth summation, $V_2$ into the third summation and $V_1$ into the second summation of sub-power series $y_3(z)$, etc.\footnote{$y_1(z)$ means the sub-power series in (\ref{eq:70014}) contains one term of $B_n's$, $y_2(z)$ means the sub-power series in (\ref{eq:70014}) contains two terms of $B_n's$, $y_3(z)$ means the sub-power series in (\ref{eq:70014}) contains three terms of $B_n's$, etc.}
\begin{theorem}
The general expression of an integral form of Lame equation in the algebraic form for infinite series about $x=0$ using R3TRF is given by
\begin{eqnarray}
 y(z)&=& \sum_{n=0}^{\infty } y_{n}(z)= y_0(z)+ y_1(z)+ y_2(z)+y_3(z)+\cdots \nonumber\\
&=& c_0 z^{\lambda } \left\{ \sum_{i_0=0}^{\infty }\frac{\left( -\varphi +\lambda\right)_{i_0}\left(\varphi +\lambda \right)_{i_0}}{(1+\lambda )_{i_0}\left(\frac{1}{2} +\lambda \right)_{i_0}} \rho ^{i_0}\right.  \nonumber\\
&+& \sum_{n=1}^{\infty } \left\{\prod _{k=0}^{n-1} \Bigg\{ \int_{0}^{1} dt_{n-k}\;t_{n-k}^{2(n-k)-1+\lambda } \int_{0}^{1} du_{n-k}\;u_{n-k}^{2(n-k)- \frac{3}{2} +\lambda }\right. \nonumber\\
&\times&  \frac{1}{2\pi i}  \oint dv_{n-k} \frac{1}{v_{n-k}} \left( \frac{v_{n-k}-1}{v_{n-k}}\right)^{\varphi -2(n-k)-\lambda} \nonumber\\
&\times& \left( 1- \overleftrightarrow {w}_{n-k+1,n}(1-t_{n-k})(1-u_{n-k})v_{n-k}\right)^{-\varphi -2(n-k)-\lambda}  \nonumber\\
&\times& \overleftrightarrow {w}_{n-k,n}^{-(2(n-k)- \frac{3}{2}+\frac{\alpha }{2} +\lambda )}\left(  \overleftrightarrow {w}_{n-k,n} \partial _{ \overleftrightarrow {w}_{n-k,n}}\right) \overleftrightarrow {w}_{n-k,n}^{\alpha +\frac{1}{2}}\left(  \overleftrightarrow {w}_{n-k,n} \partial _{ \overleftrightarrow {w}_{n-k,n}}\right) \overleftrightarrow {w}_{n-k,n}^{2(n-k)-2-\frac{\alpha}{2} +\lambda } \Bigg\} \nonumber\\
&\times& \left.\left.\sum_{i_0=0}^{\infty }\frac{\left( -\varphi +\lambda\right)_{i_0}\left(\varphi +\lambda \right)_{i_0}}{(1+\lambda )_{i_0}\left(\frac{1}{2} +\lambda \right)_{i_0}} \overleftrightarrow {w}_{1,n}^{i_0}\right\} \eta ^n \right\}  \label{eq:70032}
\end{eqnarray}
In the above, the first sub-integral form contains one term of $B_n's$, the second one contains two terms of $B_n$'s, the third one contains three terms of $B_n$'s, etc.
\end{theorem}
\begin{proof} 
In (\ref{eq:70014}) sub-power series  $y_0(z) $, $y_1(z)$, $y_2(z)$ and $y_3(z)$ of the Lame equation for infinite series about $x=a$ using R3TRF  are given by
\begin{subequations}
\begin{equation}
 y_0(z)= c_0 z^{\lambda } \sum_{i_0=0}^{\infty } \frac{(-\varphi +\lambda )_{i_0} \left( \varphi  +\lambda \right)_{i_0}}{\left( 1+\lambda \right)_{i_0}\left( \frac{1}{2} +\lambda \right)_{i_0}} \rho ^{i_0} \label{er:70021a}
\end{equation}
\begin{eqnarray}
 y_1(z) &=& c_0 z^{\lambda } \left\{\sum_{i_0=0}^{\infty}\frac{\left( i_0 +\frac{\alpha}{2}+ \frac{1}{2}+\lambda  \right) \left( i_0-\frac{\alpha}{2}+\lambda  \right)}{\left(i_0 +2+ \lambda\right) \left( i_0+ \frac{3}{2}+\lambda  \right)}\frac{(-\varphi +\lambda)_{i_0} \left( \varphi +\lambda \right)_{i_0}}{\left( 1+\lambda \right)_{i_0}\left( \frac{1}{2} +\lambda \right)_{i_0}}\right. \nonumber\\
&&\times \left. \sum_{i_1=i_0}^{\infty} \frac{(-\varphi +2 +\lambda)_{i_1}\left(  \varphi +2 +\lambda \right)_{i_1}\left( 3+\lambda \right)_{i_0}\left( \frac{5}{2} +\lambda \right)_{i_0}}{(-\varphi +2 +\lambda)_{i_0}\left(  \varphi +2 +\lambda \right)_{i_0}\left( 3+\lambda \right)_{i_1}\left( \frac{5}{2} +\lambda \right)_{i_1}} \rho ^{i_1} \right\}\eta  \label{er:70021b}
\end{eqnarray}
\begin{eqnarray}
 y_2(z) &=& c_0 z^{\lambda } \left\{\sum_{i_0=0}^{\infty}\frac{\left( i_0+\frac{\alpha}{2}+ \frac{1}{2} +\lambda  \right) \left( i_0-\frac{\alpha}{2}+\lambda  \right)}{\left(i_0 +2+ \lambda\right) \left( i_0+ \frac{3}{2}+\lambda  \right)}\frac{(-\varphi +\lambda)_{i_0} \left( \varphi +\lambda \right)_{i_0}}{\left( 1+\lambda \right)_{i_0}\left( \frac{1}{2} +\lambda \right)_{i_0}}\right. \nonumber\\
&&\times  \sum_{i_1=i_0}^{\infty} \frac{\left( i_1 +\frac{\alpha}{2}+ \frac{5}{2}+\lambda  \right) \left( i_1-\frac{\alpha}{2}+2+\lambda  \right)}{\left(i_1 +4+ \lambda\right) \left( i_1+ \frac{7}{2}+\lambda  \right)} \frac{(-\varphi +2 +\lambda)_{i_1}\left( \varphi +2 +\lambda \right)_{i_1}\left( 3+\lambda \right)_{i_0}\left( \frac{5}{2} +\lambda \right)_{i_0}}{(-\varphi +2 +\lambda)_{i_0}\left( \varphi +2 +\lambda \right)_{i_0}\left( 3+\lambda \right)_{i_1}\left( \frac{5}{2} +\lambda \right)_{i_1}} \nonumber\\
&&\times \left.\sum_{i_2=i_1}^{\infty} \frac{(-\varphi +4+\lambda)_{i_2}\left( \varphi +4+\lambda \right)_{i_2}\left( 5+\lambda \right)_{i_1}\left(\frac{9}{2}+\lambda \right)_{i_1}}{(-\varphi +4+\lambda)_{i_1}\left(  \varphi +4+\lambda \right)_{i_1}\left( 5+\lambda \right)_{i_2}\left(\frac{9}{2}+\lambda \right)_{i_2}} \rho ^{i_2} \right\} \eta ^2  \label{er:70021c}
\end{eqnarray}
\begin{eqnarray}
 y_3(z) &=&  c_0 z^{\lambda } \left\{\sum_{i_0=0}^{\infty}\frac{\left( i_0+\frac{\alpha}{2}+ \frac{1}{2} +\lambda  \right) \left( i_0-\frac{\alpha}{2}+\lambda  \right)}{\left(i_0 +2+ \lambda\right) \left( i_0+ \frac{3}{2}+\lambda  \right)}\frac{(-\varphi  +\lambda)_{i_0} \left( \varphi +\lambda\right)_{i_0}}{\left( 1+\lambda \right)_{i_0}\left( \frac{1}{2} +\lambda \right)_{i_0}} \right. \nonumber\\
&&\times  \sum_{i_1=i_0}^{\infty} \frac{\left( i_1 +\frac{\alpha}{2}+ \frac{5}{2}+\lambda  \right) \left( i_1-\frac{\alpha}{2}+2+\lambda  \right)}{\left(i_1 +4+ \lambda\right) \left( i_1+ \frac{7}{2}+\lambda  \right)} \frac{(-\varphi +2 +\lambda)_{i_1}\left( \varphi +2 +\lambda\right)_{i_1}\left( 3+\lambda \right)_{i_0}\left( \frac{5}{2} +\lambda \right)_{i_0}}{(-\varphi +2 +\lambda)_{i_0}\left( \varphi +2 +\lambda\right)_{i_0}\left( 3+\lambda \right)_{i_1}\left( \frac{5}{2} +\lambda \right)_{i_1}} \nonumber\\
&&\times \sum_{i_2=i_1}^{\infty} \frac{\left( i_2 +\frac{\alpha}{2}+ \frac{9}{2}+\lambda  \right) \left( i_2-\frac{\alpha}{2}+4+\lambda  \right)}{\left(i_2 +6+ \lambda\right) \left( i_2+ \frac{11}{2}+\lambda  \right)}  \frac{(-\varphi +4+\lambda)_{i_2}\left( \varphi +4+\lambda\right)_{i_2}\left( 5+\lambda \right)_{i_1}\left(\frac{9}{2}+\lambda \right)_{i_1}}{(-\varphi +4+\lambda)_{i_1}\left( \varphi +4+\lambda \right)_{i_1}\left( 5+\lambda \right)_{i_2}\left(\frac{9}{2}+\lambda \right)_{i_2}} \nonumber\\
&&\times \left.\sum_{i_3=i_2}^{\infty}\frac{( -\varphi +6+\lambda)_{i_3}\left(  \varphi +6+\lambda \right)_{i_3}\left( 7+\lambda \right)_{i_2}\left(\frac{13}{2}+\lambda \right)_{i_2}}{(-\varphi +6+\lambda)_{i_2}\left( \varphi +6+\lambda \right)_{i_3}\left( 7+\lambda \right)_{i_3}\left(\frac{13}{2}+\lambda \right)_{i_3}}\rho ^{i_3} \right\} \eta ^3  \label{er:70021d} 
\end{eqnarray}
\end{subequations}
Put $l=1$ in (\ref{er:70020}). Take the new (\ref{er:70020}) into (\ref{er:70021b}).
\begin{eqnarray}
y_1(z) &=& \int_{0}^{1} dt_1\;t_1^{1+\lambda } \int_{0}^{1} du_1\;u_1^{\frac{1}{2} +\lambda } \frac{1}{2\pi i} \oint dv_1 \;\frac{1}{v_1} 
\left( \frac{v_1-1}{v_1}  \right)^{\varphi -(2+\lambda )}  \nonumber\\
&&\times (1-\rho (1-t_1)(1-u_1)v_1)^{-\varphi -(2+\lambda )} \overleftrightarrow {w}_{1,1}^{-(\frac{1}{2}+\frac{\alpha }{2}+\lambda )} \left(\overleftrightarrow {w}_{1,1} \partial_{\overleftrightarrow {w}_{1,1}} \right) \overleftrightarrow {w}_{1,1}^{\alpha +\frac{1}{2} } \left(\overleftrightarrow {w}_{1,1} \partial_{\overleftrightarrow {w}_{1,1}} \right) \overleftrightarrow {w}_{1,1}^{- \frac{\alpha }{2}+\lambda }\nonumber\\
&&\times \left\{ c_0 z^{\lambda } \sum_{i_0=0}^{\infty } \frac{(-\varphi +\lambda )_{i_0} \left( \varphi  +\lambda \right)_{i_0}}{\left( 1+\lambda \right)_{i_0}\left( \frac{1}{2} +\lambda \right)_{i_0}} \overleftrightarrow {w}_{1,1} ^{i_0}\right\}\eta  \label{er:70022}
\end{eqnarray}
where
\begin{equation}
\overleftrightarrow {w}_{1,1} = \frac{v_1}{v_1-1} \frac{\rho  t_1 u_1}{1-\rho (1-t_1)(1-u_1)v_1}\nonumber
\end{equation}
Put $l=2$ in (\ref{er:70020}). Take the new (\ref{er:70020}) into (\ref{er:70021c}).
\begin{eqnarray}
y_2(z) &=& c_0 z^{\lambda } \int_{0}^{1} dt_2\;t_2^{3+\lambda } \int_{0}^{1} du_2\;u_2^{\frac{5}{2} +\lambda } \frac{1}{2\pi i} \oint dv_2 \;\frac{1}{v_2} 
\left( \frac{v_2-1}{v_2} \right)^{\varphi -(4+\lambda )}  \nonumber\\
&&\times (1-\rho (1-t_2)(1-u_2)v_2)^{-\varphi -(4+\lambda )} 
 \overleftrightarrow {w}_{2,2}^{-(\frac{5}{2}+\frac{\alpha }{2}+\lambda) } \left(\overleftrightarrow {w}_{2,2} \partial_{\overleftrightarrow {w}_{2,2}} \right) \overleftrightarrow {w}_{2,2}^{ \alpha +\frac{1}{2}}\left(\overleftrightarrow {w}_{2,2} \partial_{\overleftrightarrow {w}_{2,2}} \right) \overleftrightarrow {w}_{2,2}^{ 2-\frac{\alpha}{2} +\lambda } \nonumber\\
&&\times \left\{\sum_{i_0=0}^{\infty}\frac{\left( i_0 +\frac{\alpha}{2}+ \frac{1}{2}+\lambda  \right) \left( i_0-\frac{\alpha}{2}+\lambda  \right)}{\left(i_0 +2+ \lambda\right) \left( i_0+ \frac{3}{2}+\lambda  \right)}\frac{(-\varphi +\lambda)_{i_0} \left( \varphi +\lambda \right)_{i_0}}{\left( 1+\lambda \right)_{i_0}\left( \frac{1}{2} +\lambda \right)_{i_0}}\right. \nonumber\\
&&\times \left. \sum_{i_1=i_0}^{\infty} \frac{(-\varphi +2 +\lambda)_{i_1}\left(  \varphi +2 +\lambda \right)_{i_1}\left( 3+\lambda \right)_{i_0}\left( \frac{5}{2} +\lambda \right)_{i_0}}{(-\varphi +2 +\lambda)_{i_0}\left(  \varphi +2 +\lambda \right)_{i_0}\left( 3+\lambda \right)_{i_1}\left( \frac{5}{2} +\lambda \right)_{i_1}} \overleftrightarrow {w}_{2,2}^{i_1} \right\} \eta ^2 \label{er:70023}
\end{eqnarray}
where
\begin{equation}
\overleftrightarrow {w}_{2,2} = \frac{v_2}{v_2-1} \frac{\rho  t_2 u_2}{1-\rho (1-t_2)(1-u_2)v_2}\nonumber
\end{equation}
Put $l=1$ and $\rho  = \overleftrightarrow {w}_{2,2}$ in (\ref{er:70020}). Take the new (\ref{er:70020}) into (\ref{er:70023}).
\begin{eqnarray}
y_2(z) &=& c_0 z^{\lambda } \int_{0}^{1} dt_2\;t_2^{3+\lambda } \int_{0}^{1} du_2\;u_2^{\frac{5}{2} +\lambda } \frac{1}{2\pi i} \oint dv_2 \;\frac{1}{v_2} 
\left( \frac{v_2-1}{v_2} \right)^{\varphi -(4+\lambda )}  \nonumber\\
&&\times (1-\rho (1-t_2)(1-u_2)v_2)^{-\varphi -(4+\lambda )} 
 \overleftrightarrow {w}_{2,2}^{-(\frac{5}{2}+\frac{\alpha }{2}+\lambda) } \left(\overleftrightarrow {w}_{2,2} \partial_{\overleftrightarrow {w}_{2,2}} \right) \overleftrightarrow {w}_{2,2}^{ \alpha +\frac{1}{2}}\left(\overleftrightarrow {w}_{2,2} \partial_{\overleftrightarrow {w}_{2,2}} \right) \overleftrightarrow {w}_{2,2}^{ 2-\frac{\alpha}{2} +\lambda } \nonumber\\
&&\times \int_{0}^{1} dt_1\;t_1^{1+\lambda } \int_{0}^{1} du_1\;u_1^{\frac{1}{2} +\lambda } \frac{1}{2\pi i} \oint dv_1 \;\frac{1}{v_1} 
\left( \frac{v_1-1}{v_1} \right)^{\varphi -(2+\lambda )}  \nonumber\\
&&\times (1- \overleftrightarrow {w}_{2,2} (1-t_1)(1-u_1)v_1)^{-\varphi -(2+\lambda )} \overleftrightarrow {w}_{1,2}^{-(\frac{1}{2}+\frac{\alpha }{2}+\lambda) } \left(\overleftrightarrow {w}_{1,2} \partial_{\overleftrightarrow {w}_{1,2}} \right) \overleftrightarrow {w}_{1,2}^{ \alpha +\frac{1}{2}}\left(\overleftrightarrow {w}_{1,2} \partial_{\overleftrightarrow {w}_{1,2}} \right) \overleftrightarrow {w}_{1,2}^{-\frac{\alpha }{2}+\lambda } \nonumber\\
&&\times \left\{ c_0 z^{\lambda } \sum_{i_0=0}^{\infty } \frac{(-\varphi +\lambda )_{i_0} \left( \varphi  +\lambda \right)_{i_0}}{\left( 1+\lambda \right)_{i_0}\left( \frac{1}{2} +\lambda \right)_{i_0}} \overleftrightarrow {w}_{1,2} ^{i_0}\right\} \eta ^2 \label{er:70024}
\end{eqnarray}
where
\begin{equation}
\overleftrightarrow {w}_{1,2} = \frac{v_1}{v_1-1} \frac{\overleftrightarrow {w}_{2,2} t_1 u_1}{1- \overleftrightarrow {w}_{2,2}(1-t_1)(1-u_1)v_1}\nonumber
\end{equation}
By using similar process for the previous cases of integral forms of $y_1(z)$ and $y_2(z)$, the integral form of sub-power series expansion of $y_3(z)$ is given by
\begin{eqnarray}
y_3(z) &=& \int_{0}^{1} dt_3\;t_3^{5+\lambda } \int_{0}^{1} du_3\;u_3^{\frac{9}{2} +\lambda } \frac{1}{2\pi i} \oint dv_3 \;\frac{1}{v_3} 
\left( \frac{v_3-1}{v_3} \right)^{\varphi -(6+\lambda )}  \nonumber\\
&&\times (1-\rho (1-t_3)(1-u_3)v_3)^{-\varphi -(6+\lambda )}\nonumber\\
&&\times  \overleftrightarrow {w}_{3,3}^{-(\frac{9}{2}+\frac{\alpha }{2}+\lambda) } \left(\overleftrightarrow {w}_{3,3} \partial_{\overleftrightarrow {w}_{3,3}} \right) \overleftrightarrow {w}_{3,3}^{ \alpha +\frac{1}{2}}\left(\overleftrightarrow {w}_{3,3} \partial_{\overleftrightarrow {w}_{3,3}} \right) \overleftrightarrow {w}_{3,3}^{ 4-\frac{\alpha}{2} +\lambda } \nonumber\\
&&\times \int_{0}^{1} dt_2\;t_2^{3+\lambda } \int_{0}^{1} du_2\;u_2^{\frac{5}{2} +\lambda } \frac{1}{2\pi i} \oint dv_2 \;\frac{1}{v_2} 
\left( \frac{v_2-1}{v_2} \right)^{\varphi -(4+\lambda )}  \nonumber\\
&&\times (1- \overleftrightarrow {w}_{3,3} (1-t_2)(1-u_2)v_2)^{-\varphi -(4+\lambda )} \nonumber\\
&&\times  \overleftrightarrow {w}_{2,3}^{-(\frac{5}{2}+\frac{\alpha }{2}+\lambda) } \left(\overleftrightarrow {w}_{2,3} \partial_{\overleftrightarrow {w}_{2,3}} \right) \overleftrightarrow {w}_{2,3}^{ \alpha +\frac{1}{2}}\left(\overleftrightarrow {w}_{2,3} \partial_{\overleftrightarrow {w}_{2,3}} \right) \overleftrightarrow {w}_{2,3}^{ 2-\frac{\alpha}{2} +\lambda } \nonumber\\
&&\times \int_{0}^{1} dt_1\;t_1^{1+\lambda } \int_{0}^{1} du_1\;u_1^{\frac{1}{2} +\lambda } \frac{1}{2\pi i} \oint dv_1 \;\frac{1}{v_1} 
\left( \frac{v_1-1}{v_1}  \right)^{\varphi -(2+\lambda )}  \nonumber\\
&&\times (1- \overleftrightarrow {w}_{2,3} (1-t_1)(1-u_1)v_1)^{-\varphi -(2+\lambda )} \nonumber\\
&&\times \overleftrightarrow {w}_{1,3}^{-(\frac{1}{2}+\frac{\alpha }{2}+\lambda) } \left(\overleftrightarrow {w}_{1,3} \partial_{\overleftrightarrow {w}_{1,3}} \right) \overleftrightarrow {w}_{1,3}^{ \alpha +\frac{1}{2}}\left(\overleftrightarrow {w}_{1,3} \partial_{\overleftrightarrow {w}_{1,3}} \right) \overleftrightarrow {w}_{1,3}^{ -\frac{\alpha}{2} +\lambda } \nonumber\\
&&\times \left\{ c_0 z^{\lambda } \sum_{i_0=0}^{\infty } \frac{(-\varphi +\lambda )_{i_0} \left( \varphi  +\lambda \right)_{i_0}}{\left( 1+\lambda \right)_{i_0}\left( \frac{1}{2} +\lambda \right)_{i_0}} \overleftrightarrow {w}_{1,3} ^{i_0}\right\} \eta ^3  \label{er:70025}
\end{eqnarray}
where
\begin{equation}
\begin{cases} \overleftrightarrow {w}_{3,3} = \frac{v_3}{v_3-1} \frac{ \rho  t_3 u_3}{1- \rho (1-t_3)(1-u_3)v_3} \cr
\overleftrightarrow {w}_{2,3} = \frac{v_2}{v_2-1} \frac{\overleftrightarrow {w}_{3,3} t_2 u_2}{1- \overleftrightarrow {w}_{3,3}(1-t_2)(1-u_2)v_2} \cr
\overleftrightarrow {w}_{1,3} = \frac{v_1}{v_1-1} \frac{\overleftrightarrow {w}_{2,3} t_1 u_1}{1- \overleftrightarrow {w}_{2,3}(1-t_1)(1-u_1)v_1}
\end{cases}
\nonumber
\end{equation}
By repeating this process for all higher terms of integral forms of sub-summation $y_m(x)$ terms where $m \geq 4$, we obtain every integral forms of $y_m(x)$ terms. 
Since we substitute (\ref{er:70021a}), (\ref{er:70022}), (\ref{er:70024}), (\ref{er:70025}) and including all integral forms of $y_m(x)$ terms where $m \geq 4$ into (\ref{eq:70014}), we obtain (\ref{eq:70032}).\footnote{Or replace the finite summation with an interval $[0, q_0]$ by infinite summation with an interval  $[0,\infty ]$ in (\ref{eq:70023}). Replace $q_0$ and $q_{n-k}$ by $\varphi -\lambda $ and $\varphi -2(n-k)-\lambda$ into the new (\ref{eq:70023}). Its solution is also equivalent to (\ref{eq:70032})}
\end{proof}
Put $c_0$= 1 as $\lambda =0$  for the first independent solution of Lame equation and $\lambda =\frac{1}{2}$ for the second one in (\ref{eq:70032}). 
\begin{remark}
The integral representation of Lame equation in the algebraic form of the first kind for infinite series about $x=a$ using R3TRF is
\begin{eqnarray}
 y(z)&=& LF^R\left( a, b, c, \alpha, q, \varphi = \sqrt{\frac{\alpha (\alpha +1)a-q}{4(2a-b-c)}}; z= x-a, \rho = -\frac{2a-b-c}{(a-b)(a-c)} z, \eta = \frac{-z^2}{(a-b)(a-c)} \right) \nonumber\\
&=& \; _2F_1 \left(  -\varphi, \varphi; \frac{1}{2}; \rho \right) + \sum_{n=1}^{\infty } \left\{\prod _{k=0}^{n-1} \Bigg\{ \int_{0}^{1} dt_{n-k}\;t_{n-k}^{2(n-k)-1  } \int_{0}^{1} du_{n-k}\;u_{n-k}^{2(n-k)- \frac{3}{2} }\right. \nonumber\\
&\times&  \frac{1}{2\pi i}  \oint dv_{n-k} \frac{1}{v_{n-k}} \left( \frac{v_{n-k}-1}{v_{n-k}}\right)^{\varphi -2(n-k) }  \left( 1- \overleftrightarrow {w}_{n-k+1,n}(1-t_{n-k})(1-u_{n-k})v_{n-k}\right)^{-\varphi -2(n-k) }  \nonumber\\
&\times& \overleftrightarrow {w}_{n-k,n}^{-(2(n-k)- \frac{3}{2}+\frac{\alpha }{2}  )}\left(  \overleftrightarrow {w}_{n-k,n} \partial _{ \overleftrightarrow {w}_{n-k,n}}\right) \overleftrightarrow {w}_{n-k,n}^{\alpha +\frac{1}{2}}\left(  \overleftrightarrow {w}_{n-k,n} \partial _{ \overleftrightarrow {w}_{n-k,n}}\right) \overleftrightarrow {w}_{n-k,n}^{2(n-k)-2-\frac{\alpha}{2} } \Bigg\} \nonumber\\
&\times& \left.\; _2F_1 \left(  -\varphi, \varphi; \frac{1}{2}; \overleftrightarrow {w}_{1,n} \right) \right\} \eta ^n   \label{eq:70033}
\end{eqnarray}
\end{remark}
\begin{remark}
The integral representation of Lame equation in the algebraic form of the second kind for infinite series about $x=a$ using R3TRF is
\begin{eqnarray}
 y(z)&=& LS^R\left( a, b, c, \alpha, q, \varphi = \sqrt{\frac{\alpha (\alpha +1)a-q}{4(2a-b-c)}}; z= x-a, \rho = -\frac{2a-b-c}{(a-b)(a-c)} z, \eta = \frac{-z^2}{(a-b)(a-c)} \right) \nonumber\\
&=& z^{\frac{1}{2} } \left\{ \; _2F_1 \left(  -\varphi +\frac{1}{2}, \varphi +\frac{1}{2}; \frac{3}{2}; \rho \right)  \right. + \sum_{n=1}^{\infty } \left\{\prod _{k=0}^{n-1} \Bigg\{ \int_{0}^{1} dt_{n-k}\;t_{n-k}^{2(n-k)-\frac{1}{2}} \int_{0}^{1} du_{n-k}\;u_{n-k}^{2(n-k)-1 }\right. \nonumber\\
&\times&  \frac{1}{2\pi i}  \oint dv_{n-k} \frac{1}{v_{n-k}} \left( \frac{v_{n-k}-1}{v_{n-k}}\right)^{\varphi -2(n-k)-\frac{1}{2}}  \left( 1- \overleftrightarrow {w}_{n-k+1,n}(1-t_{n-k})(1-u_{n-k})v_{n-k}\right)^{-\varphi -2(n-k)-\frac{1}{2}}  \nonumber\\
&\times& \overleftrightarrow {w}_{n-k,n}^{-(2(n-k)-1 +\frac{\alpha }{2} )}\left(  \overleftrightarrow {w}_{n-k,n} \partial _{ \overleftrightarrow {w}_{n-k,n}}\right) \overleftrightarrow {w}_{n-k,n}^{\alpha +\frac{1}{2}}\left(  \overleftrightarrow {w}_{n-k,n} \partial _{ \overleftrightarrow {w}_{n-k,n}}\right) \overleftrightarrow {w}_{n-k,n}^{2(n-k)-\frac{3}{2}+ \frac{\alpha}{2}} \Bigg\} \nonumber\\
&\times& \left.\left.\; _2F_1 \left(  -\varphi +\frac{1}{2}, \varphi +\frac{1}{2}; \frac{3}{2}; \overleftrightarrow {w}_{1,n} \right) \right\} \eta ^n \right\}  \label{eq:70034}
\end{eqnarray}
\end{remark}
\section[Generating function for the Lame polynomial of type 2]{Generating function for the Lame polynomial of type 2
\sectionmark{Generating function for the Lame polynomial of type 2}}
\sectionmark{Generating function for the Lame polynomial of type 2}
I consider the generating function for the Lame polynomial of type 2 in the algebriac form by applying R3TRF. Since the generating function for the Lame polynomial is derived, we might be possible to construct orthogonal relations of it.\footnote{In the next series, I will construct the generating function for the type 3 Lame polynomial. For the type 3 polynomial, I treat $\alpha $ and $q$ as fixed values.} 
\begin{lemma}
The generating function for the Jacobi polynomial using hypergeometric functions is given by
\begin{eqnarray}
&&\sum_{q_0=0}^{\infty }\frac{(\gamma )_{q_0}}{ q_0!} w^{q_0} \;_2F_1(-q_0, q_0+A; \gamma; x) \label{eq:70035}\\
&&= 2^{A -1}\frac{\left(1-w+\sqrt{w^2-2(1-2x)w+1}\right)^{1-\gamma } \left(1+w+\sqrt{w^2-2(1-2x)w+1}\right)^{\gamma -A}}{\sqrt{w^2-2(1-2x)w+1}} \nonumber\\
&& \hspace{.5cm} \mbox{where}\;|w|<1 \nonumber
\end{eqnarray}
\end{lemma}
\begin{proof}
The proof of this lemma is given in Lemma 3.2.1
\end{proof}

\begin{definition}
I define that
\begin{equation}
\begin{cases}
\displaystyle { s_{a,b}} = \begin{cases} \displaystyle {  s_a\cdot s_{a+1}\cdot s_{a+2}\cdots s_{b-2}\cdot s_{b-1}\cdot s_b}\;\;\mbox{if}\;a>b \cr
s_a \;\;\mbox{if}\;a=b\end{cases}
\cr
\cr
\displaystyle { \widetilde{w}_{i,j}}  = 
\begin{cases} \displaystyle { \frac{ \widetilde{w}_{i+1,j}\; t_i u_i \left\{ 1+ (s_i+2\widetilde{w}_{i+1,j}(1-t_i)(1-u_i))s_i\right\}}{2(1-\widetilde{w}_{i+1,j}(1-t_i)(1-u_i))^2 s_i}} \cr
\displaystyle {-\frac{\widetilde{w}_{i+1,j}\; t_i u_i (1+s_i)\sqrt{s_i^2-2(1-2\widetilde{w}_{i+1,j}(1-t_i)(1-u_i))s_i+1}}{2(1-\widetilde{w}_{i+1,j}(1-t_i)(1-u_i))^2 s_i}} \;\;\mbox{where}\;i<j \cr
\cr
\displaystyle { \frac{\rho  t_i u_i \left\{ 1+ (s_{i,\infty }+2\rho (1-t_i)(1-u_i))s_{i,\infty }\right\}}{2(1-\rho (1-t_i)(1-u_i))^2 s_{i,\infty }}} \cr
\displaystyle {-\frac{\rho t_i u_i(1+s_{i,\infty })\sqrt{s_{i,\infty }^2-2(1-2\rho (1-t_i)(1-u_i))s_{i,\infty }+1}}{2(1-\rho (1-t_i)(1-u_i))^2 s_{i,\infty }}} \;\;\mbox{where}\;i=j 
\end{cases}
\end{cases}\label{eq:70040}
\end{equation}
where
\begin{equation}
a,b,i,j\in \mathbb{N}_{0} \nonumber
\end{equation}
\end{definition}
And we have
\begin{equation}
\sum_{q_i = q_j}^{\infty } r_i^{q_i} = \frac{r_i^{q_j}}{(1-r_i)}\label{eq:70041}
\end{equation}
Acting the summation operator $\displaystyle{ \sum_{q_0 =0}^{\infty } \frac{(\gamma')_{q_0}}{q_0!} s_0^{q_0} \prod _{n=1}^{\infty } \left\{ \sum_{ q_n = q_{n-1}}^{\infty } s_n^{q_n }\right\}}$ on (\ref{eq:70023}) where $|s_i|<1$ as $i=0,1,2,\cdots$ by using (\ref{eq:70040}) and (\ref{eq:70041}),
\begin{theorem} 
The general expression of the generating function for the Lame polynomial of type 2 in the algebraic form is given by
\begin{eqnarray}
&&\sum_{q_0 =0}^{\infty } \frac{(\gamma')_{q_0}}{q_0!} s_0^{q_0} \prod _{n=1}^{\infty } \left\{ \sum_{ q_n = q_{n-1}}^{\infty } s_n^{q_n }\right\} y(z) \nonumber\\
&&= \prod_{l=1}^{\infty } \frac{1}{(1-s_{l,\infty })} \mathbf{\Upsilon}(\lambda; s_{0,\infty } ;\rho ) \nonumber\\
&&+ \Bigg\{ \prod_{l=2}^{\infty } \frac{1}{(1-s_{l,\infty })} \int_{0}^{1} dt_1\;t_1^{1+\lambda} \int_{0}^{1} du_1\;u_1^{\frac{1}{2} +\lambda} \left(s_{1,\infty }^2-2(1-2 \rho (1-t_1)(1-u_1))s_{1,\infty }+1\right)^{-\frac{1}{2}}\nonumber\\
&&\times \left( \frac{ 1+s_{1,\infty } +\sqrt{s_{1,\infty }^2-2(1-2\rho (1-t_1)(1-u_1))s_{1,\infty }+1}}{2}\right)^{-\left( 3+2\lambda \right)} \nonumber\\
&&\times \widetilde{w}_{1,1}^{-(\frac{1}{2}+\frac{\alpha }{2}+\lambda )}\left( \widetilde{w}_{1,1} \partial _{ \widetilde{w}_{1,1}}\right) \widetilde{w}_{1,1}^{\alpha +\frac{1}{2}} \left( \widetilde{w}_{1,1} \partial _{ \widetilde{w}_{1,1}}\right) \widetilde{w}_{1,1}^{-\frac{\alpha }{2}+\lambda } \;\mathbf{\Upsilon}(\lambda ; s_0;\widetilde{w}_{1,1}) \Bigg\} \eta \nonumber\\
&&+ \sum_{n=2}^{\infty } \Bigg\{ \prod_{l=n+1}^{\infty } \frac{1}{(1-s_{l,\infty })} \int_{0}^{1} dt_n\;t_n^{2n-1+\lambda } \int_{0}^{1} du_n\;u_n^{2n-\frac{3}{2} +\lambda}  \nonumber\\
&&\times   \left( s_{n,\infty }^2-2(1-2\rho (1-t_n)(1-u_n))s_{n,\infty }+1\right)^{-\frac{1}{2}} \nonumber\\
&&\times  \left( \frac{ 1+s_{n,\infty } +\sqrt{s_{n,\infty }^2-2(1-2\rho (1-t_n)(1-u_n))s_{n,\infty }+1}}{2}\right)^{-\left( 4n-1 + 2\lambda \right)} \nonumber\\
&&\times \widetilde{w}_{n,n}^{-(2n-\frac{3}{2}+\frac{\alpha }{2}+\lambda )}\left(  \widetilde{w}_{n,n} \partial _{ \widetilde{w}_{n,n}}\right)  \widetilde{w}_{n,n}^{\alpha +\frac{1}{2}} \left(  \widetilde{w}_{n,n} \partial _{ \widetilde{w}_{n,n}}\right)\widetilde{w}_{n,n}^{ 2(n-1)-\frac{\alpha }{2} +\lambda }\nonumber\\
&&\times \prod_{k=1}^{n-1} \Bigg\{ \int_{0}^{1} dt_{n-k}\;t_{n-k}^{2(n-k)-1+\lambda } \int_{0}^{1} du_{n-k} \;u_{n-k}^{2(n-k)-\frac{3}{2} +\lambda }\label{eq:70042}\\
&&\times \left( s_{n-k}^2-2(1-2\widetilde{w}_{n-k+1,n} (1-t_{n-k})(1-u_{n-k}))s_{n-k}+1 \right)^{-\frac{1}{2}} \nonumber\\
&&\times \left( \frac{ 1+s_{n-k} +\sqrt{s_{n-k}^2-2(1-2\widetilde{w}_{n-k+1,n} (1-t_{n-k})(1-u_{n-k}))s_{n-k}+1}}{2}\right)^{-\left( 4(n-k)-1 +2\lambda \right)} \nonumber\\
&&\times  \widetilde{w}_{n-k,n}^{-(2(n-k)-\frac{3}{2}+\frac{\alpha }{2} +\lambda )}\left(  \widetilde{w}_{n-k,n} \partial _{ \widetilde{w}_{n-k,n}}\right) \widetilde{w}_{n-k,n}^{\alpha +\frac{1}{2}}\left( \widetilde{w}_{n-k,n} \partial _{ \widetilde{w}_{n-k,n}}\right) \widetilde{w}_{n-k,n}^{ 2(n-k-1)-\frac{\alpha }{2} +\lambda }\Bigg\} \mathbf{\Upsilon}(\lambda ; s_0;\widetilde{w}_{1,n}) \Bigg\} \eta ^n \nonumber
\end{eqnarray}
where
\begin{equation}
\begin{cases} 
{ \displaystyle \mathbf{\Upsilon}(\lambda; s_{0,\infty } ;\rho )= \sum_{q_0 =0}^{\infty } \frac{(\gamma')_{q_0}}{q_0!} s_{0,\infty }^{q_0} \left( c_0 z^{\lambda } \sum_{i_0=0}^{q_0} \frac{(-q_0)_{i_0} \left( q_0 +2 \lambda \right)_{i_0}}{(1+\lambda )_{i_0}\left(\frac{1}{2} +\lambda \right)_{i_0}} \rho ^{i_0} \right) }\cr
{ \displaystyle \mathbf{\Upsilon}(\lambda ; s_0;\widetilde{w}_{1,1}) = \sum_{q_0 =0}^{\infty } \frac{(\gamma')_{q_0}}{q_0!} s_0^{q_0}\left(c_0 z^{\lambda} \sum_{i_0=0}^{q_0} \frac{(-q_0)_{i_0} \left( q_0 +2 \lambda \right)_{i_0}}{(1+\lambda )_{i_0}\left( \frac{1}{2} +\lambda \right)_{i_0}} \widetilde{w}_{1,1} ^{i_0} \right) }\cr
{ \displaystyle \mathbf{\Upsilon}(\lambda; s_0 ;\widetilde{w}_{1,n}) = \sum_{q_0 =0}^{\infty } \frac{(\gamma')_{q_0}}{q_0!} s_0^{q_0}\left(c_0 z^{\lambda} \sum_{i_0=0}^{q_0} \frac{(-q_0)_{i_0} \left( q_0 +2 \lambda \right)_{i_0}}{(1+\lambda )_{i_0}\left(\frac{1}{2} +\lambda \right)_{i_0}} \widetilde{w}_{1,n} ^{i_0} \right)}
\end{cases}\nonumber 
\end{equation}
\end{theorem}
\begin{proof} 
Acting the summation operator $\displaystyle{ \sum_{q_0 =0}^{\infty } \frac{(\gamma')_{q_0}}{q_0!} s_0^{q_0} \prod _{n=1}^{\infty } \left\{ \sum_{ q_n = q_{n-1}}^{\infty } s_n^{q_n }\right\}}$ on the form of integrals of Lame equation  of type 2 $y(z)$,
\begin{eqnarray}
&&\sum_{\alpha _0 =0}^{\infty } \frac{(\gamma')_{q_0}}{q_0!} s_0^{q_0} \prod _{n=1}^{\infty } \left\{ \sum_{ q_n = q_{n-1}}^{\infty } s_n^{q_n }\right\} y(z) \label{eq:70043}\\
&&= \sum_{q_0 =0}^{\infty } \frac{(\gamma')_{q_0}}{q_0!} s_0^{q_0} \prod _{n=1}^{\infty } \left\{ \sum_{ q_n = q_{n-1}}^{\infty } s_n^{q_n }\right\} \Big\{ y_0(z)+y_1(z)+y_2(z)+y_3(z)+\cdots\Big\} \nonumber
\end{eqnarray}
Acting the summation operator $\displaystyle{ \sum_{q_0 =0}^{\infty } \frac{(\gamma')_{q_0}}{q_0!} s_0^{q_0} \prod _{n=1}^{\infty } \left\{ \sum_{ q_n = q_{n-1}}^{\infty } s_n^{q_n }\right\}}$ on (\ref{eq:70025a}),
\begin{eqnarray}
&&\sum_{q_0 =0}^{\infty } \frac{(\gamma')_{q_0}}{q_0!} s_0^{q_0} \prod _{n=1}^{\infty } \left\{ \sum_{ q_n = q_{n-1}}^{\infty } s_n^{q_n }\right\} y_0(z) \nonumber\\
&&= \prod_{l=1}^{\infty } \frac{1}{(1-s_{l,\infty })} \sum_{q_0 =0}^{\infty } \frac{(\gamma')_{q_0}}{q_0!} s_{0,\infty }^{q_0} \left(c_0 z^{\lambda } \sum_{i_0=0}^{q_0} \frac{(-q_0)_{i_0} \left(q_0 +2\lambda \right)_{i_0}}{\left( 1+\lambda \right)_{i_0}\left( \frac{1}{2} +\lambda \right)_{i_0}} \rho ^{i_0} \right) \hspace{2cm}\label{eq:70044}
\end{eqnarray}
Acting the summation operator $\displaystyle{ \sum_{q_0 =0}^{\infty } \frac{(\gamma')_{q_0}}{q_0!} s_0^{q_0} \prod _{n=1}^{\infty } \left\{ \sum_{ q_n = q_{n-1}}^{\infty } s_n^{q_n }\right\}}$ on (\ref{eq:70026}),
\begin{eqnarray}
&&\sum_{q_0 =0}^{\infty } \frac{(\gamma')_{q_0}}{q_0!} s_0^{q_0} \prod _{n=1}^{\infty } \left\{ \sum_{ q_n = q_{n-1}}^{\infty } s_n^{q_n }\right\} y_1(z) \nonumber\\
&&= \prod_{l=2}^{\infty } \frac{1}{(1-s_{l,\infty })} \int_{0}^{1} dt_1\;t_1^{1+\lambda } \int_{0}^{1} du_1\;u_1^{\frac{1}{2} +\lambda}
 \frac{1}{2\pi i} \oint dv_1 \;\frac{1}{v_1} (1-\rho (1-t_1)(1-u_1)v_1)^{-\left( 4 +2\lambda \right)} \nonumber\\
&&\times \sum_{q_1 =q_0}^{\infty }\left( \frac{v_1-1}{v_1} \frac{s_{1,\infty }}{1-\rho (1-t_1)(1-u_1)v_1}\right)^{q_1}  \overleftrightarrow {w}_{1,1}^{-(\frac{1}{2}+\frac{\alpha }{2}+\lambda )} \left(\overleftrightarrow {w}_{1,1} \partial_{\overleftrightarrow {w}_{1,1}} \right) \overleftrightarrow {w}_{1,1}^{\alpha +\frac{1}{2} } \left(\overleftrightarrow {w}_{1,1} \partial_{\overleftrightarrow {w}_{1,1}} \right) \overleftrightarrow {w}_{1,1}^{- \frac{ \alpha }{2}+\lambda }  \nonumber\\
&&\times  \sum_{q_0 =0}^{\infty } \frac{(\gamma' )_{q_0}}{q_0!}s_0^{q_0}\left( c_0 z^{\lambda } \sum_{i_0=0}^{q_0} \frac{(-q_0)_{i_0} \left(q_0 +2 \lambda \right)_{i_0}}{(1+\lambda )_{i_0}\left(\frac{1}{2} +\lambda \right)_{i_0}} \overleftrightarrow {w}_{1,1} ^{i_0} \right) \eta  \label{eq:70045}
\end{eqnarray}
Replace $q_i$, $q_j$ and $r_i$ by $q_1$, $q_0$ and ${ \displaystyle \frac{v_1-1}{v_1} \frac{s_{1,\infty }}{1-\rho (1-t_1)(1-u_1)v_1}}$ in (\ref{eq:70041}). Take the new (\ref{eq:70041}) into (\ref{eq:70045}).
\begin{eqnarray}
&&\sum_{q_0 =0}^{\infty } \frac{(\gamma')_{q_0}}{q_0!} s_0^{q_0} \prod _{n=1}^{\infty } \left\{ \sum_{ q_n = q_{n-1}}^{\infty } s_n^{q_n }\right\} y_1(z) \nonumber\\
&&= \prod_{l=2}^{\infty } \frac{1}{(1-s_{l,\infty })} \int_{0}^{1} dt_1\;t_1^{1+\lambda } \int_{0}^{1} du_1\;u_1^{\frac{1}{2} +\lambda}
 \frac{1}{2\pi i} \oint dv_1 \;\frac{(1-\rho (1-t_1)(1-u_1)v_1)^{-\left(3 +2 \lambda \right)} }{-\rho (1-t_1)(1-u_1)v_1^2+ (1-s_{1,\infty })v_1+s_{1,\infty } } \nonumber\\
&&\times \overleftrightarrow {w}_{1,1}^{-(\frac{1}{2}+\frac{\alpha }{2}+\lambda )}\left(  \overleftrightarrow {w}_{1,1} \partial _{ \overleftrightarrow {w}_{1,1}}\right)  \overleftrightarrow {w}_{1,1}^{\alpha +\frac{1}{2}} \left(  \overleftrightarrow {w}_{1,1} \partial _{ \overleftrightarrow {w}_{1,1}}\right) \overleftrightarrow {w}_{1,1}^{-\frac{\alpha }{2}+\lambda }  \nonumber\\
&&\times \sum_{q_0 =0}^{\infty } \frac{(\gamma' )_{q_0}}{q_0!} \left( \frac{v_1-1}{v_1} \frac{s_{0,\infty }}{1-\rho (1-t_1)(1-u_1)v_1}\right)^{q_0} \left( c_0 z^{\lambda } \sum_{i_0=0}^{q_0} \frac{(-q_0)_{i_0} \left( q_0 +2 \lambda \right)_{i_0}}{(1+\lambda )_{i_0}\left(\frac{1}{2} +\lambda \right)_{i_0}} \overleftrightarrow {w}_{1,1} ^{i_0} \right) \eta   \label{eq:70046}
\end{eqnarray}
By using Cauchy's integral formula, the contour integrand has poles at\\
 ${\displaystyle v_1= \frac{1-s_{1,\infty }-\sqrt{(1-s_{1,\infty })^2+4\rho (1-t_1)(1-u_1)s_{1,\infty }}}{2\rho (1-t_1)(1-u_1)} }$\\  or ${\displaystyle \frac{1-s_{1,\infty }+\sqrt{(1-s_{1,\infty })^2+4\rho (1-t_1)(1-u_1)s_{1,\infty }}}{2\rho (1-t_1)(1-u_1)} }$ and ${ \displaystyle \frac{1-s_{1,\infty }-\sqrt{(1-s_{1,\infty })^2+4\rho (1-t_1)(1-u_1)s_{1,\infty }}}{2\rho (1-t_1)(1-u_1)}}$ is only inside the unit circle. As we compute the residue there in (\ref{eq:70046}) we obtain
\begin{eqnarray}
&&\sum_{q_0 =0}^{\infty } \frac{(\gamma')_{q_0}}{q_0!} s_0^{q_0} \prod _{n=1}^{\infty } \left\{ \sum_{ q_n = q_{n-1}}^{\infty } s_n^{q_n }\right\} y_1(z) \nonumber\\
&&= \prod_{l=2}^{\infty } \frac{1}{(1-s_{l,\infty })} \int_{0}^{1} dt_1\;t_1^{1+\lambda} \int_{0}^{1} du_1\;u_1^{\frac{1}{2} +\lambda}
 \left( s_{1,\infty }^2-2(1-2\rho (1-t_1)(1-u_1))s_{1,\infty }+1 \right)^{-\frac{1}{2}} \nonumber\\
&&\times \left(\frac{1+s_{1,\infty }+\sqrt{s_{1,\infty }^2-2(1-2\rho (1-t_1)(1-u_1))s_{1,\infty }+1}}{2}\right)^{-\left(3 +2 \lambda \right)}\nonumber \\
&&\times \widetilde{w}_{1,1}^{-(\frac{1}{2}+\frac{\alpha }{2}+\lambda )}\left(  \widetilde{w}_{1,1} \partial _{ \widetilde{w}_{1,1}}\right) \widetilde{w}_{1,1}^{\alpha +\frac{1}{2}} \left(  \widetilde{w}_{1,1} \partial _{ \widetilde{w}_{1,1}}\right)\widetilde{w}_{1,1}^{-\frac{\alpha}{2} +\lambda }\nonumber \\
&&\times \sum_{q_0 =0}^{\infty } \frac{(\gamma' )_{q_0}}{q_0!} s_0^{q_0}\left( c_0 z^{\lambda } \sum_{i_0=0}^{q_0} \frac{(-q_0)_{i_0} \left( q_0 +2 \lambda \right)_{i_0}}{(1+\lambda )_{i_0}\left(\frac{1}{2} +\lambda \right)_{i_0}} \widetilde{w}_{1,1} ^{i_0} \right) \eta  \label{eq:70047}
\end{eqnarray}
where
\begin{eqnarray}
\widetilde{w}_{1,1} &=& \frac{v_1}{(v_1-1)}\; \frac{\rho  t_1 u_1}{1- \rho (1-t_1)(1-u_1)v_1}\Bigg|_{\Large v_1=\frac{1-s_{1,\infty }-\sqrt{(1-s_{1,\infty })^2+4\rho (1-t_1)(1-u_1)s_{1,\infty }}}{2\rho (1-t_1)(1-u_1)}\normalsize}\nonumber\\
&=& \frac{\rho  t_1 u_1 \left\{ 1+ (s_{1,\infty }+2\rho (1-t_1)(1-u_1) )s_{1,\infty }\right\}}{2(1-\rho (1-t_1)(1-u_1))^2 s_{1,\infty }}\nonumber\\
&&-\frac{\rho t_1 u_1(1+s_{1,\infty })\sqrt{s_{1,\infty }^2-2(1-2\rho (1-t_1)(1-u_1))s_{1,\infty }+1}}{2(1-\rho (1-t_1)(1-u_1))^2 s_{1,\infty }}\nonumber
\end{eqnarray}
Acting the summation operator $\displaystyle{ \sum_{q_0 =0}^{\infty } \frac{(\gamma')_{q_0}}{q_0!} s_0^{q_0} \prod _{n=1}^{\infty } \left\{ \sum_{ q_n = q_{n-1}}^{\infty } s_n^{q_n }\right\}}$ on (\ref{eq:70028}),
\begin{eqnarray}
&&\sum_{q_0 =0}^{\infty } \frac{(\gamma')_{q_0}}{q_0!} s_0^{q_0} \prod _{n=1}^{\infty } \left\{ \sum_{ q_n = q_{n-1}}^{\infty } s_n^{q_n }\right\} y_2(z) \nonumber\\
&&= \prod_{l=3}^{\infty } \frac{1}{(1-s_{l,\infty })} \int_{0}^{1} dt_2\;t_2^{3+\lambda} \int_{0}^{1} du_2\;u_2^{\frac{5}{2} +\lambda }
 \frac{1}{2\pi i} \oint dv_2 \;\frac{1}{v_2} (1-\rho (1-t_2)(1-u_2)v_2)^{-\left( 8 +2 \lambda \right)} \nonumber\\
&&\times \sum_{q_2 =q_1}^{\infty }\left( \frac{v_2-1}{v_2} \frac{s_{2,\infty }}{1-\rho (1-t_2)(1-u_2)v_2}\right)^{q_2} \overleftrightarrow {w}_{2,2}^{-(\frac{5}{2}+\frac{\alpha}{2} +\lambda )}\left( \overleftrightarrow {w}_{2,2} \partial _{ \overleftrightarrow {w}_{2,2}}\right) \overleftrightarrow {w}_{2,2}^{\alpha +\frac{1}{2}}\left( \overleftrightarrow {w}_{2,2} \partial _{ \overleftrightarrow {w}_{2,2}}\right) \overleftrightarrow {w}_{2,2}^{2-\frac{\alpha}{2} +\lambda }\nonumber\\
&&\times \int_{0}^{1} dt_1\;t_1^{1+\lambda} \int_{0}^{1} du_1\;u_1^{\frac{1}{2} +\lambda}
 \frac{1}{2\pi i} \oint dv_1 \;\frac{1}{v_1} (1-\overleftrightarrow {w}_{2,2} (1-t_1)(1-u_1)v_1)^{-\left( 4 +2 \lambda \right)} \nonumber\\
&&\times \sum_{q_1 =q_0}^{\infty }\left( \frac{v_1-1}{v_1} \frac{s_1}{1-\overleftrightarrow {w}_{2,2}(1-t_1)(1-u_1)v_1}\right)^{q_1}\nonumber\\
&&\times \overleftrightarrow {w}_{1,2}^{-(\frac{1}{2}+\frac{\alpha }{2}+\lambda )}\left(  \overleftrightarrow {w}_{1,2} \partial _{ \overleftrightarrow {w}_{1,2}}\right)  \overleftrightarrow {w}_{1,2}^{\alpha +\frac{1}{2}} \left(  \overleftrightarrow {w}_{1,2} \partial _{ \overleftrightarrow {w}_{1,2}}\right)\overleftrightarrow {w}_{1,2}^{-\frac{\alpha }{2}+\lambda }\nonumber\\
&&\times  \sum_{q_0 =0}^{\infty } \frac{(\gamma' )_{q_0}}{q_0!} s_0^{q_0}\left( c_0 z^{\lambda } \sum_{i_0=0}^{q_0} \frac{(-q_0)_{i_0} \left( q_0 +2 \lambda \right)_{i_0}}{(1+\lambda )_{i_0}\left(\frac{1}{2} +\lambda \right)_{i_0}} \overleftrightarrow {w}_{1,2} ^{i_0} \right)\eta ^2 \label{eq:70048}
\end{eqnarray}
Replace $q_i$, $q_j$ and $r_i$ by $q_2$, $q_1$ and ${ \displaystyle \frac{v_2-1}{v_2} \frac{s_{2,\infty }}{1-\rho (1-t_2)(1-u_2)v_2}}$ in (\ref{eq:70041}). Take the new (\ref{eq:70041}) into (\ref{eq:70048}).
\begin{eqnarray}
&&\sum_{q_0 =0}^{\infty } \frac{(\gamma')_{q_0}}{q_0!} s_0^{q_0} \prod _{n=1}^{\infty } \left\{ \sum_{ q_n = q_{n-1}}^{\infty } s_n^{q_n }\right\} y_2(z) \nonumber\\
&&= \prod_{l=3}^{\infty } \frac{1}{(1-s_{l,\infty })} \int_{0}^{1} dt_2\;t_2^{3+\lambda } \int_{0}^{1} du_2\;u_2^{\frac{5}{2} +\lambda }
 \frac{1}{2\pi i} \oint dv_2 \;\frac{\left( 1-\rho (1-t_2)(1-u_2)v_2\right)^{-\left( 7 +2 \lambda \right)}}{-\rho (1-t_2)(1-u_2)v_2^2+ (1-s_{2,\infty })v_2+s_{2,\infty } } \nonumber\\
&&\times \overleftrightarrow {w}_{2,2}^{-(\frac{5}{2}+\frac{\alpha}{2} +\lambda )}\left( \overleftrightarrow {w}_{2,2} \partial _{ \overleftrightarrow {w}_{2,2}}\right) \overleftrightarrow {w}_{2,2}^{\alpha +\frac{1}{2}}\left( \overleftrightarrow {w}_{2,2} \partial _{ \overleftrightarrow {w}_{2,2}}\right) \overleftrightarrow {w}_{2,2}^{2-\frac{\alpha}{2} +\lambda } \nonumber\\
&&\times \int_{0}^{1} dt_1\;t_1^{1+\lambda} \int_{0}^{1} du_1\;u_1^{\frac{1}{2} +\lambda}
 \frac{1}{2\pi i} \oint dv_1 \;\frac{1}{v_1} \left( 1-\overleftrightarrow {w}_{2,2} (1-t_1)(1-u_1)v_1\right)^{-\left( 4 +2 \lambda \right)} \nonumber\\
&&\times \sum_{q_1 =q_0}^{\infty }\left( \frac{v_2-1}{v_2} \frac{s_{1,\infty }}{1-\rho (1-t_2)(1-u_2)v_2} \frac{v_1-1}{v_1}\frac{1}{1-\overleftrightarrow {w}_{2,2}(1-t_1)(1-u_1)v_1}\right)^{q_1} \nonumber\\
&&\times  \overleftrightarrow {w}_{1,2}^{-(\frac{1}{2}+\frac{\alpha }{2}+\lambda )}\left(  \overleftrightarrow {w}_{1,2} \partial _{ \overleftrightarrow {w}_{1,2}}\right)  \overleftrightarrow {w}_{1,2}^{\alpha +\frac{1}{2}} \left(  \overleftrightarrow {w}_{1,2} \partial _{ \overleftrightarrow {w}_{1,2}}\right)\overleftrightarrow {w}_{1,2}^{-\frac{\alpha }{2}+\lambda }  \nonumber\\
&&\times \sum_{q_0 =0}^{\infty } \frac{(\gamma' )_{q_0}}{q_0!} s_0^{q_0}\left( c_0 z^{\lambda } \sum_{i_0=0}^{q_0} \frac{(-q_0)_{i_0} \left( q_0 +2 \lambda \right)_{i_0}}{(1+\lambda )_{i_0}\left(\frac{1}{2} +\lambda \right)_{i_0}} \overleftrightarrow {w}_{1,2} ^{i_0} \right) \eta ^2 \hspace{1cm}\label{eq:70049}
\end{eqnarray}
By using Cauchy's integral formula, the contour integrand has poles at\\
 ${\displaystyle
 v_2= \frac{1-s_{2,\infty }-\sqrt{(1-s_{2,\infty })^2+4\rho (1-t_2)(1-u_2)s_{2,\infty }}}{2\rho(1-t_2)(1-u_2)}}$ \\or${\displaystyle\frac{1-s_{2,\infty }+\sqrt{(1-s_{2,\infty })^2+4\rho(1-t_2)(1-u_2)s_{2,\infty }}}{2\rho(1-t_2)(1-u_2)} }$ 
and ${ \displaystyle\frac{1-s_{2,\infty }-\sqrt{(1-s_{2,\infty })^2+4\rho(1-t_2)(1-u_2)s_{2,\infty }}}{2\rho(1-t_2)(1-u_2)}}$ is only inside the unit circle. As we compute the residue there in (\ref{eq:70049}) we obtain
\begin{eqnarray}
&&\sum_{q_0 =0}^{\infty } \frac{(\gamma')_{q_0}}{q_0!} s_0^{q_0} \prod _{n=1}^{\infty } \left\{ \sum_{ q_n = q_{n-1}}^{\infty } s_n^{q_n }\right\} y_2(z) \nonumber\\
&&= \prod_{l=3}^{\infty } \frac{1}{(1-s_{l,\infty })} \int_{0}^{1} dt_2\;t_2^{3+\lambda} \int_{0}^{1} du_2\;u_2^{\frac{5}{2} +\lambda }
 \left( s_{2,\infty }^2-2(1-2\rho (1-t_2)(1-u_2))s_{2,\infty }+1\right)^{-\frac{1}{2}}\nonumber\\
&&\times \left(\frac{1+s_{2,\infty }+\sqrt{s_{2,\infty }^2-2(1-2\rho (1-t_2)(1-u_2))s_{2,\infty }+1}}{2}\right)^{-\left( 7 +2\lambda \right)}  \nonumber\\
&&\times \widetilde{w}_{2,2}^{-(\frac{5}{2}+\frac{\alpha }{2} +\lambda )}\left( \widetilde{w}_{2,2} \partial _{ \widetilde{w}_{2,2}}\right) \widetilde{w}_{2,2}^{\alpha+\frac{1}{2}} \left( \widetilde{w}_{2,2} \partial _{ \widetilde{w}_{2,2}}\right) \widetilde{w}_{2,2}^{2-\frac{\alpha}{2} +\lambda}\nonumber\\
&&\times \int_{0}^{1} dt_1\;t_1^{1+\lambda} \int_{0}^{1} du_1\;u_1^{\frac{1}{2} +\lambda}
 \frac{1}{2\pi i} \oint dv_1 \;\frac{1}{v_1} \left( 1-\widetilde{w}_{2,2} (1-t_1)(1-u_1)v_1\right)^{-\left( 4 +2 \lambda \right)} \nonumber\\
&&\times \sum_{q_1 =q_0}^{\infty }\left( \frac{v_1-1}{v_1}\frac{s_1}{1-\widetilde{w}_{2,2}(1-t_1)(1-u_1)v_1}\right)^{q_1}
 \ddot{w}_{1,2}^{-(\frac{1}{2}+\frac{\alpha }{2}+\lambda )}\left( \ddot{w}_{1,2} \partial _{ \ddot{w}_{1,2}}\right) \ddot{w}_{1,2}^{\alpha +\frac{1}{2}} \left( \ddot{w}_{1,2} \partial _{ \ddot{w}_{1,2}}\right) \ddot{w}_{1,2}^{-\frac{\alpha}{2} +\lambda }\nonumber\\
&&\times  \sum_{q_0 =0}^{\infty } \frac{(\gamma' )_{q_0}}{q_0!} s_0^{q_0}\left( c_0 z^{\lambda } \sum_{i_0=0}^{q_0}  \frac{(-q_0)_{i_0} \left( q_0 +2 \lambda \right)_{i_0}}{(1+\lambda )_{i_0}\left(\frac{1}{2} +\lambda \right)_{i_0}} \ddot{w}_{1,2} ^{i_0} \right) \eta ^2 \label{eq:70050}
\end{eqnarray}
where
\begin{eqnarray}
\widetilde{w}_{2,2} &=& \frac{v_2}{(v_2-1)}\; \frac{\rho t_2 u_2}{1- \rho (1-t_2)(1-u_2)v_2}\Bigg|_{\Large v_2=\frac{1-s_{2,\infty }-\sqrt{(1-s_{2,\infty })^2+4\rho (1-t_2)(1-u_2)s_{2,\infty }}}{2\rho (1-t_2)(1-u_2)}\normalsize}\nonumber\\
&=& \frac{\rho t_2 u_2 \left\{ 1+ (s_{2,\infty }+2\rho (1-t_2)(1-u_2) )s_{2,\infty }\right\}}{2(1-\rho (1-t_2)(1-u_2))^2 s_{2,\infty }} \nonumber\\
&-& \frac{\rho t_2 u_2 (1+s_{2,\infty })\sqrt{s_{2,\infty }^2-2(1-2\rho (1-t_2)(1-u_2))s_{2,\infty }+1}}{2(1-\rho (1-t_2)(1-u_2))^2 s_{2,\infty }}  \nonumber
\end{eqnarray}
and
\begin{equation}
\ddot{w}_{1,2} = \frac{v_1}{(v_1-1)}\; \frac{\widetilde{w}_{2,2} t_1 u_1}{1- \widetilde{w}_{2,2}(1-t_1)(1-u_1)v_1}\nonumber
\end{equation}
Replace $q_i$, $q_j$ and $r_i$ by $q_1$, $q_0$ and ${ \displaystyle \frac{v_1-1}{v_1}\frac{s_1}{1-\widetilde{w}_{2,2}(1-t_1)(1-u_1)v_1}}$ in (\ref{eq:70041}). Take the new (\ref{eq:70041}) into (\ref{eq:70050}).
\begin{eqnarray}
&&\sum_{q_0 =0}^{\infty } \frac{(\gamma')_{q_0}}{q_0!} s_0^{q_0} \prod _{n=1}^{\infty } \left\{ \sum_{ q_n = q_{n-1}}^{\infty } s_n^{q_n }\right\} y_2(z) \nonumber\\
&&= \prod_{l=3}^{\infty } \frac{1}{(1-s_{l,\infty })} \int_{0}^{1} dt_2\;t_2^{3+\lambda} \int_{0}^{1} du_2\;u_2^{\frac{5}{2} +\lambda }
 \left( s_{2,\infty }^2-2(1-2\rho (1-t_2)(1-u_2))s_{2,\infty }+1\right)^{-\frac{1}{2}}\nonumber\\
&&\times \left(\frac{1+s_{2,\infty }+\sqrt{s_{2,\infty }^2-2(1-2\rho (1-t_2)(1-u_2))s_{2,\infty }+1}}{2}\right)^{-\left( 7 +2\lambda \right)} \nonumber\\
&&\times \widetilde{w}_{2,2}^{-(\frac{5}{2}+\frac{\alpha }{2} +\lambda )}\left( \widetilde{w}_{2,2} \partial _{ \widetilde{w}_{2,2}}\right) \widetilde{w}_{2,2}^{ \alpha+\frac{1}{2}} \left( \widetilde{w}_{2,2} \partial _{ \widetilde{w}_{2,2}}\right) \widetilde{w}_{2,2}^{2-\frac{\alpha}{2} +\lambda}\nonumber\\
&&\times \int_{0}^{1} dt_1\;t_1^{1+\lambda} \int_{0}^{1} du_1\;u_1^{\frac{1}{2} +\lambda}
 \frac{1}{2\pi i} \oint dv_1 \;\frac{\left(1-\widetilde{w}_{2,2} (1-t_1)(1-u_1)v_1\right)^{-\left( 3 +2 \lambda \right)}}{-\widetilde{w}_{2,2} (1-t_1)(1-u_1)v_1^2 +(1-s_1)v_1 +s_1} \nonumber\\
&&\times  \ddot{w}_{1,2}^{-(\frac{1}{2}+\frac{\alpha }{2}+\lambda )}\left( \ddot{w}_{1,2} \partial _{ \ddot{w}_{1,2}}\right) \ddot{w}_{1,2}^{\alpha +\frac{1}{2}} \left( \ddot{w}_{1,2} \partial _{ \ddot{w}_{1,2}}\right) \ddot{w}_{1,2}^{-\frac{\alpha }{2} +\lambda } \sum_{q_0 =0}^{\infty } \frac{(\gamma')_{q_0}}{q_0!} \left( \frac{v_1-1}{v_1}\frac{s_{0,1}}{1-\widetilde{w}_{2,2}(1-t_1)(1-u_1)v_1}\right)^{q_0} \nonumber\\
&&\times \left( c_0 z^{\lambda } \sum_{i_0=0}^{q_0} \frac{(-q_0)_{i_0} \left( q_0 +2 \lambda \right)_{i_0}}{(1+\lambda )_{i_0}\left(\frac{1}{2} +\lambda \right)_{i_0}} \ddot{w}_{1,2} ^{i_0} \right) \eta ^2  \label{eq:70051}
\end{eqnarray}
By using Cauchy's integral formula, the contour integrand has poles at\\ ${\displaystyle
 v_1= \frac{1-s_1-\sqrt{(1-s_1)^2+4\widetilde{w}_{2,2} (1-t_1)(1-u_1)s_1}}{2\widetilde{w}_{2,2} (1-t_1)(1-u_1)}}$\\ or ${\displaystyle\frac{1-s_1+\sqrt{(1-s_1)^2+4\widetilde{w}_{2,2} (1-t_1)(1-u_1)s_1}}{2\widetilde{w}_{2,2} (1-t_1)(1-u_1)}}$ 
and ${ \displaystyle \frac{1-s_1-\sqrt{(1-s_1)^2+4\widetilde{w}_{2,2} (1-t_1)(1-u_1)s_1}}{2\widetilde{w}_{2,2} (1-t_1)(1-u_1)}}$ is only inside the unit circle. As we compute the residue there in (\ref{eq:70051}) we obtain
\begin{eqnarray}
&&\sum_{q_0 =0}^{\infty } \frac{(\gamma')_{q_0}}{q_0!} s_0^{q_0} \prod _{n=1}^{\infty } \left\{ \sum_{ q_n = q_{n-1}}^{\infty } s_n^{q_n }\right\} y_2(z) \nonumber\\
&&= \prod_{l=3}^{\infty } \frac{1}{(1-s_{l,\infty })} \int_{0}^{1} dt_2\;t_2^{3+\lambda} \int_{0}^{1} du_2\;u_2^{\frac{5}{2} +\lambda }
 \left( s_{2,\infty }^2-2(1-2\rho (1-t_2)(1-u_2))s_{2,\infty }+1\right)^{-\frac{1}{2}}\nonumber\\
&&\times \left(\frac{1+s_{2,\infty }+\sqrt{s_{2,\infty }^2-2(1-2\rho (1-t_2)(1-u_2))s_{2,\infty }+1}}{2}\right)^{-\left( 7 +2\lambda \right)} \nonumber\\
&&\times \widetilde{w}_{2,2}^{-(\frac{5}{2}+\frac{\alpha }{2} +\lambda )}\left( \widetilde{w}_{2,2} \partial _{ \widetilde{w}_{2,2}}\right) \widetilde{w}_{2,2}^{ \alpha+\frac{1}{2}} \left( \widetilde{w}_{2,2} \partial _{ \widetilde{w}_{2,2}}\right) \widetilde{w}_{2,2}^{2-\frac{\alpha}{2} +\lambda}\nonumber\\
&&\times \int_{0}^{1} dt_1\;t_1^{1+\lambda} \int_{0}^{1} du_1\;u_1^{\frac{1}{2} +\lambda}
 \left( s_1^2-2(1-2\widetilde{w}_{2,2} (1-t_1)(1-u_1))s_1+1\right)^{-\frac{1}{2}}\nonumber\\
&&\times \left(\frac{1+s_1+\sqrt{s_1^2-2(1-2\widetilde{w}_{2,2}(1-t_1)(1-u_1))s_1+1}}{2}\right)^{-\left(3 +2 \lambda \right)} \nonumber\\
&&\times \widetilde{w}_{1,2}^{-(\frac{1}{2}+\frac{\alpha }{2}+\lambda )}\left( \widetilde{w}_{1,2} \partial _{ \widetilde{w}_{1,2}}\right) \widetilde{w}_{1,2}^{\alpha +\frac{1}{2}} \left( \widetilde{w}_{1,2} \partial _{ \widetilde{w}_{1,2}}\right) \widetilde{w}_{1,2}^{-\frac{\alpha }{2}+\lambda } \nonumber\\
&&\times \sum_{q_0 =0}^{\infty } \frac{(\gamma' )_{q_0}}{q_0!} s_0^{q_0}\left( c_0 z^{\lambda } \sum_{i_0=0}^{q_0} \frac{(-q_0)_{i_0} \left( q_0 +2 \lambda \right)_{i_0}}{(1+\lambda )_{i_0}\left(\frac{1}{2} +\lambda \right)_{i_0}} \widetilde{w}_{1,2} ^{i_0} \right) \eta ^2  \label{eq:70052}
\end{eqnarray}
where
\begin{eqnarray}
\widetilde{w}_{1,2} &=& \frac{v_1}{(v_1-1)}\; \frac{\widetilde{w}_{2,2} t_1 u_1}{1- \widetilde{w}_{2,2} (1-t_1)(1-u_1)v_1}\Bigg|_{\Large v_1=\frac{1-s_1-\sqrt{(1-s_1)^2+4\widetilde{w}_{2,2} (1-t_1)(1-u_1)s_1}}{2\widetilde{w}_{2,2} (1-t_1)(1-u_1)}\normalsize}\nonumber\\
&=& \frac{\widetilde{w}_{2,2} t_1 u_1 \left\{ 1+ (s_1+2\widetilde{w}_{2,2}(1-t_1)(1-u_1) )s_1 \right\}}{2(1-\widetilde{w}_{2,2}(1-t_1)(1-u_1))^2 s_1}\nonumber\\
&-& \frac{\widetilde{w}_{2,2} t_1 u_1 (1+s_1)\sqrt{s_1^2-2(1-2\widetilde{w}_{2,2} (1-t_1)(1-u_1))s_1+1}}{2(1-\widetilde{w}_{2,2}(1-t_1)(1-u_1))^2 s_1} \nonumber
\end{eqnarray}
Acting the summation operator $\displaystyle{ \sum_{q_0 =0}^{\infty } \frac{(\gamma')_{q_0}}{q_0!} s_0^{q_0} \prod _{n=1}^{\infty } \left\{ \sum_{ q_n = q_{n-1}}^{\infty } s_n^{q_n }\right\}}$ on (\ref{eq:70029}),
\begin{eqnarray}
&&\sum_{q_0 =0}^{\infty } \frac{(\gamma')_{q_0}}{q_0!} s_0^{q_0} \prod _{n=1}^{\infty } \left\{ \sum_{ q_n = q_{n-1}}^{\infty } s_n^{q_n }\right\} y_3(z) \nonumber\\
&&= \prod_{l=4}^{\infty } \frac{1}{(1-s_{l,\infty })} \int_{0}^{1} dt_3\;t_3^{5+\lambda} \int_{0}^{1} du_3\;u_3^{\frac{9}{2} +\lambda}\left( s_{3,\infty }^2-2(1-2\rho (1-t_3)(1-u_3))s_{3,\infty }+1\right)^{-\frac{1}{2}}\nonumber\\
&&\times \left(\frac{1+s_{3,\infty }+\sqrt{s_{3,\infty }^2-2(1-2\rho (1-t_3)(1-u_3))s_{3,\infty }+1}}{2}\right)^{-\left(11 +2 \lambda   \right)}\nonumber\\
&&\times  \widetilde{w}_{3,3}^{-(\frac{9}{2}+\frac{\alpha }{2}+\lambda )}\left( \widetilde{w}_{3,3} \partial _{ \widetilde{w}_{3,3}}\right) \widetilde{w}_{3,3}^{ \alpha +\frac{1}{2}} \left( \widetilde{w}_{3,3} \partial _{ \widetilde{w}_{3,3}}\right) \widetilde{w}_{3,3}^{ 4-\frac{\alpha }{2}+\lambda } \nonumber\\
&&\times \int_{0}^{1} dt_2\;t_2^{3+\lambda} \int_{0}^{1} du_2\;u_2^{\frac{5}{2} +\lambda }\left( s_2^2-2(1-2\widetilde{w}_{3,3}(1-t_2)(1-u_2))s_2+1\right)^{-\frac{1}{2}}\nonumber\\
&&\times \left(\frac{1+s_2+\sqrt{s_2^2-2(1-2\widetilde{w}_{3,3}(1-t_2)(1-u_2))s_2+1}}{2}\right)^{-\left( 7 +2 \lambda \right)} \nonumber\\
&&\times \widetilde{w}_{2,3}^{-(\frac{5}{2}+\frac{\alpha }{2}+\lambda )}\left( \widetilde{w}_{2,3} \partial _{ \widetilde{w}_{2,3}}\right) \widetilde{w}_{2,3}^{ \alpha +\frac{1}{2}} \left( \widetilde{w}_{2,3} \partial _{ \widetilde{w}_{2,3}}\right) \widetilde{w}_{2,3}^{ 2-\frac{\alpha }{2}+\lambda } \nonumber\\
&&\times \int_{0}^{1} dt_1\;t_1^{1+\lambda } \int_{0}^{1} du_1\;u_1^{\frac{1}{2} +\lambda}\left( s_1^2-2(1-2\widetilde{w}_{2,3} (1-t_1)(1-u_1))s_1+1\right)^{-\frac{1}{2}}\nonumber\\
&&\times \left(\frac{1+s_1+\sqrt{s_1^2-2(1-2\widetilde{w}_{2,3}(1-t_1)(1-u_1))s_1+1}}{2}\right)^{-\left(3 +2 \lambda \right)}\nonumber\\
&&\times \widetilde{w}_{1,3}^{-(\frac{1}{2}+\frac{\alpha }{2}+\lambda )}\left( \widetilde{w}_{1,3} \partial _{ \widetilde{w}_{1,3}}\right) \widetilde{w}_{1,3}^{ \alpha +\frac{1}{2}} \left( \widetilde{w}_{1,3} \partial _{ \widetilde{w}_{1,3}}\right) \widetilde{w}_{1,3}^{ -\frac{\alpha }{2}+\lambda } \nonumber\\
&&\times \sum_{q_0 =0}^{\infty } \frac{(\gamma' )_{q_0}}{q_0!} s_0^{q_0}\left( c_0 z^{\lambda } \sum_{i_0=0}^{q_0} \frac{(-q_0)_{i_0} \left( q_0 +2 \lambda \right)_{i_0}}{(1+\lambda )_{i_0}\left(\frac{1}{2} +\lambda \right)_{i_0}} \widetilde{w}_{1,3} ^{i_0} \right) \eta ^3  \label{eq:70053}
\end{eqnarray}

\vspace{1cm}
where
\begin{eqnarray}
\widetilde{w}_{3,3} &=& \frac{v_3}{(v_3-1)}\; \frac{\rho  t_3 u_3}{1- \rho (1-t_3)(1-u_3)v_3}\Bigg|_{\Large v_3=\frac{1-s_{3,\infty }-\sqrt{(1-s_{3,\infty })^2+4\rho (1-t_3)(1-u_3)s_{3,\infty }}}{2\rho (1-t_3)(1-u_3)}\normalsize}\nonumber\\
&=& \frac{\rho  t_3 u_3 \left\{ 1+ (s_{3,\infty }+2\rho (1-t_3)(1-u_3) )s_{3,\infty } \right\}}{2(1-\rho (1-t_3)(1-u_3))^2 s_{3,\infty }}\nonumber\\
&-& \frac{\rho  t_3 u_3  (1+s_{3,\infty })\sqrt{s_{3,\infty }^2-2(1-2\rho  (1-t_3)(1-u_3))s_{3,\infty }+1} }{2(1-\rho (1-t_3)(1-u_3))^2 s_{3,\infty }}\nonumber
\end{eqnarray}
\begin{eqnarray}
\widetilde{w}_{2,3} &=& \frac{v_2}{(v_2-1)}\; \frac{\widetilde{w}_{3,3} t_2 u_2}{1- \widetilde{w}_{3,3} (1-t_2)(1-u_2)v_2 }\Bigg|_{\Large v_2=\frac{1-s_2-\sqrt{(1-s_2)^2+4\widetilde{w}_{3,3} (1-t_2)(1-u_2)s_2}}{2\widetilde{w}_{3,3} (1-t_2)(1-u_2)}\normalsize}\nonumber\\
&=& \frac{\widetilde{w}_{3,3} t_2 u_2 \left\{ 1+ (s_2+2\widetilde{w}_{3,3}(1-t_2)(1-u_2) )s_2 \right\}}{2(1-\widetilde{w}_{3,3}(1-t_2)(1-u_2))^2 s_2}\nonumber\\
&-& \frac{\widetilde{w}_{3,3} t_2 u_2  (1+s_2)\sqrt{s_2^2-2(1-2\widetilde{w}_{3,3} (1-t_2)(1-u_2))s_2+1} }{2(1-\widetilde{w}_{3,3}(1-t_2)(1-u_2))^2 s_2}\nonumber
\end{eqnarray}
\begin{eqnarray}
\widetilde{w}_{1,3} &=& \frac{v_1}{(v_1-1)}\; \frac{\widetilde{w}_{2,3} t_1 u_1}{1- \widetilde{w}_{2,3} (1-t_1)(1-u_1)v_1 }\Bigg|_{\Large v_1=\frac{1-s_1-\sqrt{(1-s_1)^2+4\widetilde{w}_{2,3} (1-t_1)(1-u_1)s_1}}{2\widetilde{w}_{2,3} (1-t_1)(1-u_1)}\normalsize}\nonumber\\
&=& \frac{\widetilde{w}_{2,3} t_1 u_1 \left\{ 1+ (s_1+2\widetilde{w}_{2,3}(1-t_1)(1-u_1) )s_1 \right\}}{2(1-\widetilde{w}_{2,3}(1-t_1)(1-u_1))^2 s_1}\nonumber\\
&-& \frac{\widetilde{w}_{2,3} t_1 u_1  (1+s_1)\sqrt{s_1^2-2(1-2\widetilde{w}_{2,3} (1-t_1)(1-u_1))s_1+1} }{2(1-\widetilde{w}_{2,3}(1-t_1)(1-u_1))^2 s_1} \nonumber
\end{eqnarray}
By repeating this process for all higher terms of integral forms of sub-summation $y_m(z)$ terms where $m \geq 4$, I obtain every  $\displaystyle{ \sum_{q_0 =0}^{\infty } \frac{(\gamma')_{q_0}}{q_0!} s_0^{q_0} \prod _{n=1}^{\infty } \left\{ \sum_{ q_n = q_{n-1}}^{\infty } s_n^{q_n }\right\}}  y_m(z)$ terms. 
Substitute (\ref{eq:70044}), (\ref{eq:70047}), (\ref{eq:70052}), (\ref{eq:70053}) and including all $\displaystyle{ \sum_{q_0 =0}^{\infty } \frac{(\gamma')_{q_0}}{q_0!} s_0^{q_0} \prod _{n=1}^{\infty } \left\{ \sum_{ q_n = q_{n-1}}^{\infty } s_n^{q_n }\right\}}  y_m(z)$ terms where $m > 3$ into (\ref{eq:70043}). 
\qed
\end{proof}
\begin{remark}
The generating function for the Lame polynomial of type 2 in the algebraic form of the first kind about $x=a$ as $q= \alpha (\alpha +1)a- 4(2a-b-c)(q_j+2j )^2 $ where $j,q_j \in \mathbb{N}_{0}$ is
\begin{eqnarray}
&&\sum_{q_0 =0}^{\infty } \frac{ (\frac{1}{2} )_{q_0}}{q_0!} s_0^{q_0} \prod _{n=1}^{\infty } \left\{ \sum_{ q_n = q_{n-1}}^{\infty } s_n^{q_n }\right\} LF_{q_j}^R\Bigg( a, b, c, \alpha, q= \alpha (\alpha +1)a- 4(2a-b-c)(q_j+2j )^2 \nonumber\\
&&; z= x-a, \rho = -\frac{2a-b-c}{(a-b)(a-c)} z, \eta = \frac{-z^2}{(a-b)(a-c)} \Bigg) \nonumber\\
&&=2^{-1}\Bigg\{ \prod_{l=1}^{\infty } \frac{1}{(1-s_{l,\infty })}  \mathbf{A}\left( s_{0,\infty } ;\rho \right) + \Bigg\{ \prod_{l=2}^{\infty } \frac{1}{(1-s_{l,\infty })} \int_{0}^{1} dt_1\;t_1 \int_{0}^{1} du_1\;u_1^{\frac{1}{2}} \overleftrightarrow {\mathbf{\Gamma}}_1 \left(s_{1,\infty };t_1,u_1,\rho \right)\nonumber\\
&&\times \widetilde{w}_{1,1}^{-\frac{1}{2}(1+\alpha )}\left( \widetilde{w}_{1,1} \partial _{ \widetilde{w}_{1,1}}\right) \widetilde{w}_{1,1}^{\alpha +\frac{1}{2}} \left( \widetilde{w}_{1,1} \partial _{ \widetilde{w}_{1,1}}\right) \widetilde{w}_{1,1}^{-\frac{\alpha }{2}}\mathbf{A}\left( s_{0} ;\widetilde{w}_{1,1}\right) \Bigg\} \eta  \nonumber\\
&&+ \sum_{n=2}^{\infty } \Bigg\{ \prod_{l=n+1}^{\infty } \frac{1}{(1-s_{l,\infty })} \int_{0}^{1} dt_n\;t_n^{2n-1} \int_{0}^{1} du_n\;u_n^{2n-\frac{3}{2} } \overleftrightarrow {\mathbf{\Gamma}}_n \left(s_{n,\infty };t_n,u_n,\rho \right) \nonumber\\
&&\times \widetilde{w}_{n,n}^{-\frac{1}{2}(4n-3+\alpha )}\left( \widetilde{w}_{n,n} \partial _{ \widetilde{w}_{n,n}}\right) \widetilde{w}_{n,n}^{ \alpha +\frac{1}{2}} \left( \widetilde{w}_{n,n} \partial _{ \widetilde{w}_{n,n}}\right) \widetilde{w}_{n,n}^{ \frac{1}{2}(4(n-1)-\alpha )} \nonumber\\
&&\times \prod_{k=1}^{n-1} \Bigg\{ \int_{0}^{1} dt_{n-k}\;t_{n-k}^{2(n-k)-1} \int_{0}^{1} du_{n-k} \;u_{n-k}^{2(n-k)-\frac{3}{2}}\overleftrightarrow {\mathbf{\Gamma}}_{n-k} \left(s_{n-k};t_{n-k},u_{n-k},\widetilde{w}_{n-k+1,n} \right)\label{eq:70054}\\
&&\times \widetilde{w}_{n-k,n}^{-\frac{1}{2}(4(n-k)-3+\alpha)}\left( \widetilde{w}_{n-k,n} \partial _{ \widetilde{w}_{n-k,n}}\right)  \widetilde{w}_{n-k,n}^{ \alpha +\frac{1}{2}} \left( \widetilde{w}_{n-k,n} \partial _{ \widetilde{w}_{n-k,n}}\right)  \widetilde{w}_{n-k,n}^{ \frac{1}{2}(4(n-k-1)-\alpha )} \Bigg\} \mathbf{A} \left( s_{0} ;\widetilde{w}_{1,n}\right) \Bigg\} \eta ^n \Bigg\}  \nonumber   
\end{eqnarray}
where
\begin{equation}
\begin{cases} 
{ \displaystyle \overleftrightarrow {\mathbf{\Gamma}}_1 \left(s_{1,\infty };t_1,u_1,\rho \right)= \frac{\left( \frac{1+s_{1,\infty }+\sqrt{s_{1,\infty }^2-2(1-2\rho (1-t_1)(1-u_1))s_{1,\infty }+1}}{2}\right)^{-3}}{\sqrt{s_{1,\infty }^2-2(1-2\rho (1-t_1)(1-u_1))s_{1,\infty }+1}}}\cr
{ \displaystyle  \overleftrightarrow {\mathbf{\Gamma}}_n \left(s_{n,\infty };t_n,u_n,\rho \right) =\frac{\left( \frac{1+s_{n,\infty }+\sqrt{s_{n,\infty }^2-2(1-2\rho (1-t_n)(1-u_n))s_{n,\infty }+1}}{2}\right)^{-\left( 4n-1 \right)}}{\sqrt{ s_{n,\infty }^2-2(1-2\rho (1-t_n)(1-u_n))s_{n,\infty }+1}}}\cr
{ \displaystyle \overleftrightarrow {\mathbf{\Gamma}}_{n-k} \left( s_{n-k};t_{n-k},u_{n-k},\widetilde{w}_{n-k+1,n} \right) = \frac{ \left( \frac{1+s_{n-k}+\sqrt{s_{n-k}^2-2(1-2\widetilde{w}_{n-k+1,n} (1-t_{n-k})(1-u_{n-k}))s_{n-k}+1}}{2}\right)^{-\left( 4(n-k)-1 \right)}}{\sqrt{ s_{n-k}^2-2(1-2\widetilde{w}_{n-k+1,n} (1-t_{n-k})(1-u_{n-k}))s_{n-k}+1}}}
\end{cases}\nonumber 
\end{equation}
and
\begin{equation}
\begin{cases} 
{ \displaystyle \mathbf{A} \left( s_{0,\infty } ;\rho \right)= \frac{\left(1- s_{0,\infty }+\sqrt{s_{0,\infty }^2-2(1-2\rho )s_{0,\infty }+1}\right)^{\frac{1}{2}} \left(1+s_{0,\infty }+\sqrt{s_{0,\infty }^2-2(1-2\rho  )s_{0,\infty }+1}\right)^{\frac{1}{2}}}{\sqrt{s_{0,\infty }^2-2(1-2\rho )s_{0,\infty }+1}}}\cr
{ \displaystyle  \mathbf{A} \left( s_{0} ;\widetilde{w}_{1,1}\right) = \frac{\left( 1-s_0+\sqrt{s_0^2-2(1-2\widetilde{w}_{1,1})s_0+1}\right)^{\frac{1}{2}} \left( 1+s_0+\sqrt{s_0^2-2(1-2\widetilde{w}_{1,1} )s_0+1}\right)^{\frac{1}{2}}}{\sqrt{s_0^2-2(1-2\widetilde{w}_{1,1})s_0+1}}} \cr
{ \displaystyle \mathbf{A} \left( s_{0} ;\widetilde{w}_{1,n}\right) = \frac{\left( 1-s_0+\sqrt{s_0^2-2(1-2\widetilde{w}_{1,n})s_0+1}\right)^{\frac{1}{2}} \left(1+s_0+\sqrt{s_0^2-2(1-2\widetilde{w}_{1,n} )s_0+1}\right)^{\frac{1}{2}}}{\sqrt{s_0^2-2(1-2\widetilde{w}_{1,n})s_0+1}}}
\end{cases}\nonumber 
\end{equation}
\end{remark}
\begin{proof}
Replace $w$, $\gamma $, $A$ and $x$ by $s_{0,\infty }$, $\frac{1}{2}$, 0 and $\rho $ in (\ref{eq:70035}). 
\begin{eqnarray}
&&\sum_{q_0=0}^{\infty }\frac{(\frac{1}{2})_{q_0}}{q_0!} s_{0,\infty }^{q_0} \;_2F_1\left( -q_0, q_0; \frac{1}{2} ; \rho \right) \label{eq:70055}\\
&&=  \frac{\left(1- s_{0,\infty }+\sqrt{s_{0,\infty }^2-2(1-2\rho )s_{0,\infty }+1}\right)^{\frac{1}{2}} \left(1+s_{0,\infty }+\sqrt{s_{0,\infty }^2-2(1-2\rho )s_{0,\infty }+1}\right)^{\frac{1}{2}}}{2\sqrt{s_{0,\infty }^2-2(1-2\rho )s_{0,\infty }+1}} \nonumber
\end{eqnarray} 
Replace $s_{0,\infty }$ and $\rho $  by  $s_0$ and $\widetilde{w}_{1,1}$ in (\ref{eq:70055}). 
\begin{eqnarray}
&&\sum_{q_0=0}^{\infty }\frac{(\frac{1}{2} )_{q_0}}{q_0!} s_0^{q_0} \;_2F_1\left(-q_0, q_0; \frac{1}{2}; \widetilde{w}_{1,1} \right) \label{eq:70056}\\
&&= \frac{\left( 1-s_0+\sqrt{s_0^2-2(1-2\widetilde{w}_{1,1})s_0+1}\right)^{\frac{1}{2}} \left( 1+s_0+\sqrt{s_0^2-2(1-2\widetilde{w}_{1,1} )s_0+1}\right)^{\frac{1}{2}}}{2\sqrt{s_0^2-2(1-2\widetilde{w}_{1,1})s_0+1}} \nonumber
\end{eqnarray} 
Replace  $\widetilde{w}_{1,1}$  by $\widetilde{w}_{1,n}$  in (\ref{eq:70056}). 
\begin{eqnarray}
&&\sum_{q_0=0}^{\infty }\frac{(\frac{1}{2} )_{q_0}}{q_0!} s_0^{q_0} \;_2F_1\left( -q_0, q_0; \frac{1}{2}; \widetilde{w}_{1,n} \right) \label{eq:70057}\\
&&=  \frac{ \left( 1-s_0+\sqrt{s_0^2-2(1-2\widetilde{w}_{1,n})s_0+1}\right)^{\frac{1}{2}} \left( 1+s_0+\sqrt{s_0^2-2(1-2\widetilde{w}_{1,n} )s_0+1}\right)^{\frac{1}{2}}}{2\sqrt{s_0^2-2(1-2\widetilde{w}_{1,n})s_0+1}} \nonumber
\end{eqnarray} 
Put $c_0$= 1, $\lambda =0$ and $\gamma' =\frac{1}{2} $ in (\ref{eq:70042}). Substitute (\ref{eq:70055}), (\ref{eq:70056}) and (\ref{eq:70057}) into the new (\ref{eq:70042}).\qed
\end{proof}

\begin{remark}
The generating function for the Lame polynomial of type 2 in the algebraic form of the second kind about $x=a$ as $q= \alpha (\alpha +1)a- 4(2a-b-c)\left( q_j+2j+\frac{1}{2} \right)^2$ where $j,\alpha _j \in \mathbb{N}_{0}$ is
\begin{eqnarray}
&&\sum_{q_0 =0}^{\infty } \frac{(\frac{3}{2})_{q_0}}{q_0!} s_0^{q_0} \prod _{n=1}^{\infty } \left\{ \sum_{ q_n = q_{n-1}}^{\infty } s_n^{q_n }\right\}  LS_{q_j}^R\left( a, b, c, \alpha, q= \alpha (\alpha +1)a- 4(2a-b-c)\left( q_j+2j+\frac{1}{2} \right)^2 \right. \nonumber\\
&&; z= x-a, \rho = -\frac{2a-b-c}{(a-b)(a-c)} z, \left. \eta = \frac{-z^2}{(a-b)(a-c)} \right) \nonumber\\
&&= z^{\frac{1}{2}} \Bigg\{ \prod_{l=1}^{\infty } \frac{1}{(1-s_{l,\infty })} \mathbf{B}\left( s_{0,\infty } ;\rho \right)  + \Bigg\{ \prod_{l=2}^{\infty } \frac{1}{(1-s_{l,\infty })} \int_{0}^{1} dt_1\;t_1^{\frac{3}{2}} \int_{0}^{1} du_1\;u_1 \overleftrightarrow {\mathbf{\Psi}}_1 \left(s_{1,\infty };t_1,u_1,\rho \right)\nonumber\\
&& \times \widetilde{w}_{1,1}^{-\frac{1}{2}(2+\alpha )}\left( \widetilde{w}_{1,1} \partial _{ \widetilde{w}_{1,1}}\right) \widetilde{w}_{1,1}^{\alpha +\frac{1}{2}} \left( \widetilde{w}_{1,1} \partial _{ \widetilde{w}_{1,1}}\right) \widetilde{w}_{1,1}^{ \frac{1}{2}(1-\alpha )}\mathbf{B}\left( s_{0} ;\widetilde{w}_{1,1}\right)\Bigg\} \eta \nonumber\\
&&+ \sum_{n=2}^{\infty } \Bigg\{ \prod_{l=n+1}^{\infty } \frac{1}{(1-s_{l,\infty })} \int_{0}^{1} dt_n\;t_n^{2n-\frac{1}{2} } \int_{0}^{1} du_n\;u_n^{2n-1} \overleftrightarrow {\mathbf{\Psi}}_n \left(s_{n,\infty };t_n,u_n,\rho \right) \nonumber\\
&&\times \widetilde{w}_{n,n}^{-\frac{1}{2}(4n-2+\alpha )}\left( \widetilde{w}_{n,n} \partial _{ \widetilde{w}_{n,n}}\right) \widetilde{w}_{n,n}^{ \alpha +\frac{1}{2} } \left( \widetilde{w}_{n,n} \partial _{ \widetilde{w}_{n,n}}\right) \widetilde{w}_{n,n}^{ \frac{1}{2}(4n-3-\alpha ) }\nonumber\\
&&\times \prod_{k=1}^{n-1} \Bigg\{ \int_{0}^{1} dt_{n-k}\;t_{n-k}^{2(n-k)-\frac{1}{2}} \int_{0}^{1} du_{n-k} \;u_{n-k}^{2(n-k)-1} \overleftrightarrow {\mathbf{\Psi}}_{n-k} \left( s_{n-k};t_{n-k},u_{n-k},\widetilde{w}_{n-k+1,n} \right)\label{eq:70058}\\
&&\times \widetilde{w}_{n-k,n}^{-\frac{1}{2}(4(n-k)-2+\alpha )}\left( \widetilde{w}_{n-k,n} \partial _{ \widetilde{w}_{n-k,n}}\right) \widetilde{w}_{n-k,n}^{\alpha +\frac{1}{2}} \left( \widetilde{w}_{n-k,n} \partial _{ \widetilde{w}_{n-k,n}}\right) \widetilde{w}_{n-k,n}^{ \frac{1}{2}(4(n-k)-3-\alpha )} \Bigg\}\left. \mathbf{B}\left( s_{0} ;\widetilde{w}_{1,n}\right)\Bigg\} \eta ^n  \right\} \nonumber
\end{eqnarray}
where
\begin{equation}
\begin{cases} 
{ \displaystyle \overleftrightarrow {\mathbf{\Psi}}_1 \left(s_{1,\infty };t_1,u_1,\rho \right)= \frac{\left( \frac{1+s_{1,\infty }+\sqrt{s_{1,\infty }^2-2(1-2\rho (1-t_1)(1-u_1))s_{1,\infty }+1}}{2}\right)^{-4}}{\sqrt{s_{1,\infty }^2-2(1-2\rho (1-t_1)(1-u_1))s_{1,\infty }+1}}}\cr
{ \displaystyle  \overleftrightarrow {\mathbf{\Psi}}_n \left(s_{n,\infty };t_n,u_n,\rho \right) = \frac{\left( \frac{1+s_{n,\infty }+\sqrt{s_{n,\infty }^2-2(1-2\rho (1-t_n)(1-u_n))s_{n,\infty }+1}}{2}\right)^{-4n}}{\sqrt{s_{n,\infty }^2-2(1-2\rho (1-t_n)(1-u_n))s_{n,\infty }+1}}}\cr
{ \displaystyle \overleftrightarrow {\mathbf{\Psi}}_{n-k} \left(s_{n-k};t_{n-k},u_{n-k},\widetilde{w}_{n-k+1,n} \right) = \frac{\left( \frac{ 1+s_{n-k} +\sqrt{s_{n-k}^2-2(1-2\widetilde{w}_{n-k+1,n} (1-t_{n-k})(1-u_{n-k}))s_{n-k}+1}}{2}\right)^{-4(n-k)}}{\sqrt{s_{n-k}^2-2(1-2\widetilde{w}_{n-k+1,n} (1-t_{n-k})(1-u_{n-k}))s_{n-k}+1}}}
\end{cases}\nonumber 
\end{equation}
and
\begin{equation}
\begin{cases} 
{ \displaystyle \mathbf{B} \left(  s_{0,\infty } ;\rho \right)= \frac{\left(1- s_{0,\infty }+\sqrt{s_{0,\infty }^2-2(1-2\rho) s_{0,\infty }+1}\right)^{-\frac{1}{2}} \left(1+s_{0,\infty }+\sqrt{s_{0,\infty }^2-2(1-2\rho )s_{0,\infty }+1}\right)^{\frac{1}{2}}}{\sqrt{s_{0,\infty }^2-2(1-2\rho )s_{0,\infty }+1}}}\cr
{ \displaystyle  \mathbf{B} \left( s_{0} ;\widetilde{w}_{1,1}\right) = \frac{\left(1- s_0+\sqrt{s_0^2-2(1-2\widetilde{w}_{1,1})s_0+1}\right)^{-\frac{1}{2}} \left( 1+s_0+\sqrt{s_0^2-2(1-2\widetilde{w}_{1,1} )s_0+1}\right)^{\frac{1}{2}}}{\sqrt{s_0^2-2(1-2\widetilde{w}_{1,1})s_0+1}}} \cr
{ \displaystyle \mathbf{B} \left( s_{0} ;\widetilde{w}_{1,n}\right) = \frac{\left( 1-s_0+\sqrt{s_0^2-2(1-2\widetilde{w}_{1,n})s_0+1}\right)^{-\frac{1}{2}} \left( 1+s_0+\sqrt{s_0^2-2(1-2\widetilde{w}_{1,n} )s_0+1}\right)^{\frac{1}{2}}}{\sqrt{s_0^2-2(1-2\widetilde{w}_{1,n})s_0+1}}}
\end{cases}\nonumber 
\end{equation}
\end{remark}
\begin{proof}
Replace $w$, $\gamma $, $A$ and $x$ by $s_{0,\infty }$, $\frac{3}{2}$, 1 and $\rho $ in (\ref{eq:70035}). 
\begin{eqnarray}
&&\sum_{q_0=0}^{\infty }\frac{(\frac{3}{2})_{q_0}}{q_0!} s_{0,\infty }^{q_0} \;_2F_1\left( -q_0, q_0+1; \frac{3}{2}; \rho \right) \label{eq:70059}\\
&&= \frac{\left(1- s_{0,\infty }+\sqrt{s_{0,\infty }^2-2(1-2\rho )s_{0,\infty }+1}\right)^{-\frac{1}{2}} \left(1+s_{0,\infty }+\sqrt{s_{0,\infty }^2-2(1-2\rho )s_{0,\infty }+1}\right)^{\frac{1}{2}}}{\sqrt{s_{0,\infty }^2-2(1-2\rho )s_{0,\infty }+1}} \nonumber
\end{eqnarray} 
Replace $s_{0,\infty }$ and $\rho $  by  $s_0$ and $\widetilde{w}_{1,1}$ in (\ref{eq:70059}). 
\begin{eqnarray}
&&\sum_{q_0=0}^{\infty }\frac{(\frac{3}{2})_{q_0}}{q_0!} s_0^{q_0} \;_2F_1\left( -q_0, q_0+1; \frac{3}{2}; \widetilde{w}_{1,1} \right) \label{eq:70060}\\
&&= \frac{\left(1- s_0+\sqrt{s_0^2-2(1-2\widetilde{w}_{1,1})s_0+1}\right)^{-\frac{1}{2}} \left( 1+s_0+\sqrt{s_0^2-2(1-2\widetilde{w}_{1,1} )s_0+1}\right)^{\frac{1}{2}}}{\sqrt{s_0^2-2(1-2\widetilde{w}_{1,1})s_0+1}} \nonumber
\end{eqnarray} 
Replace  $\widetilde{w}_{1,1}$  by $\widetilde{w}_{1,n}$  in (\ref{eq:70060}). 
\begin{eqnarray}
&&\sum_{q_0=0}^{\infty }\frac{(\frac{3}{2})_{q_0}}{q_0!} s_0^{q_0} \;_2F_1\left( -q_0, q_0+1; \frac{3}{2}; \widetilde{w}_{1,n} \right) \label{eq:70061}\\
&&= \frac{\left(1- s_0+\sqrt{s_0^2-2(1-2\widetilde{w}_{1,n})s_0+1}\right)^{-\frac{1}{2}} \left( 1+s_0+\sqrt{s_0^2-2(1-2\widetilde{w}_{1,n} )s_0+1}\right)^{\frac{1}{2}}}{\sqrt{s_0^2-2(1-2\widetilde{w}_{1,n})s_0+1}} \nonumber
\end{eqnarray} 
Put $ c_0= 1 $, $\lambda =\frac{1}{2} $ and $\gamma' =\frac{3}{2} $ in (\ref{eq:70042}). Substitute (\ref{eq:70059}), (\ref{eq:70060}) and (\ref{eq:70061}) into the new (\ref{eq:70042}).\qed
\end{proof}
\section{Summary}

Lame equation represented either in the algebraic form or in Weierstrass's form has a recursive relation between a 3-term in the power series of its equation as we see (\ref{eq:7006})--(\ref{eq:7007c}). 
Lame equation could not be described in the form of a definite or contour integral of any elementary functions because of its three term recurrence relation in its power series. The 3-term recursive relation between successive coefficients creates complex mathematical calculations in order to obtain analytic solutions of it in closed forms.

In Ref.\cite{zChou2012f} I show how to obtain the power series expansion in closed forms and its integral form of Lame equation in the algebraic form (for infinite series and polynomial of type 1 including all higher terms of $A_n$'s) by applying 3TRF. This was done by letting $A_n$ in sequence $c_n$ is the leading term in the analytic function $y(z)$: the sequence $c_n$ consists of combinations $A_n$ and $B_n$. For polynomial of type 1, I treat $q$ as a free variable and a fixed value of $\alpha $.

In this chapter I show how to construct the power series expansion and its integral form of Lame equation in the algebraic form for infinite series and polynomial of type 2 including all higher terms of $B_n$'s by applying R3TRF. This is done by letting $B_n$ in sequence $c_n$ is the leading term in a analytic function $y(z)$. For polynomial of type 2, I treat $\alpha $ as a free variable and a fixed value of $q$. 

The Frobenius solution and its integral form of Lame equation (ellipsoidal harmonics equation) for infinite series about $x=a$ in this chapter are equivalent to infinite series of Lame equation in Ref.\cite{zChou2012f}. In this chapter $B_n$ is the leading term in sequence $c_n$ in an analytic function $y(z)$. In Ref.\cite{zChou2012f} $A_n$ is the leading term in sequence $c_n$ in the analytic function $y(z)$.
 
In Ref.\cite{zChou2012f} and this chapter, as we see the power series expansions of Lame equation about $x=a$ for either polynomial or infinite series, the denominators and numerators in all $A_n$ or $B_n$ terms of $y(z)$ arise with Pochhammer symbol. Since we construct the power series expansions with Pochhammer symbols in numerators and denominators, we are able to describe integral representations of Lame equation analytically. As we observe representations in closed form integrals of Lame equation by applying either 3TRF or R3TRF, a $_2F_1$ function recurs in each of sub-integral forms of the $y(z)$ (each sub-integral $y_m(z)$ where $m=0,1,2,\cdots$ is composed of $2m$ terms of definite integrals and $m$ terms of contour integrals).
We are able to transform the Lame function into any well-known special functions having two recurrence relation in its power series of any linear differential equations because of a $_2F_1$ function in each of sub-integral forms of Lame function. After we replace $_2F_1$ function  in its integral forms to other special functions (such as  Bessel function, Kummer function, hypergeometric function, Laguerre function and etc), we are able to rebuild the Frobenius solution of Lame equation in a backward.    
 
In this chapter I show how to construct the generating function for the type 2 Lame polynomial from its integral representation.  We are able to derive orthogonal relations, recursion relations and expectation values of physical quantities from the generating function: the processes in order to obtain orthogonal and recursion relations of the Lame polynomials are similar as the case of a normalized wave function for the hydrogen-like atoms.  

Mathematical structure of the power series expansions, integral representations and generating functions for Lame equation closely resembles the case of Heun equation using 3TRF and R3TRF. (1) If $\gamma ,\delta ,\epsilon \rightarrow \frac{1}{2}$, $1 \rightarrow b-a$, $a \rightarrow c-a$, $\alpha \rightarrow \frac{1}{2}(\alpha +1)$, $\beta \rightarrow -\frac{1}{2}\alpha$, $q \rightarrow  -q + \frac{1}{4}\alpha (\alpha +1)a$ and $x \rightarrow z$ in the power series expansions and its integral forms of Heun equation using 3TRF in Ref.\cite{zchou2012c,zChou2012d}, its analytic solutions are correspondent to the Frobenius solutions and integral representations of Lame equation in the algebraic form in Ref.\cite{zChou2012f}; compare (10), (11), (18) and (19) in Ref.\cite{zchou2012c} with remarks 1--4 in Ref.\cite{zChou2012f}. And make a comparison between remarks 1, 2, 5 and 6 in Ref.\cite{zChou2012d} and remarks 5--8 in Ref.\cite{zChou2012f}. 
(2) By similar reason, putting (\ref{eq:7004}) into the power series, integral representations and generating functions of Heun equation using R3TRF in chapters 2 and 3, its analytic solutions are equivalent to the Frobenius solutions, integral forms and generating functions of Lame equation using R3TRF in this chapter; compare remarks 2.2.1--2.2.4, 2.3.2--2.3.4 in chapter 2 and remarks 3.2.5--3.2.6 in chapter 3 with remarks 8.2.1--8.2.4, 8.3.2--8.3.5, 8.4.4-8.4.5 in this chapter.

\addcontentsline{toc}{section}{Bibliography}
\bibliographystyle{model1a-num-names}
\bibliography{<your-bib-database>}

\chapter{Lame function in Weierstrass's form using reversible three-term recurrence formula}
\chaptermark{Lame function in Weierstrass's form using R3TRF} 
Lame ordinary differential equation in Weierstrass's form and Heun equation are of Fuchsian types with the four regular singularities. Lame equation in Weierstrass's form is derived from Heun equation by changing all coefficients $\gamma =\delta =\epsilon =\frac{1}{2}$, $a =\rho ^{-2}$, $\alpha = \frac{1}{2}(\alpha +1) $, $\beta = -\frac{1}{2}\alpha  $, $q=-\frac{1}{4}h\rho ^{-2}$ and an independent variable $x=sn^2(z,\rho)$. \cite{aHeun1889,aRonv1995} 

In Ref.\cite{aChou2012g,aChou2012h} I construct the power series expansion in closed forms, its integral form and the  generating function of Lame equation in Weierstrass's form by applying three term recurrence formula (3TRF)\cite{achou2012b};  analytic solutions of Lame equation is derived for infinite series and polynomial of type 1\footnote{polynomial of type 1 is a polynomial which makes $B_n$ term terminated in three term recursion relation of the power series in a linear differential equation.} including all higher terms of $A_n$'s.\footnote{`` higher terms of $A_n$'s'' means at least two terms of $A_n$'s.}

In this chapter I will apply reversible three term recurrence formula (R3TRF) in chapter 1 to (1) the Frobenius solution in closed forms, (2) its integral form of Lame equation in Weierstrass's form for infinite series and polynomial of type 2\footnote{polynomial of type 2 is a polynomial which makes $A_n$ term terminated in three term recursion relation of the power series in a linear differential equation.} including all higher terms of $B_n$'s\footnote{`` higher terms of $B_n$'s'' means at least two terms of $B_n$'s.}, (3) the generating function for the Lame polynomial of type 2. 

Nine examples of 192 local solutions of the Heun equation (Maier, 2007) are provided in the appendix.  For each example, I show how to convert local solutions of Heun equation by applying R3TRF to analytic solutions of Lame equation in Weierstrass's form.
 
\section{Introduction}
In general, there are two types of Lame equation which are described in the algebraic form and in Weierstrass's form. For computational practice, Jacobian (Weierstrass's) form is more convenient. For the general mathematical properties of analytic solutions of Lame equation, the algebraic form is better; especially for the asymptotic expansion of Lame equation, the algebraic form is much better and its asymptotic solutions are more various because of three different singularity parameters. In contrast, the Jacobian form only has two singularity parameters which are 1 and the modulus of the elliptic function $sn\; z$.  Various authors have worked on the boundary value problems in ellipsoidal geometry with Lame equation in Weierstrass's form rather than the algebraic form of it. 

Many great mathematicians/physicists attempted to describe Lame equation in a definite or contour integral form of any well-known simple functions such as Gauss hypergeometric, Kummer functions, and etc. Due to its mathematical complexity there are no analytic solutions in closed forms of the Lame function\cite{aErde1955,aHobs1931,aWhit1952}. Because Lame equation represented either in the algebraic form or in Weierstrass's form is a form of the power series that is expressed as three term recurrence relation\cite{aHobs1931,aWhit1952}. In contrast, most of well-known special functions consist of two term recursion relation (Hypergeometric, Bessel, Legendre, Kummer functions, etc). 
 They just left analytic solutions of Lame equation as solutions of recurrences because of a 3-term recursive relation between successive coefficients in its power series expansion. Three or more terms recursion relation in the power series of any linear ordinary differential equations creates the complex mathematical calculation. 

In Ref.\cite{aChou2012f} I construct the power series expansion in closed forms and its integral representation of Lame equation in the algebraic form by applying 3TRF \cite{achou2012b}:  analytic solutions of Lame equation are derived for infinite series and polynomial of type 1 including all higher terms of $A_n$'s. 
One interesting observation resulting from calculations is the fact that a $_2F_1$ function recurs in each of sub-integral forms: the first sub-integral form contains zero term of $A_n's$, the second one contains one term of $A_n$'s, the third one contains two terms of $A_n$'s, etc. Also asymptotic expansions of the Lame function for infinite series are derived analytically including for special cases as  $|2a-b-c|\gg 1$ and $2a-b-c\approx  0$.

In Ref.\cite{aChou2012g}, by changing all coefficients in the general expression of the power series and its integral representation of Lame equation in the algebraic form, I obtain the Frobenius solution and its integral form of Lame equation in Weierstrass's form for infinite series and polynomial of type 1. 
In Ref.\cite{aChou2012h} I construct the generating function for the Lame polynomial of type 1 in Weierstrass's form from the general expression of an integral representation of Lame equation in Weierstrass's form by applying the generating function for the Jacobi polynomial using hypergeometric functions.

In chapter 8, by applying R3TRF, I derive (1) the power series expansion, (2) its integral form of Lame equation in the algebraic form for infinite series and polynomial of type 2 including all higher terms of $B_n$'s and (3) the generating function for the Lame polynomial of type 2. 
  
In this chapter, by changing all coefficients in the general expression of power series, its integral representation and the generating function of Heun equation in chapter 2 and 3, I consider (1) the power series expansion in closed forms of Lame equation in Weierstrass's form (for infinite series and polynomial of type 2 including all higher terms of $B_n$'s), (2) its integral representation and (3) the generating function for the Lame polynomial of type 2 analytically.\footnote{Or we can obtain analytic solutions of Lame equation in Weierstrass's form by applying R3TRF directly.} 
Especially we are able to obtain a mathematical construction for an orthogonal relation of Lame equation in Weierstrass's form for polynomial of type 1 or 2 from the generating function for the Lame polynomial. After that, we might be possible to build recurrence relations of the Lame polynomial, a normalized constant of any wave functions for Laplace equation in ellipsoidal coordinates, an expectation value of its physical quantities and etc. 
For an example, a normalized wave function of hydrogen-like atoms is derived from the generating function for associated Laguerre polynomials in explicit mathematical calculations. Also the expectation value of physical quantities such as position and momentum is constructed by applying the recursive relation of associated Laguerre polynomials.

I show how to transform nine examples of 192 local solutions of Heun equation \cite{aMaie2007} to analytic solutions of Lame equation in Weierstrass's form in the appendix. For each example, I convert local solutions of Heun equation by applying R3TRF to analytic solutions of Lame equation in Weierstrass's form; the Frobenius solution in closed form of Lame equation, its integral representation and the generating function of it are derived from changing all coefficients and an independent variable in Heun equation.

Lame equation in Weierstrass's form is given by
\begin{equation}
\frac{d^2{y}}{d{z}^2} = \{ \alpha (\alpha +1)\rho^2\;sn^2(z,\rho )-h\} y(z)\label{eq:8001}
\end{equation}
where $\rho$, $\alpha $  and h are real parameters such that $0<\rho <1$ and $\alpha \geq -\frac{1}{2}$.
If we take $sn^2(z,\rho)=\xi $ as an independent variable, Lame equation becomes
\begin{equation}
\frac{d^2{y}}{d{\xi }^2} + \frac{1}{2}\left(\frac{1}{\xi } +\frac{1}{\xi -1} + \frac{1}{\xi -\rho ^{-2}}\right) \frac{d{y}}{d{\xi }} +  \frac{-\alpha (\alpha +1) \xi +h\rho ^{-2}}{4 \xi (\xi -1)(\xi -\rho ^{-2})} y(\xi ) = 0\label{eq:8002}
\end{equation}
This is an equation of Fuchsian type with the four regular singularities: $\xi=0, 1, \rho ^{-2}, \infty $. The first three, namely $0, 1, \rho ^{-2}$, have the property that the corresponding exponents are $0, \frac{1}{2}$ which is the same as the case of Lame equation in the algebraic form.

Lame equation is a second-order linear ordinary differential equation of the algebraic form\cite{aLame1837,aStie1885,aZwil1997}
\begin{equation}
\frac{d^2{y}}{d{x}^2} + \frac{1}{2}\left(\frac{1}{x-a} +\frac{1}{x-b} + \frac{1}{x-c}\right) \frac{d{y}}{d{x}} +  \frac{-\alpha (\alpha +1) x+q}{4 (x-a)(x-b)(x-c)} y = 0\label{eq:8003}
\end{equation}
Lame equation has four regular singular points: a, b, c and $\infty $; the exponents at the first three are all $0$ and $\frac{1}{2}$, and those at infinity are $-\frac{1}{2}\alpha $ and $\frac{1}{2}(\alpha +1)$.\cite{aMoon1961,aBoch1894}

As we compare (\ref{eq:8002}) with (\ref{eq:8003}), all coefficients on the above are correspondent to the following way.
\begin{equation}
\begin{split}
& a \longrightarrow   0 \\ & b \longrightarrow  1 \\ & c \longrightarrow  \rho ^{-2} \\
& q \longrightarrow  h \rho ^{-2} \\ & x \longrightarrow \xi = sn^2(z,\rho ) 
\end{split}\label{eq:8004}   
\end{equation}
Also (\ref{eq:8002}) is a special case of Heun's equation. Heun equation is a second-order linear ordinary differential equation of the form \cite{aHeun1889}.
\begin{equation}
\frac{d^2{y}}{d{x}^2} + \left(\frac{\gamma }{x} +\frac{\delta }{x-1} + \frac{\epsilon }{x-a}\right) \frac{d{y}}{d{x}} +  \frac{\alpha \beta x-q}{x(x-1)(x-a)} y = 0 \label{eq:8005}
\end{equation}
With the condition $\epsilon = \alpha +\beta -\gamma -\delta +1$. The parameters play different roles: $a \ne 0 $ is the singularity parameter, $\alpha $, $\beta $, $\gamma $, $\delta $, $\epsilon $ are exponent parameters, $q$ is the accessory parameter which in many physical applications appears as a spectral parameter. Also, $\alpha $ and $\beta $ are identical to each other. The total number of free parameters is six. It has four regular singular points which are 0, 1, $a$ and $\infty $ with exponents $\{ 0, 1-\gamma \}$, $\{ 0, 1-\delta \}$, $\{ 0, 1-\epsilon \}$ and $\{ \alpha, \beta \}$.

As we compare (\ref{eq:8002}) with (\ref{eq:8005}), all coefficients on the above are correspondent to the following way.
\begin{equation}
\begin{split}
& \gamma ,\delta ,\epsilon  \longleftrightarrow   \frac{1}{2} \\ & a\longleftrightarrow  \rho ^{-2} \\ & \alpha  \longleftrightarrow \frac{1}{2}(\alpha +1) \\
& \beta   \longleftrightarrow -\frac{1}{2} \alpha \\
& q \longleftrightarrow  -\frac{1}{4}h \rho ^{-2} \\ & x \longleftrightarrow \xi = sn^2(z,\rho ) 
\end{split}\label{eq:8006}   
\end{equation}
\section{Power series}
\subsection{Polynomial of type 2}
In general Lame (spectral) polynomial represented either in the algebraic form or in Weierstrass's form is defined as polynomial of type 3.\footnote{polynomial of type 3 is a polynomial which makes $A_n$ and $B_n$ terms terminated at the same time in three term recursion relation of the power series in a linear differential equation. If $A_n$ and $B_n$ terms are not terminated, it turns to be infinite series.} Lame polynomial comes from a Lame equation that has a fixed integer value of $\alpha $, just as it has a fixed value of $h$. In three-term recurrence formula, polynomial of type 3 I categorize as complete polynomial. In future papers I will derive the type 3 Lame polynomial. In Ref.\cite{aChou2012g,aChou2012h} I derive the Frobenius solution, its integral form and the generating function for the Lame polynomial of type 1: I treat $h$ as a free variable  and $\alpha $ as a fixed value. In this chapter I construct the power series expansion and its integral representation of the Lame polynomial of type 2:  I treat $\alpha $ as a free variable and $h$ as a fixed value. Also the generating function for the lame polynomial of type 2 is constructed analytically.   

In chapter 2 the general expression of power series of Heun equation about $x=0$ for polynomial of type 2 is given by 
\begin{eqnarray}
 y(x)&=&  \sum_{n=0}^{\infty } y_{n}(x) = y_0(x)+ y_1(x)+ y_2(x)+y_3(x)+\cdots \nonumber\\ 
&=& c_0 x^{\lambda } \left\{\sum_{i_0=0}^{q_0} \frac{(-q_0)_{i_0} \left(q_0+ \frac{\varphi +2(1+a)\lambda }{(1+a)}\right)_{i_0}}{(1+\lambda )_{i_0}(\gamma +\lambda )_{i_0}} \eta ^{i_0}\right.\nonumber\\
&&+ \left\{ \sum_{i_0=0}^{q_0}\frac{(i_0+ \lambda +\alpha ) (i_0+ \lambda +\beta )}{(i_0+ \lambda +2)(i_0+ \lambda +1+\gamma )}\frac{(-q_0)_{i_0} \left(q_0+ \frac{\varphi +2(1+a)\lambda }{(1+a)}\right)_{i_0}}{(1+\lambda )_{i_0}(\gamma +\lambda )_{i_0}} \right.\nonumber\\
&&\times  \left. \sum_{i_1=i_0}^{q_1} \frac{(-q_1)_{i_1}\left(q_1+4+ \frac{\varphi +2(1+a)\lambda }{(1+a)}\right)_{i_1}(3+\lambda )_{i_0}(2+\gamma +\lambda )_{i_0}}{(-q_1)_{i_0}\left(q_1+4+ \frac{\varphi +2(1+a)\lambda }{(1+a)}\right)_{i_0}(3+\lambda )_{i_1}(2+\gamma +\lambda )_{i_1}} \eta ^{i_1}\right\} z\nonumber\\
&&+ \sum_{n=2}^{\infty } \left\{ \sum_{i_0=0}^{q_0} \frac{(i_0+ \lambda +\alpha ) (i_0+ \lambda +\beta )}{(i_0+ \lambda +2)(i_0+ \lambda +1+\gamma )}\frac{(-q_0)_{i_0} \left(q_0+ \frac{\varphi +2(1+a)\lambda }{(1+a)}\right)_{i_0}}{(1+\lambda )_{i_0}(\gamma +\lambda )_{i_0}}\right.\nonumber\\
&&\times \prod _{k=1}^{n-1} \left\{ \sum_{i_k=i_{k-1}}^{q_k} \frac{(i_k+ 2k+\lambda +\alpha ) (i_k+ 2k+\lambda +\beta )}{(i_k+ 2(k+1)+\lambda )(i_k+ 2k+1+\gamma +\lambda )}\right. \nonumber\\
&&\times \left.\frac{(-q_k)_{i_k}\left(q_k+4k+ \frac{\varphi +2(1+a)\lambda }{(1+a)}\right)_{i_k}(2k+1+\lambda )_{i_{k-1}}(2k+\gamma +\lambda )_{i_{k-1}}}{(-q_k)_{i_{k-1}}\left(q_k+4k+ \frac{\varphi +2(1+a)\lambda }{(1+a)}\right)_{i_{k-1}}(2k+1+\lambda )_{i_k}(2k+\gamma +\lambda )_{i_k}}\right\} \label{eq:8007}\\
&&\times \left. \left.\sum_{i_n= i_{n-1}}^{q_n} \frac{(-q_n)_{i_n}\left(q_n+4n+ \frac{\varphi +2(1+a)\lambda }{(1+a)}\right)_{i_n}(2n+1+\lambda )_{i_{n-1}}(2n+\gamma +\lambda )_{i_{n-1}}}{(-q_n)_{i_{n-1}}\left(q_n+4n+ \frac{\varphi +2(1+a)\lambda }{(1+a)}\right)_{i_{n-1}}(2n+1+\lambda )_{i_n}(2n+\gamma +\lambda )_{i_n}} \eta ^{i_n} \right\} z^n \right\}  \nonumber
\end{eqnarray}
where
\begin{equation}
\begin{cases} z = -\frac{1}{a}x^2 \cr
\eta = \frac{(1+a)}{a} x \cr
\varphi = \alpha +\beta -\delta +a(\delta +\gamma -1) \cr
q= -(q_j+2j+\lambda )\{\varphi +(1+a)(q_j+2j+\lambda ) \} \;\;\mbox{as}\;j,q_j\in \mathbb{N}_{0} \cr
q_i\leq q_j \;\;\mbox{only}\;\mbox{if}\;i\leq j\;\;\mbox{where}\;i,j\in \mathbb{N}_{0} 
\end{cases}\nonumber 
\end{equation}
In the above, We have two indicial roots which are $\lambda = 0$ and $ 1-\gamma $.

Put (\ref{eq:8006}) in (\ref{eq:8007}) with replacing $q_i$ by $h_i$ where $h_i\in \mathbb{N}_{0}$. Take $c_0$= 1 as $\lambda =0$  for the first independent solution of Lame equation and $\lambda =\frac{1}{2}$ for the second one into the new (\ref{eq:8007}).
\begin{remark}
The power series expansion of Lame equation in Weierstrass's form of the first kind for polynomial of type 2  about $\xi =0$ as $h= 4(1+\rho ^2)(h_j+2j )^2 $ where $j,h_j \in \mathbb{N}_{0}$ is
\begin{eqnarray}
y(\xi )&=& LF_{h_j}^R\left( \rho , \alpha, h= 4(1+\rho ^2)(h_j+2j )^2; \xi = sn^2(z,\rho ), \eta = (1+\rho ^2)\xi, z=-\rho^2 \xi^2 \right) \nonumber\\
&=&  \sum_{i_0=0}^{h_0} \frac{(-h_0)_{i_0} \left( h_0  \right)_{i_0}}{\left( 1 \right)_{i_0}\left(\frac{1}{2} \right)_{i_0}} \eta ^{i_0} \nonumber\\
&&+ \left\{ \sum_{i_0=0}^{h_0}\frac{\left( i_0+\frac{\alpha}{2}+\frac{1}{2} \right)\left( i_0-\frac{\alpha}{2} \right)}{ \left(i_0+ 2 \right) \left( i_0+ \frac{3}{2} \right)}\frac{(-h_0)_{i_0} \left( h_0 \right)_{i_0}}{\left( 1 \right)_{i_0}\left(\frac{1}{2} \right)_{i_0}} \right. \left. \sum_{i_1=i_0}^{h_1} \frac{(-h_1)_{i_1}\left( h_1 +4 \right)_{i_1}\left( 3 \right)_{i_0}\left( \frac{5}{2} \right)_{i_0}}{(-h_1)_{i_0}\left( h_1 +4 \right)_{i_0}\left( 3 \right)_{i_1}\left( \frac{5}{2} \right)_{i_1}} \eta ^{i_1}\right\} z \nonumber\\
&&+ \sum_{n=2}^{\infty } \left\{ \sum_{i_0=0}^{h_0}\frac{\left( i_0+\frac{\alpha}{2}+\frac{1}{2} \right)\left( i_0-\frac{\alpha}{2} \right)}{ \left(i_0+ 2 \right) \left( i_0+ \frac{3}{2} \right)}\frac{(-h_0)_{i_0} \left( h_0 \right)_{i_0}}{\left( 1 \right)_{i_0}\left(\frac{1}{2} \right)_{i_0}} \right.\nonumber\\
&&\times \prod _{k=1}^{n-1} \left\{ \sum_{i_k=i_{k-1}}^{h_k} \frac{\left( i_k+2k+\frac{\alpha}{2}+\frac{1}{2} \right)\left( i_k+2k-\frac{\alpha}{2} \right)}{ \left( i_k+ 2k+2 \right) \left( i_k+2k +\frac{3}{2} \right)}\right.  \left.\frac{(-h_k)_{i_k}\left( h_k+4k \right)_{i_k}\left( 2k+1 \right)_{i_{k-1}}\left( 2k+\frac{1}{2} \right)_{i_{k-1}}}{(-h_k)_{i_{k-1}}\left( h_k+4k \right)_{i_{k-1}}\left( 2k+1 \right)_{i_k}\left( 2k+\frac{1}{2} \right)_{i_k}}\right\} \nonumber\\
&&\times  \left.\sum_{i_n= i_{n-1}}^{h_n} \frac{(-h_n)_{i_n}\left( h_n+4n \right)_{i_n}\left( 2n +1 \right)_{i_{n-1}}\left( 2n+\frac{1}{2} \right)_{i_{n-1}}}{(-h_n)_{i_{n-1}}\left( h_n+4n  \right)_{i_{n-1}}\left( 2n +1 \right)_{i_n}\left( 2n +\frac{1}{2} \right)_{i_n}} \eta ^{i_n} \right\} z ^n  \label{eq:8008}
\end{eqnarray}
\end{remark}
For the minimum value of Lame equation in Weierstrass's form of the first kind for polynomial of type 2 about $\xi =0$, put $h_0=h_1=h_2=\cdots=0$ in (\ref{eq:8008}).
\begin{eqnarray}
y(\xi )&=& LF_{0}^R\left( \rho , \alpha, h= 16(1+\rho ^2)j^2; \xi = sn^2(z,\rho ), \eta = (1+\rho ^2)\xi, z=-\rho^2 \xi^2 \right) \nonumber\\ 
&=& \; _2F_1\left( -\frac{\alpha }{4},\frac{\alpha }{4}+\frac{1}{4},\frac{3}{4},z \right) \hspace{1cm}\mbox{where}\;\;|z|< 1 \label{ccc:8001}
\end{eqnarray} 
\begin{remark}
The power series expansion of Lame equation in Weierstrass's form of the second kind for polynomial of type 2  about $\xi =0$ as $h= 4(1+\rho ^2)\left(h_j+2j +\frac{1}{2}\right)^2 $ where $j,h_j \in \mathbb{N}_{0}$ is
\begin{eqnarray}
y(\xi ) &=& LS_{h_j}^R\left(\rho , \alpha, h= 4(1+\rho ^2)\left( h_j+2j +\frac{1}{2}\right)^2; \xi = sn^2(z,\rho ), \eta = (1+\rho ^2)\xi, z=-\rho^2 \xi^2 \right) \nonumber\\
&=& \xi^{\frac{1}{2}} \left\{\sum_{i_0=0}^{h_0} \frac{(-h_0)_{i_0} \left( h_0 +1 \right)_{i_0}}{\left(\frac{3}{2} \right)_{i_0}\left(1 \right)_{i_0}} \eta ^{i_0}\right.\nonumber\\
&&+ \left\{ \sum_{i_0=0}^{h_0}\frac{\left( i_0+\frac{\alpha}{2}+1 \right)\left( i_0-\frac{\alpha}{2} + \frac{1}{2} \right)}{ \left(i_0+ \frac{5}{2} \right) \left( i_0+ 2 \right)}\frac{(-h_0)_{i_0} \left( h_0+ 1 \right)_{i_0}}{\left(\frac{3}{2} \right)_{i_0}\left( 1\right)_{i_0}} \right.   \left. \sum_{i_1=i_0}^{h_1} \frac{(-h_1)_{i_1}\left( h_1 +5 \right)_{i_1}\left( \frac{7}{2} \right)_{i_0}\left( 3 \right)_{i_0}}{(-h_1)_{i_0}\left( h_1 +5 \right)_{i_0}\left( \frac{7}{2} \right)_{i_1}\left( 3 \right)_{i_1}} \eta ^{i_1}\right\} z \nonumber\\
&&+ \sum_{n=2}^{\infty } \left\{ \sum_{i_0=0}^{h_0}\frac{\left( i_0+\frac{\alpha}{2}+1 \right)\left( i_0-\frac{\alpha}{2} +\frac{1}{2} \right)}{ \left(i_0+ \frac{5}{2} \right) \left( i_0+ 2 \right)}\frac{(-h_0)_{i_0} \left( h_0+ 1 \right)_{i_0}}{\left(\frac{3}{2} \right)_{i_0}\left(1\right)_{i_0}} \right.\nonumber\\
&&\times \prod _{k=1}^{n-1} \left\{ \sum_{i_k=i_{k-1}}^{h_k} \frac{\left( i_k+2k+\frac{\alpha}{2}+1\right)\left( i_k+2k-\frac{\alpha}{2} + \frac{1}{2} \right)}{ \left( i_k+ 2k+\frac{5}{2} \right) \left( i_k+2k +2 \right)}\right. \nonumber\\
&&\times  \left.\frac{(-h_k)_{i_k}\left( h_k+4k + 1 \right)_{i_k}\left( 2k+\frac{3}{2} \right)_{i_{k-1}}\left( 2k+1 \right)_{i_{k-1}}}{(-h_k)_{i_{k-1}}\left( h_k+4k +1 \right)_{i_{k-1}}\left( 2k+\frac{3}{2} \right)_{i_k}\left( 2k+1\right)_{i_k}}\right\} \nonumber\\
&&\times \left. \left.\sum_{i_n= i_{n-1}}^{h_n} \frac{(-h_n)_{i_n}\left( h_n+4n + 1 \right)_{i_n}\left( 2n +\frac{3}{2} \right)_{i_{n-1}}\left( 2n+1 \right)_{i_{n-1}}}{(-h_n)_{i_{n-1}}\left( h_n+4n +1 \right)_{i_{n-1}}\left( 2n +\frac{3}{2} \right)_{i_n}\left( 2n +1 \right)_{i_n}} \eta ^{i_n} \right\} z^n \right\} \label{eq:8009}
\end{eqnarray}
\end{remark}
For the minimum value of Lame equation in Weierstrass's form of the second kind for polynomial of type 2  about $\xi =0$, put $h_0=h_1=h_2=\cdots=0$ in (\ref{eq:8009}).
\begin{eqnarray}
y(\xi ) &=& LS_{0}^R\left(\rho , \alpha, h= 4(1+\rho ^2)\left( 2j +\frac{1}{2}\right)^2; \xi = sn^2(z,\rho ), \eta = (1+\rho ^2)\xi, z=-\rho^2 \xi^2 \right)\nonumber\\
&=&\xi^{\frac{1}{2}}\;  _2F_1\left( -\frac{\alpha }{4}+\frac{1}{4},\frac{\alpha }{4}+\frac{1}{2},\frac{5}{4},z \right) \hspace{1cm}\mbox{where}\;\;|z|<1 \label{ccc:8003}
\end{eqnarray} 
In (\ref{ccc:8001}) and (\ref{ccc:8003}), a polynomial of type 2 requires $\left| z=-\rho^2 \xi^2\right|<1$ for the convergence of the radius. For more details about this issue, it is explained in chapter 3  and 4 \cite{aChoun2014}. 

In Ref.\cite{aChou2012g,aChou2012h} I treat $\alpha $ as a fixed value and  $h$ as a free variable to construct Lame polynomial of type 1 in Weierstrass's: (1) if $\alpha = 2(2\alpha_j +j)$ or  $-2(2\alpha_j +j)-1$ where $j, \alpha_j \in \mathbb{N}_{0}$, an analytic solution of Lame equation turns to be the first kind of independent solutions of the Lame polynomial of type 1. (2) if  $\alpha = 2(2\alpha_j +j)+1$ or $-2(2\alpha_j +j+1)$, an analytic solution of Lame equation turns to be the second kind of independent solutions of the Lame polynomial of type 1.\footnote{If we take $\alpha \geq -\frac{1}{2}$, $\alpha =  -2(2\alpha_j +j)-1 $ and $-2(2\alpha_j +j+1)$  are not available any more. I consider $\alpha $ as arbitrary.} 

In this chapter I treat $h$ as a fixed value and $\alpha$ as a free variable to construct the Lame polynomial of type 2 in Weierstrass's form : (1) if $h= 4(1+\rho ^2)(h_j+2j )^2 $ where $j,h_j \in \mathbb{N}_{0}$, an analytic solution of Lame equation turns to be the first kind of independent solutions of the Lame polynomial of type 2. (2) if $h= 4(1+\rho ^2)\left( h_j+2j +\frac{1}{2}\right)^2$, an analytic solution of the Lame equation turns to be the second kind of independent solutions of the Lame polynomial of type 2.
\subsection{Infinite series}
In chapter 2 the general expression of power series of Heun equation for infinite series about $x=0$  by applying R3TRF is given by
\begin{eqnarray}
 y(x)&=& \sum_{n=0}^{\infty } y_n(x)= y_0(x)+ y_1(x)+ y_2(x)+ y_3(x)+\cdots \nonumber\\
&=& c_0 x^{\lambda } \left\{\sum_{i_0=0}^{\infty } \frac{\left(\Delta_0^{-}\right)_{i_0} \left(\Delta_0^{+}\right)_{i_0}}{(1+\lambda )_{i_0}(\gamma +\lambda )_{i_0}} \eta ^{i_0}\right.\nonumber\\
&+& \left\{ \sum_{i_0=0}^{\infty }\frac{(i_0+ \lambda +\alpha ) (i_0+ \lambda +\beta )}{(i_0+ \lambda +2)(i_0+ \lambda +1+\gamma )}\frac{\left(\Delta_0^{-}\right)_{i_0} \left(\Delta_0^{+}\right)_{i_0}}{(1+\lambda )_{i_0}(\gamma +\lambda )_{i_0}} \right.\nonumber\\
&&\times \left. \sum_{i_1=i_0}^{\infty } \frac{\left(\Delta_1^{-}\right)_{i_1} \left(\Delta_1^{+}\right)_{i_1}(3+\lambda )_{i_0}(2+\gamma +\lambda )_{i_0}}{\left(\Delta_1^{-}\right)_{i_0}  \left(\Delta_1^{+}\right)_{i_0}(3+\lambda )_{i_1}(2+\gamma +\lambda )_{i_1}}\eta ^{i_1}\right\} z\nonumber\\
&+& \sum_{n=2}^{\infty } \left\{ \sum_{i_0=0}^{\infty } \frac{(i_0+ \lambda +\alpha ) (i_0+ \lambda +\beta )}{(i_0+ \lambda +2)(i_0+ \lambda +1+\gamma )}\frac{\left(\Delta_0^{-}\right)_{i_0} \left(\Delta_0^{+}\right)_{i_0}}{(1+\lambda )_{i_0}(\gamma +\lambda )_{i_0}}\right.\nonumber\\
&\times& \prod _{k=1}^{n-1} \left\{ \sum_{i_k=i_{k-1}}^{\infty } \frac{(i_k+ 2k+\lambda +\alpha ) (i_k+ 2k+\lambda +\beta )}{(i_k+ 2(k+1)+\lambda )(i_k+ 2k+1+\gamma +\lambda )} \right. \nonumber\\
&\times& \left. \frac{ \left(\Delta_k^{-}\right)_{i_k} \left(\Delta_k^{+} \right)_{i_k}(2k+1+\lambda )_{i_{k-1}}(2k+\gamma +\lambda )_{i_{k-1}}}{\left(\Delta_k^{-}\right)_{i_{k-1}} \left(\Delta_k^{+} \right)_{i_{k-1}}(2k+1+\lambda )_{i_k}(2k+\gamma +\lambda )_{i_k}}\right\}\nonumber\\
&\times& \left.\left.\sum_{i_n= i_{n-1}}^{\infty } \frac{\left(\Delta_n^{-}\right)_{i_n}\left( \Delta_n^{+} \right)_{i_n}(2n+1+\lambda )_{i_{n-1}}(2n+\gamma +\lambda )_{i_{n-1}}}{\left(\Delta_n^{-}\right)_{i_{n-1}}\left(\Delta_n^{+} \right)_{i_{n-1}}(2n+1+\lambda )_{i_n}(2n+\gamma +\lambda )_{i_n}} \eta ^{i_n} \right\} z^n \right\} \label{eq:80010}
\end{eqnarray}
where
\begin{equation}
\begin{cases} 
\Delta_0^{\pm}= \frac{\{\varphi +2(1+a)\lambda \} \pm\sqrt{\varphi ^2-4(1+a)q}}{2(1+a)} \cr
\Delta_1^{\pm}=  \frac{\{\varphi +2(1+a)(\lambda+2 ) \} \pm\sqrt{\varphi ^2-4(1+a)q}}{2(1+a)} \cr
\Delta_k^{\pm}=  \frac{\{\varphi +2(1+a)(\lambda+2k )\} \pm\sqrt{\varphi ^2-4(1+a)q}}{2(1+a)} \cr
\Delta_n^{\pm}=   \frac{\{\varphi +2(1+a)(\lambda+2n )\} \pm\sqrt{\varphi ^2-4(1+a)q}}{2(1+a)}
\end{cases}\nonumber 
\end{equation}
Put (\ref{eq:8006}) in (\ref{eq:80010}). Take $c_0$= 1 as $\lambda =0$  for the first independent solution of Lame equation and $\lambda =\frac{1}{2}$ for the second one into the new (\ref{eq:80010}).
\begin{remark}
The power series expansion of Lame equation in Weierstrass's form of the first kind for infinite series about $\xi =0$ using R3TRF is given by
\begin{eqnarray}
 y(\xi )&=& LF^R\left( \rho , \alpha, h; \xi = sn^2(z,\rho ), \eta = (1+\rho ^2)\xi, z=-\rho^2 \xi^2 \right) \nonumber\\
&=& \sum_{i_0=0}^{\infty } \frac{\left(\Delta_0^{-}\right)_{i_0} \left(\Delta_0^{+}\right)_{i_0}}{(1)_{i_0}\left( \frac{1}{2}\right)_{i_0}} \eta ^{i_0} \nonumber\\
&+& \left\{ \sum_{i_0=0}^{\infty }\frac{\left( i_0+ \frac{\alpha }{2} +\frac{1}{2} \right) \left( i_0 -\frac{\alpha }{2} \right)}{(i_0+2)\left( i_0+ \frac{3}{2} \right)}\frac{\left(\Delta_0^{-}\right)_{i_0} \left(\Delta_0^{+}\right)_{i_0}}{(1 )_{i_0}\left( \frac{1}{2} \right)_{i_0}} \sum_{i_1=i_0}^{\infty } \frac{\left(\Delta_1^{-}\right)_{i_1} \left(\Delta_1^{+}\right)_{i_1}(3 )_{i_0}\left( \frac{5}{2} \right)_{i_0}}{\left(\Delta_1^{-}\right)_{i_0}  \left(\Delta_1^{+}\right)_{i_0}(3 )_{i_1}\left( \frac{5}{2} \right)_{i_1}}\eta ^{i_1}\right\} z\nonumber\\
&+& \sum_{n=2}^{\infty } \left\{ \sum_{i_0=0}^{\infty } \frac{\left( i_0+ \frac{\alpha }{2} +\frac{1}{2} \right) \left( i_0 -\frac{\alpha }{2} \right)}{(i_0+2)\left( i_0+ \frac{3}{2} \right)}\frac{\left(\Delta_0^{-}\right)_{i_0} \left(\Delta_0^{+}\right)_{i_0}}{(1 )_{i_0}\left( \frac{1}{2} \right)_{i_0}}\right.\nonumber\\
&\times& \prod _{k=1}^{n-1} \left\{ \sum_{i_k=i_{k-1}}^{\infty } \frac{\left( i_k+ 2k+\frac{\alpha }{2} +\frac{1}{2} \right) \left( i_k+ 2k-\frac{\alpha }{2} \right)}{(i_k+ 2k+2 )\left( i_k+ 2k +\frac{3}{2} \right)} \frac{ \left(\Delta_k^{-}\right)_{i_k} \left(\Delta_k^{+} \right)_{i_k}(2k+1)_{i_{k-1}}\left(2k +\frac{1}{2} \right)_{i_{k-1}}}{\left(\Delta_k^{-}\right)_{i_{k-1}} \left(\Delta_k^{+} \right)_{i_{k-1}}(2k+1 )_{i_k}\left( 2k +\frac{1}{2} \right)_{i_k}}\right\}\nonumber\\
&\times& \left.\sum_{i_n= i_{n-1}}^{\infty } \frac{\left(\Delta_n^{-}\right)_{i_n}\left( \Delta_n^{+} \right)_{i_n}\left( 2n+1 \right)_{i_{n-1}}\left( 2n+ \frac{1}{2} \right)_{i_{n-1}}}{\left(\Delta_n^{-}\right)_{i_{n-1}}\left(\Delta_n^{+} \right)_{i_{n-1}}(2n+1 )_{i_n}\left( 2n +\frac{1}{2} \right)_{i_n}} \eta ^{i_n} \right\} z^n  \label{eq:80011}
\end{eqnarray}
where
\begin{equation}
\begin{cases} 
\Delta_0^{\pm}= \pm\sqrt{\frac{h}{4(1+\rho ^2)}}  \cr
\Delta_1^{\pm}= 2 \pm\sqrt{\frac{h}{4(1+\rho ^2)}} \cr
\Delta_k^{\pm}= 2k \pm\sqrt{\frac{h}{4(1+\rho ^2)}} \cr
\Delta_n^{\pm}= 2n \pm\sqrt{\frac{h}{4(1+\rho ^2)}}
\end{cases}\nonumber 
\end{equation}
\end{remark}
\begin{remark}
The power series expansion of Lame equation in Weierstrass's form of the second kind for infinite series about $\xi =0$ using R3TRF is given by
\begin{eqnarray}
y(\xi )&=& LS^R\left( \rho , \alpha, h; \xi = sn^2(z,\rho ), \eta = (1+\rho ^2)\xi, z=-\rho^2 \xi^2 \right) \nonumber\\
&=& \xi ^{\frac{1}{2} } \left\{\sum_{i_0=0}^{\infty } \frac{\left(\Delta_0^{-}\right)_{i_0} \left(\Delta_0^{+}\right)_{i_0}}{\left(\frac{3}{2} \right)_{i_0}(1)_{i_0}} \eta ^{i_0}\right.\nonumber\\
&+& \left\{ \sum_{i_0=0}^{\infty }\frac{\left( i_0+ \frac{\alpha }{2} +1\right) \left( i_0-\frac{\alpha }{2}+\frac{1}{2} \right)}{\left( i_0+ \frac{5}{2} \right) (i_0+ 2)}\frac{\left(\Delta_0^{-}\right)_{i_0} \left(\Delta_0^{+}\right)_{i_0}}{\left( \frac{3}{2}\right)_{i_0}(1 )_{i_0}} \sum_{i_1=i_0}^{\infty } \frac{\left(\Delta_1^{-}\right)_{i_1} \left(\Delta_1^{+}\right)_{i_1}\left( \frac{7}{2}\right)_{i_0}(3)_{i_0}}{\left(\Delta_1^{-}\right)_{i_0}  \left(\Delta_1^{+}\right)_{i_0}\left( \frac{7}{2}\right)_{i_1}(3)_{i_1}}\eta ^{i_1}\right\} z\nonumber\\
&+& \sum_{n=2}^{\infty } \left\{ \sum_{i_0=0}^{\infty } \frac{\left( i_0+ \frac{\alpha }{2} +1\right) \left( i_0-\frac{\alpha }{2}+\frac{1}{2} \right)}{\left( i_0+ \frac{5}{2} \right) (i_0+ 2)}\frac{\left(\Delta_0^{-}\right)_{i_0} \left(\Delta_0^{+}\right)_{i_0}}{\left( \frac{3}{2}\right)_{i_0}(1 )_{i_0}}\right.\nonumber\\
&\times& \prod _{k=1}^{n-1} \left\{ \sum_{i_k=i_{k-1}}^{\infty } \frac{\left( i_k+ 2k+\frac{\alpha }{2} +1\right) (i_k+ 2k-\frac{\alpha }{2} +\frac{1}{2} )}{\left( i_k+ 2k +\frac{5}{2} \right) ( i_k+ 2k+2 )} \frac{ \left(\Delta_k^{-}\right)_{i_k} \left(\Delta_k^{+} \right)_{i_k}\left( 2k+\frac{3}{2} \right)_{i_{k-1}}(2k+1)_{i_{k-1}}}{\left(\Delta_k^{-}\right)_{i_{k-1}} \left(\Delta_k^{+} \right)_{i_{k-1}}\left( 2k +\frac{3}{2} \right)_{i_k} (2k+1)_{i_k}}\right\} \nonumber\\
&\times& \left.\left.\sum_{i_n= i_{n-1}}^{\infty } \frac{\left(\Delta_n^{-}\right)_{i_n}\left( \Delta_n^{+} \right)_{i_n}\left( 2n +\frac{3}{2} \right)_{i_{n-1}}(2n+1)_{i_{n-1}}}{\left(\Delta_n^{-}\right)_{i_{n-1}}\left(\Delta_n^{+} \right)_{i_{n-1}}\left( 2n +\frac{3}{2} \right)_{i_n}(2n+1)_{i_n}} \eta ^{i_n} \right\} z^n \right\} \label{eq:80012}
\end{eqnarray}
where
\begin{equation}
\begin{cases} 
\Delta_0^{\pm}= \frac{1}{2}\pm\sqrt{\frac{h}{4(1+\rho ^2)}}  \cr
\Delta_1^{\pm}= \frac{5}{2} \pm\sqrt{\frac{h}{4(1+\rho ^2)}} \cr
\Delta_k^{\pm}= 2k +\frac{1}{2} \pm\sqrt{\frac{h}{4(1+\rho ^2)}} \cr
\Delta_n^{\pm}= 2n +\frac{1}{2} \pm\sqrt{\frac{h}{4(1+\rho ^2)}}
\end{cases}\nonumber 
\end{equation}
\end{remark}
The infinite series in this chapter are equivalent to the infinite series in Ref.\cite{aChou2012g}. In this chapter $B_n$ is the leading term in sequence $c_n$ of the analytic function $y(\xi )$. In Ref.\cite{aChou2012g} $A_n$ is the leading term in sequence $c_n$ of the analytic function $y(\xi )$.

\section{Integral formalism}
Lame equation could not be built in a definite or contour integral form of any well-known simple functions such as Gauss hypergeometric, Kummer functions and etc.
According to Erdelyi (1940\cite{aErde1940}), ``there is no corresponding representation of simple integral formalisms of the solutions in ordinary linear differential equations with four regular singularities; Heun equation, Lame equation and Mathieu equation. It appears that the theory of integral equations connected with periodic solutions of Lame equation is not as complete as the corresponding theory of integral representations of, say, Legendre functions.''
Because Heun, Lame and Mathieu equations are composed of a 3-term recursive relation in the power series expansion. The three term recurrence relation in their power series creates complicated mathematical calculations to be constructed into its direct or contour integrals.
 
Instead of analyzing Lame equation into its integral representation of any simple functions, in earlier literature the integral equations of Lame equation were constructed by using simple kernels involving Legendre functions of Jacobian elliptic function or ellipsoidal harmonic functions: such integral relationships express one analytic solution in terms of another analytic solution involving Jacobian elliptic functions. There are many other forms of integral relations in the Lame equation.\cite{aErde1955,aArsc1964a,aShai1980,aSlee1968a,aVolk1982,aVolk1983,aVolk1984,aWang1989,aWhit1996}

In Ref.\cite{aChou2012g}, I show integral forms (each sub-integral $y_m(\xi )$ where $m=0,1,2,\cdots$ is composed of $2m$ terms of definite integrals and $m$ terms of contour integrals) of Lame equation in Weierstrass's form by changing all coefficients and an independent variable in the general expression of the integral form of Lame equation in the algebraic form in Ref.\cite{aChou2012f}: a $_2F_1$ function recurs in each of sub-integral forms of Lame equation.
Now I consider integral representations of Lame equation by changing all coefficients and an independent variable in the general expression of Heun equation in chapter 2.
\subsection{Polynomial of type 2}
In chapter 2 the general representation in the form of integral of Heun polynomial of type 2 about $x=0$ is given by
\begin{eqnarray}
 y(x)&=& \sum_{n=0}^{\infty } y_{n}(x)= y_0(x)+ y_1(x)+ y_2(x)+y_3(x)+\cdots \nonumber\\
&=& c_0 x^{\lambda } \Bigg\{ \sum_{i_0=0}^{q_0}\frac{(-q_0)_{i_0}\left(q_0+\frac{\varphi +2(1+a)\lambda }{(1+a)}\right)_{i_0}}{(1+\lambda )_{i_0}(\gamma +\lambda )_{i_0}}  \eta ^{i_0}\nonumber\\
&&+ \sum_{n=1}^{\infty } \Bigg\{\prod _{k=0}^{n-1} \Bigg\{ \int_{0}^{1} dt_{n-k}\;t_{n-k}^{2(n-k)-1+\lambda } \int_{0}^{1} du_{n-k}\;u_{n-k}^{2(n-k-1)+\gamma +\lambda } \nonumber\\
&&\times  \frac{1}{2\pi i}  \oint dv_{n-k} \frac{1}{v_{n-k}} \left( \frac{v_{n-k}-1}{v_{n-k}} \frac{1}{1-\overleftrightarrow {w}_{n-k+1,n}(1-t_{n-k})(1-u_{n-k})v_{n-k}}\right)^{q_{n-k}} \nonumber\\
&&\times \left( 1- \overleftrightarrow {w}_{n-k+1,n}(1-t_{n-k})(1-u_{n-k})v_{n-k}\right)^{-\left(4(n-k)+\frac{\varphi +2(1+a)\lambda }{(1+a)}\right)}\nonumber\\
&&\times \overleftrightarrow {w}_{n-k,n}^{-(2(n-k-1)+\alpha +\lambda )}\left(  \overleftrightarrow {w}_{n-k,n} \partial _{ \overleftrightarrow {w}_{n-k,n}}\right) \overleftrightarrow {w}_{n-k,n}^{\alpha -\beta} \left(  \overleftrightarrow {w}_{n-k,n} \partial _{ \overleftrightarrow {w}_{n-k,n}}\right) \overleftrightarrow {w}_{n-k,n}^{2(n-k-1)+\beta +\lambda } \Bigg\}\nonumber\\
&&\times \sum_{i_0=0}^{q_0}\frac{(-q_0)_{i_0}\left(q_0+\frac{\varphi +2(1+a)\lambda }{(1+a)}\right)_{i_0}}{(1+\lambda )_{i_0}(\gamma +\lambda )_{i_0}} \overleftrightarrow {w}_{1,n}^{i_0}\Bigg\} z^n \Bigg\} \label{eq:80013}
\end{eqnarray}
where
\begin{equation}
\begin{cases} z = -\frac{1}{a}x^2 \cr
\eta = \frac{(1+a)}{a} x \cr
\varphi = \alpha +\beta -\delta +a(\delta +\gamma -1) \cr
q= -(q_j+2j+\lambda )\{\varphi +(1+a)(q_j+2j+\lambda ) \} \;\;\mbox{as}\;j,q_j\in \mathbb{N}_{0} \cr
q_i\leq q_j \;\;\mbox{only}\;\mbox{if}\;i\leq j\;\;\mbox{where}\;i,j\in \mathbb{N}_{0} 
\end{cases}\nonumber 
\end{equation}
and
\begin{equation}\overleftrightarrow {w}_{i,j}=
\begin{cases} \displaystyle {\frac{v_i}{(v_i-1)}\; \frac{\overleftrightarrow w_{i+1,j} t_i u_i}{1- \overleftrightarrow w_{i+1,j} v_i (1-t_i)(1-u_i)}} \;\;\mbox{where}\; i\leq j\cr
\eta \;\;\mbox{only}\;\mbox{if}\; i>j
\end{cases}\nonumber 
\end{equation}
In the above, the first sub-integral form contains one term of $B_n's$, the second one contains two terms of $B_n$'s, the third one contains three terms of $B_n$'s, etc.

Put (\ref{eq:8006}) in (\ref{eq:80013}) with replacing $q_i$ by $h_i$ where $h_i\in \mathbb{N}_{0}$. Take $c_0$= 1 as $\lambda =0$  for the first independent solution of Lame equation and $\lambda =\frac{1}{2}$ for the second one into the new (\ref{eq:80013}).
\begin{remark}
The integral representation of Lame equation in Weierstrass's form of the first kind for polynomial of type 2
 about $\xi =0$ as $h= 4(1+\rho ^2)(h_j+2j )^2 $ where $j,h_j \in \mathbb{N}_{0}$ is
\begin{eqnarray}
y(\xi )&=& LF_{h_j}^R\left( \rho , \alpha, h= 4(1+\rho ^2)(h_j+2j )^2; \xi = sn^2(z,\rho ), \eta = (1+\rho ^2)\xi, z=-\rho^2 \xi^2 \right) \nonumber\\
&=&\; _2F_1 \left( -h_0, h_0; \frac{1}{2}; \eta  \right) + \sum_{n=1}^{\infty } \left\{\prod _{k=0}^{n-1} \Bigg\{ \int_{0}^{1} dt_{n-k}\;t_{n-k}^{2(n-k)-1 } \int_{0}^{1} du_{n-k}\;u_{n-k}^{2(n-k)- \frac{3}{2} }\right. \nonumber\\
&&\times  \frac{1}{2\pi i}  \oint dv_{n-k} \frac{1}{v_{n-k}} \left( \frac{v_{n-k}-1}{v_{n-k}} \frac{1}{1-\overleftrightarrow {w}_{n-k+1,n}(1-t_{n-k})(1-u_{n-k})v_{n-k}}\right)^{h_{n-k}} \nonumber\\
&&\times \left( 1- \overleftrightarrow {w}_{n-k+1,n}(1-t_{n-k})(1-u_{n-k})v_{n-k}\right)^{-4(n-k)}  \nonumber\\
&&\times \overleftrightarrow {w}_{n-k,n}^{-(2(n-k)-\frac{3}{2}+\frac{\alpha }{2} )}\left(  \overleftrightarrow {w}_{n-k,n} \partial _{ \overleftrightarrow {w}_{n-k,n}}\right) \overleftrightarrow {w}_{n-k,n}^{\alpha +\frac{1}{2}}\left(  \overleftrightarrow {w}_{n-k,n} \partial _{ \overleftrightarrow {w}_{n-k,n}}\right) \overleftrightarrow {w}_{n-k,n}^{2(n-k)-2-\frac{\alpha}{2}} \Bigg\} \nonumber\\
&&\times \left. \; _2F_1 \left( -h_0, h_0; \frac{1}{2}; \overleftrightarrow {w}_{1,n} \right) \right\} z^n  \label{eq:80014}
\end{eqnarray}
\end{remark} 
\begin{remark}
The integral representation of Lame equation in Weierstrass's form of the second kind for polynomial of type 2 about  $\xi =0$ as $h= 4(1+\rho ^2)\left( h_j+2j +\frac{1}{2}\right)^2$ where $j,h_j \in \mathbb{N}_{0}$ is
\begin{eqnarray}
y(\xi ) &=& LS_{h_j}^R\left(\rho , \alpha, h= 4(1+\rho ^2)\left( h_j+2j +\frac{1}{2}\right)^2; \xi = sn^2(z,\rho ), \eta = (1+\rho ^2)\xi, z=-\rho^2 \xi^2 \right) \nonumber\\
&=& \xi ^{\frac{1}{2} } \left\{ \; _2F_1 \left( -h_0, h_0 +1; \frac{3}{2}; \eta \right) \right. + \sum_{n=1}^{\infty } \left\{\prod _{k=0}^{n-1} \Bigg\{ \int_{0}^{1} dt_{n-k}\;t_{n-k}^{2(n-k)- \frac{1}{2} } \int_{0}^{1} du_{n-k}\;u_{n-k}^{2(n-k)-1 }\right. \nonumber\\
&&\times  \frac{1}{2\pi i}  \oint dv_{n-k} \frac{1}{v_{n-k}} \left( \frac{v_{n-k}-1}{v_{n-k}} \frac{1}{1-\overleftrightarrow {w}_{n-k+1,n}(1-t_{n-k})(1-u_{n-k})v_{n-k}}\right)^{h_{n-k}} \nonumber\\
&&\times \left( 1- \overleftrightarrow {w}_{n-k+1,n}(1-t_{n-k})(1-u_{n-k})v_{n-k}\right)^{-\left( 4(n-k) +1 \right)}  \nonumber\\
&&\times \overleftrightarrow {w}_{n-k,n}^{-(2(n-k)-1 +\frac{\alpha }{2} )}\left(  \overleftrightarrow {w}_{n-k,n} \partial _{ \overleftrightarrow {w}_{n-k,n}}\right) \overleftrightarrow {w}_{n-k,n}^{\alpha +\frac{1}{2}}\left(  \overleftrightarrow {w}_{n-k,n} \partial _{ \overleftrightarrow {w}_{n-k,n}}\right) \overleftrightarrow {w}_{n-k,n}^{2(n-k)- \frac{3}{2}-\frac{\alpha}{2}  } \Bigg\} \nonumber\\
&&\times \left.\left.\; _2F_1 \left( -h_0, h_0 +1; \frac{3}{2}; \overleftrightarrow {w}_{1,n} \right) \right\} z^n \right\}  \label{eq:80015}
\end{eqnarray}
\end{remark}
\subsection{Infinite series}
In chapter 2 the general representation in the form of integral of Heun equation for infinite series about $x=0$ using R3TRF is given by 
\begin{eqnarray}
 y(x)&=& \sum_{n=0}^{\infty } y_{n}(x)= y_0(x)+ y_1(x)+ y_2(x)+y_3(x)+\cdots \nonumber\\
&=& c_0 x^{\lambda } \left\{ \sum_{i_0=0}^{\infty }\frac{\left( \Delta_0^{-}\right)_{i_0}\left(\Delta_0^{+}\right)_{i_0}}{(1+\lambda )_{i_0}(\gamma +\lambda )_{i_0}}  \eta ^{i_0}\right. + \sum_{n=1}^{\infty } \left\{\prod _{k=0}^{n-1} \Bigg\{ \int_{0}^{1} dt_{n-k}\;t_{n-k}^{2(n-k)-1+\lambda } \int_{0}^{1} du_{n-k}\;u_{n-k}^{2(n-k-1)+\gamma +\lambda } \right.\nonumber\\
&\times& \frac{1}{2\pi i}  \oint dv_{n-k} \frac{1}{v_{n-k}}\left( \frac{v_{n-k}-1}{v_{n-k}}\right)^{-\Delta_{n-k}^{-}}  \left( 1- \overleftrightarrow {w}_{n-k+1,n}(1-t_{n-k})(1-u_{n-k})v_{n-k}\right)^{-\Delta_{n-k}^{+}}\nonumber\\
&\times& \overleftrightarrow {w}_{n-k,n}^{-(2(n-k-1)+\alpha +\lambda )}\left(  \overleftrightarrow {w}_{n-k,n} \partial _{ \overleftrightarrow {w}_{n-k,n}}\right) \overleftrightarrow {w}_{n-k,n}^{\alpha -\beta} \left(  \overleftrightarrow {w}_{n-k,n} \partial _{ \overleftrightarrow {w}_{n-k,n}}\right) \overleftrightarrow {w}_{n-k,n}^{2(n-k-1)+\beta +\lambda } \Bigg\}\nonumber\\
&\times& \left.\left.\sum_{i_0=0}^{\infty }\frac{\left( \Delta_0^{-}\right)_{i_0}\left(\Delta_0^{+}\right)_{i_0}}{(1+\lambda )_{i_0}(\gamma +\lambda )_{i_0}} \overleftrightarrow {w}_{1,n}^{i_0}\right\} z^n \right\} \label{eq:80016}
\end{eqnarray}
where
\begin{equation}
\begin{cases} 
\Delta_0^{\pm}=  \frac{ \varphi +2(1+a)\lambda \pm\sqrt{\varphi ^2-4(1+a)q}}{2(1+a)} \cr
\Delta_{n-k}^{\pm}=  \frac{ \varphi +2(1+a)(\lambda+2(n-k) ) \pm\sqrt{\varphi ^2-4(1+a)q}}{2(1+a)}
\end{cases}\nonumber 
\end{equation}
Put (\ref{eq:8006}) in (\ref{eq:80016}). Take $c_0$= 1 as $\lambda =0$  for the first independent solution of Lame equation and $\lambda =\frac{1}{2}$ for the second one into the new (\ref{eq:80016}).
\begin{remark}
The integral representation of Lame equation in Weierstrass's form of the first kind for infinite series about $\xi =0$ using R3TRF is
\begin{eqnarray}
 y(\xi )&=& LF^R\left( \rho , \alpha, h; \xi = sn^2(z,\rho ), \eta = (1+\rho ^2)\xi, z=-\rho^2 \xi^2 \right) \nonumber\\
&=& \; _2F_1 \left(  \Delta_0^{-}, \Delta_0^{+}; \frac{1}{2}; \eta \right) + \sum_{n=1}^{\infty } \left\{\prod _{k=0}^{n-1} \Bigg\{ \int_{0}^{1} dt_{n-k}\;t_{n-k}^{2(n-k)-1 } \int_{0}^{1} du_{n-k}\;u_{n-k}^{2(n-k)-\frac{3}{2} } \right.\nonumber\\
&\times& \frac{1}{2\pi i}  \oint dv_{n-k} \frac{1}{v_{n-k}}\left( \frac{v_{n-k}-1}{v_{n-k}}\right)^{-\Delta_{n-k}^{-}}  \left( 1- \overleftrightarrow {w}_{n-k+1,n}(1-t_{n-k})(1-u_{n-k})v_{n-k}\right)^{-\Delta_{n-k}^{+}}\nonumber\\
&\times& \overleftrightarrow {w}_{n-k,n}^{-(2(n-k)-\frac{3}{2}+\frac{\alpha }{2})}\left( \overleftrightarrow {w}_{n-k,n} \partial _{ \overleftrightarrow {w}_{n-k,n}}\right) \overleftrightarrow {w}_{n-k,n}^{\alpha+\frac{1}{2}} \left( \overleftrightarrow {w}_{n-k,n} \partial _{ \overleftrightarrow {w}_{n-k,n}}\right) \overleftrightarrow {w}_{n-k,n}^{2(n-k)-2 -\frac{\alpha }{2} } \Bigg\}\nonumber\\
&\times& \left.  \; _2F_1 \left( \Delta_0^{-}, \Delta_0^{+}; \frac{1}{2}; \overleftrightarrow {w}_{1,n} \right) \right\} z^n \label{eq:80017}
\end{eqnarray}
where
\begin{equation}
\begin{cases} 
\Delta_0^{\pm}=  \pm\sqrt{\frac{h}{4(1+\rho ^2)}} \cr
\Delta_{n-k}^{\pm}= 2(n-k)\pm\sqrt{\frac{h}{4(1+\rho ^2)}} 
\end{cases}\nonumber 
\end{equation}
\end{remark}
\begin{remark}
The integral representation of Lame equation in Weierstrass's form of the second kind for infinite series about $\xi =0$ using R3TRF is
\begin{eqnarray}
y(\xi )&=& LS^R\left( \rho , \alpha, h; \xi = sn^2(z,\rho ), \eta = (1+\rho ^2)\xi, z=-\rho^2 \xi^2 \right) \nonumber\\
&=& \xi ^{\frac{1}{2}} \left\{ \; _2F_1 \left( \Delta_0^{-}, \Delta_0^{+}; \frac{3}{2}; \eta \right) \right. + \sum_{n=1}^{\infty } \left\{\prod _{k=0}^{n-1} \Bigg\{ \int_{0}^{1} dt_{n-k}\;t_{n-k}^{2(n-k)-\frac{1}{2} } \int_{0}^{1} du_{n-k}\;u_{n-k}^{2(n-k)-1} \right.\nonumber\\
&\times& \frac{1}{2\pi i}  \oint dv_{n-k} \frac{1}{v_{n-k}}\left( \frac{v_{n-k}-1}{v_{n-k}}\right)^{-\Delta_{n-k}^{-}}  \left( 1- \overleftrightarrow {w}_{n-k+1,n}(1-t_{n-k})(1-u_{n-k})v_{n-k}\right)^{-\Delta_{n-k}^{+}}\nonumber\\
&\times& \overleftrightarrow {w}_{n-k,n}^{-(2(n-k)-1 +\frac{\alpha }{2} )}\left(  \overleftrightarrow {w}_{n-k,n} \partial _{ \overleftrightarrow {w}_{n-k,n}}\right) \overleftrightarrow {w}_{n-k,n}^{\alpha +\frac{1}{2}} \left(  \overleftrightarrow {w}_{n-k,n} \partial _{ \overleftrightarrow {w}_{n-k,n}}\right) \overleftrightarrow {w}_{n-k,n}^{2(n-k)-\frac{3}{2} -\frac{\alpha }{2} } \Bigg\}\nonumber\\
&\times& \left.\left. \; _2F_1 \left( \Delta_0^{-}, \Delta_0^{+}; \frac{3}{2}; \overleftrightarrow {w}_{1,n} \right)\right\} z^n \right\} \label{eq:80018}
\end{eqnarray}
where
\begin{equation}
\begin{cases} 
\Delta_0^{\pm}=  \frac{1}{2} \pm\sqrt{\frac{h}{4(1+\rho ^2)}} \cr
\Delta_{n-k}^{\pm}=  2(n-k)+\frac{1}{2} \pm\sqrt{\frac{h}{4(1+\rho ^2)}} 
\end{cases}\nonumber 
\end{equation}
\end{remark}
\section[Generating function for the Lame polynomial of type 2]{Generating function for the Lame polynomial of type 2
\sectionmark{Generating function for the Lame polynomial of type 2}}
\sectionmark{Generating function for the Lame polynomial of type 2}
\begin{definition}
I define that
\begin{equation}
\begin{cases}
\displaystyle { s_{a,b}} = \begin{cases} \displaystyle {  s_a\cdot s_{a+1}\cdot s_{a+2}\cdots s_{b-2}\cdot s_{b-1}\cdot s_b}\;\;\mbox{if}\;a>b \cr
s_a \;\;\mbox{if}\;a=b\end{cases}
\cr
\cr
\displaystyle { \widetilde{w}_{i,j}}  = 
\begin{cases} \displaystyle { \frac{ \widetilde{w}_{i+1,j}\; t_i u_i \left\{ 1+ (s_i+2\widetilde{w}_{i+1,j}(1-t_i)(1-u_i))s_i\right\}}{2(1-\widetilde{w}_{i+1,j}(1-t_i)(1-u_i))^2 s_i}} \cr
\displaystyle {-\frac{\widetilde{w}_{i+1,j}\; t_i u_i (1+s_i)\sqrt{s_i^2-2(1-2\widetilde{w}_{i+1,j}(1-t_i)(1-u_i))s_i+1}}{2(1-\widetilde{w}_{i+1,j}(1-t_i)(1-u_i))^2 s_i}} \;\;\mbox{where}\;i<j \cr
\cr
\displaystyle { \frac{\eta t_i u_i \left\{ 1+ (s_{i,\infty }+2\eta(1-t_i)(1-u_i))s_{i,\infty }\right\}}{2(1-\eta (1-t_i)(1-u_i))^2 s_{i,\infty }}} \cr
\displaystyle {-\frac{\eta t_i u_i(1+s_{i,\infty })\sqrt{s_{i,\infty }^2-2(1-2\eta (1-t_i)(1-u_i))s_{i,\infty }+1}}{2(1-\eta (1-t_i)(1-u_i))^2 s_{i,\infty }}} \;\;\mbox{where}\;i=j 
\end{cases}
\end{cases}\nonumber
\end{equation}
where
\begin{equation}
a,b,i,j\in \mathbb{N}_{0} \nonumber
\end{equation}
\end{definition}
In chapter 3 the generating function for the Heun polynomial of type 2 of the first kind about $x=0$ as $q =-(q_j+2j)\{\alpha +\beta -\delta +a(\delta +\gamma -1)+(1+a)(q_j+2j)\} $ where $j,q_j \in \mathbb{N}_{0}$ is given by
\begin{eqnarray}
&&\sum_{q_0 =0}^{\infty } \frac{(\gamma)_{q_0}}{q_0!} s_0^{q_0} \prod _{n=1}^{\infty } \left\{ \sum_{ q_n = q_{n-1}}^{\infty } s_n^{q_n }\right\} HF_{q_j}^R \Bigg( q_j =\frac{-\varphi \pm \sqrt{\varphi ^2-4(1+a)q}}{2(1+a)}-2j\nonumber\\
&&, \varphi =\alpha +\beta -\delta +a(\delta +\gamma -1), \Omega _1=\frac{\varphi }{(1+a)}; \eta = \frac{(1+a)}{a} x ; z= -\frac{1}{a} x^2 \Bigg) \nonumber\\
&&=2^{\frac{\varphi }{(1+a)}-1}\Bigg\{ \prod_{l=1}^{\infty } \frac{1}{(1-s_{l,\infty })}  \mathbf{A}\left( s_{0,\infty } ;\eta\right) + \Bigg\{ \prod_{l=2}^{\infty } \frac{1}{(1-s_{l,\infty })} \int_{0}^{1} dt_1\;t_1 \int_{0}^{1} du_1\;u_1^{\gamma} \overleftrightarrow {\mathbf{\Gamma}}_1 \left( s_{1,\infty };t_1,u_1,\eta\right)\nonumber\\
&&\times \widetilde{w}_{1,1}^{-\alpha}\left( \widetilde{w}_{1,1} \partial _{ \widetilde{w}_{1,1}}\right) \widetilde{w}_{1,1}^{\alpha -\beta } \left( \widetilde{w}_{1,1} \partial _{ \widetilde{w}_{1,1}}\right)\widetilde{w}_{1,1}^{\beta } \mathbf{A}\left( s_{0} ;\widetilde{w}_{1,1}\right)\Bigg\} z \nonumber\\
&&+ \sum_{n=2}^{\infty } \Bigg\{ \prod_{l=n+1}^{\infty } \frac{1}{(1-s_{l,\infty })} \int_{0}^{1} dt_n\;t_n^{2n-1} \int_{0}^{1} du_n\;u_n^{2(n-1)+\gamma } \overleftrightarrow {\mathbf{\Gamma}}_n \left( s_{n,\infty };t_n,u_n,\eta \right)\nonumber\\
&&\times \widetilde{w}_{n,n}^{-(2(n-1)+\alpha)}\left( \widetilde{w}_{n,n} \partial _{ \widetilde{w}_{n,n}}\right) \widetilde{w}_{n,n}^{\alpha -\beta } \left( \widetilde{w}_{n,n} \partial _{ \widetilde{w}_{n,n}}\right)\widetilde{w}_{n,n}^{2(n-1)+\beta }  \nonumber\\
&&\times \prod_{k=1}^{n-1} \Bigg\{ \int_{0}^{1} dt_{n-k}\;t_{n-k}^{2(n-k)-1} \int_{0}^{1} du_{n-k} \;u_{n-k}^{2(n-k-1)+\gamma}\overleftrightarrow {\mathbf{\Gamma}}_{n-k} \left( s_{n-k};t_{n-k},u_{n-k},\widetilde{w}_{n-k+1,n} \right)\label{eq:80019}\\
&&\times \widetilde{w}_{n-k,n}^{-(2(n-k-1)+\alpha )}\left( \widetilde{w}_{n-k,n} \partial _{ \widetilde{w}_{n-k,n}}\right) \widetilde{w}_{n-k,n}^{\alpha -\beta } \left( \widetilde{w}_{n-k,n} \partial _{ \widetilde{w}_{n-k,n}}\right)\widetilde{w}_{n-k,n}^{2(n-k-1)+\beta } \Bigg\} \mathbf{A} \left( s_{0} ;\widetilde{w}_{1,n}\right) \Bigg\} z^n \Bigg\}   \nonumber 
\end{eqnarray}
where
\begin{equation}
\begin{cases} 
{ \displaystyle \overleftrightarrow {\mathbf{\Gamma}}_1 \left( s_{1,\infty };t_1,u_1,\eta\right)= \frac{\left( \frac{1+s_{1,\infty }+\sqrt{s_{1,\infty }^2-2(1-2\eta (1-t_1)(1-u_1))s_{1,\infty }+1}}{2}\right)^{-\left(3+\Omega _1\right)}}{\sqrt{s_{1,\infty }^2-2(1-2\eta (1-t_1)(1-u_1))s_{1,\infty }+1}}}\cr
{ \displaystyle  \overleftrightarrow {\mathbf{\Gamma}}_n \left( s_{n,\infty };t_n,u_n,\eta \right) =\frac{\left( \frac{1+s_{n,\infty }+\sqrt{s_{n,\infty }^2-2(1-2\eta (1-t_n)(1-u_n))s_{n,\infty }+1}}{2}\right)^{-\left( 4n-1+\Omega _1\right)}}{\sqrt{ s_{n,\infty }^2-2(1-2\eta (1-t_n)(1-u_n))s_{n,\infty }+1}}}\cr
{ \displaystyle \overleftrightarrow {\mathbf{\Gamma}}_{n-k} \left( s_{n-k};t_{n-k},u_{n-k},\widetilde{w}_{n-k+1,n} \right) }\cr
{ \displaystyle = \frac{ \left( \frac{1+s_{n-k}+\sqrt{s_{n-k}^2-2(1-2\widetilde{w}_{n-k+1,n} (1-t_{n-k})(1-u_{n-k}))s_{n-k}+1}}{2}\right)^{-\left(4(n-k)-1+\Omega _1\right)}}{\sqrt{ s_{n-k}^2-2(1-2\widetilde{w}_{n-k+1,n} (1-t_{n-k})(1-u_{n-k}))s_{n-k}+1}}}
\end{cases}\nonumber 
\end{equation}
and
\begin{equation}
\begin{cases} 
{ \displaystyle \mathbf{A} \left( s_{0,\infty } ;\eta\right)}\cr
{ \displaystyle = \frac{\left(1- s_{0,\infty }+\sqrt{s_{0,\infty }^2-2(1-2\eta )s_{0,\infty }+1}\right)^{1-\gamma } \left(1+s_{0,\infty }+\sqrt{s_{0,\infty }^2-2(1-2\eta )s_{0,\infty }+1}\right)^{\gamma -\Omega _1}}{\sqrt{s_{0,\infty }^2-2(1-2\eta )s_{0,\infty }+1}}}\cr
{ \displaystyle  \mathbf{A} \left( s_{0} ;\widetilde{w}_{1,1}\right) = \frac{\left(1- s_0+\sqrt{s_0^2-2(1-2\widetilde{w}_{1,1})s_0+1}\right)^{1-\gamma} \left(1+s_0+\sqrt{s_0^2-2(1-2\widetilde{w}_{1,1} )s_0+1}\right)^{\gamma -\Omega _1}}{\sqrt{s_0^2-2(1-2\widetilde{w}_{1,1})s_0+1}}} \cr
{ \displaystyle \mathbf{A} \left(  s_{0} ;\widetilde{w}_{1,n}\right) = \frac{\left(1- s_0+\sqrt{s_0^2-2(1-2\widetilde{w}_{1,n})s_0+1}\right)^{1-\gamma} \left(1+s_0+\sqrt{s_0^2-2(1-2\widetilde{w}_{1,n} )s_0+1}\right)^{\gamma -\Omega _1}}{\sqrt{s_0^2-2(1-2\widetilde{w}_{1,n})s_0+1}}}
\end{cases}\nonumber 
\end{equation}
In chapter 3 the generating function for the Heun polynomial of type 2 of the second kind about $x=0$ as $q =-(q_j+2j+1-\gamma )\{\alpha +\beta +1-\gamma -(1-a)\delta +(1+a)(q_j+2j)\} $ where $j, q_j \in \mathbb{N}_{0}$ with replacing $c_0=\left( \frac{1+a}{a}\right)^{1-\gamma }$ by 1 is 
\begin{eqnarray}
&&\sum_{q_0 =0}^{\infty } \frac{(2-\gamma )_{q_0}}{q_0!} s_0^{q_0} \prod _{n=1}^{\infty } \left\{ \sum_{ q_n = q_{n-1}}^{\infty } s_n^{q_n }\right\} HS_{q_j}^R \Bigg( q_j =\frac{-\{\varphi +2(1+a)(1-\gamma )\} \pm \sqrt{\varphi ^2-4(1+a)q}}{2(1+a)}-2j\nonumber\\
&&, \varphi =\alpha +\beta -\delta +a(\delta +\gamma -1), \Omega _2= \frac{\varphi +2(1+a)(1-\gamma )}{(1+a)}; \eta = \frac{(1+a)}{a} x ; z= -\frac{1}{a} x^2 \Bigg) \nonumber\\
&&= 2^{\frac{\varphi +2(1+a)(1/2-\gamma )}{(1+a)}} x^{1-\gamma }\left\{ \prod_{l=1}^{\infty } \frac{1}{(1-s_{l,\infty })} \mathbf{B}\left( s_{0,\infty } ;\eta\right) \right. \nonumber\\
&&+ \Bigg\{ \prod_{l=2}^{\infty } \frac{1}{(1-s_{l,\infty })} \int_{0}^{1} dt_1\;t_1^{2-\gamma } \int_{0}^{1} du_1\;u_1 \overleftrightarrow {\mathbf{\Psi}}_1 \left( s_{1,\infty };t_1,u_1,\eta\right)\nonumber\\
&& \times \widetilde{w}_{1,1}^{-(\alpha-\gamma +1)}\left( \widetilde{w}_{1,1} \partial _{ \widetilde{w}_{1,1}}\right) \widetilde{w}_{1,1}^{\alpha -\beta } \left( \widetilde{w}_{1,1} \partial _{ \widetilde{w}_{1,1}}\right)\widetilde{w}_{1,1}^{\beta -\gamma +1} \mathbf{B}\left( s_{0} ;\widetilde{w}_{1,1}\right)\Bigg\} z \nonumber\\
&&+ \sum_{n=2}^{\infty } \Bigg\{ \prod_{l=n+1}^{\infty } \frac{1}{(1-s_{l,\infty })} \int_{0}^{1} dt_n\;t_n^{2n-\gamma } \int_{0}^{1} du_n\;u_n^{2n-1} \overleftrightarrow {\mathbf{\Psi}}_n \left( s_{n,\infty };t_n,u_n,\eta \right)\nonumber\\
&&\times \widetilde{w}_{n,n}^{-(2n-1+\alpha-\gamma )}\left( \widetilde{w}_{n,n} \partial _{ \widetilde{w}_{n,n}}\right) \widetilde{w}_{n,n}^{\alpha -\beta } \left( \widetilde{w}_{n,n} \partial _{ \widetilde{w}_{n,n}}\right)\widetilde{w}_{n,n}^{2n-1+\beta -\gamma } \nonumber\\
&&\times \prod_{k=1}^{n-1} \Bigg\{ \int_{0}^{1} dt_{n-k}\;t_{n-k}^{2(n-k)-\gamma} \int_{0}^{1} du_{n-k} \;u_{n-k}^{2(n-k)-1} \overleftrightarrow {\mathbf{\Psi}}_{n-k} \left( s_{n-k};t_{n-k},u_{n-k},\widetilde{w}_{n-k+1,n} \right) \label{eq:80020}\\
&&\times \widetilde{w}_{n-k,n}^{-(2(n-k)-1+\alpha-\gamma )}\left( \widetilde{w}_{n-k,n} \partial _{ \widetilde{w}_{n-k,n}}\right) \widetilde{w}_{n-k,n}^{\alpha -\beta } \left( \widetilde{w}_{n-k,n} \partial _{ \widetilde{w}_{n-k,n}}\right)\widetilde{w}_{n-k,n}^{2(n-k)-1+\beta -\gamma } \Bigg\} \left. \mathbf{B}\left( s_{0} ;\widetilde{w}_{1,n}\right)\Bigg\} z^n  \right\}  \nonumber
\end{eqnarray}
where
\begin{equation}
\begin{cases} 
{ \displaystyle \overleftrightarrow {\mathbf{\Psi}}_1 \left( s_{1,\infty };t_1,u_1,\eta\right)= \frac{\left( \frac{1+s_{1,\infty }+\sqrt{s_{1,\infty }^2-2(1-2\eta (1-t_1)(1-u_1))s_{1,\infty }+1}}{2}\right)^{-\left( 3+\Omega _2\right)}}{\sqrt{s_{1,\infty }^2-2(1-2\eta (1-t_1)(1-u_1))s_{1,\infty }+1}} }\cr
{ \displaystyle  \overleftrightarrow {\mathbf{\Psi}}_n \left( s_{n,\infty };t_n,u_n,\eta \right) = \frac{\left( \frac{1+s_{n,\infty }+\sqrt{s_{n,\infty }^2-2(1-2\eta (1-t_n)(1-u_n))s_{n,\infty }+1}}{2}\right)^{-\left(4n-1+\Omega _2\right)}}{\sqrt{s_{n,\infty }^2-2(1-2\eta (1-t_n)(1-u_n))s_{n,\infty }+1}}}\cr
{ \displaystyle \overleftrightarrow {\mathbf{\Psi}}_{n-k} \left( s_{n-k};t_{n-k},u_{n-k},\widetilde{w}_{n-k+1,n} \right) } \cr
{ \displaystyle = \frac{\left( \frac{(1+s_{n-k})+\sqrt{s_{n-k}^2-2(1-2\widetilde{w}_{n-k+1,n} (1-t_{n-k})(1-u_{n-k}))s_{n-k}+1}}{2}\right)^{-\left( 4(n-k)-1+\Omega _2\right)}}{\sqrt{s_{n-k}^2-2(1-2\widetilde{w}_{n-k+1,n} (1-t_{n-k})(1-u_{n-k}))s_{n-k}+1}}}
\end{cases}\nonumber 
\end{equation}
and
\begin{equation}
\begin{cases} 
{ \displaystyle \mathbf{B} \left( s_{0,\infty } ;\eta\right)= \frac{ \left(1+s_{0,\infty }+\sqrt{s_{0,\infty }^2-2(1-2\eta )s_{0,\infty }+1}\right)^{2-\gamma -\Omega _2}}{\left(1- s_{0,\infty }+\sqrt{s_{0,\infty }^2-2(1-2\eta )s_{0,\infty }+1}\right)^{1-\gamma}\sqrt{s_{0,\infty }^2-2(1-2\eta )s_{0,\infty }+1}}}\cr
{ \displaystyle  \mathbf{B} \left( s_{0} ;\widetilde{w}_{1,1}\right) = \frac{\left(1+s_0+\sqrt{s_0^2-2(1-2\widetilde{w}_{1,1} )s_0+1}\right)^{2-\gamma -\Omega _2}}{\left(1- s_0+\sqrt{s_0^2-2(1-2\widetilde{w}_{1,1})s_0+1}\right)^{1-\gamma}\sqrt{s_0^2-2(1-2\widetilde{w}_{1,1})s_0+1}}} \cr
{ \displaystyle \mathbf{B} \left( s_{0} ;\widetilde{w}_{1,n}\right) = \frac{\left(1+s_0+\sqrt{s_0^2-2(1-2\widetilde{w}_{1,n} )s_0+1}\right)^{2-\gamma -\Omega _2}}{\left(1- s_0+\sqrt{s_0^2-2(1-2\widetilde{w}_{1,n})s_0+1}\right)^{1-\gamma}\sqrt{s_0^2-2(1-2\widetilde{w}_{1,n})s_0+1}}}
\end{cases}\nonumber 
\end{equation}
Put (\ref{eq:8006}) in (\ref{eq:80019}) and (\ref{eq:80020}) with replacing $q_i$ by $h_i$ where $h_i\in \mathbb{N}_{0}$ for the first and second independent solutions of Lame equation.  
\begin{remark}
The generating function for the Lame polynomial of type 2 in Weierstrass's form of the first kind about $\xi =0$ as $h= 4(1+\rho ^2)(h_j+2j )^2 $ where $j,h_j \in \mathbb{N}_{0}$ is given by
\begin{eqnarray}
&&\sum_{h_0 =0}^{\infty } \frac{(\frac{1}{2})_{h_0}}{h_0!} s_0^{h_0} \prod _{n=1}^{\infty } \left\{ \sum_{ h_n = h_{n-1}}^{\infty } s_n^{h_n }\right\} LF_{h_j}^R\left( \rho , \alpha, h= 4(1+\rho ^2)(h_j+2j )^2; \xi = sn^2(z,\rho )\right. \nonumber\\
&&,\left. \eta = (1+\rho ^2)\xi, z=-\rho^2 \xi^2 \right)  \nonumber\\
&&=2^{-1}\Bigg\{ \prod_{l=1}^{\infty } \frac{1}{(1-s_{l,\infty })}  \mathbf{A}\left( s_{0,\infty } ;\eta\right) + \Bigg\{ \prod_{l=2}^{\infty } \frac{1}{(1-s_{l,\infty })} \int_{0}^{1} dt_1\;t_1 \int_{0}^{1} du_1\;u_1^{\frac{1}{2}} \overleftrightarrow {\mathbf{\Gamma}}_1 \left( s_{1,\infty };t_1,u_1,\eta\right)\nonumber\\
&&\times \widetilde{w}_{1,1}^{-\frac{1}{2}(1+\alpha )}\left( \widetilde{w}_{1,1} \partial _{ \widetilde{w}_{1,1}}\right) \widetilde{w}_{1,1}^{\alpha +\frac{1}{2} } \left( \widetilde{w}_{1,1} \partial _{ \widetilde{w}_{1,1}}\right)\widetilde{w}_{1,1}^{-\frac{\alpha }{2} } \mathbf{A}\left( s_{0} ;\widetilde{w}_{1,1}\right) \Bigg\} z \nonumber\\
&&+ \sum_{n=2}^{\infty } \Bigg\{ \prod_{l=n+1}^{\infty } \frac{1}{(1-s_{l,\infty })} \int_{0}^{1} dt_n\;t_n^{2n-1} \int_{0}^{1} du_n\;u_n^{2n-\frac{3}{2} } \overleftrightarrow {\mathbf{\Gamma}}_n \left( s_{n,\infty };t_n,u_n,\eta \right)\nonumber\\
&&\times \widetilde{w}_{n,n}^{-\frac{1}{2}(4n-3+\alpha )}\left( \widetilde{w}_{n,n} \partial _{ \widetilde{w}_{n,n}}\right) \widetilde{w}_{n,n}^{\alpha +\frac{1}{2} } \left( \widetilde{w}_{n,n} \partial _{ \widetilde{w}_{n,n}}\right)\widetilde{w}_{n,n}^{\frac{1}{2}(4(n-1)-\alpha ) }  \nonumber\\
&&\times \prod_{k=1}^{n-1} \Bigg\{ \int_{0}^{1} dt_{n-k}\;t_{n-k}^{2(n-k)-1} \int_{0}^{1} du_{n-k} \;u_{n-k}^{2(n-k)-\frac{3}{2}}\overleftrightarrow {\mathbf{\Gamma}}_{n-k} \left( s_{n-k};t_{n-k},u_{n-k},\widetilde{w}_{n-k+1,n} \right)\label{eq:80021}\\
&&\times \widetilde{w}_{n-k,n}^{-\frac{1}{2}(4(n-k)-3+\alpha )}\left( \widetilde{w}_{n-k,n} \partial _{ \widetilde{w}_{n-k,n}}\right) \widetilde{w}_{n-k,n}^{\alpha +\frac{1}{2}} \left( \widetilde{w}_{n-k,n} \partial _{ \widetilde{w}_{n-k,n}}\right)\widetilde{w}_{n-k,n}^{\frac{1}{2}(4(n-k-1)-\alpha )} \Bigg\} \mathbf{A} \left( s_{0} ;\widetilde{w}_{1,n}\right) \Bigg\} z^n \Bigg\} \nonumber  
\end{eqnarray}
where
\begin{equation}
\begin{cases} 
{ \displaystyle \overleftrightarrow {\mathbf{\Gamma}}_1 \left( s_{1,\infty };t_1,u_1,\eta\right)= \frac{\left( \frac{1+s_{1,\infty }+\sqrt{s_{1,\infty }^2-2(1-2\eta (1-t_1)(1-u_1))s_{1,\infty }+1}}{2}\right)^{-3}}{\sqrt{s_{1,\infty }^2-2(1-2\eta (1-t_1)(1-u_1))s_{1,\infty }+1}}}\cr
{ \displaystyle  \overleftrightarrow {\mathbf{\Gamma}}_n \left( s_{n,\infty };t_n,u_n,\eta \right) =\frac{\left( \frac{1+s_{n,\infty }+\sqrt{s_{n,\infty }^2-2(1-2\eta (1-t_n)(1-u_n))s_{n,\infty }+1}}{2}\right)^{-4n+1 }}{\sqrt{ s_{n,\infty }^2-2(1-2\eta (1-t_n)(1-u_n))s_{n,\infty }+1}}}\cr
{ \displaystyle \overleftrightarrow {\mathbf{\Gamma}}_{n-k} \left( s_{n-k};t_{n-k},u_{n-k},\widetilde{w}_{n-k+1,n} \right) = \frac{ \left( \frac{1+s_{n-k}+\sqrt{s_{n-k}^2-2(1-2\widetilde{w}_{n-k+1,n} (1-t_{n-k})(1-u_{n-k}))s_{n-k}+1}}{2}\right)^{-4(n-k)+1 }}{\sqrt{ s_{n-k}^2-2(1-2\widetilde{w}_{n-k+1,n} (1-t_{n-k})(1-u_{n-k}))s_{n-k}+1}}}
\end{cases}\nonumber 
\end{equation}
and
\begin{equation}
\begin{cases} 
{ \displaystyle \mathbf{A} \left( s_{0,\infty } ;\eta\right)= \frac{\left(1- s_{0,\infty }+\sqrt{s_{0,\infty }^2-2(1-2\eta )s_{0,\infty }+1}\right)^{\frac{1}{2}} \left( 1+s_{0,\infty }+\sqrt{s_{0,\infty }^2-2(1-2\eta )s_{0,\infty }+1}\right)^{\frac{1}{2}}}{\sqrt{s_{0,\infty }^2-2(1-2\eta )s_{0,\infty }+1}}}\cr
{ \displaystyle  \mathbf{A} \left( s_{0} ;\widetilde{w}_{1,1}\right) = \frac{\left( 1- s_0+\sqrt{s_0^2-2(1-2\widetilde{w}_{1,1})s_0+1}\right)^{\frac{1}{2}} \left( 1+s_0+\sqrt{s_0^2-2(1-2\widetilde{w}_{1,1} )s_0+1}\right)^{\frac{1}{2}}}{\sqrt{s_0^2-2(1-2\widetilde{w}_{1,1})s_0+1}}} \cr
{ \displaystyle \mathbf{A} \left(  s_{0} ;\widetilde{w}_{1,n}\right) = \frac{\left( 1- s_0+\sqrt{s_0^2-2(1-2\widetilde{w}_{1,n})s_0+1}\right)^{\frac{1}{2}} \left( 1+s_0+\sqrt{s_0^2-2(1-2\widetilde{w}_{1,n} )s_0+1}\right)^{\frac{1}{2}}}{\sqrt{s_0^2-2(1-2\widetilde{w}_{1,n})s_0+1}}}
\end{cases}\nonumber 
\end{equation}
\end{remark}
\begin{remark}
The generating function for the Lame polynomial of type 2 in Weierstrass's form of the second kind about $\xi =0$ as $h= 4(1+\rho ^2)\left( h_j+2j +\frac{1}{2}\right)^2$ where $j,h_j \in \mathbb{N}_{0}$ is given by 
\begin{eqnarray}
&&\sum_{h_0 =0}^{\infty } \frac{(\frac{3}{2})_{h_0}}{h_0!} s_0^{h_0} \prod _{n=1}^{\infty } \left\{ \sum_{ h_n = h_{n-1}}^{\infty } s_n^{h_n }\right\} LS_{h_j}^R\left(\rho , \alpha, h= 4(1+\rho ^2)\left( h_j+2j +\frac{1}{2}\right)^2; \xi = sn^2(z,\rho )\right. \nonumber\\
&&,\left. \eta = (1+\rho ^2)\xi, z=-\rho^2 \xi^2 \right) \nonumber\\
&&= \xi^{\frac{1}{2}} \left\{ \prod_{l=1}^{\infty } \frac{1}{(1-s_{l,\infty })} \mathbf{B}\left( s_{0,\infty } ;\eta\right) \right. + \Bigg\{ \prod_{l=2}^{\infty } \frac{1}{(1-s_{l,\infty })} \int_{0}^{1} dt_1\;t_1^{\frac{3}{2} } \int_{0}^{1} du_1\;u_1 \overleftrightarrow {\mathbf{\Psi}}_1 \left( s_{1,\infty };t_1,u_1,\eta\right)\nonumber\\
&& \times \widetilde{w}_{1,1}^{-\frac{1}{2}(2+\alpha )}\left( \widetilde{w}_{1,1} \partial _{ \widetilde{w}_{1,1}}\right) \widetilde{w}_{1,1}^{\alpha +\frac{1}{2}} \left( \widetilde{w}_{1,1} \partial _{ \widetilde{w}_{1,1}}\right)\widetilde{w}_{1,1}^{\frac{1}{2}(1-\alpha )} \mathbf{B}\left( s_{0} ;\widetilde{w}_{1,1}\right)\Bigg\} z \nonumber\\
&&+ \sum_{n=2}^{\infty } \Bigg\{ \prod_{l=n+1}^{\infty } \frac{1}{(1-s_{l,\infty })} \int_{0}^{1} dt_n\;t_n^{2n-\frac{1}{2} } \int_{0}^{1} du_n\;u_n^{2n-1} \overleftrightarrow {\mathbf{\Psi}}_n \left( s_{n,\infty };t_n,u_n,\eta \right)\nonumber\\
&&\times \widetilde{w}_{n,n}^{-\frac{1}{2}(4n-2+\alpha)}\left( \widetilde{w}_{n,n} \partial _{ \widetilde{w}_{n,n}}\right) \widetilde{w}_{n,n}^{\alpha +\frac{1}{2}} \left( \widetilde{w}_{n,n} \partial _{ \widetilde{w}_{n,n}}\right)\widetilde{w}_{n,n}^{\frac{1}{2}(4n-3-\alpha ) } \nonumber\\
&&\times \prod_{k=1}^{n-1} \Bigg\{ \int_{0}^{1} dt_{n-k}\;t_{n-k}^{2(n-k)-\frac{1}{2}} \int_{0}^{1} du_{n-k} \;u_{n-k}^{2(n-k)-1} \overleftrightarrow {\mathbf{\Psi}}_{n-k} \left( s_{n-k};t_{n-k},u_{n-k},\widetilde{w}_{n-k+1,n} \right) \label{eq:80022}\\
&&\times \widetilde{w}_{n-k,n}^{-\frac{1}{2}(4(n-k)-2+\alpha )}\left( \widetilde{w}_{n-k,n} \partial _{ \widetilde{w}_{n-k,n}}\right) \widetilde{w}_{n-k,n}^{\alpha +\frac{1}{2}} \left( \widetilde{w}_{n-k,n} \partial _{ \widetilde{w}_{n-k,n}}\right)\widetilde{w}_{n-k,n}^{\frac{1}{2}(4(n-k)-3-\alpha )} \Bigg\} \left. \mathbf{B}\left( s_{0} ;\widetilde{w}_{1,n}\right)\Bigg\} z^n  \right\}  \nonumber
\end{eqnarray}
where
\begin{equation}
\begin{cases} 
{ \displaystyle \overleftrightarrow {\mathbf{\Psi}}_1 \left( s_{1,\infty };t_1,u_1,\eta\right)= \frac{\left( \frac{1+s_{1,\infty }+\sqrt{s_{1,\infty }^2-2(1-2\eta (1-t_1)(1-u_1))s_{1,\infty }+1}}{2}\right)^{-4}}{\sqrt{s_{1,\infty }^2-2(1-2\eta (1-t_1)(1-u_1))s_{1,\infty }+1}} }\cr
{ \displaystyle  \overleftrightarrow {\mathbf{\Psi}}_n \left( s_{n,\infty };t_n,u_n,\eta \right) = \frac{\left( \frac{1+s_{n,\infty }+\sqrt{s_{n,\infty }^2-2(1-2\eta (1-t_n)(1-u_n))s_{n,\infty }+1}}{2}\right)^{-4n }}{\sqrt{s_{n,\infty }^2-2(1-2\eta (1-t_n)(1-u_n))s_{n,\infty }+1}}}\cr
{ \displaystyle \overleftrightarrow {\mathbf{\Psi}}_{n-k} \left( s_{n-k};t_{n-k},u_{n-k},\widetilde{w}_{n-k+1,n} \right) = \frac{\left( \frac{(1+s_{n-k})+\sqrt{s_{n-k}^2-2(1-2\widetilde{w}_{n-k+1,n} (1-t_{n-k})(1-u_{n-k}))s_{n-k}+1}}{2}\right)^{-4(n-k)}}{\sqrt{s_{n-k}^2-2(1-2\widetilde{w}_{n-k+1,n} (1-t_{n-k})(1-u_{n-k}))s_{n-k}+1}}}
\end{cases}\nonumber 
\end{equation}
and
\begin{equation}
\begin{cases} 
{ \displaystyle \mathbf{B} \left( s_{0,\infty } ;\eta\right)= \frac{ \left(1+s_{0,\infty }+\sqrt{s_{0,\infty }^2-2(1-2\eta )s_{0,\infty }+1}\right)^{-\frac{1}{2}}\left(1- s_{0,\infty }+\sqrt{s_{0,\infty }^2-2(1-2\eta )s_{0,\infty }+1}\right)^{\frac{1}{2}}}{\sqrt{s_{0,\infty }^2-2(1-2\eta )s_{0,\infty }+1}}}\cr
{ \displaystyle  \mathbf{B} \left( s_{0} ;\widetilde{w}_{1,1}\right) = \frac{\left(1+s_0+\sqrt{s_0^2-2(1-2\widetilde{w}_{1,1} )s_0+1}\right)^{-\frac{1}{2}}\left(1- s_0+\sqrt{s_0^2-2(1-2\widetilde{w}_{1,1})s_0+1}\right)^{\frac{1}{2}}}{\sqrt{s_0^2-2(1-2\widetilde{w}_{1,1})s_0+1}}} \cr
{ \displaystyle \mathbf{B} \left( s_{0} ;\widetilde{w}_{1,n}\right) = \frac{\left(1+s_0+\sqrt{s_0^2-2(1-2\widetilde{w}_{1,n} )s_0+1}\right)^{-\frac{1}{2}}\left(1- s_0+\sqrt{s_0^2-2(1-2\widetilde{w}_{1,n})s_0+1}\right)^{\frac{1}{2}}}{\sqrt{s_0^2-2(1-2\widetilde{w}_{1,n})s_0+1}}}
\end{cases}\nonumber 
\end{equation}
\end{remark}
\section{Summary}
Lame equation represented either in the algebraic form or in Weierstrass's form is the special case of Heun's differential equation. By changing all coefficients including an independent variable in Heun equation, we obtain the Lame differential equation (ellipsoidal harmonics equation).

In Ref.\cite{aChou2012f} I construct the power series expansion of Lame equation in the algebraic form (for infinite series and polynomial of type 1) and its integral representation by applying 3TRF. In Ref.\cite{aChou2012g} I derive the Frobenius solution and its integral form of Lame equation in Weierstrass's form for infinite series and polynomial of type 1 by changing all coefficients and a variable of Lame equation in the algebraic form in Ref.\cite{aChou2012f}. In Ref.\cite{aChou2012h} I obtain the generating function for the Lame polynomial of type 1 in Weierstrass's form by applying the generating function for the Jacobi polynomial using hypergeometric functions into the general representation in closed form integrals of Lame equation in Ref.\cite{aChou2012g}.

In chapter 8  I build a mathematical structure of the power series expansion of Lame equation in the algebraic form (for infinite series and polynomial of type 2), its integral form and the generating function for the type 2 Lame polynomial by applying R3TRF. 
In this chapter I construct analytic solutions of Lame equation in Weierstrass's form for infinite series and polynomial of type 2 by changing all coefficients and a variable in the power series of Heun equation and its integral form in chapter 2. Also generating functions for the type 2 Lame polynomial in Weierstrass's form of the first and second kinds is obtained by changing all coefficients and a variable of generating functions for the type 2 Heun polynomial of the first and second kinds in chapter 3. 

The Frobenius solutions and its integral forms of Lame equation in Weierstrass's form for infinite series about $\xi =0$ in this chapter are equivalent to the infinite series of Lame equation in Ref.\cite{aChou2012g}. In this chapter $B_n$ is the leading term in sequence $c_n$ in an analytic function $y(\xi )$. In Ref.\cite{aChou2012g} $A_n$ is the leading term in sequence $c_n$ in an analytic function $y(\xi )$.
 
As we see the power series expansions of Heun, Confluent Heun, GCH, Mathieu, Lame equations for either polynomial or infinite series in series ``Special functions and three term recurrence formula (3TRF)''\footnote{It is available as arXiv.} and every chapters in this series, the denominators and numerators in all $A_n$ or $B_n$ terms in these analytic solutions arise with Pochhammer symbol. Since we construct the power series expansions with Pochhammer symbols in numerators and denominators, we are able to describe integral representations of all these equations analytically. The mathematical technique to obtain integral forms of these equations is also available for any linear ordinary differential equations having a 3-term recursive relation in its power series. As we observe representations in closed form integrals of these equations by applying either 3TRF or R3TRF, a $_1F_1$ or $_2F_1$ function recurs in each of sub-integral forms of these equations (every integral forms of all these equations are composed of the definite and contour integrals). 
We might be possible to convert these equations into all other well-known special functions such as Legendre, Bessel, Laguerre, Kummer functions and etc because a $_1F_1$ or $_2F_1$ function recurs in each of sub-integral forms of these equations. 
We can rebuild other confluent forms of these equations after we replace $_1F_1$ or $_2F_1$ function in integral forms of them to other special functions. After that, we can reconstruct the power series expansion of confluent forms of these equations in a backward from those remodeled integral forms. 
Indeed, generating functions of all these equations are constructed analytically by applying the generating function for the Jacobi polynomial using hypergeometric functions or for confluent hypergeometric (Kummer's) polynomial of the first kind into the general expression of the integral representation of these equations. We can build orthogonal relations, recursion relations and expectation values of any physical quantities from these generating functions; i.e. the normalized wave function of hydrogen-like atoms and an expectation value of measurable parameters in quantum physics are constructed by applying the generating function for associated Laguerre polynomials.

\begin{appendices}
\section*{Appendix. Conversion from 9 out of 192 local solutions of Heun equation to 9 local solutions of Lam\'{e} equation in Weierstrass's form for an infinite series and a polynomial of type 2}
  A machine-generated list of 192 (isomorphic to the Coxeter group of the Coxeter diagram $D_4$) local solutions of the Heun equation was obtained by Robert S. Maier(2007) \cite{aMaie2007}. 
In appendix of Ref.\cite{aChou2012d}, I apply 3TRF to series solutions and integrals of Heun equation (for an infinite series and a polynomial of type 1) of nine out of the 192 local solution of Heun equation in Table 2 \cite{aMaie2007}. 

In appendix of chapters 2 and 3, by applying R3TRF, I construct the fundamental power series solutions and integrals of Heun equation (for an infinite series and a polynomial of type 2) of the previous nine out of the 192 local solution of Heun equation in Table 2 \cite{aMaie2007}, including generating functions for the Heun polynomial of type 2. 

In appendix of Ref.\cite{aChou2012g}, by applying 3TRF, 9 local solutions of Lame equation in Weierstrass's form out of 192 local solutions of Heun equation are constructed: Power series solutions and integrals of Lame equation are derived for an infinite series and a polynomial of type 1 analytically.

In this appendix, by changing all coefficients and independent variables of the previous nine examples of 192 local solutions of Heun equation into the first kind of independent solutions of Heun equation by applying R3TRF in chapters 2 and 3, Frobenius solutions and integrals of 9 local solutions of Lame equation in Weierstrass's form are constructed for an infinite series and a polynomial of type 2, including its generating functions for the type 2 polynomial.

Infinite series by applying either 3TRF or R3TRF are equivalent to each other. For infinite series by applying 3TRF, $A_n$ is the leading term in sequence $c_n$ in all 9 local solutions of Lame equation in appendix of Ref.\cite{aChou2012g}. For infinite series by applying R3TRF, $B_n$ is the leading term in sequence $c_n$ in all 9 local solutions of it in this appendix.\footnote{In appendix of Ref.\cite{aChou2012g}, I treat $h$ as a free variable and a fixed value of $\alpha $ to construct polynomials of type 1 for all 9 local solutions of Lame equation. In this appendix, I treat $\alpha $ as a free variable and a fixed value of $h$ to construct polynomials of type 2 for all 9 local solutions of it. And an independent variable $sn^2(z,\rho )$ is denoted by $\xi$.}
\section{Power series}
\addtocontents{toc}{\protect\setcounter{tocdepth}{1}}
\subsection{ ${\displaystyle (1-x)^{1-\delta } Hl(a, q - (\delta  - 1)\gamma a; \alpha - \delta  + 1, \beta - \delta + 1, \gamma ,2 - \delta ; x)}$ }
\subsubsection{Polynomial of type 2}
Replace coefficients $q$, $\alpha$, $\beta$, $\delta$, $c_0$, $\lambda $ and $q_j$ where $j, q_j \in \mathbb{N}_{0}$ by $q - (\delta - 1)\gamma a $, $\alpha - \delta  + 1 $, $\beta - \delta + 1$, $2 - \delta$, 1, zero and $h_j$ where $h_j \in \mathbb{N}_{0}$ into (\ref{eq:8007}). Multiply $(1-x)^{1-\delta }$ and (\ref{eq:8007}) together. Put (\ref{eq:8006}) into the new (\ref{eq:8007}). \footnote{For all 9 local solutions of Lame equation for polynomial of type 2 in this appendix, $h_i\leq h_j$  only if $i\leq j$ where $i,j,h_i,h_j \in \mathbb{N}_{0}$.}
\begin{eqnarray}
&& (1-\xi )^{\frac{1}{2}} y(\xi )\nonumber\\
&=& (1-\xi )^{\frac{1}{2}} Hl\left(\rho ^{-2}, -\rho ^{-2}(h_j+2j)[1+(1+\rho ^2)(h_j+2j)]; \frac{\alpha }{2}+1, -\frac{\alpha }{2}+\frac{1}{2}, \frac{1}{2},\frac{3}{2}; \xi \right) \nonumber\\
&=& (1-\xi )^{\frac{1}{2}} \left\{\sum_{i_0=0}^{h_0} \frac{(-h_0)_{i_0} \left( h_0+ \Omega _{\rho }\right)_{i_0}}{(1)_{i_0}\left(\frac{1}{2} \right)_{i_0}} \eta ^{i_0} \right.  \nonumber\\ 
&+& \left\{ \sum_{i_0=0}^{h_0}\frac{ \left( i_0+1+\frac{\alpha }{2} \right) \left( i_0+\frac{1}{2}-\frac{\alpha }{2} \right)}{(i_0+2)\left( i_0+ \frac{3}{2} \right)}\frac{(-h_0)_{i_0} \left( h_0+ \Omega _{\rho }\right)_{i_0}}{(1)_{i_0}\left( \frac{1}{2}\right)_{i_0}} \right. \left. \sum_{i_1=i_0}^{h_1} \frac{(-h_1)_{i_1}\left( h_1+4+ \Omega _{\rho }\right)_{i_1}(3)_{i_0}\left( \frac{5}{2}\right)_{i_0}}{(-h_1)_{i_0}\left( h_1+4+ \Omega _{\rho }\right)_{i_0}(3)_{i_1} \left( \frac{5}{2}\right)_{i_1}} \eta ^{i_1}\right\} z\nonumber\\
&+& \sum_{n=2}^{\infty } \left\{ \sum_{i_0=0}^{h_0} \frac{\left( i_0+1+\frac{\alpha }{2} \right) \left( i_0+\frac{1}{2}-\frac{\alpha }{2} \right)}{(i_0+ 2)\left( i_0+\frac{3}{2} \right)}\frac{(-h_0)_{i_0} \left( h_0+ \Omega _{\rho }\right)_{i_0}}{(1)_{i_0}\left(\frac{1}{2} \right)_{i_0}}\right.\nonumber\\
&\times& \prod _{k=1}^{n-1} \left\{ \sum_{i_k=i_{k-1}}^{h_k} \frac{\left( i_k+ 2k+1+\frac{\alpha }{2} \right) \left( i_k+ 2k+\frac{1}{2}-\frac{\alpha }{2} \right)}{(i_k+ 2k+2) \left( i_k+ 2k+\frac{3}{2} \right)} \right. \nonumber\\
&\times& \left. \frac{(-h_k)_{i_k}\left( h_k+4k+ \Omega _{\rho }\right)_{i_k}(2k+1)_{i_{k-1}}\left( 2k+\frac{1}{2} \right)_{i_{k-1}}}{(-h_k)_{i_{k-1}}\left(  h_k+4k+ \Omega _{\rho }\right)_{i_{k-1}}(2k+1)_{i_k}\left( 2k+\frac{1}{2} \right)_{i_k}}\right\} \nonumber\\
&\times& \left.\left. \sum_{i_n= i_{n-1}}^{h_n} \frac{(-h_n)_{i_n}\left( h_n+4n+ \Omega _{\rho }\right)_{i_n}(2n+1)_{i_{n-1}}\left( 2n+\frac{1}{2} \right)_{i_{n-1}}}{(-h_n)_{i_{n-1}}\left(  h_n+4n+ \Omega _{\rho }\right)_{i_{n-1}}(2n+1)_{i_n}\left( 2n+\frac{1}{2} \right)_{i_n}} \eta ^{i_n} \right\} z^n \right\} \label{eq:80023}
\end{eqnarray}
where
\begin{equation}
\begin{cases} 
\Omega _{\rho } = \frac{1}{1+\rho ^2}\cr
h= 4(h_j+2j)[1+(1+\rho ^2)(h_j+2j)]+1  
\end{cases}\nonumber 
\end{equation}
\subsubsection{Infinite series}
Replace coefficients $q$, $\alpha$, $\beta$, $\delta$, $c_0$ and $\lambda $ by $q - (\delta - 1)\gamma a $, $\alpha - \delta  + 1 $, $\beta - \delta + 1$, $2-\delta$, 1 and zero into (\ref{eq:80010}). Multiply $(1-x)^{1-\delta }$ and (\ref{eq:80010}) together. Put (\ref{eq:8006}) into the new (\ref{eq:80010}).
\begin{eqnarray}
&& (1-\xi )^{\frac{1}{2}} y(\xi )\nonumber\\
&=& (1-\xi )^{\frac{1}{2}} Hl\left(\rho ^{-2}, -\frac{1}{4}(h-1)\rho ^{-2}; \frac{\alpha }{2}+1, -\frac{\alpha }{2}+\frac{1}{2}, \frac{1}{2},\frac{3}{2}; \xi \right) \nonumber\\
&=& (1-\xi )^{\frac{1}{2}} \left\{\sum_{i_0=0}^{\infty } \frac{\left(\Delta_0^{-}\right)_{i_0} \left(\Delta_0^{+}\right)_{i_0}}{(1 )_{i_0}\left( \frac{1}{2} \right)_{i_0}} \eta ^{i_0}\right.\nonumber\\
&+& \left\{ \sum_{i_0=0}^{\infty }\frac{\left( i_0 +1+\frac{\alpha }{2} \right) \left( i_0 +\frac{1}{2}-\frac{\alpha }{2} \right)}{(i_0 +2)\left( i_0 +\frac{3}{2} \right)}\frac{\left(\Delta_0^{-}\right)_{i_0} \left(\Delta_0^{+}\right)_{i_0}}{(1)_{i_0}\left(\frac{1}{2} \right)_{i_0}} \sum_{i_1=i_0}^{\infty } \frac{\left(\Delta_1^{-}\right)_{i_1} \left(\Delta_1^{+}\right)_{i_1}(3 )_{i_0}\left( \frac{5}{2} \right)_{i_0}}{\left(\Delta_1^{-}\right)_{i_0}  \left(\Delta_1^{+}\right)_{i_0}(3 )_{i_1}\left( \frac{5}{2} \right)_{i_1}}\eta ^{i_1}\right\} z\nonumber\\
&+& \sum_{n=2}^{\infty } \left\{ \sum_{i_0=0}^{\infty } \frac{\left( i_0 +1+\frac{\alpha }{2} \right) \left( i_0 +\frac{1}{2}-\frac{\alpha }{2} \right)}{(i_0 +2)\left( i_0 +\frac{3}{2} \right)}\frac{\left(\Delta_0^{-}\right)_{i_0} \left(\Delta_0^{+}\right)_{i_0}}{(1 )_{i_0}\left(\frac{1}{2} \right)_{i_0}}\right.\nonumber\\
&\times& \prod _{k=1}^{n-1} \left\{ \sum_{i_k=i_{k-1}}^{\infty } \frac{\left( i_k+ 2k+1+\frac{\alpha }{2} \right) \left( i_k+ 2k+\frac{1}{2}-\frac{\alpha }{2} \right)}{(i_k+ 2k+2 )\left( i_k+ 2k+\frac{3}{2} \right)} \frac{ \left(\Delta_k^{-}\right)_{i_k} \left(\Delta_k^{+} \right)_{i_k}(2k+1 )_{i_{k-1}}\left( 2k+\frac{1}{2} \right)_{i_{k-1}}}{\left(\Delta_k^{-}\right)_{i_{k-1}} \left(\Delta_k^{+} \right)_{i_{k-1}}(2k+1 )_{i_k}\left( 2k+\frac{1}{2} \right)_{i_k}}\right\}\nonumber\\
&\times& \left.\left.\sum_{i_n= i_{n-1}}^{\infty } \frac{\left(\Delta_n^{-}\right)_{i_n}\left( \Delta_n^{+} \right)_{i_n}(2n+1 )_{i_{n-1}}\left( 2n+\frac{1}{2} \right)_{i_{n-1}}}{\left(\Delta_n^{-}\right)_{i_{n-1}}\left(\Delta_n^{+} \right)_{i_{n-1}}(2n+1 )_{i_n}\left( 2n+\frac{1}{2} \right)_{i_n}} \eta ^{i_n} \right\} z^n \right\} \label{eq:80024}
\end{eqnarray}
where
\begin{equation}
\begin{cases} 
\Delta_0^{\pm}= \frac{1 \pm\sqrt{1+(1+\rho ^2)(h-1)}}{2(1+\rho ^2)} \cr
\Delta_1^{\pm}= 2 + \frac{1 \pm\sqrt{1+(1+\rho ^2)(h-1)}}{2(1+\rho ^2)} \cr
\Delta_k^{\pm}= 2k + \frac{1 \pm\sqrt{1+(1+\rho ^2)(h-1)}}{2(1+\rho ^2)} \cr
\Delta_n^{\pm}= 2n +  \frac{1 \pm\sqrt{1+(1+\rho ^2)(h-1)}}{2(1+\rho ^2)}
\end{cases}\nonumber 
\end{equation}
On (\ref{eq:80023}) and (\ref{eq:80024}),
\begin{equation}
\begin{cases} 
\eta =(1+\rho ^2)\xi \cr
z=-\rho ^2\xi^2 
\end{cases}\nonumber 
\end{equation}
\subsection{ \footnotesize ${\displaystyle x^{1-\gamma } (1-x)^{1-\delta }Hl(a, q-(\gamma +\delta -2)a -(\gamma -1)(\alpha +\beta -\gamma -\delta +1), \alpha - \gamma -\delta +2, \beta - \gamma -\delta +2, 2-\gamma, 2 - \delta ; x)}$ \normalsize}
\subsubsection{Polynomial of type 2}
Replace coefficients $q$, $\alpha$, $\beta$, $\gamma $, $\delta$, $c_0$, $\lambda $ and $q_j$ where $j, q_j \in \mathbb{N}_{0}$ by $q-(\gamma +\delta -2)a-(\gamma -1)(\alpha +\beta -\gamma -\delta +1)$, $\alpha - \gamma -\delta +2$, $\beta - \gamma -\delta +2, 2-\gamma$, $2 - \delta$, 1, zero and $h_j$ where $h_j \in \mathbb{N}_{0}$ into (\ref{eq:8007}). Multiply $x^{1-\gamma } (1-x)^{1-\delta }$ and (\ref{eq:8007}) together. Put (\ref{eq:8006}) into the new (\ref{eq:8007}). 
\begin{eqnarray}
&&\xi ^{\frac{1}{2}} (1-\xi )^{\frac{1}{2}} y(\xi )\nonumber\\
&=&\xi ^{\frac{1}{2}} (1-\xi )^{\frac{1}{2}} Hl\left(\rho ^{-2}, -\rho ^{-2}(h_j+2j)[2+\rho ^2 +(1+\rho ^2)(h_j+2j)]; \frac{\alpha }{2}+\frac{3}{2}, -\frac{\alpha }{2}+1, \frac{3}{2},\frac{3}{2}; \xi \right)\nonumber\\
&=& \xi ^{\frac{1}{2}} (1-\xi )^{\frac{1}{2}} \left\{ \sum_{i_0=0}^{h_0} \frac{(-h_0)_{i_0} \left( h_0+ \Omega _{\rho }\right)_{i_0}}{(1)_{i_0}\left(\frac{3}{2} \right)_{i_0}} \eta ^{i_0} \right.  \nonumber\\ 
&+& \left\{ \sum_{i_0=0}^{h_0}\frac{ \left( i_0 +\frac{3}{2}+\frac{\alpha }{2} \right) \left( i_0+1-\frac{\alpha }{2} \right)}{(i_0+2)\left( i_0+ \frac{5}{2} \right)}\frac{(-h_0)_{i_0} \left( h_0+ \Omega _{\rho }\right)_{i_0}}{(1)_{i_0}\left( \frac{3}{2}\right)_{i_0}} \right. \left. \sum_{i_1=i_0}^{h_1} \frac{(-h_1)_{i_1}\left( h_1+4+ \Omega _{\rho }\right)_{i_1}(3)_{i_0}\left( \frac{7}{2}\right)_{i_0}}{(-h_1)_{i_0}\left( h_1+4+ \Omega _{\rho }\right)_{i_0}(3)_{i_1} \left( \frac{7}{2}\right)_{i_1}} \eta ^{i_1}\right\} z\nonumber\\
&+& \sum_{n=2}^{\infty } \left\{ \sum_{i_0=0}^{h_0} \frac{\left( i_0+\frac{3}{2}+\frac{\alpha }{2} \right) \left( i_0+1-\frac{\alpha }{2} \right)}{(i_0+ 2)\left( i_0+\frac{5}{2} \right)}\frac{(-h_0)_{i_0} \left( h_0+ \Omega _{\rho } \right)_{i_0}}{(1)_{i_0}\left(\frac{3}{2} \right)_{i_0}}\right.\nonumber\\
&\times& \prod _{k=1}^{n-1} \left\{ \sum_{i_k=i_{k-1}}^{h_k} \frac{\left( i_k+ 2k+\frac{3}{2}+\frac{\alpha }{2} \right) \left( i_k+ 2k+1-\frac{\alpha }{2} \right)}{(i_k+ 2k+2) \left( i_k+ 2k+\frac{5}{2} \right)} \right. \nonumber\\
&\times& \left. \frac{(-h_k)_{i_k}\left( h_k+4k+ \Omega _{\rho }\right)_{i_k}(2k+1)_{i_{k-1}}\left( 2k+\frac{3}{2} \right)_{i_{k-1}}}{(-h_k)_{i_{k-1}}\left(  h_k+4k+ \Omega _{\rho }\right)_{i_{k-1}}(2k+1)_{i_k}\left( 2k+\frac{3}{2} \right)_{i_k}}\right\} \nonumber\\
&\times& \left.\left. \sum_{i_n= i_{n-1}}^{h_n} \frac{(-h_n)_{i_n}\left( h_n+4n+ \Omega _{\rho } \right)_{i_n}(2n+1)_{i_{n-1}}\left( 2n+\frac{3}{2} \right)_{i_{n-1}}}{(-h_n)_{i_{n-1}}\left(  h_n +4n+ \Omega _{\rho }\right)_{i_{n-1}}(2n+1)_{i_n}\left( 2n+\frac{3}{2} \right)_{i_n}} \eta ^{i_n} \right\} z^n \right\} \label{eq:80025}
\end{eqnarray}
where
\begin{equation}
\begin{cases} 
\Omega _{\rho } = \frac{2+\rho ^2}{1+\rho ^2}\cr
h= 4+\rho ^2+ 4(h_j+2j)[2+\rho ^2+(1+\rho ^2)(h_j+2j)]  
\end{cases}\nonumber 
\end{equation}
\subsubsection{Infinite series}
Replace coefficients $q$, $\alpha$, $\beta$, $\gamma $, $\delta$, $c_0$ and $\lambda $ by $q-(\gamma +\delta -2)a-(\gamma -1)(\alpha +\beta -\gamma -\delta +1)$, $\alpha - \gamma -\delta +2$, $\beta - \gamma -\delta +2, 2-\gamma$, $2 - \delta$, 1 and zero into (\ref{eq:80010}). Multiply $x^{1-\gamma } (1-x)^{1-\delta }$ and (\ref{eq:80010}) together. Put (\ref{eq:8006}) into the new (\ref{eq:80010}). 
\begin{eqnarray}
&&\xi ^{\frac{1}{2}} (1-\xi )^{\frac{1}{2}} y(\xi )\nonumber\\
&=&\xi ^{\frac{1}{2}} (1-\xi )^{\frac{1}{2}} Hl\left(\rho ^{-2}, -\frac{1}{4}\left( (h-4)\rho ^{-2}-1\right); \frac{\alpha }{2}+\frac{3}{2}, -\frac{\alpha }{2}+1, \frac{3}{2},\frac{3}{2}; \xi \right)\nonumber\\
&=& \xi ^{\frac{1}{2}} (1-\xi )^{\frac{1}{2}} \left\{ \sum_{i_0=0}^{\infty } \frac{\left(\Delta_0^{-}\right)_{i_0} \left(\Delta_0^{+}\right)_{i_0}}{(1 )_{i_0}\left( \frac{3}{2} \right)_{i_0}} \eta ^{i_0}\right.\nonumber\\
&+& \left\{ \sum_{i_0=0}^{\infty }\frac{\left( i_0 +\frac{3}{2}+\frac{\alpha }{2} \right) \left( i_0 +1-\frac{\alpha }{2} \right)}{(i_0 +2)\left( i_0 +\frac{5}{2} \right)}\frac{\left(\Delta_0^{-}\right)_{i_0} \left(\Delta_0^{+}\right)_{i_0}}{(1)_{i_0}\left(\frac{3}{2} \right)_{i_0}} \sum_{i_1=i_0}^{\infty } \frac{\left(\Delta_1^{-}\right)_{i_1} \left(\Delta_1^{+}\right)_{i_1}(3 )_{i_0}\left( \frac{7}{2} \right)_{i_0}}{\left(\Delta_1^{-}\right)_{i_0}  \left(\Delta_1^{+}\right)_{i_0}(3 )_{i_1}\left( \frac{7}{2} \right)_{i_1}}\eta ^{i_1}\right\} z\nonumber\\
&+& \sum_{n=2}^{\infty } \left\{ \sum_{i_0=0}^{\infty } \frac{\left( i_0 +\frac{3}{2}+\frac{\alpha }{2} \right) \left( i_0 +1-\frac{\alpha }{2} \right)}{(i_0 +2)\left( i_0 +\frac{5}{2} \right)}\frac{\left(\Delta_0^{-}\right)_{i_0} \left(\Delta_0^{+}\right)_{i_0}}{(1 )_{i_0}\left(\frac{3}{2} \right)_{i_0}}\right.\nonumber\\
&\times& \prod _{k=1}^{n-1} \left\{ \sum_{i_k=i_{k-1}}^{\infty } \frac{\left( i_k+ 2k+\frac{3}{2}+\frac{\alpha }{2} \right) \left( i_k+ 2k+1-\frac{\alpha }{2} \right)}{(i_k+ 2k+2 )\left( i_k+ 2k+\frac{5}{2} \right)} \frac{ \left(\Delta_k^{-}\right)_{i_k} \left(\Delta_k^{+} \right)_{i_k}(2k+1 )_{i_{k-1}}\left( 2k+\frac{3}{2} \right)_{i_{k-1}}}{\left(\Delta_k^{-}\right)_{i_{k-1}} \left(\Delta_k^{+} \right)_{i_{k-1}}(2k+1 )_{i_k}\left( 2k+\frac{3}{2} \right)_{i_k}}\right\}\nonumber\\
&\times& \left.\left.\sum_{i_n= i_{n-1}}^{\infty } \frac{\left(\Delta_n^{-}\right)_{i_n}\left( \Delta_n^{+} \right)_{i_n}(2n+1 )_{i_{n-1}}\left( 2n+\frac{3}{2} \right)_{i_{n-1}}}{\left(\Delta_n^{-}\right)_{i_{n-1}}\left(\Delta_n^{+} \right)_{i_{n-1}}(2n+1 )_{i_n}\left( 2n+\frac{3}{2} \right)_{i_n}} \eta ^{i_n} \right\} z^n \right\} \label{eq:80026}
\end{eqnarray}
where
\begin{equation}
\begin{cases} 
\Delta_0^{\pm}= \frac{2+\rho ^2 \pm\sqrt{ (h-1)\rho ^2+h }}{2(1+\rho ^2)} \cr
\Delta_1^{\pm}= 2 + \frac{2+\rho ^2 \pm\sqrt{ (h-1)\rho ^2+h }}{2(1+\rho ^2)} \cr
\Delta_k^{\pm}= 2k + \frac{2+\rho ^2 \pm\sqrt{ (h-1)\rho ^2+h }}{2(1+\rho ^2)} \cr
\Delta_n^{\pm}= 2n +  \frac{2+\rho ^2 \pm\sqrt{ (h-1)\rho ^2+h }}{2(1+\rho ^2)}
\end{cases}\nonumber 
\end{equation}
On (\ref{eq:80025}) and (\ref{eq:80026}),
\begin{equation}
\begin{cases} 
\eta =(1+\rho ^2)\xi \cr
z=-\rho ^2\xi^2 
\end{cases}\nonumber 
\end{equation}
\subsection{ ${\displaystyle  Hl(1-a,-q+\alpha \beta; \alpha,\beta, \delta, \gamma; 1-x)}$} 
\subsubsection{Polynomial of type 2}
Replace coefficients $a$, $q$, $\gamma $, $\delta$, $x$, $c_0$, $\lambda $ and $q_j$ where $j, q_j \in \mathbb{N}_{0}$ by $1-a$, $-q +\alpha \beta $, $\delta $, $\gamma $, $1-x$, 1, zero and $h_j$ where $h_j \in \mathbb{N}_{0}$ into (\ref{eq:8007}). Put (\ref{eq:8006}) into the new (\ref{eq:8007}). 
\begin{eqnarray}
y(\varsigma )&=&  Hl\left( 1-\rho ^{-2}, -(2-\rho ^{-2})(h_j +2j)^2; \frac{\alpha }{2}+\frac{1}{2}, -\frac{\alpha }{2}, \frac{1}{2}, \frac{1}{2}; \varsigma \right)\nonumber\\
&=&  \sum_{i_0=0}^{h_0} \frac{(-h_0)_{i_0} \left( h_0 \right)_{i_0}}{(1)_{i_0}\left(\frac{1}{2} \right)_{i_0}} \eta ^{i_0}  \nonumber\\ 
&+& \left\{ \sum_{i_0=0}^{h_0}\frac{ \left( i_0 +\frac{1}{2}+\frac{\alpha }{2} \right) \left( i_0 -\frac{\alpha }{2} \right)}{(i_0+2)\left( i_0+ \frac{3}{2} \right)}\frac{(-h_0)_{i_0} \left( h_0 \right)_{i_0}}{(1)_{i_0}\left( \frac{1}{2}\right)_{i_0}} \right. \left. \sum_{i_1=i_0}^{h_1} \frac{(-h_1)_{i_1}\left( h_1+4 \right)_{i_1}(3)_{i_0}\left( \frac{5}{2}\right)_{i_0}}{(-h_1)_{i_0}\left( h_1+4 \right)_{i_0}(3)_{i_1} \left( \frac{5}{2}\right)_{i_1}} \eta ^{i_1}\right\} z\nonumber\\
&+& \sum_{n=2}^{\infty } \left\{ \sum_{i_0=0}^{h_0} \frac{\left( i_0+\frac{1}{2}+\frac{\alpha }{2} \right) \left( i_0 -\frac{\alpha }{2} \right)}{(i_0+ 2)\left( i_0+\frac{3}{2} \right)}\frac{(-h_0)_{i_0} \left( h_0 \right)_{i_0}}{(1)_{i_0}\left(\frac{1}{2} \right)_{i_0}}\right.\nonumber\\
&\times& \prod _{k=1}^{n-1} \left\{ \sum_{i_k=i_{k-1}}^{h_k} \frac{\left( i_k+ 2k+\frac{1}{2}+\frac{\alpha }{2} \right) \left( i_k+ 2k -\frac{\alpha }{2} \right)}{(i_k+ 2k+2) \left( i_k+ 2k+\frac{3}{2} \right)} 
 \frac{(-h_k)_{i_k}\left( h_k+4k \right)_{i_k}(2k+1)_{i_{k-1}}\left( 2k+\frac{1}{2} \right)_{i_{k-1}}}{(-h_k)_{i_{k-1}}\left(  h_k+4k \right)_{i_{k-1}}(2k+1)_{i_k}\left( 2k+\frac{1}{2} \right)_{i_k}}\right\} \nonumber\\
&\times& \left. \sum_{i_n= i_{n-1}}^{h_n} \frac{(-h_n)_{i_n}\left( h_n+4n \right)_{i_n}(2n+1)_{i_{n-1}}\left( 2n+\frac{1}{2} \right)_{i_{n-1}}}{(-h_n)_{i_{n-1}}\left(  h_n +4n \right)_{i_{n-1}}(2n+1)_{i_n}\left( 2n+\frac{1}{2} \right)_{i_n}} \eta ^{i_n} \right\} z^n  \label{eq:80027}
\end{eqnarray}
where
\begin{equation}
h= \rho ^2 \left( \alpha (\alpha +1)- 4(2-\rho ^{-2})(h_j+2j)^2 \right)  \nonumber 
\end{equation}
\subsubsection{Infinite series}
Replace coefficients $a$, $q$, $\gamma $, $\delta$, $x$, $c_0$ and $\lambda $ by $1-a$, $-q +\alpha \beta $, $\delta $, $\gamma $, $1-x$, 1 and zero into (\ref{eq:80010}). Put (\ref{eq:8006}) into the new (\ref{eq:80010}). 
\begin{eqnarray}
y(\varsigma )&=&  Hl\left( 1-\rho ^{-2}, \frac{1}{4}\left( h\rho ^{-2}- \alpha (\alpha +1)\right); \frac{\alpha }{2}+\frac{1}{2}, -\frac{\alpha }{2}, \frac{1}{2}, \frac{1}{2}; \varsigma \right)\nonumber\\
&=& \sum_{i_0=0}^{\infty } \frac{\left(\Delta_0^{-}\right)_{i_0} \left(\Delta_0^{+}\right)_{i_0}}{(1 )_{i_0}\left( \frac{1}{2} \right)_{i_0}} \eta ^{i_0} \nonumber\\
&+& \left\{ \sum_{i_0=0}^{\infty }\frac{\left( i_0 +\frac{1}{2}+\frac{\alpha }{2} \right) \left( i_0 -\frac{\alpha }{2} \right)}{(i_0 +2)\left( i_0 +\frac{3}{2} \right)}\frac{\left(\Delta_0^{-}\right)_{i_0} \left(\Delta_0^{+}\right)_{i_0}}{(1)_{i_0}\left(\frac{1}{2} \right)_{i_0}} \sum_{i_1=i_0}^{\infty } \frac{\left(\Delta_1^{-}\right)_{i_1} \left(\Delta_1^{+}\right)_{i_1}(3 )_{i_0}\left( \frac{5}{2} \right)_{i_0}}{\left(\Delta_1^{-}\right)_{i_0}  \left(\Delta_1^{+}\right)_{i_0}(3 )_{i_1}\left( \frac{5}{2} \right)_{i_1}}\eta ^{i_1}\right\} z\nonumber\\
&+& \sum_{n=2}^{\infty } \left\{ \sum_{i_0=0}^{\infty } \frac{\left( i_0 +\frac{1}{2}+\frac{\alpha }{2} \right) \left( i_0 -\frac{\alpha }{2} \right)}{(i_0 +2)\left( i_0 +\frac{3}{2} \right)}\frac{\left(\Delta_0^{-}\right)_{i_0} \left(\Delta_0^{+}\right)_{i_0}}{(1 )_{i_0}\left(\frac{1}{2} \right)_{i_0}}\right.\nonumber\\
&\times& \prod _{k=1}^{n-1} \left\{ \sum_{i_k=i_{k-1}}^{\infty } \frac{\left( i_k+ 2k+\frac{1}{2}+\frac{\alpha }{2} \right) \left( i_k+ 2k -\frac{\alpha }{2} \right)}{(i_k+ 2k+2 )\left( i_k+ 2k+\frac{3}{2} \right)} \frac{ \left(\Delta_k^{-}\right)_{i_k} \left(\Delta_k^{+} \right)_{i_k}(2k+1 )_{i_{k-1}}\left( 2k+\frac{1}{2} \right)_{i_{k-1}}}{\left(\Delta_k^{-}\right)_{i_{k-1}} \left(\Delta_k^{+} \right)_{i_{k-1}}(2k+1 )_{i_k}\left( 2k+\frac{1}{2} \right)_{i_k}}\right\}\nonumber\\
&\times& \left. \sum_{i_n= i_{n-1}}^{\infty } \frac{\left(\Delta_n^{-}\right)_{i_n}\left( \Delta_n^{+} \right)_{i_n}(2n+1 )_{i_{n-1}}\left( 2n+\frac{1}{2} \right)_{i_{n-1}}}{\left(\Delta_n^{-}\right)_{i_{n-1}}\left(\Delta_n^{+} \right)_{i_{n-1}}(2n+1 )_{i_n}\left( 2n+\frac{1}{2} \right)_{i_n}} \eta ^{i_n} \right\} z^n \label{eq:80028}
\end{eqnarray}
where
\begin{equation}
\begin{cases} 
\Delta_0^{\pm}= \pm \frac{1}{2}\sqrt{\frac{ h\rho ^{-2} +\alpha (\alpha +1) }{ 2-\rho ^{-2} }} \cr
\Delta_1^{\pm}= 2 \pm \frac{1}{2}\sqrt{\frac{ h\rho ^{-2} +\alpha (\alpha +1) }{ 2-\rho ^{-2} }} \cr
\Delta_k^{\pm}= 2k \pm \frac{1}{2}\sqrt{\frac{ h\rho ^{-2} +\alpha (\alpha +1) }{ 2-\rho ^{-2} }} \cr
\Delta_n^{\pm}= 2n \pm \frac{1}{2}\sqrt{\frac{ h\rho ^{-2} +\alpha (\alpha +1) }{ 2-\rho ^{-2} }}
\end{cases}\nonumber 
\end{equation}
On (\ref{eq:80027}) and (\ref{eq:80028}),
\begin{equation}
\begin{cases} 
\varsigma= 1-\xi \cr
\eta =\frac{2-\rho ^{-2}}{1-\rho ^{-2}}\varsigma \cr
z=\frac{-1}{1-\rho ^{-2}}\varsigma ^2  
\end{cases}\nonumber 
\end{equation}
\subsection{ \footnotesize ${\displaystyle (1-x)^{1-\delta } Hl(1-a,-q+(\delta -1)\gamma a+(\alpha -\delta +1)(\beta -\delta +1); \alpha-\delta +1,\beta-\delta +1, 2-\delta, \gamma; 1-x)}$ \normalsize}
\subsubsection{Polynomial of type 2}
Replace coefficients $a$, $q$, $\alpha $, $\beta $, $\gamma $, $\delta$, $x$, $c_0$, $\lambda $ and $q_j$ where $j, q_j \in \mathbb{N}_{0}$ by $1-a$, $-q+(\delta -1)\gamma a+(\alpha -\delta +1)(\beta -\delta +1)$, $\alpha-\delta +1 $, $\beta-\delta +1 $, $2-\delta$, $\gamma $, $1-x$, 1, zero and $h_j$ where $h_j \in \mathbb{N}_{0}$ into (\ref{eq:8007}). Multiply $(1-x)^{1-\delta }$ and (\ref{eq:8007}) together. Put (\ref{eq:8006}) into the new (\ref{eq:8007}). 
\begin{eqnarray}
&&\varsigma ^{\frac{1}{2}}y(\varsigma )\nonumber\\
&=&\varsigma ^{\frac{1}{2}} Hl\left( 1-\rho ^{-2}, -(2-\rho ^{-2})(h_j+2j)(h_j+2j+1); \frac{\alpha }{2}+1, -\frac{\alpha }{2}+\frac{1}{2}, \frac{3}{2}, \frac{1}{2}; \varsigma \right)\nonumber\\
&=& \varsigma ^{\frac{1}{2}} \left\{ \sum_{i_0=0}^{h_0} \frac{(-h_0)_{i_0} \left( h_0+ 1\right)_{i_0}}{(1)_{i_0}\left(\frac{3}{2} \right)_{i_0}} \eta ^{i_0} \right. \nonumber\\
&+& \left\{ \sum_{i_0=0}^{h_0}\frac{ \left( i_0 +1+\frac{\alpha }{2} \right) \left( i_0+\frac{1}{2}-\frac{\alpha }{2} \right)}{(i_0+2)\left( i_0+ \frac{5}{2} \right)}\frac{(-h_0)_{i_0} \left( h_0+ 1\right)_{i_0}}{(1)_{i_0}\left( \frac{3}{2}\right)_{i_0}} \right. \left. \sum_{i_1=i_0}^{h_1} \frac{(-h_1)_{i_1}\left( h_1+5\right)_{i_1}(3)_{i_0}\left( \frac{7}{2}\right)_{i_0}}{(-h_1)_{i_0}\left( h_1+5\right)_{i_0}(3)_{i_1} \left( \frac{7}{2}\right)_{i_1}} \eta ^{i_1}\right\} z\nonumber\\
&+& \sum_{n=2}^{\infty } \left\{ \sum_{i_0=0}^{h_0} \frac{\left( i_0+1+\frac{\alpha }{2} \right) \left( i_0+\frac{1}{2}-\frac{\alpha }{2} \right)}{(i_0+ 2)\left( i_0+\frac{5}{2} \right)}\frac{(-h_0)_{i_0} \left( h_0+1 \right)_{i_0}}{(1)_{i_0}\left(\frac{3}{2} \right)_{i_0}}\right.\nonumber\\
&\times& \prod _{k=1}^{n-1} \left\{ \sum_{i_k=i_{k-1}}^{h_k} \frac{\left( i_k+ 2k+1+\frac{\alpha }{2} \right) \left( i_k+ 2k+\frac{1}{2}-\frac{\alpha }{2} \right)}{(i_k+ 2k+2) \left( i_k+ 2k+\frac{5}{2} \right)}\right. \nonumber\\ 
&\times& \left. \frac{(-h_k)_{i_k}\left( h_k+4k+ 1\right)_{i_k}(2k+1)_{i_{k-1}}\left( 2k+\frac{3}{2} \right)_{i_{k-1}}}{(-h_k)_{i_{k-1}}\left(  h_k+4k+1\right)_{i_{k-1}}(2k+1)_{i_k}\left( 2k+\frac{3}{2} \right)_{i_k}}\right\} \nonumber\\
&\times& \left.\left. \sum_{i_n= i_{n-1}}^{h_n} \frac{(-h_n)_{i_n}\left( h_n+4n+1 \right)_{i_n}(2n+1)_{i_{n-1}}\left( 2n+\frac{3}{2} \right)_{i_{n-1}}}{(-h_n)_{i_{n-1}}\left(  h_n +4n+ 1\right)_{i_{n-1}}(2n+1)_{i_n}\left( 2n+\frac{3}{2} \right)_{i_n}} \eta ^{i_n} \right\} z^n \right\} \label{eq:80029}
\end{eqnarray}
where
\begin{equation} 
h= 4(2\rho ^2-1)(h_j +2j)(h_j +2j+1)-\rho ^2(\alpha -1)(\alpha +2) -1  \nonumber 
\end{equation}
\subsubsection{Infinite series}
Replace coefficients $a$, $q$, $\alpha $, $\beta $, $\gamma $, $\delta$, $x$, $c_0$ and $\lambda $ by $1-a$, $-q+(\delta -1)\gamma a+(\alpha -\delta +1)(\beta -\delta +1)$, $\alpha-\delta +1 $, $\beta-\delta +1 $, $2-\delta$, $\gamma $, $1-x$, 1 and zero into (\ref{eq:80010}). Multiply $(1-x)^{1-\delta }$ and (\ref{eq:80010}) together. Put (\ref{eq:8006}) into the new (\ref{eq:80010}).
\begin{eqnarray}
&&\varsigma ^{\frac{1}{2}}y(\varsigma )\nonumber\\
&=& \varsigma ^{\frac{1}{2}} Hl\left( 1-\rho ^{-2}, -\frac{1}{4}\left( (1+h)\rho ^{-2}+(\alpha -1)(\alpha +2)\right); \frac{\alpha }{2}+1, -\frac{\alpha }{2}+\frac{1}{2}, \frac{3}{2}, \frac{1}{2}; \varsigma \right)\nonumber\\
&=& \varsigma ^{\frac{1}{2}} \left\{ \sum_{i_0=0}^{\infty } \frac{\left(\Delta_0^{-}\right)_{i_0} \left(\Delta_0^{+}\right)_{i_0}}{(1 )_{i_0}\left( \frac{3}{2} \right)_{i_0}} \eta ^{i_0} \right. \nonumber\\
&+& \left\{ \sum_{i_0=0}^{\infty }\frac{\left( i_0 +1+\frac{\alpha }{2} \right) \left( i_0 +\frac{1}{2}-\frac{\alpha }{2} \right)}{(i_0 +2)\left( i_0 +\frac{5}{2} \right)}\frac{\left(\Delta_0^{-}\right)_{i_0} \left(\Delta_0^{+}\right)_{i_0}}{(1)_{i_0}\left(\frac{3}{2} \right)_{i_0}} \sum_{i_1=i_0}^{\infty } \frac{\left(\Delta_1^{-}\right)_{i_1} \left(\Delta_1^{+}\right)_{i_1}(3 )_{i_0}\left( \frac{7}{2} \right)_{i_0}}{\left(\Delta_1^{-}\right)_{i_0}  \left(\Delta_1^{+}\right)_{i_0}(3 )_{i_1}\left( \frac{7}{2} \right)_{i_1}}\eta ^{i_1}\right\} z\nonumber\\
&+& \sum_{n=2}^{\infty } \left\{ \sum_{i_0=0}^{\infty } \frac{\left( i_0 +1+\frac{\alpha }{2} \right) \left( i_0 +\frac{1}{2}-\frac{\alpha }{2} \right)}{(i_0 +2)\left( i_0 +\frac{5}{2} \right)}\frac{\left(\Delta_0^{-}\right)_{i_0} \left(\Delta_0^{+}\right)_{i_0}}{(1 )_{i_0}\left(\frac{3}{2} \right)_{i_0}}\right.\nonumber\\
&\times& \prod _{k=1}^{n-1} \left\{ \sum_{i_k=i_{k-1}}^{\infty } \frac{\left( i_k+ 2k+1+\frac{\alpha }{2} \right) \left( i_k+ 2k +\frac{1}{2}-\frac{\alpha }{2} \right)}{(i_k+ 2k+2 )\left( i_k+ 2k+\frac{5}{2} \right)} \frac{ \left(\Delta_k^{-}\right)_{i_k} \left(\Delta_k^{+} \right)_{i_k}(2k+1 )_{i_{k-1}}\left( 2k+\frac{3}{2} \right)_{i_{k-1}}}{\left(\Delta_k^{-}\right)_{i_{k-1}} \left(\Delta_k^{+} \right)_{i_{k-1}}(2k+1 )_{i_k}\left( 2k+\frac{3}{2} \right)_{i_k}}\right\}\nonumber\\
&\times& \left. \left. \sum_{i_n= i_{n-1}}^{\infty } \frac{\left(\Delta_n^{-}\right)_{i_n}\left( \Delta_n^{+} \right)_{i_n}(2n+1 )_{i_{n-1}}\left( 2n+\frac{3}{2} \right)_{i_{n-1}}}{\left(\Delta_n^{-}\right)_{i_{n-1}}\left(\Delta_n^{+} \right)_{i_{n-1}}(2n+1 )_{i_n}\left( 2n+\frac{3}{2} \right)_{i_n}} \eta ^{i_n} \right\} z^n \right\}\label{eq:80030}
\end{eqnarray}
where
\begin{equation}
\begin{cases} 
\Delta_0^{\pm}= \frac{1}{2} \pm \frac{1}{2}\sqrt{\frac{ h\rho ^{-2} +\alpha (\alpha +1) }{ 2-\rho ^{-2} }} \cr
\Delta_1^{\pm}= \frac{5}{2} \pm \frac{1}{2}\sqrt{\frac{ h\rho ^{-2} +\alpha (\alpha +1) }{ 2-\rho ^{-2} }} \cr
\Delta_k^{\pm}= 2k + \frac{1}{2} \pm \frac{1}{2}\sqrt{\frac{ h\rho ^{-2} +\alpha (\alpha +1) }{ 2-\rho ^{-2} }} \cr
\Delta_n^{\pm}= 2n +\frac{1}{2} \pm \frac{1}{2}\sqrt{\frac{ h\rho ^{-2} +\alpha (\alpha +1) }{ 2-\rho ^{-2} }}
\end{cases}\nonumber 
\end{equation}
On (\ref{eq:80029}) and (\ref{eq:80030}),
\begin{equation}
\begin{cases} 
\varsigma= 1-\xi \cr
\eta =\frac{2-\rho ^{-2}}{1-\rho ^{-2}}\varsigma \cr
z=\frac{-1}{1-\rho ^{-2}}\varsigma ^2  
\end{cases}\nonumber 
\end{equation}
\subsection{ ${\displaystyle x^{-\alpha } Hl\left(\frac{1}{a},\frac{q+\alpha [(\alpha -\gamma -\delta +1)a-\beta +\delta ]}{a}; \alpha , \alpha -\gamma +1, \alpha -\beta +1,\delta ;\frac{1}{x}\right)}$}
\subsubsection{Polynomial of type 2}
Replace coefficients $a$, $q$, $\beta $, $\gamma $, $x$, $c_0$, $\lambda $ and $q_j$ where $j, q_j \in \mathbb{N}_{0}$ by $\frac{1}{a}$, $\frac{q+\alpha [(\alpha -\gamma -\delta +1)a-\beta +\delta ]}{a}$, $\alpha-\gamma +1 $, $\alpha -\beta +1 $, $\frac{1}{x}$, 1, zero and  $h_j$ where $h_j \in \mathbb{N}_{0}$ into (\ref{eq:8007}). Multiply $x^{-\alpha }$ and (\ref{eq:8007}) together. Put (\ref{eq:8006}) into the new (\ref{eq:8007}). 
\begin{eqnarray}
&&\varsigma ^{\frac{1}{2}(\alpha +1)} y(\varsigma )\nonumber\\
&=& \varsigma ^{\frac{1}{2}(\alpha +1)} Hl\left(\rho ^2, -(1+\rho ^2)(h_j+2j)(h_j+2j+1+\alpha ); \frac{\alpha}{2} +\frac{1}{2}, \frac{\alpha}{2} +1, \alpha +\frac{3}{2}, \frac{1}{2}; \varsigma \right) \nonumber\\
&=& \varsigma ^{\frac{1}{2}(\alpha +1)} \left\{ \sum_{i_0=0}^{h_0} \frac{(-h_0)_{i_0} \left( h_0+ 1+\alpha \right)_{i_0}}{(1)_{i_0}\left(\alpha +\frac{3}{2} \right)_{i_0}} \eta ^{i_0} \right. \nonumber\\
&+& \left\{ \sum_{i_0=0}^{h_0}\frac{ \left( i_0 +\frac{1}{2}+\frac{\alpha }{2} \right) \left( i_0+1+\frac{\alpha }{2} \right)}{(i_0+2)\left( i_0+ \frac{5}{2}+ \alpha \right)}\frac{(-h_0)_{i_0} \left( h_0+ 1+\alpha \right)_{i_0}}{(1)_{i_0}\left( \alpha +\frac{3}{2}\right)_{i_0}} \right. \nonumber\\
&\times& \left. \sum_{i_1=i_0}^{h_1} \frac{(-h_1)_{i_1}\left( h_1+5+\alpha \right)_{i_1}(3)_{i_0}\left( \alpha +\frac{7}{2}\right)_{i_0}}{(-h_1)_{i_0}\left( h_1+5+\alpha \right)_{i_0}(3)_{i_1} \left( \alpha +\frac{7}{2}\right)_{i_1}} \eta ^{i_1}\right\} z\nonumber\\
&+& \sum_{n=2}^{\infty } \left\{ \sum_{i_0=0}^{h_0} \frac{\left( i_0+\frac{1}{2}+\frac{\alpha }{2} \right) \left( i_0+1+\frac{\alpha }{2} \right)}{(i_0+ 2)\left( i_0+\frac{5}{2} +\alpha \right)}\frac{(-h_0)_{i_0} \left( h_0+1 +\alpha \right)_{i_0}}{(1)_{i_0}\left(\alpha +\frac{3}{2} \right)_{i_0}}\right.\nonumber\\
&\times& \prod _{k=1}^{n-1} \left\{ \sum_{i_k=i_{k-1}}^{h_k} \frac{\left( i_k+ 2k+\frac{1}{2}+\frac{\alpha }{2} \right) \left( i_k+ 2k+1-\frac{\alpha }{2} \right)}{(i_k+ 2k+2) \left( i_k+ 2k+\frac{5}{2}+\alpha  \right)}\right.\nonumber\\
&\times& \left. \frac{(-h_k)_{i_k}\left( h_k+4k+ 1+\alpha \right)_{i_k}(2k+1)_{i_{k-1}}\left( 2k+\frac{3}{2} +\alpha \right)_{i_{k-1}}}{(-h_k)_{i_{k-1}}\left(  h_k+4k+1+\alpha \right)_{i_{k-1}}(2k+1)_{i_k}\left( 2k+\frac{3}{2} +\alpha \right)_{i_k}}\right\} \nonumber\\
&\times& \left.\left. \sum_{i_n= i_{n-1}}^{h_n} \frac{(-h_n)_{i_n}\left( h_n+4n+1+\alpha  \right)_{i_n}(2n+1)_{i_{n-1}}\left( 2n+\frac{3}{2}+\alpha  \right)_{i_{n-1}}}{(-h_n)_{i_{n-1}}\left(  h_n +4n+ 1+\alpha \right)_{i_{n-1}}(2n+1)_{i_n}\left( 2n+\frac{3}{2} +\alpha \right)_{i_n}} \eta ^{i_n} \right\} z^n \right\} \label{eq:80031}
\end{eqnarray}
where
\begin{equation}
h= (1+ \rho ^2)\left( 2(h_j +2j)+1+\alpha \right)^2  \nonumber 
\end{equation}
\subsubsection{Infinite series}
Replace coefficients $a$, $q$, $\beta $, $\gamma $, $x$, $c_0$ and $\lambda $ by $\frac{1}{a}$, $\frac{q+\alpha [(\alpha -\gamma -\delta +1)a-\beta +\delta ]}{a}$, $\alpha-\gamma +1 $, $\alpha -\beta +1 $, $\frac{1}{x}$, 1 and zero into (\ref{eq:80010}). Multiply $x^{-\alpha }$ and (\ref{eq:80010}) together. Put (\ref{eq:8006}) into the new (\ref{eq:80010}). 
\begin{eqnarray}
&&\varsigma ^{\frac{1}{2}(\alpha +1)} y(\varsigma )\nonumber\\
&=& \varsigma ^{\frac{1}{2}(\alpha +1)} Hl\left(\rho ^2,-\frac{1}{4}\left( h-(1+\rho ^2)(\alpha +1)^2\right); \frac{\alpha}{2} +\frac{1}{2}, \frac{\alpha}{2} +1, \alpha +\frac{3}{2}, \frac{1}{2}; \varsigma \right) \nonumber\\
&=&\varsigma ^{\frac{1}{2}(\alpha +1)} \left\{ \sum_{i_0=0}^{\infty } \frac{\left(\Delta_0^{-}\right)_{i_0} \left(\Delta_0^{+}\right)_{i_0}}{(1 )_{i_0}\left( \alpha +\frac{3}{2} \right)_{i_0}} \eta ^{i_0} \right. \nonumber\\
&+& \left\{ \sum_{i_0=0}^{\infty }\frac{\left( i_0 +\frac{1}{2}+\frac{\alpha }{2} \right) \left( i_0 +1+\frac{\alpha }{2} \right)}{(i_0 +2)\left( i_0 +\frac{5}{2}+\alpha \right)}\frac{\left(\Delta_0^{-}\right)_{i_0} \left(\Delta_0^{+}\right)_{i_0}}{(1)_{i_0}\left( \alpha+\frac{3}{2} \right)_{i_0}} \sum_{i_1=i_0}^{\infty } \frac{\left(\Delta_1^{-}\right)_{i_1} \left(\Delta_1^{+}\right)_{i_1}(3 )_{i_0}\left( \alpha + \frac{7}{2} \right)_{i_0}}{\left(\Delta_1^{-}\right)_{i_0}  \left(\Delta_1^{+}\right)_{i_0}(3 )_{i_1}\left( \alpha +\frac{7}{2} \right)_{i_1}}\eta ^{i_1}\right\} z\nonumber\\
&+& \sum_{n=2}^{\infty } \left\{ \sum_{i_0=0}^{\infty } \frac{\left( i_0 +\frac{1}{2}+\frac{\alpha }{2} \right) \left( i_0 +1+\frac{\alpha }{2} \right)}{(i_0 +2)\left( i_0 +\frac{5}{2}+\alpha \right)}\frac{\left(\Delta_0^{-}\right)_{i_0} \left(\Delta_0^{+}\right)_{i_0}}{(1 )_{i_0}\left( \alpha+\frac{3}{2} \right)_{i_0}}\right.\nonumber\\
&\times& \prod _{k=1}^{n-1} \left\{ \sum_{i_k=i_{k-1}}^{\infty } \frac{\left( i_k+ 2k+\frac{1}{2}+\frac{\alpha }{2} \right) \left( i_k+ 2k +1+\frac{\alpha }{2} \right)}{(i_k+ 2k+2 )\left( i_k+ 2k+\frac{5}{2} +\alpha \right)} \frac{ \left(\Delta_k^{-}\right)_{i_k} \left(\Delta_k^{+} \right)_{i_k}(2k+1 )_{i_{k-1}}\left( 2k+\frac{3}{2}+\alpha \right)_{i_{k-1}}}{\left(\Delta_k^{-}\right)_{i_{k-1}} \left(\Delta_k^{+} \right)_{i_{k-1}}(2k+1 )_{i_k}\left( 2k+\frac{3}{2}+\alpha \right)_{i_k}}\right\}\nonumber\\
&\times& \left. \left. \sum_{i_n= i_{n-1}}^{\infty } \frac{\left(\Delta_n^{-}\right)_{i_n}\left( \Delta_n^{+} \right)_{i_n}(2n+1 )_{i_{n-1}}\left( 2n+\frac{3}{2} +\alpha \right)_{i_{n-1}}}{\left(\Delta_n^{-}\right)_{i_{n-1}}\left(\Delta_n^{+} \right)_{i_{n-1}}(2n+1 )_{i_n}\left( 2n+\frac{3}{2}+\alpha \right)_{i_n}} \eta ^{i_n} \right\} z^n \right\}\label{eq:80032}
\end{eqnarray}
where
\begin{equation}
\begin{cases} 
\Delta_0^{\pm}= \frac{1}{2}(\alpha +1) \pm \frac{1}{2}\sqrt{\frac{ h}{ 1+\rho^2 }} \cr
\Delta_1^{\pm}= \frac{1}{2}(\alpha +5) \pm \frac{1}{2}\sqrt{\frac{ h}{ 1+\rho^2 }} \cr
\Delta_k^{\pm}= \frac{1}{2}(\alpha +4k+1) \pm \frac{1}{2}\sqrt{\frac{ h}{ 1+\rho^2 }} \cr
\Delta_n^{\pm}= \frac{1}{2}(\alpha +4n+1) \pm \frac{1}{2}\sqrt{\frac{ h}{ 1+\rho^2 }}
\end{cases}\nonumber 
\end{equation}
On (\ref{eq:80031}) and (\ref{eq:80032}),
\begin{equation}
\begin{cases} 
\varsigma= \xi^{-1} \cr
\eta =(1+\rho ^{-2})\varsigma \cr
z= -\rho ^{-2} \varsigma ^2  
\end{cases}\nonumber 
\end{equation}
\subsection{ ${\displaystyle \left(1-\frac{x}{a} \right)^{-\beta } Hl\left(1-a, -q+\gamma \beta; -\alpha +\gamma +\delta, \beta, \gamma, \delta; \frac{(1-a)x}{x-a} \right)}$}
\subsubsection{Polynomial of type 2}
Replace coefficients $a$, $q$, $\alpha $, $x$, $c_0$, $\lambda $ and $q_j$ where $j, q_j \in \mathbb{N}_{0}$ by $1-a$, $-q+\gamma \beta $, $-\alpha+\gamma +\delta $, $\frac{(1-a)x}{x-a}$, 1, zero and $h_j$ where $h_j \in \mathbb{N}_{0}$ into (\ref{eq:8007}). Multiply $\left(1-\frac{x}{a} \right)^{-\beta }$ and (\ref{eq:8007}) together. Put (\ref{eq:8006}) into the new (\ref{eq:8007}). 
\begin{eqnarray} 
&&(1-\rho ^2 \xi)^{\frac{\alpha }{2}} y(\varsigma )\nonumber\\
&=& (1-\rho ^2 \xi)^{\frac{\alpha }{2}} Hl\left( 1-\rho ^{-2}, (h_j+2j)\left(\alpha -(2-\rho ^{-2})(h_j+2j) \right); -\frac{\alpha }{2}+\frac{1}{2}, -\frac{\alpha }{2}, \frac{1}{2}, \frac{1}{2}; \varsigma \right) \nonumber\\
&=& (1-\rho ^2 \xi)^{\frac{\alpha }{2}} \left\{\sum_{i_0=0}^{h_0} \frac{(-h_0)_{i_0} \left( h_0- \Omega _{\rho }\right)_{i_0}}{(1)_{i_0}\left(\frac{1}{2} \right)_{i_0}} \eta ^{i_0} \right.  \nonumber\\ 
&+& \left\{ \sum_{i_0=0}^{h_0}\frac{ \left( i_0+\frac{1}{2}-\frac{\alpha }{2} \right) \left( i_0 -\frac{\alpha }{2} \right)}{(i_0+2)\left( i_0+ \frac{3}{2} \right)}\frac{(-h_0)_{i_0} \left( h_0- \Omega _{\rho }\right)_{i_0}}{(1)_{i_0}\left( \frac{1}{2}\right)_{i_0}} \right. \left. \sum_{i_1=i_0}^{h_1} \frac{(-h_1)_{i_1}\left( h_1+4 - \Omega _{\rho }\right)_{i_1}(3)_{i_0}\left( \frac{5}{2}\right)_{i_0}}{(-h_1)_{i_0}\left( h_1+4- \Omega _{\rho }\right)_{i_0}(3)_{i_1} \left( \frac{5}{2}\right)_{i_1}} \eta ^{i_1}\right\} z\nonumber\\
&+& \sum_{n=2}^{\infty } \left\{ \sum_{i_0=0}^{h_0} \frac{\left( i_0+\frac{1}{2}-\frac{\alpha }{2} \right) \left( i_0 -\frac{\alpha }{2} \right)}{(i_0+ 2)\left( i_0+\frac{3}{2} \right)}\frac{(-h_0)_{i_0} \left( h_0 - \Omega _{\rho }\right)_{i_0}}{(1)_{i_0}\left(\frac{1}{2} \right)_{i_0}}\right.\nonumber\\
&\times& \prod _{k=1}^{n-1} \left\{ \sum_{i_k=i_{k-1}}^{h_k} \frac{\left( i_k+ 2k+\frac{1}{2}-\frac{\alpha }{2} \right) \left( i_k+ 2k -\frac{\alpha }{2} \right)}{(i_k+ 2k+2) \left( i_k+ 2k+\frac{3}{2} \right)} 
\frac{(-h_k)_{i_k}\left( h_k+4k- \Omega _{\rho } \right)_{i_k}(2k+1)_{i_{k-1}}\left( 2k+\frac{1}{2} \right)_{i_{k-1}}}{(-h_k)_{i_{k-1}}\left(  h_k+4k- \Omega _{\rho }\right)_{i_{k-1}}(2k+1)_{i_k}\left( 2k+\frac{1}{2} \right)_{i_k}}\right\} \nonumber\\
&\times& \left.\left. \sum_{i_n= i_{n-1}}^{h_n} \frac{(-h_n)_{i_n}\left( h_n+4n- \Omega _{\rho }\right)_{i_n}(2n+1)_{i_{n-1}}\left( 2n+\frac{1}{2} \right)_{i_{n-1}}}{(-h_n)_{i_{n-1}}\left(  h_n+4n- \Omega _{\rho }\right)_{i_{n-1}}(2n+1)_{i_n}\left( 2n+\frac{1}{2} \right)_{i_n}} \eta ^{i_n} \right\} z^n \right\} \label{eq:80033}
\end{eqnarray}
where
\begin{equation}
\begin{cases} 
\Omega _{\rho } = \frac{\alpha }{2-\rho ^{-2}} \cr
h= 4\rho ^2 \left( \alpha \left( h_j+2j+\frac{1}{4} \right) -(2-\rho ^{-2})(h_j+2j)^2\right)  
\end{cases}\nonumber 
\end{equation}
\subsubsection{Infinite series}
Replace coefficients $a$, $q$, $\alpha $, $x$, $c_0$ and $\lambda $ by $1-a$, $-q+\gamma \beta $, $-\alpha+\gamma +\delta $, $\frac{(1-a)x}{x-a}$, 1 and zero into (\ref{eq:80010}). Multiply $\left(1-\frac{x}{a} \right)^{-\beta }$ and (\ref{eq:80010}) together. Put (\ref{eq:8006}) into the new (\ref{eq:80010}). 
\begin{eqnarray}
&&(1-\rho ^2 \xi)^{\frac{\alpha }{2}} y(\varsigma )\nonumber\\
&=& (1-\rho ^2 \xi)^{\frac{\alpha }{2}} Hl\left( 1-\rho ^{-2}, \frac{1}{4}\left( h\rho ^{-2} - \alpha \right); -\frac{\alpha }{2}+\frac{1}{2}, -\frac{\alpha }{2}, \frac{1}{2}, \frac{1}{2}; \varsigma \right) \nonumber\\
&=& (1-\rho ^2 \xi)^{\frac{\alpha }{2}}  \left\{\sum_{i_0=0}^{\infty } \frac{\left(\Delta_0^{-}\right)_{i_0} \left(\Delta_0^{+}\right)_{i_0}}{(1 )_{i_0}\left( \frac{1}{2} \right)_{i_0}} \eta ^{i_0}\right.\nonumber\\
&+& \left\{ \sum_{i_0=0}^{\infty }\frac{\left( i_0 +\frac{1}{2}-\frac{\alpha }{2} \right) \left( i_0 -\frac{\alpha }{2} \right)}{(i_0 +2)\left( i_0 +\frac{3}{2} \right)}\frac{\left(\Delta_0^{-}\right)_{i_0} \left(\Delta_0^{+}\right)_{i_0}}{(1)_{i_0}\left(\frac{1}{2} \right)_{i_0}} \sum_{i_1=i_0}^{\infty } \frac{\left(\Delta_1^{-}\right)_{i_1} \left(\Delta_1^{+}\right)_{i_1}(3 )_{i_0}\left( \frac{5}{2} \right)_{i_0}}{\left(\Delta_1^{-}\right)_{i_0}  \left(\Delta_1^{+}\right)_{i_0}(3 )_{i_1}\left( \frac{5}{2} \right)_{i_1}}\eta ^{i_1}\right\} z\nonumber\\
&+& \sum_{n=2}^{\infty } \left\{ \sum_{i_0=0}^{\infty } \frac{\left( i_0 +\frac{1}{2}-\frac{\alpha }{2} \right) \left( i_0 -\frac{\alpha }{2} \right)}{(i_0 +2)\left( i_0 +\frac{3}{2} \right)}\frac{\left(\Delta_0^{-}\right)_{i_0} \left(\Delta_0^{+}\right)_{i_0}}{(1 )_{i_0}\left(\frac{1}{2} \right)_{i_0}}\right.\nonumber\\
&\times& \prod _{k=1}^{n-1} \left\{ \sum_{i_k=i_{k-1}}^{\infty } \frac{\left( i_k+ 2k+\frac{1}{2}-\frac{\alpha }{2} \right) \left( i_k+ 2k -\frac{\alpha }{2} \right)}{(i_k+ 2k+2 )\left( i_k+ 2k+\frac{3}{2} \right)} \frac{ \left(\Delta_k^{-}\right)_{i_k} \left(\Delta_k^{+} \right)_{i_k}(2k+1 )_{i_{k-1}}\left( 2k+\frac{1}{2} \right)_{i_{k-1}}}{\left(\Delta_k^{-}\right)_{i_{k-1}} \left(\Delta_k^{+} \right)_{i_{k-1}}(2k+1 )_{i_k}\left( 2k+\frac{1}{2} \right)_{i_k}}\right\}\nonumber\\
&\times& \left.\left.\sum_{i_n= i_{n-1}}^{\infty } \frac{\left(\Delta_n^{-}\right)_{i_n}\left( \Delta_n^{+} \right)_{i_n}(2n+1 )_{i_{n-1}}\left( 2n+\frac{1}{2} \right)_{i_{n-1}}}{\left(\Delta_n^{-}\right)_{i_{n-1}}\left(\Delta_n^{+} \right)_{i_{n-1}}(2n+1 )_{i_n}\left( 2n+\frac{1}{2} \right)_{i_n}} \eta ^{i_n} \right\} z^n \right\} \label{eq:80034}
\end{eqnarray}
where
\begin{equation}
\begin{cases} 
\Delta_0^{\pm}=  \frac{-\alpha \pm\sqrt{ \alpha ^2-(2-\rho ^{-2})(h\rho ^{-2}-\alpha )}}{2(2-\rho ^{-2})}   \cr
\Delta_1^{\pm}= 2 + \frac{-\alpha \pm\sqrt{ \alpha ^2-(2-\rho ^{-2})(h\rho ^{-2}-\alpha )}}{2(2-\rho ^{-2})} \cr
\Delta_k^{\pm}= 2k + \frac{-\alpha \pm\sqrt{ \alpha ^2-(2-\rho ^{-2})(h\rho ^{-2}-\alpha )}}{2(2-\rho ^{-2})} \cr
\Delta_n^{\pm}= 2n +  \frac{-\alpha \pm\sqrt{ \alpha ^2-(2-\rho ^{-2})(h\rho ^{-2}-\alpha )}}{2(2-\rho ^{-2})}
\end{cases}\nonumber 
\end{equation}
On (\ref{eq:80033}) and (\ref{eq:80034}),
\begin{equation}
\begin{cases} 
\varsigma= \frac{(1-\rho ^{-2})\xi}{\xi-\rho ^{-2}} \cr
\eta =\frac{2-\rho ^{-2}}{1-\rho ^{-2}}\varsigma \cr
z= -\frac{1}{1-\rho ^{-2}} \varsigma ^2 
\end{cases}\nonumber 
\end{equation}
\subsection{ \footnotesize ${\displaystyle (1-x)^{1-\delta }\left(1-\frac{x}{a} \right)^{-\beta+\delta -1} Hl\left(1-a, -q+\gamma [(\delta -1)a+\beta -\delta +1]; -\alpha +\gamma +1, \beta -\delta+1, \gamma, 2-\delta; \frac{(1-a)x}{x-a} \right)}$ \normalsize}
\subsubsection{Polynomial of type 2}
Replace coefficients $a$, $q$, $\alpha $, $\beta $, $\delta $, $x$, $c_0$, $\lambda $ and $q_j$ where $j, q_j \in \mathbb{N}_{0}$ by $1-a$, $-q+\gamma [(\delta -1)a+\beta -\delta +1]$, $-\alpha +\gamma +1$, $\beta -\delta+1$, $2-\delta $, $\frac{(1-a)x}{x-a}$, 1, zero and $h_j$ where $h_j \in \mathbb{N}_{0}$ into (\ref{eq:8007}). Multiply $(1-x)^{1-\delta }\left(1-\frac{x}{a} \right)^{-\beta+\delta -1}$ and (\ref{eq:8007}) together. Put (\ref{eq:8006}) into the new (\ref{eq:8007}).
\begin{eqnarray} 
&&(1-\xi )^{\frac{1}{2}}(1-\rho ^2 \xi)^{\frac{1}{2}(\alpha -1)} y(\varsigma )\nonumber\\
&=& (1-\xi )^{\frac{1}{2}}(1-\rho ^2 \xi)^{\frac{1}{2}(\alpha -1)} Hl\bigg( 1-\rho ^{-2}, (h_j+ 2j)\left( \alpha +1 -(2-\rho ^{-2})(h_j+ 2j+1)\right); -\frac{\alpha }{2}+1, \nonumber\\
&& -\frac{\alpha }{2}+\frac{1}{2}, \frac{1}{2}, \frac{3}{2}; \varsigma \bigg) \nonumber\\
&=& (1-\xi )^{\frac{1}{2}}(1-\rho ^2 \xi)^{\frac{1}{2}(\alpha -1)} \left\{\sum_{i_0=0}^{h_0} \frac{(-h_0)_{i_0} \left( h_0+ \Omega _{\rho }\right)_{i_0}}{(1)_{i_0}\left(\frac{1}{2} \right)_{i_0}} \eta ^{i_0} \right. \nonumber\\ 
&+& \left\{ \sum_{i_0=0}^{h_0}\frac{ \left( i_0+1-\frac{\alpha }{2} \right) \left( i_0 +\frac{1}{2}-\frac{\alpha }{2} \right)}{(i_0+2)\left( i_0+ \frac{3}{2} \right)}\frac{(-h_0)_{i_0} \left( h_0+ \Omega _{\rho }\right)_{i_0}}{(1)_{i_0}\left( \frac{1}{2}\right)_{i_0}} \right. \nonumber\\
&\times& \left. \sum_{i_1=i_0}^{h_1} \frac{(-h_1)_{i_1}\left( h_1+4 + \Omega _{\rho }\right)_{i_1}(3)_{i_0}\left( \frac{5}{2}\right)_{i_0}}{(-h_1)_{i_0}\left( h_1+4+ \Omega _{\rho }\right)_{i_0}(3)_{i_1} \left( \frac{5}{2}\right)_{i_1}} \eta ^{i_1}\right\} z\nonumber\\
&+& \sum_{n=2}^{\infty } \left\{ \sum_{i_0=0}^{h_0} \frac{\left( i_0+1-\frac{\alpha }{2} \right) \left( i_0 +\frac{1}{2}-\frac{\alpha }{2} \right)}{(i_0+ 2)\left( i_0+\frac{3}{2} \right)}\frac{(-h_0)_{i_0} \left( h_0 + \Omega _{\rho }\right)_{i_0}}{(1)_{i_0}\left(\frac{1}{2} \right)_{i_0}}\right.\nonumber\\
&\times& \prod _{k=1}^{n-1} \left\{ \sum_{i_k=i_{k-1}}^{h_k} \frac{\left( i_k+ 2k+1-\frac{\alpha }{2} \right) \left( i_k+ 2k +\frac{1}{2}-\frac{\alpha }{2} \right)}{(i_k+ 2k+2) \left( i_k+ 2k+\frac{3}{2} \right)}\right.\nonumber\\
&\times& \left. \frac{(-h_k)_{i_k}\left( h_k+4k+ \Omega _{\rho } \right)_{i_k}(2k+1)_{i_{k-1}}\left( 2k+\frac{1}{2} \right)_{i_{k-1}}}{(-h_k)_{i_{k-1}}\left(  h_k+4k+ \Omega _{\rho }\right)_{i_{k-1}}(2k+1)_{i_k}\left( 2k+\frac{1}{2} \right)_{i_k}}\right\} \nonumber\\
&\times& \left.\left. \sum_{i_n= i_{n-1}}^{h_n} \frac{(-h_n)_{i_n}\left( h_n+4n+ \Omega _{\rho }\right)_{i_n}(2n+1)_{i_{n-1}}\left( 2n+\frac{1}{2} \right)_{i_{n-1}}}{(-h_n)_{i_{n-1}}\left(  h_n+4n+ \Omega _{\rho }\right)_{i_{n-1}}(2n+1)_{i_n}\left( 2n+\frac{1}{2} \right)_{i_n}} \eta ^{i_n} \right\} z^n \right\} \label{eq:80035}
\end{eqnarray}
where
\begin{equation}
\begin{cases} 
\Omega _{\rho } = \frac{-\alpha +1-\rho ^{-2}}{2-\rho ^{-2}} \cr
h= 4(h_j+ 2j)\left( (\alpha +1)\rho ^2 +(1-2\rho ^2)(h_j+2j+1)\right) +(\alpha -1)\rho ^2 +1  
\end{cases}\nonumber 
\end{equation}
\subsubsection{Infinite series}
Replace coefficients $a$, $q$, $\alpha $, $\beta $, $\delta $, $x$, $c_0$ and $\lambda $ by $1-a$, $-q+\gamma [(\delta -1)a+\beta -\delta +1]$, $-\alpha +\gamma +1$, $\beta -\delta+1$, $2-\delta $, $\frac{(1-a)x}{x-a}$, 1 and zero into (\ref{eq:80010}). Multiply $(1-x)^{1-\delta }\left(1-\frac{x}{a} \right)^{-\beta+\delta -1}$ and (\ref{eq:80010}) together. Put (\ref{eq:8006}) into the new (\ref{eq:80010}).
\begin{eqnarray}
&&(1-\xi )^{\frac{1}{2}}(1-\rho ^2 \xi)^{\frac{1}{2}(\alpha -1)} y(\varsigma )\nonumber\\
&=& (1-\xi )^{\frac{1}{2}}(1-\rho ^2 \xi)^{\frac{1}{2}(\alpha -1)} Hl\left( 1-\rho ^{-2}, \frac{1}{4}\left( (h-1)\rho ^{-2} +1- \alpha \right); -\frac{\alpha }{2}+1, -\frac{\alpha }{2}+\frac{1}{2}, \frac{1}{2}, \frac{3}{2}; \varsigma \right) \nonumber\\
&=& (1-\xi )^{\frac{1}{2}}(1-\rho ^2 \xi)^{\frac{1}{2}(\alpha -1)} \left\{\sum_{i_0=0}^{\infty } \frac{\left(\Delta_0^{-}\right)_{i_0} \left(\Delta_0^{+}\right)_{i_0}}{(1 )_{i_0}\left( \frac{1}{2} \right)_{i_0}} \eta ^{i_0}\right.\nonumber\\
&+& \left\{ \sum_{i_0=0}^{\infty }\frac{\left( i_0 +1-\frac{\alpha }{2} \right) \left( i_0 +\frac{1}{2}-\frac{\alpha }{2} \right)}{(i_0 +2)\left( i_0 +\frac{3}{2} \right)}\frac{\left(\Delta_0^{-}\right)_{i_0} \left(\Delta_0^{+}\right)_{i_0}}{(1)_{i_0}\left(\frac{1}{2} \right)_{i_0}} \sum_{i_1=i_0}^{\infty } \frac{\left(\Delta_1^{-}\right)_{i_1} \left(\Delta_1^{+}\right)_{i_1}(3 )_{i_0}\left( \frac{5}{2} \right)_{i_0}}{\left(\Delta_1^{-}\right)_{i_0}  \left(\Delta_1^{+}\right)_{i_0}(3 )_{i_1}\left( \frac{5}{2} \right)_{i_1}}\eta ^{i_1}\right\} z\nonumber\\
&+& \sum_{n=2}^{\infty } \left\{ \sum_{i_0=0}^{\infty } \frac{\left( i_0 +1-\frac{\alpha }{2} \right) \left( i_0 +\frac{1}{2}-\frac{\alpha }{2} \right)}{(i_0 +2)\left( i_0 +\frac{3}{2} \right)}\frac{\left(\Delta_0^{-}\right)_{i_0} \left(\Delta_0^{+}\right)_{i_0}}{(1 )_{i_0}\left(\frac{1}{2} \right)_{i_0}}\right.\nonumber\\
&\times& \prod _{k=1}^{n-1} \left\{ \sum_{i_k=i_{k-1}}^{\infty } \frac{\left( i_k+ 2k+1-\frac{\alpha }{2} \right) \left( i_k+ 2k +\frac{1}{2}-\frac{\alpha }{2} \right)}{(i_k+ 2k+2 )\left( i_k+ 2k+\frac{3}{2} \right)} \frac{ \left(\Delta_k^{-}\right)_{i_k} \left(\Delta_k^{+} \right)_{i_k}(2k+1 )_{i_{k-1}}\left( 2k+\frac{1}{2} \right)_{i_{k-1}}}{\left(\Delta_k^{-}\right)_{i_{k-1}} \left(\Delta_k^{+} \right)_{i_{k-1}}(2k+1 )_{i_k}\left( 2k+\frac{1}{2} \right)_{i_k}}\right\}\nonumber\\
&\times& \left.\left.\sum_{i_n= i_{n-1}}^{\infty } \frac{\left(\Delta_n^{-}\right)_{i_n}\left( \Delta_n^{+} \right)_{i_n}(2n+1 )_{i_{n-1}}\left( 2n+\frac{1}{2} \right)_{i_{n-1}}}{\left(\Delta_n^{-}\right)_{i_{n-1}}\left(\Delta_n^{+} \right)_{i_{n-1}}(2n+1 )_{i_n}\left( 2n+\frac{1}{2} \right)_{i_n}} \eta ^{i_n} \right\} z^n \right\} \label{eq:80036}
\end{eqnarray}
where
\begin{equation}
\begin{cases} 
\Delta_0^{\pm}=  \frac{-\alpha +1-\rho ^{-2} \pm\sqrt{ (-\alpha +1-\rho ^{-2})^2-(2-\rho ^{-2})((h-1)\rho ^{-2}+1-\alpha )}}{2(2-\rho ^{-2})}   \cr
\Delta_1^{\pm}= 2 +  \frac{-\alpha +1-\rho ^{-2} \pm\sqrt{ (-\alpha +1-\rho ^{-2})^2-(2-\rho ^{-2})((h-1)\rho ^{-2}+1-\alpha )}}{2(2-\rho ^{-2})} \cr
\Delta_k^{\pm}= 2k +  \frac{-\alpha +1-\rho ^{-2} \pm\sqrt{ (-\alpha +1-\rho ^{-2})^2-(2-\rho ^{-2})((h-1)\rho ^{-2}+1-\alpha )}}{2(2-\rho ^{-2})} \cr
\Delta_n^{\pm}= 2n +   \frac{-\alpha +1-\rho ^{-2} \pm\sqrt{ (-\alpha +1-\rho ^{-2})^2-(2-\rho ^{-2})((h-1)\rho ^{-2}+1-\alpha )}}{2(2-\rho ^{-2})}
\end{cases}\nonumber 
\end{equation}
On (\ref{eq:80035}) and (\ref{eq:80036}),
\begin{equation}
\begin{cases} 
\varsigma= \frac{(1-\rho ^{-2})\xi}{\xi-\rho ^{-2}} \cr
\eta =\frac{2-\rho ^{-2}}{1-\rho ^{-2}}\varsigma \cr
z= -\frac{1}{1-\rho ^{-2}} \varsigma ^2 
\end{cases}\nonumber 
\end{equation}
\subsection{ ${\displaystyle x^{-\alpha } Hl\left(\frac{a-1}{a}, \frac{-q+\alpha (\delta a+\beta -\delta )}{a}; \alpha, \alpha -\gamma +1, \delta , \alpha -\beta +1; \frac{x-1}{x} \right)}$}
\subsubsection{Polynomial of type 2}
Replace coefficients $a$, $q$, $\beta $, $\gamma $, $\delta $, $x$, $c_0$, $\lambda $ and $q_j$ where $j, q_j \in \mathbb{N}_{0}$ by $\frac{a-1}{a}$, $\frac{-q+\alpha (\delta a+\beta -\delta )}{a}$, $\alpha -\gamma +1$, $\delta $, $\alpha -\beta +1$, $\frac{x-1}{x}$, 1, zero and $h_j$ where $h_j \in \mathbb{N}_{0}$ into (\ref{eq:8007}). Multiply $x^{-\alpha }$ and (\ref{eq:8007}) together. Put (\ref{eq:8006}) into the new (\ref{eq:8007}).
\begin{eqnarray} 
&&\xi ^{-\frac{1}{2}(\alpha +1)} y(\varsigma )\nonumber\\
&=& \xi ^{-\frac{1}{2}(\alpha +1)} Hl\left( 1-\rho ^2, -(h_j+2j)\left( (1-\rho ^2)(\alpha +1)+(2-\rho ^2)(h_j+2j) \right); \frac{\alpha}{2} + \frac{1}{2}, \frac{\alpha}{2} +1,\right. \nonumber\\
&&\left. \frac{1}{2}, \alpha +\frac{3}{2}; \varsigma \right) \nonumber\\
&=& \xi ^{-\frac{1}{2}(\alpha +1)} \left\{ \sum_{i_0=0}^{h_0} \frac{(-h_0)_{i_0} \left( h_0+ \Omega _{\rho }\right)_{i_0}}{(1)_{i_0}\left(\frac{1}{2} \right)_{i_0}} \eta ^{i_0} \right. \nonumber\\ 
&+& \left\{ \sum_{i_0=0}^{h_0}\frac{ \left( i_0+\frac{1}{2}+\frac{\alpha }{2} \right) \left( i_0 +1+\frac{\alpha }{2} \right)}{(i_0+2)\left( i_0+ \frac{3}{2} \right)}\frac{(-h_0)_{i_0} \left( h_0+ \Omega _{\rho }\right)_{i_0}}{(1)_{i_0}\left( \frac{1}{2}\right)_{i_0}} \right.  \left. \sum_{i_1=i_0}^{h_1} \frac{(-h_1)_{i_1}\left( h_1+4 + \Omega _{\rho }\right)_{i_1}(3)_{i_0}\left( \frac{5}{2}\right)_{i_0}}{(-h_1)_{i_0}\left( h_1+4+ \Omega _{\rho }\right)_{i_0}(3)_{i_1} \left( \frac{5}{2}\right)_{i_1}} \eta ^{i_1}\right\} z\nonumber\\
&+& \sum_{n=2}^{\infty } \left\{ \sum_{i_0=0}^{h_0} \frac{\left( i_0+\frac{1}{2}+\frac{\alpha }{2} \right) \left( i_0 +1+\frac{\alpha }{2} \right)}{(i_0+ 2)\left( i_0+\frac{3}{2} \right)}\frac{(-h_0)_{i_0} \left( h_0 + \Omega _{\rho }\right)_{i_0}}{(1)_{i_0}\left(\frac{1}{2} \right)_{i_0}}\right.\nonumber\\
&\times& \prod _{k=1}^{n-1} \left\{ \sum_{i_k=i_{k-1}}^{h_k} \frac{\left( i_k+ 2k+\frac{1}{2}+\frac{\alpha }{2} \right) \left( i_k+ 2k +1+\frac{\alpha }{2} \right)}{(i_k+ 2k+2) \left( i_k+ 2k+\frac{3}{2} \right)} \right.\nonumber\\
&\times& \left. \frac{(-h_k)_{i_k}\left( h_k+4k+ \Omega _{\rho }\right)_{i_k}(2k+1)_{i_{k-1}}\left( 2k+\frac{1}{2} \right)_{i_{k-1}}}{(-h_k)_{i_{k-1}}\left(  h_k+4k+ \Omega _{\rho }\right)_{i_{k-1}}(2k+1)_{i_k}\left( 2k+\frac{1}{2} \right)_{i_k}}\right\} \nonumber\\
&\times& \left.\left. \sum_{i_n= i_{n-1}}^{h_n} \frac{(-h_n)_{i_n}\left( h_n+4n+ \Omega _{\rho }\right)_{i_n}(2n+1)_{i_{n-1}}\left( 2n+\frac{1}{2} \right)_{i_{n-1}}}{(-h_n)_{i_{n-1}}\left(  h_n+4n+ \Omega _{\rho }\right)_{i_{n-1}}(2n+1)_{i_n}\left( 2n+\frac{1}{2} \right)_{i_n}} \eta ^{i_n} \right\} z^n \right\} \label{eq:80037}
\end{eqnarray}
where
\begin{equation}
\begin{cases} 
\Omega _{\rho } = \frac{(\alpha +1)(1-\rho ^2)}{2-\rho ^2} \cr
h= -4(h_j+ 2j)\left( (1- \rho ^2)(\alpha +1) +(2-\rho ^2)(h_j+2j)\right) -(\alpha +1)\left( 1-\rho ^2(\alpha +1)\right)  
\end{cases}\nonumber 
\end{equation}
\subsubsection{Infinite series}
Replace coefficients $a$, $q$, $\beta $, $\gamma $, $\delta $, $x$, $c_0$ and $\lambda $ by $\frac{a-1}{a}$, $\frac{-q+\alpha (\delta a+\beta -\delta )}{a}$, $\alpha -\gamma +1$, $\delta $, $\alpha -\beta +1$, $\frac{x-1}{x}$, 1 and zero into (\ref{eq:80010}). Multiply $x^{-\alpha }$ and (\ref{eq:80010}) together. Put (\ref{eq:8006}) into the new (\ref{eq:80010}).
\begin{eqnarray}
&&\xi ^{-\frac{1}{2}(\alpha +1)} y(\varsigma )\nonumber\\
&=& \xi ^{-\frac{1}{2}(\alpha +1)} Hl\left( 1-\rho ^2, \frac{1}{4}\left[ h+ (\alpha +1)\left( 1-\rho ^2(\alpha +1)\right) \right]; \frac{\alpha}{2} + \frac{1}{2}, \frac{\alpha}{2} +1, \frac{1}{2}, \alpha +\frac{3}{2}; \varsigma \right) \nonumber\\
&=& \xi ^{-\frac{1}{2}(\alpha +1)} \left\{\sum_{i_0=0}^{\infty } \frac{\left(\Delta_0^{-}\right)_{i_0} \left(\Delta_0^{+}\right)_{i_0}}{(1 )_{i_0}\left( \frac{1}{2} \right)_{i_0}} \eta ^{i_0}\right.\nonumber\\
&+& \left\{ \sum_{i_0=0}^{\infty }\frac{\left( i_0 +\frac{1}{2}+\frac{\alpha }{2} \right) \left( i_0 +1+\frac{\alpha }{2} \right)}{(i_0 +2)\left( i_0 +\frac{3}{2} \right)}\frac{\left(\Delta_0^{-}\right)_{i_0} \left(\Delta_0^{+}\right)_{i_0}}{(1)_{i_0}\left(\frac{1}{2} \right)_{i_0}} \sum_{i_1=i_0}^{\infty } \frac{\left(\Delta_1^{-}\right)_{i_1} \left(\Delta_1^{+}\right)_{i_1}(3 )_{i_0}\left( \frac{5}{2} \right)_{i_0}}{\left(\Delta_1^{-}\right)_{i_0}  \left(\Delta_1^{+}\right)_{i_0}(3 )_{i_1}\left( \frac{5}{2} \right)_{i_1}}\eta ^{i_1}\right\} z\nonumber\\
&+& \sum_{n=2}^{\infty } \left\{ \sum_{i_0=0}^{\infty } \frac{\left( i_0 +\frac{1}{2}+\frac{\alpha }{2} \right) \left( i_0 +1+\frac{\alpha }{2} \right)}{(i_0 +2)\left( i_0 +\frac{3}{2} \right)}\frac{\left(\Delta_0^{-}\right)_{i_0} \left(\Delta_0^{+}\right)_{i_0}}{(1 )_{i_0}\left(\frac{1}{2} \right)_{i_0}}\right.\nonumber\\
&\times& \prod _{k=1}^{n-1} \left\{ \sum_{i_k=i_{k-1}}^{\infty } \frac{\left( i_k+ 2k+\frac{1}{2}+\frac{\alpha }{2} \right) \left( i_k+ 2k +1+\frac{\alpha }{2} \right)}{(i_k+ 2k+2 )\left( i_k+ 2k+\frac{3}{2} \right)} \frac{ \left(\Delta_k^{-}\right)_{i_k} \left(\Delta_k^{+} \right)_{i_k}(2k+1 )_{i_{k-1}}\left( 2k+\frac{1}{2} \right)_{i_{k-1}}}{\left(\Delta_k^{-}\right)_{i_{k-1}} \left(\Delta_k^{+} \right)_{i_{k-1}}(2k+1 )_{i_k}\left( 2k+\frac{1}{2} \right)_{i_k}}\right\}\nonumber\\
&\times& \left.\left.\sum_{i_n= i_{n-1}}^{\infty } \frac{\left(\Delta_n^{-}\right)_{i_n}\left( \Delta_n^{+} \right)_{i_n}(2n+1 )_{i_{n-1}}\left( 2n+\frac{1}{2} \right)_{i_{n-1}}}{\left(\Delta_n^{-}\right)_{i_{n-1}}\left(\Delta_n^{+} \right)_{i_{n-1}}(2n+1 )_{i_n}\left( 2n+\frac{1}{2} \right)_{i_n}} \eta ^{i_n} \right\} z^n \right\} \label{eq:80038}
\end{eqnarray}
where
\begin{equation}
\begin{cases} 
\Delta_0^{\pm}=  \frac{(1-\rho ^2)(\alpha +1) \pm\sqrt{ ( \alpha +1 )^2-(2-\rho ^{-2})(h+1+\alpha )}}{2(2-\rho ^2)}   \cr
\Delta_1^{\pm}= 2 +  \frac{(1-\rho ^2)(\alpha +1) \pm\sqrt{ ( \alpha +1 )^2-(2-\rho ^{-2})(h+1+\alpha )}}{2(2-\rho ^2)} \cr
\Delta_k^{\pm}= 2k +  \frac{(1-\rho ^2)(\alpha +1) \pm\sqrt{ ( \alpha +1 )^2-(2-\rho ^{-2})(h+1+\alpha )}}{2(2-\rho ^2)} \cr
\Delta_n^{\pm}= 2n +   \frac{(1-\rho ^2)(\alpha +1) \pm\sqrt{ ( \alpha +1 )^2-(2-\rho ^{-2})(h+1+\alpha )}}{2(2-\rho ^2)}
\end{cases}\nonumber 
\end{equation}
On (\ref{eq:80037}) and (\ref{eq:80038}),
\begin{equation}
\begin{cases} 
\varsigma= 1-\xi ^{-1} \cr
\eta =\frac{2-\rho ^2}{1-\rho ^2}\varsigma \cr
z= -\frac{1}{1-\rho ^2} \varsigma ^2 
\end{cases}\nonumber 
\end{equation}
\subsection{ ${\displaystyle \left(\frac{x-a}{1-a} \right)^{-\alpha } Hl\left(a, q-(\beta -\delta )\alpha ; \alpha , -\beta+\gamma +\delta , \delta , \gamma; \frac{a(x-1)}{x-a} \right)}$}
\subsubsection{Polynomial of type 2}
Replace coefficients $q$, $\beta $, $\gamma $, $\delta $, $x$, $c_0$, $\lambda $ and $q_j$ where $j, q_j \in \mathbb{N}_{0}$ by $q-(\beta -\delta )\alpha $, $-\beta+\gamma +\delta $, $\delta $,  $\gamma $, $\frac{a(x-1)}{x-a}$, 1, zero and $h_j$ where $h_j \in \mathbb{N}_{0}$ into (\ref{eq:8007}). Multiply $\left(\frac{x-a}{1-a} \right)^{-\alpha }$ and (\ref{eq:8007}) together. Put (\ref{eq:8006}) into the new (\ref{eq:8007}).
\begin{eqnarray} 
&&\left(\frac{\xi-\rho ^{-2}}{1-\rho ^{-2}} \right)^{-\frac{1}{2}(\alpha +1)} y(\varsigma )\nonumber\\
&=& \left(\frac{\xi-\rho ^{-2}}{1-\rho ^{-2}} \right)^{-\frac{1}{2}(\alpha +1)} Hl\left( \rho ^{-2}, -(h_j+2j)\left( (1+\rho ^{-2})(h_j+2j)+\alpha +1\right); \frac{\alpha }{2} +\frac{1}{2}, \frac{\alpha }{2} +1,\right. \nonumber\\
&&\left. \frac{1}{2}, \frac{1}{2}; \varsigma \right) \nonumber\\
&=& \left(\frac{\xi-\rho ^{-2}}{1-\rho ^{-2}} \right)^{-\frac{1}{2}(\alpha +1)} \left\{ \sum_{i_0=0}^{h_0} \frac{(-h_0)_{i_0} \left( h_0+ \Omega _{\rho }\right)_{i_0}}{(1)_{i_0}\left(\frac{1}{2} \right)_{i_0}} \eta ^{i_0} \right. \nonumber\\ 
&+& \left\{ \sum_{i_0=0}^{h_0}\frac{ \left( i_0+\frac{1}{2}+\frac{\alpha }{2} \right) \left( i_0 +1+\frac{\alpha }{2} \right)}{(i_0+2)\left( i_0+ \frac{3}{2} \right)}\frac{(-h_0)_{i_0} \left( h_0+ \Omega _{\rho }\right)_{i_0}}{(1)_{i_0}\left( \frac{1}{2}\right)_{i_0}} \right.  \left. \sum_{i_1=i_0}^{h_1} \frac{(-h_1)_{i_1}\left( h_1+4 + \Omega _{\rho } \right)_{i_1}(3)_{i_0}\left( \frac{5}{2}\right)_{i_0}}{(-h_1)_{i_0}\left( h_1+4+ \Omega _{\rho }\right)_{i_0}(3)_{i_1} \left( \frac{5}{2}\right)_{i_1}} \eta ^{i_1}\right\} z\nonumber\\
&+& \sum_{n=2}^{\infty } \left\{ \sum_{i_0=0}^{h_0} \frac{\left( i_0+\frac{1}{2}+\frac{\alpha }{2} \right) \left( i_0 +1+\frac{\alpha }{2} \right)}{(i_0+ 2)\left( i_0+\frac{3}{2} \right)}\frac{(-h_0)_{i_0} \left( h_0 + \Omega _{\rho }\right)_{i_0}}{(1)_{i_0}\left(\frac{1}{2} \right)_{i_0}}\right.\nonumber\\
&\times& \prod _{k=1}^{n-1} \left\{ \sum_{i_k=i_{k-1}}^{h_k} \frac{\left( i_k+ 2k+\frac{1}{2}+\frac{\alpha }{2} \right) \left( i_k+ 2k +1+\frac{\alpha }{2} \right)}{(i_k+ 2k+2) \left( i_k+ 2k+\frac{3}{2} \right)} \right.\nonumber\\
&\times& \left. \frac{(-h_k)_{i_k}\left( h_k+4k+ \Omega _{\rho } \right)_{i_k}(2k+1)_{i_{k-1}}\left( 2k+\frac{1}{2} \right)_{i_{k-1}}}{(-h_k)_{i_{k-1}}\left(  h_k+4k+ \Omega _{\rho } \right)_{i_{k-1}}(2k+1)_{i_k}\left( 2k+\frac{1}{2} \right)_{i_k}}\right\} \nonumber\\
&\times& \left.\left. \sum_{i_n= i_{n-1}}^{h_n} \frac{(-h_n)_{i_n}\left( h_n+4n+ \Omega _{\rho } \right)_{i_n}(2n+1)_{i_{n-1}}\left( 2n+\frac{1}{2} \right)_{i_{n-1}}}{(-h_n)_{i_{n-1}}\left(  h_n+4n+ \Omega _{\rho }\right)_{i_{n-1}}(2n+1)_{i_n}\left( 2n+\frac{1}{2} \right)_{i_n}} \eta ^{i_n} \right\} z^n \right\} \label{eq:80039}
\end{eqnarray}
where
\begin{equation}
\begin{cases} 
\Omega _{\rho } = \frac{ \alpha +1 }{1+\rho ^{-2}} \cr
h= 4(h_j+ 2j)\left( (1+ \rho ^2)(h_j+ 2j) +\rho ^2(\alpha +1) \right) +\rho ^2(\alpha +1)^2 
\end{cases}\nonumber 
\end{equation}
\subsubsection{Infinite series}
Replace coefficients $q$, $\beta $, $\gamma $, $\delta $, $x$, $c_0$ and $\lambda $ by $q-(\beta -\delta )\alpha $, $-\beta+\gamma +\delta $, $\delta $,  $\gamma $, $\frac{a(x-1)}{x-a}$, 1 and zero into (\ref{eq:80010}). Multiply $\left(\frac{x-a}{1-a} \right)^{-\alpha }$ and (\ref{eq:80010}) together. Put (\ref{eq:8006}) into the new (\ref{eq:80010}).
\begin{eqnarray}
&&\left(\frac{\xi-\rho ^{-2}}{1-\rho ^{-2}} \right)^{-\frac{1}{2}(\alpha +1)} y(\varsigma )\nonumber\\
&=& \left(\frac{\xi-\rho ^{-2}}{1-\rho ^{-2}} \right)^{-\frac{1}{2}(\alpha +1)} Hl\left( \rho ^{-2}, -\frac{1}{4}\left( h\rho ^{-2}- (\alpha +1)^2\right); \frac{\alpha }{2} +\frac{1}{2}, \frac{\alpha }{2} +1, \frac{1}{2}, \frac{1}{2}; \varsigma \right) \nonumber\\
&=& \left(\frac{\xi-\rho ^{-2}}{1-\rho ^{-2}} \right)^{-\frac{1}{2}(\alpha +1)} \left\{\sum_{i_0=0}^{\infty } \frac{\left(\Delta_0^{-}\right)_{i_0} \left(\Delta_0^{+}\right)_{i_0}}{(1 )_{i_0}\left( \frac{1}{2} \right)_{i_0}} \eta ^{i_0}\right.\nonumber\\
&+& \left\{ \sum_{i_0=0}^{\infty }\frac{\left( i_0 +\frac{1}{2}+\frac{\alpha }{2} \right) \left( i_0 +1+\frac{\alpha }{2} \right)}{(i_0 +2)\left( i_0 +\frac{3}{2} \right)}\frac{\left(\Delta_0^{-}\right)_{i_0} \left(\Delta_0^{+}\right)_{i_0}}{(1)_{i_0}\left(\frac{1}{2} \right)_{i_0}} \sum_{i_1=i_0}^{\infty } \frac{\left(\Delta_1^{-}\right)_{i_1} \left(\Delta_1^{+}\right)_{i_1}(3 )_{i_0}\left( \frac{5}{2} \right)_{i_0}}{\left(\Delta_1^{-}\right)_{i_0}  \left(\Delta_1^{+}\right)_{i_0}(3 )_{i_1}\left( \frac{5}{2} \right)_{i_1}}\eta ^{i_1}\right\} z\nonumber\\
&+& \sum_{n=2}^{\infty } \left\{ \sum_{i_0=0}^{\infty } \frac{\left( i_0 +\frac{1}{2}+\frac{\alpha }{2} \right) \left( i_0 +1+\frac{\alpha }{2} \right)}{(i_0 +2)\left( i_0 +\frac{3}{2} \right)}\frac{\left(\Delta_0^{-}\right)_{i_0} \left(\Delta_0^{+}\right)_{i_0}}{(1 )_{i_0}\left(\frac{1}{2} \right)_{i_0}}\right.\nonumber\\
&\times& \prod _{k=1}^{n-1} \left\{ \sum_{i_k=i_{k-1}}^{\infty } \frac{\left( i_k+ 2k+\frac{1}{2}+\frac{\alpha }{2} \right) \left( i_k+ 2k +1+\frac{\alpha }{2} \right)}{(i_k+ 2k+2 )\left( i_k+ 2k+\frac{3}{2} \right)} \frac{ \left(\Delta_k^{-}\right)_{i_k} \left(\Delta_k^{+} \right)_{i_k}(2k+1 )_{i_{k-1}}\left( 2k+\frac{1}{2} \right)_{i_{k-1}}}{\left(\Delta_k^{-}\right)_{i_{k-1}} \left(\Delta_k^{+} \right)_{i_{k-1}}(2k+1 )_{i_k}\left( 2k+\frac{1}{2} \right)_{i_k}}\right\}\nonumber\\
&\times& \left.\left.\sum_{i_n= i_{n-1}}^{\infty } \frac{\left(\Delta_n^{-}\right)_{i_n}\left( \Delta_n^{+} \right)_{i_n}(2n+1 )_{i_{n-1}}\left( 2n+\frac{1}{2} \right)_{i_{n-1}}}{\left(\Delta_n^{-}\right)_{i_{n-1}}\left(\Delta_n^{+} \right)_{i_{n-1}}(2n+1 )_{i_n}\left( 2n+\frac{1}{2} \right)_{i_n}} \eta ^{i_n} \right\} z^n \right\} \label{eq:80040}
\end{eqnarray}
where
\begin{equation}
\begin{cases} 
\Delta_0^{\pm}=  \frac{ \rho ^2(\alpha +1) \pm\sqrt{ h( 1 +\rho ^2) - \rho ^2( \alpha +1)^2}}{2(1+\rho ^2)}   \cr
\Delta_1^{\pm}= 2 + \frac{ \rho ^2(\alpha +1) \pm\sqrt{ h( 1 +\rho ^2) - \rho ^2( \alpha +1)^2}}{2(1+\rho ^2)} \cr
\Delta_k^{\pm}= 2k + \frac{ \rho ^2(\alpha +1) \pm\sqrt{ h( 1 +\rho ^2) - \rho ^2( \alpha +1)^2}}{2(1+\rho ^2)} \cr
\Delta_n^{\pm}= 2n + \frac{ \rho ^2(\alpha +1) \pm\sqrt{ h( 1 +\rho ^2) - \rho ^2( \alpha +1)^2}}{2(1+\rho ^2)}
\end{cases}\nonumber 
\end{equation}
On (\ref{eq:80039}) and (\ref{eq:80040}),
\begin{equation}
\begin{cases} 
\varsigma= \frac{\xi -1}{\rho ^2(\xi -\rho ^{-2})} \cr
\eta = (1+\rho ^2) \varsigma \cr
z= -\rho ^2 \varsigma ^2 
\end{cases}\nonumber 
\end{equation}
\section{Integral representation}
\subsection{ ${\displaystyle (1-x)^{1-\delta } Hl(a, q - (\delta  - 1)\gamma a; \alpha - \delta  + 1, \beta - \delta + 1, \gamma ,2 - \delta ; x)}$ }
\subsubsection{Polynomial of type 2}
Replace coefficients $q$, $\alpha$, $\beta$, $\delta$, $c_0$, $\lambda $ and $q_j$ where $j, q_j \in \mathbb{N}_{0}$ by $q - (\delta - 1)\gamma a $, $\alpha - \delta  + 1 $, $\beta - \delta + 1$, $2 - \delta$, 1, zero and $h_j$ where $h_j \in \mathbb{N}_{0}$ into (\ref{eq:80013}). Multiply $(1-x)^{1-\delta }$ and (\ref{eq:80013}) together. Put (\ref{eq:8006}) into the new (\ref{eq:80013}).
\begin{eqnarray}
&& (1-\xi )^{\frac{1}{2}} y(\xi )\nonumber\\
&=& (1-\xi )^{\frac{1}{2}} Hl\left(\rho ^{-2}, -\rho ^{-2}(h_j+2j)[1+(1+\rho ^2)(h_j+2j)]; \frac{\alpha }{2}+1, -\frac{\alpha }{2}+\frac{1}{2}, \frac{1}{2},\frac{3}{2}; \xi \right) \nonumber\\
&=& (1-\xi )^{\frac{1}{2}} \Bigg\{\;_2F_1\left( -h_0,h_0+ \Omega _{\rho }; \frac{1}{2}; \eta\right)   + \sum_{n=1}^{\infty } \Bigg\{\prod _{k=0}^{n-1} \Bigg\{ \int_{0}^{1} dt_{n-k}\;t_{n-k}^{2(n-k)-1 } \int_{0}^{1} du_{n-k}\;u_{n-k}^{2(n-k)-\frac{3}{2} } \nonumber\\
&&\times  \frac{1}{2\pi i}  \oint dv_{n-k} \frac{1}{v_{n-k}} \left( \frac{v_{n-k}-1}{v_{n-k}} \frac{1}{1-\overleftrightarrow {w}_{n-k+1,n}(1-t_{n-k})(1-u_{n-k})v_{n-k}}\right)^{h_{n-k}} \nonumber\\
&&\times \left( 1- \overleftrightarrow {w}_{n-k+1,n}(1-t_{n-k})(1-u_{n-k})v_{n-k}\right)^{-\left( 4(n-k)+\Omega _{\rho }\right)}\nonumber\\
&&\times \overleftrightarrow {w}_{n-k,n}^{-(2(n-k)-1+\frac{\alpha }{2} )}\left(  \overleftrightarrow {w}_{n-k,n} \partial _{ \overleftrightarrow {w}_{n-k,n}}\right) \overleftrightarrow {w}_{n-k,n}^{\alpha +\frac{1}{2}} \left(  \overleftrightarrow {w}_{n-k,n} \partial _{ \overleftrightarrow {w}_{n-k,n}}\right) \overleftrightarrow {w}_{n-k,n}^{2(n-k)-\frac{3}{2}-\frac{\alpha }{2} } \Bigg\}\nonumber\\
&&\times \;_2F_1\left( -h_0, h_0+\Omega _{\rho }; \frac{1}{2}; \overleftrightarrow {w}_{1,n}\right) \Bigg\} z^n \Bigg\} \label{eq:80041}
\end{eqnarray}
where
\begin{equation}
\begin{cases} 
\Omega _{\rho } = \frac{1}{1+\rho ^2}\cr
h= 4(h_j+2j)[1+(1+\rho ^2)(h_j+2j)]+1  
\end{cases}\nonumber 
\end{equation}
\subsubsection{Infinite series}
Replace coefficients $q$, $\alpha$, $\beta$, $\delta$, $c_0$ and $\lambda $ by $q - (\delta - 1)\gamma a $, $\alpha - \delta  + 1 $, $\beta - \delta + 1$, $2-\delta$, 1 and zero into (\ref{eq:80016}). Multiply $(1-x)^{1-\delta }$ and (\ref{eq:80016}) together. Put (\ref{eq:8006}) into the new (\ref{eq:80016}).
\begin{eqnarray}
&& (1-\xi )^{\frac{1}{2}} y(\xi )\nonumber\\
&=& (1-\xi )^{\frac{1}{2}} Hl\left(\rho ^{-2}, -\frac{1}{4}(h-1)\rho ^{-2}; \frac{\alpha }{2}+1, -\frac{\alpha }{2}+\frac{1}{2}, \frac{1}{2},\frac{3}{2}; \xi \right) \nonumber\\
&=& (1-\xi )^{\frac{1}{2}}  \left\{ \;_2F_1\left( -\Pi_0 ^{+}, -\Pi_0 ^{-}; \frac{1}{2}; \eta\right)  \right. + \sum_{n=1}^{\infty } \left\{\prod _{k=0}^{n-1} \Bigg\{ \int_{0}^{1} dt_{n-k}\;t_{n-k}^{2(n-k)-1 } \int_{0}^{1} du_{n-k}\;u_{n-k}^{2(n-k)-\frac{3}{2} } \right.\nonumber\\
&\times& \frac{1}{2\pi i}  \oint dv_{n-k} \frac{1}{v_{n-k}}\left( \frac{v_{n-k}-1}{v_{n-k}}\right)^{\Pi_{n-k} ^{+}}  \left( 1- \overleftrightarrow {w}_{n-k+1,n}(1-t_{n-k})(1-u_{n-k})v_{n-k}\right)^{\Pi_{n-k} ^{-}}\nonumber\\
&\times& \overleftrightarrow {w}_{n-k,n}^{-(2(n-k)-1+\frac{\alpha }{2})}\left(  \overleftrightarrow {w}_{n-k,n} \partial _{ \overleftrightarrow {w}_{n-k,n}}\right) \overleftrightarrow {w}_{n-k,n}^{\alpha +\frac{1}{2}} \left(  \overleftrightarrow {w}_{n-k,n} \partial _{ \overleftrightarrow {w}_{n-k,n}}\right) \overleftrightarrow {w}_{n-k,n}^{2(n-k)-\frac{3}{2}-\frac{\alpha }{2}  } \Bigg\}\nonumber\\
&\times& \left.\left. \;_2F_1\left( -\Pi_0 ^{+}, -\Pi_0 ^{-}; \frac{1}{2}; \overleftrightarrow {w}_{1,n}\right) \right\} z^n \right\} \label{eq:80042}
\end{eqnarray}
where
 \begin{equation}
\begin{cases} 
\Pi_0 ^{\pm}= \frac{-1\pm\sqrt{1+(1+\rho ^2)(h-1) }}{2(1+\rho ^2)} \cr
\Pi_{n-k} ^{\pm}= -2(n-k) + \frac{-1\pm\sqrt{1+(1+\rho ^2)(h-1) }}{2(1+\rho ^2)} 
\end{cases}\nonumber 
\end{equation}
On (\ref{eq:80041}) and (\ref{eq:80042}),
\begin{equation}
\begin{cases} 
\eta =(1+\rho ^2)\xi \cr
z=-\rho ^2\xi^2 
\end{cases}\nonumber 
\end{equation}
\subsection{ \footnotesize ${\displaystyle x^{1-\gamma } (1-x)^{1-\delta }Hl(a, q-(\gamma +\delta -2)a -(\gamma -1)(\alpha +\beta -\gamma -\delta +1), \alpha - \gamma -\delta +2, \beta - \gamma -\delta +2, 2-\gamma, 2 - \delta ; x)}$ \normalsize}
\subsubsection{Polynomial of type 2}
Replace coefficients $q$, $\alpha$, $\beta$, $\gamma $, $\delta$, $c_0$, $\lambda $ and $q_j$ where $j, q_j \in \mathbb{N}_{0}$ by $q-(\gamma +\delta -2)a-(\gamma -1)(\alpha +\beta -\gamma -\delta +1)$, $\alpha - \gamma -\delta +2$, $\beta - \gamma -\delta +2, 2-\gamma$, $2 - \delta$, 1, zero and $h_j$ where $h_j \in \mathbb{N}_{0}$ into (\ref{eq:80013}). Multiply $x^{1-\gamma } (1-x)^{1-\delta }$ and (\ref{eq:80013}) together. Put (\ref{eq:8006}) into the new (\ref{eq:80013}). 
\begin{eqnarray}
&&\xi ^{\frac{1}{2}} (1-\xi )^{\frac{1}{2}} y(\xi )\nonumber\\
&=&\xi ^{\frac{1}{2}} (1-\xi )^{\frac{1}{2}} Hl\left(\rho ^{-2}, -\rho ^{-2}(h_j+2j)[2+\rho ^2 +(1+\rho ^2)(h_j+2j)]; \frac{\alpha }{2}+\frac{3}{2}, -\frac{\alpha }{2}+1, \frac{3}{2},\frac{3}{2}; \xi \right)\nonumber\\
&=& \xi ^{\frac{1}{2}} (1-\xi )^{\frac{1}{2}} \left\{ \;_2F_1\left( -h_0, h_0+ \Omega _{\rho }; \frac{3}{2}; \eta\right) \right.  + \sum_{n=1}^{\infty } \Bigg\{\prod _{k=0}^{n-1} \Bigg\{ \int_{0}^{1} dt_{n-k}\;t_{n-k}^{2(n-k)-1 } \int_{0}^{1} du_{n-k}\;u_{n-k}^{2(n-k)-\frac{1}{2} } \nonumber\\
&&\times  \frac{1}{2\pi i}  \oint dv_{n-k} \frac{1}{v_{n-k}} \left( \frac{v_{n-k}-1}{v_{n-k}} \frac{1}{1-\overleftrightarrow {w}_{n-k+1,n}(1-t_{n-k})(1-u_{n-k})v_{n-k}}\right)^{h_{n-k}} \nonumber\\
&&\times \left( 1- \overleftrightarrow {w}_{n-k+1,n}(1-t_{n-k})(1-u_{n-k})v_{n-k}\right)^{-\left( 4(n-k)+\Omega _{\rho }\right)}\nonumber\\
&&\times \overleftrightarrow {w}_{n-k,n}^{-(2(n-k)-\frac{1}{2}+\frac{\alpha }{2} )}\left(  \overleftrightarrow {w}_{n-k,n} \partial _{ \overleftrightarrow {w}_{n-k,n}}\right) \overleftrightarrow {w}_{n-k,n}^{\alpha +\frac{1}{2}} \left(  \overleftrightarrow {w}_{n-k,n} \partial _{ \overleftrightarrow {w}_{n-k,n}}\right) \overleftrightarrow {w}_{n-k,n}^{2(n-k)-1-\frac{\alpha }{2} } \Bigg\}\nonumber\\
&&\times \;_2F_1\left( -h_0, h_0+\Omega _{\rho }; \frac{3}{2}; \overleftrightarrow {w}_{1,n}\right) \Bigg\} z^n \Bigg\} \label{eq:80043}
\end{eqnarray}
where
\begin{equation}
\begin{cases} 
\Omega _{\rho } = \frac{2+\rho ^2}{1+\rho ^2}\cr
h= 4+\rho ^2+ 4(h_j+2j)[2+\rho ^2+(1+\rho ^2)(h_j+2j)]  
\end{cases}\nonumber 
\end{equation}
\subsubsection{Infinite series}
Replace coefficients $q$, $\alpha$, $\beta$, $\gamma $, $\delta$, $c_0$ and $\lambda $ by $q-(\gamma +\delta -2)a-(\gamma -1)(\alpha +\beta -\gamma -\delta +1)$, $\alpha - \gamma -\delta +2$, $\beta - \gamma -\delta +2, 2-\gamma$, $2 - \delta$, 1 and zero into (\ref{eq:80016}). Multiply $x^{1-\gamma } (1-x)^{1-\delta }$ and (\ref{eq:80016}) together. Put (\ref{eq:8006}) into the new (\ref{eq:80016}). 
\begin{eqnarray}
&&\xi ^{\frac{1}{2}} (1-\xi )^{\frac{1}{2}} y(\xi )\nonumber\\
&=&\xi ^{\frac{1}{2}} (1-\xi )^{\frac{1}{2}} Hl\left(\rho ^{-2}, -\frac{1}{4}\left( (h-4)\rho ^{-2}-1\right); \frac{\alpha }{2}+\frac{3}{2}, -\frac{\alpha }{2}+1, \frac{3}{2},\frac{3}{2}; \xi \right)\nonumber\\
&=& \xi ^{\frac{1}{2}} (1-\xi )^{\frac{1}{2}} \left\{ \;_2F_1\left( -\Pi_0 ^{+}, -\Pi_0 ^{-}; \frac{3}{2}; \eta\right)  \right. + \sum_{n=1}^{\infty } \left\{\prod _{k=0}^{n-1} \Bigg\{ \int_{0}^{1} dt_{n-k}\;t_{n-k}^{2(n-k)-1 } \int_{0}^{1} du_{n-k}\;u_{n-k}^{2(n-k)-\frac{1}{2} } \right.\nonumber\\
&\times& \frac{1}{2\pi i}  \oint dv_{n-k} \frac{1}{v_{n-k}}\left( \frac{v_{n-k}-1}{v_{n-k}}\right)^{\Pi_{n-k} ^{+}}  \left( 1- \overleftrightarrow {w}_{n-k+1,n}(1-t_{n-k})(1-u_{n-k})v_{n-k}\right)^{\Pi_{n-k} ^{-}}\nonumber\\
&\times& \overleftrightarrow {w}_{n-k,n}^{-(2(n-k)-\frac{1}{2}+\frac{\alpha }{2})}\left(  \overleftrightarrow {w}_{n-k,n} \partial _{ \overleftrightarrow {w}_{n-k,n}}\right) \overleftrightarrow {w}_{n-k,n}^{\alpha +\frac{1}{2}} \left(  \overleftrightarrow {w}_{n-k,n} \partial _{ \overleftrightarrow {w}_{n-k,n}}\right) \overleftrightarrow {w}_{n-k,n}^{2(n-k)-1-\frac{\alpha }{2}  } \Bigg\}\nonumber\\
&\times& \left.\left. \;_2F_1\left( -\Pi_0 ^{+}, -\Pi_0 ^{-}; \frac{3}{2}; \overleftrightarrow {w}_{1,n}\right) \right\} z^n \right\} \label{eq:80044}
\end{eqnarray}
where
\begin{equation}
\begin{cases} 
\Pi_0 ^{\pm}= \frac{-(2+\rho ^2) \pm\sqrt{\rho ^2(h-1) +h }}{2(1+\rho ^2)} \cr
\Pi_{n-k} ^{\pm}=  -2(n-k) + \frac{-(2+\rho ^2) \pm\sqrt{\rho ^2(h-1) +h }}{2(1+\rho ^2)} 
\end{cases}\nonumber 
\end{equation}
On (\ref{eq:80043}) and (\ref{eq:80044}),
\begin{equation}
\begin{cases} 
\eta =(1+\rho ^2)\xi \cr
z=-\rho ^2\xi^2 
\end{cases}\nonumber 
\end{equation}
\subsection{ ${\displaystyle  Hl(1-a,-q+\alpha \beta; \alpha,\beta, \delta, \gamma; 1-x)}$} 
\subsubsection{Polynomial of type 2}
Replace coefficients $a$, $q$, $\gamma $, $\delta$, $x$, $c_0$, $\lambda $ and $q_j$ where $j, q_j \in \mathbb{N}_{0}$ by $1-a$, $-q +\alpha \beta $, $\delta $, $\gamma $, $1-x$, 1, zero and $h_j$ where $h_j \in \mathbb{N}_{0}$ into (\ref{eq:80013}). Put (\ref{eq:8006}) into the new (\ref{eq:80013}). 
\begin{eqnarray}
y(\varsigma )&=&  Hl\left( 1-\rho ^{-2}, -(2-\rho ^{-2})(h_j +2j)^2; \frac{\alpha }{2}+\frac{1}{2}, -\frac{\alpha }{2}, \frac{1}{2}, \frac{1}{2}; \varsigma \right)\nonumber\\
&=& \;_2F_1\left( -h_0, h_0; \frac{1}{2}; \eta\right) + \sum_{n=1}^{\infty } \Bigg\{\prod _{k=0}^{n-1} \Bigg\{ \int_{0}^{1} dt_{n-k}\;t_{n-k}^{2(n-k)-1 } \int_{0}^{1} du_{n-k}\;u_{n-k}^{2(n-k)-\frac{3}{2} } \nonumber\\
&&\times  \frac{1}{2\pi i}  \oint dv_{n-k} \frac{1}{v_{n-k}} \left( \frac{v_{n-k}-1}{v_{n-k}} \frac{1}{1-\overleftrightarrow {w}_{n-k+1,n}(1-t_{n-k})(1-u_{n-k})v_{n-k}}\right)^{h_{n-k}} \nonumber\\
&&\times \left( 1- \overleftrightarrow {w}_{n-k+1,n}(1-t_{n-k})(1-u_{n-k})v_{n-k}\right)^{-4(n-k)}\nonumber\\
&&\times \overleftrightarrow {w}_{n-k,n}^{-(2(n-k)-\frac{3}{2}+\frac{\alpha }{2} )}\left(  \overleftrightarrow {w}_{n-k,n} \partial _{ \overleftrightarrow {w}_{n-k,n}}\right) \overleftrightarrow {w}_{n-k,n}^{\alpha +\frac{1}{2}} \left(  \overleftrightarrow {w}_{n-k,n} \partial _{ \overleftrightarrow {w}_{n-k,n}}\right) \overleftrightarrow {w}_{n-k,n}^{2(n-k-1)-\frac{\alpha }{2} } \Bigg\}\nonumber\\
&&\times \;_2F_1\left( -h_0, h_0; \frac{1}{2}; \overleftrightarrow {w}_{1,n}\right) \Bigg\} z^n  \label{eq:80045}
\end{eqnarray}
where
\begin{equation}
h=  \rho ^2 \left( \alpha (\alpha +1)- 4(2-\rho ^{-2})(h_j+2j)^2 \right)  \nonumber 
\end{equation}
\subsubsection{Infinite series}
Replace coefficients $a$, $q$, $\gamma $, $\delta$, $x$, $c_0$ and $\lambda $ by $1-a$, $-q +\alpha \beta $, $\delta $, $\gamma $, $1-x$, 1 and zero into (\ref{eq:80016}). Put (\ref{eq:8006}) into the new (\ref{eq:80016}). 
\begin{eqnarray}
y(\varsigma )&=&  Hl\left( 1-\rho ^{-2}, \frac{1}{4}\left( h\rho ^{-2}- \alpha (\alpha +1)\right); \frac{\alpha }{2}+\frac{1}{2}, -\frac{\alpha }{2}, \frac{1}{2}, \frac{1}{2}; \varsigma \right)\nonumber\\
&=&  \;_2F_1\left( -\Pi_0 ^{+}, -\Pi_0 ^{-}; \frac{1}{2}; \eta\right) + \sum_{n=1}^{\infty } \left\{\prod _{k=0}^{n-1} \Bigg\{ \int_{0}^{1} dt_{n-k}\;t_{n-k}^{2(n-k)-1 } \int_{0}^{1} du_{n-k}\;u_{n-k}^{2(n-k)-\frac{3}{2} } \right.\nonumber\\
&\times& \frac{1}{2\pi i}  \oint dv_{n-k} \frac{1}{v_{n-k}}\left( \frac{v_{n-k}-1}{v_{n-k}}\right)^{\Pi_{n-k} ^{+}}  \left( 1- \overleftrightarrow {w}_{n-k+1,n}(1-t_{n-k})(1-u_{n-k})v_{n-k}\right)^{\Pi_{n-k} ^{-}}\nonumber\\
&\times& \overleftrightarrow {w}_{n-k,n}^{-(2(n-k)-\frac{3}{2}+\frac{\alpha }{2})}\left(  \overleftrightarrow {w}_{n-k,n} \partial _{ \overleftrightarrow {w}_{n-k,n}}\right) \overleftrightarrow {w}_{n-k,n}^{\alpha +\frac{1}{2}} \left(  \overleftrightarrow {w}_{n-k,n} \partial _{ \overleftrightarrow {w}_{n-k,n}}\right) \overleftrightarrow {w}_{n-k,n}^{2(n-k-1) -\frac{\alpha }{2}  } \Bigg\}\nonumber\\
&\times& \left.\left. \;_2F_1\left( -\Pi_0 ^{+}, -\Pi_0 ^{-}; \frac{1}{2}; \overleftrightarrow {w}_{1,n}\right) \right\} z^n \right\} \label{eq:80046}
\end{eqnarray}
where
\begin{equation}
\begin{cases} 
\Pi_0 ^{\pm}=  \pm\frac{1}{2}\sqrt{\frac{h\rho ^{-2}+\alpha (\alpha +1)}{ 1-\rho ^{-2}}}  \cr
\Pi_{n-k} ^{\pm}=  -2(n-k) \pm\frac{1}{2}\sqrt{\frac{h\rho ^{-2}+\alpha (\alpha +1)}{ 1-\rho ^{-2}}}
\end{cases}\nonumber 
\end{equation}
On (\ref{eq:80045}) and (\ref{eq:80046}),
\begin{equation}
\begin{cases} 
\varsigma= 1-\xi \cr
\eta =\frac{2-\rho ^{-2}}{1-\rho ^{-2}}\varsigma \cr
z=\frac{-1}{1-\rho ^{-2}}\varsigma ^2  
\end{cases}\nonumber 
\end{equation}
\subsection{ \footnotesize ${\displaystyle (1-x)^{1-\delta } Hl(1-a,-q+(\delta -1)\gamma a+(\alpha -\delta +1)(\beta -\delta +1); \alpha-\delta +1,\beta-\delta +1, 2-\delta, \gamma; 1-x)}$ \normalsize}
\subsubsection{Polynomial of type 2}
Replace coefficients $a$, $q$, $\alpha $, $\beta $, $\gamma $, $\delta$, $x$, $c_0$, $\lambda $ and $q_j$ where $j, q_j \in \mathbb{N}_{0}$ by $1-a$, $-q+(\delta -1)\gamma a+(\alpha -\delta +1)(\beta -\delta +1)$, $\alpha-\delta +1 $, $\beta-\delta +1 $, $2-\delta$, $\gamma $, $1-x$, 1, zero and $h_j$ where $h_j \in \mathbb{N}_{0}$ into (\ref{eq:80013}). Multiply $(1-x)^{1-\delta }$ and (\ref{eq:80013}) together. Put (\ref{eq:8006}) into the new (\ref{eq:80013}). 
\begin{eqnarray}
&&\varsigma ^{\frac{1}{2}}y(\varsigma )\nonumber\\
&=& \varsigma ^{\frac{1}{2}} Hl\left( 1-\rho ^{-2}, -(2-\rho ^{-2})(h_j+2j)(h_j+2j+1); \frac{\alpha }{2}+1, -\frac{\alpha }{2}+\frac{1}{2}, \frac{3}{2}, \frac{1}{2}; \varsigma \right)\nonumber\\
&=& \varsigma ^{\frac{1}{2}} \left\{  \;_2F_1\left( -h_0, h_0+1; \frac{3}{2}; \eta\right) \right. + \sum_{n=1}^{\infty } \Bigg\{\prod _{k=0}^{n-1} \Bigg\{ \int_{0}^{1} dt_{n-k}\;t_{n-k}^{2(n-k)-1 } \int_{0}^{1} du_{n-k}\;u_{n-k}^{2(n-k)-\frac{1}{2} } \nonumber\\
&&\times  \frac{1}{2\pi i}  \oint dv_{n-k} \frac{1}{v_{n-k}} \left( \frac{v_{n-k}-1}{v_{n-k}} \frac{1}{1-\overleftrightarrow {w}_{n-k+1,n}(1-t_{n-k})(1-u_{n-k})v_{n-k}}\right)^{h_{n-k}} \nonumber\\
&&\times \left( 1- \overleftrightarrow {w}_{n-k+1,n}(1-t_{n-k})(1-u_{n-k})v_{n-k}\right)^{-4(n-k)-1}\nonumber\\
&&\times \overleftrightarrow {w}_{n-k,n}^{-(2(n-k)-1+\frac{\alpha }{2} )}\left(  \overleftrightarrow {w}_{n-k,n} \partial _{ \overleftrightarrow {w}_{n-k,n}}\right) \overleftrightarrow {w}_{n-k,n}^{\alpha +\frac{1}{2}} \left(  \overleftrightarrow {w}_{n-k,n} \partial _{ \overleftrightarrow {w}_{n-k,n}}\right) \overleftrightarrow {w}_{n-k,n}^{2(n-k)-\frac{3}{2}-\frac{\alpha }{2} } \Bigg\}\nonumber\\
&&\times \;_2F_1\left( -h_0, h_0+1; \frac{3}{2}; \overleftrightarrow {w}_{1,n}\right) \Bigg\} z^n \Bigg\} \label{eq:80047}
\end{eqnarray}
where
\begin{equation} 
h= 4(2\rho ^2-1)(h_j +2j)(h_j +2j+1)-\rho ^2(\alpha -1)(\alpha +2) -1  \nonumber 
\end{equation}
\subsubsection{Infinite series}
Replace coefficients $a$, $q$, $\alpha $, $\beta $, $\gamma $, $\delta$, $x$, $c_0$ and $\lambda $ by $1-a$, $-q+(\delta -1)\gamma a+(\alpha -\delta +1)(\beta -\delta +1)$, $\alpha-\delta +1 $, $\beta-\delta +1 $, $2-\delta$, $\gamma $, $1-x$, 1 and zero into (\ref{eq:80016}). Multiply $(1-x)^{1-\delta }$ and (\ref{eq:80016}) together. Put (\ref{eq:8006}) into the new (\ref{eq:80016}).
\begin{eqnarray}
&&\varsigma ^{\frac{1}{2}}y(\varsigma )\nonumber\\
&=&\varsigma ^{\frac{1}{2}} Hl\left( 1-\rho ^{-2}, -\frac{1}{4}\left( (1+h)\rho ^{-2}+(\alpha -1)(\alpha +2)\right); \frac{\alpha }{2}+1, -\frac{\alpha }{2}+\frac{1}{2}, \frac{3}{2}, \frac{1}{2}; \varsigma \right)\nonumber\\
&=& \varsigma ^{\frac{1}{2}} \left\{ \;_2F_1\left( -\Pi_0 ^{+}, -\Pi_0 ^{-}; \frac{3}{2}; \eta\right)\right. + \sum_{n=1}^{\infty } \left\{\prod _{k=0}^{n-1} \Bigg\{ \int_{0}^{1} dt_{n-k}\;t_{n-k}^{2(n-k)-1 } \int_{0}^{1} du_{n-k}\;u_{n-k}^{2(n-k)-\frac{1}{2} } \right.\nonumber\\
&\times& \frac{1}{2\pi i}  \oint dv_{n-k} \frac{1}{v_{n-k}}\left( \frac{v_{n-k}-1}{v_{n-k}}\right)^{\Pi_{n-k} ^{+}}  \left( 1- \overleftrightarrow {w}_{n-k+1,n}(1-t_{n-k})(1-u_{n-k})v_{n-k}\right)^{\Pi_{n-k} ^{-}}\nonumber\\
&\times& \overleftrightarrow {w}_{n-k,n}^{-(2(n-k)-1+\frac{\alpha }{2})}\left(  \overleftrightarrow {w}_{n-k,n} \partial _{ \overleftrightarrow {w}_{n-k,n}}\right) \overleftrightarrow {w}_{n-k,n}^{\alpha +\frac{1}{2}} \left(  \overleftrightarrow {w}_{n-k,n} \partial _{ \overleftrightarrow {w}_{n-k,n}}\right) \overleftrightarrow {w}_{n-k,n}^{2(n-k)-\frac{3}{2} -\frac{\alpha }{2}  } \Bigg\}\nonumber\\
&\times& \left.\left. \;_2F_1\left( -\Pi_0 ^{+}, -\Pi_0 ^{-}; \frac{3}{2}; \overleftrightarrow {w}_{1,n}\right) \right\} z^n \right\} \label{eq:80048}
\end{eqnarray}
where
\begin{equation}
\begin{cases} 
\Pi_0 ^{\pm}= -\frac{1}{2} \pm\frac{1}{2}\sqrt{\frac{h\rho ^{-2}+\alpha (\alpha +1)}{ 2-\rho ^{-2}}}  \cr
\Pi_{n-k} ^{\pm}=  -2(n-k) -\frac{1}{2} \pm\frac{1}{2}\sqrt{\frac{h\rho ^{-2}+\alpha (\alpha +1)}{ 2-\rho ^{-2}}}
\end{cases}\nonumber 
\end{equation}
On (\ref{eq:80047}) and (\ref{eq:80048}),
\begin{equation}
\begin{cases} 
\varsigma= 1-\xi \cr
\eta =\frac{2-\rho ^{-2}}{1-\rho ^{-2}}\varsigma \cr
z=\frac{-1}{1-\rho ^{-2}}\varsigma ^2  
\end{cases}\nonumber 
\end{equation}
\subsection{ ${\displaystyle x^{-\alpha } Hl\left(\frac{1}{a},\frac{q+\alpha [(\alpha -\gamma -\delta +1)a-\beta +\delta ]}{a}; \alpha , \alpha -\gamma +1, \alpha -\beta +1,\delta ;\frac{1}{x}\right)}$}
\subsubsection{Polynomial of type 2}
Replace coefficients $a$, $q$, $\beta $, $\gamma $, $x$, $c_0$, $\lambda $ and $q_j$ where $j, q_j \in \mathbb{N}_{0}$ by $\frac{1}{a}$, $\frac{q+\alpha [(\alpha -\gamma -\delta +1)a-\beta +\delta ]}{a}$, $\alpha-\gamma +1 $, $\alpha -\beta +1 $, $\frac{1}{x}$, 1, zero and  $h_j$ where $h_j \in \mathbb{N}_{0}$ into (\ref{eq:80013}). Multiply $x^{-\alpha }$ and (\ref{eq:80013}) together. Put (\ref{eq:8006}) into the new (\ref{eq:80013}). 
\begin{eqnarray}
&&\varsigma ^{\frac{1}{2}(\alpha +1)} y(\varsigma )\nonumber\\
&=& \varsigma ^{\frac{1}{2}(\alpha +1)} Hl\left(\rho ^2, -(1+\rho ^2)(h_j+2j)(h_j+2j+1+\alpha ); \frac{\alpha}{2} +\frac{1}{2}, \frac{\alpha}{2} +1, \alpha +\frac{3}{2}, \frac{1}{2}; \varsigma \right) \nonumber\\
&=& \varsigma ^{\frac{1}{2}(\alpha +1)} \left\{  \;_2F_1\left( -h_0, h_0+\alpha+1; \alpha+\frac{3}{2}; \eta\right) \right. \nonumber\\
&&+ \sum_{n=1}^{\infty } \Bigg\{\prod _{k=0}^{n-1} \Bigg\{ \int_{0}^{1} dt_{n-k}\;t_{n-k}^{2(n-k)-1 } \int_{0}^{1} du_{n-k}\;u_{n-k}^{2(n-k)-\frac{1}{2} +\alpha } \nonumber\\
&&\times  \frac{1}{2\pi i}  \oint dv_{n-k} \frac{1}{v_{n-k}} \left( \frac{v_{n-k}-1}{v_{n-k}} \frac{1}{1-\overleftrightarrow {w}_{n-k+1,n}(1-t_{n-k})(1-u_{n-k})v_{n-k}}\right)^{h_{n-k}} \nonumber\\
&&\times \left( 1- \overleftrightarrow {w}_{n-k+1,n}(1-t_{n-k})(1-u_{n-k})v_{n-k}\right)^{-(4(n-k)+1+\alpha) }\nonumber\\
&&\times \overleftrightarrow {w}_{n-k,n}^{-(2(n-k)-\frac{3}{2}+\frac{\alpha }{2} )}\left(  \overleftrightarrow {w}_{n-k,n} \partial _{ \overleftrightarrow {w}_{n-k,n}}\right) \overleftrightarrow {w}_{n-k,n}^{-\frac{1}{2}} \left(  \overleftrightarrow {w}_{n-k,n} \partial _{ \overleftrightarrow {w}_{n-k,n}}\right) \overleftrightarrow {w}_{n-k,n}^{2(n-k)-1+\frac{\alpha }{2} } \Bigg\}\nonumber\\
&&\times \;_2F_1\left( -h_0, h_0+\alpha+1; \alpha+\frac{3}{2}; \overleftrightarrow {w}_{1,n}\right) \Bigg\} z^n \Bigg\} \label{eq:80049}
\end{eqnarray}
where
\begin{equation}
h= (1+ \rho ^2)\left( 2(h_j +2j)+1+\alpha \right)^2  \nonumber 
\end{equation}
\subsubsection{Infinite series}
Replace coefficients $a$, $q$, $\beta $, $\gamma $, $x$, $c_0$ and $\lambda $ by $\frac{1}{a}$, $\frac{q+\alpha [(\alpha -\gamma -\delta +1)a-\beta +\delta ]}{a}$, $\alpha-\gamma +1 $, $\alpha -\beta +1 $, $\frac{1}{x}$, 1 and zero into (\ref{eq:80016}). Multiply $x^{-\alpha }$ and (\ref{eq:80016}) together. Put (\ref{eq:8006}) into the new (\ref{eq:80016}). 
\begin{eqnarray}
&&\varsigma ^{\frac{1}{2}(\alpha +1)} y(\varsigma )\nonumber\\
&=& \varsigma ^{\frac{1}{2}(\alpha +1)} Hl\left(\rho ^2,-\frac{1}{4}\left( h-(1+\rho ^2)(\alpha +1)^2\right); \frac{\alpha}{2} +\frac{1}{2}, \frac{\alpha}{2} +1, \alpha +\frac{3}{2}, \frac{1}{2}; \varsigma \right) \nonumber\\
&=& \varsigma ^{\frac{1}{2}(\alpha +1)} \left\{ \;_2F_1\left( -\Pi_0 ^{+}, -\Pi_0 ^{-}; \alpha +\frac{3}{2}; \eta\right)\right. + \sum_{n=1}^{\infty } \left\{\prod _{k=0}^{n-1} \Bigg\{ \int_{0}^{1} dt_{n-k}\;t_{n-k}^{2(n-k)-1 } \int_{0}^{1} du_{n-k}\;u_{n-k}^{2(n-k)-\frac{1}{2} +\alpha } \right.\nonumber\\
&\times& \frac{1}{2\pi i}  \oint dv_{n-k} \frac{1}{v_{n-k}}\left( \frac{v_{n-k}-1}{v_{n-k}}\right)^{\Pi_{n-k} ^{+}}  \left( 1- \overleftrightarrow {w}_{n-k+1,n}(1-t_{n-k})(1-u_{n-k})v_{n-k}\right)^{\Pi_{n-k} ^{-}}\nonumber\\
&\times& \overleftrightarrow {w}_{n-k,n}^{-(2(n-k)-\frac{3}{2}+\frac{\alpha }{2})}\left(  \overleftrightarrow {w}_{n-k,n} \partial _{ \overleftrightarrow {w}_{n-k,n}}\right) \overleftrightarrow {w}_{n-k,n}^{-\frac{1}{2}} \left(  \overleftrightarrow {w}_{n-k,n} \partial _{ \overleftrightarrow {w}_{n-k,n}}\right) \overleftrightarrow {w}_{n-k,n}^{2(n-k)-1 +\frac{\alpha }{2}  } \Bigg\}\nonumber\\
&\times& \left.\left. \;_2F_1\left( -\Pi_0 ^{+}, -\Pi_0 ^{-}; \alpha +\frac{3}{2}; \overleftrightarrow {w}_{1,n}\right) \right\} z^n \right\} \label{eq:80050}
\end{eqnarray}
where
\begin{equation}
\begin{cases} 
\Pi_0 ^{\pm}= -\frac{1}{2}(\alpha +1) \pm\frac{1}{2}\sqrt{\frac{h }{ 1+\rho ^2}}  \cr
\Pi_{n-k} ^{\pm}=  -2(n-k) -\frac{1}{2}(\alpha +1) \pm\frac{1}{2}\sqrt{\frac{h }{ 1+\rho ^2}}
\end{cases}\nonumber 
\end{equation}
On (\ref{eq:80049}) and (\ref{eq:80050}),
\begin{equation}
\begin{cases} 
\varsigma= \xi^{-1} \cr
\eta =(1+\rho ^{-2})\varsigma \cr
z= -\rho ^{-2} \varsigma ^2  
\end{cases}\nonumber 
\end{equation}
\subsection{ ${\displaystyle \left(1-\frac{x}{a} \right)^{-\beta } Hl\left(1-a, -q+\gamma \beta; -\alpha +\gamma +\delta, \beta, \gamma, \delta; \frac{(1-a)x}{x-a} \right)}$}
\subsubsection{Polynomial of type 2}
Replace coefficients $a$, $q$, $\alpha $, $x$, $c_0$, $\lambda $ and $q_j$ where $j, q_j \in \mathbb{N}_{0}$ by $1-a$, $-q+\gamma \beta $, $-\alpha+\gamma +\delta $, $\frac{(1-a)x}{x-a}$, 1, zero and $h_j$ where $h_j \in \mathbb{N}_{0}$ into (\ref{eq:80013}). Multiply $\left(1-\frac{x}{a} \right)^{-\beta }$ and (\ref{eq:80013}) together. Put (\ref{eq:8006}) into the new (\ref{eq:80013}). 
\begin{eqnarray} 
&&(1-\rho ^2 \xi)^{\frac{\alpha }{2}} y(\varsigma )\nonumber\\
&=& (1-\rho ^2 \xi)^{\frac{\alpha }{2}} Hl\left( 1-\rho ^{-2}, (h_j+2j)\left(\alpha -(2-\rho ^{-2})(h_j+2j) \right); -\frac{\alpha }{2}+\frac{1}{2}, -\frac{\alpha }{2}, \frac{1}{2}, \frac{1}{2}; \varsigma \right) \nonumber\\
&=& (1-\rho ^2 \xi)^{\frac{\alpha }{2}} \left\{\;_2F_1\left( -h_0,h_0- \Omega _{\rho }; \frac{1}{2}; \eta\right)\right.   + \sum_{n=1}^{\infty } \Bigg\{\prod _{k=0}^{n-1} \Bigg\{ \int_{0}^{1} dt_{n-k}\;t_{n-k}^{2(n-k)-1 } \int_{0}^{1} du_{n-k}\;u_{n-k}^{2(n-k)-\frac{3}{2} } \nonumber\\
&&\times  \frac{1}{2\pi i}  \oint dv_{n-k} \frac{1}{v_{n-k}} \left( \frac{v_{n-k}-1}{v_{n-k}} \frac{1}{1-\overleftrightarrow {w}_{n-k+1,n}(1-t_{n-k})(1-u_{n-k})v_{n-k}}\right)^{h_{n-k}} \nonumber\\
&&\times \left( 1- \overleftrightarrow {w}_{n-k+1,n}(1-t_{n-k})(1-u_{n-k})v_{n-k}\right)^{-\left( 4(n-k)-\Omega _{\rho }\right)}\nonumber\\
&&\times \overleftrightarrow {w}_{n-k,n}^{-(2(n-k)-\frac{3}{2}-\frac{\alpha }{2} )}\left(  \overleftrightarrow {w}_{n-k,n} \partial _{ \overleftrightarrow {w}_{n-k,n}}\right) \overleftrightarrow {w}_{n-k,n}^{ \frac{1}{2}} \left(  \overleftrightarrow {w}_{n-k,n} \partial _{ \overleftrightarrow {w}_{n-k,n}}\right) \overleftrightarrow {w}_{n-k,n}^{2(n-k-1) -\frac{\alpha }{2} } \Bigg\}\nonumber\\
&&\times \;_2F_1\left( -h_0, h_0-\Omega _{\rho }; \frac{1}{2}; \overleftrightarrow {w}_{1,n}\right) \Bigg\} z^n \Bigg\} \label{eq:80051}
\end{eqnarray}
where
\begin{equation}
\begin{cases} 
\Omega _{\rho } = \frac{\alpha }{2-\rho ^{-2}} \cr
h= 4\rho ^2 \left( \alpha \left( h_j+2j+\frac{1}{4} \right) -(2-\rho ^{-2})(h_j+2j)^2\right)  
\end{cases}\nonumber 
\end{equation}
\subsubsection{Infinite series}
Replace coefficients $a$, $q$, $\alpha $, $x$, $c_0$ and $\lambda $ by $1-a$, $-q+\gamma \beta $, $-\alpha+\gamma +\delta $, $\frac{(1-a)x}{x-a}$, 1 and zero into (\ref{eq:80016}). Multiply $\left(1-\frac{x}{a} \right)^{-\beta }$ and (\ref{eq:80016}) together. Put (\ref{eq:8006}) into the new (\ref{eq:80016}). 
\begin{eqnarray}
&&(1-\rho ^2 \xi)^{\frac{\alpha }{2}} y(\varsigma )\nonumber\\
&=& (1-\rho ^2 \xi)^{\frac{\alpha }{2}} Hl\left( 1-\rho ^{-2}, \frac{1}{4}\left( h\rho ^{-2} - \alpha \right); -\frac{\alpha }{2}+\frac{1}{2}, -\frac{\alpha }{2}, \frac{1}{2}, \frac{1}{2}; \varsigma \right) \nonumber\\
&=& (1-\rho ^2 \xi)^{\frac{\alpha }{2}}  \left\{\;_2F_1\left( -\Pi_0 ^{+}, -\Pi_0 ^{-}; \frac{1}{2}; \eta\right)  \right. + \sum_{n=1}^{\infty } \left\{\prod _{k=0}^{n-1} \Bigg\{ \int_{0}^{1} dt_{n-k}\;t_{n-k}^{2(n-k)-1 } \int_{0}^{1} du_{n-k}\;u_{n-k}^{2(n-k)-\frac{3}{2} } \right.\nonumber\\
&\times& \frac{1}{2\pi i}  \oint dv_{n-k} \frac{1}{v_{n-k}}\left( \frac{v_{n-k}-1}{v_{n-k}}\right)^{\Pi_{n-k} ^{+}}  \left( 1- \overleftrightarrow {w}_{n-k+1,n}(1-t_{n-k})(1-u_{n-k})v_{n-k}\right)^{\Pi_{n-k} ^{-}}\nonumber\\
&\times& \overleftrightarrow {w}_{n-k,n}^{-(2(n-k)-\frac{3}{2}-\frac{\alpha }{2})}\left(  \overleftrightarrow {w}_{n-k,n} \partial _{ \overleftrightarrow {w}_{n-k,n}}\right) \overleftrightarrow {w}_{n-k,n}^{ \frac{1}{2}} \left(  \overleftrightarrow {w}_{n-k,n} \partial _{ \overleftrightarrow {w}_{n-k,n}}\right) \overleftrightarrow {w}_{n-k,n}^{2(n-k-1) -\frac{\alpha }{2}  } \Bigg\}\nonumber\\
&\times& \left.\left. \;_2F_1\left( -\Pi_0 ^{+}, -\Pi_0 ^{-}; \frac{1}{2}; \overleftrightarrow {w}_{1,n}\right) \right\} z^n \right\} \label{eq:80052}
\end{eqnarray}
where
\begin{equation}
\begin{cases} 
\Pi_0 ^{\pm}=  \frac{\alpha \pm \sqrt{\alpha ^2-(2-\rho ^{-2})(h\rho ^{-2}-\alpha )}}{2(2-\rho ^{-2})}  \cr
\Pi_{n-k} ^{\pm}=  -2(n-k) +\frac{\alpha \pm \sqrt{\alpha ^2-(2-\rho ^{-2})(h\rho ^{-2}-\alpha )}}{2(2-\rho ^{-2})}
\end{cases}\nonumber 
\end{equation}
On (\ref{eq:80051}) and (\ref{eq:80052}),
\begin{equation}
\begin{cases} 
\varsigma= \frac{(1-\rho ^{-2})\xi}{\xi-\rho ^{-2}} \cr
\eta =\frac{2-\rho ^{-2}}{1-\rho ^{-2}}\varsigma \cr
z= -\frac{1}{1-\rho ^{-2}} \varsigma ^2 
\end{cases}\nonumber 
\end{equation}
\subsection{ \footnotesize ${\displaystyle (1-x)^{1-\delta }\left(1-\frac{x}{a} \right)^{-\beta+\delta -1} Hl\left(1-a, -q+\gamma [(\delta -1)a+\beta -\delta +1]; -\alpha +\gamma +1, \beta -\delta+1, \gamma, 2-\delta; \frac{(1-a)x}{x-a} \right)}$ \normalsize}
\subsubsection{Polynomial of type 2}
Replace coefficients $a$, $q$, $\alpha $, $\beta $, $\delta $, $x$, $c_0$, $\lambda $ and $q_j$ where $j, q_j \in \mathbb{N}_{0}$ by $1-a$, $-q+\gamma [(\delta -1)a+\beta -\delta +1]$, $-\alpha +\gamma +1$, $\beta -\delta+1$, $2-\delta $, $\frac{(1-a)x}{x-a}$, 1, zero and $h_j$ where $h_j \in \mathbb{N}_{0}$ into (\ref{eq:80013}). Multiply $(1-x)^{1-\delta }\left(1-\frac{x}{a} \right)^{-\beta+\delta -1}$ and (\ref{eq:80013}) together. Put (\ref{eq:8006}) into the new (\ref{eq:80013}).
\begin{eqnarray} 
&&(1-\xi )^{\frac{1}{2}}(1-\rho ^2 \xi)^{\frac{1}{2}(\alpha -1)} y(\varsigma )\nonumber\\
&=& (1-\xi )^{\frac{1}{2}}(1-\rho ^2 \xi)^{\frac{1}{2}(\alpha -1)} Hl\bigg( 1-\rho ^{-2}, (h_j+ 2j)\left( \alpha +1 -(2-\rho ^{-2})(h_j+ 2j+1)\right); -\frac{\alpha }{2}+1, \nonumber\\
&& -\frac{\alpha }{2}+\frac{1}{2}, \frac{1}{2}, \frac{3}{2}; \varsigma \bigg) \nonumber\\
&=& (1-\xi )^{\frac{1}{2}}(1-\rho ^2 \xi)^{\frac{1}{2}(\alpha -1)} \left\{_2F_1\left( -h_0,h_0+ \Omega _{\rho }; \frac{1}{2}; \eta\right)\right.  \nonumber\\
&&+ \sum_{n=1}^{\infty } \Bigg\{\prod _{k=0}^{n-1} \Bigg\{ \int_{0}^{1} dt_{n-k}\;t_{n-k}^{2(n-k)-1 } \int_{0}^{1} du_{n-k}\;u_{n-k}^{2(n-k)-\frac{3}{2} } \nonumber\\
&&\times  \frac{1}{2\pi i}  \oint dv_{n-k} \frac{1}{v_{n-k}} \left( \frac{v_{n-k}-1}{v_{n-k}} \frac{1}{1-\overleftrightarrow {w}_{n-k+1,n}(1-t_{n-k})(1-u_{n-k})v_{n-k}}\right)^{h_{n-k}} \nonumber\\
&&\times \left( 1- \overleftrightarrow {w}_{n-k+1,n}(1-t_{n-k})(1-u_{n-k})v_{n-k}\right)^{-\left( 4(n-k)+\Omega _{\rho }\right)}\nonumber\\
&&\times \overleftrightarrow {w}_{n-k,n}^{-(2(n-k)-1-\frac{\alpha }{2} )}\left(  \overleftrightarrow {w}_{n-k,n} \partial _{ \overleftrightarrow {w}_{n-k,n}}\right) \overleftrightarrow {w}_{n-k,n}^{ \frac{1}{2}} \left(  \overleftrightarrow {w}_{n-k,n} \partial _{ \overleftrightarrow {w}_{n-k,n}}\right) \overleftrightarrow {w}_{n-k,n}^{2(n-k)-\frac{3}{2} -\frac{\alpha }{2} } \Bigg\}\nonumber\\
&&\times \;_2F_1\left( -h_0, h_0+\Omega _{\rho }; \frac{1}{2}; \overleftrightarrow {w}_{1,n}\right) \Bigg\} z^n \Bigg\} \label{eq:80053}
\end{eqnarray}
where
\begin{equation}
\begin{cases} 
\Omega _{\rho } = \frac{-\alpha +1-\rho ^{-2}}{2-\rho ^{-2}} \cr
h= 4(h_j+ 2j)\left( (\alpha +1)\rho ^2 +(1-2\rho ^2)(h_j+2j+1)\right) +(\alpha -1)\rho ^2 +1  
\end{cases}\nonumber 
\end{equation}
\subsubsection{Infinite series}
Replace coefficients $a$, $q$, $\alpha $, $\beta $, $\delta $, $x$, $c_0$ and $\lambda $ by $1-a$, $-q+\gamma [(\delta -1)a+\beta -\delta +1]$, $-\alpha +\gamma +1$, $\beta -\delta+1$, $2-\delta $, $\frac{(1-a)x}{x-a}$, 1 and zero into (\ref{eq:80016}). Multiply $(1-x)^{1-\delta }\left(1-\frac{x}{a} \right)^{-\beta+\delta -1}$ and (\ref{eq:80016}) together. Put (\ref{eq:8006}) into the new (\ref{eq:80016}).
\begin{eqnarray}
&&(1-\xi )^{\frac{1}{2}}(1-\rho ^2 \xi)^{\frac{1}{2}(\alpha -1)} y(\varsigma )\nonumber\\
&=& (1-\xi )^{\frac{1}{2}}(1-\rho ^2 \xi)^{\frac{1}{2}(\alpha -1)} Hl\left( 1-\rho ^{-2}, \frac{1}{4}\left( (h-1)\rho ^{-2} +1- \alpha \right); -\frac{\alpha }{2}+1, -\frac{\alpha }{2}+\frac{1}{2}, \frac{1}{2}, \frac{3}{2}; \varsigma \right) \nonumber\\
&=& (1-\xi )^{\frac{1}{2}}(1-\rho ^2 \xi)^{\frac{1}{2}(\alpha -1)} \left\{\;_2F_1\left( -\Pi_0 ^{+}, -\Pi_0 ^{-}; \frac{1}{2}; \eta\right) \right.  \nonumber\\
&+& \sum_{n=1}^{\infty } \left\{\prod _{k=0}^{n-1} \Bigg\{ \int_{0}^{1} dt_{n-k}\;t_{n-k}^{2(n-k)-1 } \int_{0}^{1} du_{n-k}\;u_{n-k}^{2(n-k)-\frac{3}{2} } \right.\nonumber\\
&\times& \frac{1}{2\pi i}  \oint dv_{n-k} \frac{1}{v_{n-k}}\left( \frac{v_{n-k}-1}{v_{n-k}}\right)^{\Pi_{n-k} ^{+}}  \left( 1- \overleftrightarrow {w}_{n-k+1,n}(1-t_{n-k})(1-u_{n-k})v_{n-k}\right)^{\Pi_{n-k} ^{-}}\nonumber\\
&\times& \overleftrightarrow {w}_{n-k,n}^{-(2(n-k)-1-\frac{\alpha }{2})}\left(  \overleftrightarrow {w}_{n-k,n} \partial _{ \overleftrightarrow {w}_{n-k,n}}\right) \overleftrightarrow {w}_{n-k,n}^{ \frac{1}{2}} \left(  \overleftrightarrow {w}_{n-k,n} \partial _{ \overleftrightarrow {w}_{n-k,n}}\right) \overleftrightarrow {w}_{n-k,n}^{2(n-k) -\frac{3}{2}-\frac{\alpha }{2}  } \Bigg\}\nonumber\\
&\times& \left.\left. \;_2F_1\left( -\Pi_0 ^{+}, -\Pi_0 ^{-}; \frac{1}{2}; \overleftrightarrow {w}_{1,n}\right) \right\} z^n \right\} \label{eq:80054}
\end{eqnarray}
where
\begin{equation}
\begin{cases} 
\Pi_0 ^{\pm}=  \frac{\alpha -1+\rho ^{-2}\pm \sqrt{(\alpha -1+\rho ^{-2})^2 +(2-\rho ^{-2})((1-h)\rho ^{-2}+\alpha -1)}}{2(2-\rho ^{-2})}  \cr
\Pi_{n-k} ^{\pm}=  -2(n-k) + \frac{\alpha -1+\rho ^{-2}\pm \sqrt{(\alpha -1+\rho ^{-2})^2 +(2-\rho ^{-2})((1-h)\rho ^{-2}+\alpha -1 )}}{2(2-\rho ^{-2})}
\end{cases}\nonumber 
\end{equation}
On (\ref{eq:80053}) and (\ref{eq:80054}),
\begin{equation}
\begin{cases} 
\varsigma= \frac{(1-\rho ^{-2})\xi}{\xi-\rho ^{-2}} \cr
\eta =\frac{2-\rho ^{-2}}{1-\rho ^{-2}}\varsigma \cr
z= -\frac{1}{1-\rho ^{-2}} \varsigma ^2 
\end{cases}\nonumber 
\end{equation}
\subsection{ ${\displaystyle x^{-\alpha } Hl\left(\frac{a-1}{a}, \frac{-q+\alpha (\delta a+\beta -\delta )}{a}; \alpha, \alpha -\gamma +1, \delta , \alpha -\beta +1; \frac{x-1}{x} \right)}$}
\subsubsection{Polynomial of type 2}
Replace coefficients $a$, $q$, $\beta $, $\gamma $, $\delta $, $x$, $c_0$, $\lambda $ and $q_j$ where $j, q_j \in \mathbb{N}_{0}$ by $\frac{a-1}{a}$, $\frac{-q+\alpha (\delta a+\beta -\delta )}{a}$, $\alpha -\gamma +1$, $\delta $, $\alpha -\beta +1$, $\frac{x-1}{x}$, 1, zero and $h_j$ where $h_j \in \mathbb{N}_{0}$ into (\ref{eq:80013}). Multiply $x^{-\alpha }$ and (\ref{eq:80013}) together. Put (\ref{eq:8006}) into the new (\ref{eq:80013}).
\begin{eqnarray} 
&&\xi ^{-\frac{1}{2}(\alpha +1)} y(\varsigma )\nonumber\\
&=& \xi ^{-\frac{1}{2}(\alpha +1)} Hl\left( 1-\rho ^2, -(h_j+2j)\left( (1-\rho ^2)(\alpha +1)+(2-\rho ^2)(h_j+2j) \right); \frac{\alpha}{2} + \frac{1}{2}, \frac{\alpha}{2} +1, \frac{1}{2}, \alpha +\frac{3}{2}; \varsigma \right) \nonumber\\
&=& \xi ^{-\frac{1}{2}(\alpha +1)} \left\{ _2F_1\left( -h_0,h_0+ \Omega _{\rho }; \frac{1}{2}; \eta\right)\right.   + \sum_{n=1}^{\infty } \Bigg\{\prod _{k=0}^{n-1} \Bigg\{ \int_{0}^{1} dt_{n-k}\;t_{n-k}^{2(n-k)-1 } \int_{0}^{1} du_{n-k}\;u_{n-k}^{2(n-k)-\frac{3}{2} } \nonumber\\
&&\times  \frac{1}{2\pi i}  \oint dv_{n-k} \frac{1}{v_{n-k}} \left( \frac{v_{n-k}-1}{v_{n-k}} \frac{1}{1-\overleftrightarrow {w}_{n-k+1,n}(1-t_{n-k})(1-u_{n-k})v_{n-k}}\right)^{h_{n-k}} \nonumber\\
&&\times \left( 1- \overleftrightarrow {w}_{n-k+1,n}(1-t_{n-k})(1-u_{n-k})v_{n-k}\right)^{-\left( 4(n-k)+\Omega _{\rho }\right)}\nonumber\\
&&\times \overleftrightarrow {w}_{n-k,n}^{-(2(n-k)-\frac{3}{2}+\frac{\alpha }{2} )}\left(  \overleftrightarrow {w}_{n-k,n} \partial _{ \overleftrightarrow {w}_{n-k,n}}\right) \overleftrightarrow {w}_{n-k,n}^{ -\frac{1}{2}} \left(  \overleftrightarrow {w}_{n-k,n} \partial _{ \overleftrightarrow {w}_{n-k,n}}\right) \overleftrightarrow {w}_{n-k,n}^{2(n-k)-1 +\frac{\alpha }{2} } \Bigg\}\nonumber\\
&&\times \;_2F_1\left( -h_0, h_0+\Omega _{\rho }; \frac{1}{2}; \overleftrightarrow {w}_{1,n}\right) \Bigg\} z^n \Bigg\} \label{eq:80055}
\end{eqnarray}
where
\begin{equation}
\begin{cases} 
\Omega _{\rho } = \frac{(\alpha +1)(1-\rho ^2)}{2-\rho ^2} \cr
h= -4(h_j+ 2j)\left( (1- \rho ^2)(\alpha +1) +(2-\rho ^2)(h_j+2j)\right) -(\alpha +1)\left( 1-\rho ^2(\alpha +1)\right) \cr
\end{cases}\nonumber 
\end{equation}
\subsubsection{Infinite series}
Replace coefficients $a$, $q$, $\beta $, $\gamma $, $\delta $, $x$, $c_0$ and $\lambda $ by $\frac{a-1}{a}$, $\frac{-q+\alpha (\delta a+\beta -\delta )}{a}$, $\alpha -\gamma +1$, $\delta $, $\alpha -\beta +1$, $\frac{x-1}{x}$, 1 and zero into (\ref{eq:80016}). Multiply $x^{-\alpha }$ and (\ref{eq:80016}) together. Put (\ref{eq:8006}) into the new (\ref{eq:80016}).
\begin{eqnarray}
&&\xi ^{-\frac{1}{2}(\alpha +1)} y(\varsigma )\nonumber\\
&=& \xi ^{-\frac{1}{2}(\alpha +1)} Hl\left( 1-\rho ^2, \frac{1}{4}\left[ h+ (\alpha +1)\left( 1-\rho ^2(\alpha +1)\right) \right]; \frac{\alpha}{2} + \frac{1}{2}, \frac{\alpha}{2} +1, \frac{1}{2}, \alpha +\frac{3}{2}; \varsigma \right) \nonumber\\
&=& \xi ^{-\frac{1}{2}(\alpha +1)} \left\{_2F_1\left( -\Pi_0 ^{+}, -\Pi_0 ^{-}; \frac{1}{2}; \eta\right)  \right. + \sum_{n=1}^{\infty } \left\{\prod _{k=0}^{n-1} \Bigg\{ \int_{0}^{1} dt_{n-k}\;t_{n-k}^{2(n-k)-1 } \int_{0}^{1} du_{n-k}\;u_{n-k}^{2(n-k)-\frac{3}{2} } \right.\nonumber\\
&\times& \frac{1}{2\pi i}  \oint dv_{n-k} \frac{1}{v_{n-k}}\left( \frac{v_{n-k}-1}{v_{n-k}}\right)^{\Pi_{n-k} ^{+}}  \left( 1- \overleftrightarrow {w}_{n-k+1,n}(1-t_{n-k})(1-u_{n-k})v_{n-k}\right)^{\Pi_{n-k} ^{-}}\nonumber\\
&\times& \overleftrightarrow {w}_{n-k,n}^{-(2(n-k)-\frac{3}{2}+\frac{\alpha }{2})}\left(  \overleftrightarrow {w}_{n-k,n} \partial _{ \overleftrightarrow {w}_{n-k,n}}\right) \overleftrightarrow {w}_{n-k,n}^{ -\frac{1}{2}} \left(  \overleftrightarrow {w}_{n-k,n} \partial _{ \overleftrightarrow {w}_{n-k,n}}\right) \overleftrightarrow {w}_{n-k,n}^{2(n-k) -1+\frac{\alpha }{2}  } \Bigg\}\nonumber\\
&\times& \left.\left. \;_2F_1\left( -\Pi_0 ^{+}, -\Pi_0 ^{-}; \frac{1}{2}; \overleftrightarrow {w}_{1,n}\right) \right\} z^n \right\} \label{eq:80056}
\end{eqnarray}
where
\begin{equation}
\begin{cases} 
\Pi_0 ^{\pm}=  \frac{-(1-\rho ^2)(\alpha +1)\pm \sqrt{(\alpha +1)^2-(2-\rho ^2)(h+\alpha +1)}}{2(2-\rho ^2)}  \cr
\Pi_{n-k} ^{\pm}=  -2(n-k) + \frac{-(1-\rho ^2)(\alpha +1)\pm \sqrt{(\alpha +1)^2-(2-\rho ^2)(h+\alpha +1)}}{2(2-\rho ^2)}
\end{cases}\nonumber 
\end{equation}
On (\ref{eq:80055}) and (\ref{eq:80056}),
\begin{equation}
\begin{cases} 
\varsigma= 1-\xi ^{-1} \cr
\eta =\frac{2-\rho ^2}{1-\rho ^2}\varsigma \cr
z= -\frac{1}{1-\rho ^2} \varsigma ^2 
\end{cases}\nonumber 
\end{equation}
\subsection{ ${\displaystyle \left(\frac{x-a}{1-a} \right)^{-\alpha } Hl\left(a, q-(\beta -\delta )\alpha ; \alpha , -\beta+\gamma +\delta , \delta , \gamma; \frac{a(x-1)}{x-a} \right)}$}
\subsubsection{Polynomial of type 2}
Replace coefficients $q$, $\beta $, $\gamma $, $\delta $, $x$, $c_0$, $\lambda $ and $q_j$ where $j, q_j \in \mathbb{N}_{0}$ by $q-(\beta -\delta )\alpha $, $-\beta+\gamma +\delta $, $\delta $,  $\gamma $, $\frac{a(x-1)}{x-a}$, 1, zero and $h_j$ where $h_j \in \mathbb{N}_{0}$ into (\ref{eq:80013}). Multiply $\left(\frac{x-a}{1-a} \right)^{-\alpha }$ and (\ref{eq:80013}) together. Put (\ref{eq:8006}) into the new (\ref{eq:80013}).
\begin{eqnarray} 
&&\left(\frac{\xi-\rho ^{-2}}{1-\rho ^{-2}} \right)^{-\frac{1}{2}(\alpha +1)} y(\varsigma )\nonumber\\
&=& \left(\frac{\xi-\rho ^{-2}}{1-\rho ^{-2}} \right)^{-\frac{1}{2}(\alpha +1)} Hl\left( \rho ^{-2}, -(h_j+2j)\left( (1+\rho ^{-2})(h_j+2j)+\alpha +1\right); \frac{\alpha }{2} +\frac{1}{2}, \frac{\alpha }{2} +1, \frac{1}{2}, \frac{1}{2}; \varsigma \right) \nonumber\\
&=& \left(\frac{\xi-\rho ^{-2}}{1-\rho ^{-2}} \right)^{-\frac{1}{2}(\alpha +1)} \left\{ _2F_1\left( -h_0,h_0+ \Omega _{\rho }; \frac{1}{2}; \eta\right)\right.  \nonumber\\
&&+ \sum_{n=1}^{\infty } \Bigg\{\prod _{k=0}^{n-1} \Bigg\{ \int_{0}^{1} dt_{n-k}\;t_{n-k}^{2(n-k)-1 } \int_{0}^{1} du_{n-k}\;u_{n-k}^{2(n-k)-\frac{3}{2} } \nonumber\\
&&\times  \frac{1}{2\pi i}  \oint dv_{n-k} \frac{1}{v_{n-k}} \left( \frac{v_{n-k}-1}{v_{n-k}} \frac{1}{1-\overleftrightarrow {w}_{n-k+1,n}(1-t_{n-k})(1-u_{n-k})v_{n-k}}\right)^{h_{n-k}} \nonumber\\
&&\times \left( 1- \overleftrightarrow {w}_{n-k+1,n}(1-t_{n-k})(1-u_{n-k})v_{n-k}\right)^{-\left( 4(n-k)+\Omega _{\rho }\right)}\nonumber\\
&&\times \overleftrightarrow {w}_{n-k,n}^{-(2(n-k)-\frac{3}{2}+\frac{\alpha }{2} )}\left(  \overleftrightarrow {w}_{n-k,n} \partial _{ \overleftrightarrow {w}_{n-k,n}}\right) \overleftrightarrow {w}_{n-k,n}^{ -\frac{1}{2}} \left(  \overleftrightarrow {w}_{n-k,n} \partial _{ \overleftrightarrow {w}_{n-k,n}}\right) \overleftrightarrow {w}_{n-k,n}^{2(n-k)-1 +\frac{\alpha }{2} } \Bigg\}\nonumber\\
&&\times \;_2F_1\left( -h_0, h_0+\Omega _{\rho }; \frac{1}{2}; \overleftrightarrow {w}_{1,n}\right) \Bigg\} z^n \Bigg\} \label{eq:80057}
\end{eqnarray}
where
\begin{equation}
\begin{cases} 
\Omega _{\rho } = \frac{ \alpha +1 }{1+\rho ^{-2}} \cr
h= 4(h_j+ 2j)\left( (1+ \rho ^2)(h_j+ 2j) +\rho ^2(\alpha +1) \right) +\rho ^2(\alpha +1)^2 
\end{cases}\nonumber 
\end{equation}
\subsubsection{Infinite series}
Replace coefficients $q$, $\beta $, $\gamma $, $\delta $, $x$, $c_0$ and $\lambda $ by $q-(\beta -\delta )\alpha $, $-\beta+\gamma +\delta $, $\delta $,  $\gamma $, $\frac{a(x-1)}{x-a}$, 1 and zero into (\ref{eq:80016}). Multiply $\left(\frac{x-a}{1-a} \right)^{-\alpha }$ and (\ref{eq:80016}) together. Put (\ref{eq:8006}) into the new (\ref{eq:80016}).
\begin{eqnarray}
&&\left(\frac{\xi-\rho ^{-2}}{1-\rho ^{-2}} \right)^{-\frac{1}{2}(\alpha +1)} y(\varsigma )\nonumber\\
&=& \left(\frac{\xi-\rho ^{-2}}{1-\rho ^{-2}} \right)^{-\frac{1}{2}(\alpha +1)} Hl\left( \rho ^{-2}, -\frac{1}{4}\left( h\rho ^{-2}- (\alpha +1)^2\right); \frac{\alpha }{2} +\frac{1}{2}, \frac{\alpha }{2} +1, \frac{1}{2}, \frac{1}{2}; \varsigma \right) \nonumber\\
&=& \left(\frac{\xi-\rho ^{-2}}{1-\rho ^{-2}} \right)^{-\frac{1}{2}(\alpha +1)} \left\{_2F_1\left( -\Pi_0 ^{+}, -\Pi_0 ^{-}; \frac{1}{2}; \eta\right)  \right. \nonumber\\
&+& \sum_{n=1}^{\infty } \left\{\prod _{k=0}^{n-1} \Bigg\{ \int_{0}^{1} dt_{n-k}\;t_{n-k}^{2(n-k)-1 } \int_{0}^{1} du_{n-k}\;u_{n-k}^{2(n-k)-\frac{3}{2} } \right.\nonumber\\
&\times& \frac{1}{2\pi i}  \oint dv_{n-k} \frac{1}{v_{n-k}}\left( \frac{v_{n-k}-1}{v_{n-k}}\right)^{\Pi_{n-k} ^{+}}  \left( 1- \overleftrightarrow {w}_{n-k+1,n}(1-t_{n-k})(1-u_{n-k})v_{n-k}\right)^{\Pi_{n-k} ^{-}}\nonumber\\
&\times& \overleftrightarrow {w}_{n-k,n}^{-(2(n-k)-\frac{3}{2}+\frac{\alpha }{2})}\left(  \overleftrightarrow {w}_{n-k,n} \partial _{ \overleftrightarrow {w}_{n-k,n}}\right) \overleftrightarrow {w}_{n-k,n}^{ -\frac{1}{2}} \left(  \overleftrightarrow {w}_{n-k,n} \partial _{ \overleftrightarrow {w}_{n-k,n}}\right) \overleftrightarrow {w}_{n-k,n}^{2(n-k) -1+\frac{\alpha }{2}  } \Bigg\}\nonumber\\
&\times& \left.\left. \;_2F_1\left( -\Pi_0 ^{+}, -\Pi_0 ^{-}; \frac{1}{2}; \overleftrightarrow {w}_{1,n}\right) \right\} z^n \right\} \label{eq:80058}
\end{eqnarray}
where
\begin{equation}
\begin{cases} 
\Pi_0 ^{\pm}=  \frac{-\rho ^2(\alpha +1)\pm \sqrt{(1+\rho ^2)h- \rho ^2(\alpha +1)^2}}{2(1+\rho ^2)}  \cr
\Pi_{n-k} ^{\pm}=  -2(n-k) + \frac{-\rho ^2(\alpha +1)\pm \sqrt{(1+\rho ^2)h- \rho ^2(\alpha +1)^2}}{2(1+\rho ^2)}
\end{cases}\nonumber 
\end{equation}
On (\ref{eq:80057}) and (\ref{eq:80058}),
\begin{equation}
\begin{cases} 
\varsigma= \frac{\xi -1}{\rho ^2(\xi -\rho ^{-2})} \cr
\eta = (1+\rho ^2) \varsigma \cr
z= -\rho ^2 \varsigma ^2 
\end{cases}\nonumber 
\end{equation}
\section[Generating function for the type 2 polynomial]{Generating function for the type 2 polynomial
  \sectionmark{Generating function for the type 2 polynomial}}
  \sectionmark{Generating function for the type 2 polynomial}
\subsection{ ${\displaystyle (1-x)^{1-\delta } Hl(a, q - (\delta  - 1)\gamma a; \alpha - \delta  + 1, \beta - \delta + 1, \gamma ,2 - \delta ; x)}$ }
Replace coefficients $q$, $\alpha$, $\beta$, $\delta$ and $q_j$ where $j, q_j \in \mathbb{N}_{0}$ by $q - (\delta - 1)\gamma a $, $\alpha - \delta  + 1 $, $\beta - \delta + 1$, $2 - \delta$ and $h_j$ where $h_j \in \mathbb{N}_{0}$ into (\ref{eq:80019}). Put (\ref{eq:8006}) into the new (\ref{eq:80019}). 
\begin{eqnarray}
&&\sum_{h_0 =0}^{\infty } \frac{(\frac{1}{2})_{h_0}}{h_0!} s_0^{h_0} \prod _{n=1}^{\infty } \left\{ \sum_{ h_n = h_{n-1}}^{\infty } s_n^{h_n }\right\} Hl\left(\rho ^{-2}, -\rho ^{-2}(h_j+2j)[1+(1+\rho ^2)(h_j+2j)]; \frac{\alpha }{2}+1, -\frac{\alpha }{2}+\frac{1}{2}, \frac{1}{2},\frac{3}{2}; \xi \right) \nonumber\\
&&=2^{-\Omega _{\rho } }\Bigg\{ \prod_{l=1}^{\infty } \frac{1}{(1-s_{l,\infty })}  \mathbf{A}\left( s_{0,\infty } ;\eta\right) +\Bigg\{ \prod_{l=2}^{\infty } \frac{1}{(1-s_{l,\infty })} \int_{0}^{1} dt_1\;t_1 \int_{0}^{1} du_1\;u_1^{\frac{1}{2}} \overleftrightarrow {\mathbf{\Gamma}}_1 \left( s_{1,\infty };t_1,u_1,\eta\right)\nonumber\\
&&\times \widetilde{w}_{1,1}^{-\frac{1}{2}(2+ \alpha )}\left( \widetilde{w}_{1,1} \partial _{ \widetilde{w}_{1,1}}\right) \widetilde{w}_{1,1}^{\alpha +\frac{1}{2} } \left( \widetilde{w}_{1,1} \partial _{ \widetilde{w}_{1,1}}\right)\widetilde{w}_{1,1}^{\frac{1}{2}(1-\alpha ) } \mathbf{A}\left( s_{0} ;\widetilde{w}_{1,1}\right) \Bigg\}z \nonumber\\
&&+ \sum_{n=2}^{\infty } \Bigg\{ \prod_{l=n+1}^{\infty } \frac{1}{(1-s_{l,\infty })} \int_{0}^{1} dt_n\;t_n^{2n-1} \int_{0}^{1} du_n\;u_n^{2n-\frac{3}{2} } \overleftrightarrow {\mathbf{\Gamma}}_n \left( s_{n,\infty };t_n,u_n,\eta \right)\nonumber\\
&&\times \widetilde{w}_{n,n}^{-\frac{1}{2}(4n-2+\alpha )}\left( \widetilde{w}_{n,n} \partial _{ \widetilde{w}_{n,n}}\right) \widetilde{w}_{n,n}^{\alpha +\frac{1}{2} } \left( \widetilde{w}_{n,n} \partial _{ \widetilde{w}_{n,n}}\right)\widetilde{w}_{n,n}^{\frac{1}{2}(4n-3-\alpha ) }  \nonumber\\
&&\times \prod_{k=1}^{n-1} \Bigg\{ \int_{0}^{1} dt_{n-k}\;t_{n-k}^{2(n-k)-1} \int_{0}^{1} du_{n-k} \;u_{n-k}^{2(n-k)-\frac{3}{2}}\overleftrightarrow {\mathbf{\Gamma}}_{n-k} \left( s_{n-k};t_{n-k},u_{n-k},\widetilde{w}_{n-k+1,n} \right)\label{eq:80059}\\
&&\times \widetilde{w}_{n-k,n}^{-\frac{1}{2}(4(n-k)-2+\alpha )}\left( \widetilde{w}_{n-k,n} \partial _{ \widetilde{w}_{n-k,n}}\right) \widetilde{w}_{n-k,n}^{\alpha +\frac{1}{2} } \left( \widetilde{w}_{n-k,n} \partial _{ \widetilde{w}_{n-k,n}}\right)\widetilde{w}_{n-k,n}^{\frac{1}{2}(4(n-k)-3-\alpha ) } \Bigg\} \mathbf{A} \left( s_{0} ;\widetilde{w}_{1,n}\right) \Bigg\} z^n \Bigg\}  \nonumber  
\end{eqnarray}
where
\begin{equation}
\begin{cases} 
{ \displaystyle \overleftrightarrow {\mathbf{\Gamma}}_1 \left( s_{1,\infty };t_1,u_1,\eta\right)= \frac{\left( \frac{1+s_{1,\infty }+\sqrt{s_{1,\infty }^2-2(1-2\eta (1-t_1)(1-u_1))s_{1,\infty }+1}}{2}\right)^{-\left( 3+\Omega _{\rho }\right)}}{\sqrt{s_{1,\infty }^2-2(1-2\eta (1-t_1)(1-u_1))s_{1,\infty }+1}}}\cr
{ \displaystyle  \overleftrightarrow {\mathbf{\Gamma}}_n \left( s_{n,\infty };t_n,u_n,\eta \right) =\frac{\left( \frac{1+s_{n,\infty }+\sqrt{s_{n,\infty }^2-2(1-2\eta (1-t_n)(1-u_n))s_{n,\infty }+1}}{2}\right)^{-\left( 4n-1+\Omega _{\rho }\right)}}{\sqrt{ s_{n,\infty }^2-2(1-2\eta (1-t_n)(1-u_n))s_{n,\infty }+1}}}\cr
{ \displaystyle \overleftrightarrow {\mathbf{\Gamma}}_{n-k} \left( s_{n-k};t_{n-k},u_{n-k},\widetilde{w}_{n-k+1,n} \right) }\cr
{ \displaystyle = \frac{ \left( \frac{1+s_{n-k}+\sqrt{s_{n-k}^2-2(1-2\widetilde{w}_{n-k+1,n} (1-t_{n-k})(1-u_{n-k}))s_{n-k}+1}}{2}\right)^{-\left( 4(n-k)-1+\Omega _{\rho }\right)}}{\sqrt{ s_{n-k}^2-2(1-2\widetilde{w}_{n-k+1,n} (1-t_{n-k})(1-u_{n-k}))s_{n-k}+1}}}
\end{cases}\nonumber 
\end{equation}
 \begin{equation}
\begin{cases} 
{ \displaystyle \mathbf{A} \left( s_{0,\infty } ;\eta\right)= \frac{\left(1- s_{0,\infty }+\sqrt{s_{0,\infty }^2-2(1-2\eta )s_{0,\infty }+1}\right)^{\frac{1}{2}} \left(1+s_{0,\infty }+\sqrt{s_{0,\infty }^2-2(1-2\eta )s_{0,\infty }+1}\right)^{\frac{1}{2} -\Omega _{\rho }}}{\sqrt{s_{0,\infty }^2-2(1-2\eta )s_{0,\infty }+1}}}\cr
{ \displaystyle  \mathbf{A} \left( s_{0} ;\widetilde{w}_{1,1}\right) = \frac{\left(1- s_0+\sqrt{s_0^2-2(1-2\widetilde{w}_{1,1})s_0+1}\right)^{\frac{1}{2}} \left(1+s_0+\sqrt{s_0^2-2(1-2\widetilde{w}_{1,1} )s_0+1}\right)^{\frac{1}{2} -\Omega _{\rho }}}{\sqrt{s_0^2-2(1-2\widetilde{w}_{1,1})s_0+1}}} \cr
{ \displaystyle \mathbf{A} \left(  s_{0} ;\widetilde{w}_{1,n}\right) = \frac{\left(1- s_0+\sqrt{s_0^2-2(1-2\widetilde{w}_{1,n})s_0+1}\right)^{\frac{1}{2}} \left(1+s_0+\sqrt{s_0^2-2(1-2\widetilde{w}_{1,n} )s_0+1}\right)^{\frac{1}{2} -\Omega _{\rho }}}{\sqrt{s_0^2-2(1-2\widetilde{w}_{1,n})s_0+1}}}
\end{cases}\nonumber 
\end{equation}
and
\begin{equation}
\begin{cases} 
\eta =(1+\rho ^2)\xi \cr
z=-\rho ^2\xi^2 \cr
\Omega _{\rho } = \frac{1}{1+\rho ^2}\cr
h= 4(h_j+2j)[1+(1+\rho ^2)(h_j+2j)]+1  
\end{cases}\nonumber 
\end{equation}
\subsection{ \footnotesize ${\displaystyle x^{1-\gamma } (1-x)^{1-\delta }Hl(a, q-(\gamma +\delta -2)a -(\gamma -1)(\alpha +\beta -\gamma -\delta +1), \alpha - \gamma -\delta +2, \beta - \gamma -\delta +2, 2-\gamma, 2 - \delta ; x)}$ \normalsize}
Replace coefficients $q$, $\alpha$, $\beta$, $\gamma $, $\delta$ and $q_j$ where $j, q_j \in \mathbb{N}_{0}$ by $q-(\gamma +\delta -2)a-(\gamma -1)(\alpha +\beta -\gamma -\delta +1)$, $\alpha - \gamma -\delta +2$, $\beta - \gamma -\delta +2, 2-\gamma$, $2 - \delta$ and $h_j$ where $h_j \in \mathbb{N}_{0}$ into (\ref{eq:80019}). Put (\ref{eq:8006}) into the new (\ref{eq:80019}). 
\begin{eqnarray}
&&\sum_{h_0 =0}^{\infty } \frac{(\frac{3}{2})_{h_0}}{h_0!} s_0^{h_0} \prod _{n=1}^{\infty } \left\{ \sum_{ h_n = h_{n-1}}^{\infty } s_n^{h_n }\right\} Hl\left(\rho ^{-2}, -\rho ^{-2}(h_j+2j)[2+\rho ^2 +(1+\rho ^2)(h_j+2j)]; \frac{\alpha }{2}+\frac{3}{2},\right. \nonumber\\
&&\left. -\frac{\alpha }{2}+1, \frac{3}{2},\frac{3}{2}; \xi \right)\nonumber\\
&&=2^{\Omega _{\rho }}\Bigg\{ \prod_{l=1}^{\infty } \frac{1}{(1-s_{l,\infty })}  \mathbf{A}\left( s_{0,\infty } ;\eta\right) + \Bigg\{ \prod_{l=2}^{\infty } \frac{1}{(1-s_{l,\infty })} \int_{0}^{1} dt_1\;t_1 \int_{0}^{1} du_1\;u_1^{\frac{3}{2}} \overleftrightarrow {\mathbf{\Gamma}}_1 \left( s_{1,\infty };t_1,u_1,\eta\right)\nonumber\\
&&\times \widetilde{w}_{1,1}^{-\frac{1}{2}(3+ \alpha )}\left( \widetilde{w}_{1,1} \partial _{ \widetilde{w}_{1,1}}\right) \widetilde{w}_{1,1}^{\alpha +\frac{1}{2} } \left( \widetilde{w}_{1,1} \partial _{ \widetilde{w}_{1,1}}\right)\widetilde{w}_{1,1}^{\frac{1}{2}(2-\alpha ) } \mathbf{A}\left( s_{0} ;\widetilde{w}_{1,1}\right) \Bigg\}z \nonumber\\
&&+ \sum_{n=2}^{\infty } \Bigg\{ \prod_{l=n+1}^{\infty } \frac{1}{(1-s_{l,\infty })} \int_{0}^{1} dt_n\;t_n^{2n-1} \int_{0}^{1} du_n\;u_n^{2n-\frac{1}{2} } \overleftrightarrow {\mathbf{\Gamma}}_n \left( s_{n,\infty };t_n,u_n,\eta \right)\nonumber\\
&&\times \widetilde{w}_{n,n}^{-\frac{1}{2}(4n-1+\alpha )}\left( \widetilde{w}_{n,n} \partial _{ \widetilde{w}_{n,n}}\right) \widetilde{w}_{n,n}^{\alpha +\frac{1}{2} } \left( \widetilde{w}_{n,n} \partial _{ \widetilde{w}_{n,n}}\right)\widetilde{w}_{n,n}^{\frac{1}{2}(4n-2-\alpha ) }  \nonumber\\
&&\times \prod_{k=1}^{n-1} \Bigg\{ \int_{0}^{1} dt_{n-k}\;t_{n-k}^{2(n-k)-1} \int_{0}^{1} du_{n-k} \;u_{n-k}^{2(n-k)-\frac{1}{2}}\overleftrightarrow {\mathbf{\Gamma}}_{n-k} \left( s_{n-k};t_{n-k},u_{n-k},\widetilde{w}_{n-k+1,n} \right) \label{eq:80060}\\
&&\times \widetilde{w}_{n-k,n}^{-\frac{1}{2}(4(n-k)-1+\alpha )}\left( \widetilde{w}_{n-k,n} \partial _{ \widetilde{w}_{n-k,n}}\right) \widetilde{w}_{n-k,n}^{\alpha +\frac{1}{2} } \left( \widetilde{w}_{n-k,n} \partial _{ \widetilde{w}_{n-k,n}}\right)\widetilde{w}_{n-k,n}^{\frac{1}{2}(4(n-k)-2-\alpha ) } \Bigg\} \mathbf{A} \left( s_{0} ;\widetilde{w}_{1,n}\right) \Bigg\} z^n \Bigg\}  \nonumber 
\end{eqnarray}
where
\begin{equation}
\begin{cases} 
{ \displaystyle \overleftrightarrow {\mathbf{\Gamma}}_1 \left( s_{1,\infty };t_1,u_1,\eta\right)= \frac{\left( \frac{1+s_{1,\infty }+\sqrt{s_{1,\infty }^2-2(1-2\eta (1-t_1)(1-u_1))s_{1,\infty }+1}}{2}\right)^{- 4-\Omega _{\rho } }}{\sqrt{s_{1,\infty }^2-2(1-2\eta (1-t_1)(1-u_1))s_{1,\infty }+1}}}\cr
{ \displaystyle  \overleftrightarrow {\mathbf{\Gamma}}_n \left( s_{n,\infty };t_n,u_n,\eta \right) =\frac{\left( \frac{1+s_{n,\infty }+\sqrt{s_{n,\infty }^2-2(1-2\eta (1-t_n)(1-u_n))s_{n,\infty }+1}}{2}\right)^{- 4n -\Omega _{\rho } }}{\sqrt{ s_{n,\infty }^2-2(1-2\eta (1-t_n)(1-u_n))s_{n,\infty }+1}}}\cr
{ \displaystyle \overleftrightarrow {\mathbf{\Gamma}}_{n-k} \left( s_{n-k};t_{n-k},u_{n-k},\widetilde{w}_{n-k+1,n} \right) = \frac{ \left( \frac{1+s_{n-k}+\sqrt{s_{n-k}^2-2(1-2\widetilde{w}_{n-k+1,n} (1-t_{n-k})(1-u_{n-k}))s_{n-k}+1}}{2}\right)^{- 4(n-k) -\Omega _{\rho } }}{\sqrt{ s_{n-k}^2-2(1-2\widetilde{w}_{n-k+1,n} (1-t_{n-k})(1-u_{n-k}))s_{n-k}+1}}}
\end{cases}\nonumber 
\end{equation}
 
  \begin{equation}
\begin{cases} 
{ \displaystyle \mathbf{A} \left( s_{0,\infty } ;\eta\right) }\cr 
{ \displaystyle = \frac{\left(1- s_{0,\infty }+\sqrt{s_{0,\infty }^2-2(1-2\eta )s_{0,\infty }+1}\right)^{-\frac{1}{2}} \left(1+s_{0,\infty }+\sqrt{s_{0,\infty }^2-2(1-2\eta )s_{0,\infty }+1}\right)^{\frac{1}{2} -\Omega _{\rho }}}{\sqrt{s_{0,\infty }^2-2(1-2\eta )s_{0,\infty }+1}}}\cr
{ \displaystyle  \mathbf{A} \left( s_{0} ;\widetilde{w}_{1,1}\right) = \frac{\left(1- s_0+\sqrt{s_0^2-2(1-2\widetilde{w}_{1,1})s_0+1}\right)^{-\frac{1}{2}} \left(1+s_0+\sqrt{s_0^2-2(1-2\widetilde{w}_{1,1} )s_0+1}\right)^{\frac{1}{2} -\Omega _{\rho }}}{\sqrt{s_0^2-2(1-2\widetilde{w}_{1,1})s_0+1}}} \cr
{ \displaystyle \mathbf{A} \left(  s_{0} ;\widetilde{w}_{1,n}\right) = \frac{\left(1- s_0+\sqrt{s_0^2-2(1-2\widetilde{w}_{1,n})s_0+1}\right)^{-\frac{1}{2}} \left(1+s_0+\sqrt{s_0^2-2(1-2\widetilde{w}_{1,n} )s_0+1}\right)^{\frac{1}{2} -\Omega _{\rho }}}{\sqrt{s_0^2-2(1-2\widetilde{w}_{1,n})s_0+1}}}
\end{cases}\nonumber 
\end{equation}
and
\begin{equation}
\begin{cases} 
\Omega _{\rho } = \frac{1}{1+\rho ^2}\cr
h= 4+\rho ^2+ 4(h_j+2j)[2+\rho ^2+(1+\rho ^2)(h_j+2j)]  \cr
\eta =(1+\rho ^2)\xi \cr
z=-\rho ^2\xi^2 
\end{cases}\nonumber 
\end{equation}
\subsection{ ${\displaystyle  Hl(1-a,-q+\alpha \beta; \alpha,\beta, \delta, \gamma; 1-x)}$} 
Replace coefficients $a$, $q$, $\gamma $, $\delta$, $x$ and $q_j$ where $j, q_j \in \mathbb{N}_{0}$ by $1-a$, $-q +\alpha \beta $, $\delta $, $\gamma $, $1-x$ and $h_j$ where $h_j \in \mathbb{N}_{0}$ into (\ref{eq:80019}). Put (\ref{eq:8006}) into the new (\ref{eq:80019}). 
\begin{eqnarray}
&&\sum_{h_0 =0}^{\infty } \frac{(\frac{1}{2})_{h_0}}{h_0!} s_0^{h_0} \prod _{n=1}^{\infty } \left\{ \sum_{ h_n = h_{n-1}}^{\infty } s_n^{h_n }\right\} Hl\left( 1-\rho ^{-2}, -(2-\rho ^{-2})(h_j +2j)^2; \frac{\alpha }{2}+\frac{1}{2}, -\frac{\alpha }{2}, \frac{1}{2}, \frac{1}{2}; \varsigma \right)\nonumber\\
&&=2^{-1}\Bigg\{ \prod_{l=1}^{\infty } \frac{1}{(1-s_{l,\infty })}  \mathbf{A}\left( s_{0,\infty } ;\eta\right) + \Bigg\{ \prod_{l=2}^{\infty } \frac{1}{(1-s_{l,\infty })} \int_{0}^{1} dt_1\;t_1 \int_{0}^{1} du_1\;u_1^{\frac{1}{2}} \overleftrightarrow {\mathbf{\Gamma}}_1 \left( s_{1,\infty };t_1,u_1,\eta\right)\nonumber\\
&&\times \widetilde{w}_{1,1}^{-\frac{1}{2}(1+ \alpha )}\left( \widetilde{w}_{1,1} \partial _{ \widetilde{w}_{1,1}}\right) \widetilde{w}_{1,1}^{\alpha +\frac{1}{2} } \left( \widetilde{w}_{1,1} \partial _{ \widetilde{w}_{1,1}}\right)\widetilde{w}_{1,1}^{-\frac{\alpha }{2}} \mathbf{A}\left( s_{0} ;\widetilde{w}_{1,1}\right) \Bigg\}z \nonumber\\
&&+ \sum_{n=2}^{\infty } \Bigg\{ \prod_{l=n+1}^{\infty } \frac{1}{(1-s_{l,\infty })} \int_{0}^{1} dt_n\;t_n^{2n-1} \int_{0}^{1} du_n\;u_n^{2n-\frac{3}{2} } \overleftrightarrow {\mathbf{\Gamma}}_n \left( s_{n,\infty };t_n,u_n,\eta \right)\nonumber\\
&&\times \widetilde{w}_{n,n}^{-\frac{1}{2}(4n-3+\alpha )}\left( \widetilde{w}_{n,n} \partial _{ \widetilde{w}_{n,n}}\right) \widetilde{w}_{n,n}^{\alpha +\frac{1}{2} } \left( \widetilde{w}_{n,n} \partial _{ \widetilde{w}_{n,n}}\right)\widetilde{w}_{n,n}^{\frac{1}{2}(4(n-1)-\alpha ) }  \nonumber\\
&&\times \prod_{k=1}^{n-1} \Bigg\{ \int_{0}^{1} dt_{n-k}\;t_{n-k}^{2(n-k)-1} \int_{0}^{1} du_{n-k} \;u_{n-k}^{2(n-k)-\frac{3}{2}}\overleftrightarrow {\mathbf{\Gamma}}_{n-k} \left( s_{n-k};t_{n-k},u_{n-k},\widetilde{w}_{n-k+1,n} \right)\label{eq:80061}\\
&&\times \widetilde{w}_{n-k,n}^{-\frac{1}{2}(4(n-k)-3+\alpha )}\left( \widetilde{w}_{n-k,n} \partial _{ \widetilde{w}_{n-k,n}}\right) \widetilde{w}_{n-k,n}^{\alpha +\frac{1}{2} } \left( \widetilde{w}_{n-k,n} \partial _{ \widetilde{w}_{n-k,n}}\right)\widetilde{w}_{n-k,n}^{\frac{1}{2}(4(n-k-1)-\alpha ) } \Bigg\} \mathbf{A} \left( s_{0} ;\widetilde{w}_{1,n}\right) \Bigg\} z^n \Bigg\} \nonumber   
\end{eqnarray}
where
\begin{equation}
\begin{cases} 
{ \displaystyle \overleftrightarrow {\mathbf{\Gamma}}_1 \left( s_{1,\infty };t_1,u_1,\eta\right)= \frac{\left( \frac{1+s_{1,\infty }+\sqrt{s_{1,\infty }^2-2(1-2\eta (1-t_1)(1-u_1))s_{1,\infty }+1}}{2}\right)^{-3}}{\sqrt{s_{1,\infty }^2-2(1-2\eta (1-t_1)(1-u_1))s_{1,\infty }+1}}}\cr
{ \displaystyle  \overleftrightarrow {\mathbf{\Gamma}}_n \left( s_{n,\infty };t_n,u_n,\eta \right) =\frac{\left( \frac{1+s_{n,\infty }+\sqrt{s_{n,\infty }^2-2(1-2\eta (1-t_n)(1-u_n))s_{n,\infty }+1}}{2}\right)^{-4n+1 }}{\sqrt{ s_{n,\infty }^2-2(1-2\eta (1-t_n)(1-u_n))s_{n,\infty }+1}}}\cr
{ \displaystyle \overleftrightarrow {\mathbf{\Gamma}}_{n-k} \left( s_{n-k};t_{n-k},u_{n-k},\widetilde{w}_{n-k+1,n} \right) = \frac{ \left( \frac{1+s_{n-k}+\sqrt{s_{n-k}^2-2(1-2\widetilde{w}_{n-k+1,n} (1-t_{n-k})(1-u_{n-k}))s_{n-k}+1}}{2}\right)^{-4(n-k)+1}}{\sqrt{ s_{n-k}^2-2(1-2\widetilde{w}_{n-k+1,n} (1-t_{n-k})(1-u_{n-k}))s_{n-k}+1}}}
\end{cases}\nonumber 
\end{equation}
 
 \begin{equation}
\begin{cases} 
{ \displaystyle \mathbf{A} \left( s_{0,\infty } ;\eta\right)= \frac{\left(1- s_{0,\infty }+\sqrt{s_{0,\infty }^2-2(1-2\eta )s_{0,\infty }+1}\right)^{ \frac{1}{2}} \left(1+s_{0,\infty }+\sqrt{s_{0,\infty }^2-2(1-2\eta )s_{0,\infty }+1}\right)^{\frac{1}{2} }}{\sqrt{s_{0,\infty }^2-2(1-2\eta )s_{0,\infty }+1}}}\cr
{ \displaystyle  \mathbf{A} \left( s_{0} ;\widetilde{w}_{1,1}\right) = \frac{\left(1- s_0+\sqrt{s_0^2-2(1-2\widetilde{w}_{1,1})s_0+1}\right)^{ \frac{1}{2}} \left(1+s_0+\sqrt{s_0^2-2(1-2\widetilde{w}_{1,1} )s_0+1}\right)^{\frac{1}{2} }}{\sqrt{s_0^2-2(1-2\widetilde{w}_{1,1})s_0+1}}} \cr
{ \displaystyle \mathbf{A} \left(  s_{0} ;\widetilde{w}_{1,n}\right) = \frac{\left(1- s_0+\sqrt{s_0^2-2(1-2\widetilde{w}_{1,n})s_0+1}\right)^{ \frac{1}{2}} \left(1+s_0+\sqrt{s_0^2-2(1-2\widetilde{w}_{1,n} )s_0+1}\right)^{\frac{1}{2} }}{\sqrt{s_0^2-2(1-2\widetilde{w}_{1,n})s_0+1}}}
\end{cases}\nonumber 
\end{equation}
and 
\begin{equation}
\begin{cases} 
h= \rho ^2 \left( \alpha (\alpha +1)- 4(2-\rho ^{-2})(h_j+2j)^2 \right) \cr
\varsigma= 1-\xi \cr
\eta =\frac{2-\rho ^{-2}}{1-\rho ^{-2}}\varsigma \cr
z=\frac{-1}{1-\rho ^{-2}}\varsigma ^2  
\end{cases}\nonumber 
\end{equation} 
\subsection{ \footnotesize ${\displaystyle (1-x)^{1-\delta } Hl(1-a,-q+(\delta -1)\gamma a+(\alpha -\delta +1)(\beta -\delta +1); \alpha-\delta +1,\beta-\delta +1, 2-\delta, \gamma; 1-x)}$ \normalsize}
Replace coefficients $a$, $q$, $\alpha $, $\beta $, $\gamma $, $\delta$, $x$ and $q_j$ where $j, q_j \in \mathbb{N}_{0}$ by $1-a$, $-q+(\delta -1)\gamma a+(\alpha -\delta +1)(\beta -\delta +1)$, $\alpha-\delta +1 $, $\beta-\delta +1 $, $2-\delta$, $\gamma $, $1-x$ and $h_j$ where $h_j \in \mathbb{N}_{0}$ into (\ref{eq:80019}). Put (\ref{eq:8006}) into the new (\ref{eq:80019}). 
\begin{eqnarray}
&&\sum_{h_0 =0}^{\infty } \frac{(\frac{3}{2})_{h_0}}{h_0!} s_0^{h_0} \prod _{n=1}^{\infty } \left\{ \sum_{ h_n = h_{n-1}}^{\infty } s_n^{h_n }\right\} Hl\bigg( 1-\rho ^{-2}, -(2-\rho ^{-2})(h_j+2j)(h_j+2j+1) ;\frac{\alpha }{2}+1, \nonumber\\
&& -\frac{\alpha }{2}+\frac{1}{2}, \frac{3}{2}, \frac{1}{2}; \varsigma \bigg)\nonumber\\
&&= \prod_{l=1}^{\infty } \frac{1}{(1-s_{l,\infty })}  \mathbf{A}\left( s_{0,\infty } ;\eta\right) + \Bigg\{ \prod_{l=2}^{\infty } \frac{1}{(1-s_{l,\infty })} \int_{0}^{1} dt_1\;t_1 \int_{0}^{1} du_1\;u_1^{\frac{3}{2}} \overleftrightarrow {\mathbf{\Gamma}}_1 \left( s_{1,\infty };t_1,u_1,\eta\right)\nonumber\\
&&\times \widetilde{w}_{1,1}^{-\frac{1}{2}(2+ \alpha )}\left( \widetilde{w}_{1,1} \partial _{ \widetilde{w}_{1,1}}\right) \widetilde{w}_{1,1}^{\alpha +\frac{1}{2} } \left( \widetilde{w}_{1,1} \partial _{ \widetilde{w}_{1,1}}\right)\widetilde{w}_{1,1}^{ \frac{1}{2}(1-\alpha )} \mathbf{A}\left( s_{0} ;\widetilde{w}_{1,1}\right) \Bigg\}z \nonumber\\
&&+ \sum_{n=2}^{\infty } \Bigg\{ \prod_{l=n+1}^{\infty } \frac{1}{(1-s_{l,\infty })} \int_{0}^{1} dt_n\;t_n^{2n-1} \int_{0}^{1} du_n\;u_n^{2n-\frac{1}{2} } \overleftrightarrow {\mathbf{\Gamma}}_n \left( s_{n,\infty };t_n,u_n,\eta \right)\nonumber\\
&&\times \widetilde{w}_{n,n}^{-\frac{1}{2}(4n-2+\alpha )}\left( \widetilde{w}_{n,n} \partial _{ \widetilde{w}_{n,n}}\right) \widetilde{w}_{n,n}^{\alpha +\frac{1}{2} } \left( \widetilde{w}_{n,n} \partial _{ \widetilde{w}_{n,n}}\right)\widetilde{w}_{n,n}^{\frac{1}{2}(4n-3-\alpha ) }  \nonumber\\
&&\times \prod_{k=1}^{n-1} \Bigg\{ \int_{0}^{1} dt_{n-k}\;t_{n-k}^{2(n-k)-1} \int_{0}^{1} du_{n-k} \;u_{n-k}^{2(n-k)-\frac{1}{2}}\overleftrightarrow {\mathbf{\Gamma}}_{n-k} \left( s_{n-k};t_{n-k},u_{n-k},\widetilde{w}_{n-k+1,n} \right)\label{eq:80062}\\
&&\times \widetilde{w}_{n-k,n}^{-\frac{1}{2}(4(n-k)-2+\alpha )}\left( \widetilde{w}_{n-k,n} \partial _{ \widetilde{w}_{n-k,n}}\right) \widetilde{w}_{n-k,n}^{\alpha +\frac{1}{2} } \left( \widetilde{w}_{n-k,n} \partial _{ \widetilde{w}_{n-k,n}}\right)\widetilde{w}_{n-k,n}^{\frac{1}{2}(4(n-k)-3-\alpha ) } \Bigg\} \mathbf{A} \left( s_{0} ;\widetilde{w}_{1,n}\right) \Bigg\} z^n  \nonumber
\end{eqnarray}
where
\begin{equation}
\begin{cases} 
{ \displaystyle \overleftrightarrow {\mathbf{\Gamma}}_1 \left( s_{1,\infty };t_1,u_1,\eta\right)= \frac{\left( \frac{1+s_{1,\infty }+\sqrt{s_{1,\infty }^2-2(1-2\eta (1-t_1)(1-u_1))s_{1,\infty }+1}}{2}\right)^{-4}}{\sqrt{s_{1,\infty }^2-2(1-2\eta (1-t_1)(1-u_1))s_{1,\infty }+1}}}\cr
{ \displaystyle  \overleftrightarrow {\mathbf{\Gamma}}_n \left( s_{n,\infty };t_n,u_n,\eta \right) =\frac{\left( \frac{1+s_{n,\infty }+\sqrt{s_{n,\infty }^2-2(1-2\eta (1-t_n)(1-u_n))s_{n,\infty }+1}}{2}\right)^{-4n }}{\sqrt{ s_{n,\infty }^2-2(1-2\eta (1-t_n)(1-u_n))s_{n,\infty }+1}}}\cr
{ \displaystyle \overleftrightarrow {\mathbf{\Gamma}}_{n-k} \left( s_{n-k};t_{n-k},u_{n-k},\widetilde{w}_{n-k+1,n} \right) = \frac{ \left( \frac{1+s_{n-k}+\sqrt{s_{n-k}^2-2(1-2\widetilde{w}_{n-k+1,n} (1-t_{n-k})(1-u_{n-k}))s_{n-k}+1}}{2}\right)^{-4(n-k) }}{\sqrt{ s_{n-k}^2-2(1-2\widetilde{w}_{n-k+1,n} (1-t_{n-k})(1-u_{n-k}))s_{n-k}+1}}}
\end{cases}\nonumber 
\end{equation}
  
 \begin{equation}
\begin{cases} 
{ \displaystyle \mathbf{A} \left( s_{0,\infty } ;\eta\right)= \frac{\left(1- s_{0,\infty }+\sqrt{s_{0,\infty }^2-2(1-2\eta )s_{0,\infty }+1}\right)^{ -\frac{1}{2}} \left(1+s_{0,\infty }+\sqrt{s_{0,\infty }^2-2(1-2\eta )s_{0,\infty }+1}\right)^{\frac{1}{2} }}{\sqrt{s_{0,\infty }^2-2(1-2\eta )s_{0,\infty }+1}}}\cr
{ \displaystyle  \mathbf{A} \left( s_{0} ;\widetilde{w}_{1,1}\right) = \frac{\left(1- s_0+\sqrt{s_0^2-2(1-2\widetilde{w}_{1,1})s_0+1}\right)^{ -\frac{1}{2}} \left(1+s_0+\sqrt{s_0^2-2(1-2\widetilde{w}_{1,1} )s_0+1}\right)^{\frac{1}{2} }}{\sqrt{s_0^2-2(1-2\widetilde{w}_{1,1})s_0+1}}} \cr
{ \displaystyle \mathbf{A} \left(  s_{0} ;\widetilde{w}_{1,n}\right) = \frac{\left(1- s_0+\sqrt{s_0^2-2(1-2\widetilde{w}_{1,n})s_0+1}\right)^{ -\frac{1}{2}} \left(1+s_0+\sqrt{s_0^2-2(1-2\widetilde{w}_{1,n} )s_0+1}\right)^{\frac{1}{2} }}{\sqrt{s_0^2-2(1-2\widetilde{w}_{1,n})s_0+1}}}
\end{cases}\nonumber 
\end{equation}
and 
\begin{equation}
\begin{cases}
h= 4(2\rho ^2-1)(h_j +2j)(h_j +2j+1)-\rho ^2(\alpha -1)(\alpha +2) -1  \cr  
\varsigma= 1-\xi \cr
\eta =\frac{2-\rho ^{-2}}{1-\rho ^{-2}}\varsigma \cr
z=\frac{-1}{1-\rho ^{-2}}\varsigma ^2  
\end{cases}\nonumber 
\end{equation}
\subsection{ ${\displaystyle x^{-\alpha } Hl\left(\frac{1}{a},\frac{q+\alpha [(\alpha -\gamma -\delta +1)a-\beta +\delta ]}{a}; \alpha , \alpha -\gamma +1, \alpha -\beta +1,\delta ;\frac{1}{x}\right)}$}
Replace coefficients $a$, $q$, $\beta $, $\gamma $, $x$ and $q_j$ where $j, q_j \in \mathbb{N}_{0}$ by $\frac{1}{a}$, $\frac{q+\alpha [(\alpha -\gamma -\delta +1)a-\beta +\delta ]}{a}$, $\alpha-\gamma +1 $, $\alpha -\beta +1 $, $\frac{1}{x}$ and  $h_j$ where $h_j \in \mathbb{N}_{0}$ into (\ref{eq:80019}). Put (\ref{eq:8006}) into the new (\ref{eq:80019}). 
\begin{eqnarray}
&&\sum_{h_0 =0}^{\infty } \frac{(\alpha +\frac{3}{2})_{h_0}}{h_0!} s_0^{h_0} \prod _{n=1}^{\infty } \left\{ \sum_{ h_n = h_{n-1}}^{\infty } s_n^{h_n }\right\} Hl\left(\rho ^2, -(1+\rho ^2)(h_j+2j)(h_j+2j+1+\alpha ); \frac{\alpha}{2} +\frac{1}{2},\right. \nonumber\\
&&\left. \frac{\alpha}{2} +1, \alpha +\frac{3}{2}, \frac{1}{2}; \varsigma \right) \nonumber\\
&&= 2^{\alpha }\Bigg\{ \prod_{l=1}^{\infty } \frac{1}{(1-s_{l,\infty })}  \mathbf{A}\left( s_{0,\infty } ;\eta\right) + \Bigg\{ \prod_{l=2}^{\infty } \frac{1}{(1-s_{l,\infty })} \int_{0}^{1} dt_1\;t_1 \int_{0}^{1} du_1\;u_1^{\alpha +\frac{3}{2}} \overleftrightarrow {\mathbf{\Gamma}}_1 \left( s_{1,\infty };t_1,u_1,\eta\right)\nonumber\\
&&\times \widetilde{w}_{1,1}^{-\frac{1}{2}(1+ \alpha )}\left( \widetilde{w}_{1,1} \partial _{ \widetilde{w}_{1,1}}\right) \widetilde{w}_{1,1}^{-\frac{1}{2} } \left( \widetilde{w}_{1,1} \partial _{ \widetilde{w}_{1,1}}\right)\widetilde{w}_{1,1}^{\frac{1}{2}(2+\alpha )} \mathbf{A}\left( s_{0} ;\widetilde{w}_{1,1}\right)\Bigg\} z \nonumber\\
&&+ \sum_{n=2}^{\infty } \Bigg\{ \prod_{l=n+1}^{\infty } \frac{1}{(1-s_{l,\infty })} \int_{0}^{1} dt_n\;t_n^{2n-1} \int_{0}^{1} du_n\;u_n^{2n-\frac{1}{2}+\alpha  } \overleftrightarrow {\mathbf{\Gamma}}_n \left( s_{n,\infty };t_n,u_n,\eta \right)\nonumber\\
&&\times \widetilde{w}_{n,n}^{-\frac{1}{2}(4n-3+\alpha )}\left( \widetilde{w}_{n,n} \partial _{ \widetilde{w}_{n,n}}\right) \widetilde{w}_{n,n}^{-\frac{1}{2} } \left( \widetilde{w}_{n,n} \partial _{ \widetilde{w}_{n,n}}\right)\widetilde{w}_{n,n}^{\frac{1}{2}(4n-2+\alpha ) }  \nonumber\\
&&\times \prod_{k=1}^{n-1} \Bigg\{ \int_{0}^{1} dt_{n-k}\;t_{n-k}^{2(n-k)-1} \int_{0}^{1} du_{n-k} \;u_{n-k}^{2(n-k)-\frac{1}{2}+\alpha }\overleftrightarrow {\mathbf{\Gamma}}_{n-k} \left( s_{n-k};t_{n-k},u_{n-k},\widetilde{w}_{n-k+1,n} \right)\label{eq:80063}\\
&&\times \widetilde{w}_{n-k,n}^{-\frac{1}{2}(4(n-k)-3+\alpha )}\left( \widetilde{w}_{n-k,n} \partial _{ \widetilde{w}_{n-k,n}}\right) \widetilde{w}_{n-k,n}^{-\frac{1}{2} } \left( \widetilde{w}_{n-k,n} \partial _{ \widetilde{w}_{n-k,n}}\right)\widetilde{w}_{n-k,n}^{\frac{1}{2}(4(n-k)-2+\alpha ) } \Bigg\} \mathbf{A} \left( s_{0} ;\widetilde{w}_{1,n}\right) \Bigg\} z^n \Bigg\}  \nonumber
\end{eqnarray}
where
\begin{equation}
\begin{cases} 
{ \displaystyle \overleftrightarrow {\mathbf{\Gamma}}_1 \left( s_{1,\infty };t_1,u_1,\eta\right)= \frac{\left( \frac{1+s_{1,\infty }+\sqrt{s_{1,\infty }^2-2(1-2\eta (1-t_1)(1-u_1))s_{1,\infty }+1}}{2}\right)^{-(4+\alpha )}}{\sqrt{s_{1,\infty }^2-2(1-2\eta (1-t_1)(1-u_1))s_{1,\infty }+1}}}\cr
{ \displaystyle  \overleftrightarrow {\mathbf{\Gamma}}_n \left( s_{n,\infty };t_n,u_n,\eta \right) =\frac{\left( \frac{1+s_{n,\infty }+\sqrt{s_{n,\infty }^2-2(1-2\eta (1-t_n)(1-u_n))s_{n,\infty }+1}}{2}\right)^{-(4n+\alpha ) }}{\sqrt{ s_{n,\infty }^2-2(1-2\eta (1-t_n)(1-u_n))s_{n,\infty }+1}}}\cr
{ \displaystyle \overleftrightarrow {\mathbf{\Gamma}}_{n-k} \left( s_{n-k};t_{n-k},u_{n-k},\widetilde{w}_{n-k+1,n} \right) = \frac{ \left( \frac{1+s_{n-k}+\sqrt{s_{n-k}^2-2(1-2\widetilde{w}_{n-k+1,n} (1-t_{n-k})(1-u_{n-k}))s_{n-k}+1}}{2}\right)^{-(4(n-k)+\alpha )}}{\sqrt{ s_{n-k}^2-2(1-2\widetilde{w}_{n-k+1,n} (1-t_{n-k})(1-u_{n-k}))s_{n-k}+1}}}
\end{cases}\nonumber 
\end{equation}
 
 \begin{equation}
\begin{cases} 
{ \displaystyle \mathbf{A} \left( s_{0,\infty } ;\eta\right)}\cr 
{ \displaystyle = \frac{\left(1- s_{0,\infty }+\sqrt{s_{0,\infty }^2-2(1-2\eta )s_{0,\infty }+1}\right)^{ -(\alpha +\frac{1}{2})} \left(1+s_{0,\infty }+\sqrt{s_{0,\infty }^2-2(1-2\eta )s_{0,\infty }+1}\right)^{\frac{1}{2} }}{\sqrt{s_{0,\infty }^2-2(1-2\eta )s_{0,\infty }+1}}}\cr
{ \displaystyle  \mathbf{A} \left( s_{0} ;\widetilde{w}_{1,1}\right) = \frac{\left(1- s_0+\sqrt{s_0^2-2(1-2\widetilde{w}_{1,1})s_0+1}\right)^{ -(\alpha +\frac{1}{2})} \left(1+s_0+\sqrt{s_0^2-2(1-2\widetilde{w}_{1,1} )s_0+1}\right)^{\frac{1}{2} }}{\sqrt{s_0^2-2(1-2\widetilde{w}_{1,1})s_0+1}}} \cr
{ \displaystyle \mathbf{A} \left(  s_{0} ;\widetilde{w}_{1,n}\right) = \frac{\left(1- s_0+\sqrt{s_0^2-2(1-2\widetilde{w}_{1,n})s_0+1}\right)^{ -(\alpha +\frac{1}{2})} \left(1+s_0+\sqrt{s_0^2-2(1-2\widetilde{w}_{1,n} )s_0+1}\right)^{\frac{1}{2} }}{\sqrt{s_0^2-2(1-2\widetilde{w}_{1,n})s_0+1}}}
\end{cases}\nonumber 
\end{equation}
and 
\begin{equation}
\begin{cases} 
h= (1+ \rho ^2)\left( 2(h_j +2j)+1+\alpha \right)^2 \cr
\varsigma= \xi^{-1} \cr
\eta =(1+\rho ^{-2})\varsigma \cr
z= -\rho ^{-2} \varsigma ^2  
\end{cases}\nonumber 
\end{equation}
\subsection{ ${\displaystyle \left(1-\frac{x}{a} \right)^{-\beta } Hl\left(1-a, -q+\gamma \beta; -\alpha +\gamma +\delta, \beta, \gamma, \delta; \frac{(1-a)x}{x-a} \right)}$}
Replace coefficients $a$, $q$, $\alpha $, $x$ and $q_j$ where $j, q_j \in \mathbb{N}_{0}$ by $1-a$, $-q+\gamma \beta $, $-\alpha+\gamma +\delta $, $\frac{(1-a)x}{x-a}$ and $h_j$ where $h_j \in \mathbb{N}_{0}$ into (\ref{eq:80019}). Put (\ref{eq:8006}) into the new (\ref{eq:80019}). 
\begin{eqnarray}
&&\sum_{h_0 =0}^{\infty } \frac{(\frac{1}{2})_{h_0}}{h_0!} s_0^{h_0} \prod _{n=1}^{\infty } \left\{ \sum_{ h_n = h_{n-1}}^{\infty } s_n^{h_n }\right\}  Hl\bigg( 1-\rho ^{-2}, (h_j+2j)\left(\alpha -(2-\rho ^{-2})(h_j+2j) \right); \nonumber\\
&&  -\frac{\alpha }{2}+\frac{1}{2}, -\frac{\alpha }{2}, \frac{1}{2}, \frac{1}{2}; \varsigma \bigg) \nonumber\\
&&= 2^{-\left( \Omega _{\rho }+1\right) }\Bigg\{ \prod_{l=1}^{\infty } \frac{1}{(1-s_{l,\infty })}  \mathbf{A}\left( s_{0,\infty } ;\eta\right) + \Bigg\{ \prod_{l=2}^{\infty } \frac{1}{(1-s_{l,\infty })} \int_{0}^{1} dt_1\;t_1 \int_{0}^{1} du_1\;u_1^{ \frac{1}{2}} \overleftrightarrow {\mathbf{\Gamma}}_1 \left( s_{1,\infty };t_1,u_1,\eta\right)\nonumber\\
&&\times \widetilde{w}_{1,1}^{-\frac{1}{2}(1- \alpha )}\left( \widetilde{w}_{1,1} \partial _{ \widetilde{w}_{1,1}}\right) \widetilde{w}_{1,1}^{ \frac{1}{2} } \left( \widetilde{w}_{1,1} \partial _{ \widetilde{w}_{1,1}}\right)\widetilde{w}_{1,1}^{-\frac{\alpha }{2}} \mathbf{A}\left( s_{0} ;\widetilde{w}_{1,1}\right) \Bigg\} z \nonumber\\
&&+ \sum_{n=2}^{\infty } \Bigg\{ \prod_{l=n+1}^{\infty } \frac{1}{(1-s_{l,\infty })} \int_{0}^{1} dt_n\;t_n^{2n-1} \int_{0}^{1} du_n\;u_n^{2n-\frac{3}{2}} \overleftrightarrow {\mathbf{\Gamma}}_n \left( s_{n,\infty };t_n,u_n,\eta \right)\nonumber\\
&&\times \widetilde{w}_{n,n}^{-\frac{1}{2}(4n-3-\alpha )}\left( \widetilde{w}_{n,n} \partial _{ \widetilde{w}_{n,n}}\right) \widetilde{w}_{n,n}^{ \frac{1}{2} } \left( \widetilde{w}_{n,n} \partial _{ \widetilde{w}_{n,n}}\right)\widetilde{w}_{n,n}^{\frac{1}{2}(4(n-1)-\alpha ) }  \nonumber\\
&&\times \prod_{k=1}^{n-1} \Bigg\{ \int_{0}^{1} dt_{n-k}\;t_{n-k}^{2(n-k)-1} \int_{0}^{1} du_{n-k} \;u_{n-k}^{2(n-k)-\frac{3}{2}}\overleftrightarrow {\mathbf{\Gamma}}_{n-k} \left( s_{n-k};t_{n-k},u_{n-k},\widetilde{w}_{n-k+1,n} \right)\label{eq:80064}\\
&&\times \widetilde{w}_{n-k,n}^{-\frac{1}{2}(4(n-k)-3-\alpha )}\left( \widetilde{w}_{n-k,n} \partial _{ \widetilde{w}_{n-k,n}}\right) \widetilde{w}_{n-k,n}^{ \frac{1}{2} } \left( \widetilde{w}_{n-k,n} \partial _{ \widetilde{w}_{n-k,n}}\right)\widetilde{w}_{n-k,n}^{\frac{1}{2}(4(n-k-1)- \alpha ) } \Bigg\} \mathbf{A} \left( s_{0} ;\widetilde{w}_{1,n}\right) \Bigg\} z^n \Bigg\}  \nonumber
\end{eqnarray}
where
\begin{equation}
\begin{cases} 
{ \displaystyle \overleftrightarrow {\mathbf{\Gamma}}_1 \left( s_{1,\infty };t_1,u_1,\eta\right)= \frac{\left( \frac{1+s_{1,\infty }+\sqrt{s_{1,\infty }^2-2(1-2\eta (1-t_1)(1-u_1))s_{1,\infty }+1}}{2}\right)^{-(3-\Omega _{\rho })}}{\sqrt{s_{1,\infty }^2-2(1-2\eta (1-t_1)(1-u_1))s_{1,\infty }+1}}}\cr
{ \displaystyle  \overleftrightarrow {\mathbf{\Gamma}}_n \left( s_{n,\infty };t_n,u_n,\eta \right) =\frac{\left( \frac{1+s_{n,\infty }+\sqrt{s_{n,\infty }^2-2(1-2\eta (1-t_n)(1-u_n))s_{n,\infty }+1}}{2}\right)^{-(4n-1-\Omega _{\rho }) }}{\sqrt{ s_{n,\infty }^2-2(1-2\eta (1-t_n)(1-u_n))s_{n,\infty }+1}}}\cr
{ \displaystyle \overleftrightarrow {\mathbf{\Gamma}}_{n-k} \left( s_{n-k};t_{n-k},u_{n-k},\widetilde{w}_{n-k+1,n} \right)} \cr 
{ \displaystyle = \frac{ \left( \frac{1+s_{n-k}+\sqrt{s_{n-k}^2-2(1-2\widetilde{w}_{n-k+1,n} (1-t_{n-k})(1-u_{n-k}))s_{n-k}+1}}{2}\right)^{-(4(n-k)-1-\Omega _{\rho } )}}{\sqrt{ s_{n-k}^2-2(1-2\widetilde{w}_{n-k+1,n} (1-t_{n-k})(1-u_{n-k}))s_{n-k}+1}}}
\end{cases}\nonumber 
\end{equation}
  
 \begin{equation}
\begin{cases} 
{ \displaystyle \mathbf{A} \left( s_{0,\infty } ;\eta\right)= \frac{\left(1- s_{0,\infty }+\sqrt{s_{0,\infty }^2-2(1-2\eta )s_{0,\infty }+1}\right)^{ \frac{1}{2} } \left(1+s_{0,\infty }+\sqrt{s_{0,\infty }^2-2(1-2\eta )s_{0,\infty }+1}\right)^{\frac{1}{2} +\Omega _{\rho }}}{\sqrt{s_{0,\infty }^2-2(1-2\eta )s_{0,\infty }+1}}}\cr
{ \displaystyle  \mathbf{A} \left( s_{0} ;\widetilde{w}_{1,1}\right) = \frac{\left(1- s_0+\sqrt{s_0^2-2(1-2\widetilde{w}_{1,1})s_0+1}\right)^{ \frac{1}{2} } \left(1+s_0+\sqrt{s_0^2-2(1-2\widetilde{w}_{1,1} )s_0+1}\right)^{\frac{1}{2} +\Omega _{\rho }}}{\sqrt{s_0^2-2(1-2\widetilde{w}_{1,1})s_0+1}}} \cr
{ \displaystyle \mathbf{A} \left(  s_{0} ;\widetilde{w}_{1,n}\right) = \frac{\left(1- s_0+\sqrt{s_0^2-2(1-2\widetilde{w}_{1,n})s_0+1}\right)^{ \frac{1}{2} } \left(1+s_0+\sqrt{s_0^2-2(1-2\widetilde{w}_{1,n} )s_0+1}\right)^{\frac{1}{2}+\Omega _{\rho } }}{\sqrt{s_0^2-2(1-2\widetilde{w}_{1,n})s_0+1}}}
\end{cases}\nonumber 
\end{equation}
and 
\begin{equation}
\begin{cases} 
\Omega _{\rho } = \frac{\alpha }{2-\rho ^{-2}} \cr
h= 4\rho ^2 \left( \alpha \left( h_j+2j+\frac{1}{4} \right) -(2-\rho ^{-2})(h_j+2j)^2\right)  \cr
\varsigma= \frac{(1-\rho ^{-2})\xi}{\xi-\rho ^{-2}} \cr
\eta =\frac{2-\rho ^{-2}}{1-\rho ^{-2}}\varsigma \cr
z= -\frac{1}{1-\rho ^{-2}} \varsigma ^2 
\end{cases}\nonumber 
\end{equation}
\subsection{ \footnotesize ${\displaystyle (1-x)^{1-\delta }\left(1-\frac{x}{a} \right)^{-\beta+\delta -1} Hl\left(1-a, -q+\gamma [(\delta -1)a+\beta -\delta +1]; -\alpha +\gamma +1, \beta -\delta+1, \gamma, 2-\delta; \frac{(1-a)x}{x-a} \right)}$ \normalsize}
Replace coefficients $a$, $q$, $\alpha $, $\beta $, $\delta $, $x$ and $q_j$ where $j, q_j \in \mathbb{N}_{0}$ by $1-a$, $-q+\gamma [(\delta -1)a+\beta -\delta +1]$, $-\alpha +\gamma +1$, $\beta -\delta+1$, $2-\delta $, $\frac{(1-a)x}{x-a}$ and $h_j$ where $h_j \in \mathbb{N}_{0}$ into (\ref{eq:80019}). Put (\ref{eq:8006}) into the new (\ref{eq:80019}). 
\begin{eqnarray}
&&\sum_{h_0 =0}^{\infty } \frac{(\frac{1}{2})_{h_0}}{h_0!} s_0^{h_0} \prod _{n=1}^{\infty } \left\{ \sum_{ h_n = h_{n-1}}^{\infty } s_n^{h_n }\right\}  Hl\bigg( 1-\rho ^{-2}, (h_j+ 2j)\left( \alpha +1 -(2-\rho ^{-2})(h_j+ 2j+1)\right); \nonumber\\
&& -\frac{\alpha }{2}+1, -\frac{\alpha }{2}+\frac{1}{2}, \frac{1}{2}, \frac{3}{2}; \varsigma \bigg) \nonumber\\
&&= 2^{-\Omega _{\rho }}\Bigg\{ \prod_{l=1}^{\infty } \frac{1}{(1-s_{l,\infty })}  \mathbf{A}\left( s_{0,\infty } ;\eta\right) + \Bigg\{ \prod_{l=2}^{\infty } \frac{1}{(1-s_{l,\infty })} \int_{0}^{1} dt_1\;t_1 \int_{0}^{1} du_1\;u_1^{ \frac{1}{2}} \overleftrightarrow {\mathbf{\Gamma}}_1 \left( s_{1,\infty };t_1,u_1,\eta\right)\nonumber\\
&&\times \widetilde{w}_{1,1}^{-\frac{1}{2}(2- \alpha )}\left( \widetilde{w}_{1,1} \partial _{ \widetilde{w}_{1,1}}\right) \widetilde{w}_{1,1}^{ \frac{1}{2} } \left( \widetilde{w}_{1,1} \partial _{ \widetilde{w}_{1,1}}\right)\widetilde{w}_{1,1}^{ \frac{1}{2}(1-\alpha )} \mathbf{A}\left( s_{0} ;\widetilde{w}_{1,1}\right) \Bigg\}z \nonumber\\
&&+ \sum_{n=2}^{\infty } \Bigg\{ \prod_{l=n+1}^{\infty } \frac{1}{(1-s_{l,\infty })} \int_{0}^{1} dt_n\;t_n^{2n-1} \int_{0}^{1} du_n\;u_n^{2n-\frac{3}{2}} \overleftrightarrow {\mathbf{\Gamma}}_n \left( s_{n,\infty };t_n,u_n,\eta \right)\nonumber\\
&&\times \widetilde{w}_{n,n}^{-\frac{1}{2}(4n-2-\alpha )}\left( \widetilde{w}_{n,n} \partial _{ \widetilde{w}_{n,n}}\right) \widetilde{w}_{n,n}^{ \frac{1}{2} } \left( \widetilde{w}_{n,n} \partial _{ \widetilde{w}_{n,n}}\right)\widetilde{w}_{n,n}^{\frac{1}{2}(4n-3-\alpha ) }  \nonumber\\
&&\times \prod_{k=1}^{n-1} \Bigg\{ \int_{0}^{1} dt_{n-k}\;t_{n-k}^{2(n-k)-1} \int_{0}^{1} du_{n-k} \;u_{n-k}^{2(n-k)-\frac{3}{2}}\overleftrightarrow {\mathbf{\Gamma}}_{n-k} \left( s_{n-k};t_{n-k},u_{n-k},\widetilde{w}_{n-k+1,n} \right)\label{eq:80065}\\
&&\times \widetilde{w}_{n-k,n}^{-\frac{1}{2}(4(n-k)-2-\alpha )}\left( \widetilde{w}_{n-k,n} \partial _{ \widetilde{w}_{n-k,n}}\right) \widetilde{w}_{n-k,n}^{ \frac{1}{2} } \left( \widetilde{w}_{n-k,n} \partial _{ \widetilde{w}_{n-k,n}}\right)\widetilde{w}_{n-k,n}^{\frac{1}{2}(4(n-k)-3-\alpha ) } \Bigg\} \mathbf{A} \left( s_{0} ;\widetilde{w}_{1,n}\right) \Bigg\} z^n \Bigg\} \nonumber 
\end{eqnarray}
where
\begin{equation}
\begin{cases} 
{ \displaystyle \overleftrightarrow {\mathbf{\Gamma}}_1 \left( s_{1,\infty };t_1,u_1,\eta\right)= \frac{\left( \frac{1+s_{1,\infty }+\sqrt{s_{1,\infty }^2-2(1-2\eta (1-t_1)(1-u_1))s_{1,\infty }+1}}{2}\right)^{- 4+\Omega _{\rho } }}{\sqrt{s_{1,\infty }^2-2(1-2\eta (1-t_1)(1-u_1))s_{1,\infty }+1}}}\cr
{ \displaystyle  \overleftrightarrow {\mathbf{\Gamma}}_n \left( s_{n,\infty };t_n,u_n,\eta \right) =\frac{\left( \frac{1+s_{n,\infty }+\sqrt{s_{n,\infty }^2-2(1-2\eta (1-t_n)(1-u_n))s_{n,\infty }+1}}{2}\right)^{- 4n+\Omega _{\rho } }}{\sqrt{ s_{n,\infty }^2-2(1-2\eta (1-t_n)(1-u_n))s_{n,\infty }+1}}}\cr
{ \displaystyle \overleftrightarrow {\mathbf{\Gamma}}_{n-k} \left( s_{n-k};t_{n-k},u_{n-k},\widetilde{w}_{n-k+1,n} \right) = \frac{ \left( \frac{1+s_{n-k}+\sqrt{s_{n-k}^2-2(1-2\widetilde{w}_{n-k+1,n} (1-t_{n-k})(1-u_{n-k}))s_{n-k}+1}}{2}\right)^{- 4(n-k)+\Omega _{\rho } }}{\sqrt{ s_{n-k}^2-2(1-2\widetilde{w}_{n-k+1,n} (1-t_{n-k})(1-u_{n-k}))s_{n-k}+1}}}
\end{cases}\nonumber 
\end{equation}
 
 \begin{equation}
\begin{cases} 
{ \displaystyle \mathbf{A} \left( s_{0,\infty } ;\eta\right)}\cr
{ \displaystyle = \frac{\left(1- s_{0,\infty }+\sqrt{s_{0,\infty }^2-2(1-2\eta )s_{0,\infty }+1}\right)^{ \frac{1}{2} } \left(1+s_{0,\infty }+\sqrt{s_{0,\infty }^2-2(1-2\eta )s_{0,\infty }+1}\right)^{-\frac{1}{2} +\Omega _{\rho }}}{\sqrt{s_{0,\infty }^2-2(1-2\eta )s_{0,\infty }+1}}}\cr
{ \displaystyle  \mathbf{A} \left( s_{0} ;\widetilde{w}_{1,1}\right) = \frac{\left(1- s_0+\sqrt{s_0^2-2(1-2\widetilde{w}_{1,1})s_0+1}\right)^{ \frac{1}{2} } \left(1+s_0+\sqrt{s_0^2-2(1-2\widetilde{w}_{1,1} )s_0+1}\right)^{-\frac{1}{2} +\Omega _{\rho }}}{\sqrt{s_0^2-2(1-2\widetilde{w}_{1,1})s_0+1}}} \cr
{ \displaystyle \mathbf{A} \left(  s_{0} ;\widetilde{w}_{1,n}\right) = \frac{\left(1- s_0+\sqrt{s_0^2-2(1-2\widetilde{w}_{1,n})s_0+1}\right)^{ \frac{1}{2} } \left(1+s_0+\sqrt{s_0^2-2(1-2\widetilde{w}_{1,n} )s_0+1}\right)^{-\frac{1}{2}+\Omega _{\rho } }}{\sqrt{s_0^2-2(1-2\widetilde{w}_{1,n})s_0+1}}}
\end{cases}\nonumber 
\end{equation}
and  
\begin{equation}
\begin{cases} 
\Omega _{\rho } = \frac{ \alpha +1}{2-\rho ^{-2}} \cr
h= 4(h_j+ 2j)\left( (\alpha +1)\rho ^2 +(1-2\rho ^2)(h_j+2j+1)\right) +(\alpha -1)\rho ^2 +1  \cr
\varsigma= \frac{(1-\rho ^{-2})\xi}{\xi-\rho ^{-2}} \cr
\eta =\frac{2-\rho ^{-2}}{1-\rho ^{-2}}\varsigma \cr
z= -\frac{1}{1-\rho ^{-2}} \varsigma ^2 
\end{cases}\nonumber 
\end{equation}
\subsection{ ${\displaystyle x^{-\alpha } Hl\left(\frac{a-1}{a}, \frac{-q+\alpha (\delta a+\beta -\delta )}{a}; \alpha, \alpha -\gamma +1, \delta , \alpha -\beta +1; \frac{x-1}{x} \right)}$}
\subsubsection{Polynomial of type 2}
Replace coefficients $a$, $q$, $\beta $, $\gamma $, $\delta $, $x$ and $q_j$ where $j, q_j \in \mathbb{N}_{0}$ by $\frac{a-1}{a}$, $\frac{-q+\alpha (\delta a+\beta -\delta )}{a}$, $\alpha -\gamma +1$, $\delta $, $\alpha -\beta +1$, $\frac{x-1}{x}$ and $h_j$ where $h_j \in \mathbb{N}_{0}$ into (\ref{eq:80019}). Put (\ref{eq:8006}) into the new (\ref{eq:80019}).
\begin{eqnarray}
&&\sum_{h_0 =0}^{\infty } \frac{(\frac{1}{2})_{h_0}}{h_0!} s_0^{h_0} \prod _{n=1}^{\infty } \left\{ \sum_{ h_n = h_{n-1}}^{\infty } s_n^{h_n }\right\}  Hl\bigg( 1-\rho ^2, -(h_j+2j)\left( (1-\rho ^2)(\alpha +1)+(2-\rho ^2)(h_j+2j) \right);  \nonumber\\
&& \frac{\alpha}{2} + \frac{1}{2}, \frac{\alpha}{2} +1, \frac{1}{2}, \alpha +\frac{3}{2}; \varsigma \bigg) \nonumber\\
&&= 2^{ \Omega _{\rho }-1}\Bigg\{ \prod_{l=1}^{\infty } \frac{1}{(1-s_{l,\infty })}  \mathbf{A}\left( s_{0,\infty } ;\eta\right) + \Bigg\{ \prod_{l=2}^{\infty } \frac{1}{(1-s_{l,\infty })} \int_{0}^{1} dt_1\;t_1 \int_{0}^{1} du_1\;u_1^{ \frac{1}{2}} \overleftrightarrow {\mathbf{\Gamma}}_1 \left( s_{1,\infty };t_1,u_1,\eta\right)\nonumber\\
&&\times \widetilde{w}_{1,1}^{-\frac{1}{2}(1+\alpha )}\left( \widetilde{w}_{1,1} \partial _{ \widetilde{w}_{1,1}}\right) \widetilde{w}_{1,1}^{ -\frac{1}{2} } \left( \widetilde{w}_{1,1} \partial _{ \widetilde{w}_{1,1}}\right)\widetilde{w}_{1,1}^{ \frac{1}{2}(2+\alpha )} \mathbf{A}\left( s_{0} ;\widetilde{w}_{1,1}\right)\Bigg\} z \nonumber\\
&&+ \sum_{n=2}^{\infty } \Bigg\{ \prod_{l=n+1}^{\infty } \frac{1}{(1-s_{l,\infty })} \int_{0}^{1} dt_n\;t_n^{2n-1} \int_{0}^{1} du_n\;u_n^{2n-\frac{3}{2}} \overleftrightarrow {\mathbf{\Gamma}}_n \left( s_{n,\infty };t_n,u_n,\eta \right)\nonumber\\
&&\times \widetilde{w}_{n,n}^{-\frac{1}{2}(4n-3+\alpha )}\left( \widetilde{w}_{n,n} \partial _{ \widetilde{w}_{n,n}}\right) \widetilde{w}_{n,n}^{ -\frac{1}{2} } \left( \widetilde{w}_{n,n} \partial _{ \widetilde{w}_{n,n}}\right)\widetilde{w}_{n,n}^{\frac{1}{2}(4n-2+\alpha ) }  \nonumber\\
&&\times \prod_{k=1}^{n-1} \Bigg\{ \int_{0}^{1} dt_{n-k}\;t_{n-k}^{2(n-k)-1} \int_{0}^{1} du_{n-k} \;u_{n-k}^{2(n-k)-\frac{3}{2}}\overleftrightarrow {\mathbf{\Gamma}}_{n-k} \left( s_{n-k};t_{n-k},u_{n-k},\widetilde{w}_{n-k+1,n} \right)\label{eq:80066}\\
&&\times \widetilde{w}_{n-k,n}^{-\frac{1}{2}(4(n-k)-3+\alpha )}\left( \widetilde{w}_{n-k,n} \partial _{ \widetilde{w}_{n-k,n}}\right) \widetilde{w}_{n-k,n}^{ -\frac{1}{2} } \left( \widetilde{w}_{n-k,n} \partial _{ \widetilde{w}_{n-k,n}}\right)\widetilde{w}_{n-k,n}^{\frac{1}{2}(4(n-k)-2+\alpha ) } \Bigg\} \mathbf{A} \left( s_{0} ;\widetilde{w}_{1,n}\right) \Bigg\} z^n \Bigg\} \nonumber 
\end{eqnarray}
where
\begin{equation}
\begin{cases} 
{ \displaystyle \overleftrightarrow {\mathbf{\Gamma}}_1 \left( s_{1,\infty };t_1,u_1,\eta\right)= \frac{\left( \frac{1+s_{1,\infty }+\sqrt{s_{1,\infty }^2-2(1-2\eta (1-t_1)(1-u_1))s_{1,\infty }+1}}{2}\right)^{-(3+\Omega _{\rho }) }}{\sqrt{s_{1,\infty }^2-2(1-2\eta (1-t_1)(1-u_1))s_{1,\infty }+1}}}\cr
{ \displaystyle  \overleftrightarrow {\mathbf{\Gamma}}_n \left( s_{n,\infty };t_n,u_n,\eta \right) =\frac{\left( \frac{1+s_{n,\infty }+\sqrt{s_{n,\infty }^2-2(1-2\eta (1-t_n)(1-u_n))s_{n,\infty }+1}}{2}\right)^{-(4n-1+\Omega _{\rho } )}}{\sqrt{ s_{n,\infty }^2-2(1-2\eta (1-t_n)(1-u_n))s_{n,\infty }+1}}}\cr
{ \displaystyle \overleftrightarrow {\mathbf{\Gamma}}_{n-k} \left( s_{n-k};t_{n-k},u_{n-k},\widetilde{w}_{n-k+1,n} \right) } \cr
{ \displaystyle = \frac{ \left( \frac{1+s_{n-k}+\sqrt{s_{n-k}^2-2(1-2\widetilde{w}_{n-k+1,n} (1-t_{n-k})(1-u_{n-k}))s_{n-k}+1}}{2}\right)^{-(4(n-k)-1+\Omega _{\rho } )}}{\sqrt{ s_{n-k}^2-2(1-2\widetilde{w}_{n-k+1,n} (1-t_{n-k})(1-u_{n-k}))s_{n-k}+1}}}
\end{cases}\nonumber 
\end{equation}
 
 \begin{equation}
\begin{cases} 
{ \displaystyle \mathbf{A} \left( s_{0,\infty } ;\eta\right)= \frac{\left(1- s_{0,\infty }+\sqrt{s_{0,\infty }^2-2(1-2\eta )s_{0,\infty }+1}\right)^{ \frac{1}{2} } \left(1+s_{0,\infty }+\sqrt{s_{0,\infty }^2-2(1-2\eta )s_{0,\infty }+1}\right)^{ \frac{1}{2} -\Omega _{\rho }}}{\sqrt{s_{0,\infty }^2-2(1-2\eta )s_{0,\infty }+1}}}\cr
{ \displaystyle  \mathbf{A} \left( s_{0} ;\widetilde{w}_{1,1}\right) = \frac{\left(1- s_0+\sqrt{s_0^2-2(1-2\widetilde{w}_{1,1})s_0+1}\right)^{ \frac{1}{2} } \left(1+s_0+\sqrt{s_0^2-2(1-2\widetilde{w}_{1,1} )s_0+1}\right)^{ \frac{1}{2} -\Omega _{\rho }}}{\sqrt{s_0^2-2(1-2\widetilde{w}_{1,1})s_0+1}}} \cr
{ \displaystyle \mathbf{A} \left(  s_{0} ;\widetilde{w}_{1,n}\right) = \frac{\left(1- s_0+\sqrt{s_0^2-2(1-2\widetilde{w}_{1,n})s_0+1}\right)^{ \frac{1}{2} } \left(1+s_0+\sqrt{s_0^2-2(1-2\widetilde{w}_{1,n} )s_0+1}\right)^{ \frac{1}{2}-\Omega _{\rho } }}{\sqrt{s_0^2-2(1-2\widetilde{w}_{1,n})s_0+1}}}
\end{cases}\nonumber 
\end{equation}
and  
\begin{equation}
\begin{cases} 
\Omega _{\rho } = \frac{(\alpha +1)(1-\rho ^2)}{2-\rho ^2} \cr
h= -4(h_j+ 2j)\left( (1- \rho ^2)(\alpha +1) +(2-\rho ^2)(h_j+2j)\right) -(\alpha +1)\left( 1-\rho ^2(\alpha +1)\right) \cr
\varsigma= 1-\xi ^{-1} \cr
\eta =\frac{2-\rho ^2}{1-\rho ^2}\varsigma \cr
z= -\frac{1}{1-\rho ^2} \varsigma ^2 
\end{cases}\nonumber 
\end{equation}
\subsection{ ${\displaystyle \left(\frac{x-a}{1-a} \right)^{-\alpha } Hl\left(a, q-(\beta -\delta )\alpha ; \alpha , -\beta+\gamma +\delta , \delta , \gamma; \frac{a(x-1)}{x-a} \right)}$}
Replace coefficients $q$, $\beta $, $\gamma $, $\delta $, $x$ and $q_j$ where $j, q_j \in \mathbb{N}_{0}$ by $q-(\beta -\delta )\alpha $, $-\beta+\gamma +\delta $, $\delta $,  $\gamma $, $\frac{a(x-1)}{x-a}$ and $h_j$ where $h_j \in \mathbb{N}_{0}$ into (\ref{eq:80019}). Put (\ref{eq:8006}) into the new (\ref{eq:80019}).
\begin{eqnarray}
&&\sum_{h_0 =0}^{\infty } \frac{(\frac{1}{2})_{h_0}}{h_0!} s_0^{h_0} \prod _{n=1}^{\infty } \left\{ \sum_{ h_n = h_{n-1}}^{\infty } s_n^{h_n }\right\}  Hl\left( \rho ^{-2}, -(h_j+2j)\left( (1+\rho ^{-2})(h_j+2j)+\alpha +1\right); \frac{\alpha }{2} +\frac{1}{2}, \right.\nonumber\\
&&\left. \frac{\alpha }{2} +1, \frac{1}{2}, \frac{1}{2}; \varsigma \right) \nonumber\\
&&= 2^{ \Omega _{\rho }-1}\Bigg\{ \prod_{l=1}^{\infty } \frac{1}{(1-s_{l,\infty })}  \mathbf{A}\left( s_{0,\infty } ;\eta\right) + \Bigg\{ \prod_{l=2}^{\infty } \frac{1}{(1-s_{l,\infty })} \int_{0}^{1} dt_1\;t_1 \int_{0}^{1} du_1\;u_1^{ \frac{1}{2}} \overleftrightarrow {\mathbf{\Gamma}}_1 \left( s_{1,\infty };t_1,u_1,\eta\right)\nonumber\\
&&\times \widetilde{w}_{1,1}^{-\frac{1}{2}(1+\alpha )}\left( \widetilde{w}_{1,1} \partial _{ \widetilde{w}_{1,1}}\right) \widetilde{w}_{1,1}^{ -\frac{1}{2} } \left( \widetilde{w}_{1,1} \partial _{ \widetilde{w}_{1,1}}\right)\widetilde{w}_{1,1}^{ \frac{1}{2}(2+\alpha )} \mathbf{A}\left( s_{0} ;\widetilde{w}_{1,1}\right)\Bigg\} z \nonumber\\
&&+ \sum_{n=2}^{\infty } \Bigg\{ \prod_{l=n+1}^{\infty } \frac{1}{(1-s_{l,\infty })} \int_{0}^{1} dt_n\;t_n^{2n-1} \int_{0}^{1} du_n\;u_n^{2n-\frac{3}{2}} \overleftrightarrow {\mathbf{\Gamma}}_n \left( s_{n,\infty };t_n,u_n,\eta \right)\nonumber\\
&&\times \widetilde{w}_{n,n}^{-\frac{1}{2}(4n-3+\alpha )}\left( \widetilde{w}_{n,n} \partial _{ \widetilde{w}_{n,n}}\right) \widetilde{w}_{n,n}^{ -\frac{1}{2} } \left( \widetilde{w}_{n,n} \partial _{ \widetilde{w}_{n,n}}\right)\widetilde{w}_{n,n}^{\frac{1}{2}(4n-2+\alpha ) }  \nonumber\\
&&\times \prod_{k=1}^{n-1} \Bigg\{ \int_{0}^{1} dt_{n-k}\;t_{n-k}^{2(n-k)-1} \int_{0}^{1} du_{n-k} \;u_{n-k}^{2(n-k)-\frac{3}{2}}\overleftrightarrow {\mathbf{\Gamma}}_{n-k} \left( s_{n-k};t_{n-k},u_{n-k},\widetilde{w}_{n-k+1,n} \right)\label{eq:80067}\\
&&\times \widetilde{w}_{n-k,n}^{-\frac{1}{2}(4(n-k)-3+\alpha )}\left( \widetilde{w}_{n-k,n} \partial _{ \widetilde{w}_{n-k,n}}\right) \widetilde{w}_{n-k,n}^{ -\frac{1}{2} } \left( \widetilde{w}_{n-k,n} \partial _{ \widetilde{w}_{n-k,n}}\right)\widetilde{w}_{n-k,n}^{\frac{1}{2}(4(n-k)-2+\alpha ) } \Bigg\} \mathbf{A} \left( s_{0} ;\widetilde{w}_{1,n}\right) \Bigg\} z^n \Bigg\} \nonumber 
\end{eqnarray}
where
\begin{equation}
\begin{cases} 
{ \displaystyle \overleftrightarrow {\mathbf{\Gamma}}_1 \left( s_{1,\infty };t_1,u_1,\eta\right)= \frac{\left( \frac{1+s_{1,\infty }+\sqrt{s_{1,\infty }^2-2(1-2\eta (1-t_1)(1-u_1))s_{1,\infty }+1}}{2}\right)^{-(3+\Omega _{\rho }) }}{\sqrt{s_{1,\infty }^2-2(1-2\eta (1-t_1)(1-u_1))s_{1,\infty }+1}}}\cr
{ \displaystyle  \overleftrightarrow {\mathbf{\Gamma}}_n \left( s_{n,\infty };t_n,u_n,\eta \right) =\frac{\left( \frac{1+s_{n,\infty }+\sqrt{s_{n,\infty }^2-2(1-2\eta (1-t_n)(1-u_n))s_{n,\infty }+1}}{2}\right)^{-(4n-1+\Omega _{\rho }) }}{\sqrt{ s_{n,\infty }^2-2(1-2\eta (1-t_n)(1-u_n))s_{n,\infty }+1}}}\cr
{ \displaystyle \overleftrightarrow {\mathbf{\Gamma}}_{n-k} \left( s_{n-k};t_{n-k},u_{n-k},\widetilde{w}_{n-k+1,n} \right) } \cr
{ \displaystyle = \frac{ \left( \frac{1+s_{n-k}+\sqrt{s_{n-k}^2-2(1-2\widetilde{w}_{n-k+1,n} (1-t_{n-k})(1-u_{n-k}))s_{n-k}+1}}{2}\right)^{-(4(n-k)-1+\Omega _{\rho } )}}{\sqrt{ s_{n-k}^2-2(1-2\widetilde{w}_{n-k+1,n} (1-t_{n-k})(1-u_{n-k}))s_{n-k}+1}}}
\end{cases}\nonumber 
\end{equation}
  
 \begin{equation}
\begin{cases} 
{ \displaystyle \mathbf{A} \left( s_{0,\infty } ;\eta\right)= \frac{\left(1- s_{0,\infty }+\sqrt{s_{0,\infty }^2-2(1-2\eta )s_{0,\infty }+1}\right)^{ \frac{1}{2} } \left(1+s_{0,\infty }+\sqrt{s_{0,\infty }^2-2(1-2\eta )s_{0,\infty }+1}\right)^{ \frac{1}{2} -\Omega _{\rho }}}{\sqrt{s_{0,\infty }^2-2(1-2\eta )s_{0,\infty }+1}}}\cr
{ \displaystyle  \mathbf{A} \left( s_{0} ;\widetilde{w}_{1,1}\right) = \frac{\left(1- s_0+\sqrt{s_0^2-2(1-2\widetilde{w}_{1,1})s_0+1}\right)^{ \frac{1}{2} } \left(1+s_0+\sqrt{s_0^2-2(1-2\widetilde{w}_{1,1} )s_0+1}\right)^{ \frac{1}{2} -\Omega _{\rho }}}{\sqrt{s_0^2-2(1-2\widetilde{w}_{1,1})s_0+1}}} \cr
{ \displaystyle \mathbf{A} \left(  s_{0} ;\widetilde{w}_{1,n}\right) = \frac{\left(1- s_0+\sqrt{s_0^2-2(1-2\widetilde{w}_{1,n})s_0+1}\right)^{ \frac{1}{2} } \left(1+s_0+\sqrt{s_0^2-2(1-2\widetilde{w}_{1,n} )s_0+1}\right)^{ \frac{1}{2}-\Omega _{\rho } }}{\sqrt{s_0^2-2(1-2\widetilde{w}_{1,n})s_0+1}}}
\end{cases}\nonumber 
\end{equation}
and   
\begin{equation}
\begin{cases}
\Omega _{\rho } = \frac{ \alpha +1 }{1+\rho ^{-2}} \cr
h= 4(h_j+ 2j)\left( (1+ \rho ^2)(h_j+ 2j) +\rho ^2(\alpha +1) \right) +\rho ^2(\alpha +1)^2 \cr 
\varsigma= \frac{\xi -1}{\rho ^2(\xi -\rho ^{-2})} \cr
\eta = (1+\rho ^2) \varsigma \cr
z= -\rho ^2 \varsigma ^2 
\end{cases}\nonumber 
\end{equation}
\end{appendices} 

\addcontentsline{toc}{section}{Bibliography}
\bibliographystyle{model1a-num-names}
\bibliography{<your-bib-database>}


\end{document}